\documentclass{article}
\usepackage{amsmath}
\usepackage{amssymb}
\usepackage{amsthm}
\usepackage{bbm}
\usepackage{cite}
\usepackage{stmaryrd}
\SetSymbolFont{stmry}{bold}{U}{stmry}{m}{n}
\usepackage{euscript}
\usepackage[arrow,curve,matrix,arc,2cell]{xy}
\UseAllTwocells
\usepackage[utf8]{inputenc}
\usepackage[pdftex,breaklinks,unicode]{hyperref}
\usepackage{slashed}
\usepackage{tikz}
\usepackage{tikz-cd}
\usepackage[new]{old-arrows}
\usetikzlibrary{calc}
\usetikzlibrary{positioning}
\usetikzlibrary{matrix,arrows}
\usetikzlibrary{decorations.markings}
\DeclareFontFamily{U}{rsfs}{} 
\DeclareFontShape{U}{rsfs}{n}{it}{<->
rsfs10}{} \DeclareSymbolFont{mscr}{U}{rsfs}{n}{it}
\DeclareSymbolFontAlphabet{\scr}{mscr}
\def\mathscr{\scr}
\begin{document}
\def\e#1\e{\begin{equation}#1\end{equation}}
\def\ea#1\ea{\begin{align}#1\end{align}}
\def\eq#1{{\rm(\ref{#1})}}
\theoremstyle{plain}
\newtheorem{thm}{Theorem}[section]
\newtheorem{lem}[thm]{Lemma}
\newtheorem{prop}[thm]{Proposition}
\newtheorem{cor}[thm]{Corollary}
\newtheorem{quest}[thm]{Question}
\newtheorem{prob}[thm]{Problem}
\newtheorem{princ}[thm]{Principle}
\newtheorem{claim}[thm]{Claim}
\theoremstyle{definition}
\newtheorem{dfn}[thm]{Definition}
\newtheorem{ex}[thm]{Example}
\newtheorem{rem}[thm]{Remark}
\newtheorem{conj}[thm]{Conjecture}
\newtheorem{ax}[thm]{Axiom}
\newtheorem{ass}[thm]{Assumption}
\newtheorem{property}[thm]{Property}
\newtheorem{cond}[thm]{Condition}
\newtheorem{meth}[thm]{Method}
\newtheorem{alg}[thm]{Algorithm}
\numberwithin{figure}{section}
\numberwithin{equation}{section}
\def\dim{\mathop{\rm dim}\nolimits}
\def\codim{\mathop{\rm codim}\nolimits}
\def\vdim{\mathop{\rm vdim}\nolimits}
\def\sign{\mathop{\rm sign}\nolimits}
\def\Im{\mathop{\rm Im}\nolimits}
\def\det{\mathop{\rm det}\nolimits}
\def\Res{\mathop{\rm Res}\nolimits}
\def\Ker{\mathop{\rm Ker}}
\def\Kah{\mathop{\text{\rm K\"ah}}\nolimits}
\def\Coker{\mathop{\rm Coker}}
\def\Spec{\mathop{\rm Spec}}
\def\Perf{\mathop{\rm Perf}}
\def\HNT{\mathop{\rm HNT}\nolimits}
\def\HNP{\mathop{\rm HNP}\nolimits}
\def\Vect{\mathop{\rm Vect}}
\def\LCon{\mathop{\rm LCon}}
\def\Flag{\mathop{\rm Flag}\nolimits}
\def\FlagSt{\mathop{\rm FlagSt}\nolimits}
\def\Amp{\mathop{\rm Amp}\nolimits}
\def\Div{\mathop{\rm Div}\nolimits}
\def\Iso{\mathop{\rm Iso}\nolimits}
\def\lub{\mathop{\rm lub}\nolimits}
\def\glb{\mathop{\rm glb}\nolimits}
\def\Pd{\mathop{\rm Pd}}
\def\Aut{\mathop{\rm Aut}}
\def\End{\mathop{\rm End}}
\def\Ho{\mathop{\rm Ho}}
\def\Hol{\mathop{\rm Hol}}
\def\PGL{\mathop{\rm PGL}\nolimits}
\def\GL{\mathop{\rm GL}\nolimits}
\def\SL{\mathop{\rm SL}}
\def\SO{\mathop{\rm SO}}
\def\SU{\mathop{\rm SU}}
\def\Sp{\mathop{\rm Sp}}
\def\Re{\mathop{\rm Re}}
\def\Ch{\mathop{\rm Ch}\nolimits}
\def\Hilb{\mathop{\rm Hilb}\nolimits}
\def\Spin{\mathop{\rm Spin}}
\def\Spinc{\mathop{\rm Spin^c}}
\def\SF{\mathop{\rm SF}}
\def\SFa{\mathop{\rm SF_{al}}\nolimits}
\def\SFai{\mathop{\rm SF_{al}^{ind}}\nolimits}
\def\CF{\mathop{\rm CF}\nolimits}
\def\CFi{\mathop{\rm CF^{ind}}\nolimits}
\def\Tr{\mathop{\rm Tr}}
\def\U{{\mathbin{\rm U}}}
\def\vol{\mathop{\rm vol}}
\def\inc{\mathop{\rm inc}\nolimits}
\def\ev{\mathop{\rm ev}\nolimits}
\def\go{{\rm go}}
\def\na{{\rm na}}
\def\td{\mathop{\rm td}}
\def\ad{{\rm ad}}
\def\stk{{\rm stk}}
\def\ch{\mathop{\rm ch}\nolimits}
\def\ind{\mathop{\rm ind}\nolimits}
\def\tind{{\text{\rm t-ind}}}
\def\bdim{{\mathbin{\bf dim}\kern.1em}}
\def\rk{\mathop{\rm rk}}
\def\Pic{\mathop{\rm Pic}}
\def\colim{\mathop{\rm colim}\nolimits}
\def\Stab{\mathop{\rm Stab}\nolimits}
\def\Exact{\mathop{\rm Exact}\nolimits}
\def\Crit{\mathop{\rm Crit}}
\def\Cass{\mathop{\rm Cass}}
\def\supp{\mathop{\rm supp}}
\def\Or{\mathop{\rm Or}}
\def\rank{\mathop{\rm rank}\nolimits}
\def\Hom{\mathop{\rm Hom}\nolimits}
\def\bHom{\mathop{\bf Hom}\nolimits}
\def\Ext{\mathop{\rm Ext}\nolimits}
\def\cExt{\mathop{{\mathcal E}\mathit{xt}}\nolimits}
\def\acExt{\mathop{\smash{\ac{\mathcal E}}\vphantom{{\mathcal E}}\mathit{xt}}\nolimits}
\def\racExt{\mathop{\smash{\ac{\rm E}}\vphantom{{\rm E}}{\rm xt}}\nolimits}
\def\acE{\mathop{\smash{\ac{\mathcal E}}\vphantom{{\mathcal E}}}\nolimits}
\def\fExact{\mathop{\mathfrak{Exact}}}
\def\id{{\mathop{\rm id}\nolimits}}
\def\Id{{\mathop{\rm Id}\nolimits}}
\def\Sch{\mathop{\bf Sch}\nolimits}
\def\Map{{\mathop{\rm Map}\nolimits}}
\def\AlgSp{\mathop{\bf AlgSp}\nolimits}
\def\Sta{\mathop{\bf Sta}\nolimits}
\def\AlgSp{\mathop{\bf AlgSp}\nolimits}
\def\Art{\mathop{\bf Art}\nolimits}
\def\DMSta{\mathop{\bf DMSta}\nolimits}
\def\HSta{\mathop{\bf HSta}\nolimits}
\def\DSta{\mathop{\bf DSta}\nolimits}
\def\DArt{\mathop{\bf DArt}\nolimits}
\def\VertAlg{\mathop{\bf VertAlg}\nolimits}
\def\LieAlg{\mathop{\bf LieAlg}\nolimits}
\def\rgr{{\bf gr}}
\def\TopSta{{\mathop{\bf TopSta}\nolimits}}
\def\Topho{{\mathop{\bf Top^{ho}}}}
\def\Gpds{{\mathop{\bf Gpds}\nolimits}}
\def\Aff{\mathop{\bf Aff}\nolimits}
\def\SPr{\mathop{\bf SPr}\nolimits}
\def\SSet{\mathop{\bf SSet}\nolimits}
\def\Ad{\mathop{\rm Ad}}
\def\pl{{\rm pl}}
\def\rpl{{\rm rpl}}
\def\eu{{\rm eu}}
\def\mix{{\rm mix}}
\def\fd{{\rm fd}}
\def\wt{{\rm wt}}
\def\ran{{\rm an}}
\def\fpd{{\rm fpd}}
\def\pfd{{\rm pfd}}
\def\coa{{\rm coa}}
\def\rsi{{\rm si}}
\def\rst{{\rm st}}
\def\ss{{\rm ss}}
\def\vi{{\rm vi}}
\def\smq{{\rm smq}}
\def\rsm{{\rm sm}}
\def\cla{{\rm cla}}
\def\rArt{{\rm Art}}
\def\po{{\rm po}}
\def\pe{{\rm pe}}
\def\rp{{\rm rp}}
\def\spo{{\rm spo}}
\def\Kur{{\rm Kur}}
\def\dcr{{\rm dcr}}
\def\top{{\rm top}}
\def\fc{{\rm fc}}
\def\cla{{\rm cla}}
\def\num{{\rm num}}
\def\pr{{\rm pr}}
\def\irr{{\rm irr}}
\def\red{{\rm red}}
\def\sing{{\rm sing}}
\def\virt{{\rm virt}}
\def\inv{{\rm inv}}
\def\fund{{\rm fund}}
\def\qcoh{{\rm qcoh}}
\def\coh{{\rm coh}}
\def\vect{{\rm vect}}
\def\lft{{\rm lft}}
\def\lfp{{\rm lfp}}
\def\cs{{\rm cs}}
\def\dR{{\rm dR}}
\def\Obj{{\rm Obj}}
\def\cdga{{\mathop{\bf cdga}\nolimits}}
\def\Rmod{\mathop{R\text{\rm -mod}}}
\def\Top{{\mathop{\bf Top}\nolimits}}
\def\modKQ{\mathop{\text{\rm mod-}\K Q}}
\def\modKQI{\text{\rm mod-$\K Q/I$}}
\def\modCQ{\mathop{\text{\rm mod-}\C Q}}
\def\modCdQ{\mathop{\text{\rm mod-}\C\dot Q}}
\def\modCQI{\text{\rm mod-$\C Q/I$}}
\def\modCtQ{\mathop{\text{\rm mod-}\C\ti Q}}
\def\modCtQI{\mathop{\text{\rm mod-}\C\ti Q/\ti I}}
\def\modCaQI{\text{\rm mod-$\C\ac Q/\ac I$}}
\def\ul{\underline}
\def\bs{\boldsymbol}
\def\ge{\geqslant}
\def\le{\leqslant\nobreak}
\def\boo{{\mathbin{\mathbbm 1}}}
\def\ds{/\!/}
\def\O{{\mathcal O}}
\def\bA{{\mathbin{\mathbb A}}}
\def\bG{{\mathbin{\mathbb G}}}
\def\bH{{\mathbin{\mathbb H}}}
\def\bL{{\mathbin{\mathbb L}}}
\def\P{{\mathbin{\mathbb P}}}
\def\K{{\mathbin{\mathbb K}}}
\def\R{{\mathbin{\mathbb R}}}
\def\bT{{\mathbin{\mathbb T}}}
\def\Z{{\mathbin{\mathbb Z}}}
\def\bP{{\mathbin{\mathbb P}}}
\def\Q{{\mathbin{\mathbb Q}}}
\def\N{{\mathbin{\mathbb N}}}
\def\C{{\mathbin{\mathbb C}}}
\def\CP{{\mathbin{\mathbb{CP}}}}
\def\KP{{\mathbin{\mathbb{KP}}}}
\def\RP{{\mathbin{\mathbb{RP}}}}
\def\fC{{\mathbin{\mathfrak C}\kern.05em}}
\def\fD{{\mathbin{\mathfrak D}}}
\def\fE{{\mathbin{\mathfrak E}}}
\def\fF{{\mathbin{\mathfrak F}}}
\def\A{{\mathbin{\cal A}}}
\def\haA{{\mathbin{\hat{\cal A}}}}
\def\tiA{{\mathbin{\ti{\cal A}}}}
\def\dA{{\mathbin{\dot{\cal A}}}}
\def\ddA{{\mathbin{\ddot{\cal A}}}}
\def\dddA{{\mathbin{\dddot{\cal A}}}}
\def\baA{{\mathbin{\smash{\bar{\cal A}}}\vphantom{\cal A}}}
\def\grA{{\mathbin{\smash{\grave{\cal A}}}\vphantom{\cal A}}}
\def\brA{{\mathbin{\smash{\breve{\cal A}}}\vphantom{\cal A}}}
\def\acA{{\mathbin{\smash{\acute{\cal A}}}\vphantom{\cal A}}}
\def\G{{{\cal G}}}
\def\M{{\mathbin{\cal M}}}
\def\fM{{\mathbin{\mathfrak M}}}
\def\B{{\mathbin{\cal B}}}  
\def\haB{{\mathbin{\hat{\cal B}}}}
\def\tiB{{\mathbin{\ti{\cal B}}}}
\def\dB{{\mathbin{\dot{\cal B}}}}
\def\ddB{{\mathbin{\ddot{\cal B}}}}
\def\dddB{{\mathbin{\dddot{\cal B}}}}
\def\baB{{\mathbin{\smash{\bar{\cal B}}}\vphantom{\cal B}}}
\def\grB{{\mathbin{\smash{\grave{\cal B}}}\vphantom{\cal B}}}
\def\brB{{\mathbin{\smash{\breve{\cal B}}}\vphantom{\cal B}}}
\def\acB{{\mathbin{\acute{\cal B}}}}
\def\ovB{{\mathbin{\smash{\,\overline{\!\mathcal B}}}}}
\def\cC{{\mathbin{\cal C}}}
\def\cD{{\mathbin{\cal D}}}
\def\cE{{\mathbin{\cal E}}}
\def\baE{{\mathbin{\smash{\bar{\cal E}}}\vphantom{\cal E}}}
\def\dE{{\mathbin{\dot{\cal E}}}{}}
\def\cF{{\mathbin{\cal F}}}
\def\cG{{\mathbin{\cal G}}}
\def\cH{{\mathbin{\cal H}}}
\def\cI{{\mathbin{\cal I}}}
\def\cJ{{\mathbin{\cal J}}}
\def\cK{{\mathbin{\cal K}}}
\def\cL{{\mathbin{\cal L}}}
\def\bcM{{\mathbin{\bs{\cal M}}}}
\def\cN{{\cal N}}
\def\cP{{\mathbin{\cal P}}}
\def\cQ{{\mathbin{\cal Q}}}
\def\cR{{\mathbin{\cal R}}}
\def\cS{{\mathbin{\cal S}}}
\def\T{{{\cal T}\kern .04em}}
\def\cU{{\mathbin{\cal U}}}
\def\cV{{\mathbin{\cal V}}}
\def\cW{{\mathbin{\cal W}}}
\def\cX{{\cal X}}
\def\cY{{\cal Y}}
\def\cZ{{\cal Z}}
\def\oM{{\mathbin{\smash{\,\,\overline{\!\!\mathcal M\!}\,}}\vphantom{\cal M}}}
\def\tiM{{\mathbin{\smash{\ti{\mathcal M}}}\vphantom{\cal M}}}
\def\haM{{\mathbin{\smash{\hat{\mathcal M}}}\vphantom{\cal M}}}
\def\chM{{\mathbin{\smash{\check{\mathcal M}}}\vphantom{\cal M}}}
\def\acM{{\mathbin{\smash{\,\,\ac{\!\!\mathcal M\!}\,}}\vphantom{\cal M}}}
\def\acbM{{\mathbin{\smash{\acute{\bar{\mathcal M}}}}\vphantom{\mathcal M}}}
\def\dM{{\mathbin{\smash{\dot{\mathcal M}}}\vphantom{\cal M}}}
\def\ddM{{\mathbin{\smash{\ddot{\mathcal M}}}\vphantom{\cal M}}}
\def\dddM{{\mathbin{\smash{\dddot{\mathcal M}}}\vphantom{\cal M}|}}
\def\baM{{\mathbin{\smash{\bar{\mathcal M}}}\vphantom{\mathcal M}}}
\def\grM{{\mathbin{\smash{\grave{\cal M}}}\vphantom{\cal M}}}
\def\brM{{\mathbin{\smash{\breve{\cal M}}}\vphantom{\cal M}}}
\def\cV{{\cal V}}
\def\baV{{\mathbin{\smash{\bar{\cal V}}}\vphantom{\cal V}}}
\def\cW{{\cal W}}
\def\sF{{{\mathscr F}}}
\def\sG{{{\mathscr G}}}
\def\sS{{{\mathscr S}}}
\def\baS{\mathbin{\smash{\,\,\,\bar{\!\!\!\mathscr S}}\vphantom{\cal S}}}
\def\acS{\mathbin{\smash{\,\,\,\acute{\!\!\!\mathscr S}}\vphantom{\cal S}}}
\def\grS{\mathbin{\smash{\,\,\,\grave{\!\!\!\mathscr S}}\vphantom{\cal S}}}
\def\baB{{\mathbin{\smash{\bar{\cal B}}}\vphantom{\cal B}}}
\def\sT{{{\mathscr T}}}
\def\sW{{{\mathscr W}}}
\def\baW{\mathbin{\smash{\,\,\,\bar{\!\!\!\mathscr W}}\vphantom{\cal W}}}
\def\g{{\mathfrak g}}
\def\h{{\mathfrak h}}
\def\m{{\mathfrak m}}
\def\u{{\mathfrak u}}
\def\so{{\mathfrak{so}}}
\def\su{{\mathfrak{su}}}
\def\sp{{\mathfrak{sp}}}
\def\gl{{\mathfrak{gl}}}
\def\pgl{{\mathfrak{pgl}}}
\def\fW{{\mathfrak W}}
\def\fX{{\mathfrak X}}
\def\fY{{\mathfrak Y}}
\def\fZ{{\mathfrak Z}}
\def\bM{{\bs M}}
\def\bN{{\bs N}}
\def\bO{{\bs O}}
\def\bQ{{\bs Q}}
\def\bS{{\bs S}}
\def\bU{{\bs U}}
\def\bV{{\bs V}}
\def\bW{{\bs W}\kern -0.1em}
\def\bX{{\bs X}}
\def\bY{{\bs Y}\kern -0.1em}
\def\bZ{{\bs Z}}
\def\al{\alpha}
\def\be{\beta}
\def\ga{\gamma}
\def\de{\delta}
\def\io{\iota}
\def\ep{\epsilon}
\def\la{\lambda}
\def\ka{\kappa}
\def\th{\theta}
\def\ze{\zeta}
\def\up{\upsilon}
\def\vp{\varphi}
\def\si{\sigma}
\def\om{\omega}
\def\De{\Delta}
\def\Ka{{\rm K}}
\def\La{\Lambda}
\def\Om{\Omega}
\def\Ga{\Gamma}
\def\Si{\Sigma}
\def\Th{\Theta}
\def\Up{\Upsilon}
\def\Chi{{\rm X}}
\def\Tau{{\rm T}}
\def\Nu{{\rm N}}
\def\pd{\partial}
\def\ts{\textstyle}
\def\st{\scriptstyle}
\def\sst{\scriptscriptstyle}
\def\w{\wedge}
\def\sm{\setminus}
\def\lt{\ltimes}
\def\bu{\bullet}
\def\sh{\sharp}
\def\di{\diamond}
\def\he{\heartsuit}
\def\od{\odot}
\def\op{\oplus}
\def\ot{\otimes}
\def\bt{\boxtimes}
\def\bp{\boxplus}
\def\ov{\overline}
\def\bigop{\bigoplus}
\def\bigot{\bigotimes}
\def\iy{\infty}
\def\es{\emptyset}
\def\ra{\rightarrow}
\def\rra{\rightrightarrows}
\def\Ra{\Rightarrow}
\def\Longra{\Longrightarrow}
\def\ab{\allowbreak}
\def\longra{\longrightarrow}
\def\hookra{\hookrightarrow}
\def\dashra{\dashrightarrow}
\def\lb{\llbracket}
\def\rb{\rrbracket}
\def\ha{{\ts\frac{1}{2}}}
\def\t{\times}
\def\ci{\circ}
\def\ti{\tilde}
\def\ac{\acute}
\def\gr{\grave}
\def\br{\breve}
\def\d{{\rm d}}
\def\md#1{\vert #1 \vert}
\def\ms#1{\vert #1 \vert^2}
\def\bmd#1{\big\vert #1 \big\vert}
\def\bms#1{\big\vert #1 \big\vert^2}
\def\an#1{\langle #1 \rangle}
\def\ban#1{\bigl\langle #1 \bigr\rangle}
\def\nm#1{\Vert #1 \Vert}
\def\bnm#1{\bigl\Vert #1 \bigr\Vert}
\title{Enumerative invariants and wall-crossing formulae in abelian categories}
\author{Dominic Joyce}
\date{Preliminary version. Comments welcome}
\maketitle

\begin{abstract}
An {\it enumerative invariant theory\/} in Algebraic Geometry is the study of invariants which `count' $\tau$-(semi)stable objects $E$ with fixed topological invariants $\lb E\rb=\al$ in some geometric problem, by means of a {\it virtual class\/} $[\M_\al^\ss(\tau)]_\virt$ in some homology theory, for the moduli spaces $\M_\al^\rst(\tau)\subseteq\M_\al^\ss(\tau)$ of $\tau$-(semi)stable objects. We can obtain numbers by taking integrals $\int_{[\M_\al^\ss(\tau)]_\virt}\Up$ for suitable universal cohomology classes~$\Up$.

Examples include Mochizuki's invariants for coherent sheaves on surfaces \cite{Moch}, and Donaldson--Thomas type invariants for coherent sheaves on Calabi--Yau 3- and 4-folds and Fano 3-folds,~\cite{BoJo,JoSo,KoSo,OhTh,Thom1}.

Let $\A$ be a $\C$-linear abelian category coming from Algebraic Geometry. There are two moduli stacks of objects $E$ in $\A$: the usual moduli stack $\M$, and the `projective linear' moduli stack $\M^\pl$ modulo projective linear isomorphisms, that is, we quotient out by $\la\,\id_E:E\ra E$ for $\la\in\bG_m$. Both are Artin $\C$-stacks. Previous work by the author \cite{Joyc12} gives $H_*(\M)$ the structure of a {\it graded vertex algebra}, and $H_*(\M^\pl)$ a {\it graded Lie algebra}, closely related to $H_*(\M)$. Virtual classes $[\M_\al^\ss(\tau)]_\virt$ lie in~$H_*(\M^\pl)$.

We develop a universal theory of enumerative invariants in such categories $\A$, which includes and extends many cases of interest. Virtual classes $[\M_\al^\ss(\tau)]_\virt$ are only defined when $\M_\al^\rst(\tau)=\M_\al^\ss(\tau)$. We give a systematic way to define invariants $[\M_\al^\ss(\tau)]_\inv$ in $H_*(\M^\pl)$ for all classes $\al\in C(\A)$, with $[\M_\al^\ss(\tau)]_\inv=[\M_\al^\ss(\tau)]_\virt$ when $\M_\al^\rst(\tau)=\M_\al^\ss(\tau)$. If $(\tau,T,\le)$ and $(\ti\tau,\ti T,\le)$ are two suitable (weak) stability conditions on $\A$, we prove a {\it wall-crossing formula\/} which expresses $[\M_\al^\ss(\ti\tau)]_\inv$ in terms of the $[\M_\be^\ss(\tau)]_\inv$, using the Lie bracket on $H_*(\M^\pl)$.

We apply our results when $\A$ is a category $\modCQ$ or $\modCQI$ of representations of a quiver $Q$ or quiver with relations $(Q,I)$, and when $\A=\coh(X)$ for $X$ a curve, surface, or Fano 3-fold, and when $\A$ is a category of `pairs' $\rho:V\ot_\C L\ra E$ in $\coh(X)$ for $X$ a curve or surface, where $V$ is a vector space, $L\ra X$ is a fixed line bundle, and~$E\in\coh(X)$.

We also speculate on extensions of our theory to 3-Calabi--Yau categories, which would give an alternative approach to Donaldson--Thomas theory to \cite{JoSo,Thom1}, and to 4-Calabi--Yau categories, which would give a theory of Donaldson--Thomas type invariants of Calabi--Yau 4-folds.

Our results prove conjectures made in Gross--Joyce--Tanaka~\cite{GJT}.
\end{abstract}

\setcounter{tocdepth}{2}
\tableofcontents

\section{Introduction}
\label{co1}

Let $\A$ be a $\C$-linear abelian category coming from Algebraic Geometry or Representation Theory. Some examples to have in mind are the category $\coh(X)$ of coherent sheaves on a smooth projective $\C$-scheme $X$, and the category $\modCQ$ (or $\modCQI$) of representations of a quiver $Q$ (or quiver with relations $(Q,I)$). 

Fix a quotient $K_0(\A)\twoheadrightarrow K(\A)$, such that the class $\lb E\rb\in K(\A)$ of an object $E\in\A$ measures the `topological invariants' of $E$, for example for $\A=\coh(X)$ we take $\lb E\rb=\ch(E)$ to be the Chern character, and for $\A=\modCQ$ we take $\lb E\rb=\bdim E$ to be the dimension vector. If $\A$ is well-behaved we can form the moduli stack $\M$ of objects $E$ in $\A$, and the `projective linear' moduli stack $\M^\pl$ of objects $E$ in $\A$ modulo `projective linear' isomorphisms, that is, we quotient out by $\la\,\id_E:E\ra E$ for $\la\in\bG_m$. Both are Artin $\C$-stacks. There is a morphism $\Pi^\pl:\M\ra\M^\pl$ which is a principal $[*/\bG_m]$-bundle over $\M\sm\{0\}$. We have decompositions $\M=\coprod_{\al\in K(\A)}\M_\al$, $\M^\pl=\coprod_{\al\in K(\A)}\M_\al^\pl$ into the substacks $\M_\al,\M_\al^\pl$ of $E\in\A$ with~$\lb E\rb=\al$.

Given some extra data on $\M$, for which there are natural choices in the examples above, work of the author \cite{Joyc12} explained in Chapter \ref{co4} defines the structure of a {\it graded vertex algebra\/} on the Betti $\Q$-homology $H_*(\M)=H_*(\M,\Q)$ (with a shifted grading). This is a large, complicated algebraic structure. Given a graded vertex algebra $V_*$, Borcherds \cite[\S 4]{Borc} defined a graded Lie algebra structure on $V_{*+2}/D(V_*)$. It turns out that the projection $H_*(\Pi^\pl):H_*(\M)\ra H_*(\M^\pl)$ identifies $H_*(\M^\pl)$ with $H_*(\M)/D(H_*(\M))$, so $H_*(\M^\pl)$ becomes a graded Lie algebra (with a shifted grading). The Lie bracket $[\,,\,]$ on $H_*(\M^\pl)$ is not at all obvious, and difficult to write down explicitly.

Now let $(\tau,T,\le)$ be a (weak) stability condition on $\A$, e.g.\ Gieseker stability on $\coh(X)$, or slope stability on $\modCQ$. Then we can form open substacks $\M_\al^\rst(\tau)\subseteq\M_\al^\ss(\tau)\subseteq\M_\al^\pl$ of $\tau$-(semi)stable objects $E\in\A$ with $\lb E\rb=\al$. If $\M_\al^\rst(\tau)=\M_\al^\ss(\tau)$ (that is, if there are no strictly $\tau$-semistable objects in class $\al$) then often $\M_\al^\ss(\tau)$ is not just an Artin stack, but a projective $\C$-scheme.

The story so far works in great generality, e.g.\ when $\A=\coh(X)$ for any smooth projective $\C$-scheme $X$. To go further we must make a strong assumption on $\A$, roughly that $\dim\A\le 2$ (that is, the Ext groups $\Ext_\A^i(E,F)$ are zero for $E,F\in\A$ and $i>2$), or that $\dim\A=3$ and $\A$ satisfies a Fano type condition. This holds when $\A=\coh(X)$ for $X$ a curve, surface, or Fano 3-fold. 

Under this assumption, the semistable moduli stacks $\M_\al^\ss(\tau)$ have {\it perfect obstruction theories\/} in the sense of Behrend--Fantechi \cite{BeFa1}. If $\M_\al^\rst(\tau)=\M_\al^\ss(\tau)$, which under good conditions implies that $\M_\al^\ss(\tau)$ is a proper algebraic space, we get a {\it virtual class\/} $[\M_\al^\ss(\tau)]_\virt$, which we regard as an element of $H_*(\M_\al^\pl)$. Thus, if we can construct suitable universal cohomology classes $\Up$ in $H^*(\M_\al^\pl)$ we can define numbers $I_\al^\Up(\tau)=\int_{[\M_\al^\ss(\tau)]_\virt}\Up\in\Q$.

Such $I_\al^\Up(\tau)$ are {\it enumerative invariants\/} of $\A,(\tau,T,\le),\al$. They are often invariant under deformations of the underlying geometry, and have other interesting properties. Examples include Donaldson-type invariants counting coherent sheaves on surfaces \cite{FrQi,Gott,GNY1,GoYu,MPT,MaTh,Moch} and Donaldson--Thomas invariants of Calabi--Yau 3-folds \cite{FeTh,JoSo,KiLi1,KoSo,MNOP,MOOP,Thom1}. Our point of view will be to regard the virtual class $[\M_\al^\ss(\tau)]_\inv:=[\M_\al^\ss(\tau)]_\virt\in H^*(\M_\al^\pl)$ as the primary enumerative invariant, from which we can compute integrals~$I_\al^\Up(\tau)$.

One of the main goals of this book is to answer the following questions:

\begin{quest}
\label{co1quest1}
{\bf(a)} How should one define enumerative invariants $[\M_\al^\ss(\tau)]_\inv$ when $\M_\al^\rst(\tau)\ne\M_\al^\ss(\tau),$ so that\/ $\M_\al^\ss(\tau)$ may not be a proper Deligne--Mumford stack, and the Behrend--Fantechi virtual class $[\M_\al^\ss(\tau)]_\virt$ is not defined?
\smallskip

\noindent{\bf(b)} If\/ $(\tau,T,\le),(\ti\tau,\ti T,\le)$ are (weak) stability conditions on $\A,$ can we write $[\M_\al^\ss(\ti\tau)]_\inv$ in terms of the $[\M_\be^\ss(\tau)]_\inv$? Such an expression is often called a \begin{bfseries}wall-crossing formula\end{bfseries}, as it measures how $[\M_\al^\ss(\tau)]_\inv$ changes as we deform $\tau$ so that it crosses `walls' in the space of stability conditions on $\A$.
\end{quest}

It is necessary to answer {\bf(a)} before answering {\bf(b)}, as even if $\M_\al^\rst(\ti\tau)=\M_\al^\ss(\ti\tau)$, any expression for $[\M_\al^\ss(\ti\tau)]_\inv$ in terms of the $[\M_\be^\ss(\tau)]_\inv$ is likely to involve $[\M_\be^\ss(\tau)]_\inv$ for which~$\M_\be^\rst(\tau)\ne\M_\be^\ss(\tau)$.

Previous work of the author \cite{Joyc2,Joyc3,Joyc4,Joyc5,Joyc6,Joyc7} developed a general theory answering the analogue of Question \ref{co1quest1} for certain kinds of {\it motivic invariants\/} in Algebraic Geometry. The idea is to define the motivic invariants $[\M_\al^\ss(\tau)]_{\rm mot}$ as elements of a Ringel--Hall-type algebra or Lie algebra $\SF(\M)$ of `stack functions' on $\M$ \cite{Joyc3}, and show they satisfy a universal wall-crossing formula written in terms of the Lie bracket on $\SF(\M)$. After Behrend \cite{Behr2} showed that Donaldson--Thomas invariants of Calabi--Yau 3-folds are in fact motivic, Joyce--Song \cite{JoSo} used \cite{Joyc2,Joyc3,Joyc4,Joyc5,Joyc6,Joyc7} to answer Question \ref{co1quest1} for Donaldson--Thomas invariants.

Since Behrend--Fantechi virtual classes $[\M_\al^\ss(\tau)]_\virt$ are not motivic, we cannot apply \cite{Joyc2,Joyc3,Joyc4,Joyc5,Joyc6,Joyc7,JoSo} directly. However, our answer to Question \ref{co1quest1} has the same general structure: we replace the Lie algebra $(\SF(\M),[\,,\,])$ in \cite{Joyc2,Joyc3,Joyc4,Joyc5,Joyc6,Joyc7,JoSo} by the Lie algebra $(H_*(\M^\pl),[\,,\,])$ constructed in \cite{Joyc12}, and then the wall-crossing formula in Theorem \ref{co5thm3} below solving Question \ref{co1quest1}(b) is the same universal wall-crossing formula from \cite[Th.~5.2]{Joyc7}. The proofs are very different to~\cite{Joyc7}.

We begin in Chapter \ref{co2} with background material on Artin stacks and their (co)homology, and Behrend--Fantechi obstruction theories \cite{BeFa1}. Chapter \ref{co3} discusses (weak) stability conditions $(\tau,T,\le)$ on abelian categories $\A$, and Chapter \ref{co4} summarizes results on vertex algebra and Lie algebra structures on homology of moduli stacks from~\cite{Joyc12}.

Chapter \ref{co5} states our main results, Theorems \ref{co5thm1}--\ref{co5thm3} and \ref{co5thm4}, in an abstract form. The format we have chosen is to first state a long list of conditions, given in Assumptions \ref{co5ass1}--\ref{co5ass3}, on a $\C$-linear abelian category $\A$, plus a lot of extra data, including a family $\sS$ of weak stability conditions $(\tau,T,\le)$ on $\A$. Theorems \ref{co5thm1}--\ref{co5thm3} say that if Assumptions \ref{co5ass1}--\ref{co5ass3} hold then we can define enumerative invariants $[\M_\al^\ss(\tau)]_\inv\in H^*(\M_\al^\pl)$ `counting' $\tau$-semistable objects in $\A$ in class $\al\in K(\A)$, with $[\M_\al^\ss(\tau)]_\inv=[\M_\al^\ss(\tau)]_\virt$ when the Behrend--Fantechi virtual class is defined, and if $(\tau,T,\le)$ and $(\ti\tau,\ti T,\le)$ lie in $\sS$ then we give wall-crossing formulae \eq{co5eq33}--\eq{co5eq34} which write $[\M_\al^\ss(\ti\tau)]_\inv$ in terms of the $[\M_\be^\ss(\tau)]_\inv$, using the Lie bracket on $H_*(\M^\pl)$ defined in Chapter \ref{co4}. Theorem \ref{co5thm4} extends the results to equivariant homology $H_*^G(\M^\pl)$ when $G$ acts on~$\A,\M,\M^\pl$.

Our main results prove conjectures in Gross--Joyce--Tanaka \cite[Conj.~2]{GJT} for the case of Behrend--Fantechi virtual classes in Algebraic Geometry.

Chapters \ref{co6}--\ref{co8} apply our theory to particular classes of examples. To do this we have to define the required extra data on $\A,\M,\M^\pl$ and prove that Assumptions \ref{co5ass1}--\ref{co5ass3} hold, which is quite a lot of work in some cases. Then Theorems \ref{co5thm1}--\ref{co5thm3} and \ref{co5thm4} yield results about invariants in particular situations, such as Theorems \ref{co7thm7}--\ref{co7thm8} on $\A=\coh(X)$ with $X$ a projective surface.

Chapter \ref{co6} discusses the simple case when $\A=\modCQ$ or $\modCQI$ is a category of representations of a quiver $Q$ or quiver with relations $(Q,I)$, and $(\tau,T,\le)$ is a slope stability condition on $\A$. To ensure that $\M_\al^\ss(\tau)$ is proper when $\M_\al^\rst(\tau)=\M_\al^\ss(\tau)$ we take $Q$ to have no oriented cycles. When $\A=\modCQ$, our results were proved in \cite[Th.~5.8]{GJT} by a different method.

Chapter \ref{co7} considers the case $\A=\coh(X)$ when $X$ is a smooth projective $\C$-scheme of dimension $m$. Parts of the theory work for all $m\ge 0$, but the perfect obstruction theories on $\M_\al^\ss(\tau)$ can only be defined if $m\le 2$, or $m=3$ and $X$ satisfies a Fano type condition, so we define invariants $[\M_\al^\ss(\tau)]_\inv$ only when $X$ is a curve, surface, or Fano 3-fold.

Let $X$ be a smooth projective $\C$-scheme and $L\ra X$ a line bundle. There is an extensive literature (e.g.\ \cite{Brad,BDGW,DKO,Garc,HuLe1,JoSo,KoTh3,Moch,OST,PaTh1,PaTh2,Thad}) concerning moduli spaces $\acM_{(\al,1)}^\ss(\ac\mu)$ of `stable pairs' $(E,\rho)$ where $E\in\coh(X)$, and $\rho:L\ra E$ is a morphism, and $(E,\rho)$ satisfies a stability condition $\ac\mu$, and studying associated virtual classes $[\acM_{(\al,1)}^\ss(\ac\mu)]_\virt$ and invariants. 

To include this in our theory, in Chapter \ref{co8} we take $\acA$ to be an abelian category of triples $(E,V,\rho)$ for $E\in\coh(X)$, $V\in\Vect_\C$ and $\rho:V\ot_\C L\ra E$ a morphism in $\coh(X)$. To define obstruction theories on $\acM_{(\al,d)}^\ss(\ac\mu)$ we take $X$ to be a curve or surface. Our theorems give a system of invariants $[\acM_{(\al,d)}^\ss(\ac\mu)]_\inv$ counting $\ac\mu$-semistable triples $(E,V,\rho)$ for $\lb E\rb=\al\in K(\coh(X))$ and $\dim V=d$, which include the invariants $[\M_\al^\ss(\mu)]_\inv$ from Chapter \ref{co7} when $d=0$, and invariants counting stable pairs when~$d=1$.

When $X$ is a projective surface, the invariants $[\M_\al^\ss(\mu)]_\inv$ from Chapter \ref{co7} with $\rank\al>0$ are algebraic versions of Donaldson invariants of the 4-manifold $X$, and invariants $[\acM_{(\al,1)}^\ss(\ac\mu)]_\inv$ with $\rank\al=1$ and $c_2(\al)=0$ are algebraic versions of Seiberg--Witten invariants of $X$. Section \ref{co86} explains a method to use our wall-crossing formula for the $[\acM_{(\al,d)}^\ss(\ac\mu)]_\inv$ to compute the $[\M_\al^\ss(\mu)]_\inv$ for all $\rank\al>0$ in terms of the $[\acM_{(\al,d)}^\ss(\ac\mu)]_\inv$ for $\rank\al=1$ and $d=0$ or 1. This is an algebraic analogue of writing Donaldson invariants in terms of Seiberg--Witten invariants, as in \cite{FeLe1,GNY3,KrMr,MaMo,MoWi,Witt}, for example.

The proofs of Theorems \ref{co5thm1}--\ref{co5thm3} and \ref{co5thm4} are deferred to Chapters \ref{co9}--\ref{co11}.

In the rest of this introduction, \S\ref{co1.1} defines (weak) stability conditions on abelian categories as in Chapter \ref{co3}, section \ref{co1.2} summarizes parts of Chapter \ref{co4} on vertex and Lie algebra structures on homology of moduli stacks, \S\ref{co1.3} explains the main results of Chapter \ref{co5}, and \S\ref{co1.4} discusses coherent sheaves on projective surfaces, as an example application from Chapters~\ref{co6}--\ref{co8}.

This is only a preliminary version of this book, and the final version will contain more material. This version works with virtual classes for {\it Behrend--Fantechi obstruction theories\/} \cite{BeFa1}. Section \ref{co1.5} discusses extensions of our theory to 3-Calabi--Yau and 4-Calabi--Yau categories $\A$, that the author hopes will be included in the final version of this book.
\smallskip

\noindent {\it Acknowledgements.} This research was supported by the Simons Collaboration on Special Holonomy in Geometry, Analysis and Physics. The author would like to thank Arkadij Bojko, Tom Bridgeland, Paul Feehan, Jacob Gross, Daniel Halpern-Leistner, Martijn Kool, Alyosha Latyntsev, Tom Leness, Henry Liu, Hyeonjun Park, Julius Ross, Yuuji Tanaka, Richard Thomas, and Markus Upmeier for helpful conversations.

This project borrows ideas from Mochizuki's monograph \cite{Moch} on invariants and wall-crossing for coherent sheaves on surfaces at several important points, and the author would like to thank Takuro Mochizuki for these.

\subsection{(Weak) stability conditions on abelian categories}
\label{co1.1}

We briefly review some material from \S\ref{co31}, \S\ref{co63} and \S\ref{co71}.

\begin{dfn}
\label{co1def1}
Let $\A$ be an abelian category. The {\it Grothendieck group\/} $K_0(\A)$ is the abelian group generated by isomorphism classes $[E]$ of objects $E\in\A$, with a relation $[E_2]=[E_1]+[E_3]$ for each exact sequence $0\ra E_1\ra E_2\ra E_3\ra 0$ in $\A$. Suppose we are given a surjective quotient $K_0(\A)\twoheadrightarrow K(\A)$. We write $\lb E\rb\in K(\A)$ for the class of~$E\in\A$.

Suppose that $0\in\A$ is the only object in class $0\in K(\A)$. Define the {\it positive cone\/} $C(\A)\subset K(\A)\sm\{0\}$ by~$C(\A)=\bigl\{\lb E\rb:0\ne E\in\A\bigr\}$.

Let $(T,\leq)$ be a totally ordered set and $\tau:C(\A)\ra T$ be a map. We call $(\tau,T,\leq)$ a {\it weak stability condition\/} on $\A$ if for all $\al,\be,\ga \in C(\A)$ with $\be=\al+\ga$, either 
$\tau(\al) \leq \tau(\be) \leq \tau(\ga)$, or~$\tau(\al) \geq \tau(\be) \geq \tau(\ga)$.

We call $(\tau,T,\leq)$ a {\it stability condition\/} if for all such $\al,\be,\ga$, either 
$\tau(\al)<\tau(\be)<\tau(\ga)$, or $\tau(\al)>\tau(\be)>\tau(\ga)$, or~$\tau(\al)=\tau(\be)=\tau(\ga)$.

Let $(\tau,T,\leq)$ be a weak stability condition. An object $E$ of $\A$ is called:
\begin{itemize}
\setlength{\itemsep}{0pt}
\setlength{\parsep}{0pt}
\item[(i)] {\it $\tau$-stable\/} if $\tau(\lb E'\rb)<\tau(\lb E/E'\rb)$ for all subobjects $E'\subset E$ with $E'\ne 0,E$.
\item[(ii)] {\it $\tau$-semistable\/} if $\tau(\lb E'\rb)\!\leq\!\tau(\lb E/E'\rb)$ for all $E'\subset E$ with $E'\ne 0,E$.
\item[(iii)] {\it strictly $\tau$-semistable\/} if it is $\tau$-semistable but not $\tau$-stable.
\end{itemize}
\end{dfn}

\begin{ex}
\label{co1ex1}
Let $Q=(Q_0,Q_1,h,t)$ be a quiver, so that $Q_0,Q_1$ are finite sets of vertices and edges, and $h,t:Q_1\ra Q_0$ are the head and tail maps. Let $\A=\modCQ$ be the abelian category of $\C$-representations of $Q$, as in \S\ref{co61} below. Write objects of $\A$ as $(V,\rho)=((V_v)_{v\in Q_0},(\rho_e)_{e\in Q_1})$, where $V_v$ is a finite-dimensional $\C$-vector space and $\rho_e:V_{t(e)}\ra V_{h(e)}$ a linear map. Define the {\it dimension vector\/} of $(V,\rho)$ to be $\bs d=\bdim(V,\rho)\in\N^{Q_0}\subset\Z^{Q_0}$, where $\bs d(v)=\dim_\C V_v$ for $v\in Q_0$. Define $K(\A)=\Z^{Q_0}$ to be the lattice of dimension vectors, with $\lb V,\rho\rb=\bdim(V,\rho)$. Then the positive cone is~$C(\A)=\N^{Q_0}\sm\{0\}$.

Choose real numbers $\mu_v\in\R$ for all $v\in Q_0$. Define a map $\mu:C(\A)\ra\R$ by
\begin{equation*}
\mu(\bs d)=\frac{\sum_{v\in Q_0}\bs d(v)\mu_v}{\sum_{v\in Q_0}\bs d(v)}.
\end{equation*}
Then $(\mu,\R,\le)$ is a stability condition on $\modCQ$, called {\it slope stability}.
\end{ex}

\begin{ex}
\label{co1ex2}
Let $X$ be a smooth projective $\C$-scheme of dimension $m$. Write $\coh(X)$ for the abelian category of coherent sheaves on $X$, as in \cite[\S II.5]{Hart} and \cite{HuLe1}. As in \cite[\S III.6]{Hart}, for $E,F$ in $\coh(X)$ we have {\it Ext groups\/} $\Ext^i(E,F)$ for $i=0,\ldots,m$, finite-dimensional $\C$-vector spaces, with $\Ext^0(E,F)=\Hom(E,F)$. The {\it Euler form\/} is the biadditive map $\chi:K_0(\coh(X))\t K_0(\coh(X))\ra\Z$ with
\e
\chi([E],[F])=\sum_{i\ge 0}(-1)^i\dim_\C\Ext^i(E,F)
\label{co1eq1}
\e
for all $E,F\in\coh(X)$. The {\it numerical Grothendieck group\/} is 
\begin{align*}
&K^\num(\coh(X))=K_0(\coh(X))/\Ker\chi, \quad\text{where}\\
&\Ker\chi=\bigl\{\al\in K_0(\coh(X)):\text{$\chi(\al,\be)=0$ for all $\be\in K_0(\coh(X))$}\bigr\}.
\end{align*}
Then $\chi$ descends to $\chi:K^\num(\coh(X))\t K^\num(\coh(X))\ra\Z$. Equivalently, $K^\num(\coh(X))\cong\Im\bigl(\ch:K_0(\coh(X))\ra H^{\rm even}(X,\Q)\bigr)$, where $\ch$ is the {\it Chern character}, as in \cite[App.~A]{Hart}. We will take $K(\coh(X))=K^\num(\coh(X))$ in Definition \ref{co1def1}. We write $\lb E\rb\in K(\coh(X))$ for the class of~$E\in\coh(X)$.

Since $X$ is a compact complex manifold admitting K\"ahler metrics, we can use Hodge theory. The second complex de Rham cohomology group splits~as
\begin{equation*}
H^2(X,\C)=H^{2,0}(X)\op H^{1,1}(X)\op H^{0,2}(X).
\end{equation*}
Write $H^{1,1}(X,\R)=H^{1,1}(X)\cap H^2(X,\R)$ for the vector subspace of real $(1,1)$-classes in $H^2(X,\R)$. Write $\Kah(X)\subset H^{1,1}(X,\R)$ for the {\it K\"ahler cone\/} of K\"ahler classes $\om$ on $X$. Then $\Kah(X)$ is an open convex cone in~$H^{1,1}(X,\R)$.

Let $\om\in\Kah(X)$ and $\al\in K(\coh(X))$. The {\it Hilbert polynomial of\/ $\al$ with respect to\/} $\om$ is
\e
P_\al^\om(t)=\int_X\ch(\al)\cup e^{t\om}\cup \td(X).
\label{co1eq2}
\e
Here $e^{t\om}=\sum_{i=0}^mt^i\om^i/i!$, since $\om^i=0$ if $i>m$, so $P_\al^\om(t)$ is a real polynomial of degree at most $m$, and $\td(X)\in H^{\rm even}(X,\Q)$ is the Todd class of $X$. If $\al\in C(\coh(X))$ then $P_\al^\om$ is nonzero and has positive leading coefficient. The degree $\deg P_\al^\om=\dim\al$ is independent of $\om$.

Define $G$ to be the set of monic real polynomials in $t$ of degree at most $m$:
\begin{equation*}
G=\bigl\{p(t)=t^d+a_{d-1}t^{d-1}+\cdots+a_0:d=0,\ldots,m,\;
a_0,\ldots,a_{d-1}\in\R\bigr\}.
\end{equation*}
Define a total order `$\le$' on $G$ by $p\le p'$ for $p,p'\in G$ if either
\begin{itemize}
\setlength{\itemsep}{0pt}
\setlength{\parsep}{0pt}
\item[(a)] $\deg p>\deg p'$, or
\item[(b)] $\deg p=\deg p'$ and $p(t)\le p'(t)$ for all $t\gg 0$.
\end{itemize}

Fix $\om\in\Kah(X)$, and define $\tau^\om:C(\coh(X))\ra G$ by $\tau^\om(\al)=P^\om_\al/r^\om_\al$, where $r^\om_\al$ is the leading coefficient of $P^\om_\al$, which must be positive as above. Then $(\tau^\om,G,\le)$ is a stability condition on $\coh(X)$, which we call {\it Gieseker stability\/} for the real K\"ahler class $\om$.
\end{ex}

\begin{ex}
\label{co1ex3}
Continuing in the situation of Example \ref{co1ex2}, define
\begin{equation*}
M=\bigl\{p(t)=t^d+a_{d-1}t^{d-1}:d=0,\ldots,m,\;\> a_{d-1}\in\R,\;\>
a_{-1}=0\bigr\}\subset G
\end{equation*}
and restrict the total order $\le$ on $G$ to $M$. Define $\mu^\om:C(\coh(X))\ra M$ by $\mu(\al)=t^d+a_{d-1}t^{d-1}$ when $\tau^\om(\al)=P_\al/r_\al=t^d+a_{d-1}t^{d-1}+\cdots+a_0$, that is, $\mu^\om(\al)$ is the truncation of $\tau^\om(\al)$ in Example \ref{co1ex2} at its second term. Then $(\mu^\om,M,\le)$ is a weak stability condition on $\coh(X)$, called $\mu$-{\it stability}.
\end{ex}

\subsection{Vertex and Lie algebras on homology of moduli stacks}
\label{co1.2}

We now summarize parts of Chapter \ref{co4}. Here is a version of Assumption \ref{co4ass1}, which we have simplified by supposing $\chi$ is nondegenerate, enabling us to omit Assumption~\ref{co4ass1}(g).

\begin{ass}
\label{co1ass1}
Let $\A$ be a $\C$-linear additive category. Assume:
\begin{itemize}
\setlength{\itemsep}{0pt}
\setlength{\parsep}{0pt}
\item[(a)] We can form a natural moduli stack $\M$ of objects in $\A$, an Artin $\C$-stack, locally of finite type. Then $\C$-points of $\M$ are isomorphism classes $[E]$ of objects $E\in\A$, and the isotropy groups are $\Iso_\M([E])=\Aut(E)$.
\item[(b)] There is a morphism of Artin stacks $\Phi:\M\t\M\ra\M$ which on $\C$-points acts by $\Phi_*:([E],[F])\mapsto[E\op F]$, for all objects $E,F\in\A$. That is, $\Phi$ is the morphism of moduli stacks induced by direct sum in $\A$. It is associative and commutative.
\item[(c)] There is a morphism of Artin stacks $\Psi:[*/\bG_m]\t\M\ra\M$ which on $\C$-points acts by $\Psi_*:(*,[E])\mapsto[E]$, for all objects $E$ in $\A$, and on isotropy groups acts by $\Psi_*:\Iso_{[*/\bG_m]\t\M}(*,[E])\cong\bG_m\t\Aut(E)\ra \Iso_\M([E])\cong\Aut(E)$ by $(\la,\mu)\mapsto \la\mu=(\la\cdot\id_E)\ci\mu$ for $\la\in\bG_m$ and $\mu\in\Aut(E)$. This $\Psi$ is an action of the group stack $[*/\bG_m]$ on $\M$, and is distributive over~$\Phi$.
\item[(d)] We are given a surjective quotient $K_0(\A)\twoheadrightarrow K(\A)$ of the Grothendieck group $K_0(\A)$ of $\A$. We write $\lb E\rb\in K(\A)$ for the class of $E\in\A$. We suppose that if $E\in\A$ with $\lb E\rb=0$ in $K(\A)$ then $E=0$.

We require that the map $\M(\C)\ra K(\A)$ mapping $E\mapsto\lb E\rb$ should be locally constant. This gives a decomposition $\M=\coprod_{\al\in K(\A)}\M_\al$ of $\M$ into open and closed $\C$-substacks $\M_\al\subset \M$ of objects in class $\al$, where $\M_0=\{[0]\}$. We write $\Phi_{\al,\be}=\Phi\vert_{\M_\al\t\M_\be}:\M_\al\t\M_\be\ra\M_{\al+\be}$ and $\Psi_\al=\Psi\vert_{[*/\bG_m]\t\M_\al}:[*/\bG_m]\t\M_\al\ra\M_\al$. 
\item[(e)] We are given a nondegenerate biadditive form $\chi:K(\A)\t K(\A)\ra \Z$.
\item[(f)] We are given a perfect complex $\cE^\bu$ on $\M\t\M$, such that the restriction $\cE_{\al,\be}^\bu:=\cE^\bu\vert_{\M_\al\t\M_\be}$ to $\M_\al\t\M_\be$ has rank $\chi(\al,\be)$ for all $\al,\be\in K(\A)$, and there are isomorphisms of perfect complexes
\e
\begin{split}
(\Phi\t\id_\M)^*(\cE^\bu)&\cong
\Pi_{13}^*(\cE^\bu)\op \Pi_{23}^*(\cE^\bu),\\
(\id_\M\t\Phi)^*(\cE^\bu)&\cong
\Pi_{12}^*(\cE^\bu)\op \Pi_{13}^*(\cE^\bu),\\
(\Psi\t\id_\M)^*(\cE^\bu)&\cong \Pi_1^*(L_{[*/\bG_m]})\ot \Pi_{23}^*(\cE^\bu),\\
(\Pi_2,\Psi\ci\Pi_{13})^*(\cE^\bu)&\cong \Pi_1^*(L_{[*/\bG_m]}^*)\ot\Pi_{23}^*(\cE^\bu).
\end{split}
\label{co1eq3}
\e
Here $L_{[*/\bG_m]}\ra[*/\bG_m]$ is the line bundle corresponding to the weight 1 representation of $\bG_m=\C\sm\{0\}$ on~$\C$.
\end{itemize}	
\end{ass}

\begin{rem}
\label{co1rem1}
As in \cite{Joyc12} and Chapters \ref{co6}--\ref{co8}, there are natural choices for the data of Assumption \ref{co1ass1} in large classes of interesting examples.

As in Example \ref{co1ex2}, many $\C$-linear categories $\A$ in Algebraic Geometry have well behaved {\it Ext groups\/} $\Ext^i(E,F)$ for $E,F\in\A$ and $i\ge 0$, which are finite-dimensional $\C$-vector spaces with $\Hom(E,F)=\Ext^0(E,F)$ and $\Ext^i(E,F)=0$ for $i>\dim\A$. Long exact sequences for Ext groups imply there is a biadditive {\it Euler form\/} $\chi:K_0(\A)\t K_0(\A)\ra\Z$ satisfying \eq{co1eq1}. Then in Assumption \ref{co1ass1}(d) a natural choice is $K(\A)=K_0(\A)/\Ker\chi$, the {\it numerical Grothendieck group}, so $\chi$ descends to a nondegenerate $\chi:K(\A)\t K(\A)\ra \Z$ as in Assumption \ref{co1ass1}(e). Often there is an {\it Ext complex\/} $\cExt^\bu\ra\M\t\M$, perfect in $[0,\dim\A]$, with $\cExt^\bu\vert_{([E],[F])}=\Ext^*(E,F)$, so that $\rank\cExt^\bu\vert_{\M_\al\t\M_\be}=\chi(\al,\be)$. Then in Assumption \ref{co1ass1}(f) the natural choice is~$\cE^\bu=(\cExt^\bu)^\vee$.
\end{rem}

\begin{dfn}
\label{co1def2}
Suppose Assumption \ref{co1ass1} holds. Given all this data, we define a graded vertex algebra structure on the Betti $\Q$-homology $H_*(\M)=H_*(\M,\Q)$. The inclusion of the zero object $0\in\A$ gives a morphism $[0]:*\hookra\M$ inducing $\Q\cong H_*(*)
\ra H_*(\M)$, and we define $\boo\in H_*(\M)$ to be the image of $1\in \Q$ under this map. Taking homology of $\Psi$ gives a map
\begin{equation*}
\xymatrix@C=30pt{ H_*([*/ \bG_m])\!\ot_\Q\! H_*(\M)  \ar[r]^(0.57){\bt} & H_*([*/ \bG_m] \!\t\! \M) \ar[r]^(0.61){\Psi_*} & H_*(\M). }
\end{equation*}
As $\Hom_\Q(H_*([*/ \bG_m]),\Q)\cong H^*([*/\bG_m])\cong\Q[[z]]$, this is equivalent to a map $e^{zD}:H_*(\M)\ra H_*(\M)[[z]]$.

The decomposition $\M=\coprod_{\al\in K(\A)}\M_\al$ induces an identification
\e
H_*(\M)=\ts\bigop_{\al\in K(\A)}H_*(\M_\al).
\label{co1eq4}
\e
For $u\in H_*(\M_\al)\subset H_*(\M)$ and $v\in H_*(\M_\be)\subset H_*(\M)$, define
\e
\begin{split}
&Y(u,z)v=Y(z)(u\ot v)= (-1)^{\chi(\al,\be)} \sum\nolimits_{i,j\ge 0} z^{\chi(\al, \be)+\chi(\be,\al)-i+j}\cdot{}\\
&\bigl(\Phi_{\al,\be}\ci(\Psi_\al\t\id_{\M_\be})\bigr)_*\bigl(t^j\bt ((u \bt v) \cap c_i(\cE_{\al,\be}^\bu\op \si_{\al,\be}^*(\cE^\bu_{\be,\al})^\vee))\bigr),
\end{split}
\label{co1eq5}
\e
where $t^j$ is the generator of $H_{2j}([*/\bG_m])$, and $\si_{\al,\be}:\M_\al\t\M_\be\ra\M_\be\t\M_\al$ exchanges the factors. We also use the notation $Y(u,z)v=\sum_{n\in\Z}u_n(v)z^{-n-1}$, where $(u,v)\mapsto u_n(v)$ is a bilinear operation on $H_*(\M)$.

Using \eq{co1eq4}, for $n\in\Z$ and $\al\in K(\A)$ we write
\e
\hat H_n(\M_\al)= H_{n-2\chi(\al,\al)}(\M_\al),\quad \hat H_n(\M)=\ts\bigop_{\al\in K(\A)}\hat H_n(\M_\al).
\label{co1eq6}
\e
That is, $\hat H_*(\M)$ is $H_*(\M)$, but with grading shifted by $-2\chi(\al,\al)$ in the component $H_*(\M_\al)\subset H_*(\M)$. The author \cite{Joyc12} proves:	
\end{dfn}

\begin{thm}
\label{co1thm1}
$(\hat H_*(\M), \boo,e^{zD},Y)$ above is a graded vertex algebra over $\Q$.
\end{thm}

Here {\it vertex algebras\/} are defined in \S\ref{co41}. They are very complicated algebraic objects, which arise in Conformal Field Theories in Theoretical Physics. Roughly speaking, a vertex algebra has an infinite number of Lie bracket type operations $(u,v)\mapsto u_n(v)$ for $n\in\Z$, which satisfy an infinite number of Jacobi-type identities. For the purpose of this introduction, we only need to know one fact about vertex algebras, which is that given a (graded) vertex algebra one can construct a (graded) Lie algebra, as in Borcherds~\cite[\S 4]{Borc}:

\begin{prop}
\label{co1prop1}
Let\/ $(V_*,\boo,e^{zD},Y)$ be a graded vertex algebra over $\Q$. We may construct a graded Lie algebra $(\check V_*,[\,,\,])$ over $\Q$ as follows. Noting the shift in grading, define a $\Z$-graded\/ $\Q$-vector space $\check V_*$ by $\check V_n=V_{n+2}/D(V_n)$ for\/ $n\in\Z,$ so that\/ $\check V_*=V_{*+2}/D(V_*)$. If\/ $u\in V_{a+2}$ and\/ $v\in V_{b+2},$ the Lie bracket on $\check V_*$ is
\begin{equation*}
\bigl[u+D(V_a),v+D(V_b)\bigr]=u_0(v)+D(V_{a+b})\in\check V_{a+b}.
\end{equation*}
\end{prop}

\begin{dfn}
\label{co1def3}
Continue in the situation of Assumption \ref{co1ass1} and Definition \ref{co1def2}. Then $[*/\bG_m]$ is a group stack, and $\Psi:[*/\bG_m]\t\M\ra\M$ is an action of $[*/\bG_m]$ on $\M$, which is trivial on $\M_0=\{[0]\}$, and free on $\M\sm\{[0]\}$. We may take the quotient of $\M$ by $\Psi$ to get a stack $\M^\pl$, which we call the {\it projective linear moduli stack}, with projection $\Pi^\pl:\M\ra\M^\pl$, in a 2-co-Cartesian square:
\begin{equation*}
\xymatrix@C=130pt@R=15pt{ *+[r]{[*/\bG_m]\t\M} \drtwocell_{}\omit^{}\omit{^{}} \ar[r]_(0.65){\Psi} \ar[d]^{\pi_\M} & *+[l]{\M} \ar[d]_{\Pi^\pl} \\
*+[r]{\M} \ar[r]^(0.35){\Pi^\pl} & *+[l]{\M^\pl.\!} }
\end{equation*}
Then $\Pi^\pl:\M\ra\M^\pl$ is a principal $[*/\bG_m]$-bundle over $\M\sm\{0\}$. This construction of $\M^\pl$ is known as {\it rigidification}, as in \cite{AOV,Roma}, written $\M^\pl=\M\!\!\fatslash\,\bG_m$.

We regard $\M^\pl$ as the moduli stack of objects in $\A$ `up to projective linear isomorphisms', where a {\it projective linear isomorphism\/} $[\phi]:E\ra F$ is a $\sim$-equivalence class of isomorphisms $\phi:E\ra F$ in $\coh(X)$, with $\phi\sim\phi'$ if $\phi=\la\phi'$ for some $\la\in\bG_m$. The splitting $\M=\coprod_{\al\in K(\A)}\M_\al$ descends to $\M^\pl=\coprod_{\al\in K(\A)}\M_\al^\pl$. Thus as for \eq{co1eq4} we have
\e
H_*(\M^\pl)=\ts\bigop_{\al\in K(\A)}H_*(\M_\al^\pl).
\label{co1eq7}
\e
In a similar way to \eq{co1eq6}, using \eq{co1eq7}, for $n\in\Z$ and $\al\in K(\A)$ we write
\e
\check H_n(\M_\al^\pl)=H_{n+2-2\chi(\al,\al)}(\M_\al^\pl),\quad \check H_n(\M^\pl)=\ts\bigop_{\al\in K(\A)}\check H_n(\M_\al^\pl).
\label{co1eq8}
\e
That is, $\check H_*(\M^\pl)$ is $H_*(\M^\pl)$, but with grading shifted by $2-2\chi(\al,\al)$ in the component $H_*(\M_\al^\pl)\subset H_*(\M^\pl)$. The next theorem, proved in \cite{Joyc12}, gives a geometric interpretation of the graded Lie algebra~$\hat H_{*+2}(\M)/D(\hat H_*(\M))$.
\end{dfn}

\begin{thm}
\label{co1thm2}
Work in the situation of Assumption\/ {\rm\ref{co1ass1}} and Definitions\/ {\rm\ref{co1def2}} and\/ {\rm\ref{co1def3},} and consider the graded Lie algebra $\hat H_{*+2}(\M)/D(\hat H_*(\M))$ constructed by combining Theorem\/ {\rm\ref{co1thm1}} and Proposition\/ {\rm\ref{co1prop1}}. Then $\Pi^\pl:\M\ra\M^\pl$ gives a morphism $\Pi^\pl_*:H_*(\M)\ab\ra H_*(\M^\pl)$. It is surjective, with kernel\/ $D(H_*(\M))$. This induces an isomorphism
\begin{equation*}
\Pi^\pl_*:H_*(\M)/D(H_*(\M))\longra H_*(\M^\pl).
\end{equation*}
Comparing \eq{co1eq6} and\/ {\rm\eq{co1eq8},} we see this is an isomorphism for all\/~{\rm$n\in\Z$:}
\e
\hat H_{n+2}(\M)/D(\hat H_n(\M))\longra \check H_n(\M^\pl).
\label{co1eq9}
\e
Thus there is a unique Lie bracket\/ $[\,,\,]$ on $\check H_*(\M^\pl)$ making it into a \begin{bfseries}graded Lie algebra\end{bfseries}, such that\/ \eq{co1eq9} is a Lie algebra isomorphism. Hence $\check H_0(\M^\pl)$ and\/ $\check H_{\rm even}(\M^\pl)=\bigop_{n\in\Z}\check H_{2n}(\M^\pl)$ are Lie algebras.
\end{thm}

See Remark \ref{co4rem1} below for an alternative way to think about the Lie bracket on~$\check H_*(\M^\pl)$.

\subsection{Summary of the main results}
\label{co1.3}

Our main results Theorems \ref{co5thm1}--\ref{co5thm3} depend on Assumptions \ref{co5ass1}--\ref{co5ass3} in \S\ref{co51}. For this introduction we give simplified, incomplete versions Assumptions \ref{co1ass2}--\ref{co1ass4} of these. For example, \S\ref{co51} deals with an abelian category $\A$ with an exact subcategory $\B\subseteq\A$, but here we take $\B=\A$ and omit $\B$ from the notation.

\begin{ass}
\label{co1ass2}
Let $\A$ be a noetherian $\C$-linear abelian category. Assume:

\smallskip

\noindent{\bf(a)} Assumption \ref{co1ass1} holds for $\A$. We will freely use the notation $\M,\ab\M^\pl,\ab\M_\al,\ab\M^\pl_\al,\ab\cE^\bu,\ab\chi,\ab\ldots$ of Assumption \ref{co1ass1} and Definition \ref{co1def3}, and the Lie algebra $\check H_{\rm even}(\M^\pl)$ of Theorem \ref{co1thm2}.
\smallskip

\noindent{\bf(b)} Using the notation of \S\ref{co22}, there are also {\it derived Artin\/ $\C$-stacks\/} $\bs\M,\bs\M^\pl$, which are the (projective linear) derived moduli stacks of objects in $\A$. They have classical truncations $t_0(\bs\M)=\M$ and $t_0(\bs\M^\pl)=\M^\pl$, for $\M,\M^\pl$ the classical moduli stacks in {\bf(a)}. We write $i_\M:\M\hookra\bs\M$, $i_{\M^\pl}:\M^\pl\hookra\bs\M^\pl$ for the inclusions of classical truncations. All these satisfy:
\begin{itemize}
\setlength{\itemsep}{0pt}
\setlength{\parsep}{0pt}
\item[{\bf(i)}] Since open substacks of $\bs\M,\bs\M^\pl$ correspond to open substacks of their classical truncations $\M,\M^\pl$, the splittings $\M=\coprod_{\al\in K(\A)}\M_\al$, $\M^\pl=\coprod_{\al\in K(\A)}\M_\al^\pl$ lift to $\bs\M=\coprod_{\al\in K(\A)}\bs\M_\al$, $\bs\M^\pl=\coprod_{\al\in K(\A)}\bs\M_\al^\pl$.
\item[{\bf(ii)}] The morphisms $\Phi:\M\t\M\ra\M$, $\Phi_{\al,\be}:\M_\al\t\M_\be\ra\M_{\al+\be}$, $\Psi:[*/\bG_m]\t\M\ra\M$, $\Psi_\al:[*/\bG_m]\t\M_\al\ra\M_\al$, $\Pi^\pl:\M\ra\M^\pl$ and $\Pi_\al^\pl:\M_\al\ra\M_\al^\pl$ in Assumption \ref{co1ass1} and Definition \ref{co1def3} lift to derived versions 
$\bs\Phi:\bs\M\t\bs\M\ra\bs\M$, \ldots, $\bs\Pi_\al^\pl:\bs\M_\al\ra\bs\M_\al^\pl$, which have classical truncations $\Phi=t_0(\bs\Phi)$, \ldots, $\Pi_\al^\pl=t_0(\bs\Pi_\al^\pl)$, and satisfy the same identities as $\Phi,\ldots,\Pi_\al^\pl$. The projection $\bs\Pi^\pl:\bs\M\ra\bs\M^\pl$ is a principal $[*/\bG_m]$-bundle in derived stacks on $\bs\M\sm\{0\},\bs\M^\pl\sm\{0\}$.
\item[{\bf(iii)}] $\bs\M$ and $\bs\M^\pl$ are {\it locally finitely presented}. This implies that their cotangent complexes $\bL_{\bs\M},\bL_{\bs\M^\pl}$ are perfect.
\item[{\bf(iv)}] For each $\al\in K(\A)$ and $\cE_{\al,\al}^\bu$ as in Assumption \ref{co4ass1}(f), we are given a quasi-isomorphism of perfect complexes on $\M_\al$:
\begin{equation*}
\th_\al:\De_{\M_\al}^*(\cE_{\al,\al}^\bu)[-1]\,{\buildrel\cong\over\longra}\,  i_{\M_\al}^*(\bL_{\bs\M_\al}).
\end{equation*}
\item[{\bf(v)}] For $\al,\be\in K(\A)$, the following should commute in $\Perf(\M_\al\t\M_\be)$:
\begin{align*}
\xymatrix@!0@C=265pt@R=50pt{ *+[r]{\begin{subarray}{l}\ts (\De_{\M_{\al+\be}}\ci\Phi_{\al,\be})^* \\ \ts (\cE_{\al+\be,\al+\be}^\bu)[-1]\end{subarray}} 
\ar[d]^{\text{from \eq{co1eq3}}}_\cong
\ar[r]_(0.47){ \Phi_{\al,\be}^*(\th_{\al+\be})  } & *+[l]{\Phi_{\al,\be}^*(i_{\M_{\al+\be}}^*(\bL_{\bs\M_{\al+\be}}))} \ar[dd]_(0.2){i_{\M_\al\t\M_\be}^*(\bL_{\bs\Phi_{\al,\be}})} \\
*+[r]{\begin{subarray}{l}\ts (\De_{\M_\al}\!\ci\!\Pi_{\M_\al})^*(\cE_{\al,\al}^\bu)[-1]\!\op\!  (\De_{\M_\be}\!\ci\!\Pi_{\M_\be})^*(\cE_{\be,\be}^\bu)[-1] \\
\ts {}\op \cE_{\al,\be}^\bu[-1]\op \si_{\al,\be}^*(\cE_{\be,\al}^\bu)[-1] 
\end{subarray}}
\ar[d]^(0.47){\text{project to first two factors}}
\\
*+[r]{\begin{subarray}{l}\ts (\De_{\M_\al}\!\ci\!\Pi_{\M_\al})^*(\cE_{\al,\al}^\bu)[-1]\\
\ts {}\op (\De_{\M_\be}\ci\Pi_{\M_\be})^*(\cE_{\be,\be}^\bu)[-1]
\end{subarray}} \ar[r]^(0.53){\begin{subarray}{l} \Pi_{\M_\al}^*(\th_\al)\op \\ \Pi_{\M_\be}^*(\th_\be) \end{subarray}} & *+[l]{\raisebox{30pt}{$\begin{subarray}{l}\ts i_{\M_\al\t\M_\be}^*(\bL_{\bs\M_\al\t\bs\M_\be}) \\ 
\ts {}\cong\Pi_{\M_\al}^*\ci i_{\M_\al}^*(\bL_{\bs\M_\al}) \\ \ts {}\op\Pi_{\M_\be}^*\ci i_{\M_\be}^*(\bL_{\bs\M_\be}).
\end{subarray}$}\quad\;\>}   }
\\[-35pt]
\end{align*}
\end{itemize}

\noindent{\bf(c)} We are given a subset $C(\A)_\pe\subseteq C(\A)$ of {\it permissible classes}.
\smallskip

\noindent{\bf(d)} For each $\al\in C(\A)_\pe$ we are given the following data:
\begin{itemize}
\setlength{\itemsep}{0pt}
\setlength{\parsep}{0pt}
\item[{\bf(i)}] Open $\C$-substacks $\dM_\al^\pl\subseteq\M_\al^\pl$ and $\dM_\al=(\Pi_\al^\pl)^{-1}(\dM_\al^\pl)\ab\subseteq\M_\al$, where we write $\dot\Pi_\al^\pl=\Pi_\al^\pl\vert_{\dot{\mathcal M}_\al}:\dM_\al\ra\dM^\pl_\al$. We write $\bs\dM_\al^\pl\subseteq\bs\M_\al^\pl$ and $\bs\dM_\al\subseteq\bs\M_\al$ for the corresponding derived open substacks.
\item[{\bf(ii)}] We are given a (homotopy) Cartesian square of derived Artin $\C$-stacks:
\e
\begin{gathered}
\xymatrix@C=100pt@R=15pt{ *+[r]{\bs\dM^\red_\al} \ar[r]_{\bs{\dot\Pi}_\al^\rpl} \ar[d]^{\bs j_\al} & *+[l]{\bs\dM^\rpl_\al} \ar[d]_{\bs j_\al^\pl} \\
 *+[r]{\bs\dM_\al} \ar[r]^{\bs{\dot\Pi}_\al^\pl}  & *+[l]{\bs\dM^\pl_\al,\!} }
\end{gathered}
\label{co1eq10}
\e
where $\bs\dM^{\red}_\al$ and $\bs\dM^\rpl_\al$ are locally finitely presented derived Artin stacks which are {\it quasi-smooth}, that is, the cotangent complexes $\bL_{\bs\dM^{\red}_\al}$, $\bL_{\bs\dM^{\rpl}_\al}$ are perfect in the interval~$[-1,1]$. 

We require that the classical truncations $t_0(\bs j_\al),t_0(\bs j_\al^\pl)$ should be isomorphisms. Thus we may take the classical truncations to be $t_0(\bs\dM^{\red}_\al)=\dM_\al$, $t_0(\bs\dM^\rpl_\al)=\dM_\al^\pl$. Therefore $\bs\dM^\red_\al$ and $\bs\dM_\al$ are {\it two different\/} derived enhancements of the same classical stack~$\dM_\al$. 

Since $\bs\dM^{\red}_\al,\bs\dM^\rpl_\al$ are quasi-smooth, the following are perfect obstruction theories on the classical Artin stacks $\dM_\al,\dM_\al^\pl$ in the sense of~\S\ref{co24}:
\e
\begin{split}
\bL_{i_{\dM^\red_\al}}&:i_{\dM^\red_\al}^*(\bL_{\bs\dM^\red_\al})\longra \bL_{\dM_\al},\\
\bL_{i_{\dM^\rpl_\al}}&:i_{\dM^\rpl_\al}^*(\bL_{\bs\dM^\rpl_\al})\longra \bL_{\dM_\al^\pl}.
\end{split}
\label{co1eq11}
\e
\item[{\bf(iii)}] A finite-dimensional $\C$-vector space $U_\al$ with $\dim U_\al=o_\al\ge 0$, and isomorphisms of the cotangent complexes of the morphisms $\bs j_\al,\bs j_\al^\pl$ in~\eq{co1eq10}:
\e
\bL_{\bs\dM^\red_\al/\bs\dM_\al}\cong U_\al\ot_\C\O_{\bs\dM^\red_\al}[2],\quad
\bL_{\bs\dM^\rpl_\al/\bs\dM^\pl_\al}\cong U_\al\ot_\C\O_{\bs\dM^\rpl_\al}[2].
\label{co1eq12}
\e

Taking the distinguished triangle of cotangent complexes of $\bs j_\al$, pulling back to $\dM_\al$, noting that  $\bs\dM^\red_\al,\bs\dM_\al$ have the same classical truncation, and using \eq{co1eq12}, gives a distinguished triangle on $\dM_\al$:
\e
\!\!\!\!\!\!\!\!\!\!\!\!\xymatrix@C=11pt{ U_\al\ot_\C\O_{\dot{\mathcal M}_\al}[1] \ar[r] & i_{\dM_\al}^*(\bL_{\bs\dM_\al}) \ar[r] & i_{\dM^\red_\al}^*(\bL_{\bs\dM^\red_\al}) \ar[r] & U_\al\ot_\C\O_{\dot{\mathcal M}_\al}[2]. }\!\!\!\!\!
\label{co1eq13}
\e
\item[{\bf(iv)}] If $\al,\be,\al+\be\in C(\A)_\pe$ we require that $\dim U_\al+\dim U_\be\ge \dim U_{\al+\be}$, that is,~$o_\al+o_\be\ge o_{\al+\be}$.
\end{itemize}
We refer to $\bs\dM^{\red}_\al,\bs\dM^\rpl_\al$ as `reduced' versions of $\bs\dM_\al,\bs\dM^\pl_\al$, in the sense of the `reduced obstruction theories' in \cite{KoTh1,KoTh2,MaPa,MPT,Schu,STV}. If $U_\al=o_\al=0$ then 
$\bs\dM^{\red}_\al=\bs\dM_\al$ and $\bs\dM^\rpl_\al=\bs\dM^\pl_\al$, and we call \eq{co1eq11} `non-reduced' (ordinary) obstruction theories.
\smallskip

\noindent{\bf(e)} There should exist a family of functors $F_k:\A\ra\Vect_\C$ for $k\in K$ satisfying some conditions. The invariants $[\M_\al^\ss(\tau)]_\inv$ in Theorem \ref{co1thm3} are independent of the choice of $\{F_k:k\in K\}$. The $F_k$ are used to construct auxiliary abelian categories $\baA$ which roughly lie in exact sequences $0\ra\A\ra\baA\ra\modCQ\ra 0$ for $Q$ a quiver, which are used in the proof of Theorem~\ref{co1thm3}.
\end{ass}

\begin{ass}
\label{co1ass3}
Let Assumption \ref{co1ass2} hold. Assume we are given a set $\sS$ of weak stability conditions on $\A$, satisfying:
\smallskip

\noindent{\bf(a)} If $(\tau,T,\le)\in\sS$ then $\A$ is $\tau$-artinian, as in \S\ref{co31}.
\smallskip

\noindent{\bf(b)} If $(\tau,T,\le)\in\sS$, for $E\in\A$ to be $\tau$-(semi)stable are open conditions on $[E]$ in $\M,\M^\pl$, so we have open $\C$-substacks $\M_\al^\rst(\tau)\subseteq\M_\al^\ss(\tau)\subseteq\M_\al^\pl$ parametrizing $\tau$-(semi)stable objects for all $\al\in C(\A)$. 

If $(\tau,T,\le)\in\sS$ and $\al\in C(\A)_\pe$ then $\M_\al^\rst(\tau),\M_\al^\ss(\tau)$ are of finite type.
\smallskip

\noindent{\bf(c)} If $(\tau,T,\le)\in\sS$ and $E,F\in\A$ are $\tau$-semistable with $\tau(\lb E\rb)=\tau(\lb F\rb)$ in $T$, and $\lb E\op F\rb\in C(\A)_\pe$, then $\lb E\rb,\lb F\rb\in C(\A)_\pe$.
\smallskip

\noindent{\bf(d)} Suppose $(\tau,T,\le),(\ti\tau,\ti T,\le)\in\sS$, and $I\subseteq C(\A)_\pe$ is a finite subset, and $\al\in I$, satisfy $\tau(\be)=\tau(\al)$ and $\M_\be^\ss(\tau)\ne\es$ for all $\be\in I$. Then there should exist a group morphism $\la:K(\A)\ra\R$ such that $\la(\al)=0$ and $\la(\be)>0$ (or $\la(\be)<0$) if and only if $\ti\tau(\be)>\ti\tau(\al)$ (or $\ti\tau(\be)<\ti\tau(\al)$, respectively) for all~$\be\in I$.

\smallskip

\noindent{\bf(e)} If $(\tau,T,\le)\in\sS$ and $\al\in C(\A)_\pe$ then $\M_\al^\ss(\tau)\subseteq\dM_\al^\pl$, for $\dM_\al^\pl$ as in Assumption \ref{co1ass2}(d), so the obstruction theory $\bL_{i_{\dM^\rpl_\al}}:i_{\dM^\rpl_\al}^*(\bL_{\bs\dM^\rpl_\al})\ra \bL_{\dM_\al^\pl}$ on $\dM_\al^\pl$ in \eq{co1eq11} restricts to an obstruction theory on~$\M_\al^\ss(\tau)$.
\smallskip

\noindent{\bf(f)} We are given a `rank function' $\rk:C(\A)\ra\N_{>0}=\{1,2,\ldots\}$ such that if $\al,\be\in C(\A)$ and $(\tau,T,\le)\in\sS$ with $\tau(\al)=\tau(\be)$ then~$\rk(\al+\be)=\rk(\al)+\rk(\be)$.

\smallskip

\noindent{\bf(g)} If $(\tau,T,\le)\in\sS$ and $\al\in C(\A)_\pe$ with $\M_\al^\rst(\tau)=\M_\al^\ss(\tau)$ then $\M_\al^\ss(\tau)$ is a proper algebraic space. Hence as in \S\ref{co24}, using the obstruction theory from {\bf(e)} we have a Behrend--Fantechi virtual class $[\M_\al^\ss(\tau)]_\virt,$ which we regard as lying in the Betti\/ $\Q$-homology group~$H_{2+2o_\al-2\chi(\al,\al)}(\M_\al^\pl)=\check H_{2o_\al}(\M_\al^\pl)$.
\smallskip

\noindent{\bf(h)} The analogue of {\bf(g)} holds for certain auxiliary semistable moduli stacks constructed in \S\ref{co52} using the data of Assumption \ref{co1ass2}(e).
\end{ass}

\begin{ass}
\label{co1ass4}
Let Assumptions \ref{co1ass2}--\ref{co1ass3} hold. Assume:
\smallskip

\noindent{\bf(a)} If $(\tau,T,\le),(\ti\tau,\ti T,\le)\in\sS$ then there exists a continuous family of weak stability conditions $(\tau_t,T_t,\le)_{t\in[0,1]}$ in $\sS$, in the sense of \S\ref{co32}, with  $(\tau_0,T_0,\le)\!=\!(\tau,T,\le)$ and~$(\tau_1,T_1,\le)\!=\!(\ti\tau,\ti T,\le)$.
\smallskip

\noindent{\bf(b)} Suppose $(\tau,T,\le),(\ti\tau,\ti T,\le)\in\sS$ and $\al\in C(\A)_\pe$. Then there are only finitely many decompositions $\al=\al_1+\cdots+\al_n$ with $\M_{\al_i}^\ss(\tau)\ne\es$ for all $i$ for which the combinatorial coefficient $U\bigl(\al_1,\ldots,\al_n;\tau,\ti\tau)$ defined in \S\ref{co32} is nonzero, and each such has $\al_1,\ldots,\al_n\in C(\A)_\pe$. This implies that the wall-crossing formulae \eq{co1eq14}--\eq{co1eq15} below are finite sums.
\smallskip

\noindent{\bf(c)} A more complicated version of {\bf(b)} holds, in which $(\tau,T,\le),(\ti\tau,\ti T,\le)$ vary in a continuous family $(\tau_t,T_t,\le)_{t\in[0,1]}$ in $\sS$.
\end{ass}

\begin{rem}
\label{co1rem2}
{\bf(a)} Assumptions \ref{co1ass2}--\ref{co1ass4} are designed to include as many interesting cases as possible, at the cost of making them more complex.

The reason for including the derived enhancements $\bs\M,\ab\bs\M^\pl,\ab\bs\dM^\red_\al,\ab\bs\dM^\rpl_\al$ of $\M,\ldots,\dM^\pl_\al$ is to produce the Behrend--Fantechi obstruction theories \eq{co1eq11} on $\dM_\al,\dM^\pl_\al$. Readers unfamiliar with Derived Algebraic Geometry will not lose much if they pretend a derived stack $\bs\M$ is a classical stack $\M$ together with an obstruction theory $\phi:\cF^\bu\ra\bL_\M$, except that $\cF^\bu$ is only perfect in $[-1,1]$ (one of the conditions on obstruction theories) if $\bs\M$ is quasi-smooth.

The author would have preferred to write our theory in terms of stacks with obstruction theories, without mentioning derived stacks. However, the proofs of our main results rely on enumerative invariants for moduli stacks $\baM,\baM^\pl$ of objects in auxiliary categories $\baA$ defined in \S\ref{co52}. It it is easier and more functorial to construct the corresponding derived stacks $\bs\baM,\bs\baM^\pl$, than to build the obstruction theories on $\baM,\baM^\pl$ without Derived Algebraic Geometry.
\smallskip

\noindent{\bf(b)} In most of our examples we take $\dM_\al=\M_\al$, $\dM_\al^\pl=\M_\al^\pl$, $\bs\dM^{\red}_\al=\bs\M_\al$, $\bs\dM^\rpl_\al=\bs\M_\al^\pl$, $\bs j_\al=\bs\id=\bs j_\al^\pl$, and $U_\al=o_\al=0$, so Assumption \ref{co1ass2}(d) is determined by Assumption \ref{co1ass2}(a)--(b). We call this the `non-reduced' case.

However, there are some cases in which for the natural obstruction theory $\bL_{i_{\M_\al}}:i_{\M_\al}^*(\bL_{\bs\M_\al})\ra \bL_{\M_\al}$ on $\M_\al$, $h^{-1}(i_{\M_\al}^*(\bL_{\bs\M_\al}))$ contains a trivial vector bundle $U_\al\ot_\C\O_{\M_\al}$ for $\dim U_\al>0$, so any virtual class $[\M_\al^\ss(\tau)]_\virt$ defined using this or its analogue on $\M_\al^\pl$ is zero by Theorem \ref{co2thm1}(iv). Then our theory would yield invariants $[\M_\al^\ss(\tau)]_\inv=0$, which would be boring. If $\A=\coh(X)$ for $X$ a projective surface with geometric genus $p_g=\dim H^0(K_X)>0$, and $\rank\al>0$, this holds with~$U_\al=H^0(K_X)$. 

This problem is well known in the literature \cite{KoTh1,KoTh2,MaPa,MPT,Schu,STV}, and the solution is to define a modified obstruction theory by deleting the $U_\al$ factor in $h^{-1}(i_{\M_\al}^*(\bL_{\bs\M_\al}))$, as in \eq{co1eq13}. These are known as `reduced' obstruction theories. We take our `reduced' obstruction theories to come from `reduced' derived stacks $\bs\dM^{\red}_\al,\bs\dM^\rpl_\al$. With these `reduced' obstruction theories our theory will produce invariants $[\M_\al^\ss(\tau)]_\inv$ which can be nonzero.
\end{rem}

The next theorem partly summarizes our main results, Theorems~\ref{co5thm1}--\ref{co5thm3}:

\begin{thm}
\label{co1thm3}
Suppose Assumptions\/ {\rm\ref{co1ass2}--\ref{co1ass4}} hold. Then for all\/ $(\tau,T,\le)$ in $\sS$ and\/ $\al$ in $C(\A)_\pe$ there are unique classes $[\M_\al^\ss(\tau)]_\inv$ in the Betti\/ $\Q$-homology group $H_{2o_\al+2-2\chi(\al,\al)}(\M_\al^\pl)=\check H_{2o_\al}(\M_\al^\pl)$ satisfying:
\begin{itemize}
\setlength{\itemsep}{0pt}
\setlength{\parsep}{0pt}
\item[{\bf(i)}] If\/ $\M_\al^\rst(\tau)=\M_\al^\ss(\tau)$ then\/ $[\M_\al^\ss(\tau)]_\inv=[\M_\al^\ss(\tau)]_\virt,$ where $[\M_\al^\ss(\tau)]_\virt$ is the Behrend--Fantechi virtual class of the proper algebraic space $\M_\al^\ss(\tau)$ with perfect obstruction theory, as in Assumption\/~{\rm\ref{co1ass3}(g)}.
\item[{\bf(ii)}] Suppose $(\tau,T,\le),(\ti\tau,\ti T,\le)\in\sS$ and\/ $\al\in C(\A)_\pe$ with\/ $\M_\al^\ss(\tau)=\M_\al^\ss(\ti\tau)$. Then $[\M_\al^\ss(\tau)]_\inv=[\M_\al^\ss(\ti\tau)]_\inv$.
\item[{\bf(iii)}] The $[\M_\al^\ss(\tau)]_\inv$ can be characterized by induction on $\rk\al$  using virtual classes $[\baM^\ss_{(\al,1)}(\bar\tau^0_1)]_\virt$ of auxiliary moduli spaces $\baM^\ss_{(\al,1)}(\bar\tau^0_1)$ parametrizing pairs $(E,\rho)$ of\/ $[E]\in\M_\al^\ss(\tau)$ and\/ $\rho\in\P(F_k(E))$ satisfying a stability condition, for $F_k$ as in Assumption\/ {\rm\ref{co1ass2}(e)} satisfying conditions.

In particular, this means that any \begin{bfseries}deformation invariance\end{bfseries} properties of\/ $[\baM^\ss_{(\al,1)}(\bar\tau^0_1)]_\virt$ also hold for the $[\M_\al^\ss(\tau)]_\inv$. So, for example, if\/ $\A=\coh(X)$ for a smooth projective $\C$-scheme $X$ then (suitably interpreted) the $[\M_\al^\ss(\tau)]_\inv$ will be unchanged by continuous deformations of\/~$X$.
\item[{\bf(iv)}] For\/ $(\tau,T,\le),(\ti\tau,\ti T,\le)$ in $\sS$ and\/ $\al\in C(\A)_\pe$ we have
\e
\begin{gathered}
{}\!\!\!\!\!\!\!\!\!\!\![\M_\al^\ss(\ti\tau)]_\inv= \!\!\!\!\!\!\!
\sum_{\begin{subarray}{l}n\ge 1,\;\al_1,\ldots,\al_n\in
C(\A)_\pe:\\ \M_{\al_i}^\ss(\tau)\ne\es,\; \text{all\/ $i,$} \\
\al_1+\cdots+\al_n=\al, \; o_{\al_1}+\cdots+o_{\al_n}=o_\al 
\end{subarray}} \!\!\!\!\!\!\!\!\!\!\!\!\!\!\!\!\!\!\!\!\!\!\!\!\begin{aligned}[t]
\ti U(\al_1,&\ldots,\al_n;\tau,\ti\tau)\cdot\bigl[\bigl[\cdots\bigl[[\M_{\al_1}^\ss(\tau)]_\inv,\\
&
[\M_{\al_2}^\ss(\tau)]_\inv\bigr],\ldots\bigr],[\M_{\al_n}^\ss(\tau)]_\inv\bigr]
\end{aligned}
\end{gathered}
\label{co1eq14}
\e
in the Lie algebra $\check H_{\rm even}(\M^\pl)$ from Theorem\/ {\rm\ref{co1thm2}}. Here $\ti U(-;\tau,\ti\tau)$ and\/ $U(-;\tau,\ti\tau)$ are combinatorial coefficients defined in\/ {\rm\S\ref{co32},} and there are only finitely many nonzero terms in \eq{co1eq14} by Assumption\/ {\rm\ref{co1ass4}(b)}. Equivalently, in the universal enveloping algebra $\bigl(U(\check H_{\rm even}(\M^\pl)),*\bigr)$ we have:
\e
\begin{gathered}
{}\!\!\!\!\!\!\!\!\!\!\![\M_\al^\ss(\ti\tau)]_\inv= \!\!\!\!\!\!\!
\sum_{\begin{subarray}{l}n\ge 1,\;\al_1,\ldots,\al_n\in
C(\A)_\pe:\\ \M_{\al_i}^\ss(\tau)\ne\es,\; \text{all\/ $i,$} \\
\al_1+\cdots+\al_n=\al, \; o_{\al_1}+\cdots+o_{\al_n}=o_\al \end{subarray}} \!\!\!\!\!\!\!\!\!\!\!\!\!\!\!\!\!\!\!
\begin{aligned}[t]
U(\al_1,&\ldots,\al_n;\tau,\ti\tau)\cdot[\M_{\al_1}^\ss(\tau)]_\inv *{}\\
&
[\M_{\al_2}^\ss(\tau)]_\inv*\cdots *[\M_{\al_n}^\ss(\tau)]_\inv.
\end{aligned}
\end{gathered}
\label{co1eq15}
\e

We call\/ {\rm\eq{co1eq14}--\eq{co1eq15}} \begin{bfseries}wall-crossing formulae\end{bfseries} for the $[\M_\al^\ss(\tau)]_\inv$. 
\end{itemize}

We think of the $[\M_\al^\ss(\tau)]_\inv$ as \begin{bfseries}enumerative invariants\end{bfseries} which `count' the moduli spaces $\M_\al^\ss(\tau)$ for all\/ $(\tau,T,\le)\in\sS$ and\/ $\al\in C(\A)_\pe$. They have the important property that\/ $[\M_\al^\ss(\tau)]_\inv$ is well defined even when $\M_\al^\rst(\tau)\ne\M_\al^\ss(\tau),$ so that the Behrend--Fantechi virtual class $[\M_\al^\ss(\tau)]_\virt$ is not defined.
\end{thm}

When an algebraic $\C$-group $G$ acts on the data $\A,\M,\M^\pl,\ldots,$ Theorem \ref{co5thm4} extends Theorem \ref{co1thm3} to equivariant homology $H_*^G(\M^\pl)$, as in~\S\ref{co23}.

The proof of Theorem \ref{co1thm3} in Chapters \ref{co9}--\ref{co11} is the longest and most complex the author has ever written. A very brief and inadequate summary is this: using Assumption \ref{co1ass2}(e) we define auxiliary abelian categories $\baA$ in exact sequences $0\ra\A\ra\baA\ra\modCQ\ra 0$ for $Q$ a quiver. Weak stability conditions $(\tau,T,\le)$ on $\A$ lift to families of weak stability conditions $(\bar\tau^\la_{\bs\mu},\bar T,\le)$ on~$\baA$. 

Moduli stacks $\baM_{(\al,\bs d)}^\ss(\bar\tau^\la_{\bs\mu})$ in $\baA$ have forgetful morphisms to moduli stacks $\M_\al^\ss(\tau)$ in $\A$. Even when $\M_\al^\rst(\tau)\ne\M_\al^\ss(\tau)$ there are usually choices of $\la,\bs\mu$ such that $\baM_{(\al,\bs d)}^\rst(\bar\tau^\la_{\bs\mu})=\baM_{(\al,\bs d)}^\ss(\bar\tau^\la_{\bs\mu})$, so that a virtual class $[\baM_{(\al,\bs d)}^\ss(\bar\tau^\la_{\bs\mu})]_\virt$ is defined. Thus we can form the homology class
\e
(\Pi_{\M^\ss_\al(\tau)})_*\bigl([\baM_{(\al,\bs d)}^\ss(\bar\tau^\la_{\bs\mu})]_\virt\cap c_\top(\bT_{\baM_{(\al,\bs d)}^\ss(\bar\tau^\la_{\bs\mu})/\M^\ss_{\al}(\tau)})\bigr)\in \check H_{2o_\al}(\M_\al^\pl).
\label{co1eq16}
\e
In simple cases \eq{co1eq16} is the product of the invariant $[\M_\al^\ss(\tau)]_\inv$ we seek, with the Euler characteristic of the fibre of $\Pi_{\M^\ss_\al(\tau)}:\baM_{(\al,\bs d)}^\ss(\bar\tau^\la_{\bs\mu})\ra\M^\ss_{\al}(\tau)$. Otherwise we must correct this with terms in $[\M_\be^\ss(\tau)]_\inv$ for $\rk\be<\rk\al$, leading to a definition of $[\M_\al^\ss(\tau)]_\inv$ by induction on $\rk\al$ in Assumption~\ref{co1ass3}(f).

To prove the wall-crossing formulae \eq{co1eq14}--\eq{co1eq15}, we distinguish between crossing `simple walls' involving only splittings $\al=\al_1+\al_2$, and `complex walls' involving splittings $\al=\al_1+\cdots+\al_n$ for $n\ge 3$. We prove the analogues of \eq{co1eq14}--\eq{co1eq15} for `simple walls' by using (in part) known techniques, involving a `master space' and $\bG_m$-localization of virtual classes, see for example Kiem--Li \cite{KiLi1}. Then we prove \eq{co1eq14}--\eq{co1eq15} by reducing crossing complex walls in $\A$ to crossing simple walls in auxiliary categories $\baA$, where the extension $0\ra\A\ra\baA\ra\modCQ\ra 0$ is chosen to encode the structure of the complex wall.

\begin{rem}
\label{co1rem3}
As in Remark \ref{co1rem2}, in Assumption \ref{co1ass2}(d) the obstruction theories \eq{co1eq11} are either `non-reduced', with $o_\al=0$ (this is the case in most of our examples), or `reduced', with~$o_\al>0$.

Notice the conditions $o_{\al_1}+\cdots+o_{\al_n}=o_\al$ in \eq{co1eq14}--\eq{co1eq15}. In the `non-reduced' case with $o_\al=0$ for all $\al$, these hold automatically. But in the `reduced' case, they can exclude many terms. For example, if $o_\al=c>0$ for all $\al$ (this holds for $\A=\coh(X)$ for $X$ a surface with $p_g>0$, and classes $\al$ with $\rank\al>0$) then $o_{\al_1}+\cdots+o_{\al_n}=o_\al$ holds only when $n=1$, so \eq{co1eq14}--\eq{co1eq15} reduce to $[\M_\al^\ss(\ti\tau)]_\inv=[\M_\al^\ss(\tau)]_\inv$. That is, the invariants $[\M_\al^\ss(\tau)]_\inv$ are independent of the weak stability condition $(\tau,T,\le)$ in this case.	
\end{rem}

\subsection{An example: coherent sheaves on projective surfaces}
\label{co1.4}

Rather than discussing all the applications in Chapters \ref{co6}--\ref{co8}, to show how Theorem \ref{co1thm3} works out in practice, in this introduction we will concentrate on one example, the abelian category $\coh(X)$ of coherent sheaves on a connected projective complex surface $X$. The {\it geometric genus\/} of $X$ is~$p_g=\dim H^0(K_X)$.

As in Remarks \ref{co1rem2} and \ref{co1rem3}, we distinguish the cases $p_g=0$, when we use `non-reduced' obstruction theories in Assumption \ref{co1ass2}(d), and $p_g>0$, when we use `reduced' obstruction theories for $\rank\al>0$. Because of the conditions $o_{\al_1}+\cdots+o_{\al_n}=o_\al$ in \eq{co1eq14}--\eq{co1eq15}, the wall-crossing formulae look rather different in these two cases. Here are Theorems \ref{co7thm7}--\ref{co7thm8} from~\S\ref{co77}.

\begin{thm}
\label{co1thm4}
Let\/ $X$ be a connected projective complex surface with geometric genus\/ $p_g=0$. Write $\M^\pl=\coprod_{\al\in K(\coh(X))}\M_\al^\pl$ for the projective linear moduli stack of objects in $\coh(X)$. Then:
\smallskip

\noindent{\bf(a)} Let\/ $(\tau,T,\le)$ be Gieseker stability $(\tau^\om,G,\le)$ or $\mu$-stability $(\mu^\om,M,\le)$ on $\coh(X)$ for some $\om\in\Kah(X),$ as in Examples\/ {\rm\ref{co1ex2}--\ref{co1ex3}}. Then for all\/ $\al$ in $C(\coh(X))$ there are unique classes $[\M_\al^\ss(\tau)]_\inv$ in the Betti\/ $\Q$-homology group $\check H_0(\M^\pl_\al)=H_{2-2\chi(\al,\al)}(\M_\al^\pl),$ satisfying:
\begin{itemize}
\setlength{\itemsep}{0pt}
\setlength{\parsep}{0pt}
\item[{\bf(i)}] If\/ $\M_\al^\rst(\tau)=\M_\al^\ss(\tau)$ then\/ $[\M_\al^\ss(\tau)]_\inv=[\M_\al^\ss(\tau)]_\virt,$ where $[\M_\al^\ss(\tau)]_\virt$ is the Behrend--Fantechi virtual class of the proper algebraic space $\M_\al^\ss(\tau)$ with (non-reduced) perfect obstruction theory.
\item[{\bf(ii)}] Let\/ $\al\in C(\coh(X)),$ and\/ $\O_X(1)\ra X$ be an ample line bundle, and\/ $N$ be sufficiently large that\/ $E$ is $N$-regular for $\O_X(1)$ for all\/ $[E]\in\M_\al^\ss(\tau)$. Then Example\/ {\rm\ref{co5ex1}} defines a proper moduli space\/ $\baM^\ss_{(\al,1)}(\bar\tau^0_1)$ of nonzero morphisms\/ $\rho:\O_X(-N)\ra E$ for\/ $[E]\in\M_\al^\ss(\tau),$ with an obstruction theory, and a forgetful morphism\/ $\Pi_{\M^\ss_\al(\tau)}:\baM^\ss_{(\al,1)}(\bar\tau^0_1)\!\ra\!\M_\al^\ss(\tau)$. We~have
\ea
&(\Pi_{\M^\ss_\al(\tau)})_*\bigl([\baM^\ss_{(\al,1)}(\bar\tau^0_1)]_\virt\cap c_\top(\bT_{\baM^\ss_{(\al,1)}(\bar\tau^0_1)/\M^\ss_{\al}(\tau)})\bigr)
\label{co1eq17}\\
&=\!\!\!\!\!\!\sum_{\begin{subarray}{l}n\ge 1,\;\al_1,\ldots,\al_n\in
C(\coh(X)):\\ \al_1+\cdots+\al_n=\al,\\ 
\tau(\al_i)=\tau(\al), \; \M_{\al_i}^\ss(\tau)\ne\es,\; \text{all\/ $i$}\end{subarray}\!\!\!\!\!\!\!\!\!\!\!\!\!\!\!\!\!\!\!\!\!\!\!\!\!} \!\!\!\begin{aligned}[t]
\frac{(-1)^{n+1}P_{\al_1}^{c_1(\O_X(1))}(N)}{n!}\,\cdot&\bigl[\bigl[\cdots\bigl[[\M_{\al_1}^\ss(\tau)]_\inv,\\
&\; 
[\M_{\al_2}^\ss(\tau)]_\inv\bigr],\ldots\bigr],[\M_{\al_n}^\ss(\tau)]_\inv\bigr],
\end{aligned}
\nonumber
\ea
where $P_{\al_1}^{c_1(\O_X(1))}(t)$ is the Hilbert polynomial of\/ $\al_1$ from\/ {\rm\eq{co1eq2},} the Lie brackets are in the Lie algebra $\check H_0(\M^\pl)$ from Theorem\/ {\rm\ref{co1thm2},} and there are only finitely many terms in the sum.
\end{itemize}

\noindent{\bf(b)} Suppose $(\tau,T,\le),(\ti\tau,\ti T,\le)$ are two Gieseker or $\mu$-stability conditions as in part\/ {\bf(a)\rm,} from possibly different\/ $\om,\ti\om\in\Kah(X)$. Then 
\begin{itemize}
\setlength{\itemsep}{0pt}
\setlength{\parsep}{0pt}
\item[{\bf(i)}] If\/ $\al\in C(\coh(X))$ with\/ $\M_\al^\ss(\tau)=\M_\al^\ss(\ti\tau)$ then $[\M_\al^\ss(\tau)]_\inv=[\M_\al^\ss(\ti\tau)]_\inv$.
\item[{\bf(ii)}] For all\/ $\al\in C(\coh(X))$ we have
\end{itemize}
\e
\begin{gathered}
{}
[\M_\al^\ss(\ti\tau)]_\inv= \!\!\!\!\!\!\!
\sum_{\begin{subarray}{l}n\ge 1,\;\al_1,\ldots,\al_n\in
C(\coh(X)):\\ \M_{\al_i}^\ss(\tau)\ne\es,\; \text{all\/ $i,$} \\
\al_1+\cdots+\al_n=\al 
\end{subarray}} \!\!\!\!\!\!\!\!\!\!\!\begin{aligned}[t]
\ti U(\al_1,&\ldots,\al_n;\tau,\ti\tau)\cdot\bigl[\bigl[\cdots\bigl[[\M_{\al_1}^\ss(\tau)]_\inv,\\
&
[\M_{\al_2}^\ss(\tau)]_\inv\bigr],\ldots\bigr],[\M_{\al_n}^\ss(\tau)]_\inv\bigr]
\end{aligned}
\end{gathered}
\label{co1eq18}
\e
\begin{itemize}
\setlength{\itemsep}{0pt}
\setlength{\parsep}{0pt}
\item[] in the Lie algebra $\check H_0(\M^\pl)$ from {\rm\S\ref{co1.2}}. Here $\ti U(-;\tau,\ti\tau)$ is as in\/ {\rm\S\ref{co32},} and there are only finitely many nonzero terms in \eq{co1eq18}. Equivalently, in the universal enveloping algebra $\bigl(U(\check H_0(\M^\pl)),*\bigr)$ we have:
\end{itemize}
\e
\begin{gathered}
{}
[\M_\al^\ss(\ti\tau)]_\inv= \!\!\!\!\!\!\!
\sum_{\begin{subarray}{l}n\ge 1,\;\al_1,\ldots,\al_n\in
C(\coh(X)):\\ \M_{\al_i}^\ss(\tau)\ne\es,\; \text{all\/ $i,$} \\
\al_1+\cdots+\al_n=\al \end{subarray}} \!\!\!\!\!\!\!\!\!\!\!\!
\begin{aligned}[t]
U(\al_1,&\ldots,\al_n;\tau,\ti\tau)\cdot[\M_{\al_1}^\ss(\tau)]_\inv *{}\\
&
[\M_{\al_2}^\ss(\tau)]_\inv*\cdots *[\M_{\al_n}^\ss(\tau)]_\inv.
\end{aligned}
\end{gathered}
\label{co1eq19}
\e

\noindent{\bf(c)} If an algebraic $\C$-group $G$ acts on $X$ satisfying conditions Definition\/ {\rm\ref{co7def12}(i)--(ii)} then {\bf(a)\rm,\bf(b)} generalize to equivariant homology $H_*^G(\M^\pl),$ as in\/~{\rm\S\ref{co23}}.
\end{thm}

\begin{thm}
\label{co1thm5}
Let\/ $X$ be a connected projective complex surface with geometric genus\/ $p_g=\dim H^0(K_X)>0$. Write $\M^\pl=\coprod_{\al\in K(\coh(X))}\M_\al^\pl$ for the projective linear moduli stack of objects in $\coh(X)$. Then:
\smallskip

\noindent{\bf(a)} Let\/ $(\tau,T,\le)$ be Gieseker stability $(\tau^\om,G,\le)$ or $\mu$-stability $(\mu^\om,M,\le)$ on $\coh(X)$ for some $\om\in\Kah(X)$. Then for all\/ $\al\in C(\coh(X))$ with\/ $\rank\al>0$ there are unique classes $[\M_\al^\ss(\tau)]_\inv$ in the Betti\/ $\Q$-homology group $\check H_{2p_g}(\M^\pl_\al)\ab=H_{2p_g+2-2\chi(\al,\al)}(\M_\al^\pl),$ satisfying:
\begin{itemize}
\setlength{\itemsep}{0pt}
\setlength{\parsep}{0pt}
\item[{\bf(i)}] If\/ $\M_\al^\rst(\tau)=\M_\al^\ss(\tau)$ then\/ $[\M_\al^\ss(\tau)]_\inv=[\M_\al^\ss(\tau)]_\virt,$ where $[\M_\al^\ss(\tau)]_\virt$ is the Behrend--Fantechi virtual class of the proper algebraic space $\M_\al^\ss(\tau)$ with \begin{bfseries}reduced\end{bfseries} perfect obstruction theory, constructed in\/ {\rm\S\ref{co73}}.
\item[{\bf(ii)}] The invariants $[\M_\al^\ss(\tau)]_\inv$ for $\rank\al>0$ are independent of the weak stability condition\/~$(\tau,T,\le)$.
\item[{\bf(iii)}] Let\/ $\al\in C(\coh(X))$ with\/ $\rank\al>0,$ and\/ $\O_X(1)\ra X$ be an ample line bundle, and\/ $N$ be sufficiently large that\/ $E$ is $N$-regular for $\O_X(1)$ for all\/ $[E]\in\M_\al^\ss(\tau)$. Then Example\/ {\rm\ref{co5ex1}} defines a proper moduli space\/ $\baM^\ss_{(\al,1)}(\bar\tau^0_1)$ of nonzero morphisms\/ $\rho:\O_X(-N)\ra E$ for\/ $[E]\in\M_\al^\ss(\tau),$ with a reduced perfect obstruction theory as in {\bf(i)\rm,} and a forgetful morphism\/ $\Pi_{\M^\ss_\al(\tau)}:\baM^\ss_{(\al,1)}(\bar\tau^0_1)\ra\M_\al^\ss(\tau)$. We have
\e
\begin{split}
&(\Pi_{\M^\ss_\al(\tau)})_*\bigl([\baM^\ss_{(\al,1)}(\bar\tau^0_1)]_\virt\cap c_\top(\bT_{\baM^\ss_{(\al,1)}(\bar\tau^0_1)/\M^\ss_{\al}(\tau)})\bigr)\\
&\qquad =P_\al^{c_1(\O_X(1))}(N)\cdot [\M_\al^\ss(\tau)]_\inv,
\end{split}
\label{co1eq20}
\e
where $P_\al^{c_1(\O_X(1))}(t)$ is the Hilbert polynomial of\/ $\al$ from\/ \eq{co1eq2}. 
\end{itemize}

\noindent{\bf(b)} We also define invariants $[\M_\al^\ss(\tau)]_\inv$ when\/ $\rank\al=0$. However, these use the \begin{bfseries}non-reduced\end{bfseries} obstruction theory, lie in $\check H_0(\M^\pl_\al)=H_{2-2\chi(\al,\al)}(\M_\al^\pl),$ and satisfy Theorem\/ {\rm\ref{co1thm4}(a),(b)} rather than\/~{\bf(a)(i)\rm--\bf(iii)}.
\smallskip

\noindent{\bf(c)} If an algebraic $\C$-group $G$ acts on $X$ satisfying Definition\/ {\rm\ref{co7def12}(i)--(ii)} then {\bf(a)\rm,\bf(b)} generalize to equivariant homology $H_*^G(\M^\pl)$.
\end{thm}

\begin{rem}
\label{co1rem4}
{\bf(a)} Some readers may feel Theorems \ref{co1thm4}--\ref{co1thm5} are not that useful for actual computations, since we usually cannot describe $H_*(\M^\pl)$ explicitly. However, there is a method to use them to do explicit calculations.

Write $\baM$ for the moduli stack of objects in $D^b\coh(X)$, which exists as a higher stack in the sense of \cite{Toen1} by \cite{ToVa}, and has a `projective linear' version $\baM^\pl$. The inclusion $\coh(X)\hookra D^b\coh(X)$ induces open inclusions $i:\M\hookra\baM$, $i^\pl:\M^\pl\hookra\baM^\pl$. As in \cite{Joyc12} and Remark \ref{co4rem2} below, the constructions of \S\ref{co1.2} work (with small changes) for $D^b\coh(X)$ as for $\coh(X)$, and make $H_*(\baM)$ into a graded vertex algebra, and $H_*(\baM^\pl)$ into a graded Lie algebra. Then $i_*:H_*(\M)\ra H_*(\baM)$, $i_*^\pl:H_*(\M^\pl)\ra H_*(\baM^\pl)$ are vertex and Lie algebra morphisms. So we can apply $i_*^\pl$ to regard the invariants $[\M_\al^\ss(\tau)]_\inv$ in Theorems \ref{co1thm4}--\ref{co1thm5} as elements of $H_*(\baM_\al^\pl)$, and equations \eq{co1eq17}--\eq{co1eq20} as happening in the Lie algebra~$\check H_*(\baM^\pl)$.

Now the author's PhD student Jacob Gross \cite[Th.~1.1]{Gros} proved that we can describe $H_*(\baM)$, with grading shifted as in \eq{co1eq8}, very explicitly as 
\e
\begin{split}
\hat H_*(\baM)\cong \Q[K^0_{\rm sst}(X)]&\ot_\Q \mathop{\rm Sym}\nolimits^*\bigl(K^0(X^\ran)\ot_\Z t^2\Q[t^2]\bigr)\\
&\ot_\Q\bigwedge\!\!\vphantom{\bigl(}^*\bigl(K^1(X^\ran)\ot_\Z t\Q[t^2]\bigr),
\end{split}
\label{co1eq21}
\e
where $K^0_{\rm sst}(X)$ is the {\it semi-topological K-theory\/} of $X$, as in Friedlander--Walker \cite{FrWa}, with $K^0_{\rm sst}(X)\!=\!K(\coh(X))$ modulo torsion in this case, and $K^i(X^\ran)$ are the K-theory groups of the underlying complex analytic topological space $X^\ran$ of $X$. The vertex algebra structure on $\hat H_*(\baM)$ is also described very explicitly in \cite{Gros}, and is basically a super lattice vertex algebra. So the Lie algebra $\check H_*(\baM^\pl)$ is also explicit, and $\check H_0(\baM^\pl)$ is roughly a Borcherds--Kac--Moody Lie algebra.

The reason why $H_*(\baM)$ is easier to understand than $H_*(\M)$ is that $\M$ is an `abelian monoid in stacks', with operation $\Phi:\M\t\M\ra\M$ coming from direct sum $\op$ in $\coh(X)$ as in Assumption \ref{co1ass1}(b), but $\baM$ is roughly an `abelian group in stacks', as $[1]:D^b\coh(X)\ra D^b\coh(X)$ induces a homotopy inverse for $\bar\Phi$. One can deduce that $H_*(\baM)$ should be of the general form \eq{co1eq21} from the Milnor--Moore Theorem \cite{MiMo} on the $\Q$-homology of a grouplike H-space.

Gross' result \cite[Th.~1.1]{Gros} applies if $X$ is a curve, or surface, or smooth projective toric variety, and a few other cases. However, for any smooth projective $\C$-scheme $X$ we can define a vertex algebra morphism from $\hat H_*(\M)$ to an explicit super lattice vertex algebra of the form \eq{co1eq21}. Then considering invariants $[\M_\al^\ss(\tau)]_\inv$ to take values in \eq{co1eq21}, we can do explicit computations.

Many authors define enumerative invariants as $I_\al^\Up(\tau)=\int_{[\M_\al^\ss(\tau)]_\virt}\Up\in\Q$ for suitable universal cohomology classes $\Up$. Roughly speaking, our point of view here is that \eq{co1eq21} is the dual $\cR^*$ of a $\Q$-algebra $\cR$ of universal cohomology classes $\Up$, and we consider $[\M_\al^\ss(\tau)]_\inv$ to take values in $\cR^*$, not~$H_*(\M^\pl)$.
\smallskip

\noindent{\bf(b)} The invariants $[\M_\al^\ss(\mu^\om)]_\inv$ in Theorems \ref{co1thm4}--\ref{co1thm5} for $\rank\al>0$ should be thought of as algebraic analogues of $\U(r)$-Donaldson invariants of the compact oriented 4-manifold $X$, as in Donaldson--Kronheimer~\cite[\S 9]{DoKr}.

Chapter \ref{co8} studies `stable pair invariants' $[\acM_{(\al,d)}^\ss(\ac\mu^{\om,s})]_\inv$ for $\om\in\Kah(X)$ and $s>0$ for an auxiliary abelian category $\acA$ with objects triples $(E,V,\rho)$ for $E\in\coh(X)$, $V\in\Vect_\C$ and $\rho:V\ot_\C L\ra E$ a morphism in $\coh(X)$, for a fixed line bundle $L\ra X$. These $[\acM_{(\al,d)}^\ss(\ac\mu^{\om,s})]_\inv$ include algebraic analogues of Seiberg--Witten invariants, when $\rank\al=1$, $d=1$ and $c_2(\al)=0$.

In \S\ref{co86} we explain a method similar to writing Donaldson invariants in terms of Seiberg--Witten invariants, as in \cite{FeLe1,GNY3,KrMr,MaMo,MoWi,Witt}, which allows us to compute the invariants $[\M_\al^\ss(\mu^\om)]_\inv$ in Theorems \ref{co1thm4}--\ref{co1thm5} in terms of Seiberg--Witten type invariants $[\acM_{(\al,d)}^\ss(\ac\mu^{\om,s})]_\inv$ for $\rank\al=1$ and $d=0$ or 1. The sequel \cite{Joyc13}  will use this method to compute the invariants $[\M_\al^\ss(\mu^\om)]_\inv$ in Theorems \ref{co1thm4}--\ref{co1thm5} fairly explicitly.
\smallskip

\noindent{\bf(c)} Mochizuki \cite{Moch} proves rough analogues of Theorems \ref{co1thm4}--\ref{co1thm5}. This book generalizes and builds on \cite{Moch}. There are two main ways in which Theorems \ref{co1thm4}--\ref{co1thm5} improve on~\cite{Moch}:
\begin{itemize}
\setlength{\itemsep}{0pt}
\setlength{\parsep}{0pt}
\item[{\bf(i)}] Mochizuki has no analogue of the Lie bracket on $\check H_0(\M^\pl)$ in~\eq{co1eq17}--\eq{co1eq19}.
\item[{\bf(ii)}] Mochizuki has no analogue of the universal wall-crossing formulae \eq{co1eq18}--\eq{co1eq19} or the coefficients $\ti U(\al_1,\ldots,\al_n;\tau,\ti\tau),U(\al_1,\ldots,\al_n;\tau,\ti\tau)$.
\end{itemize}

In Mochizuki's set up, it is clear on general grounds that the left hand side $[\M_\al^\ss(\ti\tau)]_\inv$ of \eq{co1eq18} is some unknown function of the terms $[\M_{\al_i}^\ss(\tau)]_\inv$ appearing on the right hand side of \eq{co1eq18}, and of $\al_i,\tau,\ti\tau$ and algebraic-topological data on $X$. In simple situations this function can be computed by hand.
\end{rem}

\subsection{Future extension to 3- and 4-Calabi--Yau categories}
\label{co1.5}

This is only a preliminary version of this book, and the final version will contain more material. Also \cite{Joyc12}, the source for \S\ref{co1.2} and Chapter \ref{co4}, is only a preliminary version, and in fact some results claimed in Chapter \ref{co4} are not in the current publicly available version of \cite{Joyc12}. We now discuss some results that the author intends will appear in the final versions of \cite{Joyc12} and this book.

In brief, the main points are these:
\begin{itemize}
\setlength{\itemsep}{0pt}
\setlength{\parsep}{0pt}
\item[(a)] There are {\it many possible variations\/} on the construction of vertex and Lie algebras on $H_*(\M)$ and $H_*(\M^\pl)$ in \S\ref{co1.2} and Chapter \ref{co4}, of which we have explained only one. These variations will be explained in the final version of \cite{Joyc12}. Some yield different kinds of vertex algebra.
\item[(b)] There are also (at least) three different flavours of obstruction theory yielding virtual classes in homology, which we will call {\it Behrend--Fantechi obstruction theories\/} \cite{BeFa1}, 3-{\it Calabi--Yau obstruction theories\/} \cite{Behr1,BeFa2,JoSo}, and 4-{\it Calabi--Yau obstruction theories\/} \cite{BoJo,KiPa,OhTh}. The current version of this book discusses only Behrend--Fantechi obstruction theories.
\item[(c)] There is a notion of $m$-{\it Calabi--Yau abelian category\/} $\A$, in which the Ext groups satisfy Serre duality $\Ext^i(E,F)\cong\Ext^{m-i}(F,E)^*$. If $X$ is a Calabi--Yau $m$-fold then $\coh(X)$ is an $m$-Calabi--Yau category. The derived moduli stack $\bs\M$ of objects in $\A$ should have a $(2-m)$-{\it shifted symplectic structure\/} in the sense of Pantev--To\"en--Vaqui\'e--Vezzosi \cite{PTVV}. This should induce an $m$-{\it Calabi--Yau obstruction theory\/} on the classical moduli stack~$\M$.
\item[(d)] The vertex algebra construction in \S\ref{co1.2} and Chapter \ref{co4} is adapted to Behrend--Fantechi obstruction theories. There are two other variations adapted to 3-Calabi--Yau and 4-Calabi--Yau obstruction theories. 

The 3-Calabi--Yau variation yields a {\it graded vertex Lie algebra\/} structure on $H_*(\M)$, and a graded Lie algebra structure on $H_*(\M^\pl)$.

The 4-Calabi--Yau variation is quite like \S\ref{co1.2} and Chapter \ref{co4}.
\item[(e)] There should be variations of Theorems \ref{co5thm1}--\ref{co5thm3} and \ref{co5thm4} which replace $\A$ by a 3- or 4-Calabi--Yau abelian category, and Behrend--Fantechi obstruction theories and virtual classes by the 3- and 4-Calabi--Yau versions, and the vertex and Lie algebras in \S\ref{co1.2} and Chapter \ref{co4} by the variations in (d).
\item[(f)] We do not really need the 3-Calabi--Yau version, as Joyce--Song \cite{JoSo} already provides analogues of Theorems \ref{co5thm1}--\ref{co5thm3} in this case. So the author may not bother extending our theory to the 3-Calabi--Yau case. However, the involvement of {\it vertex Lie algebras\/} in Donaldson--Thomas theory, as in (d), seems to be new and interesting, and may be worth further study.
\item[(g)] The author hopes to extend our theory to 4-Calabi--Yau categories in the final version of this book.
\end{itemize}
We now explain these in more detail.

\subsubsection{$m$-Calabi--Yau categories and obstruction theories}
\label{co1.5.1}

The following is not intended to be precise, just to give the general idea.

\begin{dfn}
\label{co1def4}
A $\C$-linear abelian category $\A$ is $m$-{\it Calabi--Yau\/} for $m\ge 0$ if it has Ext groups $\Ext^i(E,F)$ for $E,F\in\A$ as in Remark \ref{co1rem1}, which are finite-dimensional $\C$-vector spaces with $\Ext^0(E,F)=\Hom(E,F)$ and $\Ext^i(E,F)=0$ for $i>m$, and have functorial isomorphisms $\Ext^i(E,F)\cong \Ext^{m-i}(F,E)^*$. We also suppose that there is a moduli stack $\M$ of objects in $\A$ and an {\it Ext complex\/} $\cExt^\bu\ra\M\t\M$, perfect in the interval $[0,m]$, with $\cExt^\bu\vert_{([E],[F])}=\Ext^*(E,F)$. There should be an isomorphism $\eta:(\cExt^\bu)^\vee\ra\si^*(\cExt^\bu)[m]$, where $\si:\M\t\M\ra\M\t\M$ exchanges the factors, with $\eta^\vee=\eta[-m]$.

If $X$ is a Calabi--Yau $m$-fold, that is, a smooth projective $\C$-scheme of dimension $m$ with $K_X\cong\O_X$, then $\coh(X)$ is $m$-Calabi--Yau, where the isomorphisms $\Ext^i(E,F)\cong \Ext^{m-i}(F,E)^*=\Ext^{m-i}(F,E\ot K_X)^*$ come from Serre duality.

Pantev--To\"en--Vaqui\'e--Vezzosi \cite{PTVV} introduce $k$-{\it shifted symplectic structures\/} $\om$ on a derived stack $\bS$. This induces a quasi-isomorphism $\om\cdot{}:\bT_\bS\ra\bL_\bS[2-m]$ of (shifted) (co)tangent complexes. Derived moduli stacks $\bs\M$ of objects in an $m$-Calabi--Yau category should have $(2-m)$-shifted symplectic structures.
\end{dfn}

The next definition extends {\it symmetric obstruction theories\/} in Behrend--Fantechi \cite{BeFa2}, \cite[\S 3]{Behr1} from the 3-Calabi--Yau to the $m$-Calabi--Yau case. 

\begin{dfn}
\label{co1def5}
Let $S$ be an Artin $\C$-stack, and $m\ge 0$. An {\it $m$-Calabi--Yau obstruction theory\/} on $S$ is a morphism $\phi:\cF^\bu\ra\bL_S$ in $D_\qcoh(S)$, such that $\cF^\bu$ is perfect in the interval $[1-m,1]$, and $h^k(\phi):h^k(\cF^\bu)\ra h^k(\bL_S)$ is an isomorphism for $k\ge 0$ and surjective for $k=-1$, and we are given an isomorphism $\th:\cF^\bu\ra(\cF^\bu)^\vee[m-2]$ with~$\th^\vee=\th[2-m]$.

If $S$ is an algebraic space then $h^1(\bL_S)=0$, so $h^1(\cF^\bu)=0$, and $\cF^\bu$ is perfect in the interval $[1-m,0]$. The isomorphism $\th:\cF^\bu\ra(\cF^\bu)^\vee[m-2]$ then implies that $\cF^\bu$ is perfect in~$[2-m,0]$.

If $(\bS,\om)$ is a $(2-m)$-shifted symplectic derived Artin stack, as in \cite{PTVV}, with classical truncation $S=t_0(\bS)$ then $\bL_i:i^*(\bL_\bS)\ra\bL_S$ with $\th=i^*(\om\cdot{})$ is an $m$-Calabi--Yau obstruction theory on $S$. Thus, if $\A$ is an $m$-Calabi--Yau category and $\M,\bs\M$ the classical and derived moduli stacks of objects in $\A$ then $\bL_i:i^*(\bL_{\bs\M})\ra\bL_\M$ is a natural $m$-Calabi--Yau obstruction theory on $\M$.
\end{dfn}

\begin{rem}
\label{co1rem5}
Although we define $m$-Calabi--Yau obstruction theories for all $m\ge 0$, one should only expect to be able to define virtual classes when $m=2,3,4$. In the case $m=4$ Definition \ref{co1def5} is probably too weak, we need extra assumptions on $\phi:\cF^\bu\ra\bL_S$ to define virtual classes.	
\end{rem}

\subsubsection{Properties of virtual classes used in our proofs}
\label{co1.5.2}

We are (optimistically) hoping for versions of Theorem \ref{co1thm3} in which we can replace $\A$ by a 3- or 4-Calabi--Yau abelian category, and replace Behrend--Fantechi obstruction theories and virtual classes by 3- or 4-Calabi--Yau obstruction theories and virtual classes, and have the proofs in Chapters \ref{co9}--\ref{co11} still hold, with only a few changes, so we avoid doing lots of extra work.

We highlight five important properties of obstruction theories and virtual classes used in Chapters \ref{co9}--\ref{co11}, which need to be extended to $m$-Calabi--Yau obstruction theories for $m=3,4$:
\begin{itemize}
\setlength{\itemsep}{0pt}
\setlength{\parsep}{0pt}
\item[(i)] If $X$ is a {\it proper algebraic\/ $\C$-space\/} (not just a projective $\C$-scheme) with an obstruction theory of the desired type, there should exist a virtual class $[X]_\virt\in H_*(X)$.
\item[(ii)] As in \S\ref{co25}, if $X,Y$ are proper algebraic $\C$-spaces with obstruction theories of the desired type and $f:X\ra Y$ is smooth with fibre $\CP^m$, there should be a {\it pushforward formula\/} relating $[X]_\virt$ and $[Y]_\virt$ in~$H_*(Y)$.
\item[(iii)] As in \S\ref{co26}, a $\bG_m$-{\it localization formula\/} should hold for virtual classes on proper algebraic $\C$-spaces.
\item[(iv)] The proof of Theorem \ref{co1thm3} involves auxiliary categories $\baA$ in exact sequences $0\ra\A\ra\baA\ra\modCQ\ra 0$. To make Theorem \ref{co1thm3} work in the $m$-Calabi--Yau case, ideally we need to give $\baA$ (or a subcategory $\baB$) the structure of an $m$-Calabi--Yau category, which is not obvious. This affects the definition of the derived moduli stacks $\bs\baM,\bs\baM^\pl$ of objects in~$\baA$.
\item[(v)] The Lie bracket we define on $\check H_*(\M^\pl)$ by some vertex algebra construction has to correspond to wall-crossing behaviour of the invariants $[\M_\al^\ss(\tau)]_\inv$. It is enough for this to hold for `simple walls'. The key points in the proof are Propositions \ref{co9prop5} and \ref{co10prop14}, which write the change in $[\baM_{(\al,\bs d)}^\ss(\bar\tau)]_\virt$ on crossing a `simple wall', derived using a $\bG_m$-localization formula as in (ii), in terms of the Lie bracket on~$\check H_*(\M^\pl)$.
\end{itemize}

\subsubsection{3-Calabi--Yau virtual classes}
\label{co1.5.3}

Let $X$ be a proper algebraic $\C$-space with a 3-Calabi--Yau obstruction theory $\phi:\cF^\bu\ra\bL_X$. Then Definition \ref{co1def5} shows that $\cF^\bu$ is perfect in $[-1,0]$. Therefore $\phi:\cF^\bu\ra\bL_X$ is an obstruction theory in the usual Behrend--Fantechi sense, and a virtual class $[X]_\virt\in H_*(X)$ is defined as in \S\ref{co24}. The usual $\bG_m$-localization formula from \S\ref{co26} holds. This deals with issues (i),(iii) in~\S\ref{co1.5.2}. 

It remains to deal with (ii),(iv),(v). We discuss (v) in~\S\ref{co1.5.4}.

If $X$ is a general Artin stack with a 3-Calabi--Yau obstruction theory $\phi:\cF^\bu\ra\bL_X$, then $\cF^\bu$ is only perfect in $[-2,1]$, not $[-1,1]$, and so is not an obstruction theory in the sense of~\S\ref{co24}.

Kiem--Li \cite{KiLi1} prove a wall-crossing formula for virtual classes of 3-Calabi--Yau obstruction theories on crossing `simple walls'.

\subsubsection{Vertex Lie algebras from odd Calabi--Yau categories}
\label{co1.5.4}

A ({\it graded\/}) {\it vertex Lie algebra\/} $(V_*,e^{zD},Y_{<0})$ is a weakening of the (graded) vertex algebras $(V_*,\boo,e^{zD},Y)$ defined in \S\ref{co41}. We will not give the definition here, see Kac \cite[Def.~2.7b]{Kac}, Frenkel--Ben-Zvi \cite[\S 16.1]{FrBZ}, and~\cite{Joyc12}. 

The basic idea is that a vertex algebra has bilinear operations $(u,v)\mapsto u_n(v)$ for $n\in\Z$, but a vertex Lie algebra has only operations $(u,v)\mapsto u_n(v)$ for $n\in\N$ satisfying the same axioms as in a vertex algebra, and has no identity $\boo$. Since $Y(u,z)v=\sum_{n\in\Z}u_n(v)z^{-n-1}$, in passing from a vertex algebra $(V_*,\boo,e^{zD},Y)$ to a vertex Lie algebra $(V_*,e^{zD},Y_{<0})$ we take $Y_{<0}=Y\mod O(z^{\ge 0})$, that is, we remember only the negative powers $z^{-1},z^{-2},\ldots,$ from $Y$. Proposition \ref{co1prop1} also works for vertex Lie algebras. 

Suppose $\A$ is a $\C$-linear $m$-Calabi--Yau category for $m$ odd, as in Definition \ref{co1def4}, and write $\M$ for the moduli stack of objects in $\A$. First consider the vertex algebra $(V_*,\boo,e^{zD},Y)$ in Definition \ref{co1def2} for $V_*=\hat H_*(\M)$, with $\cE^\bu=(\cExt^\bu)^\vee$ as in Remark \ref{co1rem1}. In equation \eq{co1eq5} we have
\begin{equation*}
c_i(\cE_{\al,\be}^\bu\op \si_{\al,\be}^*(\cE^\bu_{\be,\al})^\vee)=c_i(\cE_{\al,\be}^\bu\op \cE_{\al,\be}^\bu[-m])=0,
\end{equation*}
since $(\cExt^\bu)^\vee\cong\si^*(\cExt^\bu)[m]$ as in Definition \ref{co1def4} and $m$ is odd. Also $\chi(\al,\be)=-\chi(\be,\al)$ so $z^{\chi(\al, \be)+\chi(\be,\al)-i+j}=z^{-i+j}$. Together these imply that there are no negative powers of $z$ in \eq{co1eq5}, so $u_n(v)=0$ for all $u,v$ and $n\ge 0$. Thus the Lie bracket $[\,,\,]$ on $\check H_*(\M^\pl)$ in Theorem \ref{co1thm2} is zero. So the vertex and Lie algebras in \S\ref{co1.2} and Chapter \ref{co4} are essentially trivial in the odd Calabi--Yau case.

We outline the construction of a graded vertex Lie algebra structure on $H_*(\M)$, inducing a graded Lie algebra structure on $H_*(\M^\pl)$ as in Theorem \ref{co1thm2}. Details and proofs will be given in the final version of \cite{Joyc12}. The definition of $e^{zD}:H_*(\M)\ra H_*(\M)[[z]]$ in Definition \ref{co1def2} is unchanged. To define $Y_{<0}$, we replace \eq{co1eq5} by 
\e
\begin{split}
&Y_{<0}(u,z)v=Y_{<0}(z)(u\ot v)= (-1)^{\chi(\al,\be)} \sum\nolimits_{i,j\ge 0}(-1)^ii!z^{-1-i+j}\cdot{}\\
&\bigl(\Phi_{\al,\be}\ci(\Psi_\al\t\id_{\M_\be})\bigr)_*\bigl(t^j\bt ((u \bt v) \cap \ch_i(\cE^\bu_{\al,\be})\bigr) \mod O(z^{\ge 0}).
\end{split}
\label{co1eq22}
\e
We can derive \eq{co1eq22} from \eq{co1eq5} by a limiting argument, they are related.

In contrast to \eq{co1eq6}, the correct grading shift on $H_*(\M)$ is $\hat H_n(\M)=H_{n-2}(\M)$. And in contrast to \eq{co1eq8} we take $\check H_n(\M^\pl)=H_n(\M^\pl)$, that is, there is no grading shift on $H_*(\M^\pl)$ as a graded Lie algebra.

If $u\in H_0(\M_\al)$ and $v\in H_0(\M_\be)$ then as $\ch_0(\cE^\bu_{\al,\be})=\chi(\al,\be)\cdot 1_{\M_\al\t\M_\be}$, equation \eq{co1eq22} implies that
\begin{equation*}
u_0(v)=(-1)^{\chi(\al,\be)}\chi(\al,\be)\cdot (\Phi_{\al,\be})_*(u\bt v)\in H_0(\M_{\al+\be}).
\end{equation*}
This descends to a Lie bracket on $H_0(\M^\pl)$. In the 3-Calabi--Yau case, this corresponds to the Lie bracket used in Joyce--Song \cite[Def.~5.13]{JoSo} for wall-crossing formulae of Donaldson--Thomas invariants, so this solves issue (v) in \S\ref{co1.5.2} in the 3-Calabi--Yau case.

This appearance of vertex Lie algebras in Donaldson--Thomas theory of Calabi--Yau 3-folds appears to the author to be new and interesting.

\begin{quest}
\label{co1quest2}
Is there an explanation in String Theory for the vertex Lie algebras above?
\end{quest}

These ideas may be useful in the study of $G$-{\it equivariant\/} Donaldson--Thomas invariants, say for the action of an algebraic group $G$ on a noncompact Calabi--Yau 3-fold $X$, where $G$ need not preserve the holomorphic volume form of $X$. Thomas' equivariant Vafa--Witten invariants \cite{Thom2} are of this type for $G=\bG_m$ and $X=K_S$ for $S$ a projective surface. But see Remark \ref{co5rem3}(a) below for a difficulty in applying our results to~\cite{Thom2}.

The graded vertex Lie algebra and Lie algebra above should be lifted to $H_*^G(\M)$ and $H_*^G(\M^\pl)$ as in \S\ref{co45}. The Lie bracket on $H_0^G(\M^\pl)$ may be much more complex than that on~$H_0(\M^\pl)$.

A 3-Calabi--Yau version of this book might be used to extend the Donaldson--Thomas theory of Joyce--Song \cite{JoSo} to work over other fields.

\subsubsection{4-Calabi--Yau virtual classes}
\label{co1.5.5}

Suppose $X$ is a Calabi--Yau 4-fold, and write $\M$ and $\bcM$ for the classical and derived moduli stacks of objects in $\coh(X)$, with inclusion $i:\M\hookra\bcM$. Then $\bcM$ has a $-2$-{\it shifted symplectic structure\/} in the sense of Pantev--To\"en--Vaqui\'e--Vezzosi \cite{PTVV}. Also $\bL_i:i^*(\bL_{\bcM})\ra\bL_\M$ is a 4-Calabi--Yau obstruction theory on $\M$, a classical truncation of the $-2$-shifted symplectic structure on~$\bcM$.

Borisov--Joyce \cite{BoJo} defined virtual classes for proper $-2$-shifted symplectic derived $\C$-schemes, using Derived Differential Geometry \cite{Joyc8,Joyc9,Joyc10,Joyc11}. This should have started a new subject of Donaldson--Thomas type invariants of Calabi--Yau 4-folds (DT4 invariants). Unfortunately, \cite{BoJo} was too horrible for anyone to do anything with. Recently, Oh--Thomas \cite{OhTh} gave a new, algebro-geometric definition of 4-Calabi--Yau virtual classes, equivalent to \cite{BoJo}, in the style of Behrend--Fantechi \cite{BeFa1}. The author expects \cite{OhTh} will become standard.

To define a 4-Calabi--Yau virtual class we need a choice of {\it orientation}, defined in Borisov--Joyce \cite[\S 2.4]{BoJo}. They make sense for the even Calabi--Yau case.

\begin{dfn}
\label{co1def6}
Let $X$ be an Artin $\C$-stack with an $m$-Calabi--Yau obstruction theory $\phi:\cF^\bu\ra\bL_X$, $\th:\cF^\bu\ra(\cF^\bu)^\vee[m-2]$ for $m$ even. Then we have a determinant line bundle $\det(\cF^\bu)\ra X$. As $m$ is even, $\th$ induces an isomorphism $\det\th:\det\cF^\bu\ra(\det\cF^\bu)^*$. An {\it orientation\/} for $(X,\phi,\th)$ is a choice of isomorphism $\la:\det\cF^\bu\ra\O_X$ with $\la^*\ci\la=\det\th$.

Here $\la$ is basically a square root of $\det\th$. Locally on $X$ in the \'etale topology there are two choices for $\la$, and there is a principal $\Z_2$-bundle $P\ra X$ parametrizing choices of $\la$. Then $(X,\phi,\th)$ is {\it orientable\/} if $P$ is trivializable, and an orientation is a trivialization $P\cong X\t\Z_2$.
\end{dfn}

Orientations of this type were studied in Differential Geometry by Joyce--Tanaka--Upmeier \cite{JTU}. Building on \cite{JTU}, Cao--Gross--Joyce \cite[Cor.~1.17]{CGJ} prove that if $X$ is a projective Calabi--Yau 4-fold and $\M$ is the moduli stack of objects in $\coh(X)$ (or in $D^b\coh(X)$) with its 4-Calabi--Yau obstruction theory, then $\M$ is orientable. This was extended to $\coh_\cs(X)$ for $X$ a quasiprojective Calabi--Yau 4-fold by the author's PhD student Arkadij Bojko \cite{Bojk1}.

Oh--Thomas \cite{OhTh} define their virtual class $[\M]_\virt$ only when $\M$ is a projective moduli scheme of Gieseker stable sheaves on a Calabi--Yau 4-fold $X$. However, Kiem--Park \cite[\S 4]{KiPa} provide an alternative definition which works for proper algebraic spaces and Deligne--Mumford stacks, dealing with issue (i) in \S\ref{co1.5.2}. Park \cite[Prop.~4.2]{Park} proves a virtual pushforward formula relevant to \S\ref{co1.5.2}(ii). Oh--Thomas \cite[\S 7]{OhTh} prove a $\bG_m$-localization formula for their 4-Calabi--Yau virtual classes. Park \cite[Prop.~0.9]{Park} extends this to Deligne--Mumford stacks. These address \S\ref{co1.5.2}(iii). Issue (v) will be covered in~\S\ref{co1.5.6}.

Thus it appears that most of the main ingredients necessary to extend Theorem \ref{co1thm3} to the 4-Calabi--Yau case are already available, though there are probably still technical problems to overcome. 

Computations of DT4 invariants in examples (often conjectural) can be found in Bojko \cite{Bojk2}, Cao \cite{Cao1,Cao2}, Cao--Kool \cite{CaKo1,CaKo2}, Cao--Kool--Monavari \cite{CKM1}, Cao--Leung \cite{CaLe}, Cao--Maulik--Toda \cite{CMT1,CMT2}, and Cao--Toda \cite{CaTo1,CaTo2,CaTo3}.

\subsubsection{Vertex algebras from even Calabi--Yau categories}
\label{co1.5.6}

Let $\A$ be a $\C$-linear $m$-Calabi--Yau category for $m$ even, as in Definition \ref{co1def4}, and write $\M$ for the moduli stack of objects in $\A$, with its $m$-Calabi--Yau obstruction theory. Following \cite{Joyc12}, we outline the construction of a graded vertex Lie algebra structure on $H_*(\M)$ as in Theorem \ref{co1thm1}, inducing a graded Lie algebra structure on $H_*(\M^\pl)$ as in Theorem~\ref{co1thm2}. 

As in Assumption \ref{co1ass1}(e) and Remark \ref{co1rem1} we should have a biadditive Euler form $\chi:K(\A)\t K(\A)\ra \Z$ with $\chi(\al,\be)=\rank\cExt^\bu_{\al,\be}$. Serre duality and $m$ even implies that $\chi(\al,\be)=\chi(\be,\al)$.

Suppose $\M$ is orientable in the sense of Definition \ref{co1def6}, and choose an orientation on $\M$. This is possible when $\A=\coh(X)$ for $X$ a Calabi--Yau 4-fold by \cite[Cor.~1.17]{CGJ}. If the $\M_\al$ are connected there are two choices for the orientation on $\M_\al$ for each $\al\in C(\A)\amalg\{0\}$. Write $\om_\al$ for the orientation on~$\M_\al$.

By \cite[\S 2.5]{JTU} and \cite[Th.~1.15]{CGJ}, from the orientations $\om_\al,\om_\be$ on $\M_\al,\M_\be$ we can construct an orientation $\om_\al\star\om_\be$ on $\M_{\al+\be}$ using the direct sum map $\Phi_{\al,\be}:\M_\al\t\M_\be\ra\M_{\al+\be}$ (or at least, on the connected components of $\M_{\al+\be}$ that intersect $\Im\Phi_{\al,\be}$). As in \cite[\S 2.5]{JTU} these satisfy
\e
\begin{split}
\om_\be\star\om_\al&=(-1)^{\chi(\al,\be)+\chi(\al,\al)\chi(\be,\be)}\om_\al\star\om_\be,\\
(\om_\al\star\om_\be)\star\om_\ga&=\om_\al\star(\om_\be\star\om_\ga).
\end{split}
\label{co1eq23}
\e

Let us suppose for simplicity that there exist signs $\ep_{\al,\be}=\pm 1$ for all $\al,\be\in C(\A)\amalg\{0\}$ such that $\om_{\al+\be}=\ep_{\al,\be}\cdot \om_\al\star\om_\be$. (That is, the sign comparing $\om_{\al+\be}$ and $\om_\al\star\om_\be$ should not depend on the connected component of $\M_{\al+\be}$, if $\M_{\al+\be}$ is not connected.) If $\A=\coh(X)$, and the orientations $\om_\al$ are pulled back from differential-geometric orientations on connection moduli spaces using the method of \cite{CGJ,JTU}, and $K(\A)\!\subseteq\! K^0(X^\ran)$, this is automatic. Then \eq{co1eq23}~yields
\e
\begin{split}
\ep_{\al,\be}\cdot\ep_{\be,\al}&=(-1)^{\chi(\al,\be)+\chi(\al,\al)\chi(\be,\be)},\\
\ep_{\al,\be}\cdot\ep_{\al+\be,\ga}&=\ep_{\al,\be+\ga}\cdot\ep_{\be,\ga}.
\end{split}
\label{co1eq24}
\e

The definitions of $\boo\!\in\! H_0(\M_0)$ and $e^{zD}:H_*(\M)\!\ra\! H_*(\M)[[z]]$ in Definition \ref{co1def2} are unchanged. To define $Y$, we write $\cE^\bu\!=\!(\cExt^\bu)^\vee$ as in Remark \ref{co1rem1}, and for $u\in H_a(\M_\al)\subset H_*(\M)$ and $v\in H_b(\M_\be)\subset H_*(\M)$, we replace \eq{co1eq5}~by 
\e
\begin{split}
&Y(u,z)v=Y(z)(u\ot v)= \ep_{\al,\be}(-1)^{a\chi(\be,\be)}\sum\nolimits_{i,j\ge 0} z^{\chi(\al,\be)-i+j}\cdot{}\\
&\bigl(\Phi_{\al,\be}\ci(\Psi_\al\t\id_{\M_\be})\bigr)_*\bigl(t^j\bt ((u \bt v) \cap c_i(\cE_{\al,\be}^\bu)\bigr).
\end{split}
\label{co1eq25}
\e
The analogues of \eq{co1eq6} and \eq{co1eq8} are 
\begin{align*}
\hat H_n(\M_\al)&= H_{n-\chi(\al,\al)}(\M_\al),& \hat H_n(\M)&=\ts\bigop_{\al\in K(\A)}\hat H_n(\M_\al),\\
\check H_n(\M_\al^\pl)&=H_{n+2-\chi(\al,\al)}(\M_\al^\pl),& \check H_n(\M^\pl)&=\ts\bigop_{\al\in K(\A)}\check H_n(\M_\al^\pl).
\end{align*}
Then by \cite{Joyc12} the analogues of Theorems \ref{co1thm1} and \ref{co1thm2} hold, defining graded vertex algebra and graded Lie algebra structures on $\hat H_*(\M)$ and~$\check H_*(\M^\pl)$.

We claim that when $m=4$, this Lie bracket on $\check H_{\rm even}(\M^\pl)$ is what we should use for the 4-Calabi--Yau version of Theorem \ref{co1thm3}.

\begin{rem}
\label{co1rem6}
{\bf(a)} Definition \ref{co1def2} is a special case of the construction above with $\cE^\bu\op\si^*((\cE^\bu)^\vee),\chi(\al,\be)+\chi(\be,\al),(-1)^{\chi(\al,\be)}$ in place of~$\cE^\bu,\chi(\al,\be),\ep_{\al,\be}$.
\smallskip

\noindent{\bf(b)} Equation \eq{co1eq24} is needed to show that $Y$ in \eq{co1eq25} yields a graded vertex algebra. The orientations we choose on each $\M_\al$ are arbitrary, and changing them changes the sign of $[\M_\al^\ss(\tau)]_\inv$ in the analogue of Theorem \ref{co1thm3}. These arbitrary sign changes do not affect \eq{co1eq14}--\eq{co1eq15}, as the choices of orientations are built into the Lie bracket $[\,,\,]$ via the signs $\ep_{\al,\be}$ in~\eq{co1eq25}.
\end{rem}

\section{Background from Algebraic Geometry}
\label{co2}

We begin with background on Artin stacks and their homology, Derived Algebraic Geometry and derived stacks, and Behrend--Fantechi obstruction theories and virtual classes \cite{BeFa1} and their properties. The notion of equivariant homology $H^G_*(X)$ of an Artin stack $X$ that we introduce in \S\ref{co23} is new.

\subsection{Artin stacks and their homology}
\label{co21}

\subsubsection{Basics on Artin stacks}
\label{co211}

In this book we work with {\it Artin\/ $\C$-stacks}. Some good references on Artin stacks are G\'omez \cite{Gome}, Laumon and Moret-Bailly \cite{LaMo}, and Olsson \cite{Olss}. The next remark establishes notation, and recalls some facts we will need later.

\begin{rem}
\label{co2rem1}	
{\bf(a)} Artin\/ $\C$-stacks form a strict 2-category $\Art_\C$. That is, we have objects $X,Y$ (the Artin $\C$-stacks), and 1-morphisms $f,g:X\ra Y$, and 2-morphisms $\eta:f\Ra g$ in $\Art_\C$, which are all 2-isomorphisms. Objects $X$ in $\Art_\C$ are functors
\e
X:\{\text{commutative $\C$-algebras}\}\longra\{\text{groupoids}\}
\label{co2eq1}
\e
satisfying sheaf-type conditions.

We will rarely need the 2-category structure, so we generally implicitly work in the {\it homotopy category\/} $\Ho(\Art_\C)$, which has the same objects as $\Art_\C$, and has morphisms 2-isomorphism classes $[f]:X\ra Y$ of 1-morphisms $f:X\ra Y$, and will mostly not mention 2-morphisms of stacks. {\it Fibre products\/} of stacks will always mean 2-category fibre products, not fibre products in~$\Ho(\Art_\C)$.
\smallskip

\noindent{\bf(b)} We regard $\C$-{\it schemes}, {\it algebraic\/ $\C$-spaces}, and {\it Deligne--Mumford\/} $\C$-stacks as examples of Artin $\C$-stacks, forming full 2-subcategories $\Sch_\C\subset\AlgSp_\C \subset\DMSta_\C\subset\Art_\C$, where $\Sch_\C,\AlgSp_\C$ are discrete 2-categories (i.e.\ all 2-morphisms $\eta:f\Ra f$ are identities).
\smallskip

\noindent{\bf(c)} If $X$ is an Artin $\C$-stack, we write $X(\C)$ for the set of morphisms $[x]:\Spec\C\ra X$ in $\Ho(\Art_\C)$. Elements of $X(\C)$ are called $\C$-{\it points}, or {\it geometric points}, of $X$. If $\phi:X\ra Y$ is a 1-morphism then composition with $\phi$ induces a map of sets $\phi_*:X(\C)\ra Y(\C)$. If $X$ is an Artin $\C$-stack, each $\C$-point $[x]\in X(\C)$ has an {\it isotropy group\/} $\Iso_X(x)$, an algebraic $\C$-group, which is the group of 2-morphisms $\eta:x\Ra x$ in $\Art_\C$. Here $\Iso_X(s)$ is finite if $X$ is a Deligne--Mumford $\C$-stack, and trivial if $X$ is a $\C$-scheme or algebraic $\C$-space.

\smallskip

\noindent{\bf(d)} An important class of Artin $\C$-stacks are {\it quotient stacks\/} $[S/G]$, where $S$ is a $\C$-scheme and $G$ is an algebraic $\C$-group acting on $S$. When $X=[S/G]$, the $\C$-points are $G(\C)$-orbits $xG(\C)$ for $\C$-points $x\in S(\C)$, and the isotropy groups are~$\Iso_{[S/G]}(xG(\C))=\Stab_{G(\C)}(x)$.
\smallskip

\noindent{\bf(e)} On $\C$-schemes, algebraic $\C$-spaces, Deligne--Mumford $\C$-stacks, and Artin $\C$-stacks $X$, there are abelian categories of {\it coherent sheaves\/} $\coh(X)$, {\it quasicoherent sheaves\/} $\qcoh(X)$, and $\O_X$-modules $\O_X\text{-mod}$. See Hartshorne \cite[\S II.5]{Hart} and Huybrechts and Lehn \cite{HuLe2} for schemes, and Laumon and Moret-Bailly \cite{LaMo} and Olsson \cite{Olss} for the extension to stacks.

We may then consider derived categories (as triangulated categories) of $\coh(X),\qcoh(X),\O_X\text{-mod}$, as in Gelfand and Manin \cite{GeMa} and Huybrechts \cite{Huyb}. Two important categories are the bounded derived category $D^b\coh(X)$, and $D_\qcoh(X)$, which is defined to be the full triangulated subcategory of objects $\cE^\bu$ in the unbounded derived category $D(\O_X\text{-mod})$ whose cohomology sheaves $h^k(\cE^\bu)$ are quasicoherent for $k\in\Z$. 

These derived categories have the usual Grothendieck six functor formalism of inverse images $f^*$, direct images $f_*$, \ldots, as in Huybrechts \cite{Huyb} for $\C$-schemes and Liu and Zheng \cite{LiZh} for Artin $\C$-stacks, for instance.

An object $\cE^\bu$ in $D_\qcoh(X)$ is called {\it perfect\/} if (smooth) locally on $X$ it is equivalent to a complex of vector bundles $\cdots\ra\cF^k\,{\buildrel\d\over\longra}\cF^{k+1}\ra\cdots,$ with $\cF^k$ in degree $k\in\Z$, and $\cF^k=0$ for $\md{k}\gg 0$. We call $\cE^\bu$ {\it perfect in the interval\/} $[a,b]$ for integers $a\le b$ if this holds with $\cF^k=0$ unless $a\le k\le b$.
\smallskip

\noindent{\bf(f)} An Artin $\C$-stack $X$ has a {\it cotangent complex\/} $\bL_X$ in $D_\qcoh(X)$, and a 1-morphism $f:X\ra Y$ of Artin $\C$-stacks has a cotangent complex $\bL_{X/Y}$ in $D_\qcoh(X)$ constructed by Illusie \cite{Illu}, in a distinguished triangle in~$D_\qcoh(X)$:
\begin{equation*}
\xymatrix@C=35pt{ f^*(\bL_Y) \ar[r]^(0.55){\bL_f} &
\bL_X \ar[r] & \bL_{X/Y} \ar[r] & f^*(\bL_Y)[1]. }
\end{equation*}
The ({\it relative\/}) {\it tangent complex\/} is $\bT_X=(\bL_X)^\vee$, $\bT_{X/Y}=(\bL_{X/Y})^\vee$.

We have $h^i(\bL_X)=0$ if $i>1$ for general Artin $\C$-stacks $X$, and $h^i(\bL_X)=0$ if $i>0$ when $X$ is a scheme, algebraic space, or Deligne--Mumford stack. If $X$ is a {\it smooth\/} Artin $\C$-stack (e.g.\ if $X=[S/G]$ for $S$ a smooth $\C$-scheme) then $\bL_X$ is perfect in the interval $[0,1]$. But for singular $X$, $\bL_X$ is generally not perfect.
\end{rem}

\subsubsection{Homology and cohomology of Artin stacks}
\label{co212}

The next definition will be very important to us.

\begin{dfn}
\label{co2def1}
As in Simpson \cite{Simp} and Blanc \cite[\S 3.1]{Blan}, an Artin $\C$-stack $X$ has a {\it topological realization\/} $X^\top$, which is a topological space natural up to homotopy equivalence. Topological realization gives a functor $(-)^\top:\Ho(\Art_\C)\ab\ra\Top^{\bf ho}$ to the category $\Top^{\bf ho}$ of topological spaces with morphisms homotopy classes of continuous maps. 

Let $X$ be an Artin $\C$-stack. We define the {\it homology $H_*(X)$ of\/ $X$ with coefficients in\/} $\Q$ to be $H_*(X)=H_*(X,\Q)=H_*(X^\top,\Q)$, the usual homology of the topological space $X^\top$ over $\Q$. Similarly we define the cohomology $H^*(X)=H^*(X,\Q)=H^*(X^\top,\Q)$. These are sometimes called the {\it Betti (co)homology}, to distinguish them from other (co)homology theories of stacks. 

We will almost always take (co)homology of topological spaces or stacks {\it over the rationals\/} $\Q$, and when we omit the coefficient ring we mean it to be~$\Q$. 

(Co)homology has the same functorial properties as (co)homology of topological spaces. If $X$ is a quotient stack $[S/G]$, we have a homotopy equivalence
\begin{equation*}
X^\top\simeq (S^\ran\t EG^\ran)/G^\ran, 
\end{equation*}
where $EG^\ran\ra BG^\ran$ is a classifying space for the complex analytic topological group $G^\ran=G(\C)$. If $S$ is contractible (e.g.\ if $S$ is a point $*$ or an affine space $\bA^n$) this implies that $X^\top\simeq BG^\ran$.

An important example for us will be the quotient stack $[*/\bG_m]$, where $\bG_m=\C\sm\{0\}$ is the multiplicative group. As $B\bG_m\simeq \CP^\iy$, we have an isomorphism
\e
\begin{aligned}
H^*([*/\bG_m])&\cong \Q[[z]] &&\text{as $\Q$-algebras, with $\deg z=2$, so that}\\
H^{2n}([*/\bG_m])&=\an{z^n}_\Q, &&\text{normalized so that}\quad \ts\int_{\CP^n}z^n=1,
\end{aligned}
\label{co2eq2}
\e
where $\CP^n\hookra\CP^\iy$ is the standard inclusion. We will also write 
\e
\begin{aligned}
H_*([*/\bG_m])&\cong \Q[t] &&\text{with $\deg t=2$, and $[\CP^n]=t^n$, so that}\\
H_{2n}([*/\bG_m])&=\an{t^n}_\Q, &&\text{normalized so that}\quad z^n\cdot t^n=1.
\end{aligned}
\label{co2eq3}
\e
\end{dfn}

\begin{rem}
\label{co2rem2}
{\bf(a) (Borel--Moore homology.)} There are two main kinds of homology theory of topological spaces $X$, ordinary homology $H_*(X,\Q)$ and {\it Borel--Moore homology\/} $H_*^{\rm BM}(X,\Q)$ \cite{BoMo}, also called {\it locally finite homology}, or {\it homology with closed supports}, or {\it homology of the second kind}. 

Classes in $H_*(X,\Q)$ are supported on compact sets, and pushforwards $f_*:H_*(X,\Q)\ra H_*(Y,\Q)$ are defined for all continuous $f:X\ra Y$. In contrast, classes in $H_*^{\rm BM}(X,\Q)$ can be supported on non-compact sets, and pushforwards $f_*:H_*^{\rm BM}(X,\Q)\ra H_*^{\rm BM}(Y,\Q)$ are defined only for proper $f:X\ra Y$. We have $H^n(X,\Q)\cong H_n(X,\Q)^*$ and $H_n^{\rm BM}(X,\Q)\cong H^n_\cs(X,\Q)^*$, where $H^*_\cs(X,\Q)$ is {\it compactly-supported cohomology}. If $X$ is compact then~$H_*^{\rm BM}(X,\Q)=H_*(X,\Q)$.

Borel--Moore homology is not homotopy invariant. Thus we cannot define the Borel--Moore homology of an Artin $\C$-stack $X$ to be $H_*^{\rm BM}(X^\top,\Q)$, following Definition \ref{co2def1}, as $X^\top$ is only natural up to homotopy equivalence.
\smallskip

\noindent{\bf(b) (Chow homology.)} {\it Chow homology\/} $A_*(X)$ is an algebro-geometric homology theory, defined for schemes by Fulton \cite{Fult}, and for Artin stacks by Kresch \cite{Kres}. For schemes $X$, classes in $A_*(X)$ are equivalence classes of subvarieties of $X$. We will consider {\it rational Chow homology}, written~$A_*(X,\Q)$.

Chow homology is a homology theory of Borel--Moore type, as in {\bf(a)}. So classes in $A_*(X)$ need not be supported in proper subspaces of $X$, and pushforwards $f_*:A_*(X)\ra A_*(Y)$ are defined only for proper morphisms~$f:X\ra Y$.

This has two important consequences for us. Firstly, our theory of vertex and Lie algebra structures on homology of moduli stacks in Chapter \ref{co4}, the foundation of the rest of the book, {\it does not work for Chow homology}, as it needs pushforwards $f_*$ for non-proper $f$. This is an important reason for using Betti homology.

Secondly, we sometimes need to relate Chow homology and Betti homology, for example in mapping Behrend--Fantechi \cite{BeFa1} virtual cycles in Chow homology to Betti homology. When we do this we must be careful about properness issues. If $X$ is a proper algebraic $\C$-space or Deligne--Mumford $\C$-stack then there is a natural morphism
\e
A_d(X,\Q)\longra H_{2d}(X,\Q)=H_{2d}(X).
\label{co2eq4}
\e
For Artin $\C$-stacks $X$, there is no easy relation between $A_*(X,\Q)$ and~$H_*(X)$.
\end{rem}

\subsubsection{Further properties of homology of Artin stacks}
\label{co213}

We discuss two topics on homology of Artin $\C$-stacks that will only be used in the construction of equivariant homology $H_*^G(X)$ of Artin stacks in~\S\ref{co23}.

\begin{rem}
\label{co2rem3}
Let $X$ be an Artin $\C$-stack, and $Y$ a $\C$-scheme, and $\pi:X\ra Y$ a morphism. Choose a topological realization $\pi^\top:X^\top\ra Y^\ran$ of $\pi:X\ra Y$, where we take the topological realization of $Y$ to be its complex analytic topological space $Y^\ran$. Define the {\it homology of\/ $X$ with supports proper over\/} $Y$, as an inverse limit of relative homology groups of topological spaces, by
\e
H_*^{\pr/Y}(X)=\varprojlim\nolimits_{\text{$K\subseteq Y^\ran$ compact}}H_*\bigl(X^\top,(\pi^\top)^{-1}(Y^\ran\sm K);\Q\bigr).
\label{co2eq5}
\e
If $Y$ is proper, so $Y^\ran$ is compact, taking $K=Y^\ran$ gives~$H_*^{\pr/Y}(X)\cong H_*(X)$. 

For (co)homology theories of topological spaces with prescribed supports see Bredon \cite{Bred}, for example. In ordinary homology, each class $\al$ in $H_*(X)$ is supported in a compact subset of $X$. Classes $\al$ in $H_*^{\pr/Y}(X)$ are supported on subsets of $X$ which are proper over $Y$. The topological realization $X^\top$ is only natural up to homotopy equivalence, but this does not change \eq{co2eq5} up to canonical isomorphism, so $H_*^{\pr/Y}(X)$ is well defined.

When $\pi:X\ra Y$ is $\id_Y:Y\ra Y$ we have $H_*^{\pr/Y}(Y)=H_*^{\rm BM}(Y^\ran)$, the Borel--Moore homology, as in Remark~\ref{co2rem2}(a).
\end{rem}

\begin{rem}
\label{co2rem4}
Suppose we have a 2-Cartesian square of Artin $\C$-stacks:
\e
\begin{gathered}
\xymatrix@C=80pt@R=15pt{ *+[r]{W} \ar[r]_e \ar[d]^f & *+[l]{X} \ar[d]_g \\ *+[r]{Y} \ar[r]^h & *+[l]{Z,\!}}	
\end{gathered}
\label{co2eq6}
\e
where $Y,Z$ are smooth, separated $\C$-schemes of dimensions $l,m$, and $h$ is proper and either a regular embedding (if $l<m$) or smooth (if $l\ge m$). Then we may define a {\it Gysin morphism\/} or {\it umkehr map\/} mapping the wrong way on homology:
\e
e^!:H_n(X)\longra H_{n+2l-2m}(W).
\label{co2eq7}
\e
This depends on the diagram \eq{co2eq6}, not just on $e$. See Cohen--Klein \cite{CoKl} and Buoncristiano et al.\ \cite{BRS} for general theory of Gysin morphisms, and Fulton \cite[\S 6.1]{Fult} for the analogue in Chow homology.

For regular embeddings, we can interpret this as follows: there is a class $\Pd(Y)\in H^{2m-2l}(Z)$ Poincar\'e dual to $Y\hookra Z$. Then \eq{co2eq7} maps $e^!:\al\mapsto\al\cap g^*(\Pd(Y))$, with the proviso that as $\Pd(Y)$ is supported on $Y\subset Z$, so $\al\cap g^*(\Pd(Y))$ is supported in $g^{-1}(Y)=W$, and so can be understood as a class in $H_*(W)$ rather than in $H_*(X)$. For smooth morphisms, we may choose the topological realization $W^\top=X^\top\t_{Z^\ran}Y^\ran$, and  $e^\top:W^\top\ra X^\top$ is locally a fibration with fibres compact complex manifolds of dimension $l-m$. Then we construct \eq{co2eq7} using facts about homology of fibre bundles. 

Gysin morphisms are covariantly functorial, and commute with pushforwards in 2-Cartesian squares. We can also define them for homology with supports proper over a base, as in Remark \ref{co2rem3}, with different properness conditions.
\end{rem}

\subsection{Derived Algebraic Geometry and derived stacks}
\label{co22}

\subsubsection{Basics on higher and derived stacks}
\label{co221}

We will use To\"en and Vezzosi's theory of Derived Algebraic Geometry \cite{Toen1,Toen2,ToVe1,ToVe2}. See To\"en \cite{Toen1,Toen2} for a good introduction. The next remark establishes notation, and recalls some facts we will need later.

\begin{rem}
\label{co2rem5}
{\bf(a)} In \cite[\S 3]{Toen1}, \cite[\S 2.1]{ToVe2}, To\"en and Vezzosi define an $\iy$-category of {\it higher\/ $\C$-stacks}. As these include some badly behaved objects, we restrict attention in this book to what To\"en--Vezzosi call {\it higher Artin $\C$-stacks\/} \cite[\S 3]{Toen1}, \cite[\S 2.1]{ToVe2}, which we will just call {\it higher\/ $\C$-stacks}. These form an $\iy$-category $\HSta_\C$. Objects $X$ in $\HSta_\C$ are $\iy$-functors
\e
X:\{\text{commutative $\C$-algebras}\}\longra\{\text{simplicial sets}\}
\label{co2eq8}
\e
satisfying sheaf-type conditions, as for \eq{co2eq1}. 

By including $\{$groupoids$\}\hookra\{$simplicial sets$\}$ we get a full and faithful inclusion $\Art_\C\hookra\HSta_\C$, for $\Art_\C$ as in \S\ref{co211}, and we say a higher $\C$-stack {\it is an Artin\/ $\C$-stack\/} if it is equivalent to an object in the image of this.
\smallskip

\noindent{\bf(b)} As in To\"en--Vaqui\'e \cite{ToVa}, moduli stacks of objects $\M$ in many interesting $\C$-linear triangulated categories, such as $D^b\coh(X)$ for $X$ a smooth projective $\C$-scheme, exist as higher $\C$-stacks, but not as Artin $\C$-stacks.
\smallskip

\noindent{\bf(c)} The theory of homology and cohomology of Artin stacks in \S\ref{co212} also works for higher $\C$-stacks, as in \cite{Simp}, \cite[\S 3.1]{Blan}.
\smallskip

\noindent{\bf(d)} To\"en and Vezzosi also define the $\iy$-category of {\it derived\/ $\C$-stacks\/} \cite[Def.~2.2.2.14]{ToVe2}, \cite[Def.~4.2]{Toen1}. Again, as these include some badly behaved objects, we restrict attention in this book to what To\"en calls {\it derived Artin $\C$-stacks\/} \cite[\S 3.3]{Toen1}, \cite[Def.~3.2]{Toen2}, which we will just call {\it derived\/ $\C$-stacks}. These form an $\iy$-category $\DSta_\C$. Objects $\bX$ in $\DSta_\C$ are $\iy$-functors
\e
\bX:\{\text{simplicial commutative $\C$-algebras}\}\longra
\{\text{simplicial sets}\}
\label{co2eq9}
\e
satisfying sheaf-type conditions. As $\C$ is a field of characteristic zero, there is also an equivalent definition of derived stacks as $\iy$-functors
\e
\bX\!:\!\{\text{commutative differential graded $\C$-algebras}\}\!\longra\!
\{\text{simplicial sets}\}.
\label{co2eq10}
\e

\noindent{\bf(e)} There is a {\it classical truncation functor\/} $t_0:\DSta_\C\ra\HSta_\C$ from derived stacks to (higher) stacks, induced by composing $\iy$-functors $\bX$ in \eq{co2eq9} with the inclusion $\{$commutative $\C$-algebras$\}\hookra\{$simplicial commutative $\C$-algebras$\}$ to get $\iy$-functors \eq{co2eq8}. We use bold symbols $\bX,\bY,\bs f:\bX\ra\bY,\ldots$ to denote derived stacks and their morphisms, and the corresponding non-bold symbols $X,Y,f:X\ra Y,\ldots$ to denote their classical truncations $t_0(\bX),t_0(\bY),t_0(\bs f),\ldots.$

\smallskip

\noindent{\bf(f)} We call a derived $\C$-stack $\bX$ a {\it derived Artin\/ $\C$-stack\/} if its classical truncation $X=t_0(\bX)$ is an Artin $\C$-stack in the sense of {\bf(a)}. This conflicts with the notation of To\"en, who calls these {\it derived\/ $1$-Artin stacks\/} \cite[\S 3.3]{Toen2}. We write $\DArt_\C\subset\DSta_\C$ for the full $\iy$-subcategory of derived Artin $\C$-stacks. 
\smallskip

\noindent{\bf(g)} There is a fully faithful {\it inclusion functor\/} $i:\HSta_\C\hookra\DSta_\C$ left adjoint to $t_0:\DSta_\C\ra\HSta_\C$. As $i$ is fully faithful we can regard it as embedding $\HSta_\C$ as a full $\iy$-subcategory of $\DSta_\C$. Thus, we can regard Artin $\C$-stacks and higher $\C$-stacks as examples of derived $\C$-stacks.
\smallskip

\noindent{\bf(h)} Let $\bX$ be a derived stack, and $X=t_0(\bX)$ its classical truncation. The adjoint property of $i,t_0$ implies that there is a natural {\it inclusion morphism\/} $i_X:X\hookra\bX$, which will be very important later. We often write~$i_X=i$.

We can think of $\bX$ as an {\it infinitesimal formal thickening\/} of its classical substack $i_X(X)\subset\bX$, in a similar way to how a non-reduced $\C$-scheme $S$ is an infinitesimal formal thickening of its reduced subscheme~$S^\red\subset S$.
\smallskip

\noindent{\bf(i)} Just as in many moduli problems (such as moduli of objects in $\coh(X)$ for $X$ a smooth projective $\C$-scheme) we can form a classical moduli space $\M$ as an Artin $\C$-stack (or higher $\C$-stack), we can also usually form a derived moduli space $\bs\M$ as a derived Artin $\C$-stack (or derived $\C$-stack).

Here $\M$ measures families of objects in the moduli problem over a base classical $\C$-scheme $\Spec A$, for $A$ a commutative $\C$-algebra (as in \eq{co2eq8}), and $\bs\M$ measures families of objects in the moduli problem over a base derived $\C$-scheme $\mathop{\bs\Spec} A^\bu$, for $A^\bu$ a simplicial $\C$-algebra or commutative differential graded $\C$-algebra (as in \eq{co2eq9}--\eq{co2eq10}).

To\"en--Vaqui\'e \cite{ToVa} prove existence of derived moduli stacks $\bs\M$ in many important situations, including moduli of objects in $\coh(X)$ and~$D^b\coh(X)$.
\smallskip

\noindent{\bf(j)} The (co)homology $H_*(\bX),H^*(\bX)$ of a derived stack $\bX$ is defined to be the (co)homology $H_*(X),H^*(X)$ of its classical truncation $X=t_0(\bX)$.
\smallskip

\noindent{\bf(k)} Two conditions on derived stacks $\bX$ (meaning $\bX$ is particularly nice) will be important below: for $\bX$ to be {\it locally finitely presented\/} \cite[\S 2.2.3]{ToVe2} (also called {\it locally of finite presentation\/} in \cite[\S 3.1.1]{Toen1}, and {\it fp-smooth\/} in \cite[\S 4.4]{ToVe1}), and for $\bX$ to be {\it quasi-smooth\/} \cite[\S 4.4.3]{Toen1}, \cite[\S 1]{STV}. As these are connected to the cotangent complex $\bL_\bX$, we will discuss them further in \S\ref{co222} below. Quasi-smoothness will also be important in \S\ref{co24} on obstruction theories.

If $\bX$ is locally finitely presented, or quasi-smooth, then the classical truncation $X=t_0(\bX)$ is generally {\it not\/} locally finitely presented, or quasi-smooth, when considered as a derived stack.

All the derived moduli stacks $\bs\M$ in the examples in Chapters \ref{co6}--\ref{co8} will be locally finitely presented, and most will be quasi-smooth, or at least have a large quasi-smooth open substack $\bs\dM\subset\bs\M$.
\end{rem}

\subsubsection{Cotangent complexes of derived stacks}
\label{co222}

We discuss {\it cotangent complexes\/} of derived stacks, following To\"en and Vezzosi \cite[\S 1.4]{ToVe2}, \cite[\S 4.2.4--\S 4.2.5]{Toen1}, \cite[\S 4.1]{Toen2}, \cite{Vezz}, and Lurie \cite[\S 3.4]{Luri}. Note that as in Remark \ref{co2rem5}(a),(d), our notions of higher and derived $\C$-stacks correspond to `higher Artin $\C$-stacks' and `derived Artin $\C$-stacks' in To\"en \cite{Toen1,Toen2}, which are necessary to make cotangent complexes well behaved.

\begin{rem}
\label{co2rem6}
Let $\bX$ be a derived $\C$-stack. In classical Algebraic Geometry, cotangent complexes $\bL_X$ lie in $D(\qcoh(X))$, but in Derived Algebraic Geometry they lie in a triangulated category $L_\qcoh(\bX)$ defined by To\"en \cite[\S 3.1.7, \S 4.2.4]{Toen1}, which is a substitute for $D(\qcoh(X))$. These satisfy:
\begin{itemize}
\setlength{\itemsep}{0pt}
\setlength{\parsep}{0pt}
\item[(i)] If $\bX$ is a derived $\C$-stack then $L_\qcoh(\bX)$ is a triangulated category with a t-structure whose heart is the abelian category $\qcoh(\bX)$ of quasicoherent sheaves on $\bX$. If $\bX=i(X)$ is a higher $\C$-stack or Artin $\C$-stack then $L_\qcoh(X)\simeq D(\qcoh(X))$, but in general~$L_\qcoh(\bX)\not\simeq D(\qcoh(\bX))$.
\item[(ii)] If $\bs f:\bX\ra\bY$ is a morphism of derived $\C$-stacks it induces a pullback functor $\bs f^*:L_\qcoh(\bY)\ra L_\qcoh(\bX)$, analogous to the left derived pullback functor $f^*:D(\qcoh(Y))\ra D(\qcoh(X))$ in the classical case.
\item[(iii)] There is a notion of when a complex $\cE^\bu$ in $L_\qcoh(\bX)$ is {\it perfect}.
\end{itemize}

If $\bX$ is a derived $\C$-stack, then To\"en and Vezzosi \cite[\S 4.2.5]{Toen1}, \cite[\S 1.4]{ToVe2} or Lurie \cite[\S 3.2]{Luri} define an ({\it absolute\/}) {\it cotangent complex\/} $\bL_\bX$ in $L_\qcoh(\bX)$. If $\bs f:\bX\ra\bY$ is a morphism of derived stacks, they construct a morphism $\bL_{\bs f}:\bs f^*(\bL_\bY)\ra\bL_\bX$ in $L_\qcoh(\bX)$, and the ({\it relative\/}) {\it cotangent complex\/} $\bL_{\bX/\bY}$ is defined to be the cone on this, giving the distinguished triangle
\e
\xymatrix@C=30pt{ \bs f^*(\bL_\bY) \ar[r]^(0.55){\bL_{\bs f}} &
\bL_\bX \ar[r] & \bL_{\bX/\bY} \ar[r] & \bs f^*(\bL_\bY)[1]. }
\label{co2eq11}
\e
The {\it tangent complexes\/} of $\bX$ and $\bs f$ are the duals $\bT_\bX=(\bL_\bX)^\vee$ and $\bT_{\bX/\bY}=(\bL_{\bX/\bY})^\vee$. Here are some properties of these:
\begin{itemize}
\setlength{\itemsep}{0pt}
\setlength{\parsep}{0pt}
\item[(a)] If $\bX$ is a derived algebraic space then $h^k(\bL_\bX)=0$ for $k>0$. If it is a derived Artin $\C$-stack then $h^k(\bL_\bX)=0$ for $k>1$.
\item[(b)] Suppose we have a Cartesian diagram of derived $\C$-stacks:
\begin{equation*}
\xymatrix@C=90pt@R=14pt{ *+[r]{\bW} \ar[r]_(0.3){\bs f} \ar[d]^{\bs e}
& *+[l]{\bY} \ar[d]_{\bs h} \\
*+[r]{\bX} \ar[r]^(0.7){\bs g} & *+[l]{\bZ.\!\!{}} }
\end{equation*}
Then by To\"en and Vezzosi \cite[Lem.s 1.4.1.12 \& 1.4.1.16]{ToVe2} or Lurie \cite[Prop.~3.2.10]{Luri} we have {\it base change isomorphisms\/}:
\e
\bL_{\bW/\bY}\cong\bs e^*(\bL_{\bX/\bZ}),\qquad
\bL_{\bW/\bX}\cong\bs f^*(\bL_{\bY/\bZ}).
\label{co2eq12}
\e
Note that the analogous result for classical stacks requires $g$ or $h$ to be flat, but \eq{co2eq12} holds with no flatness assumption. Because of \eq{co2eq11}--\eq{co2eq12}, cotangent complexes of derived stacks are often computable in examples.
\item[(c)] Let $\bs f:\bX\ra\bY$ be a smooth morphism of derived Artin $\C$-stacks. Then $\bL_{\bX/\bY}$ is perfect in the interval $[0,1]$.
\item[(d)] Let $\bX$ be a {\it locally finitely presented\/} derived $\C$-stack, as in Remark \ref{co2rem5}(k). Then $\bL_\bX$ is a perfect complex in $L_\qcoh(\bX)$. We say that $\bX$ has {\it virtual dimension\/} $n$, written $\vdim\bX=n$, if $\rank\bL_\bX=n$.
\item[(e)] A derived $\C$-stack $\bX$ is called {\it quasi-smooth\/} \cite[\S 4.4.3]{Toen1}, \cite[\S 1]{STV} if it is locally finitely presented and $\bL_\bX$ is perfect in the interval $[-1,\iy)$.  

To prove $\bX$ is quasi-smooth, it is sufficient to show $\bX$ is locally finitely presented and $i_X^*(\bL_\bX)$ is perfect in the interval $[-1,\iy)$ on~$X=t_0(\bX)$.
\item[(f)] Similarly, a morphism $\bs f:\bX\ra\bY$ is {\it quasi-smooth\/} if it is locally finitely presented and $\bL_{\bs f}$ is perfect in the interval $[-1,\iy)$. 
\end{itemize}
\end{rem}

\subsection{Equivariant (co)homology of Artin stacks}
\label{co23}

Let $X$ be an Artin $\C$-stack and $G$ a linear algebraic $\C$-group acting on $X$. We now explain what we mean by the {\it equivariant (co)homology\/} $H^G_*(X),H^*_G(X)$ of $X$ over $\Q$. Note that our notion of equivariant homology is {\it not obvious}, and {\it not equivalent to some other definitions}, so readers already familiar with some form of equivariant (co)homology are advised not to skip this section.

Since we can form the quotient $X/G$ as an Artin stack, we could just consider the (co)homology $H_*(X/G),H^*(X/G)$. Our definition below has $H^*_G(X)=H^*(X/G)$, but in general $H^G_*(X)\not\cong H_*(X/G)$. The reasons we chose this particular version of equivariant homology are:
\begin{itemize}
\setlength{\itemsep}{0pt}
\setlength{\parsep}{0pt}
\item[(a)] If $X$ is a proper algebraic $\C$-space with a $G$-equivariant obstruction theory, the virtual class $[X]_\virt$ exists in $H^G_*(X)$, as in \S\ref{co24}.
\item[(b)] It is compatible with $\bG_m$-equivariant localization in \S\ref{co26}.
\item[(c)] As in \S\ref{co45}, it is the correct notion to generalize the geometric constructions of vertex algebras in \S\ref{co42}, and Lie algebras in \S\ref{co43}, to the equivariant case.
\item[(d)] In \S\ref{co54} we use it to extend our theory of invariants to the equivariant case.
\end{itemize}
None of (b)--(d) work with $H_*(X/G)$ in place of~$H^G_*(X)$.

The next definition may be new, but is based on similar constructions in Totaro \cite{Tota}, Edidin--Graham \cite{EdGr2}, and Heller--Malag\'on-L\'opez~\cite{HeMa}.

\begin{dfn}
\label{co2def2}
Let $G$ be a linear algebraic $\C$-group. Following \cite[Def.~10]{HeMa}, a {\it good system of representations\/} for $G$ is a set $\{(V_j,U_j):j=0,1,\ldots\}$, such that:
\begin{itemize}
\setlength{\itemsep}{0pt}
\setlength{\parsep}{0pt}
\item[(i)] $V_j$ is a finite-dimensional $G$-representation over $\C$.
\item[(ii)] $U_j\subseteq V_j$ is a $G$-invariant Zariski open set, such that $G$ acts freely on $U_j$, and the quotient $U_j/G$ exists as a smooth separated $\C$-scheme.
\item[(iii)] For all $j$ we have $V_{j+1}=V_j\op W_j$, for $G$-subrepresentations~$V_j,W_j\subseteq V_{j+1}$.
\item[(iv)] For all $j$ we have a closed inclusion $U_j\subseteq U_{j+1}$,  which factors as $U_j=U_j\t\{0\}\subseteq U_j\t W_j\subseteq U_{j+1}$ in $V_{j+1}=V_j\t W_j$.
\item[(v)] $\lim_{j\ra\iy}\dim V_j=\iy$.
\item[(vi)] If $j<k$ then $\codim_{V_j}(V_j\sm U_j)<\codim_{V_k}(V_k\sm U_k)$.
\end{itemize}

As in \cite[Rem.~11]{HeMa} and \cite[Rem.~1.4]{Tota}, such a good system exists for any $G$, and we fix one. Let $X$ be an Artin $\C$-stack with an action of $G$. Then for each $j=0,1,\ldots$ we can form the quotient $(X\t U_j)/G$, as an Artin $\C$-stack. As in \eq{co2eq6}, we have a 2-Cartesian diagram of Artin $\C$-stacks:
\begin{equation*}
\xymatrix@C=120pt@R=15pt{ *+[r]{(X\t U_j)/G\,\,} \ar@{^{(}->}[r]_{\inc_j} \ar[d]^{\pi_j} &  *+[l]{(X\t U_{j+1})/G} \ar[d]_{\pi_{j+1}} \\
*+[r]{U_j/G\,\,} \ar@{^{(}->}[r]^{\inc'_j} &  *+[l]{U_{j+1}/G,\!} }
\end{equation*}
with the horizontal morphisms inclusions, and the vertical projections.

Since $\inc_j':U_j/G\hookra U_{j+1}/G$ is a closed embedding of smooth separated $\C$-schemes, as in \eq{co2eq7} we can form a Gysin morphism on homology
\e
\inc_j^!:H^{\pr/(U_{j+1}/G)}_{n+2\dim U_{j+1}/G}((X\t U_{j+1})/G)\longra H^{\pr/(U_j/G)}_{n+2\dim U_j/G}((X\t U_j)/G),
\label{co2eq13}
\e
for $H_*^{\pr/Y}(-)$ as in \eq{co2eq5}. For each $n\in\Z$, define the {\it equivariant homology group\/}
\e
H_n^G(X)=\varprojlim\nolimits_{j\ra\iy}H^{\pr/(U_j/G)}_{n+2\dim U_j/\dim G}((X\t U_j)/G),
\label{co2eq14}
\e
where the inverse limit of $\Q$-vector spaces is taken using the maps \eq{co2eq13}. One can show using the method of \cite[Th.~1.1]{Tota} and \cite[Th.~16]{HeMa} that $H_n^G(X)$ is independent of the choice of $\{(V_j,U_j):j\ge 1\}$, up to canonical isomorphism. For each $n=0,1,2,\ldots,$ define the {\it equivariant cohomology group\/} $H^n_G(X)$ by
\begin{equation*}
H^n_G(X)=H^n(X/G).
\end{equation*}

We can define the usual operations on equivariant (co)homology:
\begin{itemize}
\setlength{\itemsep}{0pt}
\setlength{\parsep}{0pt}
\item[(a)] Let $X,Y$ be Artin $\C$-stacks with actions of $G$, and $f:X\ra Y$ be a $G$-equivariant morphism. Then $f\t\id:X\t U_j\ra Y\t U_j$ descends to $f_j:(X\t U_j)/G\ra (Y\t U_j)/G$ and induces $(f_j)_*:H_*((X\t U_j)/G)\ra H_*((Y\t U_j)/G)$. Then $(f_j)_*,(f_{j+1})_*$ commute with the Gysin morphisms \eq{co2eq13}, and so induce a {\it pushforward map\/} $f_*:H_*^G(X)\ra H_*^G(Y)$. Pushforwards are covariantly functorial.

Similarly, $f$ descends to $f'':X/G\ra Y/G$, and so induces a {\it pullback map\/} $f^*:H^*_G(Y)\ra H^*_G(X)$. Pullbacks are contravariantly functorial.
\item[(b)] For $\al\in H_n^G(X)$, $\be\in H^m_G(X)$, there is a {\it cap product\/} $\al\cap\be\in H_{n-m}^G(X)$, characterized uniquely by $\Pi^G_j(\al\cap\be)=\Pi^G_j(\al)\cap\Pi_{X/G}^*(\be)$, where
\e
\Pi^G_j:H_*^G(X)\longra H^{\pr/(U_j/G)}_{*+2\dim V_j-2\dim G}((X\t U_j)/G)
\label{co2eq15}
\e
is the projection from the inverse limit, and $\Pi_{X/G}:(X\t U_j)/G\ra X/G$ the natural morphism. We define {\it cup products\/} in $H^*_G(X)$ to be usual cup products in $H^*(X/G)$. The {\it identity\/} is $1\in H^0_G(X)=H^0(X/G)$. Cup and cap products and identities have the usual compatibilities with each other, and with pushforwards and pullbacks.
\item[(c)] When $X=*$ is the point $\Spec\C$ with trivial $G$-action, by Poincar\'e duality for the possibly noncompact complex manifold $(*\t U_j)/G=U_j/G$ we have 
\begin{equation*}
H^{\pr/(U_j/G)}_{n+2\dim U_j/\dim G}((*\t U_j)/G)\cong H^{-n}(U_j/G).
\end{equation*}
But when $j\gg 0$ we can show using (vi) above that
\begin{equation*}
H^{-n}(U_j/G)\cong H^{-n}(V_j/G)\cong H^{-n}([*/G])=H^{-n}_G(*),
\end{equation*}
since $V_j$ is contractible. Thus taking the limit $j\ra\iy$ in \eq{co2eq14} we see that
\begin{equation*}
H_n^G(*)\cong H^{-n}_G(*)=H^{-n}([*/G]).
\end{equation*}
More generally, if $G$ acts trivially on $X$ we can use the K\"unneth isomorphism $H_*((X\t U_j)/G)\cong H_*(X)\ot_\Q H_*(U_j/G)$ to show that
\e
H_n^G(X)\cong \bigop_{m\ge 0}H_{n+m}(X)\ot_\Q H^m_G(*).
\label{co2eq16}
\e
Note in particular that we can have $H_n^G(X)\ne 0$ for all $n\in\Z$.
\item[(d)] Equation \eq{co2eq16} is an example of the K\"unneth isomorphism
\begin{equation*}
\bt:\bigop_{m\in\Z}H_{n-m}^{G_1}(X_1)\ot_\Q H_m^{G_2}(X_2)\,{\buildrel\cong\over\longra}\, H_n^{G_1\t G_2}(X_1\t X_2).
\end{equation*}
This does {\it not\/} work if we replace all of $G_1,G_2,G_1\t G_2$ by $G$.
\item[(e)] For $\Q$-(co)homology of spaces or stacks $X$ we have a bilinear pairing $H_n(X)\t H^n(X)\ra\Q$ inducing an isomorphism $H^n(X)\cong H_n(X)^*$.

For equivariant (co)homology it is generally {\it not\/} true that $H^n_G(X)\cong H_n^G(X)^*$. However, we can define a pairing
\begin{align*}
H_n^G(X)\t H^m_G(X)&\longra H_{n-m}^G(*)\cong H^{m-n}_G(*)\\
\text{by}\qquad (\al,\be)&\longmapsto \pi_*(\al\cap\be),
\end{align*}
where $\pi:X\ra *$ is the projection. This induces a map
\e
H^*_G(X)\longra \Hom_{H^*_G(*)}\bigl(H_*^G(X),H^*_G(*)\bigr),
\label{co2eq17}
\e
where $\Hom_{H^*_G(*)}(-)$ means morphisms of (graded) $H^*_G(*)$-modules. In examples \eq{co2eq17} is often an isomorphism, allowing us to regard $H^*_G(X)$ as a kind of dual of $H_*^G(X)$ over the base ring~$H^*_G(*)$.
\item[(f)] Pullback by $\pi:X\ra *$ gives an algebra morphism $H^*_G(*)\ra H^*_G(X)$. Thus cap product makes $H_*^G(X)$ into a module over $H^*_G(*)=H^*([*/G])$. Pushforwards and pullbacks, cup and cap products are all linear over $H^*_G(*)$, not just over $\Q$. So $H_*^G(-),H^*_G(-)$ behave like (co)homology theories over the ring $H^*_G(*)$, not just over $\Q$.
\item[(g)] Let $G,\{(V_j,U_j):j=0,1,\ldots\},X$ be as above, and $K\subset G$ be a linear algebraic $\C$-subgroup. Observe that $K$ acts on $X$, and $\{(V_j,U_j):j=0,1,\ldots\}$ is also a good system of representations for $K$, so we can define both $H_*^G(X)$ and $H_*^K(X)$ using $\{(V_j,U_j):j=0,1,\ldots\}$. We can then show that there is a natural $\Q$-linear morphism $\La^{G,K}:H_n^G(X)\ra H_n^K(X)$ for all $n\in\Z$, which fits into a commutative diagram for all $j=0,1,\ldots$
\begin{equation*}
\xymatrix@C=140pt@R=15pt{
*+[r]{H_n^G(X)} \ar[d]^{\La^{G,K}} \ar[r]_(0.3){\Pi_j^G} & *+[l]{H^{\pr/(U_j/G)}_{n+2\dim U_j/\dim G}((X\t U_j)/G)} \ar[d]_{\pi^!_{(X\t U_j)/K}} \\
*+[r]{H_n^K(X)}  \ar[r]^(0.3){\Pi_j^K} & *+[l]{H^{\pr/(U_j/K)}_{n+2\dim U_j/\dim K}((X\t U_j)/K),}  }
\end{equation*}
where the rows are as in \eq{co2eq15} and the right hand column is the Gysin morphism \eq{co2eq7} associated to the 2-Cartesian square, as in \eq{co2eq6}
\begin{equation*}
\xymatrix@C=130pt@R=15pt{ *+[r]{(X\t U_j)/K} \ar[r]_{\pi_{(X\t U_j)/K}} \ar[d]^{\pi_{U_j/K}} & *+[l]{(X\t U_j)/G} \ar[d]_{\pi_{U_j/G}} \\ *+[r]{U_j/K} \ar[r]^{\pi_{U_j/G}} & *+[l]{U_j/G,\!}}	
\end{equation*}
with the bottom morphism smooth between smooth separated $\C$-schemes, with fibre $G/K$. The inclusion $K\hookra G$ induces a morphism $\io^{G,K}:X/K\ra X/G$. Pullback along this gives a morphism $\La_{G,K}:H^n_G(X)\ra H^n_K(X)$ for $n=0,1,\ldots.$ These $\La^{G,K},\La_{G,K}$ commute with pushforwards, pullbacks, cup and cap products, and identities, in the obvious way.
\item[(h)] When $K=\{1\}$ we have $H_n^{\{1\}}(X)\cong H_n(X)$ and $H^n_{\{1\}}(X)\cong H^n(X)$. Thus (g) gives morphisms $\La^{G,\{1\}}:H_n^G(X)\ra H_n(X)$, $\La_{G,\{1\}}:H^n_G(X)\ra H^n(X)$, which commute with pushforwards, pullbacks, cup and cap products, and identities.
\end{itemize}
\end{dfn}

\begin{ex}
\label{co2ex1}
Let $G=\bG_m$. Then in Definition \ref{co2def2} we may take $V_j=\C^j$ with $\bG_m$-action $\ep:(z_1,\ldots,z_j)\mapsto(\ep z_1,\ldots,\ep z_j)$, and $U_j=\C^j\sm\{0\}$, so that $U_j/\bG_m=\CP^{j-1}$, and $W_j=\C$ with the obvious identification $V_{j+1}=\C^{j+1}=\C^j\op\C=V_j\op W_j$. Hence~$H_n^{\bG_m}(*)\cong \varprojlim_{j\ra\iy}H^{-n}(\CP^{j-1})=H^{-n}(\CP^\iy)$.
\end{ex}

\begin{rem}
\label{co2rem7}
{\bf(i)} An important difference between $H_*(X/G)$ and $H_*^G(X)$ is that in the situation of Definition \ref{co2def2}(g), we have natural morphisms
\begin{equation*}
\io^{G,K}_*:H_n(X/K)\longra H_n(X/G),\qquad 
\La^{G,K}:H_n^G(X)\longra H_n^K(X),
\end{equation*}
mapping in opposite directions. That is, $H_*(X/G)$ is covariantly functorial, but $H_*^G(X)$ is contravariantly functorial, under inclusions $K\hookra G$. Our construction in \S\ref{co45} of vertex algebras on equivariant homology $H_*^G(\M)$ needs contravariant functoriality in $G$, and does not work for~$H_*(\M/G)$.
\smallskip

\noindent{\bf(ii)} Fulton--MacPherson \cite{FuMa} define {\it bivariant theories\/} generalizing both homology and cohomology. A bivariant theory associates groups $B^k(f:X\ra Y)$ to a morphism $f:X\ra Y$ for $k\in\Z$, with $H^k(X)=B^k(\id_X:X\ra X)$ and $H_k(X)=B^{-k}(\pi:X\ra *)$. If we could extend (co)homology of Artin stacks to a bivariant theory $B^*(-)$, we should identify
\begin{gather*}
 H_k^G(X)=B^{-k}(\pi:X/G\longra */G),\quad H_k(X/G)=B^{-k}(\pi:X/G\longra *), \\
\text{and}\qquad H^k_G(X)=H^k(X/G)=B^k(\id:X/G\longra X/G).
\end{gather*}

\noindent{\bf(iii)} Let $X$ be a proper Deligne--Mumford $\C$-stack acted on by $G$. Edidin and Graham \cite{EdGr1,EdGr2} define the {\it equivariant Chow homology\/} $A_*^G(X,\Q)$ of $X$ over $\Q$. As in \eq{co2eq14}, their definition is equivalent to
\e
A_d^G(X,\Q)=\varprojlim\nolimits_{j\ra\iy}A_{d+\dim U_j/\dim G}((X\t U_j)/G,\Q).
\label{co2eq18}
\e
As $X$ is proper, and Chow homology is a homology theory of Borel--Moore type, generalizing \eq{co2eq4} there are natural morphisms 
\e
A_{d+\dim U_j/\dim G}((X\t U_j)/G,\Q)\longra H^{\pr/(U_j/G)}_{2d+2\dim U_j/\dim G}((X\t U_j)/G)	
\label{co2eq19}
\e
mapping Chow homology to ordinary homology with prescribed supports. Comparing \eq{co2eq14}, \eq{co2eq18} and \eq{co2eq19}, we induce a natural morphism
\e
A_d^G(X,\Q)\longra H_{2d}^G(X),
\label{co2eq20}
\e
which generalizes the morphism $A_d(X,\Q)\ra H_{2d}(X)$ in~\eq{co2eq4}.
\end{rem}

\subsection{Obstruction theories and virtual classes}
\label{co24}

In \S\ref{co24}--\S\ref{co26} we will usually only state results in the minimal generality that we need them later, although stronger results are available. So, for example, Behrend--Fantechi \cite[\S 5]{BeFa1} define virtual classes for Deligne--Mumford stacks with obstruction theories, but we state Theorem \ref{co2thm1} only for proper algebraic $\C$-spaces. We do this to make clear what properties of obstruction theories will be needed in any generalization of our theory, such as those to 3- and 4-Calabi--Yau categories discussed in~\S\ref{co1.5}.

The next definition generalizes Behrend--Fantechi \cite[Def.~4.4]{BeFa1} from Deligne--Mumford $\C$-stacks to Artin $\C$-stacks.

\begin{dfn}
\label{co2def3}
Let $X$ be an Artin $\C$-stack. A ({\it perfect\/}) {\it obstruction theory\/} on $X$ is a morphism $\phi:\cE^\bu\ra\bL_X$ in $D_\qcoh(X)$, such that $\cE^\bu$ is perfect in the interval $[-1,1]$, and $h^k(\phi):h^k(\cE^\bu)\ra h^k(\bL_X)$ is an isomorphism for $k\ge 0$ and surjective for $k=-1$. If $X$ is a scheme, algebraic space, or Deligne--Mumford stack then $h^1(\bL_X)=0$, so $h^1(\cE^\bu)=0$, and $\cE^\bu$ is perfect in the interval~$[-1,0]$.

More generally, if $f:X\ra Y$ is a morphism of Artin $\C$-stacks, a {\it relative perfect obstruction theory\/} on $X$ is a morphism $\phi:\cE^\bu\ra\bL_{X/Y}$ in $D_\qcoh(X)$, such that $\cE^\bu$ is perfect in the interval $[-1,1]$, and $h^k(\phi):h^k(\cE^\bu)\ra h^k(\bL_{X/Y})$ is an isomorphism for $k\ge 0$ and surjective for $k=-1$. This is equivalent to a perfect obstruction theory on $X$ if~$Y=\Spec\C$.
\end{dfn}

The next theorem summarizes Behrend--Fantechi \cite[\S 5]{BeFa1} over the field $\K=\C$, for proper algebraic spaces, except for (iv) which follows from Kiem--Li's theory of {\it cosection localization\/} \cite{KiLi2}, and the extension to the $G$-equivariant case,  which follows Edidin--Graham \cite{EdGr1,EdGr2} and Graber--Pandharipande \cite{GrPa}. Behrend and Fantechi assume \cite[Def.~5.2]{BeFa1} that $\cE^\bu$ has a global resolution by vector bundles, but this is shown to be unnecessary by Kresch~\cite[\S 6.2]{Kres}.

\begin{thm}
\label{co2thm1}
Let\/ $X$ be a proper algebraic\/ $\C$-space and\/ $\phi:\cE^\bu\ra\bL_X$ a perfect obstruction theory on $X$ with\/ $\rank\cE^\bu=d\in\Z$. Then one can construct a \begin{bfseries}virtual class\end{bfseries} $[X]_\virt$ in the rational Chow homology $A_d(X,\Q)$. Applying the projection $A_d(X,\Q)\ra H_{2d}(X,\Q)$ in {\rm\eq{co2eq4},} we may regard\/ $[X]_\virt$ as lying in Betti homology $H_{2d}(X,\Q)$ as in Definition\/ {\rm\ref{co2def1}}. 

These have the following properties:
\begin{itemize}
\setlength{\itemsep}{0pt}
\setlength{\parsep}{0pt}
\item[{\bf(i)}] {\bf(The smooth case.)} If\/ $X$ is smooth and\/ $\cE^\bu=T^*X$ then $[X]_\virt$ is the usual fundamental class $[X]_\fund$ of\/ $X,$ as a compact complex manifold.
\item[{\bf(ii)}] {\bf(Functorial properties.)} Virtual classes have the obvious functorial behaviour under isomorphisms, disjoint unions, and products, of proper algebraic\/ $\C$-spaces with perfect obstruction theories.
\item[{\bf(iii)}] {\bf(Deformation invariance.)} Suppose $f:X\ra Y$ is a proper morphism of algebraic\/ $\C$-spaces with\/ $Y$ connected, and\/ $\phi:\cE^\bu\ra\bL_{X/Y}$ is a relative perfect obstruction theory with\/ $\rank\cE^\bu=d$. Then for each\/ $\C$-point\/ $y\in Y(\C)$ the fibre $X_y=f^{-1}(y)\subset X$ is a proper algebraic\/ $\C$-space with a perfect obstruction theory $\phi\vert_{X_y}:\cE^\bu\vert_{X_y}\ra\bL_{X/Y}\vert_{X_y}\cong\bL_{X_y}$. Also $(i_y)_*([X_y]_\virt)\in H_{2d}(X,\Q)$ is independent of\/ $y\in Y(\C),$ where $i_y:X_y\hookra X$ is the inclusion.
\item[{\bf(iv)}] {\bf(Vanishing condition.)} Suppose there exists a surjective morphism $\si:h^1((\cE^\bu)^\vee)\!\ra\!\O_X,$ where $h^1(\cdots)$ is the cohomology sheaf. Then~$[X]_\virt\!=\!0$.
\end{itemize}

Now suppose $G$ is a linear algebraic $\C$-group acting on $X,$ and\/ $\phi:\cE^\bu\ra\bL_X$ is $G$-equivariant. Then one can construct an \begin{bfseries}equivariant virtual class\end{bfseries} $[X]_\virt$ in the $G$-equivariant Chow homology $A_d^G(X,\Q)$. Applying the projection $A_d^G(X,\Q)\ra H_{2d}^G(X,\Q)$ in {\rm\eq{co2eq20},} we may regard\/ $[X]_\virt$ as lying in equivariant homology $H_{2d}^G(X,\Q)$ as in Definition\/ {\rm\ref{co2def2}}. The analogues of\/ {\bf(i)\rm--\bf(iv)} hold.
\end{thm}

Here Behrend and Fantechi construct the virtual class $[X]_\virt\in A_d(X,\Q)$ without assuming $X$ is proper, which works as Chow homology $A_*(X,\Q)$ is of Borel--Moore type, as in Remark \ref{co2rem2}. But we need $X$ proper so that the morphisms $A_d(X,\Q)\ra H_{2d}(X,\Q)$, $A_d^G(X,\Q)\ra H_{2d}^G(X,\Q)$ are well defined. 

See Li and Tian \cite{LiTi} for an alternative construction of virtual classes, and Battistella--Carocci--Manolache \cite{BCM} and Manolache \cite{Mano1,Mano2} for further properties. Virtual classes are used to define enumerative invariants in algebraic geometry, such as Gromov--Witten \cite{Behr1} or Donaldson--Thomas invariants~\cite{JoSo,Thom1}.

\begin{rem}
\label{co2rem8}
In this book we restrict to schemes and stacks $X$ {\it over the complex numbers\/} $\C$, and define virtual classes $[X]_\virt$ in {\it rational Betti homology\/} $H_*(X,\Q)$. But we could also consider schemes and stacks over other fields $\K$, and other (generalized) (co)homology theories $E_*(-),E^*(-)$ of Artin or higher $\K$-stacks. The essential properties we need of $\K,E_*(-),E^*(-)$ are:
\begin{itemize}
\setlength{\itemsep}{0pt}
\setlength{\parsep}{0pt}
\item[(a)] $E_*(-),E^*(-)$ should have the usual properties of (co)homology theories. In particular, pushforwards $f_*:E_*(X)\ra E_*(Y)$ and pullbacks $f^*:E^*(Y)\ra E^*(X)$ should be defined for {\it arbitrary\/} morphisms $f:X\ra Y$ of Artin or higher $\K$-stacks, not just proper morphisms. 
\item[(b)] Perfect complexes $\cF^\bu\ra X$ should have Chern classes $c_i(\cF^\bu)\in E^*(X)$, with the usual properties.
\item[(c)] We have an algebra isomorphism $E^*([*/\bG_m])\cong R[z]$ or $R[[z]]$, for $R=E_*(*)=E^*(*)$ a $\Q$-algebra, and $z=c_1(L)$ for $L\ra[*/\bG_m]$ the canonical line bundle. Also $E_*([*/\bG_m])\cong R[t]$ with $z^m\cdot t^n=\de_{mn}$ for~$m,n\ge 0$.
\item[(d)] A proper $\K$-scheme $X$ with perfect obstruction theory $\phi:\cE^\bu\ra\bL_X$ should have a virtual class $[X]_\virt\in E_*(X)$, with the usual properties, including Theorems \ref{co2thm1}, \ref{co2thm2} and~\ref{co2thm3}.
\end{itemize}

Parts (a)--(c) are needed for the vertex algebras $\hat H_*(\M)$ and Lie algebras $\check H_*(\M^\pl)$ in \S\ref{co42}--\S\ref{co43} to generalize to $E_*(\M),E_*(\M^\pl)$, and (d) is needed to work with enumerative invariants in $E_*(\M^\pl)$. As in Remark \ref{co2rem2}, Chow homology $A_*(X)$ is of Borel--Moore type, so it does not satisfy (a), as $f_*$ is defined only for proper morphisms $f$, and we cannot use it.
\end{rem}

\begin{rem}
\label{co2rem9}
{\bf(Interpretation of obstruction theories using Derived Algebraic Geometry.)} Most obstruction theories studied in examples can be understood as coming from Derived Algebraic Geometry (DAG) as in \S\ref{co22}, in the following way, as in Sch\"urg--To\"en--Vezzosi \cite[\S 1]{STV}. 

Let $\bX$ be a derived Artin $\C$-stack, and $X=t_0(\bX)$ be its classical truncation, an ordinary Artin $\C$-stack. As in Remark \ref{co2rem5}(h) there is an inclusion morphism $i:X\hookra\bX$. As in Remark \ref{co2rem6} the cotangent complexes of $X,\bX$ form a distinguished triangle \eq{co2eq11} in $D_\qcoh(X)$:
\e
\xymatrix@C=40pt{ i^*(\bL_\bX) \ar[r]^{\bL_i} & \bL_X \ar[r] & \bL_{X/\bX} \ar[r] & i^*(\bL_\bX)[1]. }
\label{co2eq21}
\e

As in \cite[Prop.~1.2]{STV} we have $h^k(\bL_{X/\bX})=0$ for $k\ge -1$, by properties of classical truncations. Thus the long exact sequence of cohomology sheaves of \eq{co2eq21} implies that $h^k(\phi):h^k(i^*(\bL_\bX))\ra h^k(\bL_X)$ is an isomorphism for $k\ge 0$ and surjective for $k=-1$. That is, $\bL_i:i^*(\bL_\bX)\ra\bL_X$ automatically satisfies one of the conditions for an obstruction theory in Definition~\ref{co2def3}. 

If $\bX$ is {\it quasi-smooth}, as in Remarks \ref{co2rem5}(k) and \ref{co2rem6}(e), then $\bL_\bX$ (and hence $i^*(\bL_\bX)$) is perfect in the interval $[-1,1]$, so $\bL_i:i^*(\bL_\bX)\ra\bL_X$ satisfies the remaining condition. Thus, classical truncations $X=t_0(\bX)$ of quasi-smooth derived Artin $\C$-stacks $\bX$ have natural Behrend--Fantechi obstruction theories. So if $\bX$ is a proper quasi-smooth derived algebraic $\C$-space then Theorem \ref{co2thm1} gives a virtual class $[X]_\virt$ in $A_*(X,\Q)$ or~$H_*(X,\Q)$.

Similarly, if $\bs f:\bX\ra\bY$ is a quasi-smooth morphism of derived Artin stacks with classical truncation $f:X\ra Y$, as in Remark \ref{co2rem6}(f), then $\bL_{f/\bs f}:i^*(\bL_{\bX/\bY})\ra \bL_{X/Y}$ is a relative obstruction theory for~$f:X\ra Y$.

When one can construct a classical moduli stack $\M$ in some problem (e.g.\ the moduli of objects in an abelian category $\A$ such as $\coh(X)$), very often one can also use DAG to construct a derived moduli stack $\bcM$ with $\M=t_0(\bcM)$. There may be geometrical reasons for $\bcM$ to be quasi-smooth (e.g.\ if $\Ext^{>2}(E,E)=0$ for all $E\in\A$), and then $\bL_i:i^*(\bL_\bcM)\ra\bL_\M$ is an obstruction theory on $\M$.

In fact Behrend--Fantechi obstruction theories predate DAG, and are usually constructed by hand using derived category techniques, rather than via DAG. But secretly they are classical shadows of quasi-smooth derived moduli stacks.
\end{rem}

\subsection{Pushforward along projective space bundles}
\label{co25}

We generalize Manolache \cite[Def.~4.5]{Mano1}, \cite[Def.~2.17]{Mano2} to Artin stacks. 

\begin{dfn}
\label{co2def4}
Suppose $f:X\ra Y$ is a morphism of Artin $\C$-stacks, $\phi:\cF^\bu\ra\bL_X$, $\psi:\cG^\bu\ra\bL_Y$ are obstruction theories on $X,Y$, and $\chi:\cH^\bu\ra\bL_{X/Y}$ is a relative obstruction theory for $f$. We say that $\phi,\psi,\chi$ are a {\it compatible triple\/} if there exists a commutative diagram in $D_\qcoh(X)$, with distinguished rows:
\e
\begin{gathered}
\xymatrix@C=40pt@R=15pt{
f^*(\cG^\bu) \ar[d]^{f^*(\psi)} \ar[r] & \cF^\bu \ar[d]^\phi \ar[r] & \cH^\bu \ar[d]_\chi \ar[r] & f^*(\cG^\bu)[1] \ar[d]_{f^*(\psi)[1]} \\
f^*(\bL_Y)  \ar[r]^{\bL_f} & \bL_X  \ar[r] & \bL_{X/Y}  \ar[r] & f^*(\bL_Y)[1].  }
\end{gathered}
\label{co2eq22}
\e
Here the second row is \eq{co2eq11} for~$f:X\ra Y$.

If an algebraic $\C$-group $G$ acts on $X,Y$, and $f:X\ra Y$ is $G$-equivariant, and $\phi,\psi,\chi$ are $G$-equivariant (relative) obstruction theories, we call $\phi,\psi,\chi$ a {\it $G$-equivariant compatible triple\/} if \eq{co2eq22} lifts to $G$-equivariant complexes~$D_\qcoh^G(X)$.
\end{dfn}

\begin{ex}
\label{co2ex2}
As in Remark \ref{co2rem9}, obstruction theories often arise from quasi-smooth derived Artin stacks $\bX$, taking $X=t_0(\bX)$ with inclusion $i_X:X\hookra\bX$ and $\phi:\cF^\bu\ra\bL_X$ to be~$\bL_{i_X}:i_X^*(\bL_\bX)\ra \bL_X$. 

This also gives a natural source of examples satisfying Definition \ref{co2def4}. Suppose that $\bs f:\bX\ra\bY$ is a quasi-smooth morphism of quasi-smooth derived Artin stacks, with classical truncation $f:X\ra Y$. Then we have a diagram~\eq{co2eq22}:
\begin{equation*}
\xymatrix@C=33pt@R=15pt{
{\begin{subarray}{l}\ts f^*\!\ci\! i_Y^*(\bL_{\bY})\!=\\ \ts i_X^*\!\ci\!\bs f^*(\bL_{\bY})\end{subarray}} \ar[d]^{f^*(\bL_{i_Y})} \ar[r]_(0.55){i^*(\bL_{\bs f})} & i_X^*(\bL_\bX) \ar[d]^{\bL_{i_X}} \ar[r] & i_X^*(\bL_{\bX/\bY}) \ar[d]_{\bL_{f/\bs f}} \ar[r] & {\begin{subarray}{l}\ts f^*\!\ci\! i_Y^*(\bL_{\bY})[1]\!=\\ \ts i_X^*\!\ci\!\bs f^*(\bL_{\bY})[1]\end{subarray}} \ar[d]_{f^*(\bL_{i_Y})[1]} \\
f^*(\bL_Y) \ar[r]^{\bL_f} & \bL_X \ar[r] & \bL_{X/Y} \ar[r] & f^*(\bL_Y)[1].  }
\end{equation*}
Thus, quasi-smooth morphisms $\bs f:\bX\ra\bY$ of quasi-smooth derived Artin stacks automatically yield data $(\bL_{i_X},\bL_{i_Y},\bL_{f/\bs f})$ satisfying Definition~\ref{co2def4}.
\end{ex}

Suppose $f:X\ra Y$ is a morphism of proper Deligne--Mumford $\C$-stacks, and $\phi:\cF^\bu\ra\bL_X$, $\psi:\cG^\bu\ra\bL_Y$, $\chi:\cH^\bu\ra\bL_{X/Y}$ are as in Definition \ref{co2def4}, so we have virtual classes $[X]_\virt,[Y]_\virt$ in rational Chow homology $A_*(-,\Q)$, which may be projected to Betti homology $H_*(-,\Q)$. Manolache \cite[Th.~4.8]{Mano1} shows that $[X]_\virt=f^!([Y]_\virt)$, for $f^!:A_*(Y,\Q)\ra A_*(X,\Q)$ a {\it virtual pullback map\/} along $f:X\ra Y$, defined using $\chi:\cH^\bu\ra\bL_{X/Y}$. In \cite{Mano2} she uses this to give a {\it virtual pushforward formula\/} relating $f_*([X]_\virt\cap \eta)$ to $[Y]_\virt$ for~$\eta\in H^*(X,\Q)$.

For our theory, we will only need the following special case of the virtual pushforward formula, in Betti homology, when $f:X\ra Y$ is a projective space bundle. It follows from Manolache \cite[Th.~4.7]{Mano2}, \cite[Th.~7.7]{BCM}, noting that the constant $m+1$ in \eq{co2eq23} may be identified as
\begin{equation*}
\int_{f^{-1}(y)}c_\top(\bT_{X/Y})\vert_{f^{-1}(y)}=\int_{f^{-1}(y)}c_\top(T(f^{-1}(y)))=\chi(f^{-1}(y))=\chi(\CP^m).
\end{equation*}

\begin{thm}
\label{co2thm2}
Let\/ $X,Y$ be proper algebraic $\C$-spaces and\/ $f:X\ra Y$ be a smooth morphism which is a fibre bundle with fibre $\CP^m$. Note that\/ $\id_{\bL_{X/Y}}:\bL_{X/Y}\ra \bL_{X/Y}$ is a relative obstruction theory for $f,$ as $f$ is smooth.

Suppose $\phi:\cF^\bu\ra\bL_X,$ $\psi:\cG^\bu\ra\bL_Y$ are obstruction theories on $X,Y$ with\/ $\rank\cF^\bu=m+n$ and\/ $\rank\cG^\bu=n,$ such that\/ $(\phi,\psi,\id_{\bL_{X/Y}})$ are a compatible triple as in Definition\/ {\rm\ref{co2def4}}. Then the virtual classes\/ $[X]_\virt\in H_{2m+2n}(X,\Q)$ and\/~$[Y]_\virt\in H_{2n}(Y,\Q)$ satisfy
\e
f_*\bigl([X]_\virt\cap c_\top(\bT_{X/Y})\bigr)=(m+1)\cdot[Y]_\virt,
\label{co2eq23}
\e
where $\bT_{X/Y}\ra X$ is the relative tangent complex, a vector bundle of rank\/ $m,$ and\/ $c_\top(-)=c_m(-)$ is the top Chern class. 

If the whole situation is equivariant under a linear algebraic $\C$-group $G$ acting on $X,Y,$ then the above holds in equivariant homology\/~$H_*^G(-,\Q)$.
\end{thm}

Combining Example \ref{co2ex2} and Theorem \ref{co2thm2} and noting that $\bL_{f/\bs f}:i^*(\bL_{\bX/\bY})\ab\ra \bL_{X/Y}$ is isomorphic to $\id_{\bL_{X/Y}}:\bL_{X/Y}\ra \bL_{X/Y}$ as $\bs f$ is smooth yields:

\begin{cor}
\label{co2cor1}
Let\/ $\bX,\bY$ be proper quasi-smooth derived algebraic $\C$-spaces with\/ $\vdim\bX=m+n,$ $\vdim\bY=n$ and\/ $\bs f:\bX\ra\bY$ be a smooth morphism which is a fibre bundle with fibre $\CP^m$. Write $f:X\ra Y$ for the classical truncation of $\bs f:\bX\ra\bY$. As in\/ {\rm\S\ref{co24}} we have obstruction theories $\bL_{i_X}:i_X^*(\bL_\bX)\ra \bL_X$ and\/ $\bL_{i_Y}:i_Y^*(\bL_\bY)\ra \bL_Y$ on $X,Y,$ and virtual classes\/ $[X]_\virt\in H_{2m+2n}(X,\Q)$ and\/ $[Y]_\virt\in H_{2n}(Y,\Q)$. Then in $H_{2n}(Y,\Q)$ we have
\begin{equation*}
f_*\bigl([X]_\virt\cap c_\top(\bT_{X/Y})\bigr)=(m+1)\cdot[Y]_\virt,
\end{equation*}
where $\bT_{X/Y}\ra X$ is the relative tangent complex, a vector bundle of rank\/ $m,$ and\/ $c_\top(-)=c_m(-)$ is the top Chern class. 

If the whole situation is equivariant under a linear algebraic $\C$-group $G$ acting on $\bX,\bY,$ then the above holds in equivariant homology\/~$H_*^G(-,\Q)$.	
\end{cor}

\subsection{\texorpdfstring{$\bG_m$-localization of virtual classes}{𝔾ᵐ-localization of virtual classes}}
\label{co26}

The next theorem summarizes Graber and Pandharipande's $\bG_m$-localization formula for virtual classes \cite{GrPa}, projected from Chow homology to ordinary homology. Graber and Pandharipande make an extra technical assumption, that $\cF^\bu$ may be written as a two-term complex of equivariant vector bundles $[E^{-1}\ra E^0]$. This assumption is weakened by Chang--Kiem--Li \cite{CKL}. The final part may be deduced from the non-equivariant case using~\eq{co2eq14}.

\begin{thm}
\label{co2thm3}
Let\/ $X$ be a proper algebraic\/ $\C$-space with an action of\/ $\bG_m=\C^*,$ and a $\bG_m$-equivariant perfect obstruction theory $\phi:\cF^\bu\ra\bL_X$ with\/ $\rank \cF^\bu\ab=d\in\Z,$ so that as in {\rm\S\ref{co24}} we have a Behrend--Fantechi virtual class $[X]_\virt$ in the equivariant homology group $H_{2d}^{\bG_m}(X,\Q)$. 

Note that\/ $H_*^{\bG_m}(X,\Q)$ is a module over the $\Q$-algebra $H_*^{\bG_m}(*,\Q)\cong\Q[z]$. Thus we may form the \begin{bfseries}localized equivariant homology\end{bfseries}
\begin{equation*}
H_*^{\bG_m}(X,\Q)[z^{-1}]=H_*^{\bG_m}(X,\Q)\ot_{\Q[z]}\Q[z,z^{-1}].
\end{equation*}

Write $X^{\bG_m}$ for the $\bG_m$-fixed substack of\/ $X,$ and\/ $\{X_a:a\in A\}$ for the connected components of\/ $X^{\bG_m},$ so that\/ $A$ is finite, and\/ $X_a$ is proper. Then $\phi\vert_{X^{\bG_m}}:\cF^\bu\vert_{X^{\bG_m}}\ra\bL_X\vert_{X^{\bG_m}}$ is a morphism of complexes with\/ $\bG_m$-actions on $X^{\bG_m}$. The $\bG_m$-invariant part is a perfect obstruction theory on $X^{\bG_m}$
\e
\phi^{\bG_m}=\phi\vert_{(\cF^\bu\vert_{X^{\bG_m}})^{\bG_m}}:(\cF^\bu\vert_{X^{\bG_m}})^{\bG_m}\longra(\bL_X\vert_{X^{\bG_m}})^{\bG_m}\cong\bL_{X^{\bG_m}}.
\label{co2eq24}
\e

Write $\phi_a:\cF^\bu_a\ra \bL_{X_a}$ for the restriction of\/ \eq{co2eq24} to $X_a,$ a perfect obstruction theory on $X_a,$ and set\/ $d_a=\rank\cF^\bu_a$. Then we have virtual classes $[X_a]_\virt$ in $H_{2d_a}(X_a,\Q)$ for $a\in A$. We regard these as living in the localized equivariant homology groups for the trivial\/ $\bG_m$-action on {\rm$X_a^{\bG_m}$:} 
\begin{equation*}
H_*^{\bG_m}(X_a,\Q)[z^{-1}]\cong H_*(X_a,\Q)\ot_\Q \Q[z^{\pm 1}].
\end{equation*}

Write $\cF^\bu\vert_{X_a}=\cF^\bu_a\op \cN^\bu_a,$ where\/ $\cN^\bu_a\ra X_a$ is a perfect complex with nontrivial\/ $\bG_m$-action called the \begin{bfseries}virtual conormal bundle\end{bfseries} of\/ $X_a$ in $X,$ the $\bG_m$-moving part of\/ $\cF^\bu\vert_{X_a},$ with\/ $\rank\cN_a^\bu=d-d_a$. Write\/ $i_a:X_a\hookra X$ for the inclusion. As $\cN_a^\bu$ has only nontrivial\/ $\bG_m$-representations, it has an invertible \begin{bfseries}Euler class\end{bfseries}
\begin{equation*}
e(\cN_a^\bu)\in H^*_{\bG_m}(X_a,\Q)[z^{-1}]\cong H^*(X_a,\Q)\ot_\Q \Q[z^{\pm 1}].
\end{equation*}
Then in $H_{2d}^{\bG_m}(X,\Q)[z^{-1}]$ we have
\e
[X]_\virt=\sum_{a\in A}(-1)^{\rank\cN_a^\bu}(i_a)_*\bigl([X_a]_\virt\cap e(\cN_a^\bu)^{-1}\bigr).
\label{co2eq25}
\e

More generally, suppose\/ $G$ is a linear algebraic $\C$-group, and\/ $\bG_m\t G$ acts on $X$ preserving $\phi:\cF^\bu\ra\bL_X$. Then the above holds for $[X]_\virt\in H_{2d}^{\bG_m\t G}(X,\Q),$ $[X_a]_\virt\!\in\! H_{2d_a}^G(X_a,\Q)$ and\/ $e(\cN_a^\bu)\!\in\! H^*_{\bG_m\t G}(X_a,\Q)[z^{-1}]\!\cong\! H^*_G(X_a,\Q)\ot\Q[z^{\pm 1}]$.
\end{thm}


\begin{rem}
\label{co2rem10}
Actually, we have {\it cheated\/} in Theorem \ref{co2thm3}: Chang--Kiem--Li \cite[Th.~3.5]{CKL} still have a technical assumption that we did not include in Theorem \ref{co2thm3}, namely that $\cN^\bu_a$ has a 2 term global resolution by vector bundles on $X_a$. This exists automatically if $X_a$ is a projective $\C$-scheme. It is generally hoped that this assumption can be eliminated (see \cite[discussion of Ass.~5.3]{KiSa}, for instance), but this does not seem to be available in the literature at present.

We offer three possible solutions to this problem. Firstly, we only apply Theorem \ref{co2thm3} in Chapters \ref{co9} and \ref{co10}, when $X=\baM_{(\al,\bs d)}^\rst(\bar\tau^\la_{\bs\mu})=\baM_{(\al,\bs d)}^\ss(\bar\tau^\la_{\bs\mu})$ is a moduli space of objects in an auxiliary abelian category $\baA$ defined in \S\ref{co52}. In most of our examples one can prove that $\baM_{(\al,\bs d)}^\ss(\bar\tau^\la_{\bs\mu})$ is a projective $\C$-scheme, and Theorem \ref{co2thm3} is known in this case.

Secondly, Khan \cite[\S 2]{Khan} defines Borel--Moore homology $H_*^{\text{Kh-BM}}(X,\Q)$ of Artin and higher $\C$-stacks $X$ using sheaf cohomology methods. (Actually, this is just one of a large class of generalized homology theories Khan defines, which also have virtual classes as follows.) If $\bX$ is a quasi-smooth derived Artin $\C$-stack and $X=t_0(\bX)$ then Khan \cite[\S 3]{Khan} defines a virtual class $[X]_\virt$ in $H_*^{\text{Kh-BM}}(X,\Q)$ using sheafy methods. He relates this to the Behrend--Fantechi virtual class when $\bX$ is a quasi-smooth derived Deligne--Mumford $\C$-stack, and proves analogues of parts of Theorems \ref{co2thm1} and~\ref{co2thm2}.

Latyntsev \cite[\S A]{Laty} explains how to use Khan's methods to define an ordinary homology theory $H_*^{\text{Kh}}(X,\Q)$ of Artin and higher $\C$-stacks $X$, such that $H_*^{\text{Kh-BM}}(X,\Q)=H_*^{\text{Kh}}(X,\Q)$ if $X$ is a proper algebraic space. Thus, if $\bX$ is a proper quasi-smooth algebraic $\C$-space then Khan's virtual class $[X]_\virt$ is defined in $H_*^{\text{Kh}}(X,\Q)$. He also proves \cite[\S 8]{Laty} an analogue of Theorem \ref{co2thm3} for Khan's virtual classes. So, one option would be to use Khan's virtual classes $[X]_\virt$ in $H_*^{\text{Kh}}(X,\Q)$ throughout.

Thirdly, I am confident I can prove Theorem \ref{co2thm3} in ordinary homology (not Chow homology) without this technical assumption, using Derived Differential Geometry \cite{Joyc8,Joyc9,Joyc10,Joyc11}, and this would be sufficient for our purposes. 

For this preliminary version of the book, I am just going to work with Theorem \ref{co2thm3} as above, {\it even though it is missing a technical assumption}. Please let me know if you know how to remove the assumption.  If it seems unavoidable, I will write up my own proof of Theorem \ref{co2thm3} using DDG.
\end{rem}

Here we include the factor $(-1)^{\rank\cN_a^\bu}$ in \eq{co2eq25}, which is not in \cite{GrPa}, because we take $\cN_a^\bu$ to be the virtual {\it conormal\/} bundle, not the virtual normal bundle.

\begin{cor}
\label{co2cor2}
Suppose\/ $X$ is a proper algebraic\/ $\C$-space with an action of\/ $\bG_m=\C^*,$ and a $\bG_m$-equivariant perfect obstruction theory $\phi:\cF^\bu\ra\bL_X$ with\/ $\rank \cF^\bu=d\in\Z$. Define $\{X_a:a\in A\},$ $\phi_a:\cF^\bu_a\ra \bL_{X_a},$ $d_a,$ $\cN^\bu_a$ be as in Theorem\/ {\rm\ref{co2thm3},} so we have virtual classes $[X_a]_\virt\in H_{2d_a}(X_a)$ for $a\in A,$ and Euler classes\/~$e(\cN_a^\bu)\in H^{2d-2d_a}(X_a)\ot\Q[z^{\pm 1}]$.

Let\/ $f:X\ra Y$ be a $\bG_m$-equivariant morphism of Artin $\C$-stacks for the trivial\/ $\bG_m$-action on $Y,$ and\/ $\eta\in H^n_{\bG_m}(X)$. Then in $H_{2d-2-n}(Y)$ we have
\e
\sum_{a\in A}(-1)^{\rank\cN_a^\bu}\Res_z\bigl[(f\ci i_a)_*\bigl([X_a]_\virt\cap (e(\cN_a^\bu)^{-1}\cup i_a^*(\eta))\bigr)\bigr]=0.
\label{co2eq26}
\e

More generally, let\/ $G$ be a linear algebraic $\C$-group, and\/ $\phi:\cF^\bu\ra\bL_X,$ $f:X\ra Y$ be equivariant under actions of\/ $\bG_m\t G$ on $X,Y,$ where $\bG_m\t\{1\}$ acts trivially on $Y$. Taking\/ $[X_a]_\virt\in H_{2d_a}^G(X_a),$ $e(\cN_a^\bu)\in H^{2d-2d_a}_G(X_a)\ot\Q[z^{\pm 1}],$ $\eta\in H^n_{\bG_m\t G}(X),$ equation \eq{co2eq26} holds in\/~$H_{2d-2-n}^G(Y)$.
\end{cor}

\begin{proof} Apply $f_*(-\cap\eta)$ to \eq{co2eq25} to get an equation in $H_{2d-2-n}^{\bG_m}(Y,\Q)[z^{-1}]\cong H_{2d-2-n}(Y,\Q)\ot_\Q\Q[z^{\pm 1}]$, as $\bG_m$ acts trivially on $Y$. Since $[X]_\virt\in H_{2d}^{\bG_m}(X)$ and $\eta\in H^n_{\bG_m}(X)$ we have
\begin{align*}
f_*([X]_\virt\cap \eta)\in H_{2d-2-n}^{\bG_m}(Y,\Q)&\cong H_{2d-2-n}(Y,\Q)\ot_\Q\Q[z]\\
&\subset H_{2d-2-n}(Y,\Q)\ot_\Q\Q[z^{\pm 1}].
\end{align*}
Hence $\Res_z[f_*([X]_\virt\cap \eta)]=0$, and equation \eq{co2eq26} follows. The $G$-equivariant version follows in the same way, by the corresponding part of Theorem~\ref{co2thm3}.
\end{proof}

\section{Abelian categories and stability}
\label{co3}

Next we discuss (weak) stability conditions $(\tau,T,\le)$ on abelian categories $\A$. Sections \ref{co31}--\ref{co32} mostly follow the author \cite{Joyc6,Joyc7}, with some new results. Section \ref{co33} explains results of Alper--Halpern-Leistner--Heinloth \cite[\S 7]{AHLH} giving criteria for when a moduli stack $\M_\al^\ss(\tau)$ of $\tau$-semistable objects in $\A$ admits a proper `good moduli space', as this is important for defining virtual classes~$[\M_\al^\ss(\tau)]_\virt$.

\subsection{(Weak) stability conditions on abelian categories}
\label{co31}

The next definition comes from the author \cite{Joyc6}. See also Rudakov~\cite{Ruda}.

\begin{dfn}
\label{co3def1}
Let $\A$ be an abelian category. The {\it Grothendieck group\/} $K_0(\A)$ is the abelian group generated by isomorphism classes $[E]$ of objects $E\in\A$, with a relation $[E_2]=[E_1]+[E_3]$ for each exact sequence $0\ra E_1\ra E_2\ra E_3\ra 0$ in $\A$. Often $K_0(\A)$ is inconveniently large (e.g.\ it may be uncountably generated). So we suppose we are given a surjective quotient $K_0(\A)\twoheadrightarrow K(\A)$, and work with $K(\A)$ instead of $K_0(\A)$. We write $\lb E\rb\in K(\A)$ for the class of~$E\in\A$.

Suppose that $0\in\A$ is the only object in class $0\in K(\A)$. Define the {\it positive cone\/} $C(\A)\subset K(\A)\sm\{0\}$ by~$C(\A)=\bigl\{\lb E\rb:0\ne E\in\A\bigr\}$.

Here $K(\A)$ is an abelian group, usually free of finite rank, which we can think of as parametrizing suitable topological invariants of objects in $\A$. If $\A=\coh(X)$ for a smooth projective $\C$-scheme $X$, we can take $K(\A)\subset H^{\rm even}(X,\Q)$ to be the lattice of Chern characters $\ch E=\lb E\rb$ of $E\in\coh(X)$.

Let $(T,\leq)$ be a totally ordered set and $\tau:C(\A)\ra T$ be a map. We call $(\tau,T,\leq)$ a {\it weak stability condition\/} on $\A$ if for all $\al,\be,\ga \in C(\A)$ with $\be=\al+\ga$, either 
$\tau(\al) \leq \tau(\be) \leq \tau(\ga)$, or~$\tau(\al) \geq \tau(\be) \geq \tau(\ga)$.

We call $(\tau,T,\leq)$ a {\it stability condition\/} if for all such $\al,\be,\ga$, either 
$\tau(\al)<\tau(\be)<\tau(\ga)$, or $\tau(\al)>\tau(\be)>\tau(\ga)$, or~$\tau(\al)=\tau(\be)=\tau(\ga)$.

Let $(\tau,T,\leq)$ be a weak stability condition. An object $E$ of $\A$ is called:
\begin{itemize}
\setlength{\itemsep}{0pt}
\setlength{\parsep}{0pt}
\item[(i)] {\it $\tau$-stable\/} if $\tau(\lb E'\rb)<\tau(\lb E/E'\rb)$ for all subobjects $E'\subset E$ with $E'\ne 0,E$.
\item[(ii)] {\it $\tau$-semistable\/} if $\tau(\lb E'\rb)\!\leq\!\tau(\lb E/E'\rb)$ for all $E'\subset E$ with $E'\ne 0,E$.
\item[(iii)] {\it $\tau$-unstable} if it is not $\tau$-semistable.
\item[(iv)] {\it strictly $\tau$-semistable\/} if it is $\tau$-semistable but not $\tau$-stable.
\end{itemize}

We call an abelian category $\A$ {\it noetherian} if there exist no infinite strictly ascending chains of subobjects $E_1\subsetneq E_2\subsetneq E_3\subsetneq\cdots\subseteq E$ in~$\A$.

Let $(\tau,T,\le)$ be a weak stability condition on $\A$. We call $\A$ $\tau$-{\it artinian} if there exist no infinite strictly descending chains of subobjects $\cdots\subsetneq E_3\subsetneq E_2\subsetneq
E_1=E$ in $\A$ with $\tau(\lb E_{n+1}\rb)\ge\tau(\lb E_n/E_{n+1}\rb)$ for all~$n\ge 1$.
\end{dfn}

The next proposition is proved as in Rudakov~\cite[Th.~1]{Ruda}.

\begin{prop}
\label{co3prop1}
In Definition\/ {\rm\ref{co3def1},} suppose $0\ne E,F\in\A$ are $\tau$-semistable with $\tau(\lb E\rb)>\tau(\lb F\rb)$. Then $\Hom_\A(E,F)=0$.	
\end{prop}

Here is \cite[Th.~4.4]{Joyc6}, based on Rudakov \cite[Th.~2]{Ruda}.

\begin{thm}
\label{co3thm1}
Let\/ $(\tau,T,\le)$ be a weak stability condition on an abelian category $\A$. Suppose $\A$ is noetherian and\/ $\tau$-artinian. Then each\/ $E\in\A$ has a unique filtration by subobjects\/ $0=E_0\subsetneq E_1\subsetneq\cdots\subsetneq E_n=E$ for $n\ge 0$, such that\/ $F_i=E_i/E_{i-1}$ is $\tau$-semistable for $i=1,\ldots,n,$ and\/ $\tau([F_1])>\tau([F_2])>\cdots>\tau([F_n])$ in\/ $T$. This is called the\/ {\rm($\tau$-)}\begin{bfseries}Harder--Narasimhan filtration\end{bfseries} of\/~$E$.
\end{thm}

We define some notation for moduli stacks of $\tau$-(semi)stable objects in~$\A$:

\begin{dfn}
\label{co3def2}
Continue in the situation of Definition \ref{co3def1}. Suppose also that $\A$ is $\C$-linear, and that we can form {\it moduli stacks\/} $\M,\M^\pl$ of objects in $\A$, as Artin $\C$-stacks. As explained in \S\ref{co42}--\S\ref{co43}, in this book we consider two different moduli stacks: the (usual) moduli stack $\M$, parametrizing objects $E,F$ of $\A$ up to isomorphisms $\phi:E\ra F$, and the `projective linear' moduli stack $\M^\pl$, parametrizing objects $E$ up to `projective linear' isomorphisms, that is, we identify isomorphisms $\phi,\psi:E\ra F$ if $\psi=\phi\ci\ga\,\id_E$ for $\ga\in\bG_m$. There is a morphism $\Pi^\pl:\M\ra\M^\pl$ with fibre $[*/\bG_m]$ over $\M^\pl\sm\{0\}$. We suppose there are decompositions $\M=\coprod_{\al\in K(\A)}\M_\al$, $\M^\pl=\coprod_{\al\in K(\A)}\M^\pl_\al$, where $\M_\al,\M^\pl_\al$ parametrize objects $E\in\A$ with $\lb E\rb=\al$ in $K(\A)$.

For the weak stability conditions considered in this book, for $E$ to be $\tau$-semistable or $\tau$-stable is always an open condition on $[E]$ in $\M$ and $\M^\pl$. For each $\al\in C(\A)$ we will write $\M_\al^\rst(\tau)\subseteq\M_\al^\ss(\tau)\subseteq\M_\al^\pl$ for the moduli stacks of $\tau$-(semi)stable objects in $\A$ in class $\al\in C(\A)$, as open $\C$-substacks.
\end{dfn}

Here is a simple example, which will be discussed in more detail in~\S\ref{co63}.

\begin{ex}
\label{co3ex1}
Let $Q=(Q_0,Q_1,h,t)$ be a quiver, so that $Q_0,Q_1$ are finite sets of vertices and edges, and $h,t:Q_1\ra Q_0$ are the head and tail maps. Let $\A=\modCQ$ be the abelian category of $\C$-representations of $Q$, as in \S\ref{co61} below. Write objects of $\A$ as $(V,\rho)=((V_v)_{v\in Q_0},(\rho_e)_{e\in Q_1})$, where $V_v$ is a finite-dimensional $\C$-vector space and $\rho_e:V_{t(e)}\ra V_{h(e)}$ a linear map. Define the {\it dimension vector\/} of $(V,\rho)$ to be $\bs d=\bdim(V,\rho)\in\N^{Q_0}\subset\Z^{Q_0}$, where $\bs d(v)=\dim_\C V_v$ for $v\in Q_0$. Define $K(\A)=\Z^{Q_0}$ to be the lattice of dimension vectors, with $\lb V,\rho\rb=\bdim(\bs V,\bs\rho)$. Then the positive cone is~$C(\A)=\N^{Q_0}\sm\{0\}$.

Choose real numbers $\mu_v\in\R$ for all $v\in Q_0$. Define a map $\mu:C(\A)\ra\R$ by
\begin{equation*}
\mu(\bs d)=\frac{\sum_{v\in Q_0}\bs d(v)\mu_v}{\sum_{v\in Q_0}\bs d(v)}.
\end{equation*}
Then $(\mu,\R,\le)$ is a stability condition on $\modCQ$, called {\it slope stability}.
\end{ex}

Here are two more complicated examples, which will be studied in~\S\ref{co71}.

\begin{ex}
\label{co3ex2}
Let $X$ be a smooth projective $\C$-scheme of dimension $m$. Write $\coh(X)$ for the abelian category of coherent sheaves on $X$. Some good references on coherent sheaves are Hartshorne \cite[\S II.5]{Hart}, Huybrechts and Lehn \cite{HuLe1}, and Huybrechts \cite{Huyb}. As in \cite[\S III.6]{Hart}, for $E,F$ in $\coh(X)$ we have {\it Ext groups\/} $\Ext^i(E,F)$ for $i=0,\ldots,m$, which are finite-dimensional $\C$-vector spaces, with $\Ext^0(E,F)=\Hom(E,F)$. The {\it Euler form\/} is the biadditive map $\chi:K_0(\coh(X))\t K_0(\coh(X))\ra\Z$ with
\begin{equation*}
\chi([E],[F])=\sum_{i=0}^m(-1)^i\dim_\C\Ext^i(E,F)
\end{equation*}
for all $E,F\in\coh(X)$. The {\it numerical Grothendieck group\/} is 
\begin{align*}
&K^\num(\coh(X))=K_0(\coh(X))/\Ker\chi, \quad\text{where}\\
&\Ker\chi=\bigl\{\al\in K_0(\coh(X)):\text{$\chi(\al,\be)=0$ for all $\be\in K_0(\coh(X))$}\bigr\}.
\end{align*}
Then $\chi$ descends to $\chi:K^\num(\coh(X))\t K^\num(\coh(X))\ra\Z$. Equivalently, $K^\num(\coh(X))\cong\Im\bigl(\ch:K_0(\coh(X))\ra H^{\rm even}(X,\Q)\bigr)$, where $\ch$ is the {\it Chern character}, as in \cite[App.~A]{Hart}. We will take $K(\coh(X))=K^\num(\coh(X))$ in Definition \ref{co3def1}. We write $\lb E\rb\in K(\coh(X))$ for the class of~$E\in\coh(X)$.

Define $G$ to be the set of monic rational polynomials in $t$ of degree $\le m$:
\begin{equation*}
G=\bigl\{p(t)=t^d+a_{d-1}t^{d-1}+\cdots+a_0:d=0,1,\ldots,m,\;\> a_0,\ldots,a_{d-1}\in\Q\bigr\}.
\end{equation*}
Define a total order `$\le$' on $G$ by $p\le p'$ for $p,p'\in G$ if either
\begin{itemize}
\setlength{\itemsep}{0pt}
\setlength{\parsep}{0pt}
\item[(a)] $\deg p>\deg p'$, or
\item[(b)] $\deg p=\deg p'$ and $p(t)\le p'(t)$ for all $t\gg 0$.
\end{itemize}

Fix a very ample line bundle $\O_X(1)$ on $X$. For $E\in\coh(X)$, the {\it Hilbert polynomial\/} $P_E$ is the unique polynomial in $\Q[t]$ such that $P_E(n)=\dim H^0(E(n))$ for all $n\gg 0$. Equivalently, $P_E(n)=\chi\bigl(\lb\O_X(-n)\rb,\lb E\rb\bigr)$ for all $n\in\Z$. Thus, $P_E$ depends only on the class $\al\in K(\coh(X))$ of $E$, and we may write $P_\al$ instead of $P_E$. Define $\tau:C(\coh(X))\ra G$ by $\tau(\al)=P_\al/r_\al$, where $P_\al$ is the Hilbert polynomial of $\al$, and $r_\al$ is the leading coefficient of $P_\al$, which must be positive. Then $(\tau,G,\le)$ is a stability condition on $\coh(X)$, called {\it Gieseker stability}. Gieseker stability was introduced by Gieseker \cite{Gies} for vector bundles on algebraic surfaces, and is studied in Huybrechts and Lehn~\cite[\S 1.2]{HuLe2}.

An alternative formulation of Gieseker stability that is sometimes useful is
\e
\tau(\al)\le\tau(\be)\quad\Longleftrightarrow\quad \frac{P_\al(k)}{P_\al(l)}\le \frac{P_\be(k)}{P_\be(l)}\quad\text{for $k\gg l\gg 0$.}
\label{co3eq1}
\e
\end{ex}

\begin{ex}
\label{co3ex3}
Continuing in the situation of Example \ref{co3ex2}, define
\begin{equation*}
M=\bigl\{p(t)=t^d+a_{d-1}t^{d-1}:d=0,\ldots,m,\;\> a_{d-1}\in\Q,\;\>
a_{-1}=0\bigr\}\subset G
\end{equation*}
and restrict the total order $\le$ on $G$ to $M$. Define $\mu:C(\coh(X))\ra M$ by $\mu(\al)=t^d+a_{d-1}t^{d-1}$ when $\tau(\al)=P_\al/r_\al=t^d+a_{d-1}t^{d-1}+\cdots+a_0$, that is, $\mu(\al)$ is the truncation of the polynomial $\tau(\al)$ in Example \ref{co3ex2} at its second term. Then $(\mu,M,\le)$ is a weak stability condition (but not a stability condition) on $\coh(X)$, called $\mu$-{\it stability}. It is studied in Huybrechts and Lehn \cite[\S 1.6]{HuLe2}. It dominates $(\tau,G,\le)$ in Example~\ref{co3ex2}.
\end{ex}

\begin{dfn}
\label{co3def3}
Let $\A$ be an abelian category, with quotient $K_0(\A)\twoheadrightarrow K(\A)$ as in Definition \ref{co3def1}. Following \cite[Def.~4.10]{Joyc6}, if $(\tau,T,\le),(\ti\tau,\ti T,\le)$ are weak stability conditions on $\A$, we say that $(\ti\tau,\ti T,\le)$ {\it dominates\/} $(\tau,T,\leq)$ if $\tau (\al)\leq\tau(\be)$ implies $\ti\tau(\al)\leq\ti\tau(\be)$ for all $\al,\be\in C(\A)$. Then for $0\ne A\in\A$, we see that
\e
\text{$A$ $\ti\tau$-stable} \Longra \text{$A$ $\tau$-stable} \Longra \text{$A$ $\tau$-semistable} \Longra \text{$A$ $\ti\tau$-semistable.}
\label{co3eq2}
\e
If $\A$ is $\C$-linear with moduli stacks $\M,\M^\pl$ as in Definition \ref{co3def2}, this gives
\begin{equation*}
\M_\al^\rst(\ti\tau)\subseteq \M_\al^\rst(\tau)\subseteq \M_\al^\ss(\tau)\subseteq \M_\al^\ss(\ti\tau)
\subseteq\M^\pl_\al,\qquad \al\in C(\A).
\end{equation*}

Now let $\haA$ be another abelian category, with quotient $K_0(\haA)\twoheadrightarrow K(\haA)$, and $(\hat\tau,\hat T,\le)$ be a weak stability condition on $\haA$. Suppose $F:\A\ra\haA$ is an exact functor. We say that $(\hat\tau,\hat T,\le)$ {\it dominates\/ $(\tau,T,\leq)$ over\/} $F$ if:
\begin{itemize}
\setlength{\itemsep}{0pt}
\setlength{\parsep}{0pt}
\item[(i)] The induced group morphism $F_*:K_0(\A)\ra K_0(\haA)$ descends to a group morphism $F_*:K(\A)\ra K(\haA)$. 
\item[(ii)] Suppose $A\in\A$ with $F(A)=\hat A$ in $\haA$, and $\hat A'\subseteq \hat A$ is a subobject in $\haA$. Then there exists a subobject $A'\subseteq A$ in $\A$ with $F(A')=\hat A'$.
\item[(iii)] Suppose $\al,\be\in C(\A)$ with $F_*(\al),F_*(\be)\in C(\haA)$. Then $\tau (\al)\leq\tau(\be)$ implies $\hat\tau(F_*(\al))\leq\hat\tau(F_*(\be))$.
\end{itemize}
Generalizing \eq{co3eq2}, it is now easy to prove that if $0\ne A\in\A$ with $0\ne \hat A=F(A)\in\haA$ then
$A$ $\tau$-semistable implies $\hat A$ $\hat\tau$-semistable.

If $\A,\haA$ are $\C$-linear with moduli stacks $\M,\M^\pl$ for $\A$ and $\haM,\haM^\pl$ for $\haA$, and $F$ induces morphisms $F_*:\M\ra\haM$, $F_*^\pl:\M^\pl\ra\haM^\pl$, this implies that
\begin{equation*}
F_*^\pl:\M_\al^\ss(\tau)\longra \haM_\be^\ss(\hat\tau)\subseteq\haM_\be^\pl \;\> \text{if $\al\in C(\A)$ with $\be=F_*(\al)\in C(\haA)$.}
\end{equation*}

Note that we allow the possibility that $0\ne A\in\A$ with $F(A)=0$ in $\haA$, so that we can have $\al\in C(\A)$ with $F_*(\al)=0\notin C(\haA)$.
\end{dfn}

\begin{ex}
\label{co3ex4}
Let $\A,\haA,K(\A),K(\haA),F$ be as in Definition \ref{co3def3} and satisfy (i),(ii). Suppose also $0\ne A\in\A$ implies that $F(A)\ne 0$ in $\haA$. This is equivalent to requiring that $F_*:K(\A)\ra K(\haA)$ maps~$C(\A)\ra C(\haA)$. 

It is now easy to show that if $(\hat\tau,\hat T,\le)$ is a weak stability condition on $\haA$ then $(\hat\tau\ci F_*,\hat T,\le)$ is a weak stability condition on $\A$, and $(\hat\tau,\hat T,\le)$ dominates $(\hat\tau\ci F_*,\hat T,\le)$ over $F$. If $(\tau,T,\le)$ is another weak stability condition on $\A$, then $(\hat\tau,\hat T,\le)$ dominates $(\tau,T,\leq)$ over $F$ if and only if $(\hat\tau\ci F_*,\hat T,\le)$ dominates $(\tau,T,\leq)$, as weak stability conditions on $\A$.
\end{ex}

\subsection{\texorpdfstring{Wall crossing and coefficients $U(-;\tau,\ti\tau),\ti U(-;\tau,\ti\tau)$}{Wall crossing and coefficients U(-;τ,τ*),Ŭ(-;τ,τ*)}}
\label{co32}

In the series \cite{Joyc2,Joyc3,Joyc4,Joyc5,Joyc6,Joyc7}, the author studied motivic invariants counting $\tau$-(semi)\-stable objects in abelian categories. This was applied to Donaldson--Thomas theory of Calabi--Yau 3-folds in \cite{JoSo}. The focus of \cite{Joyc7} was wall-crossing formulae for invariants under change of weak stability condition $\tau$. These involved combinatorial coefficients $S(*;\tau,\ti\tau)\in\Z$ and $U(*;\tau, \ti\tau)\in\Q$ defined in \cite[\S 4.1]{Joyc7}. We now explain these, as they also enter our wall-crossing formulae \eq{co5eq33}--\eq{co5eq34} below for virtual classes in homology. Following \cite[\S 3.3]{JoSo}, we have changed some notation
from~\cite{Joyc7}.

\begin{dfn}
\label{co3def4}
Let $\A$ be an abelian category, and choose $K_0(\A)\twoheadrightarrow K(\A)$ as in Definition \ref{co3def1}. Let
$(\tau,T,\le),(\ti\tau,\ti T,\le)$ be weak stability conditions on
$\A$.

Let $n\ge 1$ and $\al_1,\ldots,\al_n\in C(\A)$. If for all
$i=1,\ldots,n-1$ we have either
\begin{itemize}
\setlength{\itemsep}{0pt}
\setlength{\parsep}{0pt}
\item[(a)] $\tau(\al_i)\le\tau(\al_{i+1})$ and
$\ti\tau(\al_1+\cdots+\al_i)>\ti\tau(\al_{i+1}+\cdots+\al_n)$, or
\item[(b)] $\tau(\al_i)>\tau(\al_{i+1})$ and~$\ti\tau(\al_1+\cdots+\al_i)\le\ti\tau(\al_{i+1}+\cdots+\al_n)$,
\end{itemize}
then define $S(\al_1,\ldots,\al_n;\tau,\ti\tau)=(-1)^r$, where $r$ is the number of $i=1,\ldots,n-1$ satisfying (a). Otherwise define $S(\al_1,\ldots,\al_n;\tau,\ti\tau)=0$. Now
define
\ea
&U(\al_1,\ldots,\al_n;\tau,\ti\tau)=
\label{co3eq3}\\
&\sum_{\begin{subarray}{l} \phantom{wiggle}\\
1\le l\le m\le n,\;\> 0=a_0<a_1<\cdots<a_m=n,\;\>
0=b_0<b_1<\cdots<b_l=m:\\
\text{Define $\be_1,\ldots,\be_m\in C(\A)$ by
$\be_i=\al_{a_{i-1}+1}+\cdots+\al_{a_i}$.}\\
\text{Define $\ga_1,\ldots,\ga_l\in C(\A)$ by
$\ga_i=\be_{b_{i-1}+1}+\cdots+\be_{b_i}$.}\\
\text{We require $\tau(\be_i)=\tau(\al_j)$, $i=1,\ldots,m$,
$a_{i-1}<j\le a_i$,}\\
\text{and $\ti\tau(\ga_i)=\ti\tau(\al_1+\cdots+\al_n)$,
$i=1,\ldots,l$}
\end{subarray}
\!\!\!\!\!\!\!\!\!\!\!\!\!\!\!\!\!\!\!\!\!\!\!\!\!\!\!\!\!\!\!\!\!
\!\!\!\!\!\!\!\!\!\!\!\!\!\!\!\!\!\!\!\!\!\!\!\!\!\!\!\!\!\!\!\!\!
\!\!\!\!\!\!\!\!\!\!\!\!\!\!\!\!\!\!\!\!}
\begin{aligned}[t]
\frac{(-1)^{l-1}}{l}\cdot\prod\nolimits_{i=1}^lS(\be_{b_{i-1}+1},
\be_{b_{i-1}+2},\ldots,\be_{b_i}; \tau,\ti\tau)&\\
\cdot\prod_{i=1}^m\frac{1}{(a_i-a_{i-1})!}&\,.
\end{aligned}
\nonumber
\ea
\end{dfn}

Here are some properties of the coefficients $U(-)$. Equations \eq{co3eq4}--\eq{co3eq5} come from \cite[Th.~4.8]{Joyc7}, and \eq{co3eq6} follows from \eq{co3eq3} and \cite[eq.~(60)]{Joyc7}.

\begin{thm}
\label{co3thm2}
Let\/ $\A$ be an abelian category, and choose $K_0(\A)\twoheadrightarrow K(\A)$. Suppose $(\tau,T,\le),(\hat\tau,\hat T,\le),\ab(\ti\tau,\ti T,\le)$ are weak stability conditions on
$\A$. Then for all\/ $\al_1,\ldots,\al_n\in C(\A)$ we have
\ea
U(\al_1,\ldots,\al_n;\tau,\tau)
=\begin{cases} 1, & \text{$n=1$,} \\
0, & \text{otherwise,} \end{cases}
\label{co3eq4}
\\
\begin{aligned}
\sum_{\begin{subarray}{l} m,\; a_0,\ldots,a_m: \; m=1,\ldots,n,\\ 0=a_0<a_1<\cdots<a_m=n, \\
\text{set\/ $\be_j=\al_{a_{j-1}+1}+\cdots+\al_{a_j}$,} \\  
j=1,\ldots,m \end{subarray}}\,
\begin{aligned}[t]
&U(\be_1,\ldots,\be_m;\hat\tau,\ti\tau)\cdot\\
&\ts\prod_{j=1}^mU\bigl(\al_{a_{j-1}+1},\al_{a_{j-1}+2},\ldots,\al_{a_j};\tau,\hat\tau)\\
&\qquad =U(\al_1,\ldots,\al_n;\tau,\ti\tau).
\end{aligned}
\end{aligned}
\label{co3eq5}
\ea

If also\/ $(\ti\tau,\ti T,\le)$ dominates\/ $(\tau,T,\le),$ as in Definition\/ {\rm\ref{co3def3},} then
\e
U(\al_1,\ldots,\al_n;\tau,\ti\tau)\!=\!U(\al_1,\ldots,\al_n;\ti\tau,\tau)\!=\!0\;\>\text{unless\/ $\ti\tau(\al_1)\!=\!\cdots\!=\!\ti\tau(\al_n)$.}
\label{co3eq6}
\e
\end{thm}

The next theorem is proved in \cite[Th.~5.4]{Joyc7} (see also \cite[Th.~3.14]{JoSo}). It describes a property of the coefficients $U(-;\tau,\ti\tau)$, it does not matter what $\cL$ and $\ep^\al(\tau),\ep^\al(\ti\tau)$ are. We have no explicit definition for $\ti U(\al_1,\ldots,\al_n;\tau,\ti\tau)$, we only show that \eq{co3eq7} can be rewritten in the form~\eq{co3eq8}.

\begin{thm}
\label{co3thm3}
Work in the situation of Definition\/ {\rm\ref{co3def4}}. Let\/ $\cL$ be a Lie algebra over $\Q,$ and write $U(\cL)$ for its universal enveloping algebra, with product $*$. Suppose we are given elements $\ep^\al(\tau),\ep^\al(\ti\tau)\in\cL$ for $\al\in C(\A)$ satisfying
\e
\begin{gathered}
\ep^\al(\ti\tau)= \!\!\!\!\!\!\!
\sum_{\begin{subarray}{l}n\ge 1,\;\al_1,\ldots,\al_n\in
C(\A):\\ \al_1+\cdots+\al_n=\al\end{subarray}} \!\!\!\!\!\!\!
\begin{aligned}[t]
U(\al_1,&\ldots,\al_n;\tau,\ti\tau)\cdot{}\\
&\ep^{\al_1}(\tau)*\ep^{\al_2}(\tau)*\cdots* \ep^{\al_n}(\tau)
\end{aligned}
\end{gathered}
\label{co3eq7}
\e
for each\/ $\al\in C(\A),$ with only finitely many nonzero terms. Then\/ \eq{co3eq7} may be rewritten as an equation in the Lie algebra $\cL$ using the Lie bracket\/ $[\,,\,]$. That is, we may rewrite \eq{co3eq7} in the form
\e
\begin{gathered}
\ep{}^\al(\ti\tau)= \!\!\!\!\!\!\!
\sum_{\begin{subarray}{l}n\ge 1,\;\al_1,\ldots,\al_n\in
C(\A):\\ \al_1+\cdots+\al_n=\al\end{subarray}} \!\!\!\!\!\!\!
\begin{aligned}[t]
\ti U(\al_1,&\ldots,\al_n;\tau,\ti\tau)\,\cdot\\
&[[\cdots[[\ep{}^{\al_1}(\tau),\ep{}^{\al_2}(\tau)],\ep{}^{\al_3}(\tau)],
\ldots],\ep{}^{\al_n}(\tau)],
\end{aligned}
\end{gathered}
\label{co3eq8}
\e
for some system of combinatorial coefficients\/ $\ti U(\al_1,\ldots,\al_n;\tau,\ti\tau)\in\Q,$ with only finitely many nonzero terms, such that if we expand\/ $[f,g]=f*g-g*f$ then \eq{co3eq8} becomes \eq{co3eq7}.
\end{thm}

This enters the motivic invariants theory of \cite{Joyc2,Joyc3,Joyc4,Joyc5,Joyc6,Joyc7} in the following way. Writing $\M$ for the moduli stack of objects in $\A$, we define a Ringel--Hall algebra of `stack functions' $\SFa(\M)$, with associative product $*$. For each `permissible' weak stability condition $(\tau,T,\le)$ on $\A$ we define elements $\ep^\al(\tau)\in\SFa(\M)$ for $\al\in C(\A)$ \cite[Def.s 7.6 \& 8.1]{Joyc6}, which `count' $\tau$-semistable objects in $\A$ in class $\al$, and which transform according to \eq{co3eq7} under change of stability condition \cite[Th.~5.4]{Joyc7}. We show there is a Lie subalgebra $\SFai(\M)\subset\SFa(\M)$ such that the $\ep^\al(\tau)$ lie in $\SFai(\M)$, and \eq{co3eq7} may be rewritten as a Lie algebra equation \eq{co3eq8} in $\SFai(\M)$.

These results are applied in \cite[\S 6.4--\S 6.5]{Joyc7} and \cite{JoSo} in the following way. Suppose we can define a Lie algebra morphism $\Psi:\SFai(\M)\ra \cL_{K(\A)}$, where $\cL_{K(\A)}$ is an explicit Lie algebra, often of the form $\an{\la^\al:\al\in K(\A)}_R$ for some commutative ring $R$, with Lie bracket $[\la^\al,\la^\be]=c_{\al,\be}\la^{\al+\be}$ for coefficients $c_{\al,\be}$ in $R$. Then we may define motivic invariants $J^\al(\tau)\in R$ by $\Psi(\bar\ep{}^\al(\tau))=J^\al(\tau)\la^\al$. Applying $\Psi$ to \eq{co3eq7}, interpreted using Theorem \ref{co3thm3}, then gives a wall-crossing formula for the invariants $J^\al(\tau)$. This is used in \cite[Th.~3.14]{JoSo} to prove a wall-crossing formula for Donaldson--Thomas invariants of Calabi--Yau 3-folds.

The next proposition may be deduced from \cite[Th.~4.6]{Joyc7}, which is the analogous result for the~$S(\al_1,\ldots,\al_n;\tau,\ti\tau)$.

\begin{prop}
\label{co3prop2}
In Definition {\rm\ref{co3def4},} suppose $U(\al_1,\ldots,\al_n;\tau,\ti\tau)\ne 0$. Then there exist\/ $k,l=1,\ldots,n$ such that\/ $\tau(\al_k)\le\tau(\al_i)\le\tau(\al_l)$ for all\/ $i=1,\ldots,n$ and\/~$\ti\tau(\al_k)\ge\ti\tau(\al_1+\cdots+\al_n)\ge\ti\tau(\al_l)$.	
\end{prop}

Next we prove four new results.

\begin{prop}
\label{co3prop3}
In the situation of Definition {\rm\ref{co3def4},} let\/ $(\tau,T,\le),(\ti\tau,\ti T,\le)$ be weak stability conditions on $\A,$ and suppose $n\ge 2$ and for some $\es\ne I\subsetneq\{1,\ldots,n\}$ we have $\tau(\al_i)<\tau(\al_j)$ for all\/ $i\in I$ and\/ $j\in \{1,\ldots,n\}\sm I,$ and\/ $\ti\tau(\al_i)<\ti\tau(\al_1+\cdots+\al_n)$ for all\/ $i\in I$. Then
\begin{equation*}
S(\al_1,\ldots,\al_n;\tau,\ti\tau)=U(\al_1,\ldots,\al_n;\tau,\ti\tau)=\ti U(\al_1,\ldots,\al_n;\tau,\ti\tau)=0.
\end{equation*}

The same holds if instead\/ $\tau(\al_i)>\tau(\al_j)$ and\/ $\ti\tau(\al_i)>\ti\tau(\al_1+\cdots+\al_n)$ above.
\end{prop}

\begin{proof} To show $S(\al_1,\ldots,\al_n;\tau,\ti\tau)\!=\!U(\al_1,\ldots,\al_n;\tau,\ti\tau)\!=\!0$ we divide into~cases:
\begin{itemize}
\setlength{\itemsep}{0pt}
\setlength{\parsep}{0pt}
\item[(i)] For some $1\le a<n$ we have $i\in I$ for $1\le i\le a$ and $a+1\notin I$;
\item[(ii)] For some $1\le a<b<n$ we have $i\notin I$ for $1\le i\le a$, $i\in I$ for $a<i\le b$, and $b+1\notin I$, and $\ti\tau(\al_1+\cdots+\al_a)>\ti\tau(\al_1+\cdots+\al_n)$;
\item[(iii)] For some $1\le a<b<n$ we have $i\notin I$ for $1\le i\le a$, $i\in I$ for $a<i\le b$, and $b+1\notin I$, and $\ti\tau(\al_1+\cdots+\al_a)\le\ti\tau(\al_1+\cdots+\al_n)$; and 
\item[(iv)] For some $1\le a<n$ we have $I=\{a+1,a+2,\ldots,n\}$.
\end{itemize}

In case (i) we have $\tau(\al_a)<\tau(\al_{a+1})$ and $\ti\tau(\al_1+\cdots+\al_a)<\ti\tau(\al_1+\cdots+\al_n)$ since $\ti\tau(\al_i)<\ti\tau(\al_1+\cdots+\al_n)$ for $i=1,\ldots,a$. This implies that $\ti\tau(\al_1+\cdots+\al_a)<\ti\tau(\al_{a+1}+\cdots+\al_n)$, so neither Definition \ref{co3def4}(a),(b) hold when $i=a$, and $S(\al_1,\ldots,\al_n;\tau,\ti\tau)=0$. For each term $l,m,(a_i),(b_j)$ in the sum \eq{co3eq3} for $U(\al_1,\ldots,\al_n;\tau,\ti\tau)$, as $\tau(\al_a)<\tau(\al_{a+1})$ we must have $a_p=a$, $a_{p+1}=a+1$ for some $1\le p<m$, as otherwise would contradict $\tau(\be_i)=\tau(\al_j)$ for $a_{i-1}<j\le a_i$. We cannot have $p=b_q$ for some $q=1,\ldots,l$, as then $\ti\tau(\ga_1+\cdots+\ga_q)=\ti\tau(\al_1+\cdots+\al_a)<\ti\tau(\al_1+\cdots+\al_n)$ contradicts $\ti\tau(\ga_i)=\ti\tau(\al_1+\cdots+\al_n)$. Thus $b_{q-1}<p<b_q$ for some unique $q=1,\ldots,l$, and in a similar way to the first part we find that $S(\be_{b_{q-1}+1},\be_{b_{q-1}+2},\ldots,\be_{b_q};\tau,\ti\tau)=0$. Hence each term in \eq{co3eq3} is zero, and $U(\al_1,\ldots,\al_n;\tau,\ti\tau)=0$.

In case (ii) we have $\tau(\al_a)>\tau(\al_{a+1})$ and $\ti\tau(\al_1+\cdots+\al_a)>\ti\tau(\al_1+\cdots+\al_n)$, which implies that $\ti\tau(\al_1+\cdots+\al_a)>\ti\tau(\al_{a+1}+\cdots+\al_n)$, so neither Definition \ref{co3def4}(a),(b) hold when $i=a$, and $S(\al_1,\ldots,\al_n;\tau,\ti\tau)=0$. We prove $U(\al_1,\ldots,\al_n;\tau,\ti\tau)=0$ in a similar way to (i).

In case (iii) we have $\tau(\al_b)<\tau(\al_{b+1})$ and $\ti\tau(\al_1+\cdots+\al_b)\le\ti\tau(\al_1+\cdots+\al_n)$ since $\ti\tau(\al_1+\cdots+\al_a)\le\ti\tau(\al_1+\cdots+\al_n)$ and $\ti\tau(\al_i)<\ti\tau(\al_1+\cdots+\al_n)$ for $a<i\le b$. This implies that $\ti\tau(\al_1+\cdots+\al_b)\le\ti\tau(\al_{b+1}+\cdots+\al_n)$, and we prove $S(\al_1,\ldots,\al_n;\tau,\ti\tau)=U(\al_1,\ldots,\al_n;\tau,\ti\tau)=0$ in a similar way to (i).

In case (iv) we have $\tau(\al_a)>\tau(\al_{a+1})$ and $\ti\tau(\al_1+\cdots+\al_n)>\ti\tau(\al_{a+1}+\cdots+\al_n)$ since $\ti\tau(\al_i)<\ti\tau(\al_1+\cdots+\al_n)$ for $a<i\le n$. This implies that $\ti\tau(\al_1+\cdots+\al_a)>\ti\tau(\al_{a+1}+\cdots+\al_n)$, and we prove $S(\al_1,\ldots,\al_n;\tau,\ti\tau)=U(\al_1,\ldots,\al_n;\tau,\ti\tau)=0$ in a similar way to (i). We have now shown that $S(\al_1,\ldots,\al_n;\tau,\ti\tau)=U(\al_1,\ldots,\al_n;\tau,\ti\tau)=0$ in all cases. As this holds for all permutations of $\al_1,\ldots,\al_n$, Theorem \ref{co3thm3} gives $\ti U(\al_1,\ldots,\al_n;\tau,\ti\tau)=0$.

For the last part, we reverse directions of inequalities in the proof above, and also replace (ii)--(iii) above by
\begin{itemize}
\setlength{\itemsep}{0pt}
\setlength{\parsep}{0pt}
\item[(ii$)'$] For some $1\le a<b<n$ we have $i\notin I$ for $1\le i\le a$, $i\in I$ for $a<i\le b$, and $b+1\notin I$, and $\ti\tau(\al_1+\cdots+\al_a)\le\ti\tau(\al_1+\cdots+\al_n)$;
\item[(iii$)'$] For some $1\le a<b<n$ we have $i\notin I$ for $1\le i\le a$, $i\in I$ for $a<i\le b$, and $b+1\notin I$, and $\ti\tau(\al_1+\cdots+\al_a)>\ti\tau(\al_1+\cdots+\al_n)$.\qedhere
\end{itemize}
\end{proof}

\begin{dfn}
\label{co3def5}
Work in the situation of Definition \ref{co3def4}. Let $(\tau_t,T_t,\le)_{t\in[0,1]}$ be a family of weak stability conditions on $\A$ parametrized by $t\in[0,1]$. For all $\al,\be\in C(\A)$, write $S_{\al,\be}^>, S_{\al,\be}^<,S_{\al,\be}^=$ for the subsets of $t\in[0,1]$ for which $\tau_t(\al)>\tau_t(\be)$, and $\tau_t(\al)<\tau_t(\be)$, and $\tau_t(\al)=\tau_t(\be)$, respectively, in $T_t$. Then we have a disjoint union~$[0,1]=S_{\al,\be}^>\amalg S_{\al,\be}^<\amalg S_{\al,\be}^=$. 

We call $(\tau_t,T_t,\le)_{t\in[0,1]}$ a {\it continuous family\/} of weak stability conditions on $\A$ if $S_{\al,\be}^>, S_{\al,\be}^<$ are open and $S_{\al,\be}^=$ is closed, and all three have only finitely many connected components (which must be open or closed intervals in $[0,1]$), for all $\al,\be\in C(\A)$. We can replace $[0,1]$ by other intervals in $\R$ in the obvious way.	
\end{dfn}

\begin{prop}
\label{co3prop4}
Let\/ $(\tau_t,T_t,\le)_{t\in(t_0-\de,t_0+\de)}$ be a continuous family of weak stability conditions on $\A$. Suppose\/ $\al_1,\ldots,\al_n\in C(\A),$ and let\/ $t_0-\de\ab<\ab t_-\ab\le\ab t_0\ab\le\ab t_+\ab<\ab t_0+\de$ with\/ $\md{t_\pm-t_0}$ sufficiently small (depending on $\al_1,\ldots,\al_n$). Then unless $\tau_{t_0}(\al_1)=\cdots=\tau_{t_0}(\al_n)$ we have
\begin{equation*}
S(\al_1,\ldots,\al_n;\tau_{t_-},\tau_{t_+})\!=\!U(\al_1,\ldots,\al_n;\tau_{t_-},\tau_{t_+})\!=\!\ti U(\al_1,\ldots,\al_n;\tau_{t_-},\tau_{t_+})\!=\!0.
\end{equation*}
\end{prop}

\begin{proof} Suppose it is not true that $\tau_{t_0}(\al_1)=\cdots=\tau_{t_0}(\al_n)$. Then $\min_{i=1}^n\tau_{t_0}(\al_i)\ab<\max_{i=1}^n\tau_{t_0}(\al_i)$. Divide into (non-exclusive) cases
\begin{itemize}
\setlength{\itemsep}{0pt}
\setlength{\parsep}{0pt}
\item[(a)] $\min_{i=1}^n\tau_{t_0}(\al_i)<\tau_{t_0}(\al_1+\cdots+\al_n)$; and
\item[(b)] $\tau_{t_0}(\al_1+\cdots+\al_n)<\max_{i=1}^n\tau_{t_0}(\al_i)$.
\end{itemize}
Here if $(\tau_{t_0},T_{t_0},\le)$ is a stability condition then both (a),(b) hold, but if it is only a weak stability condition, then only one need hold.

In case (a), let $I$ be the set of $i=1,\ldots,n$ for which $\tau_{t_0}(\al_i)$ is minimal. Then $\es\ne I\subsetneq\{1,\ldots,n\}$. We have $\tau_{t_0}(\al_i)<\tau_{t_0}(\al_j)$ for $i\in I$ and $j\in \{1,\ldots,n\}\sm I$, and $\tau_{t_0}(\al_i)<\tau_{t_0}(\al_1+\cdots+\al_n)$ for $i\in I$. If $\md{t_\pm-t_0}$ are sufficiently small, these imply that $\tau_{t_-}(\al_i)<\tau_{t_-}(\al_j)$ for $i\in I$ and $j\in \{1,\ldots,n\}\sm I$, and $\tau_{t_+}(\al_i)<\tau_{t_+}(\al_1+\cdots+\al_n)$ for $i\in I$, by continuity of $(\tau_t,T_t,\le)_{t\in(t_0-\de,t_0+\de)}$. The result now follows from the first part of Proposition~\ref{co3prop3}. 

For case (b), we instead take $I$ to be the set of $i=1,\ldots,n$ for which $\tau_{t_0}(\al_i)$ is maximal, and use the second part of Proposition~\ref{co3prop3}.
\end{proof}

\begin{prop}
\label{co3prop5}
In the situation of Definition {\rm\ref{co3def4},} let\/ $(\tau,T,\le),(\ti\tau,\ti T,\le)$ be weak stability conditions on $\A,$ and suppose $\al_1,\ldots,\al_n\in C(\A)$ with\/ $\tau(\al_1)>\tau(\al_2)>\cdots>\tau(\al_n)$ and\/ $\ti\tau(\al_1+\cdots+\al_i)\le \ti\tau(\al_{i+1}+\cdots+\al_n)$ for\/ $1\le i<n$.

Write $\al=\al_1+\cdots+\al_n,$ and define $N$ to be the number of\/ $i=1,\ldots,n-1$ such that\/ $\ti\tau(\al_1+\cdots+\al_i)=\ti\tau(\al_{i+1}+\cdots+\al_n)=\ti\tau(\al),$ and if\/ $1\le j<i$ then $\ti\tau(\al_{j+1}+\cdots+\al_i)\ge\ti\tau(\al),$ and if\/ $i<j<n$ then $\ti\tau(\al_{i+1}+\cdots+\al_j)\le\ti\tau(\al)$. (Here if\/ $(\ti\tau,\ti T,\le)$ is a stability condition then the first condition implies the others.) Then\/~$U(\al_1,\ldots,\al_n;\tau,\ti\tau)=1/(N+1)$.	
\end{prop}

\begin{proof} First note that if $(\ti\tau,\ti T,\le)$ is a stability condition and $1\le i<n$ with 
$\ti\tau(\al_1+\cdots+\al_i)=\ti\tau(\al_{i+1}+\cdots+\al_n)=\ti\tau(\al)$, and $1\le j<i$, then $\ti\tau(\al_1+\cdots+\al_j)\le \ti\tau(\al_{j+1}+\cdots+\al_n)$ implies that $\ti\tau(\al_{j+1}+\cdots+\al_n)\ge\ti\tau(\al)$, and this and $\ti\tau(\al_{i+1}+\cdots+\al_n)=\ti\tau(\al)$ force $\ti\tau(\al_{j+1}+\cdots+\al_i)\ge\ti\tau(\al)$, as claimed. Similarly, if $i<j<n$ then~$\ti\tau(\al_{i+1}+\cdots+\al_j)\le\ti\tau(\al)$.

Write $\{i_1,\ldots,i_N\}$ for the set of $i=1,\ldots,n-1$ satisfying the conditions above. Suppose $l,m,n,a_i,b_i,\be_i,\ga_i$ are data giving a nonzero term in the sum for $U(\al_1,\ldots,\al_n;\tau,\ti\tau)$ in \eq{co3eq3}. As $\tau(\al_1)>\cdots>\tau(\al_n)$ we must have $m=n$, $a_i=i$ and $\be_i=\al_i$ for~$i=1,\ldots,n$. 

Suppose for a contradiction that for some $1\le i<l$ we have $b_i\notin \{i_1,\ldots,i_N\}$. As $\ti\tau(\ga_j)=\ti\tau(\al)$ for all $j$ we have $\ti\tau(\ga_1+\cdots+\ga_i)=\ti\tau(\ga_{i+1}+\cdots+\ga_l)$, so $\ti\tau(\al_1+\cdots+\al_{b_i})=\ti\tau(\al_{b_i+1}+\cdots+\al_n)$. Thus the definition of $\{i_1,\ldots,i_N\}$ implies that either (i) there exists $1\le j<b_i$ with $\ti\tau(\al_{j+1}+\cdots+\al_{b_i})<\ti\tau(\al)$, or (ii) there exists $b_i<j<n$ with $\ti\tau(\al_{b_i+1}+\cdots+\al_j)>\ti\tau(\al)$. 

In case (i) we cannot have $j=b_k$ for any $k$ as then $\ti\tau(\al_{j+1}+\cdots+\al_{b_i})=\ti\tau(\ga_{k+1}+\cdots+\ga_i)=\ti\tau(\al)$. So for some $1\le k\le i$ we have $b_{k-1}<j<b_k$. From $\ti\tau(\al_{j+1}+\cdots+\al_{b_i})<\ti\tau(\al)$ and $\ti\tau(\al_{b_k+1}+\cdots+\al_{b_i})=\ti\tau(\ga_{k+1}+\cdots+\ga_i)=\ti\tau(\al)$ we see that $\ti\tau(\al_{j+1}+\cdots+\al_{b_k})<\ti\tau(\al)$, so as $\ti\tau(\al_{b_{k-1}+1}+\cdots+\al_{b_k})=\ti\tau(\ga_k)=\ti\tau(\al)$ we have $\ti\tau(\al_{b_{k-1}+1}+\cdots+\al_j)>\ti\tau(\al_{j+1}+\cdots+\al_{b_k})$. This, $\tau(\al_j)>\tau(\al_{j+1})$ and $\be_i=\al_i$ show that $S(\be_{b_{k-1}+1},\ldots,\be_{b_k};\tau,\ti\tau)=0$ by Definition \ref{co3def4}, contradicting that $l,m,\ldots,\ga_i$ give a nonzero term in \eq{co3eq3}. Case (ii) gives a contradiction in a similar way.

Thus $b_1,\ldots,b_{l-1}\in \{i_1,\ldots,i_N\}$. For any choice of $b_1<b_2<\cdots<b_{l-1}$ in $\{i_1,\ldots,i_N\}$ we see from Definition \ref{co3def4} that $S(\be_{b_{i-1}+1},\be_{b_{i-1}+2},\ldots,\be_{b_i}; \tau,\ti\tau)=1$ for $i=1,\ldots,l$, since $\be_j=\al_j$, $\tau(\al_{b_{i-1}+1})>\cdots>\tau(\al_{b_i})$ and $\ti\tau(\al_{b_{i-1}+1}+\cdots+\al_j)\le\ti\tau(\al_{j+1}+\cdots+\al_{b_i})$ for $b_{i-1}<j<b_i$. As there are $\binom{N}{l-1}$ such choices of $b_1,\ldots,b_{l-1}$, equation \eq{co3eq3} implies that
\begin{align*}
U(\al_1,\ldots,\al_n;\tau,\ti\tau)&=\sum_{l=1}^{N+1}\frac{(-1)^{l-1}}{l}\binom{N}{l-1}=\frac{1}{N+1}\sum_{l=1}^{N+1}(-1)^{l-1}\binom{N+1}{l}\\
&=\frac{1}{N+1}\bigl(1-(1-1)^{N+1}\bigr)=\frac{1}{N+1}.\qedhere
\end{align*}
\end{proof}

\begin{cor}
\label{co3cor1}
In the situation of Definition {\rm\ref{co3def4},} let\/ $(\tau,T,\le),(\ti\tau,\ti T,\le)$ be weak stability conditions on $\A,$ and suppose $E\in\A$ is $\ti\tau$-semistable and has $\tau$-Harder--Narasimhan filtration $0=E_0\subsetneq E_1\subsetneq\cdots \subsetneq E_n=E$. Set\/ $F_i=E_i/E_{i-1}$ and\/ $\al_i=\lb F_i\rb\in C(\A)$ for $i=1,\ldots,n$. Then $U(\al_1,\ldots,\al_n;\tau,\ti\tau)\ne 0$.	
\end{cor}

\begin{proof} We have $\tau(\al_1)>\cdots>\tau(\al_n)$ as $0=E_0\subsetneq\cdots \subsetneq E_n=E$ is a $\tau$-Harder--Narasimhan filtration. Also $\ti\tau(\al_1+\cdots+\al_i)=\ti\tau(E_i)\le\ti\tau(E/E_i)=\ti\tau(\al_{i+1}+\cdots+\al_n)$ for $1\le i<n$ as $E$ is $\ti\tau$-semistable. The result then follows from Proposition~\ref{co3prop5}.
\end{proof}

\subsection{Criteria for proper stable=semistable moduli stacks}
\label{co33}

In this section we will consider the following situation:

\begin{dfn}
\label{co3def6}
Let $\A$ is a $\C$-linear abelian category, $K_0(\A)\twoheadrightarrow K(\A)$ be a surjective quotient, and $(\tau,T,\le)$ be a weak stability condition on $\A$, as in \S\ref{co31}. Suppose we can form a moduli stack $\M$ and a `projective linear' moduli stack $\M^\pl$ of objects in $\A$, which are Artin $\C$-stacks locally of finite type, such that $K_0(\A)\twoheadrightarrow K(\A)$ induces decompositions $\M=\coprod_{\al\in K(\A)}\M_\al$, $\M^\pl=\coprod_{\al\in K(\A)}\M_\al^\pl$ with $\M_\al,\M_\al^\pl$ open and closed in $\M,\M^\pl$, where $\Pi_\al^\pl:\M_\al\ra\M_\al^\pl$ is a principal $[*/\bG_m]$-bundle for $\al\ne 0$. Suppose that for $\al\in C(\A)$ we have open $\C$-substacks $\M_\al^\rst(\tau)\subseteq\M_\al^\ss(\tau)\subseteq\M_\al^\pl$ of $\tau$-(semi)stable objects in class~$\al$.
\end{dfn}

As in Assumption \ref{co5ass2}(g),(h) below, it will be very important in our theory that if $\al\in C(\A)$ with $\M_\al^\rst(\tau)=\M_\al^\ss(\tau)$ then $\M_\al^\ss(\tau)$ should be a {\it proper algebraic space}, since as in Theorem \ref{co2thm1} this is a necessary condition for forming a virtual class~$[\M_\al^\ss(\tau)]_\virt$.

We will now explain some results of Alper--Halpern-Leistner--Heinloth \cite[\S 7]{AHLH}, building on work of Alper \cite{Alpe}, Halpern--Leistner \cite{Halp1,Halp2}, and Artin--Zhang \cite{ArZh}, which enable us to give some general sufficient conditions on $\A,(\tau,T,\le)$, true in many interesting examples, such that if $\M_\al^\rst(\tau)=\M_\al^\ss(\tau)$ then $\M_\al^\ss(\tau)$ is a proper algebraic space. The theory is quite technical.

\subsubsection{Good moduli spaces}
\label{co331}

First we describe Alper's theory of `good moduli spaces' for Artin stacks~\cite{Alpe}.

\begin{dfn}
\label{co3def7}
Let $Y$ be an Artin $\C$-stack, and $Z$ be an algebraic $\C$-space, and $\phi:Y\ra Z$ be a morphism. We call $\phi:Y\ra Z$ a {\it good moduli space}, or (taking $\phi$ to be implicitly given) we call $Z$ a {\it good moduli space for\/} $Y$, if:  
\begin{itemize}
\setlength{\itemsep}{0pt}
\setlength{\parsep}{0pt}
\item[(i)] $\phi$ is quasi-compact.
\item[(ii)] Pushforward $\phi_*:\qcoh(Y)\ra\qcoh(Z)$ on quasicoherent sheaves is exact.
\item[(iii)] The induced morphism $\O_Z\ra\phi_*(\O_Y)$ is an isomorphism.
\end{itemize}

Good moduli spaces have the following properties. Let $\phi:Y\ra Z$ be a good moduli space, with $Y$ locally of finite type. Then:
\begin{itemize}
\setlength{\itemsep}{0pt}
\setlength{\parsep}{0pt}
\item[(a)] $\phi$ is surjective and universally closed \cite[Th.~4.16]{Alpe}, so $Z$ has the quotient topology. 
\item[(b)] Two $\C$-points $y,y'$ in $Y$ are identified in $Z$ if and only if their closures $\ov{\{y\}},\ov{\{y'\}}$ intersect in $Y$, \cite[Th.~4.16]{Alpe}.
\item[(c)] For any closed $\C$-point $z$ in $Z$ (if $Z$ is separated, this means any $\C$-point $z$ in $Z$), there is a unique closed $\C$-point $y$ in $Y$ with $\phi(y)=z$, \cite[Prop.~9.1]{Alpe}. Also $\Iso_Y(y)$ is reductive \cite[Prop.~12.14]{Alpe}, and has strictly larger dimension than $\Iso_Y(y')$ for any other $\C$-point $y'$ in $Y$ with $\phi(y')=z$.
\item[(d)] $\phi:Y\ra Z$ is universal for morphisms $Y\ra S$ to algebraic spaces $S$, \cite[Th.~6.6]{Alpe}. Thus $Z,\phi$ are determined by $Y$ up to canonical isomorphism. Also, if $Y$ is an algebraic space then $\phi$ is an isomorphism.
\item[(e)] If $Y$ is finite type then $Z$ is finite type, \cite[Th.~4.16]{Alpe}.
\end{itemize}
\end{dfn}

\begin{ex}
\label{co3ex5}
Suppose $V$ is a quasiprojective $\C$-scheme, $G$ a reductive algebraic $\C$-group acting on $V,$ and $L\ra V$ a linearization for the action of $G$ on $V$ in the sense of Geometric Invariant Theory (GIT) of Mumford--Fogarty--Kirwan \cite{MFK}. Then GIT defines $G$-invariant open subsets $V^\rst\subseteq V^\ss\subset V$ of (semi)stable points, a quasiprojective $\C$-scheme $V\ds^\ss G$ called the {\it GIT quotient\/} with an open $\C$-subscheme $V\ds^\rst G\subseteq V\ds^\ss G$, and a morphism of Artin stacks $\pi:[V^\ss/G]\ra V\ds^\ss G$ with $[V^\rst/G]=\pi^{-1}(V\ds^\rst G)$, where $[V^\rst/G]\subseteq[V^\ss/G]$ are the quotient stacks, such that $V\ds^\ss G$ is a coarse moduli space for $[V^\ss/G]$, and $V\ds^\rst G$ is a fine moduli space for $[V^\rst/G]$. 

Alper shows \cite[Th.~13.6]{Alpe} that $\pi:[V^\ss/G]\ra V\ds^\ss G$ and $\pi:[V^\rst/G]\ra V\ds^\rst G$ are good moduli spaces. Thus we can think of the theory of good moduli spaces as a generalization of Geometric Invariant Theory to Artin stacks $Y$ which are not quotient stacks~$[V/G]$.
\end{ex}

\begin{prop}
\label{co3prop6}
In the situation of Definition\/ {\rm\ref{co3def6},} suppose $\al\in C(\A)$ with\/ $\M_\al^\rst(\tau)=\M_\al^\ss(\tau),$ and\/ $\M_\al^\ss(\tau)$ has a proper good moduli space $\M_\al^\ss(\tau)\ra\fM_\al^\ss(\tau)$. Then $\M_\al^\ss(\tau)\ra\fM_\al^\ss(\tau)$ is an isomorphism, so $\M_\al^\ss(\tau)$ is a proper algebraic space. Since $\Pi_\al^\pl:\M_\al\ra\M_\al^\pl$ is a principal $[*/\bG_m]$-bundle, a good moduli space for $\M_\al^\ss(\tau)\subseteq\M_\al^\pl$ is also a good moduli space for $(\Pi_\al^\pl)^{-1}(\M_\al^\ss(\tau))\subseteq\M_\al,$ and vice versa. 
\end{prop}

\begin{proof}
As $\tau$-stable objects $E\in\A$ have $\Aut(E)=\bG_m$, we see that $\M_\al^\rst(\tau)=\M_\al^\ss(\tau)$ has trivial isotropy groups. Then as $\M_\al^\ss(\tau)$ has a good moduli space, Alper--Hall--Rydh \cite[Th.~4.12]{AHR} shows that $\M_\al^\ss(\tau)$ is \'etale locally modelled on a $\C$-scheme near each closed point. Therefore $\M_\al^\ss(\tau)$ is an algebraic space. The first part of the proposition now follows from Definition \ref{co3def7}(d). For the second part, Definition \ref{co3def7}(d) implies that any good moduli space for $(\Pi_\al^\pl)^{-1}(\M_\al^\ss(\tau))$ must factor via $\Pi_\al^\pl\vert_{\cdots}:(\Pi_\al^\pl)^{-1}(\M_\al^\ss(\tau))\ra\M_\al^\ss(\tau)$, and one can also check that $\Pi_\al^\pl\vert_{\cdots}$ satisfies Definition~\ref{co3def7}(i)--(iii).
\end{proof}

In Theorem \ref{co3thm4} we will give necessary and sufficient conditions from \cite{AHLH} for a suitable Artin stack $Y$ to admit a good moduli space. We will need the following definitions from~\cite{AHLH,Halp1}:

\begin{dfn}
\label{co3def8}
Suppose $Y$ is an Artin $\C$-stack, and let $R$ be a {\it discrete valuation ring\/} ({\it DVR\/}) over $\C$, that is, a principal ideal domain with exactly one nonzero maximal ideal $\m\subset R$. Later we will always assume that $R$ is {\it essentially of finite type over\/} $\C$, as in~\cite[Rem.s 3.48, 4.2, 5.5 \& Lem.s 7.15--7.17]{AHLH}.
\smallskip

\noindent{\bf(a)} We say that $Y$ {\it satisfies the valuative criterion for universal closedness with respect to\/} $R$ (also called {\it the existence part of the valuative criterion for properness}, see \cite[Lem.~7.17]{AHLH}) if in any commutative square of morphisms of $\C$-stacks `$\ra$', the dotted arrow `$\dashra$' can be filled in, though {\it not\/} necessarily uniquely:
\e
\begin{gathered}
\xymatrix@C=130pt@R=15pt{*+[r]{(\Spec R)\sm\m} \ar[r] \ar[d] & *+[l]{Y} \ar[d] \\
*+[r]{\Spec R} \ar@{..>}[ur] \ar[r]  & *+[l]{\Spec\C.} }	
\end{gathered}
\label{co3eq9}
\e
Here $\m$ is a (maximal) ideal in $R$, and thus a (closed) point in $\Spec R$.
\smallskip

\noindent{\bf(b)} As in \cite[Def.~4.16]{Halp1} and \cite[Def.~3.10]{AHLH}, we say that $Y$ is $\Th$-{\it reductive with respect to\/} $R$ if in any commutative square of morphisms of $\C$-stacks `$\ra$', the dotted arrow `$\dashra$' can be uniquely filled in:
\begin{equation*}
\xymatrix@C=130pt@R=15pt{*+[r]{[(\Spec R[x]\sm (\m,x))/\bG_m(R)]} \ar[r] \ar[d] & *+[l]{Y} \ar[d] \\
*+[r]{\;\>\;\begin{subarray}{l}\ts [\Spec R[x]/\bG_m(R)] \\[2pt] \ts =[\bA^1(R)/\bG_m(R)]\end{subarray}} \ar@{..>}[ur] \ar[r]  & *+[l]{\Spec\C.} }	
\end{equation*}
Here $(\m,x)$ is the (maximal) ideal in $R[x]$ generated by $\m\subset R$ and $x$, and so is a (closed) point in $\Spec R[x]$, which is invariant under $\bG_m(R)$. Halpern-Leistner \cite{Halp1,Halp2} writes $\Th$ for the stack $[\bA^1/\bG_m]$, and much of his theory of $\Th$-stratifications (see \S\ref{co333}) concerns morphisms $\Th\ra Y$.
\smallskip

\noindent{\bf(c)} Let $\pi$ be a {\it uniformizer\/} for $R$, that is a generator of the maximal ideal $\m$, so that $\m=(\pi)$. As in \cite[Def.~3.37]{AHLH}, we say that $Y$ is {\it S-complete with respect to\/} $R$ if in any commutative square of morphisms of $\C$-stacks `$\ra$', the dotted arrow `$\dashra$' can be uniquely filled in:
\begin{equation*}
\xymatrix@C=130pt@R=15pt{*+[r]{[((\Spec R[x,y]/(xy-\pi))\sm (\m,x,y))/\bG_m(R)]} \ar[r] \ar[d] & *+[l]{Y} \ar[d] \\
*+[r]{[(\Spec R[x,y]/(xy-\pi))/\bG_m(R)]} \ar@{..>}[ur] \ar[r]  & *+[l]{\Spec\C.} }	
\end{equation*}
\end{dfn}

As $(\Spec R)\sm\m$ in \eq{co3eq9} has only one point, it is easy to see that (a) is closed under finite unions of substacks:

\begin{lem}
\label{co3lem1}
Suppose $Y$ is an Artin $\C$-stack and\/ $Y_1,\ldots,Y_n\subseteq Y$ are substacks, and\/ $R$ is a DVR over\/ $\C$. If\/ $Y_1,\ldots,Y_n$ satisfy the valuative criterion for universal closedness with respect to\/ $R$ then so does\/~$Y_1\cup\cdots\cup Y_n$.
\end{lem}

Using Definition \ref{co3def8}, the next result follows from Alper--Halpern-Leistner--Heinloth~\cite[Th.s 4.1, 5.2 \& 5.4, Prop.s 3.45, 3.47, \& Rem.s 3.48 \& 5.5]{AHLH}:

\begin{thm}
\label{co3thm4}
Suppose $Y$ is a finite type Artin $\C$-stack with affine diagonal. Then $Y$ admits a proper good moduli space $\phi:Y\ra Z$ if and only if:
\begin{itemize}
\setlength{\itemsep}{0pt}
\setlength{\parsep}{0pt}
\item[{\bf(i)}] $Y$ satisfies the valuative criterion for universal closedness with respect to any discrete valuation ring (DVR)\/ $R$ essentially of finite type over\/~$\C$.
\item[{\bf(ii)}] $Y$ is $\Th$-reductive w.r.t.\ any DVR $R$ essentially of finite type over\/~$\C$.
\item[{\bf(iii)}] $Y$ is S-complete w.r.t.\ any DVR $R$ essentially of finite type over\/~$\C$.
\end{itemize}	
\end{thm}

\subsubsection{Abelian categories of compact type}
\label{co332}

We define a class of $\C$-linear abelian categories, which we call {\it of compact type\/} (a new definition not in \cite{AHLH}) for which results of Alper--Halpern-Leistner--Heinloth \cite{AHLH}, building on foundational work of Artin--Zhang \cite{ArZh}, show that the conditions of Theorem \ref{co3thm4} hold automatically for the moduli stack $\M$ of objects in $\A$, except that $\M$ is only locally of finite type, rather than of finite type.

\begin{dfn}
\label{co3def9}
Let $\tiA$ be a $\C$-linear abelian category. An object $E$ of $\tiA$ is called {\it compact\/} (or {\it finitely presentable\/}) if the functor $\Hom(E,-):\tiA\ra\mathop{\rm Sets}$ commutes with filtered colimits, and {\it noetherian\/} if every ascending chain of subobjects of $E$ terminates. 

The category $\tiA$ is called {\it cocomplete} if arbitrary small colimits exist in $\tiA$, and {\it locally noetherian\/} if it has a set of noetherian generators.

If $\tiA$ is locally noetherian then compact objects and noetherian objects in $\tiA$ coincide \cite[Prop.~B1.3]{ArZh}, and the full subcategory $\tiA{}^{\rm co}$ of compact objects in $\tiA$ is a noetherian $\C$-linear abelian category.

As in Artin--Zhang \cite{ArZh}, summarized in Alper et al.\ \cite[\S 7.1]{AHLH}, if $\tiA$ is locally noetherian then we can use $\tiA$ to define a moduli functor for objects in $\tiA{}^{\rm co}$, and so define the moduli stack $\M_{\tiA{}^{\rm co}}$ of objects in $\tiA{}^{\rm co}$, using only~$\tiA$.

Now let $\A$ be a $\C$-linear abelian category. We say that $\A$ is {\it of compact type\/} if there exists a cocomplete, locally noetherian $\C$-linear abelian category $\ti A$ and an equivalence of categories $\io:\A\ra\tiA{}^{\rm co}$, and the moduli stack $\M=\M_{\tiA{}^{\rm co}}$ of objects in $\A$ is an Artin $\C$-stack locally of finite type.
\end{dfn}

The next theorem follows from~\cite[Lem.s 7.15, 7.16, 7.17, \& 7.19]{AHLH}.

\begin{thm}
\label{co3thm5}
Let\/ $\A$ be a $\C$-linear abelian category of compact type, and\/ $\M$ be the moduli stack of objects in $\A$. Then $\M$ is locally of finite type by assumption, has affine diagonal, and satisfies Theorem\/ {\rm\ref{co3thm4}(i)--(iii)}.	
\end{thm}

Having affine diagonal and Theorem \ref{co3thm4}(i)--(iii) are all properties inherited by closed substacks. Thus as in \cite[Th.~7.21]{AHLH}, Theorems \ref{co3thm4} and \ref{co3thm5} imply:

\begin{cor}
\label{co3cor2}
Let\/ $\A$ be a $\C$-linear abelian category of compact type, and\/ $\M$ be the moduli stack of objects in $\A$. If\/ $Y\subset\M$ is a closed substack of finite type then $Y$ admits a proper good moduli space.
\end{cor}

\begin{rem}
\label{co3rem1}
In the situation of Definition \ref{co3def6}, our goal, achieved in \S\ref{co336}, is to give useful criteria for when $\M_\al^\ss(\tau)\subseteq\M_\al^\pl$ has a proper good moduli space, or equivalently $(\Pi_\al^\pl)^{-1}(\M_\al^\ss(\tau))\subseteq\M_\al$ has a proper good moduli space, as then Proposition \ref{co3prop6} says that if $\M_\al^\rst(\tau)=\M_\al^\ss(\tau)$ then $\M_\al^\ss(\tau)$ is a proper algebraic space. So we need to verify that $(\Pi_\al^\pl)^{-1}(\M_\al^\ss(\tau))$ is of finite type (which will be part of our criteria) and satisfies Theorem~\ref{co3thm4}(i)--(iii).

Theorem \ref{co3thm5} says that $\M$, and hence $\M_\al$, satisfies Theorem \ref{co3thm4}(i)--(iii). Our strategy, following Alper--Halpern-Leistner--Heinloth \cite[\S 7]{AHLH}, will be to prove that under some conditions on $(\tau,T,\le)$, the properties Theorem \ref{co3thm4}(i)--(iii) for $\M_\al$ imply the same properties for $(\Pi_\al^\pl)^{-1}(\M_\al^\ss(\tau))$.

This comes in two parts. Sections \ref{co333}--\ref{co334} show that if $(\tau,T,\le)$ induces a `(pseudo)-$\Th$-stratification' of $\M_\al$ with semistable locus $(\Pi_\al^\pl)^{-1}(\M_\al^\ss(\tau))$, then Theorem \ref{co3thm4}(i) for $\M_\al$ implies the same for~$(\Pi_\al^\pl)^{-1}(\M_\al^\ss(\tau))$. Then \S\ref{co335} explains that if $(\tau,T,\le)$ is an `additive' weak stability condition then Theorem \ref{co3thm4}(ii)--(iii) for $\M_\al$ imply the same for $(\Pi_\al^\pl)^{-1}(\M_\al^\ss(\tau))$. 
\end{rem}

\subsubsection{$\Th$-stratifications}
\label{co333}

Next we discuss the theory of `$\Th$-stratifications' of Artin stacks, due to Halpern-Leistner \cite{Halp1,Halp2}. This gives an analogue of stability conditions and Harder--Narasimhan filtrations in \S\ref{co31} on the level of moduli stacks, so that for example the $\tau$-semistable substack $\M_\al^\ss(\tau)\subseteq\M_\al^\pl$ might be defined directly using a geometric structure on $\M_\al^\pl$ or~$\M_\al$.

Our next definition follows \cite[Def.~6.1]{AHLH}, which is a minor modification of Halpern-Leistner \cite[Def.s 2.1 \& 2.2]{Halp1} and~\cite[Def.~2.20]{Halp2}.

\begin{dfn}
\label{co3def10}
Write $\Th$ for the quotient Artin $\C$-stack $[\bA^1/\bG_m]$. Let $Y$ be an Artin $\C$-stack, locally of finite type, with affine diagonal. Then by \cite[Prop.~1.2]{Halp1} or \cite[Ex.~1.2.2]{HaPr}, the mapping stack $\Map(\Th,Y)$ is also an Artin $\C$-stack, locally of finite type, with affine diagonal. There are evaluation maps $\ev_0,\ev_1:\Map(\Th,Y)\ra Y$ mapping $f:\Th\ra Y$ to $f(0),f(1)$ for $0,1\in\bA^1$. That is, $\ev_0$ 
evaluates $f$ at the special point $[\{0\}/\bG_m]$ and $\ev_1$ evaluates $f$ at the generic point $[(\bA^1\sm\{0\})/\bG_m]\cong\Spec\C$ in $\Th$. Then:
\begin{itemize}
\setlength{\itemsep}{0pt}
\setlength{\parsep}{0pt}
\item[(a)] A $\Th$-{\it stratum\/} in $Y$ consists of a union of connected components $\cS\subseteq\Map(\Th,Y)$ such that $\ev_1:\cS\ra Y$ is a closed immersion, identifying $\cS$ with a closed $\C$-substack $\ev_1(\cS)\subseteq Y$.
\item[(b)] A $\Th$-{\it stratification\/} of $Y$ indexed by a totally ordered set $(\Ga,\le)$ with a minimal element $0\in\Ga$ is a cover of $Y$ by open substacks $Y_{\le\ga}$ for $\ga\in\Ga$ such that $Y_{\le\ga_1}\subseteq Y_{\le\ga_2}$ for $\ga_1\le\ga_2$ in $\Ga$, and a particular choice of $\Th$-stratum $\cS_\ga\subseteq\Map(\Th,Y_{\le\ga})$ in $Y_{\le\ga}$ for each $\ga\in\Ga$ such that $Y_{\le\ga}\sm\ev_1(\cS_\ga)=\bigcup_{\text{$\ga'<\ga$ in $\Ga$}}Y_{\le\ga'}$. We require that for all points $y$ in $Y$, the subset $\bigl\{\ga\in\Ga:y\in Y_{\le\ga}\bigr\}\subseteq\Ga$ has a minimal element.

The {\it semistable locus\/} of the $\Th$-stratification is $Y^\ss=Y_{\le 0}$, an open substack of $Y$ isomorphic to the minimal $\Th$-stratum~$\cS_0$.
\item[(c)] A $\Th$-stratification of $Y$ is called {\it well ordered\/} if for all $y\in Y$, the subset
\e
\bigl\{\ga\in\Ga:\ov{\{y\}}\cap\bigl(Y_{\le\ga}\sm\ts\bigcup_{\text{$\ga'<\ga$ in $\Ga$}}Y_{\le\ga'}\bigr)\ne\es\bigr\}\subseteq\Ga
\label{co3eq10}
\e
is well ordered.
\item[(d)] A $\Th$-stratification of $Y$ satisfies the {\it descending chain condition\/} if for all $y\in Y$, any descending chain $\ga_0\ge\ga_1\ge\ga_2\ge\cdots$ in \eq{co3eq10} eventually stabilizes. This holds automatically if the $\Th$-stratification is well ordered.
\end{itemize}	
\end{dfn}

The next theorem follows from \cite[Cor.~6.12]{AHLH}. They assume $(Y_{\le\ga})_{\ga\in\Ga}$ is well ordered, but they only actually use the descending chain condition.

\begin{thm}
\label{co3thm6}
Suppose $Y$ is an Artin $\C$-stack locally of finite type with affine diagonal, and\/ $(\Ga,\le),(Y_{\le\ga})_{\ga\in\Ga}$ is a\/ $\Th$-stratification of\/ $Y$ satisfying the descending chain condition. If\/  $Y$ satisfies the valuative criterion for universal closedness with respect to any DVR\/ $R$ essentially of finite type over\/ $\C,$ as in Definition\/ {\rm\ref{co3def8}(a)} and Theorem\/ {\rm\ref{co3thm4}(i),} then so does $Y_{\le\ga}$ for all\/ $\ga\in\Ga$.
\end{thm}

Now as in Halpern-Leistner \cite{Halp1,Halp2} and \cite[\S 7]{AHLH}, if we take $Y$ to be a moduli stack $\M_\al^\pl$ or $\M_\al$ of objects of class $\al$ in some $\C$-linear abelian category $\A$, there is a systematic method for constructing a well ordered $\Th$-stratification of $Y$ such that the semistable locus $Y^\ss$ is the semistable substack $\M_\al^\ss(\tau)\subseteq\M_\al^\pl$ or $(\Pi_\al^\pl)^{-1}(\M_\al^\ss(\tau))\subseteq\M_\al$ of a suitable stability condition $(\tau,T,\le)$ on $\A$, and the other $\Th$-strata $Y_\ga$ roughly parametrize $E\in\A$ with $\tau$-Harder--Narasimhan filtrations $0=E_0\subset E_1\subset\cdots\subset E_n=E$ with prescribed $\tau(E_i/E_{i-1})\in T$.

To explain why, first note that by \cite[Cor.~7.12]{AHLH} we have:

\begin{prop}
\label{co3prop7}
Let\/ $\A$ be a $\C$-linear abelian category of compact type, and\/ $\M$ be the moduli stack of objects in $\A$. Then $\C$-points $f\in\Map(\Th,\M)$ are in 1--1 correspondence with isomorphism classes of objects $E\in\A$ with a filtration by subobjects $\cdots\subseteq\bar E_{k-1}\subseteq\bar E_k\subseteq\cdots\subseteq E$ for $k\in\Z,$ such that 
\begin{itemize}
\setlength{\itemsep}{0pt}
\setlength{\parsep}{0pt}
\item[{\bf(i)}] $\bar E_k=0$ for $k\ll 0$ and\/ $\bar E_k=E$ for $k\gg 0;$
\item[{\bf(ii)}] $\ev_1(f)\cong E;$ and 
\item[{\bf(iii)}] $\ev_0(f)\cong\bigop_{k\in\Z}(\bar E_k/\bar E_{k-1}),$ with\/ $\bG_m$-action $\la\mapsto\sum_{k\in\Z}\la^k\,\id_{\bar E_k/\bar E_{k-1}}$.
\end{itemize}
\end{prop}

\begin{rem}
\label{co3rem2}
There is an {\it important difference\/} between $\tau$-Harder--Narasimhan filtrations $0=E_0\subset E_1\subset\cdots\subset E_n=E$ in Theorem \ref{co3thm1} and the filtrations $\cdots\subseteq \bar E_{k-1}\subseteq\bar E_k\subseteq\cdots\subseteq E$ in Proposition \ref{co3prop7}.

Given a filtration $0=E_0\subset E_1\subset\cdots\subset E_n=E$, to construct a filtration $\cdots\subseteq \bar E_{k-1}\subseteq\bar E_k\subseteq\cdots\subseteq E$ we need to choose the extra data of integers $k_1<k_2<\cdots<k_n$, and then we define $\bar E_k=E_i$ if $k_i\le k<k_{i+1}$, with $k_0=-\iy$, $k_{n+1}=\iy$. Then $\bar E_k/\bar E_{k-1}=E_i/E_{i-1}$ if $k=k_i$ for some $i=1,\ldots,n$, and $\bar E_k/\bar E_{k-1}=0$ otherwise.

In Halpern-Leistner \cite[\S 1.3.3]{Halp1}, a point $f\in\Map(\Th,\M)$ is considered equivalent to $f\ci\Phi^a$, where $\Phi^a:\Th\ra\Th$ is induced by the map $\bA^1\ra\bA^1$ mapping $z\mapsto z^a$ for positive integers $a$. This has the effect of replacing $k_1<k_2<\cdots<k_n$ by $ak_1<ak_2<\cdots<ak_n$. So as in \cite[\S 1.1.1]{Halp2}, the $k_i$ may be replaced by rational weights $w_1<w_2<\cdots<w_n$ in~$\Q$.

Suppose $(\tau,T,\le)$ is a stability condition on $\A$ for which Theorem \ref{co3thm1} holds. One of the central ideas of \cite{Halp1,Halp2} is to try to define a $\Th$-stratification of $\M$, such that the $\Th$-strata are parametrized by $n$-tuples $(\be_1,\ldots,\be_n)$ in $C(\A)$ for $n\ge 1$ and parametrize objects $E$ whose $\tau$-Harder--Narasimhan filtrations $0=E_0\subset E_1\subset\cdots\subset E_n=E$ have $\tau(\lb E_i/E_{i-1}\rb)=\be_i$ for $i=1,\ldots,n$. Thus, $\Th$-stratifications generalize Harder--Narasimhan stratification of~$\M$.

However, there is a difficulty in doing this, which is that to promote a Harder--Narasimhan stratification to  a $\Th$-stratification, one needs a systematic way of assigning rational weights $w_1<w_2<\cdots<w_n$ to a filtration $0=E_0\subset E_1\subset\cdots\subset E_n=E$. This is usually done in \cite{Halp1,Halp2} by solving a minimization problem. But ensuring the solutions $w_i$ lie in $\Q$ rather than $\R$ leads to extra work and strong restrictions. See for example \cite[\S 5.4]{Halp1} for discussion of cases in which the $w_i$ have no natural values, or the natural values lie in~$\R$.

In \S\ref{co334} we propose a modification to the notion of $\Th$-stratification which the author hopes gets round the problem of choosing rational weights~$w_i$.	
\end{rem}

\subsubsection{Pseudo-$\Th$-stratifications}
\label{co334}

The material of this section may be new, but is only a small modification to the big machine of \cite{Halp1,Halp2,AHLH}, so does not deserve much credit.

\begin{dfn}
\label{co3def11}
Let $Y$ be an Artin $\C$-stack, locally of finite type, with affine diagonal, and write $\Th=[\bA^1/\bG_m]$. Then as in Definition \ref{co3def10}, $\Map(\Th,Y)$ is an Artin $\C$-stack, locally of finite type, with affine diagonal, with evaluation maps~$\ev_0,\ev_1:\Map(\Th,Y)\ra Y$. 

Let $(\Ga,\preceq)$ be a partially ordered set which has a minimal element $0\in\Ga$, such that any $\ga_1,\ga_2\in\Ga$ have a unique {\it greatest lower bound\/} $\glb(\ga_1,\ga_2)$ in $\Ga$.

A {\it pseudo-$\Th$-stratification\/} of $Y$ indexed by $(\Ga,\preceq)$ is a cover $(Y_{\preceq\ga})_{\ga\in\Ga}$ of $Y$ by open substacks $Y_{\preceq\ga}$ for $\ga\in\Ga$ satisfying: 
\begin{itemize}
\setlength{\itemsep}{0pt}
\setlength{\parsep}{0pt}
\item[(i)] If $\ga_1,\ga_2\in\Ga$ then $Y_{\preceq\ga_1}\cap Y_{\preceq\ga_2}=Y_{\preceq\glb(\ga_1,\ga_2)}$. In particular, this implies that if $\ga_1\preceq\ga_2$ then $Y_{\preceq\ga_1}\subseteq Y_{\preceq\ga_2}$.
\item[(ii)] For each $\ga\in\Ga$ there exists a $\Th$-stratum $\cS_\ga\subseteq\Map(\Th,Y_{\preceq\ga})$ in $Y_{\preceq\ga}$ for each $\ga\in\Ga$ such that $Y_{\preceq\ga}\sm\ev_1(\cS_\ga)=\bigcup_{\text{$\ga'\prec\ga$ in $\Ga$}}Y_{\preceq\ga'}$. 

Note that in contrast to Definition \ref{co3def10}(b) we do not choose a particular $\Th$-stratum $\cS_\ga$ for $\ga\in\Ga$ above, and $\cS_\ga$ may be nonunique.
\item[(iii)] For all $y\in Y$, the subset $\bigl\{\ga\in\Ga:y\in Y_{\preceq\ga}\bigr\}\subseteq\Ga$ has a minimal element. 
\end{itemize}

The {\it semistable locus\/} of the pseudo-$\Th$-stratification is $Y^\ss=Y_{\preceq 0}$, an open substack of $Y$. A pseudo-$\Th$-stratification of $Y$ satisfies the {\it descending chain condition\/} if for all $y\in Y$, any descending chain $\ga_0\succeq\ga_1\succeq\ga_2\succeq\cdots$ in \eq{co3eq10} eventually stabilizes.

Clearly, any $\Th$-stratification of $Y$ in Definition \ref{co3def10} induces a pseudo-$\Th$-stratification of $Y$, by forgetting the particular choices of $\Th$-strata $\cS_\ga$, and this forgetful map preserves the descending chain conditions. Note that if $(\Ga,\le)$ is a total order then $\glb(\ga_1,\ga_2)=\min(\ga_1,\ga_2)$ exists automatically.
\end{dfn}

\begin{rem}
\label{co3rem3}
{\bf(a)} The motivating idea behind Theorem \ref{co3def11} is that in Remark \ref{co3rem2}, the choices of rational weights $w_1<\cdots<w_n$ are encoded in the choices of $\Th$-strata $\cS_\ga$, so we forget the $\cS_\ga$ and remember only their images $\ev_1(\cS_\ga)=Y_{\le\ga}\sm\bigcup_{\text{$\ga'<\ga$ in $\Ga$}}Y_{\le\ga'}$. Also the $w_i$ are needed to define the total order $\le$ on $\Ga$, so we have to replace $\le$ by a partial order~$\preceq$.
\smallskip

\noindent{\bf(b)} We can modify any pseudo-$\Th$-stratification $(Y_{\preceq\ga})_{\ga\in\Ga}$ of $Y$ over a partially ordered set $(\Ga,\preceq)$ to a total order $(\Ga,\le)$ as follows: choose an arbitrary total order $\le$ on $\Ga$ compatible with $\preceq$, and replace $Y_{\preceq\ga}$ for $\ga\in\Ga$ by 
\begin{equation*}
\bar Y_{\le\ga}=\bigcup_{\ga'\in\Ga:\ga'\le\ga}Y_{\preceq\ga'}.
\end{equation*}
However, for stratifications into Harder--Narasimhan types, as in Shatz \cite{Shat} and Nitsure \cite{Nits} for instance, using a partial order seems more natural.
\end{rem}

The next example translates work of Shatz \cite{Shat}, Nitsure \cite{Nits}, Hoskins \cite{Hosk}, and Halpern-Leistner \cite{Halp1,Halp2} into our notation. See Definition \ref{co6def9} below for an example of pseudo-$\Th$-stratifications of quiver moduli stacks.

\begin{ex}
\label{co3ex6}
Let $X$ be a smooth projective $\C$-scheme of dimension $m$ and $(\tau,G,\le)$ be Gieseker stability on $\A=\coh(X)$ with respect to an ample line bundle $\O_X(1)$ on $X$, as in Example \ref{co3ex2}. Define the set $\HNT(X)$ of {\it Harder--Narasimhan types\/} for $X$ to consist of all $n$-tuples $(p_1,\ldots,p_n)$ for $n\ge 1$, where:
\begin{itemize}
\setlength{\itemsep}{0pt}
\setlength{\parsep}{0pt}
\item[(i)] $p_i(t)\in\Q[t]$ is a nonzero rational polynomial, of degree $d_i$ with $0\le d_i\le m$, such that $p_i$ maps $\Z\ra\Z$ (such $p_i$ are called {\it numerical polynomials\/}). Write $p_i(i)=r_it^{d_i}+O(t^{d_i-1})$, so that $r_i$ is the leading coefficient of $p_i$. We require that $r_i>0$. Note that $\frac{1}{r_i}p_i(t)$ lies in $G$ in Example \ref{co3ex2}.
\item[(ii)] $\frac{1}{r_1}p_1>\frac{1}{r_2}p_2>\cdots>\frac{1}{r_n}p_n$ in the order on $G$. This implies that $d_1\ab\le\ab d_2\ab\le\ab\cdots\ab\le\ab d_n$. We call $(p_1,\ldots,p _n)$ {\it pure\/} if $d_1=\cdots=d_n$.
\end{itemize} 
If $\al\in C(\coh(X))$, with Hilbert polynomial $P_\al$ as in Example \ref{co3ex2}, we write $\HNT_\al(X)$ for the subset of $(p_1,\ldots,p_n)$ in $\HNT(X)$ with~$p_1+\cdots+p_n=P_\al$. 

If $\bs p=(p_1,\ldots,p_n)\in\HNT_\al(X)$ and $k,l\in\R$, define $x_0,\ldots,x_n\in\R^2$ by $x_j=(\sum_{i=1}^jp_j(k),\sum_{i=1}^jp_j(l))$, and define the {\it Shatz polygon\/} $\mathop{\rm Sh}(\bs p,k,l)$ to be the union of the line segments in $\R^2$ joining $x_{i-1}$ to $x_i$ for $i=1,\ldots,n$. This is a piecewise-linear curve with end-points $x_0=(0,0)$ and $x_n=(P_\al(k),P_\al(l))$. If $0\ll k\ll l$ it is the graph $\Ga_{\mathop{\rm sh}(\bs p,k,l)}$ of a continuous, piecewise linear function $\mathop{\rm sh}(\bs p,k,l):[0,P_\al(k)]\ra[0,P_\al(l)]$, and using the formulation \eq{co3eq1} of Gieseker stability, condition (ii) implies that it is a convex polygon.

Following Hoskins \cite[\S 4.5]{Hosk}, which generalizes Shatz \cite[\S 3]{Shat} and Nitsure \cite[\S 2]{Nits}, define a partial order $\preceq$ on $\HNT_\al(X)$ by $\bs p\preceq\bs p'$ if $\mathop{\rm Sh}(\bs p,k,l)$ lies below $\mathop{\rm Sh}(\bs p',k,l)$ in $[0,\iy)^2$ for $0\ll k\ll l$, that is, if $\mathop{\rm sh}(\bs p,k,l)\le\mathop{\rm sh}(\bs p',k,l)$ as functions $[0,P_\al(k)]\ra[0,P_\al(l)]$. This is independent of $k,l$ provided $0\ll k\ll l$. The minimal element in $\HNT_\al(X)$ under $\preceq$ is $\bs p=(P_\al)$, so in Definition \ref{co3def11} we write~$0=(P_\al)$.

Since $\coh(X)$ is $\tau$-artinian, for each $E\in\A$ with $\lb E\rb=\al$, by Theorem \ref{co3thm1} we have a $\tau$-Harder--Narasimhan filtration $0=E_0\subsetneq E_1\subsetneq\cdots\subsetneq E_n=E$ for $n\ge 0$, such that\/ $F_i=E_i/E_{i-1}$ is $\tau$-semistable for $i=1,\ldots,n,$ and $\tau([F_1])>\tau([F_2])>\cdots>\tau([F_n])$ in $G$. Define the {\it Harder--Narasimhan type\/} of $E$ to be
\begin{equation*}
\HNT(E)=\bigl(P_{F_1},P_{F_2},\ldots,P_{F_n}\bigr),
\end{equation*}
where $P_{F_i}$ is the Hilbert polynomial of $F_i$. Then $\HNT(E)\in\HNT_\al(X)$.

As in Shatz \cite{Shat}, Nitsure \cite[Th.~5]{Nits} and Hoskins \cite[Th.~4.15]{Hosk}, there is a locally constructible function $\HNT_\al:\M_\al\ra\HNT_\al(X)$ which maps a $\C$-point $[E]$ to $\HNT(E)$. Furthermore, $\HNT_\al$ is upper semicontinuous in the partial order $\preceq$ on $\HNT_\al(X)$. Hence for each $\bs p\in\HNT_\al(X),$ we have open substacks
\e
\M_{\al,\preceq\bs p}=\bigl\{[E]\in\M_\al:\HNT_\al([E])\preceq\bs p\bigr\}.
\label{co3eq11}
\e

Nitsure \cite[Th.s 5 \& 8]{Nits} proves that for $\bs p\in\HNT_\al(X)$ there is a natural locally closed substack $\M_{\al,\bs p}\subseteq\M_{\al,\preceq\bs p}\subseteq\M_\al$, such that in open substacks
\e
\M_{\al,\preceq\bs p}\sm\M_{\al,\bs p}=\ts\bigcup_{\begin{subarray}{l} \bs p'\in\HNT_\al(X):\bs p'\prec\bs p\end{subarray}}\M_{\al,\preceq\bs p'}.
\label{co3eq12}
\e
Furthermore $\M_{\al,\bs p}$ is a moduli stack for flat families of sheaves with $\tau$-Harder--Narasimhan filtrations of type $\bs p$. Nitsure restricts to pure sheaves, but as in Hoskins \cite[Th.~4.15 \& Cor.~4.16]{Hosk} the results also hold for general sheaves.

Writing $\cU_\al\ra X\t\M_\al$ for the universal sheaf and $\io: \M_{\al,\bs p}\hookra\M_\al$ for the inclusion, the pullback $(\id_X\t\io)^*(\cU_\al)\ra X\t\M_{\al,\bs p}$ has a canonical family $\tau$-Harder--Narasimhan filtration $0=\cE_0\subsetneq\cE_1\subsetneq\cdots \subsetneq\cE_n=(\id_X\t\io)^*(\cU_\al)$ of type $\bs p=(p_1,\ldots,p_n)$. Since $\coh(X)$ is of compact type, as in Proposition \ref{co3prop7} this filtration induces a morphism $\M_{\al,\bs p}\ra\Map(\Th,\M_{\al,\preceq\bs p})$, such that composing with $\ev_1:\Map(\Th,\M_{\al,\preceq\bs p})\ra \M_{\al,\preceq\bs p}$ gives the inclusion $\M_{\al,\bs p}\hookra \M_{\al,\preceq\bs p}$. 

Then Proposition \ref{co3prop7} and the moduli stack property of $\M_{\al,\bs p}$ imply that $\M_{\al,\bs p}\ra\Map(\Th,\M_{\al,\preceq\bs p})$ is an isomorphism with a union of connected components of $\Map(\Th,\M_{\al,\preceq\bs p})$, and so defines a $\Th$-stratum of $\M_{\al,\preceq\bs p}$, as in Definition \ref{co3def10}(a). Thus \eq{co3eq12} shows Definition \ref{co3def11}(ii) holds for $(\M_{\al,\preceq\bs p})_{\bs p\in\HNT_\al(X)}$, which is therefore a pseudo-$\Th$-stratification of~$\M_\al$. 
\end{ex}

The analogue of Theorem \ref{co3thm6} holds for pseudo-$\Th$-stratifications:

\begin{thm}
\label{co3thm7}
Suppose $Y$ is an Artin $\C$-stack locally of finite type with affine diagonal, and\/ $(\Ga,\preceq),(Y_{\preceq\ga})_{\ga\in\Ga}$ is a\/ pseudo-$\Th$-stratification of\/ $Y$ satisfying the descending chain condition. If\/  $Y$ satisfies the valuative criterion for universal closedness with respect to any DVR\/ $R$ essentially of finite type over\/ $\C,$ as in Definition\/ {\rm\ref{co3def8}(a)} and Theorem\/ {\rm\ref{co3thm4}(i),} then so does $Y_{\preceq\ga}$ for all\/~$\ga\in\Ga$.
\end{thm}

\begin{proof}
We first extend the proof of \cite[Th.~6.5]{AHLH} to pseudo-$\Th$-stratifications. In our notation, their proof starts with a DVR $R$ essentially of finite type over $\C$ and a morphism $\xi:\Spec R\ra Y$ such that $\xi(K)\in Y_{\preceq\ga}\sm\bigcup_{\text{$\ga'\prec \ga$ in $\Ga$}}Y_{\preceq\ga'}$ and $\xi(\ka)\in Y_{\preceq\ga_0}\sm\bigcup_{\text{$\ga'\prec \ga_0$ in $\Ga$}}Y_{\preceq\ga'}$, where $K$ is the generic point and $\ka$ the special point in $\Spec R$. Here given $R,\xi$, Definition \ref{co3def11}(iii) implies that $\ga,\ga_0$ exist and are unique. We will show that $\ga\preceq\ga_0$ in~$\Ga$. 

Since $Y_{\preceq\ga_0}$ is open with $\xi(\ka)\in Y_{\preceq\ga_0}$, and $\xi(\ka)$ lies in the closure of $\xi(K)$ as $\ka$ lies in the closure of $K$, we see that $\xi(K)\in Y_{\preceq\ga_0}$. As $\xi(K)\in Y_{\preceq\ga}$ this implies that $\xi(K)\in Y_{\preceq\ga}\cap Y_{\preceq\ga_0}=Y_{\preceq\glb(\ga,\ga_0)}$ by Definition \ref{co3def11}(i). But $\glb(\ga,\ga_0)\preceq\ga$, so if $\glb(\ga,\ga_0)\ne\ga$ then $\xi(K)\in Y_{\preceq\glb(\ga,\ga_0)}$ contradicts $\xi(K)\notin\bigcup_{\text{$\ga'\prec \ga$ in $\Ga$}}Y_{\preceq\ga'}$. Hence $\glb(\ga,\ga_0)=\ga$, so~$\ga\preceq\ga_0$.

The aim of the proof of \cite[Th.~6.5]{AHLH} is to modify $R,\xi$ to $R',\xi'$ using \cite[Th.~6.3]{AHLH} such that $\xi'\vert_{K'}:K'\ra Y$ factors through $\xi\vert_K:K\ra Y$, giving $\xi'(K')\in Y_{\preceq\ga}\sm\bigcup_{\text{$\ga'\prec \ga$ in $\Ga$}}Y_{\preceq\ga'}$, and also $\xi'(\ka')\in Y_{\preceq\ga}\sm\bigcup_{\text{$\ga'\prec \ga$ in $\Ga$}}Y_{\preceq\ga'}$, that is, $R',\xi'$ have $\ga_0'=\ga$. So to exclude the case $\ga=\ga_0$, they begin by supposing that $\ga<\ga_0$ in the total order $(\Ga,\le)$. In our case, as $\ga\preceq\ga_0$, if $\ga\ne\ga_0$ then $\ga\prec\ga_0$. So we can replace $\ga<\ga_0$ in \cite[Proof of Th.~6.5]{AHLH} by $\ga\prec\ga_0$ in the partial order $(\Ga,\preceq)$. The same proof then shows that we can replace $R,\xi,\ga_0$ by $R',\xi',\ga_0'$ such that $\ga_0'\prec\ga_0$, since $\xi'(\ka')\in Y_{\preceq\ga_0}\sm\ev_1(\cS_{\ga_0})$, and as $\xi'(K')\in Y_{\preceq\ga}\sm\bigcup_{\text{$\ga'\prec \ga$ in $\Ga$}}Y_{\preceq\ga'}$ the proof above shows that $\ga\preceq\ga_0'$. Thus $\ga\preceq\ga_0'\prec\ga_0$. If $\ga\ne\ga_0'$ we can iterate to $\ga\preceq\ga_0''\prec\ga_0'\prec\ga_0$, and so on. The rest of the proof of \cite[Th.~6.5]{AHLH} extends in the obvious way.

Then \cite[Cor.~6.12]{AHLH}, which uses \cite[Th.~6.5]{AHLH}, also extends immediately to pseudo-$\Th$-stratifications, and the theorem follows.	
\end{proof}

Here is a useful criterion for a (pseudo)-$\Th$-stratification to satisfy the descending chain condition:

\begin{prop}
\label{co3prop8}
Let\/ $Y$ be an Artin $\C$-stack, locally of finite type, with affine diagonal, and with a pseudo-$\Th$-stratification $(\Ga,\preceq),(Y_{\preceq\ga})_{\ga\in\Ga}$. If\/ $Y_{\preceq\ga}$ is of finite type for all\/ $\ga\in\Ga$ then $(\Ga,\preceq),(Y_{\preceq\ga})_{\ga\in\Ga}$ satisfies the descending chain condition.

The analogue holds for $\Th$-stratifications.
\end{prop}

\begin{proof}
Let $y\in Y$, and write $\Ga'$ for the subset \eq{co3eq10} of $\Ga$. Suppose $\ga_0\succeq\ga_1\succeq\ga_2\succeq\cdots$ is a descending chain in $\Ga'$. Then $Y_{\preceq\ga_0}\supseteq Y_{\preceq\ga_1}\supseteq Y_{\preceq\ga_2}\supseteq\cdots$ is a descending chain of open substacks in the finite type Artin $\C$-stack $Y_{\preceq\ga_0}$.

Choose an atlas $\pi:U_{\preceq\ga_0}\ra Y_{\preceq\ga_0}$ with $U_{\preceq\ga_0}$ a finite type $\C$-scheme, and let $U_{\preceq\ga_n}=\pi^{-1}(Y_{\preceq\ga_n})$, so that $U_{\preceq\ga_0}\supseteq U_{\preceq\ga_1}\supseteq U_{\preceq\ga_2}\supseteq\cdots$ is a descending chain of open subschemes in the finite type $\C$-scheme~$U_{\preceq\ga_0}$.

Suppose for a contradiction that $U_{\preceq\ga_0}\supseteq U_{\preceq\ga_1}\supseteq U_{\preceq\ga_2}\supseteq\cdots$ does not stabilize. Then by passing to a subsequence of $\ga_i$ we can suppose that $U_{\preceq\ga_k}\ne U_{\preceq\ga_{k+1}}$ for all $k\ge 0$. Write $d_k=\dim (U_{\preceq\ga_k}\sm U_{\preceq\ga_{k+1}})$, which makes sense as $U_{\preceq\ga_k}$ is of finite type. Then $(d_k)_{k\ge 0}$ is a sequence of nonnegative integers bounded above by $\dim U_{\preceq\ga_0}$, so $d=\limsup_{k\ra\iy}d_k$ is well defined. There is $K\ge 0$ such that $d_k\le d$ for $k\ge K$, so by replacing $\ga_k$ by $\ga_k'=\ga_{k+K}$ we can suppose $d_k\le d$ for all $k\ge 0$. As $d=\limsup_{k\ra\iy}d_k$ there exist $0\le i_0<i_1<i_2<\cdots$ with $d_{i_k}=d$ for $i=0,1,\ldots.$ Then replacing $\ga'_k$ by $\ga_k''=\ga'_{i_k}$ we can suppose that $d_k=d$ for all~$k=0,1,\ldots.$ 

Write $m=\dim U_{\ga_0}$, so that $m\ge d=\dim U_{\ga_0}\sm U_{\ga_1}$. By reverse induction on $n=m,m-1,\ldots,d$, we will construct the following data:
\begin{itemize}
\setlength{\itemsep}{0pt}
\setlength{\parsep}{0pt}
\item[(i)] elements $\de_n^i\in\Ga$ for $i=1,\ldots,M_n$ for $M_n\ge 0$ with $\de_n^i\preceq\ga_0$ for all $n,i$ and $\de_n^i\preceq\ga_k$ for all $k\ge 0$ and $n,i$ with $n>d$, so that $U_{\preceq\de_n^i}\subseteq U_{\preceq\ga_k}$ if $k=0$ or $n>d$, where $U_{\preceq\de_n^i}=\pi^{-1}(Y_{\preceq\de_n^i})$.
\item[(ii)] For $n=m,m-1,\ldots,d$ we require the closed reduced subscheme of~$U_{\ga_0}$
\e
\Bigl(U_{\ga_0}\sm\bigcup_{n<l\le m,\; 1\le i\le M_l}U_{\preceq\de_l^i}\Bigr)^\red
\label{co3eq13}
\e
to have dimension $\le n$, we write $M_n\ge 0$ for the number of irreducible components of \eq{co3eq13} of dimension exactly $n$, and we number them $U_n^i$ for $i=1,\ldots,M_n$. Note that this makes sense as $U_{\ga_0}$ is a finite type $\C$-scheme.
\item[(iii)] Write $u_n^i$ for the generic point of $U_n^i$, for $n=m,m-1,\ldots,d$ and $i=1,\ldots,M_n$.
Then $\pi(u_n^i)\in Y_{\preceq\ga_0}$, so by Definition \ref{co3def11}(iii), $\bigl\{\ga\in\Ga:\pi(u^i_n)\in Y_{\preceq\ga}\bigr\}\subseteq\Ga$ has a unique minimal element, which we write~$\de_n^i$.
\end{itemize}

For the inductive step, let $l=m,m-1,\ldots,d$ and suppose we have constructed $M_n,\de_n^i,U_n^i,u_n^i$ satisfying (i)--(iii) for all $m\le n<l$ and $i=1,\ldots,M_n$, and satisfying that the dimension of \eq{co3eq13} has dimension $\le l$ when $n=l$. For the first step $l=m$ this assumption is vacuous, noting that $\dim U_{\ga_0}=m$. Then (ii) for $n=l$ defines $M_l$ and $U_l^i$ for $i=1,\ldots,M_l$, and (iii) for $n=l$ defines $u_l^i$ and $\de_l^i$ for $i=1,\ldots,M_l$, and (i) defines~$U_{\preceq\de_l^i}$.

As $\pi(u^i_l)\in Y_{\preceq\ga_0}$ we see that $\ga_0$ lies in the subset $\bigl\{\ga\in\Ga:\pi(u^i_l)\in Y_{\preceq\ga}\bigr\}$ in (iii), so $\de_l^i\preceq\ga_0$, as we have to prove, and $Y_{\preceq\de_l^i}\subseteq Y_{\preceq\ga_0}$, and $U_{\preceq\de_l^i}\subseteq U_{\preceq\ga_0}$. If $l>d$ then $u^i_l\in U_{\preceq\ga_k}$ for all $k>0$, since $U_{\preceq\ga_k}$ is obtained from $U_{\preceq\ga_0}$ by deleting a subscheme of dimension $d$, which cannot contain $u^i_m$ as $\dim U^i_l=l>d$. Thus $\ga_k$ also lies in the subset with minimum $\de_l^i$, so $\de_m^i\preceq\ga_k$. This completes (i).

It remains to prove that \eq{co3eq13} for $n=l-1$ has dimension $\le l-1$. This holds as \eq{co3eq13} for $n=l$ has dimension $\le l$, and $M_l$ irreducible components $U^1_l,\ldots,U^{M_l}_l$ of dimension $l$. But each of these $M_l$ components has been deleted in \eq{co3eq13} for $n=l-1$ by the exclusion of $U_{\preceq\de_l^i}$ for $i=1,\ldots,M_l$, which contains the generic point $u^i_l$ of $U^i_l$. Hence \eq{co3eq13} for $n=l-1$ has no components of dimension $l$, and has dimension $\le l-1$. This completes the inductive step.

Now for $k\ge 0$ consider the closed reduced subscheme of $U_{\ga_k}$
\e
\Bigl(U_{\ga_k}\sm\bigcup_{d<l\le m,\; 1\le i\le M_l}U_{\preceq\de_l^i}\Bigr)^\red.
\label{co3eq14}
\e
Here $U_{\preceq\de_l^i}\preceq U_{\ga_k}$ as $\de_l^i\preceq\ga_k$ for $l>d$ by (i). It is contained in \eq{co3eq13} for $n=d$, and so \eq{co3eq14} has dimension $\le d$ by (ii). Write $N_k$ for the number of irreducible components of \eq{co3eq14} of dimension $d$. Then $N_0=M_d$. Comparing \eq{co3eq14} for $k,k+1$ we see that \eq{co3eq14} contains $U_{\ga_k}\sm U_{\ga_{k+1}}$, which is nonempty of dimension $d$ by construction. Thus \eq{co3eq14} for $k+1$ must have at least one fewer irreducible component of dimension $d$ than \eq{co3eq14} for $k$, so $N_{k+1}<N_k$. Hence $M_d=N_0>N_1>N_2>\cdots$ is a strictly decreasing sequence of nonnegative integers, a contradiction.

This shows that the sequence $U_{\preceq\ga_0}\supseteq U_{\preceq\ga_1}\supseteq U_{\preceq\ga_2}\supseteq\cdots$ from the beginning of the proof does stabilize, so $Y_{\preceq\ga_0}\supseteq Y_{\preceq\ga_1}\supseteq Y_{\preceq\ga_2}\supseteq\cdots$ stabilizes, and for some $N\gg 0$ we have $Y_{\preceq\ga_n}=Y_{\preceq\ga_N}$ for all $n\ge N$. But as $\ga_N\in\Ga'$ we have $Y_{\preceq\ga_N}\sm\ts\bigcup_{\text{$\ga'\prec\ga_N$ in $\Ga$}}Y_{\preceq\ga'}\ne\es$, so $Y_{\preceq\ga_n}=Y_{\preceq\ga_N}$ contradicts $\ga_N\succ\ga_n$, forcing $\ga_N=\ga_n$ for $n\ge N$. Therefore $(\Ga,\preceq),(Y_{\preceq\ga})_{\ga\in\Ga}$ satisfies the descending chain condition, as we have to prove.

The analogue for $\Th$-stratifications follows by applying the first part to the associated pseudo-$\Th$-stratification, as in Definition~\ref{co3def11}.
\end{proof}

\subsubsection{Additive weak stability conditions}
\label{co335}

We explain ideas from Alper--Halpern-Leistner--Heinloth \cite[\S 7.3]{AHLH}.

\begin{dfn}
\label{co3def12}
Work in the situation of Definition \ref{co3def6}, so that $\A$ is a $\C$-linear abelian category, and $\M$ is the moduli stack of objects in $\A$, with a decomposition $\M=\coprod_{\al\in K(\A)}\M_\al$ for $K_0(\A)\twoheadrightarrow K(\A)$. Fix~$\al\in C(\A)$.

Let $V$ be an abelian group with a total order $\le$ which is compatible with addition, in the sense that $v_0\le v_1$, $v_2\le v_3$ imply $v_0+v_2\le v_1+v_3$. Suppose $\rho:K(\A)\ra V$ is a group morphism with $\rho(\al)=0$. We call $\rho$ a {\it stability function}. If $E\in\A$ with $\lb E\rb=\al$ we call $E$ $\rho$-{\it semistable\/} if for all subobjects $F\subseteq E$ we have $\rho(\lb F\rb)\le 0$, and $\rho$-{\it unstable\/} otherwise.
\end{dfn}

The next theorem is proved in \cite[Th.~7.25]{AHLH}, using \cite[Prop.~6.14 \& Rem.~6.16]{AHLH}. The assumption that $\A$ is of compact type can be weakened significantly, but we will not need to do this.

\begin{thm}
\label{co3thm8}
In Definition\/ {\rm\ref{co3def12},} suppose that\/ $\A$ is of compact type, and being $\rho$-semistable is an open condition on $[E]\in\M_\al,$ and write $\M_\al^{\ss,\rho}\subseteq\M_\al$ for the open substack of\/ $\rho$-semistable objects. Then $\M_\al^{\ss,\rho}$ is $\Th$-reductive and S-complete with respect to any DVR $R$ essentially of finite type over\/ $\C,$ as in Definition\/ {\rm\ref{co3def8}(b),(c)} and Theorem\/~{\rm\ref{co3thm4}(ii)--(iii)}.
\end{thm}

\begin{dfn}
\label{co3def13} 
Work in the situation of Definition \ref{co3def6}. We say that a (weak) stability condition $(\tau,T,\le)$ on $\A$ is {\it additive\/} if for all $\al\in C(\A)$ there exist $V,\le,\rho$ in Definition \ref{co3def12} such that being $\rho$-semistable is an open condition on $[E]\in\M_\al$, and $(\Pi_\al^\pl)^{-1}(\M_\al^\ss(\tau))=\M_\al^{\ss,\rho}$. The name {\it additive\/} is because $\M_\al^\ss(\tau)$ is defined using an additive function~$\rho:K(\A)\ra V$.
\end{dfn}

We will show in \S\ref{co644} and \S\ref{co7110} that slope stability $(\mu,\R,\le)$ on $\A=\modCQ$ in Example \ref{co3ex1}, and Gieseker stability $(\tau,G,\le)$ on $\A=\coh(X)$ in Example \ref{co3ex2}, are both additive stability conditions. However, $\mu$-stability $(\mu,M,\le)$ on $\A=\coh(X)$ in Example \ref{co3ex3} is not additive in general when~$\dim X\ge 2$.

Theorem \ref{co3thm8} immediately implies:

\begin{cor}
\label{co3cor3}
Suppose\/ $\A$ is a $\C$-linear abelian category of compact type, and\/ $(\tau,T,\le)$ is an additive (weak) stability condition on\/ $\A$. Then\/ $(\Pi_\al^\pl)^{-1}(\M_\al^\ss(\tau))$ is $\Th$-reductive and S-complete w.r.t.\ any DVR $R$ essentially of finite type over\/ $\C,$ as in Definition\/ {\rm\ref{co3def8}(b),(c)} and Theorem\/ {\rm\ref{co3thm4}(ii),(iii),} for all\/~$\al\in C(\A)$.
\end{cor}

This is important because essentially all the stability conditions $(\tau,T,\le)$ on $\A$ that we care about are additive. However, this does {\it not\/} apply to weak stability conditions such as $\mu$-stability on $\coh(X)$ in Example \ref{co3ex3} below.

\subsubsection{Criteria for proper stable=semistable moduli stacks}
\label{co336}

We can now combine the material of \S\ref{co331}--\S\ref{co335} to prove:

\begin{thm}
\label{co3thm9}
{\bf(a)} Work in the situation of Definitions\/ {\rm\ref{co3def1}} and\/ {\rm\ref{co3def2},} so that $\A$ is a $\C$-linear abelian category, $K_0(\A)\twoheadrightarrow K(\A)$ is a surjective quotient, $\M=\coprod_{\al\in K(\A)}\M_\al,$ $\M^\pl=\coprod_{\al\in K(\A)}\M_\al^\pl$ are the (projective linear) moduli stacks of objects in $\A,$ and\/ $(\tau,T,\le)$ is a weak stability condition on\/ $\A$ such that for $\al\in C(\A)$ we have open substacks $\M_\al^\rst(\tau)\subseteq \M_\al^\ss(\tau)\subseteq\M_\al^\pl$ of\/ $\tau$-(semi)stable objects in $\A$ of class $\al$. Suppose that:
\begin{itemize}
\setlength{\itemsep}{0pt}
\setlength{\parsep}{0pt}
\item[{\bf(i)}] $\A$ is of compact type, as in {\rm\S\ref{co332}}.
\item[{\bf(ii)}] For some fixed\/ $\al\in C(\A),$ there exists a $\Th$-stratification or pseudo-$\Th$-stratification $(\Ga,\preceq),(\M_{\al,\preceq\ga})_{\ga\in\Ga}$ of\/ $\M_\al,$ as in {\rm\S\ref{co333}--\S\ref{co334},} such that\/ $\M_{\al,\preceq\ga}$ is of finite type for all\/ $\ga\in\Ga,$ with semistable locus $\M_{\al,\preceq 0}=(\Pi_\al^\pl)^{-1}(\M_\al^\ss(\tau))$. In particular, this implies that\/ $\M_\al^\ss(\tau)$ is of finite type.
\item[{\bf(iii)}] $(\tau,T,\le)$ is additive, as in {\rm\S\ref{co335}}. 
\end{itemize}
Then $\M_\al^\ss(\tau)$ admits a proper good moduli space. If also\/ $\M_\al^\rst(\tau)=\M_\al^\ss(\tau)$ then $\M_\al^\ss(\tau)$ is a proper algebraic space.
\smallskip

\noindent{\bf(b)} Suppose in addition that\/ $(\ti\tau,\ti T,\le)$ is another weak stability condition on\/ $\A$ which dominates $(\tau,T,\le),$ and\/ $\M_\al^\rst(\ti\tau)=\M_\al^\ss(\ti\tau)$. Then $
\M_\al^\ss(\ti\tau)$ is a proper algebraic space.
\end{thm}

\begin{proof} Part (a) follows from Propositions \ref{co3prop6} and \ref{co3prop8}, Theorems \ref{co3thm4}, \ref{co3thm5}, and \ref{co3thm7}, and Corollary \ref{co3cor3}. For (b), equation \eq{co3eq2} gives 
$\M_\al^\rst(\ti\tau)\subseteq\M_\al^\rst(\tau)\subseteq\M_\al^\ss(\tau)\subseteq\M_\al^\ss(\ti\tau)$, so if $\M_\al^\rst(\ti\tau)=\M_\al^\ss(\ti\tau)$ then $\M_\al^\rst(\ti\tau)=\M_\al^\rst(\tau)= \M_\al^\ss(\tau)=\M_\al^\ss(\ti\tau)$, and part (a) implies part~(b).
\end{proof}

Theorem \ref{co3thm9} is helpful for verifying Assumption \ref{co5ass2}(g),(h) below.

\begin{rem}
\label{co3rem4}
{\bf(a)} In the situations we will be interested in, usually conditions (i),(iii) above are easy to check, and the difficult things to prove are the existence of the (pseudo)-$\Th$-stratification and the finite type conditions in~(ii).

We expect that for well behaved $(\tau,T,\le)$, the $\tau$-Harder--Narasimhan filtration from Theorem \ref{co3thm2} should induce a natural pseudo-$\Th$-stratification $(\Ga,\preceq\nobreak),\ab(\M_{\al,\le\ga})_{\ga\in\Ga}$ of $\M_\al$ satisfying (ii). See Example \ref{co3ex6} when~$\A=\coh(X)$.
\smallskip

\noindent{\bf(b)} Theorems for proving the existence of $\Th$-stratifications in (ii) may be found in Halpern--Leistner \cite[Th.~2.7 \& \S 4.5]{Halp1} and \cite[Prop.~2.21 \& Th.~2.22]{Halp2}. These appear to be mainly useful for $\Q$-valued slope stability conditions, rather than $\R$-valued or more complicated Gieseker-type stability conditions.
\smallskip

\noindent{\bf(c)} Halpern-Leistner \cite[\S 5]{Halp1} studies the following situation. Let $X$ be a smooth projective $\C$-scheme, and $(Z,\A)$ be a Bridgeland stability condition on the derived category $D^b\coh(X)$, so that $\A\subset D^b\coh(X)$ is the heart of a t-structure on $D^b\coh(X)$, a $\C$-linear abelian category, and $Z:K(D^b\coh(X))=K(\A)\ra\C$ is a group morphism mapping $C(\A)\ra\bigl\{re^{i\th}:r>0$, $\th\in[0,\pi)\bigr\}$. Then $\mu:C(\A)\ra[-\iy,\iy)$ defined by $\mu(\al)=-\Re Z(\al)/\Im Z(\al)$ is a slope stability condition on~$\A$.

An object $E\in\A$ is called {\it torsion-free\/} if it has no subobject $0\ne E'\subseteq E$ with $\mu(\lb E'\rb)=-\iy$. Assuming that $Z(K(\A))\subset\Q+i\Q$ and $\A$ is noetherian and satisfies the `generic flatness' and `boundedness of quotients' conditions, Halpern-Leistner \cite[Prop.~5.40]{Halp1} constructs a $\Th$-stratification of the moduli stack $\M^{\rm tf}_\al$ of torsion-free objects in $\A$ in class $\al\in C(\A)$ for $\mu(\al)\ne -\iy$, with semistable locus $(\Pi_\al^\pl)^{-1}(\M_\al^\ss(\mu))$. He explains in \cite[\S 5.4(1)]{Halp1} how to extend this to a $\Th$-stratification $(\M_{\al,\le\ga})_{\ga\in\Ga}$ of $\M_\al$.

Under a stronger boundedness assumption which implies the $\M_{\al,\le\ga}$ are of finite type, 
Alper--Halpern-Leistner--Heinloth \cite[Ex.~7.27]{AHLH} show that Theorem \ref{co3thm9}(a)(i)--(iii) hold, so that $\M_\al^\ss(\mu)$ has a proper good moduli space for all $\al\in C(\A)$. This gives a large class of examples in which Theorem \ref{co3thm9} applies.
\smallskip

\noindent{\bf(d)} Let $X$ be a smooth projective 3-fold over $\C$. Bayer--Macr\`\i--Toda \cite{BMT} propose a `Bogomolov--Gieseker Inequality Conjecture' for $X$, and show that if it holds they can construct Bridgeland stability conditions $(Z,\A)$ on $D^b\coh(X)$. The conjecture was proved for abelian and some Calabi--Yau 3-folds by Bayer--Macr\`\i--Stellari \cite{BMS}, and for Fano 3-folds by Bernardara--Macr\`\i--Schmidt--Zhao~\cite{BMSZ}.
\smallskip

\noindent{\bf(e)} Bayer--Lahoz--Macr\`\i--Nuer--Perry--Stellari \cite[Th.~1.4]{BLMN} show that if $X$ is a smooth projective $\C$-scheme and $(Z,\A)$ is a Bridgeland stability condition on $D^b\coh(X)$ satisfying some conditions \cite[Def.~1.1]{BLMN} then moduli stacks of $(Z,\A)$-semistable objects in $D^b\coh(X)$ are of finite type and have proper good moduli spaces. See also Piyaratne--Toda \cite[Th.~1.2]{PiTo}.
\end{rem}

\section{Vertex algebras and Lie algebras on homology of moduli spaces}
\label{co4}

When $\A$ is a $\C$-linear abelian category satisfying some assumptions, and $\M$ and $\M^\pl$ are the usual and `projective linear' moduli stacks of objects in $\A$, the author \cite{Joyc12} defined the structure of a graded vertex algebra on $H_*(\M)$, and used this to induce the structure of a graded Lie algebra on $H_*(\M^\pl)$. These are deep, complicated mathematical structures. The Lie bracket on $H_*(\M^\pl)$ will be used in Chapter \ref{co5} to write down wall-crossing formulae for enumerative invariants in $\A$. We now review vertex and Lie algebras, and the theory of~\cite{Joyc12}.

\subsection{Vertex algebras and Lie algebras}
\label{co41}

For background on vertex algebras, we recommend Frenkel--Ben-Zvi \cite{FrBZ}, Kac \cite{Kac}, Lepowsky--Li \cite{LeLi}, and the author \cite{Joyc12}. Throughout we take $\K$ be a field of characteristic zero, such as $\K=\Q,\R$ or $\C$. Later we will fix~$\K=\Q$.

\begin{dfn}
\label{co4def1}
Let $V_*=\bigop_{n\in\Z}V_n$ be a graded $\K$-vector space. Write $V_*((z))\ab:= V_*[[z]][z^{-1}]$ for the vector space of Laurent series in a formal variable $z$. The vector spaces $V_*[[z]],\ab V_*((z))$ are made $\Z$-graded by declaring~$\deg z=-2$. 

A {\it field\/} on $V_*$ is a $\K$-linear map $V_*\ra V_*((z))$. The set of all fields on $V_*$ is denoted $\cF(V_*)$ and is considered as a graded $\K$-vector space by declaring $\cF(V_*)_n$ to be the set of degree $n$ fields $V\ra V((z))$ for~$n\in\Z$.

A {\it graded vertex algebra\/} $(V_*,\boo,e^{zD},Y)$ over $\K$ is a $\Z$-graded $\K$-vector space $V_*$ with an identity element $\boo\in V_0$, a grading-preserving operator $e^{zD}:V\ra V[[z]]$ with $e^{zD}v=\sum_{n\ge 0}\frac{1}{n!}D^n(v)\,z^n$ for $D:V_*\ra V_{*+2}$ the {\it translation operator}, and a grading-preserving {\it state-field correspondence\/} $Y:V_*\ra\cF(V_*)_*$ written $Y(u,z)v=\sum_{n\in\Z}u_n(v)z^{-n-1}$, where $u_n$ maps $V_*\ra V_{*+a-2n-2}$ for $u\in V_a$, satisfying:
\begin{itemize}
\setlength{\itemsep}{0pt}
\setlength{\parsep}{0pt}
\item[(i)] $Y(\boo,z)v=v$ for all $v\in V$.
\item[(ii)] $Y(v,z)\boo=e^{zD}v$ for all $v\in V$.
\item[(iii)] For all $u\in V_a$ and $v\in V_b$, there exists $N\gg 0$ such for all $w\in V_*$
\end{itemize}
\begin{equation*}
(z_1-z_2)^N\bigl(Y(u,z_1)Y(v,z_2)w-(-1)^{ab}Y(v,z_2)Y(u,z_1)w\bigr)=0\;\>\text{in $V_*[[z_1^{\pm 1},z_2^{\pm 1}]]$.}
\end{equation*}
Part (iii) is called the {\it weak commutativity property}.

Let $(V_*,\boo,e^{zD},Y)$ and $(V'_*,\boo',e^{zD'},Y')$ be graded vertex algebras over $\K$. A {\it morphism\/} $\phi:(V_*,\boo,e^{zD},Y)\ra (V'_*,\boo',e^{zD'},Y')$ is a $\K$-linear map $\phi:V_*\ra V_*'$ which preserves all the structures. That is, $\phi$ maps $V_n\ra V_n'$, and $\phi(\boo)=\boo'$, and $\phi\ci D=D'\ci\phi$, and $\phi\ci Y=Y'\ci(\phi\ot\phi)$. Such morphisms make graded vertex algebras over $\K$ into a category~$\VertAlg_\K^\rgr$.
\end{dfn}

\begin{dfn}
\label{co4def2}
A {\it graded Lie algebra\/} over $\K$ is a pair $(V_*,[\,,\,])$, where $V_*=\bigop_{a\in\Z}V_a$ is a graded $\K$-vector space, and $[\,,\,]:V_*\t V_*\ra V_*$ is a $\K$-bilinear map called the {\it Lie bracket}, which is graded (that is, $[\,,\,]$ maps $V_a\t V_b\ra V_{a+b}$ for all $a,b\in\Z$), such that for all $a,b,c\in\Z$ and $u\in V_a$, $v\in V_b$ and $w\in V_c$ we have:
\begin{gather*}
[v,u]=(-1)^{ab+1}[u,v],\\
(-1)^{ca}[[u,v],w]+(-1)^{ab}[[v,w],u]+(-1)^{bc}[[w,u],v]=0.
\end{gather*}

Let $(V_*,[\,,\,])$, $(V'_*,[\,,\,])$ be graded Lie algebras over $\K$. A {\it morphism\/} $\phi:(V_*,[\,,\,])\ra (V'_*,[\,,\,])$ is a $\K$-linear map $\phi:V_*\ra V_*'$ which preserves all the structures. That is, $\phi$ maps $V_n\ra V_n'$ and $\phi\bigl([u,v]\bigr)=\bigl[\phi(u),\phi(v)\bigr]$. Such morphisms make graded Lie algebras over $\K$ into a category~$\LieAlg_\K^\rgr$.
\end{dfn}

The next proposition is due to Borcherds \cite[\S 4]{Borc}.

\begin{prop}
\label{co4prop1}
Let\/ $(V_*,\boo,e^{zD},Y)$ be a graded vertex algebra over $\K$. We may construct a graded Lie algebra $(\check V_*,[\,,\,])$ over $\K$ as follows. Noting the shift in grading, define a $\Z$-graded\/ $\K$-vector space $\check V_*$ by 
\begin{equation*}
\check V_n=V_{n+2}/D(V_n)\qquad\text{for\/ $n\in\Z,$} 
\end{equation*}
so that\/ $\check V_*=V_{*+2}/D(V_*)$. If\/ $u\in V_{a+2}$ and\/ $v\in V_{b+2},$ the Lie bracket on $\check V_*$ is
\e
\bigl[u+D(V_a),v+D(V_b)\bigr]=u_0(v)+D(V_{a+b})\in\check V_{a+b}.
\label{co4eq1}
\e

A morphism $\phi:(V_*,\boo,e^{zD},Y)\ra (V'_*,\boo',e^{zD'},Y')$ induces a morphism $\check\phi:(V_*,[\,,\,])\ra (V'_*,[\,,\,])$ by $\check\phi\bigl(u+D(V_*)\bigr)=\phi(u)+D'(V'_*)$. Mapping $(V_*,\boo,e^{zD},Y)\ab\mapsto(\check V_*,[\,,\,])$ and\/ $\phi\mapsto\check\phi$ defines a functor~$\VertAlg_\K^\rgr\ra\LieAlg_\K^\rgr$. 
\end{prop}

\subsection{A geometric construction of vertex algebras}
\label{co42}

We now explain a special case of a geometric construction of vertex algebras from \cite{Joyc12}. We have changed some notation. Throughout $H_*(-)=H_*(-,\Q)$ and $H^*(-)=H^*(-,\Q)$ denote Betti (co)homology of Artin $\C$-stacks over~$\Q$.

\begin{ass}
\label{co4ass1}
Let $\B$ be a $\C$-linear additive category (e.g.\ an abelian category) coming from Algebraic Geometry or Representation Theory. Assume:
\begin{itemize}
\setlength{\itemsep}{0pt}
\setlength{\parsep}{0pt}
\item[(a)] We can form a natural moduli stack $\M$ of objects in $\B$, an Artin $\C$-stack, locally of finite type. Then $\C$-points of $\M$ are isomorphism classes $[E]$ of objects $E\in\B$, and the isotropy groups are $\Iso_\M([E])=\Aut(E)$.
\item[(b)] There is a natural morphism of Artin stacks $\Phi:\M\t\M\ra\M$ which on $\C$-points acts by $\Phi_*:([E],[F])\mapsto[E\op F]$, for all objects $E,F\in\B$, and on isotropy groups acts by $\Phi_*:\Iso_{\M\t\M}([E],[F])\cong\Aut(E)\t\Aut(F)\ra \Iso_\M([E\op F])\cong\Aut(E\op F)$ by $(\la,\mu)\mapsto\bigl(\begin{smallmatrix}\la & 0 \\ 0 & \mu\end{smallmatrix}\bigr)$ for $\la\in\Aut(E)$ and $\mu\in\Aut(F)$, using the obvious matrix notation for $\Aut(E\op F)$. That is, $\Phi$ is the morphism of moduli stacks induced by direct sum in the additive category $\B$. It is associative and commutative in~$\Ho(\Art_\C)$.
\item[(c)] There is a natural morphism of Artin stacks $\Psi:[*/\bG_m]\t\M\ra\M$ which on $\C$-points acts by $\Psi_*:(*,[E])\mapsto[E]$, for all objects $E$ in $\B$, and on isotropy groups acts by $\Psi_*:\Iso_{[*/\bG_m]\t\M}(*,[E])\cong\bG_m\t\Aut(E)\ra \Iso_\M([E])\cong\Aut(E)$ by $(\la,\mu)\mapsto \la\mu=(\la\cdot\id_E)\ci\mu$ for $\la\in\bG_m$ and $\mu\in\Aut(E)$. Note that $\Psi$ is {\it not\/} the same as the projection $\pi_\M:[*/\bG_m]\t\M\ra\M$ from the product $[*/\bG_m]\t\M$, which acts on isotropy groups as $(\pi_\M)_*:(\la,\mu)\mapsto\mu$. We have identities in $\Ho(\Art_\C)$:
\end{itemize}
\e
\begin{split}
\Psi\ci(\id_{[*/\bG_m]}\t\Phi)&=\Phi\ci\bigl(\Psi\ci\Pi_{12},\Psi\ci\Pi_{13}\bigr):
[*/\bG_m]\t\M^2\longra\M,\\
\Psi\ci(\id_{[*/\bG_m]}\t\Psi)&=\Psi\ci(\Om\t\id_\M):
[*/\bG_m]^2\t\M\longra\M,
\end{split}
\label{co4eq2}
\e
\begin{itemize}
\setlength{\itemsep}{0pt}
\setlength{\parsep}{0pt}
\item[] where $\Pi_{ij}$ projects to the $i^{\rm th}$ and $j^{\rm th}$ factors, and
$\Om:[*/\bG_m]^2\ra[*/\bG_m]$ is induced by the morphism $\bG_m\t\bG_m\ra\bG_m$ mapping~$(\la,\mu)\mapsto\la\mu$. 
\item[(d)] As in Definition \ref{co3def1}, we are given a surjective quotient $K_0(\B)\twoheadrightarrow K(\B)$ of the Grothendieck group $K_0(\B)$ of $\B$. We write $\lb E\rb\in K(\B)$ for the class of $E\in\B$. We suppose that if $E\in\B$ with $\lb E\rb=0$ in $K(\B)$ then $E=0$. Define the {\it positive cone\/} $C(\B)\subset K(\B)$ by~$C(\B)=\bigl\{\lb E\rb:0\ne E\in\B\bigr\}$.

We require that the map $\M(\C)\ra K(\B)$ mapping $E\mapsto\lb E\rb$ should be locally constant. This gives a decomposition $\M=\coprod_{\al\in K(\B)}\M_\al$ of $\M$ into open and closed $\C$-substacks $\M_\al\subset \M$ of objects in class $\al$, where $\M_0=\{[0]\}$. We write $\Phi_{\al,\be}=\Phi\vert_{\M_\al\t\M_\be}:\M_\al\t\M_\be\ra\M_{\al+\be}$ and $\Psi_\al=\Psi\vert_{[*/\bG_m]\t\M_\al}:[*/\bG_m]\t\M_\al\ra\M_\al$. 
\item[(e)] We are given a biadditive form $\chi:K(\B)\t K(\B)\ra \Z$.
\item[(f)] We are given a perfect complex $\cE^\bu$ on $\M\t\M$, such that the restriction $\cE_{\al,\be}^\bu:=\cE^\bu\vert_{\M_\al\t\M_\be}$ to $\M_\al\t\M_\be$ has rank $\chi(\al,\be)$ for all $\al,\be\in K(\B)$, and there are isomorphisms of perfect complexes
\ea
(\Phi\t\id_\M)^*(\cE^\bu)&\cong
\Pi_{13}^*(\cE^\bu)\op \Pi_{23}^*(\cE^\bu),
\label{co4eq3}\\
(\id_\M\t\Phi)^*(\cE^\bu)&\cong
\Pi_{12}^*(\cE^\bu)\op \Pi_{13}^*(\cE^\bu),
\label{co4eq4}\\
(\Psi\t\id_\M)^*(\cE^\bu)&\cong \Pi_1^*(L_{[*/\bG_m]})\ot \Pi_{23}^*(\cE^\bu),
\label{co4eq5}\\
(\Pi_2,\Psi\ci\Pi_{13})^*(\cE^\bu)&\cong \Pi_1^*(L_{[*/\bG_m]}^*)\ot\Pi_{23}^*(\cE^\bu).
\label{co4eq6}
\ea
Here equations \eq{co4eq3}--\eq{co4eq4} are on $\M\t\M\t\M$, and \eq{co4eq5}--\eq{co4eq6} on $[*/\bG_m]\t\M\t\M$. We write $\Pi_i$ for the projection to the $i^{\rm th}$ factor, and $\Pi_{ij}$ for the projection to the product of the $i^{\rm th}$ and $j^{\rm th}$ factors. We write $L_{[*/\bG_m]}\ra[*/\bG_m]$ for the line bundle corresponding to the weight 1 representation of $\bG_m=\C\sm\{0\}$ on~$\C$.
\item[(g)] For each $\al\in C(\B)$, we suppose that there should exist a line bundle $L_\al\ra\M_\al$ such that $\Psi_\al^*(L_\al)\cong L^{N_\al}_{[*/\bG_m]}\bt L_\al$ for some $N_\al\in\Z\sm\{0\}$. This is needed to prove that \eq{co4eq15} is an isomorphism in Theorem \ref{co4thm2} below.

If $F\in\B$ with $\lb F\rb=\be\in K(\B)$ and $\chi(\al,\be)=N_\al\ne 0$ then \eq{co4eq5} implies that $L_\al=\det(\cE_{\al,\be}\vert_{\M_\al\t\{[F]\}})$ will do. Similarly, if $\chi(\be,\al)=-N_\al\ne 0$ then \eq{co4eq6} implies that $L_\al=\det(\cE_{\be,\al}\vert_{\{[F]\}\t\M_\al})$ will do. Thus, if $\chi$ in (e) is nondegenerate then the condition is automatic.
\end{itemize}	
\end{ass}

As explained in \cite{GJT} and \cite{Joyc12}, and we will show in Chapters \ref{co6}--\ref{co8}, there are natural choices for all this data in many large classes of interesting examples.

\begin{dfn}
\label{co4def3}
Suppose Assumption \ref{co4ass1} holds. Given all this data, we define a graded vertex algebra structure on the homology $H_*(\M)$ from \S\ref{co21}. The inclusion of the zero object $0\in\B$ gives a morphism $[0]:*\hookra\M$ inducing $\Q\cong H_*(*)
\ra H_*(\M)$, and we define $\boo\in H_*(\M)$ to be the image of $1\in \Q$ under this map. Taking homology of $\Psi$ gives a map
\e
\xymatrix@C=30pt{ H_*([*/ \bG_m])\!\ot_\Q\! H_*(\M)  \ar[r]^(0.57){\bt} & H_*([*/ \bG_m] \!\t\! \M) \ar[r]^(0.61){\Psi_*} & H_*(\M). }
\label{co4eq7}
\e
As $H^*([*/\bG_m])\cong\Hom_\Q(H_*([*/ \bG_m]),\Q)$, this is equivalent to a map
\e
\xymatrix@C=35pt{ H_*(\M) \ar[r] & H_*(\M) \hat\ot_\Q H^*([*/\bG_m]) \ar[r]^(0.61){\eq{co2eq2}} & H_*(\M)[[z]], }
\label{co4eq8}
\e
using equation \eq{co2eq2}, and we denote the composition $e^{zD}$.

The decomposition $\M=\coprod_{\al\in K(\B)}\M_\al$ induces an identification
\e
H_*(\M)=\bigop_{\al\in K(\B)}H_*(\M_\al).
\label{co4eq9}
\e
For $u\in H_*(\M_\al)\subset H_*(\M)$ and $v\in H_*(\M_\be)\subset H_*(\M)$, define
\e
\begin{split}
&Y(u,z)v=Y(z)(u\ot v)= (-1)^{\chi(\al,\be)} \sum\nolimits_{i,j\ge 0} z^{\chi(\al, \be)+\chi(\be,\al)-i+j}\cdot{}\\
&\bigl(\Phi_{\al,\be}\ci(\Psi_\al\t\id_{\M_\be})\bigr)_*\bigl(t^j\bt ((u \bt v) \cap c_i(\cE_{\al,\be}^\bu\op \si_{\al,\be}^*(\cE^\bu_{\be,\al})^\vee))\bigr),
\end{split}
\label{co4eq10}
\e
where $t^j\in H_{2j}([*/\bG_m])$ is as in \eq{co2eq3}, and $\si_{\al,\be}:\M_\al\t\M_\be\ra\M_\be\t\M_\al$ exchanges the factors. Using \eq{co4eq9}, for $n\in\Z$ and $\al\in K(\B)$ we write
\e
\hat H_n(\M_\al)= H_{n-2\chi(\al,\al)}(\M_\al),\quad \hat H_n(\M)=\bigop_{\al\in K(\B)}\hat H_n(\M_\al).
\label{co4eq11}
\e
That is, $\hat H_*(\M)$ is $H_*(\M)$, but with grading shifted by $-2\chi(\al,\al)$ in the component $H_*(\M_\al)\subset H_*(\M)$. The author \cite{Joyc12} proves:	
\end{dfn}

\begin{thm}
\label{co4thm1}
$(\hat H_*(\M), \boo,e^{zD},Y)$ above is a graded vertex algebra over $\Q$.
\end{thm}

\subsection{Lie algebras from the vertex algebras of \S\ref{co42}}
\label{co43}

In the situation of \S\ref{co42}, Proposition \ref{co4prop1} makes $\hat H_{*+2}(\M)/D(\hat H_*(\M))$ into a graded Lie algebra. We can interpret this as the shifted homology $\check H_*(\M^\pl)$ of a modification $\M^\pl$ of $\M$, as in Definition \ref{co3def2}, which we now describe.

\begin{dfn}
\label{co4def4}
Continue in the situation of Assumption \ref{co4ass1} and Definition \ref{co4def3}. Then $[*/\bG_m]$ is a group stack, and $\Psi:[*/\bG_m]\t\M\ra\M$ is an action of $[*/\bG_m]$ on $\M$, which is trivial on $\M_0=\{[0]\}$, and free on $\M\sm\{[0]\}$. As explained in \cite{Joyc12}, we may take the quotient of $\M$ by $\Psi$ to get a stack $\M^\pl$, which we call the {\it projective linear moduli stack}, with projection $\Pi^\pl:\M\ra\M^\pl$, in a 2-co-Cartesian square in $\Art_\C$:
\begin{equation*}
\xymatrix@C=130pt@R=15pt{ *+[r]{[*/\bG_m]\t\M} \drtwocell_{}\omit^{}\omit{^{}} \ar[r]_(0.65){\Psi} \ar[d]^{\pi_\M} & *+[l]{\M} \ar[d]_{\Pi^\pl} \\
*+[r]{\M} \ar[r]^(0.35){\Pi^\pl} & *+[l]{\M^\pl.\!} }
\end{equation*}

This construction of $\M^\pl$ is known in the literature as {\it rigidification}, as in Abramovich--Olsson--Vistoli \cite{AOV} and Romagny \cite{Roma}, written $\M^\pl=\M\!\!\fatslash\,\bG_m$ in \cite{AOV,Roma}. We regard $\M^\pl$ as the moduli stack of objects in $\B$ `up to projective linear isomorphisms', where a {\it projective linear isomorphism\/} $[\phi]:E\ra F$ is a $\sim$-equivalence class of isomorphisms $\phi:E\ra F$ in $\coh(X)$, where $\phi\sim\phi'$ if $\phi=\la\phi'$ for some $\la\in\bG_m$. Then $\Pi^\pl:\M\ra\M^\pl$ is an isomorphism on $\C$-points. Thus, $\C$-points $x\in\M^\pl(\C)$ correspond naturally to isomorphism classes $[E]$ of objects $E\in\B$. We will write points of $\M^\pl(\C)$ as $[E]$, so $\Pi^\pl(\C)$ maps~$[E]\mapsto [E]$.  

The isotropy groups of $\M^\pl$ satisfy $\Iso_{\M^\pl}([E])\cong\Iso_\M([E])/\bG_m$ for $E\ne 0$, where the $\bG_m$-subgroup of $\Iso_\M([E])$ is determined by the action of $\Psi$ on isotropy groups. Thus by Assumption \ref{co4ass1}(c) we see that 
\e
\Iso_{\M^\pl}([E])\cong\Aut(E)/(\bG_m\cdot\id_E).
\label{co4eq12}
\e
The action of $\Pi^\pl$ on isotropy groups is given by the commutative diagram
\begin{equation*}
\xymatrix@C=150pt@R=15pt{ *+[r]{\Iso_\M([E])} \ar[r]_(0.4){\Pi^\pl_*} \ar[d]^\cong & *+[l]{\Iso_{\M^\pl}([E])} \ar[d]^\cong_{\eq{co4eq12}} \\
*+[r]{\Aut(E)} \ar[r]^(0.4){\ep\,\longmapsto\, \ep\bG_m} & *+[l]{\Aut(E)/(\bG_m\cdot\id_E).\!\!} }	
\end{equation*}

The splitting $\M=\coprod_{\al\in K(\B)}\M_\al$ descends to $\M^\pl=\coprod_{\al\in K(\B)}\M_\al^\pl$, with $\M_\al^\pl=\M_\al\!\!\fatslash\,\bG_m$ for $\al\ne 0$ and $\M_0^\pl=\M_0=*$. Thus as for \eq{co4eq9} we have
\e
H_*(\M^\pl)=\bigop_{\al\in K(\B)}H_*(\M_\al^\pl).
\label{co4eq13}
\e
In a similar way to \eq{co4eq11}, using \eq{co4eq13}, for $n\in\Z$ and $\al\in K(\B)$ we write
\e
\check H_n(\M_\al^\pl)\!=\!H_{n+2-2\chi(\al,\al)}(\M_\al^\pl),\;\> \check H_n(\M^\pl)\!=\bigop_{\!\!\!\!\!\!\!\al\in K(\B)\!\!\!\!\!\!\!\!\!\!\!}\check H_n(\M_\al^\pl).
\label{co4eq14}
\e
That is, $\check H_*(\M^\pl)$ is $H_*(\M^\pl)$, but with grading shifted by $2-2\chi(\al,\al)$ in the component $H_*(\M_\al^\pl)\subset H_*(\M^\pl)$.
\end{dfn}

The next theorem, proved in \cite{Joyc12} using Assumption \ref{co4ass1}(g), gives a geometric interpretation of the graded Lie algebra~$\hat H_{*+2}(\M)/D(\hat H_*(\M))$.

\begin{thm}
\label{co4thm2}
Work in the situation of Assumption\/ {\rm\ref{co4ass1}} and Definitions\/ {\rm\ref{co4def3}} and\/ {\rm\ref{co4def4},} and consider the graded Lie algebra $\hat H_{*+2}(\M)/D(\hat H_*(\M))$ constructed by combining Proposition\/ {\rm\ref{co4prop1}} and Theorem\/ {\rm\ref{co4thm1}}. Then $\Pi^\pl:\M\ra\M^\pl$ gives a morphism $\Pi^\pl_*:H_*(\M)\ab\ra H_*(\M^\pl)$. It is surjective, with kernel $D(H_*(\M))$. This induces an isomorphism
\begin{equation*}
\Pi^\pl_*:H_*(\M)/D(H_*(\M))\longra H_*(\M^\pl).
\end{equation*}
Comparing \eq{co4eq11} and\/ {\rm\eq{co4eq14},} we see this is an isomorphism for all\/~{\rm$n\in\Z$:}
\e
\hat H_{n+2}(\M)/D(\hat H_n(\M))\longra \check H_n(\M^\pl).
\label{co4eq15}
\e
Thus there is a unique Lie bracket\/ $[\,,\,]$ on $\check H_*(\M^\pl)$ making it into a \begin{bfseries}graded Lie algebra\end{bfseries}, such that\/ \eq{co4eq15} is a Lie algebra isomorphism. Hence $\check H_0(\M^\pl)$ and\/ $\check H_{\rm even}(\M^\pl)=\bigop_{n\in\Z}\check H_{2n}(\M^\pl)$ are Lie algebras.
\end{thm}

\begin{rem}
\label{co4rem1}
{\bf(Ways to explain the Lie bracket on $H_*(\M^\pl)$.)}

\noindent{\bf(a)} Theorem \ref{co4thm2} implies that if $u\in H_*(\M_\al)$ and $v\in H_*(\M_\be)$ then
\begin{equation*}
\bigl[\Pi^\pl_*(u),\Pi^\pl_*(v)\bigr]=\Pi^\pl_*\bigl(\Res_z(Y(u,z)v)\bigr),
\end{equation*}
where $Y(u,z)v$ is as in \eq{co4eq10}. That is, after lifting along $\Pi^\pl_*:H_*(\M)\twoheadrightarrow H_*(\M^\pl)$, the Lie bracket on $H_*(\M^\pl)$ is given by the residue in $z$ of the explicit but complicated formula \eq{co4eq10}.
\smallskip

\noindent{\bf(b)} It would be desirable to construct the Lie bracket $[\,,\,]$ on $\check H_*(\M^\pl)$ in Theorem \ref{co4thm2} directly on $\M^\pl$, rather than indirectly by pushdown from $\hat H_*(\M)$. In \cite{Joyc12} the author outlined a construction of $[\,,\,]$ using a conjectural characteristic class type operation on principal $[*/\bG_m]$-bundles such as $\Pi^\pl:\M\ra\M^\pl$, called the `projective Euler class'. These conjectures have now been proved by Upmeier \cite{Upme}, giving an alternative definition of~$[\,,\,]$.

For a heuristic explanation, consider the (informal) diagram:
\e
\xymatrix@C=25pt{ 
\M_\al^\pl\!\t\!\M_\be^\pl & {\begin{subarray}{l}\ts \bP\bigl(\cE_{\al,\be}^\bu[1]\op \\ \ts \si_{\al,\be}^*(\cE^\bu_{\be,\al})^\vee[1]\bigr)\end{subarray}} \ar[l]_(.51)\Pi \ar[r]^{\ti\Pi} & (\M_\al\!\t\!\M_\be)^\pl \ar[r]^(0.6){\Phi_{\al,\be}^\pl} & \M_{\al+\be}^\pl. }\!\!\!
\label{co4eq16}
\e
Here $(\M_\al\t\M_\be)^\pl$ means the quotient of $\M_\al\t\M_\be$ by the diagonal $[*/\bG_m]$-action. Then $\Phi_{\al,\be}:\M_\al\t\M_\be\ra\M_{\al+\be}$ descends to $\Phi_{\al,\be}^\pl$ as shown. We do not claim the projective bundle $\bP(\cdots)$ is well defined in general. But if $\cE_{\al,\be}^\bu[1]\op \si_{\al,\be}^*(\cE^\bu_{\be,\al})^\vee[1]$ were a vector bundle over $\M_\al\t\M_\be$ with fibre $\C^n$ then \eq{co4eq16} would be well defined, and $\Pi$ would be a bundle with fibre $\CP^{n-1}$, and $\ti\Pi$ a bundle with fibre $\C^n\sm\{0\}$. 

Then roughly speaking, the Lie bracket on $H_*(\M^\pl)$ should be given by the commutative diagram:
\begin{equation*}
\xymatrix@C=60pt@R=10pt{ 
*+[r]{H_a(\M_\al^\pl)\ot H_b(\M_\be^\pl)} \ar[rr]^{(-1)^{\chi(\al,\be)}[\,,\,]} \ar[d]^{\bt} && *+[l]{H_{\begin{subarray}{l} a+b-2\chi(\al,\be) \\ -2\chi(\be,\al)-2\end{subarray}}(\M_{\al+\be}^\pl)} \\
*+[r]{H_{a+b}(\M_\al^\pl\!\t\!\M_\be^\pl)} \ar@/_1.5pc/[rr]_{\mathop{\rm PE}(\cE_{\al,\be}^\bu\op\si_{\al,\be}^*(\cE^\bu_{\be,\al})^\vee)} \ar@{..>}[r]^(0.6){\Pi^!} & {\raisebox{13pt}{${\begin{subarray}{l}\ts H_{\begin{subarray}{l} a+b-2\chi(\al,\be) \\ -2\chi(\be,\al)-2\end{subarray}} \\ \ts \bigl(\bP\bigl(\cE_{\al,\be}^\bu[1]\op \\ \ts \si_{\al,\be}^*(\cE^\bu_{\be,\al})^\vee[1]\bigr)\bigr)\end{subarray}}$}}   \ar[r]^{\ti\Pi_*} & *+[l]{\begin{subarray}{l}\ts H_{\begin{subarray}{l} a+b-2\chi(\al,\be) \\ -2\chi(\be,\al)-2\end{subarray}} \\ \ts  \bigl((\M_\al\!\t\!\M_\be)^\pl\bigr).\end{subarray}} \ar[u]^{(\Phi_{\al,\be}^\pl)_*}  }
\end{equation*}
Here $\Pi^!$ (if it is defined) is a Gysin map, a wrong-way map on homology, which changes degree by the real virtual dimension of the fibre of the projective bundle $\bP(\cdots)$. The `projective Euler class' $\mathop{\rm PE}(\cE_{\al,\be}^\bu\op\si_{\al,\be}^*(\cE^\bu_{\be,\al})^\vee)$ is morally the composition $\ti\Pi_*\ci\Pi^!$, but it is defined in \cite{Upme} without defining $\bP(\cdots)$ or~$\Pi^!$.
\end{rem}

\begin{rem}
\label{co4rem2}
{\bf(Extension to derived categories.)} Most of \S\ref{co42}--\S\ref{co43} also works when $\B$ is a derived category, such as $D^b\coh(X)$ for $X$ a smooth projective $\C$-scheme, or $D^b\modCQ$ for $Q$ a quiver. These are examples of $\C$-linear additive categories. Some parts of Assumption \ref{co4ass1} do not hold for derived categories, and must be modified:
\begin{itemize}
\setlength{\itemsep}{0pt}
\setlength{\parsep}{0pt}
\item[(i)] In Assumption \ref{co4ass1}(a) and Definition \ref{co4def4}, the moduli stacks of objects in derived categories like $D^b\coh(X)$ or $D^b\modCQ$ are {\it higher\/ $\C$-stacks\/} in the sense of \S\ref{co22}, not Artin $\C$-stacks. The (co)homology of higher $\C$-stacks is well behaved, so this makes no real difference. 
\item[(ii)] In Assumption \ref{co4ass1}(d) the condition that if $E\in\B$ with $\lb E\rb=0$ in $K(\B)$ then $E=0$, which forces $\M_0=\{0\}$, is false for derived categories. In general $\M_0\ne\{0\}$, and $\M_0$ is rather large.
\item[(iii)] Assumption \ref{co4ass1}(g) is generally false when $\al=0$. 
\end{itemize}	

None of (i)--(iii) affect Definition \ref{co4def3} and Theorem \ref{co4thm1}, which hold unchanged for derived categories. In Theorem \ref{co4thm2}, part (iii) means that $\Pi^\pl_*:H_k(\M_\al)/D(H_{k-2}(\M_\al))\ra H_k(\M_\al^\pl)$ need only be an isomorphism for $\al\ne 0$. In \cite{Joyc12} we prove that $\Pi^\pl_*:H_k(\M_0)/D(H_{k-2}(\M_0))\ra H_k(\M_0^\pl)$ is an isomorphism for $k\le 2$, and give examples in which it is not for~$k=3$. 

Thus, in the derived category case, we cannot use the graded vertex algebra structure on $\hat H_*(\M)$ to induce a graded Lie bracket on $\check H_*(\M^\pl)$ (though this does work on nearly all of $\check H_*(\M^\pl)$). However, the definition of $[\,,\,]$ on $\check H_*(\M^\pl)$ using projective Euler classes outlined in Remark \ref{co4rem1} works for derived categories, and \eq{co4eq15} is a Lie algebra morphism.
\end{rem}
    
\subsection{Morphisms of the vertex and Lie algebras of \S\ref{co42}--\S\ref{co43}}
\label{co44}

The material of this section will not be used later, but is included for interest.

\begin{dfn}
\label{co4def5}
Let $\B,\M,\Phi,\Psi,K(\B),\chi,\cE^\bu$ and 
$\B',\M',\Phi',\ab\Psi',\ab K(\B'),\ab\chi',\cE^{\prime\bu}$ be two sets of data satisfying Assumption \ref{co4ass1}. Write 
$(\hat H_*(\M),\ab \boo,\ab e^{zD},\ab Y)$ and $(\hat H_*(\M'), \boo',e^{zD'},Y')$ for the corresponding graded vertex algebras from Theorem \ref{co4thm1}. Suppose:
\begin{itemize}
\setlength{\itemsep}{0pt}
\setlength{\parsep}{0pt}
\item[(a)] We are given a $\C$-linear exact functor $\Si:\B\ra\B'$. This should induce a morphism $\si:\M\ra\M'$ of moduli stacks, which acts on $\C$-points as $\si([E])=[\Si(E)]$. The induced morphism $\Si_*:K_0(\B)\ra K_0(\B')$ descends to $\Si_*:K(\B)\ra K(\B')$, with $\Si_*(\lb E\rb)=\lb \Si(E)\rb$. In $\Ho(\Art_\C)$ we have
\begin{equation*}
\Phi'\ci(\si\t\si)=\si\ci\Phi,\qquad 
\Psi'\ci(\id_{[*/\bG_m]}\t\si)=\si\ci\Psi.
\end{equation*}
\item[(b)] We are given a biadditive morphism $\xi:K(\B)\t K(\B)\ra\Z$.
\item[(c)] We are given a vector bundle $F\ra\M\t\M$ of mixed rank, such that $F_{\al,\be}:=F\vert_{\M_\al\t\M_\be}$ has rank $\xi(\al,\be)$ for all $\al,\be\in K(\B)$. We also write $G\ra\M$ for the vector bundle $\De_\M^*(F^*)$, where $\De_\M:\M\ra\M\t\M$ is the diagonal morphism. Then $G_\al:=G\vert_{\M_\al}$ has rank $\xi(\al,\al)$.
\end{itemize}

All this data should satisfy:
\begin{itemize}
\setlength{\itemsep}{0pt}
\setlength{\parsep}{0pt}
\item[(i)] $\chi'(\Si_*(\al),\Si_*(\be))=\chi(\al,\be)+\xi(\al,\be)$ for all $\al,\be\in K(\B)$.
\item[(ii)] As for \eq{co4eq3}--\eq{co4eq6}, there are isomorphisms of vector bundles
\begin{align*}
(\Phi\t\id_\M)^*(F)&\cong
\Pi_{13}^*(F)\op \Pi_{23}^*(F),\\
(\id_\M\t\Phi)^*(F)&\cong
\Pi_{12}^*(F)\op \Pi_{13}^*(F),\\
(\Psi\t\id_\M)^*(F)&\cong \Pi_1^*(L_{[*/\bG_m]})\ot \Pi_{23}^*(F),\\
(\Pi_2,\Psi\ci\Pi_{13})^*(F)&\cong \Pi_1^*(L_{[*/\bG_m]}^*)\ot\Pi_{23}^*(F).
\end{align*}
\item[(iii)] In $K_0(\Perf(\M\t\M))$ we have
\begin{equation*}
(\si\t\si)^*([\cE^{\prime\bu}])=[\cE^\bu]+[F].
\end{equation*}
\item[(iv)] There is a vector bundle $G^\pl\ra\M^\pl$ with~$G\cong(\Pi^\pl)^*(G^\pl)$. 
\end{itemize}

Write $c_\top(G)\in H^*(\M)$ for the top Chern class $c_{\rank G}(G)$ of $G$. Note that as $\rank G$ depends on the component $\M_\al\subset\M$ as in (c), we have $c_\top(G)\vert_{\M_\al}=c_{\xi(\al,\al)}(G_\al)$. Define a $\Q$-linear morphism
\e
\Om:H_*(\M)\longra H_*(\M')\quad\text{by}\quad \Om(u)=H_*(\si)\bigl(u\cap c_\top(G)\bigr).
\label{co4eq17}
\e
Here if $u\in H_a(\M_\al)$ then $\Om(u)\in H_{a-2\xi(\al,\al)}(\M'_{\Si_*(\al)})$. Combining this with \eq{co4eq11} and (i), we see that $\Om:\hat H_*(\M)\ra\hat H_*(\M')$ is grading-preserving. The next theorem is proved in~\cite{Joyc12}:
\end{dfn}

\begin{thm}
\label{co4thm3}
In Definition\/ {\rm\ref{co4def5},} $\Om:(\hat H_*(\M), \boo,e^{zD},Y)\ra 
(\hat H_*(\M'),\ab \boo',\ab e^{zD'},\ab Y')$ is a morphism of graded vertex algebras.
\end{thm}

\begin{dfn}
\label{co4def6}
Work in the situation of Definition \ref{co4def5}. Then $\si:\M\ra\M'$ descends to a morphism $\si^\pl:\M^\pl\ra\M^{\prime\pl}$ with $\si^\pl\ci\Pi^\pl=\Pi^{\prime\pl}\ci\si$ in $\Ho(\Art_\C)$. As for \eq{co4eq17}, define a $\Q$-linear morphism
\e
\Om^\pl:H_*(\M^\pl)\longra H_*(\M^{\prime\pl})\;\>\text{by}\;\> \Om^\pl(u)=H_*(\si^\pl)\bigl(u\cap c_\top(G^\pl)\bigr).
\label{co4eq18}
\e
By \eq{co4eq14} we see that $\Om^\pl:\check H_*(\M^\pl)\longra\check H_*(\M^{\prime\pl})$ is grading preserving. From $\si^\pl\ci\Pi^\pl=\Pi^{\prime\pl}\ci\si$, $G\cong(\Pi^\pl)^*(G^\pl)$, and \eq{co4eq17}--\eq{co4eq18} we see that the following diagram commutes:
\begin{equation*}
\xymatrix@C=120pt@R=15pt{
*+[r]{H_*(\M)} \ar[r]_{\Om} \ar[d]^{\Pi^\pl_*} & *+[l]{H_*(\M')} \ar[d]_{H_*(\ti\Pi^\pl)} \\ 
*+[r]{H_*(\M^\pl)} \ar[r]^{\Om^\pl} & *+[l]{H_*(\M^{\prime\pl}).\!} }	
\end{equation*}
Hence Theorems \ref{co4thm1}, \ref{co4thm2}, and \ref{co4thm3} imply:
\end{dfn}

\begin{cor}
\label{co4cor1}
In Definition\/ {\rm\ref{co4def6},} $\Om^\pl:\check H_*(\M^\pl)\ra \check H_*(\M^{\prime\pl})$ is a morphism of the graded Lie algebras in Theorem\/~{\rm\ref{co4thm2}}.
\end{cor}

\subsection{Extension to equivariant (co)homology}
\label{co45}

Next we explain how to generalize the material of \S\ref{co42}--\S\ref{co44} to be equivariant under the action of a linear algebraic $\C$-group $G$, and we work with equivariant homology $H_*^G(\M),H_*^G(\M^\pl)$, as in \S\ref{co23}.

\begin{ass}
\label{co4ass2}
Let Assumption \ref{co4ass1} hold for $\B,\M,\Phi,\Psi,\M^\pl,K(\B),\ldots.$ Suppose also that we are given a linear algebraic $\C$-group $G$, and:
\smallskip

\noindent{\bf(a)} We are given an action of the group $G(\C)$ of $\C$-points of $G$ on the $\C$-linear additive category $\B$. Explicitly, this means:
\begin{itemize}
\setlength{\itemsep}{0pt}
\setlength{\parsep}{0pt}
\item[(i)] For each $g\in G(\C)$ we are given an equivalence 
$\Ga(g):\B\,{\buildrel\sim\over\longra}\,\B$ of $\C$-linear exact categories.
\item[(ii)] For all $g,h\in G(\C)$ we are given a natural isomorphism of functors
\begin{equation*}
\Ga_{g,h}:\Ga(g)\ci \Ga(h)\Longra \Ga(gh),	
\end{equation*}
such that for all $f,g,h\in G(\C)$ the following commutes:
\e
\begin{gathered}
\xymatrix@C=140pt@R=15pt{ *+[r]{\Ga(f)\ci \Ga(g)\ci \Ga(h)} \ar@{=>}[r]_{\Ga_{f,g}*\id_{\Ga(h)}} \ar@{=>}[d]^{\id_{\Ga(f)} *\Ga_{g,h}} & *+[l]{\Ga(fg)\ci \Ga(h)}  \ar@{=>}[d]_{\Ga_{f\ci g,h}} \\
*+[r]{\Ga(f)\ci \Ga(gh)} \ar@{=>}[r]^{\Ga_{f,g\ci h}} & *+[l]{\Ga(fgh).\!} }	
\end{gathered}
\label{co4eq19}
\e
\item[(iii)] We are given a natural isomorphism $\Ga_1:\Ga(1)\Ra \Id_\B$, where for all $g\in G$
\end{itemize}
\e
\begin{split}
\Ga_{g,1}&=\id_{\Ga(g)}*\Ga_1: \Ga(g)\ci \Ga(1)\Longra \Ga(g),\\ 
\Ga_{1,g}&=\Ga_1*\id_{\Ga(g)}:\Ga(1)\ci \Ga(g)\Longra \Ga(g).
\end{split}
\label{co4eq20}
\e
\noindent{\bf(b)} The $G$-action on $\B$ in {\bf(a)} is not just at the level of $\C$-points of $G$, but is algebraic, and compatible with moduli functors. In particular, it induces an algebraic action of $G$ on the Artin $\C$-stack $\M$.
\smallskip

\noindent{\bf(c)} For $g\in G(\C)$, as $\Ga(g)$ is an equivalence of exact categories, for $E,F\in\B$ there is a canonical isomorphism $\Ga(g)(E)\op \Ga(g)(F)\cong \Ga(g)(E\op F)$, where $[E\op F]=\Phi([E],[F])$. That is, $\Phi\ci(\Ga(g),\Ga(g))([E],[F])=\Ga(g)\ci\Phi([E],[F])$ for $\C$-points $([E],[F])\in\M\t\M$. Thus {\bf(b)} implies that $\Phi:\M\t\M\ra\M$ is $G$-equivariant (that is, we are given a $G$-equivariant structure on $\Phi$), for the diagonal action of $G$ on $\M\t\M$.
\smallskip

\noindent{\bf(d)} Similarly, by the $\C$-linearity of the $G$-action in {\bf(a)},{\bf(b)}, we can show that $\Psi:[*/\bG_m]\t\M\ra\M$ is $G$-equivariant (that is, we are given a $G$-equivariant structure on $\Psi$), where $G$ acts trivially on $[*/\bG_m]$. Hence the $G$-action on $\M$ in {\bf(b)} descends to a $G$-action on $\M^\pl=\M/[*/\bG_m]$, and $\Pi^\pl:\M\ra\M^\pl$ is $G$-equivariant.
\smallskip

\noindent{\bf(e)} Classes $\lb E\rb\in K(\B)$ are invariant under the action of $G$ on $E\in\B$. Thus $G$ acts on $\M_\al\subset\M$ and $\M_\al^\pl\subset\M^\pl$ for each~$\al\in K(\B)$.
\smallskip

\noindent{\bf(f)} The perfect complex $\cE^\bu$ on $\M\t\M$ in Assumption \ref{co4ass1}(f) is equivariant under the diagonal action of $G$ on $\M\t\M$. The isomorphisms \eq{co4eq3}--\eq{co4eq6} hold in $G$-equivariant perfect complexes.
\end{ass}

Here are the analogues of Definition \ref{co4def3} and Theorem \ref{co4thm1}.

\begin{dfn}
\label{co4def7}
Suppose Assumption \ref{co4ass2} holds. We will define a graded vertex algebra structure on the equivariant homology $H_*^G(\M)$ of $\M$, as in \S\ref{co23}, which is defined by Assumption \ref{co4ass2}(b). Since $\M_0\cong *$ by Assumption \ref{co4ass1}(d), Definition \ref{co2def2}(c) shows that $H_0^G(\M_0)\cong H^0([*/G])$. Let $\boo\in H_0^G(\M)$ correspond to~$1\in H^0([*/G])$.

Generalizing \eq{co4eq7}, using Definition \ref{co2def2}(a),(d) we have morphisms
\begin{equation*}
\xymatrix@C=35pt{ {\begin{subarray}{l} \ts H_*([*/ \bG_m])\ot_\Q H_*^G(\M)= \\ \ts H_*^{\{1\}}([*/ \bG_m])\ot_\Q H_*^G(\M)\end{subarray}}  \ar[r]^(0.54){\bt} & {\begin{subarray}{l} \ts H_*^G([*/\bG_m] \!\t\! \M)= \\ \ts H_*^{\{1\}\t G}([*/\bG_m] \!\t\! \M)\end{subarray}} \ar[r]^(0.61){\Psi_*} & H_*^G(\M). }
\end{equation*}
As $H^*([*/\bG_m])\cong\Hom_\Q(H_*([*/ \bG_m]),\Q)$, as in \eq{co4eq8} this is equivalent to a map
\begin{equation*}
\xymatrix@C=35pt{ H^G_*(\M) \ar[r] & H^G_*(\M) \hat\ot_\Q H^*([*/\bG_m]) \ar[r]^(0.61){\eq{co2eq2}} & H^G_*(\M)[[z]], }
\end{equation*}
using equation \eq{co2eq2}, and we denote the composition $e^{zD}$.

Let $\al,\be\in K(\B)$ and $a,b,n\in\Z$. Then Definition \ref{co2def2}(d),(g) give morphisms
\ea
\bt&:H_a^G(\M_\al)\ot_\Q H_b^G(\M_\be)\longra H_{a+b}^{G\t G}(\M_\al\t\M_\be),
\label{co4eq21}\\
\La^{G\t G,G}&:H_n^{G\t G}(\M_\al\t\M_\be)\longra H_n^G(\M_\al\t\M_\be),
\label{co4eq22}
\ea
where the inclusion $G\hookra G\t G$ is the diagonal morphism.

The decomposition $\M=\coprod_{\al\in K(\B)}\M_\al$ induces an identification
\e
H_*^G(\M)=\bigop_{\al\in K(\B)}H_*^G(\M_\al).
\label{co4eq23}
\e
For $u\in H_a^G(\M_\al)\subset H_*^G(\M)$ and $v\in H_b^G(\M_\be)\subset H_*^G(\M)$, define
\ea
&Y(u,z)v=Y(z)(u\ot v)= (-1)^{\chi(\al,\be)} \sum\nolimits_{i,j\ge 0} z^{\chi(\al, \be)+\chi(\be,\al)-i+j}\cdot{}
\label{co4eq24}\\
&\bigl(\Phi_{\al,\be}\ci(\Psi_\al\t\id_{\M_\be})\bigr)_*\bigl(t^j\bt ((\La^{G\t G,G}(u\bt v)) \cap c_i(\cE_{\al,\be}^\bu\op \si_{\al,\be}^*(\cE^\bu_{\be,\al})^\vee))\bigr),
\nonumber
\ea
where $t^j\in H_{2j}([*/\bG_m])$ is as in \eq{co2eq3}, and $\si_{\al,\be}:\M_\al\t\M_\be\ra\M_\be\t\M_\al$ exchanges the factors, and $\La^{G\t G,G}(u\bt v)\in H_{a+b}^G(\M_\al\t\M_\be)$ is defined using \eq{co4eq21}--\eq{co4eq22}, and $c_i(\cdots)\in H^{2i}_G(\M_\al\t\M_\be)$. Using \eq{co4eq23}, for $n\in\Z$ and $\al\in K(\B)$ we write
\e
\hat H_n^G(\M_\al)= H_{n-2\chi(\al,\al)}^G(\M_\al),\quad \hat H_n^G(\M)=\bigop_{\al\in K(\B)}\hat H_n^G(\M_\al).
\label{co4eq25}
\e
That is, $\hat H_*^G(\M)$ is $H_*^G(\M)$, but with grading shifted by $-2\chi(\al,\al)$ in the component $H_*^G(\M_\al)\subset H_*^G(\M)$. The author \cite{Joyc12} proves:
\end{dfn}

\begin{thm}
\label{co4thm4}
In Definition\/ {\rm\ref{co4def7},} $(\hat H^G_*(\M),\boo,e^{zD},Y)$ is a graded vertex algebra over\/ $\Q$. In fact\/ $e^{zD},Y$ are linear over $H^*_G(*),$ so $(\hat H^G_*(\M),\boo,e^{zD},Y)$ is a vertex algebra over $H^*_G(*)$. Also $\La^{G,\{1\}}:H^G_*(\M)\ra H_*(\M)$ in Definition\/ {\rm\ref{co2def2}(h)} is a morphism of graded vertex algebras from $(\hat H^G_*(\M),\boo,e^{zD},Y)$ to $(\hat H_*(\M),\boo,e^{zD},Y)$ in Theorem\/~{\rm\ref{co4thm1}}.
\end{thm}

\begin{rem}
\label{co4rem3}
{\bf(a)} Equation \eq{co4eq24}, the crux of Definition \ref{co4def7}, is the same as \eq{co4eq10}, except for the inclusion of the Gysin morphism $\La^{G\t G,G}$ in \eq{co4eq22}. The proof of Theorem \ref{co4thm4} in \cite{Joyc12} is essentially the same as that of Theorem \ref{co4thm1}, but also including some commuting squares involving morphisms~$\La^{G\t G,G}$.
\smallskip

\noindent{\bf(b)} Section \ref{co23} discussed two kinds of $G$-equivariant homology of stacks $X$, the obvious one $H_*(X/G)$, and the non-obvious $H_*^G(X)$ with a strange definition.

Definition \ref{co4def7} and Theorem \ref{co4thm4} {\it do not work\/} if we replace $H_*^G(X)$ by $H_*(X/G)$. This is because, as in Remark \ref{co2rem7}(i), the analogue of \eq{co4eq22} is
\begin{equation*}
\io^{G\t G,G}_*:H_n((\M_\al\t\M_\be)/G)\longra H_n((\M_\al\t\M_\be)/(G\t G)),
\end{equation*}
which maps the wrong way, and cannot be substituted for $\La^{G\t G,G}$ in~\eq{co4eq24}.
\smallskip

\noindent{\bf(c)} Take $G=\bG_m$, so that $H^*_{\bG_m}(*)\cong\Q[y]$. Since $e^{zD},Y$ are $\Q[y]$-linear, the localization $\hat H^{\bG_m}_*(\M)[y^{-1}]=\hat H^{\bG_m}_*(\M)\ot_{\Q[y]}\Q[y,y^{-1}]$ is a vertex algebra over $\Q$ and over $\Q[y,y^{-1}]$. This is relevant to the $\bG_m$-localization story in~\S\ref{co26}.
\end{rem}

The material of \S\ref{co43} now extends to the $G$-equivariant case with only obvious changes. As for \eq{co4eq14} we write
\e
\check H_n^G(\M_\al^\pl)\!=\!H_{n+2-2\chi(\al,\al)}^G(\M_\al^\pl),\;\> \check H_n^G(\M^\pl)\!=\bigop_{\!\!\!\!\!\!\!\al\in K(\B)\!\!\!\!\!\!\!\!\!\!\!}\check H_n^G(\M_\al^\pl).
\label{co4eq26}
\e
The analogue of Theorem \ref{co4thm2}, proved in \cite{Joyc12}, is:

\begin{thm}
\label{co4thm5}
Work in the situation of Assumption\/ {\rm\ref{co4ass2}} and Definition\/ {\rm\ref{co4def7},} and consider the graded Lie algebra $\hat H^G_{*+2}(\M)/D(\hat H^G_*(\M))$ constructed by combining Proposition\/ {\rm\ref{co4prop1}} and Theorem\/ {\rm\ref{co4thm4}}. Then $\Pi^\pl:\M\ra\M^\pl$ gives a morphism $\Pi^\pl_*:H^G_*(\M)\ab\ra H^G_*(\M^\pl)$. It is surjective, with kernel $D(H^G_*(\M))$. This induces an isomorphism
\begin{equation*}
\Pi^\pl_*:H^G_*(\M)/D(H^G_*(\M))\longra H^G_*(\M^\pl).
\end{equation*}
Comparing \eq{co4eq25} and\/ {\rm\eq{co4eq26},} we see this is an isomorphism for all\/~{\rm$n\in\Z$:}
\e
\hat H^G_{n+2}(\M)/D(\hat H^G_n(\M))\longra \check H^G_n(\M^\pl).
\label{co4eq27}
\e
Thus there is a unique Lie bracket\/ $[\,,\,]$ on $\check H^G_*(\M^\pl)$ making it into a \begin{bfseries}graded Lie algebra\end{bfseries} over $\Q,$ such that\/ \eq{co4eq27} is a Lie algebra isomorphism. Hence $\check H^G_0(\M^\pl)$ and\/ $\check H^G_{\rm even}(\M^\pl)=\bigop_{n\in\Z}\check H^G_{2n}(\M^\pl)$ are Lie algebras. 

In fact\/ $[\,,\,]$ is bilinear over $H^*_G(*)=H^*([*/G]),$ not just over $\Q$. Also $\La^{G,\{1\}}:H^G_*(\M^\pl)\ra H_*(\M^\pl)$ in Definition\/ {\rm\ref{co2def2}(h)} is a morphism of graded Lie algebras from $\check H^G_*(\M^\pl)$ to $\check H_*(\M^\pl)$ in Theorem\/~{\rm\ref{co4thm2}}.
\end{thm}

\begin{rem}
\label{co4rem4}
As in Remark \ref{co4rem3}(c), when $G=\bG_m$, so that $H^*_{\bG_m}(*)\cong\Q[y]$, the localization $\check H^{\bG_m}_*(\M^\pl)[y^{-1}]=\check H^{\bG_m}_*(\M^\pl)\ot_{\Q[y]}\Q[y,y^{-1}]$ is a Lie algebra over $\Q$ and $\Q[y,y^{-1}]$. The Lie algebra morphism $\La^{\bG_m,\{1\}}:\check H^{\bG_m}_*(\M^\pl)\ra\check H_*(\M^\pl)$ sets $y=0$, and does not extend to~$\check H^{\bG_m}_*(\M^\pl)[y^{-1}]$.
\end{rem}

Section \ref{co44} also extends easily to the $G$-equivariant case.

\section{The main results}
\label{co5}

We can now state our main results, Theorems \ref{co5thm1}--\ref{co5thm3}, in an abstract form. These say that if a $\C$-linear abelian category $\A$, plus a lot of extra data including a family $\sS$ of weak stability conditions $(\tau,T,\le)$ on $\A$, satisfies a list of conditions given in Assumptions \ref{co5ass1}--\ref{co5ass3}, then we can define enumerative invariants $[\M_\al^\ss(\tau)]_\inv$ `counting' $\tau$-semistable objects in $\A$ in class $\al\in K(\A)$. 

Here $[\M_\al^\ss(\tau)]_\inv$ lies in the homology $H_*(\M^\pl)$ of the `projective linear' moduli stack $\M^\pl$ of objects in $\A$. If $(\tau,T,\le)$ and $(\ti\tau,\ti T,\le)$ lie in $\sS$ then we give wall-crossing formulae \eq{co5eq33}--\eq{co5eq34} which write $[\M_\al^\ss(\ti\tau)]_\inv$ in terms of the $[\M_\be^\ss(\tau)]_\inv$ for $\be\in C(\A)$, using the Lie bracket on $H_*(\M^\pl)$ defined in \S\ref{co43}. Theorem \ref{co5thm4} extends the results to equivariant homology~$H_*^G(\M^\pl)$.

Chapters \ref{co6}--\ref{co8} apply Theorems \ref{co5thm1}--\ref{co5thm3} to examples in which $\A=\modCQ$ (or $\modCQI$) is a category of representations of a quiver (with relations), or $\A=\coh(X)$ is a category of coherent sheaves on a curve, surface or Fano 3-fold $X$, or $\A$ is a category of morphisms $\rho:V\ot_\C L\ra E$ in $\coh(X)$. Some readers may prefer to look briefly at Chapters \ref{co6}--\ref{co8}, to have examples in mind, before proceeding with this chapter.

The proofs of Theorems \ref{co5thm1}--\ref{co5thm3} and \ref{co5thm4} are deferred to Chapters \ref{co9}--\ref{co11}.

\subsection{Assumptions and notation}
\label{co51}

We now state some Assumptions \ref{co5ass1}--\ref{co5ass3} on a $\C$-linear abelian category $\A$ with an exact subcategory $\B\subseteq\A$, together with some extra data.

\begin{ass}
\label{co5ass1}
{\bf(General set-up.)} Let $\A$ be a noetherian $\C$-linear abelian category, and $\B\subseteq\A$ be an exact subcategory. Assume:
\smallskip

\noindent{\bf(a)} $\B$ is {\it closed under isomorphisms in\/} $\A$ (i.e.\ if $E\cong F\in\A$ with $E\in\B$ then $F\in\B$) and {\it closed under direct summands in\/} $\A$ (i.e.\ if $E,F\in\A$ with $E\op F\in\B$ then $E,F\in\B$).
\smallskip

\noindent{\bf(b)} The inclusion $i:\B\hookra\A$ induces a morphism $i_*:K_0(\B)\ra K_0(\A)$. We should be given a surjective quotient $K_0(\A)\twoheadrightarrow K(\A)$, which we use for defining (weak) stability conditions on $\A$ as in \S\ref{co31}. We suppose that if $E\in\A$ with $\lb E\rb=0$ in $K(\A)$ then $E=0$. Write $K(\B)$ for the image of $i_*(K_0(\B))$ in $K(\A)$, so we have a commutative diagram of abelian groups:
\begin{equation*}
\xymatrix@C=100pt@R=15pt{ *+[r]{K_0(\B)} \ar[r]_{i_*} \ar@{->>}[d] & *+[l]{K_0(\A)} \ar@{->>}[d] \\
*+[r]{K(\B)\,\,} \ar@{^{(}->}[r]  & *+[l]{K(\A).\!} }	
\end{equation*}
We use the quotient $K_0(\B)\twoheadrightarrow K(\B)$ in the vertex and Lie algebra theory for $\B$ in \S\ref{co42}--\S\ref{co43}. Note in particular that $C(\B)\subseteq C(\A)\subset K(\A)$, where $C(\A)=\bigl\{\lb E\rb:0\ne E\in\A\bigr\}$ and~$C(\B)=\bigl\{\lb E\rb:0\ne E\in\B\bigr\}$.
\smallskip

\noindent{\bf(c)} Assumption \ref{co4ass1} holds for $\B$, with $K_0(\B)\twoheadrightarrow K(\B)$ and $C(\B)$ as above. We will freely use the notation $\M,\M^\pl,\M_\al,\M^\pl_\al,\cE^\bu,\chi,\ldots$ of Assumption \ref{co4ass1} and Definition \ref{co4def4}, and the Lie algebra $\check H_{\rm even}(\M^\pl)$ of Theorem \ref{co4thm2}.
\smallskip

\noindent{\bf(d)} Using the notation of \S\ref{co22}, there are also derived Artin $\C$-stacks $\bs\M,\bs\M^\pl$, which are the (projective linear) derived moduli stacks of objects in $\B$. They have classical truncations $t_0(\bs\M)=\M$ and $t_0(\bs\M^\pl)=\M^\pl$, for $\M,\M^\pl$ the classical moduli stacks in {\bf(c)}. We write $i_\M:\M\hookra\bs\M$, $i_{\M^\pl}:\M^\pl\hookra\bs\M^\pl$ for the inclusions of classical truncations. All these satisfy:
\begin{itemize}
\setlength{\itemsep}{0pt}
\setlength{\parsep}{0pt}
\item[{\bf(i)}] Since open substacks of $\bs\M,\bs\M^\pl$ correspond to open substacks of their classical truncations $\M,\M^\pl$, the splittings $\M=\coprod_{\al\in K(\B)}\M_\al$, $\M^\pl=\coprod_{\al\in K(\B)}\M_\al^\pl$ lift to $\bs\M=\coprod_{\al\in K(\B)}\bs\M_\al$, $\bs\M^\pl=\coprod_{\al\in K(\B)}\bs\M_\al^\pl$.
\item[{\bf(ii)}] The morphisms $\Phi:\M\t\M\ra\M$, $\Phi_{\al,\be}:\M_\al\t\M_\be\ra\M_{\al+\be}$, $\Psi:[*/\bG_m]\t\M\ra\M$, $\Psi_\al:[*/\bG_m]\t\M_\al\ra\M_\al$, $\Pi^\pl:\M\ra\M^\pl$ and $\Pi_\al^\pl:\M_\al\ra\M_\al^\pl$ in Assumption \ref{co4ass1} and Definition \ref{co4def4} lift to derived versions 
$\bs\Phi:\bs\M\t\bs\M\ra\bs\M$, \ldots, $\bs\Pi_\al^\pl:\bs\M_\al\ra\bs\M_\al^\pl$, which have classical truncations $\Phi=t_0(\bs\Phi)$, \ldots, $\Pi_\al^\pl=t_0(\bs\Pi_\al^\pl)$, and satisfy the same identities as $\Phi,\ldots,\Pi_\al^\pl$, so that $\bs\Phi$ is associative and commutative in $\Ho(\DArt_\C)$, the analogue of \eq{co4eq2} holds for $\bs\Phi,\bs\Psi$, and so on. The projection $\bs\Pi^\pl:\bs\M\ra\bs\M^\pl$ is a principal $[*/\bG_m]$-bundle in derived stacks on~$\bs\M\sm\{0\},\bs\M^\pl\sm\{0\}$.
\item[{\bf(iii)}] $\bs\M$ and $\bs\M^\pl$ are {\it locally finitely presented}, as in \S\ref{co22}. This implies that their cotangent complexes $\bL_{\bs\M},\bL_{\bs\M^\pl}$ are perfect.
\item[{\bf(iv)}] For each $\al\in K(\B)$ and $\cE_{\al,\al}^\bu$ as in Assumption \ref{co4ass1}(f), we are given a quasi-isomorphism of perfect complexes on $\M_\al$:
\e
\th_\al:\De_{\M_\al}^*(\cE_{\al,\al}^\bu)[-1]\,{\buildrel\cong\over\longra}\,  i_{\M_\al}^*(\bL_{\bs\M_\al}).
\label{co5eq1}
\e
Then $\rank\cE_{\al,\al}^\bu=\chi(\al,\al)$ in Assumption \ref{co4ass1}(f) and $\bs\Pi_\al^\pl:\bs\M_\al\ra\bs\M^\pl_\al$ a principal $[*/\bG_m]$-bundle for $\al\ne 0$ imply that when $\al\in C(\B)$
\e
\rank i_{\M_\al}^*(\bL_{\bs\M_\al})=-\chi(\al,\al), \quad \rank i_{\M_\al^\pl}^*(\bL_{\bs\M_\al^\pl})=1-\chi(\al,\al).
\label{co5eq2}
\e
\item[{\bf(v)}] For $\al,\be\in K(\B)$, the following should commute in $\Perf(\M_\al\t\M_\be)$:
\ea
\begin{gathered}
\!\!\!\!\!\!\!\!\!\xymatrix@!0@C=265pt@R=50pt{ *+[r]{\begin{subarray}{l}\ts (\De_{\M_{\al+\be}}\ci\Phi_{\al,\be})^* \\ \ts (\cE_{\al+\be,\al+\be}^\bu)[-1]\end{subarray}} 
\ar[d]^{\text{from \eq{co4eq3}--\eq{co4eq4}}}_\cong
\ar[r]_(0.47){ \Phi_{\al,\be}^*(\th_{\al+\be})  } & *+[l]{\Phi_{\al,\be}^*(i_{\M_{\al+\be}}^*(\bL_{\bs\M_{\al+\be}}))} \ar[dd]_(0.2){i_{\M_\al\t\M_\be}^*(\bL_{\bs\Phi_{\al,\be}})} \\
*+[r]{\begin{subarray}{l}\ts (\De_{\M_\al}\!\ci\!\Pi_{\M_\al})^*(\cE_{\al,\al}^\bu)[-1]\!\op\!  (\De_{\M_\be}\!\ci\!\Pi_{\M_\be})^*(\cE_{\be,\be}^\bu)[-1] \\
\ts {}\op \cE_{\al,\be}^\bu[-1]\op \si_{\al,\be}^*(\cE_{\be,\al}^\bu)[-1] 
\end{subarray}}
\ar[d]^(0.47){\text{project to first two factors}}
\\
*+[r]{\begin{subarray}{l}\ts (\De_{\M_\al}\!\ci\!\Pi_{\M_\al})^*(\cE_{\al,\al}^\bu)[-1]\\
\ts {}\op (\De_{\M_\be}\ci\Pi_{\M_\be})^*(\cE_{\be,\be}^\bu)[-1]
\end{subarray}} \ar[r]^(0.53){\begin{subarray}{l} \Pi_{\M_\al}^*(\th_\al)\op \\ \Pi_{\M_\be}^*(\th_\be) \end{subarray}} & *+[l]{\raisebox{30pt}{$\begin{subarray}{l}\ts i_{\M_\al\t\M_\be}^*(\bL_{\bs\M_\al\t\bs\M_\be}) \\ 
\ts {}\cong\Pi_{\M_\al}^*\ci i_{\M_\al}^*(\bL_{\bs\M_\al}) \\ \ts {}\op\Pi_{\M_\be}^*\ci i_{\M_\be}^*(\bL_{\bs\M_\be}).
\end{subarray}$}\quad\;\>}   }\!\!\!\!\!\!\!\!\!\!\!\!\!\!
\end{gathered}
\label{co5eq3}\\[-32pt]
\nonumber
\ea
\end{itemize}

\noindent{\bf(e)} We are given a subset $C(\B)_\pe\subseteq C(\B)$ of {\it permissible classes}.
\smallskip

\noindent{\bf(f)} For each $\al\in C(\B)_\pe$ we are given the following data:
\begin{itemize}
\setlength{\itemsep}{0pt}
\setlength{\parsep}{0pt}
\item[{\bf(i)}] Open $\C$-substacks $\dM_\al^\pl\subseteq\M_\al^\pl$ and $\dM_\al=(\Pi_\al^\pl)^{-1}(\dM_\al^\pl)\ab\subseteq\M_\al$, where we write $\dot\Pi_\al^\pl=\Pi_\al^\pl\vert_{\dot{\mathcal M}_\al}:\dM_\al\ra\dM^\pl_\al$. We write $\bs\dM_\al^\pl\subseteq\bs\M_\al^\pl$ and $\bs\dM_\al\subseteq\bs\M_\al$ for the corresponding derived open substacks.
\item[{\bf(ii)}] We are given a (homotopy) Cartesian square of derived Artin $\C$-stacks:
\e
\begin{gathered}
\xymatrix@C=100pt@R=15pt{ *+[r]{\bs\dM^\red_\al} \ar[r]_{\bs{\dot\Pi}_\al^\rpl} \ar[d]^{\bs j_\al} & *+[l]{\bs\dM^\rpl_\al} \ar[d]_{\bs j_\al^\pl} \\
 *+[r]{\bs\dM_\al} \ar[r]^{\bs{\dot\Pi}_\al^\pl}  & *+[l]{\bs\dM^\pl_\al,\!} }
\end{gathered}
\label{co5eq4}
\e
where $\bs\dM^{\red}_\al$ and $\bs\dM^\rpl_\al$ are locally finitely presented derived Artin stacks which are {\it quasi-smooth}, as in~\S\ref{co22}.

We require that the classical truncations $t_0(\bs j_\al),t_0(\bs j_\al^\pl)$ should be isomorphisms. Thus we may take the classical truncations to be $t_0(\bs\dM^{\red}_\al)=\dM_\al$, $t_0(\bs\dM^\rpl_\al)=\dM_\al^\pl$. Therefore $\bs\dM^\red_\al$ and $\bs\dM_\al$ are {\it two different\/} derived enhancements of the same classical stack $\dM_\al$. Note that this implies that $i_{\dM_\al}:\dM_\al\hookra\bs\dM_\al$ factors as the composition of $i_{\dM^\red_\al}:\dM_\al\hookra\bs\dM^\red_\al$ and $\bs j_\al:\bs\dM^\red_\al\ra\bs\dM_\al$. The analogues hold for $\bs\dM^\rpl_\al,\bs\dM_\al^\pl,\dM_\al^\pl$.

Since $\bs\dM^{\red}_\al,\bs\dM^\rpl_\al$ are quasi-smooth, the following are perfect obstruction theories on the classical Artin stacks $\dM_\al,\dM_\al^\pl$ in the sense of~\S\ref{co24}:
\e
\begin{split}
\bL_{i_{\dM^\red_\al}}&:i_{\dM^\red_\al}^*(\bL_{\bs\dM^\red_\al})\longra \bL_{\dM_\al},\\
\bL_{i_{\dM^\rpl_\al}}&:i_{\dM^\rpl_\al}^*(\bL_{\bs\dM^\rpl_\al})\longra \bL_{\dM_\al^\pl}.
\end{split}
\label{co5eq5}
\e

As in Remark \ref{co5rem1}(c) below, we consider $\bs\dM^\red_\al,\bs\dM^\rpl_\al$ to be `reduced' versions of $\bs\dM_\al,\bs\dM^\pl_\al$, in the sense of the `reduced obstruction theories' of \cite{KoTh1,KoTh2,MaPa,MPT,Schu,STV}, which is why we use the superscript~`${}^\red$'.
\item[{\bf(iii)}] A finite-dimensional $\C$-vector space $U_\al$ with $\dim U_\al=o_\al\ge 0$, and isomorphisms of the cotangent complexes of the morphisms $\bs j_\al,\bs j_\al^\pl$ in~\eq{co5eq4}:
\e
\bL_{\bs\dM^\red_\al/\bs\dM_\al}\cong U_\al\ot_\C\O_{\bs\dM^\red_\al}[2],\quad
\bL_{\bs\dM^\rpl_\al/\bs\dM^\pl_\al}\cong U_\al\ot_\C\O_{\bs\dM^\rpl_\al}[2].
\label{co5eq6}
\e
If $U_\al=0$, equivalently $o_\al=0$ (which happens in most of our examples) then $\bs j_\al,\bs j_\al^\pl$ in \eq{co5eq4} are isomorphisms, and we may take $\bs\dM^\red_\al=\bs\dM_\al$, $\bs\dM^\rpl_\al=\bs\dM_\al^\pl$, and $\bs j_\al,\bs j_\al^\pl$ to be the identities.

Taking \eq{co2eq11} for $\bs j_\al$, pulling back to $\dM_\al$, noting that  $t_0(\bs\dM^\red_\al)=t_0(\bs\dM_\al)\ab=\dM_\al$, and using \eq{co5eq6}, gives a distinguished triangle on $\dM_\al$:
\e
\!\!\!\!\!\!\xymatrix@C=12pt{ U_\al\ot_\C\O_{\dot{\mathcal M}_\al}[1] \ar[r] & i_{\dM_\al}^*(\bL_{\bs\dM_\al}) \ar[r] & i_{\dM^\red_\al}^*(\bL_{\bs\dM^\red_\al}) \ar[r] & U_\al\ot_\C\O_{\dot{\mathcal M}_\al}[2]. }
\label{co5eq7}
\e
The analogue works for $\bs j_\al^\pl$ on $\dM_\al^\pl$. From these and \eq{co5eq2} we deduce that for $\al\in C(\B)$ we have
\e
\begin{split}
\rank i_{\dM^\red_\al}^*(\bL_{\bs\dM^\red_\al})&=o_\al-\chi(\al,\al), \\ 
\rank i_{\dM_\al^\rpl}^*(\bL_{\bs\dM_\al^\rpl})&=o_\al+1-\chi(\al,\al).
\end{split}
\label{co5eq8}
\e
\item[{\bf(iv)}] If $\al,\be,\al+\be\in C(\B)_\pe$ we require that $\dim U_\al+\dim U_\be\ge \dim U_{\al+\be}$, that is,~$o_\al+o_\be\ge o_{\al+\be}$.
\end{itemize}

\noindent{\bf(g)} We are given a set $\bigl\{(\B_k,F_k,\la_k):k\in K\bigr\}$, where for each $k\in K$:
\begin{itemize}
\setlength{\itemsep}{0pt}
\setlength{\parsep}{0pt}
\item[{\bf(i)}] $\B_k\subseteq\B\subseteq\A$ is a full exact subcategory, closed under isomorphisms and direct summands in $\A$ as in {\bf(a)}, such that $E\in\B_k$ is an open condition on $\C$-points $[E]$ in $\M,\M^\pl$. Thus we have open $\C$-substacks $\M_k\subseteq\M$, $\M^\pl_k\subseteq\M^\pl$ which are the moduli stacks of objects in $\B_k$. We write $\M_{k,\al}=\M_\al\cap\M_k$, $\M^\pl_{k,\al}=\M_\al^\pl\cap\M^\pl_k$ for $\al\in K(\B_k)$. We write $\bs\M_k,\bs\M_k^\pl,\bs\M_{k,\al},\bs\M^\pl_{k,\al}$ for the corresponding derived moduli stacks, which are open substacks in $\bs\M,\ldots,\bs\M^\pl_\al$. 

We write $C(\B_k)=\bigl\{\lb E\rb:0\ne E\in\B_k\bigr\}\subseteq C(\B)$.
\item[{\bf(ii)}] $F_k:\A\ra\Vect_\C$ is a $\C$-linear functor. Its restriction $F_k\vert_{\B_k}:\B_k\ra\Vect_\C$ is an {\it exact functor\/} (i.e.\ it preserves short exact sequences), which extends to moduli functors, so it induces a morphism of moduli stacks
\begin{equation*}
f_k:\M_k\longra \M_{\Vect_\C}=\ts\coprod_{d\ge 0}[*/\GL(d,\C)].
\end{equation*}
There is a tautological vector bundle $\cV_{\Vect_\C}\ra\M_{\Vect_\C}$. We write $\cV_k= f_k^*(\cV_{\Vect_\C})$. Then $\cV_k\ra\M_k$ is a vector bundle with canonical isomorphisms $\cV_k\vert_{[E]}\cong F_k(E)$ for all $E\in \B_k$. We write $\cV_{k,\al}=\cV_k\vert_{\M_{k,\al}}$.

Since $F_k$ takes direct sums to direct sums, and is $\C$-linear, in a similar way to \eq{co4eq3}--\eq{co4eq6}, for all $\al,\be\in K(\B_k)$ we have isomorphisms
\ea
\Phi_{\al,\be}^*(\cV_{k,\al+\be})&\cong\cV_{k,\al}\op\cV_{k,\be},
\label{co5eq9}\\
\Psi_\al^*(\cV_{k,\al})&\cong L_{[*/\bG_m]}\ot \cV_{k,\al}.
\label{co5eq10}
\ea

All the above extends to the derived moduli stacks $\bs\M_k,\bs\M_{k,\al}$, so we have a morphism $\bs f_k:\bs\M_k\ra \M_{\Vect_\C}$ with $t_0(\bs f_k)=f_k$, and vector bundles $\bs\cV_k\ra\bs\M_k$, $\bs\cV_{k,\al}\ra\bs\M_{k,\al}$ satisfying the analogues of~\eq{co5eq9}--\eq{co5eq10}.
\item[{\bf(iii)}] If $E\in\B_k$ then $(F_k)_*:\Hom(E,E)\ra \Hom( F_k(E), F_k(E))$ is injective.
\item[{\bf(iv)}] $\la_k:K(\B_k)\ra\Z$ is a group morphism, for $K(\B_k)\subseteq K(\B)\subseteq K(\A)$, with $\dim F_k(E)=\la_k(\lb E\rb)$ for all $E\in\B_k$. Hence $\rank\cV_{k,\al}=\la_k(\al)$ for $\al\in C(\B_k)$. With {\bf(iii)}, this also implies that $\la_k(\al)>0$ if~$\al\in C(\B_k)$.
\end{itemize}
If $S\subseteq\M^\pl$ is a finite type $\C$-substack, we require that $S\subseteq\M_k^\pl$ for some~$k\in K$.
\end{ass}

Here the data $(\B_k,F_k,\la_k)$ in (g) will be used in \S\ref{co52} to define auxiliary categories $\baB\subseteq\baA$, in which we can define `pair invariants'. An example to have in mind is that if $X$ is a smooth projective $\C$-scheme and $\A=\coh(X)$, we can take $F_k:\A\ra\Vect_\C$ to map $E\mapsto H^0(E\ot \O_X(k))$ for some $k\gg 0$, and $\B_k$ to be the exact subcategory of $E\in\coh(X)$ with $H^i(E\ot \O_X(k))=0$ for $i>0$, and $\la_k(\al)=P_\al(k)$ for $P_\al$ the Hilbert polynomial of~$\al$.

\begin{ass}
\label{co5ass2}
{\bf(Weak stability conditions.)}
Let Assumption \ref{co5ass1} hold. Assume we are given a set $\sS$ of weak stability conditions on $\A$, satisfying:
\smallskip

\noindent{\bf(a)} If $(\tau,T,\le)\in\sS$ then $\A$ is $\tau$-artinian, as in Definition \ref{co3def1}. Hence each $E\in\A$ has a unique $\tau$-Harder--Narasimhan filtration $0=E_0\subsetneq E_1\subsetneq \cdots\subsetneq E_n=E$ in $\A$ by Theorem \ref{co3thm1}. We require that if $E\in\B$ with $\lb E\rb\in C(\B)_\pe$ then $E_i\in\B$ for~$i=0,\ldots,n$.
\smallskip

\noindent{\bf(b)} If $(\tau,T,\le)\in\sS$, for $E\in\B$ to be $\tau$-(semi)stable are open conditions on $[E]$ in $\M,\M^\pl$, so we have open $\C$-substacks $\M_\al^\rst(\tau)\subseteq\M_\al^\ss(\tau)\subseteq\M_\al^\pl$ parametrizing $\tau$-(semi)stable objects for all $\al\in C(\B)$. 

If $(\tau,T,\le)\in\sS$ and $\al\in C(\B)_\pe$ then $\M_\al^\rst(\tau),\M_\al^\ss(\tau)$ are of finite type.
\smallskip

\noindent{\bf(c)} If $(\tau,T,\le)\in\sS$ and $E,F\in\B$ are $\tau$-semistable with $\tau(\lb E\rb)=\tau(\lb F\rb)$ in $T$, and $\lb E\op F\rb\in C(\B)_\pe$, then $\lb E\rb,\lb F\rb\in C(\B)_\pe$.
\smallskip

\noindent{\bf(d)} Suppose $(\tau,T,\le),(\ti\tau,\ti T,\le)\in\sS$, and $I\subseteq C(\B)_\pe$ is a finite subset, and $\al\in I$, satisfy $\tau(\be)=\tau(\al)$ and $\M_\be^\ss(\tau)\ne\es$ for all $\be\in I$. Then there should exist a group morphism $\la:K(\A)\ra\R$ such that $\la(\al)=0$ and $\la(\be)>0$ (or $\la(\be)<0$) if and only if $\ti\tau(\be)>\ti\tau(\al)$ (or $\ti\tau(\be)<\ti\tau(\al)$, respectively) for all~$\be\in I$.

\smallskip

\noindent{\bf(e)} If $(\tau,T,\le)\in\sS$ and $\al\in C(\B)_\pe$ then $\M_\al^\ss(\tau)\subseteq\dM_\al^\pl$, for $\dM_\al^\pl$ as in Assumption \ref{co5ass1}(f), so the obstruction theory $\bL_{i_{\dM^\rpl_\al}}:i_{\dM^\rpl_\al}^*(\bL_{\bs\dM^\rpl_\al})\ra \bL_{\dM_\al^\pl}$ on $\dM_\al^\pl$ in \eq{co5eq5} restricts to an obstruction theory on~$\M_\al^\ss(\tau)$.
\smallskip

\noindent{\bf(f)} We are given a `rank function' $\rk:C(\A)\ra\N_{>0}=\{1,2,\ldots\}$ such that if $\al,\be\in C(\A)$ and $(\tau,T,\le)\in\sS$ with $\tau(\al)=\tau(\be)$ then~$\rk(\al+\be)=\rk(\al)+\rk(\be)$.

\smallskip

\noindent{\bf(g)} If $(\tau,T,\le)\in\sS$ and $\al\in C(\B)_\pe$ with $\M_\al^\rst(\tau)=\M_\al^\ss(\tau)$ then $\M_\al^\ss(\tau)$ is a proper algebraic space. Hence as in \S\ref{co24}, using the obstruction theory from {\bf(e)} we have a Behrend--Fantechi virtual class $[\M_\al^\ss(\tau)]_\virt,$ which by \eq{co5eq8} we regard as lying in the Betti\/ $\Q$-homology group~$H_{2+2o_\al-2\chi(\al,\al)}(\M_\al^\pl)=\check H_{2o_\al}(\M_\al^\pl)$.
\smallskip

\noindent{\bf(h)} In Definition \ref{co5def1} below, starting from $\B\subseteq\A,K(\A),\M,\M^\pl,\ldots$ above and some extra data including a quiver $Q$, we will define similar data $\baB\subseteq\baA,K(\baA),\baM,\baM^\pl,\ldots,$ and a class $\baS$ of weak stability conditions $(\bar\tau^\la_{\bs\mu},\bar T,\le)$ on $\baA$. As in {\bf(g)}, we require that if $(\bar\tau^\la_{\bs\mu},\bar T,\le)\in\baS$ and $(\al,\bs d)\in C(\baB)_\pe$ with $\baM_{(\al,\bs d)}^\rst(\bar\tau^\la_{\bs\mu})=\baM_{(\al,\bs d)}^\ss(\bar\tau^\la_{\bs\mu})$ then $\baM_{(\al,\bs d)}^\ss(\bar\tau^\la_{\bs\mu})$ is a proper algebraic space.
\end{ass}

Section \ref{co33} provides tools for proving Assumption \ref{co5ass2}(g),(h) in examples.

\begin{ass}
\label{co5ass3}
{\bf(Changes between weak stability conditions.)} Let Assumptions \ref{co5ass1}--\ref{co5ass2} hold. Assume:
\smallskip

\noindent{\bf(a)} If $(\tau,T,\le),(\ti\tau,\ti T,\le)\in\sS$ then there exists a continuous family of weak stability conditions $(\tau_t,T_t,\le)_{t\in[0,1]}$ on $\A$, in the sense of Definition \ref{co3def5}, with $(\tau_t,T_t,\le)\!\in\!\sS$ for all $t\!\in\![0,1]$, $(\tau_0,T_0,\le)\!=\!(\tau,T,\le)$, and~$(\tau_1,T_1,\le)\!=\!(\ti\tau,\ti T,\le)$.

When this holds we say $(\tau,T,\le),(\ti\tau,\ti T,\le)$ {\it are continuously connected in\/} $\sS$. This is an equivalence relation on $\sS$, and our assumption means that there is only one equivalence class.
\smallskip

\noindent{\bf(b)} Let $\al\in C(\B)_\pe$ and $(\tau_t,T_t,\le)_{t\in[0,1]}$ be a continuous family of weak stability conditions on $\A$ with $(\tau_t,T_t,\le)\in\sS$ for all $t\in[0,1]$. Then there are only finitely many sets of data $n\ge p\ge 1$, $\al_1,\ldots,\al_n\in C(\B)$ and $0=a_0<a_1<\cdots<a_p=n$ such that $\al_1+\cdots+\al_n=\al$ and, writing $\be_j=\al_{a_{j-1}+1}+\cdots+\al_{a_j}$ for $j=1,\ldots,p$, then there exist $s,t,u\in[0,1]$ with $\M_{\al_i}^\ss(\tau_s)\ne\es$ for $i=1,\ldots,n$, and $U\bigl(\al_{a_{j-1}+1},\ab\al_{a_{j-1}+2},\ab\ldots,\ab\al_{a_j};\tau_s,\tau_t)\ne 0$ in the notation of \S\ref{co32} for $j=1,\ldots,p$, and $\tau_u(\be_1)=\cdots=\tau_u(\be_p)$. For each such set we have $\al_{i_1}+\cdots+\al_{i_2}\in C(\B)_\pe$ for $1\le i_1\le i_2\le n$, so in particular~$\al_i,\be_j\in C(\B)_\pe$.

Note that taking $p=1$ and using {\bf(a)} implies that if $(\tau,T,\le),(\ti\tau,\ti T,\le)\in\sS$ and $\al\in C(\B)_\pe$, there are only finitely many decompositions $\al=\al_1+\cdots+\al_n$ with $\al_i\in C(\B)_\pe$ and $\M_{\al_i}^\ss(\tau)\ne\es$ for all $i$, and~$U\bigl(\al_1,\ldots,\al_n;\tau,\ti\tau)\ne 0$.
\end{ass}

\begin{rem}
\label{co5rem1}
{\bf(a)} Assumptions \ref{co5ass1}--\ref{co5ass3} are written to include as many interesting examples as practically possible, though at the cost of increased complexity. On a first reading it may help to take $\A=\B$, $\dM_\al=\M_\al$, $\dM_\al^\pl=\M_\al^\pl$, $U_\al=o_\al=0$, $\bs\dM^\red_\al=\bs\M_\al$, $\bs\dM^\rpl_\al=\bs\M_\al^\pl$, and $C(\B)_\pe=C(\B)$ for simplicity. The basic idea is that $\B$ is a subcategory of $\A$ on which everything we need for our theory is well behaved, though this good behaviour may fail in $\A\sm\B$. 

Similarly, we will only define invariants $[\M_\al^\ss(\tau)]_\inv$ for `permissible' classes $\al\in C(\B)_\pe\subseteq C(\B)$, so some good behaviour may fail for $\al\in C(\B)\sm C(\B)_\pe$, for example, $\M_\al^\ss(\tau)$ may not be proper if $\al\in C(\B)\sm C(\B)_\pe$ with $\M_\al^\rst(\tau)=\M_\al^\ss(\tau)$. In some examples the top morphism in \eq{co5eq5} is naturally defined on all of $\M_\al$, but $i_{\dM^\red_\al}^*(\bL_{\bs\dM^\red_\al})$ is only perfect in the interval $[-1,1]$ (required for \eq{co5eq5} to be an obstruction theory) on a smaller open substack~$\dM_\al\subsetneq\M_\al$.
\smallskip

\noindent{\bf(b)} The main reason for including the derived enhancements $\bs\M,\ab\bs\M^\pl,\ab\bs\dM^\red_\al,\ab\bs\dM^\rpl_\al$ of $\M,\ldots,\dM^\pl_\al$ is to produce the Behrend--Fantechi obstruction theories \eq{co5eq5} on $\dM_\al,\dM^\pl_\al$. Readers unfamiliar with Derived Algebraic Geometry will not lose much if they pretend a derived stack $\bs\M$ is a classical stack $\M$ together with an obstruction theory $\phi:\cF^\bu\ra\bL_\M$, except that $\cF^\bu$ is only perfect in $[-1,1]$ (one of the conditions on obstruction theories) if $\bs\M$ is quasi-smooth.

The author would have preferred to write our theory in terms of stacks with obstruction theories, without mentioning derived stacks, as this would have made it more accessible. But there is a problem with doing this. The proofs of our main results Theorems \ref{co5thm1}--\ref{co5thm3} use auxiliary abelian categories $\baA$ which roughly lie in exact sequences $0\ra\A\ra\baA\ra\modCQ\ra 0$ for $Q$ a quiver. Writing $\M,\M^\pl$ and $\baM,\baM^\pl$ for the moduli stacks of objects in $\B\subseteq\A$ and $\baB\subseteq\baA$, we have smooth morphisms $\Pi_\M:\baM\ra\M$ and $\Pi_{\M^\pl}:\baM^\pl\ra\M^\pl$.

We would like to pull back the obstruction theories on $\dM_\al,\dM^\pl_\al$ along $\Pi_\M,\Pi_{\M^\pl}$ to obstruction theories on $\dM_{(\al,\bs d)}\subseteq\baM$ and $\dM^\pl_{(\al,\bs d)}\subseteq\baM^\pl$, and so define invariants counting $\bar\tau$-semistable objects in $\baB$. But this does not work, as obstruction theories do not have functorial pullbacks along smooth morphisms. However, we can construct $\bs\baM,\ab\bs\baM^\pl,\ab\bs\dM^\red_{(\al,\bs d)},\ab\bs\dM^\rpl_{(\al,\bs d)}$ functorially as derived stacks, and so induce obstruction theories on~$\dM_{(\al,\bs d)},\dM^\pl_{(\al,\bs d)}$.

One reason why derived stacks are functorial here, but obstruction theories are not, is that to construct $\bs\baM,\ldots,\bs\dM^\rpl_{(\al,\bs d)}$ we need the derived lifts $\bs\cV_k\ra\bs\M_k$ in Assumption \ref{co5ass1}(g)(ii), and it would be messy to write the properties of these that we need in obstruction theory language.
\smallskip

\noindent{\bf(c)} We will explain the motivation for Assumption \ref{co5ass1}(f), and the `reduced' derived stacks $\bs\dM^{\red}_\al,\bs\dM^\rpl_\al$ and vector spaces $U_\al$. In most of our examples we will take $\dM_\al\!=\!\M_\al$, $\dM_\al^\pl\!=\!\M_\al^\pl$, $\bs\dM^{\red}_\al\!=\!\bs\M_\al$, $\bs\dM^\rpl_\al\!=\!\bs\M_\al^\pl$, $\bs j_\al\!=\!\bs\id\!=\!\bs j_\al^\pl$, and $U_\al=0$, so the data in Assumption \ref{co5ass1}(f) is determined by Assumption~\ref{co5ass1}(d). 

However, there are some cases in which for the natural obstruction theory $\bL_i:i^*(\bL_{\bs\M_\al})\ra \bL_{\M_\al}$ on $\M_\al$, $h^{-1}(i^*(\bL_{\bs\M_\al}))$ contains a trivial vector bundle $U_\al\ot_\C\O_{\M_\al}$ for $\dim U_\al>0$, so any virtual class $[\M_\al^\ss(\tau)]_\virt$ defined using $\bL_i$ or its analogue on $\M_\al^\pl$ is zero by Theorem \ref{co2thm1}(iv). Then our theory would yield invariants $[\M_\al^\ss(\tau)]_\inv=0$, which would be boring. If $\A=\coh(X)$ for $X$ a projective surface with geometric genus $p_g=\dim H^0(K_X)>0$, and $\rank\al>0$, this holds with~$U_\al=H^0(K_X)$. 

This problem is well known in the literature \cite{KoTh1,KoTh2,MaPa,MPT,Schu,STV}, and the solution is to define a modified obstruction theory by deleting the $U_\al$ factor in $h^{-1}(i^*(\bL_{\bs\M_\al}))$. These are known as `reduced' obstruction theories. We take our `reduced' obstruction theories to come from `reduced' derived stacks $\bs\dM^{\red}_\al,\bs\dM^\rpl_\al$. From \eq{co5eq7} we see that the obstruction theory $i^*(\bL_{\bs\dM^\red_\al})\ra \bL_{\dM_\al}$ on $\dM_\al$ is obtained by modifying $\bL_i:i^*(\bL_{\bs\dM_\al})\ra \bL_{\dM_\al}$ by deleting $U_\al\ot_\C\O_{\dM_\al}$ from $h^{-1}(\bL_{\bs\dM_\al})$. With these `reduced' obstruction theories our theory will produce invariants $[\M_\al^\ss(\tau)]_\inv$ which can be nonzero.
\end{rem}

\subsection{Auxiliary categories satisfying similar assumptions}
\label{co52}

Given $\B\subseteq\A,K(\A),\M,\M^\pl,\ldots$ satisfying Assumptions \ref{co5ass1}--\ref{co5ass2}, and some extra data including a quiver $Q$, the next long definition explains how to define related data $\baB\subseteq\baA,K(\baA),\ab\baM,\ab\baM^\pl,\ab\ldots$ which also satisfies parts of Assumptions \ref{co4ass1} and \ref{co5ass1}--\ref{co5ass2}. We will use some of these auxiliary categories $\baB\subseteq\baA$ as tools to prove Theorems \ref{co5thm1} and \ref{co5thm2} below about~$\B\subseteq\A$. 

\begin{dfn}
\label{co5def1}
Let Assumptions \ref{co5ass1}--\ref{co5ass2} hold for $\B\subseteq\A,\M,\M^\pl,K(\A),\ldots.$ Choose the following data:
\begin{itemize}
\setlength{\itemsep}{0pt}
\setlength{\parsep}{0pt}
\item[(a)] A quiver $Q=(Q_0,Q_1,h,t)$, where $Q_0,Q_1$ are finite sets of {\it vertices\/} and {\it edges\/} and $h,t:Q_1\ra Q_1$ are the {\it head\/} and {\it tail\/} maps. There should be {\it no oriented cycles\/} in $Q$, i.e.\ no loops of the form $\overset{v_0}{\bu}\,{\buildrel e_1\over \longra}\,\overset{v_1}{\bu}\,{\buildrel e_2\over \longra}\,\cdots \,{\buildrel e_n\over \longra}\,\overset{v_n=v_0}{\bu}$.
\item[(b)] A decomposition $Q_0=\dot Q_0\amalg\ddot Q_0$, with $\ddot Q_0\ne\es$, and a map $\ka:\ddot Q_0\ra K$, where $K$ is the indexing set of $\bigl\{(\B_k,F_k,\la_k):k\in K\bigr\}$ from Assumption \ref{co5ass1}(g). No edges of $Q$ should begin in $\ddot Q_0$, that is, $t(Q_1)\subseteq\dot Q_0\subseteq Q_0$. 
\end{itemize}
Here the condition in (a) that $Q$ has no oriented cycles will be needed to ensure properness of $\baM_{(\al,\bs d)}^\ss(\bar\tau^\la_{\bs\mu})$ in Assumption \ref{co5ass2}(h) above. 

When we draw diagrams of such quivers, in equations \eq{co5eq28}, \eq{co9eq7}, \eq{co10eq1}, and \eq{co10eq57} below, we will indicate vertices in $\dot Q_0$ as `$\bu$' and vertices in $\ddot Q_0$ as~`$\ci$'. 

Define $\B_\ka\subseteq\B\subseteq\A$ to be the full exact subcategory of objects $E\in\B$ with $E\in\B_{\ka(v)}$ for all $v\in\ddot Q_0$. Write $\M_\ka=\bigcap_{v\in\ddot Q_0}\M_{\ka(v)}\subseteq\M$, and $\M_\ka^\pl=\bigcap_{v\in\ddot Q_0}\M_{\ka(v)}^\pl\subseteq\M^\pl$, so that $\M_\ka,\M_\ka^\pl$ are moduli stacks of objects in $\B_\ka$, and write $\M_{\ka,\al}=\M_\ka\cap\M_\al$, $\M_{\ka,\al}^\pl=\M_\ka^\pl\cap\M_\al^\pl$ for all $\al\in C(\B_\ka)\subseteq C(\B)$. The functor $F_{\ka(v)}:\A\ra\Vect_\C$ and vector bundle $\cV_{\ka(v)}\ra \M_{\ka(v)}$ in Assumption \ref{co5ass1}(g)(ii) restrict to $\B_\ka$ and $\M_\ka$ for all~$v\in\ddot Q_0$.

Define an abelian category $\baA$ to have {\it objects\/} $(E,\bs V,\bs\rho)$, where $E\in\A$, and $\bs V=(V_v)_{v\in\dot Q_0}$ for $V_v$ a finite-dimensional $\C$-vector space for each $v\in\dot Q_0$, and if $v\in\ddot Q_0$ we write $V_v=F_{\ka(v)}(E)$, and $\bs\rho=(\rho_e)_{e\in Q_1}$ for $\rho_e:V_{t(e)}\ra V_{h(e)}$ a $\C$-linear map for each $e\in Q_1$. If $(E,\bs V,\bs\rho),(E',\bs V',\bs\rho')$ are objects in $\baA$, a {\it morphism\/} $(\th,\bs\phi):(E,\bs V,\bs\rho)\ra(E',\bs V',\bs\rho')$ consists of a morphism $\th:E\ra E'$ in $\A$, and a tuple $\bs\phi=(\phi_v)_{v\in\dot Q_0}$ for $\C$-linear maps $\phi_v:V_v\ra V_v'$ for $v\in\dot Q_0$, and if $v\in\ddot Q_0$ we write $\phi_v=F_{\ka(v)}(\th):V_v=F_{\ka(v)}(E)\ra V'_v=F_{\ka(v)}(E')$, such that $\rho_e'\ci\phi_{t(e)}=\phi'_{h(e)}\ci\rho_e:V_{t(e)}\ra V_{h(e)}'$ for all $e\in Q_1$. We define composition of morphisms, and identities $\bigl(\id_E,(\id_{V_v})_{v\in\dot Q_0}\bigr)$, in the obvious way. Since $\Vect_\C,\A$ are abelian categories, we see that $\baA$ is an abelian category.

Define $\baB\subseteq\baA$ to be the full exact subcategory of objects $(E,\bs V,\bs\rho)$ in $\baA$ with $E\in\B_\ka\subseteq\B\subseteq\A$. Define a functor $\Pi_\B:\baB\ra\B$ to map $\Pi_\B:(E,\bs V,\bs\rho)\mapsto E$ on objects and $\Pi_\B:(\th,\bs\phi)\mapsto\th$ on morphisms. We will verify parts of Assumptions \ref{co4ass1} and \ref{co5ass1}--\ref{co5ass2} for~$\baB\subseteq\baA$.

For Assumptions \ref{co4ass1}(d) and \ref{co5ass1}(b), set $K(\baA)=K(\A)\op\Z^{\dot Q_0}$. Write elements of $K(\baA)$ as $(\al,\bs d)$ for $\al\in K(\A)$ and $\bs d:\dot Q_0\ra\Z$, also written $\bs d=(\bs d(v))_{v\in\dot Q_0}$, with $\lb E,\bs V,\bs\rho\rb=(\lb E\rb,\bdim\bs V)$ for $\bdim\bs V=(\dim V_v)_{v\in\dot Q_0}$. Write $K(\baB)\subseteq K(\baA)$ for the image of the composition $K^0(\baB)\ra K^0(\baA)\twoheadrightarrow K(\baA)$.~Then
\begin{equation*}
C(\baB)=\bigl(C(\B_\ka)\t\N^{\dot Q_0}\bigr)\amalg \bigl(\{0\}\t(\N^{\dot Q_0}\sm\{0\})\bigr).
\end{equation*}

For Assumption \ref{co4ass1}(a), we will construct the classical moduli stack $\baM$ of objects in $\baB$. We write
\begin{equation*}
\baM=\coprod_{(\al,\bs d)\in C(\baB)\amalg\{(0,0)\}}\baM_{(\al,\bs d)},
\end{equation*}
with $\baM_{(\al,\bs d)}$ the moduli stack of $(E,\bs V,\bs\rho)$ in $\baB$ with $\lb E,\bs V,\bs\rho\rb=(\al,\bs d)$. To define $\baM_{(\al,\bs d)}$, first consider the Artin $\C$-stack
\e
\M_{\ka,\al}\t\ts\prod_{v\in\dot Q_0}[*/\GL(\bs d(v),\C)].
\label{co5eq11}
\e
Define vector bundles $\cV_v$ for $v\in Q_0$ over the stack \eq{co5eq11} such that $\cV_v$ is the pullback to \eq{co5eq11} of the tautological rank $\bs d(v)$ vector bundle
\e
[\bA^{\bs d(v)}/\GL(\bs d(v),\C)]\longra [*/\GL(\bs d(v),\C)]
\label{co5eq12}
\e
if $v\in\dot Q_0$ (this was written $\cV_{\Vect_\C}$ in Assumption \ref{co5ass1}(g)(ii)), and $\cV_v$ is the pullback to \eq{co5eq11} of $\cV_{\ka(v)}\vert_{\M_{\ka,\al}}$ if $v\in\ddot Q_0$, where $\cV_{\ka(v)}\ra\M_{\ka(v)}$ is as in Assumption \ref{co5ass1}(g)(ii). Now define $\baM_{(\al,\bs d)}$, as an Artin $\C$-stack, to be the total space of the vector bundle
\e
\baM_{(\al,\bs d)}=\ts\bigl(\bigop_{e\in Q_1}\cV_{t(e)}^*\ot \cV_{h(e)}\bigr)\longra \M_{\ka,\al}\t\ts\prod_{v\in\dot Q_0}[*/\GL(\bs d(v),\C)].
\label{co5eq13}
\e

It is now easy to see that $\baM_{(\al,\bs d)}$ is indeed a well-behaved moduli stack of objects $(E,\bs V,\bs\rho)$ in $\baB$ with $\lb E,\bs V,\bs\rho\rb=(\al,\bs d)$. Here $\M_{\ka,\al}$ parametrizes $E\in\B_\ka$ with $\lb E\rb=\al$, and $[*/\GL(\bs d(v),\C)]$ parametrizes $V_v$ for $v\in\dot Q_0$, and $\cV_{t(e)}^*\ot \cV_{h(e)}$ parametrizes $\rho_e$ for~$e\in Q_1$. 

As in Assumption \ref{co4ass1}(b),(c), there are natural morphisms $\bar\Phi:\baM\t\baM\ra\baM$, $\bar\Psi:[*/\bG_m]\t\baM\ra\baM$ with the required properties, which may be constructed explicitly from $\Phi\vert_{\M_\ka\t\M_\ka},\Psi\vert_{[*/\bG_m]\t\M_\ka}$ and the natural morphisms
\begin{align*}
[*/\GL(\bs d(v),\C)]\t[*/\GL(\bs e(v),\C)]&\longra[*/\GL((\bs d+\bs e)(v),\C)],\\ [*/\bG_m]\t[*/\GL(\bs d(v),\C)]&\longra [*/\GL(\bs d(v),\C)].
\end{align*}
Using $\bar\Psi$, as in Definition \ref{co4def4} we may take the quotients $\baM^\pl=\baM/[*/\bG_m]$ and $\baM_{(\al,\bs d)}^\pl=\baM_{(\al,\bs d)}/[*/\bG_m]$.

The functor $\Pi_\B:\baB\ra\B$ induces natural projections $\Pi_\M:\baM\ra\M$ and $\Pi_{\M_\al}:\baM_{(\al,\bs d)}\ra\M_\al$ on moduli stacks. The fibres of $\Pi_{\M_\al}$ are of the form $[\bA^N/\prod_{v\in\dot Q_0}\GL(\bs d(v),\C)]$, where $\bA^N$ is the fibre of the vector bundle \eq{co5eq13}. Thus the fibres of $\Pi_{\M_\al}$ are smooth Artin $\C$-stacks, so $\Pi_{\M_\al}:\baM_{(\al,\bs d)}\ra\M_\al$ is a smooth morphism, of relative dimension
\e
\dim\Pi_{\M_\al}=-\sum_{v\in\dot Q_0}\bs d(v)^2+\sum_{e\in Q_1}\bs d(t(e))\cdot \begin{cases}\bs d(h(e)), & h(e)\!\in\!\dot Q_0, \\ \la_{\ka(h(e))}(\al), & h(e)\!\in\!\ddot Q_0. \end{cases}
\label{co5eq14}
\e
As $\Pi_{\M_\al}:\baM_{(\al,\bs d)}\ra\M_\al$ is $[*/\bG_m]$-equivariant, it descends to $\Pi_{\M_\al^{\smash{\pl}}}:\baM_{(\al,\bs d)}^\pl\ra\M_\al^\pl$, which is also a smooth morphism, of dimension \eq{co5eq14}, since $\Pi_{\M_\al}$ is.

The relative tangent complex $\bT_{\baM^\pl_{(\al,\bs d)}/\M^\pl_\al}$ is a perfect complex on $\baM_{(\al,\bs d)}^\pl$ in the interval $[0,1]$, of rank \eq{co5eq14}. Define a $\Q$-linear map
\e
\begin{split}
\bar\Up_{(\al,\bs d)}&:H_*(\baM_{(\al,\bs d)}^\pl)\longra H_*(\M_\al^\pl),\\
\bar\Up_{(\al,\bs d)}&:\eta\longmapsto (\Pi_{\M_\al}^\pl)_*\bigl(\eta\cap c_{\rank \bT_{\baM^\pl_{(\al,\bs d)}/\M^\pl_\al}}(\bT_{\baM^\pl_{(\al,\bs d)}/\M^\pl_\al})\bigr).
\end{split}
\label{co5eq15}
\e
Define $\bar\Up:H_*(\baM^\pl)\ra H_*(\M^\pl)$ by $\bar\Up\vert_{H_*(\baM_{(\al,\bs d)}^\pl)}=\bar\Up_{(\al,\bs d)}$ for all $(\al,\bs d)\in C(\baB)$.

This map $\bar\Up$ will be important in the proofs of Theorems \ref{co5thm1} and \ref{co5thm2}. Although we will not use this explicitly below, one way of understanding the point of $\bar\Up$ is that there is an open substack $\baM^{\prime\pl}\subseteq\baM^\pl$ such that $H_*(\baM^{\prime\pl})$ is a Lie algebra by the construction of \S\ref{co43}, and $\Pi_{\M_\al}^\pl\vert_{\baM^{\prime\pl}}:\baM^{\prime\pl}\ra\M^\pl$ is representable, which implies that $\bT_{\Pi_{\M_\al}^\pl}\vert_{\baM^{\prime\pl}}$ is a vector bundle in degree 0 rather than just a complex in degrees $[0,1]$, and $\bar\Up\vert_{H_*(\baM^{\prime\pl})}$ is a Lie algebra morphism by the construction of \S\ref{co43}, so that $\bar\Up$ maps Lie algebra identities in $H_*(\baM^{\prime\pl})$ to Lie algebra identities in~$H_*(\M^\pl)$. 

For Assumption \ref{co5ass1}(d), we construct the derived moduli stacks $\bs\baM,\bs\baM^\pl$ of objects in $\baB$ as for $\baM,\baM^\pl$ above, but working in derived Artin stacks rather than classical Artin stacks, and using the derived stacks $\bs\M_{\ka,\al}\subseteq\bs\M_\al$ and derived vector bundles $\bs\cV_{\ka(v)}\ra\bs\M_{\ka(v)}$ from Assumption \ref{co5ass1}(d),(g)(ii) for $\B\subseteq\A$ in place of $\M_{\ka,\al}\subseteq\M_\al$ and $\cV_{\ka(v)}\ra\M_{\ka(v)}$. We also construct the morphisms $\bs{\bar\Phi}:\bs\baM\t\bs\baM\ra\bs\baM$, $\bs{\bar\Psi}:[*/\bG_m]\t\bs\baM\ra\bs\baM$, $\bs{\bar\Pi}{}^\pl:\bs\baM\ra\bs\baM^\pl$ as for the classical analogues $\bar\Phi,\bar\Psi,\bar\Pi{}^\pl$. Then Assumption \ref{co5ass1}(d)(i)--(iii) hold.

For Assumption \ref{co5ass1}(e), define
\e
C(\baB)_\pe=\bigl\{(\al,\bs d)\in C(\baB):\al\in C(\B)_\pe\amalg\{0\}\bigr\}. 
\label{co5eq16}
\e

For Assumption \ref{co5ass1}(f)(i), for $(\al,\bs d)\in C(\baB)_\pe$, if $\al\in C(\B)_\pe$ define open $\C$-substacks $\dM_{(\al,\bs d)}=\Pi_{\M_\al}^{-1}(\dM_\al)\subseteq\baM_{(\al,\bs d)}$ and $\dM_{(\al,\bs d)}^\pl=\Pi_{\M_\al^{\smash{\pl}}}^{-1}(\dM_\al^\pl)\subseteq\baM_{(\al,\bs d)}^\pl$, and if $\al=0$ define $\dM_{(0,\bs d)}=\baM_{(0,\bs d)}$ and $\dM_{(0,\bs d)}^\pl=\baM_{(0,\bs d)}^\pl$. Then $\dM_{(\al,\bs d)}=(\bar\Pi_{(\al,\bs d)}^\pl)^{-1}(\dM_{(\al,\bs d)}^\pl)$ as $\dM_\al=(\Pi_\al^\pl)^{-1}(\dM_\al^\pl)$ by Assumption \ref{co5ass1}(f)(i) for~$\B\subseteq\A$.

For Assumption \ref{co5ass1}(f)(ii), for $(\al,\bs d)\in C(\baB)_\pe$, if $\al\in C(\B)_\pe$ consider the (homotopy) commutative diagram of derived Artin $\C$-stacks:
\e
\begin{gathered}
\xymatrix@C=60pt@R=11pt{ 
*+[r]{\bs\dM^\red_{(\al,\bs d)}} \ar@{.>}[rr]_{\bs{\dot\Pi}_\al^\rpl} \ar@{.>}[dd]^{\bs{\bar\jmath}_{(\al,\bs d)}} \ar@{.>}[dr]_{\bs\Pi_{\bs\dM^\red_\al}} && *+[l]{\bs\dM^\rpl_{(\al,\bs d)}} \ar@{.>}[dd]_(0.25){\bs{\bar\jmath}_{(\al,\bs d)}^{\,\pl}} \ar@{.>}[dr]^{\bs\Pi_{\bs\dM^\rpl_\al}}
\\
& *+[r]{\bs\dM^\red_\al} \ar[rr]_(0.3){\bs{\dot\Pi}_\al^\rpl} \ar[dd]^(0.25){\bs j_\al} && *+[l]{\bs\dM^\rpl_\al} \ar[dd]_{\bs j_\al^\pl} 
\\
 *+[r]{\bs\dM_{(\al,\bs d)}} \ar[rr]^(0.25){\bs{\dot\Pi}_{(\al,\bs d)}^\pl} \ar[dr]_{\bs\Pi_{\bs\dM_\al}} && *+[l]{\bs\dM^\pl_{(\al,\bs d)}} \ar[dr]^{\bs\Pi_{\bs\dM_\al^\pl}} 
\\
& *+[r]{\bs\dM_\al} \ar[rr]^{\bs{\dot\Pi}_\al^\pl}  && *+[l]{\bs\dM^\pl_\al.\!}}
\end{gathered}
\label{co5eq17}
\e
Here $\bs\dM^\red_{(\al,\bs d)},\bs\dM^\rpl_{(\al,\bs d)}$ and the morphisms `$\dashra$' are not yet defined. Note that the rectangle of `$\ra$' is (homotopy) Cartesian by Assumption \ref{co5ass1}(f)(ii) for $\B\subseteq\A$, and the diamond of `$\ra$' is Cartesian by construction of $\bs\baM_{(\al,\bs d)},\bs\baM_{(\al,\bs d)}^\pl$.

Define $\bs\dM^\red_{(\al,\bs d)},\bs\dM^\rpl_{(\al,\bs d)}$ to be the fibre products of derived Artin stacks  such that the left hand and right hand diamonds in \eq{co5eq17} are Cartesian. Then by properties of Cartesian squares, we can fill in the final morphism $\bs{\dot\Pi}{}_\al^\rpl$ such that \eq{co5eq17} commutes, and all six faces of the cube in \eq{co5eq17} are Cartesian.

Then $\bs\dM^\red_{(\al,\bs d)},\bs\dM^\rpl_{(\al,\bs d)}$ are locally finitely presented as all the other stacks in \eq{co5eq17} are, and they are quasi-smooth as $\bs\dM^\red_\al,\bs\dM^\rpl_\al$ are by Assumption \ref{co5ass1}(f)(ii) for $\B\subseteq\A$, and the morphisms $\bs\Pi_{\bs\dM^\red_\al},\bs\Pi_{\bs\dM^\rpl_\al}$ are smooth, as $\bs\Pi_{\bs\dM_\al},\ab\bs\Pi_{\bs\dM^\pl_\al}$ are smooth, and the left and right diamonds in \eq{co5eq17} are Cartesian. 

The classical truncations of $\bs{\bar\jmath}_{(\al,\bs d)},\bs{\bar\jmath}^{\,\pl}_{(\al,\bs d)}$ are isomorphisms, as this holds for $\bs j_\al,\bs j_\al^\pl$ by Assumption \ref{co5ass1}(f)(ii) for $\B\subseteq\A$, and the left and right diamonds in \eq{co5eq17} are Cartesian, and classical truncation preserves Cartesian squares.

If $\al=0$ we define $\bs\dM^\red_{(0,\bs d)}=\bs\M_{(0,\bs d)}$ and $\bs\dM^\rpl_{(0,\bs d)}=\bs\M^\pl_{(0,\bs d)}$, and $\bs{\bar\jmath}_{(0,\bs d)},\bs{\bar\jmath}^{\,\pl}_{(0,\bs d)}$ to be the identities. Then $\bs\dM^\red_{(0,\bs d)},\bs\dM^\rpl_{(0,\bs d)}$ are smooth classical stacks, as they are moduli stacks of objects in $\modCQ$, so they are automatically quasi-smooth.

This completes Assumption \ref{co5ass1}(f)(ii) for $\baB\subseteq\baA$. Hence as for \eq{co5eq5} for $(\al,\bs d)\in C(\baB)_\pe$ the following are perfect obstruction theories on the classical Artin stacks $\dM_{(\al,\bs d)},\dM_{(\al,\bs d)}^\pl$ in the sense of~\S\ref{co24}:
\e
\begin{split}
\bL_{i_{\dM^\red_{(\al,\bs d)}}}&:i_{\dM^\red_{(\al,\bs d)}}^*(\bL_{\bs\dM^\red_{(\al,\bs d)}})\longra \bL_{\dM_{(\al,\bs d)}},\\
\bL_{i_{\dM^\rpl_{(\al,\bs d)}}}&:i_{\dM^\rpl_{(\al,\bs d)}}^*(\bL_{\bs\dM^\rpl_{(\al,\bs d)}})\longra \bL_{\dM_{(\al,\bs d)}^\pl}.
\end{split}
\label{co5eq18}
\e

From \eq{co5eq8} and \eq{co5eq14} and the Cartesian squares in \eq{co5eq17} we see that
\ea
&\rank i_{\dM^\red_{(\al,\bs d)}}^*(\bL_{\bs\dM^\red_{(\al,\bs d)}})=
\label{co5eq19}\\
&o_\al-\chi(\al,\al)-\sum_{v\in\dot Q_0}\bs d(v)^2+\sum_{e\in Q_1}\bs d(t(e))\cdot \begin{cases}\bs d(h(e)), & h(e)\!\in\!\dot Q_0, \\ \la_{\ka(h(e))}(\al), & h(e)\!\in\!\ddot Q_0, \end{cases}
\allowdisplaybreaks
\nonumber\\
&\rank i_{\dM^\rpl_{(\al,\bs d)}}^*(\bL_{\bs\dM^\rpl_{(\al,\bs d)}})=
\label{co5eq20}\\
&o_\al+1-\chi(\al,\al)-\sum_{v\in\dot Q_0}\bs d(v)^2+\sum_{e\in Q_1}\bs d(t(e))\cdot \begin{cases}\bs d(h(e)), & h(e)\!\in\!\dot Q_0, \\ \la_{\ka(h(e))}(\al), & h(e)\!\in\!\ddot Q_0. \end{cases}
\nonumber
\ea

For Assumption \ref{co5ass1}(f)(iii), for $(\al,\bs d)\in C(\baB)_\pe$, we define $\bar U_{(\al,\bs d)}=U_\al$ if $\al\in C(\B)_\pe$, and $\bar U_{(0,\bs d)}=0$ if $\al=0$. Then if $\al\in C(\B)_\pe$ we have
\begin{align*}
\bL_{\bs\dM^\red_{(\al,\bs d)}/\bs\dM_{(\al,\bs d)}}&\cong
\bs\Pi_{\bs\dM^\red_\al}^*(\bL_{\bs\dM^\red_\al/\bs\dM_\al})\cong \bs\Pi_{\bs\dM^\red_\al}^*(U_\al\ot_\C\O_{\bs\dM^\red_\al}[2])\\
&\cong \bar U_{(\al,\bs d)}\ot_\C\O_{\bs\dM^\red_{(\al,\bs d)}}[2],
\end{align*}
using the left hand diamond in \eq{co5eq17} Cartesian in the first step, the first equation of \eq{co5eq6} for $\B\subseteq\A$ in the second, and $\bar U_{(\al,\bs d)}=U_\al$ in the third. This proves the first equation of  \eq{co5eq6} for $\baB\subseteq\baA$, and the second is similar.

Assumption \ref{co5ass1}(f)(iv) for $\baB\subseteq\baA$ is immediate from the same for $\B\subseteq\A$ and the definition of $\bar U_{(\al,\bs d)}$. This completes Assumption \ref{co5ass1}(f). We will not verify the rest of Assumptions \ref{co4ass1} and \ref{co5ass1} for $\baB\subseteq\baA$, as we will not use them.

We will define a family of weak stability conditions on $\baA$. Fix $(\tau,T,\le)\in\sS$, from Assumption \ref{co5ass2} for $\A$. Define a total order $(\bar T,\le)$ by $\bar T=(T\amalg\{\pm\iy\})\t\R$, with order given by 
\begin{itemize}
\setlength{\itemsep}{0pt}
\setlength{\parsep}{0pt}
\item $(t_1,x_1)\le (t_2,x_2)$ if $t_1<t_2$ in $T\amalg\{\pm\iy\}$, or $t_1=t_2$ and $x_1\le x_2$ in~$\R$.
\item The order $\le$ on $T\amalg\{\pm\iy\}$ used above is the given order on $T$, extended by $-\iy\le t$ and $t\le\iy$ for any~$t\in T\amalg\{\pm\iy\}$.
\end{itemize}
Let $\la:K(\A)\ra\R$ be a group morphism  and let $\bs\mu=(\mu_v)_{v\in\dot Q_0}\in\R^{\dot Q_0}$. Define $\bar\tau^\la_{\bs\mu}:C(\baA)\ra \bar T$ by 
\ea
&\bar\tau^\la_{\bs\mu}(\al,\bs d)=
\label{co5eq21}\\
&\begin{cases}
\bigl(\tau(\al),[\la(\al)+\sum_{v\in\dot Q_0}\mu_vd_v]/\rk\al\bigr), & \al\ne 0, \\
\bigl(\iy,\bigl[\sum_{v\in\dot Q_0}\mu_v\bs d(v)\bigr]\big/\bigl[\sum_{v\in\dot Q_0}\bs d(v)\bigr]\bigr), & \al=0,\; \sum_{v\in\dot Q_0}\mu_v\bs d(v)>0, \\
\bigl(-\iy,\bigl[\sum_{v\in\dot Q_0}\mu_v\bs d(v)\bigr]\big/\bigl[\sum_{v\in\dot Q_0}\bs d(v)\bigr]\bigr), & \al=0,\; \sum_{v\in\dot Q_0}\mu_v\bs d(v)\le 0.
\end{cases}
\nonumber
\ea

As $(\tau,T,\le)$ is a weak stability condition, it is easy to show $(\bar\tau^\la_{\bs\mu},\bar T,\le)$ is a weak stability condition on $\baA$. As $\A$ is $\tau$-artinian, it is also easy to show that $\baA$ is $\bar\tau^\la_{\bs\mu}$-artinian, since $\bar\tau^\la_{\bs\mu}$-descending chains in $\baA$ map under $\Pi_\A$ to $\tau$-descending chains in $\A$. Write $\baS$ for the set of $(\bar\tau^\la_{\bs\mu},\bar T,\le)$ for all $(\tau,T,\le)$ in $\sS$ and $\la,\bs\mu$. Beware that $(\bar\tau^\la_{\bs\mu},\bar T,\le)$ {\it does not always depend continuously on\/} $\bs\mu$, in the sense of Definition \ref{co3def5}, as transitions between $\pm\iy$ in the second and third cases of \eq{co5eq21} may be discontinuous.

Let $(E,\bs V,\bs\rho)\in\baA$ with $E\ne 0$. If $E$ is not $\tau$-semistable it has a subobject $0\ne E'\subsetneq E$ with $\tau(\lb E'\rb)>\tau(\lb E/E'\rb)$. Then $0\ne (E',0,0)\subsetneq (E,\bs V,\bs\rho)$ with $\bar\tau^\la_{\bs\mu}(\lb E',0,0\rb)>\bar\tau^\la_{\bs\mu}(\lb E/E',\bs V,\bs\rho\rb)$, so $(E,\bs V,\bs\rho)$ is not $\bar\tau^\la_{\bs\mu}$-semistable. Conversely, $(E,\bs V,\bs\rho)$ $\bar\tau^\la_{\bs\mu}$-semistable implies that $E$ is $\tau$-semistable.

Using Assumption \ref{co5ass2} we find that for $(E,\bs V,\bs\rho)$ to be $\bar\tau^\la_{\bs\mu}$-(semi)stable is an open condition on $\C$-points $[E,\bs V,\bs\rho]$ in $\baM_{(\al,\bs d)}^\pl$. Thus we have open $\C$-substacks $\baM_{(\al,\bs d)}^\rst(\bar\tau^\la_{\bs\mu})\subseteq\baM_{(\al,\bs d)}^\ss(\bar\tau^\la_{\bs\mu})\subseteq\baM_{(\al,\bs d)}^\pl$ parametrizing $\bar\tau^\la_{\bs\mu}$-(semi)stable objects in $\baB$ in class $(\al,\bs d)$ in $C(\baB)$. Since $(E,\bs V,\bs\rho)$ $\bar\tau^\la_{\bs\mu}$-semistable implies $E$ is $\tau$-semistable, $\Pi_{\M_\al^{\smash{\pl}}}$ restricts to $\Pi_{\M_\al^\ss(\tau)}:\baM_{(\al,\bs d)}^\ss(\bar\tau^\la_{\bs\mu})\ra\M_\al^\ss(\tau)$. 

If $\al\in C(\B)_\pe$, so that $(\al,\bs d)\in C(\baB)_\pe$, then $\M_\al^\ss(\tau)\subseteq\dM_\al^\pl$ by Assumption \ref{co5ass2}(e). Hence $\baM_{(\al,\bs d)}^\ss(\bar\tau^\la_{\bs\mu})\subseteq\dM_{(\al,\bs d)}$, and the second obstruction theory in \eq{co5eq18} restricts to $\baM_{(\al,\bs d)}^\ss(\bar\tau^\la_{\bs\mu})$. If also $\baM_{(\al,\bs d)}^\rst(\bar\tau^\la_{\bs\mu})=\baM_{(\al,\bs d)}^\ss(\bar\tau^\la_{\bs\mu})$ then $\baM_{(\al,\bs d)}^\ss(\bar\tau^\la_{\bs\mu})$ is a proper algebraic space by Assumption \ref{co5ass2}(h), so as in \S\ref{co24} we have a Behrend--Fantechi virtual class
\e
[\baM_{(\al,\bs d)}^\ss(\bar\tau^\la_{\bs\mu})]_\virt\in H_{2\rank i_{\dM^\rpl_{(\al,\bs d)}}^*(\bL_{\bs\dM^\rpl_{(\al,\bs d)}})}(\baM_{(\al,\bs d)}^\pl),
\label{co5eq22}
\e
where $\rank i_{\dM^\rpl_{(\al,\bs d)}}^*(\bL_{\bs\dM^\rpl_{(\al,\bs d)}})$ is as in \eq{co5eq20}. Applying $\bar\Up_{(\al,\bs d)}$ in \eq{co5eq15} then gives
\e
\bar\Up_{(\al,\bs d)}\bigl([\baM_{(\al,\bs d)}^\ss(\bar\tau^\la_{\bs\mu})]_\virt\bigr)\in H_{2o_\al+2-2\chi(\al,\al)}(\M_\al^\pl)=\check H_{2o_\al}(\M_\al^\pl).
\label{co5eq23}
\e
The classes \eq{co5eq23} will be important in the proofs of Theorems \ref{co5thm1}--\ref{co5thm2}. Note that as $\baM_{(\al,\bs d)}^\ss(\bar\tau^\la_{\bs\mu})$ is an algebraic space, $\bT_{\baM^\pl_{(\al,\bs d)}/\M^\pl_\al}$ in \eq{co5eq15} is a vector bundle on~$\baM_{(\al,\bs d)}^\ss(\bar\tau^\la_{\bs\mu})$.

Suppose now that $(\al,\bs d)\in C(\baB)_\pe$ and $v\in\dot Q_0$ with $\bs d(v)=1$. Then we have $\baM_{(\al,\bs d)}^\pl=\baM_{(\al,\bs d)}\!\!\fatslash\,\bG_m$ with a principal $[*/\bG_m]$-bundle $\bar\Pi^\pl_{(\al,\bs d)}:\baM_{(\al,\bs d)}\ra\baM_{(\al,\bs d)}^\pl$. At a $\C$-point $[E,\bs V,\bs\rho]$ in $\baM_{(\al,\bs d)}$, quotienting $\fatslash\,\,\bG_m$ involves dividing out by a $\bG_m$ acting diagonally in $\Aut(E)$ and $\Aut(V_{v'})$ for all $v'\in\dot Q_0$. In particular, since $V_v\cong\C$ as $\bs d(v)=1$, the map $\bG_m\ra\Aut(V_v)$ is an isomorphism, so quotienting $\fatslash\,\,\bG_m$ is equivalent to fixing an isomorphism~$V_v\cong\C$. 

Because of this, $\C$-points of $\baM_{(\al,\bs d)}^\pl$ are equivalent to isomorphism classes $[(E,\bs V,\bs\rho),\io_v]$ of pairs $((E,\bs V,\bs\rho),\io_v)$ where $(E,\bs V,\bs\rho)\in\baB$ with $\lb E,\bs V,\bs\rho\rb=(\al,\bs d)$ and $\io_v:\C\ra V_v$ is an isomorphism. Thus we can define a morphism $\bar I^v_{(\al,\bs d)}:\baM_{(\al,\bs d)}^\pl\ra\baM_{(\al,\bs d)}$ mapping $[(E,\bs V,\bs\rho),\io_v]\mapsto[E,\bs V,\bs\rho]$. We have $\bar\Pi^\pl_{(\al,\bs d)}\ci\bar I^v_{(\al,\bs d)}=\id_{\baM_{(\al,\bs d)}^\pl}$. In fact there is an isomorphism $\baM_{(\al,\bs d)}\cong\baM_{(\al,\bs d)}^\pl\t[*/\bG_m]$ such that $\bar\Pi^\pl_{(\al,\bs d)}$ projects to the first factor, and $\bar I^v_{(\al,\bs d)}$ is $\id_{\bar\Pi^\pl_{(\al,\bs d)}}\t(*\ra[*/\bG_m])$. Note that $\bar I^v_{(\al,\bs d)}$ depends on the choice of $v\in\dot Q_0$ with $\bs d(v)=1$. If no such $v$ exists then in general~$\baM_{(\al,\bs d)}\not\cong\baM_{(\al,\bs d)}^\pl\t[*/\bG_m]$.

We can now define a morphism $\bar\Pi^v_{\M_{\ka,\al}}=\Pi_{\M_{\ka,\al}}\ci\bar I^v_{(\al,\bs d)}:\baM_{(\al,\bs d)}^\pl\ra\M_{\ka,\al}$, in a commutative diagram
\e
\begin{gathered}
\xymatrix@C=170pt@R=15pt{ 
*+[r]{\baM_{(\al,\bs d)}} \ar[r]_(0.3){\Pi_{\M_{\ka,\al}}} & *+[l]{\M_{\ka,\al}} \ar[d]_{\Pi^\pl_{\ka,\al}} \\
*+[r]{\baM_{(\al,\bs d)}^\pl} \ar[r]^(0.8){\Pi_{\M_{\ka,\al}^\pl}} \ar[ur]_(0.6){\bar\Pi^v_{\M_{\ka,\al}}} \ar[u]_{\bar I^v_{(\al,\bs d)}} & *+[l]{\M_{\ka,\al}^\pl.\!} }
\end{gathered}
\label{co5eq24}
\e
Generalizing \eq{co5eq13}, we identify $\baM_{(\al,\bs d)}^\pl$ with the total space of the vector bundle
\e
\baM_{(\al,\bs d)}^\pl=\ts\bigl(\bigop\limits_{e\in Q_1\!\!\!}\cV_{t(e)}^*\ot \cV_{h(e)}\bigr)\longra \M_{\ka,\al}\t\ts\prod\limits_{v'\in\dot Q_0:v'\ne v\!\!\!\!}[*/\GL(\bs d(v'),\C)],
\label{co5eq25}
\e
where $\cV_v$ is now the trivial line bundle, and $\bar\Pi^v_{\M_{\ka,\al}}$ is the projection to $\M_{\ka,\al}$. The fibres of $\Pi_{\M_\al}^v$ are of the form $[\bA^N/\prod_{v'\in\dot Q_0:v'\ne v}\GL(\bs d(v'),\C)]$, where $\bA^N$ is the fibre of \eq{co5eq25}. Hence $\bar\Pi^v_{\M_{\ka,\al}}$ is a smooth morphism (although $\bar I^v_{(\al,\bs d)}$ is not). Restricting to the preimages of $\dM_\al^\pl$, we now have a commutative diagram of smooth morphisms:
\e
\begin{gathered}
\xymatrix@C=170pt@R=15pt{ 
*+[r]{\dM_{(\al,\bs d)}} \ar[r]_(0.3){\Pi_{\dot{\mathcal M}_\al}} \ar[d]^{\dot\Pi{}^\pl_{(\al,\bs d)}} & *+[l]{\dM_\al} \ar[d]_{\dot\Pi^\pl_\al} \\
*+[r]{\dM_{(\al,\bs d)}^\pl} \ar[r]^(0.8){\Pi_{\dot{\mathcal M}_\al^\pl}} \ar[ur]_(0.6){\bar\Pi^v_{\dot{\mathcal M}_\al}}  & *+[l]{\dM_\al^\pl.\!} }
\end{gathered}
\label{co5eq26}
\e

We can also generalize \eq{co5eq24}--\eq{co5eq26} to derived moduli stacks, where for \eq{co5eq24}--\eq{co5eq25} we use the `non-reduced' versions $\bs\baM_{(\al,\bs d)},\ldots,$ and for \eq{co5eq26} we take the Cartesian product of the derived \eq{co5eq24} with $\bs\dM_\al^\rpl\ra\bs\M_\al^\pl$, and so use the `reduced' versions $\bs\dM_{(\al,\bs d)}^\red,\ldots.$ The derived version of \eq{co5eq26} is
\e
\begin{gathered}
\xymatrix@C=170pt@R=15pt{ 
*+[r]{\bs\dM_{(\al,\bs d)}^\red} \ar[r]_(0.3){\bs\Pi_{\bs\dM^\red_\al}} \ar[d]^{\bs{\dot\Pi}{}^\pl_{(\al,\bs d)}} & *+[l]{\bs\dM^\red_\al} \ar[d]_{\bs{\dot\Pi}{}^\pl_\al} \\
*+[r]{\bs\dM_{(\al,\bs d)}^\rpl} \ar[r]^(0.8){\bs\Pi_{\bs\dM_\al^\rpl}} \ar[ur]_(0.6){\bs{\bar\Pi}^v_{\bs\dM^\red_\al}}  & *+[l]{\bs\dM_\al^\rpl,\!} }
\end{gathered}
\label{co5eq27}
\e
which is again commutative with all morphisms smooth.
\end{dfn}

The next example will be used in Theorem \ref{co5thm1}. It is similar to the auxiliary categories $\B_p$ and stability conditions used in Joyce--Song \cite[\S 13.1--\S 13.2]{JoSo} to define `pair invariants' in Donaldson--Thomas theory. 

\begin{ex}
\label{co5ex1}
Let Assumptions \ref{co5ass1}--\ref{co5ass2} hold for $\B\subseteq\A,\M,\M^\pl,K(\A),\ldots.$ Fix $k\in K$, where $K$ is the indexing set of $\bigl\{(\B_k,F_k,\la_k):k\in K\bigr\}$ from Assumption \ref{co5ass1}(g). Define data $\baB\subseteq\baA,K(\baA),\ab\baM,\ab\baM^\pl,\ab\ldots$ as in Definition \ref{co5def1} using the quiver $Q$ with $\dot Q_0=\{v\}$, $\ddot Q_0=\{w\}$, $Q_1=\{e\}$, $t(e)=v$, and $h(e)=w$, and with $\ka:\ddot Q_0\ra K$ given by $\ka(w)=k$. We illustrate this by:
\e
\begin{xy}
0;<1mm,0mm>:
,(-15,3)*{e}
,(-30,3)*{v}
,(0,3)*{w}
,(10,0)*{\ka(w)=k.}
,(-30,0)*+{\bu} ; (0,0)*+{\circ} **@{-} ?>*\dir{>}
\end{xy}
\label{co5eq28}
\e

For brevity we will write objects of $\baA$ as $(E,V,\rho)$ rather than $(E,V_v,\rho_e)$, and write $K(\baA)=K(\A)\t\Z$ rather than $K(\A)\t\Z^{\{v\}}$, and write stability conditions as $\bar\tau^\la_\mu$ with $\mu\in\R$ rather than $\bs\mu=(\mu_v)$, and so on.

Suppose $(\tau,T,\le)\in\sS$ and $\al\in C(\B)_\pe$ with $\M_\al^\ss(\tau)\subseteq\M_{k,\al}^\pl$. Then Definition \ref{co5def1} defines a weak stability condition $(\bar\tau^0_1,\bar T,\le)$ on $\baA$ with $\la=0$ and $\mu_v=1$. The proofs of Theorems \ref{co5thm1}--\ref{co5thm2} will involve the moduli spaces $\baM_{(\al,1)}^\ss(\bar\tau^0_1)$. Let $(E,V,\rho)\in\baB$ with $\lb E,V,\rho\rb\ab=(\al,1)$ in $C(\baB)$. We find that $(E,V,\rho)$ is $\bar\tau^0_1$-stable if and only if it is $\bar\tau^0_1$-semistable, if and only if
\begin{itemize}
\setlength{\itemsep}{0pt}
\setlength{\parsep}{0pt}
\item[(i)] $E$ is $\tau$-semistable;
\item[(ii)] $\rho\ne 0$; and
\item[(iii)] There does not exist $0\ne E'\subsetneq E$ with $\tau(\lb E'\rb)=\tau(\lb E/E'\rb)$ and~$\rho(V)\subseteq F_k(E')\subseteq F_k(E)$.
\end{itemize}

Thus $\baM_{(\al,1)}^\rst(\bar\tau^0_1)=\baM_{(\al,1)}^\ss(\bar\tau^0_1)$, so $\baM_{(\al,1)}^\ss(\bar\tau^0_1)$ is a proper algebraic space by Assumption \ref{co5ass2}(h). Hence as in \eq{co5eq23} we may form the homology class
\ea
&\bar\Up\bigl([\baM_{(\al,1)}^\ss(\bar\tau^0_1)]_\virt\bigr)=(\Pi_{\M^\ss_\al(\tau)})_*\bigl([\baM^\ss_{(\al,1)}(\bar\tau^0_1)]_\virt\cap c_\top(\bT_{\baM^\ss_{(\al,1)}(\bar\tau^0_1)/\M^\ss_{\al}(\tau)})\bigr)
\nonumber\\
&\text{in}\qquad H_{2o_\al+2-2\chi(\al,\al)}(\M_\al^\pl)=\check H_{2o_\al}(\M_\al^\pl),
\label{co5eq29}
\ea
where $\bT_{\baM^\ss_{(\al,1)}(\bar\tau^0_1)/\M^\ss_{\al}(\tau)}$ is a vector bundle and $c_\top(\cdots)$ is its top Chern class. 
\end{ex}

\subsection{Statement of the main theorems}
\label{co53}

The next three theorems will be proved in Chapters \ref{co9}--\ref{co11} respectively. 

\begin{thm}
\label{co5thm1}
Suppose Assumptions\/ {\rm\ref{co5ass1}--\ref{co5ass2}} hold. Then for all\/ $(\tau,T,\le)\in\sS$ and\/ $\al\in C(\B)_\pe$ there are unique classes $[\M_\al^\ss(\tau)]_\inv$ in the Betti\/ $\Q$-homology group $H_{2o_\al+2-2\chi(\al,\al)}(\M_\al^\pl)=\check H_{2o_\al}(\M_\al^\pl)$ satisfying:
\begin{itemize}
\setlength{\itemsep}{0pt}
\setlength{\parsep}{0pt}
\item[{\bf(i)}] If\/ $\M_\al^\rst(\tau)=\M_\al^\ss(\tau)$ then\/ $[\M_\al^\ss(\tau)]_\inv=[\M_\al^\ss(\tau)]_\virt,$ where $[\M_\al^\ss(\tau)]_\virt$ is the Behrend--Fantechi virtual class of the proper algebraic space $\M_\al^\ss(\tau)$ with perfect obstruction theory, as in Assumption\/~{\rm\ref{co5ass2}(g)}.

Note in particular that if\/ $\M_\al^\ss(\tau)=\es$ then\/ $[\M_\al^\ss(\tau)]_\inv=0$.
\item[{\bf(ii)}] Suppose $(\tau,T,\le),(\ti\tau,\ti T,\le)\in\sS$ and\/ $\al\in C(\B)_\pe$ with\/ $\M_\al^\ss(\tau)=\M_\al^\ss(\ti\tau)$. Then $[\M_\al^\ss(\tau)]_\inv=[\M_\al^\ss(\ti\tau)]_\inv$.
\item[{\bf(iii)}] Let\/ $(\tau,T,\le)\in\sS$ and\/ $\al\in C(\B)_\pe,$ and suppose\/ $k\in K$ with\/ $\M_\al^\ss(\tau)\subseteq\M_{k,\al}^\pl\subseteq\M_\al^\pl,$ using the notation of Assumption\/ {\rm\ref{co5ass1}(g)}. Then Example\/ {\rm\ref{co5ex1}} defines a homology class\/ \eq{co5eq29} in $\check H_{2o_\al}(\M^\pl),$ which satisfies
\end{itemize}
\ea
&(\Pi_{\M^\ss_\al(\tau)})_*\bigl([\baM^\ss_{(\al,1)}(\bar\tau^0_1)]_\virt\cap c_\top(\bT_{\baM^\ss_{(\al,1)}(\bar\tau^0_1)/\M^\ss_{\al}(\tau)})\bigr)
\label{co5eq30}\\
&=\sum_{\begin{subarray}{l}n\ge 1,\;\al_1,\ldots,\al_n\in
C(\B)_\pe:\\ \al_1+\cdots+\al_n=\al,\; o_{\al_1}+\cdots+o_{\al_n}=o_\al,\\ 
\tau(\al_i)=\tau(\al), \; \M_{\al_i}^\ss(\tau)\ne\es,\; \text{all\/ $i$}\end{subarray}} \!\!\!\!\!\!\!\!\!\!\!\!\!\!\!\!\!\!\!\!\!\!\!
\begin{aligned}[t]
\frac{(-1)^{n+1}\la_k(\al_1)}{n!}\,\cdot&\bigl[\bigl[\cdots\bigl[[\M_{\al_1}^\ss(\tau)]_\inv,\\
&\; 
[\M_{\al_2}^\ss(\tau)]_\inv\bigr],\ldots\bigr],[\M_{\al_n}^\ss(\tau)]_\inv\bigr],
\end{aligned}
\nonumber
\ea	
\begin{itemize}
\setlength{\itemsep}{0pt}
\setlength{\parsep}{0pt}
\item[] where the Lie brackets are in the Lie algebra $\check H_{\rm even}(\M^\pl)$ from Theorem\/ {\rm\ref{co4thm2},} and there are only finitely many terms in the sum.
\end{itemize}

We think of the $[\M_\al^\ss(\tau)]_\inv$ as \begin{bfseries}enumerative invariants\end{bfseries} which `count' the moduli spaces $\M_\al^\ss(\tau)$ for all\/ $\al\in C(\B)_\pe$ and\/ $(\tau,T,\le)\in\sS$. They have the important property that\/ $[\M_\al^\ss(\tau)]_\inv$ is well defined even when $\M_\al^\rst(\tau)\ne\M_\al^\ss(\tau),$ so that the Behrend--Fantechi virtual class $[\M_\al^\ss(\tau)]_\virt$ is not defined.

The invariants $[\M_\be^\ss(\tau)]_\inv$ are determined uniquely by the virtual classes $[\baM^\ss_{(\al,1)}(\bar\tau^0_1)]_\virt$ by \eq{co5eq30} and an inductive argument. In particular, this means that any \begin{bfseries}deformation invariance\end{bfseries} properties of\/ $[\baM^\ss_{(\al,1)}(\bar\tau^0_1)]_\virt$ in Theorem\/ {\rm\ref{co2thm1}(iv)} also hold for the $[\M_\be^\ss(\tau)]_\inv$. So, for example, if\/ $\A=\coh(X)$ for a smooth projective $\C$-scheme $X$ then (suitably interpreted, e.g.\ by evaluating $[\M_\be^\ss(\tau)]_\inv$ on deformation-independent universal cohomology classes) the $[\M_\be^\ss(\tau)]_\inv$ will be unchanged by continuous deformations of\/~$X$.
\end{thm}

\begin{thm}
\label{co5thm2}
Let Assumptions\/ {\rm\ref{co5ass1}--\ref{co5ass2}} hold, and\/ $(\tau,T,\le),(\ti\tau,\ti T,\le)$ lie in\/ $\sS,$ and\/ $\al\in C(\B)_\pe$. Suppose that if\/ $\al_1,\ldots,\al_n\in C(\B)$ with\/ $\al=\al_1+\cdots+\al_n,$ and either 
\begin{itemize}
\setlength{\itemsep}{0pt}
\setlength{\parsep}{0pt}
\item[{\bf(i)}] $U(\al_1,\ldots,\al_n;\tau,\ti\tau)\ne 0$ and\/ $\M_{\al_i}^\ss(\tau)\ne\es$ for\/ {\rm $i=1,\ldots,n$;} or 
\item[{\bf(ii)}] $U(\al_1,\ldots,\al_n;\ti\tau,\tau)\ne 0$ and\/ $\M_{\al_i}^\ss(\ti\tau)\ne\es$ for\/ $i=1,\ldots,n,$
\end{itemize}
then\/ $\tau(\al_i)=\tau(\al)$ for all\/ $i=1,\ldots,n,$ and in case\/ {\bf(ii)} we also have $\M_{\al_i}^\ss(\ti\tau)\subseteq\M_{\al_i}^\ss(\tau)$ for\/ $i=1,\ldots,n$.

Then the\/ $[\M_\be^\ss(\tau)]_\inv,[\M_\be^\ss(\ti\tau)]_\inv$ from Theorem\/ {\rm\ref{co5thm1}} satisfy
\e
\!\!\begin{gathered}
{}
[\M_\al^\ss(\ti\tau)]_\inv= \!\!\!\!\!\!\!
\sum_{\begin{subarray}{l}n\ge 1,\;\al_1,\ldots,\al_n\in
C(\B)_\pe:\\ 
\tau(\al_i)=\tau(\al), \; \M_{\al_i}^\ss(\tau)\ne\es,\; \text{all\/ $i,$}\\ 
\al_1+\cdots+\al_n=\al,\; o_{\al_1}+\cdots+o_{\al_n}=o_\al\end{subarray}} \!\!\!\!\!\!\!\!\!\!\!\!\!\!\!\!\!\begin{aligned}[t]
\ti U(\al_1,&\ldots,\al_n;\tau,\ti\tau)\cdot\bigl[\bigl[\cdots\bigl[[\M_{\al_1}^\ss(\tau)]_\inv,\\
&\;\;
[\M_{\al_2}^\ss(\tau)]_\inv\bigr],\ldots\bigr],[\M_{\al_n}^\ss(\tau)]_\inv\bigr]
\end{aligned}
\end{gathered}\!\!
\label{co5eq31}
\e
in the Lie algebra\/ $\check H_{\rm even}(\M^\pl)$ from\/ {\rm\S\ref{co43}}. Here\/ $\ti U(-;\tau,\ti\tau)$ is as in Theorem\/ {\rm\ref{co3thm3},} and there are only finitely many terms in \eq{co5eq31}. Equivalently, by Theorem\/ {\rm\ref{co3thm3},} in the universal enveloping algebra $\bigl(U(\check H_{\rm even}(\M^\pl)),*\bigr)$ we have:
\e
\!\!\begin{gathered}
{}
[\M_\al^\ss(\ti\tau)]_\inv= \!\!\!\!\!\!\!
\sum_{\begin{subarray}{l}n\ge 1,\;\al_1,\ldots,\al_n\in
C(\B)_\pe:\\ 
\tau(\al_i)=\tau(\al), \; \M_{\al_i}^\ss(\tau)\ne\es,\; \text{all\/ $i,$}\\ 
\al_1+\cdots+\al_n=\al,\; o_{\al_1}+\cdots+o_{\al_n}=o_\al\end{subarray}} \!\!\!\!\!\!\!\!\!\!\!\!\!\!\!\!
\begin{aligned}[t]
U(\al_1,&\ldots,\al_n;\tau,\ti\tau)\cdot[\M_{\al_1}^\ss(\tau)]_\inv *{}\\
&
[\M_{\al_2}^\ss(\tau)]_\inv*\cdots *[\M_{\al_n}^\ss(\tau)]_\inv.
\end{aligned}
\end{gathered}\!\!
\label{co5eq32}
\e

We call\/ {\rm\eq{co5eq31}--\eq{co5eq32}} \begin{bfseries}dominant wall-crossing formulae\end{bfseries} for the invariants $[\M_\al^\ss(\tau)]_\inv,$ as if\/ $(\tau,T,\le)$ dominates\/ $(\ti\tau,\ti T,\le)$ then the conditions above involving {\bf(i)\rm,\bf(ii)} hold automatically.
\end{thm}

\begin{thm}
\label{co5thm3}
Let Assumptions\/ {\rm\ref{co5ass1}--\ref{co5ass3}} hold. Then for\/ $(\tau,T,\le),(\ti\tau,\ti T,\le)$ in $\sS$ and\/ $\al$ in $C(\B)_\pe,$ the invariants\/ $[\M_\be^\ss(\tau)]_\inv$ from Theorem\/ {\rm\ref{co5thm1}} satisfy
\e
\begin{gathered}
{}
[\M_\al^\ss(\ti\tau)]_\inv= \!\!\!\!\!\!\!
\sum_{\begin{subarray}{l}n\ge 1,\;\al_1,\ldots,\al_n\in
C(\B)_\pe:\\ \M_{\al_i}^\ss(\tau)\ne\es,\; \text{all\/ $i,$} \\
\al_1+\cdots+\al_n=\al, \; o_{\al_1}+\cdots+o_{\al_n}=o_\al 
\end{subarray}} \!\!\!\!\!\!\!\!\!\!\!\!\!\!\!\!\!\!\!\!\!\!\!\!\begin{aligned}[t]
\ti U(\al_1,&\ldots,\al_n;\tau,\ti\tau)\cdot\bigl[\bigl[\cdots\bigl[[\M_{\al_1}^\ss(\tau)]_\inv,\\
&
[\M_{\al_2}^\ss(\tau)]_\inv\bigr],\ldots\bigr],[\M_{\al_n}^\ss(\tau)]_\inv\bigr]
\end{aligned}
\end{gathered}
\label{co5eq33}
\e
in the Lie algebra $\check H_{\rm even}(\M^\pl)$ from {\rm\S\ref{co43}}. Here\/ $\ti U(-;\tau,\ti\tau)$ is as in Theorem\/ {\rm\ref{co3thm3},} and there are only finitely many nonzero terms in \eq{co5eq33} by Assumption\/ {\rm\ref{co5ass3}(b)}. Equivalently, in the universal enveloping algebra $\bigl(U(\check H_{\rm even}(\M^\pl)),*\bigr)$ we have:
\e
\begin{gathered}
{}
[\M_\al^\ss(\ti\tau)]_\inv= \!\!\!\!\!\!\!
\sum_{\begin{subarray}{l}n\ge 1,\;\al_1,\ldots,\al_n\in
C(\B)_\pe:\\ \M_{\al_i}^\ss(\tau)\ne\es,\; \text{all\/ $i,$} \\
\al_1+\cdots+\al_n=\al, \; o_{\al_1}+\cdots+o_{\al_n}=o_\al \end{subarray}} \!\!\!\!\!\!\!\!\!\!\!\!\!\!\!\!\!\!\!
\begin{aligned}[t]
U(\al_1,&\ldots,\al_n;\tau,\ti\tau)\cdot[\M_{\al_1}^\ss(\tau)]_\inv *{}\\
&
[\M_{\al_2}^\ss(\tau)]_\inv*\cdots *[\M_{\al_n}^\ss(\tau)]_\inv.
\end{aligned}
\end{gathered}
\label{co5eq34}
\e

We call\/ {\rm\eq{co5eq33}--\eq{co5eq34}} \begin{bfseries}(general) wall-crossing formulae\end{bfseries} for the invariants $[\M_\al^\ss(\tau)]_\inv$. Note that\/ {\rm\eq{co5eq31}--\eq{co5eq32}} are special cases of\/~{\rm\eq{co5eq33}--\eq{co5eq34}}.
\end{thm} 

\begin{rem}
\label{co5rem2}
{\bf(a)} Theorems \ref{co5thm1}--\ref{co5thm3} prove versions of the conjectures of Gross--Joyce--Tanaka \cite[\S 4.1]{GJT} for categories $\B\subseteq\A$ satisfying Assumptions \ref{co5ass1}--\ref{co5ass3}.

\smallskip

\noindent{\bf(b)} Theorem \ref{co5thm2} is a special case of Theorem \ref{co5thm3}, under weaker assumptions, and is used to prove Theorem \ref{co5thm3}. Theorem \ref{co5thm2} is not that interesting on its~own.

\smallskip

\noindent{\bf(c)} Theorem \ref{co5thm1}(iii) is a version of the `method of pair invariants' for counting strictly semistable moduli spaces, as used by Joyce--Song \cite[\S 5.4 \& \S 13.1]{JoSo}, Mochizuki \cite[\S 7.3]{Moch}, and Tanaka--Thomas~\cite{TaTh2}.
\smallskip

\noindent{\bf(d)} To construct the invariants $[\M_\al^\ss(\tau)]_\inv$ in Theorem \ref{co5thm1}, we need there to exist $\{(\B_k,F_k,\la_k):k\in K\}$ satisfying Assumption \ref{co5ass1}(g). However, the $[\M_\al^\ss(\tau)]_\inv$ are {\it independent of the choice of\/} $\{(\B_k,F_k,\la_k):k\in K\}$. To see this, let $\{(\B'_{k'},F'_{k'},\la'_{k'}):k'\in K'\}$ be an alternative choice giving invariants $[\M_\al^\ss(\tau)]_\inv'$. Then the disjoint union $\{(\B_k,F_k,\la_k):k\in K\}\amalg\{(\B'_{k'},F'_{k'},\la'_{k'}):k'\in K'\}$ also satisfies Assumption \ref{co5ass1}(g). Let $[\M_\al^\ss(\tau)]_\inv''$ be the invariants defined in Theorem \ref{co5thm1} using this disjoint union. Then using uniqueness in Theorem \ref{co5thm1}(iii) for suitable $k\in K$ shows that $[\M_\al^\ss(\tau)]_\inv=[\M_\al^\ss(\tau)]_\inv''$, and for $k'\in K'$ shows that $[\M_\al^\ss(\tau)]_\inv'=[\M_\al^\ss(\tau)]_\inv''$. Hence $[\M_\al^\ss(\tau)]_\inv=[\M_\al^\ss(\tau)]_\inv'$, and the $[\M_\al^\ss(\tau)]_\inv$ are independent of~$\{(\B_k,F_k,\la_k):k\in K\}$.
\smallskip

\noindent{\bf(e)} Wall-crossing formulae like \eq{co5eq33}--\eq{co5eq34} can be powerful tools for studying enumerative invariants, as we hope to show in the sequels to this book.

If we were only able to define $[\M_\al^\ss(\tau)]_\inv$ when $\M_\al^\rst(\tau)=\M_\al^\ss(\tau)$ then \eq{co5eq33}--\eq{co5eq34} would only make sense in rare special cases, since even if $\M_\al^\rst(\ti\tau)=\M_\al^\ss(\ti\tau)$ so that the left hand sides of \eq{co5eq33}--\eq{co5eq34} make sense, it is likely that some term on the right hand sides would involve $\al_i$ with $\M_{\al_i}^\rst(\tau)\ne\M_{\al_i}^\ss(\tau)$. Because of this, extending the definition of $[\M_\al^\ss(\tau)]_\inv$ to all $\al\in C(\B)_\pe$, as in Theorem \ref{co5thm1}, is an important prerequisite for the wall-crossing in Theorem~\ref{co5thm3}.  

\smallskip

\noindent{\bf(f)} There may be {\it many different ways\/} to extend virtual classes $[\M_\al^\ss(\tau)]_\virt$ when $\M_\al^\rst(\tau)=\M_\al^\ss(\tau)$ to invariants $[\M_\al^\ss(\tau)]_\inv$ for all $\al\in C(\B)_\pe$. For example, Mochizuki's invariants counting coherent sheaves on surfaces \cite[\S 7]{Moch} do this, but we show in \S\ref{co771}(A) that they are different from ours. 

For categories $\B$ for which $\M^\pl$ is a smooth Artin stack, we could use Kirwan's algorithm \cite{Kirw} when $\M_\al^\rst(\tau)\ne\M_\al^\ss(\tau)$, applying repeated blow ups and deletions to $\M_\al^\ss(\tau)$ until it becomes a smooth, proper Deligne--Mumford stack $\haM_\al^\ss(\tau)$, and then define $[\M_\al^\ss(\tau)]_\inv=[\haM_\al^\ss(\tau)]_\fund$ to be its fundamental class. (See Edidin--More--Rydh \cite{EdMo,EdRy} for discussion, and progress in extending this algorithm to non-smooth stacks.) Again, calculation of examples by the author shows that this gives a different answer to Theorem~\ref{co5thm1}.

The main defining feature of our invariants $[\M_\al^\ss(\tau)]_\inv$ is that they satisfy the wall-crossing formulae \eq{co5eq33}--\eq{co5eq34}. This is ensured by defining them using \eq{co5eq30}, which is derived from a special case of~\eq{co5eq33}.
\smallskip

\noindent{\bf(g)} Note that Assumption \ref{co5ass1}(f)(iv) says that $o_\al+o_\be\ge o_{\al+\be}$ if $\al,\be,\al+\be\in C(\B)_\pe$, and \eq{co5eq31}--\eq{co5eq34} involve sums over $\al_1,\ldots,\al_n$ with $\al_1+\cdots+\al_n=\al$ and $o_{\al_1}+\cdots+o_{\al_n}=o_\al$. In most of our examples we have $o_\al=0$ for all $\al\in C(\B)_\pe$, so $o_{\al_1}+\cdots+o_{\al_n}=o_\al$ holds automatically. But when $o_\al\ne 0$ for some $\al$, the condition $o_{\al_1}+\cdots+o_{\al_n}=o_\al$ has a {\it huge effect\/} in \eq{co5eq31}--\eq{co5eq34}. For example, if $\A=\B=\coh(X)$ for a projective surface $X$ with geometric genus $p_g>0$, we will take $o_\al=p_g$ for all $\al\in C(\B)_\pe$ with $\rank\al>0$. Then $o_{\al_1}+\cdots+o_{\al_n}=o_\al$ means the sums \eq{co5eq31}--\eq{co5eq34} have {\it only one term}, with $n=1$ and~$\al_1=\al$.

The restriction to $o_{\al_1}+\cdots+o_{\al_n}=o_\al$ occurs because in Chapters \ref{co9}--\ref{co10}, we will use Assumption \ref{co5ass1}(d),(f) and Theorem \ref{co2thm1}(iv) to show that if $o_{\al_1}+o_{\al_2}>o_{\al_1+\al_2}$ then contributions to wall-crossing formulae from classes $\al_1,\al_2$ are zero.
\end{rem}

\subsection{Extension to equivariant (co)homology}
\label{co54}

In \S\ref{co23} we defined the equivariant homology $H_*^G(X)$ of an Artin $\C$-stack with an action of a linear algebraic $\C$-group $G$, and in \S\ref{co45} we generalized the constructions of vertex algebras $\hat H_*(\M)$ and Lie algebras $\check H_*(\M^\pl)$ to equivariant homology $\hat H_*^G(\M)$, $\check H_*^G(\M^\pl)$. We now explain how to generalize the material of \S\ref{co51}--\S\ref{co53} to $G$-equivariant homology. Here are our additional assumptions:

\begin{ass}
\label{co5ass4}
Let Assumption \ref{co5ass1} hold for $\B\subseteq\A,\M,\M^\pl,\ab K(\A),\ab\ldots.$ This includes  Assumption \ref{co4ass1} as in Assumption \ref{co5ass1}(c). Suppose Assumption \ref{co4ass2} holds, extending Assumption \ref{co4ass1} to the $G$-invariant case, and:
\smallskip

\noindent{\bf(a)} We are given an action of the group $G(\C)$ of $\C$-points of $G$ on the $\C$-linear abelian category $\A$, in the sense of Assumption \ref{co4ass2}(a)(i)--(iii). The subcategory $\B\subseteq\A$ is $G$-invariant, and the $G$-action on $\A$ restricts to the $G$-action on $\B$ given in Assumption \ref{co4ass2}(a). As in Assumption \ref{co4ass2}(e), classes $\lb E\rb\in K(\A)$ are invariant under the action of $G$ on~$E\in\A$.
\smallskip

\noindent{\bf(b)} The $G$-actions on $\M,\M^\pl$ and the $G$-equivariance of $\Phi:\M\t\M\ra\M$, $\Psi:[*/\bG_m]\t\M\ra\M$ and $\Pi^\pl:\M\ra\M^\pl$ in Assumption \ref{co4ass2}(b)--(d) lift to the derived versions $\bs\M,\bs\M^\pl$, $\bs\Phi:\bs\M\t\bs\M\ra\bs\M$, $\bs\Psi:[*/\bG_m]\t\bs\M\ra\bs\M$ and $\bs\Pi^\pl:\bs\M\ra\bs \M^\pl$ in Assumption \ref{co5ass1}(d). Equations \eq{co5eq1} and \eq{co5eq3} in Assumption \ref{co5ass1}(d)(iv),(v) should hold in $G$-equivariant perfect complexes.
\smallskip

\noindent{\bf(c)} In Assumption \ref{co5ass1}(f), $\dM_\al\subseteq\M_\al$ and $\dM_\al^\pl\subseteq\M_\al^\pl$ are $G$-invariant open substacks, so that $G$ acts on $\dM_\al,\dM_\al^\pl$. Thus $\bs\dM_\al\subseteq\bs\M_\al$ and $\bs\dM_\al^\pl\subseteq\bs\M_\al^\pl$ are also $G$-invariant derived open substacks. The $G$-actions on $\dM_\al,\dM_\al^\pl$ should lift to $G$-actions on $\bs\dM_\al^\red$ and $\bs\dM_\al^\rpl$ such that \eq{co5eq4} is a $G$-equivariant diagram of derived stacks. This implies that \eq{co5eq5} are $G$-equivariant obstruction theories on $\dM_\al,\dM_\al^\pl$. We are given a representation of $G$ on $U_\al$ such that \eq{co5eq6} holds in $G$-equivariant perfect complexes. Hence \eq{co5eq7} also holds $G$-equivariantly.
\smallskip

\noindent{\bf(d)} In Assumption \ref{co5ass1}(g), for each $k\in K$, we are given a $G(\C)$-equivariant structure on the functor $F_k:\A\ra\Vect_\C$, where $G(\C)$ acts on $\A$ as in {\bf(a)}, and trivially on $\Vect_\C$. In the notation of Assumption \ref{co4ass2}(a), this means:
\begin{itemize}
\setlength{\itemsep}{0pt}
\setlength{\parsep}{0pt}
\item[{\bf(i)}] For all $g\in G(\C)$, we are given a natural isomorphism $\Ga^k_g:F_k\ci \Ga(g)\Ra F_k$ such that for all $g,h\in G$ the following commutes:
\begin{equation*}
\xymatrix@C=140pt@R=15pt{ *+[r]{F_k\ci \Ga(g)\ci \Ga(h)} \ar@{=>}[r]_{\Ga^k_g*\id_{\Ga(h)}} \ar@{=>}[d]^{\id_{F_k} *\Ga_{g,h}} & *+[l]{F_k\ci \Ga(h)}  \ar@{=>}[d]_{\Ga^k_h} \\
*+[r]{F_k\ci \Ga(g\ci h)} \ar@{=>}[r]^{\Ga^k_{g\ci h}} & *+[l]{F_k.\!} }	
\end{equation*}
\item[{\bf(ii)}] $\Ga^k_1=\id_{F_k}*\Ga_1:F_k\ci \Ga(1)\Longra F_k$.
\end{itemize}

The $G$-action on $\B$ preserves the subcategory $\B_k\subseteq\B$. Hence the open substacks $\M_k\subseteq\M$, $\M_k^\pl\subseteq\M^\pl$ are $G$-invariant. On $\B_k$, the $G(\C)$-equivariant structure on $F_k$ above is algebraic, and compatible with moduli functors. So the $G$-action on $\M_k$ lifts to the vector bundle $\cV_k\ra\M_k$, preserving the vector bundle structure. Equations \eq{co5eq9}--\eq{co5eq10} hold in $G$-equivariant vector bundles.

The above also extends to the derived moduli stacks $\bs\M_k,\bs\M_{k,\al}$, so that $\bs\cV_k\ra\bs\M_k$, $\bs\cV_{k,\al}\ra\bs\M_{k,\al}$ are $G$-equivariant vector bundles satisfying the analogues of \eq{co5eq9}--\eq{co5eq10} $G$-equivariantly.
\end{ass}

We explain how to extend the rest of \S\ref{co51}--\S\ref{co52} to the equivariant case:

\begin{dfn}
\label{co5def2}
Suppose Assumption \ref{co5ass4} holds. This includes Assumption \ref{co4ass2}, and generalizes Assumptions \ref{co4ass1} and \ref{co5ass1} to the $G$-equivariant case. Suppose also that Assumption \ref{co5ass2} holds. 

If $(\tau,T,\le)\in\sS$, Assumptions \ref{co4ass2}(e) and \ref{co5ass4}(a) imply that $\tau$-(semi)-stability is $G$-invariant, so the open substacks $\M_\al^\rst(\tau)\subseteq\M_\al^\ss(\tau)\subseteq\M^\pl_\al$ are $G$-invariant for $\al\in C(\B)$. If $\al\in C(\B)_\pe$ with $\M_\al^\rst(\tau)=\M_\al^\ss(\tau)$ then Assumption \ref{co5ass2}(g) says $\M_\al^\ss(\tau)$ is a proper algebraic space, which has a $G$-action. 

As $\M_\al^\ss(\tau)\subseteq\dM_\al^\pl$ by Assumption \ref{co5ass2}(f), the obstruction theory $\bL_{i_{\dM^\rpl_\al}}:i_{\dM^\rpl_\al}^*(\bL_{\bs\dM^\rpl_\al})\ra \bL_{\dM_\al^\pl}$ on $\dM_\al^\pl$ in \eq{co5eq5}, which is $G$-equivariant by Assumption \ref{co5ass4}(c), restricts to a $G$-equivariant obstruction theory on $\M_\al^\ss(\tau)$. Hence $[\M_\al^\ss(\tau)]_\virt$ in Assumption \ref{co5ass2}(g) is defined in $G$-equivariant homology by Theorem \ref{co2thm1}, and we take $[\M_\al^\ss(\tau)]_\virt\in\check H_{2o_\al}^G(\M_\al^\pl)$. As $\tau$-(semi)-stability is $G$-invariant, and $[\M_\al^\ss(\tau)]_\virt$ exists in $\check H_{2o_\al}^G(\M_\al^\pl)$, Assumptions \ref{co5ass2}--\ref{co5ass3} need {\it no changes or additions\/} for the $G$-equivariant case. 

To extend \S\ref{co52} to the $G$-equivariant case, the important point is to extend the $G$-actions to the moduli spaces $\dM_{(\al,\bs d)}\subseteq\baM_{(\al,\bs d)},\dM_{(\al,\bs d)}^\pl\subseteq\baM_{(\al,\bs d)}^\pl$ and their derived versions $\bs\dM_{(\al,\bs d)}\subseteq\bs\baM_{(\al,\bs d)}$, $\bs\dM_{(\al,\bs d)}^\pl\subseteq\bs\baM_{(\al,\bs d)}^\pl$ and $\bs\dM_{(\al,\bs d)}^\red,\bs\dM_{(\al,\bs d)}^\rpl$ in Definition \ref{co5def1}. Now $\M_{\ka,\al}=\M_\al\cap\bigcap_{v\in\ddot Q_0}\M_{\ka(v)}$ is a $G$-invariant open substack of $\M$, as $\M_\al,\M_k$ are. Take $G$ to act on the Artin stack \eq{co5eq11} by this action on $\M_{\ka,\al}$, and the trivial action on $[*/\GL(\bs d(v),\C)]$. Lift this action to the vector bundles $\cV_v$ on \eq{co5eq11} for $v\in Q_0$ by the pullback of the trivial $G$-action on \eq{co5eq12} if $v\in\dot Q_0$, and the lift of the $G$-action on $\cV_{\ka(v)}\ra\M_{\ka(v)}$ in Assumption \ref{co5ass4}(d) if $v\in\ddot Q_0$. Then $G$ acts on every term in \eq{co5eq13}, so $G$ acts on $\baM_{(\al,\bs d)}$ for each $(\al,\bs d)$, and thus on~$\baM$.

The morphisms $\bar\Phi:\baM\t\baM\ra\baM$, $\bar\Psi:[*/\bG_m]\t\baM\ra\baM$ in Definition \ref{co5def1} are constructed from all $G$-equivariant ingredients, including $\Phi,\Psi$ in Assumption \ref{co4ass2}(c),(d), and \eq{co5eq9}--\eq{co5eq10}, which hold in $G$-equivariant vector bundles by Assumption \ref{co5ass4}(d). Thus $\bar\Phi,\bar\Psi$ are $G$-equivariant. Since $\baM^\pl$ is the quotient of $\baM$ by the $[*/\bG_m]$-action $\bar\Psi$, the $G$-action on $\baM$ descends to $\baM^\pl$. The projections $\Pi_{\M_\al}:\baM_{(\al,\bs d)}\ra\M_\al$, $\Pi_{\M_\al^\pl}:\baM_{(\al,\bs d)}^\pl\ra\M_\al^\pl$ are $G$-equivariant. 

Hence $\bar\Up_{(\al,\bs d)}$ in \eq{co5eq15} extends to $\bar\Up_{(\al,\bs d)}^G:H_*^G(\baM_{(\al,\bs d)}^\pl)\longra H_*^G(\M_\al^\pl)$. As $\dM_\al,\dM_\al^\pl$ are $G$-invariant by Assumption \ref{co5ass4}(c), we see that $\dM_{(\al,\bs d)},\dM_{(\al,\bs d)}^\pl$ are $G$-invariant. 

Assumption \ref{co5ass4}(b)--(d) imply that the constructions of $\bs\baM_{(\al,\bs d)},\bs\baM_{(\al,\bs d)}^\pl$ and $\bs\dM_{(\al,\bs d)}^\red,\bs\dM_{(\al,\bs d)}^\rpl$ in Definition \ref{co5def1}, and the morphisms between them, all make sense in $G$-equivariant derived stacks and $G$-equivariant morphisms. Hence the obstruction theories \eq{co5eq18} on $\dM_{(\al,\bs d)},\dM_{(\al,\bs d)}^\pl$ are $G$-equivariant.

The argument for $\M_\al^\rst(\tau),\M_\al^\ss(\tau)$ above implies that for $(\bar\tau^\la_{\bs\mu},\bar T,\le)\in\baS$, if $(\al,\bs d)\in C(\baB)$ then $\baM_{(\al,\bs d)}^\rst(\bar\tau^\la_{\bs\mu})\subseteq\baM_{(\al,\bs d)}^\ss(\bar\tau^\la_{\bs\mu})\subseteq\baM_{(\al,\bs d)}^\pl$ are $G$-invariant open substacks. Thus if $\al\in C(\B)_\pe$ and $\baM_{(\al,\bs d)}^\rst(\bar\tau^\la_{\bs\mu}) =\baM_{(\al,\bs d)}^\ss(\bar\tau^\la_{\bs\mu})$ then the virtual class $[\baM_{(\al,\bs d)}^\ss(\bar\tau^\la_{\bs\mu})]_\virt$ in \eq{co5eq22} exists in $H_*^G(\baM_{(\al,\bs d)}^\pl)$, so \eq{co5eq23} becomes
\begin{equation*}
\bar\Up_{(\al,\bs d)}^G\bigl([\baM_{(\al,\bs d)}^\ss(\bar\tau^\la_{\bs\mu})]_\virt\bigr)\in H_{2o_\al+2-2\chi(\al,\al)}^G(\M_\al^\pl)=\check H_{2o_\al}^G(\M_\al^\pl).
\end{equation*}
\end{dfn}

The next theorem extends Theorems \ref{co5thm1}--\ref{co5thm3} to $G$-equivariant (co)homology. We describe the modifications needed to the proofs of Theorems \ref{co5thm1}--\ref{co5thm2} in \S\ref{co94} and \S\ref{co107}. Apart from replacing $\check H_*(\M^\pl)$ by $\check H_*^G(\M^\pl)$, the proof of Theorem \ref{co5thm3} needs no changes.

\begin{thm}
\label{co5thm4}
Suppose Assumption\/ {\rm\ref{co5ass4}} (which includes Assumptions {\rm\ref{co4ass1}, \ref{co4ass2}} and\/ {\rm\ref{co5ass1}}), Assumption\/ {\rm\ref{co5ass2},} and (for Theorem\/ {\rm\ref{co5thm3}}) Assumption {\rm\ref{co5ass3}} hold. 

Then Theorems\/ {\rm\ref{co5thm1}--\ref{co5thm3}} hold, replacing $\check H_*(\M^\pl)$ by $\check H_*^G(\M^\pl)$ throughout. Here $[\M_\al^\ss(\tau)]_\virt$ and\/ \eq{co5eq30} in Theorem\/ {\rm\ref{co5thm1}(i),(iii)} are defined in $\check H_{2o_\al}^G(\M_\al^\pl)$ as in Definition\/ {\rm\ref{co5def2},} and the Lie bracket on $\check H^G_*(\M^\pl)$ is given in Theorem\/ {\rm\ref{co4thm5}}. Write $[\M_\al^\ss(\tau)]_\inv^G\in\check H_{2o_\al}^G(\M_\al^\pl)$ for the corresponding invariants. The Lie algebra morphism $\La^{G,\{1\}}:H^G_*(\M^\pl)\ra H_*(\M^\pl)$ in Theorem\/ {\rm\ref{co4thm5}} maps $[\M_\al^\ss(\tau)]_\inv^G\mapsto[\M_\al^\ss(\tau)]_\inv,$ with\/ $[\M_\al^\ss(\tau)]_\inv\in\check H_{2o_\al}(\M_\al^\pl)$ as in Theorem\/~{\rm\ref{co5thm1}}.
\end{thm}

Readers are warned that in the context of $\bG_m$-localization, it is tempting to try to apply Theorem \ref{co5thm4} incorrectly, as we explain next.

\begin{rem}
\label{co5rem3}
{\bf(a)} Take $G\cong(\bG_m)^n$ to be an algebraic torus. Given a moduli stack $\M^\rst_\al(\tau)=\M^\ss_\al(\tau)$ which is a proper algebraic space with a $G$-action and a $G$-equivariant obstruction theory, then as in Theorem \ref{co2thm3} (which covers the case $G=\bG_m$), the virtual class $[\M^\ss_\al(\tau)]_\virt\in H_*^G(\M_\al^\pl)$ may be written (at least in {\it localized\/} $G$-equivariant homology $H_*^G(\M_\al^\pl)_{\rm loc}$) in terms of the virtual class $[\M^\ss_\al(\tau)^G]_\virt$ of the $G$-fixed substack $\M^\ss_\al(\tau)^G$, as in~\eq{co2eq25}.

There are important examples in which $\M^\ss_\al(\tau)$ is not proper, so $[\M^\ss_\al(\tau)]_\virt$ does not make sense, but $\M^\ss_\al(\tau)^G$ is proper. Then it is a well known technique to simply {\it define\/} $[\M^\ss_\al(\tau)]_\virt\in H_*^G(\M_\al^\pl)_{\rm loc}$ in terms of $[\M^\ss_\al(\tau)^G]_\virt$ by \eq{co2eq25}. This happens in Thomas et al.\ \cite{KoTh3,MaTh,PaTh2,TaTh1,TaTh2,Thom2}, on Vafa--Witten invariants and related topics, when $X$ is a projective surface, and $G=\bG_m$ acts on the canonical bundle $K_X$ by scaling the fibres, and $\A=\B=\coh_\cs(K_X)$ is the category of compactly-supported coherent sheaves on~$K_X$.

Readers are warned that Theorem \ref{co5thm4} requires $\M^\ss_\al(\tau)$ to be proper when $\al\in C(\B)_\pe$ and $\M^\rst_\al(\tau)=\M^\ss_\al(\tau)$, as in Assumption \ref{co5ass2}(g). {\it It is not sufficient for $\M^\ss_\al(\tau)^G$ to be proper}. Although Theorem \ref{co5thm4} would make sense for $[\M^\ss_\al(\tau)]_\virt\in H_*^G(\M_\al^\pl)_{\rm loc}$ defined by \eq{co2eq25} (in particular, $\check H_*^G(\M_\al^\pl)_{\rm loc}$ is a Lie algebra, as in Remark \ref{co4rem4}), {\it it may be false}.

The point where the proof breaks down in this case is subtle. It relates to commuting $-\cap e(\cN_a^\bu)^{-1}$ in \eq{co2eq25} past $\sum_{j\ge 0}z^jt^j\bt-$ in \eq{co4eq10}. This can involve powers $(y+z)^n$, where $y\in H^*_G(*)$ is a $G$-weight in $\cN_a^\bu$, and $z$ is the formal variable in the vertex algebra. We can expand  $(y+z)^n$ using the binomial theorem in nonnegative powers of either $y$ or $z$, but if $n<0$ these are different, which can cause the proof to fail. A very similar issue occurs in the proof of Theorem \ref{co4thm3}, and is the reason we need $F,G,G^\pl$ in Definition \ref{co4def5} to be vector bundles, and not just complexes.

\smallskip

\noindent{\bf(b)} Despite {\bf(a)}, there {\it is\/} a way to apply Theorems \ref{co5thm1}--\ref{co5thm3} to examples in which $\M^\ss_\al(\tau)^G$ is proper, but $\M^\ss_\al(\tau)$ is not. We illustrate it in Example \ref{co6ex4} below.

Instead of considering categories $\B\subseteq\A$ with a $G$-action, as in Assumption \ref{co5ass4}, we can define categories $\B^G\subseteq\A^G$ of $G$-{\it equivariant objects\/} in $\B,\A$. That is, an object $(E,\rho)$ of $\B^G$ is an object $E$ of $\B$ with a morphism $\rho:G\ra\Aut(E)$ compatible with the $G$-action on $\B$. There is a forgetful functor $\B^G\ra\B$ mapping $(E,\rho)\mapsto E$. Then $G$-fixed substacks $\M^G,\M^\ss_\al(\tau)^G$ are isomorphic to moduli stacks of all objects, and $\tau$-semistable objects, in~$\B^G$.

Under suitable assumptions, we should apply Theorems \ref{co5thm1}--\ref{co5thm3}, in non-equivariant homology, to the categories $\B^G\subseteq\A^G$. This gives a system of invariants $[\M^\ss_\al(\tau)^G]_\inv$ and wall-crossing formulae in $\check H_*((\M^G)_\pl)$, which include virtual classes $[\M^\ss_\al(\tau)^G]_\virt$ when $\M^\rst_\al(\tau)=\M^\ss_\al(\tau)$. We can then push these invariants forward to $H_*^G(\M_\al^\pl)_{\rm loc}$ as in \eq{co2eq25}, by mapping
\e
[\M^\ss_\al(\tau)^G]_\inv\mapsto (-1)^{\rank\cN^\bu}i_*\bigl([\M^\ss_\al(\tau)^G]_\inv\cap e(\cN^\bu)^{-1}\bigr)=:[\M^\ss_\al(\tau)]_\inv.
\label{co5eq35}
\e
The problem in {\bf(a)} is that if \eq{co5eq35} is not a Lie algebra morphism $\check H_*((\M^G)_\pl)\ab\ra \check H_*^G(\M_\pl)$, which requires extra conditions on $\cN^\bu$ similar to those in Theorem \ref{co4thm3}, then the $[\M^\ss_\al(\tau)]_\inv$ may not satisfy the identities~\eq{co5eq31}--\eq{co5eq34}. 	
\end{rem}

\subsection{Extension to triangulated categories}
\label{co55}

Finally we discuss extending the programme of \S\ref{co51}--\S\ref{co53} from abelian categories $\A$ to triangulated categories $\T$. We do not prove any results. See \cite[\S 7]{Joyc7} and \cite[\S 4.3.2]{GJT} for similar ideas. Throughout $\T$ is a $\C$-linear triangulated category coming from Algebraic Geometry, such as $D^b\coh(X)$ for $X$ a smooth projective surface, or $D^b\modCQI$ for a quiver with relations~$\br Q$.
\smallskip

\noindent{\bf(a)} The (projective linear) moduli stacks $\M,\M^\pl$ of objects in $\T$ are in general {\it higher stacks}, as in \S\ref{co22}, rather than Artin stacks. As explained in Remark \ref{co4rem2}, the vertex algebra structure on $\hat H_*(\M)$ in \S\ref{co42} works for triangulated categories with essentially no change. For the Lie algebra $\check H_*(\M^\pl)$, Assumption \ref{co4ass1}(g) fails when $\al=0$, so \eq{co4eq15} may not be an isomorphism on $\check H_*(\M_0^\pl)$, but $\check H_*(\M^\pl)$ still has a natural graded Lie algebra structure for which \eq{co4eq15} is a graded Lie algebra morphism.
\smallskip

\noindent{\bf(b)} As in Remark \ref{co1rem4}(a), surprisingly, as proved by Gross~\cite[Th.~1.1]{Gros}, the vertex algebra $\hat H_*(\M)$ is usually {\it easier\/} to compute explicitly for triangulated categories than for abelian categories. This is because in the abelian case $\M$ is a like an `abelian monoid in stacks', with addition $\op$, but in the triangulated case $\M$ is a like an `abelian group in stacks', since $[1]:\T\ra\T$ acts like an inverse to $\op$ at the homotopy level; and groups are simpler than monoids. 
\smallskip

\noindent{\bf(c)} Given a surjection $K_0(\T)\twoheadrightarrow K(\T)$ as in Assumption \ref{co4ass1}(d), a {\it Bridgeland stability condition\/} \cite{Brid} $(Z,\cP)$ on $\T$ is a group morphism $Z:K(\T)\ra\C$ and a family of abelian subcategories $\cP(\phi)\subset\T$ for $\phi\in\R$, such that if $0\ne E\in\cP(\phi)$ then $Z(\lb E\rb)=re^{i\pi\phi}$ for $r>0$, and $\cP(\phi)[k]=\cP(\phi+k)$ for $k\in\Z$, and other conditions hold. For $\phi\in\R$ we write $\cP([\phi,\phi+1))$ for the subcategory of $\T$ generated by $\cP(\psi)$ for $\psi\in[\phi,\phi+1)$. Then $\cP([\phi,\phi+1))$ is an abelian subcategory of $\T$, the heart of a t-structure on $\T$. Bridgeland \cite[Th.~1.2]{Brid} proves that under some conditions, the moduli space $\Stab(\T)$ of such $(Z,\cP)$ is a topological space, such that $(Z,\cP)\mapsto Z$ is a local homeomorphism~$\Stab(\T)\ra \Hom(K(\T),\C)$. 
\smallskip

\noindent{\bf(d)} Here is our approach to extending \S\ref{co51}--\S\ref{co53} to triangulated categories:

\begin{princ}
\label{co5princ1}
{\bf(i)} It is unnecessary to prove new versions of Theorems\/ {\rm\ref{co5thm1}--\ref{co5thm3}} for triangulated categories. A triangulated category $\T$ usually contains many abelian subcategories $\A\subset\T,$ such as $\cP(\phi),\cP([\phi,\phi+1))$ in {\bf(c)\rm,} and we can try to apply Theorems\/ {\rm\ref{co5thm1}--\ref{co5thm3}} to such\/ $\A$ to deduce what we want. 

However, there may be difficult issues in verifying Assumptions\/ {\rm\ref{co5ass1}--\ref{co5ass3}} for suitable $\A\subset\T,$ in particular, quasi-smoothness of derived moduli stacks in Assumption\/ {\rm\ref{co5ass1}(f),} and properness of moduli spaces in Assumption\/~{\rm\ref{co5ass2}(g),(h)}.
\smallskip

\noindent{\bf(ii)} For the analogue of Theorem\/ {\rm\ref{co5thm1}} for Bridgeland stability conditions $(Z,\cP),$ to define $[\M_\al^\ss(Z,\cP)]_\inv$ we 
apply Theorem\/ {\rm\ref{co5thm1}} to the trivial stability condition on the abelian category $\cP(\phi)$ for $\phi\in[0,2)$ with\/ $Z(\al)=re^{i\pi\phi}$ for\/ $r>0$. 
\smallskip

\noindent{\bf(iii)} We should only expect {\rm\eq{co5eq31}--\eq{co5eq34}} to make sense if\/ $(Z,\cP),(\ti Z,\ti{\mathcal P})$ are \begin{bfseries}sufficiently close\end{bfseries} in\/~$\Stab(\T)$.

If\/ $(Z,\cP),(\ti Z,\ti{\mathcal P})$ are not close then we cannot define $U(\cdots;(Z,\cP),(\ti Z,\ti{\mathcal P})),$ $\ti U(\cdots;(Z,\cP),(\ti Z,\ti{\mathcal P}))$ in\/ {\rm\eq{co5eq31}--\eq{co5eq34},} which would be infinite sums anyway.
\smallskip

\noindent{\bf(iv)} If\/ $(Z,\cP),(\ti Z,\ti{\mathcal P})$ are sufficiently close then we can construct $(\hat Z,\hat{\mathcal P})$ in $\Stab(\T)$ such that $\Im\hat Z=\Im Z,$ $\hat{\mathcal P}([0,1))=\cP([0,1)),$ $\Re\hat Z=\Re\ti Z,$ and\/ $\hat{\mathcal P}([-\ha,\ha))=\ti{\mathcal P}([-\ha,\ha))$. We can then try to wall-cross $(Z,\cP)\Ra(\hat Z,\hat{\mathcal P})$ by applying Theorem\/ {\rm\ref{co5thm3}} to the abelian category $\cP([0,1)),$ and wall-cross $(\hat Z,\hat{\mathcal P})\Ra(\ti Z,\ti{\mathcal P})$ by applying Theorem\/ {\rm\ref{co5thm3}} to the abelian category~$\ti{\mathcal P}([-\ha,\ha))$. 
\smallskip

\noindent{\bf(v)} If\/ $(Z,\cP),(\ti Z,\ti{\mathcal P})$ lie in the same connected component of\/ $\Stab(\T),$
 we should relate invariants $[\M_\al^\ss(Z,\cP)]_\inv,[\M_\be^\ss(\ti Z,\ti{\mathcal P})]_\inv$ by choosing a continuous path\/ $(Z_t,\cP_t)_{t\in[0,1]}$ from $(Z,\cP)$ to $(\ti Z,\ti{\mathcal P}),$ and\/ $N\gg 0$ such that\/ $(Z_{(i-1)/N},\cP_{(i-1)/N})$ and\/ $(Z_{i/N},\cP_{i/N})$ are sufficiently close for $i=1,\ldots,N,$ and then wall-crossing $(Z_{(i-1)/N},\cP_{(i-1)/N})\Ra(Z_{i/N},\cP_{i/N})$ for\/ $i=1,\ldots,N$ as in\/ {\bf(iii)\rm--\bf(iv)}. Thus we relate $[\M_\al^\ss(Z,\cP)]_\inv,[\M_\be^\ss(\ti Z,\ti{\mathcal P})]_\inv$ in finitely many steps.
\end{princ}

\noindent{\bf(e)} Kontsevich and Soibelman \cite{KoSo} wrote down a wall-crossing formula for motivic Donaldson--Thomas invariants of Bridgeland stability conditions on $D^b\coh(X)$ for a Calabi--Yau 3-fold $X$. We explained in \cite[\S 4.2]{GJT} how to rewrite the wall-crossing formulae \eq{co5eq33}--\eq{co5eq34} in the style of Kontsevich--Soibelman.
\smallskip

\noindent{\bf(f)} Bayer--Lahoz--Macr\`\i--Nuer--Perry--Stellari \cite[Th.~1.4]{BLMN} use Alper et al.\ \cite{AHLH} to prove properness of moduli spaces of $(Z,\cP)$-semistable objects in $D^b\coh(X)$, as required in Assumption \ref{co5ass2}(g), provided $(Z,\cP)$ satisfies suitable conditions.

\section{Applications to quivers}
\label{co6}

We now apply the results of Chapter \ref{co5} to abelian categories $\modCQ$ of a quiver $Q$, and $\modCQI$ of a quiver with relations $(Q,I)$. To get properness of semistable moduli stacks in Assumption \ref{co5ass2}(g),(h) we restrict to $Q$ with no oriented cycles. A different proof of Theorems \ref{co5thm1}--\ref{co5thm3} for $\A=\modCQ$ was already given in Gross--Joyce--Tanaka~\cite[Th.~5.8]{GJT}.

\subsection{Quivers, and quivers with relations}
\label{co61}

Here are the basic definitions in quiver theory, as in Benson~\cite[\S 4.1]{Bens}.

\begin{dfn}
\label{co6def1}
A {\it quiver\/} $Q$ is a finite directed graph. That is, $Q$ is a quadruple $(Q_0,Q_1,h,t)$, where $Q_0$ is a finite set of {\it vertices}, $Q_1$ is a finite set of {\it edges}, and $h,t:Q_1\ra Q_0$ are maps giving the {\it head\/} and {\it tail\/} of each edge.

A closed loop of directed edges $\overset{v_0}{\bu}\,{\buildrel e_1\over \longra}\,\overset{v_1}{\bu}\,{\buildrel e_2\over \longra}\,\cdots \,{\buildrel e_n\over \longra}\,\overset{v_n=v_0}{\bu}$ is an {\it oriented cycle\/} in $Q$. Sometimes we restrict to quivers with {\it no oriented cycles}.

The {\it path algebra\/} $\C Q$ is the associative algebra over $\C$ with basis all {\it paths of length\/}
$k\ge 0$, that is, sequences of the form
\e
v_0\,{\buildrel e_1\over\longra}\, v_1\ra\cdots\ra
v_{k-1}\,{\buildrel e_k\over\longra}\,v_k,
\label{co6eq1}
\e
where $v_0,\ldots,v_k\in Q_0$, $e_1,\ldots,e_k\in Q_1$,
$t(a_i)=v_{i-1}$ and $h(a_i)=v_i$. Multiplication is given by
composition of composable paths in reverse order. That~is,
\begin{align*}
&\bigl(w_0\,{\buildrel f_1\over\longra}\, w_1\ra\cdots\ra
w_{l-1}\,{\buildrel f_l\over\longra}\,w_l\bigr)\cdot\bigl(v_0\,{\buildrel e_1\over\longra}\, v_1\ra\cdots\ra
v_{k-1}\,{\buildrel e_k\over\longra}\,v_k\bigr)\\
&=\begin{cases}
v_0\,{\buildrel e_1\over\longra}\,\cdots\,{\buildrel e_k\over\longra}\, v_k=w_0\,{\buildrel f_1\over\longra}\,\cdots\,{\buildrel f_l\over\longra}\, w_l, & v_k=w_0, \\
0, & \text{otherwise.}
\end{cases}
\end{align*}

Each $v\in Q_0$ gives a basis element $\ga_v\in\C Q$ corresponding to the path of length 0 at $v$, with $\ga_v^2=\ga_v$, and the identity in $\C Q$ is $1=\sum_{v\in Q_0}\ga_v$. Also each $e\in Q_1$ gives a basis element $\ga_e\in\C Q$ corresponding to $e$ as a path of length~1.

For $n\ge 0$, write $\C Q_{(n)}$ for the vector subspace of $\C Q$ with basis all paths of length $k\ge n$. It is an ideal in $\C Q$. For $v,w\in Q_0$, write $\C Q_{vw}\subseteq\C Q$ for the vector subspace of paths from $v$ to $w$, of the form
\begin{equation*}
v=v_0\,{\buildrel e_1\over\longra}\, v_1\ra\cdots\ra
v_{k-1}\,{\buildrel e_k\over\longra}\,v_k=w.
\end{equation*}
Then $\C Q=\bigop_{v,w\in Q_0}\C Q_{vw}$, and $\C Q_{vw}=\ga_w\cdot\C Q\cdot\ga_v$, and multiplication in $\C Q$ maps~$\C Q_{vw}\t\C Q_{uv}\ra\C Q_{uw}$. 

For all $u,v,w\in Q_0$ and $e\in Q_1$, define $\C$-linear maps
\begin{align*} 
&\frac{\pd}{\pd e}:\C Q_{uw}\longra \C Q_{ut(e)}\ot_\C\C Q_{t(e)w},\quad
\frac{\pd}{\pd e}:\bigl(u=v_0\,{\buildrel e_1\over\longra}\, \cdots
\,{\buildrel e_k\over\longra}\,v_k=w\bigr)\longmapsto\\
&\sum_{i=1,\ldots,k:\; e_i=e\!\!\!\!\!\!\!\!\!\!\!\!\!\!\!\!\!}\bigl(u=v_0\,{\buildrel e_1\over\longra}\, \cdots\,{\buildrel e_{i-1}\over\longra}\,v_i=v_{t(e)}\bigr)\ot \bigl(v_{h(e)}=v_{i+1}\,{\buildrel e_{i+1}\over\longra}\, \cdots\,{\buildrel e_k\over\longra}\,v_k=w\bigr),\\
&\frac{\pd}{\pd v}:\C Q_{uw}\longra \C Q_{uv}\ot_\C\C Q_{vw},\quad
\frac{\pd}{\pd v}:\bigl(u=v_0\,{\buildrel e_1\over\longra}\, \cdots
\,{\buildrel e_k\over\longra}\,v_k=w\bigr)\longmapsto\\
&\sum_{i=0,\ldots,k:\; v_i=v\!\!\!\!\!\!\!\!\!\!\!\!\!\!\!\!\!}\bigl(u=v_0\,{\buildrel e_1\over\longra}\, \cdots\,{\buildrel e_{i-1}\over\longra}\,v_i=v\bigr)\ot \bigl(v=v_i\,{\buildrel e_i\over\longra}\, \cdots\,{\buildrel e_k\over\longra}\,v_k=w\bigr).
\end{align*}
It is then easy to show that in maps $\C Q_{uw}\longra \C Q_{uv}\ot_\C\C Q_{vw}$ we have
\e
\begin{split}
\frac{\pd}{\pd v}&=\de_{uv}(\ga_u\ot \id_{\C Q_{uw}})+\sum_{e\in Q_1:h(e)=v}((e\cdot-)\ot\id_{\C Q_{vw}})\ci\frac{\pd}{\pd e}\\
&=\de_{vw}(\id_{\C Q_{uw}}\ot\ga_w)+\sum_{e\in Q_1:t(e)=v}(\id_{\C Q_{uv}}\ot(-\cdot e))\ci\frac{\pd}{\pd e}.
\end{split}
\label{co6eq2}
\e

A {\it representation} $((V_v)_{v\in Q_0},\ab (\rho_e)_{e\in Q_1})$ of $Q$ consists of finite-dimensional $\C$-vector spaces $V_v$ for each $v\in Q_0$, and linear maps $\rho_e:V_{t(e)}\ra V_{h(e)}$ for each $e\in Q_1$. Representations of $Q$ are in 1-1 correspondence with {\it finite-dimensional left\/ $\C Q$-modules\/} $(V,\rho)$, as follows.

Given $((V_v)_{v\in Q_0},\ab (\rho_e)_{e\in Q_1})$, define $V=\bigop_{v\in Q_0}V_v$, and a linear map $\rho:\C Q\ab\ra\End(V)$ taking \eq{co6eq1} to the linear map $V\ra V$ acting as $\rho_{e_k}\ci \rho_{e_{k-1}}\ci\cdots\ci\rho_{e_1}$ on $V_{v_0}$, and 0 on $V_v$ for $v\ne v_0$. Then $(V,\rho)$ is a left $\C Q$-module. Conversely, any such $(V,\rho)$ comes from a unique representation $((V_v)_{v\in Q_0},\ab (\rho_e)_{e\in Q_1})$ of $Q$, with $V_v=\ga_v(V)\subseteq V$ and $\rho_e=\rho\vert_{V_{h(e)}}$. We will use the notation $((V_v)_{v\in Q_0},\ab (\rho_e)_{e\in Q_1})$ and $(V,\rho)$ interchangeably.

A {\it morphism of representations\/} $\phi:(V,\rho)\ra(W,\si)$ is a linear map $\phi:V\ra W$ with $\phi\ci\rho(\ga)=\si(\ga)\ci\phi$ for all $\ga\in\C Q$. Equivalently, $\phi$ defines linear maps $\phi_v:V_v\ra W_v$ for all $v\in Q_0$ with $\phi_{h(e)}\ci\rho_e=\si_e\ci\phi_{t(e)}$ for all $e\in Q_1$. Write $\modCQ$ for the $\C$-linear abelian category of representations of~$Q$.

Write $\Z^{Q_0}$ for the abelian group of maps $Q_0\ra\Z$, and $\N^{Q_0}\subset\Z^{Q_0}$ for the submonoid of maps $Q_0\ra\N$. The {\it dimension vector\/} $\bdim(V,\rho)\in\N^{Q_0}$ of a representation $(V,\rho)$ is $\bdim(V,\rho):v\mapsto\dim_\C V_v$. This induces a surjective morphism $\bdim:K_0(\modCQ)\ra\Z^{Q_0}$, where $K_0(\modCQ)$ is the Grothendieck group of $\modCQ$. We will take $K(\modCQ)=\Z^{Q_0}$ in Assumption~\ref{co4ass1}(d). 
\end{dfn}

\begin{dfn}
\label{co6def2}
A {\it quiver with relations\/} $\br Q$ is an octuple $(Q_0,Q_1,h,t,Q_2,b,e,r)$, where $Q=(Q_0,Q_1,h,t)$ is a quiver, and $Q_2$ is a finite set, and $b,e:Q_2\ra Q_0$ are maps, and $r:Q_2\ra\C Q_{(2)}\subset\C Q$ is a map such that $r(q)\in \C Q_{b(q),e(q)}$ for each $q\in Q_2$. Write $I=\an{r(q):q\in Q_2}\subseteq\C Q_{(2)}\subset\C Q$ for the two-sided ideal generated by $r(q)$ for all $q\in Q_2$. Then $\C Q/I$ is an associative $\C$-algebra, obtained by imposing the relations $r(q)=0$ in $\C Q$ for all~$q\in Q_2$.

We say that $\br Q$ {\it has no single-vertex relations\/} if $b(q)\ne e(q)$ for all~$q\in Q_2$.

We say that $\br Q$ {\it has homogeneous relations\/} if $r(q)$ is homogeneous (that is, a $\C$-linear combination of paths \eq{co6eq1} of a fixed length $k\ge 2$) for each~$q\in Q_2$.

A {\it representation\/} $(V,\rho)=((V_v)_{v\in Q_0},(\rho_e)_{e\in Q_1})$ of $\br Q$ is a representation of $Q$ such that $\rho(I)=0$ in $\Hom(V,V)$, or equivalently, such that $\rho(r(q))=0$ as a map $V_{b(q)}\ra V_{e(q)}$ for all $q\in Q_2$. 

If $(V,\rho),(W,\si)$ are representations of $\br Q$, a {\it  morphism of representations\/} $\phi:(V,\rho)\ra(W,\si)$ of $\br Q$ is just a morphism of $Q$-representations. Write $\modCQI$ for the category of representations of $\br Q$. It is a $\C$-linear abelian category, a full subcategory $\modCQI\subseteq\modCQ$. 
\end{dfn}

\begin{rem}
\label{co6rem1}
The category $\modCQI$ depends only on the quiver $Q$ and the two-sided ideal $I\subset\C Q_{(2)}$, not on the data $(Q_2,b,e,r)$ used to define $I$. A {\it quiver with relations\/} is usually defined to be a pair $(Q,I)$, rather than the octuple $\br Q=(Q_0,Q_1,h,t,Q_2,b,e,r)$. Since $I$ does not determine $(Q_2,b,e,r)$, our notion of quiver with relations is {\it more data\/} than the usual notion. 

We do this because the extra data $(Q_2,b,e,r)$ is essential for the following:
\begin{itemize}
\setlength{\itemsep}{0pt}
\setlength{\parsep}{0pt}
\item[(i)] To define $\cE^\bu\ra\M\t\M$ in Assumption \ref{co4ass1}(f), and hence the vertex and Lie algebra structures on $\hat H_*(\M),\check H_*(\M^\pl)$ in \S\ref{co42}--\S\ref{co43}.
\item[(ii)] To define the derived moduli stacks $\bs\M_{\bs d},\bs\M_{\bs d}^\pl,\bs\dM^\red_{\bs d},\bs\dM_{\bs d}^\rpl$ in Assumption \ref{co5ass1}(d),(f), and hence the obstruction theories \eq{co5eq5} on $\dM_{\bs d},\dM_{\bs d}^\pl$ and the virtual classes $[\M_{\bs d}^\ss(\mu)]_\virt$ when~$\M_{\bs d}^\rst(\mu)=\M_{\bs d}^\ss(\mu)$.
\end{itemize}	
\end{rem}

\begin{rem}
\label{co6rem2}
We can easily extend Definitions \ref{co6def1}--\ref{co6def2} and the rest of Chapter \ref{co6} to allow $Q_0,Q_1,Q_2$ to be {\it infinite\/} sets, provided:
\begin{itemize}
\setlength{\itemsep}{0pt}
\setlength{\parsep}{0pt}
\item[(a)] For all $v,w\in Q_0$, there are only finitely many $e\in Q_1$ and $q\in Q_2$ with $t(e)=v$, $h(e)=w$, $b(q)=v$, $e(q)=w$.
\item[(b)] A representation $((V_v)_{v\in Q_0},(\rho_e)_{e\in Q_1})$ has $V_v\ne 0$ for only finitely many $v\in Q_0$. Equivalently, representations $(V,\rho)$ have $\dim V<\iy$. 
\item[(c)] We define $\Z^{Q_0}$ so that $\bs d\in\Z^{Q_0}$ has $\bs d(v)\ne 0$ for only finitely many $v\in Q_0$.
\end{itemize}	
Example \ref{co6ex4} below includes such an infinite quiver.
\end{rem}

\subsection{\texorpdfstring{Moduli stacks of objects in $\modCQ$ and $\modCQI$}{Moduli stacks of objects in modCQ and modCQ/I}}
\label{co62}

We describe the moduli stacks $\M,\M^\pl$ in \S\ref{co42}--\S\ref{co43} for $\modCQ,\modCQI$.

\begin{dfn}
\label{co6def3}
Let $Q=(Q_0,Q_1,h,t)$ be a quiver. Write $\M$ for the moduli stack of objects $(V,\rho)$ in $\modCQ$. Then there is a natural decomposition
\begin{equation*}
\M=\ts\coprod_{\bs d\in\N^{Q_0}}\M_{\bs d},
\end{equation*}
where $\M_{\bs d}$ is the moduli stack of $(V,\rho)$ with $\bdim(V,\rho)=\bs d$, a smooth Artin $\C$-stack. For any such $(V,\rho)$, by considering isomorphisms $V_v\cong\C^{\bs d(v)}$ for $v\in Q_0$ we see we may write $\M_{\bs d}$ explicitly as a quotient stack $\M_{\bs d}=[R_{\bs d}/\GL_{\bs d}]$, where $R_{\bs d}=\ts\prod_{e \in Q_1} \Hom(\C^{\bs d(t(e))},\C^{\bs d(h(e))})$ and $\GL_{\bs d}=\prod_{v \in Q_0} \GL(\bs d(v),\C)$ with action
\begin{equation*}
(A_v)_{v\in Q_0}:((B_e)_{e\in Q_1})\longmapsto (A_{h(e)}\ci B_e\ci A_{t(e)}^{-1})_{e\in Q_1}.
\end{equation*}

Write $\cV_v\ra\M$ for $v\in Q_0$ for the tautological vector bundle with
\begin{equation*}
\cV_v\vert_{[((V_v)_{v\in Q_0},(\rho_e)_{e\in Q_1})]}=V_v,
\end{equation*}
and write $\cV_{v,\bs d}=\cV_v\vert_{\M_{\bs d}}$ for $\bs d\in\N^{Q_0}$, so that $\rank \cV_{v,\bs d}=\bs d(v)$. As $R_{\bs d}$ is contractible, we have $\bA^1$-homotopy equivalences
\begin{equation*}
\M_{\bs d}\simeq[*/\GL_{\bs d}]=\prod_{v\in Q_0}[*/\GL(\bs d(v),\C)].	
\end{equation*}
Thus the topological realization of $\M_{\bs d}$ is
\begin{equation*}
\M_{\bs d}^\top\simeq\prod_{v\in Q_0}B\GL(\bs d(v),\C).
\end{equation*}
As $\GL(r,\C)\simeq\U(r)$, the computation of $H^*(B\U(r))$ by Milnor and Stasheff \cite[Th.~14.5]{MiSt} implies that the cohomology of $\M_{\bs d}$ over $R$ is
\e
H^*(\M_{\bs d})=H^*(\M_{\bs d}^\top,\Q)\cong \Q[c_{v,{\bs d}}^i:v\in Q_0,\; i=1,2,\ldots,\bs d(v)],
\label{co6eq3}
\e
where $c_{v,{\bs d}}^i$ is a formal variable of degree $2i$, with $c_{v,{\bs d}}^i=c_i(\cV_{v,\bs d})$. The homology $H_*(\M_{\bs d})$ is the $\Q$-linear dual of \eq{co6eq3}, and $H_*(\M)=\bigop_{\bs d\in\N^{Q_0}}H_*(\M_{\bs d})$.

Similarly, the projective linear moduli stack $\M^\pl$ from Definition \ref{co4def4} is $\M^\pl=\coprod_{\bs d\in\N^{Q_0}}\M_{\bs d}^\pl$, where $\M_{\bs d}^\pl=[R_{\bs d}/\PGL_{\bs d}]$, for $\PGL_{\bs d}=\GL_{\bs d}/\bG_m$ with $\bG_m=\bigl\{(\la\,\id_{\bs d(v)})_{v\in Q_0}:0\ne\la\in\C\bigr\}\subseteq\GL_{\bs d}$.
\end{dfn}

\begin{dfn}
\label{co6def4}
Let $\br Q=(Q_0,Q_1,h,t,Q_2,b,e,r)$ be a quiver with relations, with quiver $Q=(Q_0,Q_1,h,t)$, giving categories $\modCQI\subset\modCQ$. Let $\M,\M^\pl$ be as in Definition \ref{co6def3} for $\modCQ$. Write $\brM$ for the moduli stack of objects in $\modCQI$, and $\brM^\pl=\brM/[*/\bG_m]$, as for $\M,\M^\pl$. The inclusion $i:\modCQI\hookra\modCQ$ induces stack morphisms $i_*:\brM\hookra\M$, $i_*^\pl:\brM^\pl\hookra\M^\pl$, which are inclusions of closed substacks. Then we have splittings $\brM=\ts\coprod_{\bs d\in\N^{Q_0}}\brM_{\bs d}$, $\brM^\pl=\ts\coprod_{\bs d\in\N^{Q_0}}\brM_{\bs d}^\pl$ with $i_*:\brM_{\bs d}\hookra\M_{\bs d}$, $i_*^\pl:\brM^\pl_{\bs d}\hookra\M^\pl_{\bs d}$.

To describe $\brM_{\bs d},\brM^\pl_{\bs d}$, define a vector bundle
\begin{equation*}
E_{\bs d}=\ts\bigop_{q\in Q_2}\cV_{\bs b(q),\bs d}^*\ot \cV_{\bs e(q),\bs d}\longra\M_{\bs d},
\end{equation*}
using the notation of Definition \ref{co6def3}, with a section $s_{\bs d}:\M_{\bs d}\ra E_{\bs d}$ mapping
\begin{equation*}
s_{\bs d}:((V_v)_{v\in Q_0},\ab (\rho_e)_{e\in Q_1})\longmapsto \ts\bigop_{q\in Q_2}r(q)\bigl((\rho_e)_{e\in Q_1}\bigr).
\end{equation*}
That is, $r(q)$ is a $\C$-linear combination of paths 
\begin{equation*}
\xymatrix@C=30pt{ \bs b(q)=v_0 \ar[r]^(0.6){e_1} & v_1 \ar[r] & \cdots \ar[r] &
v_{k-1} \ar[r]^(0.4){e_k} & v_k=\bs e(q), }
\end{equation*}
and $r(q)\bigl((\rho_e)_{e\in Q_1}\bigr)$ is the corresponding $\C$-linear combination of compositions
\begin{equation*}
\xymatrix@C=30pt{ V_{\bs b(q)}=V_{v_0} \ar[r]^(0.6){\rho_{e_1}} & V_{v_1} \ar[r] & \cdots \ar[r] &
V_{v_{k-1}} \ar[r]^(0.4){\rho_{e_k}} & V_{v_k}=V_{\bs e(q)}. }
\end{equation*}

Then $\brM_{\bs d}\subseteq\M_{\bs d}$ is the zero locus $s_{\bs d}^{-1}(0)$, a closed $\C$-substack. Equivalently, using $(r(q))_{q\in Q_2}$ we may define a $\GL_{\bs d}$-equivariant polynomial function
\e
S_{\bs d}:R_{\bs d}\longra \ts\bigop_{q\in Q_2}\Hom(\C^{\bs d(b(q))},\C^{\bs d(e(q))}),
\label{co6eq4}
\e
such that $E_{\bs d}$ and $s_{\bs d}$ are the induced vector bundle and section on $\M_{\bs d}=[R_{\bs d}/\GL_{\bs d}]$, and~$\brM_{\bs d}=[S_{\bs d}^{-1}(0)/\GL_{\bs d}]$. 

As the vector bundle $E_{\bs d}$ has $\bG_m$-weight zero, it descends to $\M_{\bs d}^\pl$, giving a vector bundle $E_{\bs d}^\pl\ra\M_{\bs d}^\pl$ with section $s_{\bs d}^\pl:\M_{\bs d}^\pl\ra E_{\bs d}^\pl$, which pull back to $E_{\bs d},s_{\bs d}$ under $\Pi^\pl_{\bs d}:\M_{\bs d}\ra\M_{\bs d}^\pl$. Then $\brM_{\bs d}^\pl\subseteq\M_{\bs d}^\pl$ is the zero locus $(s_{\bs d}^\pl)^{-1}(0)$. Equivalently, $\brM_{\bs d}^\pl=[S_{\bs d}^{-1}(0)/\PGL_{\bs d}]$. The projection $\brM_{\bs d}\ra\brM_{\bs d}^\pl$ has fibre $[\PGL_{\bs d}/\GL_{\bs d}]\cong[*/\bG_m]$ for $\bs d\ne\bs 0$, and so is a principal $[*/\bG_m]$-bundle.

If $\br Q$ has {\it homogeneous relations\/} then the components of $S_{\bs d}$ in \eq{co6eq4} are homogeneous polynomials on $R_{\bs d}$, so that $S_{\bs d}^{-1}(0)$ is a closed cone in $R_{\bs d}$, and thus is contractible. In this case we have $\bA^1$-homotopies $\brM_{\bs d}\simeq [*/\GL_{\bs d}]\simeq\M_{\bs d}$, $\brM_{\bs d}^\pl\simeq [*/\PGL_{\bs d}]\simeq\M_{\bs d}^\pl$, so that $H_*(\brM_{\bs d})\cong H_*(\M_{\bs d})$ is given in~\eq{co6eq3}.
\end{dfn}

\subsection{Slope stability conditions on quiver categories}
\label{co63}

\begin{dfn}
\label{co6def5}
Let $Q$ be a quiver, or $\br Q$ a quiver with relations. Use the notation of \S\ref{co61}--\S\ref{co62}. In the situation of \S\ref{co31} with $\A=\modCQ$ or $\A=\modCQI$, take $K(\A)=\Z^{Q_0}$, with $\lb E,\rho\rb=\bdim(V,\rho)$ for $(V,\rho)\in\A$. Then $C(\A)=\N^{Q_0}\sm\{0\}$. Fix $\mu_v\in\R$ for all $v\in Q_0$. Define $\mu:C(\A)\ra\R$ by
\begin{equation*}
\mu(\bs d)=\frac{\sum_{v\in Q_0}\mu_v\bs d(v)}{\sum_{v\in Q_0}\bs d(v)}\,.
\end{equation*}
Then $(\mu,\R,\le)$ is a stability condition on $\modCQ$ in the sense of Definition \ref{co3def1}, called {\it slope stability}. We call $\mu$ a {\it slope function}.

Writing $\M$ (or $\M^\pl$) for the (projective linear) moduli stack of objects in $\A$, for an object $E$ of $\modCQ$ to be $\mu$-stable, or $\mu$-semistable, is an open condition on the point $[E]$ in $\M$ or $\M^\pl$. Write $\M^\rst_{\bs d}(\mu)\subseteq\M^\ss_{\bs d}(\mu)\subseteq\M_{\bs d}^\pl$ for the open $\C$-substacks of $\mu$-(semi)stable objects. They are quotient stacks
\begin{equation*}
\M^\rst_{\bs d}(\mu)=[R^\rst_{\bs d}(\mu)/\PGL_{\bs d}],\quad 
\M^\ss_{\bs d}(\mu)=[R^\ss_{\bs d}(\mu)/\PGL_{\bs d}],
\end{equation*}
for $\PGL_{\bs d}$-invariant open subschemes $R^\rst_{\bs d}(\mu)\subseteq R^\ss_{\bs d}(\mu)$ in $R_{\bs d}$ or~$S_{\bs d}^{-1}(0)$.
\end{dfn}

The first part of the next proposition is proved in Gross--Joyce--Tanaka \cite[Prop.~5.7]{GJT} using results of King \cite{King}. The second part follows as $\brM^\ss_{\bs d}(\mu)\subseteq\M^\ss_{\bs d}(\mu)$ is a closed substack. King \cite{King} implies that $R^\rst_{\bs d}(\mu)\subseteq R^\ss_{\bs d}(\mu)$ are the open subschemes of GIT-(semi)stable points in $R_{\bs d}$ for a suitable linearization of the $\PGL_{\bs d}$-action on $R_{\bs d}$. Thus Geometric Invariant Theory implies that if $\M^\rst_{\bs d}(\mu)=\M^\ss_{\bs d}(\mu)$ then $\M^\ss_{\bs d}(\mu)$ is a quasiprojective $\C$-scheme, constructed as a GIT quotient. If $Q$ has no oriented cycles then $\M^\ss_{\bs d}(\mu)$ is projective.

\begin{prop}
\label{co6prop1}
Let\/ $Q$ be a quiver with \begin{bfseries}no oriented cycles\end{bfseries}, and\/ $\mu$ a slope function on $\modCQ,$ and\/ $\bs d\in \N^{Q_0}\sm\{0\}$ with\/ $\M^\rst_{\bs d}(\mu)=\M^\ss_{\bs d}(\mu)$. Then $\M^\ss_{\bs d}(\mu)$ is a smooth projective $\C$-scheme, and thus a proper algebraic\/ $\C$-space.

Similarly, let\/ $\br Q$ be a quiver with relations with no oriented cycles, and\/ $\mu$ a slope function on $\modCQI,$ and\/ $\bs d\in \N^{Q_0}\sm\{0\}$ with\/ $\brM^\rst_{\bs d}(\mu)=\brM^\ss_{\bs d}(\mu)$. Then $\brM^\ss_{\bs d}(\mu)$ is a  projective $\C$-scheme, and so a proper algebraic\/~$\C$-space.
\end{prop}

\subsection{Verifying Assumptions \ref{co4ass1} and \ref{co5ass1}--\ref{co5ass3} for quivers}
\label{co64}

This section gives the data and verifies Assumptions \ref{co4ass1}, \ref{co5ass1} and \ref{co5ass2}--\ref{co5ass3} for quiver categories $\modCQ,\modCQI$. For brevity we give only the $\modCQI$ case, but this includes $\modCQ$ by taking $Q_2=\es$ in Definition \ref{co6def2}. Theorem \ref{co6thm1} gives the final result.

\subsubsection{The data of Assumption \ref{co4ass1}}
\label{co641}

\begin{dfn}
\label{co6def6}
Let $\br Q$ be a quiver with relations, and use the notation of \S\ref{co61}--\S\ref{co62}. Set $\B=\modCQI$, a $\C$-linear abelian (hence exact) category. We will define the data and verify the conditions of Assumption \ref{co4ass1} for~$\B$.

For Assumption \ref{co4ass1}(a), the moduli stack $\brM=\coprod_{\bs d\in\N^{Q_0}}\brM_{\bs d}$ of objects in $\B=\modCQI$ is constructed in Definition \ref{co6def4}. For Assumption \ref{co4ass1}(b), for $\bs c,\bs d\in\N^{Q_0}$ we define a morphism
\begin{align*}
\Phi_{\bs c,\bs d}&:\brM_{\bs c}\t \brM_{\bs d}=[S_{\bs c}^{-1}(0)/\GL_{\bs c}]\t [S_{\bs d}^{-1}(0)/\GL_{\bs d}]\\
&\qquad\qquad \qquad\longra\brM_{\bs c+\bs d}=[S_{\bs c+\bs d}^{-1}(0)/\GL_{\bs c+\bs d}],
\end{align*}
to be induced by the morphisms $S_{\bs c}^{-1}(0)\t S_{\bs d}^{-1}(0)\ra S_{\bs c+\bs d}^{-1}(0)$ mapping $(s_{\bs c},s_{\bs d})\mapsto \bigl(\begin{smallmatrix} s_{\bs c} & 0 \\ 0 & s_{\bs d} \end{smallmatrix}\bigr)$ and $\GL_{\bs c}\t\GL_{\bs d}\ra \GL_{\bs c+\bs d}$ mapping $(A_{\bs c},A_{\bs d})\mapsto \bigl(\begin{smallmatrix} A_{\bs c} & 0 \\ 0 & A_{\bs d} \end{smallmatrix}\bigr)$, in the obvious block matrix notation. Then $\Phi:\brM\t\brM\ra\brM$ is given by $\Phi\vert_{\brM_{\bs c}\t \brM_{\bs d}}=\Phi_{\bs c,\bs d}$. Associativity and commutativity of $\Phi$ are obvious.

Similarly, for Assumption \ref{co4ass1}(c), for $\bs d\in\N^{Q_0}$ we define
\begin{equation*}
\Psi_{\bs d}:[*/\bG_m]\t \brM_{\bs d}=[*/\bG_m]\t[S_{\bs d}^{-1}(0)/\GL_{\bs d}]\longra\brM_{\bs d}=[S_{\bs d}^{-1}(0)/\GL_{\bs d}]
\end{equation*}
to be induced by $\id:S_{\bs d}^{-1}(0)\ra S_{\bs d}^{-1}(0)$ and the morphism $\bG_m\t\GL_{\bs d}\ra\GL_{\bs d}$ mapping $(\la,(A_v)_{v\in Q_0})\mapsto (\la A_v)_{v\in Q_0}$. Then $\Psi:[*/\bG_m]\t\brM\ra\brM$ is given by $\Psi\vert_{[*/\bG_m]\t \brM_{\bs d}}=\Psi_{\bs d}$. Equation \eq{co4eq2} is obvious from the definitions of~$\Psi_{\bs d},\Phi_{\bs c,\bs d}$.

For Assumption \ref{co4ass1}(d), we define $K(\A)=\Z^{Q_0}$, with quotient $K_0(\A)\twoheadrightarrow K(\A)$ mapping $[V,\rho]\mapsto\bdim(V,\rho)$. It is clear that if $E\in\B$ with $\lb E\rb=0$ then $E=0$, and the map $\brM(\C)\ra K(\B)$ taking $E\mapsto\lb E\rb$ is locally constant.

For Assumption \ref{co4ass1}(e), define $\chi:\Z^{Q_0}\t\Z^{Q_0}\ra\Z$ by
\e
\chi(\bs c,\bs d)=\ts\sum\limits_{v\in Q_0}\bs c(v)\bs d(v)-\sum\limits_{e\in
Q_1}\bs c(t(e))\bs d(h(e))+\sum\limits_{q\in Q_2}\bs c(b(q))\bs d(e(q)).
\label{co6eq5}
\e

For Assumption \ref{co4ass1}(f), define a perfect complex $\cE_{\bs c,\bs d}^\bu\ra\brM_{\bs c}\t\brM_{\bs d}$ in degrees $[-2,0]$ for all $\bs c,\bs d\in\N^{Q_0}$ by
\e
\cE_{\bs c,\bs d}^\bu=\Bigl[\xymatrix@C=17pt{
\mathop{\smash{\bigop\limits_{q\in Q_2}}\cV_{b(q)}\bt \cV_{e(q)}^*\,\,\,\,\,\,\,\,}\limits_{-2}\!\!\!\!\!\!\!\! \ar[r]^{\eta} & \mathop{\smash{\bigop\limits_{e\in Q_1}} \cV_{t(e)}\bt \cV_{h(e)}^*\,\,\,\,\,\,\,\,}\limits_{-1}\!\!\!\!\!\!\!\! \ar[r]^(0.55){\ze} & \mathop{\smash{\bigop\limits_{v\in Q_0}} \cV_v\bt \cV^*_v\,\,\,\,\,\,\,\,\,\,\,\,}\limits_0\!\!\!\!\!\!\!\!\!\!\!\!  }\Bigr], \!\!\!\!
\label{co6eq6}
\e
where $\cV_v\ra \brM_{\bs c},\brM_{\bs d}$ for $v\in Q_0$ are as in Definition \ref{co6def3}, and we write
\ea
\begin{split}
\ze&=\ts\sum_{e\in Q_1}\bigl(\rho_{e,\bs c}\bt\id:\cV_{t(e)}\bt \cV_{h(e)}^*\longra\cV_{h(e)}\bt \cV_{h(e)}^*\bigr)\\
&-\ts\sum_{e\in Q_1}\bigl(\id\bt\rho_{e,\bs d}^*:\cV_{t(e)}\bt \cV_{h(e)}^*\longra\cV_{t(e)}\bt \cV_{t(e)}^*\bigr),
\end{split}
\label{co6eq7}\\
\eta&=\ts\sum_{e\in Q_1,\; r\in Q_2}\bigl(\frac{\pd r(q)}{\pd e}(\rho_{\bs c}\ot\rho_{\bs d}):\cV_{b(q)}\bt \cV_{e(q)}^*\longra\cV_{t(e)}\bt \cV_{h(e)}^*\bigr).
\label{co6eq8}
\ea
Here $\rho_{e,\bs c}:\cV_{t(e)}\ra\cV_{h(e)}$ in \eq{co6eq7} is the vector bundle morphism on $\brM_{\bs c}$ such that $\rho_{e,\bs c}\vert_{[(V_v)_{v\in Q_0},(\rho_e)_{e\in Q_1}]}=\rho_e:V_{t(e)}\ra V_{h(e)}$, and $\rho_{e,\bs d}$ is the analogue on~$\brM_{\bs d}$. 

In \eq{co6eq8}, for $q\in Q_2$ we have $r(q)\in\C Q_{b(q)e(q)}$, in the notation of Definition \ref{co6def1}, and for $e\in Q_1$ we have $\frac{\pd r(q)}{\pd e}\in \C Q_{b(q)t(e)}\ot \C Q_{h(e)e(q)}$. That is, $r(q)$ is a $\C$-linear combination of paths $b(q)\ra e(q)$ in $Q$, and $\frac{\pd r(q)}{\pd e}$ is a $\C$-linear combination of tensor products of paths $(b(q)\ra t(e))\ot(h(e)\ra e(q))$. Writing $\rho_{\bs c}=(\rho_{e,\bs c})_{e\in Q_1}$, $\rho_{\bs d}=(\rho_{e,\bs d})_{e\in Q_1}$, by $\frac{\pd r(q)}{\pd e}(\rho_{\bs c}\ot\rho_{\bs d})$ we mean that for $\frac{\pd r(q)}{\pd e}$ in $\C Q_{b(q)t(e)}\ot \C Q_{h(e)e(q)}$, we evaluate the $\C Q_{b(q)t(e)}$ factor on $\rho_{\bs c}$ to get a vector bundle morphism $\cV_{b(q)}\ra\cV_{t(e)}$ on $\brM_{\bs c}$, and evaluate the $\C Q_{h(e)e(q)}$ factor on $\rho_{\bs d}$ to get a vector bundle morphism $\cV_{h(e)}\ra\cV_{e(q)}$ on $\brM_{\bs d}$, which is dual to a morphism $\cV_{e(q)}^*\ra\cV_{h(e)}^*$, so overall we have a morphism $\cV_{b(q)}\bt \cV_{e(q)}^*\ra\cV_{t(e)}\bt \cV_{h(e)}^*$ on $\brM_{\bs c}\t\brM_{\bs d}$. 

We will show that $\ze\ci\eta=0$ in \eq{co6eq6}. Let $q\in Q_2$ and $v\in V$, and consider the component of $\ze\ci\eta$ mapping from the $q$ term to the $v$ term. Taking equation \eq{co6eq2} with $u=b(q)$, $w=e(q)$, applying it to $r(q)\in\C Q_{b(q)e(q)}$, subtracting the two lines of \eq{co6eq2} and rearranging, yields an equation in $\C Q_{b(q)v}\ot\C Q_{ve(q)}$
\ea
&\sum_{e\in Q_1:h(e)=v\!\!\!\!\!\!}((e\cdot-)\ot\id_{\C Q_{ve(q)}})\ci\frac{\pd r(q)}{\pd e}
-\sum_{e\in Q_1:t(e)=v\!\!\!\!\!\!}(\id_{\C Q_{b(q)v}}\ot(-\cdot e))\ci\frac{\pd r(q)}{\pd e}
\nonumber\\
&=\de_{ve(q)}(r(q)\ot\ga_{e(q)})-\de_{b(q)v}(\ga_{b(q)}\ot r(q)).
\label{co6eq9}
\ea
Evaluating this on $\rho_{\bs c}\ot\rho_{\bs d}$ gives
\e
(\ze\ci\eta)_{qv}=\de_{ve(q)}\,r(q)(\rho_{\bs c})\bt\id_{\cV_{e(q)}}-\de_{b(q)v}\,\id_{\cV_{b(q)}}\bt r(q)(\rho_{\bs d}),
\label{co6eq10}
\e 
as an equation in sections of $(\cV_{b(q)}^*\ot \cV_v)\bt(\cV_v^*\ot\cV_{e(q)})\ra\brM_{\bs c}\t\brM_{\bs d}$, since the left hand of \eq{co6eq9} yields $(\ze\ci\eta)_{qv}$ by \eq{co6eq7}--\eq{co6eq8}. But $\brM_{\bs c}\subseteq\M_{\bs c}$ is defined by the equations $r(q)(\rho_{\bs c})=0$ for $q\in Q_2$, and similarly $r(q)(\rho_{\bs d})=0$ on $\brM_{\bs d}$. Hence the right hand side of \eq{co6eq10} is zero, and~$\ze\ci\eta=0$. 

Thus $\cE_{\bs c,\bs d}^\bu$ in \eq{co6eq6} is a perfect complex. Comparing \eq{co6eq5}--\eq{co6eq6} and noting that $\cV_v\ra\brM_{\bs c}$ has rank $\bs c(v)$, we see that $\rank\cE_{\bs c,\bs d}^\bu=\chi(\bs c,\bs d)$. Define $\cE^\bu\ra\brM\t\brM$ by $\cE^\bu\vert_{\brM_{\bs c}\t\brM_{\bs d}}=\cE_{\bs c,\bs d}^\bu$. Equations \eq{co4eq3}--\eq{co4eq6} for $\cE^\bu$ follow from the bilinear structure of \eq{co6eq6} and the simpler equations, as in \eq{co5eq9}--\eq{co5eq10}
\e
\Phi^*(\cV_v)\cong \cV_v\bp\cV_v, \qquad \Psi^*(\cV_v)\cong L_{[*/\bG_m]}\bt\cV_v,
\label{co6eq11}
\e
which are obvious from the definitions of $\cV_v,\Phi,\Psi$.

For Assumption \ref{co4ass1}(g), if $\bs d\in C(\B)=\N^{Q_0}\sm\{0\}$, define a line bundle $L_{\bs d}\ra\brM_{\bs d}$ by $L_{\bs d}=\bigot_{v\in Q_0}\det\cV_v$. Then the second equation of \eq{co6eq11} implies that $\Psi_{\bs d}^*(L_{\bs d})\cong L^{N_{\bs d}}_{[*/\bG_m]}\bt L_{\bs d}$ with $N_{\bs d}=\sum_{v\in Q_0}\bs d(v)>0$, as we have to prove. This completes Assumption \ref{co4ass1} for~$\B=\modCQI$.
\end{dfn}

Theorems \ref{co4thm1} and \ref{co4thm2} now define a vertex algebra structure on $\hat H_*(\brM)$, and a Lie algebra structure on $\check H_*(\brM^\pl)$. Note that as in Definition \ref{co6def3}, for $\B=\modCQ$ we have $H_*(\M)=\bigop_{\bs d\in\N^{Q_0}}H_*(\M_{\bs d})$, where $H_*(\M_{\bs d})$ may be written explicitly as the $\Q$-linear dual of \eq{co6eq3}. As in \cite{Joyc12}, we can write the vertex algebra structure on $H_*(\M)$ explicitly in these coordinates (though with rather complicated formulae). Similarly we can write the Lie algebra $H_*(\M^\pl)$ explicitly. If $\br Q$ has homogeneous relations then $H_*(\brM)\cong H_*(\M)$ and $H_*(\brM^\pl)\cong H_*(\M^\pl)$, and again we can write the vertex and Lie algebra structures on $\hat H_*(\brM),\check H_*(\brM^\pl)$ for $\B=\modCQI$ explicitly in coordinates.

\subsubsection{The data of Assumption \ref{co5ass1}}
\label{co642}

\begin{dfn}
\label{co6def7}
Continue in the situation of Definition \ref{co6def6}, and set $\A=\B=\modCQI$. As $\modCQI$ is of finite length, it is noetherian and artinian. We will define the data and verify the conditions of Assumption \ref{co5ass1} for~$\A,\B$.

Assumption \ref{co5ass1}(a) is immediate as $\A=\B$. For Assumption \ref{co5ass1}(b) we take $K(\A)=\Z^{Q_0}$ as in Definition \ref{co6def6}. Assumption \ref{co5ass1}(c) holds by Definition~\ref{co6def6}. 

For Assumption \ref{co5ass1}(d), recall that in Definition \ref{co6def4} we defined $\brM_{\bs d}$ as the zero locus $s_{\bs d}^{-1}(0)$ of a section $s_{\bs d}:\M_{\bs d}\ra E_{\bs d}$ of a vector bundle $E_{\bs d}\ra\M_{\bs d}$, where $\M_{\bs d}$ is a moduli stack of objects in $\modCQ$, a smooth Artin $\C$-stack. That is, $\brM_{\bs d}$ is the fibre product $\M_{\bs d}\t_{0,E_{\bs d},s_{\bs d}}\M_{\bs d}$ in $\Art_\C$. 

We define the derived moduli stack $\bs\brM_{\bs d}$ to be the same zero locus $s_{\bs d}^{-1}(0)$, but taken in the sense of derived stacks. That is, $\bs\brM_{\bs d}$ is the same fibre product $\M_{\bs d}\t_{0,E_{\bs d},s_{\bs d}}\M_{\bs d}$, but taken in $\DArt_\C$ rather than in $\Art_\C$. Then $t_0(\bs\brM_{\bs d})=\brM_{\bs d}$ as $t_0$ preserves fibre products. Note that as $\M_{\bs d},E_{\bs d}$ are smooth they are already correct as derived stacks, that is, we take $\bs\M_{\bs d}=\M_{\bs d}$ and $\bs E_{\bs d}=E_{\bs d}$. Also $\bs\brM_{\bs d}$ is locally finitely presented and quasi-smooth, as it is a derived fibre product of smooth Artin $\C$-stacks. Similarly, we define $\bs\brM_{\bs d}^\pl$ to be the derived zero locus of $s_{\bs d}^\pl:\M_{\bs d}^\pl\ra E_{\bs d}^\pl$, which is locally finitely presented and quasi-smooth. We set $\bs\brM=\coprod_{\bs d\in\N^{Q_0}}\bs\brM_{\bs d}$ and~$\bs\brM^\pl=\coprod_{\bs d\in\N^{Q_0}}\bs\brM_{\bs d}^\pl$.

If $\A=\modCQ$ then $E_{\bs d},E_{\bs d}^\pl$ are the zero vector bundles, so $\bs\M_{\bs d}=\M_{\bs d}$ and $\bs\M_{\bs d}^\pl=\M_{\bs d}^\pl$ are smooth classical Artin $\C$-stacks, and $\bs\M=\M$, $\bs\M^\pl=\M^\pl$.

Then part (i) of Assumption \ref{co5ass1}(d) is immediate, (ii) holds as the constructions of $\Phi,\Psi$ in \S\ref{co641} and the principal $[*/\bG_m]$-bundle property in \S\ref{co62} also work for derived stacks, and (iii) holds from above. For (iv), by \eq{co6eq6}, $\De_{\brM_{\bs d}}^*(\cE_{\bs d,\bs d}^\bu)[-1]$ is the perfect complex in degrees $[-1,1]$ on $\breve{\mathcal M}_{\bs d}$:
\e
\begin{split}
&\De_{\brM_{\bs d}}^*(\cE_{\bs d,\bs d}^\bu)[-1]=\\
&\Bigl[\xymatrix@C=30pt{
\mathop{\smash{\bigop\limits_{q\in Q_2}}\cV_{b(q)}\ot \cV_{e(q)}^*\,\,\,\,\,\,\,\,}\limits_{-1}\!\!\!\!\!\!\!\! \ar[r]^{\De^*(\eta)} & \mathop{\smash{\bigop\limits_{e\in Q_1}} \cV_{t(e)}\ot \cV_{h(e)}^*\,\,\,\,\,\,\,\,}\limits_{0}\!\!\!\!\!\!\!\! \ar[r]^(0.55){\De^*(\phi)} & \mathop{\smash{\bigop\limits_{v\in Q_0}} \cV_v\ot \cV^*_v\,\,\,\,\,\,\,\,\,\,\,\,}\limits_{1}\!\!\!\!\!\!\!\!\!\!\!\!  }\Bigr]. \!\!\!\!
\end{split}
\label{co6eq12}
\e

Using the description of $\brM_{\bs d}$ in Definition \ref{co6def4} as $\brM_{\bs d}=s_{\bs d}^{-1}(0)$ for $s_{\bs d}$ a section of the vector bundle $E_{\bs d}$ on $\M_{\bs d}=[R_{\bs d}/\GL_{\bs d}]$, we may rewrite \eq{co6eq12} as
\e
\De_{\brM_{\bs d}}^*(\cE_{\bs d,\bs d}^\bu)[-1]=\Bigl[\xymatrix@C=40pt{
\mathop{E_{\bs d}^*}\limits_{-1} \ar[r]^{\d s_{\bs d}} & \mathop{T^*R_{\bs d}}\limits_{0}\ar[r]^(0.55){\vect^*} & \mathop{\gl_{\bs d}^*}\limits_{1} }\Bigr], \!\!\!\!
\label{co6eq13}
\e
where $\gl_{\bs d}$ is the Lie algebra of $\GL_{\bs d}$, and $\d s_{\bs d}\in\Ga(E_{\bs d}\ot T^*R_{\bs d})$ is the derivative of $s_{\bs d}$ (equivalently $S_{\bs d}$), and $\vect:\gl_{\bs d}\ra\Ga(TR_{\bs d})$ is the vector fields of the $\GL_{\bs d}$-action on $R_{\bs d}$.

In \eq{co6eq13}, the right hand part $\bigl[T^*R_{\bs d}\longra \gl_{\bs d}^*\bigr]$ is the cotangent complex $\bL_{\M_{\bs d}}\vert_{\brM_{\bs d}}$ of the smooth Artin stack $\M_{\bs d}=[R_{\bs d}/\GL_{\bs d}]$, restricted to $\brM_{\bs d}\subseteq\M_{\bs d}$. Since $\bs\brM_{\bs d}$ is the derived zero locus of $s_{\bs d}:\M_{\bs d}\ra E_{\bs d}$, general properties of cotangent complexes of derived stacks as in \S\ref{co222} imply that there is a natural quasi-isomorphism $\th_{\bs d}:\De_{\brM_{\bs d}}^*(\cE_{\bs d,\bs d}^\bu)[-1]\ra i_{\bs d}^*(\bL_{\bs\brM_{\bs d}})$, as in~\eq{co5eq1}.

We can prove Assumption \ref{co5ass1}(d)(v) by considering how $\bs\Phi_{\bs c,\bs d}$ acts on cotangent complexes written in the form \eq{co6eq12}--\eq{co6eq13}, and relating this to the explanation of \eq{co4eq3}--\eq{co4eq4} in Definition~\ref{co6def6}. 

For Assumption \ref{co5ass1}(e) we set $C(\B)_\pe=C(\B)$. 

For Assumption \ref{co5ass1}(f) we take $\dM_{\bs d}=\brM_{\bs d}$, $\dM_{\bs d}^\pl=\brM_{\bs d}^\pl$ in (i), $\bs\dM_{\bs d}^\red=\bs\brM_{\bs d}$, $\bs\dM_{\bs d}^\rpl=\bs\brM_{\bs d}^\pl$, and $\bs j_{\bs d}=\bs\id=\bs j_{\bs d}^\pl$ and $\bs{\dot\Pi}{}_{\bs d}^\rpl=\bs\Pi_{\bs d}^\pl$ in (ii), and $U_{\bs d}=o_{\bs d}=0$ for all $\bs d\in C(\B)$ in (iii). Then \eq{co5eq4} is trivially Cartesian, and $\bs\dM_{\bs d}^\red,\bs\dM_{\bs d}^\rpl$ are quasi-smooth from above, proving (ii). Equation \eq{co5eq6} holds as $\bs j_{\bs d},\bs j_{\bs d}^\pl$ are identities and $U_{\bs d}=0$, giving (iii), and (iv) is trivial as~$o_{\bs c}=o_{\bs d}=o_{\bs c+\bs d}=0$.

For Assumption \ref{co5ass1}(g), take $K=\{0\}$, and set $\B_0=\B=\modCQI$, write $F_0:\modCQI\ra\Vect_\C$ for the exact functor mapping $(V,\rho)\mapsto V$ on objects and $\phi\mapsto\phi$ on morphisms, and define $\la_0:K(\B_0)=\Z^{Q_0}\ra\Z$ by $\la_0(\bs d)=\sum_{v\in Q_0}\bs d(v)$. The vector bundle $\cV_0\ra\brM$ in Assumption \ref{co5ass1}(g)(ii) is $\cV_0=\bigop_{v\in Q_0}\cV_v$. Equations \eq{co5eq9}--\eq{co5eq10} follow from \eq{co6eq11}. The rest of Assumption \ref{co5ass1}(g) is obvious. This completes Assumption~\ref{co5ass1}.
\end{dfn}

\subsubsection{The data of Assumptions \ref{co5ass2}--\ref{co5ass3}}
\label{co643}

\begin{dfn}
\label{co6def8}
Continue in the situation of Definitions \ref{co6def6}--\ref{co6def7}. We will define the data and verify the conditions of Assumptions \ref{co5ass2}(a)--(f) and \ref{co5ass3}. If $\br Q$ has no oriented cycles then we will also prove Assumption \ref{co5ass2}(g),(h) using Proposition \ref{co6prop1}. Then Assumptions \ref{co5ass1}--\ref{co5ass3} hold, and Theorems \ref{co5thm1}--\ref{co5thm3} apply.

In Assumption \ref{co5ass2}, we define $\sS$ to be the set of all slope stability conditions $(\mu,\R,\le)$ on $\A=\B=\modCQI$ defined in Definition \ref{co6def5} using constants $\mu_v\in\R$ for $v\in Q_0$. For Assumption \ref{co5ass2}(a), as $\modCQI$ is artinian it is $\mu$-artinian for any $(\mu,\R,\le)$. Parts (b),(c),(e) are obvious. 

For Assumption \ref{co5ass2}(d), let $(\mu,\R,\le),(\ti\mu,\R,\le)\in\sS$, and $I\subseteq C(\B)_\pe=\N^{Q_0}\sm\{0\}$ be a finite subset, and $\al\in I$, satisfy $\mu(\be)=\mu(\al)$ and $\brM_\be^\ss(\mu)\ne\es$ for all $\be\in I$. Define $\la:K(\A)=\Z^{Q_0}\ra\R$ by
\begin{equation*}
\la(\bs d)=\sum_{v\in Q_0}\biggl(\ti\mu_v-\frac{\sum_{w\in Q_0}\ti\mu_w\al(w)}{\sum_{w\in Q_0}\al(w)}\biggr)\bs d(v),
\end{equation*}
where $(\ti\mu,\R,\le)$ is defined using $\ti\mu_v\in\R$ for $v\in Q_0$. Then $\la(\al)=0$, and
\begin{equation*}
\ti\mu(\bs d)=\frac{\la(\bs d)}{\sum_{w\in Q_0}\bs d(w)}+\frac{\sum_{w\in Q_0}\ti\mu_w\al(w)}{\sum_{w\in Q_0}\al(w)},
\end{equation*}
for $\bs d\in \N^{Q_0}\sm\{0\}$. Hence for $\be\in I$ we have 
\begin{align*}
&\la(\be)>0 \quad\Longleftrightarrow \quad	\frac{\la(\be)}{\sum_{w\in Q_0}\be(w)}>0=\frac{\la(\al)}{\sum_{w\in Q_0}\al(w)} \quad \Longleftrightarrow \\ 
&\frac{\la(\be)}{\sum_{w\in Q_0}\be(w)}+\frac{\sum_{w\in Q_0}\ti\mu_w\al(w)}{\sum_{w\in Q_0}\al(w)}>
\frac{\la(\al)}{\sum_{w\in Q_0}\al(w)}+\frac{\sum_{w\in Q_0}\ti\mu_w\al(w)}{\sum_{w\in Q_0}\al(w)} \\ 
&\Longleftrightarrow\quad  \ti\mu(\be)>\ti\mu(\al),
\end{align*}
as we have to prove. Similarly $\la(\be)<0$ $\Longleftrightarrow$ $\ti\mu(\be)<\ti\mu(\al)$.

For Assumption \ref{co5ass2}(f), define $\rk:C(\A)=\N^{Q_0}\sm\{0\}\ra\N_{>0}$ by $\rk\bs d=\sum_{v\in Q_0}\bs d(v)$. Then~$\rk(\bs c+\bs d)=\rk \bs c+\rk\bs d$.

As above, we only prove Assumption \ref{co5ass2}(g),(h) when $\br Q$ has no oriented cycles. In this case Assumption \ref{co5ass2}(g) is immediate from Proposition~\ref{co6prop1}.

For Assumption \ref{co5ass2}(h), suppose $\br Q$ has no oriented cycles. Let $\baB\subseteq\baA,\ab K(\baA),\ab\baM,\ab\baM^\pl,\ab\ldots$ be constructed in Definition \ref{co5def1} starting from $\B\subseteq\A,\ldots$ and another quiver $\bar Q=(\bar Q_0,\bar Q_1,\bar h,\bar t)$, with $\bar Q_0=\dot{\bar Q}_0\amalg\ddot{\bar Q}_0$ and the unique map $\bar\ka:\ddot{\bar Q}\ra K=\{0\}$. Unwinding the definitions, we find that $\baB=\baA=\modCaQI$ is the abelian category of representations of a quiver with relations $\ac Q=(\ac Q_0,\ac Q_1,\ac h,\ac t,\ac Q_2,\ac b,\ac e,\ac r)$ constructed from $\br Q$ and $\bar Q$ as follows:
\begin{itemize}
\setlength{\itemsep}{0pt}
\setlength{\parsep}{0pt}
\item[(a)] The vertices are $\ac Q_0=Q_0\amalg\dot{\bar Q}_0$.
\item[(b)] The edges $\ac Q_1,\ac h,\ac t$ are of three kinds:
\begin{itemize}
\setlength{\itemsep}{0pt}
\setlength{\parsep}{0pt}
\item[(i)] Edges in $\br Q$.
\item[(ii)] Edges in $\bar Q$ starting and finishing in $\dot{\bar Q}_0$.
\item[(iii)] For each edge $\overset{\bar v}{\bu}\,{\buildrel\bar e\over \longra}\,\overset{\bar w}{\bu}$ in $\bar Q$ with $\bar v\in\dot{\bar Q}_0$ and $\bar w\in\ddot{\bar Q}_0$, which has $\bar\ka(\bar w)=0$, we get an edge $\overset{\bar v}{\bu}\,{\buildrel\bar e_v\over \longra}\,\overset{v}{\bu}$ in $\ac Q$ for each $v\in Q_0$. 
\end{itemize}
\item[(c)] The relations $\ac Q_2,\ac b,\ac e,\ac r$ are $Q_2,b,e,r$, and involve only the subquiver $Q=(Q_0,Q_1,h,t)$ in $(\ac Q_0,\ac Q_1,\ac h,\ac t)$.
\end{itemize}

Here the point of (b)(iii) is that objects $((V,\rho),\bs{\bar V},\bs{\bar\rho})$ of $\baA$ include, for each $\overset{\bar v}{\bu}\,{\buildrel\bar e\over \longra}\,\overset{\bar w}{\bu}$ as in (b)(iii), a linear map $\rho_{\bar e}:\bar V_{\bar v}\ra F_0(V,\rho)=V=\bigop_{v\in Q_0}V_v$. This is equivalent to linear maps $\rho_{\bar e_v}:\bar V_{\bar v}\ra V_v$ for each $v\in Q_0$, which are the corresponding edge maps for a representation of $\ac Q$. Observe that as $\br Q$ and $\bar Q$ have no oriented cycles, $\ac Q$ also has no oriented cycles.

Let $(\mu,\R,\le)\in\sS$ be defined using $\mu_v\in\R$ for $v\in Q_0$, and $\la:K(\A)=\Z^{Q_0}\ra\R$ be a morphism, and $\bs{\bar\mu}\in\R^{\dot{\bar Q}_0}$. Then Definition \ref{co5def1} defines a weak stability condition $(\mu^\la_{\bs{\bar\mu}},[-\iy,\iy]\t\R,\le)$ on $\baA=\modCaQI$. This is not a slope stability condition on $\modCaQI$ in the sense of Definition \ref{co6def5}. However, we can write $\mu^\la_{\bs{\bar\mu}}$-(semi)stability as a limit of slope (semi)stability on $\modCaQI$. 

Fix $(\bs d,\bs{\bar d})\in C(\baA)=(\N^{Q_0}\t\N^{\dot{\bar Q}_0})\sm\{(0,0)\}$, and if $\bs d\ne 0\ne\bs{\bar d}$ define a slope stability condition $(\ac\mu_\ep,\R,\le)$ on $\modCaQI$ for $\ep>0$ by
\begin{equation*}
\ac\mu_\ep(\bs c,\bs{\bar c})=\frac{\ep\la(\bs c)+\sum_{v\in Q_0}(\mu_v-\mu(\bs d)-\ep C)\bs c(v)+\ep\sum_{\bar v\in\dot{\bar Q}_0}\bar\mu_{\bar v}\bs{\bar c}(v)}{\sum_{v\in Q_0}\bs c(v)+\sum_{\bar v\in\dot{\bar Q}_0}\bs{\bar c}(v)}\,,
\end{equation*}
where $C=\bigl(\la(\bs d)+\sum_{\bar v\in\dot{\bar Q}_0}\bar\mu_{\bar v}\bs{\bar d}(v)\bigr)\big/\bigl(\sum_{v\in Q_0}\bs d(v)\bigr)$. One can now show that if $((V,\rho),\bs{\bar V},\bs{\bar\rho})\in\baA$ lies in class $(\bs d,\bs{\bar d})$ then $((V,\rho),\bs{\bar V},\bs{\bar\rho})$ is $\mu^\la_{\bs{\bar\mu}}$-(semi)stable if and only if $((V,\rho),\bs{\bar V},\bs{\bar\rho})$ is $\ac\mu_\ep$-(semi)stable for sufficiently small $\ep>0$. If $\bs d=0$ or $\bs{\bar d}=0$ we can do the same for a different, simpler definition of~$(\ac\mu_\ep,\R,\le)$.

As there are only finitely many splittings $(\bs d,\bs{\bar d})=(\bs c,\bs{\bar c})+(\bs e,\bs{\bar e})$ in $C(\baA)$ we can choose $\ep$ uniformly on $\baM^\pl_{(\bs d,\bs{\bar d})}$. Thus $\baM^\rst_{(\bs d,\bs{\bar d})}(\mu^\la_{\bs{\bar\mu}})=\baM^\rst_{(\bs d,\bs{\bar d})}(\ac\mu_\ep)$ and $\baM^\ss_{(\bs d,\bs{\bar d})}(\mu^\la_{\bs{\bar\mu}})=\baM^\ss_{(\bs d,\bs{\bar d})}(\ac\mu_\ep)$ for small $\ep>0$. Therefore Proposition \ref{co6prop1} for $\modCaQI$, $(\ac\mu_\ep,\R,\le)$ shows that if $\baM^\rst_{(\bs d,\bs{\bar d})}(\mu^\la_{\bs{\bar\mu}})=\baM^\ss_{(\bs d,\bs{\bar d})}(\mu^\la_{\bs{\bar\mu}})$ then $\baM^\ss_{(\bs d,\bs{\bar d})}(\mu^\la_{\bs{\bar\mu}})$ is a proper algebraic space, proving 
Assumption~\ref{co5ass2}(h).

We now again allow $\ac Q$ to have oriented cycles. For Assumption \ref{co5ass3}(a), if $(\mu,\R,\le),(\ti\mu,\R,\le)\in\sS$ we define $\mu_t=(1-t)\mu+t\ti\mu:C(\baA)\ra\R$ for $t\in[0,1]$. Then $(\mu_t,\R,\le)_{t\in[0,1]}$ is a continuous family of slope stability conditions on $\baA$, so that $(\mu_t,\R,\le)\in\sS$, with $(\mu_0,\R,\le)=(\mu,\R,\le)$ and $(\mu_1,\R,\le)=(\ti\mu,\R,\le)$, as we have to prove. For Assumption \ref{co5ass3}(b), if $\al\in C(\B)_\pe=\N^{Q_0}\sm\{0\}$ there are only finitely many splittings $\al=\al_1+\cdots+\al_n$ for $\al_1,\ldots,\al_n\in C(\B)=\N^{Q_0}\sm\{0\}$, so the finiteness claim follows, and $\al_{i_1}+\cdots+\al_{i_2}\in C(\B)_\pe$ for $1\le i_1\le i_2\le n$ is automatic as $C(\B)_\pe=C(\B)$. This completes Assumption~\ref{co5ass3}.
\end{dfn}
	
\subsubsection{Alternative proof of Assumption \ref{co5ass2}(g),(h) using \S\ref{co33}}
\label{co644}

In Definition \ref{co6def8} we showed that Assumption \ref{co5ass2}(g),(h) hold for slope stability conditions on $\A=\B=\modCQI$ when $\br Q$ has no oriented cycles, using Proposition \ref{co6prop1}. We now give an alternative proof using the results of \S\ref{co33}.

\begin{prop}
\label{co6prop2}
In the situation of Definition\/ {\rm\ref{co6def6},} suppose $\br Q$ has no oriented cycles. Then $\A=\modCQI$ is \begin{bfseries}of compact type\end{bfseries} in the sense of\/~{\rm\S\ref{co332}}.
\end{prop}

\begin{proof} The category $\tiA=\modCQI{}^{\rm \,inf}$ of representations $(\bs V,\bs\rho)=((V_v)_{v\in Q_0},\ab (\rho_e)_{e\in Q_1})$ of $\ac Q$ on possibly infinite-dimensional $\C$-vector spaces $V_v$ is a cocomplete, locally noetherian $\C$-linear abelian category, and the subcategory $\modCQI\subset\modCQI{}^{\rm \,inf}$ of objects with the $V_v$ finite-dimensional is the full subcategory of compact objects in $\modCQI{}^{\rm \,inf}$. 

Here the condition that $\br Q$ has no oriented cycles is needed for compact objects (equivalently, noetherian objects in the sense of Definition \ref{co3def9}) in $\modCQI{}^{\rm \,inf}$ to coincide with finite-dimensional representations. To see this, note that in a finite-dimensional representation $(\bs V,\bs\rho)$, considering dimensions of subrepresentations shows ascending chains of subobjects must stabilize, so $(\bs V,\bs\rho)$ is noetherian. And for infinite-dimensional $(\bs V,\bs\rho)$, if $\br Q$ has no oriented cycles we can construct an ascending chain of finite-dimensional subobjects with strictly increasing dimensions, so $(\bs V,\bs\rho)$ is not noetherian.

The moduli stack $\M$ of objects in $\modCQI$ is as described in \S\ref{co62}, and is an Artin $\C$-stack locally of finite type. Hence $\A=\modCQI$ is of compact type, as in Definition~\ref{co3def9}.
\end{proof}

Given a stability condition $(\mu,\R,\le)$ on $\A=\modCQI$ and $\bs d\in C(\A)$, the next definition uses results of Reineke \cite[\S 3]{Rein} (see also Hoskins \cite[\S 3]{Hosk}) to construct a pseudo-$\Th$-stratification of $\M_{\bs d}$ with semistable locus $(\Pi_{\bs d}^\pl)^{-1}(\M_{\bs d}^\ss(\mu))$, in the sense of \S\ref{co334}, via Harder--Narasimhan filtrations from Theorem~\ref{co3thm1}.  

\begin{dfn}
\label{co6def9}
Let $\br Q$ be a quiver with relations, and $(\mu,\R,\le)$ be a slope stability condition on $\A=\modCQI$ as in Definition \ref{co6def5}, defined using constants $\mu_v\in\R$ for $v\in Q_0$. Initially, following \cite{Hosk,Rein}, we suppose that $\mu_v\in\Z$, and explain at the end of the definition how to extend to $\mu_v\in\R$. Fix~$\bs d\in C(\A)$.

Define the set $\HNT(\br Q)$ of {\it Harder--Narasimhan types\/} for $\br Q$ to consist of all $n$-tuples $\bs e=(\bs e_1,\ldots,\bs e_n)$ for $n\ge 1$, where $\bs e_i\in\N^{Q_0}\sm\{0\}$ and $\mu(\bs e_1)>\cdots>\mu(\bs e_n)$. We write $\HNT_{\bs d}(\br Q)$ for the subset of $\bs e$ in $\HNT(\br Q)$ with~$\bs e_1+\cdots+\bs e_n=\bs d$. 

If $\bs e=(\bs e_1,\ldots,\bs e_n)\in\HNT_{\bs d}(\br Q)$, define $x_0,\ldots,x_n\in\R^2$ by
\begin{equation*}
x_j=\ts\bigl(\sum_{i=1}^j\sum_{v\in Q_0}e_{v,i},\sum_{i=1}^j\sum_{v\in Q_0}\mu_ve_{v,i}\bigr).
\end{equation*}
Define the {\it Harder--Narasimhan polygon\/} $\HNP(\bs e)$ to be the union of the line segments in $\R^2$ joining $x_{i-1}$ to $x_i$ for $i=1,\ldots,n$. This is a piecewise-linear curve with end-points $x_0=(0,0)$ and $x_n=\bigl(\sum_{v\in Q_0}d_v,\sum_{v\in Q_0}\mu_vd_v\bigr)$. As the $i^{\rm th}$ line segment has slope $\mu(\bs e_i)$ with $\mu(\bs e_1)>\cdots>\mu(\bs e_n)$, it is a convex polygon. It is the graph $\Ga_{\mathop{\rm hnp}(\bs e)}$ of a continuous, piecewise linear function~$\mathop{\rm hnp}(\bs e):[0,\sum_{v\in Q_0}d_v]\ra\R$.

As in Reineke \cite[Def.~3.6]{Rein}, define a partial order $\preceq$ on $\HNT_{\bs d}(\br Q)$ by $\bs e\preceq\bs e'$ if $\HNP(\bs e)$ lies below $\HNP(\bs e')$ in $\R^2$, that is, if $\mathop{\rm hnp}(\bs e)\le\mathop{\rm hnp}(\bs e')$ as functions $[0,\sum_{v\in Q_0}d_v]\ra\R$. The minimal element in $\HNT_{\bs d}(\br Q)$ under $\preceq$ is $(\bs d)$, so we set $0=(\bs d)$ in Definition~\ref{co3def11}.

Since $\modCQI$ is artinian, for each $(\bs V,\bs\rho)\in\modCQI$ with $\lb\bs V,\bs\rho\rb=\bs d$, by Theorem \ref{co3thm1} we have a $\mu$-Harder--Narasimhan filtration $0=(\bs V_0,\bs\rho_0)\subsetneq (\bs V_1,\bs\rho_1)\subsetneq\cdots\subsetneq (\bs V_n,\bs\rho_n)=(\bs V_n,\bs\rho_n)$ for $n>0$, such that\/ $(\bs W_i,\bs\si_i)=(\bs V_i,\bs\rho_i)/(\bs V_{i-1},\bs\rho_{i-1})$ is $\mu$-semistable for $i=1,\ldots,n,$ and $\mu(\bdim(\bs W_1,\bs\si_1))>\cdots>\mu(\bdim(\bs W_n,\bs\si_n))$. Define the {\it Harder--Narasimhan type\/} of $(\bs V,\bs\rho)$ to be
\begin{equation*}
\HNT(\bs V,\bs\rho)=\bigl(\bdim(\bs W_1,\bs\si_1),\ldots,\bdim(\bs W_n,\bs\si_n)\bigr)\in\HNT_{\bs d}(\br Q).
\end{equation*}

As in Reineke \cite[Prop.s 3.4 \& 3.7]{Rein}, there is a natural constructible function $\HNT_{\bs d}:\M_{\bs d}\ra\HNT_{\bs d}(\br Q)$ which maps a $\C$-point $[\bs V,\bs\rho]$ to $\HNT(\bs V,\bs\rho)$. Furthermore, $\HNT_{\bs d}$ is upper semicontinuous in the partial order $\preceq$ on $\HNT_{\bs d}(\br Q)$. This implies that for each $\bs e\in\HNT_{\bs d}(\br Q)$, we have an open substack
\e
\M_{\bs d,\preceq\bs e}=\bigl\{[\bs V,\bs\rho]\in\M_{\bs d}:\HNT_{\bs d}([\bs V,\bs\rho])\preceq\bs e\bigr\}.
\label{co6eq14}
\e

As in \cite[Prop.~3.4]{Rein} and \cite[\S 3.2]{Hosk}, there is a natural locally closed substack $\M_{\bs d,\bs e}\subseteq\M_{\bs d,\preceq\bs e}\subseteq\M_{\bs d}$, such that in open substacks of $\M_{\bs d}$ we have
\begin{equation*}
\M_{\bs d,\preceq\bs e}\sm\M_{\bs d,\bs e}=\ts\bigcup_{\begin{subarray}{l} \bs e'\in\HNT_{\bs d}(\br Q):\bs e'\prec\bs e\end{subarray}}\M_{\bs d,\preceq\bs e'}.
\end{equation*}
Furthermore, $\M_{\bs d,\bs e}$ is a moduli stack for flat families of $\C Q/I$-representations with $\mu$-Harder--Narasimhan filtrations of type $\bs e$. We can show as in Halpern-Leistner \cite{Halp1,Halp2} that $\M_{\bs d,\bs e}$ is the image of a $\Th$-stratum of $\M_{\bs d,\preceq\bs e}$. Hence $(\M_{\bs d,\preceq\bs e})_{\bs e\in\HNT_{\bs d}(\br Q)}$ is a {\it pseudo-$\Th$-stratification\/} of $\M_{\bs d}$, in the sense of \S\ref{co334}, with semistable locus $(\Pi_{\bs d}^\pl)^{-1}(\M_{\bs d}^\ss(\mu))\subseteq\M_{\bs d}$. Note that the $\M_{\bs d,\preceq\bs e}$ are all of finite type, as $\M_{\bs d}$ is.

So far we have assumed that $\mu_v\in\Z$. Since $\mu$-semistability and $\mu$-Harder--Narasimhan filtrations are unchanged by rescaling all $\mu_v$ by a positive factor, the construction above extends immediately to $\mu_v\in\Q$. Fix $(\bs V,\bs\rho)\in\modCQI$ with $\bdim(\bs V,\bs\rho)=\bs d$. Then as $(\mu,\R,\le)$ varies in the space of slope stability conditions on $\modCQI$, the $\mu$-Harder--Narasimhan filtration of $(\bs V,\bs\rho)$ depends only on whether $\mu(\bs e)<\mu(\bs f)$ or $\mu(\bs e)=\mu(\bs f)$ or $\mu(\bs e)>\mu(\bs f)$ for a finite set of pairs $(\bs e,\bs f)$ of dimension vectors of subquotients of $(\bs V,\bs\rho)$. 

Now $\mu(\bs e)=\mu(\bs f)$ is a linear equation on $(\mu_v)_{v\in Q_0}$ in $\R^{Q_0}$ with rational coefficients. Therefore there is a {\it finite wall and chamber decomposition\/} of $\R^{Q_0}$ (see Theorem \ref{co7thm2} below for more on this), with all walls defined by rational linear equations, such that the $\mu$-Harder--Narasimhan filtration of any $(\bs V,\bs\rho)\in\modCQI$ with $\bdim(\bs V,\bs\rho)=\bs d$ depends only on $(\bs V,\bs\rho)$ and on the stratum of the decomposition of $\R^{Q_0}$ containing~$(\mu_v)_{v\in Q_0}$. 

Hence the stratification $(\M_{\bs d,\preceq\bs e})_{\bs e\in\HNT_{\bs d}(\br Q)}$ determined by $\mu$ depends only on the stratum of $\R^{Q_0}$ containing $(\mu_v)_{v\in Q_0}$. Since the walls are defined by rational linear equations, every stratum contains a rational point. That is, given $(\mu,\R,\le)$ with $(\mu_v)_{v\in Q_0}$ in $\R^{Q_0}$ we can find a nearby $(\mu',\R,\le)$ with $(\mu'_v)_{v\in Q_0}$ in $\Q^{Q_0}$ such that the $\mu$- and $\mu'$-Harder--Narasimhan filtrations of $(\bs V,\bs\rho)$ agree for all $(\bs V,\bs\rho)$ in $\modCQI$ with $\bdim(\bs V,\bs\rho)=\bs d$. Thus the pseudo-$\Th$-stratification of $\M_{\bs d}$ defined above for $\mu'$ also works for $\mu$, as we want.
\end{dfn}

\begin{prop}
\label{co6prop3}
In the situation of Definitions\/ {\rm\ref{co6def6}--\ref{co6def8},} the slope stability conditions $(\mu,\R,\le)$ from Definition\/ {\rm\ref{co6def5}} are \begin{bfseries}additive\end{bfseries} in the sense of\/~{\rm\S\ref{co335}}.
\end{prop}

\begin{proof}
Fix $\bs d\in C(\A)=C(\modCQI)$, and define $\rho:K(\A)\ra V=\R$ by
\begin{equation*}
\rho(\bs e)=\ts\bigl(\sum_{v\in Q_0}\bs d(v)\bigr)\bigl(\sum_{v\in Q_0}\mu_v\bs e(v)\bigr)-\bigl(\sum_{v\in Q_0}\bs e(v)\bigr)\bigl(\sum_{v\in Q_0}\mu_v\bs d(v)\bigr).
\end{equation*}
Then $\rho$ is a group morphism with $\rho(\bs d)=0$. It is now easy to see that $(\bs V,\bs\rho)$ in $\modCQI$ with $\lb\bs V,\bs\rho\rb=\bs d$  is $\rho$-semistable in the sense of Definition \ref{co3def12} if and only if it is $\mu$-semistable in the sense of Definition \ref{co6def5}. Thus $\rho$-semistability is an open condition on $[\bs V,\bs\rho]$, as $\mu$-semistability is, and $(\Pi_{\bs d}^\pl)^{-1}(\M_{\bs d}^\ss(\mu))=\M_{\bs d}^{\ss,\rho}$. Hence $(\mu,\R,\le)$ is additive by Definition~\ref{co3def13}.
\end{proof}

Combining Theorem \ref{co3thm9}, Propositions \ref{co6prop2}, \ref{co6prop3} and Definition \ref{co6def9} yields:

\begin{cor}
\label{co6cor1}
In the situation of Definitions\/ {\rm\ref{co6def6}--\ref{co6def8},} for\/ $\br Q$ with no oriented cycles, the\/ $\M_{\bs d}^\ss(\mu)$ admit proper good moduli spaces. If also\/ $\M_{\bs d}^\rst(\mu)=\M_{\bs d}^\ss(\mu)$ then $\M_{\bs d}^\ss(\mu)$ is a proper algebraic space.
\end{cor}

This implies that Assumption \ref{co5ass2}(g) holds, and then Assumption \ref{co5ass2}(h) holds as in Definition~\ref{co6def8}.

\subsubsection{Summary of results}
\label{co645}

The following theorem summarizes the work of \S\ref{co641}--\S\ref{co644}.

\begin{thm}
\label{co6thm1}
In the situation above, with\/ $\A=\B=\modCQ$ or $\modCQI$ for a quiver $Q$ or quiver with relations $\br Q,$ the data in Definitions\/ {\rm\ref{co6def6}--\ref{co6def8}} satisfies Assumptions\/ {\rm\ref{co5ass1}--\ref{co5ass3},} except perhaps properness in Assumption\/ {\rm\ref{co5ass2}(g),(h)}. If\/ $Q$ or $\br Q$ has \begin{bfseries}no oriented cycles\end{bfseries} then Assumption\/ {\rm\ref{co5ass2}(g),(h)} also hold. 	

Hence Theorems\/ {\rm\ref{co5thm1}--\ref{co5thm3}} apply to\/ $\A=\B=\modCQ$ or\/ $\modCQI$ when\/ $Q$ or\/ $\br Q$ has no oriented cycles, giving enumerative invariants and wall-crossing formulae for $\B$. For\/ $\modCQ$ this reproves Gross--Joyce--Tanaka\/~{\rm\cite[Th.~5.8]{GJT}}.
\end{thm}

\begin{rem}
\label{co6rem3}
Apart from requiring a quiver with relations $\br Q$ to have no oriented cycles, the author is not aware of any useful general conditions on $\br Q$ which ensure properness of moduli spaces in Assumption~\ref{co5ass2}(g),(h).

However, if for some particular quiver with relations $\br Q$ which does have oriented cycles, and some choice of $C(\B)_\pe\subseteq C(\B)$, we can prove by hand that Assumption \ref{co5ass2}(g),(h) hold, then Theorem \ref{co6thm1} shows that Theorems \ref{co5thm1}--\ref{co5thm3} apply to $\A=\B=\modCQI$. This might happen if, for example, we had an inclusion $\modCQI\hookra D^b\coh_\cs(X)$ for some smooth quasiprojective $\C$-scheme $X$, and properness of moduli spaces followed from geometry of sheaves on~$X$.
\end{rem}

\subsection{Examples}
\label{co65}

\subsubsection{Explicit invariants for `increasing' slope stability conditions}
\label{co651}

The next definition follows Gross--Joyce--Tanaka \cite[Def.~5.5]{GJT}. 

\begin{dfn}
\label{co6def10}
Let $\br Q$ be a quiver with relations, and $(\mu,\R,\le)$ a slope stability condition on $\modCQI$ defined using $\mu_v\in\R$ for $v\in Q_0$. We call $\mu$ {\it increasing\/} if for all edges $\overset{v}{\bu}\,{\buildrel e\over \longra}\,\overset{w}{\bu}$ in $Q$ we have $\mu_v<\mu_w$. Such $\mu$ exist if and only if $\br Q$ has {\it no oriented cycles}.
\end{dfn}

We generalize \cite[Prop.~5.6 \& Th.~5.8(iii)]{GJT} to quivers with relations. Here the condition that $\ac Q$ has {\it no single-vertex relations\/} is as in Definition \ref{co6def2}.

\begin{thm}
\label{co6thm2}
Let\/ $\br Q$ be a quiver with relations with \begin{bfseries}no oriented cycles\end{bfseries} and \begin{bfseries}no single-vertex relations\end{bfseries}, and\/ $(\mu,\R,\le)$ be an \begin{bfseries}increasing\end{bfseries} slope stability condition on $\modCQI$. Then for each\/ $\bs d\in\N^{Q_0}\sm\{0\},$ either:
\begin{itemize}
\setlength{\itemsep}{0pt}
\setlength{\parsep}{0pt}
\item[{\bf(a)}] $\bs d=\de_v$ for some $v\in Q_0,$ that is, $\bs d(v)=1$ and\/ $\bs d(w)=0$ for $w\ne v$. Then $\brM^\rst_{\bs d}(\mu)=\brM^\ss_{\bs d}(\mu)$ is a single point\/~$*$.
\item[{\bf(b)}] $\bs d=n\de_v$ for some $v\in Q_0$ and\/ $n>1$. Then $\brM^\rst_{\bs d}(\mu)=\es$ and\/ $\brM^\ss_{\bs d}(\mu)\cong[*/\PGL(n,\C)]$. Note that\/ $2-2\chi(\bs d,\bs d)=2-2n^2<0$ in this case.
\item[{\bf(c)}] $\bs d\ne n\de_v$ for any\/ $v\in Q_0$ and\/ $n\ge 1$. Then $\brM^\rst_{\bs d}(\mu)=\brM^\ss_{\bs d}(\mu)=\es$.
\end{itemize}

The invariants $[\brM_{\bs d}^\ss(\mu)]_\inv$ of Theorem\/ {\rm\ref{co6thm1}} for $\A=\B=\modCQI$ are given for increasing $(\mu,\R,\le)$ by
\e
[\brM_{\bs d}^\ss(\mu)]_\inv=\begin{cases} 1\in H_0(\brM_{\bs d}^\pl)\cong\Q, & \bs d=\de_v,\; v\in Q_0, \\ 0, & \text{otherwise.} 	
\end{cases}
\label{co6eq15}
\e
Note that for any other slope stability condition $(\ti\mu,\R,\le)$ on $\modCQI,$ we can compute $[\brM_{\bs d}^\ss(\ti\mu)]_\inv$ from the $[\brM_{\bs c}^\ss(\mu)]_\inv$ in \eq{co6eq15} by\/~{\rm\eq{co5eq33}--\eq{co5eq34}}.
\end{thm}

\begin{proof} Parts (a)--(c) are proved as for $\A=\modCQ$ in \cite[Prop.~5.6]{GJT}. That $\br Q$ has no single-vertex relations ensures the relations do not change the moduli stacks in (a),(b), or their virtual dimensions. As in \cite[Th.~5.8(iii)]{GJT}, equation \eq{co6eq15} follows by Theorem \ref{co5thm1}(i) in cases (a),(c), and because $[\brM_{\bs d}^\ss(\mu)]_\inv$ lies in $H_{2o_{\bs d}+2-2\chi(\bs d,\bs d)}(\brM_{\bs d}^\pl)$ with $o_{\bs d}=0$ and $2-2\chi(\bs d,\bs d)<0$ in case~(b).
\end{proof}

\subsubsection{Examples with equivalences $D^b\modCQI\simeq D^b\coh(X)$}
\label{co652}

Quivers with relations are used to give explicit descriptions of derived categories using {\it exceptional collections}.

\begin{dfn}
\label{co6def11}
Let $\T$ be a $\C$-linear triangulated category, for example, $\T=D^b\coh(X)$ for $X$ a smooth projective $\C$-scheme. A finite ordered sequence $E_1,\ldots,E_n$ of objects of $\T$ is called an {\it exceptional collection\/} if
\begin{itemize}
\setlength{\itemsep}{0pt}
\setlength{\parsep}{0pt}
\item[(i)] $\Hom(E_i,E_i)=\C$ and $\Hom(E_i,E_i[k])=0$ if $0\ne k\in\Z$, all $i=1,\ldots,n$.
\item[(ii)] $\Hom(E_i,E_j[k])=0$ if $i>j$ and $k\in\Z$.
\end{itemize}
We call $E_1,\ldots,E_n$ {\it strong\/} if also
\begin{itemize}
\setlength{\itemsep}{0pt}
\setlength{\parsep}{0pt}
\item[(iii)] $\Hom(E_i,E_j[k])=0$ if $0\ne k\in\Z$ for all $i,j=1,\ldots,n$.
\end{itemize}
We call $E_1,\ldots,E_n$ {\it full\/} if it generates $\T$ as a triangulated category.
\end{dfn}

Bondal \cite{Bond} proves:

\begin{thm}
\label{co6thm3}
Suppose $X$ is a smooth projective $\C$-scheme and\/ $E_1,\ldots,E_n$ is a full strong exceptional collection in $D^b\coh(X)$. Then there exists a quiver with relations $\br Q$ constructed from $E_1,\ldots,E_n$ as in {\rm\cite[\S 5]{Bond}} with an equivalence of triangulated categories $D^b\coh(X)\simeq D^b\modCQI$.	
\end{thm}

This suggests the following approach to enumerative invariants in $\coh(X)$:

\begin{meth}
\label{co6meth1}
Let $X$ be a smooth projective $\C$-scheme and $E_1,\ldots,E_n$ be a full strong exceptional collection in $D^b\coh(X)$, so that Theorem \ref{co6thm3} gives a quiver with relations $\br Q$ and an inclusion $\modCQI\hookra D^b\coh(X)$. 

Suppose the enumerative invariant theory of \S\ref{co51}--\S\ref{co53} works for $\coh(X)$, with $\cE^\bu\cong(\cExt^\bu)^\vee$ in Assumption \ref{co4ass1}(f) and $U_\al=o_\al=0$ for all $\al$ in Assumption \ref{co5ass1}(f). See Chapter \ref{co7} for more details.

Suppose also that when we extend the Ext complex $\cExt^\bu$ for $\coh(X)$ to $\baM\t\baM$, where $\baM$ is the moduli stack of objects in $D^b\coh(X)$, and then restrict to $\M_{\br Q}\t\M_{\br Q}$ for $\M_{\br Q}\subset\baM$ the moduli stack of objects in $\modCQI\subset D^b\coh(X)$, we get $(\cE^\bu)^\vee$ for $\cE^\bu$ as defined in \S\ref{co641} for $\br Q$. This is only possible if $\cExt^\bu\vert_{\M_{\br Q}\t\M_{\br Q}}$ is perfect in $[0,2]$, and is a strong condition.

We can then attempt to compute enumerative invariants for (say) Gieseker semistable sheaves in $\coh(X)$ as follows: we first use Theorem \ref{co6thm2} to evaluate invariants for an increasing slope stability condition $(\mu,\R,\le)$ on $\modCQI$, and then relate invariants in $\coh(X)$ and $\modCQI$ via wall-crossing in the triangulated category $D^b\coh(X)\simeq D^b\modCQI$, using the ideas of~\S\ref{co55}.
\end{meth}

The author expects Method \ref{co6meth1} to work in the following three examples:

\begin{ex}
\label{co6ex1}
By Beilinson \cite{Beil}, $\O(0),\O(1)$ is a full strong exceptional collection in $D^b\coh(\CP^1)$. Then Bondal \cite[Ex.~5.1 \& Ex.~6.3]{Bond} implies that $D^b\coh(\CP^1)\ab\simeq D^b\modCQ$, where $Q$ is the quiver without relations
\begin{equation*}
\xymatrix@C=50pt{
\overset{v_1}{\bu} \ar@/^.6pc/[r]^{e_1} \ar@/_.6pc/[r]^{e_2} & \overset{v_2}{\bu}. }
\end{equation*}
Writing $\baM$ for the moduli stack of objects in $D^b\coh(\CP^1)\ab\simeq D^b\modCQ$, the Lie algebra $\check H_0(\baM)$ should be the affine Lie algebra~$\ti A_1$.
\end{ex}

\begin{ex}
\label{co6ex2}
By Beilinson \cite{Beil}, $\O(0),\O(1),\O(2)$ is a full strong exceptional collection in $D^b\coh(\CP^2)$. Then Bondal \cite[Ex.~5.3 \& Ex.~6.4]{Bond} implies that $D^b\coh(\CP^2)\simeq D^b\modCQI$, where $\br Q$ is the quiver with relations
\begin{equation*}
\xymatrix@C=50pt{
\overset{v_1}{\bu} \ar[r]^{e_2} \ar@/^1.2pc/[r]^{e_1} \ar@/_1.2pc/[r]^{e_3} & \overset{v_2}{\bu} \ar[r]^{f_2} \ar@/^1.2pc/[r]^{f_1} \ar@/_1.2pc/[r]^{f_3} & \overset{v_3}{\bu}, }\quad\text{relations}\quad \begin{aligned}[h] 
r(q_1)&=f_2e_3-f_3e_2, \\ r(q_2)&=f_3e_1-f_1e_3, \\ r(q_3)&=f_1e_2-f_2e_1, \\
b(q_i)&=v_1,\; e(q_i)=v_3.
\end{aligned}
\end{equation*}
\end{ex}

\begin{ex}
\label{co6ex3}
As in Bondal \cite[Ex.~5.4 \& Ex.~6.5]{Bond}, $\O(0,0),\O(1,0),\O(0,1),\ab\O(1,1)$ is a full strong exceptional collection in $D^b\coh(\CP^1\t\CP^1)$. This yields $D^b\coh(\CP^1\ab\t\CP^1)\simeq D^b\modCQI$, where $\br Q$ is the quiver with relations
\begin{equation*}
\begin{gathered}
\xymatrix@C=50pt@R=11pt{
& \overset{v_2}{\bu} \ar@/_.2pc/[dr]_{g_2} \ar@/^.4pc/[dr]^{g_1} \\
\overset{v_1}{\bu} \ar@/_.2pc/[ur]_{e_2} \ar@/^.4pc/[ur]^{e_1} \ar@/_.4pc/[dr]_{f_2} \ar@/^.2pc/[dr]^{f_1} && \overset{v_4}{\bu}, 
\\
& \overset{v_3}{\bu} \ar@/_.4pc/[ur]_{h_2} \ar@/^.2pc/[ur]^{h_1} }
\end{gathered}
\quad\text{relations}\quad \begin{aligned}[h] 
r(q_1)&=g_1e_1-h_1f_1, \\ r(q_2)&=g_1e_2-h_2f_1, \\ r(q_3)&=g_2e_1-h_1f_2, \\ r(q_4)&=g_2e_2-h_2f_2, \\
b(q_i)&=v_1,\; e(q_i)=v_4.
\end{aligned}
\end{equation*}
\end{ex}

\subsubsection{$\bG_m^2$-equivariant coherent sheaves on $\C^2$}
\label{co653}

In Remark \ref{co5rem3}(b) we outlined a method of studying enumerative invariants in $G$-equivariant problems when the fixed point moduli spaces $\M_\al^\ss(\tau)^G$ are proper, but the $\M_\al^\ss(\tau)$ are not. We explain this in a simple example.

\begin{ex}
\label{co6ex4}
There is an equivalence $F:\coh_\cs(\C^2)\,{\buildrel\sim\over\longra}\,\modCQI$ from the abelian category of compactly-supported coherent sheaves on $\C^2$ to the category of representations of the quiver with relations $\br Q$ given by
\e
\xymatrix{ \overset{v}{\bu} \ar@(ul,dl)[]_e \ar@(ur,dr)[]^f }
\qquad\text{with one relation}\qquad r(q)=ef-fe,
\label{co6eq16}
\e
acting on objects by $F:E\mapsto(V_v,(\rho_e,\rho_f))$ where $V_v=H^0(E)$ and $\rho_e=x\cdot-$, $\rho_f=y\cdot-$ for $(x,y)$ the coordinates on $\C^2$, where the relation means~$xy=yx$.

As $K(\modCQI)=\Z$, every weak stability condition is equivalent to the trivial stability condition $\mu\equiv 0$, in which every nonzero object in $\modCQI$ is semistable. However, in Assumption \ref{co5ass2}(g),(h), the moduli spaces are {\it not proper}, as $\C^2$ is not projective, and $\br Q$ has oriented cycles. For example, we have $\M_1^\rst(0)=\M_1^\ss(0)=\C^2$, the moduli space of point sheaves $\O_{(x,y)}$ on~$\C^2$.

Now let $T=\bG_m^2$ act on $\C^2$ by $(\ep,\ze):(x,y)\mapsto(\ep x,\ze y)$. This induces a $T$-action on $\A=\B=\coh_\cs(\C^2)$. The corresponding action on $\modCQI$ acts on objects by $(\ep,\ze):(V_v,(\rho_e,\rho_f))\mapsto (V_v,(\ep\rho_e,\ze\rho_f))$. It also acts on the relation by $(\ep,\ze):r(q)\mapsto \ep\ze r(q)$, which is needed to define the $T$-equivariant structure on $\cE_{\bs c,\bs d}^\bu,\bs\brM_{\bs d},\bs\brM_{\bs d}^\pl$ in Definitions \ref{co6def6}--\ref{co6def7}.

Remark \ref{co5rem3}(b) says we should consider the abelian category $\A^T$ of $T$-equivariant objects in $\A$. An object $(V,\rho)$ of $\A^T$ consists of a finite-dimensional $T$-representation $V$, and morphisms $\rho_e,\rho_f:V\ra V$ which are weighted $T$-equivariant with $T$-action $(\ep,\ze):(\rho_e,\rho_f)\mapsto (\ep\rho_e,\ze\rho_f)$, and relation $\rho_e\ci\rho_f=\rho_f\ci\rho_e$. Equivalently we may split $V=\bigop_{(i,j)\in\Z^2}V_{i,j}$, where $V_{i,j}$ has $T$-weight $(i,j)$, and $V_{i,j}=0$ for all but finitely many $i,j\in\Z^2$, and then $\rho_e,\rho_f$ should map $\rho_e:V_{i,j}\ra V_{i-1,j}$, $\rho_f:V_{i,j}\ra V_{i,j-1}$. Thus we may write $\A^T\cong\modCtQI$, where $\ti Q$ is the infinite quiver with relations
\e
\begin{gathered}
\xymatrix@C=17pt@R=11pt{ 
& \vdots \ar[d]_(0.4){f_{-1,2}} & \vdots \ar[d]_(0.4){f_{0,2}} & \vdots \ar[d]_(0.4){f_{1,2}}  \\
\cdots & \overset{v_{-1,1}}{\bu} \ar[l]_{e_{-1,1}} \ar[d]_(0.4){f_{-1,1}} & \overset{v_{0,1}}{\bu} \ar[l]_{e_{0,1}} \ar[d]_(0.4){f_{0,1}} & \overset{v_{1,1}}{\bu} \ar[l]_{e_{1,1}} \ar[d]_(0.4){f_{1,1}} &  \ar[l]_{e_{2,1}} \cdots \\
\cdots & \overset{v_{-1,0}}{\bu} \ar[l]_{e_{-1,0}} \ar[d]_(0.4){f_{-1,0}} & \overset{v_{0,0}}{\bu} \ar[l]_{e_{0,0}} \ar[d]_(0.4){f_{0,0}} & \overset{v_{1,0}}{\bu} \ar[l]_{e_{1,0}} \ar[d]_(0.4){f_{1,0}} &  \ar[l]_{e_{2,0}} \cdots \\
\cdots & \overset{v_{-1,-1}}{\bu} \ar[l]_(0.55){e_{-1,-1}} \ar[d]_(0.4){f_{-1,-1}} & \overset{v_{0,-1}}{\bu} \ar[l]_{e_{0,-1}} \ar[d]_(0.4){f_{0,-1}} & \overset{v_{1,-1}}{\bu} \ar[l]_{e_{1,-1}} \ar[d]_(0.4){f_{1,-1}} &  \ar[l]_{e_{2,-1}} \cdots \\
& \vdots  & \vdots & \vdots   }
\end{gathered}
\;\>\text{relations}\;\> \begin{aligned}[h] 
r(q_{i,j})&=f_{i-1,j}e_{i,j}\\
&-e_{i,j-1}f_{i,j}, \\
b(q_{i,j})&=v_{i,j},\\ 
e(q_{i,j})&=v_{i-1,j-1},\\
i,j&\in\Z,
\end{aligned}\!\!\!
\label{co6eq17}
\e
where infinite quivers are as in Remark~\ref{co6rem2}.

Here $\ti Q$ has no oriented cycles, in contrast to $\br Q$, so Theorem \ref{co6thm1} (which also works for infinite quivers) says that Theorems \ref{co5thm1}--\ref{co5thm3} apply with $\A=\B=\modCtQI$. We could in principle use Theorem \ref{co6thm2} to compute enumerative invariants $[\tiM_{\bs d}^\ss(0)]_\inv$ for the trivial stability condition on $\modCtQI$, which will be related to the combinatorics of counting partitions.

As in Remark \ref{co5rem3}(b), we can use enumerative invariants $[\tiM_{\bs d}^\ss(0)]_\inv$ in $\check H_0(\tiM^\pl)$ for $\B^T=\modCtQI$ to {\it define\/} enumerative invariants in localized $T$-equivariant homology $\check H_0^T(\brM^\pl)_{\rm loc}$ for $\B=\modCQI$ using \eq{co2eq25}. Here we note a subtlety in doing this. We have two moduli spaces $\brM,\brM^\pl$ for $\B=\modCQI$, and $T$ acts on both. So we have $T$-fixed substacks $\brM^T,(\brM^\pl)^T$ with inclusion morphisms $i:\brM^T\ra\brM$, $i^\pl:(\brM^\pl)^T\ra\brM^\pl$. There is a natural equivalence $\brM^T\cong\tiM$. 

However, $(\brM^\pl)^T\not\cong\tiM^\pl$. Instead, $\Z^2=\Hom(T,\bG_m)$ acts on $\tiM^\pl$, freely except at 0, and $(\brM^\pl)^T\cong\tiM^\pl/\Z^2$. This is because if $0\ne (V,\rho)\in\B^T$ and $L$ is a nontrivial 1-dimensional $T$-representation, representing a point in $\Hom(T,\bG_m)$, then $[V,\rho]$ and $[(V,\rho)\ot L]$ are distinct points in $\tiM^\pl$, but the same point in $(\brM^\pl)^T$, since taking $T$-fixed substacks and quotienting by $[*/\bG_m]$ do not commute.
\end{ex}

\begin{rem}
\label{co6rem4}
Let $X$ be an Artin $\C$-stack acted on by an algebraic $\C$-group $T$. By `$T$-fixed substack' above, we mean $\Map_T(*,X)$, the moduli stack of $T$-equivariant 1-morphisms $*=\Spec\C\ra X$, where $T$ acts trivially on~$*$.

If $X=[V/G]$ is a quotient stack for $V$ a $\C$-scheme acted on by an algebraic $\C$-group $G$, and $T$ acts on $X$ via a $T$-action on $V$ commuting with the $G$-action (that is, $T\t G$ acts on $V$), then
\begin{equation*}
X^T\cong\coprod_{\begin{subarray}{l} \text{morphisms $\rho:T\ra G$ modulo}\\ \text{conjugation $\rho\mapsto \Ad(\ga)\ci\rho$ for $\ga\in G$}\end{subarray}} [V^{(\id,\rho)(T)}/C_\rho(G)].
\end{equation*}
Here $V^{(\id,\rho)(T)}$ is the fixed $\C$-subscheme for the subgroup $(\id,\rho)(T)\subset T\t G$ acting on $V$, and $C_\rho(G)=\bigl\{\ga\in G:\ga\rho(t)=\rho(t)\ga$ for all~$t\in T\bigr\}$.

Our intuition from group actions on schemes suggests that $X^T$ should be smaller than $X$, maybe a closed substack. However, because of the 2-categorical nature of Artin stacks, this is misleading, and $X^T$ may be {\it much larger than\/} $X$. For example, if $T=\bG_m$ and $X=[*/\bG_m]$ then~$X^T\cong\Z\t[*/\bG_m]=\Z\t X$. 

In Example \ref{co6ex4}, taking the subcategory of $T$-equivariant objects turns the one vertex quiver \eq{co6eq16} with $K(\modCQI)=\Z$ into the infinite quiver \eq{co6eq17} with $K(\modCtQI)=\Z^{\Z^2}$. Note that $\Z^{\Z^2}$ is huge: many different invariants $[\tiM_{\bs d}^\ss(0)]_\inv$ in $\modCtQI$ contribute to one $[\brM_d^\ss(0)^T]_\inv$ in~$\modCQI$.
\end{rem}

\section{Applications to coherent sheaves}
\label{co7}

Next we apply the results of Chapter \ref{co5} to abelian categories $\coh(X)$ of coherent sheaves on a smooth projective $\C$-scheme $X$. Although part of the story (including the vertex algebras in Chapter \ref{co4}) works for $X$ of any dimension, to get the quasi-smooth derived stacks $\bs\dM_\al^\red,\bs\dM_\al^\rpl$ and obstruction theories \eq{co5eq5} in Assumption \ref{co5ass1}(f) we must take $X$ to be a curve, a surface, or a Fano 3-fold. 

As in \S\ref{co1.5}, the author hopes in future to extend the theory to $\coh(X)$ for $X$ a Calabi--Yau 3- or 4-fold, so our verification of much of the Assumptions in \S\ref{co72} and \S\ref{co74}--\S\ref{co75} is written to include Calabi--Yau 3- and 4-folds.

We begin in \S\ref{co71} with background on coherent sheaves, including some new material on Gieseker and $\mu$-stability for real (non-rational) K\"ahler classes $\om$ in $\Kah(X)$. Sections \ref{co72}--\ref{co75} define the data in and verify Assumptions \ref{co4ass1}, \ref{co4ass2}, \ref{co5ass1}--\ref{co5ass3} and \ref{co5ass4}. The main results are stated in~\S\ref{co76}--\S\ref{co78}.

\subsection{\texorpdfstring{Coherent sheaves, Gieseker stability, and $\mu$-stability}{Coherent sheaves, Gieseker stability, and μ-stability}}
\label{co71}

We recall some background on the abelian category $\coh(X)$ of coherent sheaves on a smooth projective $\C$-scheme $X$, and Gieseker and $\mu$-stability on $\coh(X)$. Some good references on coherent sheaves are Hartshorne \cite[\S II.5]{Hart} and Huybrechts and Lehn \cite{HuLe1}. For derived categories of coherent sheaves $D^b\coh(X)$ we recommend Huybrechts~\cite{Huyb}.

In \S\ref{co712} we will define Gieseker stability and $\mu$-stability with respect to a real K\"ahler class $\om$ on $X$. Note that working with {\it real\/} K\"ahler classes is not standard --- most of the literature considers only Gieseker stability and $\mu$-stability defined using {\it integral\/} or {\it rational\/} K\"ahler classes. We need real K\"ahler classes to ensure that two weak stability conditions in $\sS$ can be connected by a continuous family $(\tau_t,T_t,\le)_{t\in[0,1]}$ in $\sS$, as in Assumption~\ref{co5ass3}(a).

In some possibly new material, we extend well known theorems for Gieseker and $\mu$-stability for integral K\"ahler classes to real K\"ahler classes, and we also apply the results of \S\ref{co33} to Gieseker and $\mu$-stability on $\coh(X)$. Longer proofs will be deferred to~\S\ref{co7112}--\S\ref{co7116}.

\subsubsection{Basics on coherent sheaves and their moduli stacks}
\label{co711}

We will use the following notation throughout Chapter~\ref{co7}:

\begin{dfn}
\label{co7def1}
Let $X$ be a smooth, connected, projective $\C$-scheme with $\dim_\C X=m$. In \S\ref{co73} we restrict to $m=1,2$ or 3, but \S\ref{co71}--\S\ref{co72} works for all $m\ge 0$. Write $\coh(X)$ for the abelian category of coherent sheaves on $X$, as in Hartshorne \cite[\S II.5]{Hart} and Huybrechts and Lehn \cite{HuLe1}, and $D^b\coh(X)$ for its derived category, as in Huybrechts~\cite{Huyb}.

Objects $E\in\coh(X)$ are sheaves of modules over the structure sheaf $\O_X$ of $X$ satisfying a coherence condition \cite[p.~111]{Hart}. We call $E$ {\it locally free}, or a {\it vector bundle}, if it is Zariski locally isomorphic to $\O_X^{\op^r}$ for some $r\ge 0$, the {\it rank\/} of $E$. More generally, as $X$ is connected, any $E\in\coh(X)$ has a {\it rank\/} $\rank E\in\N$ such that $E\vert_U\cong\O_U^{\op^{\rank E}}$ for some dense Zariski open subset $U\subseteq X$.

The {\it support\/} $\supp E$ is the closed subset of points $x\in X$ such that the stalk $E_x$ of $E$ at $x$ is nonzero. The {\it dimension\/} of $E$ is $\dim E=\dim\supp E$, so that $\dim E\in\{0,\ldots,m\}$. A coherent sheaf $E$ is called {\it torsion\/} if $\dim E<m$, and {\it torsion-free\/} if it has no nonzero subobject $0\ne E'\subseteq E$ such that $E'$ is torsion. Also $E$ is called {\it pure\/} if it has no nonzero subobject $0\ne E'\subseteq E$ with $\dim E'<\dim E$. Then torsion-free sheaves are pure sheaves of dimension~$m$.

It may help to think of a general coherent sheaf $E$ on $X$ as being like a (singular) vector bundle $E\ra Y$ over a (singular) closed submanifold~$Y\subseteq X$.

As in \cite[\S III.6]{Hart}, for $E,F$ in $\coh(X)$ we have {\it Ext groups\/} $\Ext^i(E,F)$ for $i=0,\ldots,m$, which are finite-dimensional $\C$-vector spaces, with $\Ext^0(E,F)=\Hom(E,F)$. Writing $K_X\ra X$ for the canonical bundle, as in \cite[\S III.7]{Hart}, {\it Serre duality\/} gives canonical, functorial isomorphisms 
\begin{equation*}
\Ext^i(E,F)\cong\Ext^{m-i}(F,E\ot K_X)^*.
\end{equation*}
The {\it Euler form\/} is the biadditive map $\chi:K_0(\coh(X))\t K_0(\coh(X))\ra\Z$ with
\e
\chi([E],[F])=\sum_{i=0}^m(-1)^i\dim_\C\Ext^i(E,F)
\label{co7eq1}
\e
for all $E,F\in\coh(X)$. The {\it numerical Grothendieck group\/} is 
\begin{align*}
&K^\num(\coh(X))=K_0(\coh(X))/\Ker\chi, \quad\text{where}\\
&\Ker\chi=\bigl\{\al\in K_0(\coh(X)):\text{$\chi(\al,\be)=0$ for all $\be\in K_0(\coh(X))$}\bigr\},
\end{align*}
so we have a surjective quotient $K_0(\coh(X))\twoheadrightarrow K^\num(\coh(X))$. Then $\chi$ descends to $\chi:K^\num(\coh(X))\t K^\num(\coh(X))\ra\Z$ by definition. As $X$ is connected, the map $E\mapsto\rank E$ induces a morphism $\rank:K_0(\coh(X))\ra\Z$ which descends to $\rank:K^\num(\coh(X))\ra\Z$, since $\rank\al=\chi(\al,\lb\O_x\rb)$ for $\O_x$ the skyscraper sheaf at any $\C$-point $x\in X$. We will always take $K(\coh(X))=K^\num(\coh(X))$ in Definition \ref{co3def1} and Assumptions \ref{co4ass1}(d) and \ref{co5ass1}(b). We write $\lb E\rb\in K(\coh(X))$ for the class of~$E\in\coh(X)$.

The dimension $\dim E$ also descends to $C(\coh(X))\subset K(\coh(X))$, by for $\al\in C(\coh(X))$, regarding $K(\coh(X))$ as a subgroup of $H^{\rm even}(X,\Q)$ as above 
\e
\dim\al\!=\!\max\bigl\{\text{$0\!\le\! d\!\le\! m$: $\al$ has nonzero component in $H^{2m-2d}(X,\Q)$}\bigr\},
\label{co7eq2}
\e
so that if $0\ne E\in\coh(X)$ with $\lb E\rb=\al$ then $\dim\al=\dim E=\dim\supp E$.

As $K_0(\coh(X))=K_0(D^b\coh(X))$, we can also write $\chi([E^\bu],[F^\bu])$ for $E^\bu,F^\bu$ in $D^b\coh(X)$. The Grothendieck--Riemann--Roch Theorem \cite[\S A.4]{Hart} says that
\e
\chi(\lb E^\bu\rb,\lb F^\bu\rb)=\int_X\ch(E^\bu)^\vee\cup\ch(F^\bu)\cup\td(X),
\label{co7eq3}
\e
where $(E^\bu)^\vee$ is the derived dual in $D^b\coh(X)$, not the classical dual in $\coh(X)$, and $\ch$ is the {\it Chern character\/} and $\td(X)$ the {\it Todd class}, as in~\cite[App.~A]{Hart}. 

Since \eq{co7eq3} is a nondegenerate inner product on $\ch(E^\bu),\ch(F^\bu)$, it follows that we may identify $K^\num(\coh(X))\cong\Im\bigl(\ch:K_0(\coh(X))\ra H^{\rm even}(X,\Q)\bigr)$, so that $K_0(\coh(X))\twoheadrightarrow K^\num(\coh(X))$ is the Chern character. Then identifying $K^\num(\coh(X))$ with a subgroup of $H^{\rm even}(X,\Q)$, equation \eq{co7eq3} becomes
\begin{equation*}
\chi(\al,\be)=\int_X\al^\vee\cup\be\cup\td(X).
\end{equation*}
Here if $\al=(\al_i)_{i=0}^m$ with $\al\in H^{2i}(X,\Q)$ then we define~$\al^\vee=((-1)^i\al_i)_{i=0}^m$.\end{dfn}

\begin{rem}
\label{co7rem1}
The assumption that $X$ is {\it connected\/} is not essential but included for convenience, as it simplifies notation at several points, for example, in that $E\in\coh(X)$ have a $\rank E\in\Z$. We can easily extend our results to non-connected $X$ by applying them to each connected component of $X$.	
\end{rem}

For the rest of \S\ref{co71} we work in the situation of Definition \ref{co7def1}. We discuss moduli stacks of objects in $\coh(X)$ and $D^b\coh(X)$:

\begin{dfn}
\label{co7def2}
Write $\M$ for the moduli stack of objects in $\coh(X)$. By Laumon and Moret-Bailly \cite[\S\S 2.4.4, 3.4.4 \& 4.6.2]{LaMo} using results of Grothendieck \cite{Grot1,Grot2}, this is an Artin $\C$-stack, locally of finite type, and $\C$-points of $\M$ are isomorphism classes $[E]$ of objects $E\in\coh(X)$. It has a decomposition $\M=\coprod_{\al\in K(\coh(X))}\M_\al$, where $\M_\al$ is the moduli space of $E\in\coh(X)$ with $\lb E\rb=\al$ in~$K(\coh(X))$.

Write $\baM$ for the moduli stack of objects in $D^b\coh(X)$. This is a higher $\C$-stack in the sense of \S\ref{co22}, and exists by To\"en--Vaqui\'e \cite{ToVa}. $\C$-points of $\baM$ are isomorphism classes $[E^\bu]$ of $E^\bu\in D^b\coh(X)$. The inclusion $\coh(X)\hookra D^b\coh(X)$ induces an inclusion $\M\hookra\baM$ as an open $\C$-substack. Again we have a decomposition $\baM=\coprod_{\al\in K(\coh(X))}\baM_\al$, where $\baM_\al$ is the moduli space of $E^\bu\in D^b\coh(X)$ with $\lb E^\bu\rb=\al$ in~$K(D^b\coh(X))=K(\coh(X))$.

There is a {\it universal coherent sheaf\/} $\cU\ra X\t\M$, flat over $X$, and a {\it universal perfect complex\/} $\cU^\bu\ra X\t\baM$, with $\cU^\bu\vert_{X\t\M}\cong\cU$, such that $\cU\vert_{X\t\{[E]\}}\cong E$ and $\cU^\bu\vert_{X\t\{[E^\bu]\}}\cong E^\bu$. Write $\Perf_\C$ for the higher stack which classifies perfect complexes, as in To\"en and Vezzosi \cite[Def.~1.3.7.5]{ToVe2}, which is just $\baM$ for $X=\Spec\C$ the point. Then $\cU^\bu\ra X\t\baM$ determines a morphism $u:X\t\baM\ra\Perf_\C$ in $\HSta_\C$, such that $\cU^\bu\cong u^*(\cU_{\Perf_\C}^\bu)$ with $\cU_{\Perf_\C}^\bu\ra\Perf_\C$ the tautological perfect complex. This $u$ in turn determines a morphism $\baM\ra\Map_{\HSta_\C}(X,\Perf_\C)$ to the mapping stack in the $\iy$-category $\HSta_\C$. As in (the classical truncation of) Pantev--To\"en--Vaqui\'e--Vezzosi \cite[\S 2.1 \& Cor.~2.13]{PTVV}, this is an isomorphism $\baM\cong\Map_{\HSta_\C}(X,\Perf_\C)$.
 
The {\it Ext complex\/} $\cExt^\bu$ is a perfect complex on $\M\t\M$, given by
\e
\cExt^\bu=(\Pi_{23})_*\bigl[\Pi_{12}^*(\cU^\vee)\ot\Pi_{13}^*(\cU)\bigr],
\label{co7eq4}
\e
using derived duality, pushforward, pullback and tensor product functors as in Huybrechts \cite{Huyb}, where $\Pi_{ij}$ projects to the product of the $i^{\rm th}$ and $j^{\rm th}$ factors of $X\t\M\t\M$. It has $H^i\bigl(\cExt^\bu\vert_{([E],[F])}\bigr)\cong\Ext^i(E,F)$ for $E,F\in\coh(X)$ and $i\in\Z$, and $\rank\bigl(\cExt^\bu\vert_{\M_\al\t\M_\be}\bigr)=\chi(\al,\be)$ for~$\al,\be\in K(\coh(X))$.

As in Definition \ref{co4def4} we may also form the `projective linear' moduli stack $\M^\pl=\M/[*/\bG_m]$ of objects in $\coh(X)$ moduli projective linear isomorphisms. There is a projection $\Pi^\pl:\M\ra\M^\pl$ which is a principal $[*/\bG_m]$-bundle except over $[0]\in\M$. It has a decomposition $\M^\pl=\coprod_{\al\in K(\coh(X))}\M_\al^\pl$, such that $\Pi^\pl$ restricts to $\Pi^\pl_\al:\M_\al\ra\M^\pl_\al$. 

We will also write $\bs\M,\bs\baM$ for the derived moduli stacks of objects in $\coh(X),\ab D^b\coh(X)$, as derived $\C$-stacks in the sense of \S\ref{co22}, which exist by To\"en--Vaqui\'e \cite{ToVa}. Then $\bs\M\subset\bs\baM$ is an open substack. The classical truncations are $t_0(\bs\M)=\M$ and $t_0(\bs\baM)=\baM$, with inclusions $i:\M\hookra\bs\M$, $i:\baM\hookra\bs\baM$. The previous paragraph also holds in the world of derived rather than higher stacks, so we have universal $\bs\cU\ra X\t\bs\M$, $\bs\cU^\bu\ra X\t\bs\baM$ determining a morphism $\bs u:X\t\bs\baM\ra\bs\Perf_\C$ in $\DSta_\C$, and $\bs u$ induces an isomorphism~$\bs\baM\ra\Map_{\DSta_\C}(X,\bs\Perf_\C)$.

Since $X$ is a smooth projective $\C$-scheme, $\bs\M,\bs\baM$ are locally finitely presented. This implies that the cotangent complexes $\bL_{\bs\M},\bL_{\bs\baM}$ are perfect complexes, and we will show later that $\bL_{\bs\M}$ is perfect in~$[1-m,1]$.
\end{dfn}

\begin{rem}
\label{co7rem2}
As above, we write $\M\subset\baM$ for the moduli stacks of objects in $\coh(X)\subset D^b\coh(X)$. As in Gross \cite[Th.~1.1]{Gros} and \S\ref{co55}(b), it turns out that the vertex algebra $\hat H_*(\baM)$ can often be written down very explicitly, but $\hat H_*(\M)$ is usually much more difficult to compute. So although we will apply the results of Chapters \ref{co4}--\ref{co5} to $\A=\coh(X)$ and the homology of moduli stacks $\M,\M^\pl$, for explicit calculations we usually prefer to push forward to the homology of~$\baM,\baM^\pl$.
\end{rem}

\begin{prop}
\label{co7prop1}
In Definitions\/ {\rm\ref{co7def1}} and\/ {\rm\ref{co7def2},} $\coh(X)$ is \begin{bfseries}of compact type\end{bfseries} in the sense of\/ {\rm\S\ref{co332}}. Hence $\M$ and\/ $\M_\al$ for $\al\in K(\coh(X))$ are locally of finite type, with affine diagonal, and satisfy Theorem\/~{\rm\ref{co3thm4}(i)--(iii)}.	
\end{prop}

\begin{proof}
Let $\tiA=\qcoh(X)$ be the abelian category of quasicoherent sheaves on $X$, as in Hartshorne \cite[\S II.5]{Hart}. By \cite[\S B.2--\S B.3]{ThTr} this a cocomplete, locally noetherian $\C$-linear abelian category, and the coherent sheaves $\coh(X)\subset\qcoh(X)$ are the full subcategory of compact objects in $\qcoh(X)$ (see \cite[Tag 09M8]{deJo}). As noted in \cite[Ex.~7.23]{AHLH}, the moduli stack $\M$ in Definition \ref{co3def9} is the usual moduli stack of objects in $\coh(X)$, which is an Artin $\C$-stack locally of finite type as in \cite[Th.~4.6.2.1]{LaMo}. Hence $\coh(X)$ is of compact type, as in Definition \ref{co3def9}, and the proposition follows from Theorem~\ref{co3thm5}.	
\end{proof}

The next proposition is a consequence of the existence of Quot schemes, as in Huybrechts and Lehn \cite[\S 2.2]{HuLe2}, and generalizes an argument in \cite[Prop.~2.3.1]{HuLe2} showing Gieseker semistability is an open condition on coherent sheaves.

\begin{prop}
\label{co7prop2}
In Definitions\/ {\rm\ref{co7def1}} and\/ {\rm\ref{co7def2},} let\/ $\al_1,\ldots,\al_n$ lie in $C(\coh(X))$ for $n\ge 2,$ and set\/ $\al=\al_1+\cdots+\al_n$. Then the set of\/ $\C$-points $[E]\in \M_\al$ such that there does not exist a filtration\/ $0=E_0\subsetneq E_1\subsetneq\cdots\subsetneq E_n=E$ in $\coh(X)$ with\/ $\lb E_i/E_{i-1}\rb=\al_i$ for $i=1,\ldots,n$ is open in $\M_\al$. The analogue holds for~$\M_\al^\pl$.
\end{prop}

\begin{proof} Consider first the case $n=2$. Let $\Up:U\ra\M_\al$ be an atlas for $\M_\al$, with $U$ a $\C$-scheme locally of finite type. For $\cU\ra X\t\M$ as in Definition \ref{co7def2}, write $\cE=(\id_X\t\Up)^*(\cU)$, so that $\cE\ra X\t U$ is a coherent sheaf on $X\t U$, flat over $U$, with $\cE_u=\cE\vert_{X\t\{u\}}\cong E\vert_{\pi(u)}$ for $\C$-points $u\in U$. As in Huybrechts and Lehn \cite[\S 2.3]{HuLe2}, we may form {\it Grothendieck's Quot scheme\/} $\mathop{\rm Quot}_{X\t U/U}(\cE,\al_2)$, a $\C$-scheme with a projective morphism $\pi:\mathop{\rm Quot}_{X\t U/U}(\cE,\al_2)\ra U$, which is a fine moduli scheme for quotients $\cE_u\twoheadrightarrow\cF_u$ with $\lb\cF_u\rb=\al_2$ for $u\in U$.

Thus if $u\in U$ then $\cE_u$ admits a filtration $0=E_0\subsetneq E_1\subsetneq E_2=\cE_u$ with $\lb E_2/E_1\rb=\al_2$, so that $\lb E_1/E_0\rb=\al_1$, if and only if $u$ lies in the image of $\pi:\mathop{\rm Quot}_{X\t U/U}(\cE,\al_2)\ra U$. As $\pi$ is projective, this image is closed, so the condition that no such filtration exists is an open condition on $u\in U$. As $\Up:U\ra\M_\al$ is an atlas for $\M_\al$, the condition that no such filtration exists is also an open condition on points $[E]\in\M_\al$. The analogue holds for $\M_\al^\pl$, as $\Pi_\al^\pl:\M_\al\ra\M_\al^\pl$ identifies open substacks in $\M_\al$ and $\M_\al^\pl$. For $n>2$ we can give a similar argument involving the fibre product of $n-1$ Quot schemes.	
\end{proof}

\subsubsection{Gieseker stability, $\mu$-stability, and purity}
\label{co712}

The next three definitions define (weak) stability conditions on $\coh(X)$: Gieseker stability and $\mu$-stability with respect to a real K\"ahler class $\om$ on $X$, and purity.

Most of the literature, as in Huybrechts and Lehn \cite{HuLe2} for instance, considers only Gieseker and $\mu$-stability defined using {\it integral\/} or {\it rational\/} K\"ahler classes. Greb--Ross--Toma \cite[\S 11]{GRT1}, \cite[\S 5]{GRT3}, \cite{GRT2,GrTo1,GrTo2} discuss Gieseker and $\mu$-stability for real K\"ahler or real ample classes.

\begin{dfn}
\label{co7def3}
Since $X$ is a compact complex manifold admitting K\"ahler metrics, as in \cite[\S 4]{Joyc1}, we can use Hodge theory. The second complex de Rham cohomology group decomposes as
\begin{equation*}
H^2(X,\C)=H^{2,0}(X)\op H^{1,1}(X)\op H^{0,2}(X).
\end{equation*}
Write $H^{1,1}(X,\R)=H^{1,1}(X)\cap H^2(X,\R)$ for the vector subspace of real $(1,1)$-classes in $H^2(X,\R)$. Write $\Kah(X)\subset H^{1,1}(X,\R)$ for the {\it K\"ahler cone\/} of $X$, the set of {\it K\"ahler classes\/} on $X$, that is, cohomology classes of K\"ahler forms of K\"ahler metrics on $X$. Then $\Kah(X)$ is an open convex cone in~$H^{1,1}(X,\R)$.

Write $\Kah_\Z(X)\subset\Kah_\Q(X)\subset\Kah(X)$ for the subsets of $\om\in\Kah(X)$ which lie in the images of the natural morphisms $H^2(X,\Z)\ra H^2(X,\R)$ and $H^2(X,\Q)\ra H^2(X,\R)$. We call these {\it integral\/} and {\it rational\/} K\"ahler classes respectively. Here $\om\in\Kah_\Z(X)$ if and only if $\om=c_1(L)$ for some ample line bundle $L\ra X$. Note that if $H^{2,0}(X)\ne 0$ then $H^{1,1}(X,\R)$ is a proper subspace of $H^2(X,\R)$, and then $\Kah_\Q(X)$ need not be dense in~$\Kah(X)$.

We write $\ov{\Kah_\Q}(X)$ for the closure of $\Kah_\Q(X)$ in $\Kah(X)$. If $H^{2,0}(X)=0$ then~$\ov{\Kah_\Q}(X)=\Kah(X)$.

Let $\om\!\in\!\Kah(X)$ and $\al\!\in\! K(\coh(X))$. Noting that $K(\coh(X))\!\subset\! H^{\rm even}(X,\Q)\ab\subset H^*(X,\R)$, the {\it Hilbert polynomial of\/ $\al$ with respect to\/} $\om$ is
\e
P_\al^\om(n)=\int_X\al\cup e^{n\om}\cup \td(X).
\label{co7eq5}
\e
Here $e^{n\om}=\sum_{i=0}^mn^i\om^i/i!$, since $\om^i=0$ if $i>m$, so $P_\al^\om(n)$ is a real polynomial of degree at most $m$, and $\td(X)\in H^{\rm even}(X,\Q)$ is the Todd class of~$X$. 

If $\al\in C(\coh(X))$ then $P_\al^\om$ is nonzero with degree $\dim\al$ in~\eq{co7eq2}.

If $\om=c_1(L)$ for $L\ra X$ an ample line bundle, $\al=\lb E\rb$ for $E\in\coh(X)$, and $n\in\Z$ then $P_\al^\om(n)=\chi(E\ot L^n)\in\Z$ by the Grothendieck--Riemann--Roch Theorem, and if $E\ne 0$, so that $\al\in C(\coh(X))$, then $P_\al^\om(n)=\dim H^0(E\ot L^n)>0$ for $n\gg 0$ by Serre vanishing, so the leading coefficient of $P_\al^\om$ is positive. This also holds for general $\om\in\Kah(X)$ and~$\al\in C(\coh(X))$.

Define $G$ to be the set of monic real polynomials in $t$ of degree at most $m$:
\begin{equation*}
G=\bigl\{p(t)=t^d+a_{d-1}t^{d-1}+\cdots+a_0:d=0,\ldots,m,\;
a_0,\ldots,a_{d-1}\in\R\bigr\}.
\end{equation*}
Define a total order `$\le$' on $G$ by $p\le p'$ for $p,p'\in G$ if either
\begin{itemize}
\setlength{\itemsep}{0pt}
\setlength{\parsep}{0pt}
\item[(a)] $\deg p>\deg p'$, or
\item[(b)] $\deg p=\deg p'$ and $p(n)\le p'(n)$ for all $n\gg 0$.
\end{itemize}
We write $p<q$ if $p\le q$ and $p\ne q$. Note that $\deg p>\deg p'$ in (a) implies that $p(n)>p'(n)$ for all $n\gg 0$, the opposite to $p(n)\le p'(n)$ for $n\gg 0$ in~(b).

Fix $\om\in\Kah(X)$, and define $\tau^\om:C(\coh(X))\ra G$ by $\tau^\om(\al)=P^\om_\al/r^\om_\al$, where $r^\om_\al$ is the leading coefficient of $P^\om_\al$, which must be positive as above. As in Rudakov \cite[Lem.~2.5]{Ruda}, $(\tau^\om,G,\le)$ is a stability condition on $\coh(X)$, which we call {\it Gieseker stability\/} for the real K\"ahler class $\om$. We often omit $\om$ from the notation and write~$(\tau,G,\le)=(\tau^\om,G,\le)$.

One can show that $0\ne E\in\coh(X)$ is Gieseker semistable if and only if $E$ is pure, and for all subobjects $0\ne E'\subseteq E$ we have $\tau(\lb E'\rb)(n) \le\tau(\lb E\rb)(n)$ for $n\gg 0$. This is the definition of Gieseker semistability used by Huybrechts and Lehn \cite[Def.~1.2.4]{HuLe2} and many other authors. Here $E$ must be pure as (a) above implies that $E$ cannot have a subobject $0\ne E'\subseteq E$ with~$\dim E'<\dim E$.
\end{dfn}

\begin{dfn}
\label{co7def4}
Continue in the situation of Definition \ref{co7def3}. Define
\begin{equation*}
M=\bigl\{p(t)=t^d+a_{d-1}t^{d-1}:d=0,1,\ldots,m,\;\> a_{d-1}\in\R,\;\>
a_{-1}=0\bigr\}\subset G
\end{equation*}
and restrict the total order $\le$ on $G$ to $M$. Define $\mu^\om:C(\coh(X))\ra M$ by $\mu^\om(\al)=t^d+a_{d-1}t^{d-1}$ when $\tau^\om(\al)=P_\al/r_\al=t^d+a_{d-1}t^{d-1}+\cdots+a_0$, that is, $\mu^\om(\al)$ is the truncation of the polynomial $\tau^\om(\al)$ in Definition \ref{co7def3} at its second term. Then as in \cite[Ex.~4.17]{Joyc6}, $(\mu^\om,M,\le)$ is a weak stability condition on $\coh(X)$. It is called $\mu$-{\it stability\/} for the real K\"ahler class $\om$. It is not a stability condition if $m=\dim X\ge 2$. We often omit $\om$ from the notation and write $(\mu,M,\le)=(\mu^\om,M,\le)$. Our definition of $\mu$-(semi)stability is equivalent to Greb--Ross--Toma \cite[\S 2.1]{GRT2}, but differs from Huybrechts and Lehn~\cite[\S 1.6]{HuLe2}.

If $\tau,\mu$ are defined using the same K\"ahler class $\om$ then $(\mu,M,\le)$ dominates $(\tau,G,\le)$. Thus as in \eq{co3eq2}, if $0\ne E\in\coh(X)$ then
\e
\text{$E$ $\mu$-stable} \Longra \text{$E$ $\tau$-stable} \Longra \text{$E$ $\tau$-semistable} \Longra \text{$E$ $\mu$-semistable.}
\label{co7eq6}
\e
For $\al\in C(\coh(X))$, we also define $\bar\mu(\al)=\bar\mu{}^\om(\al)\in\R$ by $\bar\mu(\al)=\bar\mu{}^\om(\al)=a_{d-1}$ when~$\mu^\om(\al)=t^d+a_{d-1}t^{d-1}$.
\end{dfn}

\begin{dfn}
\label{co7def5}
Continue in the situation of Definitions \ref{co7def3}--\ref{co7def4}. Define
\begin{equation*}
P=\bigl\{p(t)=t^d:d=0,1,\ldots,m\bigr\}\subseteq M\subseteq G,
\end{equation*}
and restrict $\le$ on $G$ to $P$, so that $t^d\le t^e$ if and only if $d\ge e$. Define $\pi:C(\coh(X))\ra P$ by $\pi(\al)=t^{\dim\al}$, for $\dim\al$ as in \eq{co7eq2}. Then as in \cite[Ex.~4.18]{Joyc6}, $(\pi,P,\le)$ is a weak stability condition on $\coh(X)$. A nonzero $E\in\coh(X)$ is $\pi$-semistable if and only if it is {\it pure\/} in the sense of Definition~\ref{co7def1}. 

If $\om\in\Kah(X)$ and $\al\in C(\coh(X))$ with $\tau^\om(\al)=t^d+a_{d-1}t^{d-1}+\cdots+a_0$ and $\mu^\om(\al)=t^d+a_{d-1}t^{d-1}$ then $\pi(\al)=t^d$, that is, $\pi(\al)$ is the truncation of the polynomials $\tau^\om(\al),\mu^\om(\al)$ in Definitions \ref{co7def3}--\ref{co7def4} at their first term. From this we see that $(\pi,P,\le)$ dominates $(\tau^\om,G,\le),(\mu^\om,M,\le)$ for all $\om\in\Kah(X)$.

Since the moduli stacks $\M_\al^\ss(\pi)$ are generally not of finite type, $(\pi,P,\le)$ is not suitable for inclusion in $\sS$ in Assumption \ref{co5ass2}, as Assumption \ref{co5ass2}(b) would not hold. But $(\pi,P,\le)$ will still be useful in some proofs.
\end{dfn}

\begin{lem}
\label{co7lem1}
In Definitions\/ {\rm\ref{co7def3}--\ref{co7def5},} $\coh(X)$ is {\rm$\tau$-, $\mu$-} and\/ $\pi$-artinian in the sense of Definition\/ {\rm\ref{co3def1}}. Hence every $E\in\coh(X)$ has unique {\rm$\tau$-, $\mu$-} and\/ $\pi$-Harder--Narasimhan filtrations by Theorem\/~{\rm\ref{co3thm1}}.
\end{lem}

\begin{proof} By \cite[Lem.~4.19]{Joyc6}, $\coh(X)$ is $\pi$-artinian. As $(\pi,P,\le)$ dominates $(\tau,\ab G,\ab\le)$ and $(\mu,M,\le)$, it follows that $\coh(X)$ is $\tau$- and $\mu$-artinian.
\end{proof}

\begin{rem}
\label{co7rem3}
It is easy to see that $(\mu^\om,M,\le)$ depends continuously on $\om$ in $\Kah(X)$ in the sense of Definition \ref{co3def5}. However, $(\tau^\om,M,\le)$ need {\it not\/} depend continuously on $\om\in\Kah(X)$ when $m=\dim X\ge 2$. This is a well known issue, and appears in Friedman--Qin \cite{FrQi}, Matsuki--Wentworth \cite{MaWe}, and Greb--Ross--Toma \cite{GRT3}, for instance.

For example, suppose $X$ is a surface, $E\in\coh(X)$ has rank 2 and $\lb E\rb=\al$, and $0\ne E'\subsetneq E$ is a subobject of rank 1 with $\lb E'\rb=\be$, and $(\om_s)_{s\in(-\ep,\ep)}$ is a continuous family with Hilbert polynomials
\begin{equation*}
P_\be^{\om_s}(t)=t^2+st-1, \quad P_{\al-\be}^{\om_s}(t)=t^2-st+1,\quad P_\al^{\om_s}(t)=2t^2. 
\end{equation*}
Then $\tau^{\om_s}(\be)<\tau^{\om_s}(\al-\be)$ if $s\le 0$ and $\tau^{\om_s}(\be)>\tau^{\om_s}(\al-\be)$ if $s>0$, so Definition \ref{co3def5} fails. This is why Theorem \ref{co7thm2}(d),(e) below do not extend to Gieseker stability.
\end{rem}

\subsubsection{The Bogomolov inequality}
\label{co713}

Following Huybrechts and Lehn \cite[\S 3.4]{HuLe2}, we define:

\begin{dfn}
\label{co7def6}
For each $E\in\coh(X)$, the {\it discriminant\/} $\De(E)$, in the image of $H^4(X,\Z)$ in $H^4(X,\Q)$ or $H^4(X,\R)$, is 
\e
\De(E)=2\rank E\cdot c_2(E)-(\rank E-1)c_1(E)^2.
\label{co7eq7}
\e
As $\ch_0(E)=\rank E$, $\ch_1(E)=c_1(E)$ and $\ch_2(E)=\ha c_1(E)^2-2c_2(E)$ we have
\e
\De(E)=\ch_1(E)^2-2\ch_0(E)\ch_2(E),
\label{co7eq8}
\e
so $\De(E)$ is homogeneous quadratic in $\ch(E)$. Note that $\De(E)$ depends only on $\al=\lb E\rb\in K(\coh(X))$, and we also write~$\De(\al)=\De(E)$.
\end{dfn}

The next theorem is known as the {\it Bogomolov inequality}. As in Huybrechts and Lehn \cite[Th.~7.3.1]{HuLe2}, it was proved by Bogomolov for $\om\in\Kah_\Z(X)$, and by Bando--Siu \cite[Cor.~3]{BaSi} and Biswas--McKay \cite[Lem.~2.1]{BiMc} for general $\om\in\Kah(X)$. Here \cite{BaSi,BiMc} work with {\it complex analytic\/} coherent sheaves on compact K\"ahler manifolds $(X,\om)$, but when $X$ is projective, algebraic and complex analytic coherent sheaves coincide by Serre's GAGA principle~\cite{Serr}.

\begin{thm}
\label{co7thm1} 
In Definition\/ {\rm\ref{co7def1},} let\/ $m\!=\!\dim X\!\ge\! 2,$ and suppose\/ $(\mu,M,\le)=(\mu^\om,M,\le)$ is as in Definition\/ {\rm\ref{co7def4}} for some $\om\in\Kah(X),$ and\/ $E\in\coh(X)$ is a torsion-free $\mu$-semistable coherent sheaf. Then
\e
\int_X\De(E)\cup \om^{m-2}\ge 0.
\label{co7eq9}
\e	
\end{thm}

Note that if $E$ is $\tau^\om$-semistable then $E$ is $\mu$-semistable by \eq{co7eq6}, so Theorem \ref{co7thm1} also applies to Gieseker semistability defined using~$\om\in\Kah(X)$.

\subsubsection{Boundedness theorems}
\label{co714}

The next definition follows Huybrechts and Lehn~\cite[Def.~1.7.5]{HuLe2}:

\begin{dfn}
\label{co7def7}
A family of isomorphism classes $\sF$ of objects in $\coh(X)$ is called {\it bounded\/} if there exists a finite type $\C$-scheme $S$ and a coherent sheaf $E\ra X\t S$ such that $\sF$ is contained in the set of isomorphism classes $[E\vert_{X\t\{s\}}]$ for all $\C$-points~$s\in S$.
\end{dfn}

Several important theorems prove that families of coherent sheaves on $X$ satisfying certain conditions are bounded. For us, the important point will be that a moduli stack such as $\M_\al^\ss(\tau)$ is of {\it finite type\/} if and only if the family of sheaves it parametrizes is bounded. We will need to prove moduli stacks are of finite type to apply Theorem \ref{co3thm9} and verify Assumption~\ref{co5ass2}(b),(g),(h). 

\begin{dfn}
\label{co7def8}
Let $(\mu,M,\le)=(\mu^\om,M,\le)$ for some $\om\in\Kah(X)$, and let $E\in\coh(X)$ be nonzero and pure of dimension $d$. Then by Lemma \ref{co7lem1}, $E$ has a unique $\mu$-Harder--Narasimhan filtration $0=E_0\subsetneq E_1\subsetneq\cdots\subsetneq E_n=E$ for $n\ge 0$, such that $F_i=E_i/E_{i-1}$ is $\mu$-semistable for $i=1,\ldots,n$. Each $F_i$ is pure of dimension $d$ as $E$ is, so $\mu(\lb F_i\rb)=t^d+\bar\mu(\lb F_i\rb)t^{d-1}$ with $\bar\mu(\lb F_1\rb)>\cdots>\bar\mu(\lb F_n\rb)$. 

Define the {\it maximal slope\/} of $E$ to be $\mu_{\max}(E)=\mu^\om_{\max}(E)=\bar\mu(\lb F_1\rb)$.
\end{dfn}

\begin{prop}
\label{co7prop3}
Let\/ $(\mu,M,\le)=\ab(\mu^\om,M,\le)$ be as in Definition\/ {\rm\ref{co7def4}} for some $\om\in\Kah(X),$ and\/ $\sF$ be a bounded family of pure $d$-dimensional sheaves $F$ in $\coh(X)$ for some $d=0,\ldots,m,$ and\/ $\sG$ be a family of pure $d$-dimensional sheaves $G$ in $\coh(X)$ which are quotients $F\twoheadrightarrow G$ for some $F$ in $\sF$. Then:
\begin{itemize}
\setlength{\itemsep}{0pt}
\setlength{\parsep}{0pt}
\item[{\bf(a)}] The set of Chern characters $\ch_d(G)\in H_{2d}(X,\Q)$ for $G\in\cG$ is finite.
\item[{\bf(b)}] $\bar\mu(G)$ is bounded below for $G$ in $\sG$.
\item[{\bf(c)}] $\sG$ is bounded if and only if\/ $\bar\mu(G)$ is bounded above for $G$ in $\sG$.
\end{itemize}
\end{prop}

\begin{proof}
When $\om\in\Kah_\Z(X)$ (and hence when $\om\in\Kah_\Q(X)$) this is a consequence of a lemma of Grothendieck, see Huybrechts and Lehn~\cite[Lem.~1.7.9]{HuLe2}.

For real $\om$, we prove it using work of Toma \cite{Toma1}. He defines a notion of boundedness for families of complex analytic coherent sheaves on compact K\"ahler manifold $(X,\om)$, and proves boundedness results involving the slopes $\bar\mu(G)$ defined using the K\"ahler form. When $X$ is a projective complex manifold, 
algebraic and complex analytic coherent sheaves coincide by Serre's GAGA principle \cite{Serr}, and \cite[Rem.~3.3]{Toma1} his notion of boundedness agrees with Definition \ref{co7def7}. The theorem then follows from~\cite[Lem.~4.3]{Toma1}.	
\end{proof}

The next condition will be important for proving finite type and properness properties of moduli stacks~$\M_\al^\ss(\tau),\M_\al^\ss(\mu)$. 

\begin{cond}
\label{co7cond1}
Let $K\subset\Kah(X)$ be compact and $\al\in C(\coh(X))$. We say that {\it Condition\/ {\rm\ref{co7cond1}} holds for\/} $K,\al$ if for each $\mu_0\in\R,$ the family of pure sheaves $E\in\coh(X)$ with $\lb E\rb=\al$ and $\mu_{\max}^\om(E)\le\mu_0$ for some $\om\in K$ is bounded.
\end{cond}

The next proposition follows from Huybrechts and Lehn \cite[Th.~3.3.7]{HuLe2}. They prove it for integral classes $\om\in\Kah_\Z(X)$, but the extension to $\om\in\Kah_\Q(X)$ is immediate as rescaling $\om\mapsto N\om$ for $N>0$ rescales~$\mu^\om_{\max}(E)\mapsto N^{-1}\mu^\om_{\max}(E)$.

\begin{prop}
\label{co7prop4}
Let\/ $\om\in\Kah_\Q(X)$. Then Condition\/ {\rm\ref{co7cond1}} holds for\/ $K=\{\om\}$ and all\/ $\al\in C(\coh(X))$.
\end{prop}

The next proposition will be proved in \S\ref{co7112}.

\begin{prop}
\label{co7prop5}
Condition\/ {\rm\ref{co7cond1}} holds for any compact subset\/ $K\subset\Kah(X)$ and all\/ $\al\in C(\coh(X))$ with\/ $\dim\al=0$ or\/~$1$.
\end{prop}

The next proposition will be proved in \S\ref{co7113}. The restriction to $\dim\al=m$ is because the proof uses the Bogomolov inequality, Theorem \ref{co7thm1}, which works only for torsion-free sheaves, and also Proposition \ref{co7prop4}. It is based on the proofs of Greb--Ross--Toma \cite[Th.~6.8]{GrTo1} and Greb--Toma \cite[Prop.~6.3]{GrTo1}, which are corollaries of the proposition.

\begin{prop}
\label{co7prop6}
Condition\/ {\rm\ref{co7cond1}} holds for any compact subset\/ $K\subset\Kah(X)$ and all\/ $\al\in C(\coh(X))$ with\/ $\dim\al=m=\dim X$.
\end{prop}

As we explain in Remark \ref{co7rem5} below, our theory currently has problems for invariants counting sheaves of dimensions $2,3,\ldots,m-1$, for example, dimension 2 sheaves on a Fano 3-fold $X$. Some of these difficulties would be resolved if the answer to the following question is yes.

\begin{quest}
\label{co7quest1}
{\bf(a)} When $m=\dim X\ge 3,$ does Condition\/ {\rm\ref{co7cond1}} hold for any compact\/ $K\subset\Kah(X)$ and all\/ $\al\in C(\coh(X))$ with\/~$\dim\al=2,3,\ldots,m-1$? 
\smallskip

\noindent{\bf(b)} Does it hold for $K=\{\om\}$ with\/ $\om\in\Kah(X)\sm\Kah_\Q(X)$ and all\/ $\al$ in\/ $C(\coh(X))$ with\/ $\dim\al=2,3,\ldots,m-1$?	
\end{quest}

\subsubsection{Openness of Gieseker and $\mu$-(semi)stability}
\label{co715}

The next proposition will be proved in \S\ref{co7114} using Propositions \ref{co7prop2} and \ref{co7prop3}. For $\om\in\Kah_\Q(X)$ and Gieseker stability it is proved in Huybrechts and Lehn \cite[Prop.~2.3.1]{HuLe2}, and for general $\om$ it can also be deduced from Toma~\cite[Cor.~6.7]{Toma3}.

\begin{prop}
\label{co7prop7}
For each\/ $\om\in\Kah(X),$ the conditions that\/ $E\in\coh(X)$ is $\tau^\om$-(semi)stable, or $\mu^\om$-(semi)stable, are open conditions on $[E]\in\M$ and on $[E]\in\M^\pl$. Thus for each\/ $\al\in C(\coh(X))$ we have open substacks $\M_\al^\rst(\mu^\om)\subseteq\M_\al^\rst(\tau^\om)\subseteq\M_\al^\ss(\tau^\om)\subseteq\M_\al^\ss(\mu^\om)\subseteq\M_\al^\pl$ parametrizing $\tau^\om$-(semi)stable, and $\mu^\om$-(semi)stable, sheaves in $\coh(X)$ in class $\al\in C(\coh(X))$.
\end{prop}

\subsubsection{Conditions for $\M_\al^\ss(\tau),\M_\al^\ss(\mu)$ to be of finite type}
\label{co716}

When $\om\in\Kah_\Q(X)$ the next proposition is included in Huybrechts and Lehn \cite[Th.~3.3.7]{HuLe2}. For general $\om$ and $\dim\al=m$ it is part of Greb--Ross--Toma \cite[Th.~6.8]{GRT1} (see also Greb--Toma \cite[Prop.~6.3]{GrTo1}). If the answer to Question \ref{co7quest1}(b) is yes then the proposition holds for all $\om\in\Kah(X)$ and~$\al\in C(\coh(X))$.

\begin{prop}
\label{co7prop8}
Suppose\/ $\om\in\Kah(X)$ and\/ $\al\in C(\coh(X))$ are such that Condition\/ {\rm\ref{co7cond1}} holds for\/ $K=\{\om\}$ and\/ $\al$. By Propositions\/ {\rm\ref{co7prop4}--\ref{co7prop6},} this holds automatically if either\/ $\om\in\Kah_\Q(X),$ or if\/ $\dim\al=0,1$ or $m$. Then the moduli stacks $\M_\al^\ss(\tau^\om),\M_\al^\ss(\mu^\om)$ in Proposition\/ {\rm\ref{co7prop7}} are of finite type.
\end{prop}

\begin{proof}
If $E\in\coh(X)$ with $\lb E\rb=\al$ is $\tau^\om$- or $\mu^\om$-semistable then $E$ is pure of dimension $d=\dim\al$ and $\mu^\om_{\max}(E)=\bar\mu^\om(\al)$. Now the family of $F\in\coh(X)$ with $\lb F\rb=\al$ and $\mu^\om_{\max}(F)\le\bar\mu^\om(\al)$ is bounded as Condition \ref{co7cond1} holds for $K=\{\om\}$ and $\al$. Hence the families of $\tau^\om$- and $\mu^\om$-semistable $E\in\coh(X)$ with $\lb E\rb=\al$ are bounded. Together with Proposition \ref{co7prop7}, this implies that $\M_\al^\ss(\tau^\om),\M_\al^\ss(\mu^\om)$ are of finite type.
\end{proof}

\subsubsection{Wall and chamber structures for Gieseker and $\mu$-stability}
\label{co717}

The following result, proved in \S\ref{co7115}, is an example of a well known picture. Similar results can be found in Schmitt \cite[Main Theorem]{Schm1} and Greb--Ross--Toma \cite[\S 4]{GRT1}, \cite[Prop.~2.7]{GRT2}, \cite[Prop.~3.3]{GRT3}, and \cite[Th.~6.6]{GrTo1}. If the answer to Question \ref{co7quest1}(a) is yes then the theorem holds for all~$\al\in C(\coh(X))$.

\begin{thm}
\label{co7thm2}
Suppose\/ $\al\in C(\coh(X))$ is such that Condition\/ {\rm\ref{co7cond1}} holds for any compact\/ $K\subseteq\Kah(X)$ and\/ $\al$. By Propositions\/ {\rm\ref{co7prop5}} and\/ {\rm\ref{co7prop6},} this is true if\/ $\dim\al=0,1$ or\/ $m$. Then the map 
\e
\begin{split}
\Phi&:\Kah(X)\longra\{\text{pairs of open substacks of\/ $\M_\al^\pl$}\},\\
\Phi&:\om\longmapsto(\M_\al^\rst(\mu^\om),\M_\al^\ss(\mu^\om))
\end{split}
\label{co7eq10}
\e
admits a \begin{bfseries}locally finite wall and chamber structure\end{bfseries}. That is, there is a decomposition $\Kah(X)=\coprod_{i\in I}S_i,$ satisfying:
\begin{itemize}
\setlength{\itemsep}{0pt}
\setlength{\parsep}{0pt}
\item[{\bf(a)}] $S_i$ is a connected, locally closed subset of\/ $\Kah(X)$ for each\/ $i\in I,$ and\/ $\Phi\vert_{S_i}$ is constant, and\/ $S_i$ is maximal under these conditions. The $S_i$ are called \begin{bfseries}chambers\end{bfseries}.
\item[{\bf(b)}] Each compact subset of\/ $\Kah(X)$ intersects only finitely many $S_i$.
\item[{\bf(c)}] In a sufficiently small open neighbourhood\/ $U$ of any point\/ $\om_0\in\Kah(X),$ whether $\om\in U$ lies in $S_i$ depends only on whether $\mu^\om(\be)<\mu^\om(\al-\be),$ or $\mu^\om(\be)=\mu^\om(\al-\be),$ or $\mu^\om(\be)>\mu^\om(\al-\be),$ for a finite set\/ $W_U$ of\/ $\be$ in $C(\coh(X))$ with\/ $\al-\be$ in $C(\coh(X))$. The local hypersurfaces of\/ $\om$ in\/ $\Kah(X)$ with\/ $\mu^\om(\be)=\mu^\om(\al-\be)$ separating chambers are called \begin{bfseries}walls\end{bfseries}.
\item[{\bf(d)}] The map\/ $\om\mapsto\M_\al^\rst(\mu^\om)$ is lower semicontinuous (i.e.\ $\M_\al^\rst(\mu^\om)$ can only get smaller when we deform from `off a wall' to `on a wall') and\/ 
the map $\om\mapsto\M_\al^\ss(\mu^\om)$ is upper semicontinuous (i.e.\ $\M_\al^\ss(\mu^\om)$ can only get larger when we deform from `off a wall' to `on a wall').
\item[{\bf(e)}] If\/ $\M_\al^\rst(\mu^\om)=\M_\al^\ss(\mu^\om)$ with\/ $\om\in S_i$ then $S_i$ is open in $\Kah(X)$.
\end{itemize}
The analogues of\/ {\bf(a)\rm--\bf(c)} also hold if\/ $\mu^\om$ is replaced by Gieseker stability~$\tau^\om$.
\end{thm}

\subsubsection{Proper moduli stacks using Geometric Invariant Theory}
\label{co718}

The next theorem (except for the part about good moduli spaces in (vi)) summarizes Huybrechts and Lehn \cite[\S 4.3]{HuLe2}, who gather together major results by Grothendieck, Mumford, Gieseker, Maruyama, Simpson, and others. An alternative proof is given by \'Alvarez-C\'onsul and King~\cite{AlKi}.

\begin{thm}
\label{co7thm3}
Let\/ $\om\in\Kah_\Q(X)$ and write $(\tau,G,\le)=(\tau^\om,G,\le)$. Then for each\/ $\al\in C(\coh(X)),$ we can construct:
\begin{itemize}
\setlength{\itemsep}{0pt}
\setlength{\parsep}{0pt}
\item[{\bf(i)}] A projective $\C$-scheme $Q,$ which is a Quot scheme {\rm\cite[\S 2.2]{HuLe2}} parametrizing quotients $\C^N\ot \O_X(-n)\twoheadrightarrow E$ for\/ $E\in\coh(X)$ with $\lb E\rb=\al,$ where $\O_X(1)\ra X$ is an ample line bundle with\/ $\om=Mc_1(\O_X(1))$ for $M>0,$ and\/ $n\gg 0$ is fixed, and\/ $N=P^{c_1(\O_X(1))}_\al(n)$.
\item[{\bf(ii)}] A $\GL(N,\C)$-action on $Q,$ induced by the obvious action on $\C^N\ot \O_X(-n)$.\item[{\bf(iii)}] A linearization $L\ra Q$ of the $\GL(N,\C)$-action on $Q$ in the sense of Geometric Invariant Theory (GIT), as in\/ {\rm\cite{MFK}}. Thus we have $\GL(N,\C)$-invariant open subschemes $Q^\rst\subseteq Q^\ss\subseteq Q$ of GIT-(semi)stable points.
\item[{\bf(iv)}] A morphism $\io:[Q/\GL(N,\C)]\ra \M_\al,$ which restricts to isomorphisms
\e
\begin{split}
&[Q^\rst/\GL(N,\C)]\longra (\Pi_\al^\pl)^{-1}(\M^\rst_\al(\tau)),\\
&[Q^\ss/\GL(N,\C)]\longra (\Pi_\al^\pl)^{-1}(\M^\ss_\al(\tau)).
\end{split}
\label{co7eq11}
\e
\item[{\bf(v)}] The subgroup $\bG_m\cdot\Id_N\subset\GL(N,\C)$ acts trivially on $Q$. The corresponding $\bG_m$-subgroups of isotropy groups of\/ $\C$-points in the Artin stack $[Q/\GL(N,\C)]$ are mapped by $\io$ to subgroups $\bG_m\cdot\id_E$ in the isotropy groups $\Aut(E)$ of $\C$-points $[E]$ in $\M_\al$. These $\bG_m\cdot\id_E$ are quotiented out in the rigidification $\Pi_\al^\pl:\M_\al\ra\M_\al^\pl$. Thus \eq{co7eq11} descends to isomorphisms
\e
\begin{split}
&[Q^\rst/\PGL(N,\C)]\longra \M^\rst_\al(\tau),\quad [Q^\ss/\PGL(N,\C)]\longra \M^\ss_\al(\tau).
\end{split}
\label{co7eq12}
\e
\item[{\bf(vi)}] The GIT quotient $Q\ds^\ss_L \PGL(N,\C)$ is a projective $\C$-scheme with a stack morphism $\pi:[Q^\ss/\PGL(N,\C)]\ra Q\ds^\ss_L\PGL(N,\C)$. As in Example\/ {\rm\ref{co3ex5},} $Q\ds^\ss_L \PGL(N,\C)$ is a \begin{bfseries}proper good moduli space\end{bfseries} for $[Q^\ss/\PGL(N,\C)],$ and hence for $ \M^\ss_\al(\tau)$ under the isomorphism \eq{co7eq12}.
\item[{\bf(vii)}] There is an open subscheme $Q\ds^\rst_L \PGL(N,\C)\subseteq Q\ds^\ss_L \PGL(N,\C)$ such that\/ $\pi$ restricts to an isomorphism $[Q^\rst/\PGL(N,\C)]\ra Q\ds^\rst_L\PGL(N,\C),$ and hence an isomorphism $\M^\rst_\al(\tau)\ra Q\ds^\rst_L\PGL(N,\C)$.
\item[{\bf(viii)}] If\/ $\M^\rst_\al(\tau)=\M^\ss_\al(\tau)$ then $\M^\ss_\al(\tau)\cong Q\ds^\ss_L \PGL(N,\C)$ as Artin stacks, and hence $\M^\ss_\al(\tau)$ is a proper algebraic space.
\item[{\bf(ix)}] Part\/ {\bf(viii)} and\/ \eq{co7eq6} imply that if\/ $\M^\rst_\al(\mu^\om)=\M^\ss_\al(\mu^\om)$ then $\M^\ss_\al(\mu^\om)$ is a proper algebraic space.
\end{itemize}
\end{thm}

Here is a criterion for when Assumption \ref{co5ass2}(g) holds for~$(\mu^\om,M,\le)$:

\begin{cor}
\label{co7cor1}
Suppose\/ $\al\in C(\coh(X))$ is such that Condition\/ {\rm\ref{co7cond1}} holds for any compact\/ $K\subseteq\Kah(X)$ and\/ $\al$. By Propositions\/ {\rm\ref{co7prop5}} and\/ {\rm\ref{co7prop6},} this is true if\/ $\dim\al=0,1$ or\/ $m$. Suppose\/ $\om\in\ov{\Kah_\Q}(X)$ as in Definition\/ {\rm\ref{co7def3}} with\/ $\M^\rst_\al(\mu^\om)=\M^\ss_\al(\mu^\om)$. Then $\M^\ss_\al(\mu^\om)$ is a proper algebraic space.
\end{cor}

\begin{proof}
As $\M^\rst_\al(\mu^\om)=\M^\ss_\al(\mu^\om)$, Theorem \ref{co7thm2}(e) implies that $\om$ lies in an open chamber $S_i$ in the wall and chamber decomposition. Since $\om\in\ov{\Kah_\Q}(X)$ we can choose $\ti\om$ in $S_i\cap\Kah_\Q(X)$. Then $\M^\rst_\al(\mu^{\ti\om})=\M^\rst_\al(\mu^\om)$ and $\M^\ss_\al(\mu^{\ti\om})=\M^\ss_\al(\mu^\om)$. Hence $\M^\rst_\al(\mu^{\ti\om})=\M^\rst_\al(\tau^{\ti\om})=\M^\ss_\al(\tau^{\ti\om})=\M^\ss_\al(\mu^{\ti\om})$ by \eq{co7eq6}, so Theorem \ref{co7thm3}(viii) shows $\M^\ss_\al(\tau^{\ti\om})=\M^\ss_\al(\mu^\om)$ is a proper algebraic space.
\end{proof}

In \S\ref{co7111} we give an alternative proof of Assumption \ref{co5ass2}(g) using~\S\ref{co33}.

We can combine Theorems \ref{co7thm2} and \ref{co7thm3} to extend Theorem \ref{co7thm3} to $\om$ in $\ov{\Kah_\Q}(X)$ when $X$ is a curve or a surface:

\begin{cor}
\label{co7cor2}
Suppose $m=\dim X$ is $1$ or\/ $2$. Then Theorem\/ {\rm\ref{co7thm3}} holds for any\/ $\om\in\ov{\Kah_\Q}(X),$ except that we do not require that\/ $\om=Mc_1(\O_X(1))$ in~{\bf(i)}.

In particular, this implies that for all $\al\in C(\coh(X))$ the moduli stack\/ $\M^\ss_\al(\tau^\om)$ has a proper good moduli space, constructed using GIT.
\end{cor}

\begin{proof}
Given $\om\in\ov{\Kah_\Q}(X)$, we wish to find nearby $\ti\om\in\Kah_\Q(X)$ such that $\M^\rst_\al(\tau^\om)=\M^\rst_\al(\tau^{\ti\om})$ and $\M^\ss_\al(\tau^\om)=\M^\ss_\al(\tau^{\ti\om})$. The corollary then follows from Theorem \ref{co7thm3} for $\ti\om$. If $m=1$ then as $X$ is connected, $\tau^\om$-(semi)stability is independent of $\om\in\Kah(X)$, so any $\ti\om\in\Kah_\Q(X)$ will do.

In Theorem \ref{co7thm2}, the wall-and-chamber structure for Gieseker stability is defined using `walls' of the form $\tau^\om(\be)=\tau^\om(\al-\be)$. When $m=2$, such equations are either true for all $\om$, or false for all $\om$, or are rational linear equations on $\om\in\ov{\Kah_\Q}(X)$. Therefore when $m=2$ the chambers $S_i$ in Theorem \ref{co7thm2} are defined by rational linear equations, and every $S_i$ contains a rational point. Thus for each $\om\in\ov{\Kah_\Q}(X)$ we can choose $\ti\om\in\Kah_\Q(X)$ with $\om,\ti\om$ in the same chamber $S_i$, and the corollary follows.
\end{proof}

If $m=\dim X\ge 3$ then an example of Schmitt \cite[Ex.~1.1.5]{Schm1} shows that walls in $\Kah(X)$ of the form $\tau^\om(\be)=\tau^\om(\al-\be)$ or $\mu^\om(\be)=\mu^\om(\al-\be)$ may contain no rational points, so there can be chambers $S_i$ in Theorem \ref{co7thm2} containing no points $\ti\om\in\Kah_\Q(X)$, and the method of Corollary \ref{co7cor2} fails.

The next proposition follows from Greb--Ross--Toma \cite[Th.~11.6]{GRT1}, who use an alternative way to present $\M^\ss_\al(\tau^\om)$ as a GIT quotient when~$m=\dim\al=3$.

\begin{prop}
\label{co7prop9}
Suppose $m=\dim X=3$. Then Theorem\/ {\rm\ref{co7thm3}} holds for any\/ $\om\in\ov{\Kah_\Q}(X)$ and\/ $\al\in C(\coh(X))$ with\/ $\dim\al=3,$ except that we do not require that\/ $\om=Mc_1(\O_X(1))$ in~{\bf(i)}.

In particular, this implies that for all $\al\in C(\coh(X))$ with $\dim\al=3$ the moduli stack\/ $\M^\ss_\al(\tau^\om)$ has a proper good moduli space, constructed using GIT.	
\end{prop}

\subsubsection{Valuative criterion for universal closedness for $\M_\al^\ss(\tau),\M_\al^\ss(\mu)$}
\label{co719}

When $\om\in\Kah_\Z(X)$, the next theorem follows from Huybrechts and Lehn \cite[Th.~2.B.1]{HuLe2}, who generalize results of Langton \cite[\S 3]{Lang} for $\mu^\om$-semistable torsion-free sheaves and Maruyama \cite[Th.~5.7]{Maru} for $\tau^\om$-semistable torsion-free sheaves. In \S\ref{co7116} we explain how to modify the proof of \cite[Th.~2.B.1]{HuLe2} to~$\om\in\Kah(X)$.

\begin{thm}
\label{co7thm4}
Suppose\/ $\om\in\Kah(X)$ and\/ $\al\in C(\coh(X))$ are such that Condition\/ {\rm\ref{co7cond1}} holds for\/ $K=\{\om\}$ and\/ $\al$. By Propositions\/ {\rm\ref{co7prop4}--\ref{co7prop6},} this is true if either\/ $\om\in\Kah_\Q(X),$ or\/ $\dim\al=0,1$ or\/ $m$. Then $(\Pi_\al^\pl)^{-1}(\M_\al^\ss(\tau^\om))$ and\/ $(\Pi_\al^\pl)^{-1}(\M_\al^\ss(\mu^\om))$ satisfy the valuative criterion for universal closedness in Definition\/ {\rm\ref{co3def8}(a)} with respect to any discrete valuation ring $R$ over\/~$\C$.
\end{thm}

\begin{rem}
\label{co7rem4}
{\bf(a)} We could also try to prove Theorem \ref{co7thm4} using the method of \S\ref{co334}, by constructing a pseudo-$\Th$-stratification of $\M_\al$ with semistable locus $(\Pi_\al^\pl)^{-1}(\M_\al^\ss(\tau^\om))$, generalizing Example \ref{co3ex6} from $\om$ in $\Kah_\Z(X)$ to $\om$ in $\Kah(X)$. But the direct proof seems simpler in this case.
\smallskip

\noindent{\bf(b)} Toma \cite{Toma2} proves a complex analytic version of Theorem \ref{co7thm4} for compact K\"ahler manifolds $(X,\om)$. It may be possible to deduce Theorem \ref{co7thm4} from \cite[Th.~3.1]{Toma2}, without assuming Condition \ref{co7cond1} holds for~$\{\om\},\al$.
\end{rem}

\subsubsection{Gieseker stability is additive}
\label{co7110}

\begin{prop}
\label{co7prop10}
For each\/ $\om\in\Kah(X),$ Gieseker stability $(\tau^\om,G,\le)$ is an \begin{bfseries}additive\end{bfseries} stability condition on\/ $\coh(X),$ in the sense of\/~{\rm\S\ref{co335}}.
\end{prop}

\begin{proof}
Fix $\al\in C(\coh(X))$ with $\dim\al=d\le m$. Define $V=\R[1,t,\ldots,t^m]$  to be the $\R$-vector space of polynomials $P(t)$ of degree at most $m$. For such $P$, write $P^{(d)}$ for the coefficient of $t^d$ in $P$. Define a total order $\le$ on $V$ by $P(t)\le Q(t)$ if $P(n)\le Q(n)$ for $n\gg 0$. Define a group morphism $\rho:K(\coh(X))\ra V$ by
\e
\rho(\be)=P_\al^{\om,(d)}P^\om_\be-P_\be^{\om,(d)}P^\om_\al.
\label{co7eq13}
\e

Suppose $E\in\coh(X)$ with $\lb E\rb=\al$, and $0\ne F\subseteq E$ with $\lb F\rb=\be$. Then $\dim\be\le\dim\al=d$. If $\dim\be<d$ then $\rho(\be)=P_\al^{\om,(d)}P^\om_\be$, so $\rho(\be)>0$ as $P_\al^{\om,(d)}>0$ and $P^\om_\be$ has positive leading coefficient. If $\dim\be=d$ then $\rho(\be)=P_\al^{\om(d)}P_\be^{\om(d)}(\tau^\om(\be)-\tau^\om(\al))$ where $P_\al^{\om,(d)},P_\be^{\om,(d)}>0$, so $\rho(\be)>0$ if and only if $\tau^\om(\be)>\tau^\om(\al)$ as $\deg\tau^\om(\be)=\deg\tau^\om(\al)=d$. Using these we see that $E$ is $\rho$-semistable if and only if $E$ is $\tau^\om$-semistable, so $(\Pi_\al^\pl)^{-1}(\M_\al^\ss(\tau^\om))=\M_\al^{\ss,\rho}$, and $(\tau^\om,G,\le)$ is additive by Definition~\ref{co3def13}.
\end{proof}

\subsubsection{Proper good moduli spaces for $\M_\al^\ss(\tau)$ using \S\ref{co33}}
\label{co7111}

The next theorem follows from Proposition \ref{co3prop6}, Theorem \ref{co3thm4}, Corollary \ref{co3cor3}, Propositions \ref{co7prop1}, \ref{co7prop8} and \ref{co7prop10}, and Theorem \ref{co7thm4}. The last part is proved as for Theorem \ref{co3thm9}(b), since $(\mu^\om,M,\le)$ dominates $(\tau^\om,G,\le)$. If $\om\in\Kah_\Q(X)$ the result follows from Theorem \ref{co7thm3}. If $\om\in\ov{\Kah_\Q}(X)$ then Corollary \ref{co7cor1} provides an alternative proof of the final part.

\begin{thm}
\label{co7thm5}
Suppose\/ $\om\in\Kah(X)$ and\/ $\al\in C(\coh(X))$ are such that Condition\/ {\rm\ref{co7cond1}} holds for\/ $K=\{\om\}$ and\/ $\al$. This is automatic if\/ $\dim\al=0,1$ or\/ $m$. Then $\M_\al^\ss(\tau^\om)$ admits a proper good moduli space. If\/ $\M^\rst_\al(\tau^\om)=\M^\ss_\al(\tau^\om)$ then $\M^\ss_\al(\tau^\om)$ is a proper algebraic space. If\/ $\M^\rst_\al(\mu^\om)=\M^\ss_\al(\mu^\om)$ then $\M^\ss_\al(\mu^\om)$ is a proper algebraic space.
\end{thm}

The fact that $\M_\al^\ss(\tau^\om)$ admits a proper good moduli space for general $\om$ in $\Kah(X)$ (at least when $\dim\al=0,1$ or $m$) seems important. As in Corollary \ref{co7cor2} and Proposition \ref{co7prop9}, this is already known when $m=1$ or 2, and when $m=\dim\al=3$. Greb--Toma \cite{GrTo2} also prove it when $\dim\al=m\ge 4$ and $\rank\al=2$, using criteria for existence of good moduli spaces in Alper--Fedorchuk--Smyth \cite{AFS}, which are weaker than those in Alper--Halpern-Leistner--Heinloth \cite{AHLH} that we use. But the case $m=\dim\al\ge 4$ and $\rank\al>2$ may be new. 

In a related result, Greb--Ross--Toma \cite[Th.s 9.4, 9.6 \& Rem.~9.7]{GRT1} prove, in our language, that if $\al\in C(\coh(X))$ is such that Condition \ref{co7cond1} holds for $\al$ and all small compact $K\subset\Kah(X)$ then $\M_\al^\ss(\tau^{\rm mG})$ admits a proper good moduli space, which is a projective GIT quotient, for all `multi-Gieseker' stability conditions $(\tau^{\rm mG},T,\le)$ on $\coh(X)$, as in \cite[Def.~2.5]{GRT1}. The rest of the section proves Propositions \ref{co7prop5}, \ref{co7prop6} and \ref{co7prop7} and Theorems \ref{co7thm2} and~\ref{co7thm4}.

\subsubsection{Proof of Proposition \ref{co7prop5}}
\label{co7112}

If $\dim\al=0$ then for any $\om\in\Kah(X)$ and $E\in\coh(X)$ with $\lb E\rb=\al$ we have $\bar\mu^\om(E)=\mu^\om_{\max}(E)=0$. So taking $\om\in\Kah_\Q(X)$, Proposition \ref{co7prop4} shows that the family of all $E\in\coh(X)$ with $\lb E\rb=\al$ is bounded, as all such $E$ are $\mu^\om$-semistable. Therefore Condition \ref{co7cond1} holds for any compact $K\subset\Kah(X)$ and $\al\in C(\coh(X))$ with~$\dim\al=0$.

Let $K\subset\Kah(X)$ be compact, and fix $\om_0\in\Kah_\Q(X)$. As $K$ is compact there exists $C>0$ such that $\om-C^{-1}\om_0,C\om_0-\om\in\Kah(X)$ for all $\om\in K$. Let $\be\in C(\coh(X))$ with $\dim\be=1$, and write $\be=\be_{m-1}+\be_m$ with $\be_k\in H^{2k}(X,\Q)$. Then for $\om\in \Kah(X)$ we have
\e
\begin{split}
P_\be^\om(t)&=t\int_X\be_{m-1}\cup\om+\int_X(\be_m+\ha\be_{m-1}\cup c_1(X)),\\
\bar\mu^\om(\be)&=\frac{\int_X(\be_m+\ha\be_{m-1}\cup c_1(X))}{\int_X\be_{m-1}\cup\om},
\end{split}
\label{co7eq14}
\e
where $\int_X\be_{m-1}\cup\om>0$. When $\om\in K$ we have $\om-C^{-1}\om_0,C\om_0-\om\in\Kah(X)$, so $\int_X\be_{m-1}\cup(\om-C^{-1}\om_0)>0$ and $\int_X\be_{m-1}\cup(C\om_0-\om)>0$, giving 
\begin{equation*}
0<C^{-1}\int_X\be_{m-1}\cup\om_0<\int_X\be_{m-1}\cup\om<C\int_X\be_{m-1}\cup\om_0.
\end{equation*}
Combined with \eq{co7eq14} this implies that
\begin{equation*}
\sign(\bar\mu^{\om_0}(\be))=\sign(\bar\mu^\om(\be)),\quad C^{-1}\bmd{\bar\mu^{\om_0}(\be)}\le \bmd{\bar\mu^\om(\be)}\le C\bmd{\bar\mu^{\om_0}(\be)}.
\end{equation*}

Now let $\al\in C(\coh(X))$ with $\dim\al=1$, and $\om\in K$, and $\mu_0\in\R$. Suppose $E\in\coh(X)$ with $\lb E\rb=\al$ and $\mu_{\max}^\om(E)\le\mu_0$ for some $\om\in K$. Then if $0\ne E'\subseteq E$ we have
\begin{equation*}
\bar\mu^{\om_0}(\lb E'\rb)\le C\max\bigl(\bar\mu^\om(\lb E'\rb),0\bigr)\le C\max\bigl(\bar\mu^\om_{\max}(E),0\bigr)\le C\max(\mu_0,0).
\end{equation*}
As this holds for all $0\ne E'\subseteq E$ we see that $\mu_{\max}^{\om_0}(E)\le C\max(\mu_0,0)$. Thus the family of all such $E$ is bounded by Proposition \ref{co7prop4}, and the proposition~follows.

\subsubsection{Proof of Proposition \ref{co7prop6}}
\label{co7113}

Let $K\subset\Kah(X)$ be compact and $\al\in C(\coh(X))$ with $\dim\al=m$, so that $\rank\al>0$. Set $r=\rank\al$, and fix $\mu_0\in\R$. Choose $\ti\om\in\Kah_\Q(X)$. By making the compact subset $K$ larger, we can suppose that if $\om\in K$ then $(1-t)\om+t\ti\om\in K$ for all $t\in[0,1]$. Define $\sF_{\al,K,\mu_0}$ to be the family of all $E\in\coh(X)$ such that $E$ is torsion-free with $\lb E\rb=\al$ and for some $\om\in K$ we have $\mu_{\max}^\om(E)\le\mu_0$. We must prove that $\sF_{\al,K,\mu_0}$ is bounded. We will do this by constructing $\ti\mu_0\in\R$ such that $\mu_{\max}^{\ti\om}(E)\le\ti\mu_0$ for all $E\in\sF_{\al,K,\mu_0}$, and then boundedness of $\sF_{\al,K,\mu_0}$ follows from Proposition \ref{co7prop4} as~$\ti\om\in\Kah_\Q(X)$.

Let $E\in\sF_{\al,K,\mu_0}$ and $\om\in K$ with $\mu_{\max}^\om(E)\le\mu_0$. Write $\om_t=(1-t)\om+t\ti\om$ for $t\in[0,1]$, so that $\om_t\in K$ with $\om_0=\om$ and $\om_1=\ti\om$. Using ideas from Greb--Toma \cite[Proof of Lem.~6.4]{GrTo1} and Greb--Ross--Toma \cite[Proof of Th.~6.8]{GRT1}, for each $t\in[0,1]$ we construct a unique filtration of $E$ in~$\coh(X)$:
\e
0=E^t_0\subsetneq E^t_1\subsetneq\cdots\subsetneq E^t_{n^t}=E
\label{co7eq15}
\e
and write $F_i^t=E_i^t/E_{i-1}^t$ and $\be_i^t=\lb F_i^t\rb\in C(\coh(X))$ for $i=1,\ldots,n^t$, with the following properties:
\begin{itemize}
\setlength{\itemsep}{0pt}
\setlength{\parsep}{0pt}
\item[(i)] When $t=0$, \eq{co7eq15} is the $\mu^\om$-Harder--Narasimhan filtration of $E$, which exists by Theorem \ref{co3thm1} and Lemma \ref{co7lem1}. As $E$ is torsion-free, this implies that $F_i^0$ is torsion-free and $\mu^\om$-semistable with $\bar\mu^\om(\be_1^0)>\cdots>\bar\mu^\om(\be_{n^0}^0)$. Also $\bar\mu^\om(\be_1^0)=\mu_{\max}^\om(E)$, so we see that $\bar\mu^\om(\be_i^0)\le\mu_0$ for $i=1,\ldots,n^0$.
\item[(ii)] For all $t\in[0,1]$ and $i=1,\ldots,n^t$, $F_i^t$ is torsion-free and $\mu^{\om_t}$-semistable. Hence $\rank F_i^t=\rank\be_i^t>0$. As $\sum_{i=1}^{n^t}\rank\be_i^t=\rank\al=r$, this implies that $n^t\le r$.
\item[(iii)] There exist $0=t_0\le t_1<t_2<\cdots<t_l=1$ such that the filtration \eq{co7eq15} is independent of $t$ for $t$ in each of the intervals $[t_0,t_1],(t_1,t_2],\ldots,(t_{l-1},t_l]$. For brevity we make the convention that when we write $(t_{k-1},t_k]$ below, when $k=1$ we mean~$[t_0,t_1]$.
\item[(iv)] For each $0<k<l$, the filtration \eq{co7eq15} for $t\in(t_k,t_{k+1}]$ is obtained from \eq{co7eq15} for $t\in(t_{k-1},t_k]$ in the following way. There should exist at least one $i=1,\ldots,n^{t_k}$ such that $F_i^{t_k}$ is $\bar\mu^{\om_t}$-semistable for $t\in(t_{k-1},t_k]$, but is not $\bar\mu^{\om_t}$-semistable for $t\in(t_k,t_{k+1}]$, and the $\mu^{\om_t}$-Harder--Narasimhan filtration of $F_i^{t_k}$ for $t\in(t_k,t_{k+1}]$ should be independent of $t$. Write 
\e
0=F^{t_k}_{i,0}\subsetneq F^{t_k}_{i,1}\subsetneq\cdots\subsetneq F^{t_k}_{i,n_i^{\smash{t_k}}}=F_i^{t_k}
\label{co7eq16}
\e
for this $\mu^{\om_t}$-Harder--Narasimhan filtration, and $G^{t_k}_{i,j}=F^{t_k}_{i,j}/F^{t_k}_{i,j-1}$ for $j=1,\ldots,n_i^{t_k}>1$. Then $G^{t_k}_{i,j}$ is torsion-free and $\bar\mu^{\om_t}$-semistable for $t\!\in\![t_k,t_{k+1}]$ (note the $[t_k,t_{k+1}]$ here), with $\bar\mu^{\om_t}(\lb G^{t_k}_{i,1}\rb)\!>\!\cdots\!>\!\bar\mu^{\om_t}(\lb G^{t_k}_{i,n_i^{\smash{t_k}}}\rb)$ for $t\!\in\!(t_k,t_{k+1}]$, and $\bar\mu^{\om_{t_k}}(\lb G^{t_k}_{i,1}\rb)=\cdots=\bar\mu^{\om_{t_k}}(\lb G^{t_k}_{i,n_i^{\smash{t_k}}}\rb)=\bar\mu^{\om_{t_k}}(\lb F^{t_k}_i\rb)$.

We make \eq{co7eq15} for $t\in(t_k,t_{k+1}]$ by refining \eq{co7eq15} for $t\in(t_{k-1},t_k]$, inserting extra subobjects in the sequence, in such a way that $F_i^{t_k}$ is replaced by its filtration \eq{co7eq16} for each $i=1,\ldots,n^{t_k}$ with $F_i^{t_k}$ not $\bar\mu^{\om_t}$-semistable for $t\in(t_k,t_{k+1}]$. Then the sequence $F_1^t,\ldots,F_{n^t}^t$ changes as $t$ crosses $t_k$, by replacing $F_i^{t_k}$ by the sequence $G^{t_k}_{i,1},\ldots,G^{t_k}_{i,n_i^{\smash{t_k}}}$. Note that $n^t$ strictly increases as $t$ crosses $t_k$. As $n^t\le r$ by (ii) we see that~$l\le r$. 

Observe that the list $\be_1^{t_{k+1}},\ldots,\be_{n^{t_{k+1}}}^{t_{k+1}}$ is obtained from $\be_1^{t_k},\ldots,\be_{n^{t_k}}^{t_k}$ by replacing $\be_i^{t_k}=\lb F_i^{t_k}\rb$ for some $i=1,\ldots,n^{t_k}$ by $\lb G^{t_k}_{i,1}\rb,\ldots,\lb G^{t_k}_{i,n_i^{\smash{t_k}}}\rb$, where $\bar\mu^{\om_{t_k}}(\lb G^{t_k}_{i,1}\rb)=\cdots=\bar\mu^{\om_{t_k}}(\lb G^{t_k}_{i,n_i^{\smash{t_k}}}\rb)=\bar\mu^{\om_{t_k}}(\lb F^{t_k}_i\rb)$. Thus we see that
\e
\bigl\{\bar\mu^{\om_{t_k}}(\be^{t_k}_i):i=1,\ldots,n^{t_k}\bigr\}=\bigl\{\bar\mu^{\om_{t_k}}(\be^{t_{k+1}}_i):i=1,\ldots,n^{t_{k+1}}\bigr\}.
\label{co7eq17}
\e
\end{itemize}

Greb--Ross--Toma \cite[Th.~6.8]{GRT1}, \cite[Lem.~6.4]{GrTo1} construct families of filtrations \eq{co7eq15} satisfying (ii)--(iv) in the same way, except that they suppose $E$ is $\mu^\om$-semistable, so they start at $t=0$ with the trivial filtration $0=E_0^0\subsetneq E_1^1=E$, which is the $\mu^\om$-Harder--Narasimhan filtration of $E$ in this special case. In particular, the existence and uniqueness of data \eq{co7eq15} satisfying (i)--(iv) follows as in \cite{GRT1,GrTo1} using~\cite[Lem.~6.2]{GrTo1}.

Fix a positive definite inner product $\an{\,,\,}$ on $H^{1,1}(X,\R)$, and let $\nm{\,.\,}$ be the induced norm. By properties of compact K\"ahler manifolds, if $\om\in\Kah(X)$ then the following quadratic form $Q_\om$ on $H^{1,1}(X,\R)$ is positive definite:
\e
Q_\om(\ga_1,\ga_2)=2\frac{\bigl(\int_X\ga_1\cup\om^{m-1}\bigr)\bigl(\int_X\ga_2\cup\om^{m-1}\bigr)}{\int_X\om^m}-\int_X\ga_1\cup\ga_2\cup\om^{m-2}.
\label{co7eq18}
\e
Hence there exists $A>0$ such that for all $\ga\in H^{1,1}(X,\R)$ we have
\e
Q_\om(\ga,\ga)\ge A\nm{\ga}^2.
\label{co7eq19}
\e
As $K$ is compact, we can choose $A>0$ such that this holds for all $\om\in K$. Also we choose $B>0$ such that for all $\om\in K$ we have
\e
\int_X\om^m\ge B.
\label{co7eq20}
\e

By induction on $k=1,\ldots,r$ we will prove the following claims:
\begin{itemize}
\setlength{\itemsep}{0pt}
\setlength{\parsep}{0pt}
\item[$(\dag)_k$] There exists $C_k\ge 0$ such that for all $E,\om,\om_t,t_i,F_i^t,\be_i^t$ as above with $l\ge k$ we have $\bar\mu^{\om_{t_{k-1}}}(\be^{t_k}_i)\le C_k$ for all $i=1,\ldots,n^{t_k}$.
\item[$(\ddag)_k$] There exists a finite subset $S_k$ in the intersection of $H^{1,1}(X,\R)$ with the image of $H^2(X,\Z)\ra H^2(X,\R)$, such that for all $E,\om,\om_t,t_i,F_i^t,\be_i^t$ as above with $l\ge k$ we have $c_1(\be^{t_k}_i)\in S_k$ for all $i=1,\ldots,n^{t_k}$.
\end{itemize}
As $t_0=0$ and $\om=\om_0$, part (i) and \eq{co7eq17} for $k=0$ implies that $(\dag)_1$ holds with $A_1=\max(\mu_0,0)$, which is the first step in our induction. 

For $\om\in\Kah(X)$ and $\ga\in C(\coh(X))$ with $\dim\ga=m$, as $\td(X)=1+\ha c_1(X)+\cdots,$ calculation using \eq{co7eq5} shows that 
\e
\bar\mu^\om(\ga)=\frac{\int_Xc_1(\ga)\cup\om^{m-1}}{\rank\ga\int_X\om^m}+\frac{\int_Xc_1(X)\cup\om^{m-1}}{2\int_X\om^m}.
\label{co7eq21}
\e
From this with $\om_{t_{k-1}},\be^{t_k}_i$ in place of $\om,\ga$, and using $\rank\be^{t_k}_i\le r$ and $\om_{t_{k-1}}\in K$ for $K$ compact, it is easy to see that $(\ddag)_k$ implies $(\dag)_k$. But we will go the other way: for the inductive step we will first show that $(\dag)_k$ implies $(\ddag)_k$, and then show that $(\ddag)_k$ implies $(\dag)_{k+1}$ if~$k<r$.

First suppose that $k=1,\ldots,r$ and $(\dag)_k$ holds for $C_k$. We will prove $(\ddag)_k$. Substituting \eq{co7eq21} for $\om_{t_{k-1}},\be^{t_k}_i$ into $\bar\mu^{\om_{t_{k-1}}}(\be^{t_k}_i)\le C_k$, multiplying by $\rank\be^{t_k}_i\int_X\om_{t_{k-1}}^m$, and using $\rank\be^{t_k}_i\le r=\rank\al$ gives
\e
\int_Xc_1(\be^{t_k}_i)\cup\om_{t_{k-1}}^{m-1}\le r\biggl(C_k\int_X\om_{t_{k-1}}^m-\ha\min\biggl(\int_Xc_1(X)\cup\om_{t_{k-1}}^{m-1},0\biggr)\biggr).
\label{co7eq22}
\e
As $c_1(\al)=c_1(\be^{t_k}_1)+\cdots+c_1(\be^{t_k}_{n^{t_k}})$ with $n^{t_k}\le r$ we deduce that
\ea
-\int_Xc_1(\be^{t_k}_i)\cup\om_{t_{k-1}}^{m-1}&\le r(r-1)\biggl(C_k\int_X\om_{t_{k-1}}^m-\ha\min\biggl(\int_Xc_1(X)\cup\om_{t_{k-1}}^{m-1},0\biggr)\biggr)
\nonumber\\
&\qquad -\int_Xc_1(\al)\cup\om_{t_{k-1}}^{m-1}.
\label{co7eq23}
\ea
Since $K$ is compact and $\om_{t_{k-1}}\in K$ we can choose an upper bound $C_k'$ for the right hand sides of \eq{co7eq22}--\eq{co7eq23} for all $\om_{t_{k-1}}$. This yields
\e
\Big\vert\int_Xc_1(\be^{t_k}_i)\cup\om_{t_{k-1}}^{m-1}\Big\vert\le C_k'.
\label{co7eq24}
\e

The Bogomolov inequality \eq{co7eq9} for $F^{t_k}_i$ (which is $\mu^{\om_{t_{k-1}}}$-semistable by (iv)) and \eq{co7eq8} imply that
\begin{equation*}
\int_X\De(\be^{t_k}_i)\cup\om_{t_{k-1}}^{m-2}=\int_X\bigl(c_1(\be^{t_k}_i)^2-2\rank\be^{t_k}_i\ch_2(\be^{t_k}_i)\bigr)\cup\om_{t_{k-1}}^{m-2}\ge 0.
\end{equation*}
Dividing this by $\rank\be^{t_k}_i>0$, summing from $i=1,\ldots,n^{t_k}$, and using that $\sum_{i=1}^{n^{t_k}}\ch_2(\be^{t_k}_i)=\ch_2(\al)$ shows that
\e
\sum_{i=1}^{n^{t_k}}\frac{1}{\rank\be^{t_k}_i}\int_Xc_1(\be^{t_k}_i)^2\cup\om_{t_{k-1}}^{m-2}\ge 2\int_X\ch_2(\al)\cup\om_{t_{k-1}}^{m-2}.
\label{co7eq25}
\e
As $K$ is compact and $\om_{t_{k-1}}\in K$ we may choose a lower bound $C''_k$ for the right hand side of \eq{co7eq25} for all $\om_{t_{k-1}}$. Then
\e
\sum_{i=1}^{n^{t_k}}\frac{1}{\rank\be^{t_k}_i}\int_Xc_1(\be^{t_k}_i)^2\cup\om_{t_{k-1}}^{m-2}\ge C''_k.
\label{co7eq26}
\e

Combining \eq{co7eq18}--\eq{co7eq20}, \eq{co7eq24} and \eq{co7eq26}, we deduce that
\begin{equation*}
A\sum_{i=1}^{n^{t_k}}\frac{1}{\rank\be^{t_k}_i}\nm{\be^{t_k}_i}^2\le \sum_{i=1}^{n^{t_k}}\frac{1}{\rank\be^{t_k}_i}Q_{\om_{t_{k-1}}}(c_1(\be^{t_k}_i),c_1(\be^{t_k}_i))\le \frac{2r(C_k')^2}{B}-C_k''.
\end{equation*}
As $\rank\be^{t_k}_i\le r$ we see that for $i=1,\ldots,n^{t_k}$
\e
\nm{c_1(\be^{t_k}_i)}^2\le rA^{-1}(2r(C_k')^2/B-C_k'').
\label{co7eq27}
\e
Now $c_1(\be^{t_k}_i)$ lies in the intersection of $H^{1,1}(X,\R)$ with the image of $H^2(X,\Z)\ra H^2(X,\R)$, a discrete lattice. Thus \eq{co7eq27} implies that there are only finitely many possibilities for $c_1(\be^{t_k}_i)$, taken over all $E\in\sF_{\al,K,\mu_0}$ as above. Write $S_k$ for the finite set of $c_1(\be^{t_k}_i)$ satisfying \eq{co7eq27}. Then $(\ddag)_k$ holds.

Next suppose that $k=1,\ldots,r-1$ and $(\ddag)_k$ holds for $S_k$. We will prove $(\dag)_{k+1}$. Define
\begin{equation*}
C_{k+1}=\max\Bigl\{0,\sup_{\om\in K}\frac{\int_Xc_1(\ga)\cup\om^{m-1}}{\int_X\om^m}:\ga\in S_k\Bigr\}.
\end{equation*}
This makes sense as $S_k$ is finite and $K$ is compact. Let $E\in\sF_{\al,K,\mu_0}$, $\om,\ab\om_t,\ab t_i,\ab F_i^t,\ab\be_i^t$ be as above with $l>k$. Then for $i=1,\ldots,n^{t_k}$ we have
\e
\bar\mu^{\om_{t_k}}(\be^{t_k}_i)=\frac{\int_Xc_1(\be^{t_k}_i)\cup\om_{t_k}^{m-1}}{\rank\be^{t_k}_i\int_X\om_{t_k}^m}\le C_{k+1},
\label{co7eq28}
\e
as $c_1(\be^{t_k}_i)\in S_k$, $\om_{t_k}\in K$ and $\rank\be^{t_k}_i\ge 1$. Then $(\dag)_{k+1}$ follows from \eq{co7eq17} and \eq{co7eq28}. This completes the inductive step. Hence $(\dag)_k,(\ddag)_k$ hold for $k=1,\ldots,r$.

Now define
\begin{equation*}
\ti\mu_0=\max\Bigl\{0,\frac{\int_Xc_1(\ga)\cup\ti\om^{m-1}}{\int_X\ti\om^m}:\ga\in S_1\cup\cdots\cup S_r\Bigr\}.
\end{equation*}
Then if $E\in\sF_{\al,K,\mu_0}$ and $\om,\om_t,t_i,F_i^t,\be_i^t$ are as above, as for \eq{co7eq28} we have
\begin{equation*}
\bar\mu^{\ti\om}(\be^{t_l}_i)\le\ti\mu_0, \quad i=1,\ldots,n^{t_l}.
\end{equation*}
But as $t_l=1$ and $\om_{t_l}=\om_1=\ti\om$, the $F_i^{t_l}=F_i^1$ for $i=1,\ldots,n^{t_l}$ are $\mu^{\ti\om}$-semistable. Hence $E$ has a filtration \eq{co7eq15} for $t=1$ whose subquotients $F_i^1$ are $\mu^{\ti\om}$-semistable with $\bar\mu^{\ti\om}(\lb F_i^1\rb)\le\ti\mu_0$. Thus by properties of stability conditions we have $\mu_{\max}^{\ti\om}(E)\le\ti\mu_0$. Proposition \ref{co7prop6} then follows from Proposition \ref{co7prop4} as in the beginning of the proof.

\subsubsection{Proof of Proposition \ref{co7prop7}}
\label{co7114}

Let $\al\in C(\coh(X))$ with $\dim\al=d$. For a sheaf $E\in\coh(X)$ to be pure of dimension $d$ is an open condition on $[E]\in\coh(X)$ by \cite[Prop.~2.3.1]{HuLe2}. Thus we have an open substack $\M_\al^{\rm pu}=\M_\al^\ss(\pi)\subseteq\M_\al^\pl$ parametrizing pure sheaves of dimension $d$ in class $\al$. Let $\Up:U\ra\M^{\rm pu}_\al$ be an atlas for $\M^{\rm pu}_\al$. Then we can take the $\C$-scheme $U$ to be locally of finite type, as $\M^{\rm pu}_\al$ is. 

Thus $U$ can be covered by finite type open subschemes $V\subseteq U$. For such $V$, the family $\sF$ of sheaves $F\in\coh(X)$ pure of dimension $d$ with $[F]=\Up(v)$ for $\C$-points $v\in V$ is bounded, by Definition \ref{co7def7}. Such $F$ is $\tau^\om$-unstable (or $\mu^\om$-unstable) if and only if there exists a quotient $\pi:F\twoheadrightarrow G$ with $G$ pure of dimension $d$ and $\tau^\om([\Ker\pi])>\tau^\om([G])$ (or $\mu^\om([\Ker\pi])>\mu^\om([G])$, respectively). It is not $\tau^\om$-stable (or not $\mu^\om$-stable) if the same hold with $\ge$ rather than $>$, with $G\ne F$. In all four cases this implies that~$\bar\mu(G)\le\bar\mu(\al)$. 

By Proposition \ref{co7prop3}, the family $\sG$ of quotients $\pi:F\twoheadrightarrow G$ with $F$ in $\sF$ and $G$ pure of dimension $d$ with $G\ne 0,F$ and $\bar\mu(G)\le\bar\mu(\al)$ is bounded. Write $S=\bigl\{\lb G\rb:G\in\sG\bigr\}$, so that $S$ is finite as $\sG$ is bounded. Now Proposition \ref{co7prop2} implies that for each $\be\in S$, the condition that there does not exist a quotient $\pi:F\twoheadrightarrow G$ with $G$ pure of dimension $d$ and $\lb G\rb=\be$ is an open condition on $[F]\in\M_\al^{\rm pu}$. As $S$ is finite, we see that for $\C$-points $v\in V$ with $[F]=\Up(v)$, the conditions that $F$ is $\tau^\om$-(semi)stable or $\mu^\om$-(semi)stable are open conditions on $v$, as they are the intersection of at most $\md{S}$ open conditions. Because we can cover $U$ by such open $V\subseteq U$, the same holds for $\C$-points~$v\in U$. 

Since $\Up:U\ra\M^{\rm pu}_\al$ is an atlas for $\M^{\rm pu}_\al$, it follows that for $[F]\in\M^{\rm pu}_\al$ to be $\tau^\om$-(semi)stable, or $\mu^\om$-(semi)stable are open conditions on $[F]$ in $\M^{\rm pu}_\al$. Because $\tau^\om$-(semi)stable and $\mu^\om$-(semi)stable sheaves are pure, and $\M^{\rm pu}_\al$ is open in $\M_\al^\pl$, the same holds with $\M_\al^\pl$ in place of $\M^{\rm pu}_\al$. The proposition follows. 

\subsubsection{Proof of Theorem \ref{co7thm2}}
\label{co7115}

Let $\Phi$ be as in \eq{co7eq10}, and define $T=\bigl\{\Phi^{-1}(\M^\rst,\M^\ss):(\M^\rst,\M^\ss)\in\Im\Phi\bigr\}$, a set of subsets of $\Kah(X)$. Then $\Kah(X)=\coprod_{S\in T}S$. Define $\{S_i:i\in I\}$ to be the set of connected components of elements of $T$. Then $\Kah(X)=\coprod_{i\in I}S_i$, and each $S_i$ is connected with $\Phi\vert_{S_i}$ constant, and is maximal under these conditions.

Let $\om_0\in\Kah(X)$, let $U$ be a small open neighbourhood of $\om_0$ in $\Kah(X)$, and let $K=\ov U$ be the compact closure of $U$ in $\Kah(X)$. Write $\sF_{\al,K}$ for the family of $\mu^\om$-semi\-stable sheaves in class $\al$ for any $\om\in K$. Then if $F\in\sF_{\al,K}$ is $\mu^\om$-semistable we have $\mu_{\max}^\om(F)=\bar\mu^\om(\al)\le\max_{\om'\in K}\bar\mu^{\om'}(\al)=:\mu_0$, so Condition \ref{co7cond1} for $K,\al$ implies that $\sF_{\al,K}$ is bounded. 

Write $\sG_{\al,K}$ for the family of $G\in\coh(X)$ which are pure of dimension $d=\dim\al$ such that $G$ is a quotient $F\twoheadrightarrow G$ for some $F\in\sF_{\al,K}$ with $0\ne G\ne F$ and $\bar\mu^\om(\lb G\rb)\le\bar\mu^\om(\lb F\rb)=\bar\mu^\om(\al)$ for some $\om\in K$. Then Proposition \ref{co7prop3} applies to $\sF_{\al,K},\sG_{\al,K}$, so Proposition \ref{co7prop3}(a) says the set of Chern characters $\ch_d(G)\in H_{2d}(X,\Q)$ for $G\in\sG_{\al,K}$ is finite.

Using the facts that $K$ is compact and there are only finitely many values for $\ch_d(G)$ for $G\in\sG_{\al,K}$, we can show that there exists a finite subset $R\subset\Kah(X)$ and $\ti\mu_0\gg 0$ such that if $G\in\sG_{\al,K}$, so we can write $F\twoheadrightarrow G$ for $F\in\sF_{\al,K}$ with $\bar\mu^\om(\lb G\rb)\le\bar\mu^\om(\lb F\rb)=\bar\mu^\om(\al)$ for some $\om\in K$, then $\bar\mu^{\ti\om}(\lb G\rb)\le\ti\mu_0$ for some $\ti\om\in R$. Proposition \ref{co7prop3} for each $\ti\om\in R$ then implies that $\sG_{\al,K}$ is bounded. Thus the set $W_{\al,K}$ of classes $\be=\lb G\rb$ of sheaves $G$ in $\sG_{\al,K}$ is finite.

Now if $\om\in U\subset K$ and $F\in\sF_{\al,K}$ is not $\mu^\om$-semistable it has a pure maximal destabilizing quotient object $F\twoheadrightarrow G$ with $0\ne G\ne F$ and $\bar\mu^\om(\lb G\rb)<\bar\mu^\om(\lb F\rb)$, so $\lb G\rb\in W_{\al,K}$. If $\om\in K$ and $F\in\sF_{\al,K}$ is $\mu^\om$-semistable but not $\mu^\om$-stable it has a pure semistabilizing quotient object $F\twoheadrightarrow G$ with $0\ne G\ne F$ and $\bar\mu^\om(\lb G\rb)=\bar\mu^\om(\lb F\rb)$, so again $\lb G\rb\in W_{\al,K}$. Thus for each $\om\in U$ we see that
\ea
\begin{split}
\M_\al^\rst(\mu^\om)=\,&\bigl\{[F]:\text{$F\in\sF_{\al,K}$, there does not exist $F\twoheadrightarrow G$ with $G$} \\ &\;\text{pure dim $d$, $\lb G\rb=\be\in W_{\al,K}$ and $\bar\mu^\om(\be)\le\bar\mu^\om(\al)$}\bigr\},
\end{split}
\label{co7eq29}\\
\begin{split}
\M_\al^\ss(\mu^\om)=\,&\bigl\{[E]:\text{$F\in\sF_{\al,K}$, there does not exist $F\twoheadrightarrow G$ with $G$} \\ &\;\text{pure dim $d$, $\lb G\rb=\be\in W_{\al,K}$ and $\bar\mu^\om(\be)<\bar\mu^\om(\al)$}\bigr\},
\end{split}
\label{co7eq30}
\ea
where both of these are open substacks of $\bigcup_{\om'\in K}\M_\al^\ss(\mu^{\om'})\subseteq\M_\al^\pl$, which is the moduli stack of sheaves in $\sF_{\al,K}$.

From \eq{co7eq29}--\eq{co7eq30} we see that after intersecting with $U\subset\Kah(X)$, the sets in $T$ depend only on whether $\bar\mu^\om(\be)<\bar\mu^\om(\al-\be)$, or $\bar\mu^\om(\be)=\bar\mu^\om(\al-\be)$, or $\bar\mu^\om(\be)>\bar\mu^\om(\al-\be)$, for a finite set $W_{\al,K}$ of $\be$ in $C(\coh(X))$ with $\al-\be$ in $C(\coh(X))$. Thus the sets in $T$ are locally closed, as in Theorem \ref{co7thm2}(a), and satisfy (b),(c). Choosing $U$ to be a sufficiently small open ball about $\om_0$ in $\Kah(X)$, we can make the intersections of the sets of $T$ with $U$ connected. Then the intersections of sets of $T$ with $U$ coincide with intersections of sets $S_i$ with $U$. Hence the sets $S_i$ also satisfy (a)--(c). Parts (d),(e) are easy consequences of \eq{co7eq29}--\eq{co7eq30}. This proves the first part of Theorem~\ref{co7thm2}.

For the second part, the argument above with $\tau^\om$ in place of $\mu^\om$ shows that
\begin{align*}
\M_\al^\rst(\tau^\om)=\,&\bigl\{[F]:\text{$F\in\sF_{\al,K}$, there does not exist $F\twoheadrightarrow G$ with $G$} \\ &\;\text{pure dim $d$, $\lb G\rb=\be\in W_{\al,K}$ and $\tau^\om(\be)\le\tau^\om(\al)$}\bigr\},\\
\M_\al^\ss(\tau^\om)=\,&\bigl\{[E]:\text{$F\in\sF_{\al,K}$, there does not exist $F\twoheadrightarrow G$ with $G$} \\ &\;\text{pure dim $d$, $\lb G\rb=\be\in W_{\al,K}$ and $\tau^\om(\be)<\tau^\om(\al)$}\bigr\}.
\end{align*}
For fixed $\al,\be$, the inequalities $\tau^\om(\be)\le\tau^\om(\al)$ and $\tau^\om(\be)<\tau^\om(\al)$ are locally closed conditions on $\om$, which change over smooth algebraic hypersurfaces in $\Kah(X)\ni\om$. So the previous argument shows that $\tau^\om$-(semi)stability has a wall-and-chamber structure satisfying Theorem \ref{co7thm2}(a)--(c), as we have to prove.

The reason (d)--(e) do not extend to $\tau^\om$-(semi)stability is that in contrast to $\mu$-(semi)stability, $\tau^\om(\be)\le\tau^\om(\al)$ need not be a closed condition on $\om$, and $\tau^\om(\be)<\tau^\om(\al)$ need not be an open condition on $\om$, since $(\tau^\om,M,\le)$ need not depend continuously on $\om\in\Kah(X)$. See Remark \ref{co7rem3} for an example.

\subsubsection{Proof of Theorem \ref{co7thm4}}
\label{co7116}

Let $R$ be a discrete valuation ring over $\C$, and $\ga$ be the generic point and $\ka$ the special point in $\Spec R$. For $\om\in\Kah_\Z(X)$, the proof of \cite[Th.~2.B.1]{HuLe2} starts with a coherent sheaf $F\ra X\t\Spec R$, flat over $\Spec R$ such that $F\vert_{X\t \ga}$ is $\tau^\om$-semistable (or $\mu^\om$-semistable). It then constructs a subsheaf $F'\subseteq F$ such that $F'\vert_{X\t \ga}=F\vert_{X\t \ga}$ and $F'\vert_{X\t\ka}$ is $\tau^\om$-semistable (or $\mu^\om$-semistable), and this implies the valuative criterion for universal closedness for $(\Pi_\al^\pl)^{-1}(\M_\al^\ss(\tau^\om))$ (or $(\Pi_\al^\pl)^{-1}(\M_\al^\ss(\mu^\om))$, respectively) with respect to~$R$.

The proof works by inductively constructing a sequence of subsheaves $F=F^0\supseteq F^1\supseteq\cdots$ on $X\t\Spec R$ as follows: if $F^n\vert_{X\t\ka}$ is $\tau^\om$- (or $\mu^\om$-\nobreak)semistable then we set $F^n=F^{n+1}=\cdots$ and $F'=F^n$, and otherwise we define $0\ne B^n\subsetneq F^n\vert_{X\t\ka}$ to be the maximal $\tau^\om$- (or $\mu^\om$-\nobreak)\-destabilizing subsheaf of $F^n\vert_{X\t\ka}$, and set $G^n=F^n\vert_{X\t\ka}/B^n$, and define $F^{n+1}\subsetneq F^n$ to be the kernel of the composition $F^n\ra F^n\vert_{X\t\ka}\ra G^n$. There are then exact sequences in $\coh(X)$:
\begin{equation*}
\xymatrix@C=13pt{ 0 \ar[r] & B^n \ar[r] & 	F^n\vert_{X\t\ka} \ar[r] & G^n \ar[r] & 0, & 0 \ar[r] & G^n \ar[r] & F^{n+1}\vert_{X\t\ka} \ar[r] & B^n \ar[r] & 0. }
\end{equation*}
We must prove $F^n$ is independent of $n$ for $n\gg 0$. Huybrechts and Lehn do this by showing that $\tau^\om(B^0)\ge \tau^\om(B^1)\ge\cdots\ge \tau^\om(F\vert_{X\t\ka})$, and noting that as $\om\in\Kah_\Z(X)$, the coefficients of $\tau^\om(B^n)$ lie in $\Z/N\subset\Q$ for some $N\gg 0$, so the $\tau^\om(B^n)$ must be independent of $n$ for~$n\gg 0$.

When $\om\in\Kah(X)$, the coefficients of $\tau^\om(B^n)$ lie in $\R$, so we cannot use this argument. Instead we first note that using \cite[Ex.~2.B.2]{HuLe2} (the analogue of \cite[Th.~2.B.1]{HuLe2} for purity) we can take all sheaves $F,F^n,B^n,G^n$ above to be pure of dimension $d=\dim\al$. We have $\bar\mu^\om(B^n)=\mu_{\max}^\om(F^n\vert_{X\t\ka})$ by definition of $B^n$, and $\bar\mu^\om(B^0)\ge\bar\mu^\om(B^1)\ge\cdots$ as $\tau^\om(B^0)\ge \tau^\om(B^1)\ge\cdots$ and the $B^n$ are pure of dimension $d$. Also~$\bar\mu^\om(G^n)\le \bar\mu^\om(F^n\vert_{X\t\ka})=\bar\mu^\om(\al)$.

Hence $\mu_{\max}^\om(F^n\vert_{X\t\ka})\le \bar\mu^\om(B^0)$ for all $n\ge 0$, and $\lb F^n\vert_{X\t\ka}\rb=\al$. As Condition \ref{co7cond1} holds for $K=\{\om\}$ and $\al$, we see that the $F^n$ for $n\ge 0$ form a bounded family. Since we have quotients $F^n\twoheadrightarrow G^n$ with $F^n,G^n$ pure of dimension $d$ and $\bar\mu^\om(G^n)\le\bar\mu^\om(\al)$, Proposition \ref{co7prop3} shows that the $G^n$ for $n\ge 0$ also form a bounded family. Therefore $\lb G^n\rb$ and $\lb B^n\rb=\al-\lb G^n\rb$ can realize only finitely many possibilities in $K(\coh(X))$. Since $\tau^\om(B^0)\ge \tau^\om(B^1)\ge\cdots$ we see that the $\tau^\om(B^n)$ must be independent of $n$ for $n\gg 0$. This fills the gap in the proof of \cite[Th.~2.B.1]{HuLe2} when $\om$ lies in $\Kah(X)$ rather than~$\Kah_\Z(X)$.

\subsection[Verifying most of Assumptions \ref{co4ass1}, \ref{co5ass1}--\ref{co5ass3} for coherent sheaves]{Verifying most of Assumptions \ref{co4ass1} and \ref{co5ass1}--\ref{co5ass3} for \\ coherent sheaves}
\label{co72}

Our goal is to apply the results of Chapter \ref{co5} to abelian categories $\coh(X)$ of coherent sheaves on a smooth, connected, projective $\C$-scheme $X$. This requires restrictions on $X$, the main cases being when $X$ is a curve, a surface, a Fano 3-fold, and (in future work) a Calabi--Yau 4-fold.

To do this, we have to give the data and verify the conditions of Assumptions \ref{co4ass1} and \ref{co5ass1}--\ref{co5ass3} for $\A=\B=\coh(X)$. This section covers the easier parts of Assumptions \ref{co4ass1} and \ref{co5ass1}--\ref{co5ass3}. This is everything except:
\begin{itemize}
\setlength{\itemsep}{0pt}
\setlength{\parsep}{0pt}
\item[(i)] Assumption \ref{co5ass2}(g),(h) on properness of stable=semistable moduli stacks;
\item[(ii)] Assumptions \ref{co5ass1}(f) and \ref{co5ass2}(e) on the quasi-smooth derived stacks $\bs\dM_\al^\red,\ab\bs\dM_\al^\rpl$ and the associated obstruction theories \eq{co5eq5}; and
\item[(iii)] The finiteness condition Assumption~\ref{co5ass3}(b).
\end{itemize}
We will treat (i)--(iii) in \S\ref{co73}--\S\ref{co75} respectively.

In \S\ref{co724} we extend to the case when an algebraic group $G$ acts on $X$, and verify Assumptions \ref{co4ass2} and \ref{co5ass4}(a),(d) for $G$-equivariant homology.

For parts of our treatment we regard the moduli stack $\M$ of objects in $\coh(X)$ as an open substack in the (higher) moduli stack $\baM$ of objects in the derived category $D^b\coh(X)$, via the inclusion $\coh(X)\hookra D^b\coh(X)$. We will define data such as $\cE^\bu$ in Assumption \ref{co4ass1}(f) first on $\baM\t\baM$, and then restrict to the open substack $\M\t\M\subset\baM\t\baM$. This is to save work in future extensions of the theory to $D^b\coh(X)$, as in \S\ref{co55}, or in applications of Theorems \ref{co5thm1}--\ref{co5thm3} to abelian subcategories $\A\subset D^b\coh(X)$ other than $\A=\coh(X)$. It also allows us to use ideas on mapping stacks, since $\baM\cong\Map_{\HSta_\C}(X,\Perf_\C)$, but $\M$ does not have an easy mapping stack interpretation.

\begin{rem}
\label{co7rem5}
{\bf(a)} For parts of Assumptions \ref{co5ass2}(b),(g),(h) we will need Condition \ref{co7cond1} to hold for $K=\{\om\}$ with $\om\in\Kah(X)$ and $\al\in C(\coh(X))$, but as in Propositions \ref{co7prop5}--\ref{co7prop6} and Question \ref{co7quest1}(b) the author can only prove this when $\dim\al=0,1$ or $m$. The author can also only prove Assumption \ref{co5ass3}(b) when $\dim\al=0,1$ or $m$. If $m=\dim X\le 2$ this is all $\al$, but otherwise it excludes sheaves of dimensions $2,3,\ldots,m-1$. We care about $m=3$ for Fano and Calabi--Yau 3-folds, and $m=4$ for Calabi--Yau 4-folds.

We deal with this in \eq{co7eq46}--\eq{co7eq47} by requiring $C(\B)_\pe$ in Assumption \ref{co5ass1}(e) to contain only $\al\in C(\B)$ with $\dim\al=0,1$ or $m$, so that we just do not define invariants counting torsion sheaves of dimensions $2,\ldots,m-1$. If $X$ is a Fano 3-fold we also exclude $\dim\al=0$, as the quasi-smooth stacks and obstruction theories in Assumption \ref{co5ass1}(f) do not work in this case. The exclusion from our theory of sheaves of dimension 2 on Fano and Calabi--Yau 3-folds, and sheaves of dimensions 2 and 3 on Calabi--Yau 4-folds, is rather unsatisfactory.
\smallskip

\noindent{\bf(b)} Part {\bf(a)} assumes our goal is to define invariants $[\M_\al^\ss(\tau^\om)]_\inv,[\M_\al^\ss(\mu^\om)]_\inv$ for $\om\in\Kah(X)$, for which the wall-crossing formulae \eq{co5eq31}--\eq{co5eq34} in Theorems \ref{co5thm2}--\ref{co5thm3} hold for transformations between two different~$\om,\ti\om\in\Kah(X)$.

If we are willing to fix a single $\om\in\Kah_\Q(X)$ then Theorem \ref{co5thm1} still works to define invariants $[\M_\al^\ss(\tau^\om)]_\inv$ when $\dim\al=2,3,\ldots,m-1$, and \eq{co5eq31}--\eq{co5eq34} hold for transforming between $\tau^\om$ and $\mu^\om$ for such $\al$. But we cannot currently prove \eq{co5eq31}--\eq{co5eq34} for transforming between two different $\om,\ti\om\in\Kah_\Q(X)$.

\end{rem}

\subsubsection{Assumption \ref{co4ass1} for $\coh(X)$}
\label{co721}

We define the data and verify Assumption \ref{co4ass1} for $\coh(X)$:

\begin{dfn}
\label{co7def9}
Work in the situation of Definition \ref{co7def1}, and write $\B=\coh(X)$. For Assumption \ref{co4ass1}(a), the moduli stack $\M$ of objects in $\B$ is described in Definition~\ref{co7def2}.

For Assumption \ref{co4ass1}(b), we have a perfect complex on $X\t\baM\t\baM$
\e
\Pi_{12}^*(\cU^\bu)\op \Pi_{13}^*(\cU^\bu)\longra X\t\baM\t\baM.
\label{co7eq31}
\e
Here and below, for a product of stacks $S_1\t\cdots\t S_k$ we write $\Pi_i$ for the projection to the $i^{\rm th}$ factor, and $\Pi_{ij}$ for the projection to the product of the $i^{\rm th}$ and $j^{\rm th}$ factors, and so on.

As in Definition \ref{co7def2}, equation \eq{co7eq31} determines a morphism $u':X\t\baM\t\baM\ra\Perf_\C$, which determines a morphism $\bar\Phi:\baM\t\baM\ra\Map(X,\Perf_\C)\cong\baM$. As \eq{co7eq31} is isomorphic to $u^{\prime *}(\cU_{\Perf_\C}^\bu)$ we have a canonical isomorphism
\e
(\id_X\t\bar\Phi)^*(\cU^\bu)\cong \Pi_{12}^*(\cU^\bu)\op \Pi_{13}^*(\cU^\bu),
\label{co7eq32}
\e
which characterizes $\bar\Phi$. As $\coh(X)$ is closed under direct sums in $D^b\coh(X)$ we see that $\bar\Phi$ maps $\M\t\M\ra\M$, and we define $\Phi=\bar\Phi\vert_{\M\t\M}:\M\t\M\ra\M$.

Since $\cU^\bu\vert_{X\t\{[F^\bu]\}}\cong F^\bu$ we see from \eq{co7eq32} that $\bar\Phi$ acts on $\C$-points by $\bar\Phi_*:([E^\bu],[F^\bu])\mapsto[E^\bu\op F^\bu]$, so $\Phi$ act on $\C$-points by $([E],[F])\mapsto[E\op F]$. The action on isotropy groups in Assumption \ref{co4ass1}(b) is clear. 

Let $\si:\baM\t\baM\ra\baM\t\baM$ exchange the factors. Since
\e
\begin{split}
&\Pi_{12}^*(\cU^\bu)\op \Pi_{13}^*(\cU^\bu)\cong \Pi_{13}^*(\cU^\bu)\op \Pi_{12}^*(\cU^\bu)\\
&\quad\cong (\id_X\t\si)^*(\Pi_{12}^*(\cU^\bu)\op \Pi_{13}^*(\cU^\bu)),
\end{split}
\label{co7eq33}
\e
and \eq{co7eq32} characterizes $\bar\Phi$, we see there is a 2-morphism $\bar\Phi\cong\bar\Phi\ci\si$. That is, $\bar\Phi$ is commutative in $\Ho(\HSta_\C)$, so $\Phi$ is commutative in $\Ho(\Art_\C)$. Associativity of $\Phi$ follows by a similar argument on $X\t\baM\t\baM\t\baM$. This proves Assumption~\ref{co4ass1}(b).

For Assumption \ref{co4ass1}(c), we have a perfect complex on $X\t[*/\bG_m]\t\baM$
\e
\Pi_2^*(L_{[*/\bG_m]})\ot\Pi_{13}^*(\cU^\bu)\longra X\t[*/\bG_m]\t\baM.
\label{co7eq34}
\e
As above, this determines a morphism $u'':X\t[*/\bG_m]\t\baM\ra\Perf_\C$, which determines a morphism $\bar\Psi:[*/\bG_m]\t\baM\ra\Map(X,\Perf_\C)\cong\baM$. As \eq{co7eq34} is isomorphic to $u^{\prime\prime *}(\cU_{\Perf_\C}^\bu)$ we have a canonical isomorphism
\e
(\id_X\t\bar\Psi)^*(\cU^\bu)\cong \Pi_2^*(L_{[*/\bG_m]})\ot \Pi_{13}^*(\cU^\bu),
\label{co7eq35}
\e
which characterizes $\bar\Psi$. Since $\coh(X)$ is closed under tensor product with 1-dimensional complex vector spaces in $D^b\coh(X)$, we see that $\bar\Psi$ maps $[*/\bG_m]\t\M\ra\M$, and we define $\Psi=\bar\Psi\vert_{[*/\bG_m]\t\M}:[*/\bG_m]\t\M\ra\M$.

On $\C$-points $\bar\Psi$ acts by $\bar\Psi_*:(*,[F^\bu])\mapsto[\C\ot F^\bu]=[F^\bu]$, so $\Psi$ acts by $(*,[F])\mapsto[F]$. The action on isotropy groups in Assumption \ref{co4ass1}(c) is clear. 

Equation \eq{co4eq2} may be deduced from the isomorphisms
\begin{align*}
&\bigl(\Pi_2^*(L_{[*/\bG_m]})\ot\Pi_{13}^*(\cU^\bu)\bigr)\op\bigl(\Pi_2^*(L_{[*/\bG_m]})\ot\Pi_{14}^*(\cU^\bu)\bigr)\\
&\quad\cong \Pi_2^*(L_{[*/\bG_m]})\ot\bigl(\Pi_{13}^*(\cU^\bu)\op\Pi_{14}^*(\cU^\bu)\bigr)\quad\text{on $X\t[*/\bG_m]\t\M\t\M$,}\\
&\Pi_2^*(L_{[*/\bG_m]})\ot\bigl(\Pi_3^*(L_{[*/\bG_m]})\ot\Pi_{14}^*(\cU^\bu)\bigr)\\
&\quad\cong\Pi_{23}^*(\Om^*(L_{[*/\bG_m]}))\ot\Pi_{14}^*(\cU^\bu)\qquad\qquad\;\text{on $X\t[*/\bG_m]\t[*/\bG_m]\t\M$,}
\end{align*}
using the method after \eq{co7eq33}. This proves Assumption~\ref{co4ass1}(c).

For Assumption \ref{co4ass1}(d), we define $K(\coh(X))=K^\num(\coh(X))$ as in Definition \ref{co7def1}. We have a decomposition $\M=\coprod_{\al\in K(\coh(X))}\M_\al$ as in Definition~\ref{co7def2}.

If $0\ne E\in\coh(X)$ then as $X$ is a smooth projective $\C$-scheme the Hilbert polynomial of $E$ (with respect to some ample line bundle on $X$) is nonzero. But this depends on the Chern character $\lb E\rb$ of $E$ by the Grothendieck--Riemann--Roch Theorem, so $\lb E\rb\ne 0$, as required.

For Assumption \ref{co4ass1}(e), we define $\chi:K(\B)\t K(\B)\ra\Z$ as in Definition \ref{co7def1}.

For Assumption \ref{co4ass1}(f), define a perfect complex $\cExt^\bu$ on $\baM\t\baM$, the {\it Ext complex}, by
\e
\cExt^\bu=(\Pi_{23})_*\bigl[\Pi_{12}^*(\cU^\bu)^\vee\ot\Pi_{13}^*(\cU^\bu)\bigr],
\label{co7eq36}
\e
using derived pushforward, pullback and tensor product functors on $X\t\baM\t\baM$ as in Huybrechts \cite{Huyb}. Then for $F^\bu,G^\bu$ in $D^b\coh(X)$ and $k\in\Z$ we have
\begin{equation*}
H^k\bigl(\cExt^\bu\vert_{([F^\bu],[G^\bu])}\bigr)\cong \Ext^k(F^\bu,G^\bu):=\Hom_{D^b\coh(X)}\bigl(F^\bu,G^\bu[k]\bigr).
\end{equation*}
Therefore \eq{co7eq1} implies that for $\al,\be\in K(\B)$ we have
\e
\rank\bigl(\cExt^\bu\vert_{\baM_\al\t\baM_\be}\bigr)=\chi(\al,\be).
\label{co7eq37}
\e
Define $\cE^\bu,\cE^\bu_{\al,\be}$ to be the dual perfect complexes
\e
\cE^\bu=(\cExt^\bu)^\vee,\qquad \cE_{\al,\be}^\bu=(\cExt^\bu)^\vee\vert_{\baM_\al\t\baM_\be},
\label{co7eq38}
\e
so that $\rank\cE^\bu_{\al,\be}=\chi(\al,\be)$ by \eq{co7eq37}. We use the same notation $\cE^\bu$ for the restriction of $\cE^\bu$ to the open substack $\M\t\M\subset\baM\t\baM$. Equations \eq{co4eq3}--\eq{co4eq4} follow easily from \eq{co7eq32}, \eq{co7eq36} and \eq{co7eq38}, and \eq{co4eq5}--\eq{co4eq6} follow from \eq{co7eq35}, \eq{co7eq36} and~\eq{co7eq38}.

As $\chi$ is nondegenerate on $K(\coh(X))=K^\num(\coh(X))$ by definition of $K^\num(\coh(X))$, Assumption \ref{co4ass1}(g) holds automatically.
\end{dfn}

\subsubsection{Assumption \ref{co5ass1}(a)--(e),(g) for $\coh(X)$}
\label{co722}

We define the data and verify Assumption \ref{co5ass1}(a)--(e),(g) for $\A=\B=\coh(X)$:

\begin{dfn}
\label{co7def10}
Work in the situation of Definitions \ref{co7def1} and \ref{co7def9}. 

Assumption \ref{co5ass1}(a) is trivial as $\A=\B$. For Assumption \ref{co5ass1}(b) we take the quotient $K_0(\A)\twoheadrightarrow K(\A)$ to be as for Assumption \ref{co4ass1}(d) in Definition \ref{co7def9}. Assumption \ref{co5ass1}(c) was proved in Definition~\ref{co7def9}. 

For Assumption \ref{co5ass1}(d), as in Definition \ref{co7def2} the derived moduli stack $\bs\M$ of objects in $\coh(X)$ exists by To\"en--Vaqui\'e \cite{ToVa}, and is a derived Artin $\C$-stack, with classical truncation $\M=t_0(\bs\M)$. It is an open substack of the derived moduli stack $\bs\baM$ of objects in $D^b\coh(X)$, which also exists by \cite{ToVa}. There is a universal coherent sheaf $\bs\cU\ra X\t\bs\M$ and a universal perfect complex $\bs\cU^\bu\ra X\t\bs\baM$. This determines a morphism $\bs u:X\t\bs\baM\ra\bs\Perf_\C$ in $\DSta_\C$, inducing an isomorphism~$\bs\baM\ra\Map_{\DSta_\C}(X,\bs\Perf_\C)$.

As in Assumption \ref{co5ass1}(d)(i), the splittings $\M=\coprod_{\al\in K(\coh(X))}\M_\al$, $\baM=\coprod_{\al\in K(\coh(X))}\baM_\al$ lift to $\bs\M=\coprod_{\al\in K(\coh(X))}\bs\M_\al$, $\bs\baM=\coprod_{\al\in K(\coh(X))}\bs\baM_\al$.

For Assumption \ref{co5ass1}(d)(ii), in Definition \ref{co7def9} we constructed $\bar\Phi:\baM\t\baM\ra\baM$ to satisfy \eq{co7eq32}, using the universal property of $\cU^\bu\ra X\t\baM$, and set $\Phi=\bar\Phi\vert_{\M\t\M}:\M\t\M\ra\M$. The same construction works for derived stacks using the universal property of $\bs\cU^\bu\ra X\t\bs\baM$, giving a morphism $\bs{\bar\Phi}:\bs\baM\t\bs\baM\ra\bs\baM$ with an isomorphism as for \eq{co7eq32}
\begin{equation*}
(\id_X\t\bs{\bar\Phi})^*(\bs\cU^\bu)\cong \Pi_{12}^*(\bs\cU^\bu)\op \Pi_{13}^*(\bs\cU^\bu),
\end{equation*}
and we may set $\bs\Phi=\bs{\bar\Phi}\vert_{\bs\M\t\bs\M}:\bs\M\t\bs\M\ra\bs\M$. Under the map $\id\t i:X\t\baM\ra X\t\bs\baM$ we have $(\id\t i)^*(\bs\cU^\bu)=\cU^\bu$, so the universal properties of $\cU^\bu,\bs\cU^\bu$ imply that $\bs{\bar\Phi}\ci(i\t i)=i\ci\bar\Phi$. Thus taking classical truncations shows that $t_0(\bs{\bar\Phi})=\bar\Phi$, so~$t_0(\bs\Phi)=\Phi$.

Similarly, the construction of $\bar\Psi:[*/\bG_m]\t\baM\ra\baM$ and $\Psi=\bar\Psi\vert_{[*/\bG_m]\t\M}:[*/\bG_m]\t\M\ra\M$ in Definition \ref{co7def9} lift to $\bs{\bar\Psi}:[*/\bG_m]\t\bs\baM\ra\bs\baM$ and $\bs\Psi=\bs{\bar\Psi}\vert_{[*/\bG_m]\t\bs\M}:[*/\bG_m]\t\bs\M\ra\bs\M$ with $t_0(\bs{\bar\Psi})=\bar\Psi$ and $t_0(\bs\Psi)=\Psi$. The proofs of the identities satisfied by $\bar\Phi,\bar\Psi$ and $\Phi,\Psi$ lift to $\bs{\bar\Phi},\bs{\bar\Psi}$ and~$\bs\Phi,\bs\Psi$.

Thus $\bs{\bar\Psi},\bs\Psi$ are actions of the group stack $[*/\bG_m]$ on the derived stacks $\bs\baM,\bs\M$, which are free except over $0$. So we can take the quotients $\bs\baM^\pl,\bs\M^\pl$, with projections $\bs{\bar\Pi}{}^\pl:\bs\baM\ra\bs\baM^\pl$, $\bs\Pi^\pl:\bs\M\ra\bs\M^\pl$ which are principal $[*/\bG_m]$-bundles except over 0, with $t_0(\bs\baM^\pl)=\baM^\pl$ and~$t_0(\bs\M^\pl)=\M^\pl$.

For Assumption \ref{co5ass1}(d)(iii), To\"en--Vaqui\'e \cite[Cor.~3.31]{ToVa} implies that $\bs\baM$ and hence $\bs\M\subset\bs\baM$ are locally finitely presented.

For Assumption \ref{co5ass1}(d)(iv), as above we have a morphism $\bs u:X\t\bs\baM\ra\bs\Perf_\C$ inducing an isomorphism $\bs\baM\ra\Map_{\DSta_\C}(X,\bs\Perf_\C)$. As in Pantev--To\"en--Vaqui\'e--Vezzosi \cite{PTVV}, since $\bs\baM$ is a mapping stack, we can give an explicit presentation for its tangent complex:
\e
\bT_{\bs\baM}\cong (\bs\Pi_{\bs\baM})_*\ci\bs u^*(\bT_{\bs\Perf_\C}).
\label{co7eq39}
\e
But as in \cite[\S 2.3]{PTVV}, the tangent complex of $\bs\Perf_\C$ is
\e
\bT_{\bs\Perf_\C}\cong \End(\bs\cU_{\bs\Perf_\C}^\bu)[1]=(\bs\cU_{\bs\Perf_\C}^\bu)^\vee\ot \bs\cU_{\bs\Perf_\C}^\bu[1].
\label{co7eq40}
\e

Since $\bs\cU^\bu\cong\bs u^*(\bs\cU_{\bs\Perf_\C}^\bu)$, combining \eq{co7eq39}--\eq{co7eq40} yields
\begin{equation*}
\bT_{\bs\baM}\cong (\bs\Pi_{\bs\baM})_*\bigl(\End(\bs\cU^\bu)\bigr)[1]=(\bs\Pi_{\bs\baM})_*\bigl((\bs\cU^\bu)^\vee\ot \bs\cU^\bu\bigr)[1],
\end{equation*}
so pulling back by $\bar\imath:\baM\hookra\bs\baM$ gives
\e
\bar\imath^*(\bT_{\bs\baM})\cong (\Pi_{\baM})_*\bigl(\End(\cU^\bu)\bigr)[1]=(\Pi_{\baM})_*\bigl((\cU^\bu)^\vee\ot \cU^\bu\bigr)[1].
\label{co7eq41}
\e
Comparing this with the definition \eq{co7eq36} of $\cExt^\bu$, we see that
\e
\bar\imath^*(\bT_{\bs\baM})\cong \De_\baM^*(\cExt^\bu)[1],
\label{co7eq42}
\e
where $\De_\baM:\baM\ra\baM\t\baM$ is the diagonal map. Noting that $\bL_{\bs\baM}\cong\bT_{\bs\baM}^\vee$, taking duals in $D_\qcoh(\baM)$ and using $\cE^\bu=(\cExt^\bu)^\vee$ in \eq{co7eq38} yields $\bar\imath^*(\bL_{\bs\baM})\cong\De_\baM^*(\cE^\bu)[-1]$. Thus restricting to $\M_\al\subset\baM$, as in \eq{co5eq1} we have an isomorphism
\e
\th_\al:\De_{\M_\al}^*(\cE_{\al,\al}^\bu)[-1]\,{\buildrel\cong\over\longra}\,i_\al^*(\bL_{\bs\M_\al}).
\label{co7eq43}
\e

For Assumption \ref{co5ass1}(d)(v), from \eq{co7eq41} we deduce that
\ea
i_{\al,\be}^*(\bT_{\bs\M_\al\t\bs\M_\be})&\cong (\Pi_{\M_\al\t\M_\be})_*\bigl(\End(\Pi_{X\t\M_\al}^*(\cU_\al^\bu))\op\End(\Pi_{X\t\M_\be}^*(\cU_\be^\bu))\bigr)[1],
\nonumber\\
i_{\al,\be}^*\ci\bs\Phi_{\al,\be}^*&(\bT_{\bs\M_{\al+\be}})\cong \Phi_{\al,\be}^*\ci(\Pi_{\M_{\al+\be}})_*\bigl(\End(\cU_{\al+\be}^\bu)\bigr)[1]
\nonumber\\
&\cong (\Pi_{\M_\al\t\M_\be})_*\ci(\id_X\t\Phi_{\al,\be}^*)(\End(\cU_{\al+\be}^\bu)\bigr)[1]
\nonumber\\
&\cong (\Pi_{\M_\al\t\M_\be})_*\End\bigl(\Pi_{X\t\M_\al}^*(\cU_\al^\bu)\op\Pi_{X\t\M_\be}^*(\cU_\be^\bu) \bigr)[1].
\label{co7eq44}
\ea
Now $i_{\al,\be}^*(\bT_{\bs\Phi_{\al,\be}}):i_{\al,\be}^*\ci\bs\Phi_{\al,\be}^*(\bT_{\bs\M_{\al+\be}})\ra i_{\al,\be}^*(\bT_{\bs\M_\al\t\bs\M_\be})$ is induced, under the identifications \eq{co7eq44}, by the inclusion 
\ea
&\End(\Pi_{X\t\M_\al}^*(\cU_\al^\bu))\op\End(\Pi_{X\t\M_\be}^*(\cU_\be^\bu))
\nonumber\\
&\hookra\End\bigl(\Pi_{X\t\M_\al}^*(\cU_\al^\bu)\op\Pi_{X\t\M_\be}^*(\cU_\be^\bu) \bigr)
\nonumber\\
&\cong \End(\Pi_{X\t\M_\al}^*(\cU_\al^\bu))\op\End(\Pi_{X\t\M_\be}^*(\cU_\be^\bu))\op
\nonumber\\
&\Pi_{X\t\M_\al}^*(\cU_\al^\bu)^\vee\ot\Pi_{X\t\M_\be}^*(\cU_\be^\bu)\op
\Pi_{X\t\M_\be}^*(\cU_\be^\bu)^\vee\ot\Pi_{X\t\M_\al}^*(\cU_\al^\bu).
\label{co7eq45}
\ea
So the dual morphism $i_{\al,\be}^*(\bL_{\bs\Phi_{\al,\be}})$, which is the right hand morphism in \eq{co5eq3}, is induced by the dual projection to \eq{co7eq45} under the dual identifications \eq{co7eq44}. From this and the definitions \eq{co7eq36}, \eq{co7eq38} and \eq{co7eq43} of $\cE^\bu_{\ga,\de}$ and $\th_\ga$, we see that \eq{co5eq3} commutes. This completes Assumption~\ref{co5ass1}(d).

For Assumption \ref{co5ass1}(e), as in Remark \ref{co7rem5}, if $X$ is not a Fano 3-fold we define
\e
C(\B)_\pe=\bigl\{\al\in C(\B):\text{$\dim\al=0,1$ or $m$}\bigr\},
\label{co7eq46}
\e
and if $X$ is a Fano 3-fold, so that $m=3$, we define
\e
C(\B)_\pe=\bigl\{\al\in C(\B):\text{$\dim\al=1$ or 3}\bigr\},
\label{co7eq47}
\e
where $\dim\al$ is as in \eq{co7eq2}.

For Assumption \ref{co5ass1}(g), fix a very ample line bundle $\O_X(1)$ on the smooth projective $\C$-scheme $X$, and as usual write $\O_X(n)=\O_X(1)^{\ot^n}$ for $n\in\Z$. For $E\in\coh(X)$ and $n\in\Z$ we write $E(n)=E\ot\O_X(n)$. Following Huybrechts and Lehn \cite[\S 1.7]{HuLe2}, we say that $E\in\coh(X)$ is $n$-{\it regular\/} for $n\in\Z$ if $H^i(E(n-i))=0$ for all $i>0$. Write $\coh_{\text{$n$-reg}}(X)$ for the full subcategory of $n$-regular objects in $\coh(X)$. Here are some elementary properties of $n$-regularity:
\begin{itemize}
\setlength{\itemsep}{0pt}
\setlength{\parsep}{0pt}
\item[(a)] By \cite[Lem.~1.7.2]{HuLe2}, if $E$ is $n$-regular then $E$ is $m$-regular for any $m\ge n$. Hence $\cdots\subset\coh_{\text{$n$-reg}}(X)\subset \coh_{\text{$(n\!+\!1)$-reg}}(X)\subset\cdots.$
\item[(b)] $\coh_{\text{$n$-reg}}(X)$ is an exact subcategory of $\coh(X)$ for all $n\in\Z$.
\item[(c)] For $E$ to be $n$-regular is an open condition on $[E]\in\M$. Write $\M_{\text{$n$-reg}}\subset\M$ for the open substack of $n$-regular sheaves, and $\bs\M_{\text{$n$-reg}}\subset\bs\M$ for the corresponding derived open substack.
\item[(d)] If $E$ is $n$-regular then $H^i(E(n))=0$ for $i>0$ by (a), so $\dim H^0(E(n))=\chi(\lb\O_X(-n)\rb,\lb E\rb)$ depends only on $n$ and $\lb E\rb\in K(\coh(X))$. There is a natural vector bundle on $\M_{\text{$n$-reg}}$, in the sense of Artin stacks, with fibre $H^0(E(n))$ at $[E]$ in $\M_{\text{$n$-reg}}$. This also works for derived moduli stacks, giving a vector bundle on~$\bs\M_{\text{$n$-reg}}$.
\item[(e)] By \cite[Lem.~1.7.2]{HuLe2}, if $E$ is $n$-regular then $E(n)$ is globally generated. In particular, this implies that the following natural map is injective:
\begin{equation*}
\Hom(E,E)\cong \Hom(E(n),E(n))\longra \Hom\bigl(H^0(E(n)),H^0(E(n))\bigr).
\end{equation*}
\item[(f)] By \cite[Lem.~1.7.6]{HuLe2}, given any bounded family of coherent sheaves on $X$, there exists $N\gg 0$ such that every sheaf in the family is $n$-regular for all $n\ge N$. Hence, if $S\subseteq\M$ is any finite type substack, there exists $N\gg 0$ such that $S\subset\M_{\text{$n$-reg}}$ for all~$n\ge N$.
\end{itemize}

We define the data $\bigl\{(\B_k,F_k,\la_k):k\in K\bigr\}$ in Assumption \ref{co5ass1}(g). Set $K=\N$, and for each $n\in\N$ define $\B_n=\coh_{\text{$n$-reg}}(X)\subset\coh(X)$, and define $F_n:\A\ra\Vect_\C$ to map $E\mapsto H^0(E(n))$ on objects, and to map $\phi:E\ra F$ to $H^0(\phi\ot\id_{\O_X(n)}):H^0(E(n))\ra H^0(F(n))$ on morphisms. On restricting $F_n$ to $\B_n\subset\A$, since $H^1(E(n))=0$, the long exact sequence in sheaf cohomology implies that an exact sequence $0\ra E\ra F\ra G\ra 0$ in $\B_n$ maps to an exact sequence $0\ra H^0(E(n))\ra H^0(F(n))\ra H^0(G(n))\ra 0$ in $\Vect_\C$, so $F_n\vert_{\B_n}$ is an exact functor. Define a group morphism $\la_n:K(\B_n)\ra\Z$ to map $\la_n:\al\mapsto \chi(\lb\O_X(-n)\rb,\al)$. The rest of Assumption \ref{co5ass1}(g)(i)--(iv) for $\A=\coh(X)$ now follow from (a)--(f) above. 
\end{dfn}

\begin{rem}
\label{co7rem6}
In Definition \ref{co7def9} we chose $\A=\B=\coh(X)$, for simplicity. However, there are other natural choices for the exact subcategory $\B\subset\coh(X)$. For example, we could take $\B=\coh_{\rm tf}(X)$ to be the full subcategory of torsion-free coherent sheaves on $X$. The important conditions on $\B$ are that to lie in $\B$ should be an open condition in the moduli stack $\M$, and restricting to $\B$ should not break properness in Assumption \ref{co5ass2}(g),(h). For instance, taking $\B=\vect(X)$ to be the exact subcategory of vector bundles would be a bad choice when $\dim X>1$, because of properness issues.	
\end{rem}

\begin{rem}
\label{co7rem7}
The following will be needed in \S\ref{co73}. Let $\al\in C(\coh(X))$. As $\cExt^\bu\ra\M\t\M$ is perfect in the interval $[0,m]$, equations \eq{co7eq38} and \eq{co7eq43} give
\e
\text{$i_\al^*(\bL_{\bs\M_\al})$ is perfect in the interval $[1-m,1]$.}
\label{co7eq48}
\e
Note that the analogue is false for $\bs\baM_\al$, the moduli stack for $D^b\coh(X)$.

As $\bs\Pi_\al^\pl\!:\!\bs\M_\al\!\ra\!\bs\M^\pl_\al$ is a principal $[*/\bG_m]$-bundle, $\bL_{\bs\M_\al/\bs\M^\pl_\al}\!\cong\!\O_{\bs\M_\al}[-1]$, so by \eq{co2eq11}
we have a distinguished triangle on $\bs\M_\al$:
\begin{equation*}
\xymatrix@C=30pt{
\O_{\bs\M_\al}[-2]\ar[r] & (\bs\Pi_\al^\pl)^*(\bL_{\bs\M_\al^\pl}) \ar[r] & \bL_{\bs\M_\al} \ar[r] & \O_{\bs\M_\al}[-1]. }
\end{equation*}
Pulling this back by $i_\al:\M_\al\hookra\bs\M_\al$ yields a distinguished triangle on $\M_\al$:
\begin{equation*}
\xymatrix@C=15pt{
\O_{\M_\al}[-2]\ar[r] & (\Pi_\al^\pl)^*((i_\al^\pl)^*(\bL_{\bs\M_\al^\pl})) \ar[r] & i_\al^*(\bL_{\bs\M_\al}) \ar[r] & \O_{\M_\al}[-1]. }
\end{equation*}
The long exact sequence of cohomology sheaves of the dual then implies that
\begin{equation*}
\begin{aligned}
h^k(i_\al^*(\bL_{\bs\M_\al})^\vee)&\cong h^k((\Pi_\al^\pl)^*((i_\al^\pl)^*(\bL_{\bs\M_\al^\pl}))^\vee)\\
&\cong (\Pi_\al^\pl)^*(h^k((i_\al^\pl)^*(\bL_{\bs\M_\al^\pl})^\vee))
\end{aligned}
\qquad\text{unless $k=-1,-2$,}
\end{equation*}
where the second isomorphism holds as $\Pi_\al^\pl$ is a principal $[*/\bG_m]$-bundle. Since $\bL_{\bs\M_\al^\pl}$ is perfect in $(-\iy,1]$ as $\bs\M_\al^\pl$ is a derived Artin $\C$-stack, we deduce that \e
\text{$i_\al^*(\bL_{\bs\M_\al})$ is perfect in $[a,1]$ if and only if $(i_\al^\pl)^*(\bL_{\bs\M_\al^\pl})$ is, for $a\le 1$.}
\label{co7eq49}
\e
Combining \eq{co7eq48}--\eq{co7eq49} we see that
\e
\text{$(i_\al^\pl)^*(\bL_{\bs\M_\al^\pl})$ is perfect in the interval $[1-m,1]$.}
\label{co7eq50}
\e
Since $\bs\M_\al,\bs\M_\al^\pl$ are locally finitely presented as in Definition \ref{co7def10}, from Remark \ref{co2rem6}(e) we deduce that $\bs\M_\al$ and $\bs\M_\al^\pl$ are quasi-smooth if~$m\le 2$.
\end{rem}

\subsubsection{Assumptions \ref{co5ass2}(a)--(d),(f) and \ref{co5ass3}(a) for $\coh(X)$}
\label{co723}

\begin{dfn}
\label{co7def11} 
Work in the situation of Definitions \ref{co7def1}, \ref{co7def9} and \ref{co7def10}, and use the notation of Definitions \ref{co7def3} and \ref{co7def4}. We will define the set $\sS$ of weak stability conditions on $\A=\coh(X)$ in Assumption~\ref{co5ass2}.

Choose a connected open subset $U\subseteq\Kah(X)$, satisfying the conditions:
\begin{itemize}
\setlength{\itemsep}{0pt}
\setlength{\parsep}{0pt}
\item[(a)] If $m=\dim X\le 2$ then $U=\Kah(X)$.
\item[(b)] If $m>2$ then for all $\om,\ti\om\in U$, if $0\ne\eta\in H^{1,1}(X,\R)$ with $\int_X\eta\cup \om^{m-1}=0$ then $\int_X\eta^2\cup\ti\om^{m-2}<0$.
\end{itemize}
These conditions will be used in the proof in \S\ref{co752} of Proposition \ref{co7prop18}, which says that Assumption \ref{co5ass3}(b) holds for $\al\in C(\coh(X))$ with $\dim\al=m$.

By properties of compact K\"ahler manifolds, if $m\ge 2$ and $\om\in\Kah(X)$ then the quadratic form $Q_\om(\eta)=\int_X\eta^2\cup\om^{m-2}$ on $H^{1,1}(X,\R)\ni\eta$ has signature $(1,h^{1,1}(X)-1)$, and $Q_\om(\om)>0$, so if $0\ne\eta\in H^{1,1}(X,\R)$ with $\int_X\eta\cup \om^{m-1}=0$ then $Q_\om(\eta)=\int_X\eta^2\cup\om^{m-2}<0$. This is an open condition on $\om\in\Kah(X)$, so $\int_X\eta^2\cup\ti\om^{m-2}<0$ if $\ti\om$ is sufficiently close to $\om$ in $\Kah(X)$. Hence (b) holds when $m>2$ for any sufficiently small connected open $U\subset\Kah(X)$. When $m=2$, (b) holds with $U=\Kah(X)$ as  $\int_X\eta^2\cup\ti\om^{m-2}$ is independent of $\ti\om$. Define
\e
\sS=\bigl\{(\tau^\om,G,\le):\om\in U\bigr\}\cup\bigl\{(\mu^\om,M,\le):\om\in U\bigr\}.
\label{co7eq51}
\e
\end{dfn}

We now prove Assumptions \ref{co5ass2}(a)--(d),(f)--(g) and \ref{co5ass3}(a).
\smallskip

\noindent{\bf Assumption \ref{co5ass2}(a).} For $(\tau,G,\le)$ and $(\mu,M,\le)$ as in Definitions \ref{co7def3}--\ref{co7def4}, $\A=\coh(X)$ is $\tau$-artinian and $\mu$-artinian by Lemma \ref{co7lem1}. This proves the first part of Assumption \ref{co5ass2}(a) for $\A$, and the second part is trivial as $\B=\A$.
\smallskip

\noindent{\bf Assumption \ref{co5ass2}(b).} This follows from Propositions \ref{co7prop7}--\ref{co7prop8}, noting as in \eq{co7eq46}--\eq{co7eq47} that if $\al\in C(\B)_\pe$ then $\dim\al=0,1$ or~$m$.
\smallskip

\noindent{\bf Assumption \ref{co5ass2}(c).} All $(\tau,T,\le)$ in $\sS$ in \eq{co7eq51} have the property that if $\al,\be\in C(\B)$ with $\tau(\al)=\tau(\be)$ then $\dim\al=\dim\be$. Thus Assumption \ref{co5ass2}(c) follows from~\eq{co7eq46}--\eq{co7eq47}.
\smallskip

\noindent{\bf Assumption \ref{co5ass2}(d).} In the situation of Assumption \ref{co5ass2}(d), let $\dim\al=d$. Then $\tau(\be)=\tau(\al)$ for $\be\in I$ implies that $\dim\be=d$ for $\be\in I$. By \eq{co7eq51}, we divide into cases (i) $(\ti\tau,\ti T,\le)=(\tau^\om,G,\le)$ and (ii) $(\ti\tau,\ti T,\le)=(\mu^\om,M,\le)$. For $\be\in K(\A)$, write $P_\be^\om(n)=\sum_{k=0}^mC_{\be,k}n^k$. Then $C_{\be,d}=r_\be^\om$ if $\be\in C(\A)$ with~$\dim\be=d$.

In case (i), for $N$ chosen later, define a group morphism $\la:K(\A)\ra\R$ by
\begin{equation*}
\la(\be)=C_{\al,d}\cdot P_\be^\om(N)-C_{\be,d}\cdot P_\al^\om(N).
\end{equation*}
Then $\la(\al)=0$, and if $\be\in I$, so that $\dim\be=d$ and $C_{\be,d}=r^\om_\be$, we have
\e
\la(\be)=r_\al^\om r_\be^\om\bigl(\tau^\om(\be)(N)-\tau^\om(\al)(N)\bigr).
\label{co7eq52}
\e 
Since $\dim\al=\dim\be$ we see from Definition \ref{co7def3} that $\tau^\om(\be)>\tau^\om(\al)$ if and only if $\tau^\om(\be)(N)>\tau^\om(\al)(N)$ for $N\gg 0$, and as $r_\al^\om, r_\be^\om>0$ this holds if and only if $\la(\be)>0$. As $I\subseteq C(\B)_\pe$ is finite, we can choose one $N\gg 0$ to define $\la$ such that this works for all $\be\in I$, as we have to prove.

In case (ii), for $N$ chosen later, define a group morphism $\la:K(\A)\ra\R$ by
\begin{equation*}
\la(\be)=N^{d-1}\bigl(C_{\al,d}\cdot C_{\be,d-1}-C_{\be,d}\cdot C_{\al,d-1}\bigr).
\end{equation*}
Then $\la(\al)=0$, and if $\be\in I$, so that $\dim\be=d$, as in \eq{co7eq52} we have
\begin{equation*}
\la(\be)=r_\al^\om r_\be^\om\bigl(\mu^\om(\be)(N)-\mu^\om(\al)(N)\bigr).
\end{equation*}
We then complete the proof as in case~(i).
\smallskip

\noindent{\bf Assumption \ref{co5ass2}(f).} As $X$ is projective it has an ample line bundle $L\ra X$. Set $\om=c_1(L)$. As in Definition \ref{co7def3}, for $\al\in C(\A)$ the Hilbert polynomial $P^\om_\al(n)$ is of degree $\dim\al$, with positive leading coefficient $r_\al^\om$, and $P^\om_\al$ maps $\Z\ra\Z$ in this case as $\om=c_1(L)$. The fact that $P^\om_\al$ maps $\Z\ra\Z$ implies that $(\dim\al)!\cdot r_\al^\om$ is an integer (for an idea of why, note that $n\mapsto\binom{n}{\dim\al}$ is of degree $\dim\al$, maps $\Z\ra\Z$, and has leading coefficient $1/(\dim\al)!$). Define $\rk:C(\A)\ra\N_{>0}$ by
\e
\rk(\al)=(\dim\al)!\cdot r_\al^\om.
\label{co7eq53}
\e
If $\al,\be\in C(\A)$ and $(\tau,T,\le)\in\sS$ with $\tau(\al)=\tau(\be)$ then $\dim\al=\dim\be$, so $P_\al^\om,P_\be^\om,P_{\al+\be}^\om$ have the same degree $\dim\al$, and thus $r_{\al+\be}^\om=r_\al^\om+r_\be^\om$ as $P_{\al+\be}^\om=P_\al^\om+P_\be^\om$, giving $\rk(\al+\be)=\rk(\al)+\rk(\be)$.
\smallskip

\noindent{\bf Assumption \ref{co5ass3}(a).} For $\sS$ as in \eq{co7eq51}, if $\om_0,\om_1\in U$, choose a continuous family $(\om_t)_{t\in[0,1]}$ in $U$ linking $\om_0,\om_1$, which is possible as $U\subseteq\Kah(X)$ is open and connected. Then $(\mu^{\om_t},G,\le)_{t\in[0,1]}$ is a continuous family of weak stability conditions on $\A$ in the sense of Definition \ref{co3def5}, so $(\mu^{\om_0},G,\le)$ and $(\mu^{\om_1},G,\le)$ are continuously connected in $\sS$ in the sense of Assumption \ref{co5ass3}(a).

Next, if $\om\in U$, define $(\tau_t,T_t,\le)_{t\in[0,1]}$ by
\begin{equation*}
(\tau_t,T_t,\le)=\begin{cases} (\tau^\om,G,\le), & t\in[0,\ha), \\ (\mu^\om,M,\le), & t\in[\ha,1]. \end{cases}
\end{equation*}
Since $(\mu^\om,G,\le)$ dominates $(\tau^\om,G,\le)$, this is a continuous family, so $(\tau^\om,G,\le)$ and $(\mu^\om,M,\le)$ are continuously connected in $\sS$. As being continuously connected in $\sS$ is an equivalence relation, Assumption \ref{co5ass3}(a) follows.

\subsubsection{The $G$-equivariant case and Assumptions \ref{co4ass2}, \ref{co5ass4}(a),(b),(d)}
\label{co724}

In \S\ref{co45} and \S\ref{co54} we explained how to extend our theory from homology to $G$-equivariant homology for an algebraic $\C$-group $G$, by adding extra Assumptions \ref{co4ass2} and \ref{co5ass4} on the action of $G$ on the data $\A,K(\A),\M,\ldots.$ We now show that if, in the situation of \S\ref{co721}--\S\ref{co723}, $G$ acts on $X$ fixing at least one ample line bundle, then Assumptions \ref{co4ass2} and \ref{co5ass4}(a),(b)(d) hold. The proofs are easy. We verify Assumption \ref{co5ass4}(c) in \S\ref{co735}, when we discuss obstruction theories.

\begin{dfn}
\label{co7def12}
Work in the situation of \S\ref{co721}--\S\ref{co723}, defining data $\A=\B=\coh(X),K(\A),\M,\ldots$ satisfying (most of) Assumptions \ref{co4ass1} and \ref{co5ass1}--\ref{co5ass3}. Suppose in addition that $G$ is an algebraic $\C$-group with an action $\phi:G\t X\ra X$ on the $\C$-scheme $X$, such that:
\begin{itemize}
\setlength{\itemsep}{0pt}
\setlength{\parsep}{0pt}
\item[(i)] There exists at least one ample line bundle $\O_X(1)\ra X$ such that the $G$-action on $X$ lifts to $\O_X(1)$.
\item[(ii)] The action of $G$ on $K(\coh(X))=K^\num(\coh(X))$ is trivial. This is automatic if $G$ is connected, or if the action of $G$ on $H^{\rm even}(X,\Q)$ is trivial.
\end{itemize}
Note that if $X$ is a smooth elliptic curve, thought of as an algebraic $\C$-group, and $G=X$ then (i) does not hold, so (i) excludes some interesting examples. We will define the data and verify the conditions of Assumptions \ref{co4ass2} and~\ref{co5ass4}(a),(b),(d).
\smallskip

\noindent{\bf Assumption \ref{co4ass2}(a).} If $g\in G(\C)$ then $g^{-1}\in G(\C)$ so $\phi(g^{-1}):X\ra X$ is an isomorphism of projective $\C$-schemes. For (a)(i), define a functor $\Ga(g)=\phi(g^{-1})^*:\coh(X)\ra\coh(X)$. Then $\Ga(g):\B\ra\B$ is an equivalence of abelian categories as $\phi(g^{-1})$ is an isomorphism. For (a)(ii), if $g,h\in G(\C)$, define
\begin{equation*}
\Ga_{g,h}:\Ga(g)\ci\Ga(h)=\phi(g^{-1})^*\ci\phi(h^{-1})^*\Longra \phi(h^{-1}\ci g^{-1})^*= \phi((gh)^{-1})^*=\Ga(gh)
\end{equation*}
to be the natural isomorphism from composition of pullbacks. Here as $g\mapsto g^{-1}$ and $\phi\mapsto\phi^*$ are both contravariant, $g\mapsto\phi(g^{-1})^*$ is covariant. Then \eq{co4eq19} commutes as $\phi:G(\C)\t X\ra X$ is a group action. For (a)(iii), writing $1\in G(\C)$ for the identity, $\phi(1)=\id_X:X\ra X$, and we define $\Ga_1:\Ga(1)=\id_X^*\Ra \Id_\coh(X)$ to be the canonical natural isomorphism. Then \eq{co4eq20} holds.
\smallskip

\noindent{\bf Assumption \ref{co4ass2}(b)--(d).} Part (b) holds as the actions $\phi(g^{-1}):X\ra X$ for $g\in G(\C)$ are induced by the scheme-theoretic group action $\phi:G\t X\ra X$. 

In more detail, let $\cU\ra X\t\M$ be the universal sheaf in Definition \ref{co7def2}, which is flat over $\M$. Then $(\phi\t\id_\M)^*(\cU)$ is a sheaf over $(G\t X)\t\M\cong X\t(G\t\M)$, which is flat over $G\t\M$ as $\phi$ is smooth. Hence by the universal property of $\cU$ there is a morphism $\psi:G\t\M\ra\M$ such that $(\phi\t\id_\M)^*(\cU)\cong(\id_X\t\psi)^*(\cU)$. Writing $\io:G\ra G$ for the inverse map, as $\phi$ is a $G$-action on $X$ we can show that $\xi=\psi\ci(\io\t\id_\M):G\t\M\ra\M$ is a $G$-action on $\M$ (again we compose with inverses $\io$ as pullback is contravariant). If $g\in G(\C)$ then
\begin{equation*}
\xi\vert_{\{g\}\t\M}=\psi\vert_{\{g^{-1}\}\t\M}=\bigl((\phi\vert_{\{g^{-1}\}\t X})^*\bigr)_*=\bigl(\phi(g^{-1})^*\bigr)_*=\Ga(g)_*,
\end{equation*}
so on the level of $\C$-points, the $G$-action $\ti\psi:G\t\M\ra\M$ is induced by the functors $\Ga(g)$ in (a). Parts (c),(d) are obvious.

Writing $\baM$ for the moduli stack of objects in $D^b\coh(X)$ as in Definition \ref{co7def2} and using the universal complex $\cU^\bu\ra X\t\baM$, $\psi,\xi$ above extend to $\bar\psi,\bar\xi:G\t\baM\ra\baM$ with $\bar\xi$ a $G$-action.
\smallskip

\noindent{\bf Assumption \ref{co4ass2}(e).} This holds by (ii) above.
\smallskip

\noindent{\bf Assumption \ref{co4ass2}(f).} We can show using \eq{co7eq36} and the defining property $(\bar\phi\t\id_\M)^*(\cU)\cong(\id_X\t\bar\psi)^*(\cU^\bu)$ of the $G$-action on $\baM$ that the Ext complex $\cExt^\bu$ in \eq{co7eq36} is equivariant under the diagonal action of $G$ on $\baM\t\baM$. Thus $\cExt^\bu$ is also $G$-equivariant on the $G$-invariant open substack $\M\t\M\subset\baM\t\baM$. So $\cE^\bu=(\cExt^\bu)^\vee$ in \eq{co7eq38} is $G$-equivariant. The proofs of \eq{co4eq3}--\eq{co4eq6} in Definition \ref{co7def9} work $G$-equivariantly.
\smallskip

\noindent{\bf Assumption \ref{co5ass4}(a).} This follows from Assumption \ref{co4ass2}(a),(e), as~$\B=\A$.
\smallskip

\noindent{\bf Assumption \ref{co5ass4}(b).} This holds as the definitions of $G$-actions on $\M,\M^\pl$ above, and the $G$-equivariance of $\Phi,\Psi,\Pi^\pl$, also work for derived stacks, and the proofs of \eq{co5eq1} and \eq{co5eq3} in Definition \ref{co7def10} work $G$-equivariantly.
\smallskip

\noindent{\bf Assumption \ref{co5ass4}(d).} Here we make a small addition to Definition \ref{co7def10}: we  require that the $G$-action on $X$ should have a choice of lift to the ample line bundle $\O_X(1)$ on $X$ used to define the data $\bigl\{(\B_k,F_k,\la_k):k\in K=\N\bigr\}$. This is possible, for some choice of $\O_X(1)$, by condition (i) above. Then the $G$-action also lifts naturally to $\O_X(k)$ for $k\in\N$, which is used to define $(\B_k,F_k,\la_k)$. It is easy to show that $(\B_k,F_k,\la_k)$ is $G$-equivariant in the sense required.
\end{dfn}

\subsection{\texorpdfstring{Quasi-smooth derived stacks, obstruction theories, \\ and Assumptions \ref{co5ass1}(f) and \ref{co5ass2}(e)}{Quasi-smooth derived stacks, obstruction theories, and Assumptions \ref{co5ass1}(f) and \ref{co5ass2}(e)}}
\label{co73}

We now explain how to construct the quasi-smooth derived stacks $\bs\dM_\al^\red,\bs\dM_\al^\rpl$ in Assumption \ref{co5ass1}(f), which yield obstruction theories \eq{co5eq5} on~$\dM_\al,\dM_\al^\pl$.

The obvious choice is to take $\dM_\al=\M_\al$, $\dM_\al^\pl=\M_\al^\pl$, $\bs\dM_\al^\red=\bs\M_\al$, $\bs\dM_\al^\rpl=\bs\M_\al^\pl$, $\bs j_\al=\bs\id=\bs j_\al^\pl$, and $U_\al=o_\al=0$. There are two issues with this:
\begin{itemize}
\setlength{\itemsep}{0pt}
\setlength{\parsep}{0pt}
\item[(i)] By \eq{co7eq48} and \eq{co7eq50}, $i^*(\bL_{\bs\M_\al}),i^*(\bL_{\bs\M_\al^\pl})$ are perfect in $[1-m,1]$. Thus $\bs\M_\al,\bs\M_\al^\pl$ are only quasi-smooth if $m\le 2$, so for $X$ a curve or surface.

If $X$ is a Fano 3-fold and $\dim\al>0$ then we can define open substacks $\dM_\al\subseteq\M_\al$ and $\dM_\al^\pl\subseteq\M_\al^\pl$ satisfying Assumption \ref{co5ass2}(e) such that $\bs\dM_\al,\bs\dM_\al^\pl$ are quasi-smooth.

Therefore we restrict to $X$ a curve, surface or Fano 3-fold in this section.
\item[(ii)] If $X$ is a surface with geometric genus $p_g=\dim H^0(K_X)>0$, and $\rank\al>0$, then $h^{-1}(i^*(\bL_{\bs\M_\al}))$ and $h^{-1}(i^*(\bL_{\bs\M_\al^\pl}))$ contain trivial vector bundles with fibre $H^0(K_X)$. As in Remark \ref{co5rem1}(c), Theorem \ref{co2thm1}(iv) implies that virtual classes defined using the obstruction theories \eq{co5eq5} on $\M_\al,\M_\al^\pl$ are zero, so our theory would yield invariants~$[\M_\al^\ss(\tau)]_\inv=0$.

In this case we define different, `reduced' versions $\bs\dM_\al^\red,\bs\dM_\al^\rpl$ of $\bs\M_\al,\ab\bs\M_\al^\pl$, which can give nonzero invariants $[\M_\al^\ss(\tau)]_\inv$.
\end{itemize}

\subsubsection{(Quasi-)smooth moduli stacks for $X$ a curve}
\label{co731}

We verify Assumptions \ref{co5ass1}(f) and \ref{co5ass2}(e) when $X$ is a curve.

\begin{dfn}
\label{co7def13}
In Definition \ref{co7def1}, suppose $m=\dim X=1$, so that $X$ is a connected, smooth, projective complex curve. Since $\Ext^i(E,E)=0$ unless $i=0,1$ for $E\in\coh(X)$, the moduli stacks $\M,\M^\pl$ are smooth Artin stacks. The derived versions $\bs\M,\bs\M^\pl$ in Definition \ref{co7def10} are $\bs\M=\M$ and~$\bs\M^\pl=\M^\pl$.

As in \eq{co7eq46} we set $C(\coh(X))_\pe=C(\coh(X))$. In Assumption \ref{co5ass1}(f), in (i)--(ii) for each $\al\in C(\coh(X))$ we define $\dM_\al=\M_\al$, $\dM_\al^\pl=\M_\al^\pl$, $\bs\dM_\al^\red=\bs\M_\al$, $\bs\dM_\al^\rpl=\bs\M_\al^\pl$, $\bs j_\al=\bs\id=\bs j_\al^\pl$, $\bs{\dot\Pi}{}_\al^\rpl=\bs\Pi_\al^\pl$, and $U_\al=o_\al=0$. Then $\bs\dM_\al^\red,\bs\dM_\al^\rpl$ are locally finitely presented as in Definition \ref{co7def10} and quasi-smooth as in Remark \ref{co7rem7}, as $m\le 2$. Also \eq{co5eq4} is trivially Cartesian, proving (ii), equation \eq{co5eq6} holds as $\bs j_\al,\bs j_\al^\pl$ are identities and $U_\al=0$, giving (iii), and (iv) is trivial as $o_\al=o_\be=o_{\al+\be}=0$. This proves Assumption~\ref{co5ass1}(f).

Assumption \ref{co5ass2}(e) is trivial as~$\dM_\al=\M_\al$.
\end{dfn}

\begin{rem}
\label{co7rem8}
When $\A=\coh(X)$ for $X$ a curve, we do not need Behrend--Fantechi obstruction theories at all, as everywhere we use virtual classes $[S]_\virt$ (e.g.\ $[\M_\al^\ss(\tau)]_\virt$ when $\M_\al^\rst(\tau)=\M_\al^\ss(\tau)$ in Theorem \ref{co5thm1}(i)) then $S$ is a smooth proper algebraic $\C$-space, and the virtual class $[S]_\virt=[S]_\fund$ is just the fundamental class of $S$ as a compact complex manifold, as in Theorem~\ref{co2thm1}(i).
\end{rem}

\subsubsection{Quasi-smooth moduli stacks for $X$ a surface: the case $p_g=0$}
\label{co732}

\begin{dfn}
\label{co7def14}
In Definition \ref{co7def2}, suppose $m=\dim X=2$, so that $X$ is a connected, smooth, projective complex surface. The {\it geometric genus\/} of $X$ is $p_g=\dim H^0(K_X)$. In terms of the Hodge numbers $h^{p,q}(X)$ of the compact K\"ahler manifold $X$ we have $p_g=h^{2,0}(X)=h^{0,2}(X)$. Considering $X$ as a compact, oriented 4-manifold, we have~$b^2_+(X)=1+2p_g$. 	
\end{dfn}

For reasons explained in Remark \ref{co7rem9} below, we will give different definitions of the quasi-smooth stacks $\bs\dM_\al^\red,\bs\dM_\al^\rpl$ in Assumption \ref{co5ass1}(f), and hence of the obstruction theories \eq{co5eq5}, in the cases $p_g=0$, and $p_g>0$. This section handles the simpler case $p_g=0$. Definition \ref{co7def15} would also be valid for $p_g>0$, but then the invariants $[\M_\al^\ss(\tau)]_\inv$ in Theorems \ref{co5thm1}--\ref{co5thm3} would be zero for~$\rank\al>0$.

\begin{dfn}
\label{co7def15}
In Definitions \ref{co7def2} and \ref{co7def14}, suppose $m=2$ and $p_g=0$. As in \eq{co7eq46} we set $C(\coh(X))_\pe=C(\coh(X))$. In Assumption \ref{co5ass1}(f), for each $\al\in C(\coh(X))$ we define $\dM_\al=\M_\al$, $\dM_\al^\pl=\M_\al^\pl$, $\bs\dM_\al^\red=\bs\M_\al$, $\bs\dM_\al^\rpl=\bs\M_\al^\pl$, $\bs j_\al=\bs\id=\bs j_\al^\pl$, $\bs{\dot\Pi}{}_\al^\rpl=\bs\Pi_\al^\pl$, and $U_\al=o_\al=0$. Then Assumption \ref{co5ass1}(f) holds as in Definition \ref{co7def13}. Assumption \ref{co5ass2}(e) is trivial as~$\dM_\al=\M_\al$.

The obstruction theories \eq{co5eq5} on $\M_\al,\M_\al^\pl$ are essentially the same as those in Mochizuki \cite[Prop.s 5.6.4 \& 6.1.1]{Moch}, who uses classical Algebraic Geometry.
\end{dfn}

\subsubsection{Quasi-smooth moduli stacks for $X$ a surface: the case $p_g>0$}
\label{co733}

We explain why we define $\bs\dM_\al^\red,\bs\dM_\al^\rpl$ differently when~$p_g>0$:

\begin{rem}
\label{co7rem9}
{\bf(a)} In Definitions \ref{co7def2} and \ref{co7def14}, suppose $m=2$ and $p_g>0$, and let $\al\in C(\coh(X))$ with $\rank\al=r>0$. Then we could define $\bs\dM_\al^\red,\bs\dM_\al^\rpl$ and the other data of Assumption \ref{co5ass1}(f) as in Definition \ref{co7def15}, which does not use $p_g=0$. So \eq{co5eq5} would give perfect obstruction theories $\bL_i:i^*(\bL_{\bs\M_\al})\ra\bL_{\M_\al}$ and $\bL_i:i^*(\bL_{\bs\M_\al^\pl})\ra\bL_{\M_\al^\pl}$ on~$\M_\al,\M_\al^\pl$. 

At a point $[E]$ in $\M_\al$ or $\M_\al^\pl$, by Serre duality we would have
\begin{equation*}
h^1\bigl(i^*(\bL_{\bs\M_\al})^\vee\vert_{[E]}\bigr),h^1\bigl(i^*(\bL_{\bs\M_\al^\pl})^\vee\vert_{[E]}\bigr)\cong \Ext^2(E,E)\cong\Hom(E,E\ot K_X)^*.
\end{equation*}
Now we have natural linear maps
\e
\xymatrix@C=45pt{ H^0(K_X)^* \ar[r]^(0.4){\Tr^*_E} & \Hom(E,E\ot K_X)^* \ar[r]^(0.59){\Id^*_E} & H^0(K_X)^*, }
\label{co7eq54}
\e
whose composition is $r\,\id_{H^0(K_X)^*}$, so in particular $\Id_E^*$ is surjective as $r>0$.

Similarly, globally on $\M_\al,\M_\al^\pl$ we have surjective morphisms
\e
\begin{split}
\Id^*&:h^1\bigl(i^*(\bL_{\bs\M_\al})^\vee\bigr)\longra H^0(K_X)^*\ot\O_{\M_\al},\\ 
\Id^*&:h^1\bigl(i^*(\bL_{\bs\M_\al^\pl})^\vee\bigr)\longra H^0(K_X)^*\ot\O_{\M_\al^\pl}. 
\end{split}
\label{co7eq55}
\e
As $\dim H^0(K_X)^*=p_g>0$, Theorem \ref{co2thm1}(iv) implies that any virtual classes defined using $\bL_i:i^*(\bL_{\bs\M_\al})\ra\bL_{\M_\al}$ or $\bL_i:i^*(\bL_{\bs\M_\al^\pl})\ra\bL_{\M_\al^\pl}$, such as $[\M_\al^\ss(\tau)]_\virt$ in Theorem \ref{co5thm1}(i), are automatically zero. So the invariants $[\M_\al^\ss(\tau)]_\inv$ in Theorems \ref{co5thm1}--\ref{co5thm3} would be always zero when $\rank\al>0$, and so would be boring.

\smallskip

\noindent{\bf(b)} Here is a way to justify {\bf(a)}: as in Theorem \ref{co5thm1}, the invariants $[\M_\al^\ss(\tau)]_\inv$ should be unchanged by definitions of the complex structure $J$ of $X$. Now we can only have $\M_\al^\ss(\tau)_J\ne\es$ on $(X,J)$ if $c_1(\al)\in H^{1,1}_J(X,\R)\subseteq H^2(X,\R)$, where $H^{1,1}_J(X,\R)$ has codimension $2p_g>0$. If $c_1(\al)\ne 0$ then for a generic deformation $J'$ of $J$ we have $c_1(\al)\notin H^{1,1}_{J'}(X,\R)$, so $\M_\al^\ss(\tau)_{J'}=\es$, and~$[\M_\al^\ss(\tau)]_\inv=0$.
\smallskip

\noindent{\bf(c)} To get nonzero invariants when $p_g>0$, we will modify the quasi-smooth derived stacks $\bs\M_\al,\bs\M_\al^\pl$ to new quasi-smooth derived stacks $\bs\dM_\al^\red,\bs\dM_\al^\rpl$, deleting the trivial factors $H^0(K_X)^*$ in $h^1(i^*(\bL_{\bs\M_\al})^\vee),h^1(i^*(\bL_{\bs\M_\al^\pl})^\vee)$. This modifies the corresponding obstruction theories \eq{co5eq5}. At the level of obstruction theories, this modification is a well-known construction, known as {\it reduced obstruction theories}, as in Okounkov--Pandharipande \cite[\S 3.4]{OkPa}, Maulik--Pandharipande \cite[\S 2.2]{MaPa}, Maulik--Pandharipande--Thomas \cite[\S A]{MPT}, Kool--Thomas \cite{KoTh1,KoTh2}, Sch\"urg--To\"en--Vezzosi \cite[\S 4]{STV}, and Sch\"urg~\cite{Schu}.

Definition \ref{co7def16} below gives a direct construction of reduced obstruction theories in our case using Derived Algebraic Geometry.

Mochizuki \cite[\S 5 \& \S 6.1]{Moch} defines moduli spaces and obstruction theories which are essentially equivalent to our reduced obstruction theories \eq{co5eq5} in the case $b^1(X)=0$ by considering moduli stacks of sheaves of rank $r>0$ with {\it fixed determinant}. See also Thomas \cite[\S 3]{Thom1} and Huybrechts--Thomas \cite[\S 4]{HuTh} for construction of obstruction theories in the fixed determinant case.
\smallskip

\noindent{\bf(d)} The vector spaces $U_\al$ in Assumption \ref{co5ass1}(f)(iii) are there precisely to include reduced obstruction theories in our framework. We will take~$U_\al=H^0(K_X)$.	
\end{rem}

We will define the data and verify Assumptions \ref{co5ass1}(f) and~\ref{co5ass2}(e):

\begin{dfn}
\label{co7def16}
In Definitions \ref{co7def2} and \ref{co7def14}, suppose $m=2$ and $p_g>0$. As in \eq{co7eq46} we set $C(\coh(X))_\pe=C(\coh(X))$. In Assumption \ref{co5ass1}(f)(i), for each $\al\in C(\coh(X))$ we define $\dM_\al^\pl=\M_\al^\pl$, $\dM_\al=\M_\al$. 

In Assumption \ref{co5ass1}(f)(ii)--(iii) for $\al\in C(\coh(X))$ we divide into cases (A) $\rank\al=0$ and (B) $\rank\al>0$. In case (A) we define $\dM_\al=\M_\al$, $\dM_\al^\pl=\M_\al^\pl$, $\bs\dM_\al^\red=\bs\M_\al$, $\bs\dM_\al^\rpl=\bs\M_\al^\pl$, $\bs j_\al=\bs\id=\bs j_\al^\pl$, $\bs{\dot\Pi}{}_\al^\rpl=\bs\Pi_\al^\pl$, and $U_\al=o_\al=0$, so that Assumption \ref{co5ass1}(f) holds as in Definition~\ref{co7def15}.

For case (B), let $\al\in C(\coh(X))$ with $\rank\al=r>0$. Define a derived stack $\bs\dM_\al^\red$ as a fibre product $\bs\dM_\al\t_{\det,\bs\Pic(X),i}\Pic(X)$, in a Cartesian square of derived stacks
\e
\begin{gathered}
\xymatrix@C=170pt@R=15pt{
*+[r]{\bs\dM_\al^\red} \ar[r]_{\Pi_{\Pic(X)}} \ar[d]^{\bs j_\al} & *+[l]{\Pic(X)_{\vphantom{(}}} \ar@{_{(}->}[d]_i \\
*+[r]{\bs\dM_\al=\bs\M_\al} \ar[r]^{\det} & *+[l]{\bs\Pic(X).\!}
}	
\end{gathered}
\label{co7eq56}
\e
Here $\Pic(X),\bs\Pic(X)$ are the classical and derived Picard stacks of $X$, as in Sch\"urg--To\"en--Vezzosi \cite[Def.~2.2]{STV}. They are moduli stacks of line bundles $L\ra X$, and are open substacks $\Pic(X)\subset\M$, $\bs\Pic(X)\subset\bs\M$. Also $\Pic(X)=t_0(\bs\Pic(X))$ is the classical truncation of $\bs\Pic(X)$, with inclusion $i:\Pic(X)\hookra\bs\Pic(X)$, and $\det:\bs\M_\al\ra\bs\Pic(X)$ is the derived determinant morphism from \cite[\S 3.1]{STV}. Similarly, define $\bs\dM_\al^\rpl$ by the Cartesian square
\e
\begin{gathered}
\xymatrix@C=170pt@R=15pt{
*+[r]{\bs\dM_\al^\rpl} \ar[r]_{\Pi_{\Pic(X)^\pl}} \ar[d]^{\bs j_\al^\pl} & *+[l]{\Pic(X)^\pl_{\vphantom{(}}} \ar@{_{(}->}[d]_{i^\pl} \\
*+[r]{\bs\dM_\al^\pl=\bs\M_\al^\pl} \ar[r]^{\det^\pl} & *+[l]{\bs\Pic(X)^\pl.\!} }	
\end{gathered}
\label{co7eq57}
\e
Then $\bs\dM_\al^\red,\bs\dM_\al^\rpl$ are locally finitely presented as the other terms in \eq{co7eq56}--\eq{co7eq57} are, noting that $\Pic(X),\Pic(X)^\pl$ are smooth.

The universal property of \eq{co7eq57} gives a unique morphism $\bs{\dot\Pi}{}_\al^\rpl:\bs\dM_\al^\red\ra\bs\dM_\al^\rpl$ in the (homotopy) commuting diagram:
\e
\begin{gathered}
\xymatrix@!0@C=90pt@R=20pt{
*+[r]{\bs\dM_\al^\red} \ar@{..>}[dr]_{\bs{\dot\Pi}{}_\al^\rpl} \ar[rr]_(0.75){\Pi_{\Pic(X)}} \ar[dd]^{\bs j_\al} && *+[l]{\Pic(X)_{\vphantom{(}}} \ar[dr]^(0.4){\Pi^\pl} \ar@{_{(}->}[dd]^(0.25)i 
\\
& *+[r]{\bs\dM_\al^\rpl} \ar[dd]^(0.3){\bs j_\al^\pl}  \ar[rr]_(0.3){\Pi_{\Pic(X)^\pl}} && *+[l]{\Pic(X)^\pl_{\vphantom{(}}} \ar@{_{(}->}[dd]_{i^\pl}
\\
*+[r]{\bs\dM_\al} \ar[dr]_{\bs{\dot\Pi}{}_\al^\pl} \ar[rr]^(0.35){\det} && *+[l]{\bs\Pic(X)} \ar[dr]^(0.4){\bs\Pi^\pl}\\
& *+[r]{\bs\dM_\al^\pl} \ar[rr]^{\det^\pl} && *+[l]{\bs\Pic(X)^\pl.\!} }	
\end{gathered}
\label{co7eq58}
\e
Here $\bs{\dot\Pi}{}_\al^\rpl$ is a principal $[*/\bG_m]$-bundle, as $\bs{\dot\Pi}{}^\pl_\al,\Pi^\pl,\bs\Pi^\pl$ are. Also all six faces of the cube \eq{co7eq58} are Cartesian, as the four faces not including  `$\dashra$' are. 

The left hand diamond in \eq{co7eq58} is the Cartesian square \eq{co5eq4} in Assumption \ref{co5ass1}(f)(ii). Since $\Pic(X)=t_0(\bs\Pic(X))$, taking classical truncations in \eq{co7eq56} (which preserves Cartesian squares) shows that $t_0(\bs j_\al)$ is an isomorphism. Similarly $t_0(\bs j_\al^\pl)$ is an isomorphism in \eq{co7eq57}, as required in Assumption~\ref{co5ass1}(f)(ii).

If $L$ is a line bundle then $\Ext^2(L,L)^*\cong \Hom(L,L\ot K_X)\cong H^0(K_X)$ by Serre duality. Because of this, $\bs\Pic(X)$ is basically $\Pic(X)$ with the trivial vector bundle $H^0(K_X)\ot\O_{\Pic(X)}$ inserted in degree $-1$ as an obstruction bundle, so the relative cotangent complexes $\bL_{\Pic/\bs\Pic}$ of $i:\Pic(X)\hookra\bs\Pic(X)$, and similarly $\bL_{\Pic^\pl/\bs\Pic^\pl}$ of $i^\pl:\Pic(X)^\pl\hookra\bs\Pic(X)^\pl$, are
\e
\begin{split}
\bL_{\Pic/\bs\Pic}&\cong H^0(K_X)\ot\O_{\Pic(X)}[2],\\ 
\bL_{\Pic^\pl/\bs\Pic^\pl}&\cong H^0(K_X)\ot\O_{\Pic(X)^\pl}[2].
\end{split}
\label{co7eq59}
\e
Thus \eq{co2eq12} for the Cartesian squares \eq{co7eq56}--\eq{co7eq57} yields isomorphisms
\e
\begin{split}
\bL_{\bs\dM^\red_\al/\bs\dM_\al}&\cong\Pi_{\Pic(X)}^*(\bL_{\Pic/\bs\Pic})\cong H^0(K_X)\ot\O_{\bs\dM^\red_\al}[2], \\
\bL_{\bs\dM^\rpl_\al/\bs\dM^\pl_\al}&\cong\Pi_{\Pic(X)^\pl}^*(\bL_{\Pic^\pl/\bs\Pic^\pl})\cong H^0(K_X)\ot\O_{\bs\dM^\rpl_\al}[2].
\end{split}
\label{co7eq60}
\e
Set $U_\al=H^0(K_X)$ and $o_\al=\dim H^0(K_X)=p_g>0$. Then \eq{co7eq60} gives the isomorphisms \eq{co5eq6}, proving Assumption \ref{co5ass1}(f)(iii) in case~(B).

We have not yet proved that $\bs\dM_\al^\red,\bs\dM_\al^\rpl$ are quasi-smooth. As they are locally finitely presented from above, by Remark \ref{co2rem6}(e) it suffices to prove that $i^*(\bL_{\bs\dM_\al^\red}),i^*(\bL_{\bs\dM_\al^\rpl})$ are perfect in the interval $[-1,1]$. In the distinguished triangle \eq{co5eq7}, $i^*(\bL_{\bs\dM_\al})$ is perfect in $[-1,1]$ as $\bs\dM_\al=\bs\M_\al$ is quasi-smooth by Remark \ref{co7rem7}, so we see that $i^*(\bL_{\bs\dM_\al^\red})$ is perfect in $[-2,1]$.

Thus $i^*(\bL_{\bs\dM_\al^\red})^\vee$ is perfect in $[-1,2]$. If we can prove that $h^2(i^*(\bL_{\bs\dM_\al^\red})^\vee)=0$ then $i^*(\bL_{\bs\dM_\al^\red})^\vee$ is perfect in $[-1,1]$, which implies that $i^*(\bL_{\bs\dM_\al^\red})$ is perfect in $[-1,1]$, as we want. (Here we pass to duals as showing that $h^{-2}(i^*(\bL_{\bs\dM_\al^\red}))=0$ need not imply that $i^*(\bL_{\bs\dM_\al^\red})$ is perfect in~$[-1,1]$.)

The long exact sequence of cohomology sheaves of the dual of \eq{co5eq7} becomes
\e
\!\!\xymatrix@C=15pt{ \cdots \ar[r] & h^1(i_\al^*(\bL_{\bs\M_\al})^\vee) \ar[r]^(0.45)u & H^0(K_X)^*\!\ot\!\O_{\M_\al} \ar[r] & h^2(i^*(\bL_{\bs\dM_\al^\red})^\vee) \ar[r] & 0, }\!\!
\label{co7eq61}
\e
where $u$ is the top line of \eq{co7eq55}, and so is surjective as $\rank\al>0$. Hence $h^2(i^*(\bL_{\bs\dM_\al^\red})^\vee)=0$, so $i^*(\bL_{\bs\dM_\al^\red})$ is perfect in $[-1,1]$, and $\bs\dM_\al^\red$ is quasi-smooth. Similarly $\bs\dM_\al^\rpl$ is. This completes Assumption \ref{co5ass1}(f)(ii) in case~(B).

For Assumption \ref{co5ass1}(f)(iv), $o_\al+o_\be\ge o_{\al+\be}$ holds as $o_\ga=0$ if $\rank\ga=0$ and $o_\ga=p_g>0$ if $\rank\ga>0$. Assumption \ref{co5ass2}(e) is trivial as~$\dM_\al=\M_\al$.
\end{dfn}

\begin{rem}
\label{co7rem10}
In Definition \ref{co7def16}, we took $\bs\dM_\al^\red=\bs\M_\al$, $\bs\dM_\al^\rpl=\bs\M_\al^\pl$ when $\rank\al=0$, giving `non-reduced' obstruction theories \eq{co5eq5}, but gave special definitions of $\bs\dM_\al^\red,\bs\dM_\al^\rpl$ with $\bs\dM_\al^\red\ne\bs\M_\al$, $\bs\dM_\al^\rpl\ne\bs\M_\al^\pl$ in \eq{co7eq56}--\eq{co7eq57} when $\rank\al>0$, giving `reduced' obstruction theories~\eq{co5eq5}.

We chose $\bs\dM_\al^\red=\bs\M_\al$, $\bs\dM_\al^\rpl=\bs\M_\al^\pl$ when $\rank\al=0$ as then \eq{co7eq55} and $u$ in \eq{co7eq61} need not be surjective, so the $\bs\dM_\al^\red,\bs\dM_\al^\rpl$ in \eq{co7eq56}--\eq{co7eq57} might only have $i^*(\bL_{\bs\dM_\al^\red}),i^*(\bL_{\bs\dM_\al^\rpl})$ perfect in $[-2,1]$, and they may not be quasi-smooth. It is not clear that $\bs\dM_\al^\red=\bs\M_\al$, $\bs\dM_\al^\rpl=\bs\M_\al^\pl$ are always the best choices when $\rank\al=0$, but the author is primarily interested in~$\rank\al>0$.

To see why \eq{co7eq55} need not be surjective when $\rank\al=0$, let $E\in\coh(X)$ with $\lb E\rb=\al$, so that $\dim E=0$ or 1, and $0\ne s\in H^0(K_X)$ with $\supp E\subseteq s^{-1}(0)$ as subschemes of $X$. Then $s\ot\id_E=0$ in $\Hom(E,E\ot K_X)$, so $\Id_E^*$ in \eq{co7eq54} is not surjective, and \eq{co7eq55} is not surjective near $[E]\in\M_\al,\M_\al^\pl$.

If $X$ is a $K3$ surface, so that $K_X\cong\O_X$, any $0\ne s\in H^0(K_X)$ has $s^{-1}(0)=\es$, and $\supp E\subseteq s^{-1}(0)$ cannot happen when $\lb E\rb=\al\ne 0$. Then \eq{co7eq55} is surjective when $\al\in C(\coh(X))$ with $\rank\al=0$, and we can modify Definition \ref{co7def16} to also define $\bs\dM_\al^\red,\bs\dM_\al^\rpl$ as in \eq{co7eq56}--\eq{co7eq57} when~$\rank\al=0$.
\end{rem}

\subsubsection{Quasi-smooth moduli stacks for $X$ a Fano 3-fold}
\label{co734}

A smooth projective $\C$-scheme $X$ is called a {\it Fano manifold\/} if $K_X^{-1}\ra X$ is ample. Motivated by Thomas \cite[Cor.~3.39]{Thom1}, who defined Donaldson--Thomas invariants $[\M_\al^\ss(\tau)]_\virt$ for moduli spaces $\M_\al^\ss(\tau)$ of Gieseker semistable sheaves on $X$ with $\M_\al^\rst(\tau)=\M_\al^\ss(\tau)$, where $X$ is either a Calabi--Yau 3-fold or a Fano 3-fold, we verify Assumptions \ref{co5ass1}(f) and \ref{co5ass2}(e) for Fano 3-folds~$X$.

\begin{dfn}
\label{co7def17}
In Definition \ref{co7def2}, suppose $m=3$ and $X$ is a Fano 3-fold. As in \eq{co7eq47} we set $C(\coh(X))_\pe$ to be the subset of $\al\in C(\coh(X))$ with $\dim\al=1$ or 3. See Remark \ref{co7rem11} below for the reasons why. For $\al\in C(\coh(X))_\pe$, in Assumption \ref{co5ass1}(f)(i) define an open substack $\dM_\al\subseteq\M_\al$ by
\e
\dM_\al=\M_\al\sm\supp h^3\bigl(\De_{\M_\al}^*(\cExt_{\al,\al}^\bu)\bigr).
\label{co7eq62}
\e
Here $h^3\bigl(\De_{\M_\al}^*(\cExt_{\al,\al}^\bu)\bigr)$ is the third cohomology sheaf of $\De_{\M_\al}^*(\cExt_{\al,\al}^\bu)$, a coherent sheaf on $\M_\al$ with fibre $\Ext^3(E,E)$ at $[E]\in\M_\al$, and $\supp h^3(\cdots)$ is its support, a Zariski closed subset of $\M_\al$. Thus $\dM_\al$ is the open substack of $[E]$ in $\M_\al$ with $\Ext^3(E,E)=0$. This condition is invariant under the $[*/\bG_m]$-action on $\M_\al$, so there is a unique open substack $\dM_\al^\pl\subseteq\M_\al^\pl$ with~$\dM_\al=(\Pi_\al^\pl)^{-1}(\dM_\al^\pl)$.

Let $\bs\dM_\al\subseteq\bs\M_\al$ and $\bs\dM_\al^\pl\subseteq\bs\M_\al^\pl$ be the derived open substacks corresponding to $\dM_\al\subseteq\M_\al$ and $\dM_\al^\pl\subseteq\M_\al^\pl$. Define $\bs\dM_\al^\red=\bs\dM_\al$, $\bs\dM_\al^\rpl=\bs\dM_\al^\pl$, $\bs j_\al=\bs\id=\bs j_\al^\pl$, $\bs{\dot\Pi}{}_\al^\rpl=\bs{\dot\Pi}{}_\al^\pl$, and $U_\al=o_\al=0$. Then $\bs\dM_\al^\red,\bs\dM_\al^\rpl$ are locally finitely presented as in Definition \ref{co7def10}. Also \eq{co5eq4} is trivially Cartesian, equation \eq{co5eq6} holds as $\bs j_\al,\bs j_\al^\pl$ are identities and $U_\al=0$, giving (iii), and (iv) is trivial as $o_\al=o_\be=o_{\al+\be}=0$. This proves all of Assumption \ref{co5ass1}(f) except that $\bs\dM_\al^\red,\bs\dM_\al^\rpl$ are quasi-smooth. To show this, by Remark \ref{co2rem6}(e) it suffices to check that $i^*(\bL_{\bs\dM_\al^\red}),i^*(\bL_{\bs\dM_\al^\rpl})$ are perfect in the interval~$[-1,1]$.

As $\bs\dM_\al^\red=\bs\dM_\al\subseteq\bs\M_\al$ and $m=3$, equation \eq{co7eq48} shows that $i^*(\bL_{\bs\dM_\al^\red})$ is perfect in $[-2,1]$. Thus the argument at the end of Definition \ref{co7def16} shows that $i^*(\bL_{\bs\dM_\al^\red})$ is perfect in $[-1,1]$ if $h^2(i^*(\bL_{\bs\dM_\al^\red})^\vee)=0$. But equations \eq{co7eq38}, \eq{co7eq43} and \eq{co7eq62} imply that
\begin{equation*}
h^2\bigl(i^*(\bL_{\bs\dM_\al^\red})^\vee\bigr)\cong h^3\bigl(\De_{\M_\al}^*(\cExt_{\al,\al}^\bu)\bigr)\vert_{\dM_\al}=0.
\end{equation*}
Thus $i^*(\bL_{\bs\dM_\al^\red})$ is perfect in $[-1,1]$. By \eq{co7eq49} restricted to $\dM_\al,\dM_\al^\pl$, $i^*(\bL_{\bs\dM_\al^\rpl})$ is also perfect in $[-1,1]$. Hence $\bs\dM_\al^\red,\bs\dM_\al^\rpl$ are quasi-smooth, completing Assumption \ref{co5ass1}(f).

For Assumption \ref{co5ass2}(e), let $(\tau,T,\le)$ lie in $\sS$ in \eq{co7eq51}, so that  $(\tau,T,\le)$ is $(\tau^\om,G,\le)$ or $(\mu^\om,M,\le)$ for some $\om\in U\subseteq\Kah(X)$, let $\al\in C(\coh(X))_\pe$, and suppose $[E]\in\M_\al^\ss(\tau)$. We will show that $\Ext^3(E,E)=0$, so that $[E]\in\dM_\al^\pl$, and thus $\M_\al^\ss(\tau)\subseteq\dM_\al^\pl$, as we have to prove.

By Serre duality $\Ext^3(E,E)\cong\Hom(E,E\ot K_X)^*$. Suppose for a contradiction that $\Ext^3(E,E)\ne 0$. Then $\Hom(E,E\ot K_X)\ne 0$, so we can choose a nonzero morphism $s:E\ra E\ot K_X$. Let $E'=\Im s$, with $E'\ne 0$ as $s\ne 0$. Then $E\twoheadrightarrow E'$ is a nonzero quotient object, so as $E$ is $\tau$-semistable we have
\e
\tau(\lb E\rb)\le \tau(\lb E'\rb).
\label{co7eq63}
\e
Also $E'\subseteq E\ot K_X$, so $E'\ot K_X^{-1}\subseteq E$ is a nonzero subobject. This implies that $\dim E'=\dim E$ as $E$ is pure, and as $E$ is $\tau$-semistable we have
\e
\tau(\lb E'\ot K_X^{-1}\rb)\le\tau(\lb E\rb).
\label{co7eq64}
\e

Write the Hilbert polynomial of $E'$ in \eq{co7eq5} as 
\begin{equation*}
P_{\lb E'\rb}^\om(t)=a_dt^d+a_{d-1}t^{d-1}+\cdots,	
\end{equation*}
where $d=\dim E'=\dim\al>0$. Then calculation shows that
\begin{equation*}
P_{\lb E'\ot K_X^{-1}\rb}^\om(t)=a_dt^d+\Bigl(a_{d-1}+\int_X\ch_{3-d}(\lb E'\rb)\cup c_1(K^{-1}_X)\cup\om^{d-1} \Bigr)t^{d-1}+\cdots.	
\end{equation*}
As $X$ is Fano, $K^{-1}_X$ is ample, so $c_1(K^{-1}_X),\om\in\Kah(X)$, and since $E'\ne 0$ with $\dim E'=d$ we see that $\int_X\ch_{3-d}(\lb E'\rb)\cup c_1(K^{-1}_X)\cup\om^{d-1}>0$. The last two equations then imply that 
\e
\tau(\lb E'\ot K_X^{-1}\rb)>\tau(\lb E'\rb),
\label{co7eq65}
\e
for both $\tau=\tau^\om$ in Definition \ref{co7def3} and $\tau=\mu^\om$ in Definition \ref{co7def4}. But \eq{co7eq65} contradicts \eq{co7eq63}--\eq{co7eq64}. This proves Assumption~\ref{co5ass2}(e).
\end{dfn}

\begin{rem}
\label{co7rem11}
When $X$ is a Fano 3-fold, in \eq{co7eq47} we set $C(\coh(X))_\pe$ to be the subset of $\al\in C(\coh(X))$ with $\dim\al=1$ or 3, that is, we exclude $\al\in C(\coh(X))$ with $\dim\al=0$ or 2. This is for two different reasons.

As in Remark \ref{co7rem5}(a), we excluded $\al$ with $\dim\al=2$ as the author cannot currently prove Condition \ref{co7cond1} holds for $K=\{\om\}$ in this case, and this is needed to verify Assumptions \ref{co5ass2}(b),(g),(h). But it is likely that this gap could be filled, and $\al$ with $\dim\al=2$ could be included.

We excluded $\al$ with $\dim\al=0$ because the proof in Definition \ref{co7def17} that if $0\ne E\in\coh(X)$ is $\tau$-semistable then $\Ext^3(E,E)=0$ requires $\dim E>0$. If $\dim E=0$ then $K_X\vert\supp E$ is (non-canonically) trivial, so $\Ext^3(E,E)\cong\Hom(E,E)^*$, and $\Ext^3(E,E)\ne 0$. Thus if $\dim\al=0$ then \eq{co7eq62} would give $\dM_\al=\dM_\al^\pl=\es$, and $\M_\al^\ss(\tau)\subseteq\dM_\al^\pl$ in Assumption \ref{co5ass2}(e) would be false.
\end{rem}

\subsubsection{The $G$-equivariant case and Assumption \ref{co5ass4}(c)}
\label{co735}

In \S\ref{co724} we explained how to extend \S\ref{co721}--\S\ref{co723} to the $G$-equivariant case when an algebraic $\C$-group $G$ acts on $X$ satisfying Definition \ref{co7def12}(i)--(ii), and verified Assumptions \ref{co4ass2} and \ref{co5ass4}(a),(b),(d). Now that we have defined the data of Assumption \ref{co5ass1}(f) in \S\ref{co731}--\S\ref{co734}, we also verify Assumption \ref{co5ass4}(c). In fact this is immediate: all the constructions of \S\ref{co731}--\S\ref{co734} are manifestly invariant or equivariant under the action of $G$ on $X,\coh(X),\M_\al,\ldots.$ In Definition \ref{co7def16}, for a surface $X$ with $p_g>0$ and $\al\in C(\coh(X))$ with $\rank\al>0$ we defined $U_\al=H^0(K_X)$. We take the $G$-action on $U_\al$ to be that induced by the $G$-action on $X$, and then \eq{co5eq6} holds in $G$-equivariant perfect complexes.

\subsection{The properness conditions Assumption \ref{co5ass2}(g),(h)}
\label{co74}

Assumption \ref{co5ass2}(g) for $\coh(X)$ follows from Theorem \ref{co7thm5}. We now prove Assumption \ref{co5ass2}(h). Let $\baB\subseteq\baA,K(\baA),\baM,\baM^\pl,\baS,\ldots,$ be constructed in \S\ref{co52} from $\B=\A=\coh(X),\ab K(\A),\ab\M,\ab\M^\pl,\ab\sS,\ldots$ defined in \S\ref{co721}--\S\ref{co723}, using a quiver $Q$ with no oriented cycles and some extra data. Use the notation of \S\ref{co52}. 

Let $\om\in\Kah(X)$, so that $(\tau^\om,G,\le)$ and $(\mu^\om,M,\le)$ lie in $\sS$ in \eq{co7eq51}. Let $\la:K(\A)\ra\R$ be a group morphism and $\bs\mu\in\R^{\dot Q_0}$, so that Definition \ref{co5def1} gives weak stability conditions $((\bar\tau^\om)^\la_{\bs\mu},\bar G,\le)$ and $((\bar\mu^\om)^\la_{\bs\mu},\bar M,\le)$. Let $(\al,\bs d)\in C(\baB)_\pe$. To prove Assumption \ref{co5ass2}(h) we must show that:
\begin{itemize}
\setlength{\itemsep}{0pt}
\setlength{\parsep}{0pt}
\item[(A)] If $\baM_{(\al,\bs d)}^\rst((\bar\tau^\om)^\la_{\bs\mu})=\baM_{(\al,\bs d)}^\ss((\bar\tau^\om)^\la_{\bs\mu})$ then $\baM_{(\al,\bs d)}^\ss((\bar\tau^\om)^\la_{\bs\mu})$ is a proper algebraic space.
\item[(B)] If $\baM_{(\al,\bs d)}^\rst((\bar\mu^\om)^\la_{\bs\mu})=\baM_{(\al,\bs d)}^\ss((\bar\mu^\om)^\la_{\bs\mu})$ then $\baM_{(\al,\bs d)}^\ss((\bar\mu^\om)^\la_{\bs\mu})$ is a proper algebraic space.
\end{itemize}

Observe that if $\al=0$ then the moduli stacks are examples of moduli stacks of quiver representations in Chapter \ref{co6}, so (A),(B) follow from Assumption \ref{co5ass2}(g) for $\A=\B=\modCQ$, as proved in \S\ref{co643}--\ref{co644}. Similarly, if $\bs d=0$ then $\baM_{(\al,0)}^{\text{st/ss}}((\bar\tau^\om)^\la_{\bs\mu})=\M^{\text{st/ss}}_\al(\tau^\om)$ and $\baM_{(\al,0)}^{\text{st/ss}}((\bar\mu^\om)^\la_{\bs\mu})=\M^{\text{st/ss}}_\al(\mu^\om)$, so (A),(B) follow from Assumption \ref{co5ass2}(g) for $\A=\B=\coh(X)$, proved above. So we suppose that $\al\ne 0\ne\bs d$. By \eq{co5eq16} this implies that $\al\in C(\B)_\pe$, so $\dim\al=0,1$ or $m$ by \eq{co7eq46}--\eq{co7eq47}. We will extend parts of \S\ref{co71} from $\A=\coh(X)$ to $\baA$. Here is the analogue of Proposition~\ref{co7prop1}:

\begin{prop}
\label{co7prop12}
$\baA$ is of compact type in the sense of\/ {\rm\S\ref{co332}}. Hence $\baM_{(\al,\bs d)}$ is locally of finite type, with affine diagonal, and satisfies Theorem\/~{\rm\ref{co3thm4}(i)--(iii)}.
\end{prop}

\begin{proof}
 This is a straightforward combination of the proofs of Proposition \ref{co6prop2} for $\modCQ$ and Proposition \ref{co7prop1} for $\coh(X)$. The category $\ti{\bar{\mathcal A}}$ is defined as for $\baA$ in Definition \ref{co5def1}, but in objects $(E,\bs V,\bs\rho)$ we take $E\in\qcoh(X)$ and allow $V_v$ to be infinite-dimensional. Then following Propositions \ref{co6prop2} and \ref{co7prop1}, we can show that $\baA$ is of compact type.
\end{proof} 
 
We explained in Definition \ref{co5def1} that $(\bar\tau^\om)^\la_{\bs\mu}$- and $(\bar\mu^\om)^\la_{\bs\mu}$-(semi)stability are open conditions on $(E,\bs V,\bs\rho)$ in $\baA$, so we have open substacks $\baM_{(\al,\bs d)}^\rst((\bar\tau^\om)^\la_{\bs\mu})\subseteq\baM_{(\al,\bs d)}^\ss((\bar\tau^\om)^\la_{\bs\mu})\subseteq\baM_{(\al,\bs d)}^\pl$ and $\baM_{(\al,\bs d)}^\rst((\bar\mu^\om)^\la_{\bs\mu})\subseteq\baM_{(\al,\bs d)}^\ss((\bar\mu^\om)^\la_{\bs\mu})\subseteq\baM_{(\al,\bs d)}^\pl$, the analogue of Proposition \ref{co7prop7}. Here is the analogue of Proposition~\ref{co7prop8}:

\begin{prop}
\label{co7prop13}
$\baM_{(\al,\bs d)}^\ss((\bar\tau^\om)^\la_{\bs\mu})$ and\/ $\baM_{(\al,\bs d)}^\ss((\bar\mu^\om)^\la_{\bs\mu})$ are of finite type.
\end{prop}

\begin{proof}
The morphism $\Pi_{\M_\al^\pl}:\baM_{(\al,\bs d)}^\pl\ra\M_\al^\pl$ in \S\ref{co52} maps $\baM_{(\al,\bs d)}^\ss((\bar\tau^\om)^\la_{\bs\mu})\ra\M^\ss_\al(\tau^\om)$ and $\baM_{(\al,\bs d)}^\ss((\bar\mu^\om)^\la_{\bs\mu})\ra\M^\ss_\al(\mu^\om)$. Here $\M^\ss_\al(\tau^\om),\M^\ss_\al(\mu^\om)$ are of finite type by Proposition \ref{co7prop8}, as $\dim\al=0,1$ or $m$. Also the morphism $\Pi_{\M_\al^\pl}$ is of finite type, since its fibres are quiver moduli stacks $\M_{\bs d}=[\bA^N/\GL_{\bs d}]$ as in \S\ref{co62}. The proposition follows.	
\end{proof}

Here is an analogue of Theorem~\ref{co7thm4}:

\begin{prop}
\label{co7prop14}
$(\bar\Pi_{(\al,\bs d)}^\pl)^{-1}(\baM_{(\al,\bs d)}^\ss((\bar\tau^\om)^\la_{\bs\mu}))$ and\/ $(\bar\Pi_{(\al,\bs d)}^\pl)^{-1}(\baM_{(\al,\bs d)}^\ss((\bar\mu^\om)^\la_{\bs\mu}))$ satisfy the valuative criterion for universal closedness in Definition\/ {\rm\ref{co3def8}(a)} with respect to any discrete valuation ring $R$ essentially of finite type over\/~$\C$.
\end{prop}

\begin{proof} We begin with two lemmas:

\begin{lem}
\label{co7lem2}
$\Pi_{\M_\al}^{-1}\ci(\Pi_\al^\pl)^{-1}(\M_\al^\ss(\tau^\om))$ and\/ $\Pi_{\M_\al}^{-1}\ci(\Pi_\al^\pl)^{-1}(\M_\al^\ss(\mu^\om))$ satisfy the valuative criterion for universal closedness with respect to any discrete valuation ring $R$ essentially of finite type over~$\C$.
\end{lem}

\begin{proof} $(\Pi_\al^\pl)^{-1}(\M_\al^\ss(\tau^\om))$ and $(\Pi_\al^\pl)^{-1}(\M_\al^\ss(\mu^\om))$ satisfy the valuative criterion w.r.t.\ $R$ by Theorem \ref{co7thm4}, as $\dim\al=0,1$ or $m$. Also the morphism $\Pi_{\M_\al}:\baM_{(\al,\bs d)}\ra\M_\al$ satisfies the valuative criterion w.r.t.\ $R$, as its fibres are moduli stacks $\M_{\bs d}$ of a quiver $Q$ without oriented cycles, which satisfy the valuative criterion w.r.t.\ $R$ by Theorem \ref{co3thm5} and Proposition \ref{co6prop2}. The lemma follows.
\end{proof}

The next lemma is based on Example \ref{co3ex6} and Definition~\ref{co6def9}.

\begin{lem}
\label{co7lem3}
There exist pseudo-$\Th$-stratifications of\/ $\Pi_{\M_\al}^{-1}\ci(\Pi_\al^\pl)^{-1}(\M_\al^\ss(\tau^\om))$ and\/ $\Pi_{\M_\al}^{-1}\ci(\Pi_\al^\pl)^{-1}(\M_\al^\ss(\mu^\om)),$ satisfying the descending chain condition, with semistable loci $(\bar\Pi_{(\al,\bs d)}^\pl)^{-1}(\baM_{(\al,\bs d)}^\ss((\bar\tau^\om)^\la_{\bs\mu}))$ and\/~$(\bar\Pi_{(\al,\bs d)}^\pl)^{-1}(\baM_{(\al,\bs d)}^\ss((\bar\mu^\om)^\la_{\bs\mu}))$.
\end{lem}

\begin{proof}
We give the proof for $(\bar\tau^\om)^\la_{\bs\mu}$. The proof for $(\bar\mu^\om)^\la_{\bs\mu}$ is the same, replacing $(\bar\tau^\om)^\la_{\bs\mu},\tau^\om$ by $(\bar\mu^\om)^\la_{\bs\mu},\mu^\om$ throughout.

As in Definition \ref{co5def1}, $\baA$ is $(\bar\tau^\om)^\la_{\bs\mu}$-artinian, so objects $(E,\bs V,\bs\rho)$ in $\baA$ have $(\bar\tau^\om)^\la_{\bs\mu}$-Harder--Narasimhan filtrations by Theorem \ref{co3thm1}. For a $\C$-point $[E,\bs V,\bs\rho]$ in $\Pi_{\M_\al}^{-1}\ci(\Pi_\al^\pl)^{-1}(\M_\al^\ss(\tau^\om))$, $E$ is $\tau^\om$-semistable with $\lb E\rb=\al$. Let the $(\bar\tau^\om)^\la_{\bs\mu}$-Harder--Narasimhan filtration of $(E,\bs V,\bs\rho)$ be
\begin{equation*}
0=(E_0,\bs V_0,\bs\rho_0)\subsetneq (E_1,\bs V_1,\bs\rho_1)	\subsetneq\cdots\subsetneq(E_n,\bs V_n,\bs\rho_n)=(E,\bs V,\bs\rho).
\end{equation*}
Write $(F_i,\bs W_i,\bs\si_i)=(E_i,\bs V_i,\bs\rho_i)/(E_{i-1},\bs V_{i-1},\bs\rho_{i-1})$ for $i=1,\ldots,n$, so that
\e
(\bar\tau^\om)^\la_{\bs\mu}(\lb F_1,\bs W_1,\bs\si_1\rb)>\cdots>(\bar\tau^\om)^\la_{\bs\mu}(\lb F_n,\bs W_n,\bs\si_n\rb).	
\label{co7eq66}
\e
Define the {\it Harder--Narasimhan type\/} $\HNT(E,\bs V,\bs\rho)$ in $\coprod_{m\ge 1}C(\baA)$ to be
\begin{equation*}
\HNT(E,\bs V,\bs\rho)=\bigl(\lb F_1,\bs W_1,\bs\si_1\rb,\ldots,\lb F_n,\bs W_n,\bs\si_n\rb\bigr).
\end{equation*}
Define $\HNT_{(\al,\bs d)}(\baA)$ to be the set of $\HNT(E,\bs V,\bs\rho)$ for all $\C$-points $[E,\bs V,\bs\rho]$ in~$\Pi_{\M_\al}^{-1}\ci(\Pi_\al^\pl)^{-1}(\M_\al^\ss(\tau^\om))$.

In the situation above, for some $0\le k\le l\le n$ we have $E_i=0$ for $i\le k$, $E_i\ne 0$ for $i>k$, $E_i\ne E$ for $i<l$, and $E_i=E$ for $i\ge l$. Then \eq{co5eq21} and \eq{co7eq66} imply that $F_{k+1},\ldots,F_l$ are nonzero with $\tau(F_{k+1})\ge\cdots\ge\tau(F_l)$. As $E$ is $\tau^\om$-semistable, this implies that $F_i$ for $k<i\le l$ and $E_i$ for $i>k$ and $E/E_i$ for $i<l$ are all $\tau^\om$-semistable with $\tau^\om(F_i)=\tau^\om(E_i)=\tau^\om(E/E_i)=\tau^\om(\al)$, and thus $\bar\mu^\om(E/E_i)=\bar\mu^\om(\al)$ for~$i<l$.

Now the family of $E$'s is bounded by Proposition \ref{co7prop8}, as $\dim\al=0,1$ or $m$, so the family of quotients $E\twoheadrightarrow E/E_i$ is also bounded by Proposition \ref{co7prop3}, as $\bar\mu^\om(E/E_i)\le\bar\mu^\om(\al)$. Thus there are only finitely many possibilities for $\lb E/E_i\rb$ in $K(\A)$, and so only finitely many possibilities for $\lb F_i\rb=\lb E/E_{i-1}\rb-\lb E/E_i\rb$. There are also only finitely many possibilities for $\bdim(\bs W_i,\bs\si_i)$, as $\bdim(\bs W_i,\bs\si_i)\le\bs d$. Therefore $\HNT_{(\al,\bs d)}(\baA)$ is a finite set.

Let $\ep>0$ be small, to be chosen later. In a similar way to Definition \ref{co6def9}, if $(\bs\be,\bs e)=((\bs\be_1,\bs e_1),\ldots,(\bs\be_n,\bs e_n))\in\HNT(\al,\bs d)$, define $x_0,\ldots,x_n\in\R^2$ by
\e
x_j=\ts\bigl(\sum_{i=1}^j(\rk\be_j+\ep^2\sum_{v\in\dot Q_0}e_{v,i}),\sum_{i=1}^j(\la(\be_j)+\sum_{v\in \dot Q_0}(\mu_v-\ep)e_{v,i})\bigr).
\label{co7eq67}
\e
Define the {\it Harder--Narasimhan polygon\/} $\HNP(\bs\be,\bs e)$ to be the union of the line segments in $\R^2$ joining $x_{j-1}$ to $x_j$ for $j=1,\ldots,n$. It is piecewise-linear with end-points $x_0=(0,0)$ and $x_n=\bigl(\rk\al+\ep^2\sum_{v\in\dot Q_0}d_v,\la(\al)+\sum_{v\in\dot Q_0}(\mu_v-\ep)d_v\bigr)$. 

Compare equations \eq{co5eq21} and \eq{co7eq67}. If $\ep>0$ is small, we see that
\begin{itemize}
\setlength{\itemsep}{0pt}
\setlength{\parsep}{0pt}
\item[(i)] If $\be_j\ne 0$ then $(\bar\tau^\om)^\la_{\bs\mu}(\be_j,\bs e_j)=\bigl(\tau(\al),[\la(\be_j)+\sum_{v\in\dot Q_0}\mu_ve_{v,j}]/\rk\be_j\bigr)$ and the slope of $[x_{j-1},x_j]$ is $[\la(\be_j)+\sum_{v\in\dot Q_0}\mu_ve_{v,j}]/\rk\be_j+O(\ep)$.
\item[(ii)] If $\be_j=0$ and $\sum_{v\in \dot Q_0}\mu_ve_{v,j}>0$ then $(\bar\tau^\om)^\la_{\bs\mu}(\be_j,\bs e_j)=\bigl(\iy,[\sum_{v\in\dot Q_0}\mu_ve_{v,j}]\ab/[\sum_{v\in\dot Q_0}e_{v,j}]\bigr)$ and the slope of $[x_{j-1},x_j]$ is $\ep^{-2}\sum_{v\in \dot Q_0}\mu_ve_{v,j}+O(\ep^{-1})$.
\item[(iii)] If $\be_j\!=\!0$ and $\sum_{v\in \dot Q_0}\mu_ve_{v,j}\!<\!0$ then $(\bar\tau^\om)^\la_{\bs\mu}(\be_j,\bs e_j)\!=\!\bigl(-\iy,[\sum_{v\in\dot Q_0}\mu_ve_{v,j}]\ab/[\sum_{v\in\dot Q_0}e_{v,j}]\bigr)$ and the slope of $[x_{j-1},x_j]$ is $\ep^{-2}\sum_{v\in \dot Q_0}\mu_ve_{v,j}+O(\ep^{-1})$.
\item[(iv)] If $\be_j=0$ and $\sum_{v\in \dot Q_0}\mu_ve_{v,j}=0$ then $(\bar\tau^\om)^\la_{\bs\mu}(\be_j,\bs e_j)=(-\iy,0)$ and the slope of $[x_{j-1},x_j]$ is $-\ep^{-1}$.
\end{itemize}
Thus we see that the slopes of $[x_{j-1},x_j]$ have the same order in $\R$ as the $(\bar\tau^\om)^\la_{\bs\mu}(\be_j,\bs e_j)$ do in $\bar T$, provided $\ep$ is sufficiently small. We use $\mu_v-\ep$ rather than $\mu_v$ in \eq{co7eq67} to make case (iv) work.

Choose $\ep$ small enough that this holds for all $(\bs\be,\bs e)$ in $\HNT_{(\al,\bs d)}(\baA)$. Then \eq{co7eq66} implies that the slopes of $[x_{i-1},x_i]$ are decreasing in $i=1,\ldots,n$, so $\HNP(\bs\be,\bs e)$ is a convex polygon. It is the graph $\Ga_{\mathop{\rm hnp}(\bs\be,\bs e)}$ of a continuous, piecewise linear function $\mathop{\rm hnp}(\bs\be,\bs e):[0,\rk\al+\ep^2\sum_{v\in\dot Q_0}d_v]\ra\R$. Define a partial order $\preceq$ on $\HNT_{(\al,\bs d)}(\baA)$ by $(\bs\be,\bs e)\preceq(\bs\be',\bs e')$ if $\HNP(\bs\be,\bs e)$ lies below $\HNP(\bs\be',\bs e')$ in $\R^2$, that is, if $\mathop{\rm hnp}(\bs\be,\bs e)\le\mathop{\rm hnp}(\bs\be',\bs e')$ as functions $[0,\rk\al+\ep^2\sum_{v\in\dot Q_0}d_v]\ra\R$. The minimal element in $\HNT_{(\al,\bs d)}(\baA)$ under $\preceq$ is $((\al,\bs d))$, so we set $0=((\al,\bs d))$ in Definition~\ref{co3def11}.

As in Example \ref{co3ex6} and Definition \ref{co6def9}, there is an upper semicontinuous constructible function $\HNT_{(\al,\bs d)}:\Pi_{\M_\al}^{-1}\ci(\Pi_\al^\pl)^{-1}(\M_\al^\ss(\tau^\om))\ra\HNT_{(\al,\bs d)}(\baA)$ mapping $\C$-points $[E,\bs V,\bs\rho]$ to $\HNT(E,\bs V,\bs\rho)$. We can then construct a pseudo-$\Th$-stratification $(\baM_{\preceq(\bs\be,\bs e)})_{(\bs\be,\bs e)\in\HNT_{(\al,\bs d)}(\baA)}$ of $\Pi_{\M_\al}^{-1}\ci(\Pi_\al^\pl)^{-1}(\M_\al^\ss(\tau^\om))$ with semistable locus $(\bar\Pi_{(\al,\bs d)}^\pl)^{-1}(\baM_{(\al,\bs d)}^\ss((\bar\tau^\om)^\la_{\bs\mu}))$, where as in \eq{co3eq11} and \eq{co6eq14}
\begin{align*}
&\baM_{\preceq(\bs\be,\bs e)}=\\
&\bigl\{[E,\bs V,\bs\rho]\in\Pi_{\M_\al}^{-1}\ci(\Pi_\al^\pl)^{-1}(\M_\al^\ss(\tau^\om)):\HNT_{(\al,\bs d)}([E,\bs V,\bs\rho])\preceq(\bs\be,\bs e)\bigr\}.
\end{align*}
It satisfies the descending chain condition trivially, as $\HNT_{(\al,\bs d)}(\baA)$ is finite. This completes the proof.
\end{proof}

Proposition \ref{co7prop14} now follows from Theorem \ref{co3thm7} and Lemmas \ref{co7lem2}--\ref{co7lem3}.
\end{proof}

Here is the analogue of Proposition \ref{co7prop10} for $((\bar\tau^\om)^\la_{\bs\mu},\bar G,\le)$.

\begin{prop}
\label{co7prop15}
$((\bar\tau^\om)^\la_{\bs\mu},\bar G,\le)$ is an additive weak stability condition on\/ $\baA$.
\end{prop}

\begin{proof}
Fix $(\al,\bs d)\in C(\baA)$ as above, and set $n=\dim\al$, so that $0\le n\le m$. Define $V=\R[t^{-2},t^{-1},1,t,\ldots,t^m]$. For $P(t)\in V$, write $P^{(n)}$ for the coefficient of $t^n$ in $P$. Define a total order $\le$ on $V$ by $P(t)\le Q(t)$ if $P(N)\le Q(N)$ for $N\gg 0$. Generalizing \eq{co7eq13}, define a group morphism $\rho:K(\baA)\ra V$ by
\ea
&\rho(\be,\bs e)=P_\al^{\om,(n)}P^\om_\be(t)
-P_\be^{\om,(n)}P^\om_\be(t)
\nonumber\\
&+t^{-1}\bigl[n!P_\al^{c_1(L),(n)}\bigl(\la(\be)+\ts\sum_{v\in\dot Q_0}\mu_ve_v\bigr)-n!P_\be^{c_1(L),(n)}\bigl(\la(\al)+\ts\sum_{v\in\dot Q_0}\mu_vd_v\bigr)\bigr]
\nonumber\\
&-t^{-2}\bigl[n!P_\al^{c_1(L),(n)}\ts\sum_{v\in\dot Q_0}e_v-n!P_\be^{c_1(L),(n)}\ts\sum_{v\in\dot Q_0}d_v\bigr].
\label{co7eq68}
\ea

Here $L\ra X$ is the ample line bundle used to define $\rk\be$ for $\coh(X)$ in \eq{co7eq53}, so that $c_1(L)\in\Kah(X)$, and $P^\om_\be(t),P^{c_1(L)}_\be(t)$ are the Hilbert polynomials of $\be$ with respect to $\om$ and $c_1(L)$, as in \eq{co7eq5}. Then $\be\mapsto P_\be^{c_1(L),(n)}$ is a group morphism $K(\coh(X))\ra\Q$, and \eq{co7eq53} implies that $\rk\be=n!P_\be^{c_1(L),(n)}$ if $\be\in C(\A)$ with $\dim\be=n=\dim\al$. Using this and \eq{co5eq21}, which involves $\rk\al$, we can show that $\rho(\al,\bs d)=0$, and if $(\be,\bs e),(\al-\be,\bs d-\bs e)\in C(\baA)$ then $\rho(\be,\bs e)>0$ if and only if $(\bar\tau^\om)^\la_{\bs\mu}(\be,\bs e)>(\bar\tau^\om)^\la_{\bs\mu}(\al-\be,\bs d-\bs e)$. As for the term $\mu_v-\ep$ in \eq{co7eq67}, the $t^{-2}$ term in \eq{co7eq68} is to deal with the case $\be=0$ and $\sum_{v\in \dot Q_0}\mu_ve_v=0$ in \eq{co5eq21}. It then follows that $((\bar\tau^\om)^\la_{\bs\mu},\bar G,\le)$ is additive as in the proof of Proposition~\ref{co7prop10}.
\end{proof}

Unfortunately $((\bar\mu^\om)^\la_{\bs\mu},\bar M,\le)$ need not be additive, as $(\mu,M,\le)$ may not be. But when $\baM_{(\al,\bs d)}^\rst((\bar\mu^\om)^\la_{\bs\mu})=\baM_{(\al,\bs d)}^\ss((\bar\mu^\om)^\la_{\bs\mu})$ we can still use the ideas of~\S\ref{co335}.

\begin{prop}
\label{co7prop16}
Suppose $\baM_{(\al,\bs d)}^\rst((\bar\mu^\om)^\la_{\bs\mu})=\baM_{(\al,\bs d)}^\ss((\bar\mu^\om)^\la_{\bs\mu})$. Then we can construct a stability function $\rho:K(\baA)\ra V$ as in Definition\/ {\rm\ref{co3def12}} such that $\rho$-semistability is an open condition on $\baM_{(\al,\bs d)}$ and $\baM_{(\al,\bs d)}^{\ss,\rho}=\baM_{(\al,\bs d)}^\ss((\bar\mu^\om)^\la_{\bs\mu})$. 

Thus\/ $\baM_{(\al,\bs d)}^\ss((\bar\mu^\om)^\la_{\bs\mu})$ is $\Th$-reductive and S-complete with respect to any DVR $R$ essentially of finite type over\/ $\C$ by Theorem\/ {\rm\ref{co3thm8}} and Proposition\/~{\rm\ref{co7prop12}}.
\end{prop}

\begin{proof}
Use the notation $n,L,P^{(n)}$ from the proof of Proposition \ref{co7prop15}. Define $V=\R[t^{n-3/2},t^{n-4/3},1,t,\ldots,t^m]$. Define a total order $\le$ on $V$ by $P(t)\le Q(t)$ if $P(n)\le Q(n)$ for $n\gg 0$. Modifying \eq{co7eq68}, define a group morphism $\rho:K(\baA)\ra V$ by
\begin{align*}
&\rho(\be,\bs e)=P_\al^{\om,(n)}P^\om_\be(t)
-P_\be^{\om,(n)}P^\om_\be(t)\\
&+t^{n-4/3}\bigl[n!P_\al^{c_1(L),(n)}\bigl(\la(\be)\!+\!\ts\sum_{v\in\dot Q_0}\mu_ve_v\bigr)\!-\!n!P_\be^{c_1(L),(n)}\bigl(\la(\al)\!+\!\ts\sum_{v\in\dot Q_0}\mu_vd_v\bigr)\bigr]\\
&-t^{n-3/2}\bigl[n!P_\al^{c_1(L),(n)}\ts\sum_{v\in\dot Q_0}e_v-n!P_\be^{c_1(L),(n)}\ts\sum_{v\in\dot Q_0}d_v\bigr].
\end{align*}

The idea here is that if $\dim\be=\dim\al=n$ then $P^\om_\be(t)=a_nt^n+a_{n-1}t^{n-1}+\cdots+a_0$ with $a_n=P_\be^{\om,(n)}$ and $\mu^\om(\be)=t^n+(a_{n-1}/a_n)t^{n-1}$. The powers $t^{n-4/3},t^{n-3/2}$ lie between $t^{n-1}$ and $t^{n-2}$, so the terms $a_nt^n+a_{n-1}t^{n-1}$ used in $\mu^\om(\be)$ dominate the $t^{n-4/3},t^{n-3/2}$ terms, but the $t^{n-4/3},t^{n-3/2}$ terms dominate the remainder $a_{n-2}t^{n-2}+\cdots+a_0$ in $P^\om_\be(t)$, which is forgotten in~$\mu^\om(\be)$.

We can prove using \eq{co5eq21} that $\rho(\al,\bs d)=0$, and if $(\be,\bs e),(\al-\be,\bs d-\bs e)\in C(\baA)$ with $(\bar\mu^\om)^\la_{\bs\mu}(\be,\bs e)\!\ne\!(\bar\mu^\om)^\la_{\bs\mu}(\al-\be,\bs d-\bs e)$ then $(\bar\mu^\om)^\la_{\bs\mu}(\be,\bs e)\!>\!(\bar\mu^\om)^\la_{\bs\mu}(\al-\be,\bs d-\bs e)$ if and only if $\rho(\be,\bs e)\!>\!0$. Here we need $(\bar\mu^\om)^\la_{\bs\mu}(\be,\bs e)\ne(\bar\mu^\om)^\la_{\bs\mu}(\al-\be,\bs d-\bs e)$ to prevent terms $a_{n-2}t^{n-2}+\cdots+a_0$ in $P^\om_\be(t)$ influencing the result. It is then easy to see that $\rho$-semistability is an open condition and $\baM_{(\al,\bs d)}^{\ss,\rho}=\baM_{(\al,\bs d)}^\ss((\bar\mu^\om)^\la_{\bs\mu})$, and the proposition follows.
\end{proof}

Finally we show that Assumption \ref{co5ass2}(h) holds for $\coh(X)$. We must prove (A),(B) at the beginning of \S\ref{co74}. Part (A) follows from Proposition \ref{co3prop6}, Theorem \ref{co3thm4}, Corollary \ref{co3cor3}, and Propositions \ref{co7prop12}--\ref{co7prop14}, and \ref{co7prop15}. Part (B) follows from  Proposition \ref{co3prop6}, Theorem \ref{co3thm4} and Propositions \ref{co7prop12}--\ref{co7prop14} and \ref{co7prop16}. This completes the proof.

\begin{rem}
\label{co7rem12}
Objects of our categories $\baB\subseteq\baA$ for $\A=\coh(X)$ can be regarded as examples of {\it twisted quiver sheaves\/} in the sense of \'Alvarez-C\'onsul and Garc\'\i a-Prada \cite{AlGa} and Schmitt \cite{Schm2,Schm3,Schm4,Schm5}. Schmitt \cite{Schm4,Schm5} gives analogues of the GIT construction of moduli spaces in $\coh(X)$ in Theorem \ref{co7thm3} for quiver sheaves. As an alternative approach, one could try to use these to prove Assumption \ref{co5ass2}(h) for $\coh(X)$, as in Corollary \ref{co7cor1}, at least if~$\om\in\Kah_\Q(X)$.
\end{rem}

\subsection{The finiteness condition Assumption \ref{co5ass3}(b)}
\label{co75}

The next proposition will be proved in \S\ref{co751}. The proof is based on that of \cite[Prop.~5.14]{Joyc7}. It is related to Proposition \ref{co7prop5}, and in fact if $\om\in\Kah(X)$ then Proposition \ref{co7prop5} for $K=\{\om\}$ follows from Propositions \ref{co7prop4} and~\ref{co7prop17}.

\begin{prop}
\label{co7prop17}
Work in the situation of\/ {\rm\S\ref{co72}}. Then Assumption\/ {\rm\ref{co5ass3}(b)} holds for all\/ $\al\in C(\coh(X))_\pe$ with\/ $\dim\al=0$ or\/ $1$. This is also true if we take $U=\Kah(X)$ in the definition of\/ $\sS$ in \eq{co7eq51}.	
\end{prop}

The next proposition will be proved in \S\ref{co752}. The proof is based on that of \cite[Th.~5.16]{Joyc7} when $m=2$. It uses the Bogomolov inequality, Theorem \ref{co7thm1}, which is why we restrict to $m\ge 2$ and $\dim\al=m$. It is related to Proposition \ref{co7prop6}, and in fact if $\om\in\Kah(X)$ is sufficiently close to $\Kah_\Q(X)$ then Proposition \ref{co7prop6} for $K=\{\om\}$ follows from Propositions \ref{co7prop4} and~\ref{co7prop18}.

Recall that in \eq{co7eq51} we defined $\sS$ to be the set of $(\tau^\om,G,\le),(\mu^\om,M,\le)$ for $U$ in a connected open subset $U\subseteq\Kah(X)$ satisfying conditions Definition \ref{co7def11}(a),(b). The proof of Proposition \ref{co7prop18} uses these conditions.

\begin{prop}
\label{co7prop18}
Work in the situation of\/ {\rm\S\ref{co72},} and suppose $m=\dim X\ge 2$. Then Assumption\/ {\rm\ref{co5ass3}(b)} holds for all\/ $\al\in C(\coh(X))_\pe$ with\/ $\dim\al=m$. 
\end{prop}

\begin{rem}
\label{co7rem13}
As in Remark \ref{co7rem5}(a), we excluded sheaves of dimensions $2,3,\ab\ldots,\ab m-1$ from our theory by the definition \eq{co7eq46}--\eq{co7eq47} of $C(\coh(X))_\pe$. One reason for this is that the author cannot currently prove Assumption \ref{co5ass3}(b) for such $\al$. This is only an issue for sheaves of dimension 2 on Fano and Calabi--Yau 3-folds, and sheaves of dimensions 2 and 3 on Calabi--Yau 4-folds.
\end{rem}

\subsubsection{Proof of Proposition \ref{co7prop17}}
\label{co751}

If $\al\!\in\! C(\coh(X))_\pe$ with $\dim\al\!=\!0$ then $\al$ is determined by $\ch_m(\al)\!\in\! H^{2m}(X,\Q)\ab\cong\Q$, and $m!\int_X\ch_m(\al)$ is a positive integer. Thus the number of ways of splitting $\al=\al_1+\cdots+\al_n$ with $\al_i\in C(\coh(X))$ is bounded by the (finite) number of ways of writing $m!\int_X\ch_m(\al)$ as a sum of positive integers. This implies the finiteness condition in Assumption \ref{co5ass3}(b), and the condition $\al_{i_1}+\cdots+\al_{i_2}\in C(\coh(X))_\pe$ follows from~\eq{co7eq46}--\eq{co7eq47}.

Next let $\al\in C(\coh(X))_\pe$ with $\dim\al=1$, and let $(\tau_t,T_t,\le)_{t\in[0,1]}$ be a continuous family of weak stability conditions on $\A$ with $(\tau_t,T_t,\le)$ in $\sS$ in \eq{co7eq51} for all $t\in[0,1]$. Then for each $t\in[0,1]$ there is $\om_t\in\Kah(X)$ such that $(\tau_t,T_t,\le)$ is $(\tau^{\om_t},G,\le)$ or $(\mu^{\om_t},M,\le)$, and $\om_t$ depends continuously on $t\in[0,1]$. Define $K=\{\om_t:t\in[0,1]\}\subset\Kah(X)$. This is the image of a continuous map $[0,1]\ra\Kah(X)$, and so is compact. 

Fix $\ti\om\in\Kah_\Q(X)$. As in the proof of Proposition \ref{co7prop5} in \S\ref{co7112}, there exists $C>0$ such that if $\be\in C(\coh(X))$ with $\dim\be=1$ and $\om_t\in K$ then
\e
\sign(\bar\mu^{\ti\om}(\be))=\sign(\bar\mu^{\om_t}(\be)),\quad C^{-1}\bmd{\bar\mu^{\ti\om}(\be)}\le \bmd{\bar\mu^{\om_t}(\be)}\le C\bmd{\bar\mu^{\ti\om}(\be)}.
\label{co7eq69}
\e

Suppose as in Assumption \ref{co5ass3}(b) that $n\ge p\ge 1$, $\al_1,\ldots,\al_n\in C(\coh(X))$ and $0=a_0<\cdots<a_p=n$ such that $\al_1+\cdots+\al_n=\al$ and, writing $\be_j=\al_{a_{j-1}+1}+\cdots+\al_{a_j}$ for $j=1,\ldots,p$, then there exist $s,t,u\in[0,1]$ with $\M_{\al_i}^\ss(\tau_s)\ne\es$ for $i=1,\ldots,n$, and $U\bigl(\al_{a_{j-1}+1},\ab\al_{a_{j-1}+2},\ab\ldots,\ab\al_{a_j};\tau_s,\tau_t)\ne 0$ for $j=1,\ldots,p$, and $\tau_u(\be_1)=\cdots=\tau_u(\be_p)$. Then $\tau_u(\be_1)=\cdots=\tau_u(\be_p)$ imply that $\dim\be_j=\dim\al=1$ for all $j$, and $U\bigl(\al_{a_{j-1}+1},\ab\al_{a_{j-1}+2},\ab\ldots,\ab\al_{a_j};\tau_s,\tau_t)\ne 0$ implies that $\dim\al_i=\dim\be_j=1$ for all~$i=1,\ldots,n$.

Combining $\be_1+\cdots+\be_p=\al$ and $\tau_u(\be_1)=\cdots=\tau_u(\be_p)$ gives $\tau_u(\be_j)=\tau_u(\al)$ for $j=1,\ldots,p$, so $\bar\mu^{\om_u}(\be_j)=\bar\mu^{\om_u}(\al)$. Hence by \eq{co7eq69} we have
\e
\bmd{\bar\mu^{\ti\om}(\be_j)}\le C\bmd{\bar\mu^{\om_u}(\be_j)}=C\bmd{\bar\mu^{\om_u}(\al)}\le C^2\bmd{\bar\mu^{\ti\om}(\al)}.
\label{co7eq70}
\e
Since $U\bigl(\al_{a_{j-1}+1},\ab\al_{a_{j-1}+2},\ab\ldots,\ab\al_{a_j};\tau_s,\tau_t)\ne 0$, Proposition \ref{co3prop2} gives $k_j,l_j=a_{j-1}+1,\ldots,a_j$ such that $\tau_s(\al_{k_j})\le\tau_s(\al_i)\le\tau_s(\al_{l_j})$ for all $i=a_{j-1}+1,\ldots,a_j$ and $\tau_t(\al_{k_j})\!\ge\!\tau_t(\be_j)\!\ge\!\tau_t(\al_{l_j})$. These imply that $\bar\mu^{\om_s}(\al_{k_j})\!\le\!\bar\mu^{\om_s}(\al_i)\!\le\!\bar\mu^{\om_s}(\al_{l_j})$ for all $i=a_{j-1}+1,\ldots,a_j$ and $\bar\mu^{\om_t}(\al_{k_j})\ge\bar\mu^{\om_t}(\be_j)\ge\bar\mu^{\om_t}(\al_{l_j})$. Hence
\e
\begin{split}
\bar\mu^{\om_s}(\al_i)&\le \bar\mu^{\om_s}(\al_{l_j})\le C\max\bigl(\bar\mu^{\ti\om}(\al_{l_j}),0\bigr)\le C^2\max\bigl(\bar\mu^{\om_t}(\al_{l_j}),0\bigr)\\
&\le C^2\max\bigl(\bar\mu^{\om_t}(\be_j),0\bigr)\le C^3\max\bigl(\bar\mu^{\ti\om}(\be_j),0\bigr)\le C^5\bmd{\bar\mu^{\ti\om}(\al)},
\end{split}
\label{co7eq71}
\e
using \eq{co7eq69}--\eq{co7eq70}. As $\M_{\al_i}^\ss(\tau_s)\ne\es$ we may choose $[E_i]\in\M_{\al_i}^\ss(\tau_s)$ for $i=1,\ldots,n$. Then $E=E_1\op\cdots\op E_n$ has $\lb E\rb=\al_1+\cdots+\al_n=\al$. As $E_i$ is $\tau_s$-semistable it is pure of dimension $\dim\al_i=1$, so $E$ is pure. Also
\begin{equation*}
\mu_{\max}^{\om_s}(E)=\max\bigl(\bar\mu^{\om_s}(\al_1),\ldots,\bar\mu^{\om_s}(\al_n)\bigr)\le C^5\bmd{\bar\mu^{\ti\om}(\al)}
\end{equation*}
by \eq{co7eq71}, as $E_1,\ldots,E_n$ are $\tau_s$-semistable. Thus \eq{co7eq69} implies that
\begin{equation*}
\mu_{\max}^{\ti\om}(E)\le C^6\bmd{\bar\mu^{\ti\om}(\al)}=:\mu_0.
\end{equation*}
 
Now Proposition \ref{co7prop4} implies that the family of pure $E\in\coh(X)$ with $\lb E\rb=\al$ and $\mu_{\max}^{\ti\om}(E)\le\mu_0$ is bounded. But in a bounded family of sheaves $E$, the decompositions $E=E_1\op\cdots\op E_n$ for all $n\ge 1$ can realize only finitely many sets of classes $(\lb E_1\rb,\ldots,\lb E_n\rb)$ in $C(\coh(X))^n$. Thus there are only finitely many possibilities for $n\ge 1$ and $(\al_1,\ldots,\al_n)$. This implies the finiteness condition in Assumption \ref{co5ass3}(b), and the condition $\al_{i_1}+\cdots+\al_{i_2}\in C(\coh(X))_\pe$ follows from \eq{co7eq46}--\eq{co7eq47} and $\dim\al_i=1$. This completes the proof.

\subsubsection{Proof of Proposition \ref{co7prop18}}
\label{co752}

Throughout the proof we fix a continuous family of weak stability conditions $(\tau_t,T_t,\le)_{t\in[0,1]}$ on $\A$ with $(\tau_t,T_t,\le)\in\sS$ in \eq{co7eq51} for all $t\in[0,1]$, which we use in Assumption \ref{co5ass3}(b). Then for each $t\in[0,1]$ there is $\om_t\in\Kah(X)$ such that $(\tau_t,T_t,\le)$ is $(\tau^{\om_t},G,\le)$ or $(\mu^{\om_t},M,\le)$, and $\om_t$ depends continuously on $t\in[0,1]$. We will use the following notation. Fix a positive definite inner product $\an{\,,\,}$ on $H^{1,1}(X,\R)$, and let $\nm{\,.\,}$ be the induced norm. By Definition \ref{co7def11}(b), if $s,u\in[0,1]$ and $0\ne\eta\in H^{1,1}(X,\R)$ with $\int_X\eta\cup \om_u^{m-1}=0$ then $\int_X\eta^2\cup\om_s^{m-2}<0$, as $\om_s,\om_u\in U\subseteq\Kah(X)$. Since $\om_s$ depends continuously on $s$ in $[0,1]$, which is compact, there exists $C>0$ depending on $(\om_s)_{s\in[0,1]}$ such that for all such $s,u\in[0,1]$ and $\eta\in H^{1,1}(X,\R)$ with $\int_X\eta\cup \om_u^{m-1}=0$ we have
\e
\int_X\eta^2\cup\om_s^{m-2}\le -C\nm{\eta}^2.
\label{co7eq72}
\e

We will first prove the following inductive hypothesis $(*)_r$ by induction on $r=1,2,\ldots.$ Then we will deduce Proposition \ref{co7prop18} using Proposition \ref{co7prop6}.
\begin{itemize}
\setlength{\itemsep}{0pt}
\setlength{\parsep}{0pt}
\item[$(*)_r$] Let $\bar\al\in H^2(X,\Z)$ and $D\ge 0$ be given. Then there are most finitely many sets of data $n\ge 1$ and $\bar\al_1,\ldots,\bar\al_n\in H^2(X,\Z)$ with $\bar\al=\bar\al_1+\cdots+\bar\al_n$, such that there exist $\al\in C(\coh(X))$ with $\dim\al=m$, $\rank\al\le r$ and $c_1(\al)=\bar\al$, and $\al_1,\ldots,\al_n\in C(\coh(X))$ with $\al_1+\cdots+\al_n=\al$ and $c_1(\al_i)=\bar\al_i$ for $i=1,\ldots,n$, and $0=a_0<\cdots<a_p=n$ for $1\le p\le n$, and $s,t,u\in[0,1]$, such that $\int_X\De(\al)\cup\om_s^{m-2}\le D$, and writing $\be_j=\al_{a_{j-1}+1}+\cdots+\al_{a_j}$ for $j=1,\ldots,p$, then $\M_{\al_i}^\ss(\tau_s)\ne\es$ for $i=1,\ldots,n$, and for each $j=1,\ldots,p$ the sum \eq{co3eq3} for $U\bigl(\al_{a_{j-1}+1},\ab\al_{a_{j-1}+2},\ab\ldots,\ab\al_{a_j};\tau_s,\tau_t)$ contains nonzero terms, and~$\tau_u(\be_1)=\cdots=\tau_u(\be_p)$. 

Furthermore, if such data $\al,\al_i,a_j,\be_j,s,t,u$ exist (though not necessarily requiring that $\int_X\De(\al)\cup\om_s^{m-2}\le D$) then~$\int_X\De(\al)\cup\om_s^{m-2}\ge 0$.
\end{itemize}

Observe that if $\al,n,p,\al_i,a_j,\be_j,\ab s,t,u$ are as in $(*)_r$ then $\tau_u(\be_1)=\cdots=\tau_u(\be_p)$ implies that $\dim\be_j=\dim\al=m$ for all $j$, and $U\bigl(\al_{a_{j-1}+1},\ab\al_{a_{j-1}+2},\ab\ldots,\ab\al_{a_j};\tau_s,\tau_t)\ne 0$ implies that $\dim\al_i=\dim\be_j=m$ for all $i=a_{j-1}+1,\ldots,a_j$, and hence $\dim\al_i=m$ for $i=1,\ldots,n$. Thus $\rank\al_i,\rank\be_j>0$. As $\rank\al=\rank\al_1+\cdots+\rank\al_n$ we see that $p\le n\le\rank\al\le r$.
 
For the first step $r=1$, when $\rank\al=1$ the only possibility in $(*)_1$ is $n=p=1$ and $\al_1=\be_1=\al$, so $\bar\al_1=\bar\al$, proving the finiteness part. Given data $\al,\ab n,\ab p,\al_i,\ab a_j,\ab\be_j,\ab s,t,u$ satisfying $(*)_1$ with $\rank\al=1$ we have $\M_{\al_1}^\ss(\tau_s)\ne\es$ with $\al_1=\al$ and $\tau_s=\tau^{\om_s}$ or $\mu^{\om_s}$ with $\M_\al^\ss(\tau^{\om_s})\subseteq\M_\al^\ss(\mu^{\om_s})$ as $(\mu^{\om_s},M,\le)$ dominates $(\tau^{\om_s},G,\le)$. Thus $\M_\al^\ss(\mu^{\om_s})\ne\es$, so Theorem \ref{co7thm1} and $\dim\al=m\ge 2$ implies that $\int_X\De(\al)\cup\om_s^{m-2}\ge 0$. Hence $(*)_1$ holds.

We divide the inductive step into three cases, proved in three lemmas:

\begin{lem}
\label{co7lem4}
Suppose $(*)_r$ holds for $r\ge 1,$ and\/ $\bar\al\in H^2(X,\Z)$ and\/ $D>0$ are given. Then there are most finitely many sets of data $n\ge 1$ and\/ $\bar\al_1,\ldots,\bar\al_n\in H^2(X,\Z)$ with\/ $\bar\al=\bar\al_1+\cdots+\bar\al_n$ such that there exist data $\al,n,p,\al_i,a_j,\be_j,\ab s,t,u$ as in $(*)_{r+1}$ with\/ $\rank\al=r+1$ and\/ $p>1,$ and if such data exist then~$\int_X\De(\al)\cup\om_s^{m-2}\ge 0$.
\end{lem}

\begin{proof} Suppose $\al,n,p,\al_i,a_j,\be_j,\ab s,t,u$ are as in $(*)_{r+1}$ with $\rank\al=r+1$ and $p>1$. As $\rank\be_1+\cdots+\rank\be_p=\rank\al=r+1$ with $p>1$ and $\rank\be_j>0$ we see that $\rank\be_j\le r$, so we can apply $(*)_r$ with $\be_j$ in place of $\al$, for $j=1,\ldots,p$. Replacing $\al,n,p,\al_1,\ldots,\al_n$ by $\be_j,a_j-a_{j-1},1,\al_{a_{j-1}+1},\ldots,\ab\al_{a_j}$, we see that the conditions of $(*)_r$ are satisfied, except perhaps for $\int_X\De(\be_j)\cup\om_s^{m-2}\le D$. Thus the last part of $(*)_r$ gives
\e
\int_X\De(\be_j)\cup\om_s^{m-2}\ge 0.
\label{co7eq73}
\e

Define $\ga_j\in C(\coh(X))$ for $1\le j\le p$ and $\de_j\in K(\coh(X))$ for $2\le j\le p$ by
\e
\ga_j=\be_1+\cdots+\be_j\quad\text{and}\quad
\de_j=\rank\ga_{j-1}\cdot \be_j-\rank\be_j\cdot \ga_{j-1},
\label{co7eq74}
\e
so that $\ga_1=\al_1$ and $\ga_p=\al$. Then calculation shows that
\e
\frac{c_1(\ga_j)^2}{\rank\ga_j}=\frac{c_1(\ga_{j-1})^2}{\rank\ga_{j-1}}+\frac{c_1(\be_j)^2}{\rank\be_j}-\frac{c_1(\de_j)^2}{\rank\ga_j\cdot\rank\ga_{j-1}\cdot\rank\be_j}.
\label{co7eq75}
\e

Also $\tau_u(\be_1)=\cdots=\tau_u(\be_p)$ implies that $\bar\mu^{\om_u}(\be_1)=\cdots=\bar\mu^{\om_u}(\be_p)=\bar\mu^{\om_u}(\al)$, and hence $\bar\mu^{\om_u}(\ga_1)=\cdots=\bar\mu^{\om_u}(\ga_p)=\bar\mu^{\om_u}(\al)$. Then $\bar\mu^{\om_u}(\be_j)=\bar\mu^{\om_u}(\ga_{j-1})$ and equations \eq{co7eq21} and \eq{co7eq74} give $\int_Xc_1(\de_j)\cup\om_u^{m-1}=0$, so \eq{co7eq72} yields
\e
\int_Xc_1(\de_j)^2\cup\om_s^{m-2}\le -C\nm{\de_j}^2.
\label{co7eq76}
\e

Now
\ea
0&\le\!\sum_{j=1}^p\frac{1}{\rank\be_j}\int_X\De(\be_j)\!\cup\!\om_s^{m-2}\!=\!\sum_{j=1}^p\int_X
\Bigl(\frac{1}{\rank\be_j}c_1(\be_j)^2-2\ch_2(\be_j)\Bigr)\!\cup\!\om_s^{m-2}
\nonumber\\
\begin{split}
&=\int_X\Bigl(\frac{1}{\rank\al}c_1(\al)^2-2\ch_2(\al)\Bigr)\cup\om_s^{m-2}\\
&\qquad+\sum_{j=2}^p\frac{1}{\rank\ga_j\cdot\rank\ga_{j-1}\cdot\rank\be_j}
\int_Xc_1(\de_j)^2\cup\om_s^{m-2}
\end{split}
\label{co7eq77}\\
&\le \frac{1}{\rank\al}\int_X\De(\al)\cup\om_s^{m-2}-\frac{C}{r^3}\sum_{j=2}^p\nm{c_1(\de_j)}^2
\le \frac{1}{r+1}D-\frac{C}{r^3}\sum_{j=2}^p\nm{c_1(\de_j)}^2,
\nonumber
\ea
using \eq{co7eq73} in the first step, \eq{co7eq7}--\eq{co7eq8} in the second and fourth, $\ga_1=\al_1$, $\ga_p=\al$, equation \eq{co7eq75} and induction on $j=2,\ldots,p$ and $\ch_2(\al)=\sum_{j=1}^p\ch_2(\be_j)$ in the third, \eq{co7eq76} and $\rank\ga_j,\rank\be_j\le r$ in the third, and $\rank\al=r+1$ and $\int_X\De(\al)\cup\om_s^{m-2}\le D$ from $(*)_{r+1}$ in the fifth. Thus
\e
\nm{c_1(\de_j)}^2\le r^3C^{-1}D.
\label{co7eq78}
\e

As $c_1(\de_j)$ lies in the discrete lattice $H^2(X,\Z)\cap H^{1,1}(X,\R)$, equation \eq{co7eq78} implies that there are only finitely many possibilities for $c_1(\de_j)$ for $j=2,\ldots,p$. From \eq{co7eq74}, $\ga_p=\al$, $c_1(\al)=\bar\al$, and $\rank\be_j,\rank\ga_j\le r$, $p\le r+1$ we deduce that there are finitely many possibilities for $p$ and $c_1(\be_1),\ldots,c_1(\be_p)$. Also \eq{co7eq77} implies that for each $j=1,\ldots,p$ we have $\int_X\De(\be_j)\cup\om_s^{m-2}\le D$. Hence $(*)_r$ with $c_1(\be_j),\be_j,a_j-a_{j-1},1,\al_{a_{j-1}+1},\ldots,\ab\al_{a_j}$ in place of $\bar\al,\al,n,p,\al_1,\ldots,\al_n$ implies that for each of the finitely many choices of $p,c_1(\be_1),\ldots,c_1(\be_p)$ and $j=1,\ldots,p$ there are only finitely many choices for $a_j-a_{j-1}$ and $\bar\al_{a_{j-1}+1}=c_1(\al_{a_{j-1}+1}),\ldots,\bar\al_{a_j}=c_1(\al_{a_j})$. Thus overall there are only finitely many choices for $n\ge 1$ and $\bar\al_1,\ldots,\bar\al_n\in H^2(X,\Z)$ satisfying the given conditions. Finally \eq{co7eq77} implies that $\int_X\De(\al)\cup\om_s^{m-2}\ge 0$. This completes the proof.  
\end{proof}

\begin{lem}
\label{co7lem5}
Suppose $(*)_r$ holds for $r\ge 1,$ and\/ $\bar\al\in H^2(X,\Z)$ and\/ $D>0$ are given. Then there are most finitely many sets of data $n\ge 1$ and\/ $\bar\al_1,\ldots,\bar\al_n\in H^2(X,\Z)$ with\/ $\bar\al=\bar\al_1+\cdots+\bar\al_n$ such that there exist data $\al,n,p,\al_i,a_j,\be_j,\ab s,t,u$ as in $(*)_{r+1}$ with\/ $\rank\al=r+1$ and\/ $p=1,$ and the sum \eq{co3eq3} for $U\bigl(\al_1,\ldots,\al_n;\tau_s,\tau_t)$ contains nonzero terms with\/ $l>1$. If such data exist then~$\int_X\De(\al)\cup\om_s^{m-2}\ge 0$.
\end{lem}

\begin{proof} Suppose $\al,n,p=1,\al_i,a_j,\be_1=\al,s,t,u$ are as in the lemma. Choose $l>1$ and $\ga_1,\ldots,\ga_l\in C(\coh(X))$ such that there exist nonzero terms in \eq{co3eq3} for $U\bigl(\al_1,\ldots,\al_n;\tau_s,\tau_t)$ with these values of $l,\ga_1,\ldots,\ga_l$. This is well defined by assumption. From \eq{co3eq3} we have $\tau_t(\ga_1)=\cdots=\tau_t(\ga_l)=\tau_t(\al)$.

There are unique $0=c_0<c_1<\cdots<c_l=n$ such that $\ga_j=\al_{c_{j-1}+1}+\cdots+\al_{c_j}$. Each nonzero term in \eq{co3eq3} for $U\bigl(\al_1,\ldots,\al_n;\tau_s,\tau_t)$ with these values of $l,\ga_1,\ldots,\ga_l$ induces nonzero terms in \eq{co3eq3} for $U\bigl(\al_{c_{j-1}+1},\ab\al_{c_{j-1}+2},\ab\ldots,\ab\al_{c_j};\tau_s,\tau_t)$ for $j=1,\ldots,l$. Therefore the conditions of $(*_{r+1})$ hold for $\al_1,\ldots,\al_n$ with $p=l>1$, $a_j=c_j$, $\be_j=\ga_j$, the given $s,t$, and $u=t$. Thus the lemma follows from Lemma~\ref{co7lem4}.
\end{proof}

\begin{lem}
\label{co7lem6}
Suppose $(*)_r$ holds for $r\ge 1,$ and\/ $\bar\al\in H^2(X,\Z)$ and\/ $D>0$ are given. Then there are most finitely many sets of data $n\ge 1$ and\/ $\bar\al_1,\ldots,\bar\al_n\in H^2(X,\Z)$ with\/ $\bar\al=\bar\al_1+\cdots+\bar\al_n$ such that there exist data $\al,n,p,\al_i,a_j,\be_j,\ab s,t,u$ as in $(*)_{r+1}$ with\/ $\rank\al=r+1$ and\/ $p=1,$ and the sum \eq{co3eq3} for $U\bigl(\al_1,\ldots,\al_n;\tau_s,\tau_t)$ contains nonzero terms with\/ $l=1,$ but no nonzero terms with\/ $l>1$. If such data exist then~$\int_X\De(\al)\cup\om_s^{m-2}\ge 0$.
\end{lem}

\begin{proof} Suppose $\al,n,p=1,\al_i,a_j,\be_1=\al,s,t,u$ are as in the lemma, and pick a nonzero term in \eq{co3eq3} for $U\bigl(\al_1,\ldots,\al_n;\tau_s,\tau_t)$ with $l=1$. To avoid notational conflict we relabel $m,a_i,\be_i$ in \eq{co3eq3} by $q,d_j,\de_j$. Then there exist $1\le q\le n$ and  $0=d_0<d_1<\cdots<d_q=n$ such that defining $\de_j=\al_{d_{j-1}+1}+\cdots+\al_{d_j}$ for $j=1,\ldots,q$, we have $\tau_s(\al_i)=\tau_s(\de_j)$ for all $1\le j\le q$ and $d_{j-1}<i\le d_j$, and $S(\de_1,\ldots,\de_q;\tau_s,\tau_t)\ne 0$, giving a nonzero term in \eq{co3eq3}. Thus Definition \ref{co3def4} implies that for all $j=1,\ldots,q-1$ we have either
\begin{itemize}
\setlength{\itemsep}{0pt}
\setlength{\parsep}{0pt}
\item[(a)] $\tau_s(\de_j)\le\tau_s(\de_{j+1})$ and $\tau_t(\de_1+\cdots+\de_j)>\tau_t(\de_{j+1}+\cdots+\de_q)$, or
\item[(b)] $\tau_s(\de_j)>\tau_s(\de_{j+1})$ and $\tau_t(\de_1+\cdots+\de_j)\le\tau_t(\de_{j+1}+\cdots+\de_q)$.
\end{itemize}
Now if $\tau_t(\de_1+\cdots+\de_j)=\tau_t(\de_{j+1}+\cdots+\de_q)$ for some $j=1,\ldots,q-1$ then we get a nonzero term in \eq{co3eq3} for $U\bigl(\al_1,\ldots,\al_n;\tau_s,\tau_t)$ with $l=2$ and $\ga_1=\de_1+\cdots+\de_j$, $\ga_2=\de_{j+1}+\cdots+\de_q$, a contradiction. Thus we can replace (b) by
\begin{itemize}
\setlength{\itemsep}{0pt}
\setlength{\parsep}{0pt}
\item[(b$)'$] $\tau_s(\de_j)>\tau_s(\de_{j+1})$ and $\tau_t(\de_1+\cdots+\de_j)<\tau_t(\de_{j+1}+\cdots+\de_q)$.
\end{itemize}

Suppose that $q>1$. As $S(\de_1,\ldots,\de_q;\tau_s,\tau_t)\ne 0$ we must have $s\ne t$ by \cite[Th.~4.5]{Joyc7}, so either $s<t$ or $s>t$. Consider $S(\de_1,\ldots,\de_q;\tau_s,\tau_{t'})$ for $t'$ in the interval $[s,t]$ or $[t,s]$. This is nonzero if and only if (a),(b) hold for all $j=1,\ldots,q-1$ with $t'$ in place of $t$, which is true for $t'=t$ and false for $t'=s$. Thus at some $t'$ in $[t,s]$ or $[s,t]$ and some $j=1,\ldots,q-1$ we must cross between $\tau_{t'}(\de_1+\cdots+\de_j)>\tau_{t'}(\de_{j+1}+\cdots+\de_q)$ and~$\tau_{t'}(\de_1+\cdots+\de_j)\le\tau_{t'}(\de_{j+1}+\cdots+\de_q)$. 

As $(\tau_t,T_t,\le)_{t\in[0,1]}$ is a continuous family in the sense of Definition \ref{co3def5}, there are only finitely many $t'$ in $[t,s]$ or $[s,t]$ where such crossings can occur, and at such $t'$ we have $\tau_{t'}(\de_1+\cdots+\de_j)=\tau_{t'}(\de_{j+1}+\cdots+\de_q)$. Define $t''$ to be the value of $t'$ where such a crossing occurs which is closest to $t$ in $[t,s]$ or $[s,t]$. Then for all $j=1,\ldots,q-1$ we have either
\begin{itemize}
\setlength{\itemsep}{0pt}
\setlength{\parsep}{0pt}
\item[(a$)''$] $\tau_s(\de_j)\le\tau_s(\de_{j+1})$ and $\tau_{t''}(\de_1+\cdots+\de_j)>\tau_{t''}(\de_{j+1}+\cdots+\de_q)$, 
\item[(b$)''$] $\tau_s(\de_j)>\tau_s(\de_{j+1})$ and $\tau_{t''}(\de_1+\cdots+\de_j)<\tau_{t''}(\de_{j+1}+\cdots+\de_q)$, or
\item[(c$)''$] $\tau_{t''}(\de_1+\cdots+\de_j)=\tau_{t''}(\de_{j+1}+\cdots+\de_q)=\tau_{t''}(\al)$, 
\end{itemize}
and (c$)''$ occurs for at least one $j$. Let there be $p-1$ values of $j$ for which (c$)''$ occurs, for $p>1$, write them as $0<c_1<\cdots<c_{p-1}<q$, and set $c_0=0$ and $c_p=q$. Define $\be_j=\de_{c_{j-1}+1}+\cdots+\de_{c_j}$ for $j=1,\ldots,p$. Then $\al=\be_1+\cdots+\be_p$, and the conditions of $(*)_{r+1}$ hold for $\al,n,\al_i$, this $p>1$ and $\be_1,\ldots,\be_p$, the given $s$, and $t=u=t''$. Hence Lemma \ref{co7lem4} 
shows there are most finitely many $n\ge 1$ and $\bar\al_1,\ldots,\bar\al_n\in H^2(X,\Z)$ satisfying the conditions above with $q>1$, and if such data $\al,n,p,\al_i,a_j,\be_j,\ab s,t,u$ exist with $q>1$ then~$\int_X\De(\al)\cup\om_s^{m-2}\ge 0$.

Next suppose $q=1$. Then $\de_1=\al$, $d_0=0$, $d_1=n$ and $\tau_s(\al_i)=\tau_s(\al)$ for $i=1,\ldots,n$. The argument of Lemma \ref{co7lem4} with $\al_1,\ldots,\al_n$ in place of $\be_1,\ldots,\be_p$ (replacing \eq{co7eq73} with $\int_X\De(\al_i)\cup\om_s^{m-2}\ge 0$, which holds by Theorem \ref{co7thm1} as $\es\ne\M_\al^\ss(\tau_s)\subseteq\M_\al^\ss(\mu^{\om_s})$) shows there are most finitely many $n\ge 1$ and $\bar\al_1,\ldots,\bar\al_n\in H^2(X,\Z)$ satisfying the conditions above with $q=1$, and if such data $\al,n,p,\al_i,a_j,\be_j,\ab s,t,u$ exist with $q>1$ then~$\int_X\De(\al)\cup\om_s^{m-2}\ge 0$.
\end{proof}

Lemmas \ref{co7lem4}--\ref{co7lem6} prove the inductive step. Hence $(*)_r$ holds for all $r\ge 1$ by induction. Now let $\al\in C(\coh(X))_\pe$ with $\dim\al=m$. We will prove Assumption \ref{co5ass3}(b) for $\al$. Set $r=\rank\al$ and $\bar\al=c_1(\al)$, and define
\e
D=\max\bigl(\ts\sup_{x\in[0,1]}\int_X\De(\al)\cup\om_x^{m-2},0\bigr).
\label{co7eq79}
\e

Suppose $n\ge p\ge 1$, $\al_1,\ldots,\al_n\in C(\coh(X))$ and $0=a_0<a_1<\cdots<a_p=n$ such that $\al_1+\cdots+\al_n=\al$ and, writing $\be_j=\al_{a_{j-1}+1}+\cdots+\al_{a_j}$ for $j=1,\ldots,p$, then there exist $s,t,u\in[0,1]$ with $\M_{\al_i}^\ss(\tau_s)\ne\es$ for $i=1,\ldots,n$, and $U\bigl(\al_{a_{j-1}+1},\ab\al_{a_{j-1}+2},\ab\ldots,\ab\al_{a_j};\tau_s,\tau_t)\ne 0$ for $j=1,\ldots,p$, and $\tau_u(\be_1)=\cdots=\tau_u(\be_p)$. Equation \eq{co7eq79} implies that $\int_X\De(\al)\cup\om_s^{m-2}\le D$. Thus $(*)_r$ implies that there are only finitely many possible sets of data $n\ge 1$ and $\bar\al_1,\ldots,\bar\al_n\in H^2(X,\Z)$ such that $c_1(\al_i)=\bar\al_i$ for $i=1,\ldots,n$. Since $\bar\mu^{\om_s}(\al_i)$ depends only on $\rank\al_i=1,\ldots,r$ and $c_1(\al_i)=\bar\al_i$, there are only finitely many possible values for $\bar\mu^{\om_s}(\al_i)$ realized by all such choices of data $\al,n,p,\al_i,a_j,\be_j,\ab s,t,u$. Write $\mu_0\in\R$ for the maximum of these finitely many values. Then $\bar\mu^{\om_s}(\al_i)\le\mu_0$ for all such $\al,n,p,\al_i,a_j,\be_j,\ab s,t,u$ and $i=1,\ldots,n$.

As $\M_{\al_i}^\ss(\tau_s)\ne\es$ we may choose $[E_i]\in\M_{\al_i}^\ss(\tau_s)$ for $i=1,\ldots,n$. Then $E=E_1\op\cdots\op E_n$ has $\lb E\rb=\al_1+\cdots+\al_n=\al$. As $E_i$ is $\tau_s$-semistable it is pure of dimension $\dim\al_i=m$, so $E$ is pure. Also
\begin{equation*}
\mu_{\max}^{\om_s}(E)=\max\bigl(\bar\mu^{\om_s}(\al_1),\ldots,\bar\mu^{\om_s}(\al_n)\bigr)\le \mu_0.
\end{equation*}

Now Proposition \ref{co7prop6} with $K=\{\om_s:s\in[0,1]\}$, which is compact as $[0,1]$ is compact, implies that the family of pure $E\in\coh(X)$ with $\lb E\rb=\al$ and $\mu_{\max}^{\om_s}(E)\le\mu_0$ for some $s\in[0,1]$ is bounded. But in a bounded family of sheaves $E$, the decompositions $E=E_1\op\cdots\op E_n$ for all $n\ge 1$ can realize only finitely many sets of classes $(\lb E_1\rb,\ldots,\lb E_n\rb)$ in $C(\coh(X))^n$. Thus there are only finitely many possibilities for $n\ge 1$ and $(\al_1,\ldots,\al_n)$. This implies the finiteness condition in Assumption \ref{co5ass3}(b), and the condition $\al_{i_1}+\cdots+\al_{i_2}\in C(\coh(X))_\pe$ follows from \eq{co7eq46}--\eq{co7eq47} and $\dim\al_i=m$. This completes the proof.

\subsection{Counting sheaves on curves}
\label{co76}

Let $X$ be a connected smooth projective complex curve. Then \S\ref{co72}--\S\ref{co75} define data and verify Assumptions \ref{co4ass1} and \ref{co5ass1}--\ref{co5ass3} for $\A=\B=\coh(X)$. Thus Theorems \ref{co5thm1}--\ref{co5thm3} apply to define invariants $[\M_\al^\ss(\tau)]_\inv\in H_*(\M_\al^\pl)$ counting $\tau$-semistable coherent sheaves on $X$. When $X=\CP^1$ we can also consider the algebraic $\C$-groups $\bG_m$ or $\PGL(2,\C)$ acting on $X$, which satisfy Definition \ref{co7def12}(i)--(ii), and \S\ref{co724} and \S\ref{co735} verify Assumptions \ref{co4ass2} and \ref{co5ass4}, so Theorem \ref{co5thm4} applies, and we can promote the invariants $[\M_\al^\ss(\tau)]_\inv$ to equivariant homology $H_*^G(\M_\al^\pl)$.

For curves $X$ the general picture simplifies in several ways:
\begin{itemize}
\setlength{\itemsep}{0pt}
\setlength{\parsep}{0pt}
\item[(a)] For $\om\in\Kah(X)$, Gieseker stability $(\tau^\om,G,\le)$ and $\mu$-stability $(\mu^\om,M,\le)$ in \S\ref{co712} coincide, $(\tau^\om,G,\le)=(\mu^\om,M,\le)$.
\item[(b)] If $\om,\ti\om\in\Kah(X)$ then as $X$ is a connected curve, Hilbert polynomials in \eq{co7eq5} satisfy $P_\al^{\ti\om}(t)=P_\al^\om(Ct)$, where $C=\int_X\ti\om/\int_X\om>0$. Thus, although $\mu^{\ti\om}\ne\mu^\om$ if $\int_X\ti\om\ne\int_X\om$, we see that $\M^\ss_\al(\mu^\om)$ is independent of $\om\in\Kah(X)$, so $[\M_\al^\ss(\mu^\om)]_\inv$ is independent of $\om$ by Theorem~\ref{co5thm1}(ii).
\item[(c)] All the moduli stacks are smooth, so by Theorem \ref{co2thm1}(ii) we can replace the Behrend--Fantechi virtual classes $[\M_\al^\ss(\tau)]_\virt,[\baM^\ss_{(\al,1)}(\bar\tau^0_1)]_\virt$ in Theorem \ref{co5thm1} by the fundamental classes $[\M_\al^\ss(\tau)]_\fund,[\baM^\ss_{(\al,1)}(\bar\tau^0_1)]_\fund$.
\end{itemize}

Part (b) means that the wall crossing in Theorems \ref{co5thm2}--\ref{co5thm3} is trivial, so most of our theory is vacuous in this example. However, we do define invariants $[\M_\al^\ss(\tau)]_\inv$ when $\M_\al^\rst(\tau)\ne\M_\al^\ss(\tau)$, which may be new and potentially interesting. The next theorem summarizes Theorems \ref{co5thm1} and \ref{co5thm4} in this case.

\begin{thm}
\label{co7thm6}
Let $X$ be a connected, smooth, projective, complex curve, and write $\M^\pl=\coprod_{\al\in K(\coh(X))}\M_\al^\pl$ for the projective linear moduli stack of objects in $\coh(X)$. Write $(\mu,M,\le)$ for $\mu$-stability on $\coh(X)$ for some $\om\in\Kah(X)$. Then for all\/ $\al$ in $C(\coh(X))$ there are unique classes $[\M_\al^\ss(\mu)]_\inv$ in the Betti\/ $\Q$-homology group $H_{2-2\chi(\al,\al)}(\M_\al^\pl)=\check H_0(\M^\pl_\al)$ satisfying:
\begin{itemize}
\setlength{\itemsep}{0pt}
\setlength{\parsep}{0pt}
\item[{\bf(i)}] If\/ $\M_\al^\rst(\mu)=\M_\al^\ss(\mu)$ then\/ $[\M_\al^\ss(\mu)]_\inv=[\M_\al^\ss(\mu)]_\fund,$ where $[\M_\al^\ss(\mu)]_\fund$ is the fundamental class of the smooth projective\/ $\C$-scheme $\M_\al^\ss(\mu)$.
\item[{\bf(ii)}] The $[\M_\al^\ss(\mu)]_\inv$ are independent of\/ $\om\in\Kah(X)$.
\item[{\bf(iii)}] Let\/ $\al\in C(\coh(X)),$ and\/ $\O_X(1)\ra X$ be an ample line bundle, and\/ $N$ be sufficiently large that\/ $E$ is $N$-regular for $\O_X(1)$ for all\/ $[E]\in\M_\al^\ss(\mu)$. Then Example\/ {\rm\ref{co5ex1}} defines a smooth projective moduli scheme\/ $\baM^\ss_{(\al,1)}(\bar\mu^0_1)$ of nonzero morphisms\/ $\rho:\O_X(-N)\ra E$ for\/ $[E]\in\M_\al^\ss(\mu),$ with a forgetful morphism\/ $\Pi_{\M^\ss_\al(\mu)}:\baM^\ss_{(\al,1)}(\bar\mu^0_1)\ra\M_\al^\ss(\mu)$. We have
\ea
&(\Pi_{\M^\ss_\al(\mu)})_*\bigl([\baM^\ss_{(\al,1)}(\bar\mu^0_1)]_\fund\cap c_\top(\bT_{\baM^\ss_{(\al,1)}(\bar\mu^0_1)/\M^\ss_{\al}(\mu)})\bigr)
\label{co7eq80}\\
&=\!\!\!\!\!\!\sum_{\begin{subarray}{l}n\ge 1,\;\al_1,\ldots,\al_n\in
C(\coh(X)):\\ \al_1+\cdots+\al_n=\al,\\ 
\mu(\al_i)=\mu(\al), \; \M_{\al_i}^\ss(\mu)\ne\es,\; \text{all\/ $i$}\end{subarray}\!\!\!\!\!\!\!\!\!\!\!\!\!\!\!\!\!\!\!} \!\!\!\!\!\!\!\!\!\!\begin{aligned}[t]
\frac{(-1)^{n+1}P_{\al_1}^{c_1(\O_X(1))}(N)}{n!}\,\cdot&\bigl[\bigl[\cdots\bigl[[\M_{\al_1}^\ss(\mu)]_\inv,\\
&\; 
[\M_{\al_2}^\ss(\mu)]_\inv\bigr],\ldots\bigr],[\M_{\al_n}^\ss(\mu)]_\inv\bigr],
\end{aligned}
\nonumber
\ea	
where $P_{\al_1}^{c_1(\O_X(1))}(t)$ is the Hilbert polynomial of $\al_1,$ the Lie brackets are in the Lie algebra $\check H_0(\M^\pl)$ from Theorem\/ {\rm\ref{co4thm2},} and there are only finitely many terms in the sum.
\end{itemize}
When $X=\CP^1,$ so that\/ $G=\bG_m$ or $\PGL(2,\C)$ acts on $X,$ the analogue holds in $G$-equivariant homology $H_*^G(\M^\pl),$ as in\/ {\rm\S\ref{co23}}.
\end{thm}

In the sequel \cite{Joyc13} we will use the method of pair invariants in Chapter \ref{co8} to compute the invariants $[\M_\al^\ss(\mu)]_\inv$ fairly explicitly, after projecting to the homology $H_*(\baM^\pl)$ of the projective linear moduli stack $\baM^\pl$ of objects in $D^b\coh(X)$, using the explicit description of $H_*(\baM^\pl)$ in Gross~\cite{Gros}.

\subsection{Counting sheaves on surfaces}
\label{co77}

Next let $X$ be a connected smooth projective complex surface. Then \S\ref{co72}--\S\ref{co75} define data and verify Assumptions \ref{co4ass1} and \ref{co5ass1}--\ref{co5ass3} for $\A=\B=\coh(X)$. Thus Theorems \ref{co5thm1}--\ref{co5thm3} apply to define invariants $[\M_\al^\ss(\tau)]_\inv\in H_*(\M_\al^\pl)$ counting $\tau$-semistable coherent sheaves on $X$. If an algebraic $\C$-group $G$ acts on $X$ satisfying Definition \ref{co7def12}(i)--(ii) then \S\ref{co724} and \S\ref{co735} verify Assumptions \ref{co4ass2} and \ref{co5ass4}, so Theorem \ref{co5thm4} applies, and we can promote the invariants $[\M_\al^\ss(\tau)]_\inv$ to equivariant homology $H_*^G(\M_\al^\pl)$. An interesting case of this is when $X$ is a toric surface and~$G=\bG_m^2$.

The results of \S\ref{co71}--\S\ref{co74} are at their most complete and powerful for surfaces. As in Definition \ref{co7def11}, the family $\sS$ of weak stability conditions is Gieseker stability $(\tau^\om,G,\le)$ and $\mu$-stability $(\mu^\om,M,\le)$ for all $\om\in\Kah(X)$, not just $\om\in\Kah_\Q(X)$, and as $m=2$ sections \ref{co74}--\ref{co75} verify Assumptions \ref{co5ass2}(g),(h) and \ref{co5ass3}(b) for all $\al\in C(\coh(X))$, without restrictions on~$\dim\al$.

As in Definition \ref{co7def14} the {\it geometric genus\/} of $X$ is $p_g=\dim H^0(K_X)$, and as in \S\ref{co732}--\S\ref{co733}, when $\rank\al>0$ we define the quasi-smooth derived stacks $\bs\dM_\al^\red,\bs\dM_\al^\rpl$ and the associated obstruction theories \eq{co5eq5} on $\M_\al,\M_\al^\pl$ differently in the cases $p_g=0$ (with {\it non-reduced\/} obstruction theories) and $p_g>0$ (with {\it reduced\/} obstruction theories), with $o_\al=p_g$ in both cases.  

This is very important in the wall-crossing formulae \eq{co5eq30}--\eq{co5eq34}, because of the conditions $o_{\al_1}+\cdots+o_{\al_n}=o_\al$ in the sums. If $\rank\al>0$ then any nonzero term $\al_1,\ldots,\al_n$ in \eq{co5eq30}--\eq{co5eq34} must have $\rank\al_i>0$ for all $i$, so that $o_\al=o_{\al_i}=p_g$. Hence if $p_g=0$ then $o_{\al_1}+\cdots+o_{\al_n}=o_\al$ holds automatically, but if $p_g>0$ then $o_{\al_1}+\cdots+o_{\al_n}=o_\al$ forces $n=1$. Thus, if $p_g>0$ and $\rank\al>0$ then all of \eq{co5eq31}--\eq{co5eq34} reduce to $[\M_\al^\ss(\ti\tau)]_\inv=[\M_\al^\ss(\tau)]_\inv$. That is, in this case the wall crossing formulae simply show that the invariants $[\M_\al^\ss(\tau)]_\inv$ are independent of the choice of weak stability condition~$(\tau,T,\le)$.

The next two theorems, which the author considers some of the main results of this book, summarize Theorems \ref{co5thm1}--\ref{co5thm3} and \ref{co5thm4} in the two cases $p_g=0$ and $p_g>0$. In Theorem \ref{co7thm8} for simplicity we restrict to $\rank\al>0$, and we use the condition $o_{\al_1}+\cdots+o_{\al_n}=o_\al$ in \eq{co5eq30}--\eq{co5eq34} to exclude all terms except those with $n=1$, radically simplifying the formulae.

\begin{thm}
\label{co7thm7}
Let\/ $X$ be a connected projective complex surface with geometric genus\/ $p_g=0$. Write $\M^\pl=\coprod_{\al\in K(\coh(X))}\M_\al^\pl$ for the projective linear moduli stack of objects in $\coh(X)$. Then:
\smallskip

\noindent{\bf(a)} Let\/ $(\tau,T,\le)$ be Gieseker stability $(\tau^\om,G,\le)$ or $\mu$-stability $(\mu^\om,M,\le)$ on $\coh(X)$ for some $\om\in\Kah(X)$. Then for all\/ $\al\in C(\coh(X))$ there are unique classes $[\M_\al^\ss(\tau)]_\inv$ in the Betti\/ $\Q$-homology group $\check H_0(\M^\pl_\al)=H_{2-2\chi(\al,\al)}(\M_\al^\pl)$ satisfying:
\begin{itemize}
\setlength{\itemsep}{0pt}
\setlength{\parsep}{0pt}
\item[{\bf(i)}] If\/ $\M_\al^\rst(\tau)=\M_\al^\ss(\tau)$ then\/ $[\M_\al^\ss(\tau)]_\inv=[\M_\al^\ss(\tau)]_\virt,$ where $[\M_\al^\ss(\tau)]_\virt$ is the Behrend--Fantechi virtual class of the proper algebraic space $\M_\al^\ss(\tau)$ with (non-reduced) perfect obstruction theory, as in\/~{\rm\S\ref{co732}}.
\item[{\bf(ii)}] Let\/ $\al\in C(\coh(X)),$ and\/ $\O_X(1)\ra X$ be an ample line bundle, and\/ $N$ be sufficiently large that\/ $E$ is $N$-regular for $\O_X(1)$ for all\/ $[E]\in\M_\al^\ss(\tau)$. Then Example\/ {\rm\ref{co5ex1}} defines a proper moduli space\/ $\baM^\ss_{(\al,1)}(\bar\tau^0_1)$ of nonzero morphisms\/ $\rho:\O_X(-N)\ra E$ for\/ $[E]\in\M_\al^\ss(\tau),$ with an obstruction theory, and a forgetful morphism\/ $\Pi_{\M^\ss_\al(\tau)}:\baM^\ss_{(\al,1)}(\bar\tau^0_1)\!\ra\!\M_\al^\ss(\tau)$. We~have
\ea
&(\Pi_{\M^\ss_\al(\tau)})_*\bigl([\baM^\ss_{(\al,1)}(\bar\tau^0_1)]_\virt\cap c_\top(\bT_{\baM^\ss_{(\al,1)}(\bar\tau^0_1)/\M^\ss_{\al}(\tau)})\bigr)
\label{co7eq81}\\
&=\!\!\!\!\!\!\sum_{\begin{subarray}{l}n\ge 1,\;\al_1,\ldots,\al_n\in
C(\coh(X)):\\ \al_1+\cdots+\al_n=\al,\\ 
\tau(\al_i)=\tau(\al), \; \M_{\al_i}^\ss(\tau)\ne\es,\; \text{all\/ $i$}\end{subarray}\!\!\!\!\!\!\!\!\!\!\!\!\!\!\!\!\!\!\!\!\!\!\!\!\!} \!\!\!\begin{aligned}[t]
\frac{(-1)^{n+1}P_{\al_1}^{c_1(\O_X(1))}(N)}{n!}\,\cdot&\bigl[\bigl[\cdots\bigl[[\M_{\al_1}^\ss(\tau)]_\inv,\\
&\; 
[\M_{\al_2}^\ss(\tau)]_\inv\bigr],\ldots\bigr],[\M_{\al_n}^\ss(\tau)]_\inv\bigr],
\end{aligned}
\nonumber
\ea
where $P_{\al_1}^{c_1(\O_X(1))}(t)$ is the Hilbert polynomial of\/ $\al_1,$ the Lie brackets are in the Lie algebra $\check H_0(\M^\pl)$ from Theorem\/ {\rm\ref{co4thm2},} and there are only finitely many terms in the sum.
\end{itemize}

\noindent{\bf(b)} Suppose $(\tau,T,\le),(\ti\tau,\ti T,\le)$ are two Gieseker or $\mu$-stability conditions as in part\/ {\bf(a)\rm,} from possibly different\/ $\om,\ti\om\in\Kah(X)$. Then 
\begin{itemize}
\setlength{\itemsep}{0pt}
\setlength{\parsep}{0pt}
\item[{\bf(i)}] If\/ $\al\in C(\coh(X))$ with\/ $\M_\al^\ss(\tau)=\M_\al^\ss(\ti\tau)$ then $[\M_\al^\ss(\tau)]_\inv=[\M_\al^\ss(\ti\tau)]_\inv$.
\item[{\bf(ii)}] For all\/ $\al\in C(\coh(X))$ we have
\end{itemize}
\e
\begin{gathered}
{}
[\M_\al^\ss(\ti\tau)]_\inv= \!\!\!\!\!\!\!
\sum_{\begin{subarray}{l}n\ge 1,\;\al_1,\ldots,\al_n\in
C(\coh(X)):\\ \M_{\al_i}^\ss(\tau)\ne\es,\; \text{all\/ $i,$} \\
\al_1+\cdots+\al_n=\al 
\end{subarray}} \!\!\!\!\!\!\!\!\!\!\!\begin{aligned}[t]
\ti U(\al_1,&\ldots,\al_n;\tau,\ti\tau)\cdot\bigl[\bigl[\cdots\bigl[[\M_{\al_1}^\ss(\tau)]_\inv,\\
&
[\M_{\al_2}^\ss(\tau)]_\inv\bigr],\ldots\bigr],[\M_{\al_n}^\ss(\tau)]_\inv\bigr]
\end{aligned}
\end{gathered}
\label{co7eq82}
\e
\begin{itemize}
\setlength{\itemsep}{0pt}
\setlength{\parsep}{0pt}
\item[] in the Lie algebra $\check H_0(\M^\pl)$ from {\rm\S\ref{co43}}. Here $\ti U(-;\tau,\ti\tau)$ is as in Theorem\/ {\rm\ref{co3thm3},} and there are only finitely many nonzero terms in \eq{co7eq82}. Equivalently, in the universal enveloping algebra $\bigl(U(\check H_0(\M^\pl)),*\bigr)$ we have:
\end{itemize}
\e
\begin{gathered}
{}
[\M_\al^\ss(\ti\tau)]_\inv= \!\!\!\!\!\!\!
\sum_{\begin{subarray}{l}n\ge 1,\;\al_1,\ldots,\al_n\in
C(\coh(X)):\\ \M_{\al_i}^\ss(\tau)\ne\es,\; \text{all\/ $i,$} \\
\al_1+\cdots+\al_n=\al \end{subarray}} \!\!\!\!\!\!\!\!\!\!\!\!
\begin{aligned}[t]
U(\al_1,&\ldots,\al_n;\tau,\ti\tau)\cdot[\M_{\al_1}^\ss(\tau)]_\inv *{}\\
&
[\M_{\al_2}^\ss(\tau)]_\inv*\cdots *[\M_{\al_n}^\ss(\tau)]_\inv.
\end{aligned}
\end{gathered}
\label{co7eq83}
\e

\noindent{\bf(c)} If an algebraic $\C$-group $G$ acts on $X$ satisfying Definition\/ {\rm\ref{co7def12}(i)--(ii)} then {\bf(a)\rm,\bf(b)} generalize to equivariant homology $H_*^G(\M^\pl)$.
\end{thm}

\begin{thm}
\label{co7thm8}
Let\/ $X$ be a connected projective complex surface with geometric genus\/ $p_g=\dim H^0(K_X)>0$. Write $\M^\pl=\coprod_{\al\in K(\coh(X))}\M_\al^\pl$ for the projective linear moduli stack of objects in $\coh(X)$. Then:
\smallskip

\noindent{\bf(a)} Let\/ $(\tau,T,\le)$ be Gieseker stability $(\tau^\om,G,\le)$ or $\mu$-stability $(\mu^\om,M,\le)$ on $\coh(X)$ for some $\om\in\Kah(X)$. Then for all\/ $\al\in C(\coh(X))$ with\/ $\rank\al>0$ there are unique classes $[\M_\al^\ss(\tau)]_\inv$ in the Betti\/ $\Q$-homology group $\check H_{2p_g}(\M^\pl_\al)\ab=H_{2p_g+2-2\chi(\al,\al)}(\M_\al^\pl)$ satisfying:
\begin{itemize}
\setlength{\itemsep}{0pt}
\setlength{\parsep}{0pt}
\item[{\bf(i)}] If\/ $\M_\al^\rst(\tau)=\M_\al^\ss(\tau)$ then\/ $[\M_\al^\ss(\tau)]_\inv=[\M_\al^\ss(\tau)]_\virt,$ where $[\M_\al^\ss(\tau)]_\virt$ is the Behrend--Fantechi virtual class of the proper algebraic space $\M_\al^\ss(\tau)$ with \begin{bfseries}reduced\end{bfseries} perfect obstruction theory, constructed in\/ {\rm\S\ref{co733}} from the natural obstruction theory by deleting $H^0(K_X)\ot\O_{\M_\al^\ss(\tau)}$ from\/~$h^{-1}(\cdots)$.
\item[{\bf(ii)}] The invariants $[\M_\al^\ss(\tau)]_\inv$ for $\rank\al>0$ are independent of the weak stability condition\/~$(\tau,T,\le)$.
\item[{\bf(iii)}] Let\/ $\al\in C(\coh(X))$ with\/ $\rank\al>0,$ and\/ $\O_X(1)\ra X$ be an ample line bundle, and\/ $N$ be sufficiently large that\/ $E$ is $N$-regular for $\O_X(1)$ for all\/ $[E]\in\M_\al^\ss(\tau)$. Then Example\/ {\rm\ref{co5ex1}} defines a proper moduli space\/ $\baM^\ss_{(\al,1)}(\bar\tau^0_1)$ of nonzero morphisms\/ $\rho:\O_X(-N)\ra E$ for\/ $[E]\in\M_\al^\ss(\tau),$ with a reduced perfect obstruction theory as in {\bf(i)\rm,} and a forgetful morphism\/ $\Pi_{\M^\ss_\al(\tau)}:\baM^\ss_{(\al,1)}(\bar\tau^0_1)\ra\M_\al^\ss(\tau)$. We have
\e
\begin{split}
&(\Pi_{\M^\ss_\al(\tau)})_*\bigl([\baM^\ss_{(\al,1)}(\bar\tau^0_1)]_\virt\cap c_\top(\bT_{\baM^\ss_{(\al,1)}(\bar\tau^0_1)/\M^\ss_{\al}(\tau)})\bigr)\\
&\qquad =P_\al^{c_1(\O_X(1))}(N)\cdot [\M_\al^\ss(\tau)]_\inv,
\end{split}
\label{co7eq84}
\e
where $P_\al^{c_1(\O_X(1))}(t)$ is the Hilbert polynomial of\/ $\al$. 
\end{itemize}

\noindent{\bf(b)} We also define invariants $[\M_\al^\ss(\tau)]_\inv$ when $\rank\al=0$. However, these use the 
\begin{bfseries}non-reduced\end{bfseries} obstruction theory, lie in $\check H_0(\M^\pl_\al)=H_{2-2\chi(\al,\al)}(\M_\al^\pl),$ and satisfy Theorem\/ {\rm\ref{co7thm7}(a),(b)} rather than\/~{\bf(a)(i)\rm--\bf(iii)}.
\smallskip

\noindent{\bf(c)} If an algebraic $\C$-group $G$ acts on $X$ satisfying Definition\/ {\rm\ref{co7def12}(i)--(ii)} then {\bf(a)\rm,\bf(b)} generalize to equivariant homology $H_*^G(\M^\pl)$.
\end{thm}

In the sequel \cite{Joyc13} we will use the method of pair invariants in Chapter \ref{co8} to compute the invariants $[\M_\al^\ss(\tau)]_\inv$ in Theorems \ref{co7thm7}--\ref{co7thm8} fairly explicitly, after projecting to the homology $H_*(\baM^\pl)$ of the projective linear moduli stack $\baM^\pl$ of objects in $D^b\coh(X)$, using the explicit description of $H_*(\baM^\pl)$ in Gross~\cite{Gros}.

In Theorem \ref{co7thm7} or \ref{co7thm8}, suppose $\al\in C(\coh(X))$ with $\rank\al>0$ and $\int_X\De(\al)<0$, where $\De(\al)$ is the discriminant as in \S\ref{co713}. Then Theorem \ref{co7thm1} implies that $\M_\al^\rst(\tau)=\M_\al^\ss(\tau)=\es$, so Theorem \ref{co7thm7}(a)(i) and \ref{co7thm8}(a)(i) imply: 

\begin{cor}
\label{co7cor3}
In Theorems\/ {\rm\ref{co7thm7}} and\/ {\rm\ref{co7thm8},} if\/ $\al\in C(\coh(X))$ with\/ $\rank\al>0$ and\/ $\int_X\De(\al)<0$ then $[\M_\al^\ss(\tau)]_\inv=0$.
\end{cor}

\subsubsection{Discussion and open problems}
\label{co771}

\noindent{\bf(A) Relation to the work of Mochizuki \cite{Moch}.}

\noindent{\bf(a)} In his monograph \cite{Moch}, Mochizuki proves rough analogues of Theorems \ref{co7thm7}--\ref{co7thm8}. This book generalizes and builds on \cite{Moch}. The author would like to acknowledge, with thanks, that ideas from \cite{Moch} were an inspiration to the author at several points, in particular the proof of Theorem \ref{co5thm2} in Chapter~\ref{co10}.
\smallskip

\noindent{\bf(b)} Mochizuki \cite[\S 7]{Moch} defines his invariants not primarily as homology classes $[\M_\al^\ss(\tau)]_{\rm Moch}$ in $H_*(\M_\al^\pl)$, as we do, but instead as integrals $\int_{[\M_\al^\ss(\tau)]_{\rm Moch}}\Up\in\Q$ for certain universal cohomology classes $\Up\in H^*(\M_\al^\pl)$, in the style of Donaldson's polynomial invariants of 4-manifolds \cite[\S 9.2]{DoKr} (see {\bf(B)} below). So Mochizuki defines invariants in the dual $\cR^*$ of a $\Q$-algebra $\cR$ of classes~$\Up$.

Mochizuki's invariants are defined for all $\al\in C(\coh(X))$, including $\al$ with $\M_\al^\rst(\tau)\ne\M_\al^\ss(\tau)$. In our language, using \cite[Prop.~7.3.1]{Moch} the homology class $[\M_\al^\ss(\tau)]_{\rm Moch}$ underlying Mochizuki's invariants may be written, in the situations of Theorems \ref{co7thm7}(a)(ii) and \ref{co7thm8}(a)(iii), by
\ea
&[\M_\al^\ss(\tau)]_{\rm Moch}=
\label{co7eq85}\\
&\lim_{N\ra\iy}\frac{1}{P_\al^{c_1(\O_X(1))}(N)}\cdot
(\Pi_{\M^\ss_\al(\tau)})_*\bigl([\baM^\ss_{(\al,1)}(\bar\tau^0_1)]_\virt\cap c_\top(\bT_{\baM^\ss_{(\al,1)}(\bar\tau^0_1)/\M^\ss_{\al}(\tau)})\bigr).
\nonumber
\ea

If $p_g>0$ and $\rank\al>0$ then \eq{co7eq84} gives $[\M_\al^\ss(\tau)]_{\rm Moch}=[\M_\al^\ss(\tau)]_\inv$, so Mochizuki's extension of Behrend--Fantechi virtual classes $[\M_\al^\ss(\tau)]_\virt$ to the case $\M_\al^\rst(\tau)\ne\M_\al^\ss(\tau)$, when $[\M_\al^\ss(\tau)]_\virt$ is not defined, agrees with ours.

If $p_g=0$ and $\rank\al>0$ then in \eq{co7eq81} we have
\e
\lim_{N\ra\iy}\frac{P_{\al_1}^{c_1(\O_X(1))}(N)}{P_\al^{c_1(\O_X(1))}(N)}=\frac{\rank\al_1}{\rank\al}.
\label{co7eq86}
\e
Combining equations \eq{co7eq81}, \eq{co7eq85} and \eq{co7eq86} yields
\begin{align*}
&[\M_\al^\ss(\tau)]_{\rm Moch}\\
&=\sum_{\begin{subarray}{l}n\ge 1,\;\al_1,\ldots,\al_n\in
C(\coh(X)):\\ \al_1+\cdots+\al_n=\al,\\ 
\tau(\al_i)=\tau(\al), \; \M_{\al_i}^\ss(\tau)\ne\es,\; \text{all\/ $i$}\end{subarray}} \!\!\!\!\!\!\!\!\!\!\!\!\begin{aligned}[t]
\frac{(-1)^{n+1}\rank\al_1}{n!\rank\al}\,\cdot&\bigl[\bigl[\cdots\bigl[[\M_{\al_1}^\ss(\tau)]_\inv,\\
&\; 
[\M_{\al_2}^\ss(\tau)]_\inv\bigr],\ldots\bigr],[\M_{\al_n}^\ss(\tau)]_\inv\bigr].
\end{aligned}
\end{align*}
Thus Mochizuki's extension to $\M_\al^\rst(\tau)\ne\M_\al^\ss(\tau)$ is {\it different\/} from ours in general, though it does agree with ours if $\rank\al=1$ or~2.

If $\rank\al=0$ the right hand side of \eq{co7eq86} does not make sense. When $\dim\al_1=\dim\al$ (which is automatic in \eq{co7eq81}) the limit in \eq{co7eq86} exists, but if $\dim\al=1$ it may depend on the choice of ample line bundle $\O_X(1)\ra X$ in Theorem \ref{co7thm7}(a)(ii). Thus $[\M_\al^\ss(\tau)]_{\rm Moch}$ may depend on the choice of $\O_X(1)$ when $\dim\al=1$, which appears to have been overlooked in~\cite[Prop.~7.3.1]{Moch}.

\smallskip

\noindent{\bf(c)} There are two main ways in which Theorems \ref{co7thm7}--\ref{co7thm8} improve on~\cite{Moch}:
\begin{itemize}
\setlength{\itemsep}{0pt}
\setlength{\parsep}{0pt}
\item[{\bf(i)}] Mochizuki has no analogue of the Lie bracket on $\check H_0(\M^\pl)$ in~\eq{co7eq81}--\eq{co7eq83}.
\item[{\bf(ii)}] Mochizuki has no analogue of the universal wall-crossing formulae \eq{co7eq82}--\eq{co7eq83} or the coefficients $\ti U(\al_1,\ldots,\al_n;\tau,\ti\tau),U(\al_1,\ldots,\al_n;\tau,\ti\tau)$.
\end{itemize}

In Mochizuki's set up, it is clear on general grounds that the left hand side $[\M_\al^\ss(\ti\tau)]_\inv$ of \eq{co7eq82} is some unknown function of the terms $[\M_{\al_i}^\ss(\tau)]_\inv$ appearing on the right hand side of \eq{co7eq82}, and of $\al_i,\tau,\ti\tau$ and algebraic-topological data on $X$ such as $\td(X)$. In simple situations this unknown function can be computed. See \cite[Prop.~7.4.8]{Moch} for $\rank\al=2$, \cite[\S 7.7.3--\S 7.7.4]{Moch} for calculations when $\rank\al=3$, and \cite[Th.s 7.6.2 \& 7.7.1]{Moch} for $\rank\al>0$. The results often involve residue formulae reminiscent of the theory of vertex algebras, which enter our theory as in Chapter~\ref{co4}.
\smallskip

\noindent{\bf(d)} As in \cite[Th.~7.5.2]{Moch}, Mochizuki uses the method of pair invariants to write his invariants for $\rank\al>1$ in terms of invariants for $\rank\al=1$ and algebraic Seiberg--Witten invariants. We explain this in Chapter \ref{co8} and~\cite{Joyc13}.	
\medskip

\noindent{\bf(B) Relation to Donaldson invariants of oriented 4-manifolds.}

\noindent{\bf(a)} Let $X$ be a compact, oriented 4-manifold with $b^2_+(X)\ge 1$ and $G$ a compact Lie group. As in Donaldson--Kronheimer \cite[\S 9]{DoKr}, the {\it Donaldson invariants\/} of $X$ are defined by choosing a Riemannian metric $g$ on $X$ and considering integrals $\int_{\oM_\al}\Up$, where $\oM_\al$ is a compactified moduli space of instantons (anti-self-dual connections) on a principal $G$-bundle $P\ra X$ with characteristic class $\ch(P)=\al$, and $\Up$ is a universal cohomology class on $\oM_\al$. If $b^2_+(X)>1$ then $\int_{\oM_\al}\Up$ is independent of $g$, but if $b^2_+(X)=1$ then $\int_{\oM_\al}\Up$ changes by a wall-crossing formula \cite{GNY1} depending on the splitting $H^2(X,\R)=H^2_+(X,\R)\op H^2_-(X,\R)$ induced by $g$. Most of the literature takes $G=\SU(2)$ or $\SO(3)$, but Kronheimer \cite{Kron} and
Mari\~{n}o--Moore \cite{MaMo} discuss $G=\SU(r)$ and $G=\mathop{\rm PSU}(r)$ for~$r\ge 2$.

If $X$ is a projective complex surface with K\"ahler metric $g$ and K\"ahler form $\om$ and $G=\U(r)$ for $r>1$ and $\int_Xc_1(\al)\cup\om=0$ then moduli spaces $\oM_\al$ are compactifications of moduli spaces $\M^\vect_\al(\tau^\om)$ of $\tau^\om$-semistable algebraic vector bundles on $X$. The moduli spaces $\M^\ss_\al(\tau^\om)$ in Theorems \ref{co7thm7}--\ref{co7thm8} are different compactifications of the same moduli spaces $\M^\vect_\al(\tau^\om)$. Because of this, the invariants $[\M_\al^\ss(\tau)]_\inv$ in Theorems \ref{co7thm7}--\ref{co7thm8} for $\rank\al=r>1$ should be thought of as algebraic analogues of $\U(r)$-Donaldson invariants.

We have $b^2_+(X)=1+2p_g$, so the cases $p_g=0$ and $p_g>0$ in Theorems \ref{co7thm7}--\ref{co7thm8} correspond to $b^2_+(X)=1$ and $b^2_+(X)>1$ in Donaldson theory. The choice of weak stability condition $(\tau,T,\le)$ corresponds to the splitting $H^2(X,\R)=H^2_+(X,\R)\op H^2_-(X,\R)$ in Donaldson theory.

\begin{prob}
\label{co7prob1}
Find a precise relationship between the invariants $[\M_\al^\ss(\tau)]_\inv$ of Theorems\/ {\rm\ref{co7thm7}--\ref{co7thm8}} for $\rank\al=r\ge 2$ and the Donaldson invariants of the underlying K\"ahler $4$-manifold for $G=\U(r)$ or\/~$\SU(r)$.
\end{prob}

This could elucidate aspects of Donaldson theory which are not yet fully understood, for example, how to count moduli spaces with reducible connections (corresponding to our $\M_\al^\ss(\tau)$ with $\M_\al^\rst(\tau)\ne\M_\al^\ss(\tau)$), and wall-crossing formulae when $b^2_+(X)=1$ for $\SU(r)$ or $\U(r)$-Donaldson invariants when~$r>2$.
\medskip

\noindent{\bf(C) Generating functions, Gieseker versus $\mu$-stability.} Let $X$ be a connected projective complex surface with $p_g=0$. Using the notation of Theorem \ref{co7thm7}, for $\om\in\Kah(X)$, $r>0$ and $c\in H^2(X,\Z)$, define a generating function
\e
Z_{r,c}^\om(q)=\sum_{\al\in C(\coh(X)): \rank\al=r,\;
c_1(\al)=c}q^{-2\int_X\ch_2(\al)}\,[\M^\ss_\al(\mu^\om)]_\inv.
\label{co7eq87}
\e
Since $-2\int_X\ch_2(\al)\in\Z$ and $[\M^\ss_\al(\mu^\om)]_\inv=0$ if $-2\int_X\ch_2(\al)<\frac{1}{r}\int_Xc^2$ by Corollary \ref{co7cor3}, we see that $Z_{r,c}^\om(q)\in\check H_0(\M^\pl)[[q]][q^{-1}]$.

For $\tau=\mu^\om$, $\ti\tau=\mu^{\ti\om}$ the coefficients $\ti U(\al_1,\ldots,\al_n;\tau,\ti\tau),U(\al_1,\ldots,\al_n;\tau,\ti\tau)$ in \eq{co7eq82}--\eq{co7eq83} depend only on $r_i=\rank\al_i$ and $c_i=c_1(\al_i)$, not on $\ch_2(\al_i)$, so we may write $\ti U(\al_1,\ldots,\al_n;\mu^\om,\mu^{\ti\om})=\ti U((r_1,c_1),\ldots,(r_n,c_n);\mu^\om,\mu^{\ti\om})$. Thus, multiplying \eq{co7eq82} by $q^{-2\int_X\ch_2(\al)}$ and summing over $\al$ as in \eq{co7eq87} yields
\begin{equation*}
Z_{r,c}^{\ti\om}(q)= \!\!\!
\sum_{\begin{subarray}{l}n\ge 1,\; r_1,\ldots,r_n>0,\\ c_1,\ldots,c_n\in H^2(X,\Z):\\
r=r_1+\cdots+r_n,\; c=c_1+\cdots+c_n 
\end{subarray}\!\!\!\!\!\!\!\!\!\!\!\!\!\!\!\!\!\!\!\!\!\!} \!\!\!\!\!\!\!\!\!\!\begin{aligned}[t]
\ti U((r_1,c_1),\ldots,&(r_n,c_n);\mu^\om,\mu^{\ti\om})\cdot\bigl[\bigl[\cdots\bigl[Z_{r_1,c_1}^\om(q),\\
&
Z_{r_2,c_2}^\om(q)\bigr],\ldots\bigr],Z_{r_n,c_n}^\om(q)\bigr]
\end{aligned}
\end{equation*}
in the Lie algebra $\check H_0(\M^\pl)[[q]][q^{-1}]$. This argument does not work for Gieseker stability, as $\ti U(\al_1,\ldots,\al_n;\tau^\om,\tau^{\ti\om})$ may also depend on the~$\ch_2(\al_i)$.

This suggests that for surfaces with $p_g=0$, $\mu$-stability may be preferable to Gieseker stability for forming nice generating functions of invariants (say with modular properties), and for comparison with Donaldson theory as in {\bf(B)}. When $p_g>0$ there is no difference, by Theorem~\ref{co7thm8}(a)(ii).
\medskip

\noindent{\bf(D) K-theoretic Donaldson invariants, and similar.} Let $X$ be a connected projective surface, $\om\in\Amp_\Q(X)$, and $\al\in C(\coh(X))$ with $\rank\al=2$ and $\M_\al^\rst(\tau^\om)=\M_\al^\ss(\tau^\om)$, so that $\M_\al^\ss(\tau^\om)$ is a projective scheme. Then G\"{o}ttsche--Nakajima--Yoshioka \cite{GNY2} define the `K-theoretic Donaldson invariants' of $X$ to be the virtual holomorphic Euler characteristics of certain universal holomorphic line bundles on $\M_\al^\ss(\tau^\om)$. (See also \cite{Gott,GoYu}.) By the Grothendieck--Riemann--Roch Theorem we may write these invariants as $\int_{[\M_\al^\ss(\tau^\om)]_\virt}\Up$ for a universal cohomology class $\Up$. Thus K-theoretic Donaldson invariants are determined by the invariants $[\M^\ss_\al(\tau)]_\inv$ in Theorems~\ref{co7thm7}--\ref{co7thm8}.

Our invariants $[\M^\ss_\al(\tau)]_\inv$ lie in $H^*(\M_\al^\pl)$, or the derived category version $H^*(\baM_\al^\pl)$, which are infinite-dimensional spaces, so $[\M^\ss_\al(\tau)]_\inv$ contains a lot of information. The author expects that most invariants in the literature counting sheaves on surfaces can be written in terms of $\int_{[\M_\al^\ss(\tau^\om)]_\inv}\Up$ for suitable $\Up$, and so can be studied using Theorems \ref{co7thm7}--\ref{co7thm8} and the results of~\cite{Joyc13}.
\medskip

\noindent{\bf(E) Duality in $D^b\coh(X)$.} Let $X$ be a connected projective surface, and consider the derived category $D^b\coh(X)$, as in Huybrechts \cite{Huyb}, and the projective linear moduli stack $\baM^\pl=\coprod_{\al\in K(\coh(X))}\baM^\pl_\al$ of objects in $D^b\coh(X)$. {\it Duality\/} is a contravariant equivalence of categories $D:D^b\coh(X)\ra D^b\coh(X)$ mapping $D:E^\bu\ra(E^\bu)^\vee$. It takes a vector bundle $E$ to the dual vector bundle $E^*$, but in general takes coherent sheaves to complexes, not coherent sheaves. Duality acts on $K(\coh(X))\subset H^{\rm even}(X,\Q)$ by $D_*:K(\coh(X))\ra K(\coh(X))$ which multiplies $H^{2i}(X,\Q)$ by $(-1)^i$, so that $D_*$ fixes $\rank\al,\ch_2(\al)$ and changes the sign of $c_1(\al)$. Duality induces an equivalence of (higher) stacks $D_*:\baM^\pl\ra\baM^\pl$ which identifies~$\baM^\pl_\al\ra\baM^\pl_{D_*(\al)}$.

Suppose $\tau$ is some kind of weak stability condition on $D^b\coh(X)$ (e.g.\ Bridgeland stability, or Gieseker or $\mu$-stability on $\coh(X)\subset D^b\coh(X)$) for which enumerative invariants $[\oM_\al^\ss(\tau)]_\inv\in H_*(\baM^\pl_\al)$ can be defined. We can define a pull-back stability condition $D^*(\tau)$ on $D^b\coh(X)$. (Here as $D$ is contravariant we must reverse the order of objects in $D^*(\tau)$. For Bridgeland stability, take the complex conjugate of the central charge.) Then for $\al\in K(\coh(X))$ we have $D_*(\oM_\al^\ss(\tau))=\oM_{D_*(\al)}^\ss(D^*(\tau))$, so we expect that $D_*([\oM_\al^\ss(\tau)]_\inv)=[\oM_{D_*(\al)}^\ss(D^*(\tau))]_\inv$. If a wall-crossing formula like \eq{co7eq82}--\eq{co7eq83} holds between $\tau$ and $D^*(\tau)$, this gives identities on the $[\oM_\al^\ss(\tau)]_\inv$.

Now suppose $p_g>0$ and $\rank\al>0$. Then the invariants $[\oM_\al^\ss(\tau)]_\inv$ of Theorem \ref{co7thm8} are independent of $\tau$. This motivates the following:

\begin{conj}
\label{co7conj1}
Let $X$ be a connected projective complex surface with $p_g>0$. For $\al\in C(\coh(X))$ with $\rank\al>0,$ consider the projections of the invariants $[\oM_\al^\ss(\tau)]_\inv\in H_*(\M^\pl_\al)$ of Theorem \ref{co7thm8} to $H_*(\oM^\pl_\al)$, where $\oM^\pl_\al$ is the projective linear moduli stack of objects of class $\al$ in $D^b\coh(X)$. Then
\begin{equation*}
D_*([\oM_\al^\ss(\tau)]_\inv)=[\oM_{D_*(\al)}^\ss(\tau)]_\inv\quad\text{in $H_*(\oM^\pl_{D_*(\al)})$.}
\end{equation*}
\end{conj}

The author hopes to prove this conjecture in~\cite{Joyc13}.
\medskip

\noindent{\bf(F) Parabolic sheaves on surfaces.} Mochizuki \cite{Moch} develops his entire theory not just for the abelian category $\coh(X)$ on a projective surface $X$, but also for the abelian category $\coh_{\rm pa}(X,D)$ of {\it parabolic sheaves\/} on a projective surface $X$ with Cartier divisor $D\subset X$. See Maruyama--Yokogawa \cite{MaYo,Yoko} for background on parabolic sheaves. In particular the constructions of perfect obstruction theories in \cite[\S 5 \& \S 6.1]{Moch} also work for moduli stacks of parabolic sheaves. So it seems very likely that Theorems \ref{co7thm7}--\ref{co7thm8} can be extended from $\coh(X)$ to $\coh_{\rm pa}(X,D)$, hopefully with little extra work.
 
\subsection{Counting sheaves on Fano 3-folds}
\label{co78}

A smooth projective $\C$-scheme $X$ of dimension $m$ is called a {\it Fano $m$-fold\/} if the anticanonical bundle $K_X^{-1}\ra X$ is ample. In the beginning of Donaldson--Thomas theory, Thomas \cite[Cor.~3.39]{Thom1} defined Donaldson--Thomas invariants $[\M_\al^\ss(\tau)]_\virt$ for moduli spaces $\M_\al^\ss(\tau)$ of semistable sheaves on $X$ with $\M_\al^\rst(\tau)=\M_\al^\ss(\tau)$, where $X$ is either (a) a Calabi--Yau 3-fold or (b) a Fano 3-fold. 

Thomas' key observation, which we used in \S\ref{co734}, was that if $X$ is a smooth projective 3-fold then obstruction theories can be defined on $\M_\al^\ss(\tau)$ in these two cases. Donaldson--Thomas invariants of Calabi--Yau 3-folds have been intensively studied \cite{JoSo,KoSo}, but those of Fano 3-folds have so far received very little attention --- in fact the author knows of no other papers on them.

Let $X$ be a connected complex Fano 3-fold. Then \S\ref{co72}--\S\ref{co75} define data and verify Assumptions \ref{co4ass1} and \ref{co5ass1}--\ref{co5ass3} for $\A=\B=\coh(X)$. Thus Theorems \ref{co5thm1}--\ref{co5thm3} apply to define invariants $[\M_\al^\ss(\tau)]_\inv\in H_*(\M_\al^\pl)$ counting $\tau$-semistable coherent sheaves on $X$. If an algebraic $\C$-group $G$ acts on $X$ satisfying Definition \ref{co7def12}(i)--(ii) then \S\ref{co724} and \S\ref{co735} verify Assumptions \ref{co4ass2} and \ref{co5ass4}, so Theorem \ref{co5thm4} applies, and we can promote the invariants $[\M_\al^\ss(\tau)]_\inv$ to equivariant homology $H_*^G(\M_\al^\pl)$. An interesting case of this is when $X$ is a toric Fano 3-fold and~$G=\bG_m^3$.

Note that the Fano 3-fold case in \S\ref{co72}--\S\ref{co75} has some important restrictions, which are included in Theorem \ref{co7thm9} below, and discussed in \S\ref{co781}. In \eq{co7eq47} we defined $C(\coh(X))_\pe$ to be the set of $\al\in C(\coh(X))$ with $\dim\al=1$ or 3, that is, we exclude sheaves of dimensions 0, 2 from our theory. Dimension 0 sheaves are excluded since we cannot construct suitable quasi-smooth $\bs\dM_\al^\red,\bs\dM_\al^\rpl$ in Assumption \ref{co5ass1}(f) or perfect obstruction theories \eq{co5eq5} in this case, as in Remark \ref{co7rem11}. Dimension 2 sheaves are excluded as in Remark \ref{co7rem5}(a), since the author cannot prove Assumptions \ref{co5ass2}(b),(g),(h) and \ref{co5ass3}(b) in this case.

As in Remark \ref{co7rem5}(b) and Theorem \ref{co7thm9}(c) below, we can allow $\dim\al=2$ if we restrict to $\om\in\Kah_\Q(X)$, since then we can prove Assumption \ref{co5ass2}(b),(g),(h) as in \S\ref{co72}, \S\ref{co74} using Proposition \ref{co7prop4}. But we can only prove the wall-crossing formulae \eq{co7eq88}--\eq{co7eq89} for $\tau,\ti\tau$ coming from the same $\om\in\Kah_\Q(X)$, as distinct $\om,\ti\om\in\Kah_\Q(X)$ cannot be connected in $\Kah_\Q(X)$, but only in~$\Kah(X)$.

Also, in Definition \ref{co7def11} when $m>2$ we restricted to $(\tau^\om,G,\le),(\mu^\om,M,\le)$ for $\om$ in a small open $U\subset\Kah(X)$. This was needed in Proposition \ref{co7prop18} (when $\dim\al=3$), but not in Proposition \ref{co7prop17} (when $\dim\al=1$). Thus when $\dim\al=3$ we only prove \eq{co7eq88}--\eq{co7eq89} when $\om,\ti\om$ are sufficiently close in $\Kah(X)$. The next theorem summarizes Theorems \ref{co5thm1}--\ref{co5thm3} and \ref{co5thm4} for Fano 3-folds. 

\begin{thm}
\label{co7thm9}
Let\/ $X$ be a connected complex Fano $3$-fold. Write $\M^\pl=\coprod_{\al\in K(\coh(X))}\M_\al^\pl$ for the projective linear moduli stack of objects in $\coh(X)$. Then:
\smallskip

\noindent{\bf(a)} Let\/ $(\tau,T,\le)$ be Gieseker stability $(\tau^\om,G,\le)$ or $\mu$-stability $(\mu^\om,M,\le)$ on $\coh(X)$ for some $\om\in\Kah(X)$. Then for all\/ $\al\in C(\coh(X))$ with\/ $\dim\al=1$ or\/ $3$ there are unique classes $[\M_\al^\ss(\tau)]_\inv$ in the Betti\/ $\Q$-homology group $\check H_0(\M^\pl_\al)=H_{2-2\chi(\al,\al)}(\M_\al^\pl)$ satisfying:
\begin{itemize}
\setlength{\itemsep}{0pt}
\setlength{\parsep}{0pt}
\item[{\bf(i)}] If\/ $\M_\al^\rst(\tau)=\M_\al^\ss(\tau)$ then\/ $[\M_\al^\ss(\tau)]_\inv=[\M_\al^\ss(\tau)]_\virt,$ where $[\M_\al^\ss(\tau)]_\virt$ is the Behrend--Fantechi virtual class of the proper algebraic space $\M_\al^\ss(\tau)$ with perfect obstruction theory, as in\/~{\rm\S\ref{co734}}.
\item[{\bf(ii)}] Let\/ $\al\in C(\coh(X))$ with\/ $\dim\al=1$ or\/ $3,$ and\/ $\O_X(1)\ra X$ be an ample line bundle, and\/ $N$ be sufficiently large that\/ $E$ is $N$-regular for $\O_X(1)$ for all\/ $[E]\in\M_\al^\ss(\tau)$. Then Example\/ {\rm\ref{co5ex1}} defines a proper moduli space\/ $\baM^\ss_{(\al,1)}(\bar\tau^0_1)$ of nonzero morphisms\/ $\rho:\O_X(-N)\ra E$ for\/ $[E]\in\M_\al^\ss(\tau),$ with an obstruction theory, and a projection\/ $\Pi_{\M^\ss_\al(\tau)}:\baM^\ss_{(\al,1)}(\bar\tau^0_1)\ab\ra\!\M_\al^\ss(\tau)$. We~have
\begin{align*}
&(\Pi_{\M^\ss_\al(\tau)})_*\bigl([\baM^\ss_{(\al,1)}(\bar\tau^0_1)]_\virt\cap c_\top(\bT_{\baM^\ss_{(\al,1)}(\bar\tau^0_1)/\M^\ss_{\al}(\tau)})\bigr)\\
&=\!\!\!\!\!\!\sum_{\begin{subarray}{l}n\ge 1,\;\al_1,\ldots,\al_n\in
C(\coh(X)):\\ \al_1+\cdots+\al_n=\al,\\ 
\tau(\al_i)=\tau(\al), \; \M_{\al_i}^\ss(\tau)\ne\es,\; \text{all\/ $i$}\end{subarray}\!\!\!\!\!\!\!\!\!\!\!\!\!\!\!\!\!\!\!\!\!\!\!\!\!} \!\!\!\begin{aligned}[t]
\frac{(-1)^{n+1}P_{\al_1}^{c_1(\O_X(1))}(N)}{n!}\,\cdot&\bigl[\bigl[\cdots\bigl[[\M_{\al_1}^\ss(\tau)]_\inv,\\
&\; 
[\M_{\al_2}^\ss(\tau)]_\inv\bigr],\ldots\bigr],[\M_{\al_n}^\ss(\tau)]_\inv\bigr],
\end{aligned}
\nonumber
\end{align*}
where $P_{\al_1}^{c_1(\O_X(1))}(t)$ is the Hilbert polynomial of $\al_1,$ the Lie brackets are in the Lie algebra $\check H_0(\M^\pl)$ from Theorem\/ {\rm\ref{co4thm2},} and there are only finitely many terms in the sum.
\end{itemize}

\noindent{\bf(b)} Suppose $(\tau,T,\le),(\ti\tau,\ti T,\le)$ are two Gieseker or $\mu$-stability conditions as in part\/ {\bf(a)\rm,} from possibly different\/ $\om,\ti\om\in\Kah(X),$ where in {\bf(ii)} below if\/ $\dim\al=3$ we require $\om,\ti\om$ to be sufficiently close in $\Kah(X)$ that there exists connected open $U\subset\Kah(X)$ containing $\om,\ti\om$ and satisfying Definition\/ {\rm\ref{co7def11}(b)}. Then 
\begin{itemize}
\setlength{\itemsep}{0pt}
\setlength{\parsep}{0pt}
\item[{\bf(i)}] If\/ $\al\in C(\coh(X))$ with\/ $\dim\al=1$ or\/ $3$ and\/ $\M_\al^\ss(\tau)=\M_\al^\ss(\ti\tau)$ then $[\M_\al^\ss(\tau)]_\inv=[\M_\al^\ss(\ti\tau)]_\inv$.
\item[{\bf(ii)}] For all\/ $\al\in C(\coh(X))$ with\/ $\dim\al=1$ or\/ $3$ we have
\end{itemize}
\e
\begin{gathered}
{}
[\M_\al^\ss(\ti\tau)]_\inv= \!\!\!\!\!\!\!
\sum_{\begin{subarray}{l}n\ge 1,\;\al_1,\ldots,\al_n\in
C(\coh(X)):\\ \M_{\al_i}^\ss(\tau)\ne\es,\; \text{all\/ $i,$} \\
\al_1+\cdots+\al_n=\al 
\end{subarray}} \!\!\!\!\!\!\!\!\!\!\!\begin{aligned}[t]
\ti U(\al_1,&\ldots,\al_n;\tau,\ti\tau)\cdot\bigl[\bigl[\cdots\bigl[[\M_{\al_1}^\ss(\tau)]_\inv,\\
&
[\M_{\al_2}^\ss(\tau)]_\inv\bigr],\ldots\bigr],[\M_{\al_n}^\ss(\tau)]_\inv\bigr]
\end{aligned}
\end{gathered}
\label{co7eq88}
\e
\begin{itemize}
\setlength{\itemsep}{0pt}
\setlength{\parsep}{0pt}
\item[] in the Lie algebra $\check H_0(\M^\pl)$ from {\rm\S\ref{co43}}. Here $\ti U(-;\tau,\ti\tau)$ is as in Theorem\/ {\rm\ref{co3thm3},} and there are only finitely many nonzero terms in \eq{co7eq88}. Equivalently, in the universal enveloping algebra $\bigl(U(\check H_0(\M^\pl)),*\bigr)$ we have:
\end{itemize}
\e
\begin{gathered}
{}
[\M_\al^\ss(\ti\tau)]_\inv= \!\!\!\!\!\!\!
\sum_{\begin{subarray}{l}n\ge 1,\;\al_1,\ldots,\al_n\in
C(\coh(X)):\\ \M_{\al_i}^\ss(\tau)\ne\es,\; \text{all\/ $i,$} \\
\al_1+\cdots+\al_n=\al \end{subarray}} \!\!\!\!\!\!\!\!\!\!\!\!
\begin{aligned}[t]
U(\al_1,&\ldots,\al_n;\tau,\ti\tau)\cdot[\M_{\al_1}^\ss(\tau)]_\inv *{}\\
&
[\M_{\al_2}^\ss(\tau)]_\inv*\cdots *[\M_{\al_n}^\ss(\tau)]_\inv.
\end{aligned}
\end{gathered}
\label{co7eq89}
\e

\noindent{\bf(c)} If\/ $\om\in\Kah_\Q(X)$ then {\bf(a)} also holds when $\dim\al=2$. Also {\bf(b)} holds when $\dim\al=2$ and\/ $\om=\ti\om\in\Kah_\Q(X),$ that is, for changing between Gieseker and\/ $\mu$-stability for a fixed\/~$\om\in\Kah_\Q(X)$.
\smallskip

\noindent{\bf(d)} If an algebraic $\C$-group $G$ acts on $X$ satisfying Definition\/ {\rm\ref{co7def12}(i)--(ii)} then {\bf(a)\rm--\bf(c)} generalize to equivariant homology $H_*^G(\M_\al^\pl)$.
\end{thm}

\subsubsection{Discussion and open problems}
\label{co781}

\noindent{\bf(A) For what classes of projective 3-folds $X$ can we count semistable coherent sheaves?} Thomas \cite[Cor.~3.39]{Thom1} defined Donaldson--Thomas invariants $[\M_\al^\ss(\tau)]_\virt$ for moduli spaces $\M_\al^\ss(\tau)$ of semistable sheaves on $X$ with $\M_\al^\rst(\tau)=\M_\al^\ss(\tau)$, where $X$ is either a Calabi--Yau 3-fold or a Fano 3-fold. 

To see why Thomas' construction works for Calabi--Yau and Fano 3-folds, let $X$ be a smooth projective 3-fold, and define $\bs\M_\al,\bs\M_\al^\pl$ as in Definition \ref{co7def10}. Then as in Remark \ref{co2rem9}, $\bL_i:i^*(\bL_{\bs\M_\al})\ra \bL_{\M_\al}$ and $\bL_i:i^*(\bL_{\bs\M_\al^\pl})\ra \bL_{\M_\al^\pl}$ satisfy the conditions for obstruction theories on $\M_\al,\M_\al^\pl$, except that $i^*(\bL_{\bs\M_\al})$ and $i^*(\bL_{\bs\M_\al^\pl})$ are perfect in $[-2,1]$ rather than $[-1,1]$ by \eq{co7eq48} and~\eq{co7eq50}.

To make them perfect in $[-1,1]$ on some open substacks $\dM_\al\subseteq\M_\al$, $\dM_\al^\pl\subseteq\M_\al^\pl$ we want $h^2(i^*(\bL_{\bs\M_\al})^\vee)$ and $h^2(i^*(\bL_{\bs\M_\al^\pl})^\vee)$ to vanish on $\dM_\al$ and $\dM_\al^\pl$. At a $\C$-point $[E]$ in $\dM_\al$ or $\dM_\al^\pl$ we have
\begin{equation*}
h^2(i^*(\bL_{\bs\M_\al})^\vee)\vert_{[E]}\cong h^2(i^*(\bL_{\bs\M_\al^\pl})^\vee)\vert_{[E]}\cong\Ext^3(E,E)\cong \Hom(E,E\ot K_X)^*.
\end{equation*}

Consider the two cases:
\begin{itemize}
\setlength{\itemsep}{0pt}
\setlength{\parsep}{0pt}
\item[(a)] If $X$ is a Calabi--Yau 3-fold and $[E]\in\M_\al^\rst(\tau)\subseteq\M^\pl_\al$ then as $E$ is $\tau$-stable $\Hom(E,E\ot K_X)\cong\Hom(E,E)\cong\C$, so $h^2(i^*(\bL_{\bs\M_\al^\pl})^\vee\vert_{\M_\al^\rst(\tau)})\cong\O_{\M_\al^\rst(\tau)}$. In this case, in a similar way to `reduced' obstruction theories in Remark \ref{co5rem1}(c), we can modify $\bL_i:i^*(\bL_{\bs\M_\al^\pl})\ra \bL_{\M_\al^\pl}$ on $\M_\al^\rst(\tau)$, deleting $\O_{\M_\al^\rst(\tau)}$ from $h^{-2}(i^*(\bL_{\bs\M_\al^\pl}))$, to make an obstruction theory on $\M_\al^\rst(\tau)\subseteq\M^\pl_\al$. If $\M_\al^\rst(\tau)=\M_\al^\ss(\tau)$ this gives a virtual class~$[\M_\al^\ss(\tau)]_\virt$.
\item[(b)] If $X$ is a Fano 3-fold and $[E]\in\M_\al^\ss(\tau)\subseteq\M^\pl_\al$ with $\dim\al>0$, Definition \ref{co7def17} showed that $\Hom(E,E\ot K_X)=0$. Thus $h^2(i^*(\bL_{\bs\M_\al^\pl})^\vee\vert_{\M_\al^\ss(\tau)})=0$, and $\bL_i:i^*(\bL_{\bs\M_\al^\pl})\ra \bL_{\M_\al^\pl}$ is an obstruction theory on $\M_\al^\ss(\tau)\subseteq\M^\pl_\al$. If $\M_\al^\rst(\tau)=\M_\al^\ss(\tau)$ this gives a virtual class~$[\M_\al^\ss(\tau)]_\virt$.
\end{itemize}

Note that $\bL_i:i^*(\bL_{\bs\M_\al^\pl})\ra \bL_{\M_\al^\pl}$ must be modified in case (a), but not in case (b). This gives Donaldson--Thomas theory for Calabi--Yau 3-folds and for Fano 3-folds a different flavour. Donaldson--Thomas theory for Fano 3-folds fits into the framework of Chapter \ref{co5}. But as we explained in \S\ref{co1.5}, we must modify Chapter \ref{co5} to include Donaldson--Thomas theory for Calabi--Yau 3-folds.

\begin{quest}
\label{co7quest2}
Can we (partially) generalize Theorem\/ {\rm\ref{co7thm9}} to smooth projective $3$-folds $X$ satisfying weaker conditions than being a Fano\/ $3$-fold? 	
\end{quest}

To build an obstruction theory on $\M_\al^\rst(\tau)\!=\!\M_\al^\ss(\tau)$, Thomas \cite[Cor.~3.39]{Thom1} only requires $K_X^{-1}$ to be effective, not ample, and Maulik--Nekrasov--Okounkov--Pandharipande \cite[Lem.~1]{MNOP} observe that the argument also works if instead $H^0(K^{-1}_X)\ne 0$ and $\rank\al>0$, as then choosing $s\in H^0(K^{-1}_X)$ gives an injection $\Hom(E,E\ot K_X)\hookra\Hom(E,E)\cong\C$. But we need an obstruction theory on $\M_\al^\ss(\tau)$ when $\M_\al^\rst(\tau)\ne\M_\al^\ss(\tau)$, which these arguments do not give.
\smallskip

\noindent{\bf(B) Including dimension 2 sheaves.} If we could solve the following problem then Theorem \ref{co7thm9}(a),(b) would also hold for $\al\in C(\coh(X))$ with $\dim\al=2$.

\begin{prob}
\label{co7prob2}
Let\/ $X$ be a connected complex Fano $3$-fold. 
\begin{itemize}
\setlength{\itemsep}{0pt}
\setlength{\parsep}{0pt}
\item[{\bf(a)}] As in Question\/ {\rm\ref{co7quest1},} show that Condition\/ {\rm\ref{co7cond1}} holds
for\/ $K=\{\om\}$ with\/ $\om\in\Kah(X)\sm\Kah_\Q(X)$ and\/ $\al\in C(\coh(X))$ with\/~$\dim\al=2$.	
\item[{\bf(b)}] Show that Assumption\/ {\rm\ref{co5ass3}(b)} holds for\/ $\al\in C(\coh(X))$ with\/~$\dim\al=2$.
\end{itemize}	
\end{prob}

\noindent{\bf(C) Not including dimension 0 sheaves.} We exclude dimension 0 sheaves in Theorem \ref{co7thm9} since we cannot define an obstruction theory on $\M_\al^\ss(\tau)$ when $\dim\al=0$, as in Remark \ref{co7rem11}. The author can see no way to fix this problem.
\smallskip

\noindent{\bf(D) On requiring $\om,\ti\om$ to be close in $\Kah(X)$.} If $\al\in C(\coh(X))$ with $\dim\al=3$ then we only prove the wall-crossing formulae \eq{co7eq88}--\eq{co7eq89} for $\tau,\ti\tau$ coming from $\om,\ti\om$ which are sufficiently close in~$\Kah(X)$.

If $\om,\ti\om$ are not sufficiently close then we can choose a finite sequence $\om=\om_0,\om_1,\ldots,\om_N=\ti\om$ such that $\om_{i-1},\om_i$ are sufficiently close in $\Kah(X)$ for $i=1,\ldots,N$, and then \eq{co7eq88}--\eq{co7eq89} hold to transform from $\om_{i-1}$ to $\om_i$, and we can get from invariants $[\M_\al^\ss(\tau)]_\inv$ to $[\M_{\ti\al}^\ss(\ti\tau)]_\inv$ in $N$ steps.

Equation \eq{co3eq4} shows that wall crossing formulae \eq{co7eq89} can be composed. Because of this, if the last part of Assumption \ref{co5ass3}(b) holds for transformations from $\om_0$ to $\om_i$ for $i=1,\ldots,N$ then \eq{co7eq88}--\eq{co7eq89} hold for $\tau,\ti\tau$. However, if Assumption \ref{co5ass3}(b) fails then \eq{co7eq88}--\eq{co7eq89} could have infinitely many nonzero terms when $\om,\ti\om$ are not sufficiently close, and would not make sense.
\smallskip

\noindent{\bf(E) Reconstructing rank $r$ invariants from rank $1$.} Let $X$ be a Calabi--Yau 3-fold satisfying the conjectural Bogomolov--Gieseker inequality of Bayer--Macr\`\i--Toda \cite{BMT}. In an important recent advance, Feyzbakhsh--Thomas \cite{FeTh} use the wall-crossing theory of Joyce--Song \cite{JoSo} between certain weak stability conditions on $D^b\coh(X)$ to prove that the Donaldson--Thomas invariants of $X$ of rank $r>1$ can be reconstructed from the rank 1 Donaldson--Thomas invariants. If $X$ also satisfies the MNOP Conjecture \cite{MNOP} (known for many Calabi--Yau 3-folds) then the Donaldson--Thomas invariants of $X$ of rank $r\ge 1$ are determined by the Gromov--Witten invariants of $X$.

Theorem \ref{co7thm9} provides an analogue for Fano 3-folds of the wall-crossing theory of Joyce--Song \cite{JoSo}, at least in $\coh(X)$ rather than $D^b\coh(X)$. The Bogomolov--Gieseker inequality for Fano 3-folds is proved by Bernardara--Macr\`\i--Schmidt--Zhao \cite{BMSZ}. The MNOP Conjecture was proved for toric 3-folds, including toric Fano 3-folds, by Maulik--Oblomkov--Okounkov--Pandharipande~\cite{MOOP}.

It seems very likely that Feyzbakhsh--Thomas' programme can be carried out for Fano 3-folds, showing that the invariants $[\M_\al^\ss(\tau)]_\inv$ in Theorem \ref{co7thm9} for $\rank\al>0$ can be determined from those for $\rank\al=1$, and hence from the Gromov--Witten invariants of $X$ when the MNOP Conjecture holds for~$X$.

\section{Pair invariants for coherent sheaves}
\label{co8}

Let $X$ be a smooth projective $\C$-scheme, and $L\ra X$ be a fixed line bundle. A {\it pair\/} $(E,\rho)$ on $X$ is a coherent sheaf $E\in\coh(X)$ and a morphism $\rho:L\ra E$. One can impose stability conditions $\ac\mu$ on $(E,\rho)$. Write $\M^{\rm stp}_\al(\ac\mu)\subseteq\M^{\rm ssp}_\al(\ac\mu)$ for the moduli spaces of $\ac\mu$-(semi)stable pairs with $\lb E\rb=\al$ in $C(\coh(X))$. For suitable $\ac\mu$ these often have the property that $\M^{\rm stp}_\al(\ac\mu)=\M^{\rm ssp}_\al(\ac\mu)$, and $\M^{\rm ssp}_\al(\ac\mu)$ is a projective $\C$-scheme (not just an Artin stack) constructed as a GIT quotient, and if $\dim X\le 3$ then we may be able to define an obstruction theory on $\M^{\rm ssp}_\al(\ac\mu)$, and thus a virtual class~$[\M^{\rm ssp}_\al(\ac\mu)]_\virt$.

There is a large literature studying moduli spaces $\M^{\rm ssp}_\al(\ac\mu)$, and enumerative invariants using virtual (or fundamental) classes $[\M^{\rm ssp}_\al(\ac\mu)]_\virt$. For example:
\begin{itemize}
\setlength{\itemsep}{0pt}
\setlength{\parsep}{0pt}
\item Bradlow \cite{Brad} proves a Hitchin--Kobayashi correspondence between stable pairs $(E,\rho)$ on a compact K\"ahler manifold $X$ with $E\ra X$ a vector bundle, and solutions of a Hermitian--Einstein type equation on $(E,\rho)$. See Bradlow--Daskalopoulos--Garc{\'{\i}}a-Prada--Wentworth \cite{BDGW} for a survey.
\item Garc\'{\i}a-Prada \cite[Th.~4.9]{Garc} relates Bradlow's stable pairs on $X$ to $\SU(2)$-equivariant $\mu$-stable sheaves on~$X\t\CP^1$.
\item Thaddeus \cite{Thad} constructs $\M^{\rm ssp}_\al(\ac\mu)$ using GIT when $X$ is a curve. 
\item D\"urr--Kabanov--Okonek \cite{DKO} define the `Poincar\'e invariants' of a surface $X$ using virtual classes $[\M^{\rm ssp}_\al(\ac\mu)]_\virt$ of moduli spaces of pairs $(E,\rho)$ in which $E$ is a line bundle on $X$, and conjecture they agree with the full Seiberg--Witten invariants of the underlying 4-manifold of $X$. This extends work by previous authors including Friedman--Morgan~\cite{FrMo}.
\item Mochizuki \cite[\S 6.3]{Moch} defines the algebraic Seiberg--Witten invariants of a surface $X$ as in \cite{DKO}, and studies more general pair invariants on~$X$.
\end{itemize}

We will include pair invariants in the framework of Chapter \ref{co5} by considering an abelian category $\acA$ with objects $(E,V,\rho)$ where $E\in\coh(X)$, $V$ is a finite-dimensional $\C$-vector space, and $\rho:V\ot_\C L\ra E$ is a morphism. Then $(E,V,\rho)$ reduces to a pair $(E,\rho)$ when $V=\C$, and to $E\in\coh(X)$ when~$V=0$.

In \S\ref{co86} we explain a method to use wall-crossing, Theorem \ref{co5thm3}, in $\acA$ to write invariants $[\M^\ss_\al(\mu)]_\inv$ in $\coh(X)$ for $\rank\al=r>1$ in terms of $[\M^\ss_\be(\mu)]_\inv$ in $\coh(X)$ for $\rank\be=r-1$ and $[\M^\ss_{(\ga,1)}(\ac\mu)]_\inv$ in $\acA$ for $\rank\ga=1$. Thus by induction we can write rank $r>1$ invariants in $\coh(X)$ in terms of rank 1 invariants in $\acA$ (pair invariants) and rank 1 invariants in~$\coh(X)$.

This method is related to Okonek--Schmitt--Teleman \cite{OST} and Zentner \cite{Zent}, and is used by Mochizuki \cite{Moch} to show that the rank $r>1$ algebraic Donaldson invariants of a surface $X$ can in principle be written in terms of algebraic Seiberg--Witten invariants and rank 1 algebraic Donaldson invariants. We will use the method in the sequel \cite{Joyc13} to (fairly) explicitly compute the invariants $[\M_\al^\ss(\mu)]_\inv$ in Theorems \ref{co7thm6}--\ref{co7thm8} when $X$ is a curve or a surface. 

\begin{rem}
\label{co8rem1}
{\bf(a)} This chapter discusses only stable pairs in the spirit of the authors above \cite{Brad,BDGW,DKO,Garc,Moch,Thad}, in which we embed $\acA\hookra\coh(X\t\CP^1)$ as in Theorem \ref{co8thm1}, and stability on $\acA$ is induced by $\mu$-stability on $\coh(X\t\CP^1)$. Then $(E,V,\rho)\in\acA$ with $V\ne 0$ can only be semistable if $E$ is torsion-free, so $\dim E=m$ and $\rank E>0$. So we consider only invariants for classes $(\al,d)\in C(\acB)$ with $\dim\al=m$, equivalently $\rank\al>0$, and~$d=0,1$.

There are other interesting kinds of pair invariants which do not fit into the set up of this chapter, for example, invariants counting Hilbert schemes of points with $\dim\al=0$, or Pandharipande--Thomas invariants of 3-folds \cite{PaTh1,PaTh2} with $\dim\al=1$, or Joyce--Song's pair invariants for Calabi--Yau 3-folds \cite[\S 12]{JoSo} with $\dim\al=0,\ldots,3$. These other kinds may need different weak stability conditions, and different definitions of obstruction theories.
\smallskip

\noindent{\bf(b)} We have already used pair invariants in Chapter \ref{co7}: Theorems \ref{co7thm6}(iii), \ref{co7thm7}(a)(ii), \ref{co7thm8}(a)(iii), and \ref{co7thm9}(a)(ii) for $\mu$-stability involved virtual classes $[\baM^\ss_{(\al,1)}(\bar\mu^0_1)]_\virt$ of moduli spaces $\baM^\ss_{(\al,1)}(\bar\mu^0_1)$ from Example \ref{co5ex1} of pairs $(E,\rho)$ on $X$ with $L=\O_X(-N)$ for $N\gg 0$. When $\rank\al>0$ these $[\baM^\ss_{(\al,1)}(\bar\mu^0_1)]_\virt$ are a special case of the pair invariants of this chapter, and will be used in~\S\ref{co86}.

\smallskip

\noindent{\bf(c)} As in Chapter \ref{co7}, we write $\M=\coprod_{\al\in K(\coh(X))}\M_\al$ for the moduli stack of objects in $\A=\coh(X)$, and $(\mu,M,\le)$ for $\mu$-stability on $\coh(X)$, and so on. Accents $\ac{}\,$ indicate data associated to $\acA$, so that $\acM,(\ac\mu,\ac M,\le)$ are the moduli stack of objects in $\acA$, and a weak stability condition on~$\acA$.
\end{rem}

\subsection{Background on stable pairs}
\label{co81}

We define an abelian category $\acA$ of morphisms $\rho:V\ot_\C L\ra E$ in $\coh(X)$.

\begin{dfn}
\label{co8def1}
Let $X$ be a smooth, connected, projective $\C$-scheme with $\dim_\C X=m$. In \S\ref{co83} we restrict to $m=1$ or 2, but \S\ref{co81}--\S\ref{co82} work for all $m\ge 0$. Let $L\ra X$ be a fixed line bundle. Define an abelian category $\acA$ to have {\it objects\/} $(E,V,\rho)$ where $E\in\coh(X)$, $V$ is a finite-dimensional $\C$-vector space, and $\rho:V\ot_\C L\ra E$ is a morphism in $\coh(X)$. If $(E,V,\rho),(E',V',\rho')$ are objects in $\acA$, a {\it morphism\/} $(\th,\phi):(E,V,\rho)\ra(E',V',\rho')$ consists of a morphism $\th:E\ra E'$ in $\coh(X)$ and a $\C$-linear map $\phi:V\ra V'$ such that the following commutes in $\coh(X)$:
\e
\begin{gathered}
\xymatrix@C=100pt@R=15pt{ *+[r]{V\ot_\C\O_X} \ar[r]_{\phi\ot\id_{\O_X}} \ar[d]^\rho & *+[l]{V'\ot_\C\O_X} \ar[d]_{\rho'} \\ *+[r]{E} \ar[r]^\th & *+[l]{E'.\!} }	
\end{gathered}
\label{co8eq1}
\e
We define composition of morphisms, and identities $(\id_E,\id_V)$, in the obvious way. Since $\coh(X),\Vect_\C$ are abelian categories, $\acA$ is too. 

Define a surjective quotient $K_0(\acA)\twoheadrightarrow K(\acA)$ by $K(\acA)=K(\coh(X))\op\Z$, for $K(\coh(X))=K^\num(\coh(X))$ as in \S\ref{co711}, such that if $(E,V,\rho)\in\acA$ then $\lb E,V,\rho\rb=(\lb E\rb,\dim V)\in K(\acA)$. We will use this to define weak stability conditions on $\acA$, and in Assumptions \ref{co4ass1}(d) and \ref{co5ass1}(b). Elements of $K(\acA)$ will be written $(\al,d)$ for $\al\in K(\coh(X))$ and~$d\in\Z$.
\end{dfn}

When $V=\C$, an object $(E,\C,\rho)$ in $\acA$ is a morphism $\rho:L\ra E$, that is, a {\it pair\/} $(E,\rho)$ in the sense of \cite{Brad,Garc,Moch,Thad}. So we can regard $\acA$ as an `abelian category of pairs', although this term is unfortunate as objects of $\acA$ are triples. For the rest of the chapter we work in the situation of Definition \ref{co8def1}.

\subsubsection{Embedding $\acute{\mathcal A}$ into $\coh(X\t\CP^1)$}
\label{co811}

Garc\'{\i}a-Prada \cite{Garc} observed that $\acA$ can be identified with a full subcategory of $\SU(2)$-equivariant coherent sheaves $\coh^{\SU(2)}(X\t\CP^1)$ on $X\t\CP^1$: 

\begin{thm}[Garc\'{\i}a-Prada \cite{Garc}]
\label{co8thm1}
There is a full and faithful embedding
\e
\smash{\xymatrix@C=50pt{ \ac I:\acA\,\, \ar@{^{(}->}[r] & \,\coh^{\SU(2)}(X\t\CP^1) }}
\label{co8eq2}
\e
from $\acA$ to the category of\/ $\SU(2)$-equivariant coherent sheaves on $X\t\CP^1,$ where $\SU(2)$ acts trivially on $X$ and in the usual way on $\CP^1$. For each\/ $(E,V,\rho)$ in $\acA$ there is an exact sequence in {\rm$\coh^{\SU(2)}(X\t\CP^1)$:} 
\e
\smash{\xymatrix@C=20pt{ 0 \ar[r] & E\bt\O_{\CP^1} \ar[r] &  \ac I(E,V,\rho) \ar[r] & V\ot_\C(L\bt\O_{\CP^1}(2)) \ar[r] & 0, }}
\label{co8eq3}
\e
functorial in\/ $(E,V,\rho)$ under $\ac I,$ such that\/ \eq{co8eq3} corresponds to the element of
\begin{align*}
&\Ext^1_{X\t\CP^1}\bigl(V\ot_\C(L\bt\O_{\CP^1}(2)),E\bt\O_{\CP^1}\bigr)\\
&=\Hom_{D^b\coh(X\t\CP^1)}\bigl(V\ot_\C(L\bt\O_{\CP^1}(2)),E\bt\O_{\CP^1}[1]\bigr)
\end{align*}
obtained as the composition
\begin{equation*}
\xymatrix@C=49pt{ {\begin{subarray}{l}\ts V\ot_\C \\ \ts 	
(L\bt\O_{\CP^1}(2))\end{subarray}} \ar[r]^(0.49){\id_V\ot (\id_L\bt\la)} &  {\begin{subarray}{l}\ts V\ot_\C \\ \ts (L\bt\O_{\CP^1})[1]\end{subarray}} \ar[r]^(0.49){(\phi\bt\id_{\O_{\CP^1}})[1]} & E\bt\O_{\CP^1}[1], }
\end{equation*}
where $\la\in \Ext^1_{\CP^1}(\O_{\CP^1}(2),\O_{\CP^1})$ corresponds to $1$ in $\Hom(\O_{\CP^1},\O_{\CP^1}(2)\ot K_{\CP^1})^*\cong\Hom(\O_{\CP^1},\O_{\CP^1})^*\cong\C$ under Serre duality.

By abuse of notation we also write $\ac I:\acA\ra\coh(X\t\CP^1)$ for the composition of\/ $\ac I$ in \eq{co8eq2} with the faithful functor~$\coh^{\SU(2)}(X\t\CP^1)\ra\coh(X\t\CP^1)$.
\end{thm}

This will be very useful, as we can deduce definitions and theorems for $\acA$ from definitions and theorems for $\coh(X\t\CP^1)$ in Chapter~\ref{co7}.

\subsubsection{Moduli stacks of objects in $\acute{\mathcal A}$}
\label{co812}

We discuss the moduli stacks $\acM,\acM^\pl$ of objects in $\acA$, extending Definition~\ref{co7def2}.

\begin{dfn}
\label{co8def2}
Work in the situation of Definition \ref{co8def1}. As in \S\ref{co711}, write $\M=\coprod_{\al\in K(\coh(X))}\M_\al$ for the moduli stack of objects in $\coh(X)$, and $\cU\ra X\t\M$ for the universal sheaf, flat over $\M$. Write $\M_{\Vect_\C}=\coprod_{n\ge 0}[*/\GL(n,\C)]$ for the moduli stack of finite-dimensional $\C$-vector spaces. Then there is a universal vector space $\cV\ra\M_{\Vect_\C}$, a locally free sheaf on~$\M_{\Vect_\C}$.

On the product $X\t\M\t\M_{\Vect_\C}$ we have a locally free coherent sheaf $\Pi_{\Vect_\C}^*(\cV)$ and a coherent sheaf $\Pi_{X\t\M}^*(\cU)$, flat over $\M\t\M_{\Vect_\C}$, so we can consider morphisms $\rho:\Pi_{\M_{\Vect_\C}}^*(\cV)\ra\Pi_{X\t\M}^*(\cU)$ of coherent sheaves on $X\t\M\t\M_{\Vect_\C}$, locally over $\M\t\M_{\Vect_\C}$. Since $X$ is projective and $\Pi_{X\t\M}^*(\cU)$ is flat over $\M\t\M_{\Vect_\C}$, results of Grothendieck \cite[\S III.7.7.8--\S III.7.7.9]{Grot1} imply that there exists a moduli stack $\acM$ parametrizing such morphisms $\rho$, with a finite type affine morphism $(\Pi_\M,\Pi_{\M_{\Vect_\C}}):\acM\ra\M\t\M_{\Vect_\C}$. Then $\acM$ is the moduli stack of objects in $\acA$, so $\C$-points of $\acM$ may be written $[E,V,\rho]$. Since $\M\t\M_{\Vect_\C}$ is an Artin stack locally of finite type and $(\Pi_\M,\Pi_{\M_{\Vect_\C}})$ is affine and of finite type, $\acM$ is also an Artin stack locally of finite type.

On $X\t\acM$ there is a universal morphism
\e
\cR:L\bt\Pi_{\M_{\Vect_\C}}^*(\cV)\longra\Pi_{X\t\M}^*(\cU),
\label{co8eq4}
\e
such that $\cR\vert_{X\t\{[E,V,\rho]\}}$ is isomorphic to $\rho:L\ot_\C V\ra E$ for each $(E,V,\rho)\in\acA$. This has the universal property that if $S$ is an Artin $\C$-stack and $\rho:L\bt V\ra\cE$ is a morphism of coherent sheaves on $X\t S$ with $V\ra S$ locally free and $\cE$ flat over $S$, there exists a morphism $\phi:S\ra\baM$ unique up to equivalence such that there are isomorphisms $\cE\cong(\id_X\t(\Pi_\M\ci\phi))^*(\cU)$, $V\cong(\Pi_{\M_{\Vect_\C}}\ci\phi)^*(\cV)$ identifying $\rho$ with~$(\id_X\t\phi)^*(\cR)$.

Alternatively, as in Proposition \ref{co7prop12} we can show that $\acA$ is of {\it compact type\/} in the sense of \S\ref{co332}, taking $\tiA$ to be the abelian category of triples $(E,V,\rho)$ as in Definition \ref{co8def1}, but with $E\in\qcoh(X)$ and $V$ possibly infinite-dimensional. Then as in Alper--Halpern-Leistner--Heinloth \cite[\S 7]{AHLH}, the moduli stack $\acM$ of objects in $\acA$ exists as an Artin $\C$-stack locally of finite type. 

Theorem \ref{co8thm1} induces morphisms of moduli stacks
\e
\xymatrix@C=32pt{ \acM\,\, \ar@{^{(}->}[rr]^(0.3){\text{closed inclusion}} && \M_{\coh^{\SU(2)}(X\t\CP^1)}=\M_{X\t\CP^1}^{\SU(2)} \ar[r] & \M_{X\t\CP^1}, }
\label{co8eq5}
\e
where $\M_{X\t\CP^1}^{\SU(2)}$ is the fixed substack of the $\SU(2)$-action on $\M_{X\t\CP^1}$, and the first morphism is a closed inclusion, as for a point of $\M_{\coh^{\SU(2)}(X\t\CP^1)}$ to be the central term in a sequence of the form \eq{co8eq3} is a closed condition. If $H^1(\O_X)=0$, so that the line bundle $L$ is rigid, this is also an open condition, and the first morphism in \eq{co8eq5} is an open and closed inclusion.

Since the maps $\M\ra K(\coh(X))$, $[E]\mapsto\lb E\rb$ and $\M_{\Vect_\C}\ra\Z$, $[V]\mapsto\dim V$ are locally constant, the map $\acM\ra K(\acB)$, $[E,V,\rho]\mapsto\lb E,V,\rho\rb$ is locally constant, so we have a decomposition $\acM=\coprod_{(\al,d)\in K(\acB)}\acM_{(\al,d)}$ of $\acM$ into open and closed $\C$-substacks $\acM_{(\al,d)}\subset \acM$ of objects in class $(\al,d)$. 

As in Definition \ref{co3def2} we may also form the `projective linear' moduli stack $\acM^\pl=\acM/[*/\bG_m]$ of objects in $\acA$ moduli projective linear isomorphisms, where the $[*/\bG_m]$-action $\ac\Psi:[*/\bG_m]\t\acM\ra\acM$ on $\acM$ is as defined in Definition \ref{co8def5} below. It has a decomposition $\acM^\pl=\coprod_{(\al,d)\in K(\acB)}\acM_{(\al,d)}^\pl$. The projection $\ac\Pi^\pl:\acM^\pl\ra\acM$ is a principal $[*/\bG_m]$-bundle except over~0.

The above arguments also make sense in Derived Algebraic Geometry, so we have a derived moduli stack $\bs\acM$ of objects in $\acA$ with classical truncation $\acM=t_0(\bs\acM)$, so that $\bs\acM$ is a derived Artin $\C$-stack. We can construct $\bs\acM$ directly using To\"en--Vaqui\'e \cite{ToVa}, starting with a dg-category of morphisms of complexes $\rho:L\ot_\C V^\bu\ra E^\bu$ in~$D^b\qcoh(X)$. 

There is a forgetful morphism $(\Pi_{\bs\M},\Pi_{\M_{\Vect_\C}}):\bs\acM\ra\bs\M\t\M_{\Vect_\C}$, where $\bs\M$ is as in \S\ref{co711}, which realizes $\bs\acM$ as the total space of a derived family of morphisms $\rho$. As for \eq{co8eq4}, on $X\t\bs\acM$ there is a universal morphism
\e
\bs\cR:L\bt\Pi_{\M_{\Vect_\C}}^*(\cV)\longra\Pi_{X\t\bs\M}^*(\bs\cU),
\label{co8eq6}
\e
with the universal property above, but for $\bS$ a derived Artin $\C$-stack. The splitting $\acM=\coprod_{(\al,d)\in K(\acB)}\acM_{(\al,d)}$ lifts to $\bs\acM=\coprod_{(\al,d)\in K(\acB)}\bs\acM_{(\al,d)}$. 

There is also a `projective linear' version $\bs\acM^\pl=\bs\acM/[*/\bG_m]$, where the $[*/\bG_m]$-action $\bs{\ac\Psi}:[*/\bG_m]\t\bs\acM\ra\bs\acM$ on $\bs\acM$ is as in Definition \ref{co8def6} below, with $\acM^\pl=t_0(\bs\acM^\pl)$ as $t_0(\bs{\ac\Psi})=\ac\Psi$. It splits as $\bs\acM^\pl=\coprod_{(\al,d)\in K(\acB)}\bs\acM_{(\al,d)}^\pl$. The projection $\bs{\ac\Pi}{}^\pl:\bs\acM^\pl\ra\bs\acM$ is a principal $[*/\bG_m]$-bundle except over~0.
\end{dfn}

\subsubsection{Some weak stability conditions on $\acute{\mathcal A}$}
\label{co813}

We generalize $\mu$-stability $(\mu^\om,M,\le)$ on $\coh(X)$ from Definition \ref{co7def4} to~$\acA$.

\begin{dfn}
\label{co8def3}
Continue in the situation of Definition \ref{co8def1}. Let $\om\in\Kah(X)$ and $s>0$. For each $(\al,d)\in K(\acA)$, define the {\it Hilbert polynomial\/} of $(\al,d)$ by
\e
P^{\om,s}_{(\al,d)}(t)=P^\om_\al(t)+\left(\frac{ds}{m!}\int_X\om^m\right)\,t^{m-1},
\label{co8eq7}
\e
where $P^\om_\al(t)$ is the Hilbert polynomial of $\al$ on $X$ with respect to $\om$ as in Definition \ref{co7def3}. If $(\al,d)\in C(\acA)$ then $P^{\om,s}_{(\al,d)}$ is nonzero with positive leading coefficient. Define the totally ordered set $(\ac M,\le)$ of real polynomials $p(t)=t^k+a_{k-1}t^{k-1}$ with $k\le m$ as in Definition \ref{co7def4} for $X$. Define $\ac\mu^{\om,s}:C(\acA)\ra\ac M$ by $\ac\mu^{\om,s}(\al,d)=t^k+\frac{a_{k-1}}{a_k}t^{k-1}$ when $P^{\om,s}_{(\al,d)}=a_kt^k+a_{k-1}t^{k-1}+\cdots+a_0$ with $\deg P^{\om,s}_{(\al,d)}=k$, so that $a_k>0$. Then $(\ac\mu^{\om,s},\ac M,\le)$ is a weak stability condition on $\acA$, by the same argument as for $\mu$-stability in Definition~\ref{co7def4}.

For $(\al,d)\in C(\acA)$, we also define $\bar\mu^{\om,s}(\al,d)\in\R$ by 
\e
\bar\mu^{\om,s}(\al,d)=a_{k-1}\quad\text{when}\quad \ac\mu^{\om,s}(\al)=t^k+a_{k-1}t^{k-1}.
\label{co8eq8}
\e
It follows from \eq{co8eq7}--\eq{co8eq8} that if $(\al,d)\in C(\acA)$ with $\rank\al>0$ then
\e
\bar\mu^{\om,s}(\al,d)=\bar\mu^\om(\al)+\frac{ds}{\rank\al}.
\label{co8eq9}
\e
We sometimes omit $\om$ and write $(\ac\mu^s,M,\le)=(\ac\mu^{\om,s},\ac M,\le)$ and~$\bar\mu^s=\bar\mu^{\om,s}$.
\end{dfn}

We also generalize the weak stability condition $(\bar\tau^0_1,\bar T,\le)$ in Example \ref{co5ex1} on $\baA$ for $\A=\coh(X)$ and $(\tau,T,\le)=(\mu,M,\le)$ to~$\acA$.

\begin{dfn}
\label{co8def4}
Let $\om\in\Kah(X)$, so we have $\mu$-stability $(\mu^\om,M,\le)$ on $\coh(X)$ as in Definition \ref{co7def4}. Define a total order $(\bar M,\le)$ from $(M,\le)$ as in Definition \ref{co5def1}. Following \eq{co5eq21} for $\la=0$, $\bs\mu=(1)$, define $\bar\mu^{\om,0}_1:C(\acA)\ra \bar M$~by 
\e
\bar\mu^{\om,0}_1(\al,d)=\begin{cases}
\bigl(\mu^\om(\al),d/\rk\al\bigr), & \al\ne 0, \\
\bigl(\iy,1\bigr), & \al=0,
\end{cases}
\label{co8eq10}
\e
where $\rk\al$ is as defined in \S\ref{co723} for Assumption \ref{co5ass2}(f) for $\coh(X)$. Then $(\bar\mu^{\om,0}_1,\bar M,\le)$ is a weak stability condition on $\acA$ as in Definition \ref{co5def1}. Note that when $L=\O_X(-N)$, $(\bar\mu^{\om,0}_1,\bar M,\le)$ is the weak stability condition on morphisms $\rho:V\ot\O_X(-N)\ra E$ used in Theorems \ref{co7thm6}(iii), \ref{co7thm7}(a)(ii), \ref{co7thm8}(a)(iii), and  \ref{co7thm9}(a)(ii) for $\mu$-stability $(\mu^\om,M,\le)$ on $\coh(X)$.
\end{dfn}

\subsubsection{Properties of $\ac\mu^{\om,s}$- and $\bar\mu^{\om,0}_1$-stability}
\label{co814}

We study $(\ac\mu^{\om,s},\ac M,\le)$ and $(\bar\mu^{\om,0}_1,\bar M,\le)$ in Definitions~\ref{co8def3}--\ref{co8def4}.

\begin{lem}
\label{co8lem1}
$\acA$ is {\rm$\ac\mu^{\om,s}$-} and\/ $\bar\mu^{\om,0}_1$-artinian in the sense of Definition\/~{\rm\ref{co3def1}}.	
\end{lem}

\begin{proof} Suppose for a contradiction that $\cdots\subsetneq (E_3,V_3,\rho_3)\subsetneq (E_2,V_2,\rho_2)\subsetneq(E_1,V_1,\rho_1)=(E,V,\rho)$ is a strictly descending chain of subobjects in $\acA$ and $\ac\mu^{\om,s}(\lb (E_n,V_n,\rho_n)/(E_{n+1},V_{n+1},\rho_{n+1})\rb)\le\ac\mu^{\om,s}(\lb E_{n+1},V_{n+1},\rho_{n+1}\rb)$ for $n\ge 1$. Then $0\le\dim V_{n+1}\le\dim V_n$, so for some $N\ge 0$ we have $V_n=V_N$ for all $n\ge N$. Thus $\ac\mu^{\om,s}(\lb E_n/E_{n+1},0,0\rb)\le\ac\mu^{\om,s}(\lb E_{n+1},V_N,\rho_{n+1}\rb)$ in $(\acM,\le)$ for all $n\ge N$. This implies that for $n\ge N$ we have
\begin{align*}
&\dim (E_n/E_{n+1})=\deg\ac\mu^{\om,s}(\lb E_n/E_{n+1},0,0\rb)
\ge\deg\ac\mu^{\om,s}(\lb E_{n+1},V_N,\rho_{n+1}\rb)\\
&=\left.\begin{cases} \dim E_{n+1}, & V_N=0 \\ \max(\dim E_{n+1},m-1), & V_N\ne 0\end{cases}\right\} \ge \dim E_{n+1}.
\end{align*}
But the existence of such a sequence $(E_n)_{n\ge N}$ contradicts $\coh(X)$ $\pi$-artinian in Lemma \ref{co7lem1}. Hence $\acA$ is $\ac\mu^{\om,s}$-artinian. Proving $\acA$ $\bar\mu^{\om,0}_1$-artinian is similar.
\end{proof}

\begin{lem}
\label{co8lem2}
Let\/ $(E,V,\rho)\in\acA$ be {\rm$\ac\mu^{\om,s}$-} or\/ $\bar\mu^{\om,0}_1$-semistable. Then $E$ is pure.
\end{lem}

\begin{proof} If $E$ is not pure there exists $0\ne E'\subsetneq E$ with $\dim E'<\dim E$. It is then easy to check that $(E',0,0)\subset(E,V,\rho)$ $\ac\mu^{\om,s}$- and $\bar\mu^{\om,0}_1$-destabilizes~$(E,V,\rho)$.
\end{proof}

\begin{lem}
\label{co8lem3}
For $(E,V,\rho)\in\acA$ to be {\rm$\ac\mu^{\om,s}$-} or\/ $\bar\mu^{\om,0}_1$-(semi)stable are open conditions on $[E,V,\rho]$ in $\acM$ and\/ $\acM^\pl$. Thus for each\/ $(\al,d)$ in $C(\acA)$ we have open substacks $\acM_{(\al,d)}^\rst(\ac\mu^{\om,s}),\ab\acM_{(\al,d)}^\ss(\ac\mu^{\om,s}),\ab\acM_{(\al,d)}^\rst(\bar\mu^{\om,0}_1),\ab\acM_{(\al,d)}^\ss(\bar\mu^{\om,0}_1)\subseteq\acM_{(\al,d)}^\pl$ parametrizing {\rm$\ac\mu^{\om,s}$-} and\/ $\bar\mu^{\om,0}_1$-(semi)stable objects in $\acA$ in class $(\al,d)$.
\end{lem}

\begin{proof} The proof follows that of Proposition \ref{co7prop7} in \S\ref{co7114}. If $\al=0$ then 
$[0,V,0]$ is always $\ac\mu^{\om,s}$- and $\bar\mu^{\om,0}_1$-semistable, and is $\ac\mu^{\om,s}$- and $\bar\mu^{\om,0}_1$-stable if and only if $d=1$, all four conditions being trivially open. So suppose $\al\ne 0$.

Write $\acM^{\rm pu}_{(\al,d)}\subseteq\acM_{(\al,d)}^\pl$ for the open substack of $[E,V,\rho]\in\acM_{(\al,d)}^\pl$ with $E$ pure of dimension $\dim\al$. Let $\smash{\ac\Up:U\ra\acM^{\rm pu}_{(\al,d)}}$ be an atlas with $U$ locally of finite type and $V\subseteq U$ be a finite type open subscheme. Then the family $\sF$ of $E\in\coh(X)$ with $[E,V,\rho]=\Up(v)$ for $\C$-points $v\in V$ is bounded. Such $(E,V,\rho)$ is $\ac\mu^{\om,s}$-unstable if and only if there exists a quotient $\pi:(E,V,\rho)\twoheadrightarrow (F,W,\si)$ with $\Ker\pi\ne 0$ and $F$ pure of dimension $\dim\al$ (as otherwise we can replace $F$ by $F/F'$ for $F'\subsetneq F$ maximal with $\dim F'<\dim F$) and $\ac\mu^{\om,s}(\lb\Ker\pi\rb)>\ac\mu^{\om,s}(\lb F,W,\si\rb)$. It is not $\ac\mu^{\om,s}$-stable if the same holds with $\ge$ rather than~$>$. 

We can use $\ac\mu^{\om,s}(\lb\Ker\pi\rb)\ge\ac\mu^{\om,s}(\lb F,W,\si\rb)$ to derive an a priori upper bound for $\bar\mu^\om(\lb F\rb)$, for instance when $\rank\al>0$ equation \eq{co8eq9} implies that $\bar\mu^\om(\lb F\rb)\le\bar\mu^\om(\al)+\frac{ds}{\rank\al}$. The argument in \S\ref{co7114} then implies that there are only finitely many possible classes $\lb F\rb\in C(\coh(X))$ in such destabilizing quotients $\pi:(E,V,\rho)\twoheadrightarrow (F,W,\si)$ with $[E,V,\rho]=\Up(v)$, $v\in V$. Hence there are only finitely many possible $\lb F,W,\si\rb\in C(\acA)$ as $\dim W\le d$. Openness of $\ac\mu^{\om,s}$-(semi)stability follows as in \S\ref{co7114}, using Lemma \ref{co8lem2}. The argument for $\bar\mu^{\om,0}_1$ is similar but easier, as the upper bound $\bar\mu^\om(\lb F\rb)\le\bar\mu^\om(\al)$ is sufficient.
\end{proof}

\begin{lem}
\label{co8lem4}
When\/ $\al\in C(\coh(X))$ and\/ $d=0$ the obvious inclusion of subcategories $\coh(X)\hookra\acA$ mapping $E\mapsto (E,0,0)$ induces isomorphisms $\M_\al^\ss(\mu^\om)\cong\acM_{(\al,0)}^\ss(\ac\mu^{\om,s})=\acM_{(\al,0)}^\ss(\bar\mu^{\om,0}_1),$ for\/ $\M_\al^\ss(\mu^\om)$ as in {\rm\S\ref{co71}} for $\coh(X)$.
\end{lem}

\begin{proof}
If $\al,\be\in C(\coh(X))$ it is easy to show that $\mu^\om(\al)\le\mu^\om(\be)$ if and only if $\ac\mu^{\om,s}(\al,0)\le\ac\mu^{\om,s}(\be,0)$ if and only if $\bar\mu^{\om,0}_1(\al,0)\le\bar\mu^{\om,0}_1(\be,0)$. Hence $E\in\coh(X)$ is $\mu^\om$-semistable if and only if $(E,0,0)$ is $\ac\mu^{\om,s}$-semistable if and only if $(E,0,0)$ is $\bar\mu^{\om,0}_1$-semistable, and the lemma follows.
\end{proof}

\begin{prop}
\label{co8prop1}
Let\/ $\al\in C(\coh(X))$ with\/ $\rank\al>0$ and\/ $\bar\mu^\om(\al)\ge\bar\mu^\om(\lb L\rb)$. For $s>0$ used to define $(\ac\mu^{\om,s},\ac M,\le)$ in Definition\/ {\rm\ref{co8def3},} define $u>0$ by
\e
u^{-1}=\ha\bigl(\bar\mu^\om(\al)-\bar\mu^\om(\lb L\rb)\bigr)+\frac{s}{2(\rank\al+1)}.
\label{co8eq11}
\e
Then we can form $\om\bp u\,\om_{\CP^1}\in\Kah(X\t\CP^1),$ where $\om_{\CP^1}\in\Kah(\CP^1)$ is normalized by $\int_{\CP^1}\om_{\CP^1}=1,$ and so Definition\/ {\rm\ref{co7def4}} gives a weak stability condition $(\mu^{\om\bp u\,\om_{\CP^1}},M,\le)$ on $\coh(X\t\CP^1)$. Let\/ $(E,V,\rho)\in\acA$ with\/ $\lb E,V,\rho\rb=(\al,1)$. Then $(E,V,\rho)$ is $\ac\mu^{\om,s}$-(semi)stable if and only if\/ $\ac I(E,V,\rho)$ is $\mu^{\om\bp u\,\om_{\CP^1}}$-(semi)stable, for $\ac I:\acA\ra\coh(X\t\CP^1)$ as in Theorem\/ {\rm\ref{co8thm1}}.
\end{prop}

\begin{proof} Let $(E,V,\rho)\in\acA$ with $\lb E,V,\rho\rb=(\al,1)$. First suppose $E$ is not torsion-free, so we have $0\ne E'\subsetneq E$ with $E'$ torsion. Then $(E,V,\rho)$ is not $\ac\mu^{\om,s}$-semistable by Lemma \ref{co8lem2}, and $\ac I(E,V,\rho)$ is not torsion-free in $\coh(X\t\CP^1)$ as $\ac I(E',0,0)$ is a torsion subobject, so $\ac I(E,V,\rho)$ is not $\mu^{\om\bp u\,\om_{\CP^1}}$-semistable. This proves the proposition if $E$ is not torsion-free.

Next suppose $E$ is torsion-free. Then all subobjects $(E',V',\rho')\subseteq(E,V,\rho)$ have $E'=0$ or $E'$ torsion-free and $\dim V'\le 1$, so $\lb E',V',\rho'\rb=(\be,e)$ with $\be=0$ or $\rank\be>0$ and $e=0$ or 1. If $0<\rank\be<\rank\al$ then by \eq{co8eq9} we may explicitly write the inequality $\bar\mu^{\om,s}(\be,e)>\bar\mu^{\om,s}(\al-\be,1-e)$ (which is equivalent to $\ac\mu^{\om,s}(\be,e)\!>\!\ac\mu^{\om,s}(\al\!-\!\be,1\!-\!e)$ as $\deg\ac\mu^{\om,s}(\be,e)\!=\!\deg\ac\mu^{\om,s}(\al\!-\!\be,1\!-\!e)=m$)~as
\e
\bar\mu^\om(\be)+\frac{es}{\rank\be}>\bar\mu^\om(\al-\be)+\frac{(1-e)s}{\rank(\al-\be)}.
\label{co8eq12}
\e
Similarly, we write $\bar\mu^{\om\bp u\,\om_{\CP^1}}\ci \ac I_*(\be,e)>\bar\mu^{\om\bp u\,\om_{\CP^1}}\ci \ac I_*(\al-\be,1-e)$ explicitly as
\e
\begin{split}
&\frac{\rank\be\,\bar\mu^\om(\be)+e\,\bar\mu^\om(\lb L\rb)+2eu}{\rank\be+e}+u\\
&>\frac{\rank(\al-\be)\,\bar\mu^\om(\al-\be)+(1-e)\,\bar\mu^\om(\lb L\rb)+2(1-e)u}{\rank(\al-\be)+1-e}+u.
\end{split}
\label{co8eq13}
\e
Using \eq{co8eq11} and $\rank\be\,\bar\mu^\om(\be)+\rank(\al-\be)\bar\mu^\om(\al-\be)=\rank\al\,\bar\mu^\om(\al)$ we can show that \eq{co8eq12} and \eq{co8eq13} are equivalent. The cases $(\be,e)=(0,1)$ and $\rank\be=\rank\al$ can be dealt with by similar arguments. 

Hence $\ac\mu^{\om,s}(\be,e)>\ac\mu^{\om,s}(\al-\be,1-e)$ if and only if $\mu^{\om\bp u\,\om_{\CP^1}}\ci \ac I_*(\be,e)>\mu^{\om\bp u\,\om_{\CP^1}}\ci \ac I_*(\al-\be,1-e)$, so $(E',V',\rho')$ $\ac\mu^{\om,s}$-destabilizes $(E,V,\rho)$ if and only if $\ac I(E',V',\rho')$ $\mu^{\om\bp u\,\om_{\CP^1}}$-destabilizes $\ac I(E,V,\rho)$. As $\ac I(E,V,\rho)$ is $\SU(2)$-invariant, if it is $\mu^{\om\bp u\,\om_{\CP^1}}$-unstable it must have a unique maximal $\mu^{\om\bp u\,\om_{\CP^1}}$-destabilizing subobject, which is also $\SU(2)$-invariant, and so is of the form $\ac I(E',V',\rho')$ for some $(E',V',\rho')\subseteq(E,V,\rho)$ by Theorem \ref{co8thm1}. Therefore $(E,V,\rho)$ is $\ac\mu^{\om,s}$-semistable if and only if $\ac I(E,V,\rho)$ is $\mu^{\om\bp u\,\om_{\CP^1}}$-semistable. The proof for stability is similar, by allowing equality in \eq{co8eq12}--\eq{co8eq13}. The lemma follows.
\end{proof}

\begin{prop}
\label{co8prop2}
Let\/ $(\al,d)\in C(\acA)$ with\/ $\rank\al>0$ and\/ $d\le 1$. Then
\begin{itemize}
\setlength{\itemsep}{0pt}
\setlength{\parsep}{0pt}
\item[{\bf(a)}] $\acM_{(\al,d)}^\ss(\ac\mu^{\om,s}),\ab \acM_{(\al,d)}^\ss(\bar\mu^{\om,0}_1)$ are of finite type.
\item[{\bf(b)}] If\/ $\acM_{(\al,d)}^\rst(\ac\mu^{\om,s})\!=\!\acM_{(\al,d)}^\ss(\ac\mu^{\om,s})$ then $\acM_{(\al,d)}^\ss(\ac\mu^{\om,s})$ is a proper algebraic space, and similarly for\/~$\bar\mu^{\om,0}_1$.
\end{itemize}	
\end{prop}

\begin{proof}
When $d=0$, parts (a),(b) follow from Lemma \ref{co8lem4}, and Proposition \ref{co7prop8} and Theorem \ref{co7thm5} for $\coh(X)$. When $d=1$, parts (a),(b) for $\ac\mu^{\om,s}$ follow from Proposition \ref{co8prop1}, and Proposition \ref{co7prop8} and Theorem \ref{co7thm5} for $\coh(X\t\CP^1)$, and the fact that the composition \eq{co8eq5} is of finite type. 

When $d=1$, part (a) for $\bar\mu^{\om,0}_1$ follows as for $\bar\mu^\la_{\bs\mu}$ on $\baA$ in \S\ref{co52}, and part (b) is proved as for Assumption \ref{co5ass2}(h) for $\coh(X)$ in \S\ref{co74}, applied to the category $\baA$ in Example \ref{co5ex1}: \S\ref{co74} considers the case $L=\O_X(-N)$ for $N\gg 0$, so that $E$ is $N$-regular for all $[E]\in\M_\al^\ss(\tau)$, but the $N$-regularity is not needed for Assumption \ref{co5ass2}(h), so the proof works for any~$L$.
\end{proof}

\begin{thm}
\label{co8thm2}
Let\/ $\al\in C(\coh(X))$ with\/ $\rank\al>0$. Then:
\begin{itemize}
\setlength{\itemsep}{0pt}
\setlength{\parsep}{0pt}
\item[{\bf(a)}] If\/ $[E,V,\rho]$ lies in $\acM_{(\al,1)}^\ss(\ac\mu^{\om,s})$ or $\acM^\ss_{(\al,1)}(\bar\mu^{\om,0}_1)$ then $\rho\ne 0$.
\item[{\bf(b)}] If\/ $\al=\lb L\rb$ then $\acM_{(\al,1)}^\ss(\ac\mu^{\om,s})=\acM^\ss_{(\al,1)}(\bar\mu^{\om,0}_1)=\{[L,\C,\id_L]\}$.
\item[{\bf(c)}] If\/ $\bar\mu^\om(\al)\le\bar\mu^\om(\lb L\rb)$ and\/ $\al\ne\lb L\rb$ then $\acM_{(\al,1)}^\ss(\ac\mu^{\om,s})=\acM^\ss_{(\al,1)}(\bar\mu^{\om,0}_1)=\es$ for all\/~$s>0$.
\item[{\bf(d)}] If\/ $\bar\mu^\om(\al)>\bar\mu^\om(\lb L\rb)$ then for small\/ $\ep>0$ and all\/ $s\in(0,\ep)$ we have
\begin{equation*}
\acM_{(\al,1)}^\rst(\ac\mu^{\om,s})=\acM_{(\al,1)}^\ss(\ac\mu^{\om,s})=\acM^\rst_{(\al,1)}(\bar\mu^{\om,0}_1)=\acM^\ss_{(\al,1)}(\bar\mu^{\om,0}_1).
\end{equation*}
\item[{\bf(e)}] If\/ $\rank\al=1$ then for all\/ $s>0$
\ea
\begin{split}
\!\!\!\!\acM_{(\al,0)}^\rst(\ac\mu^{\om,s})&=\acM_{(\al,0)}^\ss(\ac\mu^{\om,s})=\acM^\rst_{(\al,0)}(\bar\mu^{\om,0}_1)=\acM^\ss_{(\al,0)}(\bar\mu^{\om,0}_1)\\
&=\bigl\{[E,0,0]\in\acM_{(\al,1)}^\pl:\text{$E$ is torsion-free}\bigr\},
\end{split}
\label{co8eq14}\\
\begin{split}
\!\!\!\!\acM_{(\al,1)}^\rst(\ac\mu^{\om,s})&=\acM_{(\al,1)}^\ss(\ac\mu^{\om,s})=\acM^\rst_{(\al,1)}(\bar\mu^{\om,0}_1)=\acM^\ss_{(\al,1)}(\bar\mu^{\om,0}_1)\\
&=\bigl\{[E,V,\rho]\in\acM_{(\al,1)}^\pl:\text{$E$ is torsion-free, $\rho\ne 0$}\bigr\}.
\end{split}
\label{co8eq15}
\ea
\item[{\bf(f)}] If\/ $\bar\mu^\om(\al)>\bar\mu^\om(\lb L\rb)$ and\/ $\rank\al>1$ then $\acM_{(\al,1)}^\ss(\ac\mu^{\om,s})=\es$ whenever
\e
s>\frac{\rank\al}{\rank\al-1}.
\label{co8eq16}
\e
\item[{\bf(g)}] There are $0=s_0<s_1<\cdots<s_{k+1}=\iy$ such that\/ $\acM_{(\al,1)}^\rst(\ac\mu^{\om,s_i})\ne\acM_{(\al,1)}^\ss(\ac\mu^{\om,s_i})$ for\/ $i=1,\ldots,k,$ and if\/ $s\in(s_i,s_{i+1})$ for $i=0,\ldots,k$ then $\acM_{(\al,1)}^\rst(\ac\mu^{\om,s})=\acM_{(\al,1)}^\ss(\ac\mu^{\om,s})$ is independent of\/~$s\in(s_i,s_{i+1})$.
\end{itemize}
\end{thm}

\begin{proof} For (a), if $\rho=0$ then the subobject $(0,V,0)\subset(E,V,0)$ $\ac\mu^{\om,s}$- and $\bar\mu^{\om,0}_1$-destabilizes~$(E,V,\rho)$.

For (b), if $[E,V,\rho]$ lies in $\acM_{(\al,1)}^\ss(\ac\mu^{\om,s})$ or $\acM^\ss_{(\al,1)}(\bar\mu^{\om,0}_1)$ then $E$ is torsion-free by Lemma \ref{co8lem2} and $\rho\ne 0$ by (a) and $V\cong\C$, so $\rho:V\ot_\C L\cong L\ra E$ is injective, giving $\Im\rho\cong L$. Thus $\lb E/\Im\rho\rb=\al-\lb\Im\rho\rb=\al-\lb L\rb=0$, so $E/\Im\rho=0$ and $\rho$ is an isomorphism. Hence $[E,V,\rho]=[L,\C,\id_L]$. Clearly $(L,\C,\id_L)$ is $\ac\mu^{\om,s}$- and $\bar\mu^{\om,0}_1$-semistable, so (b) follows.

For (c), suppose for a contradiction that $[E,V,\rho]$ lies in $\acM_{(\al,1)}^\ss(\ac\mu^{\om,s})$ or $\acM^\ss_{(\al,1)}(\bar\mu^{\om,0}_1)$. The argument of (b) implies that $\Im\rho\cong L$, but $\Im\rho\ne E$ as $\al\ne\lb L\rb$. We cannot have $\dim (E/\Im\rho)=m-1$ as this contradicts $\bar\mu^\om(\al)\le\bar\mu^\om(\lb L\rb)$, and we cannot have $\dim (E/\Im\rho)<m-1$ as this contradicts $L$ a line bundle and $E$ torsion-free. Hence $\dim (E/\Im\rho)=m$, and~$\rank\al>1$.

As $\bar\mu^\om(\al)\le\bar\mu^\om(\lb L\rb)$ and $\rank\al>1$ we find that $\ac\mu^{\om,s}(\lb L\rb,1)>\ac\mu^{\om,s}(\al,1)$ and $\bar\mu^{\om,0}_1(\lb L\rb,1)>\bar\mu^{\om,0}_1(\al,1)$, so the subobject $(\Im\rho,V,\rho)\subsetneq(E,V,\rho)$ $\ac\mu^{\om,s}$- and $\bar\mu^{\om,0}_1$-destabilizes $(E,V,\rho)$, a contradiction.

For (d), define $\sF$ to be the family of $E\in\coh(X)$ such that $[E,V,\rho]$ is a point in $\acM_{(\al,1)}^\ss(\bar\mu^{\om,0}_1)$ or in $\acM_{(\al,1)}^\ss(\ac\mu^{\om,s})$ for some $s\in(0,1]$. Note that all such $E$ are torsion-free by Lemma \ref{co8lem2}. If $[E,V,\rho]\in\acM_{(\al,1)}^\ss(\bar\mu^{\om,0}_1)$ then $[E]\in\M_\al^\ss(\mu^\om)$, so the family of such $E$ is bounded by Proposition \ref{co7prop6}. 

Let $u_0>u_1>0$ be the values of $u$ in \eq{co8eq11} corresponding to $s=0,1$, which are well defined as $\bar\mu^\om(\al)>\bar\mu^\om(\lb L\rb)$. If $[E,V,\rho]\in\acM_{(\al,1)}^\ss(\ac\mu^{\om,s})$ then $[\ac I(E,V,\rho)]\in\M_{\ac I_*(\al,1)}^\ss(\mu^{\om\bp u\om_{\CP^1}})$ for $u$ as in \eq{co8eq11} by Proposition \ref{co8prop1}, where $u\in[u_1,u_0]$ as $s\in(0,1]$. Thus the family of such $\ac I(E,V,\rho)$ is bounded in $\coh(X\t\CP^1)$ by Proposition \ref{co7prop6} as $[u_1,u_0]$ is compact. From this we can deduce that the family of such $E$ is bounded in $\coh(X)$, so $\sF$ is bounded.

Consider subobjects $0\ne E'\subsetneq E$ for $E\in\sF$ such that $E/E'$ is torsion-free and $\bar\mu^\om(\lb E'\rb)\ge\bar\mu^\om(\al)$ (equivalently, $\bar\mu^\om(\lb E/E'\rb)\le\bar\mu^\om(\al)$). The family of such $E'$ is bounded by Proposition \ref{co7prop3}(c) as $\sF$ is bounded, so there are only finitely many possible classes $\lb E'\rb\in C(\coh(X))$. Therefore we can choose small $\ep\in(0,1]$ such that if $0\ne E'\subsetneq E$ is as above then 
\ea
\bar\mu^\om(\lb E'\rb)&>\bar\mu^\om(\al)\quad \text{if and only if}\quad \bar\mu^\om(\lb E'\rb)\ge \bar\mu^\om(\al)+\frac{\ep}{\rank\al},
\label{co8eq17}\\
\bar\mu^\om(\lb E'\rb)&\ge\bar\mu^\om(\al)\quad \text{if and only if}\quad \bar\mu^\om(\lb E'\rb)\ge \bar\mu^\om(\al)-\ep.
\label{co8eq18}
\ea

Suppose $[E,V,\rho]\in\acM^\pl_{(\al,1)}$. We will show that $[E,V,\rho]\in\acM_{(\al,1)}^\ss(\bar\mu^{\om,0}_1)$ if and only if $[E,V,\rho]\in\acM_{(\al,1)}^\ss(\ac\mu^{\om,s})$ for any (equivalently, for all) $s\in(0,\ep)$. If $E\notin\sF$ this is immediate as $[E,V,\rho]\notin\acM_{(\al,1)}^\ss(\bar\mu^{\om,0}_1)$ and $[E,V,\rho]\notin\acM_{(\al,1)}^\ss(\ac\mu^{\om,s})$ for $s\in(0,\ep)\subset(0,1]$. So suppose~$E\in\sF$.

The only possible subobjects of $(E,V,\rho)$ which could $\bar\mu^{\om,0}_1$- or $\ac\mu^{\om,s}$-destabilize $(E,V,\rho)$ are $(E',0,0)$ for $0\ne E'\subsetneq E$ (we exclude $E'=E$ since this does not destabilize, as $\bar\mu^\om(\al)>\bar\mu^\om(\lb L\rb)$) and $(E',V,\rho)$ for $0\ne E'\subsetneq E$ with $\Im\rho\subset E'$. Since $E$ is torsion-free, $E'$ is torsion-free. If $E/E'$ is not torsion-free then we may replace $E'$ by $\Ker\bigl(E\twoheadrightarrow (E/E')/$torsion$\bigr)$, enlarging $E'$. This does not decrease $\bar\mu^\om(\lb E'\rb)$, and so makes it more likely for $(E',0,0)$ or $(E',V,\rho)$ to destabilize $(E,V,\rho)$. Thus we may suppose that $E/E'$ is torsion-free.

Using \eq{co8eq9} we find that $(E',0,0)$ $\ac\mu^{\om,s}$-destabilizes $(E,V,\rho)$ if and only if
\e
\bar\mu^\om(\lb E'\rb)>\bar\mu^\om(\al)+\frac{s}{\rank\al}.
\label{co8eq19}
\e
Also $(E',0,0)$ $\bar\mu^{\om,0}_1$-destabilizes $(E,V,\rho)$ if and only~if
\e
\bar\mu^\om(\lb E'\rb)>\bar\mu^\om(\al).
\label{co8eq20}
\e
Equation \eq{co8eq17} implies that \eq{co8eq19} for any (equivalently, all) $s\in(0,\ep)$ is equivalent to \eq{co8eq20}.

Next suppose that $0\ne E'\subsetneq E$ with $\Im\rho\subseteq E'$. Then $(E',V,\rho)$ $\ac\mu^{\om,s}$-destabilizes $(E,V,\rho)$ if and only if $\rank E'<\rank\al$ and 
\e
\bar\mu^\om(\lb E'\rb)+\frac{s}{\rank\lb E'\rb}>\bar\mu^\om(\al)+\frac{s}{\rank\al}.
\label{co8eq21}
\e
Also $(E',V,\rho)$ $\bar\mu^{\om,0}_1$-destabilizes $(E,V,\rho)$ if and only if $\rank E'<\rank\al$ and
\e
\bar\mu^\om(\lb E'\rb)\ge \bar\mu^\om(\al).
\label{co8eq22}
\e
We exclude the case that $\bar\mu^\om(\lb E'\rb)>\bar\mu^\om(\al)$, as then $[E,V,\rho]\notin\acM_{(\al,1)}^\ss(\bar\mu^{\om,0}_1)$ and $[E,V,\rho]\notin\acM_{(\al,1)}^\ss(\ac\mu^{\om,s})$ from above. Then as $1\le\rank E'<\rank\al$, by \eq{co8eq18} equation \eq{co8eq21} for any (equivalently, all) $s\in(0,\ep)$ is equivalent to \eq{co8eq22}, as both are equivalent to $\bar\mu^\om(\lb E'\rb)=\bar\mu^\om(\al)$. This proves that $[E,V,\rho]\notin\acM_{(\al,1)}^\ss(\bar\mu^{\om,0}_1)$ if and only if $[E,V,\rho]\notin\acM_{(\al,1)}^\ss(\ac\mu^{\om,s})$ for any (equivalently, for all) $s\in(0,\ep)$, so $\acM_{(\al,1)}^\ss(\ac\mu^{\om,s})=\acM_{(\al,1)}^\ss(\bar\mu^{\om,0}_1)$. There cannot be strictly $\ac\mu^{\om,s}$-semistable points $[E,V,\rho]$ in $\acM_{(\al,1)}^\ss(\ac\mu^{\om,s})$ for small $s>0$, as this would require equality in \eq{co8eq19} or \eq{co8eq21}. The definition \eq{co8eq10} also prevents strictly $\bar\mu^{\om,0}_1$-semistable points $[E,V,\rho]$ in $\acM_{(\al,1)}^\ss(\bar\mu^{\om,0}_1)$. Part (d) follows.

For (e), let $\rank\al=1$. If $[E,V,\rho]$ lies in $\acM_{(\al,d)}^\ss(\ac\mu^{\om,s})$ or $\acM^\ss_{(\al,d)}(\bar\mu^{\om,0}_1)$ for $d=0$ or 1 then $E$ is torsion-free and $\rho\ne 0$ if $d=1$ by Lemma \ref{co8lem2} and (a). Conversely, if $[E,V,\rho]$ lies in $\acM_{(\al,d)}^\pl$ with $E$ torsion-free and $\rho\ne 0$ if $d=1$ then all subobjects $0\ne (E',V',\rho')\subsetneq(E,V,\rho)$ have $E'\ne 0$ (as $\rho\ne 0$ if $d=1$), so $\rank E'=\rank E=1$, and $\rank E/E'=0$. This forces $\deg\ac\mu^{\om,s}(\lb E',V',\rho'\rb)=m$ and $\deg\ac\mu^{\om,s}(\lb (E,V,\rho)/(E',V',\rho')\rb)<m$, so $\ac\mu^{\om,s}(\lb E',V',\rho'\rb)<\ac\mu^{\om,s}(\lb (E,V,\rho)/(E',V',\rho')\rb)$, and $(E,V,\rho)$ is $\ac\mu^{\om,s}$-stable. Similarly $(E,V,\rho)$ is $\bar\mu^{\om,0}_1$-stable. Equations \eq{co8eq14}--\eq{co8eq15} follow.

For (f), suppose $\bar\mu^\om(\al)>\bar\mu^\om(\lb L\rb)$, $\rank\al>1$, \eq{co8eq16} holds, and $[E,V,\rho]\in\acM_{(\al,1)}^\ss(\ac\mu^{\om,s})$. As in (b) we have $(\Im\rho,V,\rho)\subseteq(E,V,\rho)$ with $\lb\Im\rho,V,\rho\rb=(\lb L\rb,1)$. But \eq{co8eq16} is equivalent to $\bar\mu^{\om,s}(\lb L\rb,1)>\bar\mu^{\om,s}(\al,1)$ by \eq{co8eq9}, so $(\Im\rho,V,\rho)$ $\ac\mu^{\om,s}$-destabilizes $(E,V,\rho)$. Hence~$\acM_{(\al,1)}^\ss(\ac\mu^{\om,s})=\es$.

For (g), if $\bar\mu^\om(\al)\le\bar\mu^\om(\lb L\rb)$, or if $\rank\al=1$, then (g) follows from (b),(c),(e), with $k=0$. So suppose that $\bar\mu^\om(\al)>\bar\mu^\om(\lb L\rb)$ and $\rank\al>1$. Then (d),(f) show that $\acM_{(\al,1)}^\rst(\ac\mu^{\om,s})=\acM_{(\al,1)}^\ss(\ac\mu^{\om,s})$ is independent of $s$ for $s\in(0,\ep)$ and $s\in(\frac{\rank\al}{\rank\al-1},\iy)$. Arguing as in (d), but defining $\sF$ using $s\in(0,2\frac{\rank\al}{\rank\al-1}]$ rather than $s\in(0,1]$, we find that (supposing that $E\in\sF$, $E/E'$ is torsion-free, and $\bar\mu^\om(\lb E'\rb)\le\bar\mu^\om(\al)$ in the $(E',V,\rho)$ case, which as above are sufficient to determine if $[E,V,\rho]$ lies in $\acM_{(\al,1)}^\rst(\ac\mu^{\om,s})$ or $\acM_{(\al,1)}^\ss(\ac\mu^{\om,s})$), there are only finitely many $\lb E'\rb\in C(\coh(X))$ such that $(E',0,0)$ or $(E',V,\rho)$ can $\ac\mu^{\om,s}$-destabilize or $\ac\mu^{\om,s}$-semistabilize $(E,V,\rho)$, for $[E,V,\rho]\in\acM_{(\al,d)}^\pl$ and~$s\in(0,2\frac{\rank\al}{\rank\al-1}]$. 

For each such $\lb E'\rb$ there are at most two values of $s\in(0,\iy)$ where $\lb E'\rb$ can cause $\acM_{(\al,1)}^\rst(\ac\mu^{\om,s})$ or $\acM_{(\al,1)}^\ss(\ac\mu^{\om,s})$ to be non-constant as a function of $s$, those giving equality in \eq{co8eq19} or \eq{co8eq21}. Thus $\acM_{(\al,1)}^\rst(\ac\mu^{\om,s})$ or $\acM_{(\al,1)}^\ss(\ac\mu^{\om,s})$ are locally constant in $s\in(0,\iy)$ except at finitely many values of $s$, which we write as $0<s_1<\cdots<s_{k-1}<\iy$. At each $s_i$, the moduli spaces change by equality in \eq{co8eq19} for some $0\ne (E',0,0)\subsetneq(E,V,\rho)$, or in \eq{co8eq21} for some $0\ne (E',V,\rho)\subsetneq(E,V,\rho)$, and then $(E,V,\rho)$ is strictly $\ac\mu^{\om,s_i}$-semistable, so that $\acM_{(\al,1)}^\rst(\ac\mu^{\om,s_i})\ne\acM_{(\al,1)}^\ss(\ac\mu^{\om,s_i})$. If $s\in(s_i,s_{i+1})$ then $\acM_{(\al,1)}^\rst(\ac\mu^{\om,s})$ and $\acM_{(\al,1)}^\ss(\ac\mu^{\om,s})$ are independent of $s$, and as equality cannot hold in \eq{co8eq19} or \eq{co8eq21} then $\acM_{(\al,1)}^\rst(\ac\mu^{\om,s})=\acM_{(\al,1)}^\ss(\ac\mu^{\om,s})$. This completes the proof.
\end{proof}

\subsection{Verifying most of Assumptions \ref{co4ass1}, \ref{co5ass1}--\ref{co5ass3}}
\label{co82}

In this section we work in the situation of Definition \ref{co8def1}, and give the data and verify Assumptions \ref{co4ass1} and \ref{co5ass1}--\ref{co5ass3} for $\acA$, except for Assumptions \ref{co5ass1}(f) and \ref{co5ass2}(e) on obstruction theories, which are covered in~\S\ref{co83}.

\subsubsection{Assumption \ref{co4ass1} for $\acute{\mathcal A}$}
\label{co821}

Here is the analogue of \S\ref{co721} for $\acA$:

\begin{dfn}
\label{co8def5}
Work in the situation of Definition \ref{co8def1}, and write $\acB=\acA$. For Assumption \ref{co4ass1}(a), the moduli stack $\acM$ of objects in $\acB$ is described in Definition \ref{co8def2}, and is an Artin $\C$-stack, locally of finite type.

For Assumption \ref{co4ass1}(b) we have a morphism in $\coh(X\t\acM\t\acM)$:
\e
\begin{split}
&\Pi_{12}^*\bigl(\cR:L\bt\Pi_{\M_{\Vect_\C}}^*(\cV)\longra\Pi_{X\t\M}^*(\cU)\bigr)\op{}\\ 
&\Pi_{13}^*\bigl(\cR:L\bt\Pi_{\M_{\Vect_\C}}^*(\cV)\longra\Pi_{X\t\M}^*(\cU)\bigr).
\end{split}
\label{co8eq23}
\e
Here and below, for a product of stacks $S_1\t\cdots\t S_k$ we write $\Pi_i$ for the projection to the $i^{\rm th}$ factor, and $\Pi_{ij}$ for the projection to the product of the $i^{\rm th}$ and $j^{\rm th}$ factors, and so on.

Writing $S=\acM\t\acM$, equation \eq{co8eq23} is a morphism $\rho:L\bt V\ra\cE$ in $\coh(X\t S)$ with $V=(\Pi_{\M_{\Vect_\C}}\ci\Pi_1)^*(\cV)\op(\Pi_{\M_{\Vect_\C}}\ci\Pi_2)^*(\cV)$ locally free and $\cE=(\Pi_{X\t\M}\ci\Pi_{12})^*(\cU)\op(\Pi_{X\t\M}\ci\Pi_{12})^*(\cU)$ flat over $S$. Thus the universal property of $\cR$ in \eq{co8eq4} in Definition \ref{co8def2} gives a morphism $\ac\Phi:\acM\t\acM\ra\acM$, unique up to equivalence, such that 
\e
(\id_X\t\ac\Phi)^*(\cR)\cong\eq{co8eq23}.
\label{co8eq24}
\e
Restricting this identification to $\C$-points shows that $\ac\Phi$ acts on $\C$-points by $\ac\Phi_*:([E_1,V_1,\rho_1],\ab[E_2,V_2,\rho_2])\mapsto[E_1\op E_2,V_1\op V_2,\rho_1\op\rho_2]$. The action on isotropy groups in Assumption \ref{co4ass1}(b) is clear. 

Let $\ac\si:\acM\t\acM\ra\acM\t\acM$ exchange the factors. Then the pullback of \eq{co8eq23} by $\id_X\t\ac\si$ is isomorphic to \eq{co8eq23}, by swapping the two lines of \eq{co8eq23}. So by the construction of $\ac\Phi$ there is a 2-morphism $\ac\Phi\cong\ac\Phi\ci\ac\si$. That is, $\ac\Phi$ is commutative in $\Ho(\Art_\C)$. Associativity of $\ac\Phi$ follows by a similar argument on $X\t\acM\t\acM\t\acM$. This proves Assumption~\ref{co4ass1}(b).

For Assumption \ref{co4ass1}(c), we have a morphism in $\coh(X\t[*/\bG_m]\t\acM)$:
\e
\begin{split}
\Pi_2^*(\id_{L_{[*/\bG_m]}})\ot\Pi_{13}^*\bigl(\cR)
&:L\bt L_{[*/\bG_m]}\bt\Pi_{\M_{\Vect_\C}}^*(\cV)\\
&\quad\longra \Pi_2^*(L_{[*/\bG_m]})\ot\Pi_{X\t\M}^*(\cU).
\end{split}
\label{co8eq25}
\e
As for \eq{co8eq23}, this determines a morphism  such that
\e
(\id_X\t\ac\Psi)^*(\cR)\cong\eq{co8eq25}.
\label{co8eq26}
\e
On $\C$-points $\ac\Psi$ acts by $\ac\Psi_*:( *,[E,V,\rho])\mapsto[\C\ot E,\C\ot V,\id_\C\ot\rho]=[E,V,\rho]$. The action on isotropy groups in Assumption \ref{co4ass1}(c) is clear. 

Equation \eq{co4eq2} may similarly be deduced from the identifications
\begin{align*}
&\bigl(\Pi_2^*(\id_{L_{[*/\bG_m]}})\ot\Pi_{13}^*(\cR)\bigr)\op\bigl(\Pi_2^*(\id_{L_{[*/\bG_m]}})\ot\Pi_{14}^*(\cR)\bigr)\\
&\quad\cong \Pi_2^*(\id_{L_{[*/\bG_m]}})\ot\bigl(\Pi_{13}^*(\cR)\op\Pi_{14}^*(\cR)\bigr)\qquad\text{on $X\t[*/\bG_m]\t\acM\t\acM$,}\\
&\Pi_2^*(L_{[*/\bG_m]})\ot\bigl(\Pi_3^*(L_{[*/\bG_m]})\cong\Pi_{23}^*(\Om^*(L_{[*/\bG_m]}))\qquad\text{and}\\
&\Pi_2^*(\id_{L_{[*/\bG_m]}})\ot\bigl(\Pi_3^*(\id_{L_{[*/\bG_m]}})\ot\Pi_{14}^*(\cR)\bigr)\\
&\quad\cong(\Om\ci\Pi_{23})^*(\id_{L_{[*/\bG_m]}})\ot\Pi_{14}^*(\cR)\qquad\quad\;\text{on $X\t[*/\bG_m]\t[*/\bG_m]\t\acM$.}
\end{align*}
This proves Assumption~\ref{co4ass1}(c).

For Assumption \ref{co4ass1}(d), we define $K(\acA)=K(\acB)=K(\coh(X))\op\Z$ as in Definition \ref{co8def1}.
The map $\acM\ra K(\acB)$, $[E,V,\rho]\mapsto\lb E,V,\rho\rb$ is locally constant as in Definition \ref{co8def2}, so we have a decomposition $\acM=\coprod_{(\al,d)\in K(\acB)}\acM_{(\al,d)}$ of $\acM$ into open and closed $\C$-substacks $\acM_{(\al,d)}\subset \acM$ of objects in class $(\al,d)$. If $\lb E,V,\rho\rb=(0,0)$ then $\lb E\rb=0$, so $E=0$ as in \S\ref{co721}, and $\dim V=0$, so $V=0$, and $[E,V,\rho]=(0,0,0)$. Thus~$\acM_{(0,0)}=\{[0,0,0]\}$.

For Assumption \ref{co4ass1}(e), we define $\ac\chi:K(\acB)\t K(\acB)\ra\Z$ by
\e
\ac\chi\bigl((\al,d),(\be,e)\bigr)=\chi(\al,\be)+de-d\,\chi(\lb L\rb,\be),
\label{co8eq27}
\e
where $\chi:K(\coh(X))\t K(\coh(X))\ra\Z$ is as in Definition \ref{co7def1}.

For Assumption \ref{co4ass1}(f), define a perfect complex $\acExt^\bu$ on $\acM\t\acM$, the {\it Ext complex\/} in $\acA$, by the distinguished triangle on $\acM\t\acM$:
\e
\begin{gathered}
\!\!\!\!\!\!\!\!\!\xymatrix@R=18pt{
*+[r]{\acExt^\bu} \ar@<.2pc>[d] \\
*+[r]{\;\>\,\,\begin{subarray}{l}\ts 
(\Pi_{23})_*\bigl[(\Pi_1,\Pi_\M\ci\Pi_2)^*(\cU)^\vee\ot(\Pi_1,\Pi_\M\ci\Pi_3)^*(\cU)\bigr]
\op{}\\ \ts \bigl[(\Pi_{\M_{\Vect_\C}}\ci\Pi_1)^*(\cV)^\vee\ot (\Pi_{\M_{\Vect_\C}}\ci\Pi_2)^*(\cV)\bigr]\end{subarray}}
\ar@<.2pc>[d]^{(\Pi_{23})_*(\Pi_{12}^*(\cR^\vee)\ot\id_{(\Pi_1,\Pi_\M\ci\Pi_3)^*(\cU)})\op -\id_{(\Pi_{\M_{\Vect_\C}}\ci\Pi_1)^*(\cV)^\vee}\ot(\Pi_{23})_*(\Pi_{13}^*(\cR))}
\\
*+[r]{\;\>\,\,\begin{subarray}{l}\ts
(\Pi_{23})_*\bigl[\Pi_1^*(L^*)\ot(\Pi_{\M_{\Vect_\C}}\ci\Pi_2)^*(\cV)^\vee\ot(\Pi_1,\Pi_\M\ci\Pi_3)^*(\cU)\bigr]=\\
\ts
(\Pi_{\M_{\Vect_\C}}\ci\Pi_1)^*(\cV)^\vee\ot(\Pi_{23})_*\bigl[\Pi_1^*(L^*)\ot(\Pi_1,\Pi_\M\ci\Pi_3)^*(\cU)\bigr]
\end{subarray}}
\ar@<.2pc>[d]
\\
*+[r]{\vphantom{\ac{\mathcal E}}\acExt^\bu[1],}
}
\end{gathered}
\label{co8eq28}
\e
using derived duality, pushforward, pullback and tensor product functors on $X\t\acM\t\acM$ as in Huybrechts \cite{Huyb}, where $\Pi_{23}:X\t\acM\t\acM\ra\acM\t\acM$. 

Using the definition \eq{co7eq4} of $\cExt^\bu\ra\M\t\M$ we may rewrite \eq{co8eq28} as
\e
\begin{gathered}
\!\!\!\!\!\!\!\!\!\xymatrix@R=18pt{
*+[r]{\acExt^\bu} \ar@<.2pc>[d] \\
*+[r]{(\Pi_\M\t\Pi_\M)^*
(\cExt^\bu)\op
(\Pi_{\M_{\Vect_\C}}\t\Pi_{\M_{\Vect_\C}})^*(\cV^*\bt\cV)}
\ar@<.2pc>[d]^{\Pi_1^*(\cR^\vee)_*\op -\Pi_2^*(\cR)_* }
\\
*+[r]{\Pi_{\M_{\Vect_\C}}^*(\cV^*)\bt \Pi_\M^*(\cExt^\bu\vert_{\{[L]\}\t\M})}
\ar@<.2pc>[d]
\\
*+[r]{\vphantom{\ac{\mathcal E}}\acExt^\bu[1],}
}
\end{gathered}
\label{co8eq29}
\e
where $[L]\in\M$ with $\cU^\vee\vert_{X\t\{[L]\}}=L^*$, so that in the third lines of \eq{co8eq28}--\eq{co8eq29} $(\Pi_{23})_*\bigl[\Pi_1^*(L^*)\ot(\Pi_1,\Pi_\M\ci\Pi_3)^*(\cU)=\Pi_\M^*(\cExt^\bu\vert_{\{[L]\}\t\M})$, and we write $\Pi_1^*(\cR^\vee)_*,\Pi_2^*(\cR)_*$ as a shorthand for the morphisms in \eq{co8eq28}. If $(E_1,V_1,\rho_1),\ab(E_2,V_2,\rho_2)\in\acA$, restricting \eq{co8eq29} to $([E_1,V_1,\rho_1],[E_2,V_2,\rho_2])$ yields
\e
\begin{gathered}
\!\!\!\!\!\!\!\!\!\xymatrix@R=18pt{
*+[r]{\acExt^\bu\vert_{([E_1,V_1,\rho_1],[E_2,V_2,\rho_2])}} \ar@<.2pc>[d] \\
*+[r]{\cExt^\bu\vert_{([E_1],[E_2])}\op (V_1^*\ot_\C V_2)}
\ar@<.2pc>[d]^{(\rho_1^\vee\ot\id_{E_2})\op -(\id_{V_1^*}\ot\rho_2)}
\\
*+[r]{V_1^*\ot \cExt^\bu\vert_{([L],[E_2])}}
\ar@<.2pc>[d]
\\
*+[r]{\vphantom{\ac{\mathcal E}}\acExt^\bu[1]\vert_{([E_1,V_1,\rho_1],[E_2,V_2,\rho_2])},}
}
\end{gathered}
\label{co8eq30}
\e
where $\cExt^\bu\ra\M\t\M$ is as in \eq{co7eq4}, with $\rank\bigl(\cExt^\bu\vert_{\M_\al\t\M_\be}\bigr)=\chi(\al,\be)$. Writing $\lb E_1,V_1,\rho_1\rb=(\al,d)$ and $\lb E_2,V_2,\rho_2\rb=(\be,e)$, the second line of \eq{co8eq30} has rank $\chi(\al,\be)+de$ and the third line $d\,\chi(\lb L\rb,\be)$, so $\acExt^\bu\vert_{([E_1,V_1,\rho_1],[E_2,V_2,\rho_2])}$ has rank $\ac\chi\bigl((\al,d),(\be,e)\bigr)$ by \eq{co8eq27}, and 
\e
\rank\bigl(\acExt^\bu\vert_{\acM_{(\al,d)}\t\acM_{(\be,e)}}\bigr)=\ac\chi\bigl((\al,d),(\be,e)\bigr).
\label{co8eq31}
\e

Define $\acE^\bu,\acE^\bu_{(\al,d),(\be,e)}$ to be the dual perfect complexes
\e
\acE^\bu=(\acExt^\bu)^\vee,\qquad \acE_{(\al,d),(\be,e)}^\bu=(\acExt^\bu)^\vee\vert_{\acM_{(\al,d)}\t\acM_{(\be,e)}},
\label{co8eq32}
\e
so that $\rank\acE^\bu_{(\al,d),(\be,e)}=\ac\chi((\al,d),(\be,e))$ by \eq{co8eq31}. Equations \eq{co4eq3}--\eq{co4eq4} for $\acE^\bu$ follow easily from \eq{co8eq24}, \eq{co8eq28} and \eq{co8eq32}, and \eq{co4eq5}--\eq{co4eq6} follow from \eq{co8eq26}, \eq{co8eq28} and~\eq{co8eq32}.

As $\chi$ is nondegenerate on $K(\coh(X))=K^\num(\coh(X))$ by definition, $\ac\chi$ in \eq{co8eq27} is nondegenerate on $K(\acB)=K(\coh(X))\op\Z$. Thus Assumption \ref{co4ass1}(g) holds automatically. This completes Assumption \ref{co4ass1} for~$\acB$.
\end{dfn}

\subsubsection{Assumption \ref{co5ass1}(a)--(e),(g) for $\acute{\mathcal A}$}
\label{co822}

Here is the analogue of \S\ref{co722} for $\acA$:

\begin{dfn}
\label{co8def6}
Work in the situation of Definitions \ref{co8def1} and \ref{co8def5}. We will verify Assumption \ref{co5ass1}(a)--(e),(g) for $\acB=\acA$. Assumption \ref{co5ass1}(a) is trivial as $\acB=\acA$. For Assumption \ref{co5ass1}(b) we take $K(\acA)=K(\coh(X))\op\Z$ as in Definition \ref{co8def1}. Assumption \ref{co5ass1}(c) was proved in Definition~\ref{co8def5}. 

For Assumption \ref{co5ass1}(d), the (projective linear) derived moduli stacks $\bs\acM,\ab\bs\acM^\pl$ of objects in $\acA$ exist as in Definition \ref{co8def2}. They are derived Artin $\C$-stacks with classical truncations $t_0(\bs\acM)=\acM$, $t_0(\bs\acM^\pl)=\acM^\pl$, with splittings $\bs\acM=\coprod_{(\al,d)\in K(\acB)}\bs\acM_{(\al,d)}$, $\bs\acM^\pl=\coprod_{(\al,d)\in K(\acB)}\bs\acM_{(\al,d)}^\pl$ as in Assumption \ref{co5ass1}(d)(i).

For Assumption \ref{co5ass1}(d)(ii), in Definition \ref{co8def5} we constructed $\ac\Phi:\acM\t\acM\ra\acM$ to satisfy \eq{co8eq24}, using the universal property of $\cR$ in \eq{co8eq4} in Definition \ref{co8def2}. The same construction works for derived stacks using the universal property of $\bs\cR$ in \eq{co8eq6} giving a morphism $\bs{\ac\Phi}:\bs\acM\t\bs\acM\ra\bs\acM$, unique up to equivalence, with an isomorphism as for~\eq{co8eq23}--\eq{co8eq24}:
\begin{align*}
(\id_X\t\bs{\ac\Phi})^*(\bs\cR)\cong\,
&\Pi_{12}^*\bigl(\bs\cR:L\bt\Pi_{\M_{\Vect_\C}}^*(\cV)\longra\Pi_{X\t\bs\M}^*(\bs\cU)\bigr)\op{}\\ 
&\Pi_{13}^*\bigl(\bs\cR:L\bt\Pi_{\M_{\Vect_\C}}^*(\cV)\longra\Pi_{X\t\bs\M}^*(\bs\cU)\bigr).
\end{align*}
Under the map $\id\t i:X\t\acM\ra X\t\bs\acM$ we have $(\id\t i)^*(\bs\cR)=\cR$, so the universal property of $\bs\cR$ implies that $\bs{\ac\Phi}\ci(i\t i)=i\ci\ac\Phi$. Thus taking classical truncations shows that $t_0(\bs{\ac\Phi})=\ac\Phi$.

Similarly, the construction of $\ac\Psi:[*/\bG_m]\t\acM\ra\acM$ in Definition \ref{co8def5} lifts to $\bs{\ac\Psi}:[*/\bG_m]\t\bs\acM\ra\bs\acM$ with $t_0(\bs{\ac\Psi})=\ac\Psi$. The proofs of the identities satisfied by $\ac\Phi,\ac\Psi$ lift to $\bs{\ac\Phi},\bs{\ac\Psi}$. Thus $\bs{\ac\Psi}$ is an action of the group stack $[*/\bG_m]$ on the derived stack $\bs\acM$, which is free except over $0$. So, as claimed in Definition \ref{co8def2}, we can take the quotient $\bs\acM^\pl$, with projection $\bs{\ac\Pi}{}^\pl:\bs\acM\ra\bs\acM^\pl$, which is a principal $[*/\bG_m]$-bundle except over 0, with $t_0(\bs\acM^\pl)=\acM^\pl$ as~$t_0(\bs{\ac\Psi})=\ac\Psi$.

For Assumption \ref{co5ass1}(d)(iii), Definition \ref{co8def2} gives a morphism $(\Pi_{\bs\M},\Pi_{\M_{\Vect_\C}}):\bs\acM\ra\bs\M\t\M_{\Vect_\C}$. Here $\bs\M$ is locally finitely presented as in Definition \ref{co7def10}, and $\M_{\Vect_\C}$ is locally finitely presented as it is a smooth Artin $\C$-stack. This $(\Pi_{\bs\M},\Pi_{\M_{\Vect_\C}})$ realizes $\bs\acM$ as the total space of 
\e
\Pi_{\M_{\Vect_\C}}^*(\cV^*)\ot([L],\Pi_{\bs\M})^*(\bs\cExt^\bu),
\label{co8eq33}
\e
which is a perfect complex in the interval $[0,m]$ on $\bs\M\t\M_{\Vect_\C}$. It follows that the morphism $(\Pi_{\bs\M},\Pi_{\M_{\Vect_\C}})$ is locally finitely presented, so $\bs\acM$ is locally finitely presented as $\bs\M\t\M_{\Vect_\C}$ is. Hence $\bs\acM^\pl=\bs\acM/[*/\bG_m]$ is too.

For Assumption \ref{co5ass1}(d)(iv), as in \eq{co2eq11} the morphism $(\Pi_{\bs\M},\Pi_{\M_{\Vect_\C}}):\bs\acM\ra\bs\M\t\M_{\Vect_\C}$ induces a distinguished triangle of tangent complexes:
\begin{equation*}
\xymatrix@C=15pt{ \bT_{\bs\acM}  \ar[r] & \Pi_{\bs\M}^*(\bT_{\bs\M})\op \Pi_{\M_{\Vect_\C}}^*(\bT_{\M_{\Vect_\C}}) \ar[r] & \bT_{\bs\acM/\bs\M\t\M_{\Vect_\C}}[1] \ar[r]^(0.8){[+1]} & . }
\end{equation*}
Pulling this back by $\ac\imath:\acM\hookra\bs\acM$ gives
\e
\xymatrix@C=20pt{ \ac\imath^*(\bT_{\bs\acM})  \ar[r] & {\begin{subarray}{l} \ts\;\> \Pi_{\M}^*(i^*(\bT_{\bs\M}))\op \\ \ts \Pi_{\M_{\Vect_\C}}^*(\bT_{\M_{\Vect_\C}})\end{subarray}} \ar[r] &  \ac\imath^*(\bT_{\bs\acM/\bs\M\t\M_{\Vect_\C}})[1] \ar[r]^(0.8){[+1]} & . }
\label{co8eq34}
\e

Equation \eq{co7eq42} gives an isomorphism
\e
i^*(\bT_{\bs\M})\cong\De_\M^*(\cExt^\bu)[1].
\label{co8eq35}
\e
Since $\M_{\Vect_\C}=\coprod_{n\ge 0}[*/\GL(n,\C)]$ we see that
\e
\bT_{\M_{\Vect_\C}}\cong\cV^*\ot\cV[1].
\label{co8eq36}
\e
As $(\Pi_{\bs\M},\Pi_{\M_{\Vect_\C}}):\bs\acM\ra\bs\M\t\M_{\Vect_\C}$ is the derived total space of \eq{co8eq33}, we see that
\e
\ac\imath^*(\bT_{\bs\acM/\bs\M\t\M_{\Vect_\C}})\cong\Pi_{\M_{\Vect_\C}}^*(\cV^*)\ot \Pi_\M^*(\cExt^\bu\vert_{\{[L]\}\t\M}).
\label{co8eq37}
\e
Combining the last three equations with \eq{co8eq34} gives a distinguished triangle 
\e
\xymatrix@C=20pt{ \ac\imath^*(\bT_{\bs\acM})  \ar[r] & {\begin{subarray}{l} \ts\Pi_{\M}^*(\De_\M^*(\cExt^\bu))[1]\op \\ \ts \Pi_{\M_{\Vect_\C}}^*(\cV^*\ot\cV)[1]\end{subarray}} \ar[r]^\eta &  {\begin{subarray}{l} \ts \;\> \Pi_{\M_{\Vect_\C}}^*(\cV^*)\ot {} \\ \ts \Pi_\M^*(\cExt^\bu\vert_{\{[L]\}\t\M})[1]\end{subarray}} \ar[r]^(0.8){[+1]} & . }
\label{co8eq38}
\e

Let us compare \eq{co8eq38} with the shift $[1]$ of the diagonal pullback $\De_\acM^*$ of \eq{co8eq29}. The second and third terms agree, and the natural morphism $\eta$ in \eq{co8eq38} also agrees with $\De_\acM^*(\Pi_1^*(\cR^\vee)_*\op -\Pi_2^*(\cR)_*)[1]$. Thus properties of distinguished triangles give an isomorphism $\ac\imath^*(\bT_{\bs\acM})\cong \De_\acM^*(\acExt^\bu)[1]$ analogous to \eq{co7eq42}. Taking duals, restricting to $\acM_{(\al,d)}$ for $(\al,d)\in C(\acB)$, and using \eq{co8eq32} gives an isomorphism $\ac\th_{(\al,d)}$ as in \eq{co5eq1}:
\e
\ac\th_{(\al,d)}:\De_{\acM_{(\al,d)}}^*(\acE^\bu_{(\al,d),(\al,d)})[-1]\,{\buildrel\cong\over\longra}\, \ac\imath_{(\al,d)}^*(\bL_{\bs\acM_{(\al,d)}}).
\label{co8eq39}
\e  

For Assumption \ref{co5ass1}(d)(v), consider the large diagram Figure \ref{co8fig1}. Here the left hand diagonals come from pullbacks of \eq{co8eq29}, and the right hand diagonals from pullbacks of \eq{co8eq34}. The rows are all isomorphisms, and come from pullbacks of \eq{co8eq35}--\eq{co8eq37} and \eq{co8eq39}. The top left columns are also isomorphisms, and come from \eq{co4eq3}--\eq{co4eq4} for $\cE^\bu,\acE^\bu$ and~$ (\Phi_{d,e}^{\Vect_\C})^*(\cV_{d+e})\cong\cV_d\op\cV_e$.

\begin{figure}[htbp]
\text{\begin{footnotesize}$\displaystyle
\xymatrix@!0@C=36pt@R=46pt{ *+[r]{\begin{subarray}{l}\ts (\De_{\acM_{(\al+\be,d+e)}}\ci\ac\Phi_{(\al,d),(\be,e)})^* \\ \ts (\acExt_{(\al+\be,d+e),(\al+\be,d+e)}^\bu)[1]\end{subarray}} 
\ar[dr]
&&&&&& 
*+[r]{\begin{subarray}{l}\ts \ac\Phi_{(\al,d),(\be,e)}^* \\ \ts (i^*(\bT_{\bs\acM_{(\al+\be,d+e)}}))\end{subarray}} 
\ar[llllll]^(0.3){ \begin{subarray}{l} \ac\Phi_{(\al,d),(\be,e)}^* \\ (\ac\th_{(\al+\be,d+e)})^\vee \end{subarray}}
\ar[dr]
\\
& *+[r]{\!\!\!\!\!\!\!\begin{subarray}{l} \ts \ac\Phi_{(\al,d),(\be,e)}^* \\ \ts (\Pi^*(\De_{\M_{\al+\be}}^*(\cExt_{\al+\be,\al+\be}^\bu))[1] \\ \ts \op\Pi^*(\cV_{d+e}^*\ot\cV_{d+e})[1]) \end{subarray}}   \ar[dr] &&&&&& {\begin{subarray}{l} \ts \ac\Phi_{(\al,d),(\be,e)}^*(\Pi^*(i^*(\bT_{\bs\M_{\al+\be}})) \\ \ts \op\Pi^*(\bT_{\M^{d+e}_{\Vect_\C}}))) \end{subarray}} \ar[dr] \ar[llllll]^(0.4){\eq{co8eq35}\op\eq{co8eq36}}
\\
&& *+[r]{\!\!\!\!\!\!\!\!\!\!\!\!\!\!\begin{subarray}{l} \ts \ac\Phi_{(\al,d),(\be,e)}^* (\Pi^*(\cV_{d+e}^*)\bt {} \\ \ts \Pi^*(\cExt_{\lb L\rb,\al+\be}^\bu\vert_{[L]\t\M_{\al+\be}})[1]) \end{subarray}} \ar[dr]^(0.7){[+1]} &&&&&& {\begin{subarray}{l} \ts  \ac\Phi_{(\al,d),(\be,e)}^* \ci \ac\imath^* \\
\ts (\bT_{\begin{subarray}{l} \sst \bs\acM_{(\al+\be,d+e)}/ \\ \sst \bs\M_{\al+\be}\t\M^{d+e}_{\Vect_\C}\end{subarray}})[1]
\end{subarray}\qquad} \ar[dr]^(0.7){[+1]} \ar[llllll]^(0.45){\eq{co8eq37}}
\\
&&&&&&&&&
\\
*+[r]{\begin{subarray}{l}\ts (\De_{\acM_{(\al,d)}}^*(\acExt_{(\al,d),(\al,d)}^\bu)[1]\bp\De_{\acM_{(\be,e)}}^*(\acExt_{(\be,e),(\be,e)}^\bu)[1]) \\
\ts {}\op \acExt_{\al,\be}^\bu[1]\op \si_{\al,\be}^*(\acExt_{\be,\al}^\bu)[1] 
\end{subarray}}
\ar[dr]
\ar[uuuu]_{\begin{subarray}{l} \text{from} \\ \eq{co4eq3}, \\ \eq{co4eq4} \\ \text{for} \\ \acE^\bu \end{subarray}}^\cong 
\\
& *+[r]{\!\!\!\!\begin{subarray}{l}  \ts ((\Pi^*(\De_{\M_\al}^*(\cExt_{\al,\al}^\bu))[1] \!\op\!\Pi^*(\cV_d^*\ot\cV_d)[1])\bp \\ \ts ((\Pi^*(\De_{\M_\be}^*(\cExt_{\be,\be}^\bu))[1] \!\op\!\Pi^*(\cV_e^*\ot\cV_e)[1])\op  \\ 
\ts ((\Pi^*(\cExt_{\al,\be}^\bu)\!\op\!\Pi^*(\cV_d^*\!\ot\!\cV_e))[1]\!\op\!((\Pi^*(\si^*(\cExt_{\be,\al}^\bu))\!\op\!\Pi^*(\cV_e^*\!\ot\!\cV_d))[1]
\end{subarray}} \ar[dr] \ar[uuuu]_{\begin{subarray}{l} \text{from} \\ \eq{co4eq3}, \\ \eq{co4eq4} \\ \text{and} \\ (\Phi_{d,e}^{\Vect_\C})^* \\ (\cV_{d+e})\cong \\ \cV_d\op\cV_e
\end{subarray}}^\cong
\\
&& *+[r]{\!\!\!\!\!\!\!\begin{subarray}{l} \ts ((\Pi^*(\cV_d^*)\!\ot\! \Pi^*(\cExt_{\lb L\rb,\al}^\bu\vert_{[L]\t\M_\al})[1]))\bp \\  \ts ((\Pi^*(\cV_e^*)\!\ot\! \Pi^*(\cExt_{\lb L\rb,\be}^\bu\vert_{[L]\t\M_\be})[1]))\op  \\
\ts (\Pi^*(\cV_d^*)\!\bt\! \Pi^*(\cExt_{\lb L\rb,\be}^\bu\vert_{[L]\t\M_\be})[1])\! \op\!( \Pi^*(\cExt_{\lb L\rb,\al}^\bu\vert_{[L]\t\M_\al})[1])\!\bt\!\Pi^*(\cV_e^*)) \end{subarray}}  \ar[dr]^(0.7){[+1]} \ar[uuuu]_(0.71){\begin{subarray}{l} \text{from} \\ \eq{co4eq4}  \\ \text{and} \\ (\Phi_{d,e}^{\Vect_\C})^* \\ (\cV_{d+e})\cong \\ \cV_d\op\cV_e \end{subarray}}^(0.85)\cong
\\
&&&
\\
*+[r]{\begin{subarray}{l}\ts \De_{\acM_{(\al,d)}}^*(\acExt_{(\al,d),(\al,d)}^\bu)[1]\\
\ts {}\bp\! \De_{\acM_{(\be,e)}}^*(\acExt_{(\be,e),(\be,e)}^\bu)[1]
\end{subarray}} 
\ar[dr] 
\ar[uuuu]_(0.53){\begin{subarray}{l} \text{include} \\ \text{first two} \\ \text{factors} \end{subarray}}
&&&&&&
*+[r]{\begin{subarray}{l}\ts \; i^*(\bT_{\bs\acM_{(\al,d)}\t\bs\acM_{(\be,e)}})\cong \\ 
\ts i^*(\bT_{\bs\acM_{(\al,d)}})\bp i^*(\bT_{\bs\acM_{(\be,e)}})
\end{subarray}}   
\ar[dr] 
\ar[llllll]_(0.3){(\th_{(\al,d)})^\vee\bp (\th_{(\be,e)})^\vee }
\ar[uuuuuuuu]^(0.65){i^*(\bT_{\bs{\ac\Phi}_{(\al,d),(\be,e)}})\!\!\!\!\!\!}
\\
& *+[r]{\!\!\!\!\!\!\!\!\!\!\!\!\!\!\!\!\!\!\!\!\!\!\!\!\begin{subarray}{l}  \ts ((\Pi^*(\De_{\M_\al}^*(\cExt_{\al,\al}^\bu))[1] \!\op\!\Pi^*(\cV_d^*\ot\cV_d)[1]) \\ \ts \bp((\Pi^*(\De_{\M_\be}^*(\cExt_{\be,\be}^\bu))[1] \!\op\!\Pi^*(\cV_e^*\ot\cV_e)[1]) \end{subarray}}  \ar[dr] \ar[uuuu]_(0.53){\begin{subarray}{l} \text{include} \\ \text{first two} \\ \text{factors} \end{subarray}} &&&&&& {\begin{subarray}{l}  \ts 
\; (\Pi^*(i^*(\bT_{\bs\M_\al}))\op\Pi^*(\bT_{\M^d_{\Vect_\C}}))) \\ \ts \bp(\Pi^*(i^*(\bT_{\bs\M_\be}))\op\Pi^*(\bT_{\M^e_{\Vect_\C}})))\end{subarray}\!\!\!\!\!\!} 
\ar[llllll]_(0.35){\begin{subarray}{l} \eq{co8eq35}\op \\ \eq{co8eq36}\end{subarray}} \ar[dr]  
\ar[uuuuuuuu]^(0.7){\begin{subarray}{l} \Pi^*\ci i^*(\bT_{\bs\Phi_{\al,\be}}) \op \\ \Pi^*\ci i^*(\bT_{\Phi^{\Vect_\C}_{d,e}}) \\ 
\end{subarray}\!\!\!\!\!\!\!\!\!\!\!\!\!\!\!\!\!\!\!\!\!\!}
\\
&& *+[r]{\!\!\!\!\!\!\!\!\!\!\!\!\!\!\!\!\!\!\!\!\!\!\!\!\!\!\!\!\!\!\!\!\!\!\!\!\!\!\!\!\begin{subarray}{l} \ts ((\Pi^*(\cV_d^*)\!\bt\! \Pi^*(\cExt_{\lb L\rb,\al}^\bu\vert_{[L]\t\M_\al})[1])) \\  \ts \bp((\Pi^*(\cV_e^*)\!\bt\! \Pi^*(\cExt_{\lb L\rb,\be}^\bu\vert_{[L]\t\M_\be})[1])) \end{subarray}}  \ar[dr]^(0.7){[+1]} \ar[uuuu]_(0.68){\begin{subarray}{l} \text{include} \\ \text{first two} \\ \text{factors} \end{subarray}} &&&&&& {\begin{subarray}{l}\ts \; \ac\imath^*(\bT_{\bs\acM_{(\al,d)}/\bs\M_\al\t\M^d_{\Vect_\C}})[1] \\ 
\ts {}\bp\ac\imath^*(\bT_{\bs\acM_{(\be,e)}/\bs\M_\be\t\M^e_{\Vect_\C}})[1]
\end{subarray}\qquad\qquad\qquad} \ar[dr]^(0.7){[+1]} \ar[llllll]_{\eq{co8eq37}} 
\ar[uuuuuuuu]_(0.72){\!\!\!\!\!\!\!\!\!\!\!\!\!\!\!\!\!\!\!\!\!\!\!\!\! i^*(\bT_{\begin{subarray}{l} \bs{\ac\Phi}_{(\al,d),(\be,e)}/ \\ \bs\Phi_{\al,\be}\t \Phi^{\Vect_\C}_{d,e} \end{subarray}})}
\\
&&&&&&&&&
}\!\!\!\!\!\!\!\!\!\!\!\!\!\!\!\!\!\!\!\!\!\!\!\!\!\!\!\!\!\!\!\!\!\!\!\!\!\!\!\!\!\!\!\!\!\!\!\!\!\!\!\!\!\!\!\!\!\!\!\!\!\!\!\!\!\!\!\!\!\!\!\!\!\!\!\!\!\!\!\!\!\!\!\!\!\!\!$\end{footnotesize}}
\caption{Proof of Assumption \ref{co5ass1}(d)(v)}
\label{co8fig1}
\end{figure}

In Figure \ref{co8fig1} we have three large rectangles, at top left, centre and bottom right. The top left large rectangle is the dual of \eq{co5eq3} for $\acA$, so our goal is to show that this commutes. The centre large rectangle commutes by \eq{co5eq3} for $\coh(X)$, proved in Definition \ref{co7def10}, and obvious facts about $\Vect_\C$. The bottom right large rectangle commutes as \eq{co8eq37} comes from viewing $\bs\acM\ra\bs\M\t\M_{\Vect_\C}$ as the total space of \eq{co8eq33}, and pulling \eq{co8eq33} back along $\bs\Phi\t\Phi^{\Vect_\C}$ gives the exterior direct sum of \eq{co8eq33} with itself. The thin parallelograms at the top, bottom, left and right sides of Figure \ref{co8fig1} also commute, for straightforward reasons. Thus properties of triangulated categories imply that the top left large rectangle commutes, proving Assumption~\ref{co5ass1}(d)(v).

For Assumption \ref{co5ass1}(e), set
\e
C(\acB)_\pe=\bigl\{(\al,d)\in C(\acB):\text{$\rank\al>0$ and $d\le 1$}\bigr\}.
\label{co8eq40}
\e

For Assumption \ref{co5ass1}(g), let $\bigl\{(\B_k,F_k,\la_k):k\in K\bigr\}$ be as in \S\ref{co722} for $\A=\coh(X\t\CP^1)$. Define data $\bigl\{(\acB_k,\ac F_k,\ac\la_k):k\in K\bigr\}$ for $\acA$ by setting $\acB_k\subset\acA$ to be the full subcategory of $(E,V,\rho)$ with $\ac I(E,V,\rho)\in\B_k$, and $\ac F_k=F_k\ci\ac I:\acA\ra\Vect_\C$, and $\ac\la_k=\la_k\ci\ac I_*:K(\acB_k)\ra\Z$. As $\bigl\{(\B_k,F_k,\la_k):k\in K\bigr\}$ satisfies Assumption \ref{co5ass1}(g) for $\A=\coh(X\t\CP^1)$ by \S\ref{co722}, it is easy to check that $\bigl\{(\acB_k,\ac F_k,\ac\la_k):k\in K\bigr\}$ satisfies Assumption \ref{co5ass1}(g) for~$\acA$. 
\end{dfn}

\begin{rem}
\label{co8rem2}
As $\cExt^\bu$ is perfect in the interval $[0,m]$, where $m=\dim X$, equation \eq{co8eq29} implies that $\acExt^\bu$ is perfect in $[0,m+1]$. Let $(\al,d)\in C(\acA)$. Then \eq{co8eq32}, \eq{co8eq39} and $\acExt^\bu$ perfect in $[0,m+1]$ yield an analogue of~\eq{co7eq48}:
\e
\text{$i_{(\al,d)}^*(\bL_{\bs\acM_{(\al,d)}})$ is perfect in the interval $[-m,1]$.}
\label{co8eq41}
\e
Using the argument of Remark \ref{co7rem7} we can prove the analogue of \eq{co7eq49}:
\e
\text{$i_{(\al,d)}^*(\bL_{\bs\acM_{(\al,d)}})$ is perfect in $[a,1]$ if and only if $(i_{(\al,d)}^\pl)^*(\bL_{\bs\acM_{(\al,d)}^\pl})$ is.}
\label{co8eq42}
\e
Combining \eq{co8eq41}--\eq{co8eq42} yields
\begin{equation*}
\text{$(i_{(\al,d)}^\pl)^*(\bL_{\bs\acM_{(\al,d)}^\pl})$ is perfect in the interval $[-m,1]$.}
\end{equation*}
Since $\bs\acM_{(\al,d)},\bs\acM_{(\al,d)}^\pl$ are locally finitely presented as in Definition \ref{co8def6}, from Remark \ref{co2rem6}(e) we deduce that $\bs\acM_{(\al,d)}$ and $\bs\acM_{(\al,d)}^\pl$ are quasi-smooth if~$m=1$.
\end{rem}

\subsubsection{Assumptions \ref{co5ass2}(a)--(d),(f)--(h) and \ref{co5ass3}(a) for $\acute{\mathcal A}$}
\label{co823}

Here are the analogues of \S\ref{co723} and \S\ref{co74} for $\acA$. We continue in the situation of \S\ref{co821}--\S\ref{co822}. Fix $\om\in\Kah(X)$. For $(\ac\mu^{\om,s},\ac M,\le)$ as in Definition \ref{co8def3}, define the set $\acS$ of weak stability conditions on $\acA$ in Assumption \ref{co5ass2} by   
\e
\acS=\bigl\{(\ac\mu^{\om,s},\ac M,\le):s\in(0,\iy)\bigr\}.
\label{co8eq43}
\e
We now prove Assumptions \ref{co5ass2}(a)--(d),(f)--(g) and~\ref{co5ass3}.
\smallskip

\noindent{\bf Assumption \ref{co5ass2}(a).} The first part follows from Lemma \ref{co8lem1}, and the second part is trivial as $\acB=\acA$.
\smallskip

\noindent{\bf Assumption \ref{co5ass2}(b).} Follows from \eq{co8eq40}, Lemma \ref{co8lem3} and Proposition~\ref{co8prop2}(a).

\smallskip

\noindent{\bf Assumption \ref{co5ass2}(c).} If $s>0$ and $(E,V,\rho)\in\acA$ with $\lb E,V,\rho\rb=(\al,d)$ is $\ac\mu^{\om,s}$-semistable, we may divide into two cases:
\begin{itemize}
\setlength{\itemsep}{0pt}
\setlength{\parsep}{0pt}
\item[(i)] $\rank\al>0$.
\item[(ii)] $\rank\al=0$.
\end{itemize}
Then $\ac\mu^{\om,s}(\al,d)$ has degree $m$ in case (i), and degree $<m$ in case (ii).

Suppose $(E,V,\rho),(F,W,\si)\in\acA$ are $\ac\mu^{\om,s}$-semistable with $\ac\mu^{\om,s}(\lb E,V,\rho\rb)=\ac\mu^{\om,s}(\lb F,W,\si\rb)$ and $\lb(E,V,\rho)\op(F,W,\si)\rb\in C(\acB)_\pe$. Then \eq{co8eq40} implies that $\rank\lb E\op F\rb>0$, and $\dim V\op W\le 1$. Hence $(E,V,\rho)\op(F,W,\si)$ satisfies (i), so $\ac\mu^{\om,s}(\lb (E,V,\rho)\op(F,W,\si)\rb)=\ac\mu^{\om,s}(\lb E,V,\rho\rb)=\ac\mu^{\om,s}(\lb F,W,\si\rb)$ has degree $m$. Thus $(E,V,\rho)$ also satisfies (i), so $\rank\lb E\rb>0$, and $\dim V\le\dim V\op W\le 1$. Therefore $\lb E,V,\rho\rb\in C(\acB)_\pe$ by \eq{co8eq40}, and similarly~$\lb F,W,\si\rb\in C(\acB)_\pe$.
\smallskip

\noindent{\bf Assumption \ref{co5ass2}(d).} Let $(\ac\mu^{\om,s},\ac M,\le),(\ac\mu^{\om,\ti s},\ac M,\le)\in\sS$ and $(\al,d)\in C(\acB)_\pe$. Define a group morphism $\la:K(\acA)\ra\R$ by
\begin{equation*}
\la(\be,e)=(\ti s-s)\cdot\bigl(e\rank\al-d\rank\be\bigr).
\end{equation*}
Then $\la(\al,d)=0$, and we can show from \eq{co8eq9} and \eq{co8eq40} that if $(\be,e)\in C(\acB)_\pe$ with $\ac\mu^{\om,s}(\be,e)=\ac\mu^{\om,s}(\al,d)$ then $\la(\be,e)>0$ (or $\la(\be,e)<0$) if and only if $\ac\mu^{\om,\ti s}(\be,e)>\ac\mu^{\om,\ti s}(\al,d)$ (or $\ac\mu^{\om,\ti s}(\be,e)<\ac\mu^{\om,\ti s}(\al,d)$, respectively).
\smallskip

\noindent{\bf Assumption \ref{co5ass2}(f).} We adapt the treatment of Assumption \ref{co5ass2}(f) for $\coh(X)$ in \S\ref{co723}. As $X$ is projective it has an ample line bundle $M\ra X$. As in Definition \ref{co7def3}, for $\al\in C(\coh(X))$ the Hilbert polynomial $P^{c_1(M)}_\al(t)$ is of degree $\dim\al$, with positive leading coefficient $r_\al^{c_1(M)}$, and $P^{c_1(M)}_\al$ maps $\Z\ra\Z$ as $c_1(M)\in H^2(X,\Z)$. This implies that $(\dim\al)!\cdot r_\al^{c_1(M)}$ is an integer for $\al\in C(\coh(X))$. Extending \eq{co7eq53}, define $\rk:C(\acA)\ra\N_{>0}$ by
\e
\rk(\al,d)=\begin{cases} \rank\al, & \text{$\rank\al>0$,} \\
(\dim\al)!\cdot r_\al^{c_1(M)}\!+\!d, & \text{$\al\ne 0$ and $\dim\al=m-1$,} \\
d, & \text{$d\!>\!0$, and $\al\!=\!0$ or $\dim\al\!<\!m\!-\!1$,} \\
(\dim\al)!\cdot r_\al^{c_1(M)}, & \text{$\al\ne 0$, $\dim\al<m-1$ and $d=0$.} \\
\end{cases}
\label{co8eq44}
\e

Then $\deg\ac\mu^{\om,s}(\al,d)$ is $m$ in the first case of \eq{co8eq44}, and $m-1$ in the second and third cases, and is $\dim\al<m-1$ in the fourth case. Hence if $(\al,d),(\be,e)\in C(\acA)$ and $\ac\mu^{\om,s}(\al,d)=\ac\mu^{\om,s}(\be,e)$ then $\deg\ac\mu^{\om,s}(\al,d)=\deg\ac\mu^{\om,s}(\be,e)$, so $(\al,d),\ab(\be,e),\ab(\al,d)+(\be,e)$ are either all three in the first case of \eq{co8eq44}, or all three in the second and third cases, or all three in the fourth case with $\dim\al=\dim\be=\dim(\al+\be)$. From this we deduce that~$\rk((\al,d)+(\be,e))=\rk(\al,d)+\rk(\be,e)$.
\smallskip
 
\noindent{\bf Assumption \ref{co5ass2}(g).} This follows from \eq{co8eq40} and Proposition~\ref{co8prop2}(b).
\smallskip

\noindent{\bf Assumption \ref{co5ass2}(h).} Starting from $\acB\subseteq\acA,K(\acA),\acM,\acM^\pl,\ldots$ above, the weak stability condition $(\ac\mu^{\om,s},\ac M,\le)$ on $\acA$, and some extra data including a quiver $Q$, using the data $\bigl\{(\acB_k,\ac F_k,\ac\la_k):k\in K\bigr\}$ from Definition \ref{co8def6}, Definition \ref{co5def1} defines data $\bar{\acute{\mathcal B}}\subseteq\bar{\acute{\mathcal A}},\ldots$ and families of moduli stacks $\bar{\acute{\mathcal M}}_{(\ac\al,\bs d)}^\rst(\ov{{\ac\mu}^{\om,s}}{}^\la_{\bs\mu})\subseteq\bar{\acute{\mathcal M}}_{(\ac\al,\bs d)}^\ss(\ov{{\ac\mu}^{\om,s}}{}^\la_{\bs\mu})$. We must prove that if $(\ac\al,\bs d)\in C(\bar{\acute{\mathcal B}})_\pe$ and $\bar{\acute{\mathcal M}}_{(\ac\al,\bs d)}^\rst(\ov{{\ac\mu}^{\om,s}}{}^\la_{\bs\mu})=\bar{\acute{\mathcal M}}_{(\ac\al,\bs d)}^\ss(\ov{{\ac\mu}^{\om,s}}{}^\la_{\bs\mu})$ then $\bar{\acute{\mathcal M}}_{(\ac\al,\bs d)}^\ss(\ov{{\ac\mu}^{\om,s}}{}^\la_{\bs\mu})$ is a proper algebraic space.

Similarly, starting from $\B=\A=\coh(X\t\CP^1)$, the weak stability condition $(\mu^{\om\bp u\,\om_{\CP^1}},M,\le)$ on $\A$, and some extra data including a quiver $Q$, using the data $\bigl\{(\B_k,F_k,\la_k):k\in K\bigr\}$ for $\coh(X\t\CP^1)$ from \S\ref{co722}, Definition \ref{co5def1} defines data $\baB\subseteq\baA,\ldots$ and moduli stacks $\baM_{(\al,\bs d)}^\rst(\ov{\mu^{\om\bp u\,\om_{\CP^1}}}{}^\la_{\bs\mu})\ab\subseteq\ab\baM_{(\al,\bs d)}^\ss(\ov{\mu^{\om\bp u\,\om_{\CP^1}}}{}^\la_{\bs\mu})$. 

Section \ref{co74} proves that if $\baM_{(\al,\bs d)}^\rst(\ov{\mu^{\om\bp u\,\om_{\CP^1}}}{}^\la_{\bs\mu})\ab=\ab\baM_{(\al,\bs d)}^\ss(\ov{\mu^{\om\bp u\,\om_{\CP^1}}}{}^\la_{\bs\mu})$ then $\baM_{(\al,\bs d)}^\ss(\ov{\mu^{\om\bp u\,\om_{\CP^1}}}{}^\la_{\bs\mu})$ is a proper algebraic space.

Theorem \ref{co8thm1} gives functors $\ac I:\acA\ra\coh^{\SU(2)}(X\t\CP^1)$ or $\coh(X\t\CP^1)$. As $\bigl\{(\acB_k,\ac F_k,\ac\la_k):k\in K\bigr\}$ was defined by composing $\bigl\{(\B_k,F_k,\la_k):k\in K\bigr\}$ with $\ac I$, we can extend $\ac I$ to $\bar{\ac I}:\bar{\acute{\mathcal A}}\ra\baA$. This gives a full and faithful embedding of $\bar{\acute{\mathcal A}}$ in the category $\baA_{\SU(2)}$ of $\SU(2)$-equivariant objects in $\baA$, where the $\SU(2)$-action on $\baA$ is induced by the $\SU(2)$-action on $\coh(X\t\CP^1)$ from Theorem~\ref{co8thm1}. 

We claim that for all $(\ac\al,\bs d)\in C(\bar{\acute{\mathcal B}})_\pe$, for suitable $u>0$ depending on $s>0$ and $\ac\al\in C(\acA)\amalg\{0\}$, as for \eq{co8eq5}, the functor $\bar{\ac I}:\bar{\acute{\mathcal A}}\ra\baA$ induces morphisms
\e
\begin{gathered}
\xymatrix@C=180pt@R=15pt{ *+[r]{\bar{\acute{\mathcal M}}_{(\ac\al,\bs d)}^\ss(\ov{{\ac\mu}^{\om,s}}{}^\la_{\bs\mu})\,\,} \ar@{^{(}->}[r]^(0.43){\text{closed inclusion}} & *+[l]{\baM_{(\ac I_*(\ac\al),\bs d)}^\ss(\ov{\mu^{\om\bp u\,\om_{\CP^1}}}{}^\la_{\bs\mu})^{\SU(2)}} \ar[d] \\
& *+[l]{\baM_{(\ac I_*(\ac\al),\bs d)}^\ss(\ov{\mu^{\om\bp u\,\om_{\CP^1}}}{}^\la_{\bs\mu}).} }
\end{gathered}
\label{co8eq45}
\e
To see this, by \eq{co5eq16} and \eq{co8eq40} divide into cases:
\begin{itemize}
\setlength{\itemsep}{0pt}
\setlength{\parsep}{0pt}
\item[(i)] $\ac\al=(0,0)$.
\item[(ii)] $\ac\al=(\al,0)$ for $\rank\al>0$.
\item[(iii)] $\ac\al=(\al,1)$ for $\rank\al>0$ and $\bar\mu^\om(\al)<\bar\mu^\om(\lb L\rb)$.
\item[(iv)] $\ac\al=(\al,1)$ for $\rank\al>0$ and $\bar\mu^\om(\al)\ge\bar\mu^\om(\lb L\rb)$.
\end{itemize}
In cases (i) and (ii) the moduli spaces in \eq{co8eq45} are independent of $s,u$, and both morphisms in \eq{co8eq45} are isomorphisms, as all three spaces agree with the spaces $\baM_{(\al,\bs d)}^\ss(\bar\tau^\la_{\bs\mu})$ constructed in \S\ref{co52} for $\A=\coh(X)$ and $(\tau,T,\le)=(\mu^\om,M,\le)$.

In case (iii) $\acM_{(\al,1)}^\ss(\ac\mu^{\om,s})=\es$ by Theorem \ref{co8thm2}(c), so $\bar{\acute{\mathcal M}}_{(\ac\al,\bs d)}^\ss(\ov{{\ac\mu}^{\om,s}}{}^\la_{\bs\mu})=\es$ because of $\Pi_{\M_\al^\ss(\tau)}:\baM_{(\al,\bs d)}^\ss(\bar\tau^\la_{\bs\mu})\ra\M_\al^\ss(\tau)$ in \S\ref{co52}, and Assumption \ref{co5ass2}(h) is trivial in this case. In case (iv) Proposition \ref{co8prop1} shows $\ac I:\acA\ra\coh(X\t\CP^1)$ maps $\acM_{(\al,1)}^\ss(\ac\mu^{\om,s})$ to $\M_{\ac I_*(\al,1)}^\ss(\mu^{\om\bp u\,\om_{\CP^1}})^{\SU(2)}$ for $u$ as in \eq{co8eq11}, and it easily follows that $\bar{\ac I}:\bar{\acute{\mathcal A}}\ra\baA$ maps $\bar{\acute{\mathcal M}}_{(\ac\al,\bs d)}^\ss(\ov{{\ac\mu}^{\om,s}}{}^\la_{\bs\mu})$ to $\baM_{(\ac I_*(\ac\al),\bs d)}^\ss(\ov{\mu^{\om\bp u\,\om_{\CP^1}}}{}^\la_{\bs\mu})^{\SU(2)}$. This proves the claim.

The proofs in \S\ref{co74} are compatible with restricting to $\SU(2)$-fixed substacks. Thus if $\baM_{(\ac I_*(\ac\al),\bs d)}^\rst(\ov{\mu^{\om\bp u\,\om_{\CP^1}}}{}^\la_{\bs\mu})^{\SU(2)}=\baM_{(\ac I_*(\ac\al),\bs d)}^\ss(\ov{\mu^{\om\bp u\,\om_{\CP^1}}}{}^\la_{\bs\mu})^{\SU(2)}$ then this is a proper algebraic space, even if $\baM_{(\ac I_*(\ac\al),\bs d)}^\rst(\ov{\mu^{\om\bp u\,\om_{\CP^1}}}{}^\la_{\bs\mu})\ne\baM_{(\ac I_*(\ac\al),\bs d)}^\ss\ab(\ov{\mu^{\om\bp u\,\om_{\CP^1}}}{}^\la_{\bs\mu})$. The proofs are also compatible with restricting to closed substacks. Hence as $\bar{\acute{\mathcal M}}_{(\ac\al,\bs d)}^\rst(\ov{{\ac\mu}^{\om,s}}{}^\la_{\bs\mu}),\ab\bar{\acute{\mathcal M}}_{(\ac\al,\bs d)}^\ss(\ov{{\ac\mu}^{\om,s}}{}^\la_{\bs\mu})$ are intersections of a closed substack with $\baM_{(\ac I_*(\ac\al),\bs d)}^\rst(\ov{\mu^{\om\bp u\,\om_{\CP^1}}}{}^\la_{\bs\mu})^{\SU(2)}$, $\baM_{(\ac I_*(\ac\al),\bs d)}^\ss(\ov{\mu^{\om\bp u\,\om_{\CP^1}}}{}^\la_{\bs\mu})^{\SU(2)}$, if $\bar{\acute{\mathcal M}}_{(\ac\al,\bs d)}^\rst(\ov{{\ac\mu}^{\om,s}}{}^\la_{\bs\mu})\ab=\ab\bar{\acute{\mathcal M}}_{(\ac\al,\bs d)}^\ss(\ov{{\ac\mu}^{\om,s}}{}^\la_{\bs\mu})$ then this is a proper algebraic space.
\smallskip

\noindent{\bf Assumption \ref{co5ass3}(a).} The weak stability condition $(\ac\mu^{\om,s},\ac M,\le)$ in Definition \ref{co8def3} depends continuously on $s\in(0,\iy)$ in the sense of Definition \ref{co3def5}. Thus given any $(\ac\mu^{\om,s},\ac M,\le),(\ac\mu^{\om,\ti s},\ac M,\le)$ in $\acS$ in \eq{co8eq43}, $(\ac\mu^{\om,(1-t)s+t\ti s},\ac M,\le)_{t\in[0,1]}$ is a continuous family in $\acS$ connecting $(\ac\mu^{\om,s},\ac M,\le)$ and~$(\ac\mu^{\om,\ti s},\ac M,\le)$.

\subsubsection{Assumption \ref{co5ass3}(b) for $\acute{\mathcal A}$}
\label{co824}

Here is the analogue of \S\ref{co75} for $\acA$. Let $(\al,d)\in C(\acB)_\pe$ and $(\ac\mu^{\om,s_t},\ac M,\le)_{t\in[0,1]}$ be a continuous family in $\acS$. Then $t\mapsto s_t$ is a continuous map $[0,1]\ra(0,\iy)$, so there exist $0<s_0\le s_1$ such that $\{s_t:t\in[0,1]\}=[s_0,s_1]$. Thus, to verify Assumption \ref{co5ass3}(b) for $\acA$ with $(\al,d),(\ac\mu^{\om,s_t},\ac M,\le)_{t\in[0,1]}$ in place of $\al,(\tau_t,T_t,\le)_{t\in[0,1]}$, it is sufficient to replace $\tau_s,\tau_t,\tau_u$ in Assumption \ref{co5ass3}(b) by $\ac\mu^{\om,s},\ac\mu^{\om,t},\ac\mu^{\om,u}$ with $s,t,u\in[s_0,s_1]$. We will take these $(\al,d),s_0,s_1$ to be fixed throughout the proof, though $n,p,a_j,(\al_i,d_i),(\be_j,e_j),s,t,u$ below will vary.

Assumption \ref{co5ass3}(b) concerns choices of data $n\ge p\ge 1$, $(\al_1,d_1),\ldots,(\al_n,d_n)$ in $C(\acB)$, $0=a_0<\cdots<a_p=n$ and $s,t,u\in[s_0,s_1]$ such that $(\al_1,d_1)+\cdots+(\al_n,d_n)=(\al,d)$ and, writing $(\be_j,e_j)=(\al_{a_{j-1}+1},d_{a_{j-1}+1})+\cdots+(\al_{a_j},d_{a_j})$, then $\acM_{(\al_i,d_i)}^\ss(\ac\mu^{\om,s})\ne\es$ for $i=1,\ldots,n$, and $U((\al_{a_{j-1}+1},d_{a_{j-1}+1}),\ab\ldots,\ab(\al_{a_j},d_{a_j});\ab\ac\mu^{\om,s},\ab\ac\mu^{\om,t})\ne 0$ for $j=1,\ldots,p$, and $\ac\mu^{\om,u}(\be_1,e_1)=\cdots=\ac\mu^{\om,u}(\be_p,e_p)=\ac\mu^{\om,u}(\al,d)$. As $d\le 1$ we have $e_j,d_i\le 1$. We must show there are only finitely many choices of such data, and that $(\al_{i_1},d_{i_1})+\cdots+(\al_{i_2},d_{i_2})\in C(\acB)_\pe$ for~$1\le i_1\le i_2\le n$. 

First observe that $\ac\mu^{\om,u}(\be_j,e_j)=\ac\mu^{\om,u}(\al,d)$ implies that $\deg\ac\mu^{\om,u}(\be_j,e_j)=m$, so that $\rank\be_j>0$. Then $U((\al_{a_{j-1}+1},d_{a_{j-1}+1}),\ab\ldots,\ab(\al_{a_j},d_{a_j});\ab\ac\mu^{\om,s},\ab\ac\mu^{\om,t})\ne 0$ implies that $\rank\al_i>0$ for $i=a_{j-1}+1,\ldots,a_j$, and hence for all $i=1,\ldots,n$. It now follows from \eq{co8eq40} and $\sum_{i=1}^nd_i=d\le 1$ that $(\al_{i_1},d_{i_1})+\cdots+(\al_{i_2},d_{i_2})\in C(\acB)_\pe$ for $1\le i_1\le i_2\le n$, as we have to prove.

As $\rank\al_i,\rank\be_j\!>\!0$ for all $i,j$, conditions such as $\ac\mu^{\om,u}(\be_j,e_j)\!=\!\ac\mu^{\om,u}(\al,d)$ in $\ac M$ are equivalent to $\bar\mu^{\om,u}(\be_j,e_j)=\bar\mu^{\om,u}(\al,d)$ in $\R$, where $\bar\mu^{\om,u}(\cdots)$ may be computed using \eq{co8eq9}, and we will rewrite (in)equalities in $\ac\mu^{\om,s},\ab\ac\mu^{\om,t},\ab\ac\mu^{\om,u}(\cdots)$ in terms of $\bar\mu^{\om,s},\bar\mu^{\om,t},\bar\mu^{\om,u}(\cdots)$ below. It follows from \eq{co8eq9} that if $(\be,e)\in C(\acA)$ with $e\le 1$ and $v,v'\in [s_0,s_1]$ then
\e
\bar\mu^{\om,v}(\be,e)\le \bar\mu^{\om,v'}(\be,e)+(s_1-s_0).
\label{co8eq46}
\e

Let $i=1,\ldots,n$. Then for some $j=1,\ldots,p$ we have $a_{j-1}<i\le a_j$. As $U((\al_{a_{j-1}+1},0),\ldots,(\al_{a_j},0);\ac\mu^{\om,s},\ac\mu^{\om,t})\ne 0$, Proposition \ref{co3prop2} gives $q,r=a_{j-1}+1,\ldots,a_j$ such that
\e
\begin{split}
&\bar\mu^{\om,s}(\al_q,d_q)\le\bar\mu^{\om,s}(\al_i,d_i)\le\bar\mu^{\om,s}(\al_r,d_r)\quad\text{and}\\
&\bar\mu^{\om,t}(\al_q,d_q)\ge\bar\mu^{\om,t}(\be_j,e_j)\ge\bar\mu^{\om,t}(\al_r,d_r).
\end{split}
\label{co8eq47}
\e
Equations \eq{co8eq9} and \eq{co8eq46}--\eq{co8eq47} and $\bar\mu^{\om,u}(\be_j,e_j)=\bar\mu^{\om,u}(\al,d)$ imply that
\ea
\bar\mu^{\om,s}(\al_i,d_i)&\le\bar\mu^{\om,s}(\al_r,d_r)
\le \bar\mu^{\om,t}(\al_r,d_r)+(s_1-s_0)
\nonumber\\
& \le\bar\mu^{\om,t}(\be_j,e_j)+(s_1-s_0)
\le\bar\mu^{\om,u}(\be_j,e_j)+2(s_1-s_0)
\nonumber\\
&=\bar\mu^{\om,u}(\al,d)+2(s_1-s_0)
\le\bar\mu^{\om,s_0}(\al,d)+3(s_1-s_0)=:\mu_0.
\label{co8eq48}
\ea

Pick $[E_i,V_i,\rho_i]\in \acM_{(\al_i,d_i)}^\ss(\ac\mu^{\om,s})\ne\es$ for $i=1,\ldots,n$. Then $E_i$ is torsion-free by Lemma \ref{co8lem2} and $\dim\al_i=m$. Also we have $\mu_{\max}^\om(E_i)\le \bar\mu^{\om,s}(\al_i,d_i)\le\mu_0$ by \eq{co8eq48}, as if $0\ne E_i'\subseteq E_i$ then $(E_i',0,0)\subseteq(E_i,V_i,\rho)$, so $\bar\mu^\om(\lb E_i\rb)=\bar\mu^{\om,s}(\lb E_i,0,0\rb)\le\ab\bar\mu^{\om,s}(\al_i,d_i)$ as $(E_i,V_i,\rho_i)$ is $\bar\mu^{\om,s}$-semistable. Hence $\bigop_{i=1}^nE_i$ is torsion-free in class $\al$ in $\coh(X)$, with $\mu_{\max}^\om(\bigop_{i=1}^nE_i)\le\mu_0$.

Now Proposition \ref{co7prop6} shows that the family of torsion-free sheaves $E$ in $\coh(X)$ with $\lb E\rb=\al$ and $\mu_{\max}^\om(E)\le\mu_0$ is bounded. Since $E_1\op\cdots\op E_n$ lies in a bounded family which is independent of $n,p,a_j,(\al_i,d_i),(\be_j,e_j),s,t,u$, there are only finitely many possibilities for $\al_1,\ldots,\al_n=(\lb E_1\rb,\ldots,\lb E_n\rb)$, and hence only finitely many choices for $(\al_1,d_1),\ldots,(\al_n,d_n)$. This proves Assumption~\ref{co5ass3}(b).

\subsubsection{The $G$-equivariant case and Assumptions \ref{co4ass2}, \ref{co5ass4}(a),(b),(d)}
\label{co825}

Here is the analogue of \S\ref{co724} for $\acA$.

\begin{dfn}
\label{co8def7}
Work in the situation of \S\ref{co821}--\S\ref{co824}, defining data $\acA=\acB,K(\acA),\acM,\ldots$ satisfying most of Assumptions \ref{co4ass1} and \ref{co5ass1}--\ref{co5ass3}. Suppose also that $G$ is an algebraic $\C$-group with an action $\phi:G\t X\ra X$ on the $\C$-scheme $X$, such that:
\begin{itemize}
\setlength{\itemsep}{0pt}
\setlength{\parsep}{0pt}
\item[(i)] There exists at least one ample line bundle $\O_X(1)\ra X$ such that the $G$-action on $X$ lifts to $\O_X(1)$.
\item[(ii)] The action of $G$ on $K(\coh(X))=K^\num(\coh(X))$ is trivial. This is automatic if $G$ is connected, or if the action of $G$ on $H^{\rm even}(X,\Q)$ is trivial.
\item[(iii)] We are given a chosen lift $\Up$ of the action of $G$ on $X$ to $L\ra X$ in Definition \ref{co8def1}. This includes an isomorphism $\Up(g):L\ra \phi(g)^*(L)$ for each $g\in G(\C)$ satisfying $\Up(gh)=\phi(h)^*(\Up(g))\ci\Up(h)$ for all $g,h\in G(\C)$, and $\Up(1)=\id_L$.
\end{itemize}
Here (i),(ii) are as in Definition \ref{co7def12}. We can now verify Assumptions \ref{co4ass2} and \ref{co5ass4}(a),(b),(d) in a very similar way to the proofs for $\A=\coh(X)$ in \S\ref{co724}. We explain Assumption \ref{co4ass2}(a) for $\acA$, and leave the rest as an exercise.

For Assumption \ref{co4ass2}(a), if $g\in G(\C)$ then $g^{-1}\in G(\C)$ so $\phi(g^{-1}):X\ra X$ is an isomorphism. For (a)(i), define a functor $\ac\Ga(g):\acA\ra\acA$ by
\begin{align*}
\ac\Ga(g)&:(E,V,\rho)\longmapsto \bigl(\phi(g^{-1})^*(E),V,\phi(g^{-1})^*(\rho)\ci(\id_V\ot\Up(g^{-1}))\bigr),\\
\ac\Ga(g)&:(\th,\phi)\longmapsto \bigl(\phi(g^{-1})^*(\th),\phi\bigr).
\end{align*}
That is, $\ac\Ga(g)$ acts on $(E,V,\rho)\in\acA$ as $\Ga(g)$ in \S\ref{co724} on $E$, as the identity on $V$, and on $\rho$ as the composition $V\ot L\,{\buildrel\Up(g^{-1})\over\longra}\, V\ot\phi(g^{-1})^*(L)\,{\buildrel\phi(g^{-1})^*(\rho)\over\longra}\,\phi(g^{-1})^*(E)$.

The isomorphisms $\ac\Ga_{g,h}(E,V,\rho)$, $\ac\Ga_1(E,V,\rho)$ in Assumption \ref{co4ass2}(a)(ii),(iii) act as for $\Ga_{g,h},\Ga_1$ for $\coh(X)$ in \S\ref{co724} on $E$, and as $\id_V$ on $V$, since $\ac\Ga(g)$ does not change $V$. The identities $\Up(gh)=\phi(h)^*(\Up(g))\ci\Up(h)$, $\Up(1)=\id_L$ in (iii) are needed to show that these prescriptions do give morphisms $\ac\Ga_{g,h}(E,V,\rho)$, $\ac\Ga_1(E,V,\rho)$ in $\acA$, that is, that \eq{co8eq1} commutes for each. Then \eq{co4eq19}--\eq{co4eq20} for $\ac\Ga_{g,h},\ac\Ga_1$ follow from \eq{co4eq19}--\eq{co4eq20} for $\Ga_{g,h},\Ga_1$ in~\S\ref{co724}.
\end{dfn}

\subsection{\texorpdfstring{Quasi-smooth derived stacks, obstruction theories, \\ and Assumptions \ref{co5ass1}(f) and \ref{co5ass2}(e)}{Quasi-smooth derived stacks, obstruction theories, and Assumptions \ref{co5ass1}(f) and \ref{co5ass2}(e)}}
\label{co83}

We now explain how to construct the quasi-smooth derived stacks $\bs\dM_{(\al,d)}^\red$ and $\bs\dM_{(\al,d)}^\rpl$ in Assumption \ref{co5ass1}(f) for $\acA$, which yield obstruction theories \eq{co5eq5} on $\dM_{(\al,d)}\subseteq\acM_{(\al,d)}$ and~$\dM_{(\al,d)}^\pl\subseteq\acM_{(\al,d)}^\pl$.

In Remark \ref{co7rem7}, for $\A=\coh(X)$, we saw that $i^*(\bL_{\bs\M_\al}),i^*(\bL_{\bs\M_\al^\pl})$ are perfect in $[1-m,1]$. Also, in \S\ref{co734}, if $X$ is a Fano $m$-fold we found suitable open $\dM_\al\subseteq\M_\al$ and $\dM_\al^\pl\subseteq\M_\al^\pl$ on which $i^*(\bL_{\bs\M_\al}),i^*(\bL_{\bs\M_\al^\pl})$ are perfect in $[2-m,1]$. Thus to get obstruction theories perfect in $[-1,1]$, the theory worked for $X$ a curve, surface or Fano 3-fold.

For $\acA$, $i^*(\bL_{\bs\acM_{(\al,d)}}),i^*(\bL_{\bs\acM_{(\al,d)}^\pl})$ are perfect in $[-m,1]$ by Remark \ref{co8rem2}. We will show they are perfect in $[1-m,1]$ on suitable open $\dM_{(\al,d)}\subseteq\acM_{(\al,d)}$ and $\dM_{(\al,d)}^\pl\subseteq\acM_{(\al,d)}^\pl$ if $\rank\al>0$ and $d\le 1$. So we restrict to $X$ a curve or a surface only, and include the conditions on $(\al,d)$ in $C(\acB)_\pe$ in \eq{co8eq40} above. 

\subsubsection{Quasi-smooth moduli stacks for $X$ a curve}
\label{co831}

As in \S\ref{co731}, we verify Assumptions \ref{co5ass1}(f) and \ref{co5ass2}(e) when $X$ is a curve.

\begin{dfn}
\label{co8def9}
In Definition \ref{co8def1}, suppose $m=\dim X=1$, so that $X$ is a connected, smooth, projective complex curve. In Assumption \ref{co5ass1}(f), for each $(\al,d)\in C(\acB)_\pe$ in \eq{co8eq40}, in (i)--(ii) we define $\dM_{(\al,d)}=\acM_{(\al,d)}$, $\dM_{(\al,d)}^\pl=\acM_{(\al,d)}^\pl$, $\bs\dM_{(\al,d)}^\red=\bs\acM_{(\al,d)}$, $\bs\dM_{(\al,d)}^\rpl=\bs\acM_{(\al,d)}^\pl$, $\bs j_{(\al,d)}=\bs\id=\bs j_{(\al,d)}^\pl$, $\bs{\dot\Pi}{}_{(\al,d)}^\rpl=\bs{\ac\Pi}{}_{(\al,d)}^\pl$, and $U_{(\al,d)}=o_{(\al,d)}=0$. Then $\bs\dM_{(\al,d)}^\red,\bs\dM_{(\al,d)}^\rpl$ are locally finitely presented as in Definition \ref{co8def6} and quasi-smooth as in Remark \ref{co8rem2}, as $m=1$. Also \eq{co5eq4} is trivially Cartesian, proving (ii), equation \eq{co5eq6} holds as $\bs j_{(\al,d)},\bs j_{(\al,d)}^\pl$ are identities and $U_{(\al,d)}=0$, giving (iii), and (iv) is trivial as $o_{(\al,d)}=o_{(\be,e)}=o_{(\al+\be,d+e)}=0$. This proves Assumption~\ref{co5ass1}(f).

Assumption \ref{co5ass2}(e) is trivial as~$\dM_{(\al,d)}=\acM_{(\al,d)}$.
\end{dfn}

\subsubsection{Quasi-smooth moduli stacks for $X$ a surface: the case $p_g=0$}
\label{co832}

Here is the analogue of \S\ref{co732} for $\acA$.

\begin{dfn}
\label{co8def10}
In Definition \ref{co8def1}, let $m=2$, so that $X$ is a connected projective surface, and suppose $X$ has  geometric genus $p_g=0$ in Definition~\ref{co7def14}.

For $(\al,d)\in C(\acB)_\pe$ in \eq{co8eq40}, in an analogue of \eq{co7eq62} in the Fano 3-fold case, in Assumption \ref{co5ass1}(f)(i) define an open substack $\dM_{(\al,d)}\subseteq\acM_{(\al,d)}$ by
\e
\dM_{(\al,d)}=\acM_{(\al,d)}\sm\supp h^3\bigl(\De_{\acM_{(\al,d)}}^*(\acExt_{(\al,d),(\al,d)}^\bu)\bigr).
\label{co8eq49}
\e
Here $h^3\bigl(\De_{\acM_{(\al,d)}}^*(\acExt_{(\al,d),(\al,d)}^\bu)\bigr)$ is a sheaf on $\acM_{(\al,d)}$ with fibre $\acExt^3((E,V,\rho),\ab(E,V,\rho))$ at $[E,V,\rho]\in\acM_{(\al,d)}$, and $\supp h^3(\cdots)$ is its Zariski closed support. Thus $\dM_{(\al,d)}$ is the open substack of $[E,V,\rho]$ in $\acM_{(\al,d)}$ with $\racExt^3((E,V,\rho),\ab(E,V,\rho))=0$. It is invariant under the $[*/\bG_m]$-action on $\acM_{(\al,d)}$, so there is unique open $\dM_{(\al,d)}^\pl\subseteq\acM_{(\al,d)}^\pl$ with~$\dM_{(\al,d)}=(\ac\Pi_{(\al,d)}^\pl)^{-1}(\dM_{(\al,d)}^\pl)$.

Let $\bs\dM_{(\al,d)}\subseteq\bs\acM_{(\al,d)}$ and $\bs\dM_{(\al,d)}^\pl\subseteq\bs\acM_{(\al,d)}^\pl$ be the derived open substacks corresponding to $\dM_{(\al,d)}\subseteq\acM_{(\al,d)}$ and $\dM_{(\al,d)}^\pl\subseteq\acM_{(\al,d)}^\pl$. Define $\bs\dM_{(\al,d)}^\red=\bs\dM_{(\al,d)}$, $\bs\dM_{(\al,d)}^\rpl=\bs\dM_{(\al,d)}^\pl$, $\bs j_{(\al,d)}=\bs\id=\bs j_{(\al,d)}^\pl$, $\bs{\dot\Pi}{}_{(\al,d)}^\rpl=\bs{\dot\Pi}{}_{(\al,d)}^\pl$, and $U_{(\al,d)}=o_{(\al,d)}=0$. Then $\bs\dM_{(\al,d)}^\red,\bs\dM_{(\al,d)}^\rpl$ are locally finitely presented as in Definition \ref{co8def6}. Also \eq{co5eq4} is trivially Cartesian, equation \eq{co5eq6} holds as $\bs j_{(\al,d)},\bs j_{(\al,d)}^\pl$ are identities and $U_{(\al,d)}=0$, giving (iii), and (iv) is trivial as $o_{(\al,d)}=o_{(\be,e)}=o_{(\al+\be,d+e)}=0$. This proves all of Assumption \ref{co5ass1}(f) except that $\bs\dM_{(\al,d)}^\red,\bs\dM_{(\al,d)}^\rpl$ are quasi-smooth. To show this, by Remark \ref{co2rem6}(e) it suffices to check that $i^*(\bL_{\bs\dM_{(\al,d)}^\red}),i^*(\bL_{\bs\dM_{(\al,d)}^\rpl})$ are perfect in the interval~$[-1,1]$.

As $\bs\dM_{(\al,d)}^\red=\bs\dM_{(\al,d)}\subseteq\bs\acM_{(\al,d)}$ and $m=2$, equation \eq{co8eq42} shows that $i^*(\bL_{\bs\dM_{(\al,d)}^\red})$ is perfect in $[-2,1]$. Thus the argument at the end of Definition \ref{co7def16} shows that $i^*(\bL_{\bs\dM_{(\al,d)}^\red})$ is perfect in $[-1,1]$ if $h^2(i^*(\bL_{\bs\dM_{(\al,d)}^\red})^\vee)=0$. But equations \eq{co8eq32}, \eq{co8eq39} and \eq{co8eq49} imply that
\begin{equation*}
h^2\bigl(i^*(\bL_{\bs\dM_{(\al,d)}^\red})^\vee\bigr)\cong h^3\bigl(\De_{\acM_{(\al,d)}}^*(\acExt_{(\al,d),(\al,d)}^\bu)\bigr)\vert_{\dM_{(\al,d)}}=0.
\end{equation*}
Hence $i^*(\bL_{\bs\dM_{(\al,d)}^\red})$ is perfect in $[-1,1]$. By \eq{co8eq42} restricted to $\dM_{(\al,d)},\dM_{(\al,d)}^\pl$, $i^*(\bL_{\bs\dM_{(\al,d)}^\rpl})$ is also perfect in $[-1,1]$. Therefore $\bs\dM_{(\al,d)}^\red,\bs\dM_{(\al,d)}^\rpl$ are quasi-smooth, completing Assumption~\ref{co5ass1}(f).

For Assumption \ref{co5ass2}(e), let $(\al,d)\in C(\acB)_\pe$, and suppose $[E,V,\rho]$ lies in $\acM_{(\al,d)}^\ss(\ac\mu^{\om,s})$ for some $s>0$. We will show that $\racExt^3((E,V,\rho),(E,V,\rho))=0$, so that $[E,V,\rho]\in\dM_{(\al,d)}^\pl$, and thus $\acM_\al^\ss(\ac\mu^{\om,s})\subseteq\dM_{(\al,d)}^\pl$, as we have to prove. Equation \eq{co8eq30} at $[E_1,V_1,\rho_1]=[E_2,V_2,\rho_2]=[E,V,\rho]$ gives an exact sequence
\e
\smash{\xymatrix@C=17pt{
\Ext^2(E,E) \ar[r]^(0.4){-\ci\rho} &
V^*\ot \Ext^2(L,E) \ar[r] & \racExt^3((E,V,\rho),(E,V,\rho)) \ar[r] & 0. }}
\label{co8eq50}
\e
We have to show the first morphism in \eq{co8eq50} is surjective. By Serre duality, this is equivalent to the following morphism being injective:
\e
\smash{\xymatrix@C=60pt{
\Hom(E,V\ot_\C L\ot K_X) \ar[r]^(0.55){(\rho\ot\id_{K_x})\ci -} & \Hom(E,E\ot K_X). }}
\label{co8eq51}
\e

By \eq{co8eq40} we have $\rank\al>0$ and $d=0$ or 1. If $d=0$ then $V=0$ and \eq{co8eq51} is trivially injective. If $d=1$ then $E$ is torsion-free by Lemma \ref{co8lem2} and $\rho\ne 0$ by Theorem \ref{co8thm2}(a). Thus there is a dense open subscheme $X'\subseteq X$ such that $E\vert_{X'}$ is a vector bundle and $\rho:V\ot_\C L\ra E$ is a vector bundle morphism which is nonzero at each $x\in X'$, and so is injective at $x'$ as $\rank V\ot_\C L=1$. Hence the following is injective at each $x\in X'$
\e
\smash{\xymatrix@C=60pt{
\Hom(E\vert_x,V\ot_\C L\vert_x\ot K_X\vert_x) \ar[r]^(0.55){(\rho\vert_x\ot\id_{K_X}\vert_x)\ci -} & \Hom(E\vert_x,E\vert_x\ot K_X\vert_x), }}
\label{co8eq52}
\e
which implies that \eq{co8eq51} is injective. 
\end{dfn}

\begin{rem}
\label{co8rem3}
The proof of Assumption \ref{co5ass2}(e) above is based on Mochizuki \cite[Prop.~6.1.1]{Moch}. It is the reason we assume that $\rank\al>0$ and $d\le 1$ in the definition of $C(\acB)_\pe$ in \eq{co8eq40}. If we allowed $(\al,d)$ in $C(\acB)_\pe$ with either $d>1$, or $d>0$ and $\rank\al=0$, $\al\ne 0$ so $E$ is not torsion-free, then we could not prove \eq{co8eq51} is injective, so there might exist $[E,V,\rho]\in\acM_{(\al,d)}^\ss(\ac\mu^{\om,s})$ with $\racExt^3((E,V,\rho),(E,V,\rho))\ne 0$, and Assumption \ref{co5ass2}(e) would fail.
\end{rem}

\subsubsection{Quasi-smooth moduli stacks for $X$ a surface: the case $p_g>0$}
\label{co833}

In Remark \ref{co7rem9} in \S\ref{co733}, we explained that while the quasi-smooth stacks $\bs\dM_\al^\red=\bs\M_\al$, $\bs\dM_\al^\rpl=\bs\M_\al^\pl$ and associated obstruction theories \eq{co5eq5} defined in \S\ref{co732} work for any connected projective surface $X$, when $p_g>0$ they would yield invariants $[\M_\al^\ss(\tau)]_\inv=0$ in Theorems \ref{co5thm1}--\ref{co5thm3} when $\rank\al>0$, which would be boring. So in \S\ref{co733} we defined alternative `reduced' versions when $p_g>0$ and $\rank\al>0$, which generally give nonzero~$[\M_\al^\ss(\tau)]_\inv$.

For $\acA$ the situation is similar: the quasi-smooth stacks $\bs\dM_{(\al,d)}^\red=\bs\dM_{(\al,d)}$, $\bs\dM_{(\al,d)}^\rpl=\bs\dM_{(\al,d)}^\pl$ in \S\ref{co832} satisfy Assumptions \ref{co5ass1}(f) and \ref{co5ass2}(e) for any $X$, but when $p_g>0$ they would yield invariants $[\acM_{(\al,d)}^\ss(\ac\mu^{\om,s})]_\inv=0$ if either $d=0$ and $\rank\al>0$, or $d=1$ and $\rank\al>1$. So in these cases we define alternative `reduced' versions, by the same method as Definition~\ref{co7def16}.

\begin{dfn}
\label{co8def11}
In Definition \ref{co8def1}, let $m=2$, so that $X$ is a connected projective surface, and suppose $X$ has  geometric genus $p_g>0$ in Definition~\ref{co7def14}.

Let $(\al,d)\in C(\acB)$. As for \eq{co7eq56}, define a derived stack $\bs\acM_{(\al,d)}^\red$ by the Cartesian squares of derived stacks
\e
\begin{gathered}
\xymatrix@C=100pt@R=15pt{
*+[r]{\bs\acM_{(\al,d)}^\red} \ar@/^.5pc/@<.5ex>[rr]^(0.15){\ac\Pi_{\Pic(X)}} \ar[r]_{\Pi_{\bs\M_\al^\red}} \ar[d]^{\bs{\ac\jmath}_{(\al,d)}} & {\bs\M_\al^\red} \ar[r]_{\Pi_{\Pic(X)}} \ar[d]^{\Pi_{\bs\M_\al}} & *+[l]{\Pic(X)_{\vphantom{(}}} \ar@{_{(}->}[d]_i \\
*+[r]{\bs\acM_{(\al,d)}} \ar@/_.5pc/@<-.5ex>[rr]_(0.15){\rm\ac det} \ar[r]^{\Pi_{\bs\M_\al}} & {\bs\M_\al} \ar[r]^{\det} & *+[l]{\bs\Pic(X).\!}
}	
\end{gathered}
\label{co8eq53}
\e
Here the right hand square is \eq{co7eq56} for $\coh(X)$, both squares and the large rectangle are Cartesian, and ${\rm\ac det},\ac\Pi_{\Pic(X)}$ are defined by composition.

Similarly, as for \eq{co7eq57} define $\bs\acM_{(\al,d)}^\rpl$ by the Cartesian squares
\e
\begin{gathered}
\xymatrix@C=100pt@R=15pt{
*+[r]{\bs\acM_{(\al,d)}^\rpl} \ar@/^.5pc/@<.5ex>[rr]^(0.15){\ac\Pi_{\Pic(X)^\pl}} \ar[r]_{\Pi_{\bs\M_\al^\rpl}} \ar[d]^{\bs{\ac\jmath}{}^{\,\pl}_{(\al,d)}} & {\bs\M_\al^\rpl} \ar[r]_{\Pi_{\Pic(X)^\pl}} \ar[d]^{\Pi_{\bs\M_\al^\pl}} & *+[l]{\Pic(X)^\pl_{\vphantom{(}}} \ar@{_{(}->}[d]_{i^\pl} \\
*+[r]{\bs\acM_{(\al,d)}^\pl} \ar@/_.5pc/@<-.5ex>[rr]_(0.15){\rm\ac det^\pl} \ar[r]^{\Pi_{\bs\M_\al^\pl}} & {\bs\M_\al^\pl} \ar[r]^{\det^\pl} & *+[l]{\bs\Pic(X)^\pl.\!} }	
\end{gathered}
\label{co8eq54}
\e
Since $\Pic(X)=t_0(\bs\Pic(X))$, taking classical truncations in \eq{co8eq53} (which preserves Cartesian squares) shows that $t_0(\bs\acM_{(\al,d)}^\red)=t_0(\bs\acM_{(\al,d)})=\acM_{(\al,d)}$ and $t_0(\bs{\ac\jmath}_{(\al,d)})\!=\!\id$. Similarly $t_0(\bs\acM_{(\al,d)}^\rpl)\!=\!t_0(\bs\acM_{(\al,d)}^\pl)\!=\!\acM_{(\al,d)}^\pl$ and~$t_0(\bs{\ac\jmath}_{(\al,d)}^{\,\pl})\!=\!\id$. 

Definition \ref{co8def6} showed that $\bs\acM_{(\al,d)},\bs\acM_{(\al,d)}^\pl$ are locally finitely presented. Also $\Pic(X),\bs\Pic(X),\Pic(X)^\pl,\bs\Pic(X)^\pl$ are, noting that $\Pic(X),\Pic(X)^\pl$ are smooth. So the Cartesian outer rectangles of \eq{co8eq53}--\eq{co8eq54} imply that $\bs\acM_{(\al,d)}^\red,\ab\bs\acM_{(\al,d)}^\rpl$ are locally finitely presented.

As for \eq{co7eq58} we get a (homotopy) commuting diagram:
\e
\begin{gathered}
\xymatrix@!0@C=90pt@R=23pt{
*+[r]{\bs\acM_{(\al,d)}^\red} \ar@{..>}[dr]_{\bs{\ac\Pi}{}_{(\al,d)}^\rpl} \ar[rr]_(0.75){\ac\Pi_{\Pic(X)}} \ar[dd]^(0.65){\bs{\ac\jmath}_{(\al,d)}} && *+[l]{\Pic(X)_{\vphantom{(}}} \ar[dr]^(0.4){\Pi^\pl} \ar@{_{(}->}[dd]^(0.25)i 
\\
& *+[r]{\bs\acM_{(\al,d)}^\rpl} \ar[dd]^(0.33){\bs{\ac\jmath}{}^{\,\pl}_{(\al,d)}}  \ar[rr]_(0.3){\ac\Pi_{\Pic(X)^\pl}} && *+[l]{\Pic(X)^\pl_{\vphantom{(}}} \ar@{_{(}->}[dd]_{i^\pl}
\\
*+[r]{\bs\acM_{(\al,d)}} \ar[dr]_{\bs{\ac\Pi}{}_{(\al,d)}^\pl} \ar[rr]^(0.35){\rm\ac det} && *+[l]{\bs\Pic(X)} \ar[dr]^(0.4){\bs\Pi^\pl}\\
& *+[r]{\bs\acM_{(\al,d)}^\pl} \ar[rr]^{\rm\ac det^\pl} && *+[l]{\bs\Pic(X)^\pl,\!} }	
\end{gathered}
\label{co8eq55}
\e
with all six faces of the cube Cartesian. Here $\bs{\ac\Pi}{}_{(\al,d)}^\rpl$ is a principal $[*/\bG_m]$-bundle, as $\bs{\ac\Pi}{}_{(\al,d)}^\pl,\Pi^\pl,\bs\Pi^\pl$ are. As for \eq{co7eq60}, using \eq{co2eq12} for the Cartesian outer rectangles in \eq{co8eq53}--\eq{co8eq54} and equation \eq{co7eq59} yields isomorphisms
\e
\begin{split}
\bL_{\bs\acM_{(\al,d)}^\red/\bs\acM_{(\al,d)}}&\cong\Pi_{\Pic(X)}^*(\bL_{\Pic/\bs\Pic})\cong H^0(K_X)\ot\O_{\bs\acM_{(\al,d)}^\red}[2],\\
\bL_{\bs\acM_{(\al,d)}^\rpl/\bs\acM_{(\al,d)}^\pl}&\cong\Pi_{\Pic(X)^\pl}^*(\bL_{\Pic^\pl/\bs\Pic^\pl})\cong H^0(K_X)\ot\O_{\bs\acM_{(\al,d)}^\rpl}[2].
\end{split}
\label{co8eq56}
\e

Taking $i^*$ of the distinguished triangle \eq{co2eq11} for $\bs j_{(\al,d)}:\bs\acM_{(\al,d)}^\red\ra\bs\acM_{(\al,d)}$ and using \eq{co8eq56} gives a distinguished triangle on $\acM_{(\al,d)}$:
\e
\xymatrix@C=11pt{ {\begin{subarray}{l}\ts \quad H^0(K_X) \\ \ts {}\ot\!\O_{\acM_{(\al,d)}}[1] \end{subarray}} \ar[r] & i^*(\bL_{\bs\acM_{(\al,d)}}) \ar[rrr]^{i^*(\bL_{\bs j_{(\al,d)}})} &&& i^*(\bL_{\bs\acM_{(\al,d)}^\red}) \ar[r] & {\begin{subarray}{l}\ts \quad H^0(K_X) \\ \ts {}\ot\!\O_{\acM_{(\al,d)}}[2].\! \end{subarray}} }
\label{co8eq57}
\e	
\end{dfn}

Now we can define the data and verify Assumptions \ref{co5ass1}(f) and \ref{co5ass2}(e):

\begin{dfn}
\label{co8def12}
Continue in the situation of Definition \ref{co8def11}. Let $(\al,d)\in C(\acB)_\pe$ in \eq{co8eq40}, so that $\rank\al>0$ and $d\le 1$. Divide into cases:
\begin{itemize}
\setlength{\itemsep}{0pt}
\setlength{\parsep}{0pt}
\item[(A)] $d=\rank\al=1$, and 
\item[(B)] $d=0$ and $\rank\al>0$, or $d=1$ and $\rank\al>1$.
\end{itemize}

In case (A), for reasons explained in Remark \ref{co8rem4}, we define $\dM_{(\al,d)}\subseteq\acM_{(\al,d)}$, $\dM_{(\al,d)}^\pl\subseteq\acM_{(\al,d)}^\pl$, $\bs\dM_{(\al,d)}^\red=\bs\dM_{(\al,d)}$, $\bs\dM_{(\al,d)}^\rpl=\bs\dM_{(\al,d)}^\pl$, $\bs j_{(\al,d)}=\bs\id=\bs j_{(\al,d)}^\pl$, $\bs{\dot\Pi}{}_{(\al,d)}^\rpl=\bs{\dot\Pi}{}_{(\al,d)}^\pl$, and $U_{(\al,d)}=o_{(\al,d)}=0$ as in Definition \ref{co8def10}, the `unreduced' case, so that Assumptions \ref{co5ass1}(f) and \ref{co5ass2}(e) hold.

In case (B), for Assumption \ref{co5ass1}(f)(i), as for \eq{co8eq49} define open substacks $\dM_{(\al,d)}\subseteq\acM_{(\al,d)}$ and $\dM_{(\al,d)}^\pl\subseteq\acM_{(\al,d)}^\pl$ by
\e
\begin{split}
\dM_{(\al,d)}&=\acM_{(\al,d)}\sm\supp h^2\bigl((\ac\imath_{(\al,d)}^\red)^*(\bT_{\bs\acM_{(\al,d)}^\red})\bigr),\\
\dM_{(\al,d)}^\pl&=\acM_{(\al,d)}^\pl\sm\supp h^2\bigl((\ac\imath_{(\al,d)}^\red)^*(\bT_{\bs\acM_{(\al,d)}^\rpl})\bigr).
\end{split}
\label{co8eq58}
\e
As $\bs{\ac\Pi}{}_{(\al,d)}^\rpl:\bs\acM_{(\al,d)}^\red\ra\bs\acM_{(\al,d)}^\rpl$ is a principal $[*/\bG_m]$-bundle we have 
\begin{equation*}
h^2\bigl((\ac\imath_{(\al,d)}^\red)^*(\bT_{\bs\acM_{(\al,d)}^\red})\bigr)=(\Pi_\al^\pl)^*\bigl(h^2\bigl((\ac\imath_{(\al,d)}^\red)^*(\bT_{\bs\acM_{(\al,d)}^\rpl})\bigr)\bigr),
\end{equation*}
so that $\dM_{(\al,d)}=(\ac\Pi_{(\al,d)}^\pl)^{-1}(\dM_{(\al,d)}^\pl)$. 

For Assumption \ref{co5ass1}(f)(ii) define $\bs\dM_{(\al,d)}\subseteq\bs\acM_{(\al,d)}$, $\bs\dM_{(\al,d)}^\red\subseteq\bs\acM_{(\al,d)}^\red$, $\bs\dM_{(\al,d)}^\pl\subseteq\bs\acM_{(\al,d)}^\pl$, $\bs\dM_{(\al,d)}^\rpl\subseteq\bs\acM_{(\al,d)}^\rpl$ to be the derived open substacks with classical truncations $\dM_{(\al,d)}\subseteq\acM_{(\al,d)}$ and $\dM_{(\al,d)}^\pl\subseteq\acM_{(\al,d)}^\pl$. Define $\bs j_{(\al,d)},\ab\bs j^{\,\pl}_{(\al,d)},\ab\bs{\dot\Pi}{}_{(\al,d)}^\rpl,\ab\bs{\dot\Pi}{}_{(\al,d)}^\pl$ in \eq{co5eq4} to be the restrictions of $\bs{\ac\jmath}_{(\al,d)},\bs{\ac\jmath}{}^{\,\pl}_{(\al,d)},\bs{\ac\Pi}{}_{(\al,d)}^\rpl,\bs{\ac\Pi}{}_{(\al,d)}^\pl$ in \eq{co8eq53}--\eq{co8eq55} to~$\bs\dM_{(\al,d)},\bs\dM_{(\al,d)}^\red,\bs\dM_{(\al,d)}^\pl,\bs\dM_{(\al,d)}^\rpl$. 

Then \eq{co5eq4} for $\acA$ is Cartesian, as it is the restriction of the Cartesian left hand diamond in \eq{co8eq55} to open substacks. Also $t_0(\bs j_{(\al,d)}),t_0(\bs j^{\,\pl}_{(\al,d)})$ are isomorphisms as $t_0(\bs{\ac\jmath}_{(\al,d)}),t_0(\bs{\ac\jmath}_{(\al,d)}^{\,\pl})$ are. Definition \ref{co8def11} showed that $\bs\acM_{(\al,d)}^\red,\ab\bs\acM_{(\al,d)}^\rpl$ are locally finitely presented, so $\bs\dM_{(\al,d)}^\red,\ab\bs\dM_{(\al,d)}^\rpl$ are too.

To show that $\bs\dM_{(\al,d)}^\red,\bs\dM_{(\al,d)}^\rpl$ are quasi-smooth, by Remark \ref{co2rem6}(e) it suffices to check that $i^*(\bL_{\bs\dM_{(\al,d)}^\red}),i^*(\bL_{\bs\dM_{(\al,d)}^\rpl})$ are perfect in the interval $[-1,1]$. Equation \eq{co8eq41} says that $i^*(\bL_{\bs\dM_{(\al,d)}})$ is perfect in $[-2,1]$. From \eq{co8eq57} restricted to $\dM_{(\al,d)}$ we deduce that $i^*(\bL_{\bs\dM_{(\al,d)}^\red})$ is also perfect in $[-2,1]$. 

Thus the argument at the end of Definition \ref{co7def16} shows that $i^*(\bL_{\bs\dM_{(\al,d)}^\red})$ is perfect in $[-1,1]$ if $h^2(i^*(\bL_{\bs\dM_{(\al,d)}^\red})^\vee)=0$. But this holds by definition of $\dM_{(\al,d)}$ in \eq{co8eq58}. Therefore $\bs\dM_{(\al,d)}^\red$ is quasi-smooth. The proof of \eq{co7eq49} in Remark \ref{co7rem7} also works for the principal $[*/\bG_m]$-fibration $\bs{\dot\Pi}{}_{(\al,d)}^\rpl:\bs\dM_{(\al,d)}^\red\ra\bs\dM_{(\al,d)}^\rpl$. Hence $i^*(\bL_{\bs\dM_{(\al,d)}^\rpl})$ is perfect in $[-1,1]$, and $\bs\dM_{(\al,d)}^\rpl$ is quasi-smooth.

For Assumption \ref{co5ass1}(f)(iii) set $U_{(\al,d)}\!=\!H^0(K_X)$ and $o_{(\al,d)}\!=\!\dim H^0(K_X)\!=\!p_g$. The isomorphisms \eq{co5eq6} are \eq{co8eq56} restricted to~$\bs\dM_{(\al,d)}^\red,\bs\dM_{(\al,d)}^\rpl$. For Assumption \ref{co5ass1}(f)(iv), $o_{(\al,d)}+o_{(\be,e)}\ge o_{(\al+\be,d+e)}$ holds as $o_{(\ga,f)}=0$ if $\rank\ga=d=1$ and $o_{(\ga,f)}=p_g>0$ otherwise. This completes Assumption~\ref{co5ass1}(f).

For Assumption \ref{co5ass2}(e), suppose $[E,V,\rho]$ lies in $\acM_{(\al,d)}^\ss(\ac\mu^{\om,s})$ for some $s>0$. As for \eq{co8eq50}, using equations \eq{co7eq54}, \eq{co7eq59}, \eq{co7eq61}, \eq{co8eq39}, \eq{co8eq50}, \eq{co8eq53} and \eq{co8eq57} we can construct an exact sequence
\e
\!\!\xymatrix@C=13pt{
\Ext^2(E,E) \ar[rr]^(0.42){\begin{subarray}{l} -\ci\rho \\ \op\Tr_E \end{subarray}} &&
{\begin{subarray}{l}\ts V^*\ot \Ext^2(L,E)  \\ \ts \op\Ext^2(\O_X,\O_X) \end{subarray}} \ar[r] & h^2\bigl(i^*(\bT_{\bs\dM_{(\al,d)}^\red})\bigr)\big\vert_{[E,V,\rho]} \ar[r] & 0. }\!\!
\label{co8eq59}
\e
Here by Serre duality $\Ext^2(\O_X,\O_X)\cong\Hom(\O_X,K_X)^*=H^0(K_X)^*=U_{(\al,d)}^*$, and the $\Ext^2(\O_X,\O_X)$ contribution in \eq{co8eq59} (the new term not in \eq{co8eq50}) comes from $\bT_{\Pic/\bs\Pic}$ in the dual of \eq{co7eq59}, which enters because of the morphism $i:\Pic(X)\ra\bs\Pic(X)$ in \eq{co8eq53}.

We want to show that $h^2(i^*(\bT_{\bs\dM_{(\al,d)}^\red}))\vert_{[E,V,\rho]}=0$, that is, that the first morphism in \eq{co8eq59} is surjective. As for \eq{co8eq51}, by Serre duality, this is equivalent to the following morphism being injective:
\e
\xymatrix@C=60pt{
{\begin{subarray}{l}\ts \Hom(E,V\ot_\C L\ot K_X) \\ \ts\op \Hom(\O_X,K_X) \end{subarray}} \ar[r]^(0.55){\begin{subarray}{l} (\rho\ot\id_{K_x})\ci - \\ \op (\id_E\ot -) \end{subarray}} & \Hom(E,E\ot K_X). }
\label{co8eq60}
\e

Now $E$ is torsion-free by Lemma \ref{co8lem2}, and if $d=1$ then $\rho\ne 0$ by Theorem \ref{co8thm2}(a). Thus there is a dense open subscheme $X'\subseteq X$ such that $E\vert_{X'}$ is a vector bundle, and if $d=1$ then $\rho:V\ot_\C L\ra E$ is a vector bundle morphism on $X'$ which is nonzero at each $x\in X'$. As in \eq{co8eq52}, for each $x\in X'$ consider the restriction of \eq{co8eq60} to $x$, the linear map of $\C$-vector spaces
\e
\xymatrix@C=65pt{
{\begin{subarray}{l}\ts \Hom(E\vert_x,V\ot_\C L\vert_x\ot K_X\vert_x) \\ \ts \op \Hom(\O_X\vert_x,K_X\vert_x) \end{subarray}} \ar[r]^(0.55){\begin{subarray}{l} (\rho\vert_x\ot\id_{K_X}\vert_x)\ci - \\ \op (\id_{E\vert_x}\ot -) \end{subarray}} & \Hom(E\vert_x,E\vert_x\ot K_X\vert_x), }
\label{co8eq61}
\e
Here $\dim_\C E\vert_x=\rank\al$, $\dim_\C V=d=0$ or 1, and $\dim_\C L\vert_x=\dim_\C K_X\vert_x=1$. 

By (B) above we may divide into two cases:
\begin{itemize}
\setlength{\itemsep}{0pt}
\setlength{\parsep}{0pt}
\item[(B$)'\phantom{'}$] $V=0$ and $\dim_\C E\vert_x>0$, and
\item[(B$)''$] $\dim_\C V=1$, $\dim_\C E\vert_x>1$, and $\rho\vert_x\ne 0$.
\end{itemize}
We can show using elementary linear algebra that \eq{co8eq61} is injective in both cases, for all $x\in X'$. Therefore \eq{co8eq60} is injective, so $h^2(i^*(\bT_{\bs\dM_{(\al,d)}^\red}))\vert_{[E,V,\rho]}=0$ by \eq{co8eq59}, and $[E,V,\rho]\in\dM_{(\al,d)}$ by \eq{co8eq58}, so also $[E,V,\rho]\in\dM_{(\al,d)}^\pl$ as $\dM_{(\al,d)}=(\ac\Pi_{(\al,d)}^\pl)^{-1}(\dM_{(\al,d)}^\pl)$. Thus $\acM_{(\al,d)}^\ss(\ac\mu^{\om,s})\subseteq\dM_{(\al,d)}^\pl$, proving Assumption~\ref{co5ass2}(e).
\end{dfn}

\begin{rem}
\label{co8rem4}
The proof of Assumption \ref{co5ass2}(e) above is based on Mochizuki \cite[Prop.s~6.1.1 \& 6.1.2]{Moch}, which only define `reduced' obstruction theories on pair moduli spaces like $\acM_{(\al,1)}^\ss(\ac\mu^{\om,s})$ when~$\rank\al>1$. In Definition \ref{co8def12} we used the `non-reduced' quasi-smooth stacks $\bs\dM_{(\al,d)}^\red=\bs\dM_{(\al,d)}$, $\bs\dM_{(\al,d)}^\rpl=\bs\dM_{(\al,d)}^\pl$ and obstruction theories \eq{co5eq5} in case (A) and the `reduced' versions in case (B). The reason the `reduced' versions do not work in case (A) is that if $d=\rank\al=1$ then \eq{co8eq61} is not injective (it is a linear map $\C^2\ra\C$). Thus \eq{co8eq60} might not be injective, so we could have $\acM_{(\al,d)}^\ss(\ac\mu^{\om,s})\subsetneq\dM_{(\al,d)}^\pl$, and Assumption \ref{co5ass2}(e) would fail.
\end{rem}

\subsubsection{The $G$-equivariant case and Assumption \ref{co5ass4}(c)}
\label{co834}

In \S\ref{co825} we explained how to extend \S\ref{co821}--\S\ref{co824} to the $G$-equivariant case when an algebraic $\C$-group $G$ acts on $X$, under the extra conditions Definition \ref{co8def7}(i)--(iii), and verified Assumptions \ref{co4ass2} and \ref{co5ass4}(a),(b),(d). As in \S\ref{co735}, now that we have defined the data of Assumption \ref{co5ass1}(f) in \S\ref{co831}--\S\ref{co833}, we can also verify Assumption \ref{co5ass4}(c). In fact this is immediate: all the constructions of \S\ref{co831}--\S\ref{co833} are manifestly invariant or equivariant under the action of $G$ on~$X,\coh(X),\acA,\acM_{(\al,d)},H^0(K_X),\ldots.$ 

\subsection{Counting pairs on curves}
\label{co84}

Let $X$ be a connected smooth projective complex curve, and fix a line bundle $L\ra X$ and a K\"ahler class $\om\in\Kah(X)$ on $X$. Then \S\ref{co81}--\S\ref{co83} define data and verify Assumptions \ref{co4ass1} and \ref{co5ass1}--\ref{co5ass3} for the category $\acA=\acB$ of triples $(E,V,\rho)$ where $E\in\coh(X)$, $V\in\Vect_\C$ and $\rho:V\ot_\C L\ra E$ is a morphism in $\coh(X)$. Thus Theorems \ref{co5thm1}--\ref{co5thm3} apply to define invariants $[\acM_{(\al,d)}^\ss(\ac\mu^{\om,s})]_\inv$ in $H_*(\acM_{(\al,d)}^\pl)$ counting $\ac\mu^{\om,s}$-semistable triples $(E,V,\rho)\in\acA$ for all $s>0$, where as in \eq{co8eq40} we define $[\acM_{(\al,d)}^\ss(\ac\mu^{\om,s})]_\inv$ only when $\al\in C(\coh(X))$ with $\rank\al>0$ and $d=0$ or 1. The wall-crossing formulae \eq{co5eq33}--\eq{co5eq34} hold for changing between $(\ac\mu^{\om,s},\ac M,\le)$ and $(\ac\mu^{\om,\ti s},\ac M,\le)$ for $s,\ti s>0$. We will write $[\acM_{(\al,d)}^\ss(\ac\mu^{\om,s})]_\inv=[\acM_{(\al,d)}^\ss(\ac\mu^{\om,s})]_\inv^L$ to emphasize the dependence on $L\ra X$.

The inclusion $\coh(X)\hookra\acA$ mapping $E\mapsto (E,0,0)$ induces isomorphisms $\M_\al^\pl\ra\acM_{(\al,0)}^\pl$ which identify $\M_\al^\ss(\mu^\om)$ in \S\ref{co76} with $\acM_{(\al,0)}^\ss(\ac\mu^{\om,s})$ when $d=0$, and the corresponding isomorphisms $H_*(\M_\al^\pl)\ra H_*(\acM_{(\al,0)}^\pl)$ identify the invariants $[\M_\al^\ss(\mu^\om)]_\inv$ in Theorem \ref{co7thm6} when $\rank\al>0$ with $[\acM_{(\al,0)}^\ss(\ac\mu^{\om,s})]^L_\inv$. Thus, we have embedded the invariants $[\M_\al^\ss(\mu^\om)]_\inv$ in Theorem \ref{co7thm6} inside a larger system of invariants. We explain in \S\ref{co86} how to use this larger system to compute the $[\M_\al^\ss(\mu^\om)]_\inv$ for $\rank\al>1$ in terms of rank 1 invariants.

When $X=\CP^1$ the algebraic $\C$-group $G=\bG_m$ or $\PGL(2,\C)$ acts on $X$ satisfying Definition \ref{co8def7}(i)--(iii), and \S\ref{co825} and \S\ref{co834} verify Assumptions \ref{co4ass2} and \ref{co5ass4}, so Theorem \ref{co5thm4} applies, and we can promote the $[\acM_{(\al,d)}^\ss(\ac\mu^{\om,s})]^L_\inv$ to equivariant homology~$H_*^G(\acM_{(\al,d)}^\pl)$.

In \S\ref{co76}(b) we observed that the notion of $\mu^\om$-(semi)stability on $\coh(X)$ is independent of $\om\in\Kah(X)$, so the wall crossing in Theorems \ref{co5thm2}--\ref{co5thm3} is trivial for $[\M_\al^\ss(\mu^\om)]_\inv$. However, $\ac\mu^{\om,s}$-(semi)stability on $\acA$ does depend on $s\in(0,\iy)$, so we do have nontrivial wall crossing for the $[\acM_{(\al,d)}^\ss(\ac\mu^{\om,s})]^L_\inv$ when $d=1$.

In \S\ref{co76}(c) we noted that moduli stacks $\M_\al,\M_\al^\pl$ of objects in $\coh(X)$ are smooth, so we can replace virtual classes $[\M_\al^\ss(\tau)]_\virt$ by fundamental classes $[\M_\al^\ss(\tau)]_\fund$. The moduli stacks $\acM_{(\al,d)},\acM_{(\al,d)}^\pl$ of objects in $\acA$ are in general not smooth for $d>0$ (although $\acM_{(\al,d)}^\ss(\ac\mu^{\om,s})$ is smooth when $L=\O_X(-N)$ for $N\gg 0$), so we do need virtual classes in this case.

The next theorem summarizes parts of Theorems \ref{co5thm1}, \ref{co5thm3} and \ref{co5thm4} for $\acA$, and the discussion above. Parts (ii),(iii) use Theorem \ref{co8thm2}(d),(g).

\begin{thm}
\label{co8thm3}
Let $X$ be a connected, smooth, projective, complex curve, and fix a line bundle\/ $L\ra X$ and a K\"ahler class $\om\in\Kah(X)$. 

As in\/ {\rm\S\ref{co81},} write $\acA$ for the $\C$-linear abelian category of triples $(E,V,\rho)$ where $E\in\coh(X),$ $V\in\Vect_\C$ and\/ $\rho:V\ot_\C L\ra E$ is a morphism in $\coh(X),$ and write $\acM^\pl=\coprod_{(\al,d)\in K(\acA)}\acM^\pl_{(\al,d)}$ for the projective linear moduli stack of objects in $\acA,$ and define a weak stability condition $(\ac\mu^{\om,s},\ac M,\le)$ on $\acA$ for all\/ $s>0,$ giving moduli stacks $\acM_{(\al,d)}^\rst(\ac\mu^{\om,s})\subseteq\acM_{(\al,d)}^\ss(\ac\mu^{\om,s})\subseteq\acM^\pl_{(\al,d)}$ for all\/~$(\al,d)\in C(\acA)$.

Then for all\/ $\al$ in\/ $C(\coh(X))$ with\/ $\rank\al>0$ and\/ $d=0$ or $1$ there are unique classes\/ $[\acM_{(\al,d)}^\ss(\ac\mu^{\om,s})]^L_\inv$ in the Betti\/ $\Q$-homology group $\check H_0(\acM^\pl_{(\al,d)})=H_{2-2\ac\chi((\al,d),(\al,d))}(\acM^\pl_{(\al,d)})$ for\/ $\ac\chi$ as in \eq{co8eq27} satisfying:
\begin{itemize}
\setlength{\itemsep}{0pt}
\setlength{\parsep}{0pt}
\item[{\bf(i)}] The inclusion $\coh(X)\hookra\acA$ mapping $E\mapsto (E,0,0)$ induces isomorphisms $\M_\al^\pl\ra\acM_{(\al,0)}^\pl$ for $\al\in C(\coh(X))$. When $\rank\al>0,$ the corresponding isomorphism $\check H_0(\M_\al^\pl)\ra\check H_0(\acM^\pl_{(\al,0)})$ identifies $[\M_\al^\ss(\mu^\om)]_\inv$ in Theorem\/ {\rm\ref{co7thm6}} with\/ $[\acM_{(\al,0)}^\ss(\ac\mu^{\om,s})]^L_\inv$. Thus $[\acM_{(\al,0)}^\ss(\ac\mu^{\om,s})]^L_\inv$ is independent of\/ $\om,s$. 
\item[{\bf(ii)}] If\/ $\acM_{(\al,d)}^\rst(\ac\mu^{\om,s})\!=\!\acM_{(\al,d)}^\ss(\ac\mu^{\om,s})$ then\/ $[\acM_{(\al,d)}^\ss(\ac\mu^{\om,s})]^L_\inv\!=\![\acM_{(\al,d)}^\ss(\ac\mu^{\om,s})]_\virt,$ where $[\acM_{(\al,d)}^\ss(\ac\mu^{\om,s})]_\virt$ is the virtual class of\/ $\acM_{(\al,d)}^\ss(\ac\mu^{\om,s}),$ which is a proper algebraic space with the obstruction theory defined in\/~{\rm\S\ref{co83}}.

If\/ $\rank\al>0$ and\/ $d=1$ then $\acM_{(\al,1)}^\rst(\ac\mu^{\om,s})\!=\!\acM_{(\al,1)}^\ss(\ac\mu^{\om,s})$ for all but finitely many $s\in(0,\iy)$.
\item[{\bf(iii)}] The invariant\/ $[\baM^\ss_{(\al,1)}(\bar\mu^0_1)]_\fund$ used in Theorem\/ {\rm\ref{co7thm6}(iii)} coincides with\/ $[\acM_{(\al,1)}^\ss(\ac\mu^{\om,s})]^L_\inv=[\acM_{(\al,1)}^\ss(\ac\mu^{\om,s})]_\virt$ in the case when $L=\O_X(-N)$ for $N\gg 0,$ and\/ $s>0$ is sufficiently small.
\item[{\bf(iv)}] Let\/ $\al\in C(\coh(X))$ with\/ $\rank\al>0,$ and\/ $s,\ti s>0$. Then
\ea
\!\!\!\!\!\!\!\!&[\acM_{(\al,1)}^\ss(\ac\mu^{\om,\ti s})]^L_\inv=
\label{co8eq62}\\
\!\!\!\!\!\!\!\!&\sum_{\begin{subarray}{l}n\ge 1,\;\al_1,\ldots,\al_n\in C(\coh(X)), \\
d_1,\ldots,d_n\in\{0,1\}:\; \rank\al_i>0, \\
\acM_{(\al_i,d_i)}^\ss(\ac\mu^{\om,s})\ne\es,\; \text{all\/ $i,$} \\
\al_1+\cdots+\al_n=\al, \; d_1+\cdots+d_n=1 
\end{subarray}} \begin{aligned}[t]
&\ti U((\al_1,d_1),\ldots,(\al_n,d_n);\ac\mu^{\om,s},\ac\mu^{\om,\ti s})\cdot {} \\
&\bigl[\bigl[\cdots\bigl[[\acM_{(\al_1,d_1)}^\ss(\ac\mu^{\om,s})]^L_\inv,\\
&
[\acM_{(\al_2,d_2)}^\ss(\ac\mu^{\om,s})]^L_\inv\bigr],\ldots\bigr],[\acM_{(\al_n,d_n)}^\ss(\ac\mu^{\om,s})]^L_\inv\bigr]
\end{aligned}
\nonumber
\ea
in the Lie algebra $\check H_0(\acM^\pl)$ from {\rm\S\ref{co43}}. Here $\ti U(-;\ac\mu^{\om,s},\ac\mu^{\om,\ti s})$ is as in Theorem\/ {\rm\ref{co3thm3},} and there are only finitely many nonzero terms in \eq{co8eq62}. Equivalently, in the universal enveloping algebra $\bigl(U(\check H_0(\acM^\pl)),*\bigr)$ we have:
\begin{align*}
&[\acM_{(\al,1)}^\ss(\ac\mu^{\om,\ti s})]^L_\inv=\\
&\sum_{\begin{subarray}{l}n\ge 1,\;\al_1,\ldots,\al_n\in C(\coh(X)), \\
d_1,\ldots,d_n\in\{0,1\}: \; \rank\al_i>0, \\
\acM_{(\al_i,d_i)}^\ss(\ac\mu^{\om,s})\ne\es,\; \text{all\/ $i,$} \\
\al_1+\cdots+\al_n=\al, \; d_1+\cdots+d_n=1 
\end{subarray}} \begin{aligned}[t]
&U((\al_1,d_1),\ldots,(\al_n,d_n);\ac\mu^{\om,s},\ac\mu^{\om,\ti s})\cdot {} \\
&[\acM_{(\al_1,d_1)}^\ss(\ac\mu^{\om,s})]^L_\inv *[\acM_{(\al_2,d_2)}^\ss(\ac\mu^{\om,s})]^L_\inv\\
&*\cdots *[\acM_{(\al_n,d_n)}^\ss(\ac\mu^{\om,s})]^L_\inv.
\end{aligned}
\end{align*}
\end{itemize}

When $X=\CP^1,$ so that\/ $G=\bG_m$ or $\PGL(2,\C)$ acts on $X,$ the analogue holds in $G$-equivariant homology $H_*^G(\acM^\pl),$ as in\/ {\rm\S\ref{co23}}.
\end{thm}

\begin{rem}
\label{co8rem5}
When $X$ is a curve, the restriction in \eq{co8eq40} to $(\al,d)\in C(\acA)$ with $\rank\al>0$ and $d=0,1$ is probably unnecessary, Theorem \ref{co8thm3} should extend to invariants $[\acM_{(\al,d)}^\ss(\ac\mu^{\om,s})]^L_\inv$ for all $(\al,d)\in C(\acA)$.
\end{rem}

\subsection{Counting pairs on surfaces}
\label{co85}

Let $X$ be a connected smooth projective complex surface with geometric genus $p_g\ge 0$, and fix a line bundle $L\ra X$ and a K\"ahler class $\om\in\Kah(X)$ on $X$. Then \S\ref{co81}--\S\ref{co83} define data and verify Assumptions \ref{co4ass1} and \ref{co5ass1}--\ref{co5ass3} for the category $\acA=\acB$ of triples $(E,V,\rho)$ where $E\in\coh(X)$, $V\in\Vect_\C$ and $\rho:V\ot_\C L\ra E$ is a morphism in $\coh(X)$. Thus Theorems \ref{co5thm1}--\ref{co5thm3} apply to define invariants $[\acM_{(\al,d)}^\ss(\ac\mu^{\om,s})]^L_\inv=[\acM_{(\al,d)}^\ss(\ac\mu^{\om,s})]_\inv$ in $H_*(\acM_{(\al,d)}^\pl)$ counting $\ac\mu^{\om,s}$-semistable triples $(E,V,\rho)\in\acA$ for all $s>0$, where as in \eq{co8eq40} we define $[\acM_{(\al,d)}^\ss(\ac\mu^{\om,s})]^L_\inv$ only when $\al\in C(\coh(X))$ with $\rank\al>0$ and $d=0$ or 1. The wall-crossing formulae \eq{co5eq33}--\eq{co5eq34} hold for changing between $\ac\mu^{\om,s}$ and $\ac\mu^{\om,\ti s}$ for~$s,\ti s>0$. 

The inclusion $\coh(X)\hookra\acA$ mapping $E\mapsto (E,0,0)$ induces isomorphisms $\M_\al^\pl\ra\acM_{(\al,0)}^\pl$ which identify $\M_\al^\ss(\mu^\om)$ in \S\ref{co77} with $\acM_{(\al,0)}^\ss(\ac\mu^{\om,s})$ when $d=0$, and the isomorphisms $H_*(\M_\al^\pl)\ra H_*(\acM_{(\al,0)}^\pl)$ identify $[\M_\al^\ss(\mu^\om)]_\inv$ in Theorems \ref{co7thm7} and \ref{co7thm8} when $\rank\al>0$ with~$[\acM_{(\al,0)}^\ss(\ac\mu^{\om,s})]^L_\inv$. 

When an algebraic $\C$-group $G$ acts on $X$ satisfying Definition \ref{co8def7}(i)--(iii) sections \ref{co825} and \ref{co834} verify Assumptions \ref{co4ass2} and \ref{co5ass4}, so Theorem \ref{co5thm4} applies, and the $[\acM_{(\al,d)}^\ss(\ac\mu^{\om,s})]^L_\inv$ lift to equivariant homology~$H_*^G(\acM_{(\al,d)}^\pl)$.

As in \S\ref{co73}, \S\ref{co77}, and \S\ref{co83}, we distinguish the cases $p_g=0$ and $p_g>0$, when we use `non-reduced' and `reduced' stacks $\bs\dM_{(\al,d)}^\red,\bs\dM_{(\al,d)}^\rpl$ and obstruction theories respectively. The next theorem summarizes parts of Theorems \ref{co5thm1}, \ref{co5thm3} and \ref{co5thm4} for $\acA$ when $p_g=0$. Parts (ii),(iii) use Theorem~\ref{co8thm2}(d),(g).

\begin{thm}
\label{co8thm4}
Let\/ $X$ be a connected projective complex surface with geometric genus\/ $p_g=0,$ and fix a line bundle\/ $L\ra X$ and a K\"ahler class $\om\in\Kah(X)$. 

As in\/ {\rm\S\ref{co81},} write $\acA$ for the $\C$-linear abelian category of triples $(E,V,\rho)$ where $E\in\coh(X),$ $V\in\Vect_\C$ and\/ $\rho:V\ot_\C L\ra E$ is a morphism in $\coh(X),$ and write $\acM^\pl=\coprod_{(\al,d)\in K(\acA)}\acM^\pl_{(\al,d)}$ for the projective linear moduli stack of objects in $\acA,$ and define a weak stability condition $(\ac\mu^{\om,s},\ac M,\le)$ on $\acA$ for all\/ $s>0,$ giving moduli stacks $\acM_{(\al,d)}^\rst(\ac\mu^{\om,s})\subseteq\acM_{(\al,d)}^\ss(\ac\mu^{\om,s})\subseteq\acM^\pl_{(\al,d)}$ for all\/~$(\al,d)\in C(\acA)$.

Then for all\/ $\al$ in\/ $C(\coh(X))$ with\/ $\rank\al>0$ and\/ $d=0$ or $1$ there are unique classes\/ $[\acM_{(\al,d)}^\ss(\ac\mu^{\om,s})]^L_\inv$ in the Betti\/ $\Q$-homology group $\check H_0(\acM^\pl_{(\al,d)})=H_{2-2\ac\chi((\al,d),(\al,d))}(\acM^\pl_{(\al,d)})$ for\/ $\ac\chi$ as in \eq{co8eq27} satisfying:
\begin{itemize}
\setlength{\itemsep}{0pt}
\setlength{\parsep}{0pt}
\item[{\bf(i)}] The inclusion $\coh(X)\hookra\acA$ mapping $E\mapsto (E,0,0)$ induces isomorphisms $\M_\al^\pl\ra\acM_{(\al,0)}^\pl$ for $\al\in C(\coh(X))$. When $\rank\al>0,$ the corresponding isomorphism $\check H_0(\M_\al^\pl)\ra\check H_0(\acM^\pl_{(\al,0)})$ identifies $[\M_\al^\ss(\mu^\om)]_\inv$ in Theorem\/ {\rm\ref{co7thm7}} with\/ $[\acM_{(\al,0)}^\ss(\ac\mu^{\om,s})]^L_\inv$. Thus $[\acM_{(\al,0)}^\ss(\ac\mu^{\om,s})]^L_\inv$ is independent of\/ $s$. 
\item[{\bf(ii)}] If\/ $\acM_{(\al,d)}^\rst(\ac\mu^{\om,s})\!=\!\acM_{(\al,d)}^\ss(\ac\mu^{\om,s})$ then\/ $[\acM_{(\al,d)}^\ss(\ac\mu^{\om,s})]^L_\inv\!=\![\acM_{(\al,d)}^\ss(\ac\mu^{\om,s})]_\virt,$ where $[\acM_{(\al,d)}^\ss(\ac\mu^{\om,s})]_\virt$ is the virtual class of\/ $\acM_{(\al,d)}^\ss(\ac\mu^{\om,s}),$ which is a proper algebraic space with the obstruction theory defined in\/~{\rm\S\ref{co83}}.

If\/ $\rank\al>0$ and\/ $d=1$ then $\acM_{(\al,1)}^\rst(\ac\mu^{\om,s})\!=\!\acM_{(\al,1)}^\ss(\ac\mu^{\om,s})$ for all but finitely many $s\in(0,\iy)$.
\item[{\bf(iii)}] The invariant\/ $[\baM^\ss_{(\al,1)}(\bar\tau^0_1)]_\virt$ used in Theorem\/ {\rm\ref{co7thm7}(a)(ii)} when $(\tau,T,\le)\ab=(\mu^\om,M,\le)$ coincides with\/ $[\acM_{(\al,1)}^\ss(\ac\mu^{\om,s})]^L_\inv=[\acM_{(\al,1)}^\ss(\ac\mu^{\om,s})]_\virt$ in the case when $L=\O_X(-N)$ for $N\gg 0,$ and\/ $s>0$ is sufficiently small.
\item[{\bf(iv)}] Let\/ $\al\in C(\coh(X))$ with\/ $\rank\al>0,$ and\/ $s,\ti s>0$. Then
\ea
\!\!\!\!\!\!\!\!&[\acM_{(\al,1)}^\ss(\ac\mu^{\om,\ti s})]^L_\inv=
\label{co8eq63}\\
\!\!\!\!\!\!\!\!&\sum_{\begin{subarray}{l}n\ge 1,\;\al_1,\ldots,\al_n\in C(\coh(X)), \\
d_1,\ldots,d_n\in\{0,1\}:\; \rank\al_i>0, \\
\acM_{(\al_i,d_i)}^\ss(\ac\mu^{\om,s})\ne\es,\; \text{all\/ $i,$} \\
\al_1+\cdots+\al_n=\al, \; d_1+\cdots+d_n=1 
\end{subarray}} \begin{aligned}[t]
&\ti U((\al_1,d_1),\ldots,(\al_n,d_n);\ac\mu^{\om,s},\ac\mu^{\om,\ti s})\cdot {} \\
&\bigl[\bigl[\cdots\bigl[[\acM_{(\al_1,d_1)}^\ss(\ac\mu^{\om,s})]^L_\inv,\\
&
[\acM_{(\al_2,d_2)}^\ss(\ac\mu^{\om,s})]^L_\inv\bigr],\ldots\bigr],[\acM_{(\al_n,d_n)}^\ss(\ac\mu^{\om,s})]^L_\inv\bigr]
\end{aligned}
\nonumber
\ea
in the Lie algebra $\check H_0(\acM^\pl)$ from {\rm\S\ref{co43}}. Here $\ti U(-;\ac\mu^{\om,s},\ac\mu^{\om,\ti s})$ is as in Theorem\/ {\rm\ref{co3thm3},} and there are only finitely many nonzero terms in \eq{co8eq63}. Equivalently, in the universal enveloping algebra $\bigl(U(\check H_0(\acM^\pl)),*\bigr)$ we have:
\ea
\!\!\!\!\!\!\!\!&[\acM_{(\al,1)}^\ss(\ac\mu^{\om,\ti s})]^L_\inv=
\label{co8eq64}\\
\!\!\!\!\!\!\!\!&\sum_{\begin{subarray}{l}n\ge 1,\;\al_1,\ldots,\al_n\in C(\coh(X)), \\
d_1,\ldots,d_n\in\{0,1\}: \; \rank\al_i>0, \\
\acM_{(\al_i,d_i)}^\ss(\ac\mu^{\om,s})\ne\es,\; \text{all\/ $i,$} \\
\al_1+\cdots+\al_n=\al, \; d_1+\cdots+d_n=1 
\end{subarray}} \begin{aligned}[t]
&U((\al_1,d_1),\ldots,(\al_n,d_n);\ac\mu^{\om,s},\ac\mu^{\om,\ti s})\cdot {} \\
&[\acM_{(\al_1,d_1)}^\ss(\ac\mu^{\om,s})]^L_\inv *[\acM_{(\al_2,d_2)}^\ss(\ac\mu^{\om,s})]^L_\inv\\
&*\cdots *[\acM_{(\al_n,d_n)}^\ss(\ac\mu^{\om,s})]^L_\inv.
\end{aligned}
\nonumber
\ea
\end{itemize}

If an algebraic $\C$-group $G$ acts on $X$ and Definition\/ {\rm\ref{co8def7}(i)--(iii)} hold then the above generalizes to equivariant homology $H_*^G(\acM^\pl)$.
\end{thm}

Next we give the analogue of Theorem \ref{co8thm4} for $p_g>0$. There are several important differences with Theorem~\ref{co8thm4}:
\begin{itemize}
\setlength{\itemsep}{0pt}
\setlength{\parsep}{0pt}
\item[(a)] As in \S\ref{co833}, if $\rank\al=1$ and $d=1$ then we use the `non-reduced' quasi-smooth stacks $\bs\dM_{(\al,d)}^\red=\bs\dM_{(\al,d)}$, $\bs\dM_{(\al,d)}^\rpl=\bs\dM_{(\al,d)}^\pl$ and associated obstruction theories \eq{co5eq5}, and set $o_{(\al,d)}=0$. Otherwise we use different `reduced' versions, and set $o_{(\al,d)}=p_g$. The invariants $[\acM_{(\al,d)}^\ss(\ac\mu^{\om,s})]^L_\inv$ lie in $\check H_{2o_{(\al,d)}}(\acM^\pl_{(\al,d)})$ by Theorem~\ref{co5thm1}.
\item[(b)] The wall-crossing formulae \eq{co5eq33}--\eq{co5eq34} in Theorem \ref{co5thm3} for $\acA$ include the restrictions $o_{(\al_1,d_1)}+\cdots+o_{(\al_n,d_n)}=o_{(\al,d)}$. Given the values of $o_{(\al,d)}$ in (a), this radically cuts down the possibilities for $n,(\al_i,d_i)$: we can only have $n=1$ and $(\al_1,d_1)=(\al,d)$, or $d=1$, $\rank\al>1$, $n=2$ and $\{(\al_1,d_1),(\al_2,d_2)\}=\{(\be,1),(\ga,0)\}$ where $\al=\be+\ga$ with $\rank\be=1$ and $\rank\ga>0$. As $U((\be,1),(\ga,0);\ac\mu^{\om,s},\ac\mu^{\om,\ti s})=-U((\ga,0),(\be,1);\ac\mu^{\om,s},\ac\mu^{\om,\ti s})$, we can rewrite the analogue of \eq{co8eq64} in the form~\eq{co8eq65}.
\item[(c)] In (i) below, $[\acM_{(\al,0)}^\ss(\ac\mu^{\om,s})]^L_\inv$ is independent of $\om$ by Theorem \ref{co7thm8}(a)(ii).
\item[(d)] In (iii) below, if $\rank\al>1$ then $[\acM_{(\al,1)}^\ss(\ac\mu^{\om,s})]_\virt$ and $[\baM^\ss_{(\al,1)}(\bar\tau^0_1)]_\virt$ are both defined using the `reduced' obstruction theory, and so they coincide as $\acM_{(\al,1)}^\ss(\ac\mu^{\om,s})=\baM^\ss_{(\al,1)}(\bar\tau^0_1)$ in this case by Theorem~\ref{co8thm2}(d). 

However, if $\rank\al=1$ then $[\acM_{(\al,1)}^\ss(\ac\mu^{\om,s})]_\virt$ uses the `non-reduced', and $[\baM^\ss_{(\al,1)}(\bar\tau^0_1)]_\virt$ the `reduced', obstruction theories. We can use Theorem \ref{co2thm1}(iv) to show that $[\acM_{(\al,1)}^\ss(\ac\mu^{\om,s})]_\virt=0$ in this case.
\end{itemize}

\begin{thm}
\label{co8thm5}
Let\/ $X$ be a connected projective complex surface with geometric genus\/ $p_g>0,$ and fix a line bundle\/ $L\ra X$ and a K\"ahler class $\om\in\Kah(X)$. 

As in\/ {\rm\S\ref{co81},} write $\acA$ for the $\C$-linear abelian category of triples $(E,V,\rho)$ where $E\in\coh(X),$ $V\in\Vect_\C$ and\/ $\rho:V\ot_\C L\ra E$ is a morphism in $\coh(X),$ and write $\acM^\pl=\coprod_{(\al,d)\in K(\acA)}\acM^\pl_{(\al,d)}$ for the projective linear moduli stack of objects in $\acA,$ and define a weak stability condition $(\ac\mu^{\om,s},\ac M,\le)$ on $\acA$ for all\/ $s>0,$ giving moduli stacks $\acM_{(\al,d)}^\rst(\ac\mu^{\om,s})\subseteq\acM_{(\al,d)}^\ss(\ac\mu^{\om,s})\subseteq\acM^\pl_{(\al,d)}$ for all\/~$(\al,d)\in C(\acA)$.

Then for all\/ $\al$ in\/ $C(\coh(X))$ with\/ $\rank\al>0$ and\/ $d=0$ or $1$ there are unique classes\/ $[\acM_{(\al,d)}^\ss(\ac\mu^{\om,s})]^L_\inv$ in the Betti\/ $\Q$-homology group $\check H_{2o_{(\al,d)}}(\acM^\pl_{(\al,d)})=H_{2-2\ac\chi((\al,d),(\al,d))+2o_{(\al,d)}}(\acM^\pl_{(\al,d)}),$ where $o_{(\al,d)}=0$ if\/ $\rank\al=1$ and\/ $d=1,$ and\/ $o_{(\al,d)}=p_g>0$ otherwise, and\/ $\ac\chi$ is as in\/ {\rm\eq{co8eq27},} satisfying:
\begin{itemize}
\setlength{\itemsep}{0pt}
\setlength{\parsep}{0pt}
\item[{\bf(i)}] The inclusion $\coh(X)\hookra\acA$ mapping $E\mapsto (E,0,0)$ induces isomorphisms $\M_\al^\pl\ra\acM_{(\al,0)}^\pl$ for $\al\in C(\coh(X))$. When $\rank\al>0,$ the corresponding isomorphism $H_*(\M_\al^\pl)\ra H_*(\acM^\pl_{(\al,0)})$ identifies $[\M_\al^\ss(\mu^\om)]_\inv$ in Theorem\/ {\rm\ref{co7thm8}} with\/ $[\acM_{(\al,0)}^\ss(\ac\mu^{\om,s})]^L_\inv$. Thus $[\acM_{(\al,0)}^\ss(\ac\mu^{\om,s})]^L_\inv$ is independent of\/ $\om,s$. 
\item[{\bf(ii)}] If\/ $\acM_{(\al,d)}^\rst(\ac\mu^{\om,s})\!=\!\acM_{(\al,d)}^\ss(\ac\mu^{\om,s})$ then\/ $[\acM_{(\al,d)}^\ss(\ac\mu^{\om,s})]^L_\inv\!=\![\acM_{(\al,d)}^\ss(\ac\mu^{\om,s})]_\virt,$ where $[\acM_{(\al,d)}^\ss(\ac\mu^{\om,s})]_\virt$ is the virtual class of\/ $\acM_{(\al,d)}^\ss(\ac\mu^{\om,s}),$ which is a proper algebraic space with a perfect obstruction theory. 

Here if\/ $\rank\al=1$ and\/ $d=1,$ when\/ $o_{(\al,d)}=0,$ we use the \begin{bfseries}non-reduced\end{bfseries} obstruction theory defined in {\rm\S\ref{co832},} and otherwise, when\/ $o_{(\al,d)}=p_g,$ we use the \begin{bfseries}reduced\end{bfseries} obstruction theory defined in {\rm\S\ref{co833}}.

If\/ $\rank\al>0$ and\/ $d=1$ then $\acM_{(\al,1)}^\rst(\ac\mu^{\om,s})\!=\!\acM_{(\al,1)}^\ss(\ac\mu^{\om,s})$ for all but finitely many $s\in(0,\iy)$.
\item[{\bf(iii)}] If\/ $L=\O_X(-N)$ for $N\gg 0,$ and\/ $s>0$ is sufficiently small, then $[\acM_{(\al,1)}^\ss(\ac\mu^{\om,s})]^L_\inv=[\acM_{(\al,1)}^\ss(\ac\mu^{\om,s})]_\virt$ coincides with\/ $[\baM^\ss_{(\al,1)}(\bar\tau^0_1)]_\virt$ used in Theorem\/ {\rm\ref{co7thm8}(a)(iii)} for $(\tau,T,\le)\ab=(\mu^\om,M,\le)$ when $\rank\al>1,$ and\/ $[\acM_{(\al,1)}^\ss(\ac\mu^{\om,s})]^L_\inv=[\acM_{(\al,1)}^\ss(\ac\mu^{\om,s})]_\virt=0$ when\/~$\rank\al=1$.
\item[{\bf(iv)}] Let\/ $\al\in C(\coh(X))$ with\/ $\rank\al>0,$ and\/ $s,\ti s>0$. Then
\ea
\!\!\!\!\!\!\!\!&[\acM_{(\al,1)}^\ss(\ac\mu^{\om,\ti s})]^L_\inv=[\acM_{(\al,1)}^\ss(\ac\mu^{\om,s})]^L_\inv
\label{co8eq65}\\
\!\!\!\!\!\!\!\!&+\sum_{\begin{subarray}{l}\be,\ga\in C(\coh(X)):\; \al=\be+\ga, \\
\rank\be=1,\; \rank\ga>0,\\ 
\acM_{(\be,1)}^\ss(\ac\mu^{\om,s})\ne\es\ne \acM_{(\ga,0)}^\ss(\ac\mu^{\om,s})\end{subarray}}\begin{aligned}[t] & U((\be,1),(\ga,0);\ac\mu^{\om,s},\ac\mu^{\om,\ti s})\cdot{} \\
&\;\bigl[[\acM_{(\be,1)}^\ss(\ac\mu^{\om,s})]^L_\inv,
[\acM_{(\ga,0)}^\ss(\ac\mu^{\om,s})]^L_\inv\bigr]
\end{aligned}
\nonumber
\ea
in the Lie algebra $\check H_{\rm even}(\acM^\pl)$ from {\rm\S\ref{co43}}. Here $U(-;\ac\mu^{\om,s},\ac\mu^{\om,\ti s})$ is as in\/ {\rm\eq{co3eq3},} and there are only finitely many nonzero terms in\/~\eq{co8eq65}.
\end{itemize}

If an algebraic $\C$-group $G$ acts on $X$ and Definition\/ {\rm\ref{co8def7}(i)--(iii)} hold then the above generalizes to equivariant homology $H_*^G(\acM^\pl)$.
\end{thm}

\begin{rem}
\label{co8rem6}
{\bf(a)} Parts of Theorems \ref{co8thm4} and \ref{co8thm5} were inspired by related results of Mochizuki \cite{Moch}, including the idea that if $p_g>0$ one must use the non-reduced obstruction theory when $\rank\al=1$ and $d=1$, and the reduced obstruction theory otherwise. Equation \eq{co8eq65} is analogous to~\cite[Th.~7.5.1]{Moch}.
\smallskip

\noindent{\bf(b)} Using Definition \ref{co3def4} and \eq{co8eq9} we can show that in \eq{co8eq65} we have
\e
U((\be,1),(\ga,0);\ac\mu^{\om,s},\ac\mu^{\om,\ti s})\!=\!\begin{cases}
1, & \ti s<\rank\be(\bar\mu^\om(\ga)-\bar\mu^\om(\be))<s,\\
\ha, & {\begin{subarray}{l} \ts \ti s=\rank\be(\bar\mu^\om(\ga)-\bar\mu^\om(\be))<s\;\>\text{or} \\
\ts \ti s<\rank\be(\bar\mu^\om(\ga)-\bar\mu^\om(\be))=s,\end{subarray}} \\
-1\! & s<\rank\be(\bar\mu^\om(\ga)-\bar\mu^\om(\be))<\ti s, \\
-\ha,\! & {\begin{subarray}{l} \ts s=\rank\be(\bar\mu^\om(\ga)-\bar\mu^\om(\be))<\ti s\;\>\text{or} \\
\ts s<\rank\be(\bar\mu^\om(\ga)-\bar\mu^\om(\be))=\ti s,\end{subarray}} \\
0, & \text{otherwise.}
\end{cases}\!\!
\label{co8eq66}
\e
\end{rem}

\subsection{\texorpdfstring{Reducing rank $r$ invariants to rank $1$ using stable pairs}{Reducing rank r invariants to rank 1 using stable pairs}}
\label{co86}

We now explain a method to use the wall crossing formulae \eq{co5eq30} and \eq{co5eq33}--\eq{co5eq34} to compute the invariants $[\M_\al^\ss(\mu^\om)]_\inv$ in $\coh(X)$ from Chapter \ref{co7} and $[\acM_{(\al,d)}^\ss(\ac\mu^{\om,s})]^L_\inv$ above when $\rank\al\ge 1$ and $d\le 1$, in terms of the invariants $[\acM_{(\al,d)}^\ss(\ac\mu^{\om,s})]^{\O_X}_\inv$  when $\rank\al=1$ and $d\le 1$, by induction on~$r=\rank\al$.

\begin{alg}
\label{co8alg}
Work in the situation of \S\ref{co84} or \S\ref{co85}, so that $X$ is a connected projective curve or surface with $p_g=0$ or $p_g>0$. Then for $L\ra X$ a line bundle, $\om\in\Kah(X)$ and $s>0$, $\al\in C(\coh(X))$ with $\rank\al>0$ and $d=0$ or 1, Theorems \ref{co8thm3}, \ref{co8thm4} and \ref{co8thm5} define invariants $[\acM_{(\al,d)}^\ss(\ac\mu^{\om,s})]^L_\inv\in H_*(\acM^\pl_{(\al,d)})$. As in Theorems \ref{co8thm3}(i), \ref{co8thm4}(i) and \ref{co8thm5}(i), when $\rank\al>0$ and $d=0$, these may be identified with the invariants $[\M_\al^\ss(\mu^\om)]_\inv$ in Theorems~\ref{co7thm6}--\ref{co7thm8}.

We regard $\om\in\Kah(X)$ as fixed throughout. We will explain how to compute the $[\acM_{(\al,d)}^\ss(\ac\mu^{\om,s})]^L_\inv$ and $[\M_\al^\ss(\mu^\om)]_\inv$ for all $L,s$ and $\rank\al>0$ in terms of the subset of invariants $[\acM_{(\al,d)}^\ss(\ac\mu^{\om,s})]^{\O_X}_\inv$ for $L=\O_X$ and $\rank\al=1$ and $d=0$ or 1, in the following steps:
\begin{itemize}
\setlength{\itemsep}{0pt}
\setlength{\parsep}{0pt}
\item[{\bf 1.}] Start with the family of invariants $[\acM_{(\al,d)}^\ss(\ac\mu^{\om,s})]^{\O_X}_\inv$ for $L=\O_X$, and $\rank\al=1$ and $d=0$ or 1. Note that by Theorem \ref{co8thm2}(e), in this case $\acM_{(\al,d)}^\rst(\ac\mu^{\om,s})=\acM_{(\al,d)}^\ss(\ac\mu^{\om,s})$, so $[\acM_{(\al,d)}^\ss(\ac\mu^{\om,s})]^{\O_X}_\inv=[\acM_{(\al,d)}^\ss(\ac\mu^{\om,s})]_\virt$, and $\acM_{(\al,d)}^\ss(\ac\mu^{\om,s})$ and hence $[\acM_{(\al,d)}^\ss(\ac\mu^{\om,s})]_\virt$ are independent of $\om,s$.

\item[{\bf 2.}] Let $L\ra X$ be a line bundle. Then mapping $(E,V,\rho)\mapsto (E\ot L,V,\rho)$ gives an equivalence of categories $\acA_{\O_X}\ra\acA_L$ between the categories $\acA$ defined using the line bundles $\O_X\ra X$ and $L\ra X$. This induces an isomorphism $\la_L:\acM^\pl_{\O_X}\ra\acM^\pl_L$ between the projective linear moduli stacks of objects in $\acA_{\O_X},\acA_L$. Theorem \ref{co8thm2}(e) implies that $\la_L$ identifies $\acM_{(\al,d)}^\ss(\ac\mu^{\om,s})_{\O_X}\ra \acM_{(\al\ot\lb L\rb,d)}^\ss(\ac\mu^{\om,s})_L$ when $\rank\al=1$ and $d\le 1$. Hence $(\la_L)_*:H_*(\acM^\pl_{\O_X})\ra H_*(\acM^\pl_L)$ maps $[\acM_{(\al,d)}^\ss(\ac\mu^{\om,s})]^{\O_X}_\inv\mapsto [\acM_{(\al\ot\lb L\rb,d)}^\ss(\ac\mu^{\om,s})]^L_\inv$ in this case.

Thus, we can determine the invariants $[\acM_{(\al,d)}^\ss(\ac\mu^{\om,s})]^L_\inv$ for $\rank\al=1$, $d=0$ or 1 and any $L\ra X$ in terms of the invariants in Step 1.
\end{itemize}

Steps 3, 4 explain how to compute the $[\M_\al^\ss(\mu^\om)]_\inv$ and $[\acM_{(\al,d)}^\ss(\ac\mu^{\om,s})]^L_\inv$ for any $L\ra X$, $d=0$ or 1 and $s>0$ by induction on $r=\rank\al=1,2,\ldots.$ 
\begin{itemize}
\setlength{\itemsep}{0pt}
\setlength{\parsep}{0pt}
\item[{\bf 3.}] Let $r\ge 1$, and suppose by induction that we know how to compute the invariants $[\acM_{(\al,d)}^\ss(\ac\mu^{\om,s})]^L_\inv$ for all $\al\in C(\coh(X))$ with $0<\rank\al\le r$ and $d=0$ or 1 and all $L\ra X$ and $s>0$, in terms of the subset of invariants in Step 1. The first step $r=1$ is given in Step~2.

Let $\al\in C(\coh(X))$ with $\rank\al=r+1$, and $L\ra X$ be a line bundle, and $s>0$. Take $\ti s=3$, so that Theorem \ref{co8thm2}(f) implies that $\acM_{(\al,1)}^\ss(\ac\mu^{\om,3})=\es$, giving $[\acM_{(\al,1)}^\ss(\ac\mu^{\om,3})]^L_\inv=0$. Thus in equations \eq{co8eq62}, \eq{co8eq63} and \eq{co8eq65} with $\ti s=3$, the left hand side is zero, and the $n=1$ term on the right hand side is $[\acM_{(\al,1)}^\ss(\ac\mu^{\om,s})]^L_\inv$. 

When $X$ is a curve or a surface with $p_g=0$, rearranging \eq{co8eq62} or \eq{co8eq63} with $\ti s=3$ yields
\ea
\!\!\!\!\!\!\!\!&[\acM_{(\al,1)}^\ss(\ac\mu^{\om,s})]^L_\inv=
\label{co8eq67}\\
\!\!\!\!\!\!\!\!&-\!\!\!\!\!\sum_{\begin{subarray}{l}n\ge 2,\;\al_1,\ldots,\al_n\in C(\coh(X)), \\
d_1,\ldots,d_n\in\{0,1\}:\; \rank\al_i>0, \\
\acM_{(\al_i,d_i)}^\ss(\ac\mu^{\om,s})\ne\es,\; \text{all\/ $i,$} \\
\al_1+\cdots+\al_n=\al, \; d_1+\cdots+d_n=1 
\end{subarray}} \begin{aligned}[t]
&\ti U((\al_1,d_1),\ldots,(\al_n,d_n);\ac\mu^{\om,s},\ac\mu^{\om,3})\cdot {} \\
&\bigl[\bigl[\cdots\bigl[[\acM_{(\al_1,d_1)}^\ss(\ac\mu^{\om,s})]^L_\inv,\\
&
[\acM_{(\al_2,d_2)}^\ss(\ac\mu^{\om,s})]^L_\inv\bigr],\ldots\bigr],[\acM_{(\al_n,d_n)}^\ss(\ac\mu^{\om,s})]^L_\inv\bigr].
\end{aligned}
\nonumber
\ea
When $X$ is a surface with $p_g>0$, rearranging \eq{co8eq65} with $\ti s=3$ yields
\ea
\!\!\!\!\!\!\!\!&[\acM_{(\al,1)}^\ss(\ac\mu^{\om,s})]^L_\inv=
\label{co8eq68}\\
\!\!\!\!\!\!\!\!&-\!\!\!\!\!\sum_{\begin{subarray}{l}\be,\ga\in C(\coh(X)):\; \al=\be+\ga, \\
\rank\be=1,\; \rank\ga>0,\\ 
\acM_{(\be,1)}^\ss(\ac\mu^{\om,s})\ne\es\ne \acM_{(\ga,0)}^\ss(\ac\mu^{\om,s})\end{subarray}}\begin{aligned}[t] & U((\be,1),(\ga,0);\ac\mu^{\om,s},\ac\mu^{\om,3})\cdot{} \\
&\;\bigl[[\acM_{(\be,1)}^\ss(\ac\mu^{\om,s})]^L_\inv,
[\acM_{(\ga,0)}^\ss(\ac\mu^{\om,s})]^L_\inv\bigr].
\end{aligned}
\nonumber
\ea
Here $U((\be,1),(\ga,0);\ac\mu^{\om,s},\ac\mu^{\om,3})$ is given by \eq{co8eq66}. The right hand sides of \eq{co8eq67}--\eq{co8eq68} involve only $\al_i,\be,\ga$ with rank at most $r$, and so involve only invariants already known, by our inductive hypothesis. Thus $[\acM_{(\al,1)}^\ss(\ac\mu^{\om,s})]^L_\inv$ may be computed in terms of the invariants in Step~1.

\item[{\bf 4.}] Continuing Step 3, let $\al\in C(\coh(X))$ with $\rank\al=r+1$. Choose an ample line bundle $\O_X(1)\ra X$, and choose $N\gg 0$ such that $E$ is $N$-regular for $\O_X(1)$ for all $[E]\in\M_\al^\ss(\mu^\om)$. Then Theorems \ref{co7thm6}(iii), \ref{co7thm7}(a)(ii) and \ref{co7thm8}(a)(iii) involve $[\baM^\ss_{(\al,1)}(\bar\tau^0_1)]_\virt$ for $(\tau,T,\le)=(\mu^\om,M,\le)$. Noting that $\rank\al>1$, Theorems \ref{co8thm3}(iii), \ref{co8thm4}(iii) and \ref{co8thm5}(iii) identify $[\baM^\ss_{(\al,1)}(\bar\tau^0_1)]_\virt\cong[\acM_{(\al,1)}^\ss(\ac\mu^{\om,s})]^L_\inv$ when $L=\O_X(-N)$ and $s>0$ is sufficiently small, say $s=\ep/2$ for $\ep$ as in Theorem \ref{co8thm2}(d). When $X$ is a curve or a surface with $p_g=0$, rearranging \eq{co7eq80} or \eq{co7eq81} yields 
\ea
&[\M_\al^\ss(\mu^\om)]_\inv=P_\al^{c_1(\O_X(1))}(N)^{-1}\cdot{}
\label{co8eq69}\\
&\Bigl\{(\Pi_{\M^\ss_\al(\mu^\om)})_*\bigl([\acM_{(\al,1)}^\ss(\ac\mu^{\om,\ep/2})]^{\O_X(-N)}_\inv\cap c_\top(\bT_{\acM_{(\al,1)}^\ss(\ac\mu^{\om,\ep/2})/\M^\ss_{\al}(\mu^\om)})\bigr)
\nonumber\\
&+\sum_{\begin{subarray}{l}n\ge 2,\;\al_1,\ldots,\al_n\in
C(\coh(X)):\\ \al_1+\cdots+\al_n=\al,\\ 
\mu^\om(\al_i)=\mu^\om(\al), \; \M_{\al_i}^\ss(\mu^\om)\ne\es,\; \text{all\/ $i$}\end{subarray}} \!\!\!\!\!\!\!\!\!\!\!\!\!\!\!\!\!\!\!\!\!\!\!\!\begin{aligned}[t]
&\frac{(-1)^nP_{\al_1}^{c_1(\O_X(1))}(N)}{n!}\cdot\bigl[\bigl[\cdots\bigl[[\M_{\al_1}^\ss(\mu^\om)]_\inv,\\
&\qquad\qquad\;
[\M_{\al_2}^\ss(\mu^\om)]_\inv\bigr],\ldots\bigr],[\M_{\al_n}^\ss(\mu^\om)]_\inv\bigr]\Bigr\}.
\end{aligned}
\nonumber
\ea
When $X$ is a surface with $p_g>0$, rearranging \eq{co7eq84} yields 
\ea
&[\M_\al^\ss(\mu^\om)]_\inv=P_\al^{c_1(\O_X(1))}(N)^{-1}\cdot{}
\label{co8eq70}\\
&(\Pi_{\M^\ss_\al(\mu^\om)})_*\bigl([\acM_{(\al,1)}^\ss(\ac\mu^{\om,\ep/2})]^{\O_X(-N)}_\inv\cap c_\top(\bT_{\acM_{(\al,1)}^\ss(\ac\mu^{\om,\ep/2})/\M^\ss_{\al}(\mu^\om)})\bigr).
\nonumber
\ea
By our inductive hypothesis and Step 3, the right hand sides of \eq{co8eq69}--\eq{co8eq70} involve only invariants already known. Thus $[\M_\al^\ss(\mu^\om)]_\inv$ may be computed in terms of the invariants in Step 1. By Theorems \ref{co8thm3}(i), \ref{co8thm4}(i) and \ref{co8thm5}(i), $[\M_\al^\ss(\mu^\om)]_\inv$ determines $[\acM_{(\al,0)}^\ss(\ac\mu^{\om,s})]^L_\inv$ for any $L\ra X$ and $s>0$. This completes the inductive step.
\end{itemize}
\end{alg}

We have proved:

\begin{thm}
\label{co8thm6}
Let\/ $X$ be a connected projective complex curve or surface, and fix\/ $\om\in\Kah(X)$. Then the invariants  $[\M_\al^\ss(\mu^\om)]_\inv$ in Theorems\/ {\rm\ref{co7thm6}--\ref{co7thm8}} for all\/ $\al\in C(\coh(X))$ with\/ $\rank\al>0,$ and the invariants\/ $[\acM_{(\al,d)}^\ss(\ac\mu^{\om,s})]^L_\inv$ in Theorems\/ {\rm\ref{co8thm3}, \ref{co8thm4}} and\/ {\rm\ref{co8thm5}} for all\/ $\al\in C(\coh(X))$ with\/ $\rank\al>0,$ $d=0$ or\/ $1,$ line bundles $L\ra X,$ and\/ $s>0,$ may be computed in terms of the invariants\/ $[\acM_{(\al,d)}^\ss(\ac\mu^{\om,s})]^{\O_X}_\inv$ in Theorems\/ {\rm\ref{co8thm3}, \ref{co8thm4}} and\/ {\rm\ref{co8thm5}} for all\/ $\al\in C(\coh(X))$ with\/ $\rank\al=1,$ $d=0$ or\/ $1,$ and\/ $s>0$. 

The latter invariants are Behrend--Fantechi virtual classes $[\acM_{(\al,d)}^\ss(\ac\mu^{\om,s})]_\virt,$ and are independent of\/ $\om\in\Kah(X)$ and\/~$s>0$. 
\end{thm}

\begin{rem}
\label{co8rem7}
{\bf(Relation to Donaldson and Seiberg--Witten theory.)} Let $X$ be a compact oriented 4-manifold with $b^2_+(X)\ge 1$. Then one can define the Seiberg--Witten invariants of $X$, as in Morgan \cite{Morg}, and the Donaldson invariants of $X$ for a suitable compact Lie group $G$, as in Donaldson--Kronheimer \cite[\S 9]{DoKr}. Most of the literature takes $G=\SU(2)$ or $\SO(3)$, but Kronheimer \cite{Kron} defines the invariants for $G=\SU(r)$ and $G=\mathop{\rm PSU}(r)$ for~$r\ge 2$.

It was conjectured by Witten \cite{Witt} and Moore--Witten \cite{MoWi} for $G=\SU(2)$, and by Mari\~no--Moore \cite{MaMo} for $G=\SU(N)$, that Donaldson invariants can be written in terms of Seiberg--Witten and algebro-topological invariants of $X$. Witten's conjecture is proved in many cases, see Kronheimer--Mrowka \cite{KrMr}, G\"{o}ttsche--Nakajima--Yoshioka \cite{GNY3}, and Feehan--Leness \cite{FeLe1} for results on this.

The usual method for relating Donaldson and Seiberg--Witten invariants is to introduce an auxiliary moduli space $\M_{\rm mono}$ of solutions of a monopole-type equation on $X$, with a $\U(1)$-action, such that the $\U(1)$-fixed points $\M_{\rm mono}^{\U(1)}$ can be related to either instantons on $X$ (for Donaldson invariants), or to Seiberg--Witten solutions on $X$. See Feehan--Leness \cite{FeLe2} for this programme for $G=\mathop{\rm PSU}(2)$, and Zentner \cite{Zent} for~$G=\mathop{\rm PSU}(r)$. 

Following Mochizuki \cite{Moch}, our invariants $[\M_\al^\ss(\mu^\om)]_\inv$ in Theorems \ref{co7thm7}--\ref{co7thm8} for $\rank\al=r>0$ should be understood as algebraic versions of Donaldson invariants of a projective surface $X$ for $G=\U(r)$. (If $\pi_1(X)=0$ then $\U(r)$- and $\SU(r)$-Donaldson invariants are nearly the same thing.)

Following D\"urr--Kabanov--Okonek \cite{DKO} and Mochizuki \cite[\S 6.3]{Moch}, our invariants $[\acM_{(\al,1)}^\ss(\ac\mu^{\om,s})]^{\O_X}_\inv$ in Theorems \ref{co8thm4}--\ref{co8thm5} for $\rank\al=1$ and $c_2(\al)=0$ should be understood as algebraic versions of Seiberg--Witten invariants of $X$. (Here $\rank\al=1$ and $c_2(\al)=0$ imply that if $[E,V,\rho]\in \acM_{(\al,1)}^\ss(\ac\mu^{\om,s})$ then $E$ is a line bundle.) The other invariants $[\acM_{(\al,d)}^\ss(\ac\mu^{\om,s})]^{\O_X}_\inv$ for $\rank\al=1$ and $d=0,1$ depend on the Seiberg--Witten invariants, $\pi_1(X)$, and the fundamental classes of Hilbert schemes $\Hilb^n(X)$, which can be written in terms of $H^*(X)$.

Thus, Algorithm \ref{co8alg} and Theorem \ref{co8thm6} prove an algebraic version of the conjectures of \cite{MaMo,MoWi,Witt}. Note that Mochizuki \cite{Moch} in effect proves his own version of Theorem \ref{co8thm6}, but with fewer explicit formulae.

A new feature of our work is the use of the Lie brackets $[\,,\,]$ on $H_*(\M^\pl)$, $H_*(\acM^\pl)$ from \S\ref{co43} to write the wall-crossing formulae \eq{co8eq67}--\eq{co8eq69}. This may be useful in Donaldson theory of general 4-manifolds, not just surfaces.

In the sequel \cite{Joyc13} the author will use Algorithm \ref{co8alg} to compute the invariants $[\M_\al^\ss(\mu^\om)]_\inv$ fairly explicitly for curves and surfaces.
\end{rem}

\begin{rem}
\label{co8rem8}
{\bf(Relation to Feyzbakhsh--Thomas \cite{FeTh}.)} Let $X$ be a Calabi--Yau 3-fold satisfying some conditions. As in \S\ref{co781}(E), Feyzbakhsh--Thomas \cite{FeTh} use the wall-crossing theory of Joyce--Song \cite{JoSo} to prove that the rank $r\ge 1$ Donaldson--Thomas invariants of $X$ can be written in terms of the rank 1 invariants. Their method is similar in outline to Algorithm \ref{co8alg}, though the technical details are different. 

It seems likely that Feyzbakhsh--Thomas' programme can be carried out for Fano 3-folds, showing that the invariants $[\M_\al^\ss(\tau)]_\inv$ in Theorem \ref{co7thm9} for $\rank\al>0$ can be determined from those for $\rank\al=1$.
\end{rem}

\section{Proof of Theorem \ref{co5thm1}}
\label{co9}

Suppose throughout this section that $\B\subseteq\A,\M,\M^\pl,K(\A),\ldots$ are a fixed set of data satisfying Assumptions \ref{co5ass1}--\ref{co5ass2}.

\subsection{\texorpdfstring{Defining invariants $[\M_\al^\ss(\tau)]^k_\inv$}{Defining invariants [ℳₐˢˢ(τ)]ᵏᵢᵥ}}
\label{co91}

In the next lemma, $\al_1,\ldots,\al_n$ satisfy the conditions in equations \eq{co5eq30}--\eq{co5eq32} of Theorems \ref{co5thm1}--\ref{co5thm2}. It implies there are only finitely many terms in \eq{co5eq30}--\eq{co5eq32}. Also, when we prove some statement by induction on increasing $\rk\al$ from Assumption \ref{co5ass2}(f), as $\rk\al_i<\rk\al$, the statement holds for each~$\al_i$.

\begin{lem}
\label{co9lem1}
Suppose\/ $\al=\al_1+\cdots+\al_n$ for $n\ge 2$ with\/ $\al_i\in C(\B)_\pe,$ $\tau(\al_i)=\tau(\al),$ and\/ $\M_{\al_i}^\ss(\tau)\ne\es$ for $i=1,\ldots,n$. Then $\rk\al=\rk\al_1+\cdots+\rk\al_n$ with\/ $\rk\al_i>0$ for $i=1\ldots,n$. Hence $n\le\rk\al$ and\/ $\rk\al_i<\rk\al$ for all\/ $i$. There are only finitely many sets of data $n,\al_1,\ldots,\al_n$ satisfying these conditions.
\end{lem}

\begin{proof} Since $\tau(\al_i)=\tau(\al)$, we have $\tau(\al_1+\cdots+\al_{n-1})=\tau(\al_n)$, so Assumption \ref{co5ass2}(f) gives $\rk(\al)=\rk(\al_1+\cdots+\al_{n-1})+\rk(\al_n)$. Thus $\rk\al=\rk\al_1+\cdots+\rk\al_n$ by induction on $n$. As $n\ge 2$ and $\rk\al_i>0$, we see that $n\le\rk\al$ and~$\rk\al_i<\rk\al$.

Each $\C$-point $[E]$ in $\M_\al^\ss(\tau)$ has a splitting $E=F_1\op\cdots\op F_m$ into indecomposable factors in $\B$, which have classes $\lb F_1\rb,\ldots,\lb F_m\rb$ in $C(\B)_\pe$ by Assumption \ref{co5ass2}(c). These classes are independent of choices up to order, so $(\lb F_1\rb,\ldots,\lb F_m\rb)S_m$ in $\mathop{\rm Sym}^mC(\B)_\pe:=C(\B)_\pe^m/S_m$ depends only on $[E]$. 

Consider the map 
\e
\M_\al^\ss(\tau)(\C)\longra\ts\coprod_{m\ge 1}\mathop{\rm Sym}^mC(\B)_\pe,\quad [E]\longmapsto(\lb F_1\rb,\ldots,\lb F_m\rb)S_m.
\label{co9eq1}
\e
This has an upper semicontinuity property on $\M_\al^\ss(\tau)(\C)$ in the Zariski topology, as under limits in $\M_\al^\ss(\tau)(\C)$ splittings can only refine by decomposing $\lb E_i\rb$ into more factors. But $\M_\al^\ss(\tau)(\C)$ is compact in the Zariski topology since $\M_\al^\ss(\tau)$ is finite type by Assumption \ref{co5ass2}(b), so \eq{co9eq1} can take only finitely many values. 

Let $n,\al_1,\ldots,n$ be as in the lemma. Choose $\C$-points $[E_i]$ in $\M_{\al_i}^\ss(\tau)$ for each $i$ and set $E=E_1\op\cdots\op E_n$. Write $E=F_1\op\cdots\op F_m$ for the decomposition into indecomposables. Then each $E_i$ is the direct sum of a subset of the $F_j$, so $\al_i=\lb E_i\rb$ is a sum of a subset of the $\lb F_j\rb$. As there are only finitely many possibilities for $\lb F_1\rb,\ldots,\lb F_m\rb$, there are only finitely many sets~$\al_1,\ldots,\al_n$.
\end{proof}

\begin{dfn}
\label{co9def1}
Let $(\tau,T,\le)\in\sS$, $\al\in C(\B)_\pe$, and $k\in K$ with $\M_\al^\ss(\tau)\subseteq\M_{k,\al}^\pl$. Using the notation of Example \ref{co5ex1}, in particular equation \eq{co5eq29}, as $\baM^\rst_{(\al,1)}(\bar\tau^0_1)=\baM^\ss_{(\al,1)}(\bar\tau^0_1)$, define $\bar\Up_\al^k(\tau)\in\check H_{2o_\al}(\M_\al^\pl)$ by
\e
\begin{split}
&\bar\Up_\al^k(\tau)=\bar\Up\bigl([\baM_{(\al,1)}^\ss(\bar\tau^0_1)]_\virt\bigr)\\
&=(\Pi_{\M^\ss_\al(\tau)})_*\bigl([\baM^\ss_{(\al,1)}(\bar\tau^0_1)]_\virt\cap c_\top(\bT_{\baM^\ss_{(\al,1)}(\bar\tau^0_1)/\M^\pl_\al})\bigr)\\
&=(-1)^{\la_k(\al)-1}(\Pi_{\M^\ss_\al(\tau)})_*\bigl([\baM^\ss_{(\al,1)}(\bar\tau^0_1)]_\virt\cap c_\top(\bL_{\baM^\ss_{(\al,1)}(\bar\tau^0_1)/\M^\pl_\al})\bigr),
\end{split}
\label{co9eq2}
\e
where the sign is from $c_i(\cV^*)\!=\!(-1)^ic_i(\cV)$ and $\rank\bT_{\baM^\ss_{(\al,1)}(\bar\tau^0_1)/\M^\pl_\al}\!=\!\la_k(\al)\!-\!1$.
\end{dfn}

The next proposition is essentially the same as Theorem \ref{co5thm1}, except that the invariants $[\M_\al^\ss(\tau)]^k_\inv$ are allowed to depend on~$k\in K$.

\begin{prop}
\label{co9prop1}
For all\/ $(\tau,T,\le)\in\sS,$ $\al\in C(\B)_\pe,$ and\/ $k\in K$ with\/ $\M_\al^\ss(\tau)\subseteq\M_{k,\al}^\pl,$ there exist unique $[\M_\al^\ss(\tau)]^k_\inv\in\check H_{2o_\al}(\M_\al^\pl)$ satisfying:
\begin{itemize}
\setlength{\itemsep}{0pt}
\setlength{\parsep}{0pt}
\item[{\bf(i)}] If\/ $\M_\al^\rst(\tau)=\M_\al^\ss(\tau)$ then\/ $[\M_\al^\ss(\tau)]^k_\inv=[\M_\al^\ss(\tau)]_\virt,$ where $[\M_\al^\ss(\tau)]_\virt$ is the Behrend--Fantechi virtual class as in Assumption\/ {\rm\ref{co5ass2}(g)}.
\item[{\bf(ii)}] Suppose $(\tau,T,\le),(\ti\tau,\ti T,\le)\in\sS,$ $\al\in C(\B)_\pe$ and\/ $k\in K$ with\/ $\M_\al^\ss(\tau)=\M_\al^\ss(\ti\tau)\subseteq\M_{k,\al}^\pl$. Then $[\M_\al^\ss(\tau)]^k_\inv=[\M_\al^\ss(\ti\tau)]^k_\inv$.
\item[{\bf(iii)}] For all\/ $(\tau,T,\le)\in\sS,$ $\al\in C(\B)_\pe$ and\/ $k\in K$ with\/ $\M_\al^\ss(\tau)\subseteq\M_{k,\al}^\pl,$
\end{itemize}
\e
\bar\Up_\al^k(\tau)=\sum_{\begin{subarray}{l}n\ge 1,\;\al_1,\ldots,\al_n\in
C(\B)_\pe:\\  
\tau(\al_i)=\tau(\al),\; \M_{\al_i}^\ss(\tau)\ne\es,\; \text{all $i,$} \\
\al_1+\cdots+\al_n=\al,\; o_{\al_1}+\cdots+o_{\al_n}=o_\al \end{subarray}}  \!\!\!\!\!\!\!\!\!\!\!\!\!\!\!\!\!\!\!\!\!\begin{aligned}[t]
\frac{(-1)^{n+1}\la_k(\al_1)}{n!}\,\cdot&\bigl[\bigl[\cdots\bigl[[\M_{\al_1}^\ss(\tau)]^k_\inv,\\
&\; 
[\M_{\al_2}^\ss(\tau)]^k_\inv\bigr],\ldots\bigr],\\
&\;[\M_{\al_n}^\ss(\tau)]^k_\inv\bigr].
\end{aligned}
\label{co9eq3}
\e
\begin{itemize}
\setlength{\itemsep}{0pt}
\setlength{\parsep}{0pt}
\item[] Here $\bar\Up_\al^k(\tau)$ is as in {\rm\eq{co9eq2},} the\/ $[\M_{\al_i}^\ss(\tau)]^k_\inv$ are defined in \eq{co9eq3} for all choices of\/ $\al_1,\ldots,\al_n,$ and there are only finitely many terms in the sum.
\end{itemize}
\end{prop}

\begin{proof} Lemma \ref{co9lem1} implies that the sum in \eq{co9eq3} has only finitely many terms. We will define unique $[\M_\al^\ss(\tau)]^k_\inv$ for all $(\tau,T,\le)\in\sS,$ $\al\in C(\B)_\pe$ and $k\in K$ with $\M_\al^\ss(\tau)\subseteq\M_{k,\al}^\pl$ satisfying \eq{co9eq3} by induction on increasing $\rk\al=1,2,\ldots$ in Assumption \ref{co5ass2}(f). If $\M_\al^\ss(\tau)=\es$ we define $[\M_\al^\ss(\tau)]^k_\inv=0$. This is required by (i) as $[\M_\al^\ss(\tau)]_\virt=0$, and is consistent with \eq{co9eq3} since $\baM^\ss_{(\al,1)}(\bar\tau^0_1)=\es$ so both sides of \eq{co9eq3} are zero.

Suppose by induction that $l\ge 0$ is given and there exist unique $[\M_\al^\ss(\tau)]^k_\inv$ as above for all $(\tau,T,\le),\al,k$ with $\rk\al\le l$. (In the first step $l=0$ this assumption is vacuous as $\rk\al>0$.) Suppose $(\tau,T,\le)\in\sS,$ $\al\in C(\B)_\pe$ and $k\in K$ with $\es\ne\M_\al^\ss(\tau)\subseteq\M_{k,\al}^\pl$ and $\rk\al=l+1$. We may rewrite \eq{co9eq3} as
\ea
&[\M_\al^\ss(\tau)]^k_\inv=\la_k(\al)^{-1}\cdot\bar\Up_\al^k(\tau)
\label{co9eq4}\\
&+\sum_{\begin{subarray}{l}n\ge 2,\;\al_1,\ldots,\al_n\in
C(\B)_\pe:\\  
\tau(\al_i)=\tau(\al),\; \M_{\al_i}^\ss(\tau)\ne\es,\; \text{all $i,$} \\
\al_1+\cdots+\al_n=\al,\; o_{\al_1}+\cdots+o_{\al_n}=o_\al \end{subarray}} \!\!\!\!\!\!\!\!\!\!\!\!\!\!\!\begin{aligned}[t]
&\frac{(-1)^n\,\la_k(\al_1)}{n!\,\la_k(\al)}\cdot\bigl[\bigl[\cdots\bigl[[\M_{\al_1}^\ss(\tau)]^k_\inv,\\
&\qquad\qquad\qquad\qquad 
[\M_{\al_2}^\ss(\tau)]^k_\inv\bigr],\ldots\bigr],[\M_{\al_n}^\ss(\tau)]^k_\inv\bigr],
\end{aligned}
\nonumber
\ea	
where $[\M_\al^\ss(\tau)]^k_\inv$ comes from the term $n=1$, $\al_1=\al$ in the sum in \eq{co9eq3}, and $\la_k(\al)>0$ by Assumption \ref{co5ass1}(g)(iv), so $\la_k(\al)^{-1}$ is well defined. 

Lemma \ref{co9lem1} says that for all $n\ge 2$, $\al_1,\ldots,\al_n$ as in the sum in \eq{co9eq4} we have $\rk\al_i<\rk\al=l+1$, so $\rk\al_i\le l$, and $[\M_{\al_i}^\ss(\tau)]^k_\inv$ is already defined for $i=1,\ldots,n$ by induction. Thus every term on the right hand side of \eq{co9eq4} is defined. There are only finitely many terms in \eq{co9eq4} from above. Hence \eq{co9eq4} defines unique $[\M_\al^\ss(\tau)]^k_\inv$ for which \eq{co9eq3} holds. So by induction there are unique $[\M_\al^\ss(\tau)]^k_\inv$ for all allowable $(\tau,T,\le),\al,k$ which satisfy (iii) and $[\M_\al^\ss(\tau)]^k_\inv=0$ if~$\M_\al^\ss(\tau)=\es$.

It remains to prove that these $[\M_\al^\ss(\tau)]^k_\inv$ satisfy (i),(ii). For (i), suppose $\M_\al^\rst(\tau)=\M_\al^\ss(\tau)$. Then in \eq{co9eq3} there are no terms in the sum with $n\ge 2$, as choosing $\C$-points $[E_i]$ in $\M_{\al_i}^\ss(\tau)\ne\es$ for $i=1,\ldots,n$ would give a $\C$-point $[E_1\op\cdots\op E_n]$  in $\M_\al^\ss(\tau)\sm \M_\al^\rst(\tau)$. Thus \eq{co9eq2}--\eq{co9eq3} reduce to
\ea
&(\bar\Pi_{\M^\ss_{\al}(\tau)})_*([\baM^\ss_{(\al,1)}(\bar\tau^0_1)]_\virt\!\cap\! c_\top(\bT_{\baM^\ss_{(\al,1)}(\bar\tau^0_1)/\M^\ss_{\al}(\tau)}))=\bar\Up_\al^k(\tau)
\nonumber\\
&\quad =\la_k(\al)\cdot [\M_\al^\ss(\tau)]^k_\inv.
\label{co9eq5}
\ea

In this case, the smooth morphism $\Pi_{\M_\al^\ss(\tau)}:\baM^\ss_{(\al,1)}(\bar\tau^0_1)\ra\M_\al^\ss(\tau)$ is a fibre bundle with fibre $\CP^{\la_k(\al)-1}$, the pushdown to $\M_\al^\ss(\tau)$ of $\bP(\cV_{k,\al})\ra\M_{k,\al}$, where $\cV_{k,\al}\ra\M_{k,\al}$ is the rank $\la_k(\al)$ vector bundle from Assumption \ref{co5ass1}(g)(ii). This is because when $\M_\al^\rst(\tau)=\M_\al^\ss(\tau)$, the conditions in Example \ref{co5ex1}(i)--(iii) for $(E,V,\rho)$ in class $(\al,1)$ to be $\bar\tau^0_1$-semistable reduce to $E$ $\tau$-semistable and $\rho\ne 0$. By \S\ref{co52}, $\Pi_{\M_\al^\ss(\tau)}$ is the classical truncation of a smooth morphism $\bs\Pi_{\M_\al^\ss(\tau)}:\bs\baM^\ss_{(\al,1)}(\bar\tau^0_1)\ra\bs\M_\al^\ss(\tau)$ of proper quasi-smooth derived algebraic spaces, which is also a fibre bundle with fibre $\CP^{\la_k(\al)-1}$. Thus Corollary \ref{co2cor1} yields
\e
\begin{split}
&(\bar\Pi_{\M^\ss_{\al}(\tau)})_*\bigl([\baM^\ss_{(\al,1)}(\bar\tau^0_1)]_\virt\cap c_\top(\bT_{\baM^\ss_{(\al,1)}(\bar\tau^0_1)/\M^\ss_{\al}(\tau)})\bigr)\\
&\quad =\la_k(\al)\cdot [\M_\al^\ss(\tau)]_\virt.
\end{split}
\label{co9eq6}
\e
Equations \eq{co9eq5}--\eq{co9eq6} prove part~(i).

For (ii), suppose by induction that $l\ge 0$ is given such that if $(\tau,T,\le),\ab(\ti\tau,\ti T,\le)\in\sS$, $\al\in C(\B)_\pe$ and $k\in K$ with $\M_\al^\ss(\tau)=\M_\al^\ss(\ti\tau)\subseteq\M_{k,\al}^\pl$ and $\rk\al\le l$, then $[\M_\al^\ss(\tau)]^k_\inv=[\M_\al^\ss(\ti\tau)]^k_\inv$. (This is vacuous when $l=0$ as $\rk\al>0$.) Let $(\tau,T,\le),\ab(\ti\tau,\ti T,\le)\in\sS$, $\al\in C(\B)_\pe$ and $k\in K$ with $\M_\al^\ss(\tau)=\M_\al^\ss(\ti\tau)\subseteq\M_{k,\al}^\pl$ and~$\rk\al=l+1$.

Suppose $\al_1,\ldots,\al_n$ are as in the sum in \eq{co9eq3} for $\tau$. If $[E_i]$ is a $\C$-point of $\M_{\al_i}^\ss(\tau)$ for $i=1,\ldots,n$ then $[E_1\op\cdots\op E_n]$ is a $\C$-point of $\M_\al^\ss(\tau)=\M_\al^\ss(\ti\tau)$, so $E_1\op\cdots\op E_n$ is $\ti\tau$-semistable. This implies that $E_1,\ldots,E_n$ are $\ti\tau$-semistable, so $[E_i]\in \M_{\al_i}^\ss(\ti\tau)$ and $\M_{\al_i}^\ss(\tau)\subseteq \M_{\al_i}^\ss(\ti\tau)$, and also $\ti\tau(\al_1)=\cdots=\ti\tau(\al_n)=\ti\tau(\al)$. Hence $n$, $\al_1,\ldots,\al_n$ satisfy the conditions in \eq{co9eq3} for $\ti\tau$, so the same argument gives $\M_{\al_i}^\ss(\ti\tau)\subseteq \M_{\al_i}^\ss(\tau)$, and thus~$\M_{\al_i}^\ss(\tau)=\M_{\al_i}^\ss(\ti\tau)$. 

If $n\ge 2$ then $\rk\al_i<\rk\al=l+1$ by Lemma \ref{co9lem1}, so $\rk\al_i\le l$ and $[\M_{\al_i}^\ss(\tau)]^k_\inv=[\M_{\al_i}^\ss(\ti\tau)]^k_\inv$ for $i=1,\ldots,n$ by induction. Hence, comparing \eq{co9eq3} for $\tau,\al,k$ and $\ti\tau,\al,k$, the terms in the sum over $\al_1,\ldots,\al_n$ for $n\ge 2$ agree. Also $\baM^\ss_{(\al,1)}(\bar\tau^0_1)=\baM^\ss_{(\al,1)}(\bar{\ti\tau}{}^0_1)$ as $\M_\al^\ss(\tau)=\M_\al^\ss(\ti\tau)$, so the left hand sides of \eq{co9eq3} for $\tau,\al,k$ and $\ti\tau,\al,k$ agree. Thus the only remaining terms, with $n=1$ and $\al_1=\al$, agree, showing that $\la_k(\al)\cdot [\M_\al^\ss(\tau)]^k_\inv=\la_k(\al)\cdot [\M_\al^\ss(\ti\tau)]^k_\inv$. Since $\la_k(\al)>0$ by Assumption \ref{co5ass1}(g)(iv), we see that $[\M_\al^\ss(\tau)]^k_\inv=[\M_\al^\ss(\ti\tau)]^k_\inv$. Therefore (ii) holds by induction.
\end{proof}

\subsection{\texorpdfstring{A Lie bracket relation on virtual classes $[\bar{\cal M}^{\rm ss}_{(\al,1)}(\bar\tau^0_1)]_{\rm virt}$}{A Lie bracket relation on virtual classes [ℳˢˢᵦ₁(τ⁰)]ᵥᵢ}}
\label{co92}

The goal of this section is to prove Corollary \ref{co9cor1} below, which gives relations between elements $\bar\Up_{\al_1}^{k_1}(\tau),\bar\Up_{\al_2}^{k_2}(\tau)$ from \eq{co9eq2} in the Lie algebra $\check H_{\rm even}(\M^\pl)$. The next definition sets up ingredients we will need in the proof.

\begin{dfn}
\label{co9def2}
Suppose $(\tau,T,\le)\in\sS$, $\al\in C(\B)_\pe$ and $k_1,k_2\in K$ with $\M_\al^\ss(\tau)\subseteq\M_{k_1,\al}^\pl\cap\M_{k_2,\al}^\pl$, using the notation of Assumption \ref{co5ass1}(g). Define data $\acB\subseteq\acA,K(\acA),\ab\acM,\ab\acM^\pl,\ab\ldots$ as in Definition \ref{co5def1} using the quiver $\ac Q$ with $\dot Q_0=\{v_1,v_2,v_3\}$, $\ddot Q_0=\{w_1,w_2\}$, $\ac Q_1=\{e_1,e_2,e_3,e_4\}$, $t(e_1)=v_1$, $h(e_1)=w_1$, $t(e_2)=v_2$, $h(e_2)=w_2$, $t(e_3)=v_3$, $h(e_3)=w_1$, $t(e_4)=v_3$, $h(e_4)=w_2$, and with $\ka:\ddot Q_0\ra K$ given by $\ka(w_i)=k_i$ for $i=1,2$. We illustrate this by:
\e
\begin{gathered}
\begin{xy}
0;<1mm,0mm>:
,(-15,2.5)*{e_2}
,(-38,12)*{e_3}
,(-38,4)*{e_4}
,(0,3)*{w_2}
,(-56,9.5)*{v_3}
,(12,0)*{\ka(w_2)=k_2.}
,(-30,0)*+{\bu} ; (0,0)*+{\circ} **@{-} ?>*\dir{>}
,(-15,13.5)*{e_1}
,(-30,14)*{v_1}
,(-30,3)*{v_2}
,(0,14)*{w_1}
,(12,11)*{\ka(w_1)=k_1\phantom{.}}
,(-30,11)*+{\bu} ; (0,11)*+{\circ} **@{-} ?>*\dir{>}
,(-56,6)*+{\bu} ; (-30,0)*+{\bu} **@{-} ?>*\dir{>}
,(-56,6)*+{\bu} ; (-30,11)*+{\bu} **@{-} ?>*\dir{>}\end{xy}
\end{gathered}
\label{co9eq7}
\e
Here we use notation $\acA,\acB,\acM,\ldots$ rather than $\baA,\baB,\baM,\ldots$ to distinguish them from Example~\ref{co5ex1}.

For brevity we will write objects of $\acA$ as $(E,\bs V,\bs\rho)$ with $\bs V=(V_1,V_2,V_3)$, $\bs\rho=(\rho_1,\ldots,\rho_4)$ rather than $(V_{v_1},V_{v_2},V_{v_3}),(\rho_{e_1},\ldots,\rho_{e_4})$, and write $K(\acA)=K(\A)\t\Z^3$ rather than $K(\A)\t\Z^{\{v_1,v_2,v_3\}}$, and write stability conditions as $\ac\tau^\la_{\bs\mu}$ with $\bs\mu=(\mu_1,\mu_2,\mu_3)\in\R^3$ rather than $\bs\mu=(\mu_{v_1},\mu_{v_2},\mu_{v_3})$, and so on.

As in Definition \ref{co5def1}, define a weak stability condition $(\ac\tau^0_{\bs\mu},\ac T,\le)$ on $\acA$ with $\la=0$ and $\bs\mu=(\ep,\ep,1)$, where $\ep$ satisfies
\e
0<\ep<\frac{1}{2\rk\al}.
\label{co9eq8}
\e
We are interested in the moduli stacks $\acM_{(\al,\bs 1)}^\rst(\ac\tau^0_{\bs\mu})\subseteq\acM_{(\al,\bs 1)}^\ss(\ac\tau^0_{\bs\mu})\subseteq\acM_{(\al,\bs 1)}^\pl$, where $\bs 1=(1,1,1)$. For $(E,\bs V,\bs\rho)\in\acB$ with $\lb E,\bs V,\bs\rho\rb=(\al,\bs 1)$ to be strictly $\ac\tau^0_{\bs\mu}$-semistable there must exist $0\ne (E',\bs V',\bs\rho')\subsetneq(E,\bs V,\bs\rho)$ with $\ac\tau^0_{\bs\mu}(\lb E',\bs V',\bs\rho'\rb)=\ac\tau^0_{\bs\mu}(\lb E/E',\bs V/\bs V',\bs\rho''\rb)$. This implies that $0\ne E'\subsetneq E$ with $\tau(\al')=\tau(\al-\al')$, where $\al'=\lb E'\rb$, with $1\le\rk\al'<\rk\al$ as $\rk\al,\rk(\al-\al')>0$, and 
\e
\frac{\ep\dim V_1'+\ep\dim V_2'+\dim V_3'}{\rk\al'}=\frac{2\ep+1}{\rk\al}.
\label{co9eq9}
\e
But using \eq{co9eq8} we can show that if $\dim V_3'=1$ then the left hand side of \eq{co9eq9} is greater than the right, and if $\dim V_3'=0$ it is smaller, so \eq{co9eq9} does not hold.

Thus $(E,\bs V,\bs\rho)$ is $\ac\tau^0_{\bs\mu}$-stable, and $\acM_{(\al,\bs 1)}^\rst(\ac\tau^0_{\bs\mu})=\acM_{(\al,\bs 1)}^\ss(\ac\tau^0_{\bs\mu})$. Therefore Assumption \ref{co5ass2}(h) says that $\acM_{(\al,\bs 1)}^\ss(\ac\tau^0_{\bs\mu})$ is a proper algebraic space.

Write $\ac L_{a,b,c}\ra\acM_{(\al,\bs 1)}$ for the line bundle with fibre $V_1^{\ot^a}\ot V_2^{\ot^b}\ot V_3^{\ot^c}$ at $[E,\bs V,\bs\rho]$. Then $\ac L_{a,b,c}$ has weight $a+b+c$ under the $[*/\bG_m]$-action. Hence if $a+b+c=0$ then $\ac L_{a,b,c}$ descends to a line bundle $\ac L^\pl_{a,b,c}\ra\acM_{(\al,\bs 1)}^\pl$ with $(\Pi^\pl)^*(\ac L_{a,b,c}^\pl)\cong \ac L_{a,b,c}$. This restricts to~$\ac L^\pl_{a,b,c}\ra\acM_{(\al,\bs 1)}^\ss(\ac\tau^0_{\bs\mu})$.

Define a $\bG_m$-action on $\acM_{(\al,\bs 1)}^\ss(\ac\tau^0_{\bs\mu})$ by, for $\th\in\bG_m$
\e
\th:[E,\bs V,\bs\rho]\longmapsto[E,\bs V,(\rho_1,\rho_2,\rho_3,\th\,\rho_4)].
\label{co9eq10}
\e
Now $\acM_{(\al,\bs 1)}^\ss(\ac\tau^0_{\bs\mu})$ is the classical truncation of a derived stack $\bs\acM_{(\al,\bs 1)}^\ss(\ac\tau^0_{\bs\mu})$, which we regard as an open substack of $\bs\dM_{(\al,\bs 1)}^\rpl$ defined in \eq{co5eq17}. This $\bs\acM_{(\al,\bs 1)}^\ss(\ac\tau^0_{\bs\mu})$ is quasi-smooth as $\bs\dM_{(\al,\bs 1)}^\rpl$ is, and is a proper derived algebraic space as $\acM_{(\al,\bs 1)}^\ss(\ac\tau^0_{\bs\mu})$ is a proper algebraic space. The $\bG_m$-action on $\acM_{(\al,\bs 1)}^\ss(\ac\tau^0_{\bs\mu})$ lifts to $\bs\acM_{(\al,\bs 1)}^\ss(\ac\tau^0_{\bs\mu})$, as the definition of $\bs\dM_{(\al,\bs 1)}^\rpl$ is evidently $\bG_m$-equivariant. Thus the obstruction theory $\bL_i:i^*(\bL_{\bs\acM_{(\al,\bs 1)}^\ss(\ac\tau^0_{\bs\mu})})\ra \bL_{\acM_{(\al,\bs 1)}^\ss(\ac\tau^0_{\bs\mu})}$ on $\acM_{(\al,\bs 1)}^\ss(\ac\tau^0_{\bs\mu})$, which is the restriction of the second line of \eq{co5eq18} for $\bs\dM_{(\al,\bs 1)}^\rpl$, is $\bG_m$-equivariant.
\end{dfn}

As in Definition \ref{co9def2} we have a $\bG_m$-action and a $\bG_m$-equivariant perfect obstruction theory on $\acM_{(\al,\bs 1)}^\ss(\ac\tau^0_{\bs\mu})$. We will apply the equivariant localization result Corollary \ref{co2cor2} to $\acM_{(\al,\bs 1)}^\ss(\ac\tau^0_{\bs\mu})$. To do this we first describe the $\bG_m$-fixed substack $\acM_{(\al,\bs 1)}^\ss(\ac\tau^0_{\bs\mu})^{\bG_m}$, in the next proposition.

\begin{prop}
\label{co9prop2}
In the situation of Definition\/ {\rm\ref{co9def2},} the $\bG_m$-fixed substack\/ $\acM_{(\al,\bs 1)}^\ss(\ac\tau^0_{\bs\mu})^{\bG_m}$ is the disjoint union of the following pieces:
\begin{itemize}
\setlength{\itemsep}{0pt}
\setlength{\parsep}{0pt}
\item[{\bf(a)}] The substack\/ $\acM_{(\al,\bs 1)}^\ss(\ac\tau^0_{\bs\mu})_{\rho_4=0}$ of all\/ $\C$-points $[E,\bs V,\bs\rho]$ in $\acM_{(\al,\bs 1)}^\ss(\ac\tau^0_{\bs\mu})$ with\/ $\rho_4=0$. All such\/ $[E,\bs V,\bs\rho]$ have $\rho_1,\rho_2,\rho_3\ne 0$. Thus $\ac L^\pl_{1,0,-1}\vert_{\acM_{(\al,\bs 1)}^\ss(\ac\tau^0_{\bs\mu})_{\rho_4=0}}$ is trivial, with nonvanishing section $[E,\bs V,\bs\rho]\mapsto\rho_3$. There is a morphism
\e
\begin{split}
\Pi_{\rho_4=0}&:\acM_{(\al,\bs 1)}^\ss(\ac\tau^0_{\bs\mu})_{\rho_4=0}\longra \baM_{(\al,1)}^\ss(\bar\tau^0_1)_{k_1},\\
\Pi_{\rho_4=0}&:[E,\bs V,\bs\rho]\longmapsto [E,V_1,\rho_1],
\label{co9eq11}
\end{split}
\e
where $\baM_{(\al,1)}^\ss(\bar\tau^0_1)_{k_1}$ is as in Example\/ {\rm\ref{co5ex1}} with\/ $k=k_1$. Also $\Pi_{\rho_4=0}$ is smooth of relative dimension $\la_{k_2}(\al)-1,$ with fibre $\bP(F_{k_2}(E))\cong\CP^{\la_{k_2}(\al)-1}$ over $[E,V_1,\rho_1],$ by mapping $[E,\bs V,\bs\rho]$ to $\rho_2(V_2)$ in\/~$\bP(F_{k_2}(E))$.
\item[{\bf(b)}] The substack\/ $\acM_{(\al,\bs 1)}^\ss(\ac\tau^0_{\bs\mu})_{\rho_3=0}$ of all\/ $\C$-points $[E,\bs V,\bs\rho]$ in $\acM_{(\al,\bs 1)}^\ss(\ac\tau^0_{\bs\mu})$ with\/ $\rho_3=0$.  All such\/ $[E,\bs V,\bs\rho]$ have $\rho_1,\rho_2,\rho_4\ne 0$. Thus $\ac L^\pl_{0,1,-1}\vert_{\acM_{(\al,\bs 1)}^\ss(\ac\tau^0_{\bs\mu})_{\rho_3=0}}$ is trivial, with nonvanishing section $[E,\bs V,\bs\rho]\mapsto\rho_4$. There is a morphism
\e
\begin{split}
\Pi_{\rho_3=0}&:\acM_{(\al,\bs 1)}^\ss(\ac\tau^0_{\bs\mu})_{\rho_3=0}\longra \baM_{(\al,1)}^\ss(\bar\tau^0_1)_{k_2},\\
\Pi_{\rho_3=0}&:[E,\bs V,\bs\rho]\longmapsto [E,V_2,\rho_2].
\label{co9eq12}
\end{split}
\e
Here $\Pi_{\rho_3=0}$ is smooth of dimension $\la_{k_1}(\al)-1,$ with fibre $\CP^{\la_{k_1}(\al)-1}$.
\item[{\bf(c)}] For\/ $\al_1,\al_2\in C(\B)_\pe$ with\/ $\al_1+\al_2=\al,$ $\tau(\al_1)=\tau(\al_2)$ and\/ $\M_{\al_i}^\ss(\tau)\ne\es$ for $i=1,2,$ there is a substack\/ $\acM_{(\al,\bs 1)}^\ss(\ac\tau^0_{\bs\mu})_{\al_1,\al_2}$ with\/ $\C$-points $[E_1\op E_2,\bs V,\bs\rho]$ such that\/ $\lb E_i\rb=\al_i$ and $\phi_i(V_i)\subseteq F_{k_i}(E_i)\subset F_{k_i}(E_1\op E_2)$ for $i=1,2$. The splitting $E=E_1\op E_2$ is determined uniquely by the requirement that\/ $\phi_i(V_i)\subseteq F_{k_i}(E_i)$ for $i=1,2$ and\/ $\tau(\al_1)=\tau(\al_2)$. 

Taking $\th\in\bG_m$ to act by $\th\,\id_{E_1}+\id_{E_2},\th\,\id_{V_1},\id_{V_2},\th\,\id_{V_3}$ on $E=E_1\op E_2,\ab V_1,V_2,V_3,$ respectively, we see that\/ $\rho_1,\ldots,\rho_4$ are $\bG_m$-equivariant under the action {\rm\eq{co9eq10},} so $[E_1\op E_2,\bs V,\bs\rho]$ is a $\bG_m$-fixed point.

All such\/ $[E_1\op E_2,\bs V,\bs\rho]$ have $\rho_1,\ldots,\rho_4\ne 0$. Thus $[E_1\op E_2,\bs V,\bs\rho]\mapsto\rho_3^a\rho_4^b$ gives a trivialization of\/ $\ac L^\pl_{a,b,-a-b}\vert_{\acM_{(\al,\bs 1)}^\ss(\ac\tau^0_{\bs\mu})_{\al_1,\al_2}}$ for all\/ $a,b\in\Z$.

There is an isomorphism of algebraic spaces
\e
\begin{split}
\Pi_{\al_1,\al_2}&:\acM_{(\al,\bs 1)}^\ss(\ac\tau^0_{\bs\mu})_{\al_1,\al_2}\longra \baM_{(\al_1,1)}^\ss(\bar\tau^0_1)_{k_1}\t\baM_{(\al_2,1)}^\ss(\bar\tau^0_1)_{k_2},\\
\Pi_{\al_1,\al_2}&:[E,\bs V,\bs\rho]\longmapsto \bigl([E_1,V_1,\rho_1],[E_2,V_2,\rho_2]\bigr).
\label{co9eq13}
\end{split}
\e
\end{itemize}

Parts\/ {\bf(a)\rm,\bf(b)} lift to derived moduli stacks, where we regard\/ $\acM_{(\al,\bs 1)}^\ss(\ac\tau^0_{\bs\mu})$ as an open substack of\/ $\bs\dM_{(\al,\bs 1)}^\rpl$ defined in\/ {\rm\eq{co5eq17},} and similarly for $\bs\baM_{(\al,1)}^\ss(\bar\tau^0_1)_{k_1}$ and\/ $\bs\baM_{(\al,1)}^\ss(\bar\tau^0_1)_{k_2}$. However, be warned that\/ {\bf(c)} may \begin{bfseries}not\end{bfseries} lift to derived stacks. This is because the morphism $\bs\Phi_{\al_1,\al_2}:\bs\M_{\al_1}\t\bs\M_{\al_2}\ra\bs\M_\al$ need not lift to $\bs{\dot\Phi}{}^\red_{\al_1,\al_2}:\bs\dM_{\al_1}^\red\t\bs\dM_{\al_2}^\red\ra\bs\dM_\al^\red,$ and\/ $\bs\baM_{(\al_1,1)}^\ss(\bar\tau^0_1)_{k_1},\bs\baM_{(\al_2,1)}^\ss(\bar\tau^0_1)_{k_2},\bs\acM_{(\al,\bs 1)}^\ss(\ac\tau^0_{\bs\mu})$ are constructed using $\bs\dM_{\al_1}^\red,\bs\dM_{\al_2}^\red,\bs\dM_\al^\red$ respectively.
\end{prop}

\begin{proof} First observe that if $[E,\bs V,\bs\rho]\in\acM_{(\al,\bs 1)}^\ss(\ac\tau^0_{\bs\mu})^{\bG_m}$ then $\rho_1\ne 0$, as otherwise $(0,(V_1,0,0),\bs 0)$ would $\ac\tau^0_{\bs\mu}$-destabilize $(E,\bs V,\bs\rho)$. Similarly $\rho_2\ne 0$ because of $(0,(0,V_2,0),\bs 0)$, and $\rho_3,\rho_4$ cannot both be zero, because of $(0,(0,0,V_3),\bs 0)$. Hence we can divide into three cases:
\begin{itemize}
\setlength{\itemsep}{0pt}
\setlength{\parsep}{0pt}
\item[(a$)'$] $\rho_4=0$, $\rho_1,\rho_2,\rho_3\ne 0$; 
\item[(b$)'$] $\rho_3=0$, $\rho_1,\rho_2,\rho_4\ne 0$; and
\item[(c$)'$] $\rho_1,\ldots,\rho_4\ne 0$.
\end{itemize}

It is obvious from \eq{co9eq10} that any $[E,\bs V,\bs\rho]\in\acM_{(\al,\bs 1)}^\ss(\ac\tau^0_{\bs\mu})$ with $\rho_4=0$ lies in $\acM_{(\al,\bs 1)}^\ss(\ac\tau^0_{\bs\mu})^{\bG_m}$, so that $\acM_{(\al,\bs 1)}^\ss(\ac\tau^0_{\bs\mu})_{\rho_4=0}\subseteq\acM_{(\al,\bs 1)}^\ss(\ac\tau^0_{\bs\mu})^{\bG_m}$, as in (a). If $[E,\bs V,\bs\rho]\in\acM_{(\al,\bs 1)}^\ss(\ac\tau^0_{\bs\mu})$ with $\rho_3=0$ then by letting $\th\in\bG_m$ act on $E,V_1,V_2,V_3$ by $\id_E,\id_{V_1},\id_{V_2},\th\,\id_{V_3}$ shows $\acM_{(\al,\bs 1)}^\ss(\ac\tau^0_{\bs\mu})_{\rho_3=0}\subseteq\acM_{(\al,\bs 1)}^\ss(\ac\tau^0_{\bs\mu})^{\bG_m}$, as in~(b). 

If $[E,\bs V,\bs\rho]$ is a $\bG_m$-fixed point of type (c$)'$ then there must exist morphisms from $\bG_m$ to $\Aut(E)$ and $\Aut(V_i)$ for $i=1,2,3$ such that $\rho_1,\ldots,\rho_4$ are $\bG_m$-equivariant under the action \eq{co9eq10}. As $V_1,V_2,V_3\cong\C$ and we are free to change the actions on $E,V_i$ by an overall factor, we can choose the $\bG_m$-action on $V_3$ to be $\th\mapsto\th\,\id_{V_3}$. Since $\rho_3,\rho_4\ne 0$ this implies that $\bG_m$ must act on $V_1,V_2$ by $\th\,\id_{V_1},\id_{V_2}$, as in part (c) of the proposition.

The action of $\bG_m$ on $E$ corresponds to a direct sum $E=\bigop_{i\in\Z}E^i$, where $E^i\ne 0$ for only finitely many $i\in\Z$, and $\th\in\bG_m$ acts by $\sum_{i\in\Z}\th^i\,\id_{E^i}$. For $\rho_1,\rho_2$ to be $\bG_m$-equivariant we see that $\rho_1(V_1)\subseteq F_{k_1}(E^1)\subseteq F_{k_1}(E)$ and $\rho_2(V_2)\subseteq F_{k_2}(E^0)\subseteq F_{k_2}(E)$. Then $E^i=0$ for $i\ne 0,1$, as otherwise the subobject $(E^0\op E^1,\bs V,\bs\rho)$ would $\ac\tau^0_{\bs\mu}$-destabilize $(E,\bs V,\bs\rho)$. Thus writing $E_1=E^1$ and $E_2=E^0$ we see that $E=E_1\op E_2$ and $\bG_m$ acts on $E$ by $\th\,\id_{E_1}+\id_{E_2}$, as in (c). Since $E$ is $\tau$-semistable this forces $E_1,E_2$ to be $\tau$-semistable with $\tau(\lb E_1\rb)=\tau(\lb E_2\rb)$, so writing $\al_i=\lb E_i\rb$ we have $\al_1+\al_2=\al,$ $\tau(\al_1)=\tau(\al_2)$ and\/ $\M_{\al_i}^\ss(\tau)\ne\es$ for $i=1,2$, as in (c). In fact $E=E_1\op E_2$ is the unique splitting with $\rho_i(V_i)\subseteq F_{k_i}(E_i)\subseteq F_{k_i}(E)$ for $i=1,2$ and $\tau(\lb E_1\rb)=\tau(\lb E_2\rb)$, as if $E=E'_1\op E'_2$ were a distinct such splitting then the subobject $((E_1\cap E_1')\op (E_2\cap E_2'),\bs V,\bs\rho)$ would $\ac\tau^0_{\bs\mu}$-destabilize~$(E,\bs V,\bs\rho)$.

We have described $\bG_m$-fixed substacks $\acM_{(\al,\bs 1)}^\ss(\ac\tau^0_{\bs\mu})_{\rho_4=0}$, $\acM_{(\al,\bs 1)}^\ss(\ac\tau^0_{\bs\mu})_{\rho_3=0}$ and $\acM_{(\al,\bs 1)}^\ss(\ac\tau^0_{\bs\mu})_{\al_1,\al_2}$ in $\acM_{(\al,\bs 1)}^\ss(\ac\tau^0_{\bs\mu})$, which are closed in $\acM_{(\al,\bs 1)}^\ss(\ac\tau^0_{\bs\mu})$, are disjoint, and include all $\bG_m$-fixed points. There are finitely many of them, as there are only finitely many possibilities for $\al_1,\al_2$ by Lemma \ref{co9lem1}. Hence they are open and closed, and $\acM_{(\al,\bs 1)}^\ss(\ac\tau^0_{\bs\mu})^{\bG_m}$ is the disjoint union of the substacks in~(a)--(c).

The claims on triviality of line bundles in (a)--(c) are immediate from those on when $\rho_3\ne 0$ or $\rho_4\ne 0$. This completes the first parts of (a)--(c). It remains to prove that \eq{co9eq11}--\eq{co9eq13} are well defined and smooth/an isomorphism.

For \eq{co9eq11}, suppose $[E,\bs V,\bs\rho]$ is a $\C$-point in $\acM_{(\al,\bs 1)}^\ss(\ac\tau^0_{\bs\mu})_{\rho_4=0}$, so that $\rho_4=0$ and $\rho_1,\rho_2,\rho_3\ne 0$. Then $E$ is $\tau$-semistable as in Definition \ref{co9def2}. If $0\ne E'\subsetneq E$ with $\rho_1(V_1)\subseteq F_{k_1}(E')\subsetneq F_{k_1}(E)$ and $\tau(\lb E'\rb)=\tau(\lb E/E'\rb)$ then using \eq{co5eq21} we find that $(E',(V_1,0,V_3),\bs\rho')$ would $\ac\tau^0_{\bs\mu}$-destabilize $(E,\bs V,\bs\rho)$, a contradiction. So no such $E'$ exists, and $(E,V_1,\rho_1)$ is $\bar\tau^0_1$-semistable by Example \ref{co5ex1}(i)--(iii). Thus the morphism \eq{co9eq11} is well~defined.

Consider the fibre of \eq{co9eq11} over a point $[E,V_1,\rho_1]$ in $\baM_{(\al,1)}^\ss(\bar\tau^0_1)_{k_1}$. For $[E,\bs V,\bs\rho]$ in the fibre we have $\rho_3\ne 0$, so $\rho_3:V_3\ra V_1$ is an isomorphism, and we are free to take $V_3=V_1$ and $\rho_3=\id_{V_1}$. As $\rho_4=0$, the remaining data in $[E,\bs V,\bs\rho]$ not fixed in $[E,V_2,\rho_2]$ is $V_2\cong\C$ and $\rho_2:V_2\ra F_{k_2}(E_2)$. Given that $(E,V_1,\rho_1)$ is $\bar\tau^0_1$-semistable, we can show using \eq{co5eq21} and \eq{co9eq8} that $(E,(V_1,V_2,V_1),\ab (\rho_1,\rho_2,\id_{V_1},0))$ is $\ac\tau^0_{\bs\mu}$-semistable if and only if $\rho_2\ne 0$. Thus the pair $(V_2,\rho_2)$ is equivalent to a point $\rho_2(V_2)$ in $\bP(F_{k_2}(E))$, so \eq{co9eq11} has fibre $\bP(F_{k_2}(E))\cong\CP^{\la_{k_2}(\al)-1}$ over $[E,V_1,\rho_1]$. Hence \eq{co9eq11} is smooth of dimension $\la_{k_2}(\al)-1$, with fibre $\CP^{\la_{k_2}(\al)-1}$. This proves (a). The rest of (b) is similar.

For \eq{co9eq13}, suppose $[E,\bs V,\bs\rho]$ is a $\C$-point in $\acM_{(\al,\bs 1)}^\ss(\ac\tau^0_{\bs\mu})_{\al_1,\al_2}$. Then $E$ is $\tau$-semistable as in Definition \ref{co9def2}. If $0\ne E_1'\subsetneq E_1$ with $\rho_1(V_1)\subseteq F_{k_1}(E_1')\subsetneq F_{k_1}(E_1)$ and $\tau(\lb E_1'\rb)=\tau(\lb E_1/E_1'\rb)$ then $(E_1'\op E_2,\bs V,\bs\rho)$ would $\ac\tau^0_{\bs\mu}$-destabilize $(E,\bs V,\bs\rho)$, a contradiction. So no such $E_1'$ exists, and $(E_1,V_1,\rho_1)$ is $\bar\tau^0_1$-semistable by Example \ref{co5ex1}(i)--(iii). Similarly $(E_2,V_2,\rho_2)$ is $\bar\tau^0_1$-semistable. Thus the morphism \eq{co9eq13} is well defined.

If $([E_1,V_1,\rho_1],[E_2,V_2,\rho_2])$ is a $\C$-point of $\baM_{(\al_1,1)}^\ss(\bar\tau^0_1)_{k_1}\t\baM_{(\al_2,1)}^\ss(\bar\tau^0_1)_{k_2}$ then as $V_1\cong\C\cong V_2$ we may take $V_3=V_1$, and choose any isomorphism $\rho_4:V_1\ra V_2$, and show that $[E,\bs V,\bs\rho]=[E_1\op E_2,(V_1,V_2,V_1),(\rho_1,\rho_2,\id_{V_1},\rho_4)]$ is a $\C$-point of $\acM_{(\al,\bs 1)}^\ss(\ac\tau^0_{\bs\mu})_{\al_1,\al_2}$. As $[E,\bs V,\bs\rho]$ is unchanged by $\rho_4\mapsto\th\,\rho_4$ by \eq{co9eq10}, it is independent of the choice of isomorphism $\rho_4:V_1\ra V_2$. Thus we may define a stack morphism $([E_1,V_1,\rho_1],[E_2,V_2,\rho_2])\mapsto[E,\bs V,\bs\rho]$ inverse to \eq{co9eq13}, so \eq{co9eq13} is an isomorphism. This completes (c). The extension of (a),(b) to derived stacks is immediate, as the proofs above also work in derived stacks.
\end{proof}

Next we describe the obstruction theory on $\acM_{(\al,\bs 1)}^\ss(\ac\tau^0_{\bs\mu})^{\bG_m}$.

\begin{prop}
\label{co9prop3}
As in\/ {\rm\S\ref{co26},} the restriction of the $\bG_m$-equivariant obstruction theory\/ $\bL_i:i^*(\bL_{\bs\acM_{(\al,\bs 1)}^\ss(\ac\tau^0_{\bs\mu})})\ra \bL_{\acM_{(\al,\bs 1)}^\ss(\ac\tau^0_{\bs\mu})}$ on $\acM_{(\al,\bs 1)}^\ss(\ac\tau^0_{\bs\mu})$ to $\acM_{(\al,\bs 1)}^\ss(\ac\tau^0_{\bs\mu})^{\bG_m}$ splits as the direct sum of a $\bG_m$-fixed part, which is the natural obstruction theory $\bL_{i^{\bG_m}}:(i^{\bG_m})^*(\bL_{\bs\acM_{(\al,\bs 1)}^\ss(\ac\tau^0_{\bs\mu})^{\bG_m}})\ra \bL_{\acM_{(\al,\bs 1)}^\ss(\ac\tau^0_{\bs\mu})^{\bG_m}}$ associated to the proper quasi-smooth derived algebraic space $\bs\acM_{(\al,\bs 1)}^\ss(\ac\tau^0_{\bs\mu})^{\bG_m},$ and a\/ $\bG_m$-moving part\/ $\cN^\bu,$ the \begin{bfseries}virtual conormal bundle\end{bfseries} of\/ $\acM_{(\al,\bs 1)}^\ss(\ac\tau^0_{\bs\mu})^{\bG_m}$ in $\acM_{(\al,\bs 1)}^\ss(\ac\tau^0_{\bs\mu})$.

In Proposition\/ {\rm\ref{co9prop2}(a)--(c)} we may identify\/ $\cN^\bu$ as follows:
\begin{itemize}
\setlength{\itemsep}{0pt}
\setlength{\parsep}{0pt}
\item[{\bf(a)}] The virtual conormal bundle of\/ $\acM_{(\al,\bs 1)}^\ss(\ac\tau^0_{\bs\mu})_{\rho_4=0}$ in $\acM_{(\al,\bs 1)}^\ss(\ac\tau^0_{\bs\mu})$ is $\cN_{\rho_4=0}^\bu\cong \ac L^\pl_{0,-1,1}\vert_{\acM_{(\al,\bs 1)}^\ss(\ac\tau^0_{\bs\mu})_{\rho_4=0}},$ with\/ $\bG_m$-weight\/~$-1$. 
\item[{\bf(b)}] The virtual conormal bundle of\/ $\acM_{(\al,\bs 1)}^\ss(\ac\tau^0_{\bs\mu})_{\rho_3=0}$ in $\acM_{(\al,\bs 1)}^\ss(\ac\tau^0_{\bs\mu})$ is $\cN_{\rho_3=0}^\bu\cong \ac L^\pl_{-1,0,1}\vert_{\acM_{(\al,\bs 1)}^\ss(\ac\tau^0_{\bs\mu})_{\rho_3=0}},$ with\/ $\bG_m$-weight\/~$1$. 
\item[{\bf(c)}] The virtual conormal bundle of\/ $\acM_{(\al,\bs 1)}^\ss(\ac\tau^0_{\bs\mu})_{\al_1,\al_2}$ in $\acM_{(\al,\bs 1)}^\ss(\ac\tau^0_{\bs\mu})$ is $\cN_{\al_1,\al_2}^\bu\ab =\cN_{\al_1,\al_2}^{+,\bu}\op \cN_{\al_1,\al_2}^{-,\bu},$ where $\cN_{\al_1,\al_2}^{\pm,\bu}$ has\/ $\bG_m$-weight\/ $\pm 1,$ and using the identification {\rm\eq{co9eq13},} in $K^0(\Perf_{\baM_{(\al_1,1)}^\ss(\bar\tau^0_1)_{k_1}\t\baM_{(\al_2,1)}^\ss(\bar\tau^0_1)_{k_2}})$ we have
\end{itemize}	
\e
\begin{split}
\bigl[\cN_{\al_1,\al_2}^{+,\bu}\bigr]&=\bigl[(\bar\Pi^v_{\dot{\mathcal M}_{\al_2}})^*(\cV^*_{k_1,\al_2})\bigr]
-\bigl[(\bar\Pi^v_{\dot{\mathcal M}_{\al_1}}\t\bar\Pi^v_{\dot{\mathcal M}_{\al_2}})^*(\cE^\bu_{\al_1,\al_2})\bigr],
\\
\bigl[\cN_{\al_1,\al_2}^{-,\bu}\bigr]&=\bigl[(\bar\Pi^v_{\dot{\mathcal M}_{\al_1}})^*(\cV^*_{k_2,\al_1})\bigr]
-\bigl[(\bar\Pi^v_{\dot{\mathcal M}_{\al_1}}\t\bar\Pi^v_{\dot{\mathcal M}_{\al_2}})^*(\si_{\al_1,\al_2}^*(\cE_{\al_2,\al_1}^\bu)\bigr].
\end{split}
\label{co9eq14}
\e
\begin{itemize}
\setlength{\itemsep}{0pt}
\setlength{\parsep}{0pt}
\item[] Here $\bar\Pi^v_{\dot{\mathcal M}_{\al_1}}:\baM_{(\al_1,1)}^\ss(\bar\tau^0_1)_{k_1}\ra\dM_{\al_1}$ and\/ $\bar\Pi^v_{\dot{\mathcal M}_{\al_2}}:\baM_{(\al_2,1)}^\ss(\bar\tau^0_1)_{k_2}\ra\dM_{\al_2}$ are as in equation\/~\eq{co5eq26}.
\item[{\bf(d)}] For\/ $\al_1,\al_2$ as in Proposition\/ {\rm\ref{co9prop2}(c),} we have $o_{\al_1}+o_{\al_2}\ge o_\al$. Under the isomorphism {\rm\eq{co9eq13},} the virtual classes satisfy
\ea
&[\acM_{(\al,\bs 1)}^\ss(\ac\tau^0_{\bs\mu})_{\al_1,\al_2}]_\virt
\label{co9eq15}\\
&=\begin{cases}
[\baM_{(\al_1,1)}^\ss(\bar\tau^0_1)_{k_1}]_\virt\bt[\baM_{(\al_2,1)}^\ss(\bar\tau^0_1)_{k_2}]_\virt, & o_{\al_1}+o_{\al_2}=o_\al, \\
0, & o_{\al_1}+o_{\al_2}>o_\al,
\end{cases}
\nonumber
\ea
where\/ $[\baM_{(\al_1,1)}^\ss(\bar\tau^0_1)_{k_1}]_\virt,[\baM_{(\al_2,1)}^\ss(\bar\tau^0_1)_{k_2}]_\virt$ are as in Example\/~{\rm\ref{co5ex1}}. 
\end{itemize}	
Part\/ {\bf(d)} will compensate for the problem noted at the end of Proposition\/ {\rm\ref{co9prop2}}.
\end{prop}

\begin{proof} For the first part, general properties of obstruction theories imply that the $\bG_m$-fixed part of the restriction to $\acM_{(\al,\bs 1)}^\ss(\ac\tau^0_{\bs\mu})^{\bG_m}$ of an obstruction theory on $\acM_{(\al,\bs 1)}^\ss(\ac\tau^0_{\bs\mu})$ is an obstruction theory on $\acM_{(\al,\bs 1)}^\ss(\ac\tau^0_{\bs\mu})^{\bG_m}$, and as the obstruction theory comes from a $\bG_m$-equivariant quasi-smooth derived enhancement $\bs\acM_{(\al,\bs 1)}^\ss(\ac\tau^0_{\bs\mu})$, this $\bG_m$-fixed part comes from the quasi-smooth derived enhancement $\bs\acM_{(\al,\bs 1)}^\ss(\ac\tau^0_{\bs\mu})^{\bG_m}$ of~$\acM_{(\al,\bs 1)}^\ss(\ac\tau^0_{\bs\mu})^{\bG_m}$.

For (a), consider the (2-)commutative diagram of Artin $\C$-stacks
\e
\begin{gathered}
\xymatrix@!0@C=140pt@R=20pt{
& \,\,\acM_{(\al,\bs 1)}^\ss(\ac\tau^0_{\bs\mu})_{\rho_4=0}\, \ar@{_{(}->}[dl]_(0.45){i_{\rho_4=0}} \ar[dr]^(0.45){\Pi_{\rho_4=0}}_(0.4){\eq{co9eq11}} \ar[dd]_{\Pi_{\M^\ss_\al(\tau)}^{\rho_4=0}} \\
*+[r]{\acM_{(\al,\bs 1)}^\ss(\ac\tau^0_{\bs\mu})} \ar[dr]_{\Pi_{\M^\ss_\al(\tau)}} && *+[l]{\baM_{(\al,1)}^\ss(\bar\tau^0_1)_{k_1}} \ar[dl]^{\Pi_{\M^\ss_\al(\tau)}} \\
& \M^\ss_\al(\tau). }	
\end{gathered}
\label{co9eq16}
\e
All morphisms apart from $i_{\rho_4=0}$ are smooth. Equation \eq{co9eq16} is the classical truncation of a corresponding diagram of derived stacks, with all morphisms apart from $\bs i_{\rho_4=0}$ smooth.

Taking the pullback of \eq{co2eq11} for $\bs\Pi_{\M^\ss_\al(\tau)}:\bs\acM_{(\al,\bs 1)}^\ss(\ac\tau^0_{\bs\mu})\ra\bs\M^\ss_\al(\tau)$ to $\acM_{(\al,\bs 1)}^\ss(\ac\tau^0_{\bs\mu})$, where the pullback of $\bL_{\bs\acM_{(\al,\bs 1)}^\ss(\ac\tau^0_{\bs\mu})/\bs\M^\ss_\al(\tau)}$ is $\bL_{\acM_{(\al,\bs 1)}^\ss(\ac\tau^0_{\bs\mu})/\M^\ss_\al(\tau)}$ as $\bs\Pi_{\M^\ss_\al(\tau)}$ is smooth, gives a distinguished triangle on $\acM_{(\al,\bs 1)}^\ss(\ac\tau^0_{\bs\mu})$:
\e
\xymatrix@C=15pt{ \Pi_{\M^\ss_\al(\tau)}^*(i^*(\bL_{\bs\M^\ss_\al(\tau)})) \ar[r] & i^*(\bL_{\bs\acM_{(\al,\bs 1)}^\ss(\ac\tau^0_{\bs\mu})}) \ar[r] & \bL_{\acM_{(\al,\bs 1)}^\ss(\ac\tau^0_{\bs\mu})/\M^\ss_\al(\tau)} \ar[r]^(0.8){[1]} & . }
\label{co9eq17}
\e
Now $\bs\Pi_{\M^\ss_\al(\tau)}:\bs\acM_{(\al,\bs 1)}^\ss(\ac\tau^0_{\bs\mu})\ra\bs\M^\ss_\al(\tau)$ is $\bG_m$-equivariant for the trivial $\bG_m$-action on $\bs\M^\ss_\al(\tau)$. Thus the restriction of \eq{co9eq17} to $\acM_{(\al,\bs 1)}^\ss(\ac\tau^0_{\bs\mu})_{\rho_4=0}$ is equivariant under a $\bG_m$-action, where the $\bG_m$-action on $\Pi_{\M^\ss_\al(\tau)}^*(i^*(\bL_{\bs\M^\ss_\al(\tau)}))$ is trivial.

An easy calculation shows that in vector bundles on $\acM_{(\al,\bs 1)}^\ss(\ac\tau^0_{\bs\mu})_{\rho_4=0}$
\ea
&\bL_{\acM_{(\al,\bs 1)}^\ss(\ac\tau^0_{\bs\mu})/\M^\ss_\al(\tau)}\vert_{\acM_{(\al,\bs 1)}^\ss(\ac\tau^0_{\bs\mu})_{\rho_4=0}}
\nonumber\\
&\cong \biggl(\frac{\ts (\ac\Pi^{v_3}_{\dot{\mathcal M}_\al})^*(\cV_{k_1,\al})}{\ts\rho_1(\O)}\biggr)^*_{\wt=0}\!\op\! \biggl(\frac{\ts (\ac\Pi^{v_3}_{\dot{\mathcal M}_\al})^*(\cV_{k_2,\al})}{\ts\rho_2(\O)}\biggr)^*_{\wt=0}\!\op\! (\ac L_{0,-1,1}^\pl)_{\wt=-1}
\nonumber\\
&\cong\bigl(\bL_{\acM_{(\al,\bs 1)}^\ss(\ac\tau^0_{\bs\mu})_{\rho_4=0}/\M^\ss_\al(\tau)}\bigr){}_{\wt=0}\op (\ac L_{0,-1,1}^\pl)_{\wt=-1}.
\label{co9eq18}
\ea
Here $\ac\Pi^{v_3}_{\dot{\mathcal M}_\al}$ is as in \eq{co5eq26}. In the middle line, at a point $[E,\bs V,\bs\rho]\mapsto[E]$, the first term is the contribution of $(V_1,\rho_1)$, or equivalently $(V_1,V_3,\rho_1,\rho_3)$ (noting that $\rho_1,\rho_3\ne 0$), the second the contribution of $(V_2,\rho_2)$ (noting that $\rho_2\ne 0$), and the third the contribution of $\rho_4$, which gives a section of $\ac L_{0,1,-1}^\pl$ on $\acM_{(\al,\bs 1)}^\ss(\ac\tau^0_{\bs\mu})$. The subscripts `$\wt=0$', `$\wt=-1$' indicate the $\bG_m$-weights, which hold as $\bG_m$  acts on $E,V_1,V_2,V_3$ by $\id_E,\id_{V_1},\id_{V_2},\id_{V_3}$ and by \eq{co9eq10} on $\rho_1,\ldots,\rho_4$, as in the proof of Proposition~\ref{co9prop2}.

Taking $\bG_m$-moving parts in the restriction of \eq{co9eq17} to $\acM_{(\al,\bs 1)}^\ss(\ac\tau^0_{\bs\mu})_{\rho_4=0}$ gives a distinguished triangle in perfect complexes with nonzero $\bG_m$-weights, where the first term is zero as $\Pi_{\M^\ss_\al(\tau)}^*(i^*(\bL_{\bs\M^\ss_\al(\tau)}))$ has trivial $\bG_m$-action, the second term is $\cN_{\rho_4=0}^\bu$, and the third term is $(\ac L_{0,-1,1}^\pl)_{\wt=-1}$ by \eq{co9eq18}. Thus $\cN_{\rho_4=0}^\bu\cong \ac L^\pl_{0,-1,1}\vert_{\acM_{(\al,\bs 1)}^\ss(\ac\tau^0_{\bs\mu})_{\rho_4=0}}$ with $\bG_m$-weight $-1$, proving~(a).

The proof of (b) is essentially the same, except that $(\ac L_{0,-1,1}^\pl)_{\wt=-1}$ in \eq{co9eq18} is replaced by $(\ac L_{-1,0,1}^\pl)_{\wt=1}$, where the $\bG_m$-weight is $1$ as $\th\in\bG_m$ acts on $E,V_1,V_2,V_3$ by $\id_E,\id_{V_1},\id_{V_2},\th\,\id_{V_3}$ as in the proof of Proposition \ref{co9prop2}.

For (c),(d), first note that using \eq{co5eq17} we can construct a $\bG_m$-equivariant Cartesian square of derived stacks
\e
\begin{gathered}
\xymatrix@C=100pt@R=15pt{
*+[r]{\bs\dM_{(\al,\bs 1)}^\rpl} \ar[r]_{\bs{\ac\Pi}{}^{v_3}_{\bs\dM^\red_\al}} \ar[d]^{\,\bs{\ac\jmath}{}_{(\al,\bs 1)}^{\,\pl}} & *+[l]{\bs\dM_\al^\red} \ar[d]_{\bs j_\al} \\
*+[r]{\bs\dM_{(\al,\bs 1)}^\pl} \ar[r]^{\bs{\ac\Pi}{}^{v_3}_{\bs\dM_\al}}  & *+[l]{\bs\dM_\al,\!} }
\end{gathered}
\label{co9eq19}
\e
where $\bs j_\al,\bs{\ac\jmath}{}_{(\al,\bs 1)}^{\,\pl}$ are as in \eq{co5eq4} for $\B,\acB$, and $\bs{\ac\Pi}{}^{v_3}_{\bs\dM^\red_\al}$ is as in \eq{co5eq27}, and $\bs{\ac\Pi}{}^{v_3}_{\bs\dM_\al}$ is the analogue of $\bs{\ac\Pi}{}^{v_3}_{\bs\dM^\red_\al}$ for the `non-reduced' versions of the derived stacks, where $\bs{\ac\Pi}{}^{v_3}_{\bs\dM^\red_\al},\bs{\ac\Pi}{}^{v_3}_{\bs\dM_\al}$ are smooth, and the $\bG_m$-actions on $\bs\dM_\al^\red,\bs\dM_\al$ are trivial.

There is a $\bG_m$-equivariant commutative diagram in $D_\qcoh(\dM_{(\al,\bs 1)}^\pl)$, with rows and diagonals distinguished triangles:
\e
\,\,\xymatrix@!0@C=38pt@R=40pt{
{\!\!\!\!\!\begin{subarray}{l}\ts (\ac\Pi^{v_3}_{\dot{\mathcal M}_\al})^*\!\ci\! i^* \\ \ts (\bL_{\sst\bs\dM^\red_\al/\bs\dM_\al})[-1]\!\!\!\!\!\!\!\!\!\!\!\!\!\!\!\!\!\!\!\!\!\!\!\!\end{subarray}} \ar[dr] \ar@/_4pc/[ddddrr]_0 \ar@<1ex>[rr]^(0.35)\cong && {\begin{subarray}{l}\ts i^*(\bL_{\sst\bs\dM_{(\al,\bs 1)}^\rpl/\bs\dM_{(\al,\bs 1)}^\pl})\!\!\! \\ \ts\qquad\qquad [-1]\end{subarray}} \ar@/_2.2pc/[ddddrr]_(0.8)0 \ar[dr] \ar[rr] && 0 \ar[dr] \ar[rr]_(0.7){[+1]} &&
\\
& \!\!\!\!\!\!\!\!(\ac\Pi^{v_3}_{\dot{\mathcal M}_\al})^*(i^*(\bL_{\bs\dM_\al})) \ar@/_1pc/[dddr]_{\begin{subarray}{l}(\ac\Pi^{v_3}_{\dot{\mathcal M}_\al})^* \\ (i^*(\bL_i))\end{subarray}\!\!\!\!\!\!\!\!\!\!\!} \ar[dr] \ar[rr] && i^*(\bL_{\bs\dM_{(\al,\bs 1)}^\pl}) \ar@/_.4pc/[dddr]^(0.55){\bL_i} \ar[dr] \ar[rr] && i^*(\bL_{\bs\dM_{(\al,\bs 1)}^\pl/\bs\dM_\al}) \ar[dr]^(0.6)\cong \ar[rr]^(0.7){[+1]} && 
\\
&& \!\!\!\!\!\!\!\!(\ac\Pi^{v_3}_{\dot{\mathcal M}_\al})^*(i^*(\bL_{\bs\dM^\red_\al})) \ar[dr]^(0.7){[+1]} \ar[dd]^(0.55){\!\!\!\!\!\!\!\!\!\!\!\!\begin{subarray}{l}(\ac\Pi^{v_3}_{\dot{\mathcal M}_\al})^* \\ (i^*(\bL_i))\end{subarray}} \ar[rr] && i^*(\bL_{\bs\dM_{(\al,\bs 1)}^\rpl}) \ar[dd]^{\bL_i} \ar[dr]^(0.7){[+1]} \ar[rr] && i^*(\bL_{\bs\dM_{(\al,\bs 1)}^\rpl/\bs\dM_\al^\red}) \ar[dr]^(0.7){[+1]} \ar[dd]_\cong \ar[rr]^(0.7){[+1]} && 
\\
&&&&&&&
\\
&& (\ac\Pi^{v_3}_{\dot{\mathcal M}_\al})^*(\bL_{\dot{\mathcal M}_\al}) \ar[rr] && \bL_{\dM_{(\al,\bs 1)}^\pl}  \ar[rr] && \bL_{\dM_{(\al,\bs 1)}^\pl/\dM_\al} \ar[rr]^(0.7){[+1]} &&  }
\!\!\!\!\!\!\!\!\!\!\!\!\!\!\!\!\!
\label{co9eq20}
\e

Here the top three rows come from \eq{co9eq19} by properties of cotangent complexes of derived stacks in Cartesian squares, and in the bottom row we project to cotangent complexes of the underlying classical stacks. The bottom middle column $\bL_i$ is the obstruction theory on $\dM_{(\al,\bs 1)}^\pl$ which restricts to that on $\acM_{(\al,\bs 1)}^\ss(\ac\tau^0_{\bs\mu})$. The top and middle isomorphisms marked in \eq{co9eq20} hold as \eq{co9eq19} is Cartesian in derived stacks, and the bottom holds as $\bs{\ac\Pi}{}^{v_3}_{\bs\dM^\red_\al}$ is smooth.

We rewrite \eq{co9eq20} using the marked isomorphisms and \eq{co5eq1}, \eq{co5eq6} to obtain
\e
\xymatrix@!0@C=38pt@R=35pt{
U_\al\!\ot\!\O[1] \ar[dr] \ar@/_4pc/[ddddrr]_0 \ar[rr]^\id && U_\al\!\ot\!\O[1] \ar@/_2.2pc/[ddddrr]_(0.8)0 \ar[dr] \ar[rr] && 0 \ar[dr] \ar[rr]_(0.7){[+1]} &&
\\
& {\begin{subarray}{l}\ts (\De_{\dot{\mathcal M}_\al}\ci\ac\Pi^{v_3}_{\dot{\mathcal M}_\al})^* \\ \ts (\cE_{\al,\al}^\bu)[-1]\end{subarray}} \ar@/_1pc/[dddr]_(0.55){\begin{subarray}{l}(\ac\Pi^{v_3}_{\dot{\mathcal M}_\al})^* \\ (i^*(\bL_i)) \\ \quad \ci\th_\al \end{subarray}\!\!\!\!\!\!\!\!\!\!\!} \ar[dr] \ar[rr] && i^*(\bL_{\bs\dM_{(\al,\bs 1)}^\pl}) \ar@/_.4pc/[dddr]^(0.55){\bL_i} \ar[dr] \ar[rr] && \bL_{\dM_{(\al,\bs 1)}^\pl/\dM_\al} \ar[dr]^(0.6)\id \ar[rr]^(0.7){[+1]} && 
\\
&& \!\!\!\!\!\!\!\!(\ac\Pi^{v_3}_{\dot{\mathcal M}_\al})^*(i^*(\bL_{\bs\dM^\red_\al})) \ar[dr]^(0.7){[+1]} \ar[dd]^(0.55){\!\!\!\!\!\!\!\!\!\!\!\!\begin{subarray}{l}(\ac\Pi^{v_3}_{\dot{\mathcal M}_\al})^* \\ (i^*(\bL_i))\end{subarray}} \ar[rr] && i^*(\bL_{\bs\dM_{(\al,\bs 1)}^\rpl}) \ar[dd]^{\bL_i} \ar[dr]^(0.7){[+1]} \ar[rr] && \bL_{\dM_{(\al,\bs 1)}^\pl/\dM_\al} \ar[dr]^(0.7){[+1]} \ar[dd]_\id \ar[rr]^(0.7){[+1]} && 
\\
&&&&&&&
\\
&& (\ac\Pi^{v_3}_{\dot{\mathcal M}_\al})^*(\bL_{\dot{\mathcal M}_\al}) \ar[rr] && \bL_{\dM_{(\al,\bs 1)}^\pl}  \ar[rr] && \bL_{\dM_{(\al,\bs 1)}^\pl/\dM_\al} \ar[rr]^(0.7){[+1]} &&  }
\!\!\!\!\!\!\!\!\!\!\!\!\!
\label{co9eq21}
\e

Pulling back \eq{co9eq21} by $i_{\al_1,\al_2}:\acM_{(\al,\bs 1)}^\ss(\ac\tau^0_{\bs\mu})_{\al_1,\al_2}\hookra\acM_{(\al,\bs 1)}^\ss(\ac\tau^0_{\bs\mu})\subset\dM_{(\al,\bs 1)}^\pl$ gives a similar diagram in $D_\qcoh(\acM_{(\al,\bs 1)}^\ss(\ac\tau^0_{\bs\mu})_{\al_1,\al_2})$ equivariant under a $\bG_m$-action on each term, where we write $\ac\imath_{\dot{\mathcal M}_\al}=\ac\Pi^{v_3}_{\dot{\mathcal M}_\al}\ci i_{\al_1,\al_2}$ and simplify notation:
\e
\!\!\!\!\!\xymatrix@!0@C=40pt@R=40pt{
U_\al\!\ot\!\O[1] \ar@/_4pc/[ddddrr]_0 \ar[dr] \ar[rr]^\id && U_\al\!\ot\!\O[1] \ar@/_2.2pc/[ddddrr]_(0.8)0  \ar[dr] \ar[rr] && 0 \ar[dr] \ar[rr]_(0.7){[+1]} &&
\\
& {\begin{subarray}{l}\ts (\De_{\dot{\mathcal M}_\al}\ci\ac\imath_{\dot{\mathcal M}_\al})^* \\ \ts (\cE_{\al,\al}^\bu)[-1]\end{subarray}} \ar@/_1pc/[dddr]_(0.53){\begin{subarray}{l}\ac\imath_{\dot{\mathcal M}_\al}^* \\ (i^*(\bL_i))\ci{} \\  i_{\al_1,\al_2}^*(\th_\al) \end{subarray}\!\!\!\!\!\!\!\!\!\!\!} \ar[dr] \ar[rr] && {\begin{subarray}{l}\ts \quad i_{\al_1,\al_2}^* \\ \ts (\bL_{\bs\dM_{(\al,\bs 1)}^\pl})\end{subarray}} \ar@/_.4pc/[dddr]^(0.6){\!\!\!\!\!\!\!\!\!\begin{subarray}{l} i_{\al_1,\al_2}^* \\ (\bL_i)\end{subarray}} \ar[dr] \ar[rr] && {\begin{subarray}{l}\ts \quad i_{\al_1,\al_2}^* \\ \ts (\bL_{\acM_{(\al,\bs 1)}^\ss(\ac\tau^0_{\bs\mu})/\dM_\al})\end{subarray}} \ar[dr]^(0.6)\id \ar[rr]^(0.7){[+1]} && \\
&& \ac\imath_{\dot{\mathcal M}_\al}^*(i^*(\bL_{\bs\dM^\red_\al})) \ar[dr]^(0.7){[+1]} \ar[dd]^(0.6){\!\!\!\!\!\!\!\begin{subarray}{l}\ac\imath_{\dot{\mathcal M}_\al}^* \\ (i^*(\bL_i)) \end{subarray}} \ar[rr] && {\begin{subarray}{l}\ts \quad i_{\al_1,\al_2}^* \\ \ts (\bL_{\bs\acM_{(\al,\bs 1)}^\ss(\ac\tau^0_{\bs\mu})})\end{subarray}} \ar[dd]^(0.7){i_{\al_1,\al_2}^*(\bL_i)} \ar[dr]^(0.7){[+1]} \ar[rr] && {\begin{subarray}{l}\ts \quad i_{\al_1,\al_2}^* \\ \ts (\bL_{\acM_{(\al,\bs 1)}^\ss(\ac\tau^0_{\bs\mu})/\dM_\al})\end{subarray}} \ar[dr]^(0.7){[+1]} \ar[dd]_\id \ar[rr]^(0.7){[+1]} && 
\\
&&&&&&&
\\
&& \ac\imath_{\dot{\mathcal M}_\al}^*(\bL_{\dot{\mathcal M}_\al}) \ar[rr] && {\begin{subarray}{l}\ts \quad i_{\al_1,\al_2}^* \\ \ts (\bL_{\acM_{(\al,\bs 1)}^\ss(\ac\tau^0_{\bs\mu})})\end{subarray}}  \ar[rr] && {\begin{subarray}{l}\ts \quad i_{\al_1,\al_2}^* \\ \ts (\bL_{\acM_{(\al,\bs 1)}^\ss(\ac\tau^0_{\bs\mu})/\dM_\al})\end{subarray}} \ar[rr]^(0.7){[+1]} &&  }
\!\!\!\!\!\!\!\!\!\!\!\!\!\!\!\!\!\!\!\!\!\!\!\!\!
\label{co9eq22}
\e
Here we have replaced $\dM_{(\al,\bs 1)}^\pl$ and $\bs\dM_{(\al,\bs 1)}^\rpl$ by their open substacks $\acM_{(\al,\bs 1)}^\ss(\ac\tau^0_{\bs\mu})$ and $\bs\acM_{(\al,\bs 1)}^\ss(\ac\tau^0_{\bs\mu})$, as these contain $\acM_{(\al,\bs 1)}^\ss(\ac\tau^0_{\bs\mu})_{\al_1,\al_2}$.

Under the identification \eq{co9eq13}, $\ac\imath_{\dot{\mathcal M}_\al}$ becomes the composition
\e
\xymatrix@C=40pt{ \baM_{(\al_1,1)}^\ss(\bar\tau^0_1)_{k_1}\!\t\!\baM_{(\al_2,1)}^\ss(\bar\tau^0_1)_{k_2} \ar[r]^(0.63){\bar\Pi^v_{\dot{\mathcal M}_{\al_1}}\t\bar\Pi^v_{\dot{\mathcal M}_{\al_2}} } & \dM_{\al_1}\!\t\!\dM_{\al_2} \ar[r]^(0.54){\Phi_{\al_1,\al_2}\vert_{\cdots}} & \dM_\al, }	
\label{co9eq23}
\e
where $\bar\Pi^v_{\dot{\mathcal M}_{\al_1}},\bar\Pi^v_{\dot{\mathcal M}_{\al_2}}$ are as in \eq{co5eq26}. 

We can now show that
\ea
&(\De_{\dot{\mathcal M}_\al}\ci\ac\imath_{\dot{\mathcal M}_\al})^*(\cE_{\al,\al}^\bu)[-1]=(\bar\Pi^v_{\dot{\mathcal M}_{\al_1}}\t\bar\Pi^v_{\dot{\mathcal M}_{\al_2}})^*\ci\Phi_{\al_1,\al_2}^*(\De_{\M_\al}^*(\cE_{\al,\al}^\bu))[-1]
\nonumber\\
&\cong(\si_{23}\ci(\De_{\M_{\al_1}}\t\De_{\M_{\al_2}})\ci(\bar\Pi^v_{\dot{\mathcal M}_{\al_1}}\t\bar\Pi^v_{\dot{\mathcal M}_{\al_2}}))^*\ci(\Phi_{\al_1,\al_2}\t\Phi_{\al_1,\al_2})^*(\cE_{\al,\al}^\bu)[-1]
\nonumber\\
&\cong(\si_{23}\ci(\De_{\M_{\al_1}}\t\De_{\M_{\al_2}})\ci(\bar\Pi^v_{\dot{\mathcal M}_{\al_1}}\t\bar\Pi^v_{\dot{\mathcal M}_{\al_2}}))^*
\nonumber\\
&\qquad\bigl((\cE_{\al_1,\al_1}^\bu)_{13}\op(\cE_{\al_2,\al_2}^\bu)_{24}\op(\cE_{\al_1,\al_2}^\bu)_{14}\op(\cE_{\al_2,\al_1}^\bu)_{23}\bigr)[-1]
\nonumber\\
&\cong(\bar\Pi^v_{\dot{\mathcal M}_{\al_1}}\!\t\!\bar\Pi^v_{\dot{\mathcal M}_{\al_2}})^*
\bigl(\De_{\M_{\al_1}}^*(\cE_{\al_1,\al_1}^\bu)[-1]\op\De_{\M_{\al_2}}^*(\cE_{\al_2,\al_2}^\bu)[-1]\bigr)_{\wt=0}
\nonumber\\
&\qquad\op\bigl((\bar\Pi^v_{\dot{\mathcal M}_{\al_1}}\t\bar\Pi^v_{\dot{\mathcal M}_{\al_2}})^*(\cE_{\al_1,\al_2}^\bu[-1])\bigr)_{\wt=1}
\nonumber\\
&\qquad\op\bigl((\bar\Pi^v_{\dot{\mathcal M}_{\al_1}}\t\bar\Pi^v_{\dot{\mathcal M}_{\al_2}})^*(\si_{\al_1,\al_2}^*(\cE_{\al_2,\al_1}^\bu)[-1])\bigr)_{\wt=-1},
\label{co9eq24}
\ea
using \eq{co9eq23} in the first step and the following commutative diagram 
\e
\begin{gathered}
\xymatrix@C=190pt@R=15pt{ *+[r]{\M_{\al_1}\t\M_{\al_2}} \ar[r]_(0.6){\Phi_{\al_1,\al_2}} \ar[d]^{\De_{\M_{\al_1}}\t\De_{\M_{\al_2}}} & *+[l]{\M_\al} \ar[dd]_{\De_{\M_\al}} \\
*+[r]{\M_{\al_1}\t\M_{\al_1}\t\M_{\al_2}\t\M_{\al_2}} \ar[d]^{\si_{23}} \\
*+[r]{\M_{\al_1}\t\M_{\al_2}\t\M_{\al_1}\t\M_{\al_2}} \ar[r]^(0.6){\Phi_{\al_1,\al_2}\t\Phi_{\al_1,\al_2}} & *+[l]{\M_\al\t\M_\al!}
}	
\end{gathered}
\label{co9eq25}
\e
in the second, where $\si_{23}$ exchanges second and third factors. 

In the third step of \eq{co9eq24} we use \eq{co4eq3}--\eq{co4eq4}, where the subscripts `$ab$' on $(\cE_{\al_i,\al_j}^\bu)_{ab}$ refer to which factors $a,b=1,\ldots,4$ in the bottom left corner of \eq{co9eq25} the $\cE_{\al_i,\al_j}^\bu$ depends. In the fourth we work out $(\si_{23}\ci(\De_{\M_{\al_1}}\t\De_{\M_{\al_2}}))^*$ explicitly, and reorganize according to the $\bG_m$-weights of each term, indicated by the subscript `$\wt=\cdots$'. Since $\th\in\bG_m$ acts on $E$ by $\th\,\id_{E_1}+\id_{E_2}$ by Proposition \ref{co9prop2}(c), and $(\th_1\id_{E_1},\th_2\id_{E_2})$ acts on $\cE^\bu\vert_{([E_1],[E_2])}$ by $\th_1^{-1}\th_2$ by \eq{co4eq5}--\eq{co4eq6}, the $\bG_m$-weight of $\cE_{\al_i,\al_j}^\bu$ in \eq{co9eq24} is $\de_{j2}-\de_{i2}$.

Write $\acM_{(\al,\bs 1)}^\ss(\ac\tau^0_{\bs\mu})_{\rho_3,\rho_4\ne 0}\subset\acM_{(\al,\bs 1)}^\ss(\ac\tau^0_{\bs\mu})$ for the open substack of $\C$-points $[E,\bs V,\bs\rho]$ with $\rho_3,\rho_4\ne 0$. This contains $\acM_{(\al,\bs 1)}^\ss(\ac\tau^0_{\bs\mu})_{\al_1,\al_2}$ by Proposition \ref{co9prop2}(c). On $\acM_{(\al,\bs 1)}^\ss(\ac\tau^0_{\bs\mu})_{\rho_3,\rho_4\ne 0}$ we can take $V_1=V_2=V_3=\C$ and $\rho_3,\rho_4=\id_\C$, so
\e
\bL_{\acM_{(\al,\bs 1)}^\ss(\ac\tau^0_{\bs\mu})/\dM_\al}\vert_{\acM_{(\al,\bs 1)}^\ss(\ac\tau^0_{\bs\mu})_{\rho_3,\rho_4\ne 0}}\cong (\ac\Pi^{v_3}_{\dot{\mathcal M}_\al})^*(\cV_{k_1,\al}^*\op\cV_{k_2,\al}^*).
\label{co9eq26}
\e
Applying $i_{\al_1,\al_2}^*$ and using \eq{co9eq13}, \eq{co9eq23}, and $\Phi_{\al_1,\al_2}^*(\cV_{k,\al}^*)\cong \cV_{k,\al_1}^*\bp\cV_{k,\al_2}^*$, we see that
\ea
&i_{\al_1,\al_2}^*(\bL_{\acM_{(\al,\bs 1)}^\ss(\ac\tau^0_{\bs\mu})/\dM_\al})\cong
(\bar\Pi^v_{\dot{\mathcal M}_{\al_1}})^*(\cV^*_{k_1,\al_1})_{\wt=0}\op(\bar\Pi^v_{\dot{\mathcal M}_{\al_2}})^*(\cV^*_{k_2,\al_2})_{\wt=0}\nonumber\\
&\qquad\qquad \op(\bar\Pi^v_{\dot{\mathcal M}_{\al_2}})^*(\cV^*_{k_1,\al_2})_{\wt=1} \op(\bar\Pi^v_{\dot{\mathcal M}_{\al_1}})^*(\cV^*_{k_2,\al_1})_{\wt=-1}.
\label{co9eq27}
\ea
Here the $\bG_m$-weights are as shown as $\th\in\bG_m$ acts on $E,V_1,V_2$ by $\th\,\id_{E_1}+\id_{E_2},\th\,\id_{V_1},\id_{V_2}$ by Proposition \ref{co9prop2}(c), so $\cV^*_{k_i,\al_j}$ has $\bG_m$-weight $\de_{i1}-\de_{j1}$, with $\de_{i1}$ coming from $\th\,\id_{V_1}$, and $-\de_{j1}$ from~$\th\,\id_{E_1}$.

Equations \eq{co9eq24} and \eq{co9eq27} describe four terms in \eq{co9eq22}, with their $\bG_m$-weights, where the only possible $\bG_m$-weights are $0,1,-1$. The terms $U_\al\ot\O[1]$ in \eq{co9eq22} also have $\bG_m$-weight 0. Thus we can split \eq{co9eq22} into the direct sum of three diagrams with $\bG_m$-weights $0,1,-1$. For weight 0 this gives:
\e
\text{\begin{footnotesize}$\displaystyle
\!\!\!\!\!\xymatrix@!0@C=41pt@R=45pt{
U_\al\!\ot\!\O[1] \ar@/_4.6pc/[ddddrr]_(0.8)0 \ar[dr] \ar[rr]^\id && U_\al\!\ot\!\O[1] \ar@/_2.2pc/[ddddrr]_(0.8)0 \ar[dr] \ar[rr] && 0 \ar[dr] \ar[rr]_(0.6){[+1]} && \qquad\qquad
\\
& {\!\!\!\!\begin{subarray}{l}\ts (\bar\Pi^v_{\dot{\mathcal M}_{\al_1}}\!\t\!\bar\Pi^v_{\dot{\mathcal M}_{\al_2}})^* \\ \ts
\bigl(\De_{\M_{\al_1}}^*(\cE_{\al_1,\al_1}^\bu)\op \\ \ts\De_{\M_{\al_2}}^*(\cE_{\al_2,\al_2}^\bu)\bigr)[-1]\!\!\!\!\!\end{subarray}} \ar@/_1pc/[dddr]_(0.53){\begin{subarray}{l}
(\bar\Pi^v_{\dot{\mathcal M}_{\al_1}})^* \\ (i^*(\bL_i))\ci\th_{\al_1} \\ 
{}\op(\bar\Pi^v_{\dot{\mathcal M}_{\al_2}})^* \\ (i^*(\bL_i))\ci\th_{\al_2}\end{subarray}\!\!\!\!\!\!\!\!\!\!\!\!\!\!\!\!\!} \ar[dr] \ar[rr] && 
{\begin{subarray}{l}\ts \;\> i_{\al_1}^*(\bL_{\bs\baM{}_{(\al_1,1)}^\pl}) \\ \ts {}\bp i_{\al_2}^*(\bL_{\bs\baM{}_{(\al_2,1)}^\pl}) \end{subarray}} \ar@/_1pc/[dddr]^(0.65){\!\!\!\!\!\!\!\!\!\begin{subarray}{l} i_{\al_1}^*(\bL_i)\bp \\ i_{\al_2}^*(\bL_i) \end{subarray}} \ar[dr] \ar[rr] && {\begin{subarray}{l}\ts \; (\bar\Pi^v_{\dot{\mathcal M}_{\al_1}})^*(\cV^*_{k_1,\al_1}) \\ \ts \op(\bar\Pi^v_{\dot{\mathcal M}_{\al_2}})^*(\cV^*_{k_2,\al_2})\end{subarray}} \ar[dr]_(0.4)\id \ar[rr]^(0.6){[+1]} && \qquad\qquad
\\
&& {\begin{subarray}{l}\ts \ac\imath_{\dot{\mathcal M}_\al}^* \\ \ts (i^*(\bL_{\bs\dM^\red_\al}))^{\bG_m}\end{subarray}} \ar[dr]^(0.7){[+1]} \ar[dd]  \ar[rr] && {\begin{subarray}{l}\ts \quad i_{\al_1,\al_2}^* \\ \ts (\bL_{\bs\acM_{(\al,\bs 1)}^\ss(\ac\tau^0_{\bs\mu})})^{\bG_m}\end{subarray}} \ar[dr]^(0.7){[+1]} \ar[dd]^(0.7){i_{\al_1,\al_2}^*(\bL_i)^{\bG_m}} \ar[rr] && {\begin{subarray}{l}\ts \; (\bar\Pi^v_{\dot{\mathcal M}_{\al_1}})^*(\cV^*_{k_1,\al_1}) \\ \ts \op(\bar\Pi^v_{\dot{\mathcal M}_{\al_2}})^*(\cV^*_{k_2,\al_2})\end{subarray}} \ar[dr]^(0.7){[+1]} \ar[dd]_\id \ar[rr]^(0.6){[+1]} && \qquad\qquad
\\
&&&&&&&
\\
&& {\begin{subarray}{l}\ts \; (\bar\Pi^v_{\dot{\mathcal M}_{\al_1}})^*(\bL_{\dot{\mathcal M}_{\al_1}}) \\ \ts \op(\bar\Pi^v_{\dot{\mathcal M}_{\al_2}})^*(\bL_{\dot{\mathcal M}_{\al_2}})\end{subarray}} \ar[rr] && {\begin{subarray}{l}\ts \;\bL_{\baM_{(\al_1,1)}^\ss(\bar\tau^0_1)_{k_1}} \\ \ts \bp\bL_{\baM_{(\al_2,1)}^\ss(\bar\tau^0_1)_{k_2}} \end{subarray}}  \ar[rr] && {\begin{subarray}{l}\ts \; (\bar\Pi^v_{\dot{\mathcal M}_{\al_1}})^*(\cV^*_{k_1,\al_1}) \\ \ts \op(\bar\Pi^v_{\dot{\mathcal M}_{\al_2}})^*(\cV^*_{k_2,\al_2})\end{subarray}} \ar[rr]^(0.6){[+1]} && \qquad\qquad }
\!\!\!\!\!\!\!\!\!\!\!\!\!\!\!\!\!\!\!\!\!\!\!\!\!\!\!\!\!\!\!\!\!\!\!\!\!\!\!\!\!\!
$\end{footnotesize}}
\label{co9eq28}
\e

Here we have rewritten the second term on the second line using the fact that $(\bs\dM_{(\al,\bs 1)}^\pl)^{\bG_m}\cong \bs\baM{}_{(\al_1,1)}^\pl\t \bs\baM{}_{(\al_2,1)}^\pl$. Note that we cannot do the same for the second term on the third line because of the issue noted at the end of Proposition \ref{co9prop2}. We have also rewritten the morphisms from the second to the fourth row, using \eq{co5eq3} for the left hand morphism.

For $\bG_m$-weights $1,-1$ the diagrams yield distinguished triangles:
\ea
&\xymatrix@!0@C=85pt{
{\begin{subarray}{l}\ts (\bar\Pi^v_{\dot{\mathcal M}_{\al_1}}\!\t\!\bar\Pi^v_{\dot{\mathcal M}_{\al_2}})^* \\ \ts (\cE_{\al_1,\al_2}^\bu[-1]) \end{subarray}} \ar[r] & \cN_{\al_1,\al_2}^{+,\bu}  \ar[r] & (\bar\Pi^v_{\dot{\mathcal M}_{\al_2}})^*(\cV^*_{k_1,\al_2})  \ar[r]^(0.7){[+1]} &   }\!\!
\label{co9eq29}\\
&\xymatrix@!0@C=83pt{
{\begin{subarray}{l}\ts (\bar\Pi^v_{\dot{\mathcal M}_{\al_1}}\!\t\!\bar\Pi^v_{\dot{\mathcal M}_{\al_2}})^* \\ \ts (\si_{\al_1,\al_2}^*(\cE_{\al_2,\al_1}^\bu)[-1]) \end{subarray}}
\ar[r] & \cN_{\al_1,\al_2}^{-,\bu}  \ar[r] & (\bar\Pi^v_{\dot{\mathcal M}_{\al_1}})^*(\cV^*_{k_2,\al_1})  \ar[r]^(0.7){[+1]} &   }\!\!
\label{co9eq30}
\ea
Taking K-theory classes in \eq{co9eq29}--\eq{co9eq30} proves part (c).

It remains to prove (d). To do this we must compare the obstruction theories on $\acM_{(\al,\bs 1)}^\ss(\ac\tau^0_{\bs\mu})_{\al_1,\al_2}$ and $\baM_{(\al_1,1)}^\ss(\bar\tau^0_1)_{k_1}\t\baM_{(\al_2,1)}^\ss(\bar\tau^0_1)_{k_2}$ under the isomorphism \eq{co9eq13}. As $\al_1,\al_2,\al=\al_1+\al_2$ lie in $C(\B)_\pe$, Assumption \ref{co5ass1}(f)(iv) says that $o_{\al_1}+o_{\al_2}\ge o_\al$, as claimed in part (d). Thus we may choose an identification $U_{\al_1}\op U_{\al_2}=U_\al\op V$, where~$V\cong \C^{o_{\al_1}+o_{\al_2}-o_\al}$. 

Consider the following modified version of~\eq{co9eq28}:
\e
\text{\begin{footnotesize}$\displaystyle
\!\!\!\!\!\xymatrix@!0@C=41pt@R=50pt{
{\begin{subarray}{l}\ts \; U_\al\!\ot\!\O[1] \\
\ts \op V\!\ot\!\O[1]\end{subarray}} \ar@/_4.6pc/[ddddrr]_(0.8)0 \ar[dr]^(0.4){e_1} \ar[rr]^\id && {\begin{subarray}{l}\ts \; U_\al\!\ot\!\O[1] \\
\ts \op V\!\ot\!\O[1]\end{subarray}} \ar@/_2.2pc/[ddddrr]_(0.8)0 \ar[dr]^(0.4){f_1} \ar[rr] && 0 \ar[dr] \ar[rr]_(0.6){[+1]} && \qquad\qquad
\\
& {\!\!\!\!\begin{subarray}{l}\ts (\bar\Pi^v_{\dot{\mathcal M}_{\al_1}}\!\t\!\bar\Pi^v_{\dot{\mathcal M}_{\al_2}})^* \\ \ts
\bigl(\De_{\M_{\al_1}}^*(\cE_{\al_1,\al_1}^\bu)\op \\ \ts\De_{\M_{\al_2}}^*(\cE_{\al_2,\al_2}^\bu)\bigr)[-1]\!\!\!\!\! \\
\ts \qquad \op V\ot\O[1]\end{subarray}} \ar@/_1pc/[dddr]_(0.53){\begin{subarray}{l}
(\bar\Pi^v_{\dot{\mathcal M}_{\al_1}})^* \\ (i^*(\bL_i))\ci\th_{\al_1} \\ 
{}\op(\bar\Pi^v_{\dot{\mathcal M}_{\al_2}})^* \\ (i^*(\bL_i))\ci\th_{\al_2} \\ \qquad\op 0\end{subarray}\!\!\!\!\!\!\!\!\!\!\!\!\!\!\!\!\!} \ar[dr] \ar[rr] && 
{\begin{subarray}{l}\ts \;\> i_{\al_1}^*(\bL_{\bs\baM{}_{(\al_1,1)}^\pl}) \\ \ts {}\bp i_{\al_2}^*(\bL_{\bs\baM{}_{(\al_2,1)}^\pl}) \\
\ts \qquad \op V\ot\O[1] \end{subarray}} \ar@/_1pc/[dddr]^(0.65){\!\!\!\!\!\!\!\!\!\begin{subarray}{l} i_{\al_1}^*(\bL_i)\bp \\ i_{\al_2}^*(\bL_i) \\ \quad \op 0\end{subarray}} \ar[dr] \ar[rr] && {\begin{subarray}{l}\ts \; (\bar\Pi^v_{\dot{\mathcal M}_{\al_1}})^*(\cV^*_{k_1,\al_1}) \\ \ts \op(\bar\Pi^v_{\dot{\mathcal M}_{\al_2}})^*(\cV^*_{k_2,\al_2})\end{subarray}} \ar[dr]_(0.4)\id \ar[rr]^(0.6){[+1]} && \qquad\qquad
\\
&& {\!\!\!\!\!\!\!\ac\imath_{\dot{\mathcal M}_\al}^*(i^*(\bL_{\bs\dM^\red_\al}))^{\bG_m}} \ar[dr]^(0.7){[+1]} \ar[dd]  \ar[rr] && {\begin{subarray}{l}\ts \quad i_{\al_1,\al_2}^* \\ \ts (\bL_{\bs\acM_{(\al,\bs 1)}^\ss(\ac\tau^0_{\bs\mu})})^{\bG_m}\end{subarray}} \ar[dr]^(0.7){[+1]} \ar[dd]^(0.7){i_{\al_1,\al_2}^*(\bL_i)^{\bG_m}} \ar[rr] && {\begin{subarray}{l}\ts \; (\bar\Pi^v_{\dot{\mathcal M}_{\al_1}})^*(\cV^*_{k_1,\al_1}) \\ \ts \op(\bar\Pi^v_{\dot{\mathcal M}_{\al_2}})^*(\cV^*_{k_2,\al_2})\end{subarray}} \ar[dr]^(0.7){[+1]} \ar[dd]_\id \ar[rr]^(0.6){[+1]} && \qquad\qquad
\\
&&&&&&&
\\
&& {\begin{subarray}{l}\ts \; (\bar\Pi^v_{\dot{\mathcal M}_{\al_1}})^*(\bL_{\dot{\mathcal M}_{\al_1}}) \\ \ts \op(\bar\Pi^v_{\dot{\mathcal M}_{\al_2}})^*(\bL_{\dot{\mathcal M}_{\al_2}})\end{subarray}} \ar[rr] && {\begin{subarray}{l}\ts \;\bL_{\baM_{(\al_1,1)}^\ss(\bar\tau^0_1)_{k_1}} \\ \ts \bp\bL_{\baM_{(\al_2,1)}^\ss(\bar\tau^0_1)_{k_2}} \end{subarray}}  \ar[rr] && {\begin{subarray}{l}\ts \; (\bar\Pi^v_{\dot{\mathcal M}_{\al_1}})^*(\cV^*_{k_1,\al_1}) \\ \ts \op(\bar\Pi^v_{\dot{\mathcal M}_{\al_2}})^*(\cV^*_{k_2,\al_2})\end{subarray}} \ar[rr]^(0.6){[+1]} && \qquad\qquad }
\!\!\!\!\!\!\!\!\!\!\!\!\!\!\!\!\!\!\!\!\!\!\!\!\!\!\!\!\!\!\!\!\!\!\!\!\!\!\!\!\!\!
$\end{footnotesize}}
\label{co9eq31}
\e
Here we have taken the sum of \eq{co9eq28} with $V\ot\O[1],\id_{V\ot\O[1]}$ in the top left parallelogram, so \eq{co9eq31} still commutes, with distinguished rows and diagonals. The morphisms $e_1,f_1$ are labelled for later use.

Next consider the diagram:
\e
\text{\begin{footnotesize}$\displaystyle
\!\!\!\!\!\xymatrix@!0@C=42pt@R=55pt{
{\begin{subarray}{l}\ts \; U_{\al_1}\!\ot\!\O[1] \\
\ts \op U_{\al_2}\!\ot\!\O[1]\end{subarray}} \ar@/_4.6pc/[ddddrr]_(0.8)0 \ar@{.>}[dr]^(0.4){e_2} \ar[rr]^\id && {\begin{subarray}{l}\ts \; U_{\al_1}\!\ot\!\O[1] \\
\ts \op U_{\al_2}\!\ot\!\O[1]\end{subarray}} \ar@/_2.2pc/[ddddrr]_(0.8)0 \ar@{.>}[dr]^(0.4){f_2} \ar[rr] && 0 \ar[dr] \ar[rr]_(0.6){[+1]} && \qquad\qquad
\\
& {\!\!\!\!\!\!\!\!\!\!\!\!\begin{subarray}{l}\ts (\bar\Pi^v_{\dot{\mathcal M}_{\al_1}}\!\t\!\bar\Pi^v_{\dot{\mathcal M}_{\al_2}})^* \\ \ts
\bigl(\De_{\M_{\al_1}}^*(\cE_{\al_1,\al_1}^\bu)\op \\ \ts\De_{\M_{\al_2}}^*(\cE_{\al_2,\al_2}^\bu)\bigr)[-1]\!\!\!\!\! \\
\ts \qquad \op V\ot\O[1]\end{subarray}} \ar@/_1pc/[dddr]_(0.6){\begin{subarray}{l}
(\bar\Pi^v_{\dot{\mathcal M}_{\al_1}})^* \\ (i^*(\bL_i))\ci\th_{\al_1} \\ 
{}\op(\bar\Pi^v_{\dot{\mathcal M}_{\al_2}})^* \\ (i^*(\bL_i))\ci\th_{\al_2} \\ \qquad\op 0\end{subarray}\!\!\!\!\!\!\!\!\!\!\!\!\!\!\!\!\!} \ar@{.>}[dr] \ar[rr] && 
{\begin{subarray}{l}\ts \;\> i_{\al_1}^*(\bL_{\bs\baM{}_{(\al_1,1)}^\pl}) \\ \ts {}\bp i_{\al_2}^*(\bL_{\bs\baM{}_{(\al_2,1)}^\pl}) \\
\ts \qquad \op V\ot\O[1] \end{subarray}} \ar@/_1pc/[dddr]^(0.65){\!\!\!\!\!\!\!\!\!\begin{subarray}{l} i_{\al_1}^*(\bL_i)\bp \\ i_{\al_2}^*(\bL_i) \\ \quad \op 0\end{subarray}} \ar@{.>}[dr] \ar[rr] && {\begin{subarray}{l}\ts \; (\bar\Pi^v_{\dot{\mathcal M}_{\al_1}})^*(\cV^*_{k_1,\al_1}) \\ \ts \op(\bar\Pi^v_{\dot{\mathcal M}_{\al_2}})^*(\cV^*_{k_2,\al_2})\end{subarray}} \ar[dr]_(0.4)\id \ar[rr]^(0.6){[+1]} && \qquad\qquad
\\
&& {\!\!\!\!\!\!\!\!\!\!\!\!\!\!\!\!\!\!\!\!\!\!\!\!\!\!\!\!\!\!\!\!\!\!\!\!\begin{subarray}{l}\ts  (\bar\Pi^v_{\dot{\mathcal M}_{\al_1}})^*(i^*(\bL_{\bs\dM^\red_{\al_1}})) \\ \ts \op(\bar\Pi^v_{\dot{\mathcal M}_{\al_2}})^*(i^*(\bL_{\bs\dM^\red_{\al_2}})) \\ \ts \quad \op V\ot\O[1]\end{subarray}} \ar@{.>}[dr]^(0.7){[+1]} \ar@{.>}[dd]  \ar@{.>}[rr] && {\begin{subarray}{l}\ts \;\> i_{\al_1}^*(\bL_{\bs\baM_{(\al_1,1)}^\ss(\bar\tau^0_1)_{k_1}}) \\ \ts {}\bp i_{\al_2}^*(\bL_{\bs\baM_{(\al_2,1)}^\ss(\bar\tau^0_1)_{k_2}}) \\
\ts \qquad \op V\ot\O[1] \end{subarray}} \ar@{.>}[dr]^(0.7){[+1]} \ar@{.>}[dd]^(0.7){\begin{subarray}{l} i_{\al_1}^*(\bL_i)\bp \\ i_{\al_2}^*(\bL_i) \\ \quad \op 0\end{subarray}} \ar@{.>}[rr] && {\begin{subarray}{l}\ts (\bar\Pi^v_{\dot{\mathcal M}_{\al_1}})^* \\ \ts (\cV^*_{k_1,\al_1}) \\ \ts \op(\bar\Pi^v_{\dot{\mathcal M}_{\al_2}})^* \\ \ts (\cV^*_{k_2,\al_2})\end{subarray}} \ar[dr]^(0.7){[+1]} \ar[dd]_\id \ar@{.>}[rr]^(0.6){[+1]} && \qquad\qquad
\\
&&&&&&&
\\
&& {\begin{subarray}{l}\ts \; (\bar\Pi^v_{\dot{\mathcal M}_{\al_1}})^*(\bL_{\dot{\mathcal M}_{\al_1}}) \\ \ts \op(\bar\Pi^v_{\dot{\mathcal M}_{\al_2}})^*(\bL_{\dot{\mathcal M}_{\al_2}})\end{subarray}} \ar[rr] && {\begin{subarray}{l}\ts \;\bL_{\baM_{(\al_1,1)}^\ss(\bar\tau^0_1)_{k_1}} \\ \ts \bp\bL_{\baM_{(\al_2,1)}^\ss(\bar\tau^0_1)_{k_2}} \end{subarray}}  \ar[rr] && {\begin{subarray}{l}\ts \; (\bar\Pi^v_{\dot{\mathcal M}_{\al_1}})^*(\cV^*_{k_1,\al_1}) \\ \ts \op(\bar\Pi^v_{\dot{\mathcal M}_{\al_2}})^*(\cV^*_{k_2,\al_2})\end{subarray}} \ar[rr]^(0.6){[+1]} && \qquad\qquad }
\!\!\!\!\!\!\!\!\!\!\!\!\!\!\!\!\!\!\!\!\!\!\!\!\!\!\!\!\!\!\!\!\!\!\!\!\!\!\!\!\!\!\!\!\!\!\!$\end{footnotesize}}
\label{co9eq32}
\e
This is the analogue of \eq{co9eq21} for the construction of the obstruction theory on $\baM_{(\al_1,1)}^\ss(\bar\tau^0_1)_{k_1}\t\baM_{(\al_2,1)}^\ss(\bar\tau^0_1)_{k_2}$, direct summed with $V\ot\O[1],\id_{V\ot\O[1]}$ in the centre left parallelogram. The morphisms $e_2,f_2$ are labelled for later use.

Here the morphisms `$\ra$' in \eq{co9eq32} are the same as in \eq{co9eq31}, but the morphisms `$\dashra$' may be different, and the left two objects in the third rows are different. In particular, the second and fourth rows of \eq{co9eq31} and \eq{co9eq32} are the same. The first rows of \eq{co9eq31} and \eq{co9eq32} also agree under~$U_{\al_1}\op U_{\al_2}=U_\al\op V$.

Working on $\baM_{(\al_1,1)}^\ss(\bar\tau^0_1)_{k_1}\t\baM_{(\al_2,1)}^\ss(\bar\tau^0_1)_{k_2}\t\bA^1$, where we write $y$ for the coordinate on $\bA^1$, consider the commuting diagram in $D_{\qcoh(\cdots)}$ with distinguished rows and diagonals:
\e
\text{\begin{footnotesize}$\displaystyle
\!\!\!\!\!\xymatrix@!0@C=41pt@R=45pt{
(U_\al\!\op\! V)\!\ot\!\O[1] \ar@/_4.6pc/[ddddrr]_(0.8)0 \ar[dr]^(0.4){(1-y)e_1+ye_2} \ar[rr]^\id && (U_\al\!\op\! V)\!\ot\!\O[1] \ar@/_2.2pc/[ddddrr]_(0.8)0 \ar[dr]^(0.4){(1-y)f_1+yf_2} \ar[rr] && 0 \ar[dr] \ar[rr]_(0.6){[+1]} && \qquad\qquad
\\
& {\!\!\!\!\!\!\!\!\!\!\!\!\begin{subarray}{l}\ts (\bar\Pi^v_{\dot{\mathcal M}_{\al_1}}\!\t\!\bar\Pi^v_{\dot{\mathcal M}_{\al_2}})^* \\ \ts
\bigl(\De_{\M_{\al_1}}^*(\cE_{\al_1,\al_1}^\bu)\op \\ \ts\De_{\M_{\al_2}}^*(\cE_{\al_2,\al_2}^\bu)\bigr)[-1]\!\!\!\!\! \\
\ts \qquad \op V\ot\O[1]\end{subarray}} \ar@/_1pc/[dddr]_(0.52){\begin{subarray}{l}
(\bar\Pi^v_{\dot{\mathcal M}_{\al_1}})^* \\ (i^*(\bL_i))\ci\th_{\al_1} \\ 
{}\op(\bar\Pi^v_{\dot{\mathcal M}_{\al_2}})^* \\ (i^*(\bL_i))\ci\th_{\al_2} \\ \qquad\op 0\end{subarray}\!\!\!\!\!\!\!\!\!\!\!\!\!\!\!\!\!} \ar@{.>}[dr] \ar[rr] && 
{\begin{subarray}{l}\ts \;\> i_{\al_1}^*(\bL_{\bs\baM{}_{(\al_1,1)}^\pl}) \\ \ts {}\bp i_{\al_2}^*(\bL_{\bs\baM{}_{(\al_2,1)}^\pl}) \\
\ts \qquad \op V\ot\O[1] \end{subarray}} \ar@/_1pc/[dddr]^(0.65){\!\!\!\!\!\!\!\!\!\begin{subarray}{l} i_{\al_1}^*(\bL_i)\bp \\ i_{\al_2}^*(\bL_i) \\ \quad \op 0\end{subarray}} \ar@{.>}[dr] \ar[rr] && {\begin{subarray}{l}\ts \; (\bar\Pi^v_{\dot{\mathcal M}_{\al_1}})^*(\cV^*_{k_1,\al_1}) \\ \ts \op(\bar\Pi^v_{\dot{\mathcal M}_{\al_2}})^*(\cV^*_{k_2,\al_2})\end{subarray}} \ar[dr]_(0.4)\id \ar[rr]^(0.6){[+1]} && \qquad\qquad
\\
&& \cF^\bu \ar@{.>}[dr]^(0.7){[+1]} \ar@{.>}[dd]^(0.6){\phi}  \ar@{.>}[rr] && \cG^\bu \ar@{.>}[dr]^(0.7){[+1]} \ar@{.>}[dd]^(0.6){\psi} \ar@{.>}[rr] && {\begin{subarray}{l}\ts (\bar\Pi^v_{\dot{\mathcal M}_{\al_1}})^*(\cV^*_{k_1,\al_1}) \\ \ts \op(\bar\Pi^v_{\dot{\mathcal M}_{\al_2}})^*(\cV^*_{k_2,\al_2})\end{subarray}} \ar[dr]^(0.7){[+1]} \ar[dd]_\id \ar@{.>}[rr]^(0.6){[+1]} && \qquad\qquad
\\
&&&&&&&
\\
&& {\begin{subarray}{l}\ts \; (\bar\Pi^v_{\dot{\mathcal M}_{\al_1}})^*(\bL_{\dot{\mathcal M}_{\al_1}}) \\ \ts \op(\bar\Pi^v_{\dot{\mathcal M}_{\al_2}})^*(\bL_{\dot{\mathcal M}_{\al_2}})\end{subarray}} \ar[rr] && {\begin{subarray}{l}\ts \;\bL_{\baM_{(\al_1,1)}^\ss(\bar\tau^0_1)_{k_1}} \\ \ts \bp\bL_{\baM_{(\al_2,1)}^\ss(\bar\tau^0_1)_{k_2}} \end{subarray}}  \ar[rr] && {\begin{subarray}{l}\ts \; (\bar\Pi^v_{\dot{\mathcal M}_{\al_1}})^*(\cV^*_{k_1,\al_1}) \\ \ts \op(\bar\Pi^v_{\dot{\mathcal M}_{\al_2}})^*(\cV^*_{k_2,\al_2})\end{subarray}} \ar[rr]^(0.6){[+1]} && \qquad\qquad }
\!\!\!\!\!\!\!\!\!\!\!\!\!\!\!\!\!\!\!\!\!\!\!\!\!\!\!\!\!\!\!\!\!\!\!\!\!\!\!\!\!\!\!\!\!\!\!$\end{footnotesize}}
\label{co9eq33}
\e
Here the two top left diagonal morphisms $(1-y)(\cdots)+y(\cdots)$ interpolate between the two top left diagonal morphisms in \eq{co9eq31} at $y=0$ and \eq{co9eq32} at $y=1$, after the identification $U_{\al_1}\op U_{\al_2}=U_\al\op V$. We define $\cF^\bu$ and $\cG^\bu$ to be the cones of these morphisms. As the left morphism from first to fourth row is zero, the left morphism from second to fourth row factors via $\cF^\bu$, giving the morphism $\phi$. The morphism $\psi$ exists by the same argument. The remaining morphisms `$\dashra$' can be filled in by properties of triangulated categories.

We have constructed a morphism
\e
\begin{split}
\psi:\cG^\bu\longra &\,\bL_{\baM_{(\al_1,1)}^\ss(\bar\tau^0_1)_{k_1}}\bp\bL_{\baM_{(\al_2,1)}^\ss(\bar\tau^0_1)_{k_2}}\\
=&\,\bL_{\baM_{(\al_1,1)}^\ss(\bar\tau^0_1)_{k_1}\t\baM_{(\al_2,1)}^\ss(\bar\tau^0_1)_{k_2}\t\bA^1/\bA^1}.
\end{split}
\label{co9eq34}
\e
We would like to say that \eq{co9eq34} is a {\it relative obstruction theory\/} for $\pi_{\bA^1}:\baM_{(\al_1,1)}^\ss(\bar\tau^0_1)_{k_1}\t\baM_{(\al_2,1)}^\ss(\bar\tau^0_1)_{k_2}\t\bA^1\ra\bA^1$, in the sense of Definition \ref{co2def3}. In fact this may not be quite true. It is clear from \eq{co9eq33} that $\cG^\bu$ is a perfect complex, but not that $\cG^\bu$ is perfect in the interval $[-1,1]$, or that $h^i(\psi)$ is an isomorphism for $i\ge 0$ and surjective for~$i=-1$. 

As these are open conditions, they hold on the complement of a minimal closed $\C$-substack $C\subset\baM_{(\al_1,1)}^\ss(\bar\tau^0_1)_{k_1}\t\baM_{(\al_2,1)}^\ss(\bar\tau^0_1)_{k_2}\t\bA^1$. Over $y=0$ in $\bA^1$, equation \eq{co9eq33} coincides with \eq{co9eq31}, so \eq{co9eq34} coincides with $i_{\al_1,\al_2}^*(\bL_i)^{\bG_m}$ in \eq{co9eq31}, which is an obstruction theory. Hence $C$ does not intersect the hyperplane $y=0$ in $\baM_{(\al_1,1)}^\ss(\bar\tau^0_1)_{k_1}\t\baM_{(\al_2,1)}^\ss(\bar\tau^0_1)_{k_2}\t\bA^1$. Similarly, over $y=1$ in $\bA^1$, \eq{co9eq33} coincides with \eq{co9eq32}, so \eq{co9eq34} coincides with 
\ea
(i_{\al_1}^*(\bL_i)&\!\bp\! i_{\al_2}^*(\bL_i))\!\op\! 0:(i_{\al_1}^*(\bL_{\bs\baM_{(\al_1,1)}^\ss(\bar\tau^0_1)_{k_1}}) \!\bp\! i_{\al_2}^*(\bL_{\bs\baM_{(\al_2,1)}^\ss(\bar\tau^0_1)_{k_2}}))\!\op\! (V\!\ot\!\O[1]) 
\nonumber\\
&\longra\bL_{\baM_{(\al_1,1)}^\ss(\bar\tau^0_1)_{k_1}}\bp\bL_{\baM_{(\al_2,1)}^\ss(\bar\tau^0_1)_{k_2}},
\label{co9eq35}
\ea
which is an obstruction theory. Thus $C$ does not intersect~$y=1$.

As $\baM_{(\al_1,1)}^\ss(\bar\tau^0_1)_{k_1}\t\baM_{(\al_2,1)}^\ss(\bar\tau^0_1)_{k_2}$ is proper, and $C$ is closed, we see that $\pi_{\bA^1}\vert_C:C\ra\bA^1$ is proper, so its image is a Zariski closed set $D\subset\bA^1$, which does not contain 0 or 1 as $C$ does not intersect $y=0$ or $y=1$. So restricting \eq{co9eq34} to $\baM_{(\al_1,1)}^\ss(\bar\tau^0_1)_{k_1}\t\baM_{(\al_2,1)}^\ss(\bar\tau^0_1)_{k_2}\t(\bA^1\sm D)$ gives a relative obstruction theory for the proper morphism
\e
\pi_{\bA^1\sm D}:\baM_{(\al_1,1)}^\ss(\bar\tau^0_1)_{k_1}\t\baM_{(\al_2,1)}^\ss(\bar\tau^0_1)_{k_2}\t(\bA^1\sm D)\longra(\bA^1\sm D).
\label{co9eq36}
\e
For each $\C$-point $y\in\bA^1\sm D$, restricting \eq{co9eq34} gives a perfect obstruction theory on the fibre $\baM_{(\al_1,1)}^\ss(\bar\tau^0_1)_{k_1}\t\baM_{(\al_2,1)}^\ss(\bar\tau^0_1)_{k_2}$ of \eq{co9eq36} at $y$, which is isomorphic to $i_{\al_1,\al_2}^*(\ac\psi_{(\al,\bs 1)})^{\bG_m}$ when $y=0$ and to \eq{co9eq35} when $y=1$. 

As $\bA^1\sm D$ is connected, Theorem \ref{co2thm1}(iii) implies that the virtual class of the obstruction theory $i_{\al_1,\al_2}^*(\bL_i)^{\bG_m}$ on $\acM_{(\al,\bs 1)}^\ss(\ac\tau^0_{\bs\mu})_{\al_1,\al_2}$ is identified under \eq{co9eq13} with the virtual class of \eq{co9eq35} on $\baM_{(\al_1,1)}^\ss(\bar\tau^0_1)_{k_1}\t\baM_{(\al_2,1)}^\ss(\bar\tau^0_1)_{k_2}$. If $o_{\al_1}+o_{\al_2}=o_\al$, so that $V=0$, equation \eq{co9eq35} is the product obstruction theory on $\baM_{(\al_1,1)}^\ss(\bar\tau^0_1)_{k_1}\t\baM_{(\al_2,1)}^\ss(\bar\tau^0_1)_{k_2}$, and the first case of \eq{co9eq15} follows. If $o_{\al_1}+o_{\al_2}>o_\al$, so that $V\ne 0$, then by choosing a surjective linear map $V^*\twoheadrightarrow\C$ we can use Theorem \ref{co2thm1}(iv) to show that the virtual class of the obstruction theory \eq{co9eq35} is zero, proving the second case of \eq{co9eq15}, and completing the proposition.
\end{proof}

\begin{dfn}
\label{co9def3}
Work in the situation above. Define an equivariant cohomology class
\ea
\eta&\in H^{2(o_\al+\la_{k_1}(\al)+\la_{k_2}(\al)-\chi(\al,\al)-1)}_{\bG_m}(\acM_{(\al,\bs 1)}^\ss(\ac\tau^0_{\bs\mu}))\quad\text{by}
\label{co9eq37}\\
\eta&=(-1)^{\la_{k_1}(\al)+\la_{k_2}(\al)} c_{\la_{k_1}(\al)+\la_{k_2}(\al)-2}\bigl([\bL_{\acM_{(\al,\bs 1)}^\ss(\ac\tau^0_{\bs\mu})/\M^\pl_\al}]\!-\![\ac L_{-1,-1,2}^\pl]\bigr).
\nonumber
\ea
Here $\Pi_{\M_\al^{\smash{\pl}}}:\acM_{(\al,\bs 1)}^\ss(\ac\tau^0_{\bs\mu})\ra\M^\pl_\al$ is smooth of dimension $\la_{k_1}(\al)+\la_{k_2}(\al)-1,$ so its relative cotangent complex $\bL_{\acM_{(\al,\bs 1)}^\ss(\ac\tau^0_{\bs\mu})/\M^\pl_\al}$ is a vector bundle of this rank on $\acM_{(\al,\bs 1)}^\ss(\ac\tau^0_{\bs\mu})$. Both $\bL_{\acM_{(\al,\bs 1)}^\ss(\ac\tau^0_{\bs\mu})/\M^\pl_\al}$ and\/ $L_{-1,-1,2}^\pl$ are $\bG_m$-equivariant, so have classes $[\bL_{\acM_{(\al,\bs 1)}^\ss(\ac\tau^0_{\bs\mu})/\M^\pl_\al}],[L_{-1,-1,2}^\pl]$ in $K^0_{\bG_m}(\Perf_{\acM_{(\al,\bs 1)}^\ss(\ac\tau^0_{\bs\mu})})$. Thus the equivariant Chern class $c_{\cdots}(\cdots)$ makes sense. As $\Pi_{\M_\al^{\smash{\pl}}}$ is $\bG_m$-equivariant, we have a pullback map $\Pi_{\M_\al^{\smash{\pl}}}^*: H^*(\M_\al^\pl)\ra H^*_{\bG_m}(\acM_{(\al,\bs 1)}^\ss(\ac\tau^0_{\bs\mu}))$. Hence \eq{co9eq37} is well defined.

We apply Corollary \ref{co2cor2} to the proper algebraic space $\acM_{(\al,\bs 1)}^\ss(\ac\tau^0_{\bs\mu})$ with its $\bG_m$-action and $\bG_m$-equivariant obstruction theory $\bL_i:i^*(\bL_{\bs\acM_{(\al,\bs 1)}^\ss(\ac\tau^0_{\bs\mu})})\ra \bL_{\acM_{(\al,\bs 1)}^\ss(\ac\tau^0_{\bs\mu})}$, and cohomology class $\eta$, and with $\Pi_{\M_\al^{\smash{\pl}}}:\acM_{(\al,\bs 1)}^\ss(\ac\tau^0_{\bs\mu})\ra\M^\pl_\al$ in place of $f:X\ra Y$. Using the notation of Propositions \ref{co9prop2}--\ref{co9prop3}, \eq{co2eq26} becomes:
\ea
&(-1)^{\rank\cN_{\rho_4=0}^\bu}\Res_z\bigl\{
(\Pi_{\M_\al^{\smash{\pl}}}\ci i_{\rho_4=0})_*\bigl([\acM_{(\al,\bs 1)}^\ss(\ac\tau^0_{\bs\mu})_{\rho_4=0}]_\virt
\nonumber\\
&\qquad\qquad\qquad \cap(e(\cN_{\rho_4=0}^\bu)^{-1}\cup i_{\rho_4=0}^*(\eta))\bigr)\bigr\}
\nonumber\\
\begin{split}
&+(-1)^{\rank\cN_{\rho_3=0}^\bu}\Res_z\bigl\{(\Pi_{\M_\al^{\smash{\pl}}}\ci i_{\rho_3=0})_*\bigl([\acM_{(\al,\bs 1)}^\ss(\ac\tau^0_{\bs\mu})_{\rho_3=0}]_\virt\\
&\qquad\qquad\qquad\cap (e(\cN_{\rho_3=0}^\bu)^{-1}\cup i_{\rho_3=0}^*(\eta))\bigr)\bigr\}
\end{split}
\label{co9eq38}\\
&+\sum_{\!\!\begin{subarray}{l} \al_1,\al_2\in
C(\B)_\pe:\; \al_1+\al_2=\al, \\ 
\tau(\al_1)=\tau(\al_2),\;\M_{\al_i}^\ss(\tau)\ne\es,\; i=1,2\end{subarray}\!\!\!\!\!\!\!\!\!\!\!\!\!\!\!\!\!\!\!\!\!\!\!\!\!\!\!\!\!\!\!\!\!\!\!\!\!\!\!\!\!\!\!\!\!\!\!\!\!\!\!}
\begin{aligned}[t] 
&(-1)^{\rank\cN_{\al_1,\al_2}^\bu}\Res_z\bigl\{(\Pi_{\M_\al^{\smash{\pl}}}\ci i_{\al_1,\al_2})_*\bigl([\acM_{(\al,\bs 1)}^\ss(\ac\tau^0_{\bs\mu})_{\al_1,\al_2}]_\virt\\
&\qquad\qquad\qquad\qquad \cap(e(\cN_{\al_1,\al_2}^\bu)^{-1}\cup i_{\al_1,\al_2}^*(\eta))\bigr)\bigr\}=0.
\end{aligned}
\nonumber
\ea
\end{dfn}

\begin{prop}
\label{co9prop4}
In\/ {\rm\eq{co9eq38},} using the notation of\/ \eq{co9eq2} we have
\ea
\begin{split}
&(-1)^{\rank\cN_{\rho_4=0}^\bu}\Res_z\bigl\{
(\Pi_{\M_\al^{\smash{\pl}}}\ci i_{\rho_4=0})_*\bigl([\acM_{(\al,\bs 1)}^\ss(\ac\tau^0_{\bs\mu})_{\rho_4=0}]_\virt\\
&\qquad\qquad\qquad \cap(e(\cN_{\rho_4=0}^\bu)^{-1}\cup i_{\rho_4=0}^*(\eta))\bigr)\bigr\}
=\la_{k_2}(\al)\,\bar\Up_\al^{k_1}(\tau),
\end{split}
\label{co9eq39}\\
\begin{split}
&(-1)^{\rank\cN_{\rho_3=0}^\bu}\Res_z\bigl\{(\Pi_{\M_\al^{\smash{\pl}}}\ci i_{\rho_3=0})_*\bigl([\acM_{(\al,\bs 1)}^\ss(\ac\tau^0_{\bs\mu})_{\rho_3=0}]_\virt\\
&\qquad\qquad\qquad\cap (e(\cN_{\rho_3=0}^\bu)^{-1}\cup i_{\rho_3=0}^*(\eta))\bigr)\bigr\}
=-\la_{k_1}(\al)\,\bar\Up_\al^{k_2}(\tau). 
\end{split}
\label{co9eq40}
\ea
\end{prop}

\begin{proof} To prove \eq{co9eq39}, writing $c_{\la_{k_1}(\al)+\la_{k_2}(\al)-2}(\cdots)=c_\top(\cdots)$, first note that
\ea
i_{\rho_4=0}^*(\eta)&=i_{\rho_4=0}^*\bigl((-1)^{\la_{k_1}(\al)+\la_{k_2}(\al)}c_\top\bigl([\bL_{\acM_{(\al,\bs 1)}^\ss(\ac\tau^0_{\bs\mu})/\M^\pl_\al}]-[\ac L_{-1,-1,2}^\pl]\bigr)\bigr)
\nonumber\\
&=(-1)^{\la_{k_1}(\al)+\la_{k_2}(\al)}c_\top\bigl([\bL_{\acM_{(\al,\bs 1)}^\ss(\ac\tau^0_{\bs\mu})_{\rho_4=0}/\M^\pl_\al}\op \ac L_{0,-1,1}^\pl\vert_{\acM_{\rho_4=0}}]
\nonumber\\
&\qquad\qquad\qquad\qquad\qquad\qquad -[\ac L_{-1,-1,2}^\pl\vert_{\acM_{\rho_4=0}}]\bigr)
\nonumber\\
&=(-1)^{\la_{k_1}(\al)+\la_{k_2}(\al)}c_\top(\bL_{\acM_{(\al,\bs 1)}^\ss(\ac\tau^0_{\bs\mu})_{\rho_4=0}/\M^\pl_\al})
\nonumber\\
&=c_\top(\bT_{\acM_{(\al,\bs 1)}^\ss(\ac\tau^0_{\bs\mu})_{\rho_4=0}/\M^\pl_\al}),
\label{co9eq41}
\ea
using \eq{co9eq37} in the first step, \eq{co9eq16} and \eq{co9eq18} in the second, $\ac L_{0,-1,1}^\pl\vert_{\acM_{\rho_4=0}}\cong \ac L_{-1,-1,2}^\pl\vert_{\acM_{\rho_4=0}}$ as $\ac L^\pl_{1,0,-1}\vert_{\acM_{\rho_4=0}}$ is trivial by Proposition \ref{co9prop2}(a) in the third, and $\rank \bT_{\acM_{\rho_4=0}/\M^\pl_\al}=\la_{k_1}(\al)+\la_{k_2}(\al)-2$ in the fourth.

We now prove equation \eq{co9eq39} by
\ea
&(-1)^{\rank\cN_{\rho_4=0}^\bu}\Res_z\bigl\{
(\Pi_{\M_\al^{\smash{\pl}}}\ci i_{\rho_4=0})_*\bigl([\acM_{(\al,\bs 1)}^\ss(\ac\tau^0_{\bs\mu})_{\rho_4=0}]_\virt
\nonumber\\
&\qquad\qquad\qquad \cap(e(\cN_{\rho_4=0}^\bu)^{-1}\cup i_{\rho_4=0}^*(\eta))\bigr)\bigr\}
\nonumber\\
&=-\Res_z\bigl\{(\Pi_{\M_\al^{\smash{\pl}}}\ci i_{\rho_4=0})_*\bigl([\acM_{(\al,\bs 1)}^\ss(\ac\tau^0_{\bs\mu})_{\rho_4=0}]_\virt
\nonumber\\
&\qquad\qquad\qquad \cap(e(\cN_{\rho_4=0}^\bu)^{-1}\cup c_\top(\bT_{\acM_{(\al,\bs 1)}^\ss(\ac\tau^0_{\bs\mu})_{\rho_4=0}/\M^\pl_\al}))\bigr)\bigr\}
\nonumber\\
&=(\Pi_{\M_\al^{\smash{\pl}}}\ci i_{\rho_4=0})_*\bigl([\acM_{(\al,\bs 1)}^\ss(\ac\tau^0_{\bs\mu})_{\rho_4=0}]_\virt\cap c_\top(\bT_{\acM_{(\al,\bs 1)}^\ss(\ac\tau^0_{\bs\mu})_{\rho_4=0}/\M^\pl_\al})\bigr)
\allowdisplaybreaks
\nonumber \\
&=(\bar\Pi_{\M_\al^{\smash{\pl}}}\ci\Pi_{\rho_4=0})_*\bigl([\acM_{(\al,\bs 1)}^\ss(\ac\tau^0_{\bs\mu})_{\rho_4=0}]_\virt\cap \bigl[\Pi_{\rho_4=0}^*\bigl(c_\top(\bT_{\baM_{(\al,1)}^\ss(\bar\tau^0_1)_{k_1}/\M^\pl_\al})\bigr)\nonumber\\
&\qquad\qquad \cup c_\top(\bT_{\acM_{(\al,\bs 1)}^\ss(\ac\tau^0_{\bs\mu})_{\rho_4=0}/\baM_{(\al,1)}^\ss(\bar\tau^0_1)_{k_1}})\bigr]\bigr)
\allowdisplaybreaks
\nonumber \\
&=\la_{k_2}(\al)\,(\bar\Pi_{\M_\al^{\smash{\pl}}})_*\bigl([\baM_{(\al,1)}^\ss(\bar\tau^0_1)_{k_1}]_\virt\cap c_\top(\bT_{\baM_{(\al,1)}^\ss(\bar\tau^0_1)_{k_1}/\M^\pl_\al})\bigr)
\nonumber \\
&=\la_{k_2}(\al)\,\bar\Up_\al^{k_1}(\tau).
\label{co9eq42}
\ea
Here in the first step we use \eq{co9eq41} and $\rank\cN_{\rho_4=0}^\bu=1$. For the second, as $\cN_{\rho_4=0}^\bu\cong \ac L^\pl_{0,-1,1}$ with $\bG_m$-weight $-1$ by Proposition \ref{co9prop3}(a) we have  $e(\cN_{\rho_4=0}^\bu)=-1\cdot z+c_1(\ac L^\pl_{0,-1,1})$, where the coefficient $-1$ of $z$ is the weight $-1$ of the $\bG_m$-action, so $e(\cN_{\rho_4=0}^\bu)^{-1}=-z^{-1}+O(z^{<-1})$. But $[\acM_{(\al,\bs 1)}^\ss(\ac\tau^0_{\bs\mu})_{\rho_4=0}]_\virt$ and $c_\top(\bT_{\acM_{(\al,\bs 1)}^\ss(\ac\tau^0_{\bs\mu})_{\rho_4=0}/\M^\pl_\al})$ are both $\bG_m$-invariant, and so are independent of $z$. Thus taking $\Res_z\{\cdots\}$ gives a factor $-1$ from $e(\cN_{\rho_4=0}^\bu)^{-1}$, and the second step of \eq{co9eq42} follows.

In the third step we use the right hand triangle of smooth morphisms in \eq{co9eq16}, and substitute in $\Pi_{\M_\al^{\smash{\pl}}}\ci i_{\rho_4=0}=\bar\Pi_{\M_\al^{\smash{\pl}}}\ci\Pi_{\rho_4=0}$ as in \eq{co9eq16}, where $\bar\Pi_{\M_\al^{\smash{\pl}}}:\baM_{(\al,1)}^\ss(\bar\tau^0_1)_{k_1}\ra\M_\al^\pl$ is the projection. In the fourth we apply Corollary \ref{co2cor1} to the derived version $\bs\Pi_{\rho_4=0}$ of the smooth morphism $\Pi_{\rho_4=0}$ in \eq{co9eq11}, which as for $\Pi_{\rho_4=0}$ in Proposition \ref{co9prop2}(a) is a smooth fibration with fibre $\CP^{\la_{k_2}(\al)-1}$. In the fifth we use equation \eq{co9eq2}. The proof of \eq{co9eq40} is similar, using Propositions \ref{co9prop2}(b) and \ref{co9prop3}(b). The sign difference in \eq{co9eq40} comes from $e(\cN_{\rho_3=0}^\bu)=1\cdot z+c_1(\ac L^\pl_{-1,0,1})$, since $\cN_{\rho_3=0}^\bu$ has $\bG_m$-weight~$1$.
\end{proof}

The next proposition is crucial to our programme, as it shows how the Lie bracket on $\check H_{\rm even}(\M^\pl)$ from \S\ref{co43} is related to properties of Behrend--Fantechi virtual classes under $\bG_m$-localization, as in~\eq{co2eq26}.

\begin{prop}
\label{co9prop5}
In equation {\rm\eq{co9eq38},} for each\/ $\al_1,\al_2$ we have
\ea
&(-1)^{\rank\cN_{\al_1,\al_2}^\bu}\Res_z\bigl\{(\Pi_{\M_\al^{\smash{\pl}}}\ci i_{\al_1,\al_2})_*\bigl([\acM_{(\al,\bs 1)}^\ss(\ac\tau^0_{\bs\mu})_{\al_1,\al_2}]_\virt
\label{co9eq43}\\
&\cap(e(\cN_{\al_1,\al_2}^\bu)^{-1}\cup i_{\al_1,\al_2}^*(\eta))\bigr)\bigr\}
=\begin{cases} \bigl[\bar\Up_{\al_1}^{k_1}(\tau),\bar\Up_{\al_2}^{k_2}(\tau)\bigr], & o_{\al_1}+o_{\al_2}=o_\al, \\
0, & o_{\al_1}+o_{\al_2}>o_\al, \end{cases}
\nonumber
\ea
using the Lie bracket on $\check H_{\rm even}(\M^\pl)$ from\/ {\rm\S\ref{co43},} where $o_{\al_1}+o_{\al_2}\ge o_\al$ by Assumption\/ {\rm\ref{co5ass1}(f)(iv)}.
\end{prop}

\begin{proof} The second case of \eq{co9eq43}, when $o_{\al_1}+o_{\al_2}>o_\al$, is immediate from \eq{co9eq15}. So suppose that $o_{\al_1}+o_{\al_2}=o_\al$. Writing $H^*_{\bG_m}(\acM_{(\al,\bs 1)}^\ss(\ac\tau^0_{\bs\mu})_{\al_1,\al_2})=H^*_{\bG_m}(\baM_{(\al_1,1)}^\ss(\bar\tau^0_1)_{k_1}\t\baM_{(\al_2,1)}^\ss(\bar\tau^0_1)_{k_2})$ by \eq{co9eq13}, and including $\bG_m$-weights in the notation, we compute
\ea
&i_{\al_1,\al_2}^*(\eta)=(-1)^{\la_{k_1}(\al)+\la_{k_2}(\al)}c_{\la_{k_1}(\al)+\la_{k_2}(\al)-2}\bigl([i_{\al_1,\al_2}^*(\bL_{\acM_{(\al,\bs 1)}^\ss(\ac\tau^0_{\bs\mu})/\M^\pl_\al})]
\nonumber\\
&\qquad\qquad\qquad\qquad\qquad\qquad\qquad\qquad -[i_{\al_1,\al_2}^*(\ac L_{-1,-1,2}^\pl)]\bigr)
\nonumber\\
&=(-1)^{\la_{k_1}(\al)+\la_{k_2}(\al)}c_{\la_{k_1}(\al)+\la_{k_2}(\al)-2}\bigl((\bar\Pi^v_{\dot{\mathcal M}_{\al_1}})^*(\cV^*_{k_1,\al_1})_{\wt=0}
\nonumber\\
&\quad\op(\bar\Pi^v_{\dot{\mathcal M}_{\al_2}})^*(\cV^*_{k_2,\al_2})_{\wt=0}\op(\bar\Pi^v_{\dot{\mathcal M}_{\al_2}})^*(\cV^*_{k_1,\al_2})_{\wt=1}\op
(\bar\Pi^v_{\dot{\mathcal M}_{\al_1}})^*(\cV^*_{k_2,\al_1})_{\wt=-1}\bigr)
\allowdisplaybreaks
\nonumber\\
&=(-1)^{\la_{k_1}(\al)+\la_{k_2}(\al)}c_{\la_{k_1}(\al)+\la_{k_2}(\al)-2}\bigl((((\bar\Pi^v_{\dot{\mathcal M}_{\al_1}})^*(\cV_{k_1,\al_1})/\rho_1(\O))^*
\nonumber\\
&\quad \!\op\!((\bar\Pi^v_{\dot{\mathcal M}_{\al_2}})^*(\cV_{k_2,\al_2})/\rho_2(\O))^*)_{\wt=0}
\nonumber\\
&\quad\op(\bar\Pi^v_{\dot{\mathcal M}_{\al_2}})^*(\cV^*_{k_1,\al_2})_{\wt=1}\op
(\bar\Pi^v_{\dot{\mathcal M}_{\al_1}})^*(\cV^*_{k_2,\al_1})_{\wt=-1}\bigr)
\nonumber\\
&=(-1)^{\la_{k_1}(\al)+\la_{k_2}(\al)}c_\top\bigl(\bL_{\baM_{(\al_1,1)}^\ss(\bar\tau^0_1)_{k_1}/\dM^\pl_{\al_1}}{}_{\wt=0}\bigr)
\nonumber\\
&\quad\cup c_\top\bigl(\bL_{\baM_{(\al_2,1)}^\ss(\bar\tau^0_1)_{k_2}/\dM^\pl_{\al_2}}{}_{\wt=0}\bigr)
\nonumber\\
&\quad\cup c_\top\bigl((\bar\Pi^v_{\dot{\mathcal M}_{\al_2}})^*(\cV^*_{k_1,\al_2})_{\wt=1}\bigr)\cup
c_\top\bigl((\bar\Pi^v_{\dot{\mathcal M}_{\al_1}})^*(\cV^*_{k_2,\al_1})_{\wt=-1}\bigr).
\label{co9eq44}
\ea

Here in the first step we use \eq{co9eq37}, and in the second we use \eq{co5eq26}, \eq{co9eq27}, and that $i_{\al_1,\al_2}^*(\ac L_{-1,-1,2}^\pl)$ is trivial by Proposition \ref{co9prop2}(c). In the third we use that  $(\bar\Pi^v_{\dot{\mathcal M}_{\al_i}})^*(\cV_{k_i,\al_i})$ has a trivial line subbundle $\rho_i(\O)$ for $i=1,2$, so quotienting by $\rho_i(\O)$ does not change the Chern class. In the fourth step  we use that in $c_{\la_{k_1}(\al)+\la_{k_2}(\al)-2}(\cdots)$, the `$\cdots$' is the sum of four vector bundles of total rank $\la_{k_1}(\al)+\la_{k_2}(\al)-2$, and so splits as the product of the top Chern classes of the four, and the isomorphisms~$((\bar\Pi^v_{\dot{\mathcal M}_{\al_i}})^*(\cV_{k_i,\al_i})/\rho_i(\O))^*\cong\bL_{\baM_{(\al_i,1)}^\ss(\bar\tau^0_1)_{k_i}/\dM^\pl_{\al_i}}$.

We can now prove the first case of \eq{co9eq43}, in the long equation
\ea
&(-1)^{\rank\cN_{\al_1,\al_2}^\bu}\!\Res_z\bigl\{(\Pi_{\M_\al^{\smash{\pl}}}\ci i_{\al_1,\al_2})^{\bG_m}_*\bigl([\acM_{(\al,\bs 1)}^\ss(\ac\tau^0_{\bs\mu})_{\al_1,\al_2}]_\virt
\nonumber\\
&\qquad\qquad \cap (e(\cN_{\al_1,\al_2}^\bu)^{-1}\cup i_{\al_1,\al_2}^*(\eta))\bigr)\bigr\}
\nonumber\\
&=(-1)^{\la_{k_1}(\al_2)+\la_{k_2}(\al_1)+\chi(\al_1,\al_2)+\chi(\al_2,\al_1)}\Res_z\bigl\{(\Pi_{\M_\al^{\smash{\pl}}}\ci i_{\al_1,\al_2})_*\nonumber\\
&\qquad \bigl(([\baM_{(\al_1,1)}^\ss(\bar\tau^0_1)_{k_1}]_\virt\bt
 [\baM_{(\al_2,1)}^\ss(\bar\tau^0_1)_{k_2}]_\virt)\cap (e(\cN_{\al_1,\al_2}^\bu)^{-1}\!\cup\! i_{\al_1,\al_2}^*(\eta))\bigr)\bigr\}
\allowdisplaybreaks
\nonumber\\
&=(-1)^{\la_{k_1}(\al_2)+\la_{k_2}(\al_1)+\chi(\al_1,\al_2)+\chi(\al_2,\al_1)}(-1)^{\la_{k_1}(\al)+\la_{k_2}(\al)}\cdot{}
\nonumber\\
&\;\>\Res_z\bigl\{(\Pi_{\M_\al^{\smash{\pl}}}\ci i_{\al_1,\al_2})^{\bG_m}_*\bigl(([\baM_{(\al_1,1)}^\ss(\bar\tau^0_1)_{k_1}]_\virt\bt[\baM_{(\al_2,1)}^\ss(\bar\tau^0_1)_{k_2}]_\virt)
\nonumber\\
&\;\>\cap\bigl[e((\bar\Pi^v_{\dot{\mathcal M}_{\al_1}}\t\bar\Pi^v_{\dot{\mathcal M}_{\al_2}})^*(\cE^\bu_{\al_1,\al_2})_{\wt=1})\cup e((\bar\Pi^v_{\dot{\mathcal M}_{\al_2}})^*(\cV^*_{k_1,\al_2})_{\wt=1})^{-1}
\nonumber\\
&\;\>\cup e((\bar\Pi^v_{\dot{\mathcal M}_{\al_1}}\t\bar\Pi^v_{\dot{\mathcal M}_{\al_2}})^*(\si_{\al_1,\al_2}^*(\cE_{\al_2,\al_1}^\bu))_{\wt=-1})
\cup e((\bar\Pi^v_{\dot{\mathcal M}_{\al_1}})^*(\cV^*_{k_2,\al_1})_{\wt=-1})^{-1}
\nonumber\\
&\;\> \cup c_\top\bigl(\bL_{\baM_{(\al_1,1)}^\ss(\bar\tau^0_1)_{k_1}/\M^\pl_\al}{}_{\wt=0}\bigr)
\cup c_\top\bigl(\bL_{\baM_{(\al_2,1)}^\ss(\bar\tau^0_1)_{k_2}/\M^\pl_\al}{}_{\wt=0}\bigr)
\nonumber\\
&\;\>\cup c_\top\bigl((\bar\Pi^v_{\dot{\mathcal M}_{\al_2}})^*(\cV^*_{k_1,\al_2})_{\wt=1}\bigr)\cup
c_\top\bigl((\bar\Pi^v_{\dot{\mathcal M}_{\al_1}})^*(\cV^*_{k_2,\al_1})_{\wt=-1}\bigr)
\bigr]\bigr)\bigr\}
\allowdisplaybreaks
\nonumber\\
&=(-1)^{\la_{k_1}(\al_1)+\la_{k_2}(\al_2)+\chi(\al_1,\al_2)+\chi(\al_2,\al_1)}\Res_z\bigl\{(\Pi^\pl_\al\!\ci\!\Phi_{\al_1,\al_2}\!\ci\!(\bar\Pi^v_{\dot{\mathcal M}_{\al_1}}\!\!\!\t\!\bar\Pi^v_{\dot{\mathcal M}_{\al_2}}))^{\bG_m}_*
\nonumber\\
&\;\>\bigl(([\baM_{(\al_1,1)}^\ss(\bar\tau^0_1)_{k_1}]_\virt
\bt[\baM_{(\al_2,1)}^\ss(\bar\tau^0_1)_{k_2}]_\virt)
\nonumber\\
&\;\>\cap\bigl[(\bar\Pi^v_{\dot{\mathcal M}_{\al_1}}\!\!\!\t\!\bar\Pi^v_{\dot{\mathcal M}_{\al_2}})^*\bigl(e(\cE^\bu_{\al_1,\al_2}{}_{\wt=1})\cup e(\si_{\al_1,\al_2}^*(\cE^\bu_{\al_2,\al_1}){}_{\wt=-1})
\nonumber\\
&\;\>\cup c_\top\bigl(\bL_{\baM_{(\al_1,1)}^\ss(\bar\tau^0_1)_{k_1}/\dM^\pl_{\al_1}}{}_{\wt=0}\bigr)
\cup c_\top\bigl(\bL_{\baM_{(\al_2,1)}^\ss(\bar\tau^0_1)_{k_2}/\dM^\pl_{\al_2}}{}_{\wt=0}\bigr)\bigr]\bigr)\bigr\}
\allowdisplaybreaks
\nonumber\\
&=(-1)^{\chi(\al_1,\al_2)}\!\Res_z\bigl\{(\Pi^\pl_\al\!\ci\!\Phi_{\al_1,\al_2})^{\bG_m}_*\bigl(\bigl[(\bar\Pi^v_{\dot{\mathcal M}_{\al_1}})_*\bigl([\baM_{(\al_1,1)}^\ss(\bar\tau^0_1)_{k_1}]_\virt\nonumber\\
&\;\>\cap c_\top(\bT_{\baM_{(\al_1,1)}^\ss(\bar\tau^0_1)_{k_1}/\dM^\pl_{\al_1}})\bigr)
\bt(\bar\Pi^v_{\dot{\mathcal M}_{\al_2}})_*\bigl([\baM_{(\al_2,1)}^\ss(\bar\tau^0_1)_{k_2}]_\virt
\nonumber\\
&\;\>\cap c_\top(\bT_{\baM_{(\al_2,1)}^\ss(\bar\tau^0_1)_{k_2}/\dM^\pl_{\al_2}})\bigr)\bigr]
\cap \bigl[e(\cE^\bu_{\al_1,\al_2}{}_{\wt=1})\cup e(\si_{\al_1,\al_2}^*(\cE^\bu_{\al_2,\al_1})^\vee_{\wt=1})\bigr]\bigr)\bigr\}
\allowdisplaybreaks
\nonumber\\
&=\Res_z\bigl\{(-1)^{\chi(\al_1,\al_2)}\bigl((\Pi^\pl_\al\ci\Phi_{\al_1,\al_2})^{\bG_m}_*\bigl(
\bigl[\hat\Up_{\al_1}^{k_1}(\tau)\bt\hat\Up_{\al_2}^{k_2}(\tau)\bigr]
\nonumber\\
&\qquad\qquad \cap \bigl[e(\cE^\bu_{\al_1,\al_2}\op\si_{\al_1,\al_2}^*(\cE^\bu_{\al_2,\al_1})^\vee)_{\wt=1}\bigr]\bigr)\bigr\}
\allowdisplaybreaks
\nonumber\\
&=\Res_z\bigl\{(-1)^{\chi(\al_1,\al_2)}\bigl((\Pi^\pl_\al\ci\Phi_{\al_1,\al_2})^{\bG_m}_*\bigl(
\bigl[\hat\Up_{\al_1}^{k_1}(\tau)\bt\hat\Up_{\al_2}^{k_2}(\tau)\bigr]
\nonumber\\
&\qquad\qquad \cap \bigl[\ts \sum_{i\ge 0} z^{\chi(\al_1,\al_2)+\chi(\al_2,\al_1)-i}c_i(\cE^\bu_{\al_1,\al_2}\op\si_{\al_1,\al_2}^*(\cE^\bu_{\al_2,\al_1})^\vee)
\bigr]\bigr)\bigr\}
\allowdisplaybreaks
\nonumber\\
&=\Res_z\bigl\{(-1)^{\chi(\al_1,\al_2)}\bigl((\Pi^\pl_\al\ci\Phi_{\al_1,\al_2}\ci(\Psi_{\al_1}\t\id_{\M_{\al_2}}))_*
\nonumber\\
&\qquad\qquad\bigl(
\ts\sum_{j\ge 0}z^jt^j\bt\bigl(
\bigl[\hat\Up_{\al_1}^{k_1}(\tau)\bt\hat\Up_{\al_2}^{k_2}(\tau)\bigr]
\nonumber\\
&\qquad\qquad \cap \bigl[\ts \sum_{i\ge 0} z^{\chi(\al_1,\al_2)+\chi(\al_2,\al_1)-i}c_i(\cE^\bu_{\al_1,\al_2}\op\si_{\al_1,\al_2}^*(\cE^\bu_{\al_2,\al_1})^\vee)
\bigr]\bigr)\bigr)\bigr\}
\allowdisplaybreaks
\nonumber\\
&=(\Pi^\pl_\al)_*\bigl(\Res_z\bigl\{(-1)^{\chi(\al_1,\al_2)}\bigl[\ts \sum_{i,j\ge 0} z^{\chi(\al_1,\al_2)+\chi(\al_2,\al_1)-i+j}\bigl(\Phi_{\al_1,\al_2}\ci
\nonumber\\
&\;\> 
(\Psi_{\al_1}\!\t\!\id_{\M_{\al_2}})\bigr)_*\bigl[t^j\!\bt\!\bigl((\hat\Up_{\al_1}^{k_1}(\tau)\!\bt\!\hat\Up_{\al_2}^{k_2}(\tau))\!\cap\! c_i(\cE^\bu_{\al_1,\al_2}\!\op\!\si_{\al_1,\al_2}^*(\cE^\bu_{\al_2,\al_1})^\vee)\bigr)\bigr]\bigr\}\bigr)
\allowdisplaybreaks
\nonumber\\
&=(\Pi^\pl_\al)_*\bigl(\Res_z\bigl\{Y(\hat\Up_{\al_1}^{k_1}(\tau),z)\hat\Up_{\al_2}^{k_2}(\tau)\bigr\}\bigr)=(\Pi^\pl_\al)_*\bigl\{\bigl(\hat\Up_{\al_1}^{k_1}(\tau)\bigr)_0\bigl(\hat\Up_{\al_2}^{k_2}(\tau)\bigr)\bigr\}
\allowdisplaybreaks
\nonumber\\
&=\bigl[(\Pi^\pl_{\al_1})_*(\hat\Up_{\al_1}^{k_1}(\tau)),(\Pi^\pl_{\al_2})_*(\hat\Up_{\al_2}^{k_2}(\tau))\bigr]
=\bigl[\bar\Up_{\al_1}^{k_1}(\tau),\bar\Up_{\al_2}^{k_2}(\tau)\bigr].
\label{co9eq45}
\ea
Here when we write $(\cdots)_*^{\bG_m}$ it is to emphasize that we are doing a pushforward along the 1-morphism $\cdots$ in $\bG_m$-{\it equivariant\/} homology, with respect to a {\it particular\/ $\bG_m$-equivariant structure\/} on $\cdots$. 

In the first step of \eq{co9eq45} we compute $(-1)^{\rank\cN_{\al_1,\al_2}^\bu}$ using \eq{co9eq14}, make the identification \eq{co9eq13}, and use the first case of \eq{co9eq15}. In the second step we substitute in \eq{co9eq14} and \eq{co9eq44} and use the multiplicative property of Euler classes. In the third step we cancel four terms using for $i=1,2$ 
\begin{equation*}
e((\bar\Pi^v_{\dot{\mathcal M}_{\al_i}})^*(\cV^*_{k_{3-i},\al_i}))=c_{\la_{k_{3-i}}(\al_i)}((\bar\Pi^v_{\dot{\mathcal M}_{\al_i}})^*(\cV^*_{k_{3-i},\al_i})),
\end{equation*}
as $(\bar\Pi^v_{\dot{\mathcal M}_{\al_i}})^*(\cV^*_{k_{3-i},\al_i})$ is a vector bundle of rank $\la_{k_{3-i}}(\al_i)$, and we substitute in
\begin{equation*}
\Pi_{\M_\al^{\smash{\pl}}}\ci i_{\al_1,\al_2}=\Pi^\pl_\al\ci\ac\Pi^{v_3}_{\dot{\mathcal M}_\al}\ci i_{\al_1,\al_2}=\Pi^\pl_\al\ci\ac\imath_{\dot{\mathcal M}_\al}=\Pi^\pl_\al\ci\Phi_{\al_1,\al_2}\ci(\bar\Pi^v_{\dot{\mathcal M}_{\al_1}}\t\bar\Pi^v_{\dot{\mathcal M}_{\al_2}}),
\end{equation*}
using \eq{co9eq23}. In the fourth step we use $e(\cD^\bu)=(-1)^{\rank\cD^\bu}e((\cD^\bu)^\vee),$ and $(-1)^{\la_{k_i}(\al_i)-1}c_\top(\bL_{\baM_{(\al_i,1)}^\ss(\bar\tau^0_1)_{k_i}/\dM^\pl_{\al_i}})=c_\top(\bT_{\baM_{(\al_i,1)}^\ss(\bar\tau^0_1)_{k_i}/\dM^\pl_{\al_i}})$ for $i=1,2$, and properties of (co)homology. In the fifth, in a similar way to \eq{co9eq2} we write
\begin{equation*}
\hat\Up_{\al_i}^{k_i}(\tau)=(\bar\Pi^v_{\dot{\mathcal M}_{\al_i}})_*\bigl([\baM^\ss_{(\al_i,1)}(\bar\tau^0_1)_{k_i}]_\virt\cap c_\top(\bT_{\baM^\ss_{(\al_i,1)}(\bar\tau^0_1)_{k_i}/\M^\pl_{\al_i}})\bigr)
\end{equation*}
for $i=1,2$, so that by \eq{co5eq26} and \eq{co9eq2} we have
\e
(\Pi_{\al_i}^\pl)_*(\hat\Up_{\al_i}^{k_i}(\tau))=\bar\Up_{\al_i}^{k_i}(\tau).
\label{co9eq46}
\e

For the sixth step of \eq{co9eq45}, as $\bG_m$ acts trivially on $\M_{\al_1}\t\M_{\al_2}$, so that $(\M_{\al_1}\t\M_{\al_2})/\bG_m\cong [*/\bG_m]\t\M_{\al_1}\t\M_{\al_2}$, we have
\e
\begin{split}
H^*_{\bG_m}(\M_{\al_1}\t\M_{\al_2})&\cong H^*([*/\bG_m]\t\M_{\al_1}\t\M_{\al_2})\\
&\cong H^*(\M_{\al_1}\t\M_{\al_2})\ot_\Q\Q[[z]],
\end{split}
\label{co9eq47}
\e
where $H^*([*/\bG_m])\cong\Q[[z]]$ as in \eq{co2eq2}. The term $e(\cE^\bu_{\al_1,\al_2}\op\si_{\al_1,\al_2}^* (\cE^\bu_{\al_2,\al_1})^\vee_{\wt=1})$ in the fifth step lies in the localized version $H^*(\M_{\al_1}\t\M_{\al_2})\ot_\Q\Q[[z]][z^{-1}]$ of \eq{co9eq47}. We expand it using
\e
e(\cD^\bu_{\wt=d})=\sum_{i\ge 0}(dz)^{\rank\cD^\bu-i}c_i(\cD^\bu)_{\text{non-eq.}},\qquad d\ne 0,
\label{co9eq48}
\e
where $c_i(\cD^\bu)_{\text{non-eq.}}$ is the Chern class in non-equivariant cohomology. 

For the seventh step of \eq{co9eq45}, we rewrite the $\bG_m$-equivariant pushforward $(\Pi^\pl_\al\ci\Phi_{\al_1,\al_2})^{\bG_m}_*$ in the sixth step in terms of non-equivariant homology. There is a subtlety here because although $\bG_m$ acts trivially on $\M_{\al_1}\t\M_{\al_2},\M_\al^\pl$, the correct $\bG_m$-equivariant structure on the 1-morphism $\Pi^\pl_\al\ci\Phi_{\al_1,\al_2}:\M_{\al_1}\t\M_{\al_2}\ra\M_\al^\pl$ is {\it not\/} the trivial one. This is a phenomenon to do with stacks, which form a 2-category, and does not happen for group actions on schemes.

To understand the $\bG_m$-equivariant structure on $\Pi^\pl_\al\ci\Phi_{\al_1,\al_2}$, observe that it allows us to push $\Pi^\pl_\al\ci\Phi_{\al_1,\al_2}$ down to a 1-morphism $(\M_{\al_1}\t\M_{\al_2})/\bG_m\ra\M_\al^\pl/\bG_m$. We claim this is given by the commutative diagram
\e
\begin{gathered}
\xymatrix@C=180pt@R=15pt{ *+[r]{[*/\bG_m]\!\t\!\M_{\al_1}\!\t\!\M_{\al_2}} \ar[r]_{\id_{[*/\bG_m]}\t\Psi_{\al_1}\t\id_{\M_{\al_2}}} \ar@{..>}[d]^{\id_{[*/\bG_m]}\t(\Pi^\pl_\al\ci\Phi_{\al_1,\al_2})^{\bG_m}} & *+[l]{[*/\bG_m]\!\t\!\M_{\al_1}\!\t\!\M_{\al_2}} \ar[d]_{\id_{[*/\bG_m]}\t\Phi_{\al_1,\al_2}} \\
*+[r]{[*/\bG_m]\t\M_\al^\pl} & *+[l]{[*/\bG_m]\t\M_\al,\!} \ar[l]_{\id_{[*/\bG_m]}\t\Pi^\pl_\al}  }\!\!\!	
\end{gathered}
\label{co9eq49}
\e
for $\Psi_{\al_1}$ in Assumption \ref{co4ass1}(c). To see this, note that at $[E_1\op E_2,\bs V,\bs\rho]$ in $\acM_{(\al,\bs 1)}^\ss(\ac\tau^0_{\bs\mu})_{\al_1,\al_2}$, $\th\in\bG_m$ acts on $E=E_1\op E_2$ by $\th\,\id_{E_1}+\id_{E_2}$ by Proposition \ref{co9prop2}(c). The 1-morphism $\id_{[*/\bG_m]}\t\Psi_{\al_1}\t\id_{\M_{\al_2}}$ in \eq{co9eq49} encodes this action, which gives the correct $\bG_m$-equivariance.

We convert \eq{co9eq49} into a commutative diagram in equivariant homology by
\e
\begin{gathered}
\xymatrix@C=160pt@R=20pt{ *+[r]{\begin{subarray}{l}\ts H_*^{\bG_m}(\M_{\al_1}\!\t\!\M_{\al_2})\cong \\ \ts H_*(\M_{\al_1}\!\t\!\M_{\al_2})\!\ot\!\Q[z]\end{subarray}} \ar[r]_{\begin{subarray}{l} (\Psi_{\al_1}\ci\id_{\M_{\al_2}})_* \\ (\sum_{j\ge 0}z^jt^j\bt-)\end{subarray}} \ar[d]^{(\Pi^\pl_\al\ci\Phi_{\al_1,\al_2})^{\bG_m}_*} & *+[l]{\begin{subarray}{l}\ts H_*^{\bG_m}(\M_{\al_1}\!\t\!\M_{\al_2})\cong \\ \ts H_*(\M_{\al_1}\!\t\!\M_{\al_2})\!\ot\!\Q[z]\end{subarray}} \ar[d]_{(\Phi_{\al_1,\al_2})_*\ot\id_{\Q[z]}} \\
*+[r]{\begin{subarray}{l}\ts H_*^{\bG_m}(\M_\al^\pl)\cong \\ \ts H_*(\M_\al^\pl)\ot\Q[z]\end{subarray}} & *+[l]{\begin{subarray}{l}\ts H_*^{\bG_m}(\M_\al)\cong \\ \ts H_*(\M_\al)\ot\Q[z],\!\end{subarray}} \ar[l]_{(\Pi^\pl_\al)_*\ot\id_{\Q[z]}}  }\!\!\!
\end{gathered}
\label{co9eq50}
\e
where we use the dual bases $(z^j)_{j=0}^\iy,(t^j)_{j=0}^\iy$ of $H^*([*/\bG_m]),H_*([*/\bG_m])$ from \eq{co2eq2}--\eq{co2eq3}. In the seventh step of \eq{co9eq45} we substitute in \eq{co9eq50}.

In the eighth step we rearrange, and in the ninth we recognize the definition \eq{co4eq10} of the state-field correspondence $Y(-,z)$ of the vertex algebra constructed in \S\ref{co42}. In the tenth step we use Definition \ref{co4def1}, in the eleventh the definition of the Lie bracket $[\,,\,]$ on $\check H_{\rm even}(\M^\pl)$ in \eq{co4eq1} and Theorem \ref{co4thm2}, and in the twelfth \eq{co9eq46}. This finally proves \eq{co9eq45}, and the proposition.
\end{proof}

Equations \eq{co9eq38}, \eq{co9eq39}, \eq{co9eq40}, and \eq{co9eq43} now imply:

\begin{cor}
\label{co9cor1}
Suppose\/ $(\tau,T,\le)\in\sS,$ $\al\in C(\B)_\pe$ and\/ $k_1,k_2\in K$ with\/ $\M_\al^\ss(\tau)\subseteq\M_{k_1,\al}^\pl\cap\M_{k_2,\al}^\pl$. Then in the Lie algebra $\check H_{\rm even}(\M^\pl)$ we have
\e
\begin{split}
0&=\la_{k_2}(\al)\,\bar\Up_\al^{k_1}(\tau)-\la_{k_1}(\al)\,\bar\Up_\al^{k_2}(\tau)+\!\!\!\!\!\!\sum_{\!\!\begin{subarray}{l} \al_1,\al_2\in
C(\B)_\pe:\\ 
\tau(\al_1)=\tau(\al_2),\;\M_{\al_i}^\ss(\tau)\ne\es,\; i=1,2,\\
\; \al_1+\al_2=\al, \; o_{\al_1}+o_{\al_2}=o_\al\end{subarray}\!\!\!\!\!\!\!\!\!\!\!\!\!\!\!\!\!\!\!\!\!\!\!\!\!\!\!\!\!\!\!\!\!\!\!\!\!\!\!\!\!\!\!\!\!\!\!\!\!\!\!\!\!\!\!\!\!\!}\bigl[\bar\Up_{\al_1}^{k_1}(\tau),\bar\Up_{\al_2}^{k_2}(\tau)\bigr],
\end{split}
\label{co9eq51}
\e
where there are only finitely many terms in the sum by Lemma\/ {\rm\ref{co9lem1}} with\/~$n=2$.
\end{cor}

\subsection{Proof of Theorem \ref{co5thm1}}
\label{co93}

Substituting \eq{co9eq3} for $k_1,k_2$ into each term in \eq{co9eq51} now gives:

\begin{cor}
\label{co9cor2}
In the Lie algebra\/ $\check H_{\rm even}(\M^\pl)$ we have
\ea 
&0=\!\!\!\!\!\!\!\!\sum_{\begin{subarray}{l}n\ge 1,\;\al_1,\ldots,\al_n\in
C(\B)_\pe:\\  
\tau(\al_i)=\tau(\al),\; \M_{\al_i}^\ss(\tau)\ne\es,\; \text{all $i,$} \\
\al_1+\cdots+\al_n=\al,\; o_{\al_1}+\cdots+o_{\al_n}=o_\al \end{subarray}}\,\, \begin{aligned}[t]
&\frac{(-1)^{n+1}\la_{k_2}(\al)\,\la_{k_1}(\al_1)}{n!}\cdot\bigl[\bigl[\cdots\bigl[[\M_{\al_1}^\ss(\tau)]^{k_1}_\inv, \\
&[\M_{\al_2}^\ss(\tau)]^{k_1}_\inv\bigr],\ldots\bigr],[\M_{\al_n}^\ss(\tau)]^{k_1}_\inv\bigr]\nonumber\\
\end{aligned}
\\
&-\!\!\!\!\!\sum_{\begin{subarray}{l}n\ge 1,\;\al_1,\ldots,\al_n\in
C(\B)_\pe:\\  
\tau(\al_i)=\tau(\al),\; \M_{\al_i}^\ss(\tau)\ne\es,\; \text{all $i,$} \\
\al_1+\cdots+\al_n=\al,\; o_{\al_1}+\cdots+o_{\al_n}=o_\al \end{subarray}} \,\,\begin{aligned}[t]
&\frac{(-1)^{n+1}\la_{k_1}(\al)\,\la_{k_2}(\al_1)}{n!}\cdot\bigl[\bigl[\cdots\bigl[[\M_{\al_1}^\ss(\tau)]^{k_2}_\inv, \\
&[\M_{\al_2}^\ss(\tau)]^{k_2}_\inv\bigr],\ldots\bigr],[\M_{\al_n}^\ss(\tau)]^{k_2}_\inv\bigr]\!\!\!\!\!\!\\
\end{aligned}
\nonumber\\
&+\!\!\!\!\!\sum_{\begin{subarray}{l} n>m\ge 1,\\ \al_1,\ldots,\al_n\in
C(\B)_\pe:\\  
\tau(\al_i)=\tau(\al),\\ \M_{\al_i}^\ss(\tau)\ne\es,\; \text{all $i,$} \\
\al_1+\cdots+\al_n=\al,\\ o_{\al_1}+\cdots+o_{\al_n}=o_\al \end{subarray}}\, \begin{aligned}[t]
&\frac{(-1)^n\,\la_{k_1}(\al_1)\,\la_{k_2}(\al_{m+1})}{m!\,(n-m)!}\cdot{}\\
&\bigl[\bigl[\bigl[\cdots\bigl[[\M_{\al_1}^\ss(\tau)]^{k_1}_\inv,[\M_{\al_2}^\ss(\tau)]^{k_1}_\inv\bigr],\ldots\bigr],[\M_{\al_m}^\ss(\tau)]^{k_1}_\inv\bigr],\\
&\bigl[\bigl[\cdots\bigl[[\M_{\al_{m+1}}^\ss(\tau)]^{k_2}_\inv,[\M_{\al_{m+2}}^\ss(\tau)]^{k_2}_\inv\bigr],\ldots\bigr],[\M_{\al_n}^\ss(\tau)]^{k_2}_\inv\bigr]\bigr].
\end{aligned}\!\!\!\!\!\!\!\!\!\!\!\!\!\!\!\!\!\!\!\!\!\!\!\!\!\!\!\!\!\!\!\!\!\!\!\!\!\!\!\!\!\!\!\!\!\!\!\!
\label{co9eq52}
\ea
\end{cor}

Theorem \ref{co5thm1} now follows from Proposition \ref{co9prop1} and the next proposition:

\begin{prop}
\label{co9prop6}
Suppose\/ $(\tau,T,\le)\in\sS,$ $\al\in C(\B)_\pe,$ and\/ $k_1,k_2\in K$ with\/ $\M_\al^\ss(\tau)\ab\subseteq\M_{k_1,\al}^\pl\cap\M_{k_2,\al}^\pl$. Then\/ $[\M_\al^\ss(\tau)]^{k_1}_\inv=[\M_\al^\ss(\tau)]^{k_2}_\inv$.
\end{prop}

\begin{proof} We will prove that $[\M_\al^\ss(\tau)]^{k_1}_\inv=[\M_\al^\ss(\tau)]^{k_2}_\inv$ by induction on $\rk\al=1,2,\ldots.$ Let $l\ge 0$ be given, and suppose by induction that for all $(\tau,T,\le)$ in $\sS$, $\al\in C(\B)_\pe$ and $k_1,k_2\in K$ with $\rk\al\le l$ and $\M_\al^\ss(\tau)\ab\subseteq\M_{k_1,\al}^\pl\cap\M_{k_2,\al}^\pl$ we have $[\M_\al^\ss(\tau)]^{k_1}_\inv=[\M_\al^\ss(\tau)]^{k_2}_\inv$. (This is vacuous when $l=0$ as~$\rk\al>0$.)

Let $(\tau,T,\le)\in\sS$, $\al\in C(\B)_\pe$ and $k_1,k_2\in K$ with $\rk\al=l+1$ and $\M_\al^\ss(\tau)\ab\subseteq\M_{k_1,\al}^\pl\cap\M_{k_2,\al}^\pl$. Consider equation \eq{co9eq52} for these $(\tau,T,\le),\al,k_1,k_2$. For all the terms in the sums in $\al_1,\ldots,\al_n$ with $n\ge 2$, Lemma \ref{co9lem1} shows that $\rk\al_i<\rk\al=l+1$ for $i=1,\ldots,n$, so $\rk\al_i\le l$ and $[\M_{\al_i}^\ss(\tau)]^{k_1}_\inv=[\M_{\al_i}^\ss(\tau)]^{k_2}_\inv$. Thus we may rewrite \eq{co9eq52} as
\ea
&\la_{k_1}(\al)\,\la_{k_2}(\al)\bigl([\M_\al^\ss(\tau)]^{k_2}_\inv\!-\![\M_\al^\ss(\tau)]^{k_1}_\inv\bigr)=\!\!\!\!\!\!\sum_{\begin{subarray}{l}n\ge 2,\;\al_1,\ldots,\al_n\in
C(\B)_\pe:\\  
\tau(\al_i)=\tau(\al),\; \M_{\al_i}^\ss(\tau)\ne\es,\; \text{all $i,$} \\
\al_1+\cdots+\al_n=\al,\; o_{\al_1}+\cdots+o_{\al_n}=o_\al \end{subarray}}\!\!\!\!\!\!\!\!\!\!\frac{(-1)^n}{(n!)^2}
\nonumber\\[-5pt]
&\sum_{\si\in S_n}\biggl\{
\begin{aligned}[t]
& -\la_{k_2}(\al)\,\la_{k_1}(\al_{\si(1)})\cdot{}\\
&\bigl[\bigl[\cdots\bigl[[\M_{\al_{\si(1)}}^\ss(\tau)]^{k_1}_\inv,
[\M_{\al_{\si(2)}}^\ss(\tau)]^{k_1}_\inv\bigr],\ldots\bigr],[\M_{\al_{\si(n)}}^\ss(\tau)]^{k_1}_\inv\bigr]\\
&+\la_{k_1}(\al)\,\la_{k_2}(\al_{\si(1)})\cdot{} \\
&\bigl[\bigl[\cdots\bigl[[\M_{\al_{\si(1)}}^\ss(\tau)]^{k_1}_\inv,
[\M_{\al_{\si(2)}}^\ss(\tau)]^{k_1}_\inv\bigr],\ldots\bigr],[\M_{\al_{\si(n)}}^\ss(\tau)]^{k_1}_\inv\bigr]\\
&+\sum_{m=1}^{n-1}\binom{n}{m}\la_{k_1}(\al_{\si(1)})\,\la_{k_2}(\al_{\si(m+1)})\cdot{}\\
&\bigl[\bigl[\bigl[\cdots\bigl[[\M_{\al_{\si(1)}}^\ss(\tau)]^{k_1}_\inv,[\M_{\al_{\si(2)}}^\ss(\tau)]^{k_1}_\inv\bigr],\ldots\bigr],[\M_{\al_{\si(m)}}^\ss(\tau)]^{k_1}_\inv\bigr],\\
&\bigl[\bigl[\cdots\bigl[[\M_{\al_{\si(m+1)}}^\ss(\tau)]^{k_1}_\inv,[\M_{\al_{\si(m+2)}}^\ss(\tau)]^{k_1}_\inv\bigr],\ldots\bigr],[\M_{\al_{\si(n)}}^\ss(\tau)]^{k_1}_\inv\bigr]\bigr]\biggr\}.
\end{aligned}\!\!\!\!\!\!\!\!\!\!\!\!\!\!\!\!\!
\label{co9eq53}
\ea
Here the left hand side of \eq{co9eq53} is the terms in $n=1$ from the first two lines of \eq{co9eq52}. For the remaining terms, we have replaced $[\M_{\al_i}^\ss(\tau)]^{k_2}_\inv$ by $[\M_{\al_i}^\ss(\tau)]^{k_1}_\inv$, and included a symmetrization over permutations of $\al_1,\ldots,\al_n$, which is cancelled by the extra factor of $1/n!$, with~$n!=\md{S_n}$.

We claim that in the universal enveloping algebra $U(\check H_{\rm even}(\M^\pl))$ we have
\ea
&\sum_{\si\in S_n}
\begin{aligned}[t]
&\la_k(\al_{\si(1)})\cdot{} \\
&\bigl[\bigl[\cdots\bigl[[\M_{\al_{\si(1)}}^\ss(\tau)]^{k_1}_\inv,
[\M_{\al_{\si(2)}}^\ss(\tau)]^{k_1}_\inv\bigr],\ldots\bigr],[\M_{\al_{\si(n)}}^\ss(\tau)]^{k_1}_\inv\bigr]
\end{aligned}
\label{co9eq54}\\
&=\sum_{\si\in S_n}
\begin{aligned}[t]
&\sum_{p=0}^{n-1}(-1)^p\binom{n-1}{p}\la_k(\al_{\si(p+1)})\cdot{}
\\
&[\M_{\al_{\si(1)}}^\ss(\tau)]^{k_1}_\inv*
[\M_{\al_{\si(2)}}^\ss(\tau)]^{k_1}_\inv*\cdots*[\M_{\al_{\si(n)}}^\ss(\tau)]^{k_1}_\inv,
\end{aligned}
\nonumber
\ea
for $k=k_1$ or $k_2$. To see this, note that the top line of \eq{co9eq54} is a sum of $n-1$ iterated Lie brackets $[\,,\,]$ in $\check H_{\rm even}(\M^\pl)$. In $U(\check H_{\rm even}(\M^\pl))$ we may expand $[f,g]=f*g-g*f$. So expanding $n-1$ iterated Lie brackets $[\,,\,]$ gives $2^{n-1}$ terms. For each of these, let $p=0,\ldots,n-1$ be the number of Lie brackets for which we chose the $-g*f$ term from $f*g-g*f$. For fixed $p$, there are $\binom{n-1}{p}$ such terms, and the term has sign $(-1)^p$ from the $p$ minus signs in $-g*f$. Also the term $[\M_{\al_{\si(1)}}^\ss(\tau)]^{k_1}_\inv$ ends up in the $(p+1)^{\rm th}$ place in the product of $[\M_{\al_i}^\ss(\tau)]^{k_1}_\inv$, as taking $p$ `$-g*f$' terms moves $p$ factors $[\M_{\al_i}^\ss(\tau)]^{k_1}_\inv$ to the left of $[\M_{\al_{\si(1)}}^\ss(\tau)]^{k_1}_\inv$ in the product. So after reordering $1,\ldots,n$, which we can do because of the symmetrization $\sum_{\si\in S_n}\cdots$, we get the second line of \eq{co9eq54}.
 
We may now write the sum $\sum_{\si\in S_n}\{\cdots\}$ in \eq{co9eq53} as
\ea
&\sum_{\si\in S_n}
\biggl\{-\sum_{p=0}^{n-1}(-1)^p\binom{n-1}{p}\la_{k_2}(\al)\,\la_{k_1}(\al_{\si(p+1)})
\label{co9eq55}\\
&+\sum_{p=0}^{n-1}(-1)^p\binom{n-1}{p}\la_{k_1}(\al)\,\la_{k_2}(\al_{\si(p+1)})
\nonumber\\
&+\sum_{m=1}^{n-1}\binom{n}{m}\sum_{p=0}^{m-1}(-1)^p\binom{m-1}{p}\sum_{q=0}^{n-m-1}(-1)^q\binom{n-m-1}{q}
\nonumber\\
&\bigl(\la_{k_1}(\al_{\si(p+1)})\la_{k_2}(\al_{\si(m+q+1)})\!-\!\la_{k_2}(\al_{\si(p+1)})\la_{k_1}(\al_{\si(m+q+1)})\bigr)
\biggr\}\cdot{}
\nonumber\\
&[\M_{\al_{\si(1)}}^\ss(\tau)]^{k_1}_\inv*
[\M_{\al_{\si(2)}}^\ss(\tau)]^{k_1}_\inv*\cdots*[\M_{\al_{\si(n)}}^\ss(\tau)]^{k_1}_\inv.
\nonumber
\ea
Here the first term in the sum $\sum_{\si\in S_n}\{\cdots\}$ in \eq{co9eq53} equals the first term in \eq{co9eq55} by \eq{co9eq54} with $k=k_1$, the second terms agree by \eq{co9eq54} with $k=k_2$, and the third terms agree by \eq{co9eq54} with $m,k_1$ and $n-m,k_2$ in place of~$n,k$. 

Now consider the expression $\{\cdots\}$ in \eq{co9eq55}. Writing $\al=\al_{\si(1)}+\cdots+\al_{\si(n)}$, we have $\la_{k_a}(\al)=\sum_{i=1}^n\la_{k_a}(\al_{\si(i)})$. Hence we may rewrite $\{\cdots\}$ as
\ea
\sum_{i,j=0}^{n-1} &\la_{k_1}(\al_{\si(i+1)})\la_{k_2}(\al_{\si(j+1)})\cdot\biggl\{-(-1)^i\binom{n-1}{i}+(-1)^j\binom{n-1}{j}
\nonumber\\
&+\sum_{\begin{subarray}{l} m=i+1,\ldots,j\\
\text{if $i<j$}\end{subarray}}\binom{n}{m}(-1)^i\binom{m-1}{i}(-1)^{j-m}\binom{n-m-1}{j-m}
\label{co9eq56}\\
&-\sum_{\begin{subarray}{l} m=j+1,\ldots,i\\
\text{if $j<i$}\end{subarray}}\binom{n}{m}(-1)^j\binom{m-1}{j}(-1)^{i-m}\binom{n-m-1}{i-m}\biggr\},
\nonumber
\ea
substituting $p=i$ in the first term, $p=j$ in the second, $p=i$ and $q=j-m$ in the third, and $p=j$ and $q=i-m$ in the fourth.

Write $S_{ij}^n$ for the expression $\{\cdots\}$ in \eq{co9eq56}. Then we have
\begin{align*}
&\sum_{i,j=0}^{n-1}S_{ij}^nx^iy^j=-\sum_{i=0}^{n-1}(-1)^i\binom{n-1}{i}x^i\sum_{j=0}^{n-1}y^j+\sum_{i=0}^{n-1}x^i\sum_{j=0}^{n-1}(-1)^j\binom{n-1}{j}y^j
\\
&+\sum_{m=1}^{n-1}\binom{n}{m}\sum_{i=0}^{m-1}(-1)^i\binom{m-1}{i}x^i\sum_{j=m}^{n-1}(-1)^{j-m}\binom{n-m-1}{j-m}y^{j-m}y^m\\
&-\sum_{m=1}^{n-1}\binom{n}{m}\sum_{i=m}^{n-1}(-1)^{i-m}\binom{n-m-1}{i-m}x^{i-m}x^m\sum_{j=0}^{m-1}(-1)^j\binom{m-1}{j}y^j
\allowdisplaybreaks\\
&=-(1-x)^{n-1}\frac{1-y^n}{1-y}+\frac{1-x^n}{1-x}(1-y)^{n-1}\\
&+\sum_{m=1}^{n-1}\binom{n}{m}(1-x)^{m-1}(1-y)^{n-m-1}y^m\\
&-\sum_{m=1}^{n-1}\binom{n}{m}(1-x)^{n-m-1}x^m(1-y)^{m-1}
\allowdisplaybreaks\\
&=-(1-x)^{n-1}\frac{1-y^n}{1-y}+\frac{1-x^n}{1-x}(1-y)^{n-1}\\
&+\Bigl[\Bigl(1+\frac{(1-x)y}{1-y}\Bigr)^n-1-\frac{(1-x)^ny^n}{(1-y)^n}\Bigr]\frac{(1-y)^{n-1}}{1-x}\\
&-\Bigl[\Bigl(1+\frac{x(1-y)}{1-x}\Bigr)^n-1-\frac{x^n(1-y)^n}{(1-x)^n}\Bigr]\frac{(1-x)^{n-1}}{1-y}
\allowdisplaybreaks\\
&=(1-x)^{-1}(1-y)^{-1}\bigl[-(1-x)^n(1-y^n)+(1-x^n)(1-y)^n\\
&+\bigl((1-y)+(1-x)y\bigr)^n-(1-y)^n-(1-x)^ny^n\\
&-\bigl((1-x)+x(1-y)\bigr)^n+(1-x)^n+x^n(1-y)^n\bigr]=0.
\end{align*}
Hence $S_{ij}^n=0$ for all $i,j,n$, so \eq{co9eq56} is zero, and \eq{co9eq55} is zero, and the right hand side of \eq{co9eq53} is zero. As $\la_{k_1}(\al),\la_{k_2}(\al_1)>0$ by Assumption \ref{co5ass1}(g)(iv), equation \eq{co9eq53} implies that $[\M_\al^\ss(\tau)]^{k_2}_\inv-[\M_\al^\ss(\tau)]^{k_1}_\inv=0$. This proves the inductive step, and the proposition follows by induction.
\end{proof}

\subsection{\texorpdfstring{Extension to the $G$-equivariant case}{Extension to the G-equivariant case}}
\label{co94}

Theorem \ref{co5thm4} extends Theorems \ref{co5thm1}--\ref{co5thm3} to $G$-equivariant homology $H_*^G(\M^\pl)$. The proof of Theorem \ref{co5thm1} in \S\ref{co91}--\S\ref{co93} must be modified as follows:
\begin{itemize}
\setlength{\itemsep}{0pt}
\setlength{\parsep}{0pt}
\item[(a)] We replace $H_*(-),H^*(-)$ by $H^G_*(-),H_G^*(-)$ throughout.
\item[(b)] In Proposition \ref{co9prop2}(c), $G$ acts on $\acM_{(\al,\bs 1)}^\ss(\ac\tau^0_{\bs\mu})_{\al_1,\al_2}$ on the left of \eq{co9eq13}, but $G\!\t\! G$ acts on $\baM_{(\al_1,1)}^\ss(\bar\tau^0_1)_{k_1}\!\t\!\baM_{(\al_2,1)}^\ss(\bar\tau^0_1)_{k_2}$. Then \eq{co9eq13} is $G$-equivariant for the diagonal $G$-action on~$\baM_{(\al_1,1)}^\ss(\bar\tau^0_1)_{k_1}\!\t\!\baM_{(\al_2,1)}^\ss(\bar\tau^0_1)_{k_2}$.
\item[(c)] In Proposition \ref{co9prop3}(d) we must replace \eq{co9eq15} by
\end{itemize}
\ea
&[\acM_{(\al,\bs 1)}^\ss(\ac\tau^0_{\bs\mu})_{\al_1,\al_2}]_\virt
\label{co9eq57}\\
&=\begin{cases}
\La^{G\t G,G}\bigl([\baM_{(\al_1,1)}^\ss(\bar\tau^0_1)_{k_1}]_\virt\bt[\baM_{(\al_2,1)}^\ss(\bar\tau^0_1)_{k_2}]_\virt\bigr), & o_{\al_1}+o_{\al_2}=o_\al, \\
0, & o_{\al_1}+o_{\al_2}>o_\al.
\end{cases}
\nonumber
\ea
\begin{itemize}
\setlength{\itemsep}{0pt}
\setlength{\parsep}{0pt}
\item[] Here $\La^{G\t G,G}$ maps $H_*^{G\t G}(\cdots)$ to $H_*^G(\cdots)$ as in Definition~\ref{co2def2}(g).
\item[(c)] In the first step of equation \eq{co9eq45} in the proof of Proposition \ref{co9prop5}, we substitute in \eq{co9eq57} rather than \eq{co9eq15}. This gives an extra operator $\La^{G\t G,G}$, which persists through the second--eighth steps of \eq{co9eq45}. In the ninth step we use the $G$-equivariant version \eq{co4eq24} of the state-field correspondence $Y(-,z)$, which includes an extra operator $\La^{G\t G,G}$ compared to the usual definition \eq{co4eq10}. Thus the $\La^{G\t G,G}$ is absent in the ninth--twelfth steps of \eq{co9eq45}, and \eq{co9eq43} holds in $\check H_{\rm even}^G(\M^\pl)$.
\end{itemize}

\section{Proof of Theorem \ref{co5thm2}}
\label{co10}

\subsection{Set up for the proof of Theorem \ref{co5thm2}}
\label{co101}

We define notation that will be used in the proof of Theorem \ref{co5thm2}.

\begin{dfn}
\label{co10def1}
As in Theorem \ref{co5thm2}, let Assumptions \ref{co5ass1}--\ref{co5ass2} hold for $\B\subseteq\A,\M,\M^\pl,K(\A),\ldots,$ and fix $(\tau,T,\le),(\ti\tau,\ti T,\le)$ in $\sS$ and $\al\in C(\B)_\pe$. Suppose that if $\al_1,\ldots,\al_n\in C(\B)$ with $\al=\al_1+\cdots+\al_n,$ and either 
\begin{itemize}
\setlength{\itemsep}{0pt}
\setlength{\parsep}{0pt}
\item[(i)] $U(\al_1,\ldots,\al_n;\tau,\ti\tau)\ne 0$ and $\M_{\al_i}^\ss(\tau)\ne\es$ for $i=1,\ldots,n$; or 
\item[(ii)] $U(\al_1,\ldots,\al_n;\ti\tau,\tau)\ne 0$ and $\M_{\al_i}^\ss(\ti\tau)\ne\es$ for $i=1,\ldots,n,$
\end{itemize}
then $\tau(\al_i)=\tau(\al)$ for all $i=1,\ldots,n$, and in case (ii) we also have $\M_{\al_i}^\ss(\ti\tau)\subseteq\M_{\al_i}^\ss(\tau)$ for~$i=1,\ldots,n$.

Using the notation of Assumption \ref{co5ass1}(g), pick $k\in K$ such that $\M_\al^\ss(\tau)\subseteq\M_{k,\al}^\pl$, which is possible by the last part of Assumption \ref{co5ass1}(g). This $k$ will be fixed throughout Chapter \ref{co10}. Set $r=\la_k(\al)$, so that $r>0$ by Assumption \ref{co5ass1}(g)(iv). Define data $\acB\subseteq\acA,K(\acA),\ab\acM,\ab\acM^\pl,\ab\ldots$ as in Definition \ref{co5def1} using the quiver $\ac Q=(Q_0,Q_1,h,t)$ illustrated by:
\e
\begin{xy}
0;<.87mm,0mm>:
,(-90,3)*{e_1}
,(-100,3)*{v_1}
,(-72,3)*{e_2}
,(-80,3)*{v_2}  
,(-49,3)*{e_{r-2}}
,(-60,0)*+{\cdots}
,(-30,3)*{e_{r-1}}
,(-40,3)*{v_{r-1}}
,(-10,3)*{e_r}
,(-20,3)*{v_r}
,(0,3)*{w}
,(11,0)*{\ka(w)\!=\!k,}
,(-20,0)*+{\bu} ; (0,0)*+{\circ} **@{-} ?>*\dir{>}
,(-40,0)*+{\bu} ; (-20,0)*+{\bu} **@{-} ?>*\dir{>}
,(-56,0)*+{} ; (-40,0)*+{\bu} **@{-} ?>*\dir{>}
,(-80,0)*+{\bu} ; (-64,0)*+{} **@{-} ?>*\dir{>}
,(-100,0)*+{\bu} ; (-80,0)*+{\bu} **@{-} ?>*\dir{>}
\end{xy}
\label{co10eq1}
\e
with vertices $Q_0=\{v_1,\ldots,v_r,w\}$, edges $Q_1=\{e_1,\ldots,e_r\}$, and head and tail maps $h(e_i)=v_{i+1}$, $i<r$, $h(e_r)=w$, and $t(e_i)=v_i$. We take $\dot Q_0=\{v_1,\ldots,v_r\}$ and $\ddot Q_0=\{w\}$, and define $\ka:\ddot Q_0\ra K$ by $\ka(w)=k$. We write the stability conditions on $\acA$ defined in Definition \ref{co5def1} as $(\ac\tau^\la_{\bs\mu},\ac T,\le)$, which form a set $\acS$. We write $\ac\Up_{(\al,\bs d)},\ac\Up$ for the maps $\bar\Up_{(\al,\bs d)},\bar\Up$ defined in \eq{co5eq15} for $\acB\subseteq\acA,\ldots.$

Here we use notation $\acA,\acB,K(\acA),\acM,(\ac\tau^\la_{\bs\mu},\ac T,\le),\ac\Up_{(\al,\bs d)},\ldots$ rather than $\baA,\ab\baB,\ab K(\baA),\ab\baM,\ab(\bar\tau^\la_{\bs\mu},\bar T,\le),\bar\Up_{(\al,\bs d)},\ab\ldots$ to distinguish them from $\baA,\ab\baB,\ldots$ in Example \ref{co5ex1}, which will also appear in the proof of Proposition \ref{co10prop7} below.

For brevity we will write objects of $\acA$ as $(E,\bs V,\bs\rho)$ with $\bs V=(V_1,\ldots,V_r)$, $\bs\rho=(\rho_1,\ldots,\rho_r)$ rather than $(V_{v_1},\ldots,V_{v_r}),(\rho_{e_1},\ldots,\rho_{e_r})$, and write $K(\acA)=K(\A)\t\Z^r$ rather than $K(\A)\t\Z^{\{v_1,\ldots,v_r\}}$, and write stability conditions as $\ac\tau^\la_{\bs\mu}$ with $\bs\mu=(\mu_1,\ldots,\mu_r)\in\R^r$ rather than $\bs\mu=(\mu_{v_1},\ldots,\mu_{v_r})$, and so on. We will use these $\acB,\acA,\ldots$ for the whole of Chapter~\ref{co10}.

Define a biadditive map $\ac\chi:K(\acB)\t K(\acB)\ra\Z$, the analogue of $\chi:K(\B)\t K(\B)\ra\Z$ in Assumption \ref{co4ass1}(e), by
\e
\ac\chi\bigl((\be,\bs d),(\ga,\bs e)\bigr)=\ts\chi(\be,\ga)+\sum_{i=1}^rd_ie_i-\sum_{i=1}^{r-1}d_ie_{i+1}-d_r\,\la_k(\ga).
\label{co10eq2}
\e

There is also a natural inclusion functor $I_\acB:\B_\ka\hookra\acB$ acting on objects by $I_\acB:E\mapsto(E,\bs 0,\bs 0)$, with $\Pi_\B\ci I_\acB=\Id_{\B_\ka}$. It identifies $\B_\ka$ with the full subcategory of objects in $\acB$ in classes $(\be,0)$ in $K(\acB)$. It induces morphisms $I_\acM:\M_\ka\hookra\acM$, $I_{\acM^\pl}:\M_\ka^\pl\hookra\acM^\pl$ of moduli stacks, which are open inclusions, with $\Pi_\M\ci I_\acM=\id_{\M_\ka}$, $\Pi_{\M^\pl}\ci I_{\acM^\pl}=\id_{\M_\ka^\pl}$. 

The following sets $R_\al,\ti R_\al$ will be important in the proofs below:
\e
\begin{split}
R_\al&=\bigl\{\be\in C(\B):\,\text{either $\be=\al$, or $\ga=\al-\be\in C(\B)$}\\
&\qquad\qquad\qquad\text{with $\tau(\al)=\tau(\be)=\tau(\ga)$ and $\M_\be^\ss(\tau),\M_\ga^\ss(\tau)\ne\es$}\bigr\},\\
\ti R_\al&=\bigl\{\be\in C(\B):\,\text{either $\be=\al$, or $\ga=\al-\be\in C(\B)$}\\
&\qquad\qquad\qquad\text{with $\ti\tau(\al)=\ti\tau(\be)=\ti\tau(\ga)$ and $\M_\be^\ss(\ti\tau),\M_\ga^\ss(\ti\tau)\ne\es$}\bigr\}.
\end{split}
\label{co10eq3}
\e
Lemma \ref{co9lem1} with $n=2$ shows $R_\al$ is finite. As $\al\in C(\B)_\pe$, Assumption \ref{co5ass2}(c) implies that $R_\al\subseteq C(\B)_\pe$. Proposition \ref{co10prop1} shows $\ti R_\al\subseteq R_\al$. Our proofs will only involve moduli spaces $\M^\ss_\be(\tau)$ for $\be\in R_\al$ and $\M^\ss_\be(\ti\tau)$ for~$\be\in\ti R_\al$.

If $\be\in R_\al$ and $[E]$ is a $\C$-point of $\M_\be^\ss(\tau)$ then either $\be=\al$ and we set $F=0$, or for $\ga$ as in \eq{co10eq3} we can choose a $\C$-point $[F]$ of $\M_\ga^\ss(\tau)$. Then in both cases $[E\op F]$ is a $\C$-point of $\M_\al^\ss(\tau)$. Since $\M_\al^\ss(\tau)\subseteq\M_{k,\al}^\pl$, we see that $E\op F\in\B_\ka$, so $E\in\B_\ka$ as $\B_\ka$ is closed under direct summands in $\A$, and~$\be\in C(\B_\ka)$. 

We also define finite sets $S_\al,\ti S_\al$ by
\e
\begin{split}
S_\al&=\bigl\{(\be,\bs d):\be\in R_\al\amalg\{0\}, \; \bs d\in\N^r,\; (\be,\bs d)\ne (0,0),\; d_i\le i,\; \text{all $i$}\bigr\},\\
\ti S_\al&=\bigl\{(\be,\bs d):\be\in\ti R_\al\amalg\{0\}, \; \bs d\in\N^r,\; (\be,\bs d)\ne (0,0),\; d_i\le i,\; \text{all $i$}\bigr\}.
\end{split}
\label{co10eq4}
\e
Then $S_\al\subseteq C(\acB)_\pe$ by \eq{co5eq16}, as $R_\al\subseteq C(\B)_\pe$. Proposition \ref{co10prop1} will show that $\ti S_\al\subseteq S_\al$. Our proofs only involve moduli spaces $\acM^\ss_{(\be,\bs d)}(\ac\tau^\la_{\bs\mu_x})$ for~$(\be,\bs d)\in S_\al$.

Choose generic $\mu_1,\ldots,\mu_r\in\R$ satisfying
\e
1\gg\mu_1\gg \mu_2\gg \cdots\gg \mu_{r-1}\gg \mu_r>0.
\label{co10eq5}
\e
Here the `$\gg$' mean that $1/\mu_1$ and $\mu_{i-1}/\mu_i$ for $1<i\le r$ must satisfy finitely many positive lower bounds, depending on $\rk\al,S_\al$ and $\ti\la$ chosen in Proposition \ref{co10prop4}, that will appear in the proof of Proposition \ref{co10prop8} below. In particular, we will require that $\mu_i/\mu_j>\rk\al$ if $1\le i<j\le n$.

For each $x\in[-1,0]$, define $\bs\mu_x=(\mu_1+x,\mu_2+x,\ldots,\mu_r+x)$ in $\R^r$. For each group morphism $\la:K(\A)\ra\R$, define a weak stability  condition $(\ac\tau^\la_{\bs\mu_x},\ac T,\le)$ on $\acA$ as in Definition \ref{co5def1}. We will use $(\ac\tau^0_{\bs\mu_{\smash{-1}}},\ac T,\le)$ as our reference weak stability condition on~$\acA$.

As in Definition \ref{co5def1}, for such each $(\ac\tau^\la_{\bs\mu_x},\ac T,\le)$ and $(\be,\bs d)\in S_\al$ we have moduli spaces $\acM_{(\be,\bs d)}^\rst(\ac\tau^\la_{\bs\mu_x})\subseteq\acM_{(\be,\bs d)}^\ss(\ac\tau^\la_{\bs\mu_x})\subseteq\dM_{(\be,\bs d)}^\pl\subseteq\acM_{(\be,\bs d)}^\pl$. Here $\dM_{(\be,\bs d)}^\pl$ is the classical truncation of the quasi-smooth derived stack $\bs\dM_{(\be,\bs d)}^\rpl$ defined in \eq{co5eq17}, so as in \eq{co5eq18} we have an obstruction theory $\bL_i:i^*(\bL_{\bs\dM_{(\be,\bs d)}^\rpl})\ra \bL_{\dM_{(\be,\bs d)}^\pl}$ on $\dM_{(\be,\bs d)}^\pl$, which restricts to an obstruction theory on $\acM_{(\be,\bs d)}^\ss(\ac\tau^\la_{\bs\mu_x})$.

If $\acM_{(\be,\bs d)}^\rst(\ac\tau^\la_{\bs\mu_x})=\acM_{(\be,\bs d)}^\ss(\ac\tau^\la_{\bs\mu_x})$ then $\acM_{(\be,\bs d)}^\ss(\ac\tau^\la_{\bs\mu_x})$ is a proper algebraic space by Assumption \ref{co5ass2}(h), so we have a virtual class $[\acM_{(\be,\bs d)}^\ss(\ac\tau^\la_{\bs\mu_x})]_\virt$, as in \S\ref{co24}. Thus as in \eq{co5eq15}, \eq{co5eq23} and \eq{co5eq29} we may form the homology class
\e
\begin{split}
&\ac\Up_{(\be,\bs d)}(\ac\tau^\la_{\bs\mu_x}):=\ac\Up\bigl([\acM_{(\be,\bs d)}^\ss(\ac\tau^\la_{\bs\mu_x})]_\virt\bigr)\\
&=(\Pi_{\M_\be^\ss(\tau)})_*\bigl([\acM_{(\be,\bs d)}^\ss(\ac\tau^\la_{\bs\mu_x})]_\virt\cap c_\top(\bT_{\acM_{(\be,\bs d)}^\ss(\ac\tau^\la_{\bs\mu_x})/\M_\be^\ss(\tau)})\bigr)\\
&=(-1)^{\sum_{i<r}d_id_{i+1}+d_r\la_k(\be)-\sum_{i\le r}d_i^2}\cdot{} \\
&\qquad (\Pi_{\M_\be^\ss(\tau)})_*\bigl([\acM_{(\be,\bs d)}^\ss(\ac\tau^\la_{\bs\mu_x})]_\virt\cap c_\top(\bL_{\acM_{(\be,\bs d)}^\ss(\ac\tau^\la_{\bs\mu_x})/\M_\be^\ss(\tau)})\bigr)\\
&\text{in}\qquad H_{2o_\be+2-2\chi(\be,\be)}(\M_\be^\pl)=\check H_{2o_\be}(\M_\be^\pl).
\end{split}
\label{co10eq6}
\e
Our proof of Theorem \ref{co5thm2} will involve studying the classes \eq{co10eq6}.
\end{dfn}

The next result uses the assumptions involving Definition~\ref{co10def1}(i),(ii).

\begin{prop}
\label{co10prop1}
In Definition\/ {\rm\ref{co10def1}} we have $\ti R_\al\subseteq R_\al,$ and\/ $\ti S_\al\subseteq S_\al,$ and if\/ $\be\in \ti R_\al$ then $\M_\be^\ss(\ti\tau)\subseteq\M_\be^\ss(\tau)$. If\/ $\be\in\ti R_\al$ and\/ $[E]\in\M_\be^\ss(\tau)\sm\M_\be^\ss(\ti\tau)$ there exists $0\ne E'\subsetneq E$ such that\/ $E',E/E'$ are $\tau$-semistable, $\lb E'\rb,\lb E/E'\rb\in R_\al,$ and\/~$\ti\tau(\lb E'\rb)>\ti\tau(\lb E/E'\rb)$.
\end{prop}

\begin{proof} Let $[E]$ be a $\C$-point of $\M_\al^\ss(\ti\tau)$. By Assumption \ref{co5ass2}(a), $E$ has a unique $\tau$-Harder--Narasimhan filtration $0=E_0\subsetneq E_1\subsetneq \cdots\subsetneq E_n=E$ in $\B$, such that $F_i=E_i/E_{i-1}$ is $\tau$-semistable with $\tau(\al_1)>\cdots>\tau(\al_n)$, where $\al_i=\lb F_i\rb$ for $i=1,\ldots,n$. Corollary \ref{co3cor1} shows that $U(\al_1,\ldots,\al_n;\tau,\ti\tau)\ne 0$, and $[F_i]\in\M_{\al_i}^\ss(\tau)\ne\es$, so the assumption using Definition \ref{co10def1}(i) gives $\tau(\al_i)=\tau(\al)$ for $i=1,\ldots,n$. But as $\tau(\al_1)>\cdots>\tau(\al_n)$ we see that $n=1$, so $E=F_1$ is $\tau$-semistable. This proves that
\e
\M_\al^\ss(\ti\tau)\subseteq\M_\al^\ss(\tau).
\label{co10eq7}
\e

Now let $\be\in\ti R_\al$ with $\be\ne\al$, and write $\ga=\al-\be$. Then \eq{co10eq3} says that $\ti\tau(\al)=\ti\tau(\be)=\ti\tau(\ga)$ and $\M_\be^\ss(\ti\tau),\M_\ga^\ss(\ti\tau)\ne\es$. Choose $\C$-points $[E]\in\M_\be^\ss(\ti\tau)$ and $[F]\in\M_\ga^\ss(\ti\tau)$. Then $E\op F$ is $\ti\tau$-semistable as $E,F$ are with $\ti\tau(\be)=\ti\tau(\ga)$, so $[E\op F]\in \M_\al^\ss(\ti\tau)\subseteq\M_\al^\ss(\tau)$ by \eq{co10eq7}. Hence $E,F$ are $\tau$-semistable, as nonzero direct summands of $\tau$-semistable objects are $\tau$-semistable. Thus $[E]\in\M_\be^\ss(\tau)$, $[F]\in\M_\ga^\ss(\tau)$, so $\M_\be^\ss(\ti\tau)\subseteq\M_\be^\ss(\tau)$, as we have to prove, and $\M_\ga^\ss(\ti\tau)\subseteq\M_\ga^\ss(\tau)$. Therefore $\M_\be^\ss(\tau)\ne\es$ and $\M_\ga^\ss(\tau)\ne\es$, so $\be\in R_\al$ by \eq{co10eq3}. Since $\al\in R_\al$, this proves that $\ti R_\al\subseteq R_\al$, as we want, and $\ti S_\al\subseteq S_\al$ follows from \eq{co10eq4}. This completes the first part of the proposition.

For the second part, suppose $\be\in\ti R_\al$ and $[E]\in\M_\be^\ss(\tau)\sm\M_\be^\ss(\ti\tau)$. If $\be\ne\al$ then set $\ga=\al-\be$, so that $\M_\ga^\ss(\ti\tau)\ne\es$ by \eq{co10eq3}, and pick a $\C$-point $[G]$ in $\M_\ga^\ss(\ti\tau)\subseteq\M_\ga^\ss(\tau)$. If $\be=\al$ set $G=0$. Then in both cases $[E\op G]$ is a $\C$-point of $\M_\al^\ss(\tau)\sm\M_\al^\ss(\ti\tau)$. By Assumption \ref{co5ass2}(a), $E$ has a unique $\ti\tau$-Harder--Narasimhan filtration $0=E_0\subsetneq E_1\subsetneq \cdots\subsetneq E_n=E$ in $\B$, such that $F_i=E_i/E_{i-1}$ is $\ti\tau$-semistable with $\ti\tau(\al_1)>\cdots>\ti\tau(\al_n)$, where $\al_i=\lb F_i\rb$ for $i=1,\ldots,n$. Also $n\ge 2$ as $E$ is not $\ti\tau$-semistable. 

If $G\ne 0$, as $G$ is $\ti\tau$-semistable, the $\ti\tau$-Harder--Narasimhan filtration of $G$ is $0\subsetneq G$. Write $0=E'_0\subsetneq E'_1\subsetneq \cdots\subsetneq E'_{n'}=E\op G$ with $F_i'=E_i'/E_{i-1}'$ and $\al_i'=\lb F_i'\rb$ for the $\ti\tau$-Harder--Narasimhan filtration of $E\op G$. This is obtained by combining the $\ti\tau$-Harder--Narasimhan filtrations of $E,G$, in such a way that
\begin{align*}
&(F_1',\ldots,F'_{n'})=(F_1,\ldots,F_{i-1},F_i\op G,F_{i+1},\ldots,G_n) \quad\text{if $\ti\tau(\al_i)=\ti\tau(\ga)$,}\\
&(F_1',\ldots,F'_{n'})=(F_1,\ldots,F_i,G,F_{i+1},\ldots,F_n) \quad\text{if $\ti\tau(\al_i)<\ti\tau(\ga)<\ti\tau(\al_{i+1})$,}
\end{align*}
omitting conditions $\ti\tau(\al_i)<\ti\tau(\ga)$ if $i=0$ and $\ti\tau(\ga)<\ti\tau(\al_{i+1})$ if~$i=n$.

Corollary \ref{co3cor1} shows that  $U(\al'_1,\ldots,\al'_{n'};\ti\tau,\tau)\ne 0$, and $[F'_i]\in\M_{\al'_i}^\ss(\ti\tau)\ne\es$, so the assumption using Definition \ref{co10def1}(ii) gives $\tau(\al'_i)=\tau(\al)$ and $\M_{\al'_i}^\ss(\ti\tau)\subseteq\M_{\al'_i}^\ss(\tau)$ for $i=1,\ldots,n'$. Thus $F'_i$ is $\tau$-semistable for $i=1,\ldots,n'$. As the $F_i'$ are of the form $F_i,F_i\op G$, or $G$, these imply that $F_i$ is $\tau$-semistable and $\tau(\al_i)=\tau(\al)$ for~$i=1,\ldots,n$.

Set $E'=E_1\subsetneq E$. Then $E'=F_1$ is $\tau$-semistable, and $E/E'$ has a filtration $0=E_0''\subsetneq E_1''\subsetneq \cdots\subsetneq E_{n-1}''=E/E'$ with $E_i''=E_{i+1}/E_1$, and $E_i''/E_{i-1}''=F_{i+1}$ $\tau$-semistable with $\tau(\lb E_i''/E_{i-1}''\rb)=\tau(\al_{i+1})=\tau(\al)$ for $i=1,\ldots,n-1$, so $E/E'$ is $\tau$-semistable. We have $\tau(\lb E'\rb)=\tau(\al_1)=\tau(\al)$ and $\tau(\lb E/E'\rb)=\tau(\al_2+\cdots+\al_n)=\tau(\al)$. Also $G$ is $\tau$-semistable with $\tau(\lb G\rb)=\tau(\ga)=\tau(\al)$ if $G\ne 0$, so $E'\op G$ and $(E/E')\op G$ are also $\tau$-semistable with $\tau(\lb E'\op G\rb)=\tau(\lb (E/E')\op G\rb)=\tau(\al)$. By considering $[E'\op ((E/E')\op G)]$ and $[(E/E')\op (E'\op G)]$ as $\C$-points of $\M_\al^\ss(\tau)$, we see from \eq{co10eq3} that $\lb E'\rb,\lb E/E'\rb\in R_\al$. Also $\ti\tau(\lb E'\rb)>\ti\tau(\lb E/E'\rb)$ as $\ti\tau(\al_1)>\cdots>\ti\tau(\al_n)$ and $\lb E'\rb=\al_1$, $\lb E/E'\rb=\al_2+\cdots+\al_n$ with~$n\ge 2$.
\end{proof}

We prove some necessary conditions on $(\be,\bs d)$ for $\acM_{(\be,\bs d)}^\ss(\ac\tau^\la_{\bs\mu_x})\ne\es$.

\begin{prop}
\label{co10prop2}
Let\/ $(\ac\tau^\la_{\bs\mu_x},\ac T,\le)$ be as in Definition\/ {\rm\ref{co10def1}}.
\smallskip

\noindent{\bf(a)} Suppose $(E,\bs V,\bs\rho)$ is a $\ac\tau^\la_{\bs\mu_x}$-semistable object of\/ $\acB$. Then:
\begin{itemize}
\setlength{\itemsep}{0pt}
\setlength{\parsep}{0pt}
\item[{\bf(i)}] If\/ $E=0$ then for some $1\le a\le b\le r$ we have $V_i=0$ unless $a\le i\le b,$ and\/ $\rho_i:V_i\ra V_{i+1}$ is an isomorphism for $a\le i<b$. 
\item[{\bf(ii)}] If\/ $E\ne 0$ then $\rho_i$ is injective for\/ $i=1,\ldots,r$.
\item[{\bf(iii)}] If\/ $E\ne 0$ and\/ $\frac{1}{r-a}(\mu_{a+1}+\cdots+\mu_r)+x\le 0$ for some $1\le a<r$ then $\rho_i$ is an isomorphism for $a\le i<r$.
\item[{\bf(iv)}] If\/ $E\ne 0$ and\/ $\frac{1}{r}(\mu_1+\cdots+\mu_r)+x\le 0$ then\/~$V_1=\cdots=V_r=0$.
\end{itemize}

\noindent{\bf(b)} Suppose $(\be,\bs d)\in S_\al$ with\/ $\acM_{(\be,\bs d)}^\ss(\ac\tau^\la_{\bs\mu_x})\ne\es$. Then:
\begin{itemize}
\setlength{\itemsep}{0pt}
\setlength{\parsep}{0pt}
\item[{\bf(i)}] If\/ $\be=0$ then for some then for some $1\le a\le b\le r$ and\/ $n>0$ we have $d_i=n$ if\/ $a\le i\le b$ and $d_i=0$ otherwise.
\item[{\bf(ii)}] If\/ $\be\ne 0$ then $d_1\le d_2\le \cdots\le d_r\le \la_k(\be)$. 
\item[{\bf(iii)}] If\/ $\be\ne 0$ and\/ $\frac{1}{r-a}(\mu_{a+1}+\cdots+\mu_r)+x\le 0$ for some\/ $1\le a<r$ then\/~$d_a=d_{a+1}=\cdots=d_r$.
\item[{\bf(iv)}] If\/ $\be\ne 0$ and\/ $\frac{1}{r}(\mu_1+\cdots+\mu_r)+x\le 0$ then\/~$\bs d=0$.
\end{itemize}
\end{prop}

\begin{proof} For (a), suppose $(E,\bs V,\bs\rho)$ is $\ac\tau^\la_{\bs\mu_x}$-semistable. Let $b=r+1$ if $E\ne 0$, and $b=1,\ldots,r$ be largest such that $V_b\ne 0$ if $E=0$. Define vector subspaces $U_i\subseteq V_i$ for $i=1,\ldots,r$ by $U_i=\Ker\bigl(\rho_{b-1}\ci\rho_{b-2}\ci\cdots\ci\rho_i:V_i\ra V_b\bigr)$ if $i<b$, writing $V_{r+1}=F_k(E)$, and $U_i=0$ if~$i\ge b$. 

Define vector subspaces $W_i\subseteq V_i$ with $V_i=U_i\op W_i$ by $W_i=V_i$ for $i\ge b$, and by reverse induction on $i=b-1,b-2,\ldots,1$ choose $W_i$ with $V_i=U_i\op W_i$ and $\rho_i(W_i)\subseteq W_{i+1}$. Then $\rho_i(U_i)\subseteq U_{i+1}$, $\rho_i(W_i)\subseteq W_{i+1}$ for all $1\le i<r$, and we have a direct sum $(E,\bs V,\bs\rho)=(0,\bs U,\bs\rho\vert_{\bs U})\op (E,\bs W,\bs\rho\vert_{\bs W})$ in $\acB$. Here $(E,\bs W,\bs\rho\vert_{\bs W})\ne 0$ as $E\ne 0$ or~$W_b=V_b\ne 0$.

If $E\ne 0$ then we must have $(0,\bs U,\bs\rho\vert_{\bs U})=0$, as otherwise one of $(0,\bs U,\bs\rho\vert_{\bs U}),\ab(E,\bs W,\bs\rho\vert_{\bs W})$ must $\ac\tau^\la_{\bs\mu_x}$-destabilize $(E,\bs V,\bs\rho)$ by \eq{co5eq21}. If $E=0$ and $(0,\bs V,\bs\rho)$ is $\ac\tau^\la_{\bs\mu_x}$-stable then $(0,\bs U,\bs\rho\vert_{\bs U})=0$ as stable objects are indecomposable.

If $E=0$ and $(0,\bs U,\bs\rho\vert_{\bs U})\ne 0$ then $\lb 0,\bs U,\bs\rho\vert_{\bs U}\rb,\lb 0,\bs W,\bs\rho\vert_{\bs W}\rb$ are not proportional as $U_b=0$, $W_b\ne 0$, so as $\mu_1,\ldots,\mu_r$ are generic we see that one of $(0,\bs U,\bs\rho\vert_{\bs U})$ or $(0,\bs W,\bs\rho\vert_{\bs W})$ must $\ac\tau^\la_{\bs\mu_x}$-destabilize $(E,\bs V,\bs\rho)$, a contradiction. Thus in all cases $U_i=0$ for all $i=1,\ldots,r$. Therefore $\rho_i$ is injective for $i=1,\ldots,b-1$. If $E\ne 0$ then $b=r+1$ and this proves~(ii). 

If $E=0$ this shows that $V_i=0$ for $i\ge b$ and $\rho_i$ is injective for $1\le i<b$. A similar argument taking $a=1,\ldots,b$ smallest with $V_a\ne 0$, defining $U_i=\Im\bigl(\rho_{i-1}\ci\rho_{i-2}\ci\cdots\ci\rho_a:V_a\ra V_i\bigr)$, and choosing $W_i$ with $V_i=U_i\op W_i$ and $\rho_i(W_i)\subseteq W_{i+1}$, shows that $W_i=0$ for all $i$ and hence $\rho_i$ is surjective for $a\le i<r$. Therefore $V_i=0$ unless $a\le i\le b,$ and $\rho_i$ is an isomorphism for $a\le i<b$, proving~(i).

For (iii)--(iv), suppose $E\ne 0$ and $\frac{1}{r-a}(\mu_{a+1}+\cdots+\mu_r)+x\le 0$ for some $0\le a<r$. Define vector subspaces $U_i\subseteq V_i$ by $U_i=V_i$ for $i\le a$ and $U_i=\rho_{i-1}\ci\rho_{i-2}\ci \cdots\ci\rho_a(U_a)$ for $i>a>0$, and $U_i=0$ when $a=0$. Define $W_i=V_i/U_i$ for $i=1,\ldots,r$. Then we have an exact sequence in $\acB$
\e 
\xymatrix@C=15pt{
0 \ar[r] & (E,\bs U,\bs\rho\vert_{\bs U}) \ar[r] & (E,\bs V,\bs\rho) \ar[r] &
(0,\bs W,\bs\si) \ar[r] & 0.  }
\label{co10eq8}
\e
As $W_i=0$ if $i\le a$, and $\dim W_i\le\dim W_{i+1}$ for $i<r$ as $\dim U_i\le\dim U_{i+1}$ by (ii), and $\mu_1>\cdots>\mu_r$, we see that if the $W_i$ are not all zero then
\begin{equation*}
\frac{1}{\sum_{i=1}^r\dim W_i}\sum_{i=1}^r(\mu_i+x)\dim W_i\le \frac{1}{r-a}(\mu_{a+1}+\cdots+\mu_r)+x\le 0.
\end{equation*}
Therefore by \eq{co5eq21} we see that $\ac\tau^\la_{\bs\mu_x}(\lb 0,\bs W,\bs\si\rb)=(-\iy,y)$ for $y\in\R$, so equation \eq{co10eq8} $\ac\tau^\la_{\bs\mu_x}$-destabilizes $(E,\bs V,\bs\rho)$, a contradiction. Thus $W_i=0$ for all $i$, so $\rho_i$ is surjective for $a\le i<r$. Part (iii) then follows from (ii) when $a>0$, and part (iv) follows when $a=0$. This completes part (a), and (b) is immediate.
\end{proof}

We compute the moduli spaces $\acM_{(\be,\bs d)}^\ss(\ac\tau^0_{\bs\mu_{\smash{-1}}})$:

\begin{prop}
\label{co10prop3}
If\/ $(\be,\bs d)\in S_\al,$ then
\ea
\acM_{(\be,\bs d)}^\ss(\ac\tau^0_{\bs\mu_{\smash{-1}}})\cong\begin{cases} 
I_\acM^\pl(\M^\ss_\be(\tau)), & \bs d=0, \\
[*/\PGL(n,\C)], & \begin{aligned}[h]&\text{$d_i\!=\!n\!>\!0$ if\/ $a\!\le\! i\!\le\! b,$ $d_i\!=\!0$ otherwise} \\
&\text{for some $1\!\le\! a\!\le\! b\le\! r,$ and\/ $\be\!=\!0,$}
\end{aligned} \\
\es, & \text{otherwise.}
\end{cases}
\nonumber\\[-15pt]
\label{co10eq9}
\ea
\end{prop}

\begin{proof} As $\ac\tau^0_{\bs\mu_{\smash{-1}}}(\lb E,\bs 0,\bs 0\rb)=(\tau(\lb E\rb),0)$ we see that an object $(E,\bs 0,\bs 0)$ in $\acB$ is $\ac\tau^0_{\bs\mu_{\smash{-1}}}$-semistable if and only if $E$ is $\tau$-semistable, so the first case of \eq{co10eq9} follows. If $\be=0$ and $(0,\bs V,\bs\rho)$ is a $\C$-point of $\acM_{(\be,\bs d)}^\ss(\ac\tau^0_{\bs\mu_{\smash{-1}}})$ then Proposition \ref{co10prop2}(a)(i) and (b)(i) show that $d_i=n>0$ if $a\le i\le b$ and $d_i=0$ otherwise for $1\le a\le b\le r$, and $\rho_i:V_i\ra V_{i+1}$ is an isomorphism for $a\le i<b$. So up to isomorphism we may take $V_i=\C^n$ for $a\le i\le b$ and $\rho_i=\id_{\C^n}$ for $a\le i<b$. It is easy to check this does give a $\ac\tau^0_{\bs\mu_{\smash{-1}}}$-semistable object, with isotropy group $\GL(n,\C)$ in $\acM_{(0,\bs d)}$ and $\PGL(n,\C)$ in $\acM_{(0,\bs d)}^\pl$. Thus, $\acM_{(\be,\bs d)}^\ss(\ac\tau^0_{\bs\mu_{\smash{-1}}})$ is one point with isotropy group $\PGL(n,\C)$ for such $\bs d$, giving the second case of \eq{co10eq9}. If $\bs d$ is not of this form then $\acM_{(\be,\bs d)}^\ss(\ac\tau^0_{\bs\mu_{\smash{-1}}})=\es$, giving part of the third case of~\eq{co10eq9}.

The remaining case is when $\be\ne 0\ne\bs d$. Then any object $(E,\bs V,\bs\rho)$ in class $(\be,\bs d)$ in $\acB$ is $\ac\tau^0_{\bs\mu_{\smash{-1}}}$-destabilized by the exact sequence
\begin{equation*}
\xymatrix@C=18pt{ 0 \ar[r] & (E,\bs 0,\bs 0) \ar[r] & (E,\bs V,\bs\rho) \ar[r] & (0,\bs V,(\rho_1,\ldots,\rho_{r-1},0))) \ar[r] & 0, } 	
\end{equation*}
so $\acM_{(\be,\bs d)}^\ss(\ac\tau^0_{\bs\mu_{\smash{-1}}})=\es$. This completes the proof.	
\end{proof}

\subsection{\texorpdfstring{Reducing Theorem \ref{co5thm2} to a wall-crossing in $\ac{\cal B}$}{Reducing Theorem \ref{co5thm2} to a wall-crossing in ℬ}}
\label{co102}

The next proposition shows that \eq{co5eq31} is equivalent to a wall-crossing between weak stability conditions $(\ac\tau^0_{\bs\mu_{\smash{-1}}},\ac T,\le),(\ac\tau^{\ti\la}_{\bs\mu_{\smash{-1}}},\ac T,\le)$ in $\acB$. It is important that these are both defined {\it from the same\/} $(\tau,T,\le)\in\sS$. This enables us to interpolate continuously between them in the family $(\ac\tau^{s\ti\la}_{\bs\mu_{\smash{-1}}},\ac T,\le)$ for~$s\in[0,1]$.

\begin{prop}
\label{co10prop4}
In the situation of Definition\/ {\rm\ref{co10def1},} there exists a group morphism $\ti\la:K(\A)\ra\R$ with\/ $\ti\la(\al)=0,$ such that if\/ $\be\in R_\al$ with\/ $\ti\la(\be)\ne 0$ then\/ $\bmd{\ti\la(\be)}\ge\ha r(r+1)\rk\al,$ and:
\smallskip

\noindent{\bf(a)} For all\/ $\be\in R_\al,$ there are inverse isomorphisms of moduli stacks:
\e
\xymatrix@C=100pt{
*+[r]{\acM^\ss_{(\be,0)}(\ac\tau_{\bs\mu_{\smash{-1}}}^0)} \ar@<.5ex>[r]^(0.55){\Pi_{\M^\pl}} & *+[l]{\M^\ss_\be(\tau).} \ar@<.5ex>[l]^(0.45){I_{\acM^\pl}} }   
\label{co10eq10}
\e
Similarly, for all\/ $\be\in\ti R_\al$ there are inverse isomorphisms
\e
\xymatrix@C=100pt{*+[r]{\acM^\ss_{(\be,0)}(\ac\tau_{\bs\mu_{\smash{-1}}}^{\ti\la})} \ar@<.5ex>[r]^(0.55){\Pi_{\M^\pl}} & *+[l]{\M^\ss_\be(\ti\tau).} \ar@<.5ex>[l]^(0.45){I_{\acM^\pl}}   }
\label{co10eq11}
\e

\noindent{\bf(b)} If\/ $\al_1,\ldots,\al_n\in R_\al$ with\/ $\al_1+\cdots+\al_n=\al$ then
\e
U((\al_1,0),\ldots,(\al_n,0);\ac\tau_{\bs\mu_{\smash{-1}}}^0,\ac\tau_{\bs\mu_{\smash{-1}}}^{\ti\la})=U(\al_1,\ldots,\al_n;\tau,\ti\tau).
\label{co10eq12}
\e
Hence we may also take
\e
\ti U((\al_1,0),\ldots,(\al_n,0);\ac\tau_{\bs\mu_{\smash{-1}}}^0,\ac\tau_{\bs\mu_{\smash{-1}}}^{\ti\la})=\ti U(\al_1,\ldots,\al_n;\tau,\ti\tau).
\label{co10eq13}
\e
\end{prop}

\begin{proof}
By Definition \ref{co10def1}, $R_\al\subseteq C(\B)_\pe$ is a finite subset with $\al\in R_\al$ and $\tau(\be)=\tau(\al)$, $\M_\be^\ss(\tau)\ne\es$ for all $\be\in R_\al$. Thus Assumption \ref{co5ass2}(d) with $I=R_\al$ gives a group morphism $\ti\la:K(\A)\ra\R$ such that $\ti\la(\al)=0$ and $\ti\la(\be)>0$ (or $\ti\la(\be)>0$) if and only if $\ti\tau(\be)>\ti\tau(\al)$ (or $\ti\tau(\be)<\ti\tau(\al)$) for all $\be\in R_\al$. Multiplying $\ti\la$ by a large positive constant, we can also suppose that if $\be\in R_\al$ with $\ti\la(\be)\ne 0$ then~$\md{\ti\la(\be)}\ge\ha r(r+1)\rk\al$. By \eq{co10eq3} we see that
\e
\ti R_\al=\bigl\{\be\in R_\al:\ti\la(\be)=0\bigr\}.
\label{co10eq14}
\e

For (a), equation \eq{co10eq10} is obvious from the first case of \eq{co10eq9}. Suppose $\ga,\de\in R_\al$ with $\be:=\ga+\de\in\ti R_\al$. Then we have
\ea
\ti\tau(\ga)&\le\ti\tau(\de) && \Longleftrightarrow & \ti\tau(\ga)&\le\ti\tau(\be)\le\ti\tau(\de) && \Longleftrightarrow
\nonumber\\
\ti\tau(\ga)&\le\ti\tau(\al)\le\ti\tau(\de) && \Longleftrightarrow & \ti\la(\ga)&\le 0\le \ti\la(\de) && \Longleftrightarrow 
\nonumber\\
\frac{\ti\la(\ga)}{\rk\ga}&\le \frac{\ti\la(\de)}{\rk\de} && \Longleftrightarrow 
&\ac\tau_{\bs\mu_{\smash{-1}}}^{\ti\la}(\ga,0)&\le\ac\tau_{\bs\mu_{\smash{-1}}}^{\ti\la}(\de,0).
\label{co10eq15}
\ea 
The first step of \eq{co10eq15} holds as $\be=\ga+\de$, the second as $\ti\tau(\be)=\ti\tau(\al)$ since $\be\in\ti R_\al$, the third by choice of $\ti\la$ the fourth as $\ti\la(\ga)+\ti\la(\de)=\ti\la(\be)=0$ and 
$\rk\ga,\ab\rk\de>0$ by Assumption \ref{co5ass1}(g)(iv), and the fifth by \eq{co5eq21} and as $\tau(\ga)=\tau(\de)=\tau(\al)$ by~\eq{co10eq3}.

Let $\be\in\ti R_\al$. We claim that $[E]$ is a $\C$-point of $\M_\be^\ss(\ti\tau)$ if and only if $[E,\bs 0,\bs 0]=I_{\acM^\pl}([E])$ is a $\C$-point of $\acM_{(\be,0)}^\ss(\ac\tau_{\bs\mu_{\smash{-1}}}^{\ti\la})$. To see this, note that $\be\in R_\al$ and $\M_\be^\ss(\ti\tau)\subseteq\M_\be^\ss(\tau)$ by the first part of Proposition \ref{co10prop1}, so it is enough to show that if $[E]\in\M_\be^\ss(\tau)$ then $E$ is $\ti\tau$-semistable if and only if $(E,\bs 0,\bs 0)$ is $\ac\tau_{\bs\mu_{\smash{-1}}}^{\ti\la}$-semistable. By the second part of Proposition \ref{co10prop1}, it is enough to test $\ti\tau$-semistability of $E$ upon subobjects $0\ne E'\subsetneq E$ with $E',E/E'$ $\tau$-semistable and $\lb E'\rb,\lb E/E'\rb\in R_\al$. Setting $\ga=\lb E'\rb$ and $\de=\lb E/E'\rb$, equation \eq{co10eq15} shows that $\ti\tau(\lb E'\rb)\le\ti\tau(\lb E/E'\rb)$ if and only if $\ac\tau_{\bs\mu_{\smash{-1}}}^{\ti\la}(\lb E',\bs 0,\bs 0\rb)\le\ac\tau_{\bs\mu_{\smash{-1}}}^{\ti\la}(\lb E/E',\bs 0,\bs 0\rb)$, proving the claim. Equation \eq{co10eq11} follows from this and $\Pi_{\M^\pl}\ci I_{\acM^\pl}=\id$. This completes part~(a).

For (b), let $\al_1,\ldots,\al_n\in R_\al$ with $\al_1+\cdots+\al_n=\al$, and consider the definitions \eq{co3eq3} of $U((\al_1,0),\ldots,(\al_n,0);\ac\tau_{\bs\mu_{\smash{-1}}}^0,\ac\tau_{\bs\mu_{\smash{-1}}}^{\ti\la})$ and $U(\al_1,\ldots,\al_n;\tau,\ti\tau)$. We claim that the sums are equal, term by term. To see this, note that as $\ac\tau_{\bs\mu_{\smash{-1}}}^0(\al_i,0)\ab=(\tau(\al_i),0)$ and $\tau(\al_i)=\tau(\al)$ for all $i$, the conditions $\ac\tau_{\bs\mu_{\smash{-1}}}^0(\be_i,0)=\ac\tau_{\bs\mu_{\smash{-1}}}^0(\al_j,0)$ and $\tau(\be_i)=\tau(\al_j)$ in \eq{co3eq3} are trivial in both cases. Equation \eq{co10eq15} with $\ga_i,\al-\ga_i$ in place of $\be,\ga$ implies that the conditions $\ac\tau_{\bs\mu_{\smash{-1}}}^{\ti\la}(\ga_i,0)=\ac\tau_{\bs\mu_{\smash{-1}}}^{\ti\la}(\al_1+\cdots+\al_n,0)$ and $\ti\tau(\ga_i)=\ti\tau(\al_1+\cdots+\al_n)$ in \eq{co3eq3} are equivalent. Thus the sums over $l,m$, $a_0,\ldots,a_m$, $b_0,\ldots,b_l$ in \eq{co3eq3} for $U((\al_1,0),\ldots,(\al_n,0);\ac\tau_{\bs\mu_{\smash{-1}}}^0,\ac\tau_{\bs\mu_{\smash{-1}}}^{\ti\la})$ and $U(\al_1,\ldots,\al_n;\tau,\ti\tau)$ are the same. For each such term in the sum and each $i=1,\ldots,l$ we have
\begin{equation*}
S((\be_{b_{i-1}+1},0),\ldots,(\be_{b_i},0);\ac\tau_{\bs\mu_{\smash{-1}}}^0,\ac\tau_{\bs\mu_{\smash{-1}}}^{\ti\la})=S(\be_{b_{i-1}+1},\ldots,\be_{b_i}; \tau,\ti\tau)
\end{equation*}
by Definition \ref{co3def4}, as if $b_{i-1}<j<b_i$ then $\ac\tau_{\bs\mu_{\smash{-1}}}^0(\be_j,0)=\ac\tau_{\bs\mu_{\smash{-1}}}^0(\be_{j+1},0)$, $\tau(\be_j)=\tau(\be_{j+1})$, and since $\ti\tau(\be_{b_{i-1}+1}+\cdots+\be_{b_i})=\ti\tau(\ga_i)=\ti\tau(\al)$, \eq{co10eq15} implies that
\begin{align*}
\ac\tau_{\bs\mu_{\smash{-1}}}^{\ti\la}(\be_{b_{i-1}+1}+\cdots+\be_j,0)&\le \ac\tau_{\bs\mu_{\smash{-1}}}^{\ti\la}(\be_{j+1}+\cdots+\be_{b_i},0)\qquad\Longleftrightarrow\\\ti\tau(\be_{b_{i-1}+1}+\cdots+\be_j)&\le \ti\tau(\be_{j+1}+\cdots+\be_{b_i}).
\end{align*}
Thus \eq{co3eq3} implies equation \eq{co10eq12}, and \eq{co10eq13} follows from this and Theorem \ref{co3thm3}, proving~(b).
\end{proof}

For the rest of Chapter \ref{co10} we will consider only weak stability conditions $(\ac\tau^{s\ti\la}_{\bs\mu_x},\ac T,\le)$ on $\acA$ for $x\in[-1,0]$ and $s\in[0,1]$, as these will be all we need.

\subsection{\texorpdfstring{A Lie algebra ${\cal L}_{{\cal B},\ac{\cal B}}$}{A Lie algebra ℒᵦᵦ}}
\label{co103}

\begin{dfn}
\label{co10def2}
Continue in the situation of Definition \ref{co10def1}. Define a $\Q$-vector space $\cL_{\B,\acB}$ by
\begin{equation*}
\cL_{\B,\acB}=\Q[p_1,\ldots,p_r]\op \bigl(\Q[q_1,\ldots,q_r]\ot \check H_{\rm even}(\M^\pl)\bigr).
\end{equation*}
As a shorthand we write $p^{\bs d}=p_1^{d_1}\cdots p_r^{d_r}$ and $q^{\bs e}=q_1^{e_1}\cdots q_r^{e_r}$ for $\bs d,\bs e\in\N^r$. Then $\cL_{\B,\acB}$ is spanned by vectors $p^{\bs d}$ and $q^{\bs e}\ot\eta_\al$, where $\bs d,\bs e\in\N^r$ and $\eta_\al\in\check H_{\rm even}(\M_\al^\pl)\subseteq\check H_{\rm even}(\M^\pl)$ for $\al\in C(\B)\amalg\{0\}$. Define a $\Q$-bilinear Lie bracket $[\,,\,]:\cL_{\B,\acB}\t\cL_{\B,\acB}\ra\cL_{\B,\acB}$ on these generators by
\e
\begin{aligned}{}
&[p^{\bs d},p^{\bs d'}]=\bigl(\ac\chi((0,\bs d),(0,\bs d'))-\ac\chi((0,\bs d'),(0,\bs d))\bigr)p^{\bs d+\bs d'},\\
&[p^{\bs d},q^{\bs e'}\ot\eta_{\al'}]=\bigl(\ac\chi((0,\bs d),(\al',\bs e'))-\ac\chi((\al',\bs e'),(0,\bs d))\bigr)q^{\bs d+\bs e'}\ot\eta_{\al'},\\
&[q^{\bs e}\ot\eta_\al,p^{\bs d'}]=\bigl(\ac\chi((\al,\bs e),(0,\bs d'))-\ac\chi((0,\bs d'),(\al,\bs e))\bigr)q^{\bs d'+\bs e}\ot\eta_\al,\\
&[q^{\bs e}\ot\eta_\al,q^{\bs e'}\ot\eta_{\al'}]=q^{\bs e+\bs e'}\ot[\eta_\al,\eta_{\al'}],
\end{aligned}
\label{co10eq16}
\e
where $[\eta_\al,\eta_{\al'}]$ is the Lie bracket in $\check H_{\rm even}(\M^\pl)$. An easy but tedious calculation shows that $[\,,\,]$ is antisymmetric and satisfies the Jacobi identity.
\end{dfn}

Because of Lemma \ref{co9lem1}, all the wall-crossing type formulae in Chapter \ref{co10}, such as \eq{co10eq18}, \eq{co10eq21}, \ldots\ in $\check H_{\rm even}(\acM^\pl)$ or $\cL_{\B,\acB}$ below, have only finitely many terms, but for brevity we will usually not say this.

\begin{dfn}
\label{co10def3}
{\bf(a)} Motivated by \eq{co10eq9}, define elements $[\acM_{(\be,\bs d)}^\ss(\ac\tau^0_{\bs\mu_{\smash{-1}}})]_\cL$ in $\cL_{\B,\acB}$ for all $(\be,\bs d)\in S_\al$ by
\ea
[\acM_{(\be,\bs d)}^\ss(\ac\tau^0_{\bs\mu_{\smash{-1}}})]_\cL=
\begin{cases} 
q^{\bs 0}\ot[\M^\ss_\be(\tau)]_\inv, & \bs d=0, \\
p_ap_{a+1}\cdots p_b, &  \begin{aligned}[h]&\text{$d_i\!=\!1$ if\/ $a\!\le\! i\!\le\! b,$ $d_i=0$ otherwise} \\
&\text{for some $1\!\le\! a\!\le\! b\!\le\! n,$ and\/ $\be\!=\!0,$} \end{aligned}  \\
0, & \text{otherwise.}
\end{cases}
\nonumber\\[-17pt]
\label{co10eq17}
\ea
These should be thought of as analogues of invariants $[\acM_{(\be,\bs d)}^\ss(\ac\tau^0_{\bs\mu_{\smash{-1}}})]_\inv$ in Theorem \ref{co5thm1}, but defined `by hand'. Here in the second case of \eq{co10eq9} with $\acM_{(\be,\bs d)}^\ss(\ac\tau^0_{\bs\mu_{\smash{-1}}})\cong[*/\PGL(n,\C)]$, when $n>1$ the virtual dimension is $2-2n^2<0$, so we set the invariants to be zero.

Now for any $x\in[-1,0]$ and $s\in[0,1]$ define elements $[\acM_{(\be,\bs d)}^\ss(\ac\tau^{s\ti\la}_{\bs\mu_x})]_\cL$ in $\cL_{\B,\acB}$ for all $(\be,\bs d)\in S_\al$ by
\ea
&[\acM_{(\be,\bs d)}^\ss(\ac\tau^{s\ti\la}_{\bs\mu_x})]_\cL= 
\label{co10eq18}\\
&\sum_{\begin{subarray}{l}n\ge 1,\;(\al_1,\bs d_1),\ldots,(\al_n,\bs d_n)\in S_\al:\\ (\al_1,\bs d_1)+\cdots+(\al_n,\bs d_n)=(\be,\bs d),\\ 
o_{\al_1}+\cdots+o_{\al_n}=o_\al \end{subarray}\!\!\!\!\!\!\!\!\!\!\!\!\!\!\!\!\!\!\!\!\!\!\!\!\!\!\!\!} \!\!\!\!\!\!\!
\begin{aligned}[t]
\ti U((\al_1,\bs d_1),\ldots,&(\al_n,\bs d_n);\ac\tau^0_{\bs\mu_{\smash{-1}}},\ac\tau^{s\ti\la}_{\bs\mu_x})\cdot\bigl[\bigl[\cdots\bigl[[\acM_{(\al_1,\bs d_1)}^\ss(\ac\tau^0_{\bs\mu_{\smash{-1}}})]_\cL,\\
&
[\acM_{(\al_2,\bs d_2)}^\ss(\ac\tau^0_{\bs\mu_{\smash{-1}}})]_\cL\bigr],\ldots\bigr],[\acM_{(\al_n,\bs d_n)}^\ss(\ac\tau^0_{\bs\mu_{\smash{-1}}})]_\cL\bigr].
\end{aligned}
\nonumber
\ea
If $x=-1$ and $s=0$  this recovers the same elements \eq{co10eq17} by~\eq{co3eq4}.
\smallskip

\noindent{\bf(b)} Similarly, define $[\acM_{(\be,\bs d)}^\ss(\ac\tau^{\ti\la}_{\bs\mu_{\smash{-1}}})]_{\smash{\ti{\cal L}}}\in\cL_{\B,\acB}$ for all $(\be,\bs d)\in\ti S_\al$ by 
\ea
[\acM_{(\be,\bs d)}^\ss(\ac\tau^{\ti\la}_{\bs\mu_{\smash{-1}}})]_{\smash{\ti{\cal L}}}=
\begin{cases} 
q^{\bs 0}\ot[\M^\ss_\be(\ti\tau)]_\inv, & \bs d=0, \\
p_ap_{a+1}\cdots p_b, &  \begin{aligned}[h]&\text{$d_i\!=\!1$ if\/ $a\!\le\! i\!\le\! b,$ $d_i=0$ otherwise} \\
&\text{for some $1\!\le\! a\!\le\! b\!\le\! n,$ and\/ $\be\!=\!0,$} \end{aligned}  \\
0, & \text{otherwise.}
\end{cases}
\nonumber\\[-17pt]
\label{co10eq19}
\ea
For $x\in[-1,0]$, define $[\acM_{(\be,\bs d)}^\ss(\ac\tau^{\ti\la}_{\bs\mu_x})]_{\smash{\ti{\cal L}}}$ in $\cL_{\B,\acB}$ for all $(\be,\bs d)\in\ti S_\al$ by
\ea
&[\acM_{(\be,\bs d)}^\ss(\ac\tau^{\ti\la}_{\bs\mu_x})]_{\smash{\ti{\cal L}}}= 
\label{co10eq20}\\
&\sum_{\begin{subarray}{l}n\ge 1,\;(\al_1,\bs d_1),\ldots,(\al_n,\bs d_n)\in\ti S_\al:\\ (\al_1,\bs d_1)+\cdots+(\al_n,\bs d_n)=(\be,\bs d),\\ 
o_{\al_1}+\cdots+o_{\al_n}=o_\al \end{subarray}\!\!\!\!\!\!\!\!\!\!\!\!\!\!\!\!\!\!\!\!\!\!\!\!\!\!\!\!} \!\!\!\!\!\!\!
\begin{aligned}[t]
\ti U((\al_1,\bs d_1),\ldots,&(\al_n,\bs d_n);\ac\tau^{\ti\la}_{\bs\mu_{\smash{-1}}},\ac\tau^{\ti\la}_{\bs\mu_x})\cdot\bigl[\bigl[\cdots\bigl[[\acM_{(\al_1,\bs d_1)}^\ss(\ac\tau^{\ti\la}_{\bs\mu_{\smash{-1}}})]_{\smash{\ti{\cal L}}},\\
&
[\acM_{(\al_2,\bs d_2)}^\ss(\ac\tau^{\ti\la}_{\bs\mu_{\smash{-1}}})]_{\smash{\ti{\cal L}}}\bigr],\ldots\bigr],[\acM_{(\al_n,\bs d_n)}^\ss(\ac\tau^{\ti\la}_{\bs\mu_{\smash{-1}}})]_{\smash{\ti{\cal L}}}\bigr].
\end{aligned}
\nonumber
\ea
If $x=-1$ this recovers the same elements \eq{co10eq19} by~\eq{co3eq4}.
\end{dfn}

Equation \eq{co3eq5} of Theorem \ref{co3thm2} and \eq{co10eq18}, \eq{co10eq20} imply:

\begin{cor}
\label{co10cor1}
Let\/ $x,x'\in[-1,0],$ $s,s'\in[0,1]$ and\/ $(\be,\bs d)\in S_\al$. Then in the Lie algebra $\cL_{\B,\acB}$ we~have
\ea
&[\acM_{(\be,\bs d)}^\ss(\ac\tau^{s\ti\la}_{\bs\mu_x})]_\cL= 
\label{co10eq21}\\
&\sum_{\begin{subarray}{l}n\ge 1,\;(\al_1,\bs d_1),\ldots,(\al_n,\bs d_n)\in S_\al:\\ (\al_1,\bs d_1)+\cdots+(\al_n,\bs d_n)=(\be,\bs d),\\
o_{\al_1}+\cdots+o_{\al_n}=o_\be \end{subarray}\!\!\!\!\!\!\!\!\!\!\!\!\!\!\!\!\!\!\!\!\!\!\!\!\!\!\!\!} \!\!\!\!\!\!\!
\begin{aligned}[t]
\ti U((\al_1,\bs d_1),\ldots,&(\al_n,\bs d_n);\ac\tau^{s'\ti\la}_{\bs\mu_{x'}},\ac\tau^{s\ti\la}_{\bs\mu_x})\cdot\bigl[\bigl[\cdots\bigl[[\acM_{(\al_1,\bs d_1)}^\ss(\ac\tau^{s'\ti\la}_{\bs\mu_{x'}})]_\cL,\\
&
[\acM_{(\al_2,\bs d_2)}^\ss(\ac\tau^{s'\ti\la}_{\bs\mu_{x'}})]_\cL\bigr],\ldots\bigr],[\acM_{(\al_n,\bs d_n)}^\ss(\ac\tau^{s'\ti\la}_{\bs\mu_{x'}})]_\cL\bigr].
\end{aligned}
\nonumber
\ea

If also\/ $(\be,\bs d)\in\ti S_\al$ we have
\ea
&[\acM_{(\be,\bs d)}^\ss(\ac\tau^{\ti\la}_{\bs\mu_x})]_{\smash{\ti{\cal L}}}= 
\label{co10eq22}\\
&\sum_{\begin{subarray}{l}n\ge 1,\;(\al_1,\bs d_1),\ldots,(\al_n,\bs d_n)\in\ti S_\al:\\ (\al_1,\bs d_1)+\cdots+(\al_n,\bs d_n)=(\be,\bs d),\\
o_{\al_1}+\cdots+o_{\al_n}=o_\be \end{subarray}\!\!\!\!\!\!\!\!\!\!\!\!\!\!\!\!\!\!\!\!\!\!\!\!\!\!\!\!} \!\!\!\!\!\!\!
\begin{aligned}[t]
\ti U((\al_1,\bs d_1),\ldots,&(\al_n,\bs d_n);\ac\tau^{\ti\la}_{\bs\mu_{x'}},\ac\tau^{\ti\la}_{\bs\mu_x})\cdot\bigl[\bigl[\cdots\bigl[[\acM_{(\al_1,\bs d_1)}^\ss(\ac\tau^{\ti\la}_{\bs\mu_{x'}})]_{\smash{\ti{\cal L}}},\\
&
[\acM_{(\al_2,\bs d_2)}^\ss(\ac\tau^{\ti\la}_{\bs\mu_{x'}})]_{\smash{\ti{\cal L}}}\bigr],\ldots\bigr],[\acM_{(\al_n,\bs d_n)}^\ss(\ac\tau^{\ti\la}_{\bs\mu_{x'}})]_{\smash{\ti{\cal L}}}\bigr].
\end{aligned}
\nonumber
\ea
\end{cor}

The proofs below involve proving the following for some $x,s,(\be,\bs d)$.

\begin{claim}
\label{co10cla1}
{\bf(a)} Let\/ $x\in[-1,0],$ $s\in[0,1]$ and\/ $(\be,\bs d)\in S_\al$ with\/ $\acM_{(\be,\bs d)}^\rst(\ac\tau^{s\ti\la}_{\bs\mu_x})=\acM_{(\be,\bs d)}^\ss(\ac\tau^{s\ti\la}_{\bs\mu_x})\subseteq\acM^\pl_{(\be,\bs d)},$ giving a virtual class\/ $[\acM_{(\be,\bs d)}^\ss(\ac\tau^{s\ti\la}_{\bs\mu_x})]_\virt\in\check H_{\rm even}(\acM^\pl_{(\be,\bs d)})$. Then, possibly under extra assumptions, for $\ac\Up_{(\be,\bs d)}(\ac\tau^{s\ti\la}_{\bs\mu_x})$ as in {\rm\eq{co10eq6},} 
\e
q^{\bs d}\ot\ac\Up_{(\be,\bs d)}(\ac\tau^{s\ti\la}_{\bs\mu_x})=[\acM_{(\be,\bs d)}^\ss(\ac\tau^{s\ti\la}_{\bs\mu_x})]_\cL
\qquad\text{in\/ $\cL_{\B,\acB}$.}
\label{co10eq23}
\e
If this holds we will say that \begin{bfseries}Claim \ref{co10cla1}(a) holds for\end{bfseries}~$x,s,(\be,\bs d)$.
\smallskip

\noindent{\bf(b)} Let\/ $x\in[-1,0]$ and\/ $(\be,\bs d)\in\ti S_\al$ with\/ $\acM_{(\be,\bs d)}^\rst(\ac\tau^{\ti\la}_{\bs\mu_x})=\acM_{(\be,\bs d)}^\ss(\ac\tau^{\ti\la}_{\bs\mu_x})\subseteq\acM^\pl_{(\be,\bs d)}$. Then, possibly under extra assumptions
\e
q^{\bs d}\ot\ac\Up_{(\be,\bs d)}(\ac\tau^{\ti\la}_{\bs\mu_x})=[\acM_{(\be,\bs d)}^\ss(\ac\tau^{\ti\la}_{\bs\mu_x})]_{\smash{\ti{\cal L}}}
\qquad\text{in\/ $\cL_{\B,\acB}$.}
\label{co10eq24}
\e
If this holds we will say that \begin{bfseries}Claim \ref{co10cla1}(b) holds for\end{bfseries}\/~$x,(\be,\bs d)$.
\end{claim}

\begin{prop}
\label{co10prop5}
If\/ $(\be,\bs d)\in\ti S_\al$ then in $\cL_{\B,\acB}$ we have
\ea
&[\acM_{(\be,\bs d)}^\ss(\ac\tau^{\ti\la}_{\bs\mu_0})]_\cL= 
\label{co10eq25}\\
&\sum_{\begin{subarray}{l}n\ge 1,\;(\al_1,\bs d_1),\ldots,(\al_n,\bs d_n)\in\ti S_\al:\\ (\al_1,\bs d_1)+\cdots+(\al_n,\bs d_n)=(\be,\bs d),\\
o_{\al_1}+\cdots+o_{\al_n}=o_\be \end{subarray}\!\!\!\!\!\!\!\!\!\!\!\!\!\!\!\!\!\!\!\!\!\!\!\!\!\!\!\!} \!\!\!\!\!\!\!
\begin{aligned}[t]
\ti U((\al_1,\bs d_1),\ldots,&(\al_n,\bs d_n);\ac\tau^{\ti\la}_{\bs\mu_{\smash{-1}}},\ac\tau^{\ti\la}_{\bs\mu_0})\cdot\bigl[\bigl[\cdots\bigl[[\acM_{(\al_1,\bs d_1)}^\ss(\ac\tau^{\ti\la}_{\bs\mu_{\smash{-1}}})]_\cL,\\
&
[\acM_{(\al_2,\bs d_2)}^\ss(\ac\tau^{\ti\la}_{\bs\mu_{\smash{-1}}})]_\cL\bigr],\ldots\bigr],[\acM_{(\al_n,\bs d_n)}^\ss(\ac\tau^{\ti\la}_{\bs\mu_{\smash{-1}}})]_\cL\bigr].
\end{aligned}
\nonumber
\ea
\end{prop}

\begin{proof} Equation \eq{co10eq25} is \eq{co10eq21} with $x=0$, $x'=-1$ and $s=s'=1$, except that we require $(\al_i,\bs d_i)\in\ti S_\al$ rather than $(\al_i,\bs d_i)\in S_\al$. Suppose that $n\ge 1$, $(\al_1,\bs d_1),\ldots,(\al_n,\bs d_n)$ is a term in the sum \eq{co10eq21}, such that at least one $(\al_i,\bs d_i)$ lies in $S_\al\sm\ti S_\al$. By \eq{co10eq4} and \eq{co10eq14} this means that $\ti\la(\al_i)\ne 0$ for at least one $i$. As $(\be,\bs d)\in\ti S_\al$ we have $\ti\la(\al_1)+\cdots+\ti\la(\al_n)=\ti\la(\be)=0$. Thus there exist $i\ne j=1,\ldots,n$ with $\ti\la(\al_i)>0$ and~$\ti\la(\al_j)<0$.

Write $I=\bigl\{i\in\{1,\ldots,n\}:\ti\la(\al_i)>0\bigr\}$. Then $\es\ne I\subsetneq\{1,\ldots,n\}$. Using equations \eq{co5eq21}, \eq{co10eq4}, \eq{co10eq5} and the fact from Proposition \ref{co10prop4} that if $\ga\in R_\al$ with\/ $\ti\la(\ga)\ne 0$ then $\bmd{\ti\la(\ga)}\ge\ha r(r+1)\rk\al$, we can show that
\begin{align*}
\ac\tau^{\ti\la}_{\bs\mu_{\smash{-1}}}(\al_i,\bs d_i)&=\begin{cases} (t,x)\;\text{with}\; x>0 & \text{if $i\in I$,} \\
(t,x)\;\text{with}\; x<0 & \text{if $i\notin I$ and $\al_i\ne 0$,} \\
(-\iy,x) & \text{if $i\notin I$ and $\al_i=0$,}\end{cases}	
\\
\ac\tau^{\ti\la}_{\bs\mu_0}(\al_i,\bs d_i)&=(t,x)\;\text{with}\; x\ge\ha r(r+1) \; \text{if $i\in I$,}
\\
\ac\tau^{\ti\la}_{\bs\mu_0}(\be,\bs d)&=(t,x)\;\text{with}\; x<\ha r(r+1).
\end{align*}
Therefore $\ac\tau^{\ti\la}_{\bs\mu_{\smash{-1}}}(\al_i,\bs d_i)>\ac\tau^{\ti\la}_{\bs\mu_{\smash{-1}}}(\al_j,\bs d_j)$ if $i\in I$ and $j\in\{1,\ldots,n\}\sm I$, and $\ac\tau^{\ti\la}_{\bs\mu_0}(\al_i,\bs d_i)>\ac\tau^{\ti\la}_{\bs\mu_0}(\be,\bs d)$ if $i\in I$. Thus the last part of Proposition \ref{co3prop3} shows that $\ti U((\al_1,\bs d_1),\ldots,(\al_n,\bs d_n);\ac\tau^{\ti\la}_{\bs\mu_{\smash{-1}}},\ac\tau^{\ti\la}_{\bs\mu_0})=0$. Hence all terms in \eq{co10eq21} with at least one $(\al_i,\bs d_i)$ in $S_\al\sm\ti S_\al$ are zero, and \eq{co10eq25} follows.
\end{proof}

\subsection{Some vanishing results}
\label{co104}

Here is an analogue of Proposition \ref{co10prop2}(b):

\begin{prop}
\label{co10prop6}
Suppose\/ $x\in[-1,0],$ $s\in[0,1]$ and\/ $(\be,\bs d)\in S_\al$ with\/ $[\acM_{(\be,\bs d)}^\ss(\ac\tau^{s\ti\la}_{\bs\mu_x})]_\cL\ab\ne 0$. Then:
\begin{itemize}
\setlength{\itemsep}{0pt}
\setlength{\parsep}{0pt}
\item[{\bf(i)}] If\/ $\be=0$ then for some then for some $1\le a\le b\le r$ we have $d_i=1$ if\/ $a\le i\le b$ and $d_i=0$ otherwise.
\item[{\bf(ii)}] If\/ $\be\ne 0$ then $d_1\le d_2\le \cdots\le d_r$. If also $d_{i+1}\le d_i+1$ for $1\le i<r$ then~$d_r\le \la_k(\be)$.
\item[{\bf(iii)}] If\/ $\be\ne 0$ and\/ $\frac{1}{r-a}(\mu_{a+1}+\cdots+\mu_r)+x\le 0$ for some\/ $1\le a<r$ then\/~$d_a=d_{a+1}=\cdots=d_r$.
\item[{\bf(iv)}] If\/ $\be\ne 0$ and\/ $\frac{1}{r}(\mu_1+\cdots+\mu_r)+x\le 0$ then\/~$\bs d=0$.
\end{itemize}
The analogues hold if\/ $x\in[-1,0],$ $(\be,\bs d)\in\ti S_\al$ and\/ $[\acM_{(\be,\bs d)}^\ss(\ac\tau^{\ti\la}_{\bs\mu_x})]_{\smash{\ti{\cal L}}}\ne 0$.
\end{prop}

\begin{proof} For (i), we see from \eq{co5eq21} that $\ac\tau^{s\ti\la}_{\bs\mu_x}(0,\bs d)\le\ac\tau^{s\ti\la}_{\bs\mu_x}(0,\bs e)$ if and only if $\ac\tau^0_{\bs\mu_{\smash{-1}}}(0,\bs d)\le\ac\tau^0_{\bs\mu_{\smash{-1}}}(0,\bs e)$, so using Definition \ref{co3def4} we can show that
\begin{equation*}
U((0,\bs d_1),\ldots,(0,\bs d_n);\ac\tau^0_{\bs\mu_{\smash{-1}}},\ac\tau^{s\ti\la}_{\bs\mu_x})=\begin{cases} 1, & n=1, \\ 0, & n>1. \end{cases}
\end{equation*}
Thus Theorem \ref{co3thm3} and \eq{co10eq18} imply that $[\acM_{(0,\bs d)}^\ss(\ac\tau^{s\ti\la}_{\bs\mu_x})]_\cL=[\acM_{(0,\bs d)}^\ss(\ac\tau^0_{\bs\mu_{\smash{-1}}})]_\cL$, so (i) follows from \eq{co10eq17}. 

For (ii)--(iv), suppose for a contradiction that there exist $x\in[-1,0],$ $s\in[0,1]$ and $(\be,\bs d)\in S_\al$ with $[\acM_{(\be,\bs d)}^\ss(\ac\tau^{s\ti\la}_{\bs\mu_x})]_\cL\ne 0$, but at least one of (ii)--(iv) do not hold. Choose such $x,s,(\be,\bs d)$ with $\md{\bs d}=d_1+\cdots+d_r$ least, and with $\rk\be$ least for this fixed $\md{\bs d}$. We have $\be\ne 0$ and $\bs d\ne 0$ as one of (ii)--(iv) do not hold. 

Consider $[\acM_{(\be,\bs d)}^\ss(\ac\tau^{us\ti\la}_{\bs\mu_{u(x+1)-1}})]_\cL$ for $u\in[0,1]$. This is zero for $u=0$ by \eq{co10eq17}, and nonzero when $u=1$, and can be discontinuous in $u$ at only finitely many $u\in[0,1]$. Hence there exists $u_0\in[0,1]$ such that for any small $\de>0$, there exist $u',u''\in[0,1]$ with $u_0-\de<u'\le u_0\le u''<u_0+\de$ such that  
\begin{equation*}
[\acM_{(\be,\bs d)}^\ss(\ac\tau^{u's\ti\la}_{\bs\mu_{x'}})]_\cL=0,\qquad [\acM_{(\be,\bs d)}^\ss(\ac\tau^{u''s\ti\la}_{\bs\mu_{x''}})]_\cL\ne 0,
\end{equation*}
where $x_0\!=\!u_0(x\!+\!1)\!-\!1$, $x'\!=\!u'(x\!+\!1)\!-\!1$, $x''\!=\!u''(x\!+\!1)\!-\!1$. By \eq{co10eq21}  we have
\ea
&[\acM_{(\be,\bs d)}^\ss(\ac\tau^{u''s\ti\la}_{\bs\mu_{x''}})]_\cL= 
\label{co10eq26}\\
&\sum_{\begin{subarray}{l}n\ge 1,\;(\al_1,\bs d_1),\ldots,(\al_n,\bs d_n)\in
S_\al:\\ (\al_1,\bs d_1)+\cdots+(\al_n,\bs d_n)=(\be,\bs d),\\
o_{\al_1}+\cdots+o_{\al_n}=o_\be \end{subarray}\!\!\!\!\!\!\!\!\!\!\!\!\!\!\!\!\!\!\!\!\!\!\!\!\!\!\!\!} \!\!\!\!\!\!\!\!\!\!\!
\begin{aligned}[t]
\ti U((\al_1,\bs d_1),\ldots,&(\al_n,\bs d_n);\ac\tau^{u's\ti\la}_{\bs\mu_{x'}},\ac\tau^{u''s\ti\la}_{\bs\mu_{x''}})\cdot\bigl[\bigl[\cdots\bigl[[\acM_{(\al_1,\bs d_1)}^\ss(\ac\tau^{u's\ti\la}_{\bs\mu_{x'}})]_\cL,\\
&
[\acM_{(\al_2,\bs d_2)}^\ss(\ac\tau^{u's\ti\la}_{\bs\mu_{x'}})]_\cL\bigr],\ldots\bigr],[\acM_{(\al_n,\bs d_n)}^\ss(\ac\tau^{u's\ti\la}_{\bs\mu_{x'}})]_\cL\bigr].
\end{aligned}
\nonumber
\ea
As $[\acM_{(\be,\bs d)}^\ss(\ac\tau^{u''s\ti\la}_{\bs\mu_{x''}})]_\cL\ne 0$, there must be some $(\al_1,\bs d_1),\ldots,(\al_n,\bs d_n)$ giving a nonzero term in the sum \eq{co10eq26}. We cannot have $n=1$ as then $(\al_1,\bs d_1)=(\be,\bs d)$ but $[\acM_{(\be,\bs d)}^\ss(\ac\tau^{u's\ti\la}_{\bs\mu_{x'}})]_\cL=0$, so $n>1$. Thus each $(\al_i,\bs d_i)$ has $\md{\bs d_i}<\md{\bs d}$ or $\md{\bs d_i}=\md{\bs d}$ and $\rk\al_i<\rk\be$, so by choice of $(\be,\bs d)$ minimal, we see that (ii)--(iv) above hold with $(\al_i,\bs d_i)$ in place of $(\be,\bs d)$ for~$i=1,\ldots,n$. 

First suppose that $\al_i\ne 0$ for $i=1,\ldots,n$. Then (ii)--(iv) hold for $(\al_i,\bs d_i),x'$ for all $i$,  so they hold for $(\al_i,\bs d_i),x$ for all $i$ as $x'<x$. Summing the (in)equalities on $\bs d_{i,j}$ in (ii)--(iv) for $(\al_i,\bs d_i),x$ over $i=1,\ldots,n$ shows that (ii)--(iv) hold for $(\be,\bs d),x$, a contradiction. Hence $\al_i=0$ for some~$i=1,\ldots,n$.

Suppose there exists $i=1,\ldots,n$ with $\al_i=0$ and $\sum_{j=1}^r(\mu_j+x_0)d_{i,j}<0$. Write $I$ for the set of $i=1,\ldots,n$ such that $\al_i=0$ and $\sum_{j=1}^r(\mu_j+x_0)d_{i,j}$ is minimal amongst such $i$. Then $\es\ne I\subsetneq\{1,\ldots,n\}$. As $x',x''$ are close to $x_0$ for small $\de>0$ we also have $\sum_{j=1}^r(\mu_j+x')d_{i,j}<0$ and $\sum_{j=1}^r(\mu_j+x'')d_{i,j}<0$. Using \eq{co5eq21} we can now show that $\ac\tau^{u's\ti\la}_{\bs\mu_{x'}}(\al_i,\bs d_i)<\ac\tau^{u's\ti\la}_{\bs\mu_{x'}}(\al_j,\bs d_j)$ if $i\in I$ and $j\notin I$, and $\ac\tau^{u''s\ti\la}_{\bs\mu_{x''}}(\al_i,\bs d_i)<\ac\tau^{u''s\ti\la}_{\bs\mu_{x''}}(\be,\bs d)$ if $i\in I$. So the first part of Proposition \ref{co3prop3} gives $\ti U((\al_1,\bs d_1),\ldots,(\al_n,\bs d_n);\ac\tau^{u's\ti\la}_{\bs\mu_{x'}},\ac\tau^{u''s\ti\la}_{\bs\mu_{x''}})=0$, a contradiction. 

Similarly, if there exists $i=1,\ldots,n$ with $\al_i=0$ and $\sum_{j=1}^r(\mu_j+x_0)d_{i,j}>0$, we write $I$ for the set of $i=1,\ldots,n$ with $\al_i=0$ and $\sum_{j=1}^r(\mu_j+x_0)d_{i,j}$ maximal, and the second part of Proposition \ref{co3prop3} gives a contradiction. 

This proves that if $i=1,\ldots,n$ with $\al_i=0$ then $\sum_{j=1}^r(\mu_j+x_0)d_{i,j}=0$. The possibilities for $\bs d_i$ are listed in part (i). As $\mu_1,\ldots,\mu_r$ are generic, only one of these possibilities can satisfy $\sum_{j=1}^r(\mu_j+x_0)d_{i,j}=0$. Hence there are unique $1\le a\le b\le r$ such that if $\al_i=0$ for $i=1,\ldots,n$ then $d_{i,j}=1$ if $a\le j\le b$ and $d_{i,j}=0$ otherwise, and $x_0$ is determined by
\e
\ts\frac{1}{b-a+1}\sum_{j=a}^b\mu_j+x_0=0.
\label{co10eq27}
\e
This, $\mu_1>\cdots>\mu_r$ and $x'<x_0$ imply that $\frac{1}{r-a+1}(\mu_a+\cdots+\mu_r)+x'<0$. As the $(\al_i,\bs d_i),x'$ satisfy (ii)--(iv), with $\al_i\ne 0$ for some $i$, we see that $a>1$, as $a=1$ would contradict (iv), and
\e
d_{i,a-1}=d_{i,a}=\cdots=d_{i,r} \qquad\text{if $\al_i\ne 0$.}
\label{co10eq28}
\e

Divide into cases:
\begin{itemize}
\setlength{\itemsep}{0pt}
\setlength{\parsep}{0pt}
\item[(a)] $\al_1=\al_2=0$;
\item[(b)] $b<r$, $\al_1=0$ and $\al_2\ne 0$; 
\item[(c)] $b<r$, $\al_i\ne 0$ for $i=1,\ldots,m<n$ and $\al_{m+1}=0$;
\item[(d)] $b=r$ and $d_r\le\la_k(\be)$;
\item[(e)] $b=r$, $d_r>\la_k(\be)$, and $\al_l=0$ for exactly one $l=1,\ldots,n$; and 
\item[(f)] $b=r$, $d_r>\la_k(\be)$, and $\al_l=0$ for more than one $l=1,\ldots,n$,
\end{itemize}
In case (a) $\bs d_1=\bs d_2$ by the argument above, so we see that
\e
\bigl[[\acM_{(\al_1,\bs d_1)}^\ss(\ac\tau^{u's\ti\la}_{\bs\mu_{x'}})]_\cL,[\acM_{(\al_2,\bs d_2)}^\ss(\ac\tau^{u's\ti\la}_{\bs\mu_{x'}})]_\cL\bigr]=0,
\label{co10eq29}
\e
contradicting that $(\al_1,\bs d_1),\ldots,(\al_n,\bs d_n)$ give a nonzero term in~\eq{co10eq26}.

In case (b), $d_{1,j}=1$ if $a\le j\le b$ and $d_{1,j}=0$ otherwise, \eq{co10eq28} for $i=2$, and \eq{co10eq2} imply that $\ac\chi((\al_1,\bs d_1),(\al_2,\bs d_2))=\ac\chi((\al_2,\bs d_2),\ab(\al_1,\bs d_1))=0$, so from \eq{co10eq16} and \eq{co10eq19} we deduce that \eq{co10eq29} holds, a contradiction. Here we use $[\acM_{(\al_1,\bs d_1)}^\ss(\ac\tau^{u's\ti\la}_{\bs\mu_{x'}})]_\cL=p^{\bs d_1}$ and the second case of~\eq{co10eq16}.

Case (c) is proved as for (b), but with $(\al_{m+1},\bs d_{m+1}),(\al_1+\cdots+\al_m,\ab\bs d_1+\cdots+\bs d_m)$ in place of $(\al_1,\bs d_1),(\al_2,\bs d_2)$.

In case (d), if $i=1,\ldots,n$ with $\al_i\ne 0$ then $(\al_i,\bs d_i)$ satisfies (ii)--(iv) for $x'$, which imply (ii)--(iv) for $x$ as $x'<x$. And if $i=1,\ldots,n$ with $\al_i=0$ then $d_{i,j}=0$ if $j<a$, $d_{i,j}=1$ for $j\ge a$, so by \eq{co10eq27} with $b=r$ we see that $\bs d_i$ also satisfies the inequalities in (ii)--(iv) for $x$, except $d_r\le\la_k(\be)$. So summing over $i=1,\ldots,n$, we see that $\bs d$ satisfies all the inequalities in (ii)--(iv), including $d_r\le\la_k(\be)$ which we have assumed. Thus (ii)--(iv) hold for $(\be,\bs d)$, a contradiction.

In case (e), let $l$ be unique with $\al_l=0$. As $d_r>\la_k(\be)$, and $d_{l,r}=1$, and $d_{j,r}\le\la_k(\al_j)$ for $j\ne l$, we see that $d_r=\la_k(\be)+1$ and $d_{j,r}=\la_k(\al_j)$ for $j\ne l$. First suppose that $l=1$. Then $d_{1,j}=1$ if $a\le j\le r$ and $d_{1,j}=0$ otherwise, \eq{co10eq28} for $i=2$, $d_{2,r}=\la_k(\al_2)$, and \eq{co10eq2} imply that $\ac\chi((\al_1,\bs d_1),(\al_2,\bs d_2))=\ac\chi((\al_2,\bs d_2),\ab(\al_1,\bs d_1))=0$, so from \eq{co10eq16} and \eq{co10eq19} we deduce that \eq{co10eq29} holds, a contradiction. If $l>1$ then the same argument with $(\al_l,\bs d_l),(\al_1+\cdots+\al_{l-1},\bs d_1+\cdots+\bs d_{l-1})$ in place of $(\al_1,\bs d_1),(\al_2,\bs d_2)$ implies that $\bigl[\cdots\bigl[\acM_{(\al_1,\bs d_1)}^\ss(\ac\tau^{u's\ti\la}_{\bs\mu_{x'}})]_\cL,\ldots\bigr],[\acM_{(\al_l,\bs d_l)}^\ss(\ac\tau^{u's\ti\la}_{\bs\mu_{x'}})]_\cL\bigr]=0$, a contradiction.

In case (f), we have $d_{i,a-1}\le d_{i,a-1}$ if $\al_i\ne 0$ by (ii) for $(\al_i,\bs d_i)$, and $d_{l,a-1}=0$, $d_{l,a}=1$ if $\al_l=0$, and there are at least two such $l$. Summing over $i=1,\ldots,n$ shows that $d_a\ge d_{a-1}+2$. Thus (ii) for $(\be,\bs d)$ does not require $d_r\le \la_k(\be)$. All the rest of (ii)--(iv) for $(\be,\bs d)$ follow by adding (ii)--(iv) for $(\al_i,\bs d_i)$ if $\al_i\ne 0$, and using $d_{l,j}=1$ if $j\ge a$ and $d_{l,j}=0$ otherwise if $\al_l=0$. Thus (ii)--(iv) hold for $(\be,\bs d)$, a contradiction.

We have now excluded all cases, proving (ii)--(iv). For the final part, we use the same proof with $\ac\tau^{\ti\la}_{\bs\mu_{\smash{-1}}},[\acM_{(\be,\bs d)}^\ss(\ac\tau^{\ti\la}_{\bs\mu_x})]_{\smash{\ti{\cal L}}}$ in place of $\ac\tau^0_{\bs\mu_{\smash{-1}}},[\acM_{(\be,\bs d)}^\ss(\ac\tau^{s\ti\la}_{\bs\mu_x})]_\cL$, and require $(\al_i,\bs d_i)\in\ti S_\al$ for all $i$ as in \eq{co10eq20} and~\eq{co10eq22}.
\end{proof}

Combining Propositions \ref{co10prop2}(b) and \ref{co10prop6} yields:

\begin{cor}
\label{co10cor2}
Suppose $x\in[-1,0],$ $s\in[0,1]$ and\/ $(\be,\bs d)\in S_\al$ do not satisfy the conditions in Proposition\/ {\rm\ref{co10prop6}(i)--(iv)}. Then Claim\/ {\rm\ref{co10cla1}(a)} holds for $x,s,(\be,\bs d),$ as both sides of\/ \eq{co10eq23} are zero. The analogue holds for $x\in[-1,0],$ $(\be,\bs d)\in\ti S_\al,$ and Claim\/~{\rm\ref{co10cla1}(b)}.	
\end{cor}

\subsection{Completing the proof of Theorem \ref{co5thm2}}
\label{co105}

The next proposition will be proved by showing that equation \eq{co5eq30} in Theorem \ref{co5thm1} for $[\baM^\ss_{(\be,1)}(\bar\tau^0_1)]_\virt$ is equivalent to equation \eq{co10eq23} in Claim \ref{co10cla1}(a) for $[\acM_{(\be,\bs d)}^\ss(\ac\tau^0_{\bs\mu_x})]_\virt$, for certain $\bs d,x$ with $\baM^\ss_{(\be,1)}(\bar\tau^0_1)\cong \acM_{(\be,\bs d)}^\ss(\ac\tau^0_{\bs\mu_x})$. This motivates \eq{co5eq30}. Note that as $\mu_1>\cdots>\mu_r$, the left hand side of \eq{co10eq30} when $a<r$ is smaller than the right, so there exists $x\in[-1,0]$ satisfying~\eq{co10eq30}.

\begin{prop}
\label{co10prop7}
Suppose\/ $(\be,\bs d)\in S_\al$ with\/ $\be\ne 0,$ and for some\/ $a=1,\ldots,r$ we have\/ $d_i=0$ for\/ $i<a$ and\/ $d_i=1$ for\/ $i\ge a,$ and\/ $x\in[-1,0]$ satisfies
\e
\begin{aligned}
-\ts\frac{1}{r-a+1}(\mu_a+\cdots+\mu_r)&\le x<-\ts\frac{1}{r-a}(\mu_{a+1}+\cdots+\mu_r) && \text{if\/ $a<r,$} \\
-\mu_r&\le x\le 0 && \text{if\/ $a=r.$}
\end{aligned}
\label{co10eq30}
\e
Then\/ $\acM_{(\be,\bs d)}^\rst(\ac\tau^0_{\bs\mu_x})=\acM_{(\be,\bs d)}^\ss(\ac\tau^0_{\bs\mu_x}),$ and Claim\/ {\rm\ref{co10cla1}(a)} holds for $x,$ $s=0,$ $(\be,\bs d)$. If also\/ $(\be,\bs d)\in\ti S_\al$ then Claim\/ {\rm\ref{co10cla1}(b)} holds for\/ $x,(\be,\bs d)$.
\end{prop}

\begin{proof} Write $\baB$ for the category discussed in Example \ref{co5ex1}, which should be distinguished from $\acB$ in Definition \ref{co10def1}. Let $\be,\bs d,a,x$ be as in the proposition. Define a functor $F:\baB\hookra\acB$ to map $(E,V,\rho)\mapsto(E,\bs V,\bs\rho)$ on objects and $(\th,\phi)\mapsto(\th,\bs\phi)$ on morphisms, where
\begin{equation*}
V_i=\begin{cases} 0, & i<a, \\ V, & a\le i\le r, \end{cases} \quad \rho_i=\begin{cases} 0, & i<a, \\ \id_V, & a\le i<r, \\ \rho, & i=r,\end{cases}\quad \phi_i=\begin{cases} 0, & i<a, \\ \phi, & a\le i\le 
r. \end{cases}
\end{equation*}
Define a functor $G:\acB\twoheadrightarrow\baB$ to map $(E,\bs V,\bs\rho)\mapsto(E,V_r,\rho_r)$ on objects and $(\th,\bs\phi)\mapsto(\th,\rho_r)$ on morphisms. Then $G\ci F=\Id$. Here $F$ is an equivalence from $\baB$ to a full subcategory of $\acB$, with inverse~$G$.

We claim that $F_*:\baM^\pl\hookra\acM^\pl$ gives an isomorphism from the moduli space $\baM^\ss_{(\be,1)}(\bar\tau^0_1)$ in Example \ref{co5ex1}, to the moduli space $\acM_{(\be,\bs d)}^\ss(\ac\tau^0_{\bs\mu_x})$ in Definition \ref{co10def1}. Also $\baM^\rst_{(\be,1)}(\bar\tau^0_1)=\baM^\ss_{(\be,1)}(\bar\tau^0_1)$, as in Example \ref{co5ex1}, and $\acM_{(\be,\bs d)}^\rst(\ac\tau^0_{\bs\mu_x})=\acM_{(\be,\bs d)}^\ss(\ac\tau^0_{\bs\mu_x})$. To see this, let $[E,\bs V,\bs\rho]$ be a $\C$-point of $\acM_{(\be,\bs d)}^\ss(\ac\tau^0_{\bs\mu_x})$. When $a<r$, Proposition \ref{co10prop2}(a)(ii) shows $\rho_i:V_i\ra V_{i+1}$ is injective for $i=a,a+1,\ldots,r-1$, so $\rho_i$ is an isomorphism as $\dim V_i=\dim V_{i+1}=1$. Hence, for all $a$, we see that~$(E,\bs V,\bs\rho)\cong F(E,V_r,\rho_r)$.

Thus it is sufficient to show that if $(E,V,\rho)$ lies in class $(\be,1)$ in $\baB$ then $(E,V,\rho)$ is $\ac\tau^0_1$-semistable if and only if $F(E,V,\rho)$ is $\ac\tau^0_{\bs\mu_x}$-semistable, and in both cases semistable implies stable. We do this by considering the possible subobjects of $(E,V,\rho)$ in $\baB$ and $F(E,V,\rho)$ in $\acB$, in an elementary way. 

It turns out that $(E,V,\rho)$ can only be unstable if there exists $0\ne E'\subsetneq E$ with $\tau(\lb E'\rb)=\tau(\lb E/E'\rb)=t$ and $\rho(V)\subseteq F_k(E')\subset F_k(E)$, and then $F(E',V,\rho)\subset F(E,V,\rho)$ destabilizes $F(E,V,\rho)$. There are no other possible subobjects which destabilize $F(E,V,\rho)$, or which make $(E,V,\rho)$ or $F(E,V,\rho)$ strictly semistable. There is one subtle point in proving this: when $a<r$ we should consider subobjects $(E,\bs V',\bs\rho')\subsetneq F(E,V,\rho)$ such that for some $a<b\le r$ we have $V'_i=0$ for $i<b$ and $V'_i=V$ for $i\ge b$. We can show using the first inequality of \eq{co10eq30} that this does not destabilize~$F(E,V,\rho)$.

The isomorphism $F_*:\baM^\ss_{(\be,1)}(\bar\tau^0_1)\ra\acM_{(\be,\bs d)}^\ss(\ac\tau^0_{\bs\mu_x})$ identifies obstruction theories and virtual cycles. It also commutes with the linear maps $\bar\Up$ and $\ac\Up$ in \eq{co5eq15} for $\baB$ and $\acB$. So equations \eq{co5eq30} with $\be$ in place of $\al$ and \eq{co10eq6} yield
\ea
&\ac\Up_{(\be,\bs d)}(\ac\tau^0_{\bs\mu_x})=\bar\Up\bigl([\baM^\ss_{(\be,1)}(\bar\tau^0_1)]_\virt\bigr)
\label{co10eq31}\\
&=\sum_{\begin{subarray}{l}n\ge 1,\;\al_1,\ldots,\al_n\in
R_\al:\\ \al_1+\cdots+\al_n=\be, \; o_{\al_1}+\cdots+o_{\al_n}=o_\be \end{subarray}\!\!\!\!\!\!\!\!\!\!} \!\!\!\!\!\!\!\!\!\!
\begin{aligned}[t]
\frac{(-1)^{n+1}\la_k(\al_1)}{n!}\,\cdot&\bigl[\bigl[\cdots\bigl[[\M_{\al_1}^\ss(\tau)]_\inv,\\
&\; 
[\M_{\al_2}^\ss(\tau)]_\inv\bigr],\ldots\bigr],[\M_{\al_n}^\ss(\tau)]_\inv\bigr].
\end{aligned}
\nonumber
\ea	
Here we use that the conditions on $\al_1,\ldots,\al_n$ in \eq{co5eq30} imply that $\al_i\in R_\al$, and that $\al_i\in R_\al$ implies that $\tau(\al_i)=\tau(\al)$ and~$\M_{\al_i}^\ss(\tau)\ne\es$.

Let us apply $q^{\bs d}\ot-$ to \eq{co10eq31} and interpret it as an equation in $\cL_{\B,\acB}$. As $\ac\chi((0,\bs d),(\al_1,0))=-\la_k(\al_1)$ and $\ac\chi((\al_1,0),(0,\bs d))=0$ by \eq{co10eq2},
equation \eq{co10eq16} gives $[p^{\bs d},q^{\bs 0}\ot [\M_{\al_1}^\ss(\tau)]^k_\inv]=-\la_k(\al_1)\,q^{\bs d}\ot [\M_{\al_1}^\ss(\tau)]^k_\inv$. Hence by \eq{co10eq16}, equation \eq{co10eq31} is equivalent to an equation in the Lie algebra~$\cL_{\B,\acB}$:
\begin{align*}
&q^{\bs d}\ot\ac\Up_{(\be,\bs d)}(\ac\tau^0_{\bs\mu_x})=\\
&\sum_{\begin{subarray}{l}n\ge 1,\;\al_1,\ldots,\al_n\in
R_\al:\\ \al_1+\cdots+\al_n=\be, \; o_{\al_1}+\cdots+o_{\al_n}=o_\be \end{subarray}} \!\!\!\!\!\!\!\!\!\!\!\begin{aligned}[t]
\frac{(-1)^n}{n!}\,\cdot&\bigl[\bigl[\cdots\bigl[\bigl[p^{\bs d},q^{\bs 0}\ot [\M_{\al_1}^\ss(\tau)]_\inv\bigr],\\
&\; 
q^{\bs 0}\ot [\M_{\al_2}^\ss(\tau)]_\inv\bigr],\ldots\bigr],q^{\bs 0}\ot [\M_{\al_n}^\ss(\tau)]_\inv\bigr].
\end{aligned}
\end{align*}
By \eq{co10eq17}, as $d_i=0$ for $i<a$ and $d_i=1$ for $a\le i\le m$, this may be rewritten
\ea
&q^{\bs d}\ot\ac\Up_{(\be,\bs d)}(\ac\tau^0_{\bs\mu_x})=
\label{co10eq32}\\
&\sum_{\begin{subarray}{l}n\ge 1,\;\al_1,\ldots,\al_n\in
R_\al:\\ \al_1+\cdots+\al_n=\be, \; o_{\al_1}+\cdots+o_{\al_n}=o_\be \end{subarray}} \!\!\!\!\!\!\!\!\!\!\!\begin{aligned}[t]
\frac{(-1)^n}{n!}\,\cdot&\bigl[\bigl[\cdots\bigl[\bigl[[\acM_{(0,\bs d)}^\ss(\ac\tau^0_{\bs\mu_{\smash{-1}}})]_\cL,[\acM_{(\al_1,0)}^\ss(\ac\tau^0_{\bs\mu_{\smash{-1}}})]_\cL\bigr],\\
&\; 
[\acM_{(\al_2,0)}^\ss(\ac\tau^0_{\bs\mu_{\smash{-1}}})]_\cL\bigr],\ldots\bigr],[\acM_{(\al_n,0)}^\ss(\ac\tau^0_{\bs\mu_{\smash{-1}}})]_\cL\bigr].
\end{aligned}
\nonumber
\ea	

Using the facts that for $i,j=1,\ldots,n$ we have
\begin{align*}
\ac\tau^0_{\bs\mu_{\smash{-1}}}(\al_i,0)&=\ac\tau^0_{\bs\mu_{\smash{-1}}}(\al_j,0),& \ac\tau^0_{\bs\mu_x}(\al_i,0)&=\ac\tau^0_{\bs\mu_x}(\al_j,0),\\
\ac\tau^0_{\bs\mu_{\smash{-1}}}(\al_i,0)&>\ac\tau^0_{\bs\mu_{\smash{-1}}}(0,\bs d),& \ac\tau^0_{\bs\mu_x}(\al_i,0)&<\ac\tau^0_{\bs\mu_x}(0,\bs d),
\end{align*}
a calculation of Joyce and Song \cite[Prop.~13.8]{JoSo} implies that
\begin{equation*}
U\bigl((\al_1,0),\ldots,(\al_k,0),(0,\bs d),(\al_{k+1},0),\ldots,(\al_n,0);\ac\tau^0_{\bs\mu_{\smash{-1}}},\ac\tau^0_{\bs\mu_x})=\frac{(-1)^{n-k}}{k!(n-k)!}.
\end{equation*} 
Therefore as in \cite[Prop.~13.10]{JoSo}, in Theorem \ref{co3thm3} we may take
\e
\begin{split}
&\ti U\bigl((\al_1,0),\ldots,(\al_k,0),(0,\bs d),(\al_{k+1},0),\ldots,(\al_n,0);\ac\tau^0_{\bs\mu_{\smash{-1}}},\ac\tau^0_{\bs\mu_x})\\
&=\begin{cases}
\frac{(-1)^n}{n!}, & k=0, \\
0, & k>0.
\end{cases}
\end{split}
\label{co10eq33}
\e

Now compare equations \eq{co10eq18} for $s=0$, \eq{co10eq32}, and \eq{co10eq33}. These show that \eq{co10eq32} is the sum of {\it only some\/} terms in the sum \eq{co10eq18} for $[\acM_{(\be,\bs d)}^\ss(\ac\tau^0_{\bs\mu_x})]_\cL$, namely those terms in which $(\al_1,\bs d_1),\ldots,(\al_n,\bs d_n)$ consist of exactly one term $(0,\bs d)$ and the remaining terms of the form $(\al_i,0)$. We claim that {\it all other terms in the sum\/ \eq{co10eq18} are zero}, so that both sides of \eq{co10eq32} equal $[\acM_{(\be,\bs d)}^\ss(\ac\tau^0_{\bs\mu_x})]_\cL$, proving \eq{co10eq23} for $s=0$, and verifying Claim \ref{co10cla1}(a) for $x,$ $s=0$, $(\be,\bs d)$.

To see this, note that any other term in \eq{co10eq18} must contain at least two terms $(\al_i,\bs d_i)$ with $\al_i=0$ and $\bs d_i\ne 0$ of the form in the second case in \eq{co10eq17}, and these terms $(0,\bs d_i)$ must sum to $(0,\bs d)$. To get $d_r=1$, we see that there is exactly one term $(\al_i,\bs d_i)=(0,\bs d_i)$ with $d_{i,r}=1$, and this has $\bs d_{i,j}=0$ if $j<b$ and $\bs d_{i,j}=1$ if $j\ge b$, for some $b$ with $a<b\le r$. 

We can then show that for all $j=1,\ldots,n$ with $j\ne i$ we have
\begin{equation*}
\ac\tau^0_{\bs\mu_{\smash{-1}}}(\al_i,\bs d_i)<\ac\tau^0_{\bs\mu_{\smash{-1}}}(\al_j,\bs d_j),\qquad \ac\tau^0_{\bs\mu_x}(\al_i,\bs d_i)<\ac\tau^0_{\bs\mu_x}(\al_j,\bs d_j),	
\end{equation*}
where we prove the second inequality when $(\al_j,\bs d_j)=(\al_j,0)$ using the second inequality of \eq{co10eq30}. Hence Proposition \ref{co3prop3} with $I=\{i\}$ shows that
\begin{equation*}
\ti U\bigl((\al_1,\bs d_1),\ldots,(\al_n,\bs d_n);\ac\tau^0_{\bs\mu_{\smash{-1}}},\ac\tau^0_{\bs\mu_x})=0,
\end{equation*}
so the term in \eq{co10eq18} is zero. This proves the first part of the proposition. 

For the second part, let $(\be,\bs d)\in\ti S_\al$, so $\be\in\ti R_\al$, and write $(\bar{\ti\tau}{}^0_1,\bar{\ti T},\le)$ and $\baM^\ss_{(\be,1)}(\bar{\ti\tau}{}^0_1)$ for the stability condition and moduli space in Example \ref{co5ex1} constructed starting from $(\ti\tau,\ti T,\le)$ rather than $(\tau,T,\le)$. We claim that the action of $F$ on moduli spaces induces an isomorphism
\e
F_*:\baM^\ss_{(\be,1)}(\bar{\ti\tau}{}^0_1)\,{\buildrel\cong\over\longra}\,\acM_{(\be,\bs d)}^\ss(\ac\tau^{\ti\la}_{\bs\mu_x}).
\label{co10eq34}
\e
As for the proof for $F_*:\baM^\ss_{(\be,1)}(\bar\tau^0_1)\,{\buildrel\cong\over\longra}\,\acM_{(\be,\bs d)}^\ss(\ac\tau^0_{\bs\mu_x})$ above, it is sufficient to show that if $(E,V,\rho)$ lies in class $(\be,1)$ in $\baB$ then $(E,V,\rho)$ is $\bar{\ti\tau}{}^0_1$-semistable if and only if $F(E,V,\rho)$ is $\ac\tau^{\ti\la}_{\bs\mu_x}$-semistable.

Now $(E,V,\rho)$ is $\bar{\ti\tau}{}^0_1$-semistable if and only if $E$ is $\ti\tau$-semistable and there exists no subobject $E'\subsetneq E$ with $E=0$ or $\ti\tau(\lb E'\rb)=\ti\tau(\lb E/E'\rb)=\ti\tau(\be)$ with $\rho(V)\subseteq F_k(E')\subseteq F_k(E)$. This implies that $\lb E'\rb,\lb E/E'\rb\in\ti R_\al\subset R_\al$.

Similarly, $F(E,V,\rho)$ is $\ac\tau^{\ti\la}_{\bs\mu_x}$-semistable if and only if:
\begin{itemize}
\setlength{\itemsep}{0pt}
\setlength{\parsep}{0pt}
\item[(i)] There exist no subobjects $0\ne E'\subsetneq E$ with $\tau(\lb E'\rb)>\tau(\lb E/E'\rb)$, or $\tau(\lb E'\rb)=\tau(\lb E/E'\rb)$ and $\ti\la(\lb E'\rb)/\rk\lb E'\rb>\sum_{i=a}^r(\mu_i+x)/\rk\be$; and
\item[(ii)] There exist no $E'\subsetneq E$ with $E=0$ or $\tau(\lb E'\rb)=\tau(\lb E/E'\rb)$, and $\rho(V)\subseteq F_k(E')\subseteq F_k(E)$, and $\ti\la(\lb E'\rb)/\rk\lb E/E'\rb>-\sum_{i=a}^r(\mu_i+x)/\rk\be$.
\end{itemize}
Using equation \eq{co10eq30}, and $\rk\lb E'\rb,\rk\lb E/E'\rb\le\rk\be\le\rk\al$, and the condition from Proposition \ref{co10prop4} that if $\ga\in R_\al$ with $\ti\la(\ga)\ne 0$ then $\bmd{\ti\la(\ga)}\ge\ab\ha r(r+1)\rk\al$, we can show that the inequalities on $\ti\la(\lb E'\rb)$ in (i),(ii) are equivalent to $\ti\la(\lb E'\rb)>0$ and $\ti\la(\lb E'\rb)\ge 0$, respectively. Equation \eq{co10eq34} then follows as for the proof of the isomorphism \eq{co10eq11} in Proposition~\ref{co10prop4}.

As for \eq{co10eq31}, equation \eq{co5eq30} for $\be,(\ti\tau,T,\le)$ and \eq{co10eq34} yield
\ea
&\ac\Up_{(\be,\bs d)}(\ac\tau^{\ti\la}_{\bs\mu_x})=\bar\Up\bigl([\baM^\ss_{(\be,1)}(\bar{\ti\tau}{}^0_1)]_\virt\bigr)
\label{co10eq35}\\
&=\sum_{\begin{subarray}{l}n\ge 1,\;\al_1,\ldots,\al_n\in
\ti R_\al:\\ \al_1+\cdots+\al_n=\be, \; o_{\al_1}+\cdots+o_{\al_n}=o_\be \end{subarray}} \!\!\!\!\!\!\!\!\!\!
\begin{aligned}[t]
\frac{(-1)^{n+1}\la_k(\al_1)}{n!}\,\cdot&\bigl[\bigl[\cdots\bigl[[\M_{\al_1}^\ss(\ti\tau)]_\inv,\\
&\; 
[\M_{\al_2}^\ss(\ti\tau)]_\inv\bigr],\ldots\bigr],[\M_{\al_n}^\ss(\ti\tau)]_\inv\bigr].
\end{aligned}
\nonumber
\ea	
If $\al_1,\ldots,\al_n\in\ti R_\al$ then since $\ti\la(\al_i)=0$ by \eq{co10eq14} we have
$\ac\tau^{\ti\la}_{\bs\mu_{\smash{-1}}}(\al_i,0)=\ac\tau^0_{\bs\mu_{\smash{-1}}}(\al_i,0)$ and $\ac\tau^{\ti\la}_{\bs\mu_x}(\al_i,0)=\ac\tau^0_{\bs\mu_x}(\al_i,0)$, so as for \eq{co10eq33} we have
\begin{align*}
&\ti U\bigl((\al_1,0),\ldots,(\al_k,0),(0,\bs d),(\al_{k+1},0),\ldots,(\al_n,0);\ac\tau^{\ti\la}_{\bs\mu_{\smash{-1}}},\ac\tau^{\ti\la}_{\bs\mu_x})\\
&=\begin{cases}
\frac{(-1)^n}{n!}, & k=0, \\
0, & k>0.
\end{cases}
\end{align*}
The rest of the proof follows that of the first part.
\end{proof}

In the next proof, we will see that $\acM_{(\be,\bs d)}^\ss(\ac\tau^0_{\bs\mu_x})$ can change discontinuously as $x$ varies in $[-1,0]$. Thus the claim in (a) that $\acM_{(\be,\bs d)}^\rst(\ac\tau^0_{\bs\mu_x})=\acM_{(\be,\bs d)}^\ss(\ac\tau^0_{\bs\mu_x})$ for all $x$ may be surprising, as one would expect strictly semistables at points of discontinuity. This is possible as $\ac\tau^0_{\bs\mu_x}$ does not depend continuously on $x$ in the sense of Definition \ref{co3def5}, as noted after \eq{co5eq21}. To see why, note that in \eq{co10eq36} below we do not have $\ac\tau^0_{\bs\mu_x}(0,\bs e)=\ac\tau^0_{\bs\mu_x}(\be,\bs f)$ for any~$x$.

Similarly, in (b) we also show that $\acM_{(\be,\bs d)}^\rst(\ac\tau^{\ti\la}_{\bs\mu_x})=\acM_{(\be,\bs d)}^\ss(\ac\tau^{\ti\la}_{\bs\mu_x})$ for all $x$ in $[-1,0]$. This depends on $\ti\la$ being chosen sufficiently large in Proposition \ref{co10prop4}. We can have $\acM_{(\be,\bs d)}^\rst(\ac\tau^{s\ti\la}_{\bs\mu_x})\ne\acM_{(\be,\bs d)}^\ss(\ac\tau^{s\ti\la}_{\bs\mu_x})$ for $x\in[-1,0]$ and~$s\in(0,1)$.

\begin{prop}
\label{co10prop8}
Let\/ $x\in[-1,0]$ and\/ $(\be,\bs d)\in S_\al$ satisfy $\be\ne 0\ne\bs d$ and\/ $d_1\le 1,$ $d_i\ab\le \ab d_{i+1}\ab\le \ab d_i+1$ for\/ $i<r,$ and\/ $d_r\le \la_k(\be)$. Then:
\begin{itemize}
\setlength{\itemsep}{0pt}
\setlength{\parsep}{0pt}
\item[{\bf(a)}] For all\/ $x\in[-1,0]$ we have\/ $\acM_{(\be,\bs d)}^\rst(\ac\tau^0_{\bs\mu_x})=\acM_{(\be,\bs d)}^\ss(\ac\tau^0_{\bs\mu_x}),$ and Claim\/ {\rm\ref{co10cla1}(a)} holds for $x,$ $s=0,$ $(\be,\bs d)$.
\item[{\bf(b)}] Suppose $(\be,\bs d)\in\ti S_\al$. Then for all\/ $x\in[-1,0]$ we have\/ $\acM_{(\be,\bs d)}^\rst(\ac\tau^{\ti\la}_{\bs\mu_x})=\acM_{(\be,\bs d)}^\ss(\ac\tau^{\ti\la}_{\bs\mu_x}),$ and Claim\/ {\rm\ref{co10cla1}(b)} holds for $x,(\be,\bs d)$.
\end{itemize}
\end{prop}

\begin{proof} Let $(\be,\bs d)\in S_\al$ satisfy the conditions in the proposition. For the first part of the proof we suppose also that $d_r\ge 2$. Let $a=1,\ldots,r$ be least with $d_a=d_r$. Then $d_a=d_{a+1}=\cdots=d_r$ as $d_i\le d_{i+1}$ for $i<r$, and $a>1$ as $d_1\le 1$ and $d_r\ge 2$. Define $\bs e\in\N^r$ by $e_i=0$ if $i<a$ and $e_i=1$ if $i\ge a$, and set $\bs f=\bs d-\bs e$. Then $(0,\bs e),(\be,\bs f)\in S_\al$. Write $x_0=-\frac{1}{r-a+1}\sum_{i=a}^r\mu_r\in(-1,0)$. Then for $x\in[-1,0]$, equation \eq{co5eq21} implies that
\e
\begin{aligned}
\ac\tau^0_{\bs\mu_x}(0,\bs e)&<\ac\tau^0_{\bs\mu_x}(\be,\bs f)&&\text{if\/ $x\le x_0$,}\\
\ac\tau^0_{\bs\mu_x}(0,\bs e)&>\ac\tau^0_{\bs\mu_x}(\be,\bs f)&&\text{if\/ $x>x_0$.}
\end{aligned}
\label{co10eq36}
\e

We will first prove that:
\begin{itemize}
\setlength{\itemsep}{0pt}
\setlength{\parsep}{0pt}
\item[(i)] If $x\in[-1,x_0]$ then $\acM_{(\be,\bs d)}^\ss(\ac\tau^0_{\bs\mu_x})=\es$ and $[\acM_{(\be,\bs d)}^\ss(\ac\tau^0_{\bs\mu_x})]_\cL=0$. Hence Claim \ref{co10cla1}(a) holds for $x,$ $s=0,$ $(\be,\bs d)$.
\item[(ii)] Suppose that $x,x'\in[x_0,0]$ and $(\be,\bs d)=(\al_1,\bs d_1)+\cdots+(\al_n,\bs d_n)$ with $(\al_i,\bs d_i)\in S_\al$ and either $\acM_{(\al_i,\bs d_i)}^\ss(\ac\tau^0_{\bs\mu_x})\ne\es$ or $[\acM_{(\al_i,\bs d_i)}^\ss(\ac\tau^0_{\bs\mu_x})]_\cL\ne 0$ for $i=1,\ldots,n$, and $U((\al_1,\bs d_1),\ldots,(\al_n,\bs d_n);\ac\tau^0_{\bs\mu_x},\ac\tau^0_{\bs\mu_{x'}})\ne 0$.

Then either $n=1$ and $(\al_1,\bs d_1)=(\be,\bs d)$, or $n=2$ and $\bigl\{(\al_1,\bs d_1),(\al_2,\bs d_2)\bigr\}\ab=\bigl\{(0,\bs e),(\be,\bs f)\bigr\}$. If $x,x'\in(x_0,0]$ then $n=1$.
\item[(iii)] If $x\in(x_0,0]$ then $\acM_{(\be,\bs d)}^\ss(\ac\tau^0_{\bs\mu_x}),[\acM_{(\be,\bs d)}^\ss(\ac\tau^0_{\bs\mu_x})]_\cL$ are both independent of $x$, with $\acM_{(\be,\bs d)}^\rst(\ac\tau^0_{\bs\mu_x})=\acM_{(\be,\bs d)}^\ss(\ac\tau^0_{\bs\mu_x})$.
\item[(iv)] If $x\in(x_0,0]$ then 
\e
[\acM_{(\be,\bs d)}^\ss(\ac\tau^0_{\bs\mu_x})]_\cL=\bigl[[\acM_{(\be,\bs f)}^\ss(\ac\tau^0_{\bs\mu_{x_0}})]_\cL,[\acM_{(0,\bs e)}^\ss(\ac\tau^0_{\bs\mu_{x_0}})]_\cL\bigr].
\label{co10eq37}
\e
\item[(v)] If $x\in(x_0,0]$ then
\e
\ac\Up_{(\be,\bs d)}(\ac\tau^0_{\bs\mu_x})=(\la_k(\be)-d_a+1)\cdot\ac\Up_{(\be,\bs f)}(\ac\tau^0_{\bs\mu_x}).
\label{co10eq38}
\e
Here both sides are defined by (ii) for $(\be,\bs d),(\be,\bs f)$.
\item[(vi)] Suppose Claim \ref{co10cla1}(a) holds for $x,$ $s=0,$ $(\be,\bs f)$ for all $x\in[-1,0]$. Then Claim \ref{co10cla1}(a) holds for $x,$ $s=0,$ $(\be,\bs d)$ for all $x\in[-1,0]$. 
\end{itemize}
Part (a) of the proposition will then follow by induction on~$d_r=1,2,\ldots.$

In case (i), the condition $x\le x_0$ is equivalent to the inequalities on $x$ in Propositions \ref{co10prop2}(b)(iii), \ref{co10prop6}(iii) with $a-1$ in place of $a$, noting that $a>1$. But the conclusions $d_{a-1}=d_a=\cdots=d_r$ and $\bs d=0$ of these results are false. Thus Propositions \ref{co10prop2}(b)(iii)--(iv), \ref{co10prop6}(iii)--(iv) show that $\acM_{(\be,\bs d)}^\ss(\ac\tau^0_{\bs\mu_x})=\es$ and $[\acM_{(\be,\bs d)}^\ss(\ac\tau^0_{\bs\mu_x})]_\cL=0$, so Claim \ref{co10cla1}(a) holds for $x,$ $s=0,$ $(\be,\bs d)$, as in Corollary \ref{co10cor2}, proving~(i).

For (ii), suppose $x,x'$ and $(\al_1,\bs d_1),\ldots,(\al_n,\bs d_n)$ satisfy the given conditions, and write $\bs d_i=(d_{i,1},\ldots,d_{i,r})$. As $\acM_{(\al_i,\bs d_i)}^\ss(\ac\tau^0_{\bs\mu_x})\ne\es$ or $[\acM_{(\al_i,\bs d_i)}^\ss(\ac\tau^0_{\bs\mu_x})]_\cL\ne 0$, Proposition \ref{co10prop2}(b) or \ref{co10prop6} applies to $(\al_i,\bs d_i)$, so $d_{i,1}\le d_{i,2}\le\cdots\le d_{i,r}$ if $\al_i\ne 0$. We will prove:
\begin{itemize}
\setlength{\itemsep}{0pt}
\setlength{\parsep}{0pt}
\item[(A)] There does not exist $i=1,\ldots,n$ with $\al_i=0$ and $\ac\tau^0_{\bs\mu_x}(0,\bs d_i)<\ac\tau^0_{\bs\mu_x}(0,\bs e)$.
\item[(B)] There does not exist $i=1,\ldots,n$ with $\al_i=0$ and $\ac\tau^0_{\bs\mu_x}(0,\bs d_i)>\ac\tau^0_{\bs\mu_x}(0,\bs e)$.
\item[(C)] There is at most one $i=1,\ldots,n$ with $\al_i=0$, and then $\bs d_i=\bs e$.
\item[(D)] There is at most one $i=1,\ldots,n$ with $\al_i\ne 0$.
\end{itemize}

For (A), suppose for a contradiction that $\al_i=0$ and $\ac\tau^0_{\bs\mu_x}(0,\bs d_i)<\ac\tau^0_{\bs\mu_x}(0,\bs e)$. Then by Proposition \ref{co10prop2}(b)(i) or \ref{co10prop6}(i) we must have $d_{i,j}=m$ if $a'\le j\le b'$ and $d_{i,j}=0$ otherwise, for $m>0$ and $1\le a'\le b'\le r$, and from \eq{co5eq21} and $\mu_1>\cdots>\mu_r$ we deduce that $a<a'$. Choose such $i$ with $a'$ least. Then we have  $d_{j,a'-1}\le d_{j,a'}$ for all $j=1,\ldots,n$, with $d_{i,a'-1}<d_{i,a'}$. Summing over $j=1,\ldots,n$ gives $d_{a'-1}<d_{a'}$, contradicting $d_a=d_{a+1}=\cdots=d_r$ and~$a'>a$.

For (B), suppose for a contradiction that such $i$ exist. Let $I$ be the set of $i=\{1,\ldots,n\}$ with $\al_i=0$ and $\ac\tau^0_{\bs\mu_x}(0,\bs d_i)$ maximal amongst such $i$. Then $\es\ne I\subsetneq\{1,\ldots,n\}$, as $\al_j\ne 0$ for some $j$. From \eq{co5eq21}, $\ac\tau^0_{\bs\mu_x}(0,\bs d_i)>\ac\tau^0_{\bs\mu_x}(0,\bs e)$, $x,x'\ge x_0$ and the definition of $x_0$ we see that $\ac\tau^0_{\bs\mu_x}(0,\bs d_i)=(\iy,y)$, $\ac\tau^0_{\bs\mu_{x'}}(0,\bs d_i)=(\iy,y')$ for some $y,y'\in\R$ and all $i\in I$. Therefore $\ac\tau^0_{\bs\mu_x}(\al_i,\bs d_i)>\ac\tau^0_{\bs\mu_x}(\al_j,\bs d_j)$ and $\ac\tau^0_{\bs\mu_{x'}}(\al_i,\bs d_i)>\ac\tau^0_{\bs\mu_{x'}}(\be,\bs d)$ for all $i\in I$ and $j\in \{1,\ldots,n\}\sm I$, so Proposition \ref{co3prop3} shows that $U((\al_1,\bs d_1),\ldots,\ab(\al_n,\bs d_n);\ab \ac\tau^0_{\bs\mu_x},\ab\ac\tau^0_{\bs\mu_{x'}})=0$, a contradiction. 

For (C), if $\al_i=0$ then we must have $\ac\tau^0_{\bs\mu_x}(0,\bs d_i)=\ac\tau^0_{\bs\mu_x}(0,\bs e)$ by (A)--(B), and as $\mu_1,\ldots,\mu_r$ are generic this implies that $\bs d_i$ is proportional to $(0,\bs e)$. Let the sum of $\bs d_i$ for all $i=1,\ldots,n$ with $\al_i=0$ be $m\bs e$ for $m>0$. As $a>1$ we have $d_a\le d_{a-1}+1$, and $d_{j,a-1}\le d_{j,a}$ for $j=1,\ldots,n$ with $\al_j\ne 0$, and $d_{i,a-1}=0$, $\sum_id_{i,a}=m$ for $i=1,\ldots,n$ with $\al_i=0$, so $m\le 1$. Part (C) follows.

For (D), suppose for a contradiction that $\al_i\ne 0$ for at least two $i=1,\ldots,n$. Divide into subcases:
\begin{itemize}
\setlength{\itemsep}{0pt}
\setlength{\parsep}{0pt}
\item[(D1)] Either $\al_i\ne 0$ for all $i$, or $\al_i=0$ for some $i$, so that $\bs d_i=\bs e$ by (C), and $x=x_0$. This implies that $\ac\tau^0_{\bs\mu_x}(0,\bs e)=(-\iy,y)$ for $y\in\R$ by \eq{co5eq21}.
\item[(D2)] We have $\al_i=0$ for some $i$, so that $\bs d_i=\bs e$, and $x>x_0$. This implies that $\ac\tau^0_{\bs\mu_x}(0,\bs e)=(\iy,y)$ for $y\in\R$.
\end{itemize}

In case (D1) we define $I$ to be the set of $i=1,\ldots,n$ such that $\al_i\ne 0$ and $\ac\tau^0_{\bs\mu_x}(\al_i,\bs d_i)$ is maximal in $M$ for such $i$. We claim that $\es\ne I\subseteq\{1,\ldots,n\}$, and $\ac\tau^0_{\bs\mu_x}(\al_i,\bs d_i)>\ac\tau^0_{\bs\mu_x}(\al_j,\bs d_j)$ and $\ac\tau^0_{\bs\mu_{x'}}(\al_i,\bs d_i)>\ac\tau^0_{\bs\mu_{x'}}(\be,\bs d)$ for all $i\in I$ and $j\in \{1,\ldots,n\}\sm I$, so Proposition \ref{co3prop3} shows that $U((\al_1,\bs d_1),\ldots,\ab(\al_n,\bs d_n);\ab \ac\tau^0_{\bs\mu_x},\ab\ac\tau^0_{\bs\mu_{x'}})=0$, a contradiction.

Similarly, in case (D2) we define $I$ to be the set of $i=1,\ldots,n$ such that $\al_i\ne 0$ and $\ac\tau^0_{\bs\mu_x}(\al_i,\bs d_i)$ is minimal in $M$ for such $i$. We claim that $\es\ne I\subseteq\{1,\ldots,n\}$, and $\ac\tau^0_{\bs\mu_x}(\al_i,\bs d_i)<\ac\tau^0_{\bs\mu_x}(\al_j,\bs d_j)$ and $\ac\tau^0_{\bs\mu_{x'}}(\al_i,\bs d_i)<\ac\tau^0_{\bs\mu_{x'}}(\be,\bs d)$ for all $i\in I$ and $j\in \{1,\ldots,n\}\sm I$, so again Proposition \ref{co3prop3} shows that $U((\al_1,\bs d_1),\ldots,\ab(\al_n,\bs d_n);\ab \ac\tau^0_{\bs\mu_x},\ab\ac\tau^0_{\bs\mu_{x'}})=0$, a contradiction.

To prove the claims, note that $\ac\tau^0_{\bs\mu_x}(\al_i,\bs d_i)>\ac\tau^0_{\bs\mu_x}(\al_j,\bs d_j)$ in case (D1) is automatic from the definition of $I$ when $\al_j\ne 0$, and from $\ac\tau^0_{\bs\mu_x}(\al_j,\bs d_j)=(-\iy,y)$ when $\al_j=0$. Similarly $\ac\tau^0_{\bs\mu_x}(\al_i,\bs d_i)<\ac\tau^0_{\bs\mu_x}(\al_j,\bs d_j)$ in case (D2) is automatic.

Suppose $i=1,\ldots,n$ with $\al_i\ne 0$. If $\bs d_i\ne 0$, let $b=1,\ldots,n$ be least with $d_{i,b}>0$. Then $d_b=1$ as $d_1\le 1$, $d_{j+1}\le d_j+1$, so $b<a\le r$ as $d_a=d_{a+1}=\cdots=d_r$ and $d_r\ge 2$. Since $\mu_1\gg \cdots\gg\mu_r>0$, and $x_0=-\frac{1}{r-a+1}\sum_{i=a}^r\mu_r\in(-1,0)$, and $x,x'\in[x_0,0]$, we see from \eq{co5eq21} that
\e
\ac\tau^0_{\bs\mu_x}(\al_i,\bs d_i)=(\tau(\al),y),\;\> \ac\tau^0_{\bs\mu_{x'}}(\al_i,\bs d_i)=(\tau(\al),y')\;\>\text{with}\;\> y,y'\approx \frac{\mu_b}{\rk\al_i}.
\label{co10eq39}
\e
If $\bs d_i=0$ then $\ac\tau^0_{\bs\mu_x}(\al_i,\bs d_i)=\ac\tau^0_{\bs\mu_{x'}}(\al_i,\bs d_i)=(\tau(\al),0)$.

Similarly, let $c=1,\ldots,n$ be least with $d_c>0$. Then 
\e
\ac\tau^0_{\bs\mu_{x'}}(\be,\bs d)=(\tau(\al),z)\quad\text{with}\quad z\approx \frac{\mu_c}{\rk\be}.
\label{co10eq40}
\e
In case (D1) there is a unique $i=1,\ldots,n$ with $d_{i,c}>0$, and then $I=\{i\}$, so that $\es\ne I\subseteq\{1,\ldots,n\}$, and $\ac\tau^0_{\bs\mu_{x'}}(\al_i,\bs d_i)>\ac\tau^0_{\bs\mu_{x'}}(\be,\bs d)$ by \eq{co10eq39}--\eq{co10eq40} since $\mu_b/\rk\al_i>\mu_c/\rk\be$ as $\mu_c>0$ and $\rk\al_i<\rk\be$. This proves case~(D1).

In case (D2) the minimal value of $\ac\tau^0_{\bs\mu_x}(\al_i,\bs d_i)$ for $\al_i\ne 0$ must be either $(\tau(\al),y)$ for $y\approx \mu_b/\rk\al_j$ with $b>c$ when $\bs d_i\ne 0$, or $(\tau(\al),0)$ when $\bs d_i=0$. We must have $\es\ne I\subseteq\{1,\ldots,n\}$ as $d_{j,c}=1$ for some $j$ and then $j\notin I$, and $\ac\tau^0_{\bs\mu_{x'}}(\al_i,\bs d_i)<\ac\tau^0_{\bs\mu_{x'}}(\be,\bs d)$ by \eq{co10eq39}--\eq{co10eq40}, noting that $\mu_c/\mu_b>\rk\al$ as in Definition \ref{co10def1}, so we can suppose that $\mu_c/\mu_b>\rk\al\ge \rk\be\ge\rk\be/\rk\al_i$, implying that $\mu_c/\rk\be>\mu_b/\rk\al_i$. This proves case (D2), and completes~(D).

Together (C)--(D) prove the first part of (ii). For the second part, if $x,x'\in(x_0,0]$ and $\bigl\{(\al_1,\bs d_1),(\al_2,\bs d_2)\bigr\}\ab=\bigl\{(0,\bs e),(\be,\bs f)\bigr\}$ then from \eq{co10eq36} and Definition \ref{co3def4} we can show that $U((\al_1,\bs d_1),(\al_2,\bs d_2);\ab \ac\tau^0_{\bs\mu_x},\ab\ac\tau^0_{\bs\mu_{x'}})=0$, a contradiction, so the case $n=2$ is excluded, as we want.

For (iii), let $x,x'\in(x_0,0]$, and suppose $[E,\bs V,\bs\rho]\in\acM_{(\be,\bs d)}^\ss(\ac\tau^0_{\bs\mu_{x'}})$. Write $(\al_1,\bs d_1),\ldots,(\al_n,\bs d_n)$ for the classes of the $\ac\tau^0_{\bs\mu_x}$-semistable subquotients of $\ac\tau^0_{\bs\mu_x}$-Harder--Narasimhan filtration of $(E,\bs V,\bs\rho)$ in $\acA$. Then $\acM_{(\al_i,\bs d_i)}^\ss(\ac\tau^0_{\bs\mu_x})\ne\es$ for $i=1,\ldots,n$, and $U((\al_1,\bs d_1),\ldots,(\al_n,\bs d_n);\ac\tau^0_{\bs\mu_x},\ac\tau^0_{\bs\mu_{x'}})\ne 0$ by Corollary \ref{co3cor1}. Hence $n=1$ by (ii), so $(E,\bs V,\bs\rho)$ is $\ac\tau^0_{\bs\mu_x}$-semistable, and $\acM_{(\be,\bs d)}^\ss(\ac\tau^0_{\bs\mu_{x'}})\subseteq\acM_{(\be,\bs d)}^\ss(\ac\tau^0_{\bs\mu_x})$. Similarly $\acM_{(\be,\bs d)}^\ss(\ac\tau^0_{\bs\mu_x})\subseteq\acM_{(\be,\bs d)}^\ss(\ac\tau^0_{\bs\mu_{x'}})$, so~$\acM_{(\be,\bs d)}^\ss(\ac\tau^0_{\bs\mu_{x'}})=\acM_{(\be,\bs d)}^\ss(\ac\tau^0_{\bs\mu_x})$. 

Thus $\acM_{(\be,\bs d)}^\ss(\ac\tau^0_{\bs\mu_x})$ is independent of $x\in(x_0,0]$. Also \eq{co10eq21} with $x,x'$ exchanged reduces to $[\acM_{(\be,\bs d)}^\ss(\ac\tau^0_{\bs\mu_{x'}})]_\cL=[\acM_{(\be,\bs d)}^\ss(\ac\tau^0_{\bs\mu_x})]_\cL$, since $U((\al_1,\bs d_1),\ab\ldots,\ab(\al_n,\bs d_n);\ab\ac\tau^0_{\bs\mu_x},\ac\tau^0_{\bs\mu_{x'}})=0$ when $n\ge 2$ by (ii), which implies that $\ti U((\al_1,\bs d_1),\ab\ldots,\ab(\al_n,\bs d_n);\ab\ac\tau^0_{\bs\mu_x},\ac\tau^0_{\bs\mu_{x'}})=0$. Hence $[\acM_{(\be,\bs d)}^\ss(\ac\tau^0_{\bs\mu_x})]_\cL$ is independent of~$x\in(x_0,0]$.

For generic $x\in(x_0,0]$ we see that $\acM_{(\be,\bs d)}^\rst(\ac\tau^0_{\bs\mu_x})=\acM_{(\be,\bs d)}^\ss(\ac\tau^0_{\bs\mu_x})$, as we cannot split $(\be,\bs d)=(\al_1,\bs d_1)+\cdots+(\al_n,\bs d_n)$ for $n\ge 2$ with $\ac\tau^0_{\bs\mu_x}(\al_1,\bs d_1)=\cdots=\ac\tau^0_{\bs\mu_x}(\al_n,\bs d_n)$. Since $\acM_{(\be,\bs d)}^\ss(\ac\tau^0_{\bs\mu_x})$ is independent of $x\in(x_0,0]$, this holds for all $x\in(x_0,0]$, proving~(iii).

For (iv), consider equation \eq{co10eq21} with $x_0$ in place of $x'$. By (ii), there are only three possibly nonzero terms, namely $n=1$ and $(\al_1,\bs d_1)=(\be,\bs d)$, and $n=2$ and $(\al_1,\bs d_1)=(\be,\bs f)$, $(\al_2,\bs d_2)=(0,\bs e)$, and $n=2$ and $(\al_1,\bs d_1)=(0,\bs e)$, $(\al_2,\bs d_2)=(\be,\bs f)$. The first term gives zero by (i). By \eq{co10eq36} and Definition \ref{co3def4} we can show that
\begin{equation*}
U((\be,\bs f),(0,\bs e);\ac\tau^0_{\bs\mu_{x'}},\ac\tau^0_{\bs\mu_x})=1,\qquad
U((0,\bs e),(\be,\bs f);\ac\tau^0_{\bs\mu_{x'}},\ac\tau^0_{\bs\mu_x})=-1,
\end{equation*}
so by Theorem \ref{co3thm3} we may take
\e
\ti U((\be,\bs f),(0,\bs e);\ac\tau^0_{\bs\mu_{x'}},\ac\tau^0_{\bs\mu_x})=1,\qquad
\ti U((0,\bs e),(\be,\bs f);\ac\tau^0_{\bs\mu_{x'}},\ac\tau^0_{\bs\mu_x})=0.
\label{co10eq41}
\e
Thus \eq{co10eq21} reduces to \eq{co10eq37}, proving~(iv).

For (v), in \eq{co10eq38} note that $\acM_{(\be,\bs d)}^\ss(\ac\tau^0_{\bs\mu_x})$ is independent of $x\in(x_0,0]$ by (iii), and so is $\acM_{(\be,\bs d)}^\ss(\ac\tau^0_{\bs\mu_x})$ by (iii) for $(\be,\bs f)$ (note that $x_0$ for $(\be,\bs f)$ is smaller than $x_0$ for $(\be,\bs d)$). So \eq{co10eq38} is independent of $x\in(x_0,0]$. We will prove it for fixed $x\in(x_0,0]$ very close to~$x_0$. 

For each $\C$-point $[(E,\bs V,\bs\rho)]$ in $\acM_{(\be,\bs d)}^\ss(\ac\tau^0_{\bs\mu_x})$, define a subobject $(E,\bs V',\bs\rho')\subset(E,\bs V,\bs\rho)$ by $V'_i=V_i$ for $i<a$, and $V_i'=\rho_{i-1}\ci\cdots\ci\rho_{a-1}(V_{a-1})$ if $a\le i\le r$, and $\rho'_i=\rho_i\vert_{V_i'}$. As each $\rho_i$ is injective by Proposition \ref{co10prop2}(a)(ii), we have $V_i'\cong V_{a-1}$ for $i\ge a-1$, so $\lb E,\bs V',\bs\rho'\rb=(\be,\bs f)$. 

Suppose for a contradiction that $(E,\bs V',\bs\rho')$ is $\ac\tau^0_{\bs\mu_x}$-unstable. Then there exists $0\ne(E'',\bs V'',\bs\rho'')\subsetneq(E,\bs V',\bs\rho')$ with
\e
\ac\tau^0_{\bs\mu_x}(\lb E'',\bs V'',\bs\rho''\rb)>\ac\tau^0_{\bs\mu_x}(\lb E/E'',\bs V'/\bs V'',\bs\si'\rb).
\label{co10eq42}
\e
If $\tau(\lb E''\rb)>\tau(\lb E/E''\rb)$ then $(E'',\bs V'',\bs\rho'')$ $\ac\tau^0_{\bs\mu_x}$-destabilizes $(E,\bs V,\bs\rho)$, a contradiction. If $\tau(\lb E''\rb)=\tau(\lb E/E''\rb)$ then writing
\begin{align*}
\ac\tau^0_{\bs\mu_x}(\lb E,\bs V,\bs\rho\rb )&=(\tau(\al),y), &
\ac\tau^0_{\bs\mu_x}(\lb E'',\bs V'',\bs\rho''\rb)&=(\tau(\al),y''), \\
\ac\tau^0_{\bs\mu_x}(\lb E,\bs V',\bs\rho'\rb)&=(\tau(\al),y'), &
\ac\tau^0_{\bs\mu_x}(\lb E/E'',\bs V'/\bs V'',\bs\si'\rb)&=(\tau(\al),y'''),
\end{align*}
then \eq{co10eq42} implies that $y''>y'>y'''$. But $y$ is close to $y'$ as $x$ is close to $x_0=-\frac{1}{r-a+1}\sum_{i=a}^r\mu_r\in(-1,0)$, so $y''>y$, and thus $(E'',\bs V'',\bs\rho'')$ $\ac\tau^0_{\bs\mu_x}$-destabilizes $(E,\bs V,\bs\rho)$, a contradiction. Hence $(E,\bs V',\bs\rho')$ is $\ac\tau^0_{\bs\mu_x}$-semistable. 

Thus we may define a morphism of moduli stacks
\e
\Pi:\acM_{(\be,\bs d)}^\ss(\ac\tau^0_{\bs\mu_x})\longra\acM_{(\be,\bs f)}^\ss(\ac\tau^0_{\bs\mu_x}),\quad
\Pi:[E,\bs V,\bs\rho]\longmapsto[E,\bs V',\bs\rho'].
\label{co10eq43}
\e

We claim that $\Pi$ in \eq{co10eq43} is smooth with fibre $\CP^{\la_k(\be)-d_a}$. To see this, note that for $(E,\bs V,\bs\rho)$, $(E,\bs V',\bs\rho')$ as above we have vector subspaces $\rho'_r(V_r')\subseteq\rho_r(V_r)\subseteq F_k(E)$, where as $\rho'_r,\rho_r$ are injective by Proposition \ref{co10prop2}(a)(ii) we have $\dim \rho'_r(V_r')=d_r-1$, $\dim \rho_r(V_r)=d_r$ and $\dim F_k(E)=\la_k(\be)$ by Assumption \ref{co5ass1}(g)(iv). Hence $\rho_r(V_r)/\rho'_r(V_r')$ is a point of the projective space $\bP(F_k(E)/\rho'_r(V_r'))\cong \CP^{\la_k(\be)-d_a}$. It is not difficult to show that $(E,\bs V,\bs\rho)$ is determined uniquely up to isomorphism by $(E,\bs V',\bs\rho')$ and $\rho_r(V_r)/\rho'_r(V_r')$, and $\Pi$ is smooth with fibre $\bP(F_k(E)/\rho'_r(V_r'))$ over~$(E,\bs V',\bs\rho')$.

Similarly, we can show the derived version $\bs\Pi:\bs\acM_{(\be,\bs d)}^\ss(\ac\tau^0_{\bs\mu_x})\ra\bs\acM_{(\be,\bs f)}^\ss(\ac\tau^0_{\bs\mu_x})$ of \eq{co10eq43} is a smooth morphism of proper quasi-smooth derived algebraic spaces with fibre $\CP^{\la_k(\be)-d_a}$. Hence Corollary \ref{co2cor1} implies that
\e
\begin{split}
&\Pi_*\bigl([\acM_{(\be,\bs d)}^\ss(\ac\tau^0_{\bs\mu_x})]_\virt\cap c_\top(\bT_{\acM_{(\be,\bs d)}^\ss(\ac\tau^0_{\bs\mu_x})/\acM_{(\be,\bs f)}^\ss(\ac\tau^0_{\bs\mu_x})})\bigr)\\
&\quad =(\la_k(\be)-d_a+1)\cdot
[\acM_{(\be,\bs f)}^\ss(\ac\tau^0_{\bs\mu_x})]_\virt.
\end{split}
\label{co10eq44}
\e
Then \eq{co10eq38} follows from
\begin{align*}
&\ac\Up_{(\be,\bs d)}(\ac\tau^0_{\bs\mu_x})=
(\Pi_{\M_\be}^\pl)_*\bigl([\acM_{(\be,\bs d)}^\ss(\ac\tau^0_{\bs\mu_x})]_\virt\cap c_\top(\bT_{\acM_{(\be,\bs d)}^\ss(\ac\tau^0_{\bs\mu_x})/\M^\pl_\be})\bigr)\\
&\quad =(\Pi_{\M_\be}^\pl\ci\Pi)_*\bigl([\acM_{(\be,\bs d)}^\ss(\ac\tau^0_{\bs\mu_x})]_\virt\cap (c_\top(\bT_{\acM_{(\be,\bs d)}^\ss(\ac\tau^0_{\bs\mu_x})/\acM_{(\be,\bs f)}^\ss(\ac\tau^0_{\bs\mu_x})})\\
&\qquad\qquad \cup \Pi^*(c_\top(\bT_{\acM_{(\be,\bs f)}^\ss(\ac\tau^0_{\bs\mu_x})/\M^\pl_\be}))
\bigr)\\
&\quad=(\Pi_{\M_\be}^\pl)_*\bigl\{\Pi_*\bigl([\acM_{(\be,\bs d)}^\ss(\ac\tau^0_{\bs\mu_x})]_\virt\cap (c_\top(\bT_{\acM_{(\be,\bs d)}^\ss(\ac\tau^0_{\bs\mu_x})/\acM_{(\be,\bs f)}^\ss(\ac\tau^0_{\bs\mu_x})})\bigr)\\
&\qquad\qquad \cap \Pi^*(c_\top(\bT_{\acM_{(\be,\bs f)}^\ss(\ac\tau^0_{\bs\mu_x})/\M^\pl_\be}))
\bigr\}\\
&\quad=(\la_k(\be)-d_a+1)\cdot\ac\Up_{(\be,\bs f)}(\ac\tau^0_{\bs\mu_x}),
\end{align*}
where we use \eq{co10eq6} in the first step, the following commutative diagram of smooth morphisms in the second:
\begin{equation*}
\xymatrix@C=80pt@R=11pt{ 
\acM_{(\be,\bs d)}^\ss(\ac\tau^0_{\bs\mu_x}) \ar[rr]_{\Pi_{\M_\be}^\pl} \ar[dr]_\Pi && \M^\pl_\be, \\
& \acM_{(\be,\bs f)}^\ss(\ac\tau^0_{\bs\mu_x}) \ar[ur]_{\Pi_{\M_\be}^\pl} }	
\end{equation*}
properties of (co)homology in the third, and \eq{co10eq6} and \eq{co10eq44} in the fourth. This proves~(v).

For (vi), suppose Claim \ref{co10cla1}(a) holds for all $x\in[-1,0]$, $s=0,$ $(\be,\bs f)$. Then
\begin{equation*}
[\acM_{(\be,\bs f)}^\ss(\ac\tau^0_{\bs\mu_{x_0}})]_\cL=
[\acM_{(\be,\bs f)}^\ss(\ac\tau^0_{\bs\mu_x})]_\cL=q^{\bs f}\ot\ac\Up_{(\be,\bs f)}(\ac\tau^0_{\bs\mu_x}).
\end{equation*}
Here the first step follows from (iii) for $(\be,\bs f)$, noting that $x_0$ for $(\be,\bs f)$ is smaller than $x_0$ for $(\be,\bs d)$. As $[\acM_{(0,\bs e)}^\ss(\ac\tau^0_{\bs\mu_x})]_\cL$ is independent of $x$, equation \eq{co10eq17} implies that
\begin{equation*}
[\acM_{(0,\bs e)}^\ss(\ac\tau^0_{\bs\mu_{x_0}})]_\cL=p^{\bs e}=p_ap_{a+1}\cdots p_r.
\end{equation*}
Substituting the last two equations into \eq{co10eq37} and using \eq{co10eq16} and 
\begin{align*}
\ac\chi((\be,\bs f),(0,\bs e))&=f_r-f_{a-1}=0,\\
\ac\chi((0,\bs e),(\be,\bs f))&=f_a-\la_k(\be)=d_a-1-\la_k(\be),
\end{align*}
which follow from \eq{co10eq2} and the definitions of $\bs d,\bs e,\bs f$, we deduce that
\e
[\acM_{(\be,\bs d)}^\ss(\ac\tau^0_{\bs\mu_x})]_\cL=(\la_k(\be)-d_a+1)
\cdot q^{\bs d}\ot \ac\Up_{(\be,\bs f)}(\ac\tau^0_{\bs\mu_x}).
\label{co10eq45}
\e
Combining \eq{co10eq38} and \eq{co10eq45} proves \eq{co10eq23}, so Claim \ref{co10cla1}(a) holds for $x,s=0,(\be,\bs d)$. This proves~(vi).

We can now prove part (a) of the proposition. Proposition \ref{co10prop7} when $d_r=1$, and (iii) above when $d_r\ge 2$, show that $\acM_{(\be,\bs d)}^\rst(\ac\tau^0_{\bs\mu_x})=\acM_{(\be,\bs d)}^\ss(\ac\tau^0_{\bs\mu_x})$ for all $x\in[-1,0]$. We will show by induction on $d_r=1,2,\ldots$ that Claim \ref{co10cla1}(a) holds for all $x\in[-1,0],$ $s=0$, $(\be,\bs d)$. The first step $d_r=1$ holds by Proposition \ref{co10prop7}, and the inductive step holds by (vi) above, since $f_r=d_r-1$. Thus part (a) follows by induction.

The proof of (b) is similar, taking $(\be,\bs d)\in\ti S_\al$ and replacing $\ac\tau^0_{\bs\mu_x}$ by $\ac\tau^{\ti\la}_{\bs\mu_x}$ throughout. For (ii), we need to add one extra point to the proof above:

If $(\al_i,\bs d_i)\in\ti S_\al$ for some $i$, so that $\ti\la(\al_i)=0$ by \eq{co10eq4} and \eq{co10eq14}, then $\ac\tau^{\ti\la}_{\bs\mu_x}(\al_i,\bs d_i)=\ac\tau^0_{\bs\mu_x}(\al_i,\bs d_i)$, so \eq{co10eq39} holds for $\ac\tau^{\ti\la}_{\bs\mu_x}$. But if $(\al_i,\bs d_i)\in S_\al\sm\ti S_\al$, so that $\ti\la(\al_i)\ne 0,$ then by \eq{co5eq21}, equation \eq{co10eq39} must be replaced by
\begin{equation*}
\ac\tau^{\ti\la}_{\bs\mu_x}(\al_i,\bs d_i)=(\tau(\al),y),\;\> \ac\tau^{\ti\la}_{\bs\mu_{x'}}(\al_i,\bs d_i)=(\tau(\al),y')\;\>\text{with}\;\>   y,y'\approx \frac{\ti\la(\al_i)+\mu_b}{\rk\al_i}.
\end{equation*}
The rest of the proof needs only obvious changes.
\end{proof}

From the proofs of Propositions \ref{co10prop7} and \ref{co10prop8}, in particular equations \eq{co10eq31}, \eq{co10eq35}, \eq{co10eq38}, \eq{co10eq44}, and \eq{co10eq45}, using induction on $d_r$, we deduce:

\begin{cor}
\label{co10cor3}
Let\/ $(\be,\bs d)\in S_\al$ satisfy $\be\ne 0\ne\bs d$ and\/ $d_1\le 1,$ $d_i\ab\le \ab d_{i+1}\ab\le \ab d_i+1$ for\/ $i<r,$ and\/ $d_r\le \la_k(\be)$. Then
\ea
&[\acM_{(\be,\bs d)}^\ss(\ac\tau^0_{\bs\mu_0})]_\cL=\frac{(\la_k(\be)-1)!}{(\la_k(\be)-d_r)!}\cdot q^{\bs d}\ot
\label{co10eq46}\\
&\sum_{\begin{subarray}{l}n\ge 1,\;\al_1,\ldots,\al_n\in
R_\al:\\ \al_1+\cdots+\al_n=\be, \; o_{\al_1}+\cdots+o_{\al_n}=o_\be \end{subarray}} \!\!\!\!\!\!\!\!\!\!
\begin{aligned}[t]
\frac{(-1)^{n+1}\la_k(\al_1)}{n!}\,\cdot&\bigl[\bigl[\cdots\bigl[[\M_{\al_1}^\ss(\tau)]_\inv,\\
&\; 
[\M_{\al_2}^\ss(\tau)]_\inv\bigr],\ldots\bigr],[\M_{\al_n}^\ss(\tau)]_\inv\bigr],
\end{aligned}
\nonumber\\
&[\acM_{(\be,\bs d)}^\ss(\ac\tau^{\ti\la}_{\bs\mu_0})]_{\smash{\ti{\cal L}}}
=\frac{(\la_k(\be)-1)!}{(\la_k(\be)-d_r)!}\cdot q^{\bs d}\ot
\label{co10eq47}\\
&\sum_{\begin{subarray}{l}n\ge 1,\;\al_1,\ldots,\al_n\in
\ti R_\al:\\ \al_1+\cdots+\al_n=\be, \; o_{\al_1}+\cdots+o_{\al_n}=o_\be \end{subarray}} \!\!\!\!\!\!\!\!\!\!
\begin{aligned}[t]
\frac{(-1)^{n+1}\la_k(\al_1)}{n!}\,\cdot&\bigl[\bigl[\cdots\bigl[[\M_{\al_1}^\ss(\ti\tau)]_\inv,\\
&\; 
[\M_{\al_2}^\ss(\ti\tau)]_\inv\bigr],\ldots\bigr],[\M_{\al_n}^\ss(\ti\tau)]_\inv\bigr].
\end{aligned}
\nonumber
\ea
\end{cor}

\begin{prop}
\label{co10prop9}
Let\/ $s\in[0,1]$ and\/ $(\be,\bs d)\in S_\al$ satisfy $\be\ne 0\ne\bs d,$ $d_1\le 1,$ $d_j\ab\le\ab d_{j+1}\ab\le\ab d_j+1$ for\/ $j<r,$ $d_r=\la_k(\be)$. Suppose there is a splitting $(\be,\bs d)=(\al_1,\bs d_1)+\cdots+(\al_n,\bs d_n)$ for $n\ge 2,$ such~that 
\e
\ac\tau^{s\ti\la}_{\bs\mu_0}(\al_1,\bs d_1)=\ac\tau^{s\ti\la}_{\bs\mu_0}(\al_2,\bs d_2)=\cdots=\ac\tau^{s\ti\la}_{\bs\mu_0}(\al_n,\bs d_n),
\label{co10eq48}
\e
and for each $i\!=\!1,\ldots,n$ we have $(\al_i,\bs d_i)\!\in\! S_\al$ and either $\acM_{(\al_i,\bs d_i)}^\ss(\ac\tau^{s\ti\la}_{\bs\mu_0})\!\ne\!\es$ or $[\acM_{(\al_i,\bs d_i)}^\ss(\ac\tau^{s\ti\la}_{\bs\mu_0})]_\cL\ne 0$. Then $s\in(0,1),$ $n=2,$ and for $i=1,2$ we have $\al_i\ne 0\ne\bs d_i,$ $d_{i,1}\le 1,$ $d_{i,j}\ab\le\ab d_{i,j+1}\ab\le\ab d_{i,j}+1$ for\/ $j<r,$ $d_{i,r}=\la_k(\al_i),$ and\/ $\acM_{(\al_i,\bs d_i)}^\rst(\ac\tau^{s\ti\la}_{\bs\mu_0})\!=\!\acM_{(\al_i,\bs d_i)}^\ss(\ac\tau^{s\ti\la}_{\bs\mu_0})$. For given $(\be,\bs d)$ there are only finitely many possibilities for $s,(\al_1,\bs d_1),(\al_2,\bs d_2)$.	
\end{prop}

\begin{proof} Divide into cases:
\begin{itemize}
\setlength{\itemsep}{0pt}
\setlength{\parsep}{0pt}
\item[(i)] $\al_i=0$ for some $i=1,\ldots,n$.
\item[(ii)] $\al_i\ne 0$ for all $i$ and $\bs d_j=0$ for some $j=1,\ldots,n$.
\item[(iii)] $\al_i\!\ne\! 0$, $\bs d_i\!\ne\! 0$ for all $i$, and $\bs d_j,\bs d_{j'}$ are proportional for some $1\!\le\! j\!<\!j'\!\le\! n$.
\item[(iv)] $\al_i\ne 0$, $\bs d_i\ne 0$ for all $i$, and the $\bs d_i$ are pairwise not proportional.
\end{itemize}
Case (i) does not occur by \eq{co5eq21}, as $\ac\tau^{s\ti\la}_{\bs\mu_0}(\al_i,\bs d_i)=(\pm\iy,y)$ if $\al_i=0$ and $\ac\tau^{s\ti\la}_{\bs\mu_0}(\al_j,\bs d_j)=(t,z)$ if $\al_j\ne 0$, with $t\in T$, $y,z\in\R$, so $\ac\tau^{s\ti\la}_{\bs\mu_0}(\al_i,\bs d_i)\ne\ac\tau^{s\ti\la}_{\bs\mu_0}(\al_j,\bs d_j)$.

In case (ii) we have $\sum_{i\ne j}d_{i,r}>\sum_{i\ne j}\la_k(\al_i)$ as $d_r=\la_k(\be)$ and $\la_k(\al_j)>0$, so for some $i\ne j$ we have $d_{i,r}>\la_k(\al_i)$. Also $d_{i,l+1}\le d_{i,l}+1$ for $1\le l<r$ as $d_{l+1}\le d_l+1$ by assumption and $d_{m,l}\le d_{m,l+1}$ for $m\ne i$ by Proposition \ref{co10prop2}(b)(ii). Hence Propositions \ref{co10prop2}(b)(ii) and \ref{co10prop6}(ii) give $\acM_{(\al_i,\bs d_i)}^\ss(\ac\tau^{s\ti\la}_{\bs\mu_0})=\es$ and $[\acM_{(\al_i,\bs d_i)}^\ss(\ac\tau^{s\ti\la}_{\bs\mu_0})]_\cL=0$, a contradiction. So case (ii) does not occur.

In case (iii), let $a=1,\ldots,r$ be least with $d_{j,a}>0$, so that $d_{j',a}>0$ as $\bs d_j,\bs d_{j'}$ are nonzero and proportional. If $a=1$ then $d_1\ge d_{j,a}+d_{j',a}\ge 2$ contradicting $d_1\le 1$. If $a>1$ then as $d_{i,a-1}\le d_{i,a}$ for $i=1,\ldots,n$ by Propositions \ref{co10prop2}(b)(ii) and \ref{co10prop6}(ii) for $(\al_i,\bs d_i)$ and $d_{j,a-1}=d_{j',a-1}=0$, $d_{j,a},d_{j',a}>0$ we have $d_a\ge d_{a-1}+2$, contradicting $d_a\le\ab d_{a-1}+1$. So case (iii) does not occur.

In case (iv), as the $\bs d_i$ are pairwise non-proportional, the $\R$ component of \eq{co10eq48} in $M=(T\amalg\{\pm\iy\})\t\R$ is $n-1$ independent equations on one variable $s\in[0,1]$, so as $\mu_1,\ldots,\mu_r$ are generic we see that \eq{co10eq48} has no solution $s$ if $n>2$. Thus $n=2$. By \eq{co5eq21}, the $\R$ component of \eq{co10eq48} is then of the form $Ps+Q=0$ for $P,Q\in\R$, where $Q\ne 0$ as $\bs d_1,\bs d_2$ are not proportional, so there is at most one solution $s\in\R$. We exclude $s=0,1$ as $\mu_1,\ldots,\mu_r$ are generic and $\ti\la$ is integral, so $s\in(0,1)$ if \eq{co10eq48} has a solution~$s\in[0,1]$. 

As $\acM_{(\al_i,\bs d_i)}^\ss(\ac\tau^{s\ti\la}_{\bs\mu_0})\ne\es$ or $[\acM_{(\al_i,\bs d_i)}^\ss(\ac\tau^{s\ti\la}_{\bs\mu_0})]_\cL\ne 0$, Proposition \ref{co10prop2}(b)(ii) or \ref{co10prop6} and $d_{i+1}\le d_i+1$ give $d_{i,1}\le 1,$ $d_{i,j}\ab\le\ab d_{i,j+1}\ab\le\ab d_{i,j}+1$ for $j<r$ and $d_{i,r}\le\la_k(\al_i)$. Since $d_{1,r}+d_{2,r}=d_r=\la_k(\be)=\la_k(\al_1)+\la_k(\al_2)$ we see that $d_{i,r}=\la_k(\al_i)$. As $S_\al$ is finite and $(\al_1,\bs d_1),(\al_2,\bs d_2)$ determine $s$, there are only finitely many possibilities for $s,(\al_1,\bs d_1),(\al_2,\bs d_2)$.

If $\acM_{(\al_2,\bs d_2)}^\rst(\ac\tau^{s\ti\la}_{\bs\mu_0})\ne\acM_{(\al_2,\bs d_2)}^\ss(\ac\tau^{s\ti\la}_{\bs\mu_0})$ there is a $\C$-point in $\acM_{(\al_2,\bs d_2)}^\ss(\ac\tau^{s\ti\la}_{\bs\mu_0})$ representing a sum of two nonzero $\ac\tau^{s\ti\la}_{\bs\mu_0}$-semistables, say in classes $(\al'_2,\bs d'_2),(\al'_3,\bs d'_3)$ in $C(\acB)$. Then $(\be,\bs d)=(\al_1,\bs d_1)+(\al'_2,\bs d'_2)+(\al'_3,\bs d'_3)$ satisfies the conditions of the proposition, contradicting $n=2$. Hence $\acM_{(\al_2,\bs d_2)}^\rst(\ac\tau^{s\ti\la}_{\bs\mu_0})=\acM_{(\al_2,\bs d_2)}^\ss(\ac\tau^{s\ti\la}_{\bs\mu_0})$, and similarly $\acM_{(\al_1,\bs d_1)}^\rst(\ac\tau^{s\ti\la}_{\bs\mu_0})=\acM_{(\al_1,\bs d_1)}^\ss(\ac\tau^{s\ti\la}_{\bs\mu_0})$. This completes the proof.
\end{proof}

\begin{prop}
\label{co10prop10}
Let\/ $s\in[0,1]$ and\/ $(\be,\bs d)\in S_\al$ satisfy $\be\ne 0\ne\bs d$ and\/ $d_1\le 1,$ $d_i\ab\le\ab d_{i+1}\ab\le\ab d_i+1$ for\/ $i<r,$ and\/ $d_r=\la_k(\be)$. Suppose\/ $\acM_{(\be,\bs d)}^\rst(\ac\tau^{s\ti\la}_{\bs\mu_0})=\acM_{(\be,\bs d)}^\ss(\ac\tau^{s\ti\la}_{\bs\mu_0}),$ which holds except for finitely many values of\/ $s$ in $(0,1)$. Then Claim\/ {\rm\ref{co10cla1}(a)} holds for\/ $x=0,s,(\be,\bs d)$.
\end{prop}

\begin{proof} We will prove the proposition by induction on $d_r=0,1,\ldots.$ The first step $d_r=0$ is vacuous as there are no $(\be,\bs d)$ satisfying the conditions with $d_r=0$. So suppose by induction that $l\ge 0$ is given such that the proposition holds for all such $s\in[0,1]$, $(\be,\bs d)$ with $d_r\le l$, and let $\bar s\in[0,1]$, $(\be,\bs d)$ satisfy the conditions of the proposition with~$d_r=l+1$.

Proposition \ref{co10prop9} gives a finite set of triples $s\in(0,1)$, $(\al_1,\bs d_1),(\al_2,\bs d_2)$ associated to $(\be,\bs d)$. Let $0<s_1<s_2<\cdots<s_p\le\bar s$ be the possible values of $s$ in $(0,\bar s]$ for such triples, and for each $i=1,\ldots,p$, let $\bigl\{(\ga_{i,j},\bs e_{i,j}),(\de_{i,j},\bs f_{i,j})\bigr\}$ for $j=1,\ldots,q_i$ be the distinct possibilities for $\bigl\{(\al_1,\bs d_1),(\al_2,\bs d_2)\bigr\}$ in Proposition \ref{co10prop9} for $s=s_i$, ordered such that $\ac\tau^0_{\bs\mu_0}(\ga_{i,j},\bs e_{i,j})>\ac\tau^0_{\bs\mu_0}(\de_{i,j},\bs f_{i,j})$. This is possible since $s_i$ is the unique $s\in[0,1]$ with $\ac\tau^{s\ti\la}_{\bs\mu_0}(\al_1,\bs d_1)=\ac\tau^{s\ti\la}_{\bs\mu_0}(\al_2,\bs d_2)$ by the proof of Proposition \ref{co10prop9}, so $\ac\tau^0_{\bs\mu_0}(\al_1,\bs d_1)\ne\ac\tau^0_{\bs\mu_0}(\al_2,\bs d_2)$. 

Let $a=1,\ldots,r$ be least such that $d_a>0$. Then $d_a=1$, as $d_1\le 1$ and $d_{l+1}\ab\le\ab d_l+1$ for $l<r$. As $\bs d=\bs e_{i,j}+\bs f_{i,j}$ we see that either
\begin{itemize}
\setlength{\itemsep}{0pt}
\setlength{\parsep}{0pt}
\item[(i)] $e_{i,j,a}=1$ and $f_{i,j,a}=0$; or 
\item[(ii)] $e_{i,j,a}=0$ and $f_{i,j,a}=1$.
\end{itemize}
Expanding $\ac\tau^0_{\bs\mu_0}(\ga_{i,j},\bs e_{i,j})>\ac\tau^0_{\bs\mu_0}(\de_{i,j},\bs f_{i,j})$ using \eq{co5eq21}, since $\mu_1\gg\cdots\gg\mu_r>0$ by \eq{co10eq5}, we see that the dominant contributions come from $\mu_ae_{i,j,a}=\mu_a$ in case (i), and $\mu_af_{i,j,a}=\mu_a$ in case (ii), and thus that case (ii) does not occur.

Since $\mu_1,\ldots,\mu_r$ in \eq{co10eq5} were chosen generic and are independent of $\ti\la$, the equation $\ac\tau^{s_i\ti\la}_{\bs\mu_0}(\ga_{i,j},\bs e_{i,j})=\ac\tau^{s_i\ti\la}_{\bs\mu_0}(\de_{i,j},\bs f_{i,j})$ implies that for each $i$ (which fixes $s_i$) there can be only one possible pair $\{\bs e_{i,j},\bs f_{i,j}\}$, and (i) above distinguishes the order of $\bs e_{i,j},\bs f_{i,j}$. Thus $\bs e_{i,j},\bs f_{i,j}$ are independent of $j=1,\ldots,q_i$, so from now on we will write them as $\bs e_i,\bs f_i$. Define $s'_0=0$ and $s'_p=\bar s$, and choose $s'_i\in(s_i,s_{i+1})$ for $1\le i<p$. Then~$0=s_0'<s_1<s_1'<\cdots<s_{p-1}'<s_p\le s_p'=\bar s$.

Now $(\ga_{i,j},\bs e_i),(\de_{i,j},\bs f_i)$ satisfy the conditions of the proposition by Proposition \ref{co10prop9}, and $e_{i,r},f_{i,r}\le l$ as $d_{i,r}=e_{i,r}+f_{i,r}=l+1$ and $e_{i,r},f_{i,r}>0$, so by induction Claim \ref{co10cla1}(a) holds for $x=0,s_i,(\ga_{i,j},\bs e_i),(\de_{i,j},\bs f_i)$, giving
\e
\begin{split}
q^{\bs e_i}\ot\ac\Up_{(\ga_{i,j},\bs e_i)}(\ac\tau^{s_i\ti\la}_{\bs\mu_0})&=[\acM_{(\ga_{i,j},\bs e_i)}^\ss(\ac\tau^{s_i\ti\la}_{\bs\mu_0})]_\cL,\\ 
q^{\bs f_i}\ot\ac\Up_{(\de_{i,j},\bs f_i)}(\ac\tau^{s_i\ti\la}_{\bs\mu_0})&=[\acM_{(\de_{i,j},\bs f_i)}^\ss(\ac\tau^{s_i\ti\la}_{\bs\mu_0})]_\cL,
\end{split}\qquad \text{all $i,j$.}
\label{co10eq49}
\e

Proposition \ref{co10prop9} gives $\acM_{(\ga_{i,j},\bs e_i)}^\ss(\ac\tau^{s_i\ti\la}_{\bs\mu_0})\ne\es$ or $[\acM_{(\ga_{i,j},\bs e_i)}^\ss(\ac\tau^{s_i\ti\la}_{\bs\mu_0})]_\cL\ne 0$, but then \eq{co10eq49} implies that $\acM_{(\ga_{i,j},\bs e_i)}^\ss(\ac\tau^{s_i\ti\la}_{\bs\mu_0})\ne\es$ in both cases, and similarly $\acM_{(\de_{i,j},\bs f_i)}^\ss(\ac\tau^{s_i\ti\la}_{\bs\mu_0})\ne\es$. Hence $\acM_{(\be,\bs d)}^\rst(\ac\tau^{s_i\ti\la}_{\bs\mu_0})\ne\acM_{(\be,\bs d)}^\ss(\ac\tau^{s_i\ti\la}_{\bs\mu_0})$ as direct sums of semistables in classes $(\ga_{i,j},\bs e_i),(\de_{i,j},\bs f_i)$ give strictly semistables in class $(\be,\bs d)$. Therefore $s_p<s_p'=\bar s$ as $\acM_{(\be,\bs d)}^\rst(\ac\tau^{\bar s\ti\la}_{\bs\mu_0})\ne\acM_{(\be,\bs d)}^\ss(\ac\tau^{\bar s\ti\la}_{\bs\mu_0})$ by assumption.

If $s\in[0,\bar s]\sm\{s_1,\ldots,s_p\}$ then as there are no splittings $(\be,\bs d)=(\al_1,\bs d_1)+\cdots+(\al_n,\bs d_n)$ for $s$ and $n\ge 2$ satisfying the conditions of Proposition \ref{co10prop9} we have $\acM_{(\be,\bs d)}^\rst(\ac\tau^{s\ti\la}_{\bs\mu_0})=\acM_{(\be,\bs d)}^\ss(\ac\tau^{s\ti\la}_{\bs\mu_0})$, so $[\acM_{(\be,\bs d)}^\ss(\ac\tau^{s\ti\la}_{\bs\mu_0})]_\virt$ is defined. 

Observe that $(\ac\tau^{s\ti\la}_{\bs\mu_0},\ac T,\le)_{s\in[0,1]}$ is a continuous family of weak stability conditions on $\acA$ in the sense of Definition \ref{co3def5}. This holds as $x=0$ and $\mu_i>0$, so the third case of \eq{co5eq21} does not happen; when $x<0$ we can get discontinuity in $(\ac\tau^{s\ti\la}_{\bs\mu_x},\ac T,\le)$ by switching between the second and third cases of \eq{co5eq21}. 

We claim that both $\acM_{(\be,\bs d)}^\ss(\ac\tau^{s\ti\la}_{\bs\mu_0})$ and $[\acM_{(\be,\bs d)}^\ss(\ac\tau^{s\ti\la}_{\bs\mu_0})]_\cL$ are locally constant in $s$ in each connected component of $[0,\bar s]\sm\{s_1,\ldots,s_p\}$. For $\acM_{(\be,\bs d)}^\ss(\ac\tau^{s\ti\la}_{\bs\mu_0})$ this holds as $\ac\tau^{s\ti\la}_{\bs\mu_0}$ depends continuously on $s$, and so $\acM_{(\be,\bs d)}^\ss(\ac\tau^{s\ti\la}_{\bs\mu_0})$ can only change by crossing a wall on which there are strictly semistables; for $[\acM_{(\be,\bs d)}^\ss(\ac\tau^{s\ti\la}_{\bs\mu_0})]_\cL$ we can deduce it from Propositions \ref{co3prop4} and \ref{co10prop9} and equation~\eq{co10eq21}.

Claim \ref{co10cla1}(a) holds for $x=0,s'_0=0,(\be,\bs d)$ by Proposition \ref{co10prop10}, giving
\e
q^{\bs d}\ot\ac\Up_{(\be,\bs d)}(\ac\tau^{s_0'\ti\la}_{\bs\mu_0})=[\acM_{(\be,\bs d)}^\ss(\ac\tau^{s_0'\ti\la}_{\bs\mu_0})]_\cL.
\label{co10eq50}
\e
We will show below that in the Lie algebra $\check H_{\rm even}(\M^\pl)$ we have
\e
\begin{split}
&\ac\Up_{(\be,\bs d)}(\ac\tau^{s_i'\ti\la}_{\bs\mu_0})=\ac\Up_{(\be,\bs d)}(\ac\tau^{s_{i-1}'\ti\la}_{\bs\mu_0})\\
&+\sum_{\begin{subarray}{l} j=1,\ldots,q_i: \\ o_{\ga_{i,j}}+o_{\de_{i,j}}=o_\be \end{subarray}\!\!\!\!\!\!\!\!\!\!\!\!\!\!\!}\bigl[\ac\Up_{(\ga_{i,j},\bs e_i)}(\ac\tau^{s_i\ti\la}_{\bs\mu_0}),\ac\Up_{(\de_{i,j},\bs f_i)}(\ac\tau^{s_i\ti\la}_{\bs\mu_0})\bigr].
\end{split}
\label{co10eq51}
\e
To avoid interrupting the current argument, the proof of the \eq{co10eq51} will be deferred to \S\ref{co106}, as it is rather long.

We claim that \eq{co10eq21} with $x=x'=0$, $s'=s_{i-1}'$ and $s=s_i'$ may be written
\e
\begin{split}
&[\acM_{(\be,\bs d)}^\ss(\ac\tau^{s_i'\ti\la}_{\bs\mu_0})]_\cL=[\acM_{(\be,\bs d)}^\ss(\ac\tau^{s_{i-1}'\ti\la}_{\bs\mu_0})]_\cL \\
&+\sum_{\begin{subarray}{l} j=1,\ldots,q_i: \\ o_{\ga_{i,j}}+o_{\de_{i,j}}=o_\be \end{subarray}\!\!\!\!\!\!\!\!\!\!\!\!\!\!\!}\bigl[[\acM_{(\ga_{i,j},\bs e_i)}^\ss(\ac\tau^{s_i\ti\la}_{\bs\mu_0})]_\cL,[\acM_{(\de_{i,j},\bs f_i)}^\ss(\ac\tau^{s_i\ti\la}_{\bs\mu_0})]_\cL\bigr],
\end{split}
\label{co10eq52}
\e
for $i=1,\ldots,p$. To see this, note that as $[\acM_{(\be,\bs d)}^\ss(\ac\tau^{s_i'\ti\la}_{\bs\mu_0})]_\cL$ depends only on the component of $[0,\bar s]\sm\{s_1,\ldots,s_p\}$ containing $s_i'$, and similarly for $s_{i-1}'$, we may take $s_i',s_{i-1}'$ arbitrarily close to $s_i$. Thus Proposition \ref{co3prop4} shows that in \eq{co10eq21} we have $\ti U((\al_1,\bs d_1),\ldots,(\al_n,\bs d_n);\ac\tau^{s_{i-1}'\ti\la}_{\bs\mu_0},\ac\tau^{s_i'\ti\la}_{\bs\mu_0})\ne 0$ only if $\ac\tau^{s_i\ti\la}_{\bs\mu_0}(\al_1,\bs d_1)=\cdots=\ac\tau^{s_i\ti\la}_{\bs\mu_0}(\al_n,\bs d_n)$. Then Proposition \ref{co10prop9} implies that the only possibilities for $(\al_1,\bs d_1),\ldots,(\al_n,\bs d_n)$ giving a nonzero term in \eq{co10eq21} are when $n=1$ and $(\al_1,\bs d_1)=(\be,\bs d)$, or $n=2$ and $\bigl\{(\al_1,\bs d_1),(\al_2,\bs d_2)\bigr\}=\bigl\{(\ga_{i,j},\bs e_i),(\de_{i,j},\bs f_i)\bigr\}$ for $j=1,\ldots,q_i$. Computing the $\ti U(\cdots;\ac\tau^{s_{i-1}'\ti\la}_{\bs\mu_0},\ac\tau^{s_i'\ti\la}_{\bs\mu_0})$ as for \eq{co10eq41} yields~\eq{co10eq52}.

We will prove by induction on $i=0,\ldots,p$ that
\e
q^{\bs d}\ot\ac\Up_{(\be,\bs d)}(\ac\tau^{s_i'\ti\la}_{\bs\mu_0})=[\acM_{(\be,\bs d)}^\ss(\ac\tau^{s_i'\ti\la}_{\bs\mu_0})]_\cL.
\label{co10eq53}
\e
The first step $i=0$ is \eq{co10eq50}. For the inductive step, suppose \eq{co10eq53} holds for $i-1$ for $0<i\le p$. Apply $q^{\bs d}\ot-$ to \eq{co10eq51}, subtract \eq{co10eq52}, and add \eq{co10eq53} for $i-1$. Since the contributions from the sums $\sum_{j=1}^{q_i}\cdots$ in \eq{co10eq51}--\eq{co10eq52} cancel by \eq{co10eq49}, this proves \eq{co10eq53} for~$i$. 

When $i=p$, so $s_i'=\bar s$, equation \eq{co10eq53} proves Claim \ref{co10cla1}(a) holds for $x=0,\bar s,(\be,\bs d)$. This completes the inductive step in our induction on $d_r=0,1,\ldots$ above, and proves the proposition.
\end{proof}

\begin{cor}
\label{co10cor4}
Let\/ $(\be,\bs d)\in\ti S_\al$ satisfy $\be\ne 0\ne\bs d$ and\/ $d_1\le 1,$ $d_i\ab\le\ab d_{i+1}\ab\le\ab d_i+1$ for\/ $i<r,$ and\/ $d_r=\la_k(\be)$. Then in $\cL_{\B,\acB}$ we have
\e
[\acM_{(\be,\bs d)}^\ss(\ac\tau^{\ti\la}_{\bs\mu_0})]_\cL=[\acM_{(\be,\bs d)}^\ss(\ac\tau^{\ti\la}_{\bs\mu_0})]_{\smash{\ti{\cal L}}}.
\label{co10eq54}
\e
In particular, this holds when\/ $\be=\al$ and\/ $\bs d=(1,2,\ldots,r)$.
\end{cor}

\begin{proof} For $(\be,\bs d)$ satisfying the given conditions, Proposition \ref{co10prop10} proves that Claim \ref{co10cla1}(a) holds for $x=0,s=1,(\be,\bs d)$, and Proposition \ref{co10prop8}(b) shows that Claim \ref{co10cla1}(b) holds for $x=0,(\be,\bs d)$. Therefore combining equations \eq{co10eq23} and \eq{co10eq24} proves \eq{co10eq54}. For the last part, as $r=\la_k(\al)$ by Definition \ref{co10def1}, $\be=\al$ and $\bs d=(1,2,\ldots,r)$ satisfy the given conditions.
\end{proof}

We can now complete the proof of Theorem \ref{co5thm2}. Let $(\be,\bs d)\in\ti S_\al$ satisfy $\be\ne 0\ne\bs d$ and $d_1\le 1,$ $d_i\ab\le\ab d_{i+1}\ab\le\ab d_i+1$ for $i<r$, and $d_r=\la_k(\be)$. Then combining equations \eq{co10eq20} for $x=0$, \eq{co10eq25}, and \eq{co10eq54} yields
\ea
&\sum_{\begin{subarray}{l}n\ge 1,\\
(\al_i,\bs d_i)\in\ti S_\al,\\
i=1,\ldots,n:\\ \sum_{i=1}^n(\al_i,\bs d_i)\\
=(\be,\bs d), \\
o_{\al_1}+\cdots+o_{\al_n} \\ =o_\be \end{subarray}} 
\begin{aligned}[t]
&\ti U((\al_1,\bs d_1),\ldots,(\al_n,\bs d_n);\ac\tau^{\ti\la}_{\bs\mu_{\smash{-1}}},\ac\tau^{\ti\la}_{\bs\mu_0})\cdot{}\\
&\Bigl\{\bigl[\bigl[\cdots\bigl[[\acM_{(\al_1,\bs d_1)}^\ss(\ac\tau^{\ti\la}_{\bs\mu_{\smash{-1}}})]_\cL,[\acM_{(\al_2,\bs d_2)}^\ss(\ac\tau^{\ti\la}_{\bs\mu_{\smash{-1}}})]_\cL\bigr],\ldots\bigr],\\
&[\acM_{(\al_n,\bs d_n)}^\ss(\ac\tau^{\ti\la}_{\bs\mu_{\smash{-1}}})]_\cL\bigr]-\bigl[\bigl[\cdots\bigl[[\acM_{(\al_1,\bs d_1)}^\ss(\ac\tau^{\ti\la}_{\bs\mu_{\smash{-1}}})]_{\smash{\ti{\cal L}}},\\
&[\acM_{(\al_2,\bs d_2)}^\ss(\ac\tau^{\ti\la}_{\bs\mu_{\smash{-1}}})]_{\smash{\ti{\cal L}}}\bigr],\ldots\bigr],[\acM_{(\al_n,\bs d_n)}^\ss(\ac\tau^{\ti\la}_{\bs\mu_{\smash{-1}}})]_{\smash{\ti{\cal L}}}\bigr]\Bigr\}=0.
\end{aligned}
\label{co10eq55}
\ea

Observe that if $(\al_i,\bs d_i)\in\ti S_\al$ with $\al_i\ne 0\ne\bs d_i$ then $[\acM_{(\al_i,\bs d_i)}^\ss(\ac\tau^{\ti\la}_{\bs\mu_{\smash{-1}}})]_\cL=0$ by Proposition \ref{co10prop6}(iv), and $[\acM_{(\al_i,\bs d_i)}^\ss(\ac\tau^{\ti\la}_{\bs\mu_{\smash{-1}}})]_{\smash{\ti{\cal L}}}=0$ by \eq{co10eq19}. Thus all nonzero terms in \eq{co10eq55} have $\al_i=0$ or $\bs d_i=0$ for each $i=1,\ldots,n$. If $\al_i=0$ then $[\acM_{(0,\bs d_i)}^\ss(\ac\tau^{\ti\la}_{\bs\mu_{\smash{-1}}})]_\cL=[\acM_{(0,\bs d_i)}^\ss(\ac\tau^{\ti\la}_{\bs\mu_{\smash{-1}}})]_{\smash{\ti{\cal L}}}$ by \eq{co10eq17}, \eq{co10eq19} and the first part of the proof of Proposition \ref{co10prop6}.

We will prove by induction on $\la_k(\be)=1,\ldots,r$ that $[\acM_{(\be,0)}^\ss(\ac\tau^{\ti\la}_{\bs\mu_{\smash{-1}}})]_\cL=[\acM_{(\be,0)}^\ss(\ac\tau^{\ti\la}_{\bs\mu_{\smash{-1}}})]_{\smash{\ti{\cal L}}}$ if $\be\in\ti R_\al$. Suppose by induction that $0\le l<r$ is given with $[\acM_{(\be,0)}^\ss(\ac\tau^{\ti\la}_{\bs\mu_{\smash{-1}}})]_\cL=[\acM_{(\be,0)}^\ss(\ac\tau^{\ti\la}_{\bs\mu_{\smash{-1}}})]_{\smash{\ti{\cal L}}}$ if $\be\in\ti R_\al$ with $\la_k(\be)\le l$ (this assumption is vacuous in the first step $l=0$ as $\la_k(\be)>0$), and let $\be\in\ti R_\al$ with $\la_k(\be)=l+1$. Define $\bs d=(d_1,\ldots,d_r)$ by $d_i=\max(i+\la_k(\be)-r,0)$ for $i=1,\ldots,r$. Then $(\be,\bs d)$ satisfy the conditions above, so \eq{co10eq55} holds.

Suppose $n\ge 1$, $(\al_1,\bs d_1),\ldots,(\al_n,\bs d_n)$ give a nonzero term in the sum \eq{co10eq55}. Then $\al_i=0$ or $\bs d_i=0$ for each $i$ as above. Divide into cases:
\begin{itemize}
\setlength{\itemsep}{0pt}
\setlength{\parsep}{0pt}
\item[(i)] $\al_i\ne 0$ for at least two $i=1,\ldots,n$; and
\item[(ii)] $\al_i=\be$, $\bs d_i=0$ for some $i=1,\ldots,n$, and $\al_j=0$, $\bs d_j\ne 0$ for $j\ne i$.
\end{itemize}
In case (i), if $\al_i\ne 0$ we have $\la_k(\al_i)<\la_k(\be)=l+1$, so $\la_k(\al_i)\le l$, and $[\acM_{(\al_i,0)}^\ss(\ac\tau^{\ti\la}_{\bs\mu_{\smash{-1}}})]_\cL=[\acM_{(\al_i,0)}^\ss(\ac\tau^{\ti\la}_{\bs\mu_{\smash{-1}}})]_{\smash{\ti{\cal L}}}$ by induction. Hence $[\acM_{(\al_i,\bs d_i)}^\ss(\ac\tau^{\ti\la}_{\bs\mu_{\smash{-1}}})]_\cL=[\acM_{(\al_i,\bs d_i)}^\ss(\ac\tau^{\ti\la}_{\bs\mu_{\smash{-1}}})]_{\smash{\ti{\cal L}}}$ for all $i=1,\ldots,n$, and the expression $\{\cdots\}$ in \eq{co10eq55} for $(\al_1,\bs d_1),\ldots,(\al_n,\bs d_n)$ is zero.

Thus the sum of all terms of type (ii) in \eq{co10eq55} is zero. Using the explicit expressions for $[\acM_{(0,\bs d_i)}^\ss(\ac\tau^{\ti\la}_{\bs\mu_{\smash{-1}}})]_\cL=[\acM_{(0,\bs d_i)}^\ss(\ac\tau^{\ti\la}_{\bs\mu_{\smash{-1}}})]_{\smash{\ti{\cal L}}}$ in \eq{co10eq19} and the Lie bracket \eq{co10eq16}, this becomes a large combinatorial sum, linear in $[\acM_{(\be,0)}^\ss(\ac\tau^{\ti\la}_{\bs\mu_{\smash{-1}}})]_\cL-[\acM_{(\be,0)}^\ss(\ac\tau^{\ti\la}_{\bs\mu_{\smash{-1}}})]_{\smash{\ti{\cal L}}}$. We can identify the answer from Corollary \ref{co10cor3}, as the term $n=1$, $\al_1=\be$ in \eq{co10eq47}. Thus \eq{co10eq55} implies that
\begin{equation*}
\la_k(\be)!\cdot q^{\bs d}\ot
\bigl([\acM_{(\be,0)}^\ss(\ac\tau^{\ti\la}_{\bs\mu_{\smash{-1}}})]_\cL-[\acM_{(\be,0)}^\ss(\ac\tau^{\ti\la}_{\bs\mu_{\smash{-1}}})]_{\smash{\ti{\cal L}}}\bigr)=0.
\end{equation*}
Hence $[\acM_{(\be,0)}^\ss(\ac\tau^{\ti\la}_{\bs\mu_{\smash{-1}}})]_\cL=[\acM_{(\be,0)}^\ss(\ac\tau^{\ti\la}_{\bs\mu_{\smash{-1}}})]_{\smash{\ti{\cal L}}}$, proving the inductive step. 

As $\al\in\ti R_\al$ this proves that $[\acM_{(\al,0)}^\ss(\ac\tau^{\ti\la}_{\bs\mu_{\smash{-1}}})]_\cL=[\acM_{(\al,0)}^\ss(\ac\tau^{\ti\la}_{\bs\mu_{\smash{-1}}})]_{\smash{\ti{\cal L}}}$. By \eq{co10eq16} and \eq{co10eq17}--\eq{co10eq19} we may write this as
\ea
&q^0\ot [\M_\al^\ss(\ti\tau)]_\inv=
\label{co10eq56}\\
&\sum_{\begin{subarray}{l}n\ge 1,\;\al_1,\ldots,\al_n\in
R_\al:\\ \al_1+\cdots+\al_n=\al,\\ 
o_{\al_1}+\cdots+o_{\al_n}=o_\al \end{subarray}} \!\!\!\!\!\!\!
\begin{aligned}[t]
\ti U(&(\al_1,0),\ldots,(\al_n,0);\ac\tau^0_{\bs\mu_{\smash{-1}}},\ac\tau^{\ti\la}_{\bs\mu_{-1}})
\cdot  q^0\ot \\
&\bigl[\bigl[\cdots\bigl[[\M_{\al_1}^\ss(\tau)]_\inv,
[\M_{\al_2}^\ss(\tau)]_\inv\bigr],\ldots\bigr],[\M_{\al_n}^\ss(\tau)]_\inv\bigr].
\end{aligned}
\nonumber
\ea
Substituting $U((\al_1,0),\ldots,(\al_n,0);\ac\tau_{\bs\mu_{\smash{-1}}}^0,\ac\tau_{\bs\mu_{\smash{-1}}}^{\ti\la})=U(\al_1,\ldots,\al_n;\tau,\ti\tau)$ by \eq{co10eq12}, omitting $q^0\ot-$, and noting by \eq{co10eq3} that the sum over $\al_i\in R_\al$ in \eq{co10eq56} is equivalent to the sum over $\al_i\in C(\B)_\pe$ with $\tau(\al_i)=\tau(\al)$ and $\M_{\al_i}^\ss(\tau)\ne\es$ in \eq{co5eq31}, this proves \eq{co5eq31}. Equation \eq{co5eq32} is equivalent to \eq{co5eq31} by Theorem \ref{co3thm3}. This completes the proof of Theorem \ref{co5thm2}, except for the proof of \eq{co10eq51} above, which we give in the next section.

\subsection{Proof of equation \eq{co10eq51}}
\label{co106}

Finally we prove \eq{co10eq51}, which is equation \eq{co10eq84} in Corollary \ref{co10cor6} below. Equation \eq{co10eq51} is similar to \eq{co9eq51} in Corollary \ref{co9cor1}, whose proof occupied the whole of \S\ref{co92}. We will prove it by a very similar method, although the moduli spaces involved are more complicated in this case. We begin with an analogue of Definition~\ref{co9def2}:

\begin{dfn}
\label{co10def4}	
Let $(\be,\bs d),s_1,\ldots,s_p$, and $(\ga_{i,j},\bs e_i),(\de_{i,j},\bs f_i)$ for $i=1,\ab\ldots,\ab p$ and $j=1,\ldots,q_i$ be as in the proof of Proposition \ref{co10prop10}. Let $a=1,\ldots,r$ be least such that $d_a>0$. Then $d_a=1$, as $d_1\le 1$ and $d_{l+1}\ab\le\ab d_l+1$ for $l<r$, and as in the proof of Proposition \ref{co10prop10} we have $e_{i,a}=1$ and $f_{i,a}=0$. Let $b_i=1,\ldots,r$ be least such that $f_{i,b_i}>0$. Then $1\le a<b_i\le r$, and~$f_{i,b_i}=1$.

Fix $i=1,\ldots,p$. Define data $\grB\subseteq\grA,K(\grA),\ab\grM,\ab\grM^\pl,\ab\ldots$ as in Definition \ref{co5def1} using the quiver $\gr Q=(Q_0,Q_1,h,t)$ illustrated by:
\ea
\text{\begin{footnotesize}$\displaystyle
\begin{xy}
0;<.77mm,0mm>:
,(-107.5,3)*{e_1}
,(-107.5,-6)*{\gr Q}
,(-42,-4)*{\text{relation}}
,(-115,3)*{v_1}
,(-92,3)*{e_2}
,(-100,3)*{v_2}  
,(-74,3)*{e_{a-1}}
,(-82,0)*+{\cdots}
,(-49,0)*+{\cdots}
,(2,0)*+{\cdots}
,(-57,3)*{e_a}
,(-42,3)*{e_{b_i-2}}
,(-65,3)*{v_a}
,(-32,3)*{v_{b_i-1}}
,(-7.5,3)*{e_{b_i}}
,(-22,3)*{e_{b_i-1}}
,(26.5,3)*{e_r}
,(11.5,3)*{e_{r-1}}
,(-60,-8)*{e_{-1}}
,(-28.5,-8)*{e_0}
,(-14,3)*{v_{b_i}}
,(19,3)*{v_r}
,(34,4)*{w}
,(27,-6)*{\ka(w)\!=\!k,}
,(-15,0)*+{\bu} ; (-2,0)*+{} **@{-} ?>*\dir{>}
,(6,0)*+{} ; (19,0)*+{\bu} **@{-} ?>*\dir{>}
,(19,0)*+{\bu} ; (34,0)*+{\circ} **@{-} ?>*\dir{>}
,(-31,0)*+{\bu} ; (-15,0)*+{\bu} **@{-} ?>*\dir{>}
,(-65,0)*+{\bu} ; (-52,0)*+{} **@{-} ?>*\dir{>}
,(-45,0)*+{} ; (-31,0)*+{\bu} **@{-} ?>*\dir{>}
,(-79,0)*+{} ; (-65,0)*+{\bu} **@{-} ?>*\dir{>}
,(-100,0)*+{\bu} ; (-85,0)*+{} **@{-} ?>*\dir{>}
,(-115,0)*+{\bu} ; (-100,0)*+{\bu} **@{-} ?>*\dir{>}
,(-53.5,-13)*{v_0}
,(-65,0)*+{\bu} ; (-50,-12)*+{\bu} **@{-} ?>*\dir{>}
,(-15,0)*+{\bu} ; (-50,-12)*+{\bu} **@{-} ?>*\dir{>}
,(-31,0)*+{\bu} ; (-50,-12)*+{\bu} **@{..} ?>*\dir{>}
\end{xy}$\end{footnotesize}}
\nonumber \\[-13pt]
\label{co10eq57}
\ea
with vertices $Q_0=\{v_0,v_1,\ldots,v_r,w\}$, edges $Q_1=\{e_{-1},e_0,\ldots,e_r\}$, and head and tail maps $h(e{-1})=h(e_0)=v_0$, $h(e_l)=v_{l+1}$, $1\le l<r$, $h(e_r)=w$, and $t(e{-1})=v_a$, $t(e_0)=v_{b_i}$, $t(e_l)=v_l$, $1\le l\le r$. We take $\dot Q_0=\{v_0,\ldots,v_r\}$ and $\ddot Q_0=\{w\}$, and define $\ka:\ddot Q_0\ra K$ by $\ka(w)=k$. The arrow `$\dashra$' is not an edge, but indicates a relation we will impose below. We write the stability conditions on $\grA$ defined in Definition \ref{co5def1} as $(\gr\tau^\la_{\bs\mu},\gr T,\le)$, which form a set $\grS$. We write $\gr\Up_{(\be,\bs d)},\gr\Up$ for the maps $\bar\Up_{(\be,\bs d)},\bar\Up$ defined in \eq{co5eq15} for $\grB\subseteq\grA,\ldots.$ Note that $\gr Q$ is related to $\ac Q$ in \eq{co10eq1} by adding the vertex $v_0$ and edges~$e_{-1},e_0$.

Here we use notation $\grA,\grB,K(\grA),\grM,(\gr\tau^\la_{\bs\mu},\gr T,\le),\gr\Up_{(\be,\bs d)},\ldots$ rather than $\baA,\ab\baB,\ab K(\baA),\ab\baM,\ab(\bar\tau^\la_{\bs\mu},\bar T,\le),\bar\Up_{(\be,\bs d)},\ab\ldots$ to distinguish them from $\baA,\ab\baB,\ldots$ in Example \ref{co5ex1}, which will also appear in the proof of Proposition \ref{co10prop7} below.

For brevity we will write objects of $\grA$ as $(E,\bs V,\bs\rho)$ with $\bs V=(V_0,\ldots,V_r)$, $\bs\rho=(\rho_{-1},\ldots,\rho_r)$ rather than $(V_{v_0},\ldots,V_{v_r}),(\rho_{e_{-1}},\ldots,\rho_{e_r})$, and write $K(\grA)=K(\A)\t\Z^{r+1}$ rather than $K(\A)\t\Z^{\{v_0,\ldots,v_r\}}$, and write stability conditions as $\gr\tau^\la_{\bs\mu}$ with $\bs\mu=(\mu_0,\ldots,\mu_r)\in\R^{r+1}$ rather than $\bs\mu=(\mu_{v_0},\ldots,\mu_{v_r})$, and so on.

Let $\ep>0$ be very small --- satisfying smallness conditions that will appear after \eq{co10eq58} and in the proof of Proposition \ref{co10prop11}. As in Definition \ref{co5def1}, define a weak stability condition $(\gr\tau^{s_i\ti\la}_{\bs{\gr\mu}},\gr T,\le)$ on $\acA$, where $\bs{\gr\mu}=(-\ep,\mu_1,\ldots,\mu_r)$. 

Form the moduli spaces $\grM^\rst_{(\be,(1,\bs d))}(\gr\tau^{s_i\ti\la}_{\bs{\gr\mu}})\subseteq\grM^\ss_{(\be,(1,\bs d))}(\gr\tau^{s_i\ti\la}_{\bs{\gr\mu}})$. We must have $\grM^\rst_{(\be,(1,\bs d))}(\gr\tau^{s_i\ti\la}_{\bs{\gr\mu}})=\grM^\ss_{(\be,(1,\bs d))}(\gr\tau^{s_i\ti\la}_{\bs{\gr\mu}})$, as in any splitting $(\be,(1,\bs d))=(\ga,(1,\bs e))+(\de,(0,\bs f))$ in $C(\grB)$ we cannot have $\gr\tau^{s_i\ti\la}_{\bs{\gr\mu}}(\ga,(1,\bs e))\ab=\gr\tau^{s_i\ti\la}_{\bs{\gr\mu}}(\de,(0,\bs f))$, since $\ep$ is small and contributes to only one side. Therefore $\grM^\ss_{(\be,(1,\bs d))}(\gr\tau^{s_i\ti\la}_{\bs{\gr\mu}})$ is a proper algebraic space by Assumption~\ref{co5ass2}(h).

Consider the morphism
\e
\begin{split}
\ac\Pi_{(\be,\bs d)}^\pl&:\grM^\ss_{(\be,(1,\bs d))}(\gr\tau^{s_i\ti\la}_{\bs{\gr\mu}})\longra\acM_{(\be,\bs d)}^\pl,\\
\ac\Pi_{(\be,\bs d)}^\pl&:[E,\bs V,\bs\rho]\longmapsto [E,\bs V',\bs\rho']=[E,(V_1,\ldots,V_r),(\rho_1,\ldots,\rho_r)],
\end{split}
\label{co10eq58}
\e
omitting $V_0,\rho_{-1},\rho_0$. In a similar way to Proposition \ref{co10prop1}, one can show that only decompositions $(\be,(1,\bs d))=(\ga,(1,\bs e))+(\de,(0,\bs f))$ for $(\ga,\bs e),(\de,\bs f)$ in the finite set $S_\al$ in \eq{co10eq4} are relevant to determining $\gr\tau^{s_i\ti\la}_{\bs{\gr\mu}}$-semistability in class $(\be,(1,\bs d))$. Since $\gr\tau^{s_i\ti\la}_{\bs{\gr\mu}}(\ga,(l,\bs e))=\ac\tau^{s_i\ti\la}_{{\bs\mu}_0}(\ga,\bs e)+O(\ep)$, if $\ep>0$ is small enough (depending on finitely many smallness conditions involving $(\ga,\bs e),(\de,\bs f)\in S_\al$) we see that \eq{co10eq58} maps to~$\acM_{(\be,\bs d)}^\ss(\ac\tau^{s_i\ti\la}_{\bs\mu_0})$.

Define $\grM^\ss_{(\be,(1,\bs d))}(\gr\tau^{s_i\ti\la}_{\bs{\gr\mu}})_0\subset\grM^\ss_{(\be,(1,\bs d))}(\gr\tau^{s_i\ti\la}_{\bs{\gr\mu}})$ to be the closed $\C$-substack of points $[E,\bs V,\bs\rho]$ such that $\rho_0\ci\rho_{b_i-1}=0$. Then $\grM^\ss_{(\be,(1,\bs d))}(\gr\tau^{s_i\ti\la}_{\bs{\gr\mu}})_0$ is a proper algebraic space. Write 
\e
\ac\Pi^0_{(\be,\bs d)}=\ac\Pi_{(\be,\bs d)}^\pl\vert_{\grM^\ss_{(\be,(1,\bs d))}(\gr\tau^{s_i\ti\la}_{\bs{\gr\mu}})_0}:\grM^\ss_{(\be,(1,\bs d))}(\gr\tau^{s_i\ti\la}_{\bs{\gr\mu}})_0\longra\acM_{(\be,\bs d)}^\ss(\ac\tau^{s_i\ti\la}_{\bs\mu_0}).
\label{co10eq59}
\e

Note that $\rho_{b_i-1}:V_{b_i-1}\ra V_{b_i}$ is injective by Proposition \ref{co10prop2}(a)(ii) as \eq{co10eq58} maps to $\acM_{(\be,\bs d)}^\ss(\ac\tau^{s_i\ti\la}_{\bs\mu_0})$. Also $d_{b_i}=d_{b_i-1}+1$ by choice of $(\be,\bs d)$ and $b_i$, so $V_{b_i}/\rho_{b_i-1}(V_{b_i-1})\cong\C$. Thus if $[E,\bs V,\bs\rho]\in \grM^\ss_{(\be,(1,\bs d))}(\gr\tau^{s_i\ti\la}_{\bs{\gr\mu}})_0$ then $\rho_0$ factors through, and is determined by, a $\C$-linear map
\e
\ti\rho_0:V_{b_i}/\rho_{b_i-1}(V_{b_i-1})\cong\C\longra V_0\cong\C.
\label{co10eq60}
\e
If $[E,\bs V,\bs\rho]\in\grM^\ss_{(\be,(1,\bs d))}(\gr\tau^{s_i\ti\la}_{\bs{\gr\mu}})_0$ then $\rho_{-1},\rho_0$ are not both zero, since otherwise $0\ne(E,(0,V_1,\ldots,V_r),\bs\rho)\subsetneq(E,\bs V,\bs\rho)$ would $\gr\tau^{s_i\ti\la}_{\bs{\gr\mu}}$-destabilize $(E,\bs V,\bs\rho)$. Thus the fibre of \eq{co10eq59} over $[E,\bs V',\bs\rho']\in \acM_{(\be,\bs d)}^\ss(\ac\tau^{s_i\ti\la}_{\bs\mu_0})$ is an open substack of
\begin{align*}
&\bigl((\Hom(V_a,V_0)\op\Hom(V_{b_i}/\rho_{b_i-1}(V_{b_i-1}),V_0))\sm\{(0,0)\}\bigr)\big/\Aut(V_0)\\
&\qquad \cong \bigl(\C^2\sm\{(0,0)\}\bigr)\big/\bG_m=\CP^1.
\end{align*}
Hence $\ac\Pi^0_{(\be,\bs d)}$ in \eq{co10eq59} is representable and smooth of relative dimension~1.

We will define a derived enhancement $\bs\grM^\ss_{(\be,(1,\bs d))}(\gr\tau^{s_i\ti\la}_{\bs{\gr\mu}})_0$ of $\grM^\ss_{(\be,(1,\bs d))}(\gr\tau^{s_i\ti\la}_{\bs{\gr\mu}})_0$. As in Definition \ref{co5def1} we have an open substack $\grM^\ss_{(\be,(1,\bs d))}(\gr\tau^{s_i\ti\la}_{\bs{\gr\mu}})\subseteq\dM_{(\be,(1,\bs d))}^\pl$, where $\dM_{(\be,(1,\bs d))}^\pl$ is the classical truncation of the quasi-smooth derived stack $\bs\dM_{(\be,(1,\bs d))}^\rpl$ defined in \eq{co5eq17}. Write $\bs{\grM}^\ss_{(\be,(1,\bs d))}(\gr\tau^{s_i\ti\la}_{\bs{\gr\mu}})\subseteq\bs\dM_{(\be,(1,\bs d))}^\rpl$ for the corresponding derived open substack, with classical truncation $\grM^\ss_{(\be,(1,\bs d))}(\gr\tau^{s_i\ti\la}_{\bs{\gr\mu}})$. Then $\bs{\grM}^\ss_{(\be,(1,\bs d))}(\gr\tau^{s_i\ti\la}_{\bs{\gr\mu}})$ is quasi-smooth as $\bs\dM_{(\be,(1,\bs d))}^\rpl$ is, and is a proper derived algebraic space as $\grM^\ss_{(\be,(1,\bs d))}(\gr\tau^{s_i\ti\la}_{\bs{\gr\mu}})$ is a proper algebraic space.

We have vector bundles $\bs\cV_0,\bs\cV_{b_i-1}\ra\bs{\grM}^\ss_{(\be,(1,\bs d))}(\gr\tau^{s_i\ti\la}_{\bs{\gr\mu}})$ with fibres $V_0,V_{b_i-1}$ over $(E,\bs V,\bs\rho)$, the pullbacks of vector bundles $\bs\cV_v$ in Definition \ref{co5def1}, and a section $\bs s:\bs{\grM}^\ss_{(\be,(1,\bs d))}(\gr\tau^{s_i\ti\la}_{\bs{\gr\mu}})\ra \bs\cV_{b_i-1}^*\ot\bs\cV_0$ mapping $(E,\bs V,\bs\rho)\mapsto \rho_0\ci\rho_{b_i-1}$. Let $\bs\grM^\ss_{(\be,(1,\bs d))}(\gr\tau^{s_i\ti\la}_{\bs{\gr\mu}})_0=\bs s^{-1}(0)$ be the derived zero section of $\bs s$ in $\bs{\grM}^\ss_{(\be,(1,\bs d))}(\gr\tau^{s_i\ti\la}_{\bs{\gr\mu}})$. Then $\bs\grM^\ss_{(\be,(1,\bs d))}(\gr\tau^{s_i\ti\la}_{\bs{\gr\mu}})_0$ is a proper quasi-smooth derived algebraic space, as these are closed under taking derived zero sections of vector bundles, and has classical truncation $\grM^\ss_{(\be,(1,\bs d))}(\gr\tau^{s_i\ti\la}_{\bs{\gr\mu}})_0$. Thus as in \S\ref{co24} we have an obstruction theory $\bL_i:i^*(\bL_{\bs\grM^\ss_{(\be,(1,\bs d))}(\gr\tau^{s_i\ti\la}_{\bs{\gr\mu}})_0})\ra\bL_{\grM^\ss_{(\be,(1,\bs d))}(\gr\tau^{s_i\ti\la}_{\bs{\gr\mu}})_0}$ on~$\grM^\ss_{(\be,(1,\bs d))}(\gr\tau^{s_i\ti\la}_{\bs{\gr\mu}})_0$.

As in \eq{co9eq10}, define a $\bG_m$-action on $\grM^\ss_{(\be,(1,\bs d))}(\gr\tau^{s_i\ti\la}_{\bs{\gr\mu}})_0$ by, for $\th\in\bG_m$
\e
\th:[E,\bs V,\bs\rho]\longmapsto[E,\bs V,(\th\,\rho_{-1},\rho_0,\rho_1,\ldots,\rho_r)].
\label{co10eq61}
\e
This lifts to the derived enhancement $\bs\grM^\ss_{(\be,(1,\bs d))}(\gr\tau^{s_i\ti\la}_{\bs{\gr\mu}})_0$. Hence the obstruction theory on $\grM^\ss_{(\be,(1,\bs d))}(\gr\tau^{s_i\ti\la}_{\bs{\gr\mu}})_0$ above is $\bG_m$-equivariant. In \eq{co10eq80} we will apply Corollary \ref{co2cor2} to~$\grM^\ss_{(\be,(1,\bs d))}(\gr\tau^{s_i\ti\la}_{\bs{\gr\mu}})_0$.
\end{dfn}

For the rest of the section we work in the situation of Definition \ref{co10def4}. The next two results are analogues of Propositions \ref{co9prop2} and~\ref{co9prop3}:

\begin{prop}
\label{co10prop11}
In the situation of Definition\/ {\rm\ref{co10def4},} the\/ $\bG_m$-fixed substack\/ $\grM^\ss_{(\be,(1,\bs d))}(\gr\tau^{s_i\ti\la}_{\bs{\gr\mu}})_0^{\bG_m}$ is the disjoint union of the following pieces:
\begin{itemize}
\setlength{\itemsep}{0pt}
\setlength{\parsep}{0pt}
\item[{\bf(a)}] The substack\/ $\grM^\ss_{(\be,(1,\bs d))}(\gr\tau^{s_i\ti\la}_{\bs{\gr\mu}})_{\rho_{-1}=0}$ of\/ $\C$-points $[E,\bs V,\bs\rho]$ with\/ $\rho_{-1}=0$. All such\/ $[E,\bs V,\bs\rho]$ have $\rho_0\ne 0$. There is an isomorphism of stacks
\e
\begin{split}
\Pi_{\rho_{-1}=0}&:\grM^\ss_{(\be,(1,\bs d))}(\gr\tau^{s_i\ti\la}_{\bs{\gr\mu}})_{\rho_{-1}=0}\longra \acM_{(\be,\bs d)}^\ss(\ac\tau^{\smash{s_{i-1}'\ti\la}}_{\bs\mu_0}),\\
\Pi_{\rho_{-1}=0}&:[E,\bs V,\bs\rho]\longmapsto [E,(V_1,\ldots,V_r),(\rho_1,\ldots,\rho_r)].
\end{split}
\label{co10eq62}
\e
\item[{\bf(b)}] The substack\/ $\grM^\ss_{(\be,(1,\bs d))}(\gr\tau^{s_i\ti\la}_{\bs{\gr\mu}})_{\rho_0=0}$ of\/ $\C$-points $[E,\bs V,\bs\rho]$ with\/ $\rho_0=0$. All such\/ $[E,\bs V,\bs\rho]$ have $\rho_{-1}\ne 0$. There is an isomorphism of stacks
\e
\begin{split}
\Pi_{\rho_0=0}&:\grM^\ss_{(\be,(1,\bs d))}(\gr\tau^{s_i\ti\la}_{\bs{\gr\mu}})_{\rho_0=0}\longra \acM_{(\be,\bs d)}^\ss(\ac\tau^{s_i'\ti\la}_{\bs\mu_0}),\\
\Pi_{\rho_0=0}&:[E,\bs V,\bs\rho]\longmapsto [E,(V_1,\ldots,V_r),(\rho_1,\ldots,\rho_r)].
\end{split}
\label{co10eq63}
\e
\item[{\bf(c)}] For\/ $j=1,\ldots,q_i,$ there is a substack\/ $\grM^\ss_{(\be,(1,\bs d))}(\gr\tau^{s_i\ti\la}_{\bs{\gr\mu}})_j$ with\/ $\C$-points $[E,\bs V,\bs\rho]=[E'\op E'',(V_0,\bs V'\op\bs V''),(\rho_{-1},\rho_0,\bs\rho'\op\bs\rho'')],$ where\/ $(E',\bs V',\bs\rho'),\ab(E'',\bs V'',\bs\rho'')$ lie in $\acB$ with\/ $\lb E',\bs V',\bs\rho'\rb=(\ga_{i,j},\bs e_i)$ and\/ $\lb E'',\bs V'',\bs\rho''\rb=(\de_{i,j},\bs f_i),$ and\/ $V_a'=V_a,$ $V_a''=0$ with\/ $\rho_{-1}:V_a'\ra V_0$ an isomorphism, and\/ $\rho_0:V_{b_i}\ra V_0$ is zero on $V_{b_i}'$ and an isomorphism $V_{b_i}''\ra V_0$. We take $\th\in\bG_m$ to act by $\th\,\id_{E'}+\id_{E''}$ on $E=E'\op E'',$ and by\/ $\id_{V_0}$ on\/ $V_0,$ and by $\th\,\id_{V_l'}+\id_{V_l''}$ on $V_l=V_l'\op V_l''$ for $l=1,\ldots,r$. Then\/ $\rho_{-1},\ldots,\rho_r$ are $\bG_m$-equivariant under the action {\rm\eq{co10eq61},} so $[E,\bs V,\bs\rho]$ is a $\bG_m$-fixed point. The splittings $E=E'\op E'',$ $(V_1,\ldots,V_r)=\bs V'\op\bs V''$ are uniquely determined by these conditions. There is an isomorphism of stacks
\end{itemize}
\e
\begin{split}
\Pi_j&:\grM^\ss_{(\be,(1,\bs d))}(\gr\tau^{s_i\ti\la}_{\bs{\gr\mu}})_j\longra \acM_{(\ga_{i,j},\bs e_i)}^\ss(\ac\tau^{s_i\ti\la}_{\bs\mu_0})\t\acM_{(\de_{i,j},\bs f_i)}^\ss(\ac\tau^{s_i\ti\la}_{\bs\mu_0}),\\
\Pi_j&:[E,\bs V,\bs\rho]\longmapsto \bigl([E',\bs V',\bs\rho'],[E'',\bs V'',\bs\rho'']\bigr).
\end{split}
\label{co10eq64}
\e
\end{prop}

\begin{proof} Let $[E,\bs V,\bs\rho]\in\grM^\ss_{(\be,(1,\bs d))}(\gr\tau^{s_i\ti\la}_{\bs{\gr\mu}})_0^{\bG_m}$. As in Definition \ref{co10def4}, $\rho_{-1},\rho_0$ cannot both be zero. Thus we can divide into 
three cases:
\begin{itemize}
\setlength{\itemsep}{0pt}
\setlength{\parsep}{0pt}
\item[(a$)'$] $\rho_{-1}=0$, $\rho_0\ne 0$; 
\item[(b$)'$] $\rho_0=0$, $\rho_{-1}\ne 0$; and
\item[(c$)'$] $\rho_{-1},\rho_0\ne 0$.
\end{itemize}

It is clear from \eq{co10eq61} that any $[E,\bs V,\bs\rho]$ in $\grM^\ss_{(\be,(1,\bs d))}(\gr\tau^{s_i\ti\la}_{\bs{\gr\mu}})$ with $\rho_{-1}=0$ is $\bG_m$-invariant, so $\grM^\ss_{(\be,(1,\bs d))}(\gr\tau^{s_i\ti\la}_{\bs{\gr\mu}})_{\rho_{-1}=0}\subseteq\grM^\ss_{(\be,(1,\bs d))}(\gr\tau^{s_i\ti\la}_{\bs{\gr\mu}})_0^{\bG_m}$, as in (a). If $[E,\bs V,\bs\rho]\in\grM^\ss_{(\be,(1,\bs d))}(\gr\tau^{s_i\ti\la}_{\bs{\gr\mu}})_{\rho_0=0}$ then letting $\th\in\bG_m$ act on $E,V_0,\ldots,\ab V_r$ by $\id_E,\th^{-1}\,\id_{V_0},\id_{V_1},\ldots,\id_{V_r}$ shows $\grM^\ss_{(\be,(1,\bs d))}(\gr\tau^{s_i\ti\la}_{\bs{\gr\mu}})_{\rho_0=0}\subseteq\ab\grM^\ss_{(\be,(1,\bs d))}(\gr\tau^{s_i\ti\la}_{\bs{\gr\mu}})_0^{\bG_m}$, as in~(b). 

If $[E,\bs V,\bs\rho]$ is a $\bG_m$-fixed point of type (c$)'$ then there must exist morphisms from $\bG_m$ to $\Aut(E)$ and $\Aut(V_l)$ for $l=0,\ldots,r$ such that $\rho_{-1},\ldots,\rho_r$ are $\bG_m$-equivariant under \eq{co10eq61}. As $V_0\cong\C$ and we are free to change the actions on $E,V_l$ by an overall factor, we can choose the $\bG_m$-action on $V_0$ to be trivial.

The actions of $\bG_m$ on $E,V_l$ correspond to direct sums $E=\bigop_{n\in\Z}E^n$, $V_l=\bigop_{n\in\Z}V_l^n$, where $E^n,V_l^n\ne 0$ for only finitely many $n\in\Z$, and $\th\in\bG_m$ acts by $\sum_{n\in\Z}\th^n\,\id_{E^n}$, $\sum_{n\in\Z}\th^n\,\id_{V_l^n}$. For $\rho_{-1},\rho_0$ to be $\bG_m$-equivariant we see that $V_a=V_a^1$ for $\rho_{-1}:V_a\ra V_0$ to be an isomorphism, and $V_{b_i}^0\ne 0$ for $\ti\rho_0$ in \eq{co10eq60} to be an isomorphism as $\rho_0\ne 0$. Then $E^n=V_l^n=0$ for $n\ne 0,1$, as otherwise the subobject $(\hat E,\bs{\hat V},\bs{\hat\rho})\subsetneq(E,\bs V,\bs\rho)$ with $\hat E=E^0\op E^1$, $\hat V_l=V_l^0\op V_l^1$ would $\gr\tau^{s_i\ti\la}_{\bs{\gr\mu}}$-destabilize $(E,\bs V,\bs\rho)$. Thus writing $E'=E^1$, $E''=E^0$, $V_l'=V_l^1$, $V_l''=V_l^0$ we see that $E=E'\op E''$, $V_l=V_l'\op V_l''$ and $\th\in\bG_m$ acts on $V_0$ by $\id_{V_0}$, and on $E,V_l$ by $\th\,\id_{E'}+\id_{E''},\th\,\id_{V_l'}+\id_{V_l''}$ for $l>0$, and $V_a'=V_a,$ $V_a''=0$ with $\rho_{-1}:V_a'\ra V_0$ an isomorphism, and $\rho_0(V_{b_i}')=0$, as in~(c).

Thus in \eq{co10eq59} we have $\ac\Pi_{(\be,\bs d)}^0([E,\bs V,\bs\rho])=[(E',\bs V',\bs\rho')\op(E'',\bs V'',\bs\rho'')]$, so $(E',\bs V',\bs\rho')\op(E'',\bs V'',\bs\rho'')$ is strictly $\ac\tau^{s_i\ti\la}_{\bs\mu_0}$-semistable. The proof of Proposition \ref{co10prop10} showed this forces $\bigl\{\lb E',\bs V',\bs\rho'\rb,\lb E'',\bs V'',\bs\rho''\rb\}=\bigl\{(\ga_{i,j},\bs e_i),(\de_{i,j},\bs f_i)\bigr\}$ for some $j=1,\ldots,q_i$. As $V_a'=V_a$ and $f_{i,a}=0$, we see that $\lb E',\bs V',\bs\rho'\rb=(\ga_{i,j},\bs e_i)$, $\lb E'',\bs V'',\bs\rho''\rb=(\de_{i,j},\bs f_i)$, as in (c). Write $\grM^\ss_{(\be,(1,\bs d))}(\gr\tau^{s_i\ti\la}_{\bs{\gr\mu}})_j$ for the substack of points of type (c$)'$ with $\lb E',\bs V',\bs\rho'\rb=(\ga_{i,j},\bs e_i)$ and~$\lb E'',\bs V'',\bs\rho''\rb=(\de_{i,j},\bs f_i)$.

We have now defined substacks $\grM^\ss_{(\be,(1,\bs d))}(\gr\tau^{s_i\ti\la}_{\bs{\gr\mu}})_{\rho_{-1}=0},\grM^\ss_{(\be,(1,\bs d))}(\gr\tau^{s_i\ti\la}_{\bs{\gr\mu}})_{\rho_0=0}$ and $\grM^\ss_{(\be,(1,\bs d))}(\gr\tau^{s_i\ti\la}_{\bs{\gr\mu}})_j$ for $j=1,\ldots,q_i$. These are closed in $\grM^\ss_{(\be,(1,\bs d))}(\gr\tau^{s_i\ti\la}_{\bs{\gr\mu}})_0$, $\bG_m$-fixed, disjoint, and include all $\bG_m$-fixed points. Hence they are open and closed, and $\grM^\ss_{(\be,(1,\bs d))}(\gr\tau^{s_i\ti\la}_{\bs{\gr\mu}})_0^{\bG_m}$ is the disjoint union of the substacks~(a)--(c).

To prove \eq{co10eq62}--\eq{co10eq63}, we first claim that
\ea
\acM_{(\be,\bs d)}^\ss(\ac\tau^{\smash{s_{i-1}'\ti\la}}_{\bs\mu_0})&=\bigl\{[E,\bs V,\bs\rho]\in\acM_{(\be,\bs d)}^\ss(\ac\tau^{s_i\ti\la}_{\bs\mu_0}):\nexists\, (E',\bs V',\bs\rho')\subset(E,\bs V,\bs\rho),
\nonumber\\
&\qquad \lb E',\bs V',\bs\rho'\rb=(\ga_{i,j},\bs e_j),\;\text{some $j=1,\ldots,q_i$}\bigr\},
\label{co10eq65}\\
\acM_{(\be,\bs d)}^\ss(\ac\tau^{\smash{s_i'\ti\la}}_{\bs\mu_0})&=\bigl\{[E,\bs V,\bs\rho]\in\acM_{(\be,\bs d)}^\ss(\ac\tau^{s_i\ti\la}_{\bs\mu_0}):\nexists\, (E',\bs V',\bs\rho')\subset(E,\bs V,\bs\rho),
\nonumber\\
&\qquad \lb E',\bs V',\bs\rho'\rb=(\de_{i,j},\bs f_j),\;\text{some $j=1,\ldots,q_i$}\bigr\}.
\label{co10eq66}
\ea

To show these, note that as in the proof of Proposition \ref{co10prop10}, $\acM_{(\be,\bs d)}^\ss(\ac\tau^{s\ti\la}_{\bs\mu_0})$ is constant for $s\in(s_{i-1},s_i)$, where $s_{i-1}'\in (s_{i-1},s_i)$, and for $s\in(s_i,s_{i+1})$, where $s_i'\in(s_i,s_{i+1})$. As $\ac\tau^{s\ti\la}_{\bs\mu_0}$ is upper semicontinuous in $s$, it follows that
\begin{equation*}
\acM_{(\be,\bs d)}^\ss(\ac\tau^{\smash{s_{i-1}'\ti\la}}_{\bs\mu_0}),\acM_{(\be,\bs d)}^\ss(\ac\tau^{\smash{s_i'\ti\la}}_{\bs\mu_0})\subseteq\acM_{(\be,\bs d)}^\ss(\ac\tau^{s_i\ti\la}_{\bs\mu_0}).
\end{equation*}
Also, $[E,\bs V,\bs\rho]\in\acM_{(\be,\bs d)}^\ss(\ac\tau^{s_i\ti\la}_{\bs\mu_0})$ does not lie in $\acM_{(\be,\bs d)}^\ss(\ac\tau^{\smash{s_{i-1}'\ti\la}}_{\bs\mu_0})$ (or $\acM_{(\be,\bs d)}^\ss(\ac\tau^{\smash{s_i'\ti\la}}_{\bs\mu_0})$) if and only if there exists  $0\ne(E',\bs V',\bs\rho')\subsetneq(E,\bs V,\bs\rho)$ with $\ac\tau^{s_i\ti\la}_{\bs\mu_0}(\lb E',\bs V',\bs\rho'\rb)=\ac\tau^{s_i\ti\la}_{\bs\mu_0}(\lb E/E',\bs V/\bs V',\bs\si'\rb)$ and $\ac\tau^{\smash{s_{i-1}'\ti\la}}_{\bs\mu_0}(\lb E',\bs V',\bs\rho'\rb)\!>\!\ac\tau^{\smash{s_{i-1}'\ti\la}}_{\bs\mu_0}(\lb E/E',\bs V/\bs V',\bs\si'\rb)$ (or $\ac\tau^{\smash{s_i'\ti\la}}_{\bs\mu_0}(\lb E',\bs V',\bs\rho'\rb)>\ac\tau^{\smash{s_i'\ti\la}}_{\bs\mu_0}(\lb E/E',\bs V/\bs V',\bs\si'\rb)$, respectively). From the proof of Proposition \ref{co10prop10}, $\ac\tau^{s_i\ti\la}_{\bs\mu_0}(\lb E',\bs V',\bs\rho'\rb)=\ac\tau^{s_i\ti\la}_{\bs\mu_0}(\lb E/E',\ab\bs V/\bs V',\ab\bs\si'\rb)$ implies that $\bigl\{\lb E',\bs V',\bs\rho'\rb,\lb E/E',\bs V/\bs V',\bs\si'\rb\bigr\}=\bigl\{(\ga_{i,j},\bs e_j),(\de_{i,j},\bs f_j)\bigr\}$ for $j=1,\ldots,q_i$. Since 
\begin{equation*}
\ac\tau^{\smash{s_{i-1}'\ti\la}}_{\bs\mu_0}(\ga_{i,j},\bs e_j)>	\ac\tau^{\smash{s_{i-1}'\ti\la}}_{\bs\mu_0}(\de_{i,j},\bs f_j)\quad\text{and}\quad \ac\tau^{\smash{s_i'\ti\la}}_{\bs\mu_0}(\ga_{i,j},\bs e_j)<	\ac\tau^{\smash{s_i'\ti\la}}_{\bs\mu_0}(\de_{i,j},\bs f_j),
\end{equation*}
equations \eq{co10eq65}--\eq{co10eq66} follow.

For \eq{co10eq62}, let $[E,\bs V,\bs\rho]\in \grM^\ss_{(\be,(1,\bs d))}(\gr\tau^{s_i\ti\la}_{\bs{\gr\mu}})_{\rho_{-1}=0}$. Suppose there is a subobject $(E',\bs V',\bs\rho')\!\subset\!(E,(V_1,\ldots,V_r),(\rho_1,\ldots,\rho_r))$ in $\acB$ with $\lb E',\bs V',\bs\rho'\rb\!=\!(\ga_{i,j},\bs e_j)$. Now $\rho_{b_i-1}:V_{b_i-1}'\ra V_{b_i}'$ is an isomorphism, as it is injective by by Proposition \ref{co10prop2}(a)(ii), and $\dim V_{b_i-1}'=\dim V_{b_i}'$ as $\lb E',\bs V',\bs\rho'\rb=(\ga_{i,j},\bs e_j)$. But $\rho_0\ci\rho_{b_i-1}=0$ as $[E,\bs V,\bs\rho]\in \grM^\ss_{(\be,(1,\bs d))}(\gr\tau^{s_i\ti\la}_{\bs{\gr\mu}})_0$, so $\rho_0(V_{b_i}')=0$. Since $\rho_{-1}=0$, it follows that $(E',(0,\bs V'),(0,0,\bs\rho'))\subset(E,\bs V,\bs\rho)$ is a subobject in $\grB$, which $\gr\tau^{s_i\ti\la}_{\bs{\gr\mu}}$-destabilizes $(E,\bs V,\bs\rho)$, a contradiction. As no such subobject exists, equations \eq{co10eq59} and \eq{co10eq65} show that $[E,(V_1,\ldots,V_r),(\rho_1,\ldots,\rho_r)]\in \acM_{(\be,\bs d)}^\ss(\ac\tau^{\smash{s_{i-1}'\ti\la}}_{\bs\mu_0})$. 

Thus the morphism \eq{co10eq62} is well defined. To see that it is an isomorphism, note that we can define an inverse by
\begin{align*}
&\Pi_{\rho_{-1}=0}^{-1}:[E,(V_1,\ldots,V_r),(\rho_1,\ldots,\rho_r)]\longmapsto[E,\bs V,\bs\rho]\\
&\qquad =[E,(V_{b_i}/\rho_{b_i-1}(V_{b_i-1}),V_1,\ldots,V_r),(0,\pi_{V_{b_i}/\rho_{b_i-1}(V_{b_i-1})},\rho_1,\ldots,\rho_r)]
\end{align*}
where $V_0:=V_{b_i}/\rho_{b_i-1}(V_{b_i-1})\cong\C$ as in \eq{co10eq60} and $\pi_{V_{b_i}/\rho_{b_i-1}(V_{b_i-1})}:V_{b_i}\ra V_0$ is the projection. This completes part (a).

For \eq{co10eq63}, let $[E,\bs V,\bs\rho]\in \grM^\ss_{(\be,(1,\bs d))}(\gr\tau^{s_i\ti\la}_{\bs{\gr\mu}})_{\rho_0=0}$. Suppose there is a subobject $(E',\bs V',\bs\rho')\!\subset\!(E,(V_1,\ldots,V_r),(\rho_1,\ldots,\rho_r))$ with $\lb E',\bs V',\bs\rho'\rb\!=\!(\de_{i,j},\bs f_j)$. Then $V_a'\!=\!0$ as $f_{i,a}\!=\!0$, so $\rho_{-1}(V_a')\!=\!0$, and $\rho_0\!=\!0$. Thus $(E',(0,\bs V'),(0,0,\bs\rho'))\ab\subset(E,\bs V,\bs\rho)$ is a subobject, which $\gr\tau^{s_i\ti\la}_{\bs{\gr\mu}}$-destabilizes $(E,\bs V,\bs\rho)$, a contradiction. As no such subobject exists, $[E,(V_1,\ldots,V_r),(\rho_1,\ldots,\rho_r)]\in \acM_{(\be,\bs d)}^\ss(\ac\tau^{\smash{s_i'\ti\la}}_{\bs\mu_0})$ by \eq{co10eq59} and \eq{co10eq66}. Thus \eq{co10eq63} is well defined. We define an inverse by
\begin{align*}
\Pi_{\rho_0=0}^{-1}:[E,(V_1,\ldots,V_r),(\rho_1,\ldots,\rho_r)]\!\mapsto\![E,(V_a,V_1,\ldots,V_r),(\id_{V_a},0,\rho_1,\ldots,\rho_r)].
\end{align*}
Hence \eq{co10eq63} is an isomorphism, proving~(b).

For \eq{co10eq64}, if $[E,\bs V,\bs\rho]\in\grM^\ss_{(\be,(1,\bs d))}(\gr\tau^{s_i\ti\la}_{\bs{\gr\mu}})_j$ then as above we construct $(E',\bs V',\bs\rho'),(E'',\bs V'',\bs\rho'')$ with $(E',\bs V',\bs\rho')\op(E'',\bs V'',\bs\rho'')$ i$\ac\tau^{s_i\ti\la}_{\bs\mu_0}$-semistable, so that $(E',\bs V',\bs\rho'),(E'',\bs V'',\bs\rho'')$ are $\ac\tau^{s_i\ti\la}_{\bs\mu_0}$-semistable, with $\lb E',\bs V',\bs\rho'\rb=(\ga_{i,j},\bs e_i)$ and $\lb E'',\bs V'',\bs\rho''\rb=(\de_{i,j},\bs f_i)$. Therefore $\Pi_j$ in \eq{co10eq64} is well defined. 

If $([E',\bs V',\bs\rho'],[E'',\bs V'',\bs\rho''])$ is a $\C$-point of $\acM_{(\ga_{i,j},\bs e_i)}^\ss(\ac\tau^{s_i\ti\la}_{\bs\mu_0})\t\acM_{(\de_{i,j},\bs f_i)}^\ss(\ac\tau^{s_i\ti\la}_{\bs\mu_0})$ then as $V_a'\cong\C\cong V_{b_i}''$ as $e_{i,a}=f_{i,b_i}=1$ we may take $V_0=V''_{b_i}$, and choose any isomorphism $\rho_{-1}:V'_a\ra V''_{b_i}$, and show that
\begin{equation*}
[E,\bs V,\bs\rho]=[E'\op E'',(V''_{b_i},\bs V'\op\bs V''),(\rho_{-1},0\op\id_{V_{b_i}''},\bs\rho'\op\bs\rho')]
\end{equation*}
is a $\C$-point of $\grM^\ss_{(\be,(1,\bs d))}(\gr\tau^{s_i\ti\la}_{\bs{\gr\mu}})_j$. As $[E,\bs V,\bs\rho]$ is unchanged by $\rho_{-1}\mapsto\th\,\rho_{-1}$ by \eq{co10eq61}, it is independent of the choice of isomorphism $\rho_{-1}:V'_a\ra V''_{b_i}$. Thus we may define a stack morphism $([E',\bs V',\bs\rho'],[E'',\bs V'',\bs\rho''])\mapsto[E,\bs V,\bs\rho]$ inverse to \eq{co10eq64}, so \eq{co10eq64} is an isomorphism. This completes~(c).
\end{proof}

\begin{prop}
\label{co10prop12}
As in\/ {\rm\S\ref{co26},} the restriction of the $\bG_m$-equivariant obstruction theory\/ $\bL_i:i^*(\bL_{\bs\grM^\ss_{(\be,(1,\bs d))}(\gr\tau^{s_i\ti\la}_{\bs{\gr\mu}})_0})\ra\bL_{\grM^\ss_{(\be,(1,\bs d))}(\gr\tau^{s_i\ti\la}_{\bs{\gr\mu}})_0}$ on $\grM^\ss_{(\be,(1,\bs d))}(\gr\tau^{s_i\ti\la}_{\bs{\gr\mu}})_0$ to $\grM^\ss_{(\be,(1,\bs d))}(\gr\tau^{s_i\ti\la}_{\bs{\gr\mu}})_0^{\bG_m}$ splits as the direct sum of a $\bG_m$-fixed part, which is the natural obstruction theory $\bL_{i^{\bG_m}}:(i^{\bG_m})^*(\bs\grM^\ss_{(\be,(1,\bs d))}(\gr\tau^{s_i\ti\la}_{\bs{\gr\mu}})_0^{\bG_m})\ra \bL_{\grM^\ss_{(\be,(1,\bs d))}(\gr\tau^{s_i\ti\la}_{\bs{\gr\mu}})_0^{\bG_m}}$ from the proper quasi-smooth derived algebraic space $\bs\grM^\ss_{(\be,(1,\bs d))}(\gr\tau^{s_i\ti\la}_{\bs{\gr\mu}})_0^{\bG_m},$ and a\/ $\bG_m$-moving part\/ $\cN^\bu,$ the \begin{bfseries}virtual conormal bundle\end{bfseries} of\/ $\grM^\ss_{(\be,(1,\bs d))}(\gr\tau^{s_i\ti\la}_{\bs{\gr\mu}})_0^{\bG_m}$.

In Proposition\/ {\rm\ref{co10prop11}(a)--(c)} we may identify these as follows:
\begin{itemize}
\setlength{\itemsep}{0pt}
\setlength{\parsep}{0pt}
\item[{\bf(a)}] The obstruction theory on $\grM^\ss_{(\be,(1,\bs d))}(\gr\tau^{s_i\ti\la}_{\bs{\gr\mu}})_{\rho_{-1}=0}$ is identified by \eq{co10eq62} with that on $\acM_{(\be,\bs d)}^\ss(\ac\tau^{s_{i-1}'\ti\la}_{\bs\mu_0})$ in Definition\/ {\rm\ref{co10def1}}. The virtual conormal bundle\/ $\cN_{\rho_{-1}=0}^\bu$ is a line bundle with\/ $\bG_m$-weight\/~$-1$.
\item[{\bf(b)}] The obstruction theory on $\grM^\ss_{(\be,(1,\bs d))}(\gr\tau^{s_i\ti\la}_{\bs{\gr\mu}})_{\rho_0=0}$ is identified by \eq{co10eq63} with that on $\acM_{(\be,\bs d)}^\ss(\ac\tau^{s_i'\ti\la}_{\bs\mu_0})$ in Definition\/ {\rm\ref{co10def1}}. The virtual conormal bundle $\cN_{\rho_0=0}^\bu$ is a line bundle with\/ $\bG_m$-weight\/ $1$.
\item[{\bf(c)}] The virtual conormal bundle of\/ $\grM^\ss_{(\be,(1,\bs d))}(\gr\tau^{s_i\ti\la}_{\bs{\gr\mu}})_j$ in $\grM^\ss_{(\be,(1,\bs d))}(\gr\tau^{s_i\ti\la}_{\bs{\gr\mu}})_0$ is $\cN_j^\bu\ab =\cN_j^{+,\bu}\op \cN_j^{-,\bu},$ where $\cN_j^{\pm,\bu}$ has\/ $\bG_m$-weight\/ $\pm 1,$ and using the identification {\rm\eq{co10eq64},} in $K^0(\Perf_{\acM_{(\ga_{i,j},\bs e_i)}^\ss(\ac\tau^{s_i\ti\la}_{\bs\mu_0})\t\acM_{(\de_{i,j},\bs f_i)}^\ss(\ac\tau^{s_i\ti\la}_{\bs\mu_0})})$ we have
\end{itemize}	
\ea
\bigl[\cN_j^{+,\bu}\bigr]&=\sum_{1\le l<r\!\!\!}\bigl[[\cV'_l\bt(\cV_{l+1}'')^*]^\pl\bigr]+\bigl[[\cV'_r\bt(\cV_{k,\de_{i,j}}'')^*]^\pl\bigr]-\sum_{1\le l\le r\!\!\!}\bigl[[\cV'_l\bt(\cV_l'')^*]^\pl\bigr]
\nonumber\\
&-\bigl[(\ac\Pi^{v_a}_{\dot{\mathcal M}_{\ga_{i,j}}}\t\ac\Pi^{v_{b_i}}_{\dot{\mathcal M}_{\de_{i,j}}})^*(\cE^\bu_{\ga_{i,j},\de_{i,j}})\bigr],
\label{co10eq67}\\
\bigl[\cN_j^{-,\bu}\bigr]&=\sum_{1\le l<r\!\!\!}\bigl[[(\cV_{l+1}')^*\bt\cV''_l]^\pl\bigr]+\bigl[(\cV_{k,\ga_{i,j}}')^*\bt\cV''_r]^\pl\bigr]
-\sum_{1\le l\le r\!\!\!}\bigl[[(\cV_l')^*\bt\cV''_l]^\pl\bigr]
\nonumber\\
&-\bigl[(\ac\Pi^{v_a}_{\dot{\mathcal M}_{\ga_{i,j}}}\t\ac\Pi^{v_{b_i}}_{\dot{\mathcal M}_{\de_{i,j}}})^*(\si^*(\cE_{\de_{i,j},\ga_{i,j}}^\bu)\bigr].
\label{co10eq68}
\ea
\begin{itemize}
\setlength{\itemsep}{0pt}
\setlength{\parsep}{0pt}
\item[] Here $\cV',\cV''$ denote vector bundles on $\acM_{(\ga_{i,j},\bs e_i)}^\ss(\ac\tau^{s_i\ti\la}_{\bs\mu_0}),\acM_{(\de_{i,j},\bs f_i)}^\ss(\ac\tau^{s_i\ti\la}_{\bs\mu_0}),$ respectively, and\/ $[\cdots]^\pl$ indicates a vector bundle descending from $\acM_{(\be,\bs d)}$ to $\acM^\pl_{(\be,\bs d)},$ and\/ $\ac\Pi^{v_a}_{\dot{\mathcal M}_{\ga_{i,j}}}:\acM_{(\ga_{i,j},\bs e_i)}^\ss(\ac\tau^{s_i\ti\la}_{\bs\mu_0})\ra\dM_{\ga_{i,j}},$ $\ac\Pi^{v_{b_i}}_{\dot{\mathcal M}_{\de_{i,j}}}:\acM_{(\de_{i,j},\bs f_i)}^\ss(\ac\tau^{s_i\ti\la}_{\bs\mu_0})\ab\ra\dM_{\de_{i,j}}$ are as in equation\/~\eq{co5eq26}.
\item[{\bf(d)}] For\/ $j=1,\ldots,q_i,$ we have $o_{\ga_{i,j}}+o_{\de_{i,j}}\ge o_\be$. Under the isomorphism {\rm\eq{co10eq64},} the virtual classes in Proposition\/ {\rm\ref{co10prop11}(c)} satisfy
\ea
&[\grM^\ss_{(\be,(1,\bs d))}(\gr\tau^{s_i\ti\la}_{\bs{\gr\mu}})_j]_\virt
\label{co10eq69}\\
&=\begin{cases}
[\acM_{(\ga_{i,j},\bs e_i)}^\ss(\ac\tau^{s_i\ti\la}_{\bs\mu_0})]_\virt\bt[\acM_{(\de_{i,j},\bs f_i)}^\ss(\ac\tau^{s_i\ti\la}_{\bs\mu_0})]_\virt, & o_{\ga_{i,j}}+o_{\de_{i,j}}=o_\be, \\
0, & o_{\ga_{i,j}}+o_{\de_{i,j}}>o_\be,
\end{cases}
\nonumber
\ea
for\/ $[\acM_{(\ga_{i,j},\bs e_i)}^\ss(\ac\tau^{s_i\ti\la}_{\bs\mu_0})]_\virt,[\acM_{(\de_{i,j},\bs f_i)}^\ss(\ac\tau^{s_i\ti\la}_{\bs\mu_0})]_\virt$ as in Definition\/~{\rm\ref{co10def1}}. 
\end{itemize}	
\end{prop}

\begin{proof} 
The decomposition of $\grM^\ss_{(\be,(1,\bs d))}(\gr\tau^{s_i\ti\la}_{\bs{\gr\mu}})_0^{\bG_m}$ into a disjoint union of pieces $\grM^\ss_{(\be,(1,\bs d))}(\gr\tau^{s_i\ti\la}_{\bs{\gr\mu}})_{\rho_{-1}=0},\ldots,\grM^\ss_{(\be,(1,\bs d))}(\gr\tau^{s_i\ti\la}_{\bs{\gr\mu}})_j$ in Proposition \ref{co10prop11} lifts to split $\bs\grM^\ss_{(\be,(1,\bs d))}(\gr\tau^{s_i\ti\la}_{\bs{\gr\mu}})_0^{\bG_m}$ into $\bs\grM^\ss_{(\be,(1,\bs d))}(\gr\tau^{s_i\ti\la}_{\bs{\gr\mu}})_{\rho_{-1}=0},\ab\ldots,\ab\bs\grM^\ss_{(\be,(1,\bs d))}(\gr\tau^{s_i\ti\la}_{\bs{\gr\mu}})_j$. 

For the first part of (a), we claim that the isomorphism $\Pi_{\rho_{-1}=0}$ in \eq{co10eq62} lifts to an isomorphism $\bs\Pi_{\rho_{-1}=0}:\bs\grM^\ss_{(\be,(1,\bs d))}(\gr\tau^{s_i\ti\la}_{\bs{\gr\mu}})_{\rho_{-1}=0}\ra \bs\acM_{(\be,\bs d)}^\ss(\ac\tau^{\smash{s_{i-1}'\ti\la}}_{\bs\mu_0})$ of proper quasi-smooth derived algebraic spaces. Then the associated obstruction theories are identified by~\eq{co10eq62}.

To see this, note that $\Pi_{\rho_{-1}=0}$ is the restriction to $\grM^\ss_{(\be,(1,\bs d))}(\gr\tau^{s_i\ti\la}_{\bs{\gr\mu}})_{\rho_{-1}=0}$ of a morphism $\Pi_{\rho_{-1}=0}':\grM^\ss_{(\be,(1,\bs d))}(\gr\tau^{s_i\ti\la}_{\bs{\gr\mu}})\ra \acM_{(\be,\bs d)}^\ss(\ac\tau^{\smash{s_{i-1}'\ti\la}}_{\bs\mu_0})$ defined by the same formula. This $\Pi_{\rho_{-1}=0}'$ is smooth, as it is a forgetful morphism forgetting $V_0,\rho_{-1},\rho_0$. This works on the derived level too, giving a smooth morphism $\bs\Pi_{\rho_{-1}=0}':\bs\grM^\ss_{(\be,(1,\bs d))}(\gr\tau^{s_i\ti\la}_{\bs{\gr\mu}})\ra \bs\acM_{(\be,\bs d)}^\ss(\ac\tau^{\smash{s_{i-1}'\ti\la}}_{\bs\mu_0})$. We define $\bs\Pi_{\rho_{-1}=0}$ to be the restriction of $\bs\Pi_{\rho_{-1}=0}'$ to $\bs\grM^\ss_{(\be,(1,\bs d))}(\gr\tau^{s_i\ti\la}_{\bs{\gr\mu}})_{\rho_{-1}=0}$. 

Now $\bs\grM^\ss_{(\be,(1,\bs d))}(\gr\tau^{s_i\ti\la}_{\bs{\gr\mu}})_{\rho_{-1}=0}$ is defined as a substack of $\bs\grM^\ss_{(\be,(1,\bs d))}(\gr\tau^{s_i\ti\la}_{\bs{\gr\mu}})$ by equations setting $\rho_0\ci\rho_{b_i-1}=0$ and $\rho_{-1}=0$. Noting that at each point $\rho_{b_i-1}$ is injective and $\ti\rho_0$ in \eq{co10eq60} is an isomorphism as $\rho_0\ne 0$, we can show that these equations are transverse on the smooth fibres of $\bs\Pi_{\rho_{-1}=0}'$ and turn these fibres into points, so $\bs\Pi_{\rho_{-1}=0}$ is an \'etale morphism of  derived stacks. Thus $\bs\Pi_{\rho_{-1}=0}$ is an isomorphism, as its classical truncation $\Pi_{\rho_{-1}=0}$ is by Proposition~\ref{co10prop11}(a).

Similarly, there is a derived enhancement of $\ac\Pi^0_{(\be,\bs d)}$ in \eq{co10eq59}:
\e
\bs{\ac\Pi}{}^0_{(\be,\bs d)}:\bs\grM^\ss_{(\be,(1,\bs d))}(\gr\tau^{s_i\ti\la}_{\bs{\gr\mu}})_0\longra\bs\acM_{(\be,\bs d)}^\ss(\ac\tau^{s_i\ti\la}_{\bs\mu_0}),
\label{co10eq70}
\e
which is the restriction of $\bs\Pi_{\rho_{-1}=0}'$ to $\bs\grM^\ss_{(\be,(1,\bs d))}(\gr\tau^{s_i\ti\la}_{\bs{\gr\mu}})_0$, and is smooth (of dimension 1) as the restrictions of the equation $\rho_0\ci\rho_{b_i-1}=0$ to the smooth fibres of $\bs\Pi_{\rho_{-1}=0}'$ are transverse. Taking the distinguished triangle \eq{co2eq11} of cotangent complexes of \eq{co10eq70}, pulling back to $\grM^\ss_{(\be,(1,\bs d))}(\gr\tau^{s_i\ti\la}_{\bs{\gr\mu}})_0$, and noting that $i^*(\bL_{\bs\grM/\bs\acM})=\bL_{\grM/\acM}$ as $\bs{\ac\Pi}{}^0_{(\be,\bs d)}$ is smooth, gives a distinguished triangle:
\e
\begin{gathered}
\xymatrix@R=10pt{ (\ac\Pi^0_{(\be,\bs d)})^*(i^*(\bL_{\bs\acM{}_{(\be,\bs d)}^\ss(\ac\tau^{s_i\ti\la}_{\bs\mu_0})})) \ar[d] \\ i^*(\bL_{\bs{\grM}{}^\ss_{(\be,(1,\bs d))}(\gr\tau^{s_i\ti\la}_{\bs{\gr\mu}})_0}) \ar[d] \\ \bL_{\grM/\acM}=\bL_{\grM^\ss_{(\be,(1,\bs d))}(\gr\tau^{s_i\ti\la}_{\bs{\gr\mu}})_0/\acM_{(\be,\bs d)}^\ss(\ac\tau^{s_i\ti\la}_{\bs\mu_0})} \ar[d]^(0.85){[+1]} \\ {} }
\end{gathered}
\label{co10eq71}
\e

This is $\bG_m$-equivariant, so restricting it to $\grM^\ss_{(\be,(1,\bs d))}(\gr\tau^{s_i\ti\la}_{\bs{\gr\mu}})_{\rho_{-1}=0}$, it is equivariant under a $\bG_m$-action on each term. As $\ac\Pi^0_{(\be,\bs d)}$ is $\bG_m$-invariant, the $\bG_m$-action on $(\ac\Pi^0_{(\be,\bs d)})^*(i^*(\bL_{\bs\acM{}_{(\be,\bs d)}^\ss(\ac\tau^{s_i\ti\la}_{\bs\mu_0})}))$ has $\bG_m$-weight 0. Also $\bL_{\grM/\acM}$ in \eq{co10eq71} is a line bundle, as $\ac\Pi^0_{(\be,\bs d)}$ is representable and smooth of dimension 1. It has $\bG_m$-weight $-1$, as the dual $\bL^*_{\grM/\acM}\vert_{[E,\bs V,\bs\rho]}$ has fibre $\Hom(V_a,V_0)$ which parametrizes $\rho_{-1}$, which has $\bG_m$-weight 1 by \eq{co10eq61}, and $\bG_m$ acts trivially on $V_0,V_a$. Thus the $\bG_m$-moving part of $i^*(\bL_{\bs{\grM}{}^\ss_{(\be,(1,\bs d))}(\gr\tau^{s_i\ti\la}_{\bs{\gr\mu}})_0})\vert_{\grM^\ss_{(\be,(1,\bs d))}(\gr\tau^{s_i\ti\la}_{\bs{\gr\mu}})_{\rho_{-1}=0}}$ is $\cN_{\rho_{-1}=0}^\bu=\bL_{\grM/\acM}$, a line bundle with $\bG_m$-weight $-1$. This proves~(a).

The proof of (b) is essentially the same, except that at $[E,\bs V,\bs\rho]$ the fibre of the line bundle $\bL_{\grM/\acM}$ is $\Hom(V_{b_i}/\rho_{b_i-1}(V_{b_i-1}),V_0)^*$, the dual of the space containing $\ti\rho_0$ in \eq{co10eq60}, and has weight $1$ as $\th\in\bG_m$ acts on $E,V_0,\ldots,\ab V_r$ by $\id_E,\th^{-1}\,\id_{V_0},\id_{V_1},\ldots,\id_{V_r}$ as in the proof of Proposition~\ref{co10prop11}.

The proofs of (c),(d) follow the (long) proofs of Proposition \ref{co9prop3}(c),(d) closely, so we just explain the differences. The analogue of \eq{co9eq19} is
\begin{equation*}
\xymatrix@C=120pt@R=15pt{
*+[r]{(\bs\dM_{(\be,(1,\bs d))}^\rpl)_0} \ar[r]_(0.55){\bs{\gr\Pi}{}^{v_0}_{\bs\dM^\red_\be}} \ar[d]^{\,\bs{\gr\jmath}{}_{(\be,(1,\bs d))}^{\,\pl}} & *+[l]{\bs\dM_\be^\red} \ar[d]_{\bs j_\be} \\
*+[r]{(\bs\dM_{(\be,(1,\bs d))}^\pl)_0} \ar[r]^(0.55){\bs{\gr\Pi}{}^{v_0}_{\bs\dM_\be}}  & *+[l]{\bs\dM_\be,\!} }
\end{equation*}
where $\bs j_\be,\bs{\gr\jmath}{}_{(\be,(1,\bs d))}^{\,\pl}$ are as in \eq{co5eq4} for $\B,\grB$, and $(\cdots)_0$ means we take the derived zero section of $\bs s:\cdots\ra \bs\cV_{b_i-1}^*\ot\bs\cV_0$ mapping $(E,\bs V,\bs\rho)\mapsto \rho_0\ci\rho_{b_i-1}$ as in Definition \ref{co10def4}, and $\bs{\gr\Pi}{}^{v_0}_{\bs\dM^\red_\be}$ is as in \eq{co5eq27} restricted to $(\cdots)_0$, and $\bs{\gr\Pi}{}^{v_0}_{\bs\dM_\be}$ is the analogue of $\bs{\gr\Pi}{}^{v_0}_{\bs\dM^\red_\be}$ for the `non-reduced' versions of the derived stacks. Then $\bs\grM^\ss_{(\be,(1,\bs d))}(\gr\tau^{s_i\ti\la}_{\bs{\gr\mu}})\subseteq\bs\dM_{(\be,(1,\bs d))}^\rpl$ is open, so $\bs\grM^\ss_{(\be,(1,\bs d))}(\gr\tau^{s_i\ti\la}_{\bs{\gr\mu}})_0\subseteq(\bs\dM_{(\be,(1,\bs d))}^\rpl)_0$ is open. The argument above for $\bs\Pi_{\rho_{-1}=0}$ and $\bs{\ac\Pi}{}^0_{(\be,\bs d)}$ shows that $\bs{\gr\Pi}{}^{v_0}_{\bs\dM^\red_\be},\bs{\gr\Pi}{}^{v_0}_{\bs\dM_\be}$ are smooth over $\bs\grM^\ss_{(\be,(1,\bs d))}(\gr\tau^{s_i\ti\la}_{\bs{\gr\mu}})_0$ and its image in $(\bs\dM_{(\be,(1,\bs d))}^\pl)_0$. They are also $\bG_m$-equivariant for the trivial $\bG_m$-actions on $\bs\dM_\be^\red,\bs\dM_\be$.

The main non-cosmetic changes in the proofs are that the analogue of the cotangent complex \eq{co9eq26} is more complicated, and this affects~\eq{co9eq27}--\eq{co9eq33}. The analogue of \eq{co9eq26} is that the relative cotangent complex of $\gr\Pi^{v_0}_{\dot{\mathcal M}_\be}:\grM^\ss_{(\be,(1,\bs d))}(\gr\tau^{s_i\ti\la}_{\bs{\gr\mu}})_0\ra\dM_\be$ fits into an exact sequence
\e
\begin{gathered}
\xymatrix@R=.5pt{
0 \ar[dd] \\ \\
\bL_{\grM^\ss_{(\be,(1,\bs d))}(\gr\tau^{s_i\ti\la}_{\bs{\gr\mu}})_0/\dM_\be}
\ar[ddd] \\ \\ \\
{\begin{subarray}{l} \ts
\O\op\bigl[(\cV_a\ot\cV_0^*)\op((\cV_{d_i}/\rho_{d_i-1}(\cV_{d_i-1}))\ot\cV_0^*) \\
\ts \qquad \op 
\bigop_{l=1}^{r-1}(\cV_l\ot\cV_{l+1}^*)\op (\cV_r\ot\cV_{k,\be}^*)\bigr]^\pl
\end{subarray}} \ar[ddd] \\ \\ \\
\bigl[\bigop_{l=0}^r(\cV_l\ot\cV_l^*)\bigr]^\pl \ar[dd] \\ \\ 0.\!\!
}
\end{gathered}
\label{co10eq72}
\e
Here the third term parametrizes the dual spaces to the data of the edges in \eq{co10eq57}, and the fourth term the dual spaces to the automorphism Lie algebras of the vertices, and we add $\O$ in the third term because of the quotient $\fatslash\,\,\bG_m$ used to form~$\grM^\pl_{(\be,(1,\bs d))}$.

For the analogue of \eq{co9eq27}, we pull equation \eq{co10eq72} back by the inclusion $i_j:\grM^\ss_{(\be,(1,\bs d))}(\gr\tau^{s_i\ti\la}_{\bs{\gr\mu}})_j\hookra\grM^\ss_{(\be,(1,\bs d))}(\gr\tau^{s_i\ti\la}_{\bs{\gr\mu}})_0$ and rewrite to give
\e
\begin{gathered}
\xymatrix@R=.5pt{
0 \ar[dd] \\ \\
i_j^*(\bL_{\grM^\ss_{(\be,(1,\bs d))}(\gr\tau^{s_i\ti\la}_{\bs{\gr\mu}})_0/\dM_\be})
\ar[ddd] \\ \\ \\
{\begin{subarray}{l} \ts
\O^3\op\bigl[\bigop_{l=1}^{r-1}((\cV'_l\bp\cV_l'')\ot(\cV'_{l+1}\bp\cV_{l+1}'')^*)\\
\ts \quad \op ((\cV'_r\bp\cV_r'')\ot(\cV_{k,\ga_{i,j}}'\bp\cV_{k,\de_{i,j}}'')^*)\bigr]^\pl
\end{subarray}} \ar[ddd] \\ \\ \\
\O\op\bigl[\bigop_{l=1}^r((\cV'_l\bp\cV_l'')\ot(\cV'_l\bp\cV_l'')^*)\bigr]^\pl \ar[dd] \\ \\ 0.\!\!
}
\end{gathered}
\label{co10eq73}
\e
Here we identify the pullbacks of $\cV_a\ot\cV_0^*,(\cV_{d_i}/\rho_{d_i-1}(\cV_{d_i-1}))\ot\cV_0^*$ and $\cV_0\ot\cV_0^*$ with $\O$, and use the isomorphism \eq{co10eq64}, and the notation of \eq{co10eq67}--\eq{co10eq68}.

Now $\bG_m$ acts on each term in \eq{co10eq73}. By Proposition \ref{co10prop11}(c), $\th\in\bG_m$ acts by $\th\,\id_{E'}+\id_{E''}$ on $E=E'\op E'',$ and by\/ $\id_{V_0}$ on\/ $V_0,$ and by $\th\,\id_{V_l'}+\id_{V_l''}$ on $V_l=V_l'\op V_l''$. Thus the vector bundles $\cV'$ in \eq{co10eq73} have $\bG_m$-weight 1, so that $(\cV')^*$ has $\bG_m$-weight $-1$, and $\cV'',(\cV'')^*,\O$ have $\bG_m$-weight 0. Hence we can split \eq{co10eq73} into the direct sum of three exact sequences, with $\bG_m$-weights $0,1,-1$. After cancelling $\O$ from the third and fourth lines, the weight 0 sequence is the direct sum of the analogues of \eq{co10eq72} for $\ac\Pi^{v_a}_{\dot{\mathcal M}_{\ga_{i,j}}}:\acM_{(\ga_{i,j},\bs e_i)}^\ss(\ac\tau^{s_i\ti\la}_{\bs\mu_0})\ra\dM_{\ga_{i,j}}$ and $\ac\Pi^{v_{b_i}}_{\dot{\mathcal M}_{\de_{i,j}}}:\acM_{(\de_{i,j},\bs f_i)}^\ss(\ac\tau^{s_i\ti\la}_{\bs\mu_0})\ra\dM_{\de_{i,j}}$. Thus we have
\e
\begin{split}
&i_j^*(\bL_{\grM^\ss_{(\be,(1,\bs d))}(\gr\tau^{s_i\ti\la}_{\bs{\gr\mu}})_0/\dM_\be})_{\wt=0}\\
&\quad \cong \bL_{\acM_{(\ga_{i,j},\bs e_i)}^\ss(\ac\tau^{s_i\ti\la}_{\bs\mu_0})/\dM_{\ga_{i,j}}}\bp
\bL_{\acM_{(\de_{i,j},\bs f_i)}^\ss(\ac\tau^{s_i\ti\la}_{\bs\mu_0})/\dM_{\de_{i,j}}}.
\end{split}
\label{co10eq74}
\e
The $\bG_m$-weight $1,-1$ components of \eq{co10eq73} give exact sequences
\e
\begin{gathered}
\xymatrix@R=.5pt{
0 \ar[dd] \\ \\
i_j^*(\bL_{\grM^\ss_{(\be,(1,\bs d))}(\gr\tau^{s_i\ti\la}_{\bs{\gr\mu}})_0/\dM_\be})_{\wt=1}
\ar[ddd] \\ \\ \\
{\begin{subarray}{l} \ts
\bigl[\bigop_{l=1}^{r-1}\cV'_l\bt(\cV_{l+1}'')^*\\
\ts \quad \op \cV'_r\bt(\cV_{k,\de_{i,j}}'')^*\bigr]^\pl
\end{subarray}} \ar[ddd] \\ \\ \\
\bigl[\bigop_{l=1}^r\cV'_l\bt(\cV_l'')^*\bigr]^\pl \ar[dd] \\ \\ 0,\!\!
}
\quad
\xymatrix@R=.5pt{
0 \ar[dd] \\ \\
i_j^*(\bL_{\grM^\ss_{(\be,(1,\bs d))}(\gr\tau^{s_i\ti\la}_{\bs{\gr\mu}})_0/\dM_\be})_{\wt=-1}
\ar[ddd] \\ \\ \\
{\begin{subarray}{l} \ts
\bigl[\bigop_{l=1}^{r-1}(\cV_{l+1}')^*\bt\cV''_l\\
\ts \quad \op (\cV_{k,\ga_{i,j}}')^*\bt\cV''_r\bigr]^\pl
\end{subarray}} \ar[ddd] \\ \\ \\
\bigl[\bigop_{l=1}^r(\cV_l')^*\bt\cV''_l\bigr]^\pl \ar[dd] \\ \\ 0.\!\!
}
\end{gathered}
\label{co10eq75}
\e
Equations \eq{co10eq74}--\eq{co10eq75} are the analogue of \eq{co9eq27}.

As in the proof of Proposition \ref{co9prop3}(c),(d), the analogue of \eq{co9eq22} holds in our situation, with the three occurrences of $i_{\al_1,\al_2}^*(\bL_{\acM^\pl_{(\al,\bs 1)}/\dM_\al})$ replaced by $i_j^*(\bL_{\grM^\ss_{(\be,(1,\bs d))}(\gr\tau^{s_i\ti\la}_{\bs{\gr\mu}})_0/\dM_\be})$ in \eq{co10eq73}. The $\bG_m$-weight $0,1,-1$ parts of the analogue of \eq{co9eq22} give the analogues of \eq{co9eq28}, which we do not write out, and of \eq{co9eq29}--\eq{co9eq30}, which are
\ea
&\xymatrix@C=3pt{
{\begin{subarray}{l}\ts (\bar\Pi^{v_a}_{\dot{\mathcal M}_{\ga_{i,j}}}\!\!\t\!\bar\Pi^{v_{b_i}}_{\dot{\mathcal M}_{\de_{i,j}}})^* \\ \ts (\cE_{\ga_{i,j},\de_{i,j}}^\bu[-1]) \end{subarray}} \ar[rr] && \cN_j^{+,\bu}  \ar[rrr] &&& i_j^*(\bL_{\grM^\ss_{(\be,(1,\bs d))}(\gr\tau^{s_i\ti\la}_{\bs{\gr\mu}})_0/\dM_\be})_{\wt=1}  \ar[rr]^(0.9){[+1]} &&  , }\!\!
\label{co10eq76}\\
&\xymatrix@C=2.5pt{
{\begin{subarray}{l}\ts (\bar\Pi^{v_a}_{\dot{\mathcal M}_{\ga_{i,j}}}\!\!\t\!\bar\Pi^{v_{b_i}}_{\dot{\mathcal M}_{\de_{i,j}}})^* \\ \ts (\si^*(\cE_{\de_{i,j},\ga_{i,j}}^\bu)[-1])\!\!\! \end{subarray}}
\ar[rr] && \cN_j^{-,\bu}  \ar[rrr] &&& i_j^*(\bL_{\grM^\ss_{(\be,(1,\bs d))}(\gr\tau^{s_i\ti\la}_{\bs{\gr\mu}})_0/\dM_\be})_{\wt=-1}  \ar[rr]^(0.9){[+1]} && .  }\!\!
\label{co10eq77}
\ea
Equations \eq{co10eq67}--\eq{co10eq68} now follow from \eq{co10eq75}--\eq{co10eq77}, proving~(c).

The proof of (d) follows that of Proposition \ref{co9prop3}(d) over \eq{co9eq31}--\eq{co9eq36}. An important point in the proof is that in \eq{co9eq31}--\eq{co9eq32}, the right hand terms $(\bar\Pi^v_{\dot{\mathcal M}_{\al_1}})^*(\cV^*_{k_1,\al_1})\op(\bar\Pi^v_{\dot{\mathcal M}_{\al_2}})^*(\cV^*_{k_2,\al_2})$ are the same. The analogue of this is \eq{co10eq74}, which equates the corresponding terms in the analogues of \eq{co9eq31}--\eq{co9eq32} for the obstruction theories on $\grM^\ss_{(\be,(1,\bs d))}(\gr\tau^{s_i\ti\la}_{\bs{\gr\mu}})_j$ and $\acM_{(\ga_{i,j},\bs e_i)}^\ss(\ac\tau^{s_i\ti\la}_{\bs\mu_0})\t\acM_{(\de_{i,j},\bs f_i)}^\ss(\ac\tau^{s_i\ti\la}_{\bs\mu_0})$. The rest of the proof needs only obvious changes.
\end{proof}

Here is the analogue of Definition \ref{co9def3}.

\begin{dfn}
\label{co10def5}
Work in the situation above. First observe that if $(\be,\bs d)$ satisfies the conditions of Proposition \ref{co10prop10}, one can show that
\e
\rank\bT_{\acM^\pl_{(\be,\bs d)}/\M^\pl_\be}\!=\!\sum_{i=1}^{r-1}d_id_{i+1}\!+\!d_r\la_k(\be)\!-\!\sum_{i=1}^rd_i^2\!=\!\ha\la_k(\be)(\la_k(\be)\!-\!1).
\label{co10eq78}
\e
This allows us to simplify the sign appearing in \eq{co10eq6}.

Define an equivariant cohomology class
\e
\begin{split}
\eta&\in H^{\la_k(\be)(\la_k(\be)-1)}_{\bG_m}(\grM^\ss_{(\be,(1,\bs d))}(\gr\tau^{s_i\ti\la}_{\bs{\gr\mu}})_0)\quad\text{by}
\\
\eta&=(\ac\Pi^0_{(\be,\bs d)})^*\bigl(c_{\la_k(\be)(\la_k(\be)-1)/2}(\bT_{\acM^\pl_{(\be,\bs d)}/\M^\pl_\be})\bigr)\\
&=(-1)^{\la_k(\be)(\la_k(\be)-1)/2}(\ac\Pi^0_{(\be,\bs d)})^*\bigl(c_\top(\bL_{\acM^\pl_{(\be,\bs d)}/\M^\pl_\be})\bigr).
\end{split}
\label{co10eq79}
\e
Here one can show using the conditions on $(\be,\bs d)$ in Proposition \ref{co10prop10} that $\Pi_{\M^\ss_\be(\tau)}:\acM_{(\be,\bs d)}^\ss(\ac\tau^{s_i\ti\la}_{\bs\mu_0})\ra\M^\ss_\be(\tau)$ is smooth of  dimension \eq{co10eq78}, so its relative cotangent complex $\bL_{\acM_{(\be,\bs d)}^\ss(\ac\tau^{s_i\ti\la}_{\bs\mu_0})/\M^\ss_\be(\tau)}=\bL_{\acM^\pl_{(\be,\bs d)}/\M^\pl_\be}$ is a vector bundle of this rank on $\grM^\ss_{(\be,(1,\bs d))}(\gr\tau^{s_i\ti\la}_{\bs{\gr\mu}})$. As $\ac\Pi^0_{(\be,\bs d)}$ is $\bG_m$-equivariant for the trivial $\bG_m$-action on $\acM_{(\be,\bs d)}^\ss(\ac\tau^{s_i\ti\la}_{\bs\mu_0})$, the pullback in \eq{co10eq79} makes sense in equivariant cohomology $H^*_{\bG_m}(\grM^\ss_{(\be,(1,\bs d))}(\gr\tau^{s_i\ti\la}_{\bs{\gr\mu}})_0)$. Hence \eq{co10eq79} is well defined.

Apply Corollary \ref{co2cor2} to the proper algebraic space $\grM^\ss_{(\be,(1,\bs d))}(\gr\tau^{s_i\ti\la}_{\bs{\gr\mu}})_0$ with its $\bG_m$-action, $\bG_m$-equivariant obstruction theory $\bL_i:i^*(\bL_{\bs\grM^\ss_{(\be,(1,\bs d))}(\gr\tau^{s_i\ti\la}_{\bs{\gr\mu}})_0})\ab\ra \bL_{\grM^\ss_{(\be,(1,\bs d))}(\gr\tau^{s_i\ti\la}_{\bs{\gr\mu}})_0}$, and cohomology class $\eta$, and with $\Pi_{\M_\be^{\smash{\pl}}}:\grM^\ss_{(\be,(1,\bs d))}(\gr\tau^{s_i\ti\la}_{\bs{\gr\mu}})_0\ab\ra\M^\pl_\be$ in place of $f:X\ra Y$. Using the notation of Propositions \ref{co10prop11}--\ref{co10prop12}, and writing $i_{\rho_{-1}=0}:\grM^\ss_{(\be,(1,\bs d))}(\gr\tau^{s_i\ti\la}_{\bs{\gr\mu}})_{\rho_{-1}=0}\hookra\grM^\ss_{(\be,(1,\bs d))}(\gr\tau^{s_i\ti\la}_{\bs{\gr\mu}})_0$ for the inclusion, and similarly for $i_{\rho_0=0},i_j$, equation \eq{co2eq26} for $\grM^\ss_{(\be,(1,\bs d))}(\gr\tau^{s_i\ti\la}_{\bs{\gr\mu}})_0$ becomes 
\ea
&0=(-1)^{\rank\cN_{\rho_{-1}=0}^\bu}\Res_z\bigl\{
(\Pi_{\M_\be^{\smash{\pl}}}\ci i_{\rho_{-1}=0})_*\bigl([\grM^\ss_{(\be,(1,\bs d))}(\gr\tau^{s_i\ti\la}_{\bs{\gr\mu}})_{\rho_{-1}=0}]_\virt
\nonumber\\
&\qquad\qquad\qquad \cap (e(\cN_{\rho_{-1}=0}^\bu)^{-1}\!\cup\! i_{\rho_{-1}=0}^*(\eta))\bigr)\bigr\}
\nonumber\\
\begin{split}
&+(-1)^{\rank\cN_{\rho_0=0}^\bu}\Res_z\bigl\{
(\Pi_{\M_\be^{\smash{\pl}}}\ci i_{\rho_0=0})_*\bigl([\grM^\ss_{(\be,(1,\bs d))}(\gr\tau^{s_i\ti\la}_{\bs{\gr\mu}_-})_{\rho_0=0}]_\virt\\
&\qquad\qquad\qquad
\cap (e(\cN_{\rho_0=0}^\bu)^{-1}\cup i_{\rho_0=0}^*(\eta))\bigr)\bigr\}
\end{split}
\label{co10eq80}\\
&+\sum_{j=1,\ldots,q_i\!\!\!\!\!\!\!\!\!}\begin{aligned}[t] &(-1)^{\rank\cN_j^\bu}\Res_z\bigl\{
(\Pi_{\M_\be^{\smash{\pl}}}\ci i_j)_*\bigl([\grM^\ss_{(\be,(1,\bs d))}(\gr\tau^{s_i\ti\la}_{\bs{\gr\mu}})_j]_\virt\\
&\qquad\qquad\qquad\qquad
\cap(e(\cN_j^\bu)^{-1}\cup i_j^*(\eta))\bigr)\bigr\}.
\end{aligned}
\nonumber
\ea
\end{dfn}

The next two results are the analogues of Propositions \ref{co9prop4} and~\ref{co9prop5}.

\begin{prop}
\label{co10prop13}
In\/ {\rm\eq{co10eq80},} using the notation of\/ \eq{co10eq6} we have
\ea
&(-1)^{\rank\cN_{\rho_{-1}=0}^\bu}\Res_z\bigl\{
(\Pi_{\M_\be^{\smash{\pl}}}\ci i_{\rho_{-1}=0})_*\bigl([\grM^\ss_{(\be,(1,\bs d))}(\gr\tau^{s_i\ti\la}_{\bs{\gr\mu}})_{\rho_{-1}=0}]_\virt
\nonumber\\
&\qquad\qquad\qquad \cap (e(\cN_{\rho_{-1}=0}^\bu)^{-1}\!\cup\! i_{\rho_{-1}=0}^*(\eta))\bigr)\bigr\}=\ac\Up_{(\be,\bs d)}(\ac\tau^{s_{i-1}'\ti\la}_{\bs\mu_0}),
\label{co10eq81}\\
&(-1)^{\rank\cN_{\rho_0=0}^\bu}\Res_z\bigl\{
(\Pi_{\M_\be^{\smash{\pl}}}\ci i_{\rho_0=0})_*\bigl([\grM^\ss_{(\be,(1,\bs d))}(\gr\tau^{s_i\ti\la}_{\bs{\gr\mu}_-})_{\rho_0=0}]_\virt
\nonumber\\
&\qquad\qquad\qquad
\cap (e(\cN_{\rho_0=0}^\bu)^{-1}\cup i_{\rho_0=0}^*(\eta))\bigr)\bigr\}=-\ac\Up_{(\be,\bs d)}(\ac\tau^{s_i'\ti\la}_{\bs\mu_0}). 
\label{co10eq82}
\ea
\end{prop}

\begin{proof} 
To prove \eq{co10eq81}, we will show that
\ea
&(-1)^{\rank\cN_{\rho_{-1}=0}^\bu}\Res_z\bigl\{
(\Pi_{\M_\be^{\smash{\pl}}}\ci i_{\rho_{-1}=0})_*\bigl([\grM^\ss_{(\be,(1,\bs d))}(\gr\tau^{s_i\ti\la}_{\bs{\gr\mu}})_{\rho_{-1}=0}]_\virt
\nonumber\\
&\qquad\qquad\qquad \cap (e(\cN_{\rho_{-1}=0}^\bu)^{-1}\!\cup\! i_{\rho_{-1}=0}^*(\eta))\bigr)\bigr\}
\nonumber\\
\begin{split}
&=-\Res_z\bigl\{
(\Pi_{\M_\be^{\smash{\pl}}}\ci i_{\rho_{-1}=0})_*\bigl([\grM^\ss_{(\be,(1,\bs d))}(\gr\tau^{s_i\ti\la}_{\bs{\gr\mu}})_{\rho_{-1}=0}]_\virt\\
&\qquad\qquad\qquad \cap (e(\cN_{\rho_{-1}=0}^\bu)^{-1}\!\cup\!  \Pi_{\rho_{-1}=0}^*(c_\top(\bT_{\acM^\pl_{(\be,\bs d)}/\M^\pl_\be}))\bigr)\bigr\}
\end{split}
\label{co10eq83}\\
&=(\Pi_{\M_\be^{\smash{\pl}}}\ci i_{\rho_{-1}=0})_*\bigl([\grM^\ss_{(\be,(1,\bs d))}(\gr\tau^{s_i\ti\la}_{\bs{\gr\mu}})_{\rho_{-1}=0}]_\virt\cap \Pi_{\rho_{-1}=0}^*(c_\top(\bT_{\acM^\pl_{(\be,\bs d)}/\M^\pl_\be}))\bigr)
\nonumber\\
&=(\Pi_{\M_\be^{\smash{\pl}}})_*\bigl([\acM_{(\be,\bs d)}^\ss(\ac\tau^{\smash{s_{i-1}'\ti\la}}_{\bs\mu_0})]_\virt\cap c_\top(\bT_{\acM^\pl_{(\be,\bs d)}/\M^\pl_\be})\bigr)=\ac\Up_{(\be,\bs d)}(\ac\tau^{s_{i-1}'\ti\la}_{\bs\mu_0}).
\nonumber
\ea
Here in the first step we use $\rank\cN_{\rho_{-1}=0}^\bu=1$, equation \eq{co10eq79}, and $\Pi_{\rho_{-1}=0}=\ac\Pi^0_{(\be,\bs d)}\ci i_{\rho_{-1}=0}$. As $\cN_{\rho_{-1}=0}^\bu$ is a line bundle $\bL_{\grM/\acM}$ with $\bG_m$-weight $-1$ by Proposition \ref{co10prop12}(a) we have $e(\cN_{\rho_{-1}=0}^\bu)=-1\cdot z+c_1(\bL_{\grM/\acM})$, where the coefficient $-1$ of $z$ is the weight $-1$ of the $\bG_m$-action, so $e(\cN_{\rho_{-1}=0}^\bu)^{-1}\ab =-z^{-1}+O(z^{<-1})$. But $[\grM^\ss_{(\be,(1,\bs d))}(\gr\tau^{s_i\ti\la}_{\bs{\gr\mu}})_{\rho_{-1}=0}]_\virt$ and $\Pi_{\rho_{-1}=0}^*(c_\top(\bT_{\acM^\pl_{(\be,\bs d)}/\M^\pl_\be}))$ are both $\bG_m$-invariant, and so are independent of $z$. Thus taking $\Res_z\{\cdots\}$ gives a factor $-1$ from $e(\cN_{\rho_{-1}=0}^\bu)^{-1}$, and the second step of \eq{co10eq83} follows. In the third step we use Proposition \ref{co10prop12}(a), and in the fourth \eq{co10eq6}. The proof of \eq{co10eq82} is similar, using Proposition \ref{co10prop12}(b). The sign difference in \eq{co10eq82} comes from $e(\cN_{\rho_0=0}^\bu)=1\cdot z+c_1(\bL_{\grM/\acM})$ since $\cN_{\rho_0=0}^\bu$ has $\bG_m$-weight~$1$.
\end{proof}

\begin{prop}
\label{co10prop14}
In equation {\rm\eq{co10eq80},} for each\/ $j=1,\ldots,q_i$ we have
\ea
&(-1)^{\rank\cN_j^\bu}\Res_z\bigl\{
(\Pi_{\M_\be^{\smash{\pl}}}\ci i_j)_*\bigl([\grM^\ss_{(\be,(1,\bs d))}(\gr\tau^{s_i\ti\la}_{\bs{\gr\mu}})_j]_\virt\cap(e(\cN_j^\bu)^{-1}\cup i_j^*(\eta))\bigr)\bigr\}
\nonumber\\
&=\begin{cases} \bigl[\ac\Up_{(\ga_{i,j},\bs e_i)}(\ac\tau^{s_i\ti\la}_{\bs\mu_0}),\ac\Up_{(\de_{i,j},\bs f_i)}(\ac\tau^{s_i\ti\la}_{\bs\mu_0})\bigr], & o_{\ga_{i,j}}+o_{\de_{i,j}}=o_\be, \\
0, & o_{\ga_{i,j}}+o_{\de_{i,j}}>o_\be, \end{cases}
\label{co10eq84}
\ea
using the Lie bracket on $\check H_{\rm even}(\M^\pl)$ from\/ {\rm\S\ref{co43},} where $o_{\ga_{i,j}}+o_{\de_{i,j}}\ge o_\be$ by Assumption\/~{\rm\ref{co5ass1}(f)(iv)}.
\end{prop}

\begin{proof} The second case of \eq{co10eq84}, when $o_{\ga_{i,j}}+o_{\de_{i,j}}>o_\be$, is immediate from \eq{co10eq69}. So suppose that $o_{\ga_{i,j}}+o_{\de_{i,j}}=o_\be$. By a very similar computation to \eq{co10eq72}--\eq{co10eq75}, we deduce that $(\ac\Pi^0_{(\be,\bs d)}\ci i_j)^*(\bL_{\acM^\pl_{(\be,\bs d)}/\M^\pl_\be})$ has components of $\bG_m$-weights $0,1,-1$, with an isomorphism and exact sequences
\begin{gather}
(\ac\Pi^0_{(\be,\bs d)}\!\ci\! i_j)^*(\bL_{\acM^\pl_{(\be,\bs d)}/\M^\pl_\be})_{\wt=0}\!\cong\! \bL_{\acM_{(\ga_{i,j},\bs e_i)}^\pl/\M_{\ga_{i,j}}^\pl}\bp
\bL_{\acM_{(\de_{i,j},\bs f_i)}^\pl/\M_{\de_{i,j}}^\pl},
\label{co10eq85}\\
\begin{gathered}
\xymatrix@R=.5pt{
0 \ar[dd] \\ \\
(\ac\Pi^0_{(\be,\bs d)}\ci i_j)^*(\bL_{\acM^\pl_{(\be,\bs d)}/\M^\pl_\be})_{\wt=1}
\ar[ddd] \\ \\ \\
{\begin{subarray}{l} \ts
\bigl[\bigop_{l=1}^{r-1}\cV'_l\bt(\cV_{l+1}'')^*\\
\ts \quad \op \cV'_r\bt(\cV_{k,\de_{i,j}}'')^*\bigr]^\pl
\end{subarray}} \ar[ddd] \\ \\ \\
\bigl[\bigop_{l=1}^r\cV'_l\bt(\cV_l'')^*\bigr]^\pl \ar[dd] \\ \\ 0,\!\!
}
\quad
\xymatrix@R=.5pt{
0 \ar[dd] \\ \\
(\ac\Pi^0_{(\be,\bs d)}\ci i_j)^*(\bL_{\acM^\pl_{(\be,\bs d)}/\M^\pl_\be})_{\wt=-1}
\ar[ddd] \\ \\ \\
{\begin{subarray}{l} \ts
\bigl[\bigop_{l=1}^{r-1}(\cV_{l+1}')^*\bt\cV''_l\\
\ts \quad \op (\cV_{k,\ga_{i,j}}')^*\bt\cV''_r\bigr]^\pl
\end{subarray}} \ar[ddd] \\ \\ \\
\bigl[\bigop_{l=1}^r(\cV_l')^*\bt\cV''_l\bigr]^\pl \ar[dd] \\ \\ 0.\!\!
}
\end{gathered}
\label{co10eq86}
\end{gather}
Note that the third and fourth rows of \eq{co10eq75} and \eq{co10eq86} agree.

Writing $H^*_{\bG_m}(\grM^\ss_{(\be,(1,\bs d))}(\gr\tau^{s_i\ti\la}_{\bs{\gr\mu}})_j)=H^*_{\bG_m}(\acM_{(\ga_{i,j},\bs e_i)}^\ss(\ac\tau^{s_i\ti\la}_{\bs\mu_0})\t\acM_{(\de_{i,j},\bs f_i)}^\ss(\ac\tau^{s_i\ti\la}_{\bs\mu_0}))$ by \eq{co10eq64}, and including $\bG_m$-weights in the notation, we compute
\ea
&i_j^*(\eta)=(-1)^{\la_k(\be)(\la_k(\be)-1)/2}(\ac\Pi^0_{(\be,\bs d)}\ci i_j)^*(c_\top(\bL_{\acM^\pl_{(\be,\bs d)}/\M^\pl_\be}))
\label{co10eq87}\\
&=(-1)^{\la_k(\be)(\la_k(\be)-1)/2}c_\top\bigl((\ac\Pi^0_{(\be,\bs d)}\ci i_j)^*(\bL_{\acM^\pl_{(\be,\bs d)}/\M^\pl_\be})_{\wt=0}\bigr)
\nonumber\\
&\cup c_\top\bigl((\ac\Pi^0_{(\be,\bs d)}\!\ci\! i_j)^*(\bL_{\acM^\pl_{(\be,\bs d)}/\M^\pl_\be})_{\wt=1}\bigr)
\!\cup\! c_\top\bigl((\ac\Pi^0_{(\be,\bs d)}\!\ci\! i_j)^*(\bL_{\acM^\pl_{(\be,\bs d)}/\M^\pl_\be})_{\wt=-1}\bigr)
\nonumber\\
&=(-1)^{\la_k(\be)(\la_k(\be)-1)/2}c_\top\bigl(
\ts\sum\limits_{1\le l<r\!\!\!\!}[[\cV'_l\bt(\cV_{l+1}'')^*]^\pl]+[[\cV'_r\bt(\cV_{k,\de_{i,j}}'')^*]^\pl]
\nonumber\\
&-\ts\sum\limits_{1\le l\le r\!\!\!\!}[[\cV'_l\bt(\cV_l'')^*]^\pl]\bigr)_{\wt=1}\cup c_\top\bigl(
\ts\sum\limits_{1\le l<r\!\!\!\!}[[(\cV_{l+1}')^*\bt\cV''_l]^\pl]+[[(\cV_{k,\ga_{i,j}}')^*\bt\cV''_r]^\pl]
\nonumber\\
&-\ts\sum\limits_{1\le l\le r\!\!\!\!}[[(\cV_l')^*\bt\cV''_l]^\pl]\bigr)_{\wt=-1}\cup c_\top\bigl(\bL_{\acM_{(\ga_{i,j},\bs e_i)}^\pl/\M_{\ga_{i,j}}^\pl}\bp\bL_{\acM_{(\de_{i,j},\bs f_i)}^\pl/\M_{\de_{i,j}}^\pl}\bigr){}_{\wt=0}.
\nonumber
\ea
Here the first step uses \eq{co10eq79}, the second that $(\ac\Pi^0_{(\be,\bs d)}\ci i_j)^*(\bL_{\acM^\pl_{(\be,\bs d)}/\M^\pl_\be})$ has $\bG_m$-weights $0,1,-1$, and the third \eq{co10eq85}--\eq{co10eq86}.

We can now prove the first case of \eq{co10eq84}, in the long equation
\ea
&(-1)^{\rank\cN_j^\bu}\Res_z\bigl\{
(\Pi_{\M_\be^{\smash{\pl}}}\ci i_j)_*\bigl([\grM^\ss_{(\be,(1,\bs d))}(\gr\tau^{s_i\ti\la}_{\bs{\gr\mu}})_j]_\virt\cap(e(\cN_j^\bu)^{-1}\cup i_j^*(\eta))\bigr)\bigr\}
\nonumber\\
&=(-1)^{\la_k(\ga_{i,j})\la_k(\de_{i,j})+\chi(\ga_{i,j},\de_{i,j})+\chi(\de_{i,j},\ga_{i,j})}
\nonumber\\
&\;\>\Res_z\bigl\{
(\Pi^\pl_\be\ci\Phi_{\ga_{i,j},\de_{i,j}}\ci(\ac\Pi^{v_a}_{\dot{\mathcal M}_{\ga_{i,j}}}\t\ac\Pi^{v_{b_i}}_{\dot{\mathcal M}_{\de_{i,j}}}))^{\bG_m}_*
\nonumber\\
&\;\>\bigl(([\acM_{(\ga_{i,j},\bs e_i)}^\ss(\ac\tau^{s_i\ti\la}_{\bs\mu_0})]_\virt\!\bt\![\acM_{(\de_{i,j},\bs f_i)}^\ss(\ac\tau^{s_i\ti\la}_{\bs\mu_0})]_\virt)\cap
(e(\cN_j^\bu)^{-1}\!\cup\! i_j^*(\eta))\bigr)\bigr\}
\allowdisplaybreaks
\nonumber\\
&=(-1)^{\la_k(\ga_{i,j})\la_k(\de_{i,j})+\chi(\ga_{i,j},\de_{i,j})+\chi(\de_{i,j},\ga_{i,j})}(-1)^{\la_k(\be)(\la_k(\be)-1)/2}\cdot{}
\nonumber\\
&\;\>\Res_z\bigl\{(\Pi^\pl_\be\ci\Phi_{\ga_{i,j},\de_{i,j}}\ci(\ac\Pi^{v_a}_{\dot{\mathcal M}_{\ga_{i,j}}}\t\ac\Pi^{v_{b_i}}_{\dot{\mathcal M}_{\de_{i,j}}}))^{\bG_m}_*
\nonumber\\
&\;\>\bigl(([\acM_{(\ga_{i,j},\bs e_i)}^\ss(\ac\tau^{s_i\ti\la}_{\bs\mu_0})]_\virt\!\bt\![\acM_{(\de_{i,j},\bs f_i)}^\ss(\ac\tau^{s_i\ti\la}_{\bs\mu_0})]_\virt)\cap
\nonumber\\
&\;\>\bigl[e\bigl(\ts\sum_{1\le l<r}[[\cV'_l\bt(\cV_{l+1}'')^*]^\pl]+[[\cV'_r\bt(\cV_{k,\de_{i,j}}'')^*]^\pl]
\nonumber\\
&\;\>-\ts\sum_{1\le l\le r}[[\cV'_l\bt(\cV_l'')^*]^\pl]\bigr)_{\wt=1}^{-1}
\cup e\bigl((\ac\Pi^{v_a}_{\dot{\mathcal M}_{\ga_{i,j}}}\t\ac\Pi^{v_{b_i}}_{\dot{\mathcal M}_{\de_{i,j}}})^*(\cE^\bu_{\ga_{i,j},\de_{i,j}})\bigr)_{\wt=1}
\nonumber\\
&\;\>\cup e\bigl(\ts\sum_{1\le l<r}[[(\cV_{l+1}')^*\bt\cV''_l]^\pl]+[(\cV_{k,\ga_{i,j}}')^*\bt\cV''_r]^\pl]
\nonumber\\
&\;\>-\ts \sum_{1\le l\le r}[[(\cV_l')^*\bt\cV''_l]^\pl]\bigr)_{\wt=-1}^{-1}\cup e\bigl((\ac\Pi^{v_a}_{\dot{\mathcal M}_{\ga_{i,j}}}\!\!\t\!\ac\Pi^{v_{b_i}}_{\dot{\mathcal M}_{\de_{i,j}}})^*(\si^*(\cE_{\de_{i,j},\ga_{i,j}}^\bu))\bigr)_{\wt=-1}
\nonumber\\
&\;\>\cup c_\top\bigl(
\ts\sum\limits_{1\le l<r\!\!\!\!}[[\cV'_l\bt(\cV_{l+1}'')^*]^\pl]+[[\cV'_r\bt(\cV_{k,\de_{i,j}}'')^*]^\pl]
-\sum\limits_{1\le l\le r\!\!\!\!}[[\cV'_l\bt(\cV_l'')^*]^\pl]\bigr)_{\wt=1}
\nonumber\\
&\;\>\cup c_\top\bigl(
\ts\sum\limits_{1\le l<r\!\!\!\!}[[(\cV_{l+1}')^*\bt\cV''_l]^\pl]\!+\![[(\cV_{k,\ga_{i,j}}')^*\bt\cV''_r]^\pl]
\!-\!\ts\sum\limits_{1\le l\le r\!\!\!\!}[[(\cV_l')^*\bt\cV''_l]^\pl]\bigr)_{\wt=-1}
\nonumber\\
&\;\>\cup c_\top\bigl(\bL_{\acM_{(\ga_{i,j},\bs e_i)}^\pl/\M_{\ga_{i,j}}^\pl}\bp\bL_{\acM_{(\de_{i,j},\bs f_i)}^\pl/\M_{\de_{i,j}}^\pl}\bigr){}_{\wt=0}\bigr]\bigr)\bigr\}
\allowdisplaybreaks
\nonumber\\
&=(-1)^{\la_k(\ga_{i,j})(\la_k(\ga_{i,j})-1)/2+\la_k(\de_{i,j})(\la_k(\de_{i,j})-1)/2+\chi(\ga_{i,j},\de_{i,j})+\chi(\de_{i,j},\ga_{i,j})}\cdot{}
\nonumber\\
&\;\>\Res_z\bigl\{(\Pi^\pl_\be\ci\Phi_{\ga_{i,j},\de_{i,j}}\ci(\ac\Pi^{v_a}_{\dot{\mathcal M}_{\ga_{i,j}}}\t\ac\Pi^{v_{b_i}}_{\dot{\mathcal M}_{\de_{i,j}}}))^{\bG_m}_*
\nonumber\\
&\;\>\bigl(([\acM_{(\ga_{i,j},\bs e_i)}^\ss(\ac\tau^{s_i\ti\la}_{\bs\mu_0})]_\virt\!\bt\![\acM_{(\de_{i,j},\bs f_i)}^\ss(\ac\tau^{s_i\ti\la}_{\bs\mu_0})]_\virt)\nonumber\\
&\;\>\cap\bigl[
e\bigl((\ac\Pi^{v_a}_{\dot{\mathcal M}_{\ga_{i,j}}}\t\ac\Pi^{v_{b_i}}_{\dot{\mathcal M}_{\de_{i,j}}})^*(\cE^\bu_{\ga_{i,j},\de_{i,j}})\bigr)_{\wt=1}
\nonumber\\
&\;\>\cup (-1)^{\chi(\de_{i,j},\ga_{i,j})}e\bigl((\ac\Pi^{v_a}_{\dot{\mathcal M}_{\ga_{i,j}}}\!\!\t\!\ac\Pi^{v_{b_i}}_{\dot{\mathcal M}_{\de_{i,j}}})^*(\si^*(\cE_{\de_{i,j},\ga_{i,j}}^\bu)^\vee)\bigr)_{\wt=1}
\nonumber\\
&\;\>\cup c_\top\bigl(\bL_{\acM_{(\ga_{i,j},\bs e_i)}^\pl/\M_{\ga_{i,j}}^\pl}\bp\bL_{\acM_{(\de_{i,j},\bs f_i)}^\pl/\M_{\de_{i,j}}^\pl}\bigr){}_{\wt=0}\bigr]\bigr)\bigr\}
\allowdisplaybreaks
\nonumber\\
&=(-1)^{\la_k(\ga_{i,j})(\la_k(\ga_{i,j})-1)/2+\la_k(\de_{i,j})(\la_k(\de_{i,j})-1)/2+\chi(\ga_{i,j},\de_{i,j})}\cdot\Res_z\bigl\{
\nonumber\\
&\;\>(\Pi^\pl_\be\ci\Phi_{\ga_{i,j},\de_{i,j}})^{\bG_m}_*\bigl(\bigl((\ac\Pi^{v_a}_{\dot{\mathcal M}_{\ga_{i,j}}})_*([\acM_{(\ga_{i,j},\bs e_i)}^\ss(\ac\tau^{s_i\ti\la}_{\bs\mu_0})]_\virt\cap c_\top(\bL_{\acM_{(\ga_{i,j},\bs e_i)}^\pl/\M_{\ga_{i,j}}^\pl}))
\nonumber\\
&\;\>\bt(\ac\Pi^{v_{b_i}}_{\dot{\mathcal M}_{\de_{i,j}}})_*([\acM_{(\de_{i,j},\bs f_i)}^\ss(\ac\tau^{s_i\ti\la}_{\bs\mu_0})]_\virt\cap c_\top(\bL_{\acM_{(\de_{i,j},\bs f_i)}^\pl/\M_{\de_{i,j}}^\pl}))\bigr)
\nonumber\\
&\;\>
\cap\bigl[e(\cE^\bu_{\ga_{i,j},\de_{i,j}}\op\si^*(\cE^\bu_{\de_{i,j},\ga_{i,j}})^\vee)_{\wt=1}\bigr]\bigr)\bigr\}
\allowdisplaybreaks
\nonumber\\
&=(-1)^{\chi(\ga_{i,j},\de_{i,j})}\cdot{}
\Res_z\bigl\{(\Pi^\pl_\be\ci\Phi_{\ga_{i,j},\de_{i,j}})^{\bG_m}_*\bigl(
\bigl(\hat\Up_{(\ga_{i,j},\bs e_i)}(\ac\tau^{s_i\ti\la}_{\bs\mu_0})\bt
\hat\Up_{(\de_{i,j},\bs f_i)}(\ac\tau^{s_i\ti\la}_{\bs\mu_0})\bigr)
\nonumber\\
&\qquad\qquad\cap \bigl[e(\cE^\bu_{\ga_{i,j},\de_{i,j}}\op\si^*(\cE^\bu_{\de_{i,j},\ga_{i,j}})^\vee)_{\wt=1}  \bigr]\bigr)\bigr\}
\allowdisplaybreaks
\nonumber\\
&=\Res_z\bigl\{(-1)^{\chi(\ga_{i,j},\de_{i,j})}(\Pi^\pl_\be\ci\Phi_{\ga_{i,j},\de_{i,j}})^{\bG_m}_*\bigl(
\bigl(\hat\Up_{(\ga_{i,j},\bs e_i)}(\ac\tau^{s_i\ti\la}_{\bs\mu_0})\bt
\hat\Up_{(\de_{i,j},\bs f_i)}(\ac\tau^{s_i\ti\la}_{\bs\mu_0})\bigr)
\nonumber\\
&\qquad\qquad \cap \bigl[\ts \sum_{p\ge 0} z^{\chi(\ga_{i,j},\de_{i,j})+\chi(\de_{i,j},\ga_{i,j})-p}
c_p(\cE^\bu_{\ga_{i,j},\de_{i,j}}\op\si^*(\cE^\bu_{\de_{i,j},\ga_{i,j}})^\vee)
\bigr]\bigr)\bigr\}
\allowdisplaybreaks
\nonumber\\
&=\Res_z\bigl\{(-1)^{\chi(\ga_{i,j},\de_{i,j})}
\bigl((\Pi^\pl_\be\ci\Phi_{\ga_{i,j},\de_{i,j}}\ci(\Psi_{\ga_{i,j}}\t\id_{\M_{\de_{i,j}}}))_*
\nonumber\\
&\qquad\qquad\bigl(
\ts\sum_{q\ge 0}z^qt^q\bt\bigl(
\bigl[\hat\Up_{(\ga_{i,j},\bs e_i)}(\ac\tau^{s_i\ti\la}_{\bs\mu_0})\bt
\hat\Up_{(\de_{i,j},\bs f_i)}(\ac\tau^{s_i\ti\la}_{\bs\mu_0})\bigr]
\nonumber\\
&\qquad\qquad \cap \bigl[\ts \sum_{p\ge 0} z^{\chi(\ga_{i,j},\de_{i,j})+\chi(\de_{i,j},\ga_{i,j})-p}
c_p(\cE^\bu_{\ga_{i,j},\de_{i,j}}\op\si^*(\cE^\bu_{\de_{i,j},\ga_{i,j}})^\vee)
\bigr]\bigr)\bigr\}
\allowdisplaybreaks
\nonumber\\ 
&=(\Pi^\pl_\be)_*\Res_z\bigl\{(-1)^{\chi(\ga_{i,j},\de_{i,j})}\ts \sum_{p,q\ge 0} z^{\chi(\ga_{i,j},\de_{i,j})+\chi(\de_{i,j},\ga_{i,j})-p+q}
\nonumber\\
&\;\>\bigl(\Phi_{\ga_{i,j},\de_{i,j}}\ci(\Psi_{\ga_{i,j}}\!\t\!\id_{\M_{\de_{i,j}}})\bigr)_*\bigl[t^q\!\bt\!\bigl((\hat\Up_{(\ga_{i,j},\bs e_i)}(\ac\tau^{s_i\ti\la}_{\bs\mu_0})\!\bt\!\hat\Up_{(\de_{i,j},\bs f_i)}(\ac\tau^{s_i\ti\la}_{\bs\mu_0}))
\nonumber\\
&\;\>\cap c_p(\cE^\bu_{\ga_{i,j},\de_{i,j}}\op\si^*(\cE^\bu_{\de_{i,j},\ga_{i,j}})^\vee)\bigr)\bigr]\bigr\}
\allowdisplaybreaks
\nonumber\\
&=(\Pi^\pl_\be)_*\Res_z\bigl\{Y(\hat\Up_{(\ga_{i,j},\bs e_i)}(\ac\tau^{s_i\ti\la}_{\bs\mu_0}),z)\hat\Up_{(\de_{i,j},\bs f_i)}(\ac\tau^{s_i\ti\la}_{\bs\mu_0})\bigr\}
\allowdisplaybreaks
\nonumber\\
&=(\Pi^\pl_\be)_*\bigl\{\bigl(\hat\Up_{(\ga_{i,j},\bs e_i)}(\ac\tau^{s_i\ti\la}_{\bs\mu_0})\bigr)_0\bigl(\hat\Up_{(\de_{i,j},\bs f_i)}(\ac\tau^{s_i\ti\la}_{\bs\mu_0})\bigr)\bigr\}
\allowdisplaybreaks
\nonumber\\
&=\bigl[(\Pi^\pl_{\ga_{i,j}})_*(\hat\Up_{(\ga_{i,j},\bs e_i)}(\ac\tau^{s_i\ti\la}_{\bs\mu_0})),(\Pi^\pl_{\de_{i,j}})_*(\hat\Up_{(\de_{i,j},\bs f_i)}(\ac\tau^{s_i\ti\la}_{\bs\mu_0}))\bigr]
\nonumber\\
&=\bigl[\ac\Up_{(\ga_{i,j},\bs e_i)}(\ac\tau^{s_i\ti\la}_{\bs\mu_0}),\ac\Up_{(\de_{i,j},\bs f_i)}(\ac\tau^{s_i\ti\la}_{\bs\mu_0})\bigr].
\label{co10eq88}
\ea
Here as in \eq{co9eq45}, when we write $(\cdots)_*^{\bG_m}$ it is to emphasize that we are doing a pushforward along the 1-morphism $\cdots$ in $\bG_m$-{\it equivariant\/} homology, with respect to a {\it particular\/ $\bG_m$-equivariant structure\/} on~$\cdots$. 

In the first step of \eq{co10eq88} we make the identification \eq{co10eq64} and use
\begin{equation*}
\Pi_{\M^\pl_\be}\ci i_j=\Pi^\pl_\be\ci\Phi_{\ga_{i,j},\de_{i,j}}\ci(\ac\Pi^{v_a}_{\dot{\mathcal M}_{\ga_{i,j}}}\t\ac\Pi^{v_{b_i}}_{\dot{\mathcal M}_{\de_{i,j}}})\ci\Pi_j
\end{equation*}
where $\Pi^\pl_\be:\M_\be\ra\M^\pl_\be$ is the projection, $\Phi_{\ga_{i,j},\de_{i,j}}$ is as in Assumption \ref{co4ass1}(d), $\ac\Pi^{v_a}_{\dot{\mathcal M}_{\ga_{i,j}}},\ac\Pi^{v_{b_i}}_{\dot{\mathcal M}_{\de_{i,j}}}$ are as in Definition \ref{co5def1}, and $\Pi_j$ is as in \eq{co10eq64}, but we omit $\Pi_j$ in \eq{co10eq88} as we consider \eq{co10eq64} an identification. We also use the first case of \eq{co10eq69}, and we compute $\rank\cN_j^\bu$ using \eq{co10eq67}--\eq{co10eq68}, noting that $\dim\cV'_l=e_{i,l}$, $\dim\cV''_l=f_{i,l}$, and that $(\ga_{i,j},\bs e_i),(\de_{i,j},\bs f_i)$ satisfy the conditions of Proposition \ref{co10prop10}, so we can simplify as in~\eq{co10eq78}. 

In the second step  of \eq{co10eq88} we substitute in \eq{co10eq67}--\eq{co10eq68} and \eq{co10eq87} and use the multiplicative property of Euler classes. In the third step we cancel four terms using $e(\cV)=c_\top(\cV)$ for a vector bundle $\cV$, noting that the terms concerned, which come from \eq{co10eq75} and \eq{co10eq86}, are the K-theory classes of vector bundles. We also substitute in~$e(\si^*(\cE^\bu_{\de_{i,j},\ga_{i,j}}))=(-1)^{\chi(\de_{i,j},\ga_{i,j})}e(\si^*(\cE^\bu_{\de_{i,j},\ga_{i,j}})^\vee)$. In the fourth we use properties of (co)homology and the multiplicative property of Euler classes. In the fifth, in a similar way to \eq{co9eq2} we write
\begin{align*}
&\hat\Up_{(\ga_{i,j},\bs e_i)}(\ac\tau^{s_i\ti\la}_{\bs\mu_0})=(\ac\Pi^{v_a}_{\dot{\mathcal M}_{\ga_{i,j}}})_*([\acM_{(\ga_{i,j},\bs e_i)}^\ss(\ac\tau^{s_i\ti\la}_{\bs\mu_0})]_\virt\cap c_\top(\bT_{\acM_{(\ga_{i,j},\bs e_i)}^\pl/\M_{\ga_{i,j}}^\pl})=\\
&(-1)^{\la_k(\ga_{i,j})(\la_k(\ga_{i,j})-1)/2}(\ac\Pi^{v_a}_{\dot{\mathcal M}_{\ga_{i,j}}})_*([\acM_{(\ga_{i,j},\bs e_i)}^\ss(\ac\tau^{s_i\ti\la}_{\bs\mu_0})]_\virt\!\cap\! c_\top(\bL_{\acM_{(\ga_{i,j},\bs e_i)}^\pl/\M_{\ga_{i,j}}^\pl}),\\
&\hat\Up_{(\de_{i,j},\bs f_i)}(\ac\tau^{s_i\ti\la}_{\bs\mu_0})=(\ac\Pi^{v_{b_i}}_{\dot{\mathcal M}_{\de_{i,j}}})_*([\acM_{(\de_{i,j},\bs f_i)}^\ss(\ac\tau^{s_i\ti\la}_{\bs\mu_0})]_\virt\cap c_\top(\bT_{\acM_{(\de_{i,j},\bs f_i)}^\pl/\M_{\de_{i,j}}^\pl})=\\
&(-1)^{\la_k(\de_{i,j})(\la_k(\de_{i,j})-1)/2}(\ac\Pi^{v_{b_i}}_{\dot{\mathcal M}_{\de_{i,j}}})_*([\acM_{(\de_{i,j},\bs f_i)}^\ss(\ac\tau^{s_i\ti\la}_{\bs\mu_0})]_\virt\!\cap\! c_\top(\bL_{\acM_{(\de_{i,j},\bs f_i)}^\pl/\M_{\de_{i,j}}^\pl}),
\end{align*}
computing the signs as in \eq{co10eq78}, so that by \eq{co10eq6} we have
\e
\begin{split}
(\Pi_{\ga_{i,j}}^\pl)_*(\hat\Up_{(\ga_{i,j},\bs e_i)}(\ac\tau^{s_i\ti\la}_{\bs\mu_0}))&=\ac\Up_{(\ga_{i,j},\bs e_i)}(\ac\tau^{s_i\ti\la}_{\bs\mu_0}),\\
(\Pi_{\de_{i,j}}^\pl)_*(\hat\Up_{(\de_{i,j},\bs f_i)}(\ac\tau^{s_i\ti\la}_{\bs\mu_0}))&=\ac\Up_{(\de_{i,j},\bs f_i)}(\ac\tau^{s_i\ti\la}_{\bs\mu_0}).
\end{split}
\label{co10eq89}
\e

For the sixth step of \eq{co10eq88}, we expand $e(\cE^\bu_{\ga_{i,j},\de_{i,j}}\op\si^* (\cE^\bu_{\de_{i,j},\ga_{i,j}})^\vee)_{\wt=1}$ as in \eq{co9eq47}--\eq{co9eq48}. For the seventh, we rewrite the $\bG_m$-equivariant pushforward $(\Pi^\pl_\be\ci\Phi_{\ga_{i,j},\de_{i,j}})^{\bG_m}_*$ in the sixth step in terms of non-equivariant homology as in \eq{co9eq49}--\eq{co9eq50}. In the eighth we rearrange, and in the ninth we recognize the definition \eq{co4eq10} of the state-field correspondence $Y(-,z)$ of the vertex algebra constructed in \S\ref{co42}. In the tenth step we use Definition \ref{co4def1}, in the eleventh the definition of the Lie bracket $[\,,\,]$ on $\check H_{\rm even}(\M^\pl)$ in \eq{co4eq1} and Theorem \ref{co4thm2}, and in the twelfth \eq{co10eq89}. This proves \eq{co10eq88}, and the proposition.
\end{proof}

Finally, here is the analogue of Corollary \ref{co9cor1}, which proves \eq{co10eq51}, and completes the proof of Theorem \ref{co5thm2}. It follows from \eq{co10eq80}--\eq{co10eq82} and~\eq{co10eq84}.

\begin{cor}
\label{co10cor6}
Equation \eq{co10eq51} holds in the Lie algebra $\check H_{\rm even}(\M^\pl),$ that~is,
\begin{equation*}
\ac\Up_{(\be,\bs d)}(\ac\tau^{s_i'\ti\la}_{\bs\mu_0})=\ac\Up_{(\be,\bs d)}(\ac\tau^{s_{i-1}'\ti\la}_{\bs\mu_0})+\sum_{\begin{subarray}{l} j=1,\ldots,q_i: \\ o_{\ga_{i,j}}+o_{\de_{i,j}}=o_\be \end{subarray}\!\!\!\!\!\!\!\!\!\!\!\!\!\!\!}\bigl[\ac\Up_{(\ga_{i,j},\bs e_i)}(\ac\tau^{s_i\ti\la}_{\bs\mu_0}),\ac\Up_{(\de_{i,j},\bs f_i)}(\ac\tau^{s_i\ti\la}_{\bs\mu_0})\bigr].
\end{equation*}
\end{cor}

\subsection{\texorpdfstring{Extension to the $G$-equivariant case}{Extension to the G-equivariant case}}
\label{co107}

Theorem \ref{co5thm4} extends Theorems \ref{co5thm1}--\ref{co5thm3} to $G$-equivariant homology $H_*^G(\M^\pl)$. As in \S\ref{co94}, the proof of Theorem \ref{co5thm2} in \S\ref{co101}--\S\ref{co106} must be modified as follows:
\begin{itemize}
\setlength{\itemsep}{0pt}
\setlength{\parsep}{0pt}
\item[(a)] We replace $H_*(-),H^*(-)$ by $H^G_*(-),H_G^*(-)$ throughout.
\item[(b)] $G$ acts naturally on all the classical moduli stacks, and all the derived moduli stacks, in the proofs. Therefore the obstruction theories used to define virtual classes, which come from cotangent complexes of quasi-smooth derived moduli stacks, are $G$-equivariant, and define virtual classes in equivariant homology.
\item[(c)] In Proposition \ref{co10prop11}(c), $G$ acts on $\grM^\ss_{(\be,(1,\bs d))}(\gr\tau^{s_i\ti\la}_{\bs{\gr\mu}})_j$ on the left of \eq{co10eq64}, but $G\t G$ acts on $\acM_{(\ga_{i,j},\bs e_i)}^\ss(\ac\tau^{s_i\ti\la}_{\bs\mu_0})\t\acM_{(\de_{i,j},\bs f_i)}^\ss(\ac\tau^{s_i\ti\la}_{\bs\mu_0})$. Then \eq{co10eq64} is $G$-equivariant for the diagonal $G$-action on $\acM_{(\ga_{i,j},\bs e_i)}^\ss(\ac\tau^{s_i\ti\la}_{\bs\mu_0})\t\acM_{(\de_{i,j},\bs f_i)}^\ss(\ac\tau^{s_i\ti\la}_{\bs\mu_0})$.
\item[(d)] In Proposition \ref{co10prop12}(d) we must replace \eq{co10eq69} by
\end{itemize}
\ea
&[\grM^\ss_{(\be,(1,\bs d))}(\gr\tau^{s_i\ti\la}_{\bs{\gr\mu}})_j]_\virt
\label{co10eq90}\\
&=\begin{cases}
\La^{G\t G,G}\bigl([\acM_{(\ga_{i,j},\bs e_i)}^\ss(\ac\tau^{s_i\ti\la}_{\bs\mu_0})]_\virt\bt[\acM_{(\de_{i,j},\bs f_i)}^\ss(\ac\tau^{s_i\ti\la}_{\bs\mu_0})]_\virt\bigr), & o_{\ga_{i,j}}+o_{\de_{i,j}}=o_\be, \\
0, & o_{\ga_{i,j}}+o_{\de_{i,j}}>o_\be.
\end{cases}
\nonumber
\ea
\begin{itemize}
\setlength{\itemsep}{0pt}
\setlength{\parsep}{0pt}
\item[] Here $\La^{G\t G,G}$ maps $H_*^{G\t G}(\cdots)$ to $H_*^G(\cdots)$ as in Definition~\ref{co2def2}(g).
\item[(e)] In the first step of equation \eq{co10eq88} in the proof of Proposition \ref{co10prop14}, we substitute in \eq{co10eq90} rather than \eq{co10eq69}. This gives an extra operator $\La^{G\t G,G}$, which persists through the second--eighth steps of \eq{co10eq88}. In the ninth step we use the $G$-equivariant version \eq{co4eq24} of the state-field correspondence $Y(-,z)$, which includes an extra operator $\La^{G\t G,G}$ compared to the usual definition \eq{co4eq10}. Thus the $\La^{G\t G,G}$ is absent in the ninth--twelfth steps of \eq{co10eq88}, and \eq{co10eq84} holds in~$\check H_{\rm even}^G(\M^\pl)$.
\end{itemize}

\section{Proof of Theorem \ref{co5thm3}}
\label{co11}

Let Assumptions \ref{co5ass1}--\ref{co5ass3} hold. We will prove Theorem \ref{co5thm3} by induction on $\rk\al=1,2,\ldots$ in Assumption \ref{co5ass2}(f). Let $l\ge 0$ be given, and suppose by induction that \eq{co5eq34} holds for all $(\tau,T,\le),(\ti\tau,\ti T,\le)\in\sS$ and $\al\in C(\B)_\pe$ with $\rk\al\le l$ (this is vacuous when $l=0$). Suppose that $(\tau,T,\le),(\ti\tau,\ti T,\le)\in\sS$ and $\al\in C(\B)_\pe$ with $\rk\al=l+1$. Assumption \ref{co5ass3}(a) gives a continuous family $(\tau_t,T_t,\le)_{t\in[0,1]}$ in $\sS$ with $(\tau_0,T_0,\le)=(\tau,T,\le)$ and $(\tau_1,T_1,\le)=(\ti\tau,\ti T,\le)$. These $(\tau,T,\le),(\ti\tau,\ti T,\le),\ab\al,\ab(\tau_t,T_t,\le)_{t\in[0,1]}$ will be fixed throughout the proof.

\begin{dfn}
\label{co11def1}
Let $\al\in C(\B)_\pe$ be fixed as above. Consider sets of data $n\ge p\ge 1$, $\al_1,\ldots,\al_n\in C(\B)$ with $\al_1+\cdots+\al_n=\al$, and $0=a_0<a_1<\cdots<a_p=n$, and set $\be_j=\al_{a_{j-1}+1}+\cdots+\al_{a_j}$ for $j=1,\ldots,p$. Write $\bs\al=(\al_i)_{i=1}^n$ and $\bs a=(a_j)_{j=0}^p$. For such $\bs\al,\bs a$, define the condition:
\begin{itemize}
\setlength{\itemsep}{0pt}
\setlength{\parsep}{0pt}
\setlength{\itemindent}{5pt}
\item[$(*)_{\bs\al,\bs a}$] There exist $s,t,u\in[0,1]$ such that $\M_{\al_i}^\ss(\tau_s)\ne\es$ for all $i$, and $U(\al_{a_{j-1}+1},\ab\ldots,\ab\al_{a_j};\tau_s,\tau_t)\ne 0$ for all $j=1,\ldots,p$ and~$\tau_u(\be_1)=\cdots=\tau_u(\be_p)$.
\end{itemize}
Assumption \ref{co5ass3}(b) says there are only finitely many $\bs\al,\bs a$ satisfying~$(*)_{\bs\al,\bs a}$.

For each $\bs\al,\bs a$ satisfying $(*)_{\bs\al,\bs a}$ and for $t\in[0,1]$, we have (in)equalities
\e
\begin{aligned}
\tau_t(\al_1+\cdots+\al_i)&<	\tau_t(\al_{i+1}+\cdots+\al_n),&& \tau_t(\al_i)\le\tau_t(\al_{i+1}),\\
\tau_t(\al_1+\cdots+\al_i)&=\tau_t(\al_{i+1}+\cdots+\al_n), && \tau_t(\al_i)>\tau_t(\al_{i+1}), \\
\tau_t(\al_1+\cdots+\al_i)&>\tau_t(\al_{i+1}+\cdots+\al_n), && \text{for all $1\le i<n$.}
\end{aligned}
\label{co11eq1}
\e
As $(\tau_t,T_t,\le)_{t\in[0,1]}$ is a continuous family, Definition \ref{co3def5} implies that the set of $t\in[0,1]$ satisfying any (in)equality in \eq{co11eq1} is the disjoint union of a finite set of intervals in $[0,1]$. Write $0=t_0<t_1<\cdots<t_N=1$ for the finite set of $t\in[0,1]$ which are the end points of a connected component of the set of $t\in[0,1]$ satisfying one of the (in)equalities in \eq{co11eq1} for some $\bs\al,\bs a$ satisfying $(*)_{\bs\al,\bs a}$, together with 0,1 if these are not already end points. There are only finitely many $t_a$ as there are only finitely many choices of $\bs\al,\bs a$, and for each $\bs\al,\bs a$ only finitely many (in)equalities \eq{co11eq1}. Choose $s_a\in(t_{a-1},t_a)$ for $a=1,\ldots,N$, so that~$0=t_0<s_1<t_1<\cdots <t_{N-1}<s_N<t_N=1$.

Let $\bs\al,\bs a$ satisfy $(*)_{\bs\al,\bs a}$, and $s,t\in[0,1]$. Observe that the conditions that $\tau_t(\be_1)=\cdots=\tau_t(\be_p)$, and the values of $U(\al_{a_{j-1}+1},\ab\ldots,\ab\al_{a_j};\tau_s,\tau_t)$ for $j=1,\ldots,p$, depend only on whether certain inequalities from \eq{co11eq1}, for $s$ and for $t$, hold or not. Thus we see that:
\begin{itemize}
\setlength{\itemsep}{0pt}
\setlength{\parsep}{0pt}
\item[(i)] The values of $U(\al_{a_{j-1}+1},\ab\ldots,\ab\al_{a_j};\tau_s,\tau_t)$ for $j=1,\ldots,p$ depend only on whether $s=t_a$ for some $0\le a\le N$, or $s\in(t_{a-1},t_a)$ for $a=1,\ldots,N$, and $t=t_b$ for some $0\le b\le N$, or $t\in(t_{b-1},t_b)$ for $b=1,\ldots,N$, and not otherwise on the values of $s,t$.

In particular, for fixed $t$ the function $[0,1]\ra\Q$ mapping $s\mapsto U(\al_{a_{j-1}+1},\ab\ldots,\ab\al_{a_j};\tau_s,\tau_t)$ can be discontinuous only at $s=t_a$ for $a=0,\ldots,N$, and is constant on $(t_{a-1},t_a)$ for each $a=1,\ldots,N$. Similarly, for fixed $s$ the function $[0,1]\ra\Q$ mapping $t\mapsto U(\al_{a_{j-1}+1},\ab\ldots,\ab\al_{a_j};\tau_s,\tau_t)$ is constant on $(t_{a-1},t_a)$ for each $a=1,\ldots,N$.
\item[(ii)] Whether or not $\tau_u(\be_1)=\cdots=\tau_u(\be_p)$ holds depends only on whether $u=t_a$ for some $0\le a\le N$, or $u\in(t_{a-1},t_a)$ for $a=1,\ldots,N$, and not otherwise on the value of~$u$.
\end{itemize}
\end{dfn}

\begin{lem}
\label{co11lem1}
If\/ $s\in[0,1]$ and\/ $\al_1,\ldots,\al_n\in C(\B)$ with\/ $\al_1+\cdots+\al_n=\al$ in $C(\B)_\pe$ and\/ $\M_{\al_i}^\ss(\tau_s)\ne\es$ for all\/ $i$ then the function $[0,1]\ra\Q$ mapping $t\mapsto U(\al_1,\ldots,\al_n;\tau_s,\tau_t)$ can be discontinuous only at\/ $t=t_a$ for some $a=0,\ldots,N,$ so it is constant in\/ $(t_{a-1},t_a)$ for\/~$a=1,\ldots,N$.	
\end{lem}

\begin{proof} This follows from the case $p=1$, $\bs a=(0,1)$ above.	
\end{proof}

\begin{lem}
\label{co11lem2}
The moduli spaces $\M_\al^\ss(\tau_t)$ and classes $[\M_\al^\ss(\tau_t)]_\inv$ in Theorem\/ {\rm\ref{co5thm1}} are constant functions of\/ $t\in(t_{a-1},t_a)$ for $a=1,\ldots,N$.
\end{lem}

\begin{proof} Let $a=1,\ldots,N$ and $s,t\in(t_{a-1},t_a)$. Suppose $[E]$ is a $\C$-point of $\M_\al^\ss(\tau_t)$. By Assumption \ref{co5ass2}(a), $E$ has a unique $\tau_s$-Harder--Narasimhan filtration $0=E_0\subsetneq E_1\subsetneq \cdots\subsetneq E_n=E$ in $\B$, such that $F_i=E_i/E_{i-1}$ is $\tau_s$-semistable (so that $\M_{\al_i}^\ss(\tau_s)\ne\es$), with $\tau_s(\al_1)>\cdots>\tau_s(\al_n)$, where $\al_i=\lb F_i\rb$ for $i=1,\ldots,n$. Corollary \ref{co3cor1} shows that $U(\al_1,\ldots,\al_n;\tau_s,\tau_t)\ne 0$.

Setting $p=1$, $a_0=0$, $a_1=n$ we see that $\bs\al,\bs a$ satisfy $(*)_{\bs\al,\bs a}$ in Definition \ref{co11def1}. Therefore the condition $U(\al_1,\ldots,\al_n;\tau_s,\tau_t)\ne 0$ depends only on the fact that $t\in (t_{a-1},t_a)$, not on the actual value of $t$, so we also have $U(\al_1,\ldots,\al_n;\tau_s,\tau_s)\ne 0$ as $s\in(t_{a-1},t_a)$, and thus $n=1$ by equation \eq{co3eq4}. Hence the $\tau_s$-Harder--Narasimhan filtration of $E$ is $0=E_0\subset E_1=E$, so $E$ is $\tau_s$-semistable, and $\M_\al^\ss(\tau_t)\subseteq\M_\al^\ss(\tau_s)$. Exchanging $s,t$ gives $\M_\al^\ss(\tau_s)\subseteq\M_\al^\ss(\tau_t)$, so $\M_\al^\ss(\tau_s)=\M_\al^\ss(\tau_t)$. This implies that $[\M_\al^\ss(\tau_s)]_\inv =[\M_\al^\ss(\tau_t)]_\inv$ by Theorem \ref{co5thm1}(ii), and the lemma follows.
\end{proof}

\begin{lem}
\label{co11lem3}
For each\/ $a=1,\ldots,N,$ if\/ $\al_1,\ldots,\al_n\in C(\B)$ with\/ $\al_1+\cdots+\al_n=\al,$ and either 
\begin{itemize}
\setlength{\itemsep}{0pt}
\setlength{\parsep}{0pt}
\item[{\bf(i)}] $U(\al_1,\ldots,\al_n;\tau_{t_a},\tau_{s_a})\ne 0$ and\/ $\M_{\al_i}^\ss(\tau_{t_a})\ne\es$ for\/ {\rm $i=1,\ldots,n$;} or 
\item[{\bf(ii)}] $U(\al_1,\ldots,\al_n;\tau_{s_a},\tau_{t_a})\ne 0$ and\/ $\M_{\al_i}^\ss(\tau_{s_a})\ne\es$ for\/ $i=1,\ldots,n,$
\end{itemize}
then\/ $\tau_{t_a}(\al_i)=\tau_{t_a}(\al)$ for all\/ $i=1,\ldots,n$. Also in case {\bf(ii)} we have $\M_{\al_i}^\ss(\tau_{s_a})\subseteq\M_{\al_i}^\ss(\tau_{t_a})$ for\/ $i=1,\ldots,n$. Hence Theorem\/ {\rm\ref{co5thm2}} shows that\/ {\rm\eq{co5eq31}--\eq{co5eq32}} hold with $\tau_{t_a},\tau_{s_a}$ in place of\/ $\tau,\ti\tau$.

The analogue also holds with $\tau_{t_{a-1}}$ in place of\/~$\tau_{t_a}$.
\end{lem}

\begin{proof} Let $a=1,\ldots,N$ and $\al_1,\ldots,\al_n\in C(\B)$ with $\al_1+\cdots+\al_n=\al$, and set $p=1$, $a_0=0$, $a_1=1$ and $\be_1=\al$. If (i) holds then $(*)_{\bs\al,\bs a}$ holds with $t_a,s_a$ in place of $s,t$, and Definition \ref{co11def1} implies that $U(\al_1,\ldots,\al_n;\tau_{t_a},\tau_{s_a})=U(\al_1,\ldots,\al_n;\tau_{t_a},\tau_s)\ne 0$ for any $s\in (t_{a-1},t_a)$. Applying Proposition \ref{co3prop4} with $t_a,t_a,s$ in place of $t_-,t_0,t_+$ (here we deal with $t_->t_+$ by changing the sign of the $t$ variable in Proposition \ref{co3prop4}), and taking $\md{s-t_a}$ sufficiently small, shows that $\tau_{t_a}(\al_1)=\cdots=\tau_{t_a}(\al_n)$, so $\tau_{t_a}(\al_i)=\tau_{t_a}(\al)$ for all $i$, as we want.

Similarly, if (ii) holds then $(*)_{\bs\al,\bs a}$ holds with $s_a,t_a$ in place of $s,t$, and Definition \ref{co11def1} implies that $U(\al_1,\ldots,\al_n;\tau_{s_a},\tau_{t_a})=U(\al_1,\ldots,\al_n;\tau_s,\tau_{t_a})\ne 0$ for any $s\in (t_{a-1},t_a)$. Applying Proposition \ref{co3prop4} with $s,t_a,t_a$ in place of $t_-,t_0,t_+$ and taking $\md{s-t_a}$ sufficiently small gives $\tau_{t_a}(\al_i)=\tau_{t_a}(\al)$ for all~$i$.

Suppose $a,\al_1,\ldots,\al_n$ satisfy (ii), so that $\tau_{t_a}(\al_i)=\tau_{t_a}(\al)$ for all $i$ from above, and let $r\in(t_{a-1},t_a)$. Choose a $\C$-point $[E_i]$ in $\M_{\al_i}^\ss(\tau_{s_a})$ for $i=1,\ldots,n$. By Assumption \ref{co5ass2}(a) $E_i$ has a unique $\tau_r$-Harder--Narasimhan filtration $0=E_i^0\subsetneq E_i^1\subsetneq \cdots\subsetneq E_i^{m_i}=E_i$ in $\B$, such that $F_i^j=E_i^j/E_i^{j-1}$ is $\tau_r$-semistable with $\tau_r(\be_i^1)>\cdots>\tau_r(\be_i^{m_i})$, where $\be_i^j=\lb F_i^j\rb$. Corollary \ref{co3cor1} shows that $U(\be_i^1,\ldots,\be_i^{m_i};\tau_r,\tau_{s_a})\ne 0$ for~$i=1,\ldots,n$. 

Writing $\bs\be=(\be_1^1,\ldots,\be_1^{m_1},\be_2^1,\ldots,\be_n^{m_n})$ and $\bs a=(0,m_1,m_1+m_2,\ab\ldots,\ab m_1+\cdots+m_n)$ we see that $(*)_{\bs\be,\bs a}$ in Definition \ref{co11def1} holds with $r,s_a,t_a,n,\sum_im_i$ in place of $s,t,u,p,n$. Therefore $U(\be_i^1,\ldots,\be_i^{m_i};\tau_r,\tau_{s_a})\ne 0$ is independent of $r\in(t_{a-1},t_a)$. Putting $r=s_a$ and using \eq{co3eq4} shows that $m_i=1$. Thus the $\tau_r$-Harder--Narasimhan filtration of $E_i$ is $0\subsetneq E_i$, and $E_i$ is $\tau_r$-semistable. Suppose $0\ne E_i'\subsetneq E_i$ is a proper subobject of $E_i$. Then $\tau_r(\lb E_i'\rb)\le \tau_r(\lb E_i/E_i'\rb)$ for $r\in(t_{a-1},t_a)$ as $E_i$ is $\tau_r$-semistable. Taking the limit $r\ra t_a$ and using the continuity of $(\tau_t,T_t,\le)_{t\in[0,1]}$ shows that $\tau_{t_a}(\lb E_i'\rb)\le \tau_{t_a}(\lb E_i/E_i'\rb)$. Thus $E_i$ is $\tau_{t_a}$-semistable, so $\M_{\al_i}^\ss(\tau_{s_a})\subseteq\M_{\al_i}^\ss(\tau_{t_a})$, as we have to prove.

This proves the hypotheses of Theorem \ref{co5thm2} with $\tau_{t_a},\tau_{s_a}$ in place of $\tau,\ti\tau$, so \eq{co5eq31}--\eq{co5eq32} hold with $\tau_{t_a},\tau_{s_a}$ in place of $\tau,\ti\tau$. The last part is proved in an almost identical way, by replacing $t_a$ by $t_{a-1}$ above, and changing the sign of the $t$ variable in Proposition \ref{co3prop4} for (ii) rather than for~(i).	
\end{proof}

\begin{lem}
\label{co11lem4}
Let\/ $a=1,\ldots,N$. Then \eq{co5eq34} holds with $\tau_{t_0},\tau_{t_a}$ in place of\/ $\tau,\ti\tau$ if and only if\/ \eq{co5eq34} holds with $\tau_{t_0},\tau_{s_a}$ in place of\/ $\tau,\ti\tau$.	

The same holds with $\tau_{t_{a-1}}$ in place of\/ $\tau_{t_a}$.
\end{lem}

\begin{proof} First suppose that \eq{co5eq34} holds with $\tau_{t_0},\tau_{t_a}$ in place of $\tau,\ti\tau$. Recall the inductive hypothesis at the beginning of Chapter \ref{co11}. From these we see that if either $\be\in C(\B)_\pe$ with $\rk\be\le l$, or $\be=\al$, then
\e
\begin{gathered}
{}
[\M_\be^\ss(\tau_{t_a})]_\inv= \!\!\!\!\!\!\!
\sum_{\begin{subarray}{l}n\ge 1,\;\al_1,\ldots,\al_n\in
C(\B)_\pe:\\ \M_{\al_i}^\ss(\tau_{t_0})\ne\es,\; \text{all $i,$} \\
\al_1+\cdots+\al_n=\be, \; o_{\al_1}+\cdots+o_{\al_n}=o_\be 
\end{subarray}\!\!\!\!\!\!\!\!\!\!\!\!\!\!\!\!\!\!\!\!\!\!\!\!\!\!\!\!} \!\!\!\!\!\!\!
\begin{aligned}[t]
U(\al_1,&\ldots,\al_n;\tau_{t_0},\tau_{t_a})\cdot [\M_{\al_1}^\ss(\tau_{t_0})]_\inv*{}\\
&\,\,\,
[\M_{\al_2}^\ss(\tau_{t_0})]_\inv*\cdots *[\M_{\al_n}^\ss(\tau_{t_0})]_\inv.
\end{aligned}
\end{gathered}
\label{co11eq2}
\e
with only finitely many possibilities for $n\ge 1$, $\al_1,\ldots,\al_n$ in the sum.

Lemma \ref{co11lem3} says \eq{co5eq32} holds with $\tau_{t_a},\tau_{s_a}$ in place of $\tau,\ti\tau$, giving 
\e
\begin{gathered}
{}
[\M_\al^\ss(\tau_{s_a})]_\inv= \!\!\!\!\!\!\!
\sum_{\begin{subarray}{l}p\ge 1,\;\be_1,\ldots,\be_p\in
C(\B)_\pe:\\ \M_{\be_j}^\ss(\tau_{t_a})\ne\es,\; \text{all $j$} \\
\be_1+\cdots+\be_p=\al, \; o_{\be_1}+\cdots+o_{\be_p}=o_\al
\end{subarray}\!\!\!\!\!\!\!\!\!\!\!\!\!\!\!\!\!\!\!\!\!\!\!\!\!\!\!\!} \!\!\!\!\!\!
\begin{aligned}[t]
U(\be_1,\ldots,\be_p;\tau_{t_a},\tau_{s_a})\cdot [\M_{\be_1}^\ss(\tau_{t_a})]_\inv*{}&\\
[\M_{\be_2}^\ss(\tau_{t_a})]_\inv*\cdots *[\M_{\be_p}^\ss(\tau_{t_a})]_\inv.&
\end{aligned}
\end{gathered}
\label{co11eq3}
\e
We claim that we can replace \eq{co11eq3} by the sum
\ea
&[\M_\al^\ss(\tau_{s_a})]_\inv= 
\label{co11eq4}\\
&\sum_{\begin{subarray}{l}p\ge 1,\;\be_1,\ldots,\be_p\in
C(\B)_\pe: \\
\be_1+\cdots+\be_p=\al, \; o_{\be_1}+\cdots+o_{\be_p}=o_\be, \\ 
U(\be_1,\ldots,\be_p;\tau_{t_a},\tau_{s_a})\ne 0, \; \tau_{t_a}(\be_1)=\cdots=\tau_{t_a}(\be_p),\\ 
\exists \,n\ge p,\;\al_1,\ldots,\al_n\in C(\B)_\pe\; \text{and}\; 0=a_0<a_1<\cdots<a_p=n\; \text{with}\; \M_{\al_i}^\ss(\tau_{t_0})\ne\es,\; \text{all $i$}\\
\al_{a_{j-1}+1}\!+\cdots+\al_{a_j}=\be_j,\, o_{\al_{a_{j-1}+1}}+\cdots+o_{\al_{a_j}}=o_{\be_j},\, U(\al_{a_{j-1}+1},\ldots,\al_{a_j};\tau_{t_0},\tau_{t_a})\ne 0,\, \text{all $j$}
\end{subarray}\!\!\!\!\!\!\!\!\!\!\!\!\!\!\!\!\!\!\!\!\!\!\!\!\!\!\!\!\!\!\!\!\!\!\!\!\!\!\!\!\!\!\!\!\!\!\!\!\!\!\!\!\!\!\!\!\!\!\!\!\!\!\!\!\!\!\!\!\!\!\!\!} \!\!\!\!\!\!\!\!\!\!\!\!\!\!\!\!\!\!\!\!\!\!\!\!\!\!\!\!\!\!\!\!\!\!\!\!\!\!\!\!\!\!\begin{aligned}[t]
U(\be_1,&\ldots,\be_p;\tau_{t_a},\tau_{s_a})\cdot [\M_{\be_1}^\ss(\tau_{t_a})]_\inv*{}\\
&
[\M_{\be_2}^\ss(\tau_{t_a})]_\inv*\cdots *[\M_{\be_p}^\ss(\tau_{t_a})]_\inv.
\end{aligned}
\nonumber
\ea
Here we have replaced the condition $\M_{\be_j}^\ss(\tau_{t_a})\ne\es$ in \eq{co11eq3} by much more complicated conditions on $\be_1,\ldots,\be_p$. To see that \eq{co11eq4} makes sense, and \eq{co11eq3} and \eq{co11eq4} are equivalent, we will prove:
\begin{itemize}
\setlength{\itemsep}{0pt}
\setlength{\parsep}{0pt}
\item[(i)] There are only finitely many terms in the sum \eq{co11eq4}.
\item[(ii)] Every nonzero term in \eq{co11eq3} also appears in \eq{co11eq4}.
\item[(iii)] Every term in \eq{co11eq4} which does not appear in \eq{co11eq3} is zero.
\end{itemize}

For (i), Assumption \ref{co5ass3}(b) says that there are only finitely many sets of data
$n\ge p\ge 1$, $\al_1,\ldots,\al_n\in C(\B)_\pe$ and $0=a_0<a_1<\cdots<a_p=n$ such that $\al_1+\cdots+\al_n=\al$, $\M_{\al_i}^\ss(\tau_{t_0})\ne\es$ for all $i$, $U\bigl(\al_{a_{j-1}+1},\ldots,\al_{a_j};\tau_{t_0},\tau_{t_a})\ne 0$ for all $j$, and $\tau_{t_a}(\be_1)=\cdots=\tau_{t_a}(\be_p)$, where $\be_j=\al_{a_{j-1}+1}+\cdots+\al_{a_j}$ for all $j$. Thus the sum in \eq{co11eq4} is finite.

For (ii), if $p\ge 1$, $\be_1,\ldots,\be_p$ give a nonzero term in \eq{co11eq3} then $U(\be_1,\ab\ldots,\ab\be_p;\ab\tau_{t_a},\ab\tau_{s_a})\ne 0$ and $\be_j\in C(\B)_\pe$ with $o_{\be_1}+\cdots+o_{\be_p}=o_\al$, $\M_{\be_j}^\ss(\tau_{t_a})\ne\es$ and $[\M_{\be_j}^\ss(\tau_{t_a})]_\inv\ne 0$ for all $j$. Thus Lemma \ref{co11lem3}(ii) gives $\tau_{t_a}(\be_1)=\cdots=\tau_{t_a}(\be_p)$. Also in \eq{co11eq2} for $\be=\be_j$ there is at least one nonzero term, say from $n^j\ge 1$, $\al_1^j,\ldots,\al_{n^j}^j$ with $\al_1^j+\cdots+\al_{n^j}^j=\be_j$, $o_{\al_1^j}+\cdots+o_{\al_{n^j}^j}=o_{\be_j}$, $\M_{\al_i^j}^\ss(\tau_{t_0})\ne\es$ and $U(\al_1^j,\ldots,\al_{n^j}^j;\tau_{t_0},\tau_{t_a})\ne 0$. Then by taking $n=n^1+\cdots+n^p$, $a_j=n^1+\cdots+n^j$ and $\al_1,\ldots,\al_n=\al_1^1,\ldots,\al_{n^1}^1,\al_2^1,\ldots,\al_{n^p}^p$, we obtain data satisfying the conditions in the sum \eq{co11eq4} for $\be_1,\ldots,\be_p$.

For (iii), suppose $p\ge 1$, $\be_1,\ldots,\be_p$ satisfy the conditions for the sum in \eq{co11eq4}, but not for the sum in \eq{co11eq3}. Then we must have $\M_{\be_j}^\ss(\tau_{t_a})=\es$ for some $j=1,\ldots,p$. But then $[\M_{\be_j}^\ss(\tau_{t_a})]_\inv=0$ by Theorem \ref{co5thm1}(i), so the term in \eq{co11eq4} is zero. This proves~\eq{co11eq4}.

Observe that for each term $p\ge 1$, $\be_1,\ldots,\be_p$ in \eq{co11eq4}, either $p=1$ and $\be_1=\al$, or $p>1$ and $\rk\be_j<\rk\al=l+1$ by Assumption \ref{co5ass2}(f) since $\al=\be_1+\cdots+\be_p$ and $\tau_{t_a}(\be_1)=\cdots=\tau_{t_a}(\be_p)$, so \eq{co11eq2} holds with $\be=\be_j$. Substituting \eq{co11eq2} with $\be=\be_j$ into \eq{co11eq4} for all $j$ and rearranging yields 
\ea
&[\M_\al^\ss(\tau_{s_a})]_\inv= 
\label{co11eq5}\\
&\sum_{\begin{subarray}{l}n\ge 1,\\ \al_1,\ldots,\al_n\in
C(\B)_\pe:\\ \al_1+\cdots+\al_n=\al,\\ o_{\al_1}+\cdots+o_{\al_n}=o_\al, \\ 
\M_{\al_i}^\ss(\tau_{t_0})\ne\es,\; \text{all\/ $i$}\end{subarray}}
\begin{aligned}[t]
&\raisebox{-7pt}{$\Biggl\{$}\sum_{\begin{subarray}{l} p,\; a_0,\ldots,a_p: \; p=1,\ldots,n,\\ 0=a_0<a_1<\cdots<a_p=n, \\
\text{set\/ $\be_j=\al_{a_{j-1}+1}+\cdots+\al_{a_j}$,} \\  
\text{$1\!\le\! j\!\le\! p$. Suppose}\; o_{\al_{a_{j-1}+1}}+\cdots+o_{\al_{a_j}}=o_{\be_j},\\
\be_j\in C(\B)_\pe,\; \text{all $j$, and}\; \tau_{t_a}(\be_1)=\cdots=\tau_{t_a}(\be_p)
\end{subarray}\!\!\!\!\!\!\!\!\!\!\!\!\!\!\!\!}\!\!\!\!\!\!
\begin{aligned}[t]
&U(\be_1,\ldots,\be_p;\tau_{t_a},\tau_{s_a})\cdot\\
&\ts\prod_{j=1}^pU(\al_{a_{j-1}+1},\al_{a_{j-1}+2},\\
&\qquad\qquad \ldots,\al_{a_j};\tau_{t_0},\tau_{t_a})
\end{aligned}\raisebox{-25pt}{$\Biggr\}$}\\
&\qquad\cdot [\M_{\al_1}^\ss(\tau_{t_0})]_\inv*[\M_{\al_2}^\ss(\tau_{t_0})]_\inv*\cdots *[\M_{\al_n}^\ss(\tau_{t_0})]_\inv.
\end{aligned}
\nonumber
\ea

We would like to drop the assumptions that $o_{\al_{a_{j-1}+1}}+\cdots+o_{\al_{a_j}}=o_{\be_j}$, $\be_j\in C(\B)_\pe$, and $\tau_{t_a}(\be_1)=\cdots=\tau_{t_a}(\be_p)$ in the inner sum $\{\cdots\}$ in \eq{co11eq5}. If this was allowed without changing the sum, then the inner sum $\{\cdots\}$ would equal $U(\al_1,\ldots,\al_n;\tau_{t_0},\tau_{s_a})$ by equation \eq{co3eq5} of Theorem \ref{co3thm2}, so \eq{co11eq5} would imply \eq{co5eq34} with $\tau_{t_0},\tau_{s_a}$ in place of $\tau,\ti\tau$, as we have to prove.

Assumption \ref{co5ass3}(b) gives $\al_{i_1}+\cdots+\al_{i_2}\in C(\B)_\pe$ for $1\le i_1\le i_2\le n$. Thus Assumption \ref{co5ass1}(f)(iv) implies that $o_{\al_{i_1}+\cdots+\al_{i_2-1}}+o_{\al_{i_2}}\ge o_{\al_{i_1}+\cdots+\al_{i_2}}$ for $1\le i_1<i_2\le n$. Using this, $o_{\al_1}+\cdots+o_{\al_n}=o_\al$ and inductive arguments we can show that $o_{\al_{i_1}}+\cdots+o_{\al_{i_2}}=o_{\al_{i_1}+\cdots+\al_{i_2}}$ for all $1\le i_1<i_2\le n$. Hence $o_{\al_1^j}+\cdots+o_{\al_{n^j}^j}=o_{\be_j}$ for all $j$. Also Assumption \ref{co5ass3}(b) implies that $\be_1,\ldots,\be_p\in C(\B)_\pe$. Thus we can omit the first two  assumptions.

Note that each nonzero term $\bs\al=(\al_i)_{i=1}^n$, $\bs a=(a_j)_{j=0}^p$ in \eq{co11eq5} satisfies $(*)_{\bs\al,\bs a}$ in Definition \ref{co11def1}, since~$\M_{\al_i}^\ss(\tau_{t_0})\ne\es$.

Now let $\bs\al,\bs a$ be given, which need not satisfy the conditions in \eq{co11eq5}. Writing $\be_j=\al_{a_{j-1}+1}+\cdots+\al_{a_j}$ as usual, suppose that $U(\be_1,\ldots,\be_p;\tau_{t_a},\tau_{s_a})\ne 0$. Suppose also that $(*)_{\bs\al,\bs a'}$ holds for some $\bs a'=(a'_j)_{j=0}^{p'}$, though $(*)_{\bs\al,\bs a}$ need not hold. Then $U(\be_1,\ldots,\be_p;\tau_{t_a},\tau_s)=U(\be_1,\ldots,\be_p;\tau_{t_a},\tau_{s_a})\ne 0$ for any $s\in(t_{a-1},t_a)$ by Definition \ref{co11def1}, as $(*)_{\bs\al,\bs a'}$ holds, so $U(\be_1,\ldots,\be_p;\tau_{t_a},\tau_s)$ depends only on the (in)equalities \eq{co11eq1}. Thus Proposition \ref{co3prop4} with $t_a,t_a,s$ in place of $t_-,t_0,t_+$ implies that $\tau_{t_a}(\be_1)=\cdots=\tau_{t_a}(\be_p)$, since $U(\be_1,\ab\ldots,\ab\be_p;\ab\tau_{t_a},\ab\tau_s)\ne 0$ for $s\in(t_{a-1},t_a)$ with $\md{s-t_a}$ arbitrarily small.

This implies that if we fix $n,\al_1,\ldots,\al_n$ in the outer sum in \eq{co11eq5} such that $(*)_{\bs\al,\bs a'}$ holds for some $\bs a'$, then in the inner sum $\{\cdots\}$ over $p,a_0,\ldots,a_p$, if $U(\be_1,\ldots,\be_p;\tau_{t_a},\tau_{s_a})\ne 0$ then the condition $\tau_{t_a}(\be_1)=\cdots=\tau_{t_a}(\be_p)$ holds automatically, and can be omitted. Thus, for such $n,\al_1,\ldots,\al_n$ the inner sum equals $U(\al_1,\ldots,\al_n;\tau_{t_0},\tau_{s_a})$ by \eq{co3eq5}. Since $(*)_{\bs\al,\bs a}$ holds for each nonzero term in \eq{co11eq5}, we deduce that
\ea
&[\M_\al^\ss(\tau_{s_a})]_\inv= \!\!\!\!\!\!\sum_{\begin{subarray}{l}n\ge 1,\; \al_1,\ldots,\al_n\in
C(\B)_\pe:\\ \al_1+\cdots+\al_n=\al,\; o_{\al_1}+\cdots+o_{\al_n}=o_\al, \\
\M_{\al_i}^\ss(\tau_{t_0})\ne\es,\; \text{all\/ $i$, $\exists\,\,\bs a'$ such that $(*)_{\bs\al,\bs a'}$ holds}\end{subarray}\!\!\!\!\!\!\!\!\!\!\!\!\!\!\!\!\!\!\!\!\!\!\!\!\!\!\!\!\!\!\!\!\!\!\!\!\!\!\!\!\!\!\!\!\!\!\!\!\!\!\!\!\!\!\!\!\!\!\!\!\!\!\!}
\begin{aligned}[t]
U(\al_1,\ldots,\al_n;\tau_{t_0},\tau_{s_a})\cdot{} &[\M_{\al_1}^\ss(\tau_{t_0})]_\inv\\
*\cdots *&[\M_{\al_n}^\ss(\tau_{t_0})]_\inv.
\end{aligned}
\label{co11eq6}
\ea

If $\al_1,\ldots,\al_n\in C(\B)_\pe$ with $\al_1+\cdots+\al_n=\al$, $\M_{\al_i}^\ss(\tau_{t_0})\ne\es$ for all $i$, and $U(\al_1,\ldots,\al_n;\tau_{t_0},\tau_{s_a})\ne 0$, then $(*)_{\bs\al,\bs a'}$ holds with $p'=1$, $a_0'=0$, $a_1'=n$. Thus we can omit the condition on $(*)_{\bs\al,\bs a'}$ in \eq{co11eq6}, as it only excludes terms with $U(\al_1,\ldots,\al_n;\tau_{t_0},\tau_{s_a})=0$. This proves \eq{co5eq34} with $\tau_{t_0},\tau_{s_a}$ in place of $\tau,\ti\tau$, giving the `only if' part of the lemma.

To extend this to `if and only if', consider the argument above, but without first supposing that \eq{co5eq34} holds with $\tau_{t_0},\tau_{t_a}$ in place of $\tau,\ti\tau$. Then \eq{co11eq2} may not hold when $\be=\al$. Since \eq{co11eq2} with $\be=\al$ appears with coefficient $U(\al;\hat\tau,\ti\tau)=1$ in \eq{co11eq5}, we see that the error in \eq{co5eq34} with $\tau_{t_0},\tau_{t_a}$ in place of $\tau,\ti\tau$ is the same as the error in \eq{co5eq34} with $\tau_{t_0},\tau_{s_a}$ in place of $\tau,\ti\tau$, so the former holds if and only if the latter does.
\end{proof}

By equation \eq{co3eq4} of Theorem \ref{co3thm2}, as $t_0=0$, equation \eq{co5eq34} with $\tau_{t_0},\tau_{t_0}$ in place of $\tau,\ti\tau$ reduces to $[\M_\al^\ss(\tau_{t_0})]_\inv=[\M_\al^\ss(\tau_{t_0})]_\inv$, and so is true. By induction on $a=1,\ldots,N$, we suppose by induction that \eq{co5eq34} holds with $\tau_{t_0},\tau_{s_b}$ in place of $\tau,\ti\tau$ for $1\le b<a$, and that \eq{co5eq34} holds with $\tau_{t_0},\tau_{t_b}$ in place of $\tau,\ti\tau$ for $0\le b<a$. In the inductive step, we apply Lemma \ref{co11lem4} to the second assumption with $b=a-1$ to show that \eq{co5eq34} holds with $\tau_{t_0},\tau_{s_a}$ in place of $\tau,\ti\tau$, and apply it again to show that \eq{co5eq34} holds with $\tau_{t_0},\tau_{t_a}$ in place of $\tau,\ti\tau$, completing the induction.

Therefore \eq{co5eq34} holds with $\tau_{t_0},\tau_{t_N}$ in place of $\tau,\ti\tau$. As $\tau_{t_0}=\tau_0=\tau$ and $\tau_{t_N}=\tau_1=\ti\tau$, this proves \eq{co5eq34} for $\al$. This completes the inductive step that we began at the beginning of Chapter \ref{co11}. Hence by induction \eq{co5eq34} holds for all $(\tau,T,\le),(\ti\tau,\ti T,\le)\in\sS$ and $\al\in C(\B)_\pe$. By Theorem \ref{co3thm3}, this implies that \eq{co5eq33} also holds for all $(\tau,T,\le),(\ti\tau,\ti T,\le)\in\sS$ and $\al\in C(\B)_\pe$, and the proof of Theorem \ref{co5thm3} is complete.

\begin{rem}
\label{co11rem1}
The argument above when we replace \eq{co11eq3} by \eq{co11eq4}, and the finiteness condition Assumption \ref{co5ass3}(b), may appear unnecessarily complicated. Na\"\i vely one might think we could just substitute \eq{co11eq2} into \eq{co11eq3}, and rearrange using \eq{co3eq5} to get \eq{co5eq34} with $\tau_{t_0},\tau_{s_a}$ in place of $\tau,\ti\tau$. 

The issue is that we would first have to drop the condition $\M_{\be_j}^\ss(\tau_{t_a})\ne\es$ in \eq{co11eq3}, as if $\M_{\be_j}^\ss(\tau_{t_a})=\es$ there may still be nonzero terms in \eq{co11eq2} for $\be=\be_j$ which are needed in the subsequent argument. But then \eq{co11eq3} without $\M_{\be_j}^\ss(\tau_{t_a})\ne\es$ could become an infinite sum, and the rearrangements to prove \eq{co5eq34} might involve infinitely many cancellations, as in the fallacious argument
\begin{equation*}
0=(1-1)+(1-1)+\cdots=1+(-1+1)+(-1+1)+\cdots=1.	
\end{equation*}
So we are careful to work only with finite sums at each step. 
\end{rem}

\medskip

\noindent{\small\sc The Mathematical Institute, Radcliffe
Observatory Quarter, Woodstock Road, Oxford, OX2 6GG, U.K.

\noindent E-mail: {\tt joyce@maths.ox.ac.uk.}}

\end{document}